\pdfoutput=1
\RequirePackage{ifpdf}
\ifpdf 
\documentclass[pdftex]{sigma}
\else
\documentclass{sigma}
\fi

\numberwithin{equation}{section}
\numberwithin{figure}{section}
\numberwithin{footnote}{section}

\newtheorem{thm}{Theorem}[section]
\newtheorem*{Theorem*}{Theorem}
\newtheorem{cor}[thm]{Corollary}
\newtheorem{lem}[thm]{Lemma}
\newtheorem{sublem}[thm]{Sublemma}
\newtheorem{prop}[thm]{Proposition}

\newtheorem{conj}[thm]{Conjecture}

\newtheorem{lemdef}[thm]{Lemma--Definition}
\newtheorem{sumary}[thm]{Informal Summary}

 { \theoremstyle{definition}
\newtheorem{defn}[thm]{Definition}

\newtheorem{exm}[thm]{Example}
\newtheorem{rem}[thm]{Remark}
\newtheorem{conds}[thm]{Condition}

\newtheorem{defnlem}[thm]{Definition--Lemma}
\newtheorem{assump}[thm]{Assumption}

\newtheorem{situ}[thm]{Situation}
\newtheorem{choice}[thm]{Choice}
\newtheorem{notation}[thm]{Notation} }

\usepackage{amscd}
\usepackage{amstext}
\usepackage{amsmath}
\usepackage[all]{xy}
\usepackage{yhmath}
\usepackage{mathrsfs}
\usepackage{bbding}

\usepackage[original]{imakeidx}
\makeindex 
\makeindex[name=syindex, title=List of notations]

\def\E{\ifmmode{\mathbb E}\else{$\mathbb E$}\fi} 
\def\N{\ifmmode{\mathbb N}\else{$\mathbb N$}\fi} 
\def\R{\ifmmode{\mathbb R}\else{$\mathbb R$}\fi} 
\def\Q{\ifmmode{\mathbb Q}\else{$\mathbb Q$}\fi} 
\def\C{\ifmmode{\mathbb C}\else{$\mathbb C$}\fi} 
\def\H{\ifmmode{\mathbb H}\else{$\mathbb H$}\fi} 
\def\Z{\ifmmode{\mathbb Z}\else{$\mathbb Z$}\fi} 
\def\P{\ifmmode{\mathbb P}\else{$\mathbb P$}\fi} 
\def\T{\ifmmode{\mathbb T}\else{$\mathbb T$}\fi} 
\def\SS{\ifmmode{\mathbb S}\else{$\mathbb S$}\fi} 
\def\DD{\ifmmode{\mathbb D}\else{$\mathbb D$}\fi} 
\def\K{\ifmmode{\mathbb K}\else{$\mathbb K$}\fi}

\begin{document}

\allowdisplaybreaks

\newcommand{\arXivNumber}{1706.02131}

\renewcommand{\thefootnote}{}

\renewcommand{\PaperNumber}{031}

\FirstPageHeading

\ShortArticleName{Unobstructed Immersed Lagrangian Correspondence and Filtered $A_{\infty}$ Functor}

\ArticleName{Unobstructed Immersed Lagrangian Correspondence \\and Filtered $\boldsymbol{A_{\infty}}$ Functor\footnote{This paper is a~contribution to the Special Issue on Integrability, Geometry, Moduli in honor of Motohico Mulase for his 65th birthday. The~full collection is available at \href{https://www.emis.de/journals/SIGMA/Mulase.html}{https://www.emis.de/journals/SIGMA/Mulase.html}}}

\Author{Kenji FUKAYA}

\AuthorNameForHeading{K.~Fukaya}

\Address{Yau Mathematical Sciences Center, Jingzhai, Tsinghua University,\\ Haidian District, Beijing, 100084, P.R.~China}
\Email{\href{mailto:fukayakenji@mail.tsinghua.edu.cn}{fukayakenji@mail.tsinghua.edu.cn}}
\URLaddress{\url{https://www.math.kyoto-u.ac.jp/~fukaya/fukaya.html}}

\ArticleDates{Received October 11, 2019, in final form March 06, 2025; Published online April 29, 2025}

\Abstract{In this paper, we `construct' a 2-functor from the unobstructed immersed Weinstein category to the category of all filtered $A_{\infty}$ categories. We consider arbitrary (compact) symplectic manifolds and its arbitrary (relatively spin) immersed Lagrangian submanifolds. The filtered $A_{\infty}$ category associated to $(X,\omega)$ is defined by using Lagrangian Floer theory in such generality, see Akaho--Joyce (2010) and Fukaya--Oh--Ohta--Ono (2009). The morphism of unobstructed immersed Weinstein category (from~$(X_1,\omega_1)$ to $(X_2,\omega_2)$) is by definition a~pair of an immersed Lagrangian submanifold of the direct product and its bounding cochain (in the sense of Akaho--Joyce (2010) and Fukaya--Oh--Ohta--Ono (2009)). Such a~morphism transforms an (immersed) Lagrangian submanifold of $(X_1,\omega_1)$ to one of~$(X_2,\omega_2)$. The key new result proved in this paper shows that this geometric transformation preserves unobstructedness of the Lagrangian Floer theory. Thus, this paper generalizes earlier results by Wehrheim--Woodward and Mau's--Wehrheim--Woodward so that it works in complete generality in the compact case. The main idea of the proofs are based on Lekili--Lipyanskiy's~Y diagram and a lemma from homological algebra, together with systematic use of Yoneda functor. In other words, the proofs are based on a different idea from those which are studied by Bottmann--Mau's--Wehrheim--Woodward, where strip shrinking and figure 8 bubble plays the central role.}


\Classification{53D35; 53D40; 57R56; 53D12; 53D37; 57R17}

\renewcommand{\thefootnote}{\arabic{section}.\arabic{footnote}}
\setcounter{footnote}{0}

\tableofcontents

\section{Introduction}\label{sec:intro}

The purpose of this paper is to `construct' a 2-functor
from the `unobstructed immersed Weinstein category'
to the `category of all filtered $A_{\infty}$ categories'.
The next definition is a variation of
one proposed by Weinstein \cite{Wi}.

\begin{defn}[informal definition]\label{def11}
The {\it unobstructed immersed Weinstein category}
\index{Unobstructed immersed Weinstein category}
is a category whose object is a compact symplectic
manifold and a morphism from $(X_1,\omega_1)$
to $(X_2,\omega_2)$ is a pair of an immersed Lagrangian
submanifold $L_{12}$ of $(X_1 \times X_2,-\pi_1^*(\omega_1)+\pi_2^*(\omega_2)$)
and a~bounding cochain $b_{12}$ on it.
\end{defn}

The notion of a bounding cochain is introduced in \cite{fooobook} and
its generalization to the immersed case is by \cite{AJ}.
We emphasis that Definition~\ref{def11} is an informal definition.
Various issues which will appear when one tries to define such a
2-category literary are discussed in
Section~\ref{issueformal}.

We consider a 2-category whose objects are (strict and unital)
filtered $A_{\infty}$ categories (see \mbox{\cite{fu4,ancher}}
and Section~\ref{sec:HFIm} for its definition) and
morphisms are (strict and unital) filtered $A_{\infty}$
functors. See Section~\ref{2-category formulation}, Theorem~\ref{thm94} and Section~\ref{sec:catofAinfcat} for a version of the
construction of such a 2-category.

The main result of this paper could be summarized as follows.

\begin{sumary}\label{thm02}
There exists a $2$-functor from the unobstructed immersed
Weinstein category to the $2$-category of all filtered $A_{\infty}$
categories.
\end{sumary}

This statement is informal and the author does {\it not} claim that its
proof is in this paper. The precise
statements which we prove in this paper will be given in this introduction
and the main body of the paper.
The relation between those results (proved in this paper) and
the results which would literary prove Informal Summary~\ref{thm02}
is discussed in Section~\ref{issueformal}.

The idea to associate an $A_{\infty}$ category (whose object is a Lagrangian
submanifold and whose morphisms are Floer cohomology) is
started by the author's paper \cite{fu0}
(inspired by a S. Donaldson's talk at University of Warwick 1992).
The most essential step to make this construction rigorous
was carried out in \cite{fooobook}, based on the virtual fundamental
chain technique (see~\cite{FO}).
The work~\cite{fooobook} contains the detailed proof of the cases of a single Lagrangian submanifold
and a~pair of Lagrangian submanifolds.
The construction of a (unital and strict) filtered $A_{\infty}$ category
based on the Lagrangian Floer theory in \cite{fooobook} along the same line as \cite{fooobook}
was written in \cite{AFOOO,fu4,ancher}.
Akaho and Joyce \cite{AJ} generalized this story and include
Lagrangian submanifolds which are not necessary embedded
but are immersed.
Thus we obtain the next theorem.

\begin{thm}\label{thm03}
Let $(X,\omega)$ be a compact symplectic manifold
and $\mathbb L$ a finite set of its spin immersed Lagrangian
submanifolds.\footnote{In the introduction, we assume
spinness of Lagrangian submanifolds rather than relatively-spinness,
for simplicity. The statement in the relatively spin case will be
given in the main body of the paper.}
We assume that the self intersection of elements of $\mathbb L$ and
intersection between two elements of $\mathbb L$ are transversal.
Then there exists a $($strict and unital$)$ filtered $A_{\infty}$ category,
$\mathfrak{Fukst}((X,\omega),\mathbb L)$, such that
\begin{enumerate}\itemsep=0pt
\item[$(1)$]
An object of $\mathfrak{Fukst}((X,\omega),\mathbb L)$ is a pair
$(L,b)$ where $L$ is an element of $\mathbb L$ and
$b$ is a bounding cochain of $L$ in the sense of {\rm\cite{AJ,fooobook}}.
\item[$(2)$]
The module of morphisms $CF((L_1,b_1),(L_2,b_2))$ from $(L_1,b_1)$ to
$(L_2,b_2)$ is given as follows:
\begin{enumerate}\itemsep=0pt
\item[$(a)$] If $L_1 \ne L_2$, then $CF((L_1,b_1),(L_2,b_2))$ is the free $\Lambda_0$ module whose basis is identified with
the intersection $L_1 \cap L_2$.
Here the universal
Novikov ring $\Lambda_0$ is defined in Definition {\rm\ref{def21}}.
\item[$(b)$]
If $L_1 = L_2 = L$, then $CF((L_1,b_1),(L_2,b_2))$ is the
completion of the tensor product~\smash{$\Omega\bigl(\tilde L \times_X \tilde L\bigr) \otimes \Lambda_0$} of the de Rham complex
\smash{$\Omega\bigl(\tilde L \times_X \tilde L\bigr)$}
and $\Lambda_0$.
Here our immersed Lagrangian submanifold $L$ is given by an immersion
$\tilde L \to X$ and $\tilde L \times_X \tilde L$ is the fiber product
of $\tilde L$ with itself.
\end{enumerate}
\item[$(3)$]
The cohomology group of $CF((L_1,b_1),(L_2,b_2))$ is the Floer cohomology
$HF((L_1,b_1),(L_2,\allowbreak b_2))$ defined in {\rm\cite{AJ,fooobook}}.
\end{enumerate}
\end{thm}
Theorem~\ref{thm03}
is Theorem~\ref{AJtheorem} in Section~\ref{sec:HFIm}, which is slightly more general.
\begin{rem}
In item (2b), we may also take \smash{$H\bigl(\tilde L \times_X \tilde L;\Lambda_0\bigr)$}
(the cohomology group with $\Lambda_0$ coefficient)
instead of \smash{$\Omega\bigl(\tilde L \times_X \tilde L\bigr)\,\widehat{\otimes}\,\Lambda_0$}.
The process to produce a structure on
\smash{$H\bigl(\tilde L \times_X \tilde L;\Lambda_0\bigr)$}
from one on \smash{$\Omega\bigl(\tilde L \times_X \tilde L\bigr) \otimes \Lambda_0$}
is purely algebraic and automatic. See \cite[Theorem 5.4.2']{fooobook}
for example.
\end{rem}
Theorem~\ref{thm03} is not new and the most essential
part of its proof had been given in \cite{AJ,fooobook}.
(We use the de Rham version, while \cite{AJ,fooobook} uses the
singular homology version. This difference however is
not important but is rather of technical nature.)
In the de Rham version, it is also written and proved in
\cite{AFOOO}.

Theorem~\ref{thm03} is the object part of `2-functor'
mentioned in Informal Summary \ref{thm02}.
The main new point of this paper is the
morphism part of the `2 functor' mentioned in Informal Summary~\ref{thm02}.
The next theorem is the key new result.
Let $(X_i,\omega_i)$ be a compact symplectic manifold
for $i=1,2$. We assume they are spin.\footnote{The case
when $X_1$ or $X_2$ is not spin is included in the main
body of the paper.}
Let $L_1$, $L_{12}$ be spin immersed Lagrangian
submanifolds of $(X_1,\omega_1)$ and $(X_1 \times X_2,-\pi_1^*(\omega_1)
+\pi_2^*(\omega_2))$, respectively.
We say that they are transversal if the
fiber product
$\tilde L_1 \times_{X_1} \tilde L_{12}$ is transversal.
In that case, the map $\tilde L_1 \times_{X_1} \tilde L_{12} \to X_2$
defines an immersed Lagrangian submanifold which we
write $L_1 \times_{X_1} L_{12}$.
\begin{thm}\label{them4}
If $L_1$ and $L_{12}$ are unobstructed\footnote{A Lagrangian
submanifold is said to be unobstructed if
there exists a bounding cochain of it.} and the immersion
$\tilde L_1 \times_{X_1} \tilde L_{12} \to X_2$ is self-clean, then
$L_1 \times_{X_1} L_{12}$ is also unobstructed.
There exists a way to obtain a bounding cochain
of~${L_1 \times_{X_1} L_{12}}$ from bounding cochains
of $L_1$ and of $L_{12}$, which is independent of the choices up to gauge equivalence.
\end{thm}
Theorem~\ref{them4} is Theorems~\ref{thm61} and \ref{th72}, which are
proved in Sections \ref{sec:Unobstructedness}
and \ref{sec:represent}.
Note that for generic (embedded) Lagrangian submanifolds $L_1$, $L_{12}$
of $(X_1,\omega_1)$ and $(X_1 \times X_2,-\pi_1^*(\omega_1)
+\pi_2^*(\omega_2))$,
the fiber product $L_1 \times_{X_1} L_{12}$ is an
immersed Lagrangian submanifold of $X_2$.
However, it is not necessary embedded.
Therefore, including immersed Lagrangian submanifolds
is inevitable.
\begin{rem}\label{rem16}
\quad
\begin{enumerate}\itemsep=0pt
\item[(1)]
The relation between a Lagrangian correspondence and an $A_{\infty}$ functor
was studied in the earlier works by Wehrheim--Woodward, Ma'u--Wehrheim--Woodward
(see \cite{MWW,WW2} etc.). Note that, in their situation where all the Lagrangian submanifolds
involved are embedded and monotone, the statement corresponding to
Theorem~\ref{them4} is classical (due to Oh), since we can take $0$ as the bounding cochain.
\item[(2)] To include more general objects than embedded and monotone
Lagrangian submanifolds,
Wehrheim--Woodward proceeds
as follows. They first consider embedded and monotone
Lagrangian submanifolds (with bounding cochain $0$).
They then enhance the set of such Lagrangian submanifolds
so that a sequence of Lagrangian correspondences
\[
L_1 \times_{X_1} L_{12} \times_{X_2} L_{23} \times_{X_3} \dots \times_{X_{k-1}} L_{(k-1)k}
\]
is regarded as an object of $\mathfrak{Fuk}^{\#}(X_k,\omega_k)$, the extended version of
$\mathfrak{Fuk}(X_k,\omega_k)$.
Theorem~\ref{them4} enables us to work with genuine geometric
Lagrangian submanifolds rather than extended objects.
We will discuss the relation between our results and one by \cite{MWW,WW2}
more in Section~\ref{relWWM}
\item[(3)]
The statement of Theorem~\ref{them4} was known as a conjecture for a while.
For example, the author discussed this conjecture with several mathematicians
during the years
2008--2015.
It was also mentioned by K.~Wehrheim's talk in 2012 \cite{Weh}
and is written as a `conjecture' in \cite{bww}.
More precisely, it had been conjectured that
the virtual fundamental chain of an appropriate
moduli space of Figure 8 bubbles gives the bounding
cochain in Theorem~\ref{them4}. The conjecture of this
form is still open.
It is the opinion of the author that to prove this version of the conjecture is a
very interesting analytic problem.
If this conjecture is proved and a bounding cochain is obtained
as the virtual fundamental chain of the moduli space
of Figure 8 bubbles, the author has no doubt that such
a bounding cochain is gauge equivalent to the bounding cochain we obtained
in Theorem~\ref{them4}.
We will discuss this point more in Section~\ref{reltoBW}.
\item[(4)]
Until 2015, the author did not have an idea to prove Theorem~\ref{them4}
other than those by using `strip shrinking' and `Figure 8 bubble',
which are emphasized in \cite{bww}.\footnote{See \cite{bww}
or Sections~\ref{reltoBW} and \ref{relWWM} for `strip shrinking' and `Figure 8 bubble'.
Studying them certainly are interesting in its own and potentially can
be applied to various geometric problems.}
By this reason the author did not have a plan to study this conjecture until 2015.
In May 2015, the author realized that using the method of
Lekili--Lipyanskiy \cite{LL} and homological algebra we can
prove Theorem~\ref{them4} much easier than the idea using
`strip shrinking' or `Figure 8 bubble'.
He then started working on Lagrangian correspondence and
its relation to Lagrangian Floer theory.
(The main motivation of the author's study
is its application to the gauge theory (see \cite{df,fu9,takagi})).
This paper is an outcome of that study.
\end{enumerate}
\end{rem}
We also remark that to define Floer cohomology of a Lagrangian
submanifold (beyond the monotone or exact cases) we need a
bounding cochain. So proving the existence of a bounding cochain
is the key step for applications of Lagrangian Floer theory.
In general, proving the existence of a bounding cochain is not easy.
Theorem~\ref{them4} provides a useful tool to prove it.

The next theorem is a more functorial version of Theorem~\ref{them4}.
Let $\mathbb L_1$, $\mathbb L_2$ and $\mathbb L_{12}$ be finite sets of
spin immersed Lagrangian
submanifolds of $(X_1,\omega_1)$, $(X_2,\omega_2)$ and $(X_1 \times X_2,-\pi_1^*(\omega_1)
+\pi_2^* (\omega_2))$, respectively.
We assume each of them satisfies the transversality conditions in Theorem~\ref{thm03}.
Moreover, we assume that for each $L_1 \in \mathbb L_1$ and $L_{12} \in \mathbb L_{12}$
the fiber product $L_1 \times_{X_1} L_{12}$ is transversal and is an element
of $\mathbb L_2$.
\begin{thm}\label{thm5}
In the situation of Theorem {\rm\ref{them4}},
there exists a filtered $A_{\infty}$ bi-functor
\[
\mathfrak{Fukst}((X_1,\omega_1),\mathbb L_1)
\times
\mathfrak{Fukst}((X_1 \times X_2,-\pi_1^*(\omega_1)
+\pi_2^* (\omega_2)),\mathbb L_{12})
\to
\mathfrak{Fukst}((X_2,\omega_2),\mathbb L_2)
\]
such that it sends the pair of objects $(L_1,b_1)$,
$(L_{12},b_{12})$ to $L_1 \times_{X_1} L_{12}$ equipped
with the bounding cochain given in Theorem {\rm\ref{them4}}.
\end{thm}
See Definition~\ref{bi-functor} for the definition of a filtered $A_{\infty}$ bi-functor.
Theorem~\ref{thm5} provides the morphism part of the `2-functor' mentioned in Informal Summary \ref{thm02}.
Theorem~\ref{thm5} is
Corollary~\ref{cor73} and is proved in Section~\ref{sec:represent}.
We call the bi-functor in Theorem~\ref{thm5} the correspondence
bi-functor.
We like to mention that in the situation when all the Lagrangian
submanifolds involved are embedded and monotone,
Theorem~\ref{thm5} was proved by Ma'u--Wehrheim--Woodward
in \cite{MWW}.
See Section~\ref{relWWM} for more explanation on the
relation of Theorem~\ref{thm5} to \cite{MWW}.

The next theorem gives a definition of the composition of morphisms
in unobstructed immersed Weinstein category.
In other words, Theorem~\ref{thm05} could be used to give a definition of
unobstructed immersed Weinstein category as a
(topological) 2-category.\footnote{We say `topological' 2-category since to
compose two (unobstructed immersed) Lagrangian correspondences
we need to assume transversality. Therefore, morphisms can be composed only
generically.}

Let $(X_i,\omega_i)$, $i=1,2,3$, be compact symplectic manifolds which are spin.
Let $\mathbb L_{ij}$, $(ij) = (12)$, $(23)$ or $(13)$, be
finite sets of spin Lagrangian submanifolds of
$(X_i \times X_j,-\pi_1^*(\omega_i)
+\pi_2^* (\omega_j))$.
We assume that for any $L_{12} \in \mathbb L_{12}$, $L_{23} \in \mathbb L_{23}$
the fiber product $L_{12} \times_{X_2} L_{23}$ is transversal
and becomes an element of $\mathbb L_{13}$.
\begin{thm}\label{thm05}
There exists a filtered $A_{\infty}$ bi-functor
\begin{align*}
\mathfrak{comp}\colon\ &
\mathfrak{Fukst}((X_1\! \times X_2,-\pi_1^*(\omega_1)
+\pi_2^* \omega_2),\mathbb L_{12})\times
\mathfrak{Fukst}((X_2\! \times X_3,-\pi_1^*(\omega_2)
+\pi_2^* (\omega_3)),\mathbb L_{23}) \\
 & \to
\mathfrak{Fukst}((X_1 \times X_3,-\pi_1^*(\omega_1)
+\pi_2^* (\omega_3)),\mathbb L_{13})
\end{align*}
such that it sends a pair of objects
$(L_{12},b_{12})$, $(L_{23},b_{23})$ to
$(L_{13},b_{13})$, where $L_{13} = L_{12} \times_{X_2} L_{23}$
and $b_{13}$ is a bounding cochain of $L_{13}$ which is determined
from $b_{12}$ and $b_{23}$ in a way independent of the choices up to gauge equivalence.
We call this functor the composition functor.

The composition functor is associative, in the sense
that the next diagram commutes up to homotopy equivalence, as $A_{\infty}$ tri-functors
\begin{equation}\label{Diagram1}
\begin{CD}
\displaystyle{\mathfrak{Fukst}(X_1 \times X_2)
\times \mathfrak{Fukst}(X_2 \times X_{3})
\atop \!\!\!\!\!\!\!\!\!\!\!\!\!\!\!\!\!\!\!\!\!
\!\!\!\!\!\!\!\!\!\!\!\!\!\!\!\!\!\!\!\times
\mathfrak{Fukst}(X_3 \times X_{4})} @ >>>
\displaystyle{\mathfrak{Fukst}(X_1 \times X_{3})
\atop\times\mathfrak{Fukst}(X_3 \times X_{4})} \\
@ VVV @ VVV\\
\mathfrak{Fukst}(X_1 \times X_{2})
\times \mathfrak{Fukst}(X_2 \times X_{4}) @ > >> \mathfrak{Fukst}(X_1 \times X_{4}).
\end{CD}
\end{equation}

\end{thm}
The first half of Theorem~\ref{thm05} is Theorems~\ref{comp} and \ref{comp2}
which are proved in Section~\ref{sec:comp}.
The second half of Theorem~\ref{thm05} is Theorem~\ref{asscompmain} proved in
Section~\ref{sec:associ}.
\begin{rem}
Actually Theorem~\ref{thm5} follows from
Theorem~\ref{thm05} by putting $X_1$ to be a point.
\end{rem}

In the situation when all the Lagrangian
submanifolds involved are embedded and monotone,
Theorem~\ref{thm05} (and Theorem~\ref{thm18} below) were also proved by Ma'u--Wehrheim--Woodward
in~\cite{MWW}.
We also remark that Wehrheim--Woodward and Ma'u--Wehrheim--Woodward
studied a~fiber product of Lagrangian correspondences
(under the assumption that it becomes embedded Lagrangian correspondence)
in their study of the composition of Lagrangian correspondences.
See Section~\ref{relWWM} for more explanation on the
relation of Theorems~\ref{thm5}, \ref{thm05} and \ref{thm18} to \cite{MWW}.

The next theorem says that the correspondence bi-functor in Theorem~\ref{thm5}
is compatible with the composition functor in Theorem~\ref{thm05}.
To state it, we need some digression.
Let $\mathscr C_i$ be a~strict and unital filtered $A_{\infty}$
category for $i=1,2$. Then we can define a filtered
$A_{\infty}$
category~${\mathcal{FUNC}(\mathscr C_1,\mathscr C_2)}$
whose object is a strict and unital filtered $A_{\infty}$
functor. (This is the unital and strict version of
Theorem~\ref{th21010} whose proof is the same as Theorem~\ref{th21010}.)
For three strict and unital filtered $A_{\infty}$ categories
$\mathscr C_i$, $i=1,2,3$, we can define a filtered $A_{\infty}$
bi-functor
\begin{equation}\label{11comp}
\mathcal{FUNC}(\mathscr C_1,\mathscr C_2)
\times
\mathcal{FUNC}(\mathscr C_2,\mathscr C_3)
\to \mathcal{FUNC}(\mathscr C_1,\mathscr C_3),
\end{equation}
which gives a composition of filtered $A_{\infty}$
functors among objects.
(See Theorem~\ref{thm94}.)
The bi-functor \eqref{11comp} is associative.
Roughly speaking, \eqref{11comp} is defined as follows.
We first define a~homotopy equivalence from functor category
$\mathcal{FUNC}(\mathscr C_1,\mathscr C_2)$ to a
full subcategory of the DG-category of left $\mathscr C_1$ and right $\mathscr C_2$ bi-modules.
(This is a version of Yoneda's lemma.)
We also prove that the composition of $A_{\infty}$ functors
corresponds to the tensor product of the bi-modules.
Then using the fact that tensor product of bi-modules is an
object part of the DG-bi-functor,
we obtain \eqref{11comp}.
(See Section~\ref{sec:catofAinfcat}.)

On the other hand, the correspondence bi-functor in Theorem~\ref{thm5} can be reinterpreted
as a~filtered $A_{\infty}$
functor
\begin{gather}
\mathfrak{Fukst}((X_1 \times X_2,-\pi_1^*(\omega_1)
+\pi_2^*(\omega_2)),\mathbb L_{12})\nonumber \\
\qquad\to
\mathcal{FUNC}(\mathfrak{Fukst}((X_1,\omega_1),\mathbb L_1),
\mathfrak{Fukst}((X_2,\omega_2),\mathbb L_2))\label{map11}
\end{gather}
to the functor category.
\begin{thm}\label{thm18}
The next diagram commutes up to homotopy equivalence
\begin{equation}\label{thm110dia}
\begin{CD}
\mathfrak{Fukst}(X_1 \times X_2)
\times \mathfrak{Fukst}(X_2 \times X_3) @ >>>
\mathfrak{Fukst}(X_1 \times X_3) \\
@ VVV @ VVV\\
\displaystyle{\mathcal{FUNC}(\mathfrak{Fukst}(X_1),\mathfrak{Fukst}(X_2))
\atop \times \mathcal{FUNC}(\mathfrak{Fukst}(X_2),\mathfrak{Fukst}(X_3))} @ > >> \mathcal{FUNC}(\mathfrak{Fukst}(X_1),\mathfrak{Fukst}(X_3)).
\end{CD}
\end{equation}
Here the vertical arrows are functors \eqref{map11}, the
upper horizontal arrow is the composition functor in Theorem~{\rm\ref{thm05}} and
lower horizontal arrow is the functor~\eqref{11comp}.\footnote{We take appropriate finite sets $\mathbb L_{ij}$ of Lagrangian submanifolds of $X_i \times X_j$.}
\end{thm}
Theorem~\ref{thm18} is Theorems~\ref{thm93} and \ref{thm109}
which are proved in Sections~\ref{sec:comptibility} and~\ref{2-category formulation}.

\begin{rem}
The object part of Theorem~\ref{thm18}, that
is, the commutativity of the diagram~\eqref{thm110dia} as the maps between the sets of objects, implies Theorem~\ref{thm5} by putting $X_1$ to be a point.
The (homotopy) commutativity of diagram~\eqref{thm110dia} as $A_{\infty}$ bi-functors is
more involved.
\end{rem}

All the constructions of this paper are based on a
study of moduli spaces of pseudo-holomor\-phic curves.
Even though we use the moduli space of pseudo-holomorphic quilts
in the sense of~\cite{WW3} we do {\it not} use
the most difficult part of the analytic study of
the moduli space of pseudo-holomorphic quilts.
Especially we do {\it not} study `strip shrinking' and
`Figure 8 bubble'.
Our proof relies much on the cobordism argument
which was initiated by Y.~Lekili and M.~Lipyanskiy
\cite{LL} and various technique
from homological algebra.
By this reason, we do not need new analytic detail to carry out
in this paper,
except we need to take a slightly different compactification
of the moduli space of pseudo-holomorphic disks
bounding a Lagrangian submanifold $L_{12}$ of
the product. This is because otherwise the
moduli space of pseudo-holomorphic quilts
would not carry a Kuranishi structure.
We will explain this point in Section~\ref{sec:directcomp} and also
provide the detail of this different compactification.

In Sections \ref{sec:homotopyfafunc}--\ref{sec:independence3},
we show that various filtered $A_{\infty}$ (bi)-functors
we construct in this paper are independent of the
choices involved and also of the Hamiltonian isotopies of the
Lagrangian submanifolds involved.

In Section~\ref{sec:kuneth}, we show that by a similar
method used in the other part of this paper,
we can show K\"unneth theorem in Lagrangian Floer theory.
(We remark that K\"unneth theorem in Lagrangian Floer theory
is also proved by \cite{Limo1,Limo}.)

Section~\ref{sec:orient} is devoted to the discussion
of sign and orientation. More arguments on sign and orientation are given
in the paper \cite{ono2} written by K.~Ono.

Section~\ref{sec:remark} is a brief discussion
on two points. One is the relation of this paper
to the works by Wehrheim--Woodwards--Ma'u--Bottman.
The other is an issue which will appear to define/prove
`Definition~\ref{def11}'/`Informal Summary~\ref{thm02}' literary.

We expect that there are various applications of the whole
construction (especially the part to construct a
filtered $A_{\infty}$ functor from a Lagrangian correspondence and
several compatibility statements about it, which is new in this paper).
Some of the applications are now on the way being worked out and being written
or already available as a preprint. (See \cite{df,elekili,fu9,takagi}
and etc.)
In this paper, we concentrate in defining the basic objects in as much
general form as possible, leaving applications to other
papers.
A generalization of the story to the case of non-compact Lagrangian
submanifolds is now studied by Yuan--Gao \cite{yu}.

The construction of this paper is based on various earlier works.
The author tried to make this paper independent from various
earlier papers, except the detail of the proofs, as much as possible.
By this reason, this paper contains several review sections.
Another reason why the review sections are included is that
we need to rewrite some of the earlier results to those
based on the de Rham version of virtual fundamental chain
technique, which we use systematically in this paper.
We refer \cite{foootech2, foootech22,fooonewbook} for the most detailed exposition
of the de Rham version of virtual fundamental chain
technique (Kuranishi structures and CF-perturbations).
If the reader wants to know definitions and
 statements of the theory in \cite{foootech2, foootech22, fooonewbook} only
(and not its proof),
there is a summary in \cite[Part~7]{fooospectr}.

The construction of Kuranishi structures on the moduli spaces
of pseudo-holomorphic curves are written in detail in
\cite[Part~4]{foootech}, \cite{fooo:const1, fooo:const2,foooanalysis}.
It is written also in Section~\ref{sec:directcomp} of this paper
emphasizing the part where the construction we need in this
paper is (slightly) different from the other papers.

The results of this paper were announced in \cite{fu9,takagi} together
with the main idea of its proof.

\section[Filtered $A_{\infty}$ category: Review]{Filtered $\boldsymbol{A_{\infty}}$ category: Review}
\label{sec:Ainfcat}

This section is a review of the homological algebra of filtered $A_{\infty}$ categories.
There is nothing really new in this section.
Our purpose here is to provide the precise definitions of the various
notions we use in this paper.
We give proofs only in the case when the author is unable to find
an appropriate reference in the literature.
Our main reference in this section is \cite{fu4}.
There are other references such as \cite{BLM, DL, fu2, Ke, Ke1, KS, Lef, Se}.
In this section, we will discuss the algebraic side of the story only.
In the case when the reader has certain knowledge of $A_{\infty}$
categories, the reader can skip this section and comes back
when it is used in later sections.

\subsection[$A_\infty$ category]{$\boldsymbol{A_{\infty}}$ category}
\label{subsec:Ainfcat}

We first recall certain notations.
\begin{defn}\label{def21}\quad
\begin{enumerate}\itemsep=0pt
\item[(1)]
Let $R$ be a commutative ring with unit.
We denote by
$
\Lambda_0^R
$
the set of all the formal sums
\begin{equation}\label{form21}
\sum_{i=0}^{\infty} a_i T^{\lambda_i},
\end{equation}
where $a_i \in R$, $\lambda_i \in \R$ and
$0 = \lambda_0 < \lambda_1 < \dots < \lambda_i < \lambda_{i+1}
< \cdots$ with
$\lim_{i\to \infty} \lambda_i = +\infty$.
We can define a ring structure on $
\Lambda_0^R
$ in an obvious way.

We call $\Lambda_0^R$ the {\it universal Novikov ring}. \index{universal Novikov ring}
In the case when $R$ is a field, its maximal ideal is the set of formal sums \eqref{form21}
with $a_0 = 0$. We write it as $\Lambda_+^R$.
In the case when $R$ is a~field,
the field of fractions of $\Lambda_0^R$
is the set of the formal sums of the form~\eqref{form21}
such that~${\lambda_0 < \lambda_1 < \dots < \lambda_i < \lambda_{i+1}
< \cdots}$ with $\lambda_i \in \R$ and
$\lim_{i\to \infty} \lambda_i = +\infty$.
We denote it by~$\Lambda^R$ and call it
the {\it universal Novikov field}.\index{universal Novikov field}
We use the same notation $\Lambda^R_+$ \big(resp.\ $\Lambda^R$\big) for this ideal (resp.\ ring)
in case $R$ is a ring but is not a field.
\index[syindex]{Lambdaplus@$\Lambda_+$}\index[syindex]{Lambda@$\Lambda$}
We call $R$ the {\it ground ring}.\index{ground ring}
Sometimes we omit $R$ from the notation and write
$\Lambda_0$ etc.\ \index[syindex]{Lambda0@$\Lambda_0$} in place of $\Lambda^R_0$ etc.
In the geometric applications in this paper,
we use $R = \R$, since we use the de Rham model
for homology theory of spaces.
\item[(2)]
We define a filtration $\bigl\{\mathfrak F^{\lambda}\Lambda_0
\mid \lambda \ge 0\bigr\}$ as follows.
The subset $\mathfrak F^{\lambda}\Lambda_0$ of $\Lambda_0$ consists of elements~\eqref{form21}
such that $\lambda_i < \lambda$ implies $a_i = 0$.
We call it the {\it energy filtration}.
\index{energy filtration}
It induces a filtration on $\Lambda$ and $\Lambda_+$ in an obvious way.
The energy filtration defines a metric on $\Lambda_0$,
$\Lambda$, $\Lambda_+$
such that
$\mathfrak F^{\lambda}\Lambda_0$ is the $e^{-\lambda}$-neighborhood
of $0$.
The rings
$\Lambda_0$,
$\Lambda$, $\Lambda_+$
are complete with respect to this metric.
We call this metric the {\it $T$-adic metric}.\index{$T$-adic metric}
We use also the name energy filtration for the filtration
of various $\Lambda_0$ (or $\Lambda$) modules induced by this filtration of~$\Lambda_0$.\looseness=-1

\item[(3)]
A {\it discrete monoid}\index{discrete monoid} $G$ is a discrete subset of $\R_{\ge 0}$
such that $0 \in G$ and $g_1,g_2 \in G \Rightarrow g_1+g_2 \allowbreak\in G$.
\item[(4)]
For a discrete monoid $G$, we define a subring $\Lambda_G$ of $\Lambda_0$,
where a formal sum \eqref{form21} is an element of $\Lambda_G$ if and only if
$\lambda_i \in G$ for all $i$ with $a_i \ne 0$.
The $T$-adic metric of $\Lambda_0$ induces one on $\Lambda_G$.
\item[(5)]
Let $\overline C$ be a free $R$ module.
We denote by $C$ the completion of $\overline C \otimes_R \Lambda_0^R$.
Here the $T$ adic metric on $\overline C \otimes_R \Lambda_0^R$
is induced from one on $\Lambda_0^R$ in an obvious way and
the completion is taken with respect to this metric.
We call such $C$ a {\it completed free $\Lambda_0$ module}.
\index{completed free $\Lambda_0$ module}
We write~${
\mathfrak F^{\lambda}C =\big\{x \in C \mid x \equiv 0 \mod T^{\lambda}\big\}}$.
An element of $C$ is identified with an
(infinite) sum
$
\sum_{i=0}^{\infty} T^{\lambda_i}x_i
$
such that $x_i \in \overline C$ and $\lambda_i \in \R_{\ge 0}$
with $\lim_{i\to \infty}\lambda_{i}=+\infty$.
\item[(6)]
For two completed free $\Lambda_0$ modules $C_1$, $C_2$, we denote by
$
C_1 \,\widehat\otimes\, C_2
$
the $T$-adic completion of the algebraic tensor product
over $\Lambda_0$.
When $C_i$ is the completion of $\overline C_i \otimes_R \Lambda_0$,
for $i=1,2$, $
C_1 \,\widehat\otimes\, C_2
$ \index[syindex]{tzzensor@$ \widehat\otimes$}is the completion of $\overline C_1 \otimes_R \overline C_2 \otimes_R \Lambda_0$.
An element of $C_1 \,\widehat\otimes_{\Lambda_0}\, C_2$ is identified with an
(infinite) sum
$
\sum_{i=0}^{\infty} x_i \otimes y_i
$
such that $x_i \in \mathfrak F^{\lambda_{i,1}}C_1$, $y_i \in \mathfrak F^{\lambda_{i,2}}C_2$
with $\lim_{i\to \infty}\lambda_{i,1} + \lambda_{i,2}=+\infty$.
\item[(7)]
If $\overline C$ is graded, then $C$ is graded.
(Here we consider either $\Z$ grading or $\Z_{2N}$ grading.
In our geometric application, we mostly use $\Z_2$ grading, for the sake of simplicity.)
Suppose~$C$ is graded. We define its {\it degree shift} $C[1]$ as follows.
$C[1]^m= C^{m+1}$. Here $C^m$ is degree~$m$ part.\looseness=-1\index{degree shift}
\item[(8)]
An element $x$ of a completed free $\Lambda_0$ module $C$ is said to be {\it $G$-gapped}
\index{$G$-gapped} if
\smash{$
x = \sum_{g \in G}T^g x_g
$}
where $x_g \in \overline C$.
\item[(9)]
A $\Lambda_0$ module homomorphism $\varphi$ between completed free $\Lambda_0$ modules
$C_1$, $C_2$
are said to be $G$-{\it gapped} \index{$G$-gapped} if it sends an arbitrary $G$-gapped element to
a $G$-gapped element.
This condition is equivalent to the condition that
\begin{equation}\label{phiexp}
\varphi = \sum _{g \in G}T^g \varphi_g,
\end{equation}
where $\varphi_g \colon \overline C_1 \to \overline C_2$ are $R$ module homomorphisms.
\item[(10)]
For a $G$-gapped homomorphism $\varphi$ as in \eqref{phiexp},
we write
$
\overline{\varphi} = \varphi_0 \colon \overline{C}_1 \to \overline{C}_2
$
and call it the {\it $R$-reduction} \index{$R$-reduction} of $\varphi$.
\end{enumerate}
\end{defn}

\begin{defn}\label{defn22}
A {\it non-unital curved filtered $A_{\infty}$ category}
\index{filtered $A_{\infty}$ category}\index{non-unital}
\index{curved} $\mathscr C$ is a collection of the set
$\mathfrak{Ob}(\mathscr C)$, \index[syindex]{OB@$\mathfrak{Ob}$} the set of objects, a graded completed free $\Lambda_0$ module $\mathscr
C(c_1,c_2)$ for each
$c_1,c_2 \in \mathfrak{Ob}(\mathscr C)$, and the operations
\[
{\mathfrak m}_k \colon\ \mathscr C[1](c_0,c_1) \, \widehat\otimes \cdots \widehat\otimes \, \mathscr C[1](c_{k-1},c_k)
\to \mathscr C[1](c_0,c_k),
\]
of degree $+1$ for $k=0,1,2,\dots$ and $c_i \in \mathfrak{Ob}(\mathscr C)$.
(Note that in the case when $k=0$ the domain is
0 if $c_0 \ne c_1$ and is
$\Lambda_0$ if $c_0 = c_1$.)

We call $\mathscr C[1](c_0,c_1)$ the {\it module of morphisms}
\index{module of morphisms}
and $\mathfrak m_k$ the {\it structure operations}.
\index{structure operations}

We assume the following three conditions:
\begin{enumerate}\itemsep=0pt
\item[(1)]
We require $\mathfrak m_k$ to satisfy the {\it $A_{\infty}$ formula}
\index{$A_{\infty}$ formula} \eqref{form22}
described below.
\item[(2)]
The operations $\mathfrak m_k$ preserves the filtration.\footnote{Actually
this condition follows automatically from $\Lambda_0$ linearity.} Namely,
\[
\mathfrak m_k\bigl(\mathfrak F^{\lambda}
\bigl(
 \mathscr C[1](c_0,c_1) \, \widehat\otimes \cdots \widehat\otimes\, \mathscr C[1](c_{k-1},c_k)
\bigr)\bigr)
\subseteq
\mathfrak F^{\lambda}
 (
\mathscr C[1](c_0,c_k) ).
\]
\item[(3)]
We have
$
\mathfrak m_0(1) \equiv 0 \mod T^{\varepsilon}$,
for some $\varepsilon > 0$.
\end{enumerate}
\end{defn}
To describe the $A_{\infty}$ formula, we introduce notations.
Let $a,b\in \mathfrak{Ob}(\mathscr C)$. We put
\begin{equation}\label{form24new}
B_k\mathscr C[1](a,b) := \underset{a=c_0,c_1,\dots,c_{k-1},c_k=b}{\widehat{\bigoplus}}
\mathscr C[1](c_0,c_1) \,\widehat\otimes \cdots \widehat\otimes\, \mathscr C[1](c_{k-1},c_k).
\end{equation}
\index[syindex]{BkC@$B_k\mathscr C[1](a,b)$}
(Here and hereafter $\widehat\oplus$ denotes the $T$-adic completion of
the direct sum.)
\index[syindex]{pzzlus@$\widehat\oplus$}

Note that in the case when $k=0$
\begin{equation}\label{form23}
B_0\mathscr C[1](c_0,c_1)
:=
\begin{cases}
0 &\text{if $c_0 \ne c_1$},
\\
\Lambda_0 &\text{if $c_0 = c_1$}.
\end{cases}
\end{equation}
We denote
\[
B\mathscr C[1](a,b) = \underset{k=0,1,2,\dots}{\widehat{\bigoplus}} B_k\mathscr C[1](a,b),
\qquad
B\mathscr C[1] := \widehat{\underset{a,b}{{\bigoplus}}}B\mathscr C[1](a,b).
\]
We define a homomorphism
\begin{equation*}
\Delta \colon\ B_k\mathscr C[1](a,b) \to \underset{k_1+k_2=k}{\widehat{\bigoplus}}
\underset{c}{\widehat{\bigoplus}}
B_{k_1}\mathscr C[1](a,c) \,\widehat\otimes\, B_{k_2}\mathscr C[1](c,b)
\end{equation*}
\index[syindex]{Delta@$\Delta$} by
\[
\Delta(x_1\otimes\cdots\otimes x_k) := \sum_{k_1=0}^k(x_1 \otimes\cdots\otimes x_{k_1})
\otimes
(x_{k_1+1} \otimes\cdots\otimes x_{k}).
\]

It induces maps
\[
\Delta \colon\ B_k\mathscr C[1] \to \underset{k_1+k_2=k,
\atop k_1=0,\dots,k}{\widehat{\bigoplus}}
B_{k_1}\mathscr C[1]\,\widehat\otimes\, B_{k_2}\mathscr C[1],\qquad k=0,1,2,\dots,
\]
and
$\Delta \colon B\mathscr C[1] \to B\mathscr C[1]\,\widehat\otimes\, B\mathscr C[1]$.
Then $(B\mathscr C[1](a,a),\Delta)$ and $(B\mathscr C[1],\Delta)$ are graded formal
coalgebras.\footnote{The coalgebra structure is defined by a map
$\Delta \colon C \to C\otimes C$. Here the target of our $\Delta$ is
the completion $C \,\widehat\otimes\, C$. In such a case it is called a {\it
formal coalgebra}. Such a notion appears in the theory of formal groups.}
\index{formal coalgebra}

Operations $\mathfrak m_k$ define homomorphisms:
$ B_k\mathscr C[1](a,b) \to \mathscr C[1](a,b)$.
It can be extended uniquely to coderivations
$
\hat d_k \colon\ B\mathscr C[1] \to B\mathscr C[1]$, $
\hat d_k \colon\ B\mathscr C[1](a,b) \to B\mathscr C[1](a,b)
$
by
\[
\hat d_k(x_1\otimes\cdots\otimes x_n)
:=
\sum_{\ell}
(-1)^{*} x_1 \otimes\cdots \otimes
{\mathfrak m}_k(x_{\ell},\dots,x_{\ell+k-1}) \otimes
\cdots \otimes x_n,
\]
where $* = (\deg x_1 +1) + \cdots +
(\deg x_{\ell-1} +1)$.
We put
\begin{equation}\label{hatddd}
\hat d := \sum_k \hat d_k.
\end{equation}
\index[syindex]{dhat@$\hat d$}
Now the $A_{\infty}$ formula is
\begin{equation}\label{form22}
\hat d\circ \hat d= 0.
\end{equation}
Note that \eqref{form22} is equivalent to the equality
\begin{equation}\label{formula25}
0=\sum_{k_1+k_2=k+1}\sum_{i=0}^{k_1-1}
(-1)^* \mathfrak m_{k_1}(x_1,\dots,x_i,\mathfrak m_{k_2}(x_{i+1},\dots,x_{k_2}),
\dots,x_k),
\end{equation}
where
$* = i +\sum_{j=1}^i \deg x_j$.

We use the notation
\index[syindex]{degprime@$\deg'$}
$
\deg' x := \deg x - 1
$
then $* = \sum _{j=1}^i \deg' x_j$.
\begin{rem}
The sign convention in \eqref{formula25} is the same as
\cite{fooobook} but is different from \cite{Se}.
It seems that two different conventions are related
to each other by the process to take opposite category
(see Definition~\ref{opcate}).
\end{rem}
\begin{rem}\label{rem23}
We can define the notion of a non-unital $A_{\infty}$ category
over a ring $R$ (which is not filtered) in the same way except the following:
\begin{enumerate}\itemsep=0pt
\item[(1)]
We do not require the structure operations $\mathfrak m_k$ to preserve
the filtration.
\item[(2)]
We require $\mathfrak m_0 = 0$. In other words, we require
strictness, in the sense of Definition~\ref{defn2333}\,(2).
\end{enumerate}
Note that item (2) is our convention.
At this stage this is only a matter of
convention. Namely, we may include the curved case
over (unfiltered) ring.
It may be natural to do so in the case when we study the situation where
structure operations are converging (in the Lagrangian Floer theory)
and the version over $\C$.
Also in the case of monotone Lagrangian submanifolds
with minimal Maslov number $2$ such a situation appears
naturally.

Since we required $\mathfrak m_0 \equiv 0 \mod T^{\varepsilon}$,
we include this condition.

There will appear more serious reasons related to item (2), as the story
goes on. See Remarks~\ref{nonfil1} and \ref{rem210}.
\end{rem}
\begin{defn}\label{defn2333}
Let $\mathscr C$ be a non-unital curved filtered $A_{\infty}$ category.
\begin{enumerate}\itemsep=0pt
\item[(1)]
We say $\mathscr C$ is {\it $G$-gapped}
\index{$G$-gapped} if all the operations $\mathfrak m_k$ are
$G$-gapped.
\item[(2)]
We say $\mathscr C$ is a
{\it non-unital filtered $A_{\infty}$ category} if $\mathfrak m_0 =0$.
We also say that $\mathscr C$ is {\it strict}
\index{strict} instead.
\item[(3)]
If $\mathscr C$ is $G$-gapped, we define
{\it $R$-reduction}
\index{$R$-reduction} $\overline{\mathscr C}$ of our filtered $A_{\infty}$
category as follows. It is an~$A_{\infty}$ category over $R$ in the sense of Remark~\ref{rem23}.
\begin{enumerate}\itemsep=0pt
\item[(a)]
$\mathfrak{OB}(\overline{\mathscr C}) = \mathfrak{OB}(\mathscr C)$.
\item[(b)]
For $c,c' \in \mathfrak{OB}(\overline{\mathscr C})$,
there is a free $R$ module $\overline{\mathscr C}(c,c')$
such that
$
{\mathscr C}(c,c')$
is a completion of $\overline{\mathscr C}(c,c') \otimes_R {\Lambda_0}$,
by the definition of a completed free $\Lambda_0$ module.
We take this $R$ module $\overline{\mathscr C}(c,c')$ as the module of morphisms
of $\overline{\mathscr C}$.
\item[(c)]
The structure morphisms $\overline{\mathfrak m}_k$ are the
$R$-reductions of $\mathfrak m_k$.
\end{enumerate}
Note that $\overline{\mathfrak m}_0 = 0$ by
Definition~\ref{defn22}\,(3).
Other conditions for $\overline{\mathscr C}$ to be an $A_{\infty}$ category
follow from the corresponding properties of ${\mathscr C}$.
\item[(4)]
We say $\mathscr C$ is {\it unital}\index{unital}
(or $\mathscr C$ is a curved filtered $A_{\infty}$ category) if there exists ${\bf e}_c \in \mathscr C^0(c,c)$
for each $c \in \mathfrak{Ob}(\mathscr C)$ such that the following holds:
\begin{enumerate}\itemsep=0pt
\item[(a)] If $x_1 \in \mathscr C(c,c')$, $x_2 \in \mathscr C(c',c)$ then
$
\mathfrak m_2(\text{\bf e}_c,x_1) = x_1$,
$ \mathfrak m_2(x_2,\text{\bf e}_c) = (-1)^{\deg x_2}x_2$.
\item[(b)]
If $k+\ell \ne 1$,
$x_1 \otimes \dots \otimes x_{\ell}
\in B_{\ell}\mathscr C[1](a,c)$,
$y_1 \otimes \dots \otimes y_{k}
\in B_k\mathscr C[1](c,b)$ then
\begin{equation}\label{formidenty}
\mathfrak m_{k+\ell+1}(x_1,\dots,x_{\ell},\text{\bf e}_c,y_1,\dots,y_{k}) = 0.
\end{equation}
\end{enumerate}
\item[(5)]
A {\it filtered $A_{\infty}$ algebra} is a
\index{filtered $A_{\infty}$ algebra}
non-unital curved filtered $A_{\infty}$ category with one object.
Its unitality and strictness is defined as its unitality and strictness
as a non-unital curved filtered $A_{\infty}$ category.
\item[(6)] Let $\mathfrak C = (C,\{\mathfrak m_k\})$ be an
$A_{\infty}$ algebra. We define $\widetilde{\mathfrak M}(C;\Lambda_+)$,
\index[syindex]{MClambdaplus@$\widetilde{\mathfrak M}(C;\Lambda_+)$}
the {\it Maurer--Cartan solution space}
\index{Maurer--Cartan solution space} of $C$, as the set of all
elements $b \in C^1$ such that
\begin{enumerate}\itemsep=0pt
\item
$b \equiv 0 \mod \Lambda_+$.\footnote{We study
the case when this condition is not satisfied sometimes and
define $\widetilde{\mathfrak M}(\mathfrak C;\Lambda_0)$.
In such a case, the equation \eqref{MCeq} is more delicate to define
since the left-hand side may not converge in $T$-adic topology.
We do not discuss this generalization in this paper. See, for example, \cite{fooo092}.}
\item
\begin{equation}\label{MCeq}
\sum_{k=0}^{\infty} \mathfrak m_k(b,\dots,b) = 0.
\end{equation}
We remark that the left-hand side is an infinite sum, which converges in
$T$-adic topology.\footnote{In case $\mathscr C$ is unital, we sometimes
study the weaker equation which replaces the right-hand side by $C {\bf e}$
for some $C \in \Lambda_+$. 
See \cite[Section 4.3]{fooobook}.}\,\footnote{We can define
an equivalence relation called gauge equivalence so that
the $A_{\infty}$ structure defined by the deformed operators $\mathfrak m_k^{b}$
depends only on the gauge equivalence class of $b$.
See \cite[Section 3.6.3]{fooobook}.}
An element of $\widetilde{\mathfrak M}(\mathfrak C;\Lambda_+)$ is called a
{\it bounding cochain}\index{bounding cochain}.
\end{enumerate}
\item[(7)]
Let $\mathscr C$ be a non-unital curved filtered $A_{\infty}$
category. We define a non-unital filtered $A_{\infty}$ category
$\mathscr C^s$ as follows:\index[syindex]{Cscrs@$\mathscr C^s$}
\begin{enumerate}\itemsep=0pt
\item
An object of $\mathscr C^s$ is a pair $(c,b)$, where $c \in
\mathfrak{OB}(\mathscr C)$
and $b \in \widetilde{\mathfrak M}(\mathscr C(c,c);\Lambda_+)$.
\item
If $(c,b),(c',b')$ are objects of $\mathscr C^s$, then
$
\mathscr C^s((c,b),(c',b')) = \mathscr C(c,c')
$
by definition.
\item
If $(c_i,b_i) \in \mathfrak{OB}(\mathscr C')$ for $i=0,\dots,k$
and $x_i \in \mathscr C^s((c_{i-1},b_{i-1}),(c_i,b_i)) = \mathscr C(c_{i-1},c_i)$
for~${i=1,\dots,k}$. Then we define the structure operations~\smash{$\mathfrak m_k^{(b_0,\dots,b_k)}$}
of $\mathscr C^s$ as follows:
\[
\mathfrak m_k^{(b_0,\dots,b_k)}(x_1,\dots,x_k)
= \sum_{\ell_0,\dots,\ell_k}
\mathfrak m_{k+\ell_0+\cdots+\ell_k}\bigl(b_0^{\ell_0},x_1,b_1^{\ell_1},\dots,
b_{k-1}^{\ell_{k-1}},x_k,b_k^{\ell_k}\bigr).
\]
\end{enumerate}
The proof of the fact that this formula defines a non-unital filtered $A_{\infty}$ category
is similar to the corresponding result in the case of an $A_{\infty}$ algebra,
which is proved as \cite[Proposition~3.6.10]{fooobook}.
We call $\mathscr C^s$ the {\it associated strict category}
\index{associated strict category} to $\mathscr C$.
If $\mathscr C$ is unital, then $\mathscr C^s$ is also unital.
\end{enumerate}
\end{defn}

\begin{rem}\label{nonfil1}
Note that \eqref{MCeq} does not make sense in the case
of an (unfiltered) $A_{\infty}$ category.
In fact, the left-hand side is an infinite sum.
(This is one reason why we assume strictness
for (unfiltered) $A_{\infty}$ category.)

There are several ways to go around this point.
We will not discuss it here.
\end{rem}
\begin{rem}
In Definition~\ref{defn2333}\,(4), we required strict unitality.
In \cite[Definition 3.3.2]{fooobook}, we defined the notion of
a homotopy unit of a filtered $A_{\infty}$ {\it algebra}.
We can define the notion of a~homotopy unit of a filtered $A_{\infty}$ {\it category}
in the same way.
We do not discuss it in this paper, since in our geometric application
we can construct a strict unit by using the de Rham model.
\end{rem}
\begin{defn}[Bondal and Kapranov \cite{bondkap}]
\label{defnDGcate}
An $A_{\infty}$ category is said to be a {\it differential graded category}
\index{differential graded category}
or a {\it DG-category}
\index{DG category} if $\mathfrak m_k = 0$ for $k\ne 1,2$.
\end{defn}
\begin{rem}\label{remDGcat}
In the usual definition of a DG-category,
the space of morphisms $\mathscr C(c_1,c_2)$ is a~chain complex
with boundary operator $d$ and the composition map
\[
\circ \colon\ \mathscr C(c_1,c_2) \otimes \mathscr C(c_2,c_3)
\to \mathscr C(c_1,c_3)
\]
is assumed to be a chain map and the compositions are assumed to be
associative (strictly). We change the sign and define
$
\mathfrak m_1(x) := (-1)^{\deg x+1} d(x)$,
$
\mathfrak m_2(x,y) := (-1)^{\deg x(\deg y + 1)}y\circ x$.
Then it satisfies $A_{\infty}$ relation \eqref{formula25}.
(See \cite[Example--Lemma 1.7]{fu4}.)
There is an alternative choice of the sign,
that is,
$
\mathfrak m_1(x) := d(x)$,
$
\mathfrak m_2(x,y) := (-1)^{\deg x}x\circ y$.
This is the choice in
\cite[Definition~21.21\,(2)(3)]{fooonewbook}.
\end{rem}

\subsection[$A_\infty$ functor]{$\boldsymbol{A_{\infty}}$ functor}
\label{subsec:Ainffun}
\begin{defn}\label{defn224}
Let $\mathscr C_i$, $i=1,2$, be non-unital curved filtered $A_{\infty}$ categories.
A {\it filtered $A_{\infty}$ functor} $\mathscr F \colon \mathscr C_1 \to \mathscr C_2$
\index{filtered $A_{\infty}$ functor} is a collection of $\mathscr F_{\rm ob}$, $\mathscr F_k$, $k = 0,1,2,\dots$, such that
\begin{enumerate}\itemsep=0pt
\item[(1)] We are given a set theoretical map $\mathscr F_{\rm ob} \colon \mathfrak{Ob}(\mathscr C_1)
\to \mathfrak{Ob}(\mathscr C_2)$, which we call the
{\it object part} of $\mathscr F$.

\item[(2)] For $c_1,c_2 \in \mathfrak{Ob}(\mathscr C_1)$,
$\mathscr F_k(c_1,c_2) \colon B_k\mathscr C_1[1](c_1,c_2) \to \mathscr C_2[1](\mathscr F_{\rm ob}(c_1),\mathscr F_{\rm ob}(c_2))$
is a $\Lambda_0$ module homomorphism of degree 0.
It preserves filtration in a similar sense as Definition~\ref{defn22}\,(2).
We write $\mathscr F_k$ in place of $\mathscr F_k(c_1,c_2)$ sometimes.
\item[(3)] We require that $\mathscr F_{0} \equiv 0 \mod T^{\varepsilon}$, $\varepsilon > 0$.
Note that $\mathscr F_{0}$ consists of maps
$\mathscr F_{0}(c) \colon
\Lambda_0 \to \mathscr C_2[1](\mathscr F_{\rm ob}(c),\mathscr F_{\rm ob}(c))$
for each $c \in \mathfrak{OB}(\mathscr C_1)$.
\item[(4)] We extend $\mathscr F_k(c_1,c_2)$ to a formal
coalgebra homomorphism
\[
\widehat{\mathscr F}(c_1,c_2) \colon\ B\mathscr C_1[1](c_1,c_2)
\to B\mathscr C_2[1](\mathscr F_{\rm ob}(c_1),\mathscr F_{\rm ob}(c_2)).
\]
Then $\widehat{\mathscr F}(c_1,c_2)$ is a chain map with respect to the boundary operator
$\hat d$ in \eqref{hatddd}.
\end{enumerate}
\end{defn}
\begin{rem}
In Definition~\ref{defn224}, we include the case $\mathscr F_0 \ne 0$,
that is, a `curved' filtered~$A_{\infty}$ functor.
(In \cite{fu4}, we did not include it. However,
the definition of filtered $A_{\infty}$ {\it algebra homomorphism}
\index{filtered $A_{\infty}$ algebra homomorphism} in \cite[Definition 3.2.29]{fooobook}
includes the case $\mathfrak f_0 \ne 0$.)

The map $\widehat{\mathscr F}$ on
$B_k\mathscr C_1[1](c_1,c_2)$
is defined by
\begin{equation}\label{form29}
\widehat{\mathscr F}(x_1,\dots,x_k) :=
\sum_{\ell=1}^{\infty}\sum_{k_1,\dots,k_{\ell} \atop
k_1 +\dots +k_{\ell} = k}
{\mathscr F}_{k_1}(x_1,\dots,x_{k_1})
\otimes \dots \otimes
{\mathscr F}_{k_{\ell}}(x_{k-k_{\ell}+1},\dots,x_{k}),
\end{equation}
for $k\ge 1$. For $k=0$, it is
\[
\widehat{\mathscr F}(1)
:= 1 + \sum_{\ell=1}^{\infty} {\mathscr F}_{0}(1)^{\otimes \ell},
\]
$\widehat{\mathscr F}$ is a formal coalgebra homomorphism.

\end{rem}
\begin{rem}\label{rem210}
We define an $A_{\infty}$ functor between (unfiltered) $A_{\infty}$
categories in the same way.
We require ${\mathscr F}_0 = 0$ in the unfiltered situation.
There is more serious reason to require it compared to Remark~\ref{rem23}\,(2).
We remark that in our situation where ${\mathscr F}_0 \ne 0$,
the right-hand side of~\eqref{form29} is an {\it infinite} sum.
It converges in $T$-adic topology thanks to Definition~\ref{defn224}\,(3).
In the case when we work over the ground ring, the unfiltered case,
the right-hand side of~\eqref{form29} should be a finite sum.
\end{rem}

\begin{defn}
Let $\mathscr F \colon \mathscr C_1 \to \mathscr C_2$
be a filtered $A_{\infty}$ functor between
non-unital curved filtered~$A_{\infty}$ categories.
\begin{enumerate}\itemsep=0pt
\item[(1)]
We say $\mathscr F$ is {\it strict}
\index{strict} if $\mathscr F_0 =0$.
\item[(2)]
Suppose $\mathscr C_1$, $\mathscr C_2$ are $G$-gapped.
We say $\mathscr F$ is $G$-gapped if all
\index{$G$-gapped}
the maps $\mathscr F_k$ are $G$-gapped for $k=0,1,2, \dots$.
\item[(3)]
A $G$-gapped filtered $A_{\infty}$
functor between
non-unital curved filtered $A_{\infty}$ categories
induce an $A_{\infty}$
functor between their $R$-reductions.
\item[(4)]
Suppose $\mathscr C_1$ and $\mathscr C_2$ are
unital in addition.
We say $\mathscr F$ is (strictly) {\it unital} if
\index{unital}
the following two conditions are satisfied:
\begin{enumerate}\itemsep=0pt
\item
$\mathscr F_{1}({\bf e}_c)
= {\bf e}_{\mathscr F_{\rm ob}(c)}$.
\item
$
\mathscr F_{k+\ell+1}(x_1,\dots,x_{\ell},{\bf e}_c,y_1,\dots,y_{\ell})
= 0
$
for $k+\ell >0$.
\end{enumerate}
\item[(5)]
If $\mathscr F \colon \mathscr C_1 \to \mathscr C_2$
is a filtered $A_{\infty}$ functor between
non-unital curved filtered $A_{\infty}$ categories,
then we obtain a strict
filtered $A_{\infty}$ functor
$\mathscr F^s \colon \mathscr C^s_1 \to \mathscr C^s_2$
between their associated strict categories as follows.
\begin{enumerate}\itemsep=0pt
\item
Let $c \in \mathfrak{OB}(\mathscr C_1)$ and
$b \in \widetilde{\mathfrak M}(\mathscr C(c,c);\Lambda_+)$.
We put
\[
\mathscr F_*(b)
:= \sum_{k=0}^{\infty} \mathscr F_k(b,\dots,b).
\]
We can prove $\mathscr F_*(b) \in \widetilde{\mathfrak M}
(\mathscr C(\mathscr F_{\rm ob}(c),\mathscr F_{\rm ob}(c));\Lambda_+)$.
We define
\[
\mathscr F^s_{\rm ob}(c,b) := (\mathscr F_{\rm ob}(c),\mathscr F_*(b)).
\]
\item
Let $(c_i,b_i) \in \mathfrak{OB}(\mathscr C'_1)$,
$i=0,\dots,k$, and
$x_i \in \mathscr C^s_1((c_{i-1},b_{i-1}),(c_i,b_i))
= \mathscr C_1(c_{i-1},c_i)$, $i=1,\dots,k$.
We put\index[syindex]{Fsk@$\mathscr F^s_{k}$}
\[
\mathscr F^s_{k}(x_1,\dots,x_k)
:= \sum_{\ell_0,\dots,\ell_k}
\mathscr F_{k+\ell_0+\cdots+\ell_k}\big(b_0^{\ell_0},x_1,b_1^{\ell_1},\dots,
b_{k-1}^{\ell_{k-1}},x_k,b_k^{\ell_k}\big).
\]
\end{enumerate}
We also put $\mathscr F^s_{0} = 0$.
We can show that $\mathscr F^s_{\rm ob}$ and $\mathscr F^s_{k}$ define a
strict filtered $A_{\infty}$ functor
$\mathscr F^s \colon \mathscr C^s_1 \to \mathscr C^s_2$,
in the same way as \cite[Lemma 3.6.36, Definition--Lemma 5.2.15, Lemma~5.2.16]{fooobook}.
(They discuss the case of $A_{\infty}$ algebra.)
We say $\mathscr F^s$ is the {\it associated strict functor}
\index{associated strict functor} to $\mathscr F$.
If $\mathscr F$ is unital (resp.\ $G$-gapped), then so is $\mathscr F^s$.
\item[(6)]
The {\it identity functor} \index{identity functor} \index[syindex]{ID@$\mathscr {ID}$}
$\mathscr {ID} \colon \mathscr C \to \mathscr C$ is
defined by
\begin{enumerate}\itemsep=0pt
\item
$\mathscr {ID}_{\rm ob} =$ the identity map: $\mathfrak{OB}(\mathscr C) \to \mathfrak{OB}(\mathscr C)$.
\item
$\mathscr {ID}_{1}(c_1,c_2) \colon \mathscr C(c_1,c_2) \to \mathscr C(c_1,c_2)$
is the identity map.
\item
$\mathscr {ID}_{k} = 0$ for $k\ne {\rm ob}, 1$.
\end{enumerate}
$\mathscr {ID}$ is unital (resp.\ $G$-gapped)
if so is $\mathscr C$.
\end{enumerate}
\end{defn}
\begin{defnlem}
Let $\mathscr F^1 \colon \mathscr C_1 \to \mathscr C_2$,
$\mathscr F^2 \colon \mathscr C_2 \to \mathscr C_3$
be filtered $A_{\infty}$ functors.
\begin{enumerate}\itemsep=0pt
\item[(1)]
We define their {\it composition}
\index{composition}
$\mathscr F = \mathscr F^2 \circ \mathscr F^1 \colon \mathscr C_1 \to \mathscr C_3$
as follows:
\begin{gather*}
\mathscr F_{\rm ob} = \mathscr F^2_{\rm ob} \circ \mathscr F^1_{\rm ob},\\
\widehat{\mathscr F}(c_1,c_2) = \widehat{ \mathscr F}^2\bigl(\mathscr F^1_{\rm ob}(c_1),\mathscr F^1_{\rm ob}(c_2)\bigr)
\circ\widehat{\mathscr F^1}(c_1,c_2)\colon\\
\hphantom{\widehat{\mathscr F}(c_1,c_2) =} B\mathscr C_1(c_1,c_2) \to B\mathscr C_3
(\mathscr F_{\rm ob}(c_1),\mathscr F_{\rm ob}(c_2)).
\end{gather*}
\item[(2)]
If $\mathscr F^1$, $\mathscr F^2$ are strict (resp.\ unital, $G$-gapped), then
$\mathscr F = \mathscr F^2 \circ \mathscr F^1$ is strict (resp.\ unital, $G$-gapped).
\item[(3)]
If $\mathscr F^{1 s}$, $\mathscr F^{2 s}$
are strict functors associated to $\mathscr F^1$, $\mathscr F^2$,
respectively, then
$\mathscr F^{1 s} \circ \mathscr F^{2 s}$
is the strict functor associated to
$\mathscr F^1 \circ \mathscr F^2$.
\item[(4)]
$\mathscr F \circ \mathscr{ID} = \mathscr{ID} \circ \mathscr F
= \mathscr F$.
\end{enumerate}
\end{defnlem}
The proof is easy and is omitted.
\subsection{Functor category}
\label{subsec:funcat}

\begin{defn}[{\cite[Definition 7.49]{fu4}}]\label{defn215555}
Let $\mathscr F, \mathscr G \colon \mathscr C_1 \to \mathscr C_2$
be two curved filtered $A_{\infty}$ functors between
non-unital curved filtered $A_{\infty}$ categories.

A {\it pre-natural transformation}
\index{pre-natural transformation} from $\mathscr F$ to $\mathscr G$
of degree $d$
is a family of operators $\mathcal T = \{\mathcal T_k(a,b)\}$
\[
\mathcal T_k(a,b) \colon\ B_k\mathscr C_1[1](a,b)
\to \mathscr C_2[1](\mathscr F_{\rm ob}(a),\mathscr G_{\rm ob}(b))
\]
of degree $d$ for $k=0,1,2,\dots$ and
$a,b \in \mathfrak{OB}(\mathscr C_1)$, which preserves filtration in the same sense as
Definition~\ref{defn22}\,(2).\footnote{It means that
$\deg' \mathcal T_k(a,b)({\bf x}) = \deg'{\bf x} + d$, where
$\deg'x_1\otimes\dots\otimes x_k = \sum \deg' x_i$.}
We require that the image of $\mathcal T_0$ has strictly positive energy.

We write $\mathfrak{deg}\,{\mathcal T} := d+1$
\index[syindex]{deg@$\mathfrak{deg}$}
and $\mathfrak{deg}' := \mathfrak{deg} - 1 =d$.\footnote{This convention is different from \cite[two lines below equation~(7.44)]{fu4}.
Actually in \cite[equation~(7.12.2)]{fu4} $\deg'$ is defined as $\deg + 1$,
which is different from our convention here.
The convention of this paper coincides with \cite[equation~(3.2.2)]{fooobook},
which seems better than one in \cite{fu4}.
(See a line after Definition~\ref{defnidentittrans}, for example.)
This inconsistency does not affect to the calculation of the sign in \cite{fu4}
which we use much in this paper.}
We say that $\mathcal T$ is $G$-gapped if each of $\mathcal T_k(a,b)$ is $G$-gapped.
We denote by $\mathcal{FUNC}(\mathscr F,\mathscr G)$ the set of all
\index[syindex]{FUNCF@$\mathcal{FUNC}(\mathscr F,\mathscr G)$}
pre-natural transformations from $\mathscr F$ to $\mathscr G$.
It is a~completed free $\Lambda_{0}$ module and is graded.
We denote by $\mathcal{FUNC}^d(\mathscr F,\mathscr G)$ the degree $d$ part.
In other words, if $\mathcal T \in \mathcal{FUNC}^d(\mathscr F,\mathscr G)$,
then $\mathfrak{deg}\,{\mathcal T} = d+1$ and $\mathfrak{deg}'{\mathcal T} = d$.
\end{defn}
\begin{rem}
We remark that $\mathcal T_{0}(a,b) = 0$
if $a \ne b$ and $\mathcal T_{0}(a,a)$ is a
$\Lambda_0$ module homomorphism
\[
\mathcal T_{0}(a,a) \colon\
B_0\mathscr C_1[1](a,a) = \Lambda_0
\to \mathscr C_2[1](\mathscr F_{\rm ob}(a),\mathscr G_{\rm ob}(a)).
\]
We denote by $\mathcal T_{0}(a) \in
\mathscr C_2[1](\mathscr F_{\rm ob}(a),\mathscr G_{\rm ob}(a))$
the element $\mathcal T_{0}(a,a)(1)$.
\end{rem}
For $a',b' \in \mathfrak{Ob}(\mathscr
C_2)$, let
$
\pi_{a',b'} \colon B\mathscr C_2[1](a',b') \to \mathscr C_2[1](a',b')
$
be the projection.
\begin{lemdef}\label{defn280}\quad
\begin{enumerate}\itemsep=0pt
\item[$(1)$]
For each $\mathcal T =\{\mathcal T_k(a,b)\}
\in \mathcal{FUNC}^d(\mathscr F,\mathscr G)$, there exists uniquely a family
\[
\widehat{\mathcal T}(a,b)\colon\ B\mathscr C_1[1](a,b) \to B\mathscr C_2[1](\mathscr F_{\rm ob}(a),\mathscr G_{\rm ob}(b)),
\]
of $\Lambda_0$ module homomorphisms with the following properties:
\begin{gather}
\pi_{\mathscr F_{\rm ob}(a),\mathscr G_{\rm ob}(b)}\circ\widehat{\mathcal T}(a,b)
= \mathcal T_k(a,b) \qquad \text{on}\  B_k\mathscr C_1[1](a,b), \nonumber\\
\Delta\circ\widehat{\mathcal T}(a,b) = \sum_c\bigl(\widehat{\mathscr F} \otimes_s \widehat{\mathcal T}(c,b)
+ \widehat{\mathcal T}(a,c) \otimes_s \widehat{\mathscr G}\bigr) \circ \Delta.\label{form2626}
\end{gather}
Here $\otimes_s$ is defined by $(A \otimes_s B)(\text{\bf x},\text{\bf y})
= (-1)^{\deg B\deg' \text{\bf x}}A(\text{\bf x}) \otimes B(\text{\bf y})$.
\item[$(2)$]
There exists
$\delta \mathcal T = \{(\delta \mathcal T)_k(a,b)\}
\in \mathcal{FUNC}^{\mathfrak{deg}\,{\mathcal T}+1}(\mathscr F,\mathscr G)$
uniquely such that
\[
\widehat{\delta \mathcal T} = \hat d \circ \widehat{\mathcal T} + (-1)^{\mathfrak{deg}\,{\mathcal T}+1} \widehat{\mathcal T}\circ\hat d.
\]
\item[$(3)$]
$\delta (\delta\mathcal T) =0$.
\item[$(4)$]
A pre-natural transformation $\mathcal T$ is said to be a
{\it natural transformation}
\index{natural transformation} if $\delta\mathcal T = 0$.
\end{enumerate}
\end{lemdef}
(1) is \cite[Lemma 7.45]{fu4}. (2) is \cite[Lemma 7.48]{fu4}.
(3) is \cite[Corollary 7.50]{fu4}.

\begin{defn}\label{defn2929}
Let $\mathscr F^{(i)}\colon \mathscr C_1 \to \mathscr C_2$,
$i=0,\dots,k$, be curved filtered $A_{\infty}$ functors between
non-unital curved $A_{\infty}$ categories
and
$\mathcal T^{(i)} \in \mathcal{FUNC}^{d_i}\bigl(\mathscr F^{(i-1)},\mathscr F^{(i)}\bigr)$
for $i=1,\dots,k$
(here $k=1,2,\dots$).
We define
$\mathfrak m_k\big(\mathcal T^{(1)},\dots,\mathcal T^{(k)}\big) = \mathcal T
\in \mathcal{FUNC}^{d}\bigl(\mathscr F^{(0)},\mathscr F^{(k)}\bigr)$
\index[syindex]{m2kT1@$\mathfrak m_k\big(\mathcal T^{(1)},\dots,\mathcal T^{(k)}\big)$}
as follows
($d = d_1 + \dots + d_k +1$).

If $k=1$, then $\mathfrak m_1\big(\mathcal T^{(1)}\big) = -\delta \big(\mathcal T^{(1)}\big)$,
where $\delta$ is as in Definition~\ref{defn280}\,(2).

Suppose $k\ge 2$.
Let
$\text{\bf x} \in B(\mathscr C_1[1])$. We consider
\[
\Delta^{2k}\text{\bf x} = \sum_a \text{\bf x}_a^{(1)} \otimes
\cdots \otimes \text{\bf x}_a^{(2k+1)}.
\]
Here \smash{$\Delta^{m} \colon B(\mathscr C_1[1]) \to \underbrace{B(\mathscr C_1[1])
\otimes \dots \otimes B(\mathscr C_1[1])}_{m+1}$}
is defined inductively by $\Delta^{m} := (\Delta \otimes {\rm id}) \circ \Delta^{m-1}$, $\Delta^{1} = \Delta$.
We put
\[
{\mathcal T}(\text{\bf x}) := -\sum_a (-1)^{*_a}
\mathfrak m\bigl(
\widehat{\mathscr F^{(0)}}\bigl(\text{\bf x}_a^{(1)}\bigr),\mathcal T^{(1)}
\bigl(\text{\bf x}_a^{(2)}\bigr),
\dots,
\mathcal T^{(k)}\bigl(\text{\bf x}_a^{(2k)}\bigr),\widehat{\mathscr F^{(k)}}\bigl(\text{\bf x}_a^{(2k+1)}\bigr)
\bigr),
\]
where
\smash{$
*_a = \sum_{j=1}^k\sum_{i=1}^{2j-1}
d_j\deg'\text{\bf x}_a^{(i)}.
$}
Note that
\[
\deg' {\mathcal T}(\text{\bf x})
=
\sum_{i=1}^{k+1}\deg' \text{\bf x}_a^{(i)} +
\sum_{i=1}^{k} \mathfrak{deg}'\mathcal T^{(i)}
+ 1.
\]
Therefore,
$
\mathfrak{deg}'{\mathcal T} = \sum_{i=1}^{k} \mathfrak{deg}'\mathcal T^{(i)} + 1.
$
This is consistent with Definition~\ref{defn22}.

We consider the case when $k=0$.
We will define $\mathfrak m_0^{\mathscr{F}}(1)
\in \mathcal{FUNC}(\mathscr{F},\mathscr{F})$
(the $\mathfrak m_0$ operator of the functor category).
For $c \in \mathfrak{OB}(\mathscr C_2)$, the $\mathfrak m_0$ operator of $\mathscr C_2$
determine an element
$\mathfrak m_0(1)_c \in \mathscr C[1](\mathscr{F}_{\rm ob}(c),\mathscr{F}_{\rm ob}(c))$.
We put
$
\mathfrak m_0^{\mathscr{F}}(1)_c := -\mathfrak m_0(1)_c$.

\end{defn}
\begin{thm}\label{th21010}
Let $\mathscr C_1$, $\mathscr C_2$ be curved filtered $A_{\infty}$ categories.
Then,
there exists a non-unital curved filtered $A_{\infty}$ category
$\mathcal{FUNCC}(\mathscr C_1,\mathscr C_2)$ \index[syindex]{FUNCC@$\mathcal{FUNCC}(\mathscr C_1,\mathscr C_2)$}
\index[syindex]{FUNC@$\mathcal{FUNC}(\mathscr C_1,\mathscr C_2)$} such that
\begin{enumerate}\itemsep=0pt
\item[$(1)$]
The set of its objects
$\mathfrak{OB}(\mathcal{FUNCC}(\mathscr C_1,\mathscr C_2))$
consists of filtered $A_{\infty}$ functors
$\mathscr F \colon \mathscr C_1 \to \mathscr C_2$.
\item[$(2)$]
For $\mathscr F, \mathscr G \in \mathfrak{OB}(\mathcal{FUNCC}(\mathscr C_1,\mathscr C_2))$,
$\mathcal{FUNCC}(\mathscr F, \mathscr G)$ is the module of morphisms from $\mathscr F$ to~$\mathscr G$.
\item[$(3)$]
The structure operations
\[
\mathfrak m_k \colon\ B_k\mathcal{FUNCC}(\mathscr F,\mathscr G) \to \mathcal{FUNCC}(\mathscr F,\mathscr G)
\]
are as in Definition~{\rm\ref{defn2929}}.
\end{enumerate}
We denote by $\mathcal{FUNC}(\mathscr C_1,\mathscr C_2)$
the full subcategory of $\mathcal{FUNCC}(\mathscr C_1,\mathscr C_2)$
the set of whose objects are strict filtered $A_{\infty}$ functors.

If $\mathscr C_2$ is strict, then $\mathcal{FUNCC}(\mathscr C_1,\mathscr C_2)$
and $\mathcal{FUNC}(\mathscr C_1,\mathscr C_2)$ are strict.

In case $\mathscr C_1$, $\mathscr C_2$ are unital and/or strict,
we consider only unital and/or strict filtered $A_{\infty}$
functors as objects of $\mathcal{FUNC}(\mathscr C_1,\mathscr C_2)$.
In that way, we obtain strict and unital filtered $A_{\infty}$ category.

\end{thm}
This is \cite[Theorem--Definition 7.55]{fu4}.
(Note that only the strict case is proved in
\cite[Theorem--Definition 7.55]{fu4}. However, the proof there can be
applied without change in our case.
The functor category in the curved case is also studied in \cite[Section~3.4]{DL}.)
We call $\mathcal{FUNCC}(\mathscr C_1,\mathscr C_2)$, $\mathcal{FUNC}(\mathscr C_1,\mathscr C_2)$ the {\it functor category}.
\index{functor category}
\begin{prop}\label{prop211}
A strict $A_{\infty}$ functor $\mathscr F \colon \mathscr C_1 \to \mathscr C_2$ induces
strict $A_{\infty}$ functors
$\mathscr F_* \colon \mathcal{FUNC}(\mathscr C,\allowbreak\mathscr C_1) \to \mathcal{FUNC}(\mathscr C,\mathscr C_2)$,
$\mathscr F^* \colon \mathcal{FUNC}(\mathscr C_2,\mathscr C) \to \mathcal{FUNC}(\mathscr C_1,\mathscr C)$
such that
$(\mathscr F_*)_{\rm ob}(\mathscr G) = \mathscr F\circ\mathscr G$,
$(\mathscr F^*)_{\rm ob}(\mathscr G) = \mathscr G\circ\mathscr F$.
The same is true if we replace $\mathcal{FUNC}$
by $\mathcal{FUNCC}$.
$($In that case, we do not need to assume~$\mathscr F$ to be strict.$)$

\end{prop}
This is \cite[Proposition--Definition 8.41]{fu4}.
\begin{defn}
In the situation of Theorem~\ref{th21010}, we assume that
$\mathscr C_1$, $\mathscr C_2$ are $G$-gapped.
We define a $G$-gapped filtered $A_{\infty}$ category
$\mathcal{FUNC}^G(\mathscr C_1,\mathscr C_2)$ as follows.
Its object is a $G$-gapped
filtered $A_{\infty}$ functors
$\mathscr F \colon \mathscr C_1 \to \mathscr C_2$.
The morphisms and the structure maps are the same as
Theorem~\ref{th21010}\,(2)(3).
It is easy to see that the structure maps are $G$-gapped.
\index[syindex]{FUNCG@$\mathcal{FUNC}^G(\mathscr C_1,\mathscr C_2)$}
The $G$-gapped version of Proposition~\ref{prop211} holds.
We may replace $\mathcal{FUNC}$
by $\mathcal{FUNCC}$.
\end{defn}

Hereafter, in the case of $G$-gapped category, we omit $G$ and
write $\mathcal{FUNC}(\mathscr C_1,\mathscr C_2)$
in place of~$\mathcal{FUNC}^G(\mathscr C_1,\mathscr C_2)$.

\begin{defn}[{\cite[Definition~8.2]{fu4}}]
\label{defnidentittrans}
Let $\mathscr F \colon \mathscr C_1 \to \mathscr C_2$
be a curved filtered $A_{\infty}$ functors between
non-unital curved filtered $A_{\infty}$ categories.
We assume $\mathscr C_2$ is unital in addition.
We define the {\it identity natural transformation}
\index{identity natural transformation}
$\operatorname{Id}^{\mathscr F}$ as follows.
\index[syindex]{idF@$\operatorname{Id}^{\mathscr F}$}
Let $\text{\bf e}_c \in \mathscr C_2^0(c,c)$ be the unit
in $\mathscr C_2$. We put
\[
\operatorname{Id}^{\mathscr F}_{\rm ob}(a)
= -\text{\bf e}_{\mathscr F_{\rm ob}(a)}
\in \mathscr C_2^{pb}(\mathscr F_0(a),\mathscr F_0(a)),\qquad
\operatorname{Id}^{\mathscr F}_k= 0 \qquad
\text{for}\  k \ge 1.
\]
\end{defn}

Note that \smash{$\deg' \text{\bf e}_{\mathscr F_{\rm ob}(a)}
= \deg \text{\bf e}_{\mathscr F_{\rm ob}(a)} - 1 = - 1$}.
Therefore, $\mathfrak{deg}\operatorname{Id}^{\mathscr F} = 0$.

It is easy to see from definition that $\operatorname{Id}^{\mathscr F}$
satisfies \eqref{formidenty} for the structure
operations $\mathfrak m_k$
of~$\mathcal{FUNCC}(\mathscr C_1,\mathscr C_2)$.
Therefore, we have the following.

\begin{lem}\label{lem214}
If $\mathscr C_2$ is unital, then $\mathcal{FUNCC}(\mathscr C_1,\mathscr C_2)$
is unital.
\end{lem}

\subsection[$A_\infty$-Whitehead theorem]{$\boldsymbol{A_{\infty}}$-Whitehead theorem}
\label{subsec:Whitehead}

In this subsection, all filtered $A_{\infty}$ categories are assumed to be
strict (except in Remark~\ref{rem2230}).

\begin{defn}\label{6.22}
Let $\mathscr C$ be a strict filtered $A_{\infty}$ category and $c,c' \in \mathfrak{Ob}(\mathscr C)$.
Let $x \in \mathscr C^0(c,c')$. We say that $x$ is a {\it homotopy equivalence}
\index{homotopy equivalence}
if there exists $y \in \mathscr C^0(c',c)$ such that
\begin{enumerate}\itemsep=0pt
\item[(1)] $\mathfrak m_1(x) =\mathfrak m_1(y) = 0$,
\item[(2)] $\mathfrak m_2(y,x) - \text{\bf e}_c \in
\operatorname{Im}\mathfrak m_1$, $\mathfrak m_2(x,y) - \text{\bf e}_{c'} \in
\operatorname{Im}\mathfrak m_1$.
\end{enumerate}
Two objects $c,c' \in \mathfrak{Ob}(\mathscr C)$ are said to be
{\it homotopy equivalent}
\index{homotopy equivalent} to each other
if there exists a~homotopy equivalence
between them.

Homotopy equivalence is an equivalence relation by
\cite[Lemma 6.24]{fu4}.
\end{defn}

\begin{defn}\label{defn216}
Suppose that we are in the situation of Lemma~\ref{lem214} and we assume that $\mathscr C_2$
is strict.
Two strict filtered $A_{\infty}$ functors $\mathscr F,\mathscr G \colon \mathscr C_1 \to \mathscr C_2$
are said to be {\it homotopy equivalent}
\index{homotopy equivalent} to each other
if they are homotopy equivalent as
objects of $\mathcal{FUNC}(\mathscr C_1,\mathscr C_2)$ in the sense of
Definition~\ref{6.22}.
(Note that $\mathcal{FUNC}(\mathscr C_1,\mathscr C_2)$ is strict if $\mathscr C_2$ is strict.)
The homotopy equivalence among strict filtered $A_{\infty}$ functors is an equivalence relation.
We can define the notion two $G$-gapped strict filtered $A_{\infty}$ functors
to be homotopy equivalent (as $G$-gapped strict filtered $A_{\infty}$ functors)
in a similar way.
\end{defn}

\begin{rem}
We consider the case when $\mathscr C_1$, $\mathscr C_2$ have
only one object. In this case, curved
filtered $A_{\infty}$ functors $\mathscr F,\mathscr G \colon \mathscr C_1 \to \mathscr C_2$
are nothing but filtered $A_{\infty}$ homomorphisms.
The notion two (curved) filtered $A_{\infty}$ homomorphisms to be homotopic
is defined in \cite[Definition~4.2.35]{fooobook}.
We will define its category version in Definition~\ref{defn135}.
To distinguish one, we defined here from one in
Definition~\ref{defn135} we will use the
terminology `homotopy equivalent' in place of `homotopic'
in Definition~\ref{defn216}.\footnote{This notation
is different from \cite{fu4} at this point.}
We will prove in Section~\ref{sec:homotopyfafunc} that `homotopic'
implies 'homotopy equivalent' (see Proposition~\ref{prop1313}).
The converse is not correct (see Example~\ref{example1315}).
\end{rem}
\begin{defn}\label{defn225}
Let $\mathscr C_1$, $\mathscr C_2$ be strict filtered $A_{\infty}$ categories.
We assume that they are unital.
A strict $A_{\infty}$ functor $\mathscr F \colon \mathscr C_1 \to \mathscr C_2$ is
said to be a {\it homotopy equivalence}
\index{homotopy equivalence} if there exists
a strict $A_{\infty}$ functor $\mathscr G \colon \mathscr C_2 \to \mathscr C_1$
such that the composition $\mathscr F\circ \mathscr G$ is homotopy equivalent
to the identity functor $\mathscr{ID}^{\mathscr C_2}$ and that
$\mathscr G\circ \mathscr F$ is homotopy equivalent
to the identity functor $\mathscr{ID}^{\mathscr C_1}$.
We say $\mathscr G$ a {\it homotopy inverse}
\index{homotopy inverse} to $\mathscr F$.
Two strict
$A_{\infty}$ categories are said to be {\it homotopy equivalent}
\index{homotopy equivalent}
to each other if
there exists a~homotopy equivalence between them.
\end{defn}

We assume that the ground ring $R$ is a field in the next theorem.
\begin{thm}\label{white}
Let $\mathscr C_1$, $\mathscr C_2$ be filtered $A_{\infty}$ categories.
We assume they are unital, strict and gapped.
Let $\mathscr F \colon \mathscr C_1 \to \mathscr C_2$ be a strict
and gapped $A_{\infty}$
functor such that
\begin{enumerate}\itemsep=0pt
\item[$(1)$] ${\mathscr F}_1 \colon {\mathscr C}_1(c_1,c'_1)
\to {\mathscr C}_2(\mathscr F_{\rm ob}(c_1),
\mathscr F_{\rm ob}(c'_1))$ induces an isomorphism on
${\mathfrak m}_1$ homology.
\item[$(2)$] For any $c_2
\in \mathfrak{Ob}(\mathscr C'_2)$, there exists $c_1
\in \mathfrak{Ob}(\mathscr C'_1)$ such that $\mathscr F_{\rm ob}(c_1)$ is homotopy
equivalent to~$c_2$.
\end{enumerate}
Then $\mathscr F$ is a homotopy equivalence.

If $\mathscr C_1$, $\mathscr C_2$, $\mathscr F$ are $G$-gapped,
we may take homotopy inverse which is $G$-gapped also.
Moreover, homotopy equivalence in $(2)$ is taken to be $G$-gapped.
\end{thm}
The non-filtered version of this theorem is \cite[Theorem 8.6]{fu4}.
We can prove Theorem~\ref{white} in the same way.
\begin{rem}\label{rem2230}
Note that we assumed strictness of $\mathscr C$ here.
Actually Theorem~\ref{white}\,(1) does not make sense
in case $\mathfrak m_0 \ne 0$.
In a slightly different way, we can define homotopy equivalence of filtered
$A_{\infty}$ categories in the curved case
and Theorem~\ref{white} holds in that generality.
See Section~\ref{sec:homotopyfafunc} Theorem~\ref{white2}.
(We remark that the assumption~(1) of Theorem~\ref{white2} does make sense
in the curved case since $\overline{\mathscr C}$
is strict.
In the curved case, we replace (2) by the condition that
$\mathscr F_{\rm ob}$ is a bijection.)
\end{rem}

\subsection[$A_\infty$-Yoneda embedding]{$\boldsymbol{A_{\infty}}$-Yoneda embedding}
\label{subsec:Yoneda}

\begin{defn}[{\cite[Definition 7.8]{fu4}}]\label{opcate}
Let $\mathscr C$ be a non-unital curved filtered $A_{\infty}$ category. We define
its {\it opposite $A_{\infty}$ category}
\index{opposite $A_{\infty}$ category} $\mathscr C^{\rm op}$
\index[syindex]{cscrop@$\mathscr C^{\rm op}$} as follows:
\begin{enumerate}\itemsep=0pt
\item[(1)] $\mathfrak{Ob}(\mathscr C^{\rm op}) = \mathfrak{Ob}(\mathscr C)$.
\item[(2)] Let $c,c' \in
\mathfrak{Ob}(\mathscr C^{\rm op}) = \mathfrak{Ob}(\mathscr C)$. We put
$\mathscr C^{\rm op}(c,c') = \mathscr C(c',c)$.
\item[(3)] We define structure operations
$\mathfrak m_k^{\rm op}$ of $\mathscr C^{\rm op}$ by
$
\mathfrak m_k^{\rm op}(x_1,\dots,x_k)
= (-1)^{*}\mathfrak m_k(x_k,\dots,x_1)$,
where
\smash{$
* = \sum_{1\le i<j\le k} (\deg x_i+1)(\deg x_j+1) + 1$}.
\end{enumerate}
\end{defn}

\begin{lem}\label{lem232}
\quad
\begin{enumerate}\itemsep=0pt
\item[$(1)$]
$\mathscr C^{\rm op}$ is a non-unital curved filtered $A_{\infty}$ category.
\item[$(2)$]
If $\mathscr C$ is unital $($resp.\ strict, $G$-gapped$)$, then so is $\mathscr C^{\rm op}$.
\end{enumerate}

\end{lem}
(1) is \cite[Lemma 7.10]{fu4}.
(2) is immediate from the definition.
(Definition~\ref{defn2333}\,(3)\,(a).)
\begin{defn}\label{oppfunctordef}
Let $\mathscr C_1$, $\mathscr C_2$ be non-unital curved filtered $A_{\infty}$ categories.
For a filtered $A_{\infty}$ functor $\mathscr F \colon \mathscr C_1 \to \mathscr C_2$,
we can construct its {\it opposite $A_{\infty}$ functor}
\index{opposite $A_{\infty}$ functor}
$\mathscr F^{\rm op} \colon \mathscr C^{\rm op}_1 \to \mathscr C^{\rm op}_2$ as follows.
\begin{enumerate}\itemsep=0pt
\item[(1)] $\mathscr F^{\rm op}_{\rm ob} := \mathscr F_{\rm ob}$.
\item[(2)]
$\mathscr F_k^{\rm op}(\text{\bf x})
:= (-1)^{\varepsilon(\text{\bf x})}\mathscr F_k(\text{\bf x}^{\rm op})$. Here we put
\begin{equation}\label{form214}
\text{\bf x}^{\rm op} := x_k \otimes \cdots\otimes x_1,
\end{equation}
and
\begin{gather}
\varepsilon(\text{\bf x}) := \sum_{1\le i < j \le k}
(\deg x_i + 1)(\deg x_j + 1).\label{form215}
\end{gather}
\end{enumerate}
It is checked in
\cite[Definition--Lemma 7.23]{fu4} that
$\mathscr F^{\rm op}$ is a filtered $A_{\infty}$ functor.

It is easy to see that
if $\mathscr F$ is unital (resp.\ $G$-gapped), then
the functor $\mathscr F^{\rm op}$ is also unital (resp.\ $G$-gapped).
\end{defn}

The next lemma is easy to show.
\begin{lem}\label{opopoplemma}
If $\mathscr C_1$, $\mathscr C_2$ are non-unital strict filtered $A_{\infty}$ categories,
then the set theoretical map~${\mathscr F \mapsto \mathscr F^{\rm op}}$
is the object part of the isomorphism
$
\mathcal{FUNC}(\mathscr C_1,\mathscr C_2)^{\rm op}
\cong
\mathcal{FUNC}\bigl(\mathscr C^{\rm op}_1,\mathscr C^{\rm op}_2\bigr)$.
The same holds in the curved case.
\end{lem}
The proof is easy and is omitted.

\begin{defn}
We define a filtered $A_{\infty}$ category $\mathcal{CH}$ as
\index[syindex]{CH@$\mathcal{CH}$}
follows.
$\mathfrak{Ob}(\mathcal{CH})$ is the set of (all) chain complexes
of completed free $\Lambda_0$ modules.\footnote{To avoid Russell paradox in set theory, we fix a sufficiently
large set (a universe) and consider only
completed free $\Lambda_0$ modules contained in this set.}
Let $(C,d),(C',d) \in \mathfrak{Ob}(\mathcal{CH})$. Then
\[
\mathcal{CH}^k((C,d),(C',d)) = \bigoplus_{\ell} \operatorname{Hom}_{R}\big(C^{\ell},C^{\prime\ell+k}\big).
\]
We define
\begin{equation}\label{form213}
\mathfrak m_1(x) := d\circ x + (-1)^{\deg x+1} x \circ d,
\qquad
\mathfrak m_2(x,y) := (-1)^{\deg x(\deg y + 1)}y\circ x,
\end{equation}
where $\circ$ is the composition.
We put $\mathfrak m_k = 0$ for $k\ge 3$ and $k=0$.
It is checked in \cite[Proposition~7.7]{fu4}
that $\mathcal{CH}$ is a filtered $A_{\infty}$ category.
It is strict and unital.
\end{defn}

Suppose $\mathscr C$ is a non-unital {\it strict} filtered $A_{\infty}$ category.
We will define the following four functors
 \index[syindex]{Opyon@$\mathfrak{OpYon}$} \index[syindex]{Yonop@$\mathfrak{Yon}^{\rm op}$} \index[syindex]{OpYonop@$\mathfrak{OpYon}^{\rm op}$}
\begin{alignat*}{3}
& \mathfrak{Yon} \colon\ \mathscr C \to \mathcal{FUNC}(\mathscr C^{\rm op},\mathcal{CH}), \qquad&&
\mathfrak{OpYon}\colon\ \mathscr C^{\rm op} \to \mathcal{FUNC}(\mathscr C,\mathcal{CH}),&\\
& \mathfrak{Yon}^{\rm op}\colon\ \mathscr C^{\rm op} \to \mathcal{FUNC}(\mathscr C,\mathcal{CH}^{\rm op}), \qquad&&
\mathfrak{OpYon}^{\rm op}\colon\
\mathscr C \to
\mathcal{FUNC}(\mathscr C^{\rm op},\mathcal{CH}^{\rm op}).&
\end{alignat*}
The object parts of them are defined by
\begin{alignat*}{3}
& (\mathfrak{Yon}_{\rm ob}(c)_{\rm ob})(b) := \mathscr C(b,c), \qquad&&
(\mathfrak{OpYon}_{\rm ob}(c))_{\rm ob}(b) := \mathscr C(c,b), & \\
& (\mathfrak{Yon}^{\rm op}_{\rm ob}(c))_{\rm ob}(b) := \mathscr C(b,c), \qquad&&
(\mathfrak{OpYon}^{\rm op}_{\rm ob}(c))_{\rm ob}(b) := \mathscr C(c,b).&
\end{alignat*}

\begin{defn}\label{defn26}
Let $\mathscr C$ be a strict filtered $A_{\infty}$ category.
We define a filtered
$A_{\infty}$ functor $\mathfrak{Yon}_{\rm ob}(c)
\colon \mathscr C^{\rm op} \to \mathcal{CH}$ as follows.
$
\mathfrak{Yon}_{\rm ob}(c)_{\rm ob}(b_0) := \mathscr C(b_0,c)$.
Let
$\text{\bf x} \in B_k\mathscr C^{\rm op}(b_0,b_k) = B_k\mathscr C(b_k,b_0)$,
$y \in \mathfrak{Yon}_0(c)_{\rm ob}(b_0) = \mathscr C(c,b_0)$,
$k=1,2,\dots$. Then
\[
\mathfrak{Yon}_{\rm ob}(c)_k(\text{\bf x})(y)
:= (-1)^{\varepsilon(\text{\bf x})}\mathfrak m(\text{\bf x}^{\rm op},y).
\]
\end{defn}
See \cite[Definitions 7.28]{fu4}, where
$\mathfrak{Rep}$ is used instead of $\mathfrak{YON}$.
We apply the construction of~${\mathfrak{Yon}_{\rm ob}(c)}$ to the
opposite filtered $A_{\infty}$ category $\mathscr C^{\rm op}$ and define
$\mathfrak{OpYon}_{\rm ob}(c)$ as follows.
\begin{defn}\label{defn27}
$\mathfrak{OpYon}_{\rm ob}(c) \colon \mathscr C^{\rm op} \to \mathcal{CH}$ is defined by
\begin{enumerate}\itemsep=0pt
\item[(1)] $\mathfrak{OpYon}_{\rm ob}(c)(b_0) := \mathscr C(c,b_0)$,
\item[(2)] $\mathfrak{OpYon}_k(\text{\bf x})(y)
:= -(-1)^{\deg'y\deg'\text{\bf x}}\mathfrak m_{k+1}(y,\text{\bf x})$.\footnote{There are errors in \cite{fu4} on the
corresponding statements. It is corrected here.
It does not affect other parts of \cite{fu4}
since the functor $\mathfrak{OpYon}$ is not used in \cite{fu4}.
It will be used in this paper.}
\end{enumerate}
\end{defn}
\begin{lem}
There exist a filtered $A_{\infty}$ functors
$\mathfrak{Yon} \colon \mathscr C \to \mathcal{FUNC}(\mathscr C^{\rm op},\mathcal{CH})$
and $\mathfrak{OpYon}\colon\allowbreak \mathscr C \to
\mathcal{FUNC}(\mathscr C^{\rm op},\mathcal{CH}^{\rm op})$
such that its object part is given by Definitions {\rm\ref{defn26}} and {\rm\ref{defn27}}.
\end{lem}
The case of $\mathfrak{Yon}$ is \cite[Definitions 9.6 and Lemma 9.8]{fu4}.
The case of $\mathfrak{OpYon}$ is obtained by applying the case of $\mathfrak{Yon}$
to the opposite category $\mathscr C^{\rm op}$.

\begin{defn}
We define $\mathfrak{Yon}^{\rm op}\colon \mathscr C^{\rm op}\! \to\! \mathcal{FUNC}(\mathscr C,\mathcal{CH}^{\rm op})$
and
$\mathfrak{OpYon}^{\rm op}\colon \mathscr C^{\rm op} \!\to \!\mathcal{FUNC}(\mathscr C^{\rm op},\allowbreak\mathcal{CH}^{\rm op})$
to be the opposite functors of
$\mathfrak{Yon}$ and $\mathfrak{OpYon}$, respectively.
\end{defn}
\begin{rem}
The functors
$\mathfrak{Yon}^{\rm op}_{\rm ob}(c)$ and $\mathfrak{OpYon}^{\rm op}_{\rm ob}(c)$
are written as $\mathscr F^c$, ${}^c\mathscr F$ respectively
in~\cite[Section~7]{fu4}.
\end{rem}

\begin{defn}\label{defn28}
We say strict filtered $A_{\infty}$ functors: $ \mathscr C \to \mathcal{CH}^{\rm op}$,
$\mathscr C \to \mathcal{CH}$, $\mathscr C^{\rm op} \to \mathcal{CH}$,
$\mathscr C^{\rm op} \to \mathcal{CH}^{\rm op}$
are {\it representable}
\index{representable} if they are homotopy equivalent
to $\mathfrak{Yon}_{\rm ob}(c)$, $\mathfrak{OpYon}_{\rm ob}(c)$ and
$\mathfrak{Yon}^{\rm op}_{\rm ob}(c)$, $\mathfrak{OpYon}^{\rm op}_{\rm ob}(c)$,
for some $c \in \mathfrak{Ob}(\mathscr C) = \mathfrak{Ob}(\mathscr C^{\rm op})$,
respectively.
\end{defn}
The next lemma is easy to show.
\begin{lem}
The unitality and $G$-gappedness are preserved by
Definitions {\rm\ref{defn26}}, {\rm\ref{defn27}}, {\rm\ref{defn28}}.
\end{lem}
\begin{defn}\label{defn241}
We denote by $\mathfrak{Rep}(\mathscr C^{\rm op},\mathcal{CH})$
\index[syindex]{RepCop@$\mathfrak{Rep}(\mathscr C^{\rm op},\mathcal{CH})$} the full subcategory of
$\mathcal{FUNC}(\mathscr C^{\rm op},\mathcal{CH})$ such that
$\mathfrak{Ob}(\mathfrak{Rep}(\mathscr C^{\rm op},\mathcal{CH}))$ is the
set of all filtered representable $A_{\infty}$ functors.
\end{defn}
We denote by $\mathfrak{Rep}^G(\mathscr C^{\rm op},\mathcal{CH})$
the full subcategory of $\mathcal{FUNC}(\mathscr C^{\rm op},\mathcal{CH})$
whose objects consist of the $G$-gapped filtered representable $A_{\infty}$ functors.
The filtered $A_{\infty}$ category $\mathfrak{Rep}^G(\mathscr C^{\rm op},\mathcal{CH})$ is $G$-gapped.

We next define a filtered $A_{\infty}$ functor
$\mathfrak{Yon} \colon \mathscr C \cong \mathfrak{Rep}^G(\mathscr C^{\rm op},\mathcal{CH})$.
\begin{defn}
For an object $c$ of $\mathscr C$, the object $\mathfrak{Yon}_{\rm ob}(c)$ of
$\mathfrak{Rep}^G(\mathscr C^{\rm op},\mathcal{CH})$ is defined by Definition~\ref{defn26}.
\end{defn}
\begin{thm}[Yoneda's lemma]\label{Yoneda}
Let $\mathscr C$ be a $G$-gapped strict and unital filtered $A_{\infty}$
category.\index{Yoneda's lemma}
Then, there exists a homotopy equivalences of $G$-gapped filtered
$A_{\infty}$ categories
$\mathfrak{Yon} \colon \mathscr C \cong \mathfrak{Rep}^G(\mathscr C^{\rm op},\mathcal{CH})$, such that
$\mathfrak{Yon}_{\rm ob}(c)$ is as in Definition {\rm\ref{defn26}}.

\end{thm}
This is a filtered version of \cite[Theorem 9.1]{fu4}.
Using Theorem~\ref{white} instead of \cite[Theorem~8.6]{fu4},
the proof of Theorem~\ref{Yoneda} is the same as the proof of
\cite[Theorem 9.1]{fu4}.
\begin{defn}
We call $\mathfrak{Yon} \colon \mathscr C \cong \mathfrak{Rep}^G(\mathscr C^{\rm op},\mathcal{CH})$
the {\it $A_{\infty}$ Yoneda functor}.
\index{$A_{\infty}$ Yoneda functor}\index[syindex]{Yon@$\mathfrak{Yon}$}
\end{defn}

\begin{rem}\label{remark244}
In Section~\ref{sec:Ainfcat}, we describe the result over $\Lambda_0$ coefficient.
In most of the places we can use $\Lambda$ coefficient and forget
the filtration.
However, we then need to assume that our filtered~$A_{\infty}$ category
is strict.
So to work over $\Lambda$ coefficient
in our geometric application, a~natural way is
to proceed as follows.
We first define a curved filtered $A_{\infty}$ category over $\Lambda_0$.
Take its associated strict category.
Change the coefficient ring from $\Lambda_0$ to $\Lambda$.
This is the way taken in~\cite{AFOOO}.

We call a filtered $A_{\infty}$ category, {\it $\Lambda_0$ linear}
\index{$\Lambda_0$ linear} if its module of morphisms are $\Lambda_0$ module
and its structure equations are $\Lambda_0$ linear.
An $A_{\infty}$ category over $\Lambda$ (resp.\ $R$) is called also to be {\it $\Lambda$ linear}
(resp.\ {\it $R$ linear}).
\index{$\Lambda$ linear}\index{$R$ linear}

There are certain cases where it is better to work
over $\Lambda_0$. For example, reduction to $R$ works only for
$\Lambda_0$ linear category.
\end{rem}
Yoneda's lemma in the case of curved filtered $A_{\infty}$ category
is discussed in \cite[Section~4]{DL}.

\section{Floer theory of immersed Lagrangian submanifolds: Review}
\label{sec:HFIm}

This section is a review of Floer theory of immersed Lagrangian submanifolds.
Our main purpose here is to provide the precise definitions of various
notions we use in this paper.
We also include certain discussions on orientation in the Morse--Bott case,
which we use in later sections.
If the reader has certain knowledge on Lagrangian Floer theory
and its immersed version, the reader may skip this section
and comes back when it is quoted in later sections.

The Floer theory of immersed Lagrangian submanifolds is developed by
Akaho--Joyce in \cite{AJ}, generalizing the case of embedded Lagrangian
submanifolds
in \cite{fooobook,fooobook2}.
Here we rewrite the story by using the de Rham model.
The main reference we use on the virtual fundamental chain technique in the de Rham model is
\cite{foootech2,foootech22,fooonewbook}.
Note that \cite{fooobook,fooobook2,foootech2,foootech22,fooonewbook} do not discuss
the construction of filtered $A_{\infty}$-{\it categories}
but focus on filtered $A_{\infty}$ {\it algebras}.
The references on the category case are \cite{AFOOO, fu4, ancher}.

\subsection{Immersed Lagrangian submanifold}
\label{subsec:immlag}

Let $(X,\omega)$ be a symplectic manifold of real dimension $2n$.
We assume it is either compact or tame.
We sometimes say that $X$ is a symplectic manifold for simplicity.

\begin{notation}\label{not3131}
For a symplectic manifolds $(X,\omega_X)$, $(Y,\omega_Y)$,
we denote $(X\times Y,\pi_1^*(\omega_X)+\pi_2^*(\omega_Y))$
by $(X,\omega_X) \times(Y,\omega_Y)$.
Sometimes we denote $(X,-\omega_X)$ by $-X$ by an abuse of notation.
We also denote $(X,\omega_X) \times(Y,\omega_Y)$ by $X\times Y$ sometimes.
Moreover, we write $-X\times Y$ instead of $(X,-\omega_X) \times (Y,\omega_Y)$\index[syindex]{Xtimes@$-X\times Y$}
sometimes.
\end{notation}

\begin{defn}\label{def3131}
\quad
\begin{enumerate}\itemsep=0pt
\item[(1)]
An {\it immersed Lagrangian submanifold}
\index{immersed Lagrangian submanifold} $L$ of $(X,\omega)$ is a pair
$\bigl(\tilde L,i_L\bigr)$ where $\tilde L$ is a smooth manifold of dimension $n$
and $i_L$ is a smooth map \smash{$i_L \colon \tilde L \to X$} such that its
derivative~$d_pi_L \colon \allowbreak T_p\tilde L \to T_{i_L(p)}X$ is injective
and that $i_L^*\omega = 0$.
\item[(2)]
Sometimes we denote by $L$ the image of \smash{$i_L \colon \tilde L \to X$}
by an abuse of notation.
\item[(3)]
In this paper, all immersed Lagrangian submanifolds are assumed
to be compact and oriented unless otherwise mentioned.
\item[(4)]
We say $L = \bigl(\tilde L,i_L\bigr)$ has {\it clean self-intersection}
\index{clean self-intersection} if the
following holds.
\begin{enumerate}\itemsep=0pt
\item
The fiber product
\[
\tilde L \times_{X} \tilde L
:= \bigl\{(p,q) \in \tilde L \times \tilde L \mid i_L(p) = i_L(q)\bigr\}
\]
is a smooth submanifold of $\tilde L \times \tilde L$.
\item
For \smash{$(p,q) \in \tilde L \times_{X} \tilde L$}, we have
\[
T_{(p,q)}(\tilde L \times_{X} \tilde L)
=
\bigl\{(V,W) \in T_p\tilde L \times T_q\tilde L
\mid (d_pi_L)(V) = (d_qi_L)(W)\bigr\}.
\]
We remark that the left-hand side is automatically contained in the
right-hand side. The condition is that the right-hand side is contained in
the left-hand side.
Hereafter, we put
\smash{$
L(+) = \tilde L \times_{X} \tilde L$}.
\index[syindex]{Lplus@$L(+)$}
\end{enumerate}
\item[(5)]
We decompose $L(+)$ into the disjoint union of
finitely many connected components
as
\begin{equation}\label{form3333}
L(+)= \tilde L \sqcup \coprod_{a \in \mathcal A_L} L(a),
\end{equation}
\index[syindex]{$ \mathcal A_L$}
where $\tilde L$ is identified with the intersection of $L(+)$
and the diagonal.
We say $\tilde L$ the {\it diagonal component}
\index{diagonal component} of $L(+)$
and other $L(a)$'s the {\it switching components}.
\index{switching components}
We put $\mathcal A^+_L = \{o\} \cup \mathcal A_L$,
and $L(o) = \tilde L$.
\item[(6)]
We say that $L$ has {\it transversal self-intersection}
\index{transversal self-intersection}
when it has clean self-intersection and all the switching components
are zero-dimensional.
\end{enumerate}

\end{defn}
\begin{rem}
We consider sometimes the case when $\tilde L$ is not
connected. In such a case, the diagonal component is not actually a
connected component. We however call
it the diagonal component by an abuse of notation.

\end{rem}
We next define the notion of a relative spin structure of an immersed Lagrangian submanifold,
following \cite[Definition 8.1.2]{fooobook2}.
\begin{defn}\label{relspin}
Let $L$ be an immersed Lagrangian submanifold which has
clean self-intersection.
We fix a triangulation of $X$ such that $L$ is a subcomplex.
It induces a triangulation of $\tilde L$ such that
$i_L$ sends each simplex of $\tilde L$ to a simplex of $X$
by a diffeomorphism.

A {\it relative spin structure}
\index{relative spin structure} of $L$ is the following objects.
\begin{enumerate}\itemsep=0pt
\item[(1)]
A real and oriented vector bundle $V$ on the 3 skeleton $X_{[3]}$ of $X$.
\item[(2)]
A spin structure $\sigma$ of the bundle $i_L^*(V) \oplus T\tilde L$ on the
3 skeleton $\tilde L_{[3]}$ of $\tilde L$.
\end{enumerate}
We call $V$ the {\it background datum}
\index{background data} of our relative
spin structure.
We say also $\sigma$ is a {\it $V$-relative spin structure}.

\end{defn}
\begin{rem}
Let us put $[{\rm st}] = w^2(V) \in H^2(X;\Z_2)$.
Then a spin structure $\sigma$ of the bundle~${i_L^*(V) \oplus T\tilde L}$ on the
3 skeleton $\tilde L_{[3]}$ of $\tilde L$ exists if and only if
$w^2(L) = i_L^*([{\rm st}])$. Sometimes~${[{\rm st}]}$ is called the {\it background class}.
We use $V$ rather than $[{\rm st}]$ since to define the notion of a~relative spin structure
it is more precise when we use it.
(We may say $L$ is $[{\rm st}]$-relatively spin if~${w^2(L) = i_L^*([{\rm st}])}$. We need to
be more precise to define the notion of a relative spin structure of~$L$.)

\end{rem}
Note that the notion of a relative spin structure in Definition~\ref{relspin}
depends on the choice of a~triangulation of $X$.
We can however show that this notion is independent of such a~choice
in a similar way as \cite[Proposition 8.1.6]{fooobook2}.

The immersed Lagrangian Floer theory associates a
filtered $A_{\infty}$ algebra to an immersed
Lagrangian submanifold (which is relatively spin
and has clean self-intersection).
The underlying vector space of filtered $A_{\infty}$ algebra
is the vector space of differential forms on $\tilde L \times_X \tilde L$.
(More precisely, the completion of its tensor product with $\Lambda_0$.)
By the same reason as the Floer cohomology of a pair of
Lagrangian submanifolds (with clean intersection),
we need to use a certain principal ${\rm O}(1)$ bundle, which is
equivalent to a $\Z_2$-local system,
on the switching components.
(It is unnecessary in the self-transversal case which
was the case of \cite{AJ}.)
We next discuss this point following \cite[Section 3.7.5]{fooobook},
\cite[Section 8.8]{fooobook2} and will define $\Theta_a^-$.

\begin{defn}\label{defubution36}
\quad
\begin{enumerate}\itemsep=0pt
\item[(1)]
Let $L(a)$ be one of the switching components of $L(+)$.
$L(a)$ is a submanifold of $\bigl(\tilde L \times \tilde L\bigr) \setminus$
diagonal.
We compose \smash{$L(a) \to \tilde L \times \tilde L$} with the projection
to the first factor to obtain $i_{a,l} \colon L(a) \to \tilde L$.
\index[syindex]{Ial@$i_{a,l}$} This is a smooth immersion.
Using the projection to the second factor, we obtain $i_{a,r} \colon L(a) \to \tilde L$.
\index[syindex]{iar@$i_{a,r}$}
\item[(2)]
For $x \in X$, we denote by
$\mathcal{LGR}_x$
\index[syindex]{LGRx@$\mathcal{LGR}_x$} the set of all the oriented $n$-dimensional
subspaces $V$ of~$T_xX$ such that $\omega = 0$ on $V$.
\smash{$\bigcup_{x\in X} \mathcal{LGR}_x$} is a fiber bundle over $X$
which we \index[syindex]{LGR@$\mathcal{LGR}$} write $\mathcal{LGR}$.
\end{enumerate}

\end{defn}
Below we assume
\begin{equation}\label{codimensioncond}
\dim L - \dim L(a) \ge 2
\end{equation}
for switching components $L(a)$.
The orientation problem of the general case can be reduced to this
case by the following trick.
Let $u \colon (\Sigma, \partial\Sigma) \to (X,L)$ be a
pseudo-holomorphic map (see Definition~\ref{def3737}).
For the orientation problem, it suffices to consider the case
$\Sigma \subset \C$.
When we replace $X$, $L$, $u$ by
$X \times \C$, $L \times \partial\Sigma$, $u\times {\rm identity}$,
the moduli spaces of pseudo-holomorphic maps (together with their
Kuranishi structures), do not change by this process.
Therefore, we may assume~\eqref{codimensioncond} without loss of generality.

We take and fix a Riemannian
metric on $\tilde L$. This is nothing but the reduction of the structure
group of its tangent bundle to ${\rm SO}(n)$.
\begin{defn}[{see \cite[p.~687 and p.~721]{fooobook2}}]\label{lem34}
Let $x \in L(a)$.
\begin{enumerate}\itemsep=0pt
\item[(1)]
We denote by $\mathcal P^a_x$
\index[syindex]{Pax@$\mathcal P^a_x$}
the set of all smooth maps
$\lambda_x \colon [0,1] \to \mathcal{LGR}_x$
such that
\begin{enumerate}\itemsep=0pt
\item
$\lambda_{x}(0) = (d_xi_{a,l})\bigl(T_{i_{a,l}(x)}{\tilde L}\bigr)$,
\item
$\lambda_{x}(1) = (d_xi_{a,r})\bigl(T_{i_{a,r}(x)}{\tilde L}\bigr)$,
\item
$\lambda_{x}(t) \supseteq
(d_xi_{a,l})\bigl(T_{i_{a,l}(x)}{\tilde L}\bigr)
\cap (d_xi_{a,r})\bigl(T_{i_{a,r}(x)}{\tilde L}\bigr)$.
\end{enumerate}
\item[(2)]
For $\lambda_x \in \mathcal P^a_x$, we define the space
${\bf I}_{\lambda_x}$
\index[syindex]{Ilambdax@${\bf I}_{\lambda_x}$} as the set of all smooth maps
$\sigma \colon [0,1] \times \R^n \to TX$
such that $\sigma(t;\cdot) \colon \R^n \to T_xX$ is a
linear isometry between $\R^n$ and the linear subspace $\lambda_x(t)$
of $T_{x}X$.
\item[(3)]
Let $P_{\rm SO}L$ be the principal ${\rm SO}(n)$ bundle
associated to the tangent bundle.
\index[syindex]{PSOL@$P_{\rm SO}L$}
We may identify its fiber at $p \in \tilde L$ with the set of all
orientation preserving isometries $\R^n \to T_p\tilde L$.
For $x = (p,q) \in L(a)$, we consider
\[
P_x =
\frac{{(P_{\rm Spin}L)_p} \times {(P_{\rm Spin}L)_q}}
{\{-1,+1\}}.
\]
Here ${(P_{\rm Spin}L)_p}$ is the double cover of
the fiber of $P_{\rm SO}L$ at $p$
and can be identified with ${\rm Spin}(n)$.
The denominator $\{-1,+1\}$ is the group ${\rm O}(1)$ consisting
of $(1,1) \in {\rm Spin}(n)\times {\rm Spin}(n)$ and
$(-1,-1) \in {\rm Spin}(n)\times {\rm Spin}(n)$.
\item[(4)]
For $x = (p,q) \in L(a)$,
we define a map ${\bf I}_{\lambda_x} \to (P_{\rm SO} L)_p
\times (P_{\rm SO} L)_q$
by restricting $\sigma \in {\bf I}_{\lambda_x}$ to~${t=0,1}$.
We also have a double cover $P_x \to (P_{\rm SO} L)_p
\times (P_{\rm SO} L)_q$.
We define the space~$\widetilde {\bf I}_{\lambda_x}$ by the
fiber product\index[syindex]{Ilambdax@$\widetilde {\bf I}_{\lambda_x}$}
\smash{$
\widetilde {\bf I}_{\lambda_x}
= {\bf I}_{\lambda_x}
\times_{(P_{\rm SO} L)_p
\times (P_{\rm SO} L)_q}
P_x$}.
\item[(5)]
We put
\[
\mathcal I_x = \bigcup_{\lambda_x \in \mathcal P^a_x}
{\bf I}_{\lambda_x},
\qquad
\widetilde{\mathcal I}_x = \bigcup_{\lambda_x \in \mathcal P^a_x}
\widetilde {\bf I}_{\lambda_x}.
\]
The projection $\mathcal I_x \to \mathcal P^a_x$ is a fiber \index[syindex]{Ix@$\mathcal I_x$}
bundle.
\end{enumerate}

\end{defn}
We next want to regard $\bigcup_{x \in L(a)}\widetilde{\mathcal I}_{x}$
as a fiber bundle over $L(a)_{[3]}$.
We use a relative spin structure for this purpose.
Let $V$ be a real and oriented vector bundle on the 3 skeleton $X_{[3]}$.
We fix a metric on it.
We may assume that $L(a)_{[3]}$ is contained in $X_{[3]}$.
Let $x = (p,q) \in L(a)_{[3]}$.
We denote by $P_{\rm SO}(TL\oplus V)$
the principal ${\rm SO}$ bundle on \smash{$\tilde L_{[3]}$}
whose fiber at $p$ is the set of linear isometries
$\R^{n+m} \to T_pL \oplus V_y$. (Here $y = i_L(p)$.)
The spin structure of $TL \oplus V$ defines a
fiber-wise double cover $P_{\rm Spin}(TL\oplus V)$
of $P_{\rm SO}(TL\oplus V)$ on $L(a)_{[3]}$.
(Note that such a double cover
may not exist for $P_{\rm SO}L$.)

We {\it choose} an orientation preserving isometry $I_y \colon V_y \cong \R^m$.
It induces an embedding
$
(P_{\rm SO}L)_p \to (P_{\rm SO}(TL\oplus V))_p
$.
By taking a double cover, we have
\begin{equation}\label{spindoubleform}
(P_{\rm Spin}L)_p \to (P_{\rm Spin}(TL\oplus V))_p.
\end{equation}
We put
\[
P_x(V) := \frac{(P_{\rm Spin}(TL\oplus V))_p \times {(P_{\rm Spin}(TL\oplus V))_q}}
{\{-1,+1\}}.
\]
Then it is a double cover of
$
{(P_{\rm SO}(TL\oplus V))_p} \times(P_{\rm SO}(TL\oplus V))_q.
$
By using \eqref{spindoubleform},
\[
\widetilde {\mathcal I}_{x}
\cong \mathcal I_{x}
\times_{(P_{\rm SO}(TL\oplus V))_p
\times(P_{\rm SO}(TL\oplus V))_q}
P_x(V)
\]
for $(p,q) \in L(a)_{[3]}$.
Note that this identification is independent of the choice of
$I_y \colon V_y \cong \R^m$.
This is because we use the same identification
for the first factor and the second factor of the
numerator.

We remark again that we are given a spin structure of the vector bundle
$T\tilde L \oplus i_L^*V$. Therefore, the unions of
${(P_{\rm Spin} TL\oplus V)_p}$
(resp.\ ${(P_{\rm Spin} TL\oplus V)_q}$) for $p$ (resp.\ $q$)
becomes a principal bundle over $L(a)_{[3]}$.
We thus obtain a fiber bundle
\smash{$\widetilde {\mathcal I} \to L(a)_{[3]}$} whose fiber at $x$
is \smash{$\widetilde {\mathcal I}_{x}$}.

We remark that
$\mathcal P^a_x$
is homotopy equivalent to the loop space
of the oriented Lagrangian Grassmannian, $\Omega\mathcal{LGR}(n-d)$,
where $d = \dim L(a)$.
It is well know that
\[
\mathcal{LGR}(n-d) = {\rm U}(n-d)/{\rm SO}(n-d).
\]
In fact, ${\rm U}(n-d)$ acts transitively to the set of all
the oriented Lagrangian linear subspaces of~$\C^{n-d}$ and the isotropy group of this action at
$\R^{n-d}$ is ${\rm SO}(n-d)$.
Therefore, we have an exact sequence
\begin{align*}
1= \pi_2({\rm U}(n-d))& \to \pi_2(\mathcal{LGR}(n-d))
\to \pi_1({\rm SO}(n-d)) \\
&\to \pi_1({\rm U}(n-d))
\to \pi_1(\mathcal{LGR}(n-d)) \to 1.
\end{align*}
We assumed \eqref{codimensioncond}, that is, $n-d \ge 2$. Therefore, $\pi_1({\rm SO}(n-d)) = \Z_2$ and $\pi_1({\rm U}(n-d)) = \Z$.
Therefore,
$
\pi_0(\Omega\mathcal{LGR}(n-d)) = \Z$,
$
\pi_1(\Omega\mathcal{LGR}(n-d)) = \Z_2$.
Moreover, the map
$
\pi_1(\Omega\mathcal{LGR}(n-d)) \to \pi_1({\rm SO}(n-d))
$
is an isomorphism.
It implies the next lemma.
\begin{lem}\label{lem37}
 The double cover $\widetilde{\mathcal I}_x \to \mathcal I_x$ is nontrivial.

\end{lem}
Using $\lambda_x$ as in Definition~\ref{lem34},
we define a Fredholm operator as follows.
(We follow \cite[Section 8.1.3]{fooobook2} here.)
We put
\index[syindex]{Zminus@$Z_-$}
\begin{gather*}
Z_-= \{z \in \C \mid \vert z\vert \le 1\} \cup
\{z \in \C \mid \operatorname{Re}z \ge 0, \, \vert\operatorname{Im}z\vert \le 1\},
\\
Z_+ = \{ -x + \sqrt{-1}y \mid x,y\in \R, \, x + \sqrt{-1}y \in Z_-\}.
\end{gather*}

\begin{figure}[ht]
\centering
\includegraphics[scale=0.4]{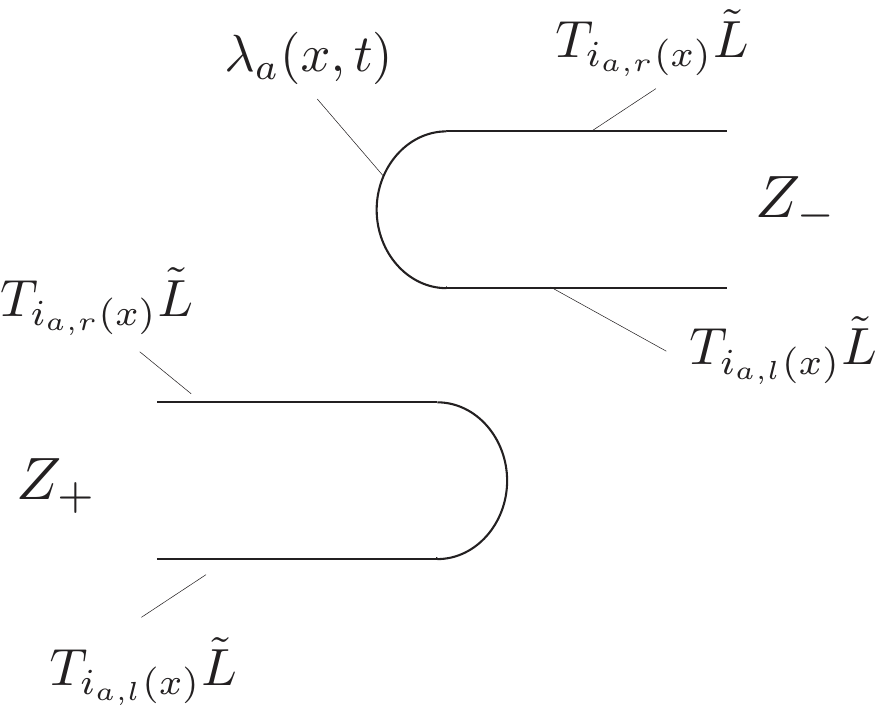}
\caption{Domains $Z_+$, $Z_-$.}
\label{Figure31}
\end{figure}

Let $k$ be a sufficiently large integer. (For example we may take
$k=100$.)
We consider
the set of locally $L^2_k$ maps $u \colon Z_- \to T_xX$
(resp.\ $u \colon Z_+ \to T_xX$)
with the following properties:
\begin{enumerate}\itemsep=0pt
\item[(1)]
$u\bigl(t + \sqrt{-1}\bigr) \in (di_L)\bigl(T_{i_{a,r}(x)}{\tilde L}\bigr)$
for $t \in \R_{\ge 0}$
(resp.\ $u\bigl(t + \sqrt{-1}\bigr) \in (di_L)\bigl(T_{i_{a,r}(x)}{\tilde L}\bigr)$
for $t \in \R_{\le 0}$),
\item[(2)]
$u\bigl(t - \sqrt{-1}\bigr) \in (di_L)\bigl(T_{i_{a,l}(x)}{\tilde L}\bigr)$
for $t \in \R_{\ge 0}$
(resp.\ $u\bigl(t - \sqrt{-1}\bigr) \in (di_L)\bigl(T_{i_{a,l}(x)}{\tilde L}\bigr)$
for $t \in \R_{\le 0}$),
\item[(3)]
$u\bigl(\exp\bigl(\pi\sqrt{-1}(3/2-t)\bigr)\bigr) \in \lambda_a(x,t)$
(resp.\ $u\bigl(\exp\bigl(\pi\sqrt{-1}(t-1/2)\bigr)\bigr) \in \lambda_a(x,t)$),
\item[(4)]
\begin{equation}\label{form34}
\sum_{\ell=0}^k\int_{W} e^{\delta\vert{\operatorname{Re}z}\vert } \Vert\nabla^{\ell} u\Vert^2  d z d \overline z
< \infty.
\end{equation}
Here $\delta > 0$ is a fixed small number.
See Figure~\ref{Figure31}.
\end{enumerate}
We consider the totality of such maps $u$ and use the
left-hand side of \eqref{form34} as its norm. We denote it by
$L^2_k(Z_-;T_xX;\lambda_a;\delta)$
(resp.\ $L^2_k(Z_+;T_xX;\lambda_a;\delta)$).
This is a Hilbert space.
We consider the set of all the locally $L^2_{k-1}$
maps $u \colon Z_- \to T_xX$ (resp.\ $u \colon Z_+ \to T_xX$)
which satisfies \eqref{form34} with $k$ replaced by $k-1$
and denote it by
$L^2_{k-1}(Z_-;T_xX;\delta)$
(resp.\ $L^2_{k-1}(Z_+;T_xX;\delta)$).
(In other words, we do not require (1), (2), (3) for this function space.)

We use the Cauchy--Riemann operator
to define the operators
$\overline\partial_{Z_-,\lambda_x}$,
$\overline\partial_{Z_+,\lambda_x}$
\begin{gather}
\overline\partial_{Z_-,\lambda_x} \colon\ L^2_k(Z_-;T_xX;\lambda_a;\delta)
\to L^2_{k-1}(Z_-;T_xX;\delta),\nonumber
\\
\overline\partial_{Z_+,\lambda_x} \colon\ L^2_k(Z_+;T_xX;\lambda_a;\delta)
\to L^2_{k-1}(Z_+;T_xX;\delta).\label{3535}
\end{gather}
\index[syindex]{doverlineZ-x@$\overline\partial_{Z_-,\lambda_x}$}
The next lemma is now standard.
\begin{lem}
The operators $\overline\partial_{Z_-,\lambda_x}$,
$\overline\partial_{Z_+,\lambda_x}$ are Fredholm.
\end{lem}

By moving $x$ and $\lambda_x$, we obtain the
family index bundles
${\rm Ind}\bigl(\overline\partial_{Z_-,\lambda_x}\bigr)$,
${\rm Ind}\bigl(\overline\partial_{Z_+,\lambda_x}\bigr)$
and their determinant real line bundles
$\operatorname{Det} {\rm Ind}\bigl(\overline\partial_{Z_-,\lambda_x}\bigr)$,
$\operatorname{Det} {\rm Ind}\bigl(\overline\partial_{Z_+,\lambda_x}\bigr)$.
They are bundles over $\mathcal I$.
\begin{lemdef}\label{lemdef39}
\quad
\begin{enumerate}\itemsep=0pt
\item[$(1)$]
The restriction of the pullback of
$\operatorname{Det} {\rm Ind}\bigl(\overline\partial_{Z_-,\lambda_x}\bigr)$,
$\operatorname{Det} {\rm Ind}\bigl(\overline\partial_{Z_+,\lambda_x}\bigr)$
to $\widetilde{\mathcal I}_x$ is trivial.
\item[$(2)$]
Moreover, we can define a real line bundles on $L(a)$ in a canonical way,
which pulls back to $\operatorname{Det} {\rm Ind}\bigl(\overline\partial_{Z_-,\lambda_x}\bigr)$,
$\operatorname{Det} {\rm Ind}\bigl(\overline\partial_{Z_+,\lambda_x}\bigr)$.
\item[$(3)$]
We denote by $\Theta^-_a$, $\Theta^+_a$
\index[syindex]{thetaminusa@$\Theta^-_a$},
\index[syindex]{thetaplusa@$\Theta^+_a$ }
the principal ${\rm O}(1)$ bundles which correspond to the real line bundles on $L(a)$ in item $(2)$.
\item[$(4)$]
There exists an isomorphism
$
\Theta_a^- \otimes \Theta_a^+ \cong \operatorname{Det} TL(a)$.
\end{enumerate}

\end{lemdef}
\begin{proof}
This is \cite[Proposition 8.8.1]{fooobook2}.
We refer to \cite[pp.~721--722]{fooobook2}
for the proof of item~(1).

We provide a bit more detail of the proof of item (2) than
\cite{fooobook2} here, since we use a certain part of the construction
in the proof of Lemma~\ref{lem310}.
Recall $\pi_0({\mathcal I}_x)$
corresponds one to one to integers $k$.
Let ${\mathcal I}_{x,k}$ be the corresponding connected component.
Its union for \smash{$x \in L(a)_{[3]}$} is denoted by
${\mathcal I}_{k}$.
Its pullback to $\tilde{\mathcal I}$
is denoted by $\tilde{\mathcal I}_{k}$.

The double cover $\widetilde{\mathcal I}_{x,k} \to
{\mathcal I}_{x,k}$ is nontrivial by Lemma~\ref{lem37}. Therefore,
$
\pi_0\bigl(\tilde{\mathcal I}_{x,k}\bigr) = \pi_0({\mathcal I}_{x,k})
$
is trivial. We then have an exact sequence
\[
\pi_1\bigl(\widetilde{\mathcal I}_{x,k}\bigr)
\to \pi_1\bigl(\widetilde{\mathcal I}_k\bigr)
\to \pi_1(L(a)_{[3]})
\to 1.
\]

Therefore, \smash{$\pi_1\bigl(\widetilde{\mathcal I}_k\bigr)
\to \pi_1(L(a)_{[3]})$} is surjective.
Thus item (1) implies that
we can define a~group homomorphisms $\pi_1(L(a)_{[3]}) \to \Z_2$
which pulls back to the homomorphism
$\pi_1\bigl(\widetilde{\mathcal I}_k\bigr) \to \Z_2$
given by $\operatorname{Det} {\rm Ind}\bigl(\overline\partial_{Z_-,\lambda_x}\bigr)$,
$\operatorname{Det} {\rm Ind}\bigl(\overline\partial_{Z_+,\lambda_x}\bigr)$.
Thus we obtain real vector bundles
on \smash{$L(a)_{[3]}$} which pull back to
$\operatorname{Det} {\rm Ind}\bigl(\overline\partial_{Z_-,\lambda_x}\bigr)$ and
$\operatorname{Det} {\rm Ind}\bigl(\overline\partial_{Z_+,\lambda_x}\bigr)$
on $\widetilde{\mathcal I}_k$, respectively.

See the proof of \cite[Proposition 8.8.1]{fooobook2}
for the proof of the
fact that this line bundle is independent of $k$.
Since any real line bundles on the 3-skeleton uniquely
extend to the whole space,
we obtain a real line bundles on $L(a)$.

We next discuss item (4). This is proved in \cite[Section 8.8]{fooobook2}.
We sketch its proof below since we use a similar
argument in the proof of Lemma~\ref{lem310}.
It suffices to show that the isomorphism for
an arbitrary loop $\gamma$ in $L(a)_{[3]}$.
We choose a loop $\gamma$ and fix a trivialization of $V$
on $\gamma$.

We will prove the isomorphism of family indices
\begin{equation}\label{form382}
{\rm Ind}\bigl(\overline\partial_{Z_-,\lambda_x}\bigr) \oplus
{\rm Ind}\bigl(\overline\partial_{Z_+,\lambda_x}\bigr) \oplus
T_x L(a)
\cong
{\rm Ind}\bigl(\overline\partial_{Z,\lambda^2_x}\bigr),
\end{equation}
where the right-hand side is defined as follows.
We put
\index[syindex]{ZR@$Z(R)$}
\[
Z(R) =
\{z \in \C \mid \vert z - R\vert \le 1\}
\cup
\{z \in \C \mid \vert z + R\vert \le 1\}
\cup
\{z \in \C \mid \vert \operatorname{Im} z \vert \le 1,\,
 \vert{\operatorname{Re} z}\vert \le R\}.
\]
See Figure~\ref{Figure32}.
We use $\lambda_x$ on
$\partial Z(R) \cap \partial \{z \in \C \mid \vert z - R\vert \le 1\}$
and on $\partial Z(R) \cap \partial \{z \in \C \mid \vert z + R\vert \le 1\}$
and $d_x{i_{a,L}}\bigl(T_{i_{a,l}(x)}{\tilde L}\bigr)$
\big(resp.\ $d_x{i_{a,r}}\bigl(T_{i_{a,r}(x)}{\tilde L}\bigr)$\big)
on $\partial Z(R) \cap
\{z \mid \operatorname{Im} z = -1\}$
(resp.\ $\partial Z(R) \cap
\{z \mid \operatorname{Im} z = 1\}$)
to define the boundary condition
\smash{$\lambda_x^2$} on $\partial Z(R)$.
We then
obtain\index[syindex]{doverlineZR@$\overline\partial_{Z(R),\lambda^2_x}$}
\[
\overline\partial_{Z(R),\lambda^2_x} \colon\ L^2_k\bigl(Z(R);T_xX;\lambda^2_x\bigr)
\to L^2_{k-1}(Z(R);T_xX)
\]
in the same way as
$\overline\partial_{Z_-,\lambda_x}$.
(Since $Z(R)$ is compact we do not use weighted Sobolev
space but use usual Sobolev space.)

\begin{figure}[ht]
\centering
\includegraphics[scale=0.5]{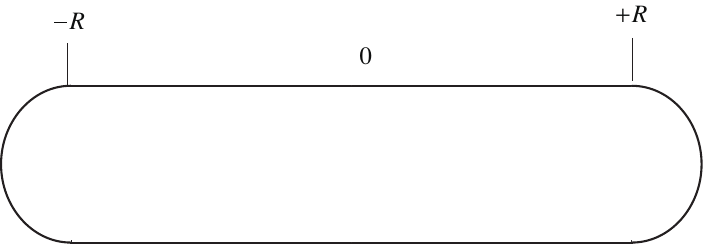}
\caption{Domains $Z(R)$.}
\label{Figure32}
\end{figure}
We glue two index problems
\smash{$\overline\partial_{Z_-,\lambda_x}$}
and \smash{$\overline\partial_{Z_+,\lambda_x}$}
at their ends and the result is
\smash{$\overline\partial_{Z(R),\lambda^2_x}$}.
Note that, however, there is a degeneration of the
operators \smash{$\overline\partial_{Z_-,\lambda_x}$}
and \smash{$\overline\partial_{Z_+,\lambda_x}$}
at the end.
The eigenspace of $0$ of this degeneration is
exactly $T_xL(a)$.
(See Definition~\ref{defubution36}\,(1c).)
Therefore, the standard family index sum formula
(see, for example, \cite[Theorem 4.9]{fu-1})
gives~\eqref{form382}.

Now we consider the family of indices of \smash{$\overline\partial_{Z(R),\lambda^2_x}$},
where we move $x$ and $\lambda^2_x$, and regard
it as a bundle on $\widetilde{\mathcal I}_k$.
Then since we are working on $\widetilde{\mathcal I}_k$,
the boundary has a canonical spin structure as family
on the boundary condition $\lambda^2_x$.
In fact, we fixed a trivialization of $V$.
So the spin structure of $\lambda^2_x(t) \oplus V$
corresponds one to one to the spin structure of $\lambda^2_x(t)$.
Therefore, the determinant line bundle of the family
\smash{$\overline\partial_{Z(R),\lambda^2_x}$}
on $\widetilde{\mathcal I}_k$ is trivial.
We thus proved item (4).
\end{proof}

We use the next lemma in the later sections.
We consider two $V$-relative spin structures~$\sigma_1$ and~$\sigma_2$
of $L$. Then the difference $\sigma_1 - \sigma_2$ is
regarded as an element of $H^1\bigl(\tilde L;\Z_2\bigr)$.
Using Lemma--Definition~\ref{lemdef39},
we obtain a line bundle $\Theta^-_a$ for each of
the $V$-relative spin structures $\sigma_1$ and $\sigma_2$.
In the next lemma, we write them as
$\Theta^-_{a,\sigma_1}$ and $\Theta^-_{a,\sigma_2}$,
respectively.
\begin{lem}\label{lem310}
The tensor product\footnote{See also \cite[Proposition~3.10]{fooo:inv}.} $\Theta^-_{a,\sigma_1}\otimes \Theta^-_{a,\sigma_2}$
is the principal ${\rm O}(1)$ bundle corresponding~to
\[
i_{a,l}^*(\sigma_1 - \sigma_2)
- i_{a,r}^*(\sigma_1 - \sigma_2)
\in H^1(L(a);\Z_2),
\]
where
$\sigma_1 - \sigma_2 \in H^1\bigl(\tilde L;\Z_2\bigr)$
is as above.
\end{lem}
\begin{proof}
We consider a loop $\gamma \colon S^1 \to L(a)_{[3]}$
and fix a trivialization of $V$ on $\gamma$.
In the case when
\begin{gather}\label{form3900}
(i_{a,l}^*(\sigma_1 - \sigma_2)
- i_{a,r}^*(\sigma_1 - \sigma_2))
\cap [\gamma] = 0,
\end{gather}
we will prove that
$\Theta^-_{a,\sigma_1}\otimes \Theta^-_{a,\sigma_2}$
is trivial on $\gamma$.
In the case when
\begin{equation}\label{form39}
(i_{a,l}^*(\sigma_1 - \sigma_2)
- i_{a,r}^*(\sigma_1 - \sigma_2))
\cap [\gamma] \ne 0,
\end{equation}
we will prove that
$\Theta^-_{a,\sigma_1}\otimes \Theta^-_{a,\sigma_2}$
is nontrivial on $\gamma$.

Let $\widetilde{\mathcal I}^{
\sigma_1}_k$ \big(resp.\ $\widetilde{\mathcal I}^{
\sigma_2}_k$\big)
be the space $\widetilde{\mathcal I}_k$
we obtain as in the proof of Lemma--Definition~\ref{lemdef39}
using relative spin structure $\sigma_1$
(resp.\ $\sigma_2$).
As we proved during the proof of Lemma--Definition~\ref{lemdef39}\,(2), the loop $\gamma$ lifts to $\widetilde{\mathcal I}^{
\sigma_1}_k$
\big(resp.\ $\widetilde{\mathcal I}^{
\sigma_2}_k$\big).

We take the lift $\tilde\gamma^{\sigma_1} \colon S^1 \to \widetilde{\mathcal I}^{
\sigma_1}_k$
\big(resp.
$\tilde\gamma^{\sigma_2} \colon S^1 \to \widetilde{\mathcal I}^{
\sigma_2}_k$\big).
We compose it with the projection to obtain
$\gamma^{\sigma_1} \colon S^1 \to {\mathcal I}^{
\sigma_1}_k$
\big(resp.
$\gamma^{\sigma_2} \colon S^1 \to {\mathcal I}^{
\sigma_2}_k$\big).
(As we can show from the discussion below, $\gamma^{\sigma_1}$~is not
homotopic to $\gamma^{\sigma_2}$
if \eqref{form39} holds.)

For each $s \in S^1$, the element $\gamma^{\sigma_1}(s)$
defines a path $\lambda^{\sigma_1}_{s}(\cdot) \colon [0,1] \to \mathcal{LGR}_{\gamma(s)}$
satisfying Definition~\ref{lem34}\,(1)\,(a)\,(b)\,(c) for
$x =\gamma(s)$.
We obtain $\lambda^{\sigma_2}_{s}(\cdot)$ in the same way.
We use them to obtain Fredholm operators
$\overline\partial_{Z_-,\lambda^{\sigma_1}_s}$,
$\overline\partial_{Z_+,\lambda^{\sigma_1}_s}$,
$\overline\partial_{Z_-,\lambda^{\sigma_2}_s}$,
$\overline\partial_{Z_+,\lambda^{\sigma_2}_s}$.
by \eqref{3535}.
It suffices to show that
\begin{equation}\label{form392}
\operatorname{Det} \operatorname{Index}\bigl(\overline\partial_{Z_-,\lambda^{\sigma_1}_s}\bigr)
\otimes
\operatorname{Det} \operatorname{Index}\bigl(\overline\partial_{Z_-,\lambda^{\sigma_2}_s}\bigr)
\end{equation}
is a nontrivial real line bundle as a family index bundles over $S^1
= \gamma$, if and only if \eqref{form39} holds.

By Lemma--Definition~\ref{lemdef39},
\begin{align*}
\operatorname{Det} \operatorname{Index}\bigl(\overline\partial_{Z_-,\lambda^{\sigma_1}_s}\bigr)
\otimes
\operatorname{Det} \operatorname{Index}\bigl(\overline\partial_{Z_+,\lambda^{\sigma_1}_s}\bigr)
&\cong
\operatorname{Det} TL(a)\\
&
\cong
\operatorname{Det} \operatorname{Index}\bigl(\overline\partial_{Z_-,\lambda^{\sigma_2}_s}\bigr)
\otimes
\operatorname{Det} \operatorname{Index}\bigl(\overline\partial_{Z_+,\lambda^{\sigma_2}_s}\bigr).
\end{align*}
Therefore, \eqref{form392} is isomorphic to
\begin{equation*}
\operatorname{Det} \operatorname{Index}\bigl(\overline\partial_{Z_-,\lambda^{\sigma_1}_s}\bigr)
\otimes
\operatorname{Det} \operatorname{Index}\bigl(\overline\partial_{Z_+,\lambda^{\sigma_2}_s}\bigr)
\otimes \operatorname{Det} TL(a).
\end{equation*}
We define $\lambda_s^{\sigma_1,\sigma_2}(x)
\subset T_xX$
for $z \in \partial Z(R)$ as follows.

We use $\lambda^{\sigma_1}_s$ on $\partial Z(R) \cap \partial \{z \in \C \mid \vert z + R\vert \le 1\}$,
$\lambda^{\sigma_2}_s$ on $\partial Z(R) \cap \partial \{z \in \C \mid \vert z - R\vert \le 1\}$,
and $d_{i_L}\bigl(T_{i_{a,l}(\gamma(s))}{\tilde L}\bigr)$
\big(resp.\ $d_{i_r}\bigl(T_{i_{a,l}(\gamma(s))}{\tilde L}\bigr)$\big)
on $\partial Z(R) \cap
\{z \mid \operatorname{Im} z = -1\}$
(resp.\ $\partial Z(R) \cap
\{z \mid \operatorname{Im} z = 1\}$).
We then obtain
\smash{$\overline\partial_{Z(R),\lambda_s^{\sigma_1,\sigma_2}}$}
in the same way as \smash{$\lambda_s^{\sigma_1,\sigma_2}$}.
See Figure~\ref{Figure33}.
\begin{figure}[ht]
\centering
\includegraphics[scale=0.5]{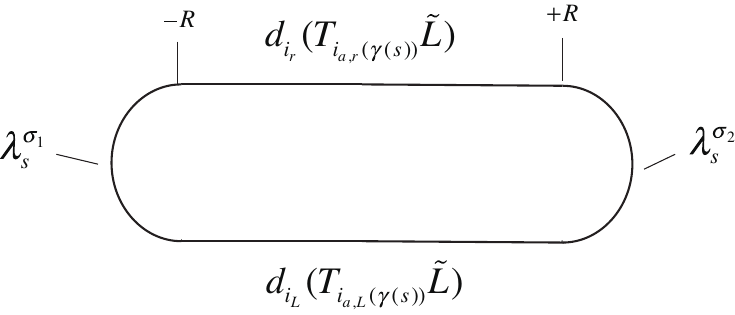}
\caption{Domains $Z(R)$.}
\label{Figure33}
\end{figure}

In the same way as the proof of
Lemma--Definition~\ref{lemdef39}\,(4),
the bundle \eqref{form392} is isomorphic to%
\begin{equation}\label{form311}
\operatorname{Det} \operatorname{Index}\bigl(\overline\partial_{Z(R),\lambda_s^{\sigma_1,\sigma_2}}\bigr).
\end{equation}
Note that we consider the family of $n$-dimensional
real vector spaces $\lambda_s^{\sigma_1,\sigma_2}(z)$
parametrized by~${(s,z) \in S^1 \times \partial Z(R)
\cong S^1 \times S^1}$.
\eqref{form39} implies that this bundle has nontrivial second
Stiefel--Whitney class.
Therefore, the example given in the proof of \cite[Proposition 8.1.7]{fooobook2}
implies that the family index bundle
\eqref{form311} is nontrivial on $S^1$.

\eqref{form3900} implies that the boundary
condition corresponds to a trivial bundle.
Therefore, we can show that the family of index bundle
\eqref{form311} is trivial in this case.

The proof of Lemma~\ref{lem310} is complete.
\end{proof}

\begin{exm}
Let $L_0$ be an embedded and spin Lagrangian submanifold
of $X$,
$\tilde L$ a disjoint union of two copies of
$L_0$, and $i_L \colon \tilde L \to X$ the identity
maps on each components.
The fiber product $\tilde L \times_X \tilde L$ is disjoint
union of 4 copies of $L_0$, where two are diagonal components and
two are switching components.
We take two different spin structures $\sigma_1$ and
$\sigma_2$ on $L_0$ and use them for the two
connected components of $\tilde L$ and obtain
a relative spin structure.

Then $\Theta^-$ is the trivial bundle on the diagonal components
and is the line bundle corresponding to
$\sigma_1 - \sigma_2 \in H^1(L_0;\Z_2)$ on the
switching components.
\end{exm}

Let $\Theta$ be a principal ${\rm O}(1)$ bundle on a manifold $M$.
We denote by $\Omega(M;\Theta)$ the $\R$ vector space of smooth differential forms on $M$
with coefficient $\Theta$, that is, the set of smooth sections of~${\Omega^M \otimes \Theta}$.
Here $\Omega^M$ is the real vector bundle of differential forms on $M$ and
$\Theta$ is the real line bundle corresponding to the principal
${\rm O}(1)$ bundle
$\Theta$.
\begin{defn}\label{defn313}
Suppose $R = \R$ or $\C$.
We put
\begin{align}
CF\bigl(L;\Lambda^R_0\bigr)
&=
\Omega(L(+),\Theta^-) \,\widehat\otimes_{R}\, \Lambda^R_0 \nonumber\\
&=
\bigl(\Omega\bigl(\tilde L\bigr) \,\widehat\otimes_{R}\, \Lambda^R_0\bigr)
\oplus
\bigoplus_{a \in \mathcal A_L}
\bigl(\Omega(L(a),\Theta_a^-) \,\widehat\otimes_{R}\, \Lambda^R_0\bigr).\label{form315}
\end{align}
\index[syindex]{CFL@$CF\bigl(L;\Lambda^R_0\bigr)$}
Here $\widehat\otimes_{R}$ is the $T$-adic completion of the algebraic tensor
product.

We remark that $CF\bigl(L;\Lambda^R_0\bigr)$ is a completed free $\Lambda^R_0$ module.

We also denote
\[
CF(L;\R)
=
\Omega(L(+),\Theta^-) =
\Omega\bigl(\tilde L\bigr)
\oplus
\bigoplus_{a \in \mathcal A_L} \Omega(L(a),\Theta_a^-).
\]

\end{defn}

\subsection{Moduli space of pseudo-holomorphic polygons}
\label{subsec:modpolygon}

The purpose of this subsection is to prove the next theorem.

\begin{thm}\label{AJtheorem}
Let $L$ be a relatively spin immersed Lagrangian submanifold of $(X,\omega)$.
We assume that $L$ has clean self-intersection.
Then we can define a structure of
filtered $A_{\infty}$ algebra on the
completed free graded $\Lambda^R_0$ module $CF\bigl(L;\Lambda^R_0\bigr)$.
It is unital and is $G$-gapped for some discrete submonoid $G$.
\end{thm}

\begin{rem}
Theorem~\ref{AJtheorem} is proved
by Akaho--Joyce in \cite{AJ} except the following points.
Those points are of technical nature.
\begin{enumerate}\itemsep=0pt
\item[(1)]
We include the case of clean self-intersection.
Akaho--Joyce \cite{AJ} restrict themselves to the case of transversal self-intersection.
This difference is not essential.
In fact, Lagrangian Floer theory in the Morse--Bott situation
is fully worked out in~\cite{fooobook}.
I think \cite{AJ} restricted themselves to the transversal case only for the sake of simplicity of notations.
We include it, since we need to use the clean self-intersection case in Section~\ref{sec:Unobstructedness}.
\item[(2)]
We use the de Rham model to work out the transversality issue, while
\cite{AJ} used the singular homology.
The author of this paper together with joint authors has completed detailed account explaining the way
to use the de Rham model in the virtual
fundamental chain technique after~\cite{AJ} was published (see \cite{foootech,foootech2,foootech22,fooonewbook}).
In his opinion using the de Rham model is the shortest way to
work out the virtual fundamental chain technique in the chain level
in detail and rigorously
when we include Morse--Bott situation.
\item[(3)]
We prove (exact) unitality of the algebra.
In fact, using the de Rham model we can obtain an exact unit
(see \cite{fooo091}).
When using singular homology, we obtain
a homotopy unit but it is hard to obtain an exact unit
(see \cite[Section~3.3]{fooobook} and \cite[Section~7.3]{fooobook2}.)
\end{enumerate}
\end{rem}
\begin{rem}
The filtered $A_{\infty}$ algebra $\bigl(CF\bigl(L;\Lambda^R_0\bigr),
\{\mathfrak m_k\}\bigr)$ in Theorem~\ref{AJtheorem} depends on various choices but
is independent of the choices up to homotopy equivalence.
See Remark~\ref{rem340} and Section~\ref{sec:independence2}.
(In Section~\ref{sec:independence2}, we will prove the case of filtered $A_{\infty}$ category,
which implies the case of filtered $A_{\infty}$ algebra.)

\end{rem}
The proof of Theorem~\ref{AJtheorem} will be completed in Section
\ref{subsec:Ainfalgim}.
In this subsection, we describe the moduli spaces of pseudo-holomorphic
polygons, which are used to define the structure operations~$\mathfrak m_k$ of our filtered $A_{\infty}$ algebra.

Let $L$ be as in Theorem~\ref{AJtheorem} and
$\vec a = (a_0,\dots,a_k)$, $a_i \in \mathcal A^+_L$.
(Here $\mathcal A^+_L$ is as in Definition~\ref{def3131}\,(5).)
We fix a compatible almost complex structure $J_X$ on $X$.

\begin{defn}\label{def3737}
Let $E \in \R_{\ge 0}$.
We define the set
\smash{$\mathring{\widetilde{\mathcal M}}(L;\vec a;E)$}
\index[syindex]{M1LAE@$\mathring{\widetilde{\mathcal M}}(L;\vec a;E)$}
as the totality of all the objects~${(\Sigma;u;\vec z;\gamma)}$ with the following
properties.
\begin{enumerate}\itemsep=0pt
\item[(1)]
The space $\Sigma$ is a union of a disk plus
a finite number of trees of sphere components
attached to the interior of the disk.
$\Sigma$ is connected, simply connected and has
at worst double points as singularities.
\big(In particular, $\partial\Sigma = S^1$.\big)
(See Figure~\ref{Figure34}.)
\begin{figure}[ht]
\centering
\includegraphics[scale=0.37]{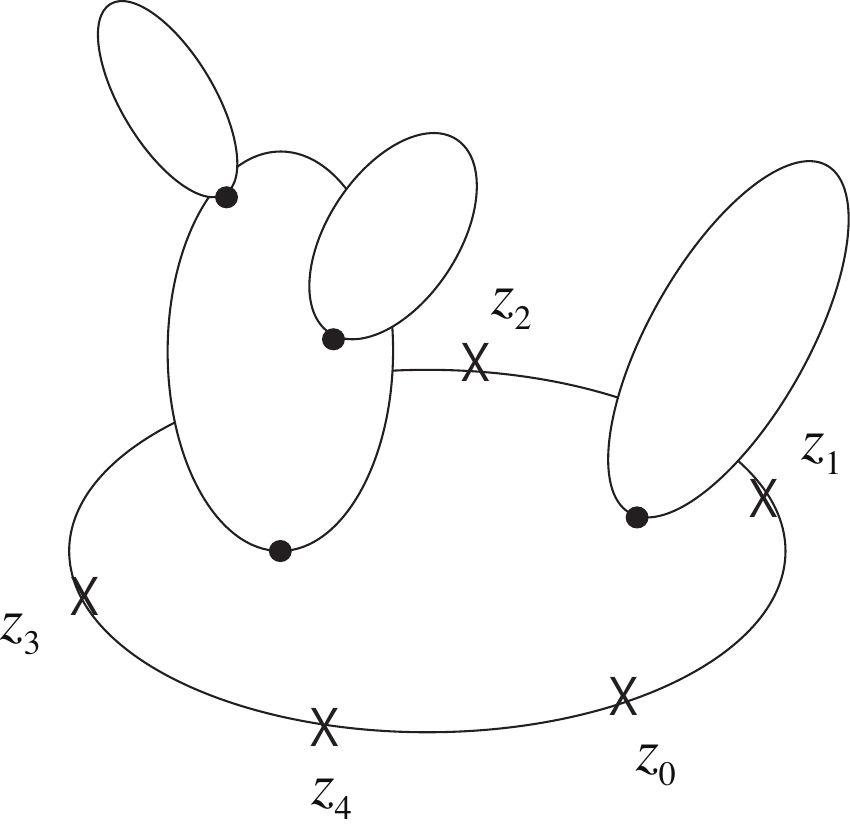}
\caption{Domain $\Sigma$.}
\label{Figure34}
\end{figure}
\item[(2)]
The map $u \colon \Sigma \to X$ is $J_X$-holomorphic.
\item[(3)]
We put $\vec z = (z_0,\dots,z_k)$.
Then, the points $z_i \in \partial\Sigma = S^1$
are mutually distinct.
The $k+1$ tuple of points $(z_0,\dots,z_k)$ respects the
counter clockwise cyclic order of $S^1$.
\item[(4)]
The map $\gamma \colon S^1 \setminus
\{z_0,\dots,z_k\} \to \tilde L$ is smooth and
satisfies
$
i_L(\gamma(z)) = u(z)
$
for $z \in S^1 \setminus
\{z_0,\dots,z_k\}$.
\item[(5)]
For $i=0,\dots,k$, we have
\smash{$
(\lim_{z \uparrow z_i}\gamma(z),\lim_{z \downarrow z_i}\gamma(z))
\in L(a_i)$}. \index[syindex]{limuparrow@$\lim_{z \uparrow z_i}$}
Here $L(a_i)$ is as in \eqref{form3333}.
The limit in the left-hand side is defined as follows.
Let \smash{$x_m = e^{t_m \sqrt{-1}} \in S^1$}, where
$t_m$ is an increasing sequence of real numbers converging to
$t$ with \smash{$e^{t \sqrt{-1}} = z_i$}. We say ${\lim_{z \uparrow z_i}\gamma(z) = y}$ if
$
\lim_{m\to\infty}\gamma(x_m) = y
$
for any such sequence $x_m$.
The definition of $\lim_{z \downarrow z_i}\gamma(z)$ is similar.
(See Figure~\ref{Figure35}.)
\item[(6)]
$
\int_{D^2} u^*\omega = E$.
\item[(7)]
(Stability)
The set of the maps $v \colon \Sigma \to \Sigma$ with the following
properties is a finite set:
\begin{enumerate}\itemsep=0pt
\item
$u\circ v = u$,
\item
$v$ is biholomorphic,
\item
$v(z_i) = z_i$,
\item
$\gamma \circ v = \gamma$.\footnote{Actually
this condition follows from~(a).}
\end{enumerate}
\end{enumerate}
We write ${\rm Aut}(\Sigma;u;\vec z;\gamma)$
the finite group consisting of the maps $v$ satisfying (a), (b), (c), (d)
above. We call it the {\it group of automorphisms} \index{group of automorphisms} of
$(\Sigma;u;\vec z;\gamma)$.

\end{defn}
\begin{figure}[ht]
\centering
\includegraphics[scale=0.4]{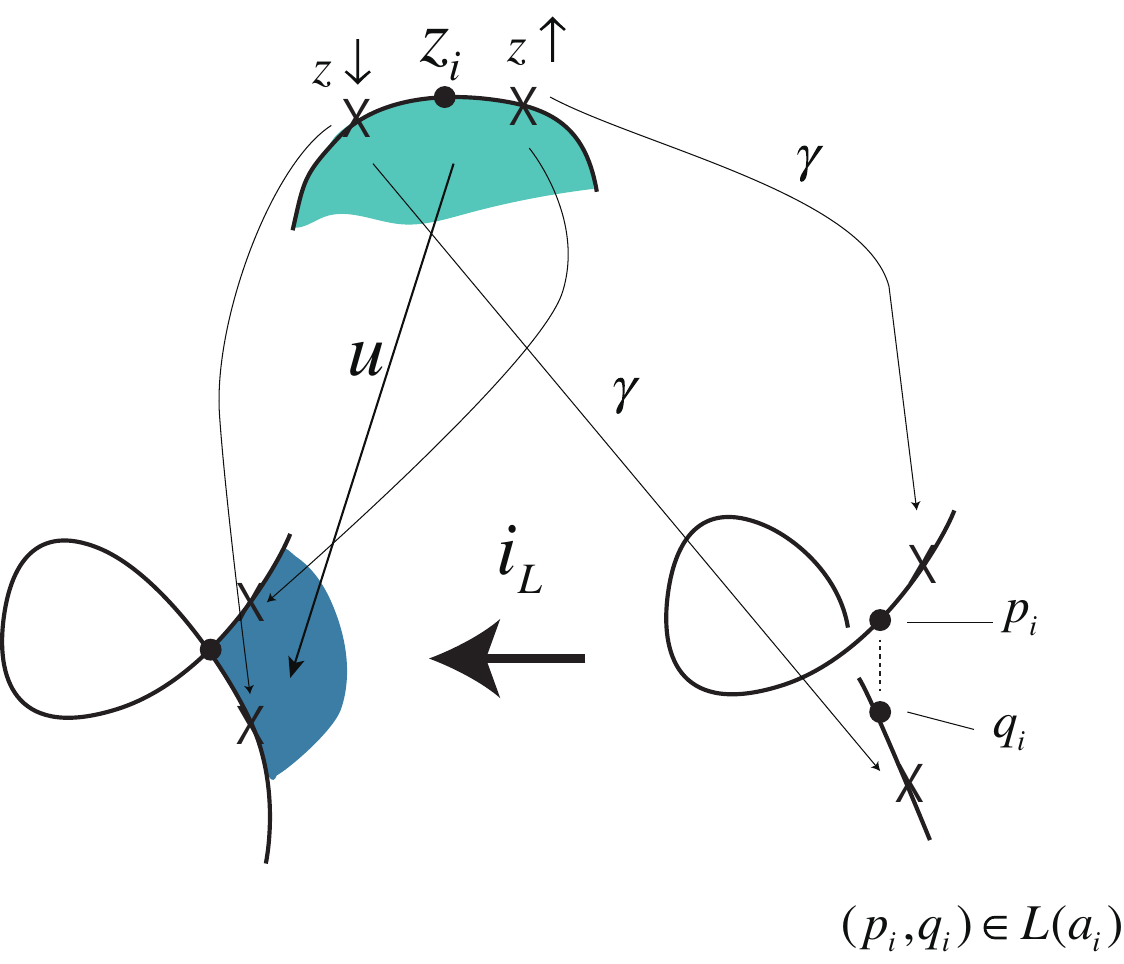}
\caption{$\lim_{z \uparrow z_i}\gamma(z)$.}
\label{Figure35}
\end{figure}
\begin{rem}
Item (4) implies that $\gamma$ extends continuously to $z_i$
if $L(a_i)$ is the diagonal component and that
$\gamma$ does not extend continuously to $z_i$
if $L(a_i)$ is a switching component.

\end{rem}
\begin{defn}\label{defn314}
Let \smash{$(\Sigma;u;\vec z;\gamma),
(\Sigma';u';\vec z^{\,\prime};\gamma') \in \mathring{\widetilde{\mathcal M}}(L;\vec a;E)$}.
We say that they are equivalent and write
$(\Sigma;u;\vec z;\gamma) \sim
(\Sigma';u';\vec z^{\,\prime};\gamma')$ if there exists a map
$v \colon D^2 \to D^2$ such that
\begin{enumerate}\itemsep=0pt
\item[(1)]
the map $v$ is biholomorphic,
\item[(2)]
$u = u'\circ v$,
\item[(3)]
$z'_i = v(z_i)$,
\item[(4)]
$\gamma = \gamma' \circ v$ on
$\partial D^2 \setminus
\{z_0,\dots,z_k\}$.
\end{enumerate}
We denote by \smash{$\mathring{\mathcal M}(L;\vec a;E)$}
\index[syindex]{M1LaE@$\mathring{\mathcal M}(L;\vec a;E)$}
the set of all the equivalence classes of this equivalence relation
$\sim$.

We define \index{evaluation map} {\it evaluation maps}
\begin{equation}\label{35form}
{\rm ev} = ({\rm ev}_0,\dots,{\rm ev}_k) \colon\
\mathring{\mathcal M}(L;\vec a;E)
\to \prod_{i=0}^k L(a_i)
\end{equation}
by
\begin{equation}\label{form3188}
{\rm ev}_i(u;\vec z;\gamma)
:=(\lim_{z \uparrow z_i}\gamma(z),\lim_{z \downarrow z_i}\gamma(z)).
\end{equation}
Here the right-hand side is as in
Definition~\ref{def3737}\,(4).

Gromov compactness implies that the set
\begin{equation}\label{defnG0(L)}
G_0(L) = \bigl\{E \in \R_{\ge 0} \mid
\exists \vec a \, \mathring{\mathcal M}(L;\vec a;E)
\ne \varnothing \bigr\}
\end{equation}
is discrete.
We define $G(L)$
\index[syindex]{GL@$G(L)$} to be the monoid generated by $G_0(L)$.
In other words, $G(L)$ is the set of all nonnegative numbers
which are sums of finitely many elements in $G_0(L)$.
The subset~${G(L) \subset \R}$ is discrete since $G_0(L)$ is discrete.

\end{defn}
We next define a compactification of \smash{$\mathring{\mathcal M}(L;\vec a;E)$}.
We first describe a combinatorial or a~topological
structure of an element in the compactification
by a tree
with additional data.
(This is a standard method used by various people in
various related situations.)

\begin{figure}[ht]
\centering
\includegraphics[scale=0.3]{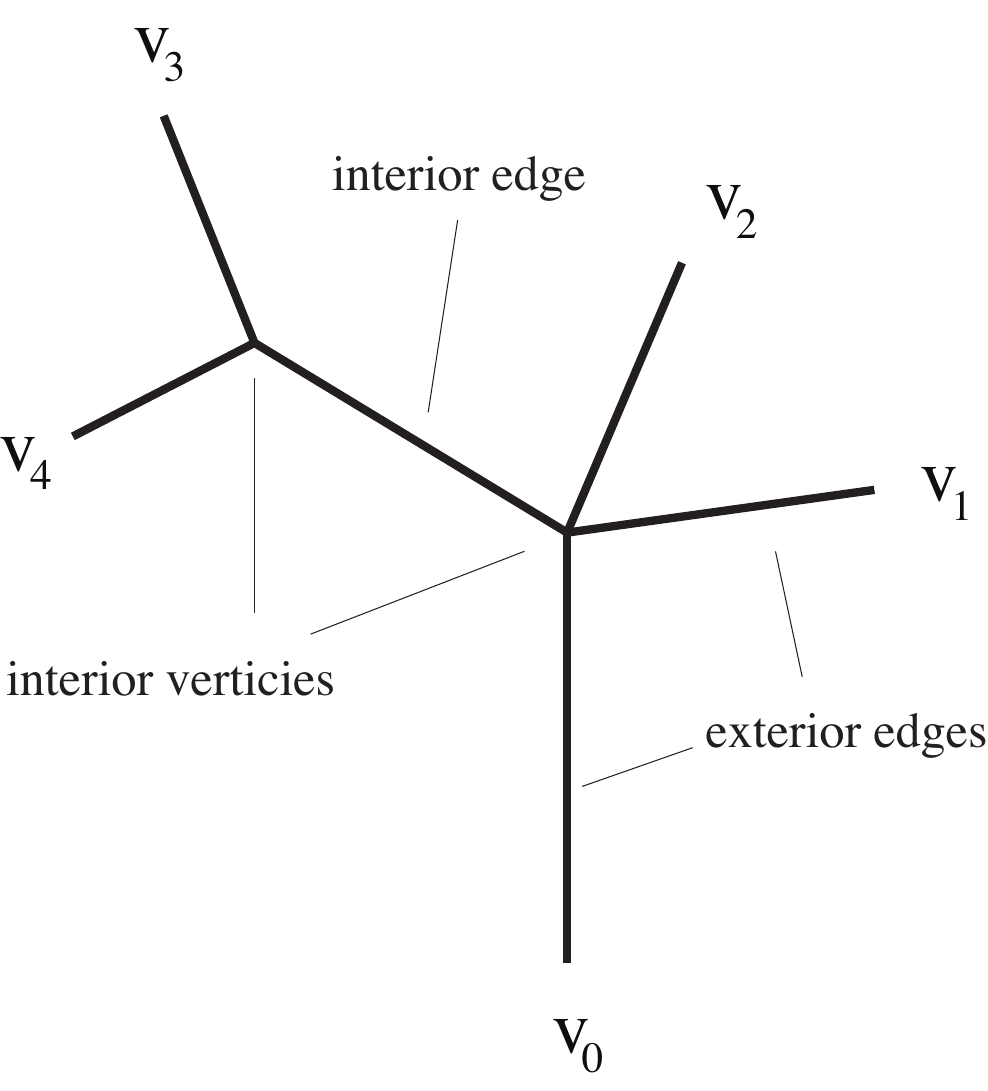}
\caption{Tree $\Gamma$.}
\label{Figure36}
\end{figure}

\begin{defn}
A {\it stable decorated ribbon tree}
\index{stable decorated ribbon tree} with $k+1$ exterior vertices and
energy $E$, which we denote by $(\Gamma,E(),a(),{\rm v}_0)$,
is a~connected tree $\Gamma$ with additional data described below.
Let~$C_0(\Gamma)$ be the set of vertices
and $C_1(\Gamma)$ the set of edges.
\begin{enumerate}\itemsep=0pt
\item[(1)]
The set $C_0(\Gamma)$ is decomposed as
$C_0(\Gamma) = C_0^{\rm int}(\Gamma)
\sqcup C_0^{\rm ext}(\Gamma)$.
\index[syindex]{C0int@$C_0^{\rm int}(\Gamma)$}
\index[syindex]{C0etx@$C_0^{\rm ext}(\Gamma)$}
We call an element of $C_0^{\rm int}(\Gamma)$
\big(resp.\ $C_0^{\rm ext}(\Gamma)$\big) \index{interior vertex}\index{exterior vertex}
an {\it interior vertex} (resp.\ {\it an exterior vertex}).
\item[(2)]
All the vertices in $C_0^{\rm ext}(\Gamma)$ have
exactly one edge containing it.
\item[(3)]
A ribbon structure of our tree $\Gamma$ is given.
In other words, an embedding $\Gamma \to \R^2$
is given up to isotopy.
\item[(4)]
The set $C_0^{\rm ext}(\Gamma)$ contains exactly $k+1$
elements.
The choice of $0$-th vertex ${\rm v}_0 \in C_0^{\rm ext}(\Gamma)$
\index{$0$-th vertex}
is given.
\item[(5)]
A map $E(\cdot) \colon C_0^{\rm int}(\Gamma) \to \R_{\ge 0}$
is given and
\[
E = \sum_{{\rm v} \in C_0^{\rm int}(\Gamma)} E({\rm v}).
\]
\item[(6)]
A map $a(\cdot) \colon C_1(\Gamma) \to \mathcal A_L$ is given.
\item[(7)] (Stability)
For each ${\rm v} \in C_0^{\rm int}(\Gamma)$,
one of the following holds:
\begin{enumerate}\itemsep=0pt
\item
$E({\rm v}) > 0$.
\item
The number of edges containing ${\rm v}$ is not smaller than $3$.
\end{enumerate}
\item[(8)]
$E({\rm v}) \in G(L)$ for any ${\rm v} \in C_0^{\rm int}(\Gamma)$.
\end{enumerate}
We denote by
$\mathscr{TR}_{k+1,E}$
\index[syindex]{TRk+1E@$\mathscr{TR}_{k+1,E}$} the set of all such $(\Gamma,E(),a(),{\rm v}_0)$.
We remark that we do not include the data (1), (3) in the notation
$(\Gamma,E(\cdot),a(\cdot),{\rm v}_0)$. However, they are included as a part of the
data which an element of $\mathscr{TR}_{k+1,E}$ comprises.
\par
We remark that $\mathscr{TR}_{k+1,E} = \varnothing$ unless $E \in G(L)$.
\par
Note that $C^{\rm ext}_0(\Gamma)$ consists of $k+1$ elements.
We enumerate them as
$
{\rm v}_0, {\rm v}_1, \dots, {\rm v}_{k}
$
so that ${\rm v}_0$ is one determined by item (4),
and the order respects the counter clockwise
orientation of $\R^2$ (into which $\Gamma$ is embedded by
using ribbon structure).
We call ${\rm v}_i$ the $i$-{\it th exterior vertex}.
\index{$i$-th exterior vertex}
(See Figure~\ref{Figure36}.)

\end{defn}
Note that we have a decomposition
$
C_1(\Gamma)
= C_1^{\rm int}(\Gamma) \sqcup C_1^{\rm ext}(\Gamma)$,
where $C_1^{\rm ext}(\Gamma)$ is the set of~${k+1}$ edges which contain one of the exterior vertices.
We call an element of $C_1^{\rm ext}(\Gamma)$
\index[syindex]{C1etx@$C_1^{\rm ext}(\Gamma)$}
an {\it exterior edge} \index{exterior edge}
and an element of $C_1^{\rm int}(\Gamma)$
an {\it interior edge}.\index{interior edge}
\index[syindex]{C1int@$C_1^{\rm int}(\Gamma)$}

We next associate a fiber product of the spaces
\smash{$\mathring{\mathcal M}(L;\vec a;E)$} to each
element of $\mathscr{TR}_{k+1,E}$.
\begin{defn}
Let $\hat\Gamma = (\Gamma,E(\cdot),a(\cdot),{\rm v}_0) \in \mathscr{TR}_{k+1,E}$.
Suppose ${\rm v} \in C^{\rm int}_0(\Gamma)$.
There exists a~unique edge ${\rm e}_0({\rm v})$
such that ${\rm e}_0({\rm v})$ lies in the same connected
component as ${\rm v}_0$ in $\Gamma \setminus {\rm v}$.
Thus, using the ribbon structure we enumerate the edges containing
${\rm v}$ as
\begin{equation}\label{edgename37}
{\rm e}_0({\rm v}), {\rm e}_1({\rm v}),\dots,
{\rm e}_{k_{{\rm v}}}({\rm v}),
\end{equation}
so they respect the counter clockwise cyclic ordering.
(See Figure~\ref{Figure37}.)
We put
\begin{equation}\label{form3838}
\vec a({\rm v})
= (a({\rm e}_0({\rm v})),\dots, a({\rm e}_{k_{{\rm v}}}({\rm v}))).
\end{equation}

\begin{figure}[ht]
\centering
\includegraphics[scale=0.3]{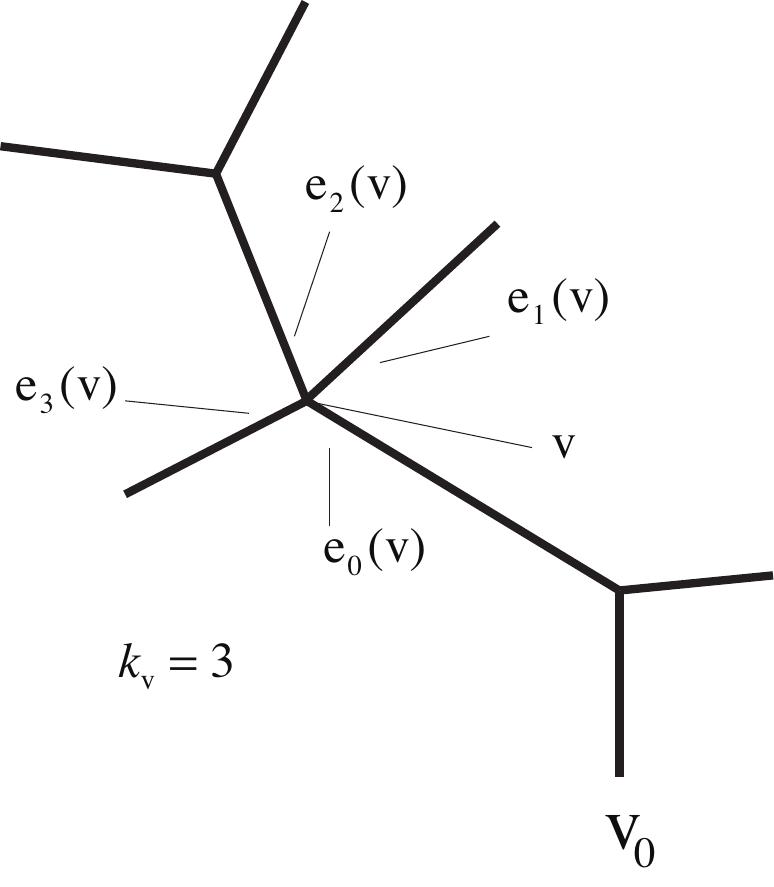}
\caption{$e_i({\rm v})$.}
\label{Figure37}
\end{figure}

We take the direct product
\begin{equation}\label{eq3888}
\prod_{{\rm v} \in C^{\rm int}_0(\Gamma)}
\mathring{\mathcal M}(L;\vec a({\rm v});E({\rm v})).
\end{equation}
We will define a map
\begin{equation}\label{form310}
\mathscr{EV} \colon\
\prod_{{\rm v} \in C^{\rm int}_0(\Gamma)}
\mathring{\mathcal M}(L;\vec a({\rm v});E({\rm v}))
\to
\prod_{{\rm e} \in C^{\rm int}_1(\Gamma)} L(a({\rm e})) \times
L(a({\rm e}))
\end{equation}
as follows.\index[syindex]{evscr@$\mathscr{EV}$}
Let ${\rm e} \in C^{\rm int}_1(\Gamma)$.
Suppose $\partial ({\rm e}) = \{{\rm v}, {\rm v}'\}$.
If ${\rm v}$ lies in the same connected component
as ${\rm v}_0$ in $\Gamma \setminus {\rm Int} {\rm e}$,
then we put $v_{\rm t}({\rm e}) = {\rm v}$.
Otherwise, ${\rm v}'$ lies in the same connected component
as~${\rm v}_0$ in $\Gamma \setminus {\rm Int} {\rm e}$.
We put $v_{\rm t}({\rm e}) = {\rm v}'$ in the latter case.

We define $v_{\rm s}({\rm e})$ such that
$\partial ({\rm e}) = \{v_{\rm s}({\rm e}), v_{\rm t}({\rm e})\}$.
(See Figure~\ref{Figure38}.)
\begin{figure}[ht]
\centering
\includegraphics[scale=0.3]{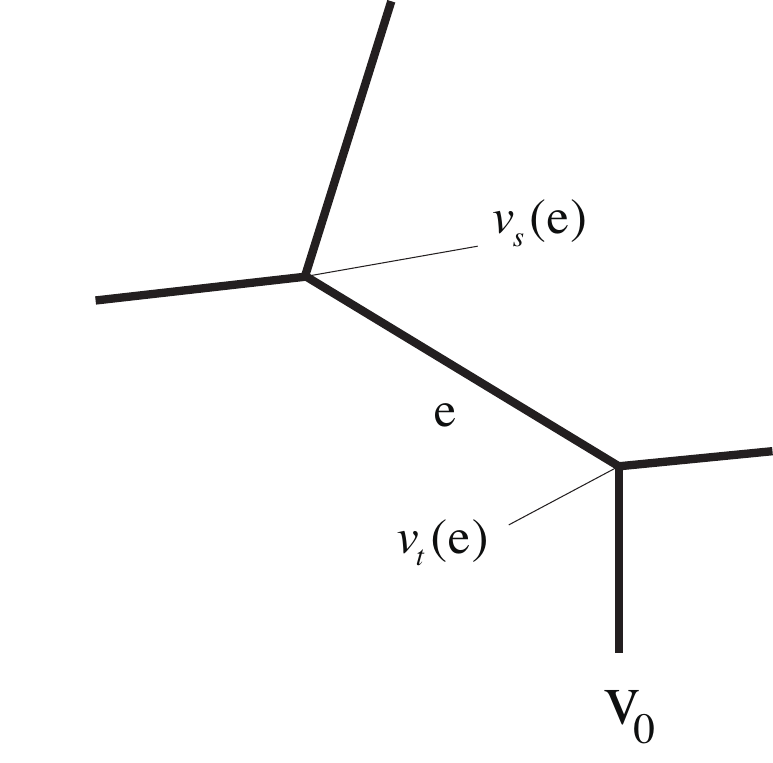}
\caption{$v_s({\rm e})$, $v_t({\rm e})$.}
\label{Figure38}
\end{figure}

Now let $\vec{{\bf x}} = \big({\bf x}_{\rm v} : {\rm v} \in C^{\rm int}_0(\Gamma)\big)$
be an element of \eqref{eq3888}.
We will define
\[
\mathscr{EV}(\vec{{\bf x}})
=
(\mathscr{EV}_{\rm e}(\vec{{\bf x}}):{\rm e} \in C^{\rm int}_1(\Gamma)).
\]
Here
$\mathscr{EV}_{\rm e}(\vec{{\bf x}})
\in L(a({\rm e})) \times L(a({\rm e}))$.
Let ${\rm e} \in C_1(\Gamma)$.
Then there exists $k_s$ and $k_t$ such that
\[
{\rm e} = {\rm e}_{k_s}(v_{\rm s}({\rm e}))
= {\rm e}_{k_t}(v_{\rm t}({\rm e})).
\]
(Actually $k_s = 0$.)
We define
$
\mathscr{EV}_{\rm e}(\vec{{\bf x}})
:=
(
{\rm ev}_{k_s}({\bf x}_{{\rm v}_s}),
{\rm ev}_{k_t}({\bf x}_{{\rm v}_t})
)$,
where ${\rm ev}_{k_s}$, ${\rm ev}_{k_t}$ are the evaluation maps \eqref{35form}.

Now we define
\[
\mathring{\mathcal M}\big(L;\hat\Gamma\big)
:=
\prod_{{\rm v} \in C^{\rm int}_0(\Gamma)}
\mathring{\mathcal M}(L;\vec a({\rm v});E({\rm v}))
{}_{\mathscr{EV}}\times
\prod_{{\rm e} \in C^{\rm int}_1(\Gamma)} \Delta_{L(a({\rm e}))}.
\]
Here\index[syindex]{M1Lgamma@$\mathring{\mathcal M}\big(L;\hat\Gamma\big)$}
$\Delta_{L(a({\rm e}))} \cong L(a({\rm e}))
\subset L(a({\rm e})) \times L(a({\rm e}))$
is the diagonal and the fiber product is taken over~${\prod_{{\rm e} \in C_1(\Gamma)} L(a({\rm e})) \times
L(a({\rm e}))}$.
See Figures \ref{Figure39-1} and \ref{Figure39-2}.

Let ${\rm e}_i$ be the unique edge containing ${\rm v}_i$.
We then put
\[
a_i\bigl(\hat\Gamma\bigr) := a({\rm e}_i),
\qquad
\vec a\bigl(\hat\Gamma\bigr)
:= \bigl(a_0\bigl(\hat\Gamma\bigr), a_1\bigl(\hat\Gamma\bigr), \dots, a_k\bigl(\hat\Gamma\bigr)\bigr).
\]

We put
$\mathscr{TR}_{E,\vec a} = \bigl\{\hat\Gamma \in \mathscr{TR}_{k+1,E}
\mid \vec a\bigl(\hat\Gamma\bigr) = \vec a\bigr\}$ and denote
\index[syindex]{TREa@$\mathscr{TR}_{E,\vec a}$}
\begin{equation}\label{def33314}
{\mathcal M}(L;\vec a;E)
:=
\coprod_{\hat\Gamma \in \mathscr{TR}_{E,\vec a}
}
\mathring{\mathcal M}\big(L;\hat\Gamma\big).
\end{equation}
\index[syindex]{M1LaE@${\mathcal M}(L;\vec a;E)$}
Moreover, we put
\begin{equation}\label{form317}
{\mathcal M}_{k+1}(L;E)
:=
\coprod_{\vec a \in (\mathcal A_L)^{k+1}}
{\mathcal M}(L;\vec a;E).
\end{equation}
\end{defn}
\index[syindex]{M1k+1LE@${\mathcal M}_{k+1}(L;E)$}
\begin{figure}[ht]
\centering
\includegraphics[scale=0.3]{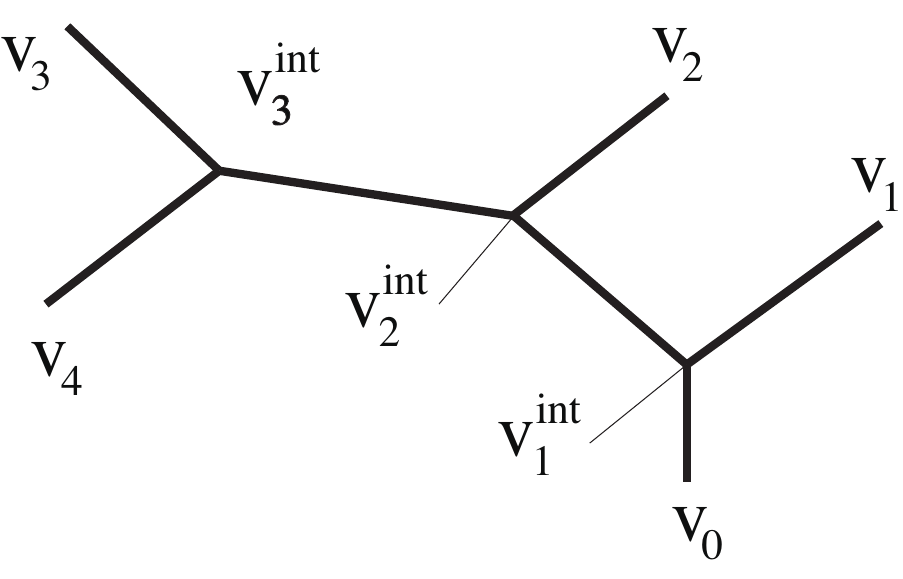}
\caption{The graph $\Gamma$.}
\label{Figure39-1}
\end{figure}
\begin{figure}[ht]
\centering
\includegraphics[scale=0.3]{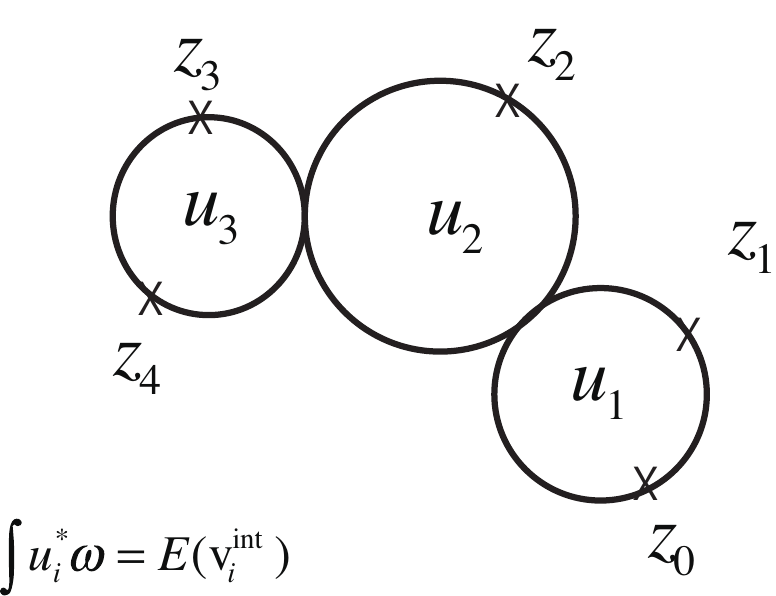}
\caption{An element of $\mathring{\mathcal M}\big(L;\hat\Gamma\big)$.}
\label{Figure39-2}
\end{figure}
\begin{defn}\label{defn32222}
Let $\vec a = (a_0,\dots,a_k) \in (\mathcal A_L)^{k+1}$.
We put $L(\vec a) = L(a_0) \times \dots \times L(a_k)$.
We define an {\it evaluation map} \index{evaluation map}
\begin{equation}\label{form3666}
{\rm ev} = ({\rm ev}_0,\dots,{\rm ev}_k)\colon\
{\mathcal M}(L;\vec a;E) \to L(\vec a)
\end{equation}
as follows.
Let $\hat\Gamma \in \mathscr{TR}_{E,\vec a}$,
${\rm v}_i$ its $i$-th exterior vertex and ${\rm e}_i$
the edge containing ${\rm v}_i$.
In other words, ${\rm e}_i$ is the $i$-th exterior edge.
Let ${\rm v}'_i$ be the other vertex of ${\rm e}_i$.
There exists $j_i$ such that ${\rm e}_i$ is the $j_i$-th
edge of ${\rm v}'_i$.
(Here we enumerate the edges of ${\rm v}'_i$ as in \eqref{edgename37}.)
For $\vec{{\bf x}} = \bigl({\bf x}_{\rm v} : {\rm v} \in C^{\rm int}_0(\Gamma)\bigr)$,
we put
$
{\rm ev}_i(\vec{{\bf x}}) := {\rm ev}_{j_i}({\bf x}_{{\rm v}'_i})$.

Using \eqref{form3666}, we define
\[
\mathscr{EV}\colon\
\prod_{{\rm v} \in C^{\rm int}_0(\Gamma)}
{\mathcal M}(L;\vec a({\rm v});E({\rm v}))
\to
\prod_{{\rm e} \in C^{\rm int}_1(\Gamma)} L(a({\rm e})) \times
L(a({\rm e}))
\]
in the same way as \eqref{form310}.
We put
\begin{equation}\label{form3737}
{\mathcal M}\big(L;\hat\Gamma\big)
:=
\prod_{{\rm v} \in C^{\rm int}_0(\Gamma)}
{\mathcal M}(L;\vec a({\rm v});E({\rm v}))
{}_{\mathscr{EV}}\times
\prod_{{\rm e} \in C^{\rm int}_1(\Gamma)} \Delta_{L(a({\rm e}))}.
\end{equation}
This is a compactification\index[syindex]{M1LHATGAMMA@${\mathcal M}\big(L;\hat\Gamma\big)$}
 of \smash{$\mathring{\mathcal M}\big(L;\hat\Gamma\big)$}.

\end{defn}
We remark that we can also define ${\mathcal M}_{k+1}(L;E)$
or ${\mathcal M}(L;\vec a;E)$
as the set of the stable maps~$(\Sigma,u,\vec z,\gamma)$
with certain properties similar to Definition~\ref{def3737},
which we omit.
(See \cite[Definition~4.2]{AJ}.)
Then we can define a stable map topology on it in the same way as
\cite[Definitions~7.1.39 and 7.1.42]{fooobook2}
and \cite[Definition~10.3]{FO}.
(See also Section~\ref{sec:difcompex}.)
\begin{thm}\label{thm323}
The spaces
${\mathcal M}_{k+1}(L;E)$ and ${\mathcal M}(L;\vec a;E)$ are compact and Hausdorff.

\end{thm}
The proof is the same as the proof of \cite[Lemma 10.4 and Theorem 11.1]{FO},
\cite[Theorem~7.1.43]{fooobook2} and is now standard.
\par
\begin{thm}\label{thekuraexist}
The spaces ${\mathcal M}(L;\vec a;E)$ for various $\vec a$, $E$ have Kuranishi structures with
corners, which enjoy the following properties:
\begin{enumerate}\itemsep=0pt
\item[$(1)$]
The codimension $m$ normalized corner, which we denote by
$S_m{\mathcal M}(L;\vec a;E)$,
of ${\mathcal M}(L;\vec a;E)$ is
identified with the disjoint union of
${\mathcal M}\big(L;\hat\Gamma\big)$, where $\hat \Gamma$
is an element of $\mathscr{TR}_{E,\vec a}$ such that~${\# C_0^{\rm int}(\Gamma) = m+1}$.
\item[$(2)$]
The map \eqref{form3666} is the underlying continuous map
of a strongly smooth map.\footnote{See \cite[Definition 3.35\,(4)]{fooonewbook}.} Moreover, ${\rm ev}_0$~is weakly submersive.\footnote{See \cite[Definition 3.35\,(5)]{fooonewbook}.}
\item[$(3)$]
The induced Kuranishi structure on
${\mathcal M}\big(L;\hat\Gamma\big)
\subseteq S_m{\mathcal M}(L;\vec a;E)$\footnote{See \cite[Proposition 24.16]{fooonewbook}.}
is isomorphic to the fiber product Kuranishi structure
\eqref{form3737}.
\item[$(4)$]
The isomorphism in item $(3)$ satisfies the corner compatibility
conditions, Condition {\rm\ref{cccond}}, below.
\item[$(5)$]
The Kuranishi structures are compatible with the forgetful maps of marked points
corresponding to the diagonal component, in the sense of {\rm\cite[\emph{Definition} 3.1]{fooo091}}.
\end{enumerate}

\end{thm}
\begin{rem}
The notion of a normalized corner is defined in \cite{foootech22,fooonewbook}.
See also~\cite{jo}.
For example, the normalized boundary of $[0,\infty)^2$
is the disjoint union of two copies of $[0,\infty)$.
Note that the two elements $0$ of the two copies of $[0,\infty)$ correspond to the
same point $(0,0)$ in $[0,\infty)^2$
but are different in the normalized corner (boundary). This is the point
where the notion of a~{\it normalized} corner (boundary)
\index{normalized corner}\index{normalized boundary} is different from
the notion of a corner (boundary).
See Figure~\ref{Figure310}.
\begin{figure}[ht]
\centering
\includegraphics[scale=0.35]{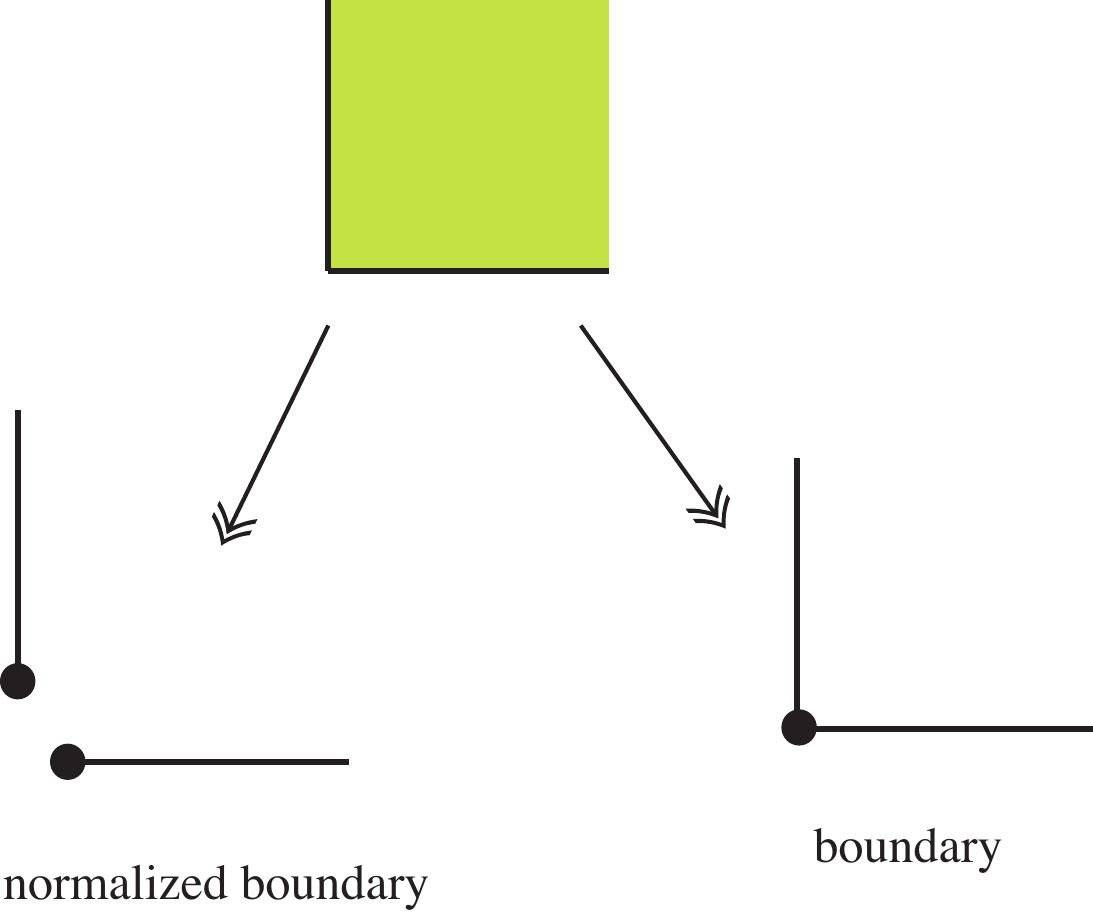}
\caption{Normalized boundary.}
\label{Figure310}
\end{figure}

\end{rem}
We describe the corner compatibility conditions.
We need some digression and discuss graph insertion.
\begin{defn}
Let $\hat\Gamma = (\Gamma,E(\cdot),a(\cdot),{\rm v}_0) \in \mathscr{TR}_{E,\vec a}$.
We assume that for each ${\rm v} \in C^{\rm int}_0(\Gamma)$
we have an element
\smash{$\hat\Gamma_{\rm v} = (\Gamma_{\rm v},
E_{\rm v}(\cdot),a_{\rm v}(\cdot),({\rm v}_{\rm v})_0) \in
\mathscr{TR}_{E({\rm v}),\vec a({\rm v})}$}.
Here $\vec a({\rm v})$ is as in \eqref{form3838}.

We define
\[
\hat\Gamma^{\bullet}
=
\hat\Gamma \# \bigl(\hat\Gamma_{\rm v}: {\rm v} \in C^{\rm int}_0(\Gamma)\bigr)
=
(\Gamma^{\bullet},E^{\bullet}(\cdot),a^{\bullet}(\cdot),{\rm v}^{\bullet}_0)
 \in \mathscr{TR}_{E,\vec a}\]
as follows:
\begin{enumerate}\itemsep=0pt
\item[(1)]
We put the tree $\Gamma_{\rm v}$ at the position of the vertex ${\rm v}$ of $\Gamma$.
We join $i$-th exterior edge of~$\Gamma_{\rm v}$ with the $i$-th edge of
$\Gamma$ containing $\rm v$. We perform this construction to all the
interior vertices ${\rm v}$ of $\Gamma$. We thus obtain $\Gamma^{\bullet}$.
(See Figure~\ref{Figure311}.)
\begin{figure}[ht]
\centering
\includegraphics[scale=0.34]{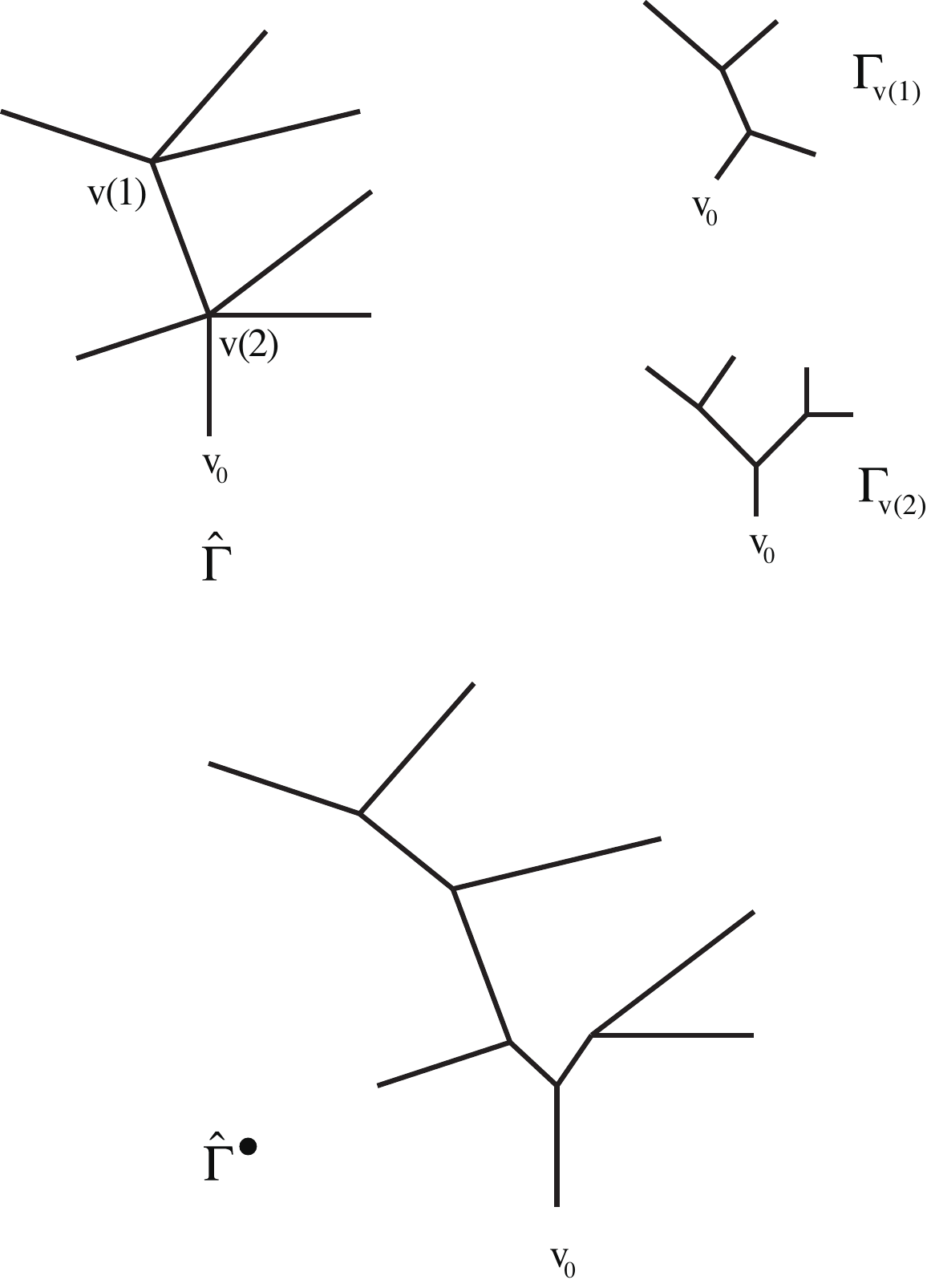}
\caption{$\hat\Gamma^{\bullet}$.}
\label{Figure311}
\end{figure}
\item[(2)] The decomposition
$C_0(\Gamma) = C_0^{\rm int}(\Gamma)
\sqcup C_0^{\rm ext}(\Gamma)$ induces
$C_0(\Gamma^{\bullet}) = C_0^{\rm int}(\Gamma^{\bullet})
\sqcup C_0^{\rm ext}(\Gamma^{\bullet})$ by
\begin{equation}\label{form380}
C^{\rm int}_0(\Gamma^{\bullet}) = \coprod_{{\rm v} \in C^{\rm int}_0(\Gamma)}
C^{\rm int}_0(\Gamma_{\rm v}).
\end{equation}
\item[(3)]
In view of \eqref{form380}, $E_{\rm v}(\cdot)$ and $a_{\rm v}(\cdot)$
induce $E^{\bullet}(\cdot)$ and $a^{\bullet}(\cdot)$, respectively.
\item[(4)]
Let ${\rm e}_0$ be the $0$-th exterior edge and ${\rm v}'_0$
the vertex of ${\rm e}_0$ such that ${\rm v}'_0 \ne {\rm v}_0$.
The $0$-th exterior vertex ${\rm v}^{\bullet}_0$ of
$\hat\Gamma^{\bullet}$ is by definition
the $0$-th exterior vertex $({\rm v}_{{\rm v}'_0})_0$ of $\hat\Gamma_{{\rm v}'_0}$
\end{enumerate}

\end{defn}

Let \smash{$\vec{{\bf x}} = \bigl({\bf x}_{\rm v} : {\rm v} \in C^{\rm int}_0(\Gamma^{\bullet})\bigr)$}
be an element of ${\mathcal M}\bigl(L;\hat\Gamma^{\bullet}\bigr)$.
(Here ${\bf x}_{\rm v} \in {\mathcal M}(L;\vec a^{\bullet}({\rm v});E^{\bullet}({\rm v}))$.)
For ${\rm v} \in C^{\rm int}_0(\Gamma)$, we use \eqref{form380}
to obtain $\vec{{\bf x}}({\rm v})$ from $\vec{{\bf x}}$.
It is easy to see that
\smash{$\vec{{\bf x}}({\rm v}) \in {\mathcal M}\bigl(L;\hat\Gamma_{\rm v}\bigr)$}.
Furthermore,
$({\bf x}({\rm v})\colon{\rm v} \in C^{\rm int}_0(\Gamma))$
is an element of ${\mathcal M}\big(L;\hat\Gamma\big)$.

Suppose $\#C_0^{\rm int}(\Gamma) = m+1$,
$\#C_0^{\rm int}(\Gamma_{\rm v}) = \ell_{\rm v}+1$
and \smash{$\ell = \sum _{\rm v \in C_0^{\rm int}(\Gamma)}\ell_{\rm v}$}.
Then
\[
\ell + m + 1 = \#C_0^{\rm int}(\Gamma^{\bullet}).
\]
Theorem~\ref{thekuraexist}\,(1) then claims
\begin{gather}
\vec{{\bf x}} \in {\mathcal M}\bigl(L;\hat\Gamma^{\bullet}\bigr)
\subseteq S_{\ell+m}({\mathcal M}(L;\vec a,E)),\label{form319}
\\
\vec{{\bf x}}({\rm v}) \in {\mathcal M}\bigl(L;\hat\Gamma_{\rm v}\bigr)
\subseteq S_{\ell_{\rm v}}({\mathcal M}(L;\vec a_{\rm v},E({\rm v}))).\label{form32020}
\end{gather}
Note that ${\mathcal M}\big(L;\hat\Gamma\big)$ is obtained as the
fiber product of ${\mathcal M}(L;\vec a_{\rm v},E({\rm v}))$.
Therefore, \eqref{form32020} implies
\begin{gather}\label{eq321}
\vec{{\bf x}} = ({\bf x}({\rm v}) : {\rm v} \in C^{\rm int}_0(\Gamma))
\in S_{\ell}\big({\mathcal M}\big(L;\hat\Gamma\big)\big).
\end{gather}
On the other hand, Theorem~\ref{thekuraexist}\,(1) claims
\begin{equation}\label{eq322}
{\mathcal M}\big(L;\hat\Gamma\big) \subseteq S_m({\mathcal M}(L;\vec a,E)).
\end{equation}
Combining \eqref{eq321} and \eqref{eq322}, we obtain
\begin{equation}\label{eq323}
{\mathcal M}\bigl(L;\hat\Gamma^{\bullet}\bigr)
\subseteq S_{\ell}(S_m({\mathcal M}(L;\vec a,E))).
\end{equation}
We have an $(\ell+m)!/\ell!m!$ fold covering map of spaces
with Kuranishi structures,
\begin{equation}\label{eq324}
S_{\ell}(S_m({\mathcal M}(L;\vec a,E)))
\to
S_{\ell+m}({\mathcal M}(L;\vec a,E))).
\end{equation}
(See \cite{foootech22}, \cite[Proposition 24.16]{fooonewbook}.)
By restricting to ${\mathcal M}\bigl(L;\hat\Gamma^{\bullet}\bigr)
\subseteq S_{\ell}(S_m({\mathcal M}(L;\vec a,E)))$
(see equation~\eqref{eq323}),
this map is a homeomorphism to its image.
Now the corner compatibility condition is stated as follows.
\begin{conds}[corner compatibility condition]\label{cccond}
\index{corner compatibility condition}
We consider two Kuranishi structures on ${\mathcal M}\bigl(L;\hat\Gamma^{\bullet}\bigr)$.
One (which we call the fiber product Kuranishi structure)
is obtained as the fiber product~\eqref{form3737}.
The other (which we call the induced Kuranishi structure) is induced from the Kuranishi structure on
${\mathcal M}(L;\vec a,E)$ by the open inclusion
\eqref{form319}.
We consider two isomorphisms between them:
\begin{enumerate}\itemsep=0pt
\item[(1)]
The isomorphism required in Theorem~\ref{thekuraexist}\,(3).
\item[(2)]
Applying Theorem~\ref{thekuraexist}\,(3) to each of $\hat\Gamma_{\rm v}$,
the inclusion in \eqref{form32020} is extended to an isomorphism between
the induced Kuranishi structure and the fiber product Kuranishi structure.
It then induces an isomorphism between the induced Kuranishi structure and the fiber product Kuranishi structure
on the space \eqref{eq323}.
By \eqref{eq324} \big(which is an isomorphism on~\smash{${\mathcal M}\bigl(L;\hat\Gamma^{\bullet}\bigr)$}\big),
it induces an isomorphism between the
induced Kuranishi structure and the fiber product Kuranishi structure
on \smash{${\mathcal M}\bigl(L;\hat\Gamma^{\bullet}\bigr)$}.
\end{enumerate}
We require that the two isomorphisms (1), (2) above coincide with each other.
\end{conds}

\begin{rem}
Condition \ref{cccond} looks rather complicated. Actually, in our geometric
situation, it is rather trivial that Condition \ref{cccond} is satisfied.
Corner compatibility conditions such as Condition \ref{cccond}
are spelled out in \cite{foootech22}, \cite[Chapters 16 and 21]{fooonewbook} for the purpose of
axiomatizing the construction of a compatible
system of perturbations of the compatible system of
Kuranishi structures.
In other words, we spelled out the properties we need to
construct a compatible system of perturbations
in a way independent of the geometric origin of the
system of Kuranishi structures.

\end{rem}

In the case when $L$ is an embedded Lagrangian submanifold,
Theorem~\ref{thekuraexist} is \cite[Propositions~7.1.1 and 7.1.2]{fooobook2}.\footnote{Note that item (4) is not stated in \cite[Propositions~7.1.1 and 7.1.2]{fooobook2}.
However, this compatibility is fairly obvious from the construction.}
Its generalization is in \cite{AJ} in the case when $L$ has transversal self-intersection.
In the general case, we can use Morse--Bott gluing which can be worked out in the
same way as \cite[Section 7.1.3]{fooobook2}.
(See also \cite{fu15,FO}.)
In fact, the analytic detail of \cite[Section~7.1.3]{fooobook2} is designed
so that it works also in the Morse--Bott case in general.
The detail of the analysis to prove Theorem~\ref{thekuraexist} is given also
in \cite[Parts 2 and 3]{foootech} and in \cite{fooo:const1, fooo:const2,foooanalysis}.

We next discuss the orientation.
We first recall the definition of orientation local systems of spaces
with Kuranishi structure.
Let $\widehat{\mathcal U} =\{(U_p,E_p,s_p,\psi_p) \mid p \in X\}$ be a Kuranishi structure of~$X$.
(We use the definition of \cite[Definition 3.8]{foootech2}. So it has a tangent bundle in the sense of~\cite[Definition A1.4]{fooobook2}.)
We obtain a principal ${\rm O}(1)$ bundle
$
O_p = \operatorname{Det} TU_p \otimes \operatorname{Det} TE_p
$
on~$U_p$.
\index[syindex]{Ozp@$O_p$} By the condition of the coordinate change in \cite[Definition 3.2\,(8)]{foootech2},
$O_q \cong \varphi_{pq}^*O_p$
and this isomorphism is compatible in the sense that the map
$
O_r \cong \varphi_{qr}^*O_q \cong
\varphi^*_{qr}\varphi_{pq}^* O_p
$
coinsides with~${
O_r \cong \varphi_{pr}^*O_p}$.
We call such collection of $\{O_p \mid p \in X\}$
together with isomorphisms~${O_q \cong \varphi_{pq}^*O_p}$
satisfying the above explained compatibility
conditions,
the {\it orientation local system}
\index{orientation local system} of our space with
Kuranishi structure \smash{$\bigl(X,\widehat{\mathcal U}\bigr)$}
and write it as \smash{$O_{(X,\widehat{\mathcal U})}$}.
\index[syindex]{OXU@$O_{\bigl(X,\widehat{\mathcal U}\bigr)}$}
We write it also as $O_X$ by an abuse of notation.

If we construct a compatible good coordinate system
$\{(U_{\mathfrak p},E_{\mathfrak p},s_{\mathfrak p},\psi_{\mathfrak p}) \mid {\mathfrak p} \in {\mathfrak P}\}$
then $\{O_p \mid p \in X\}$ induces a system of principal ${\rm O}(1)$ bundles
$\{O_{\mathfrak p} \mid {\mathfrak p} \in {\mathfrak P}\}$ which is compatible
with the coordinate change in a similar sense as above.
(Here $O_{\mathfrak p}$ is a principal ${\rm O}(1)$ bundle on $U_{\mathfrak p}$.)
We can use it to define and study integration along the fiber
in a similar way as the case of manifolds.
(See \cite[Chapter 27]{foootech2} and \cite{fooonewbook}.)

Suppose that $\mathfrak f = \{f_p \mid p \in X\}$ is a weakly submersive
strongly smooth map $\bigl(X,\widehat{\mathcal U}\bigr) \to N$ to a smooth manifold $N$.
If $\Theta$ is a principal ${\rm O}(1)$ bundle on $N$, we pull it back to each $U_p$ to obtain
$f_p^*\Theta$. They are compatible with the coordinate change
in a similar sense as above. We denote the system $\{f_p^*\Theta \mid p \in X\}$
by $\mathfrak f^*\Theta$.
We can define a tensor product of several systems
in an obvious way.

An isomorphism between $\mathfrak f^*\Theta$ and $O$ is
a system of isomorphisms of real line bundles $f_p^*\Theta
\cong O_p$ which commute with coordinate changes.

Let $\vec a = (a_0,\dots,a_k)$, $a_i \in \mathcal A^+_L$.
We use the principal ${\rm O}(1)$ bundles $\Theta^-_{a_i}$ and
$\Theta^+_{a_i}$, which are defined in Lemma--Definition~\ref{lemdef39}.

\begin{prop}\label{prop329}
The $V$-relative spin structure of $L$ canonically induces an isomorphism of
principal ${\rm O}(1)$ bundles
\begin{equation}\label{oriisom}
O_{{\mathcal M}(L;\vec a;E)}
\cong
\bigotimes_{i=0}^k {\rm ev}_i^*\Theta^-_{a_i}.
\end{equation}

\end{prop}
\begin{proof}
This is a straightforward generalization of
\cite[Proposition 8.8.6]{fooobook2}.
We provide the proof below for completeness.

Let $\mathfrak x = (\Sigma;u;\vec z;\gamma)$ be an element of
\smash{$\mathring{\widetilde{\mathcal M}}(L;\vec a;E)$}.
If suffices to consider the case when $\Sigma = D^2$.
We may also assume that the image of $\gamma$ lies in
$L_{[3]}$.

We put ${\rm ev}_i(\mathfrak x) = (p_i,q_i) \in L(a_i)$
and $i_L(p_i) = i_L(q_i) = x_i$.
We write $x_i = (p_i,q_i) \in L(a_i)$ by an abuse of notation.
We fix a trivialization $V_{x_i} \cong \R^m$
and take an element $\lambda_{x_i} \in \mathcal P^{a_i}_{x_i}$
for each~$i$. (See Definition~\ref{lem34}.)

We show that those choices together with the $V$-relative spin
structure of $L$ determine an isomorphism of the
principal ${\rm O}(1)$ bundles of the left and right-hand sides of
\eqref{oriisom} at $\mathfrak x$.

For each $i$, we use $\lambda_{x_i}$ to define an
elliptic operator
\[
\overline\partial_{Z_-,\lambda_{x_i}} \colon\ L^2_k(Z_-;T_{x_i}X;\lambda_{a_i};\delta)
\to L^2_{k-1}(Z_-;T_{x_i}X;\delta)
\]
as in \eqref{3535}.
We can glue it with the linearized operator of the defining equation of
\smash{$\mathring{\widetilde{\mathcal M}}(L;\vec a;E)$}
at~$\mathfrak x$ to obtain an elliptic operator
$P$ on $\Sigma=D^2$ whose symbol is the same as one
of the Cauchy--Riemann operator with $u^*TX$ coefficient.
Its boundary condition is given by
concatenating the family $z \in \partial \Sigma \mapsto (di_{L})(T_{\gamma(z)}\tilde L)$
with $\lambda_{x_i}$'s.
(See Figure~\ref{Figure312}.)
We denote this family of Lagrangian subspaces (the boundary condition)
by $\lambda$.
\begin{figure}[ht]
\centering
\includegraphics[scale=0.3]{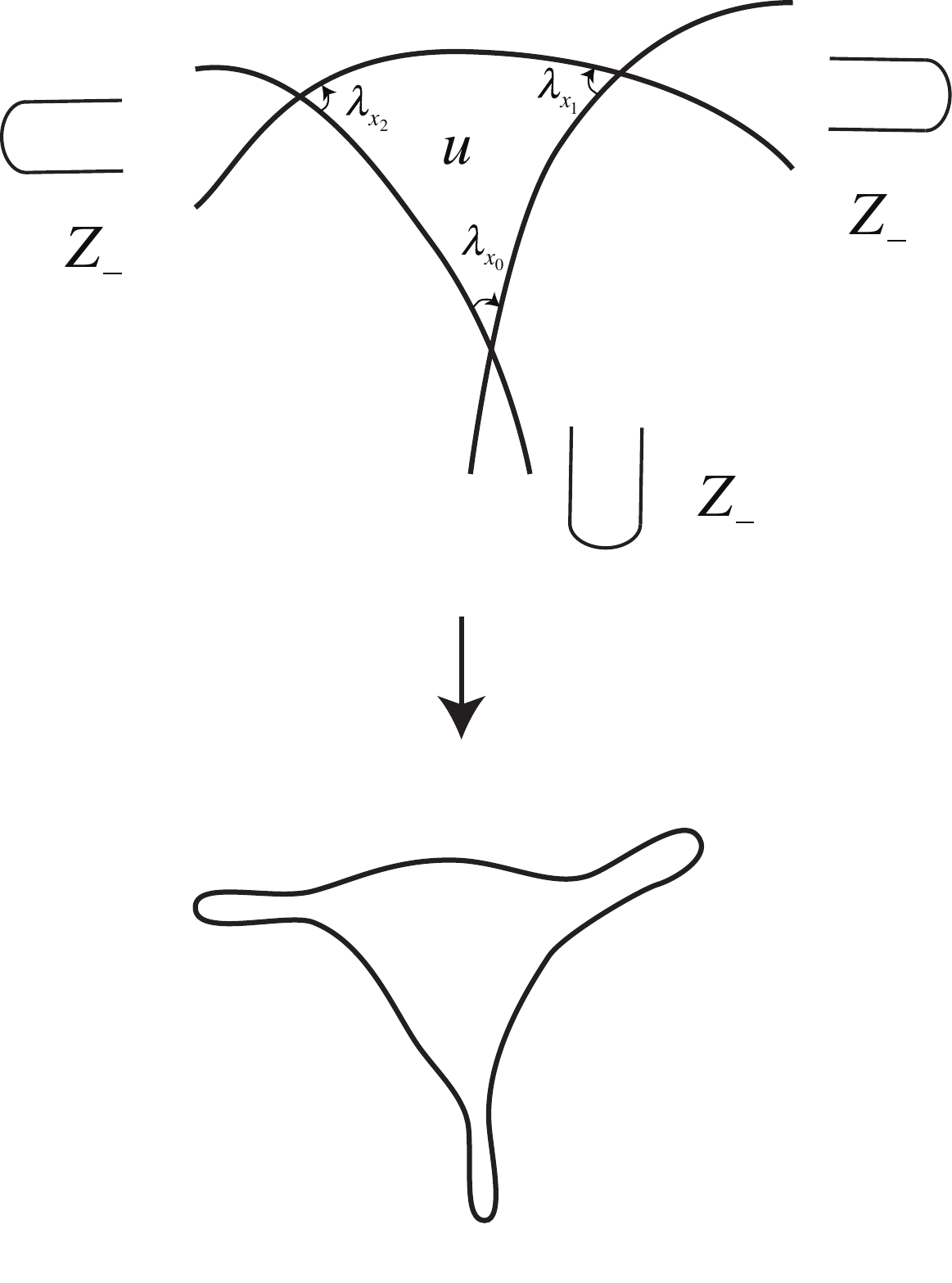}
\caption{Family of Lagrangian subspaces $\lambda$.}
\label{Figure312}
\end{figure}

We claim there is a canonical orientation of the determinant line bundle of
the index of $P$.
We prove it below.
Using the isomorphism $u^*TX \cong \C^n \times \Sigma$ \big(note that $\Sigma = D^2$\big),
we may regard $\lambda$ as an $S^1$ parametrized family of Lagrangian subspaces of $\C^n$.
The trivialization of $V$ and relative spin structure
determine a trivialization of this family of subspaces as an abstract vector bundle.
We thus have a trivial complex vector bundle $\xi_0 = \C^n \times \Sigma$ on $\Sigma = D^2$.

On the other hand, by
an identification $\xi_{0,\R} \vert_{\partial\Sigma} = \R^n \times \partial\Sigma$
with $\lambda$, (which may not be
consistent with the trivialization $\xi_0 = \C^n \times \Sigma$). This identification induces $\xi_0\vert_{\partial\Sigma} \cong u^*TX\vert_{\partial\Sigma}$.
In the same way as the proof of \cite[Theorem 8.1.1]{fooobook2},
we can show that the difference of the index of $P$ and the Cauchy--Riemann operator
of the bundle $E$ with $\xi_{0,\R}$ boundary condition has a~canonical orientation.
(This is based on the fact that this difference can be identified with an index of
a certain family of
operators on $\C P^1$ with complex linear symbols.)

The index of the Cauchy--Riemann operator
of the bundle $E$ with $E_{\R}$ as a boundary condition is canonically identified with $\R^n$.
Therefore, the determinant bundle of the index bundle of $P$
is canonically trivialized.

On the other hand, we find that
\begin{equation}\label{aaa337}
TU_{\mathfrak x} \ominus E_{\mathfrak x} \oplus \bigoplus_{i=0}^k
\operatorname{Index}\overline\partial_{Z_-,\lambda_{x_i}}
\cong
\operatorname{Index} P
\end{equation}
as a virtual vector space.
(Here $U_{\xi}$ is a Kuranishi neighborhood of $\xi$ and $E_{\xi}$
is an obstruction bundle.)

We remark \eqref{aaa337}
induces an isomorphism \eqref{oriisom} at a point $\mathfrak x$.
In fact, $\operatorname{Det} \operatorname{Index}  P$ is trivial.
Therefore, we obtain an isomorphism $\operatorname{Det}  \operatorname{Index}\overline\partial_{Z_-,\lambda_{x_i}}
\cong \Theta_{x_i}^-$ by Lemma--Definition~\ref{lemdef39}\,(3), and an isomorphism $\operatorname{Det} TU_{\mathfrak x} \otimes \operatorname{Det} E_{\mathfrak x}^*
\cong O_{{\mathcal M}(L;\vec a;E)}$ by definition.
Thus we obtain \eqref{oriisom}.

We next explain the way to obtain a family of isomorphisms $V_{x_i} \cong \R^m$ and of $\lambda_{x_i}$
and so that the above isomorphisms induce a global isomorphism \eqref{oriisom}.

In fact, the independence of the choice of $\lambda_{x_i}$
is the consequence of Lemma--Definition~\ref{lemdef39}\,(2).

Let us discuss the dependence of the identification $V_{x_i} \cong \R^m$.
We first remark that to prove Proposition~\ref{prop329}
it suffices to prove this isomorphism on each loop of the domain,
since both sides are principal ${\rm O}(1)$ bundles.
Let $S^1 \to C^{\infty}\big(\big(D^2,\partial D^2\big),(X,L)\big)$ be a smooth map.
It induces a map
$\big(S^1\times D^2,S^1\times \partial D^2\big) \to (X,L)$.
The pullback of $V$ by this map is a trivial bundle since it is
an oriented real bundle on $S^1\times D^2$.
Therefore, we have a continuous family of isomorphisms~${V_{x_i} \cong \R^m}$
on this $S^1$ parametrized family.

The proof of Proposition~\ref{prop329} is complete.
\end{proof}

\subsection[The filtered $A_\infty$ algebra associated to an immersed
Lagrangian submanifold]{The filtered $\boldsymbol{A_{\infty}}$ algebra associated\\ to an immersed
Lagrangian submanifold}
\label{subsec:Ainfalgim}

We now use Theorem~\ref{thekuraexist} to prove Theorem~\ref{AJtheorem}.
We refer \cite[Definition 9.1]{foootech2} and \cite{fooonewbook} for the definition of CF-perturbations
\index{CF-perturbations} on
Kuranishi structures.

\begin{prop}\label{prop330}
Let $E_0 > 0$. Then there exists a system of CF-perturbations
\smash{${\widehat{\mathfrak S}}$}
\index[syindex]{Sfrak@${\widehat{\mathfrak S}}$} on the moduli spaces
${\mathcal M}(L;\vec a;E)$ with Kuranishi structures
which are outer collarings\footnote{See \cite[Chapter 17]{fooonewbook} for the definition of an outer collaring
(it was called $\tau$-collaring in \cite{foootech22}).}\index{outer collaring} \index{thickening}
of thickenings of the structures given in Theorem {\rm\ref{thekuraexist}}, for various $\vec a$, $E$ with
$E< E_0$.
It enjoys the following properties
$($see {\rm\cite[\emph{Definition} 5.3]{foootech2} \emph{and} \cite{fooonewbook}} for the definition of a thickening$)$:
\index{thickening}
\begin{enumerate}\itemsep=0pt
\item[$(1)$]
Each of \smash{${\widehat{\mathfrak S}}$} is transversal to zero.
\item[$(2)$]
The evaluation map ${\rm ev}_0$ is strongly submersive \index{strongly submersive} with respect to
this CF-perturbation
$($see {\rm\cite[\emph{Definition} 9.2]{foootech2} \emph{and} \cite{fooonewbook}} for the definition of strong submersivity$)$.
\item[$(3)$]
They are compatible at the corners in the following sense.
We consider the left-hand side~${\mathcal M}\big(L;\hat\Gamma\big)$ of \eqref{form3737}
and require that the following two CF-perturbations on it coincide each other.
\begin{enumerate}\itemsep=0pt
\item[$(a)$]
The space ${\mathcal M}\big(L;\hat\Gamma\big)$ is a stratum of
${\mathcal M}(L;\vec a;E)$
with respect to its corner structure stratification \index{corner structure stratification}
$($see {\rm\cite[\emph{Definition} 4.15]{foootech2} \emph{and} \cite{fooonewbook}} for the definition of the corner structure stratification$)$.
We restrict CF-perturbation \smash{${\widehat{\mathfrak S}}$} on ${\mathcal M}(L;\vec a;E)$
to ${\mathcal M}\big(L;\hat\Gamma\big)$ and obtain a CF-perturbation on it.
\item[$(b)$]
The right-hand side of \eqref{form3737} is a fiber product of
various connected components of ${\mathcal M}(L;\vec a;E)$.
We take the restriction of \smash{${\widehat{\mathfrak S}}$} to the moduli spaces
appearing as the fiber product factors of the right-hand side of \eqref{form3737}
and take the fiber product CF-perturbation, in the sense of {\rm\cite[\emph{Definition} 10.13]{foootech2} \emph{and} \cite{fooonewbook}}.
Since ${\rm ev}_0$ is strongly submersive, we can take the
fiber product CF-perturbation $(${\rm \cite[\emph{Lemma--Defini\-tion}~10.12]{foootech2} \emph{and}~\cite{fooonewbook}}$)$.
\end{enumerate}
\item[$(4)$]
They are compatible with the forgetful map of the marked points which
correspond to the diagonal component other than $0$-th one. The
precise definition of compatibility is written in {\rm\cite[\emph{Definition} 5.1]{fooo091}}.
\end{enumerate}

\end{prop}
\begin{proof}
The proof is mostly the same as the proof of \cite[Section 5]{fooo091},
\cite{foootech2}, \cite[Chapter~12]{fooonewbook} and~\cite{foootech22}, \cite[Chapter 17]{fooonewbook}.
We explain only the points where the discussion is slightly different.

We first observe that it suffices to define a CF-perturbation of
the space ${\mathcal M}(L;\vec a;E)$ such that~$a_i$ is not the diagonal component for $i\ne 0$.
In fact, then the CF-perturbation in the general case is
automatically determined by item (4).

We then remark the following.
\begin{lem}\label{Lema31}
There exits $k_0$ $($depending on $E_0)$ such that
the following holds. Let $\vec a = (a_0,\dots,a_k)$.
Suppose ${\mathcal M}(L;\vec a;E) \ne \varnothing$,
$a_i$ is not the diagonal component for $i\ne 0$, and $E < E_0$ then $k \le k_0$.
\end{lem}
\begin{proof}
This is a direct consequence of Gromov compactness.
\end{proof}

Thus we need to construct CF-perturbations on only finitely many
spaces with Kuranishi structures.

The rest of the proof is the same as
\cite[Section 5]{fooo091}, \cite{foootech2, foootech22, fooonewbook}.
The construction is by induction on $E$.
Suppose we have constructed CF-perturbations
with the required properties for ${\mathcal M}(L;\vec a;E)$
with $E < E_1 <E_0$.
We will construct one for
${\mathcal M}(L;\vec a;E_1)$.
By the induction hypothesis Proposition~\ref{prop330}\,(3),
the boundary and corners of ${\mathcal M}(L;\vec a;E_1)$
are fiber products of the moduli spaces for which
CF-perturbations are already defined by the induction hypothesis.
We take the fiber product of those CF-perturbations
to obtain CF-perturbations of the boundary
and corners of ${\mathcal M}(L;\vec a;E_1)$
that are compatible with each other.
Therefore, by using the relative version of the existence theorem
of CF-perturbations (see \cite[Proposition 17.65 or 15.7]{foootech22} or \cite[Proposition 17.81 or
15.7]{fooonewbook}),
we can extend it to ${\mathcal M}(L;\vec a;E_1)$.
The proof is now complete by induction.
\end{proof}

We use the CF-perturbations obtained in Proposition~\ref{prop330} to
define the structure operations of our filtered $A_{\infty}$ structure.

\begin{defn}\label{def32}
\quad
\begin{enumerate}\itemsep=0pt
\item[(1)]
Let $E < E_0$, $E \in G(L)$. For $(E,k) \ne (0,1)$, we define multi-linear maps\index[syindex]{m2kEepsilon@$\mathfrak m_{k}^{E,\varepsilon}$}
\[
\mathfrak m_{k}^{E,\varepsilon}\colon\ CF(L;\R)^{\otimes k} \to CF(L;\R)
\]
by
\begin{equation}\label{form38}
\mathfrak m_{k}^{E,\varepsilon}(h_1,\dots,h_k)
:=
(-1)^*{\rm ev}_0! \big({\rm ev}_1^* h_1 \times \dots \times {\rm ev}_k^* h_k ; {\widehat{\mathfrak S^{\varepsilon}}}\big).
\end{equation}
Here $CF(L;\R) = \Omega(L;\Theta^-)$ (see Definition~\ref{defn313} and \eqref{form315}).
Note that $h_1,\dots,h_k \in CF(L;\R)$ and
${\rm ev}_i^* h_i$ are the pullbacks of differential forms with respect to the
strongly smooth map~${\rm ev}_i$ (see \cite[Definition 7.70]{foootech2} and~\cite{fooonewbook}).
\smash{${\rm ev}_0! \bigl(*; {\widehat{\mathfrak S^{\varepsilon}}}\bigr)$} is the integration along the
fiber \index{integration along the
fiber} of the differential form with respect to the CF-perturbation
(see \cite[Definition 9.13]{foootech2} and \cite{fooonewbook}). It depends on a
positive number $\varepsilon$
(see Remark~\ref{newnewrem} below).
Here we consider the moduli spaces ${\mathcal M}(L;\vec a;E)$
and their Kuranishi structures and the CF-perturbations
obtained in Theorem~\ref{thekuraexist} and Proposition~\ref{prop330}
to define them.

The sign $*$ in \eqref{form38} is
\smash{$
* = \sum_{i=1}^k i(\deg h_i + 1)+1$} when we take the convention of \cite[p.~552]{fooonewbook}. (The
same correction term also appears in \cite[Section~2.2.2]{ST}.)
(However, in this paper we do not use this particular formula of $*$ as we will mention in Remark~\ref{rem172}.)

We remark that by Lemma~\ref{lem333} below the right-hand side of equation~\eqref{form38} is an element of~$CF(L;\R)$.

We also define
\begin{equation}\label{form3420000}
\mathfrak m_{1}^0(h) := dh.
\end{equation}
There is no sign when we use the convention of \cite[Definition 21.29\,(4)]{fooonewbook}.
\item[(2)]
We define
$
\mathfrak m_{k}^{<E_0,\varepsilon} \colon CF(L;\Lambda_0)^{\otimes k} \to CF(L;\Lambda_0)
$
by\index[syindex]{m2k<Eepsilon@$\mathfrak m_{k}^{<E_0,\varepsilon}$}
\[
\mathfrak m_{k}^{<E_0,\varepsilon}
:=
\sum_{E < E_0, E \in G(L)} T^{E}\mathfrak m_{k}^{E,\varepsilon}.
\]
\end{enumerate}

\end{defn}
\begin{rem}\label{newnewrem}
Here $\varepsilon$ is a sufficiently small positive number.
It is proved in \cite{foootech2}, \cite[Theorem~9.15]{fooonewbook}
that the integration along the
fiber \smash{${\rm ev}_0!\big(\cdots; {\widehat{\mathfrak S^{\varepsilon}}}\big)$}
is independent of various choices such as partition of unity,
if $\varepsilon$ is sufficiently small. (The integration along the
fiber depends on $\varepsilon$ and the CF-perturbation.)
How much $\varepsilon$ should be small for this well-definedness to hold
is also CF-perturbation dependent. Note that, however, for a fixed $E_0$, we have only
finitely many moduli spaces to perturb. Therefore, we can take $\varepsilon_0$,
which is $E_0$ dependent, so that the integration along the
fiber is well-defined
if $\varepsilon < \varepsilon_0$ for those finitely many moduli spaces and their CF-perturbations.

\end{rem}
\begin{lem}\label{lem333}
The right-hand side of
\eqref{form38} is an element of $CF(L;\R)$.

\end{lem}
\begin{proof}
Note that we may decompose $\mathfrak m_{k}^{E,\varepsilon}$ to the sum
$
\mathfrak m_{k}^{E,\varepsilon} = \sum_{\vec a} \mathfrak m_{\vec a}^{E,\varepsilon}$
where $\mathfrak m_{\vec a}^{E,\varepsilon}$ is defined by
${\mathcal M}(L;\vec a;E)$.

We may assume $h_i \in \Omega(L(a_i);\Theta^-_{a_i})$.
Then $\mathfrak m_{\vec a}^{E,\varepsilon}(h_1,\dots,h_k)$
is nonzero for $\vec a = (a_0,a_1,\dots,a_k)$ with $a_0 \in \mathcal A(L)$.
We consider the following two cases separately.

Case 1: $L(a_0)$ is a switching component.

We define an involution $\tau \colon \tilde L \times_X \tilde L \to \tilde L \times_X \tilde L$
by $\tau(p,q) = (q,p)$.
We take $a'_0 \in \mathcal A(L)$ such that
$\tau(L(a_0)) = L(a'_0)$.
By definition it is easy to see that
\smash{$
\tau^*(\Theta^-_{a'_0}) = \Theta^+_{a_0}.
$}
Therefore, by Lemma--Definition~\ref{lemdef39} we have
\begin{equation}\label{form341}
\tau^*(\Theta^-_{a'_0}) = \Theta^-_{a_0} \otimes \operatorname{Det} TL(a_0).
\end{equation}
Proposition~\ref{prop329} implies that for $h_0 \in \Omega(L(a_0);\Theta^-_{a_0})$
we can define the integration
\[
\int_{({\mathcal M}(L;\vec a;E),{\widehat{\mathfrak S^{\varepsilon}}})}
{\rm ev}_1^*h_1 \times \dots \times {\rm ev}_k^*h_k
\times {\rm ev}_0^*h_0 \in \R.
\]
In other words, we may regard
\[
\mathfrak m_{\vec a}^{E,\varepsilon}(h_1,\dots,h_k)
\in \Omega(L(a_0);\Theta^-_{a_0} \otimes \operatorname{Det} TL(a_0)).
\]
Therefore, by \eqref{form341} we may regard
\begin{equation}\label{form342}
\mathfrak m_{\vec a}^{E,\varepsilon}(h_1,\dots,h_k)
\in \Omega(L(a'_0);\Theta^-_{a'_0}),
\end{equation}
as required.

Case 2: $L(a_0)$ is the diagonal component.

In this case, $L(a_0) \cong \tilde L$ is oriented.
By definition $a'_0 = a_0$.
Moreover, $\Theta^-_{a_0}$ and $\Theta^+_{a_0}$ are both
trivial bundles. We can prove \eqref{form342} easily
in this case, by using
Proposition~\ref{prop329}.
\end{proof}

\begin{prop}\label{prop334}
\smash{$\mathfrak m_{k}^{<E_0,\varepsilon}$, $k=0,1,\dots,$}
defines a filtered $A_{\infty}$ structure modulo $T^{E_0}$.
\index{filtered $A_{\infty}$ structure modulo $T^{E_0}$}
Namely, we have
\begin{gather*}
0\equiv\sum_{k_1+k_2=k+1}\sum_{i=0}^{k_1-1}
(-1)^{*_i} \mathfrak m^{<E_0,\varepsilon}_{k_1}(x_1,\dots,x_i,
\mathfrak m^{<E_0,\varepsilon}_{k_2}(x_{i+1},\dots,x_{i+k_2}),
\dots,x_k)\quad
\mod T^{E_0}
\end{gather*}
for sufficiently small $\varepsilon >0$.
Here ${*_i} = \deg' x_1 + \dots + \deg' x_i$.
Moreover, $1 = [L(a_0)] \in CF(L(a_0),\R)$
$($the differential form $($function$)$ $1$ on the diagonal component$)$
is a unit.
\end{prop}
\begin{proof}
The proof is now a routine using Proposition~\ref{prop330},
Stokes' formula \cite[Proposition~9.26]{foootech2}, \cite{fooonewbook} and
the composition formula
\cite[Theorem 10.20]{foootech2}, \cite{fooonewbook}
and proceed as follows.

It suffices to show
\begin{equation}\label{form345}
\sum_{E_1+E_2 = E} \sum_{k_1+k_2=k+1}\sum_{i=1,\dots,k_1} (-1)^*
\mathfrak m_{k_1}^{E_1,\varepsilon}\big(h_1,\dots, \mathfrak m_{k_2}^{E_2,\varepsilon}(h_{i+1},\dots,
h_{i+k_2}),\dots,h_k\big) = 0,
\end{equation}
with $* = \deg'h_1 + \dots + \deg' h_i$.
We denote by ${\rm ev}_0!\big(*;\big(\mathcal M,\widehat{\mathfrak S^{\varepsilon}}\big)\big)$
the integration along the fiber of the differential form
defined by a CF-perturbation \smash{$\widehat{\mathfrak S^{\varepsilon}}$}
of the space $\mathcal M$ with Kuranishi structure.
Now by Stokes' theorem (see \cite[Proposition 9.26]{foootech2} and \cite{fooonewbook}) and the definition, we have
\begin{gather}
\big(d \circ \mathfrak m_{\vec a}^{E,\varepsilon}\big)(h_1,\dots,h_k)
+ \sum_{i=1}^k (-1)^*\mathfrak m_{\vec a}^{E,\varepsilon}(h_1,\dots,d h_i,\dots,h_k) \nonumber\\
\qquad=
{\rm ev}_0! \big({\rm ev}_1^* h_1 \times \dots \times {\rm ev}_k^* h_k ;
\big(\partial {\mathcal M}(L;\vec a;E),{\widehat{\mathfrak S^{\varepsilon}}}\big)\big).\label{form346}
\end{gather}
Here $* = \deg'h_1 + \dots + \deg' h_{i-1} + 1$.
Let $b \in \mathcal A(L)$ and $1\le i \le j \le k$.
We define
\[
\vec a(b,i,j,1) := (a_0,\dots,a_{i},b,a_{j+1},\dots,a_k),
\qquad
\vec a(b,i,j,2) := (b,a_{i+1},\dots,a_j).
\]
Then by Theorem~\ref{thekuraexist}\,(3), we have
\begin{equation}\label{form347}
\partial {\mathcal M}(L;\vec a;E) =
\coprod_{b,i,j
\atop E_1 + E_2 = E, (*)}
{\mathcal M}(L;\vec a(b,i,j,2);E_2) {}_{{\rm ev}_0} \times_{{\rm ev}_{i+1}} {\mathcal M}(L;\vec a(b,i,j,1);E_1).
\end{equation}
Here the condition $(*)$ in the notation of direct sum is
\begin{enumerate}\itemsep=0pt
\item[(*.1)]
$E_1 = 0$ and $(a_0,\dots,a_{i},b,a_{j+1},\dots,a_k) = (b,b)$ does not hold.
\item[(*.2)]
$E_2 = 0$ and $(b,a_{i+1},\dots,a_j) = (b,b)$ does not hold.
\end{enumerate}
See Figure~\ref{Figure3page35}.
\begin{figure}[ht]
\centering
\includegraphics[scale=0.3]{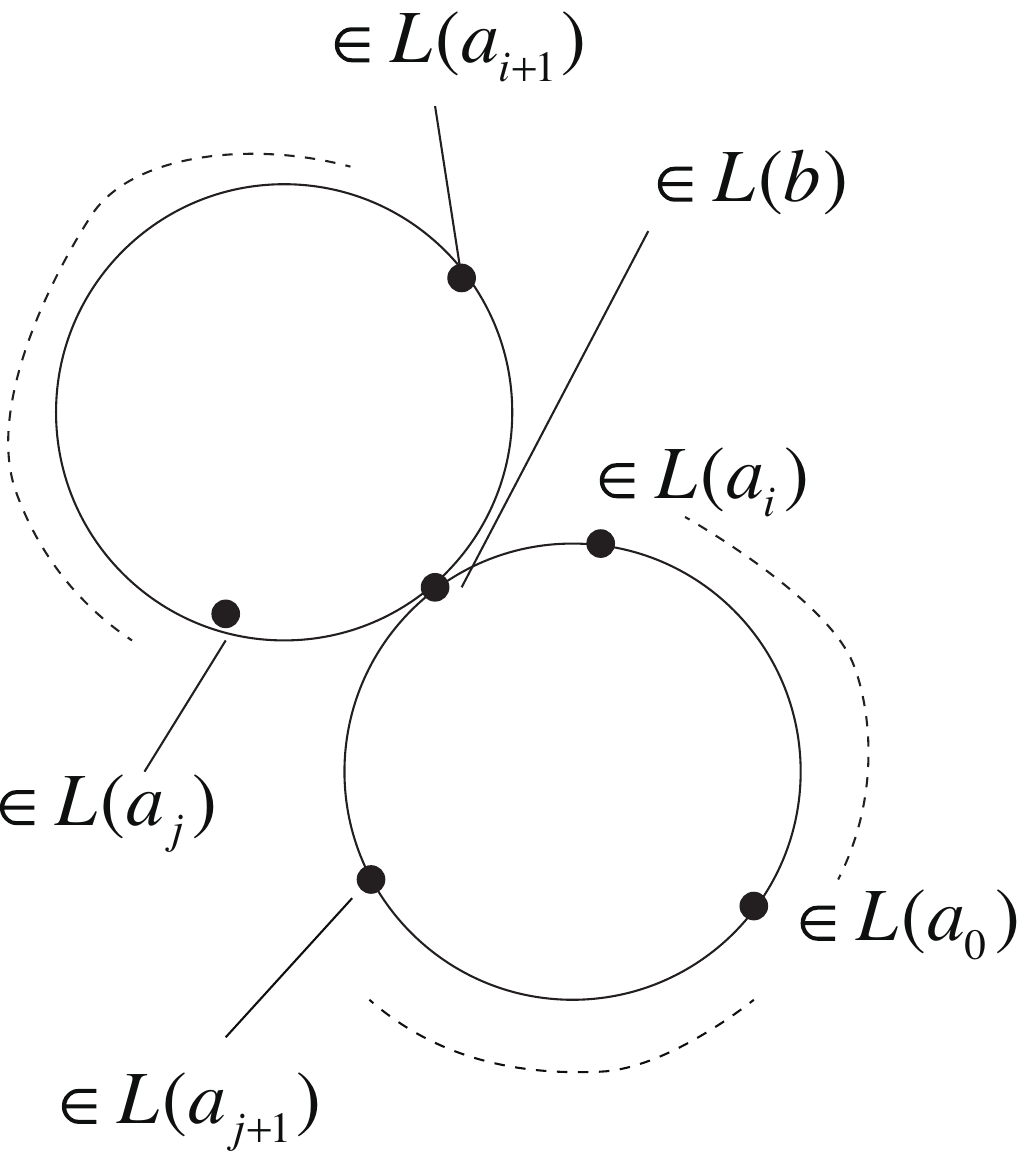}
\caption{Boundary of ${\mathcal M}(L;\vec a;E)$.}
\label{Figure3page35}
\end{figure}

By Proposition~\ref{prop330}\,(3), our CF-perturbations are
compatible with the isomorphism \eqref{form347}.
Therefore, by the composition formula (see \cite[Theorem 10.20]{foootech2} and \cite{fooonewbook}), the right-hand side of~\eqref{form346} is equal to the sum
\begin{gather*}
\sum_{b,i,j
\atop E_1 + E_2 = E, E_1,E_2>0}
{\rm ev}_0!\big({\rm ev}_1^*h_1 \times \dots \times {\rm ev}_{i+1}^*{\rm ev}_0!({\rm ev}_1^*h_{i+1}\\
\qquad \times
\cdots {\rm ev}_{j-i}^*h_{j});
\big({\mathcal M}(L;\vec a(b,i,j,2);E_2),{\widehat{\mathfrak S^{\varepsilon}}}\big)
 \times {\rm ev}_{i+2}^*h_{j+1}\\
 \qquad \times \cdots {\rm ev}_{j-i}^*h_{k} ;
\big({\mathcal M}(L;\vec a(b,i,j,1);E_1),{\widehat{\mathfrak S^{\varepsilon}}}\big)\big).
\end{gather*}
By definition, this sum is \eqref{form345} minus left-hand side of \eqref{form346}
up to sign.
This proves \eqref{form345} up to sign. See \cite[Chapter 8]{fooobook}
and \cite{fooonewbook} for the sign in the case of an embedded Lagrangian submanifold $L$.
In the case $L$ is immersed and has transversal self-intersection
see \cite{AJ}. The way to generalize them to the case
of an immersed Lagrangian submanifold which has clean self-intersection is explained
in Section~\ref{oriAinfMB} and in the paper \cite{ono2} by Kaoru Ono.

The unitality is a consequence of Proposition~\ref{prop330}\,(4).
\end{proof}

Note that one of the reasons why we stop our construction at $E=E_0$ is the
running out problem, which is explained in detail in \cite[Section 7.2.3]{fooobook2}.
(See also \cite[Section 14]{fooo091}, \cite{foootech22,fooonewbook}.)

The other reason why we need to fix $E_0$ and stop the
construction at $E = E_0$ appears in Remark~\ref{newnewrem}.
The well-definedness of the integration along the fiber
(as well as Stokes' theorem and the composition formula)
holds only for $\varepsilon< \varepsilon_0$, where $\varepsilon$ is the parameter of
our CF-perturbation, and $\varepsilon_0$ is
dependent on moduli spaces (spaces with Kuranishi structures) which we work with.\footnote{It might be possible to see carefully the moduli space itself
and obtain a certain estimate of this number. However, to include
such explicit estimate to the whole story of virtual fundamental
chain (such as CF-perturbation and de Rham theory) is rather cumbersome
and so to use only the fact that there exists such $\varepsilon_0$ for each
individual moduli space seems to be a better choice.}
As far as we consider the construction up to energy $E_0$ and $k < k_0$
($k$ is the number of
input),
we need to use only a finite number of moduli spaces so we can
take the same~$\varepsilon_0$ for all of them.

On the other hand, the CF-perturbation we obtain
this way is actually $E_0$ dependent.

Note that we require the compatibility of
CF-perturbations with forgetful maps of the marked
points corresponding to the diagonal component.
Therefore, we only need finiteness of the number of input which
does not correspond to the diagonal component.
The number of such inputs can be estimated by the energy
because of Lemma~\ref{Lema31}.

Even though we need to stop at $E=E_0$ and so can define
only a filtered $A_{\infty}$ structure modulo $T^{E_0}$,
we can still use it to define a filtered $A_{\infty}$ structure
as follows.
The method is the same as \cite[Section 7.2]{fooobook2},
\cite[Section 14]{fooo091} and \cite{foootech22,fooonewbook}.
(Our discussion here is slightly sketchy since it is the same as the papers
quoted above. More detail is given in \cite{AFOOO}.)

\begin{defn}[{\cite[Definition 8.5]{fooo091}}]\label{pisotopydef}
We consider $t \in [0,1]$ dependent families of operations\index[syindex]{m2kt@$\mathfrak m_k^t$} \index[syindex]{ckt@$\mathfrak c_k^t$}
\[
\mathfrak m_k^t \colon\ CF(L;\Lambda_0)^{\otimes k} \to CF(L;\Lambda_0),
\qquad
\mathfrak c_k^t\colon\ CF(L;\Lambda_0)^{\otimes k} \to CF(L;\Lambda_0)
\]
which are $G(L)$-gapped.
We say $\bigl(\bigl\{\mathfrak m^t_{k}\bigr\},\bigl\{\mathfrak c^t_{k}\bigr\}\bigr)$
is a {\it pseudo-isotopy} modulo $T^E$ \index{pseudo-isotopy modulo $T^E$} of $G$-gapped filtered~$A_{\infty}$ algebra
structures modulo $T^E$ on $CF(L)$
between $\bigl\{\mathfrak m^0_{k}\bigr\}$ and $\bigl\{\mathfrak m^1_{k}\bigr\}$ if the following holds:
\begin{enumerate}\itemsep=0pt
\item[(1)]
The operations $\mathfrak m^t_{k}$ and $\mathfrak c^t_{k}$ are
continuous in $C^{\infty}$ topology.
The map which sends $t$ to $\mathfrak m^t_{k}$ or $\mathfrak c^t_{k}$
is smooth. Here we use operator topology with respect to the $C^{\infty}$
topology for $\mathfrak m^t_{k}$ or $\mathfrak c^t_{k}$ to define
this smoothness.
\item[(2)] For each (but fixed) $t$, the set of operators $\bigl\{\mathfrak m^t_{k}\bigr\}$ defines a $G$-gapped filtered $A_{\infty}$ algebra
structures modulo $T^E$ on $CF(L;\Lambda_0)$.
\item[(3)] For each $h_i \in CF(L;\Lambda_0)$,
\begin{gather}
\frac{d}{dt} \mathfrak m_{k}^t(h_1,\ldots,h_k) + \sum_{k_1+k_2=k+1}\sum_{i=1}^{k-k_2+1}
(-1)^{*}\mathfrak c^t_{k_1}\bigl(h_1,\ldots, \mathfrak m_{k_2}^t(h_i,\ldots),\ldots,h_k\bigr)\nonumber \\
\qquad- \sum_{k_1+k_2=k+1}\sum_{i=1}^{k-k_2+1}
\mathfrak m^t_{k_1}\bigl(h_1,\ldots, \mathfrak c_{k_2}^t(h_i,\ldots),\ldots,h_k\bigr)\nonumber\\
\qquad\qquad\equiv 0 \mod T^E.\label{isotopymaineq}
\end{gather}
Here $* = \deg' h_1 + \dots + \deg'h_{i-1}$.
\item[(4)]
$\mathfrak c_{k}^t \equiv 0 \mod \Lambda_+$.
\end{enumerate}

\end{defn}
We put $G(L) = \{E_0, E_1,\dots, E_k, \dots\}$
with $0 = E_0 < E_1 < E_2 < \cdots$.
By Proposition~\ref{prop334}, we obtain
$\bigl\{\mathfrak m_{k}^{<E_i,\varepsilon}\bigr\}$ which defines a
$G(L)$-gapped filtered $A_{\infty}$ algebra modulo $E_i$ for each~${i =1,2,\dots}$.

We may regard \smash{$\bigl\{\mathfrak m_{k}^{<E_{i+1},\varepsilon}\bigr\}$} as a
$G(L)$ gapped filtered $A_{\infty}$ algebra modulo $E_i$
by forgetting the terms involving $T^{E_{i+1}}$.
We write it \smash{$\bigl\{\mathfrak m_{k}^{<E_{i+1},\varepsilon}\vert_{E_i}\bigr\}$}.
\begin{prop}\label{prop336}
There exits $\varepsilon_i$ such that if
$\varepsilon < \varepsilon_i$
then there exists $\bigl(\bigl\{\mathfrak m^{t,i,\varepsilon}_{k}\bigr\},\bigl\{\mathfrak c^{t,i,\varepsilon}_{k}\bigr\}\bigr)$
which is a pseudo-isotopy modulo $T^{E_i}$ of $G(L)$-gapped filtered $A_{\infty}$ algebra
structures modulo~$T^{E_i}$ on~$CF(L)$
between \smash{$\bigl\{\mathfrak m_{k}^{<E_i,\varepsilon}\bigr\}$} and \smash{$\bigl\{\mathfrak m_{k}^{<E_{i+1},\varepsilon}\vert_{E_i}\bigr\}$}.

\end{prop}
\begin{proof}
We remark that both $\bigl\{\mathfrak m_{k}^{<E_i,\varepsilon}\bigr\}$ and $\bigl\{\mathfrak m_{k}^{<E_{i+1},\varepsilon}\vert_{E_i}\bigr\}$
are defined as in Definition~\ref{def32}.
The only difference is we use different CF-perturbations to define them.
We use homotopy between those two different CF-perturbations.
We consider Kuranishi structures on ${\mathcal M}(L;\vec a;E) \times [0,1]$
which is a direct product with one on ${\mathcal M}(L;\vec a;E)$
given in Theorem~\ref{thekuraexist} and
the trivial Kuranishi structure on $[0,1]$.
\begin{lem}\label{lem337}
There exists a system of CF-perturbations \smash{${\widehat{\mathfrak S}}_{\rm para}$} of \index[syindex]{Spara@${\widehat{\mathfrak S}}_{\rm para}$}
outer collarings of thickenings of ${\mathcal M}(L;\vec a;E) \times [0,1]$ for
various $\vec a$ and $E<E_1$ with the following properties:
\begin{enumerate}\itemsep=0pt
\item[$(1)$]
Each of \smash{${\widehat{\mathfrak S}}_{\rm para} $} is transversal to zero.
\item[$(2)$]
${\rm ev}_0 \times \pi \colon {\mathcal M}(L;\vec a;E) \times [0,1] \to L \times [0,1]$ is strongly submersive with respect to this CF-perturbation.
\item[$(3)$]
They are compatible in a similar sense as Proposition {\rm\ref{prop330}\,(3)}.
\item[$(4)$]
Its restriction to ${\mathcal M}(L;\vec a;E) \times \{0\}$
coincides with the CF-perturbation we used to define~\smash{$\bigl\{\mathfrak m_{k}^{<E_i,\varepsilon}\bigr\}$}.
Its restriction to ${\mathcal M}(L;\vec a;E) \times \{1\}$
coincides with the CF-perturbation we used to define \smash{$\bigl\{\mathfrak m_{k}^{<E_{i+1},\varepsilon}\vert_{E_i}\bigr\}$}.
\item[$(5)$]
They are compatible with the forgetful maps of the marked points which
corresponds to the diagonal component other than $0$-th one. The
precise definition of compatibility is written in {\rm\cite[\emph{Definition} 5.1]{fooo091}}.
\end{enumerate}

\end{lem}
See \cite[Section 21]{foootech22}, \cite[Chapter 21]{fooonewbook} for the precise meaning of the compatibility in
item~(3).
The proof of Lemma~\ref{lem337} is mostly the same as Proposition~\ref{prop330}
and is omitted. See \cite[Section~21]{foootech22}, \cite[Chapter 21]{fooonewbook}.
\begin{rem}
Note that \smash{$\mathfrak m_k^{<E_i,\varepsilon}$} is {\it different} from \smash{$\mathfrak m_{k}^{<E_{i+1},\varepsilon}\vert_{E_i}$}
even for sufficiently small $\varepsilon$.
One of the reasons why it is difficult to take them to be the same
is explained in \cite[Section 7.2.3]{fooobook}.
Another reason appears in Remark~\ref{newnewrem}.
It is an opinion of the author that it is safer (if not inevitable) to use
``homotopy inductive limit'' than working out infinitely many
moduli spaces simultaneously and check that we can take the same $\varepsilon$ independent of
them.
\end{rem}
\begin{rem}
We mention thickenings and outer collarings in Lemma~\ref{lem337}.
This is the way to construct a CF-perturbation with appropriate properties
taken in
\cite{foootech22,fooonewbook}.
As far as applications concern, a CF-perturbation constructed
on an outer collaring of a thickening of the original Kuranishi structure
can be used in the same way as a CF-perturbation constructed on
the original Kuranishi structure.
(See \cite{foootech22,fooonewbook} for its reason.)

\end{rem}

Now we put
\begin{gather}
\mathfrak m_{E,k}^{i,t,\varepsilon}(h_1,\dots,h_k) + \mathfrak c_{E,k}^{i,t,\varepsilon}(h_1,\dots,h_k) \wedge  d t :=
({\rm ev}_0 \times \pi)! \bigl({\rm ev}_1^* h_1 \times \dots \times {\rm ev}_k^* h_k ; {\widehat{\mathfrak S^{\varepsilon}}}_{\rm para}\bigr).\!\!\label{form3833}
\end{gather}
Here we use the space ${\mathcal M}(L;\vec a;E) \times [0,1]$ with a Kuranishi structure
and its CF-perturbation to define the right-hand side.
The variable $t$ is the coordinate of $[0,1]$.
Note that \smash{$\mathfrak m_{E,k}^{i,t,\varepsilon}(h_1,\dots,h_k)$}
and \smash{$\mathfrak c_{E,k}^{i,t,\varepsilon}(h_1,\dots,h_k)$} are $t$-parametrized families
of elements of $CF(L;\R)$ which may be regarded as smooth forms on
$\tilde L \times_X \tilde L \times [0,1]$ that do not contain $ d t$.
(See \cite[Section 22.4]{fooonewbook} and \cite[Section 4.1]{ST} for the sign.)
We put
\[
\mathfrak m_{k}^{i,t,\varepsilon} := \sum_{E < E_i} T^E \mathfrak m_{E,k}^{i,t,\varepsilon},
\qquad
\mathfrak c_{k}^{i,t,\varepsilon} := \sum_{E < E_i} T^E \mathfrak c_{E,k}^{i,t,\varepsilon}.
\]
Using Lemma~\ref{lem337} in place of Proposition~\ref{prop330},
we can apply Stokes' formula and the composition formula in the same way as the
proof of Proposition~\ref{prop334} and obtain \eqref{isotopymaineq}.
\end{proof}

\begin{prop}\label{prop338}
There exits a positive number $\varepsilon_i$ such that if
$\varepsilon, \varepsilon' < \varepsilon_i$,
then there exists $\bigl(\bigl\{\mathfrak m^{\prime,t,i,\varepsilon}_{k}\},\{\mathfrak c^{\prime,t,i,\varepsilon}_{k}\bigr\}\bigr)$
which is a pseudo-isotopy modulo $T^{E_i}$ of $G(L)$-gapped filtered $A_{\infty}$ algebra
structures modulo $T^{E_i}$ on $CF(L)$
between \smash{$\bigl\{\mathfrak m_{k}^{<E_i,\varepsilon}\bigr\}$} and \smash{$\bigl\{\mathfrak m_{k}^{<E_i,\varepsilon'}\bigr\}$}.
\end{prop}
The proof is the same as the proof of Proposition~\ref{prop336} and so is omitted.
\par
We also use the next algebraic result.
\begin{lem}\label{lem339}
Let $E < E'$ and $\bigl\{\mathfrak m_{k}^{0}\bigr\}$ $\big($resp.\ $\bigl\{\mathfrak m_{k}^{1}\bigr\}\big)$ be
$G$-gapped filtered $A_{\infty}$ algebra modulo $T^E$ $\big($resp.\ \smash{$T^{E'}\big)$}
on $C(L;\Lambda_0)$.
We regard $\bigl\{\mathfrak m_{k}^{1}\bigr\}$ as a $G$-gapped filtered $A_{\infty}$ algebra modulo $T^E$
and denote it by $\bigl\{\mathfrak m_{k}^{1}\vert_{T^E}\bigr\}$
Let $\bigl\{\mathfrak c^{t}_{k}\bigr\}$ be a
pseudo-isotopy modulo $T^E$ of $G$-gapped filtered $A_{\infty}$ algebra
between $\bigl\{\mathfrak m_{k}^{0}\bigr\}$ and $\bigl\{\mathfrak m_{k}^{1}\vert_{T^E}\bigr\}$.
Then there exists $\bigl\{\mathfrak m_{k}^{0 +}\bigr\}$ and $\bigl\{\mathfrak c^{t,+}_{k}\bigr\}$
such that
\begin{enumerate}\itemsep=0pt
\item[$(1)$]
 $\bigl\{\mathfrak m_{k}^{0 +}\bigr\}$ is a $G$-gapped filtered $A_{\infty}$ algebra modulo $T^{E'}$.
 \item[$(2)$]
 If we regard $\bigl\{\mathfrak m_{k}^{0 +}\bigr\}$ as a
 $G$-gapped filtered $A_{\infty}$ algebra modulo $T^{E}$, then
 it coincides with~$\bigl\{\mathfrak m_{k}^{0}\bigr\}$.
 \item[$(3)$]
 $\bigl(\bigl\{\mathfrak m^{t,+}_{k}\bigr\},\bigl\{\mathfrak c^{t,+}_{k}\bigr\}\bigr)$ is a pseudo-isotopy modulo $T^{E'}$ of $G$-gapped filtered $A_{\infty}$ algebras
 between~$\bigl\{\mathfrak m_{k}^{0 +}\bigr\}$ and~$\bigl\{\mathfrak m_{k}^{1}\bigr\}$.
 \item[$(4)$]
 If we regard $\bigl\{\mathfrak c^{t,+}_{k}\bigr\}$ as a
pseudo-isotopy modulo $T^{E}$ of $G$-gapped filtered $A_{\infty}$ algebras, then
 it coincides with $\bigl\{\mathfrak c^{t}_{k}\bigr\}$.
\end{enumerate}

\end{lem}
\begin{proof}
We may assume $G(L) \cap [E,E'] = \{E'\}$.
We put $\bigl\{\mathfrak c^{t,+}_{k}\bigr\} := \bigl\{\mathfrak c^{t}_{k}\bigr\}$.
They we can solve differential equation \eqref{isotopymaineq}
to obtain a coefficient of \smash{$T^{E'}$} of $\{\mathfrak m_{k}^{0 +}\}$.
See \cite[Section 7--1]{fooobook2}, \cite[Section 14]{fooo091} and \cite{foootech22,fooonewbook} for the proof of a similar but more difficult
result.
\end{proof}

We take $\varepsilon_i$ which is smaller than the constants in
Propositions \ref{prop336} and \ref{prop338}.
Then we use Propositions \ref{prop336} and \ref{prop338}
and Lemma~\ref{lem339} inductively to find systems of operations
\smash{$\bigl\{\mathfrak m_{k}^{<E_i,j,\varepsilon_i}\bigr\}$}, $\bigl\{\mathfrak m_{k}^{<E_i,j,\varepsilon_{i+1}}\bigr\}$,
$\bigl(\bigl\{\mathfrak m^{t,i,j,\varepsilon}_{k}\bigr\},\bigl\{\mathfrak c^{t,i,j,\varepsilon}_{k}\bigr\}\bigr)$
$\bigl(\bigl\{\mathfrak m^{\prime,t,i,j,\varepsilon}_{k}\bigr\},\bigl\{\mathfrak c^{
\prime,t,i,j,\varepsilon}_{k}\bigr\}\bigr)$
for $j>i$ with the following properties:
\begin{enumerate}\itemsep=0pt
\item[(1)]
The operators \smash{$\bigl\{\mathfrak m_{k}^{<E_i,j,\varepsilon_i}\bigr\}$},
\smash{$\bigl\{\mathfrak m_{k}^{<E_i,j,\varepsilon_{i+1}}\bigr\}$} define structures of $G(L)$
gapped filtered $A_{\infty}$ algebras modulo $T^{E_j}$ on $CF(L;\Lambda_0)$.
\item[(2)]
The pair \smash{$\bigl(\bigl\{\mathfrak m^{t,i,j,\varepsilon}_{k}\bigr\},\bigl\{\mathfrak c^{t,i,j,\varepsilon}_{k}\bigr\}\bigr)$}
is a pseudo-isotopy modulo $T^{E_j}$ of $G$-gapped filtered $A_{\infty}$ algebras
 between \smash{$\bigl\{\mathfrak m_{k}^{<E_i,j,\varepsilon_{i+1}}\bigr\}$}
 and \smash{$\bigl\{\mathfrak m_{k}^{<E_{i+1},j,\varepsilon_{i+1}}\bigr\}$}.
\item[(3)]
The pair \smash{$\bigl(\bigl\{\mathfrak m^{\prime,t,i,j,\varepsilon}_{k}\bigr\},\bigl\{\mathfrak c^{\prime,t,i,j,\varepsilon}_{k}\bigr\}\bigr)$}
is a pseudo-isotopy modulo $T^{E_j}$ of $G$-gapped filtered $A_{\infty}$ algebras
 between \smash{$\bigl\{\mathfrak m_{k}^{<E_i,j,\varepsilon_{i}}\bigr\}$}
 and \smash{$\bigl\{\mathfrak m_{k}^{<E_{i},j,\varepsilon_{i+1}}\bigr\}$}.
 \item[(4)]
 If $j' < j$, then the system of structures
 \smash{$\bigl\{\mathfrak m_{k}^{<E_i,j,\varepsilon_i}\bigr\}$}, \smash{$\bigl\{\mathfrak m_{k}^{<E_i,j,\varepsilon_{i+1}}\bigr\}$},
\smash{$\bigl(\bigl\{\mathfrak m^{t,i,j,\varepsilon}_{k}\bigr\},\bigl\{\mathfrak c^{t,i,j,\varepsilon}_{k}\bigr\}\bigr)$},
\smash{$\bigl(\bigl\{\mathfrak m^{\prime,t,i,j,\varepsilon}_{k}\bigr\},\bigl\{\mathfrak c^{
\prime,t,i,j,\varepsilon}_{k}\bigr\}\bigr)$}
coincide with the system of structures
\smash{$\bigl\{\mathfrak m_{k}^{<E_i,j',\varepsilon_i}\bigr\}$}, \smash{$\bigl\{\mathfrak m_{k}^{<E_i,j',\varepsilon_{i+1}}\bigr\}$},
\smash{$\bigl(\bigl\{\mathfrak m^{t,i,j',\varepsilon}_{k}\bigr\},\bigl\{\mathfrak c^{t,i,j',\varepsilon}_{k}\bigr\}\bigr)$},
\smash{$\bigl(\bigl\{\mathfrak m^{\prime,t,i,j',\varepsilon}_{k}\bigr\},\bigl\{\mathfrak c^{
\prime,t,i,j',\varepsilon}_{k}\bigr\}\bigr)$}
as filtered $A_{\infty}$ structures modulo $T^{E_{j'}}$
or as pseudo-isotopies modulo $T^{E_{j'}}$.
\item[(5)]
If $j=i$, then
 \smash{$\bigl\{\mathfrak m_{k}^{<E_i,i,\varepsilon_i}\bigr\}$}, \smash{$\bigl\{\mathfrak m_{k}^{<E_i,i,\varepsilon_{i+1}}\bigr\}$},
\smash{$\bigl(\bigl\{\mathfrak m^{t,i,i,\varepsilon}_{k}\bigr\},\bigl\{\mathfrak c^{t,i,i,\varepsilon}_{k}\bigr\}\bigr)$},
\smash{$\bigl(\bigl\{\mathfrak m^{\prime,t,i,i,\varepsilon}_{k}\bigr\},\bigl\{\mathfrak c^{
\prime,t,i,j,\varepsilon}_{k}\bigr\}\bigr)$}
coincide with
 \smash{$\bigl\{\mathfrak m_{k}^{<E_i,\varepsilon_i}\bigr\}$}, \smash{$\bigl\{\mathfrak m_{k}^{<E_i,\varepsilon_{i+1}}\bigr\}$},
\smash{$\bigl(\bigl\{\mathfrak m^{t,i,\varepsilon}_{k}\bigr\},\bigl\{\mathfrak c^{t,i,\varepsilon}_{k}\bigr\}\bigr)$},
\smash{$\bigl(\bigl\{\mathfrak m^{\prime,t,i,\varepsilon}_{k}\bigr\},\bigl\{\mathfrak c^{
\prime,t,i,j,\varepsilon}_{k}\bigr\}\bigr)$}, respectively.
Note that \smash{$\bigl(\bigl\{\mathfrak m^{t,i,\varepsilon}_{k}\bigr\},\bigl\{\mathfrak c^{t,i,\varepsilon}_{k}\bigr\}\bigr)$}
is obtained by Proposition~\ref{prop336}
and \smash{$\bigl(\bigl\{\mathfrak m^{\prime,t,i,\varepsilon}_{k}\bigr\},\bigl\{\mathfrak c^{
\prime,t,i,j,\varepsilon}_{k}\bigr\}\bigr)$} is obtained by
Proposition~\ref{prop338}.
\end{enumerate}
Now we put
$
\mathfrak m_k = \lim_{j\to \infty}\mathfrak m_{k}^{<E_i,j,\varepsilon_i}$.
Note that the right-hand side converges in $T$ adic topology by item (4).
This is the required filtered $A_{\infty}$ structure.
The proof of Theorem~\ref{AJtheorem} is now complete.

\begin{rem}\label{rem340}
The filtered $A_{\infty}$ structure obtained by Theorem~\ref{AJtheorem}
is independent of the choices up to pseudo-isotopy.
We can prove it as follows.
We can prove that for each $E_i$ the structure $\bigl\{\mathfrak m_{k}^{<E_i,j,\varepsilon_i}\bigr\}$
is independent of the choices up to pseudo-isotopy modulo $T^{E_i}$
in the same way as Proposition~\ref{prop336}.
We can next show that this pseudo-isotopy modulo $T^{E_i}$
is independent of the choices up to
pseudo-isotopy of pseudo-isotopies modulo $T^{E_i}$.
We can use it in the same way as above to show
the required independence of the filtered $A_{\infty}$ structure
up to pseudo-isotopy.
See \cite[Section 21.3]{foootech22,fooonewbook}.
We will discuss this point more in Section~\ref{sec:independence2}.

\end{rem}

\begin{rem}\label{rem134141}
Let $h_1$, $h_2$ be differential forms on a connected component $L(a)$
of $\tilde L \times_X \tilde L$.
We can choose $\mathfrak m_{2,\beta_0}$ so that\index[syindex]{m22beta0@$\mathfrak m_{2,\beta_0}$}
$
\mathfrak m_{2,\beta_0}(h_1,h_2) = (-1)^{\deg h_1} h_1 \wedge h_2$.
(Here we use the sign convention of~\cite[Definition 21.29\,(5)]{fooonewbook}.)
In fact, the right-hand side is induced by the
moduli space of constant maps to $L(a)$ with three marked
points. This moduli space is transversal and we do {\it not}
perturb~it.

To define $\mathfrak m_{k,\beta_0}(h_1,h_2)$ for $k\ge 3$, we use
the moduli space of constant maps to $L(a)$ with more than three marked
points, which may be obstructed. Such a moduli space may be nonempty after
perturbation.
In the situation when $L(a)$ is zero-dimensional (which was the situation
of~\cite{AJ}), except the case of $\dim L = 1$
this moduli space has negative dimension and so we may assume~$\mathfrak m_{k,\beta_0}$ for $k\ge 3$ to be zero.

\end{rem}

\subsection[Filtered $A_\infty$ categories of immersed Lagrangian
Floer theory]{Filtered $\boldsymbol{A_{\infty}}$ categories of immersed Lagrangian
Floer theory}
\label{subsec:Ainfcatim}

\begin{situ}\label{situ320}
Let $(X,\omega)$ be a symplectic manifold
which is compact or convex at infinity.
We take $V$ a real oriented vector bundle on the 3-skeleton $X_{[3]}$.
We consider a finite set \index[syindex]{Lmathbb@$\mathbb L$}
\[
\mathbb L = \{(L_c,\sigma_c) \mid c \in \mathfrak O\}
\]
of pairs
$\mathcal L_c = (L_c,\sigma_c)$ of immersed Lagrangian submanifolds
$L_c$ and their $V$-relatively spin structure $\sigma_c$.
We assume the next two conditions:
\begin{enumerate}\itemsep=0pt
\item[(1)]
The self-intersection of each $L_c$ is clean.
\item[(2)]
The submanifold $L_c$ has clean intersection with $L_{c'}$
for any $c,c' \in \mathfrak O$.
\end{enumerate}
We call such $\mathbb L$ a {\it clean collection of $V$-relatively spin
immersed Lagrangian submanifolds}.
\index{clean collection of $V$-relatively spin
immersed Lagrangian submanifolds}
\end{situ}
The purpose of this subsection is to associate a filtered
$A_{\infty}$ category
$\mathfrak{Fuk}((X,\omega);V;\mathbb L)$
\index[syindex]{FukXomega@$\mathfrak{Fuk}((X,\omega);V;\mathbb L)$} to a clean collection $\mathbb L$ of $V$-relatively spin
immersed Lagrangian submanifolds.
The actual work to carry out for this purpose is in fact completed in the
last subsection and we only need to rephrase the outcome of the
last subsection.

Let $L_c = \bigl(\tilde L_c,i_{L_c}\bigr)$, where
$i_{L_c} \colon \tilde L_c \to X$ is a Lagrangian immersion.
We consider the disjoint union
\smash{$
\tilde L = \bigcup_{c \in \mathfrak O} \tilde L_c$},
and use \smash{$i_{L_a}$} to obtain a Lagrangian immersion
\smash{$i_L \colon \tilde L \to X$}.
We put~${L = \bigl(\tilde L,i_L\bigr)}$ and apply Theorem~\ref{AJtheorem}.
For $c,c' \in \mathfrak O$, we decompose \smash{$\tilde L_c \times_X \tilde L_{c'}$} to
connected components as
\begin{gather*}
\tilde L_c \times_X \tilde L_{c'}
=
\begin{cases}
\displaystyle
\bigcup_{a \in \mathcal A_{c,c'}} L_{c,c'}(a) &\text{if $c \ne c'$},\\
\displaystyle
\tilde L_c \cup \bigcup_{a \in \mathcal A_{c,c'}} \tilde L_{c,c'}(a) &\text{if $c = c'$}.
\end{cases}
\end{gather*}
We then put
\[
L(+) =
\tilde L \times_X \widetilde L
=
\bigcup_{c \in \mathfrak O} \widetilde L_c
\cup \bigcup_{c, c' \in \mathfrak O \atop a \in \mathcal A_{c,c'}}
L_{c,c'}(a).
\]
By Lemma--Definition~\ref{lemdef39}, we obtain a principal ${\rm O}(1)$ bundle ($\Z_2$ local system)
$\Theta^-$ on $\tilde L \times_X \widetilde L$.
We denote its restriction to
$L_{c,c'}(a)$ by $\Theta^-_{c,c';a}$.
\index[syindex]{Lccprimea@$L_{c,c'}(a)$} \index[syindex]{Thetaminus@$\Theta^-_{c,c';a}$}
We also remark that $\Theta^-$ is a trivial bundle
on the diagonal component.

According to Definition~\ref{defn313}, we have
\[
CF(L)
=
\bigoplus_{c \in \mathfrak O}
\Omega\bigl(\widetilde L_c\bigr)
\,\widehat\otimes\, \Lambda_0
\oplus \bigoplus_{c, c' \in \mathfrak O \atop a \in \mathcal A_{c,c'}}
\Omega(L_{c,c'}(a);\Theta^-_{c,c';a})
 \,\widehat\otimes\, \Lambda_0.
\]
\begin{defn}\label{330}
\quad
\begin{enumerate}\itemsep=0pt
\item[(1)]
The set of objects $\mathfrak{OB}(\mathfrak{Fuk}((X,\omega);V;\mathbb L))$
consists of the pairs $(L_c,\sigma_c)$ for $c \in \mathfrak O$.
\item[(2)]
If $c,c' \in \mathfrak O$ with $c \ne c'$, then the module of morphisms
from $(L_c,\sigma_c)$ to $(L_{c'},\sigma_{c'})$,
which we denote by \index[syindex]{FukXomega@$\mathfrak{Fuk}((X,\omega);V;\mathbb L)((L_c,\sigma_c),(L_{c'},\sigma_{c'}))$}
$\mathfrak{Fuk}((X,\omega);V;\mathbb L)((L_c,\sigma_c),(L_{c'},\sigma_{c'}))$,
is
\[
\bigoplus_{a \in \mathcal A_{c,c'}}\Omega(L_{c,c'}(a);\Theta^-_{c,c';a})
 \,\widehat\otimes\, \Lambda_0.
\]
\item[(3)]
In case $c=c'$,
the module of morphisms
from $(L_c,\sigma_c)$ to $(L_{c},\sigma_{c})$,
which we denote by
$\mathfrak{Fuk}((X,\omega);V;\mathbb L)((L_c,\sigma_c),(L_{c},\sigma_{c}))$
is
\[
\Omega\bigl(\widetilde L_c\bigr)
 \,\widehat\otimes\, \Lambda_0
\oplus
\bigoplus_{a \in \mathcal A_{c,c}}\Omega(L_{c,c}(a);\Theta^-_{c,c;a})
 \,\widehat\otimes\, \Lambda_0.
\]
\end{enumerate}
Hereafter, we write
$CF((L_c,\sigma_c),(L_{c'},\sigma_{c'}))$
in place of $\mathfrak{Fuk}((X,\omega);V;\mathbb L)((L_c,\sigma_c),(L_{c'},\sigma_{c'}))$.

\end{defn}
In Theorem~\ref{AJtheorem}, we obtained the structure operation of
our filtered $A_{\infty}$ algebra $\bigl(\Omega\bigl(\widetilde L_c\bigr),\allowbreak\{\mathfrak m_k\}\bigr)$
where
\begin{equation}\label{eq338}
\mathfrak m_k \colon\ CF(L)^{\otimes k} \to CF(L).
\end{equation}
\begin{defn}\label{defn331}
Let $c_0,\dots,c_k \in \mathfrak O$ and
$(L_{c_i},\sigma_{c_i})$, $i=0,\dots,k$, be corresponding objects.
We define
the structure operations
\begin{equation}\label{form354}
\mathfrak m_k \colon\ \bigotimes_{i=1}^k
CF((L_{c_{i-1}},\sigma_{c_{i-1}}),(L_{c_{i}},\sigma_{c_{i}}))
\to
CF((L_{c_{0}},\sigma_{c_{0}}),(L_{c_{k}},\sigma_{c_{k}}))
\end{equation}
as the corresponding component of \eqref{eq338}.
\end{defn}
\begin{rem}
\quad
\begin{enumerate}\itemsep=0pt
\item[(1)]
In Definition~\ref{defn22} (see \eqref{form23}), we required that the map $\mathfrak m_0$ of filtered $A_{\infty}$
category is
\[
\mathfrak m_0 \colon\ \Lambda_0 \to CF((L_c,\sigma_c),(L_{c'},\sigma_{c'}))
\]
and is nonzero only when $c=c'$.
We can check that our structure morphism
is zero in case $k=0$ and $c_0 \ne c_1$ as follows.

By definition, $\mathfrak m_0$ is defined by using the
moduli space of pseudo-holomorphic disks with one boundary
marked point. It consists of $(\Sigma,z_0,u,\gamma)$
where $\Sigma$ is bordered Riemann surface with one boundary
component and of genus $0$, $z_0 \in \partial \Sigma$,
$u \colon (\Sigma,\partial \Sigma) \to (X,L)$ and~${\gamma \colon \partial \Sigma \setminus \{z_0\}\to \tilde L}$.
We require $u = i_: \circ \gamma$ on
$\partial \Sigma \setminus \{z_0\}$.
(See Definition~\ref{def3737}\,(4).) Since~${\partial \Sigma \setminus \{z_0\}}$
is connected, the image of~$\gamma$ is contained in
one of the connected components of~$\tilde L$, say $\tilde L_c$.
In that case ${\rm ev}$ of this element
goes to $\tilde L_c \times_X \tilde L_c$.
So $\mathfrak m_0(1)$ is contained in the subspace
mentioned above.
\item[(2)]
It is also clear from the definition that the structure operation
$\mathfrak m_k$ of $L$ is decomposed as~\eqref{form354}.
Namely, the $CF((L_{c_{0}},\sigma_{c_{0}}),(L_{c_{k}},\sigma_{c_{k}}))$
component of $\mathfrak m_k(x_1,\dots,x_k)$
depends only on the component $(x_1,\dots,x_k)$ of $CF(L)^{\otimes k}$ such that
$x_1 \in CF((L_{c_{0}},\sigma_{c_{0}}),(L_{c_{1}},\sigma_{c_{1}}))$,
$x_i \in CF((L_{c_{i-1}},\sigma_{c_{i-1}}),(L_{c_{i}},\sigma_{c_{i}}))$
and
$x_k \in CF((L_{c_{k-1}},\sigma_{c_{k-1}}),(L_{c_{k}},\sigma_{c_{k}}))$
for some $c_1,\dots,c_{k-1}$.
\end{enumerate}
\end{rem}
\begin{thm}\label{prop333}
Definitions {\rm\ref{330}} and {\rm\ref{defn331}} define a curved filtered $A_{\infty}$
category.
$1 \in \Omega\bigl(\tilde L_c\bigr)$ becomes its unity.
\end{thm}
\begin{proof}
This is immediate from Theorem~\ref{AJtheorem} and Definition~\ref{defn22}.
\end{proof}

\begin{defn}\label{definition450}
Let $\mathscr C$ be a filtered $A_{\infty}$
category and $c$ its object. Then
$\mathscr C(c,c)$ together with restrictions of
structure operations define a structure of a filtered $A_{\infty}$
algebra.
Let $c$, $c'$ be two objects.
The restriction of structure operations
define a map
\[
\mathfrak n\colon\ B\mathscr C(c,c)[1] \otimes \mathscr C(c,c')
\otimes B\mathscr C(c',c')[1]
\to \mathscr C(c,c'),
\]
where $\mathscr C(c,c')$ is
the space of morphisms and is a completed free $\Lambda_0$ module.
We denote the restriction of $\mathfrak n$ to
$B_k\mathscr C(c,c)[1] \otimes \mathscr C(c,c')
\otimes B_{\ell}\mathscr C(c',c')[1]$
by $\mathfrak n_{k,\ell}$, $k,\ell =0,1,2,\dots$.
They define a~structure of filtered $A_{\infty}$
bi-module
\index{filtered $A_{\infty}$
bi-module} on $\mathscr C(c,c')$ over
$\mathscr C(c,c)$--$\mathscr C(c',c')$
in the sense of \cite[Definition~3.75]{fooobook}.
(See also Section~\ref{subsec:bi-functoralg}.)

\end{defn}
In the case of a filtered $A_{\infty}$ category,
 $\mathscr C(c,c)$
is nothing but the filtered~$A_{\infty}$ algebra associated to a~single (immersed) Lagrangian submanifold
$(L_c,\sigma_c)$. Moreover, $\mathscr C(c,c')$
is nothing but the
filtered~$A_{\infty}$ bi-module
associated to a pair of (immersed) Lagrangian submanifolds~${(L_c,\sigma_c)}$, $(L_{c'},\sigma_{c'})$.
Thus Theorem~\ref{prop333} reproduces
many of the constructions in \cite{fooobook}.
However, by this trick to include the immersed case
to reduce the construction of a filtered $A_{\infty}$ category
to one of a filtered $A_{\infty}$ algebra,
one aspect which we mention below is lost.
 Let $\varphi \colon X \to X$ be a~Hamiltonian diffeomorphism.
 As we is proved in \cite[Theorem 4.1.5]{fooobook},
 we have an homotopy equivalence
\begin{equation}\label{eq340}
CF((L_c,\sigma_c), (L_{c'},\sigma_{c'})) \otimes_{\Lambda_0}
\Lambda
 \cong
CF((L_c,\sigma_c), (\varphi(L_{c'}),\varphi(\sigma_{c'})))\otimes_{\Lambda_0}
\Lambda
\end{equation}
of filtered $A_{\infty}$ bi-module, in the case when $L_c$ and $L_{c'}$ are
embedded. See Section~\ref{sec:independence3} and Theorem~\ref{thm154}
for the immersed case. (The proof of Theorem~\ref{thm154} is actually
the same as the embedded case.)
Note that here we move $L_{c'}$ by $\varphi$ but do not move $L_c$.
It is difficult to see what is the corresponding construction
in the case of a single immersed Lagrangian submanifold other than the
obvious one. Namely, we move various connected components by
 different Hamiltonian diffeomorphisms.
However, it is rather hard to see in which sense
filtered $A_{\infty}$ algebra $CF(L)$ of an immersed Lagrangian submanifold
(with many components) is invariant.\footnote{Provably the unobstructed
(immersed) Lagrangian cobordism is the correct formulation to work with, see~\cite{BC}.}
One big reason for it is in \eqref{eq340} we have to use
$\Lambda$ coefficient rather than $\Lambda_0$ coefficient.
We will discuss related issue in Section~\ref{sec:independence3} more.
Note that \eqref{eq340} is the most important
property of Lagrangian Floer homology for applications.
In fact, the motivation of Floer to define Lagrangian Floer homology
is to study intersection of a pair of Lagrangian submanifolds
and the most important property of Floer homology
for that purpose is \eqref{eq340}.

The invariance of Floer homology (of $\Lambda$ coefficient) of a pair
under the Hamiltonian diffeomorphisms in the sense of
\cite[Theorem 4.1.5]{fooobook} will be discussed
in Section~\ref{sec:equivalencehammilton} in a slightly more sophisticated form.

\begin{rem}
In this subsection and in this paper,
we take and fix a finite set of Lagrangian submanifolds and
define our category by using those finitely many
Lagrangian submanifolds only.
It is more canonical to use all the Lagrangian submanifolds
and construct a single big filtered $A_{\infty}$ category.
We do not try to do so in this paper since for
the purpose of most of the applications
choosing an appropriate finite set of Lagrangian submanifolds
and using only those Lagrangian submanifolds
are good enough and since it is simpler to write the detail
in the case when we work on a finite set of Lagrangian submanifolds.
See Section~\ref{issueformal} for more discussion on this point.
\end{rem}

\subsection[Opposite $A_\infty$ category and $\omega \mapsto -\omega$]{Opposite $ \boldsymbol{A_{\infty}}$ category and $\boldsymbol{\omega \mapsto -\omega}$}
\label{subsec:Opposite}

In this subsection, we explain how the $A_{\infty}$ category
$\mathfrak{Fuk}((X,\omega);V;\mathbb L)$ behaves when we
replace the symplectic form $\omega$ by $-\omega$.
We use this relationship when we study Lagrangian correspondences.

Let $\mathbb L = \{(L_c,\sigma_c) \mid c \in \mathfrak O\}$ be
a clean collection of $V$-relatively spin
immersed Lagrangian submanifolds as in Situation \ref{situ320}.

\begin{lem}
We can regard $\mathbb L$ as a clean collection of $V$-relatively spin
immersed Lagrangian submanifolds of $(X,-\omega)$.
\end{lem}
The proof is obvious.
\begin{lem}\label{lem348}
There exists $V \oplus TX$-relatively spin structure $\sigma'_c$ of $L_c$ such that
$\mathbb L' = \{(L_c,\sigma'_c) \mid c \in \mathfrak O\}$
is a clean collection of $V \oplus TX$-relatively spin
immersed Lagrangian submanifolds of~$(X,\omega)$.
\end{lem}
\begin{proof}
We remark that $TX\vert_L = TL \oplus TL$.
Therefore, the lemma follows from the well known fact that
for any oriented real vector bundle $W$ there exists a
canonical spin structure on the bundle $W \oplus W$.
\end{proof}

From now on, we frequently identify the set $\mathbb L$ and $\mathbb L'$.
Now the main result of this subsection is the following.
\begin{thm}\label{opthere}
We may take the various choices made in the definitions
so that we have the next isomorphism of filtered $A_{\infty}$
categories
\[
\mathfrak{Fuk}((X,\omega);V;\mathbb L)
\cong
\mathfrak{Fuk}((X,-\omega);V\oplus TX;\mathbb L')^{\rm op}.
\]

\end{thm}
\begin{rem}
See also \cite[Remark 5.3.3]{WW5}.
\end{rem}
\begin{proof}
The proof is similar to the proof of \cite[Theorem 1.3]{fooo:inv}.
It is obvious that the set of objects and the modules of morphisms
are the same.
We can show that the local system $\Theta_-$ we use
to define module of morphisms does not change when we
replace background datum $V$ by~${V \oplus TX}$.
This is because $TX\vert_L$ is spin and its spin
structure is canonical.
We need to study a~certain sign issue which will be discussed during
the proof of Lemma~\ref{lem3.38} below.

Thus it remains to check that the structure operations coincide with each other.
By the argument of Section~\ref{subsec:Ainfcatim},
it suffices to consider the case when $\mathbb L$ consists of a single
$V$-relatively spin~immersed Lagrangian submanifold $(L,\sigma)$.

We consider the moduli space
\smash{$\raisebox{-1pt}{$\mathring{\widetilde{\mathcal M}}$}(L;\vec a;E)$} in
Definition~\ref{def3737}.
To specify the almost complex structure and the
symplectic form, we denote this moduli space as
\smash{$\raisebox{-1pt}{$\mathring{\widetilde{\mathcal M}}$}((X,\omega,J_X);L;\vec a;E)$}.
For~${\vec a = (a_0,\dots,a_k)}$, we put
$\vec a^{\rm op} = (a_k,\dots,a_{0})$
and define a map
\[
I \colon\ \mathring{\widetilde{\mathcal M}}((X,\omega,J_X);L;\vec a;E)
\to \mathring{\widetilde{\mathcal M}}((X,-\omega,-J_X);L;\vec a^{\rm op};E)
\]
as follows.
Let $(\Sigma;u;\vec z;\gamma) \in \mathring{\widetilde{\mathcal M}}((X,\omega,J_X);L;\vec a;E)$.
For simplicity, we assume $\Sigma = D^2$.
Then we put
$\vec z^{\, \prime} := (\overline z_0,\overline z_k,\dots,\overline z_1)
$, where $\vec z = (z_0,\dots,x_k)$, and
$u'(z) := u(\overline z)$, $ \gamma'(z) := \gamma(\overline z)$.
It is easy to see that
\[
I\big(D^2;u;\vec z;\gamma\big) := \big(D^2;u';\vec z^{\,\prime};\gamma'\big)
\in
\mathring{\widetilde{\mathcal M}}((X,-\omega,-J_X);L;\vec a^{\rm op};E).
\]
It is easy to see that we can extend $I$ to a homeomorphism
\begin{equation}\label{map342}
I \colon\ {{\mathcal M}}((X,\omega,J_X);L;\vec a;E)
\to {{\mathcal M}}((X,-\omega,-J_X);L;\vec a^{\rm op};E).
\end{equation}
\begin{lem}\label{lem3.38}
The map {\rm\eqref{map342}} is a underlying continuous map
of an isomorphism of Kuranishi structures.
The next diagram commutes:
\[
\begin{CD}
{{\mathcal M}}((X,\omega,J_X);L;\vec a;E) @ >{I}>>
 {{\mathcal M}}((X,-\omega,-J_X);L;\vec a^{\rm op};E) \\
@ V{{\rm ev}}VV @ VV{{\rm ev}}V\\
L(+)^{k+1} @ > {} >>
L(+)^{k+1},
\end{CD}
\]
where the map in the second horizontal arrow is
$
(x_0,x_1,\dots,x_k) \mapsto (x_0,x_k,\dots,x_1)$.
\end{lem}
\begin{proof}
The commutativity of the diagram is obvious from the definition.
The proof of the first half is the same as the proof of \cite[Proposition 4.5]{fooo:inv}.
\end{proof}

We need to study the orientation carefully to complete the
proof of Theorem~\ref{opthere}.
We decompose
\[
{{\mathcal M}}((X,\omega,J_X);L;\vec a;E)
=
\bigcup_{d}{{\mathcal M}}((X,\omega,J_X);L;\vec a;E;d),
\]
where ${{\mathcal M}}((X,\omega,J_X);L;\vec a;E;d)$
is the compactification of the
moduli space which consists of the elements $\big(D^2;\vec z;u,\gamma\big)$
with virtual dimension $d+\sum_{i=1}^k\dim L(a_i)$.
We define the moduli space~${{\mathcal M}((X,-\omega,-J_X);L;\vec a^{\rm op};E;d)}$
in the same way.

Let $h_0 \in \Omega^{d_0}(L(a_0);\Theta^-_{a_0}),
\dots, h_k \in \Omega^{d_k}(L(a_k);\Theta^-_{a_k})$.
We take a CF-perturbation to integrate
differential forms on the space
${{\mathcal M}}((X,\omega,J_X);L;\vec a;E;d)$
with Kuranishi structure (see \cite[Definition 10.22]{foootech2}).
By Lemma~\ref{lem3.38},
it induces a CF-perturbation
on ${{\mathcal M}}((X,\omega,J_X);L;\vec a;E;d)$.

We compare the integrations
\begin{equation}\label{form343}
\int_{{{\mathcal M}}((X,\omega,J_X);L;\vec a;E;d)}
{\rm ev}^*(h_0 \times h_1 \times \dots \times h_k)
\end{equation}
and
\begin{equation}\label{form344}
\int_{{{\mathcal M}}((X,-\omega,-J_X);L;\vec a^{\rm op};E;d)}
{\rm ev}^*(h_0 \times h_k \times \dots \times h_1).
\end{equation}
Here integrations are defined by using CF-perturbations (see \cite[Definition 10.22]{foootech2}).

We consider the case
\begin{equation}\label{form345-a}
d = \sum \deg h_i.
\end{equation}
\begin{lem}\label{lem339new}
We use the $V$-relative spin structure to define the
orientation of the moduli space which we use for
integration. Then
$\eqref{form343}
= (-1)^*\times \eqref{form344}$,
where
\[
* = 1 + \sum_{1\le i<j\le k}\deg' h_i\deg' h_j + \varepsilon.
\]
Here $\varepsilon =0$ if and only if
$
d - (k-2)
$
is divisible by $4$. Otherwise, ${\varepsilon} = 1$.

\end{lem}
\begin{proof}
The proof is mostly the same as \cite[Proposition 4.9]{fooo:inv}.

The sign
$\sum_{1\le i<j\le k}\deg' h_i\deg' h_j$
is induced by the fact that we exchange the order of $i$-th and $j$-th marked points.
Here $\deg'$ rather than $\deg$ appears since the moduli parameter
which moves those marked points are exchanged also.
The first term $1$ appears since the moduli parameter to move $0$-th marked point
is reversed.
See the proof of \cite[Proposition 4.9]{fooo:inv} for the detail of the argument
on those points.
We finally explain the reason why $\varepsilon$ appears.
During the proof of Proposition~\ref{prop329},
we use the fact that the index $\operatorname{Index} P$
appearing \eqref{aaa337} is isomorphic to a~complex vector space.
This is because $P$ is an operator on $S^2$ whose
symbol is the same as the Cauchy--Riemann operator.
It implies that its (real) determinant bundle
is trivial.

Since our isomorphism $I$ in Lemma~\ref{lem3.38}
sends a $J_X$-holomorphic map $u$ to a
$-J_X$-holomorphic map $\overline u$,
the map which is induced to $\operatorname{Index}  P$ by $I$
is not complex linear.
It is actually anti complex linear.
Therefore, it induces an orientation preserving map on
$\operatorname{Det} P$
if and only if the numerical index
(the complex dimension) of $P$ is even.
Note that $P$ is homotopic to the Cauchy--Riemann operator
on $S^2$ of a bundle with Chern number $m$,
where $m$ is the half of the Maslov index.
Therefore, this map is orientation preserving
if and only if the Maslov index~${
d - (k-2)}
$
is divisible by $4$.
This is the reason why $\varepsilon$ appears.
(This point is also similar to the proof of~\cite[Theorem~4.6]{fooo:inv}.)

The rest of the proof is entirely similar to the proof of \cite[Proposition 4.9]{fooo:inv}.
\end{proof}

We next show the following lemma.

\begin{lem}\label{lem3392}
Suppose \eqref{form345-a} holds.
The orientation which we obtain when using $V\oplus TX$-relative spin structure
is different from one we obtain when using $V$-relative spin structure
if and only if $(-1)^{\varepsilon} =-1$.
\end{lem}
\begin{proof}
We consider a map $u \colon \big(D^2,\partial D^2\big) \to (X,L)$.
It induces a trivialization of
$u\vert_{\partial D^2}^*(TX)$ since $D^2$ is contractible.
On the other hand since $TX\vert_L = TL \oplus TL$,
we have another trivialization of $u\vert_{\partial D^2}^*(TX)$.

It is easy to see that these two trivializations are homotopic
each other if and only if $(-1)^{\varepsilon} = 1$.

We can use this fact to prove the lemma as follows.
Let $\lambda \colon S^1 \to {\rm SO}(n)$ be a loop representing the
generator of $\pi_1({\rm SO}(n)) = \Z_2$.
Since the map $\pi_1({\rm SO}(n)) \to \pi_1({\rm U}(n)) = \Z$ induced by the inclusion is
trivial, we have a map $\lambda_+ \colon D^2 \to {\rm U}(n)$ which coincides with $\lambda$
on the boundary. We identify $\{0\} \times D^2 \times \C^n$ with
$\{1\} \times D^2 \times \C^n$ by using $\lambda_+$ and obtain a
rank $n$ complex vector bundle $E$ on $D^2 \times S^1$.
By construction, $E\vert_{\partial D^2 \times S^1}$ has a
real $n$-dimensional subbundle which is obtained by gluing
$\{0\} \times \partial D^2 \times \R^n$ with
$\{1\} \times \partial D^2 \times \R^n$ using $\lambda$.
We denote it by $F$. Note that the 2nd Stiefel--Whitney class of $F$ is nonzero
by the choice of $\lambda$.

Using the pair $(E,F)$, we obtain an $S^1$-parametrized family of
Cauchy--Riemann operators with boundary condition.
Namely, for $t \in S^1$ we consider
\[
\overline \partial \colon\ L^{2}_1\big(D^2;E\vert_{\{t\} \times D^2},F\vert_{\{t\} \times S^1}\big) \to L^{2}\big(D^2;E\vert_{\{t\} \times D^2}\big)
\]
on $E\vert_{\{t\} \times D^2}$
with boundary condition determined by $F$ and
obtain family of index bundle that is a real vector bundle over $S^1$.
Using the fact that the 2nd Stiefel--Whitney class of $F$ is nonzero,
the calculation in the proof of \cite[Proposition 8.1.7]{fooobook2} shows that
this bundle is unoriented.

This implies that the two orientations obtained by different trivializations
of $u\vert_{\partial D^2}^*(TX)$ are different.
This implies Lemma~\ref{lem3392}.
\end{proof}

Theorem~\ref{opthere} follows from Lemmas \ref{lem339new}
and \ref{lem3392} and the definition of opposite category
(see Definition~\ref{opcate} especially its item (3)).
\end{proof}

\section{Preliminary on Lagrangian correspondence}
\label{sec:lagcorr}

The review of the theory of filtered $A_{\infty}$ categories and
the construction of the filtered $A_{\infty}$ category associated to
a symplectic manifold is completed in the previous sections.
In this section, we start studying the relationship between
Lagrangian correspondences
and filtered $A_{\infty}$ functors, which is the main subject of
this paper.
This section is rather formal. We introduce certain notations
which we will use in later sections.

\begin{defnlem}
Let $L_1$ (resp.\ $L_{12}$) be an immersed Lagrangian submanifold
of $(X_1,\allowbreak\omega_1)$ (resp.\ $(X_{1} \times X_2, - \pi_1^*(\omega_1)
+ \pi_2^*(\omega_2))$).
\begin{enumerate}\itemsep=0pt
\item[(1)]
We say $L_1$ is transversal to $L_{12}$ if
the fiber product
$\tilde L_1 \times_{X_1} \tilde L_{12}$ is transversal.
\item[(2)]
Assume $L_1$ is transversal to $L_{12}$.
We put $\tilde L_2 = \tilde L_1 \times_{X_1} \tilde L_{12}$.
The composition $i_{L_2} \colon \tilde L_2 \to X_1 \times X_2 \to X_2$
is a Lagrangian immersion.
\item[(3)]
We call $L_2 = \bigl(\tilde L_2,i_{L_2}\bigr)$ the {\it geometric
transformation} \index{geometric transformation} of $L_1$ by $L_{12}$.
\end{enumerate}

\end{defnlem}
\begin{proof}
We prove item (2).
Let $x = (y,(p,q)) \in \tilde L_2$ and $V \in \operatorname{Ker}(d_xi_{L_2})$.
Then $V = (w,v)$ where $w \in T_y\tilde L_1, v \in T_p\tilde L_{12}$.
$(d_yi_{L_1})(w) = (d_pi_{L_{12}})(v)$ and $(d_pi_{L_{12}})(v) \in TX_1 \oplus 0$.
Since $\tilde L_1 \times_{X_1} \tilde L_{12}$ is transversal, there exists
$v' \in T_p\tilde L_{12}$ such that
$\omega_1((d_pi_{L_{12}})(v),(d_pi_{L_{12}})(v')) \ne 0$.
Since $(d_pi_{L_{12}})(v) \in TX_1 \oplus 0$ this implies
$
\omega((d_pi_{L_{12}})(v),(d_pi_{L_{12}})(v')) \ne 0$.
This contradicts the assumption that $L_{12}$ is an immersed
Lagrangian submanifold.
We have proved that $L_2$ is an immersed submanifold.

Let $(v_1,w_1), (v_2,w_2) \in T_x\tilde L_2$
where $v_i \in T_y\tilde L_1$, $w_i \in T_{(p,q)}\tilde L_{12}$.
Then we have $(d_xi_{L_1})(v_i) = (\pi_1(d_{(p,q)}i_{L_{12}}))(w_i)$.
Hence
$\omega_1(w_1,w_2) = 0$.
Since $\omega(w_1,w_2) = 0$, it follows that $\omega_2(w_1,w_2)
= 0$. We proved that $L_2$ is an immersed Lagrangian submanifold.
\end{proof}

It is not in general correct that
the geometric transformation of an embedded Lagrangian submanifold
by an embedded Lagrangian correspondence has clean self-intersection.

\begin{exm}
Let $X = (-1,1) \times S^1$, $L_1 = \{0\} \times S^1$.
We take a symplectic diffeomorphism~$\varphi$
which is a composition of $(s,t) \to (s,t+1/2)$ and a $C^1$ small
Hamiltonian diffeomorphism. \big(Here we identify $[0,1]/0\sim1 = S^1$.\big)
We can choose $\varphi$ such that $L_1 \cap \varphi(L_1)$ is
not clean.
Let~${L_{12} \subset - X_1 \times X_1}$ be the disjoint union
of the diagonal and the graph of $\varphi$.
The geometric transformation of $L_1$ by $L_{12}$ is not clean.
\end{exm}

\begin{defn}\label{defn43}
Let $L_1 \subset X_1$ and $L_{12} \subset X_{12}$ be
immersed Lagrangian submanifolds.
We say~$L_1$ has {\it clean transformation} by $L_{12}$
\index{clean transformation}
if
\begin{enumerate}\itemsep=0pt\samepage
\item[(1)]
The fiber product $\tilde L_1 \times_{X_1} \tilde L_{12}$
is transversal.
\item[(2)]
The geometric transformation
$L_2$ has clean self-intersection.
\end{enumerate}

\end{defn}

\begin{lem}\label{lem44}
Suppose $L_1$ has clean transformation by $L_{12}$ and let $L_2$ be its
geometric transformation.
\begin{enumerate}\itemsep=0pt
\item[$(1)$] If $L_1$ and $L_{12}$ are oriented so is $L_2$.
\item[$(2)$]
If $L_1$ $($resp.\ $L_{12})$ has $V_1$-relative spin structure
$($resp.\ $\pi_1^*(TX_1 \oplus V_1) \oplus \pi_2^*(V_2)$-relative spin structure$)$,
then $L_2$ has $V_2$-relative spin structure.
\end{enumerate}
\end{lem}
\begin{proof}
Let $x = (y,z) \in \tilde L_2$.
Then there exists a canonical isomorphism of vector spaces
\begin{equation}\label{form41}
T_{x} L_2 \oplus T_yX_1
\cong T_yL_1\oplus T_zL_{12}.
\end{equation}
This implies (1).

To prove (2), we first remark the following.
Suppose we have a transversal fiber product $\mathcal X \times_{\mathcal Y} \mathcal Z$.
Then we can choose smooth triangulations of $\mathcal X$, $\mathcal Y$, $\mathcal Z$,
 $\mathcal X \times_{\mathcal Y} \mathcal Z$ such that
\begin{enumerate}\itemsep=0pt
\item[(1)]
The maps $\mathcal X \to \mathcal Y$ and $\mathcal Z \to \mathcal Y$
send $2$ skeletons $\mathcal X_{[2]}$, $\mathcal Z_{[2]}$ to the $2$-skeleton.
\item[(2)]
The $2$-skeleton of $\mathcal X \times_{\mathcal Y} \mathcal Z$
is contained in
 $\mathcal X_{[2]} \times_{\mathcal Y_{[2]}} \mathcal Z_{[2]}$.
\end{enumerate}
To find such a triangulation, we first
take a triangulation of $\mathcal X \times_{\mathcal Y} \mathcal Z$.
We then can take enough many vertices of $\mathcal X$, $\mathcal Y$, $\mathcal Z$
such that
 $\mathcal X_{[0]} \times_{\mathcal Y_{[0]}} \mathcal Z_{[0]}$
 contains the $0$ skeleton of $\mathcal X \times_{\mathcal Y} \mathcal Z$.
 We can then take
$\mathcal X_{[1]}$, $\mathcal Y_{[1]}$, $\mathcal Z_{[1]}$
(subdividing $0$ skeleton if necessary), such that
$\mathcal X_{[1]} \times_{\mathcal Y_{[1]}} \mathcal Z_{[1]}$
contains the $1$ skeleton of the fiber product.
We then can find a required triangulation.

On the other hand, the trivialization on 2 skeleton $(L_{12})_{[2]}$ of
$\pi_1^*(TX_1 \oplus V_1) \oplus \pi_2^*(V_2) \oplus T L_{12}$
(that is nothing but the $\pi_1^*(TX_1 \oplus V_1) \oplus \pi_2^*(V_2)$-relative spin structure)
and the trivialization on~2 skeleton $(L_{1})_{[2]}$ of $TL_1 \oplus i_{L_1}^*V_1$
induce a trivialization of
\[
T_y X_1 \oplus (V_1)_y \oplus (V_2)_z \oplus T_zL_{12}
\oplus T_y L_1 \oplus (V_1)_y
\]
on the fiber product $(L_1)_{[2]} \times_{(X_1)_{[2]}} (L_1)_{[2]}$.
In view of \eqref{form41}, it induces a trivialization of
\begin{equation}\label{form4242}
T_{x} L_2 \oplus T_yX_1 \oplus T_yX_1
 \oplus (V_1)_y
\oplus (V_1)_y \oplus (V_2)_z
\end{equation}
on $(L_{2})_{[2]}$. (Note that we use our choice of triangulation and item (2) here.)

We remark that if $E$ is an oriented vector bundle then $E \oplus E$ is spin.
In fact,
\[
\sum_k w_k(E\oplus E) = \bigg(\sum_k w_k(E\oplus E)\bigg)^2.
\]
Hence $w_2(E\oplus E) = w_1(E)\cup w_1(E) = 0$ since $E$ is oriented.

Therefore, the existence of a trivialization of \eqref{form4242} on 2 skeleton $(L_2)_{[2]}$
implies the existence of a trivialization of
$
T L_2 \oplus V_2
$
on $(L_2)_{[2]}$. (Note that the trivialization and the spin structure are identical
notions on the 2 skeleton.)
Therefore, $L_2$ is $V_2$-relatively spin as required.
\end{proof}

\begin{rem}
The proof of Lemma~\ref{lem44} gives some particular relative spin
structure of $L_2$. However, in this paper we use the existence of
relative spin
structure of $L_2$ only. We make the choice of its relative
spin structure later during the proof of Theorem~\ref{trimain}
(see Lemma~\ref{exirespi}).
This relative spin structure seems to be related to one
obtained from the proof of Lemma~\ref{lem44}.
We however do not try to clarify the relationship
between those two relative spin structures in this paper.
\end{rem}

The next lemma will be used in later sections.
\begin{lem}\label{lem45}
Let $L_1$, $L_{12}$ be immersed submanifolds of
$X_1$ and $- X_1 \times X_2$ respectively.\footnote{See
Notation \ref{not3131} for $-X_1\times X_2$.}
We assume that $L_1 = \bigl(\tilde L_1,i_{L_1}\bigr)$ has clean transformation by $L_{12}
= \bigl(\tilde L_{12},i_{L_{12}}\bigr)$ and denote
by \smash{$L_2 = \bigl(\tilde L_2,i_{L_2}\bigr)$} the geometric transformation.
Then
\begin{equation}\label{fiberprod1}
\bigl(\tilde L_1 \times \tilde L_2\bigr) \times_{X_1 \times X_2}
\bigl(\tilde L_{12}\bigr)
\end{equation}
is diffeomorphic to
\begin{equation}\label{fiberprod2}
\tilde L_2 \times_{X_2} \tilde L_2.
\end{equation}

\end{lem}
\begin{proof}
\eqref{fiberprod1} is the left-hand side of
\begin{equation}\label{eq43}
\bigl(\tilde L_1 \times \bigl(\tilde L_1 \times_{X_1} \tilde L_{12}\bigr)\bigr) \times_{X_1 \times X_2}
\tilde L_{12}
=
\bigl(\tilde L_1 \times_{X_1} \tilde L_{12}\bigr) \times_{X_2}
\bigl(\tilde L_1 \times_{X_1} \tilde L_{12}\bigr).
\end{equation}
 On the other hand, \eqref{fiberprod2} is the right-hand side of \eqref{eq43}.
Note that the equality \eqref{eq43} is given by
\[
((x_1,(x_2,y_1)),y_2) \mapsto ((x_2,y_1),(x_1,y_2)).\tag*{\qed}
\]\renewcommand{\qed}{}
\end{proof}

We can generalize the definitions and lemmas of this section to the case
when we have
three symplectic manifolds, as follows.

\begin{defnlem}\label{defnlen47}
Let $(X_i,\omega_i)$ be a compact symplectic manifolds and $V_i$ its
background datum, for $i=1,2,3$.
Let $L_{12}$, $L_{23}$ be Lagrangian submanifolds of
$-X_{1} \times X_2$, $-X_{2} \times X_3$, respectively.
\begin{enumerate}\itemsep=0pt
\item[(1)]
If the fiber product $\tilde L_{13} = \tilde L_{12} \times_{X_2} \tilde L_{23}$
is transversal, then the map $\tilde L_{13} \to -X_1 \times X_3$
induced by $\tilde L_{13} \to -X_1 \times X_2\times -X_2 \times X_3 \to -X_1 \times X_3$ is a Lagrangian
immersion.
We assume that $L_{13}$ is self clean.
In such situation, we call $L_{13}$ the {\it geometric composition} of~$L_{12}$ and~$L_{23}$.\index{geometric composition}
\item[(2)]
If $L_{12}$ and $L_{23}$ are oriented, then so is the geometric composition $L_{13}$.
\item[(3)]
If $L_{12}$ has $\pi_1^*(TX_1 \oplus V_1) \oplus \pi_2^*(V_2)$-relative spin structure
and
$L_{23}$ has $\pi_1^*(TX_2 \oplus V_2) \oplus \pi_2^*(V_3)$-relative spin structure,
then the geometric composition
$L_{13}$ has $\pi_1^*(TX_1 \oplus V_1) \oplus \pi_2^*(V_3)$-relative spin structure.
\end{enumerate}
\end{defnlem}
\begin{proof}
The case when $X_1$ is a point is proved already.
The proof of the general case is the same and so is omitted.
\end{proof}

\section{The K\"unneth bi-functor in Lagrangian Floer theory}
\label{sec:Kunneth}

\subsection[Algebraic framework of $A_\infty$ bi-functors and tri-functors]{Algebraic framework of $\boldsymbol{A_{\infty}}$ bi-functors and tri-functors}
\label{subsec:bi-functoralg}

To define the notion of filtered $A_{\infty}$
bi-functor, we recall the following.
Let $(B_1,\Delta_1)$, $(B_1,\Delta_2)$ be
graded coalgebras.
We define graded coalgebra structure\index[syindex]{Delta@$\Delta$}
\[
\Delta \colon\ B_1 \otimes B_2
\to (B_1 \otimes B_2) \otimes (B_1 \otimes B_2)
\]
of
$B_1 \otimes B_2$ by the next formula
\begin{equation}\label{form5151}
\Delta(x \otimes y)
= \mathcal S(\Delta_1(x) \otimes \Delta_2(y)),
\end{equation}
where \index[syindex]{Scal@$\mathcal S$}
\[
\mathcal S((x_1 \otimes x_2) \otimes (y_1 \otimes y_2))
= (-1)^{\deg' y_1\deg' x_2}
((x_1 \otimes y_1) \otimes (x_2 \otimes y_2)).
\]
The case of completed tensor product of formal coalgebra is the same.
Note that in Definition~\ref{bi-functor} etc.\ we use the
shifted degree. So we used $\deg'$ in the above formula
instead of $\deg$.
\begin{defn}\label{bi-functor}
Let $\mathscr{C}_1$, $\mathscr{C}_2$, $\mathscr{C}_3$ be
non-unital curved filtered $A_{\infty}$ categories.
 A {\it filtered $A_{\infty}$
bi-functor}
\index{filtered $A_{\infty}$
bi-functor}
$
\mathscr{F}\colon \mathscr{C}_1 \times \mathscr{C}_2 \to \mathscr{C}_3
$
consists of $\mathscr{F}_{\rm ob}$ and $\mathscr{F}_{k_1,k_2}$,
$k_1,k_2 = 0,1,2,3,\dots$, of degree $0$ with the following properties (1), (2),(3), (4):
\begin{enumerate}\itemsep=0pt
\item[(1)]
$\mathscr{F}_{\rm ob} \colon \mathfrak{OB}(\mathscr{C}_1) \times
\mathfrak{OB}(\mathscr{C}_2) \to \mathfrak{OB}(\mathscr{C}_3)$
is a map between sets of objects.
\index[syindex]{Fob@$\mathscr{F}_{\rm ob}$}
\item[(2)]
For each $c_{1,1},c_{1,2} \in \mathfrak{OB}(\mathscr{C}_1)$
and $c_{2,1},c_{2,2} \in \mathfrak{OB}(\mathscr{C}_2)$,
the bi-functor $\mathscr{F}_{k_1,k_2}$ associates a $\Lambda_0$
linear map
\index[syindex]{Fk1k2@$\mathscr{F}_{k_1,k_2}$}
\begin{gather*}
\mathscr F_{k_1,k_2}(c_{1,1},c_{1,2};c_{2,1},c_{2,2})
\colon\ B_{k_1}\mathscr C_1[1](c_{1,1},c_{1,2})
\,\widehat\otimes\, B_{k_2}\mathscr C_2[1](c_{2,1},c_{2,2}) \\
\qquad\to
\mathscr C_3[1](\mathscr F_{\rm ob}(c_{1,1},c_{2,1}),\mathscr F_{\rm ob}(c_{1,2},c_{2,2})).
\end{gather*}
\item[(3)]
We require $\mathscr F_{k_1,k_2}(c_{1,1},c_{1,2};c_{2,1},c_{2,2})$
to preserve the filtration in a similar sense as
Definition~\ref{defn22}\,(2).
\end{enumerate}
Note that the symbol $\mathscr{C}_1 \times \mathscr{C}_2$
is used here. However, we do not define the product $\mathscr{C}_1 \times \mathscr{C}_2$
of two $A_{\infty}$ categories in this paper.
In other words, $\mathscr{C}_1 \times \mathscr{C}_2$ is simply a notation.

To describe the most important condition, we introduce certain notations.

Let $\Delta_i \colon B\mathscr C_i[1]((c_{i,1},c_{i,2})
\to B\mathscr C_i[1](c_{i,1},c_{i,2}) \,\widehat\otimes\,
B\mathscr C_i[1](c_{i,1},c_{i,2})$ be the formal coalgebra structure
for $i=1,2,3$. We define the formal coalgebra structure
$\Delta$ on the completed tensor product~${ B\mathscr C_1[1](c_{1,1},c_{1,2})
\,\widehat{\otimes}\,
 B\mathscr C_2[1](c_{2,1},c_{2,2})}$
 by \eqref{form5151}.

The system of maps
$\{\mathscr F_{k_1,k_2}\}$ induces uniquely a formal coalgebra
homomorphism
\begin{gather*}
\widehat{\mathscr F}(c_{1,1},c_{1,2};c_{2,1},c_{2,2})
\colon\ B\mathscr C_1[1](c_{1,1},c_{1,2})
\,\widehat{\otimes}\, B\mathscr C_2[1](c_{2,1},c_{2,2}) \\
\qquad\to
B\mathscr C_3[1](\mathscr F_{\rm ob}(c_{1,1},c_{2,1}),\mathscr F_{\rm ob}(
c_{1,2},c_{2,2})).
\end{gather*}
Note that the structure operations of $\mathscr C_i$
induce a coderivation \[
\hat d_i \colon\ B\mathscr C_i[1](c_{i,1},c_{i,2})
\to B\mathscr C_i[1](c_{i,1},c_{i,2}).
\]
\begin{enumerate}\itemsep=0pt
\item[(4)]
We regard \smash{$\widehat{\mathscr F}(c_{1,1},c_{1,2};c_{2,1},c_{2,2})$} as a chain map.
Namely, we require
\[
\hat d_3 \circ \widehat{\mathscr F}(c_{1,1},c_{1,2};c_{2,1},c_{2,2})
=
\widehat{\mathscr F}(c_{1,1},c_{1,2};c_{2,1},c_{2,2})
\circ \big(\hat d_1 \,\widehat{\otimes}\, {\rm id} + {\rm id}\, \widehat{\otimes}\, \hat d_2\big),
\]
\end{enumerate}
where $\widehat{\otimes}$ is as in Definition~\ref{def21}\,(6).

\end{defn}
\begin{defn}
Let $\mathscr{C}_1$, $\mathscr{C}_2$, $\mathscr{C}_3$ be
non-unital curved filtered $A_{\infty}$ categories and $
\mathscr{F}\colon \mathscr{C}_1 \times \mathscr{C}_2\allowbreak \to \mathscr{C}_3
$
a filtered $A_{\infty}$
bi-functor.
\begin{enumerate}\itemsep=0pt
\item[(1)] We say $\mathscr{F}$ is {\it strict} \index{strict} if $\mathscr{F}_{0,0} = 0$.
\item[(2)]
Suppose $\mathscr{C}_1$, $\mathscr{C}_2$, $\mathscr{C}_3$ are $G$-gapped.
We say $\mathscr{F}$ is {\it $G$-gapped} \index{$G$-gapped} if $\mathscr{F}_{k_1,k_2}$ are all
$G$-gapped.
\item[(3)]
Suppose $\mathscr{C}_1$, $\mathscr{C}_2$, $\mathscr{C}_3$ are unital.
We say $\mathscr{F}$ is {\it unital} \index{unital} if the following holds:
\begin{enumerate}\itemsep=0pt
\item
$\mathscr F_{1,0}({\bf e}_{c_1} \otimes 1) =
\mathscr F_{0,1}(1 \otimes {\bf e}_{c_2})
= {\bf e}_{\mathscr F_{\rm ob}(c_1,c_2)}$,
\item
$
\mathscr F_{k_1+\ell_1+1,k_2}\bigl(x^1_1,\dots,x^1_{k_1},{\bf e}_{c_1},y^1_1,\dots,y^1_{\ell_1}
;x^2_1,\dots,x^2_{k_2}\bigr)
= 0
$
for $k_1+k_2+\ell_1 >0$,
\item
$
\mathscr F_{k_1,k_2+\ell_2+1}\bigl(
x^1_1,\dots,x^1_{k_1};x^2_1,\dots,x^2_{k_2},{\bf e}_{c_2},y^1_1,\dots,y^1_{\ell_2}\bigr)
= 0
$
for $k_1+k_2+\ell_2 >0$.
\end{enumerate}
\end{enumerate}
\end{defn}

\begin{exm}\label{ex54}
Suppose $\mathscr{C}_1$, $\mathscr{C}_2$ have only one object.
We also assume that they are strict.
Then we may regard them as filtered $A_{\infty}$ algebras, which we denote by
$(C_1,\{\mathfrak m_k\})$, $(C_2,\{\mathfrak m_k\})$.
We call
a strict filtered $A_{\infty}$ bi-functor
\[
\mathscr F \colon\ \bigl(C^{\rm op}_1,\{\mathfrak m_k\}\bigr)\times (C_2,\{\mathfrak m_k\})
\to \mathcal{CH}
\]
a {\it filtered $A_{\infty}$ bi-module}
\index{filtered $A_{\infty}$ bi-module} over $(C_1,\{\mathfrak m_k\})$-$(C_2,\{\mathfrak m_k\})$.
We also say left $C_1$ and right $C_2$ filtered~$A_{\infty}$ bi-module.

The notion of filtered $A_{\infty}$ bi-module is introduced in
\cite[Definition 3.7.5]{fooobook}.
Below we will check that Definition~\ref{bi-functor} coincides with the definition in
\cite{fooobook} in this case.

Since $\mathscr{C}_1$, $\mathscr{C}_2$ have unique objects,
$\mathscr F_{\rm ob}$ determines a chain complex,
which we write $(D,d)$.
(Here $d$ is the boundary operator of this chain complex.)
$\mathscr F_{k_1,k_2}$ becomes a map
\[
\mathscr F_{k_1,k_2}\colon\ B_{k_1}C_1[1] \otimes B_{k_2}C_2[1]
\to \operatorname{Hom}(D,D)[1].
\]
of degree one. We will define
\[
\mathfrak n_{k_1,k_2}\colon\
B_{k_1}C_1[1] \otimes D \otimes B_{k_2}C_2[1] \to D.
\]
We first define
$
{\rm OP} \colon B_kC_1[1] \to B_kC_1[1]
$
\index[syindex]{O2P@${\rm OP}$}
by
$
{\rm OP}({\bf x}) = (-1)^{\varepsilon({\bf x})}{\bf x}^{\rm op}$,
where $\varepsilon({\bf x})$ and ${\bf x}^{\rm op}$ are
\eqref{form214}, \eqref{form215}, respectively.
We remark
\begin{equation}\label{form54new}
\mathfrak m_k^{\rm op}({\bf x}) = -\mathfrak m_k({\rm OP}({\bf x})).
\end{equation}
We now put
\[
\mathfrak n_{k_1,k_2}(\text{\bf x};y;\text{\bf z})
:= (-1)^{\deg' y \deg'\text{\bf z}}(\mathscr F_{k_1,k_2}({\rm OP}({\bf x});\text{\bf z}))(y)
\]
for $(k_1,k_2) \ne (0,0)$
(note that $\deg'(x_1\otimes\dots\otimes x_{k}) = k+ \sum \deg x_i$) and
\begin{equation}\label{form56}
\mathfrak n_{0,0}(y) := (-1)^{\deg y} dy.
\end{equation}
We call $\mathfrak n_{k_1,k_2}$ the structure operations of filtered $A_{\infty}$
bi-module
(compare Definition~\ref{definition450}).

We will prove that
Definition~\ref{bi-functor}\,(4) becomes the following equality:
\begin{align}
0 ={} &\sum (-1)^{\deg'\text{\bf x}^{(2;1)}_{c_x}}\mathfrak n\bigl(\text{\bf x}^{(2;1)}_{c_x};
\mathfrak n(\text{\bf x}^{(2;2)}_{c_x};y;\text{\bf z}^{(2;1)}_{c_z});
\text{\bf z}^{(2;2)}_{c_z}\bigr)\nonumber
\\
&+
\mathfrak n\bigl(\hat d_1(\text{\bf x});y;\text{\bf z}\bigr)
+
(-1)^{\deg'\text{\bf x} + \deg y}\mathfrak n\bigl(\text{\bf x};y;\hat d_2(\text{\bf z})\bigr).\label{form51}
\end{align}
The notation is as follows.
The symbol $\hat d_1$ \big(resp.\ $\hat d_2$\big) is the
coderivation induced by the $A_{\infty}$ operations
on $C_1$ (resp.\ $C_2$).
We put
\[
\Delta_1(\text{\bf x}) =
\sum_{c_x} \text{\bf x}^{(2;1)}_{c_x} \otimes \text{\bf x}^{(2;2)}_{c_x},
\qquad
\Delta_2(\text{\bf z}) =
\sum_{c_z} \text{\bf z}^{(2;1)}_{c_z} \otimes \text{\bf z}^{(2;2)}_{c_z}.
\]
The formula \eqref{form51} is the defining relation of a filtered $A_{\infty}$ bi-module
in \cite[Definition 3.7.5]{fooobook}.
We also call it the $A_{\infty}$ relation.
\begin{rem}
In \cite[Definition 3.7.5]{fooobook}, the sign in the third term of
right-hand side is \[
(-1)^{\deg'{\bf z} + \deg' y}.
\]\ So it is slightly
different. In \cite{fooobook}, the bi-module is written $D(1)$.
Here we use the notation $D$ for a bi-module. So the definitions of this paper and of \cite{fooobook}
are consistent. We discuss this point more in Remarks \ref{rem55111} and \ref{lem5757}.
\end{rem}

Let us prove that Definition~\ref{bi-functor}\,(4) implies
\eqref{form51}. (The main part of the proof is the check of the
sign.)
Definition~\ref{bi-functor}\,(4) becomes the following identity
in $\operatorname{Hom}(D,D)[1]$:
\begin{gather}
\sum(-1)^{\deg'\text{\bf x}^{(2;2)}_{c_x}\deg'\text{\bf z}^{(2;1)}_{c_z}}
\mathfrak m_2\bigl(\mathscr F\bigl(\text{\bf x}^{(2;1)}_{c_x},\text{\bf z}^{(2;1)}_{c_z}\bigr),
\mathscr F\bigl(\text{\bf x}^{(2;2)}_{c_x},\text{\bf z}^{(2;2)}_{c_z}\bigr)\bigr)+
\mathfrak m_1(\mathscr F(\text{\bf x},\text{\bf z}))\nonumber
\\
\qquad
= \mathscr F\bigl(\hat d_1(\text{\bf x}),\text{\bf z}\bigr)
+ (-1)^{\deg'\text{\bf x}}\mathscr F\bigl(\text{\bf x},\hat d_2(\text{\bf z})\bigr).\label{form57}
\end{gather}
Here $\mathfrak m_1$, $\mathfrak m_2$ in the left-hand side are the
structure operations of $\mathcal{CH}[1]$.
They are related to the composition and the differential
by \eqref{form213}.

We plug in $y \in D$ in the first term of the left-hand side
and obtain
\begin{gather}
\sum(-1)^{*_1}
\mathscr F\bigl(\text{\bf x}^{(2;2)}_{c_x},\text{\bf z}^{(2;2)}_{c_z}\bigr)
\bigl(\mathscr F\bigl(\text{\bf x}^{(2;1)}_{c_x},\text{\bf z}^{(2;1)}_{c_z}\bigr)(y)\bigr)\nonumber
\\
\qquad=\sum (-1)^{*_2}
\mathfrak n\bigl({\rm OP}\bigl(\text{\bf x}^{(2;2)}_{c_x}\bigr)
;\mathfrak n
\bigl({\rm OP}(\text{\bf x}^{(2;1)}_{c_x});y;\text{\bf z}^{(2;1)}_{c_z}\bigr);
\text{\bf z}^{(2;2)}_{c_z}\bigr).\label{form58}
\end{gather}
Here
\begin{align*}
*_1 &= \bigl(\deg'\text{\bf x}^{(2;2)}_{c_x} +
\deg'\text{\bf z}^{(2;2)}_{c_z}\bigr)
\bigl(\deg'\text{\bf x}^{(2;1)}_{c_x} +
\deg'\text{\bf z}^{(2;1)}_{c_z} + 1\bigr)
+
\deg'\text{\bf x}^{(2;2)}_{c_x}\deg'\text{\bf z}^{(2;1)}_{c_z}, \\
*_2
&= *_1 + \deg'\text{\bf z}^{(2;1)}_{c_z}\deg' y
+ \deg'\text{\bf z}^{(2;2)}_{c_z}\bigl(\deg'\text{\bf x}^{(2;1)}_{c_x}
+ \deg'y + \deg'\text{\bf z}^{(2;1)}_{c_z} + 1\bigr) \\
&= \deg'\text{\bf x}^{(2;1)}_{c_x}\deg'\text{\bf x}^{(2;2)}_{c_x}
+\deg'\text{\bf x}^{(2;2)}_{c_x} + \deg'y\bigl(\deg'\text{\bf z}^{(2;1)}_{c_z}
+ \deg'\text{\bf z}^{(2;2)}_{c_z}\bigr).
\end{align*}
Note that in the sum \eqref{form58} the case
\begin{gather*}
\bigl(\text{\bf x}^{(2;1)}_{c_x},\text{\bf z}^{(2;1)}_{c_z}\bigr)
= 1 \otimes 1 \in B_{0}C_1[1] \otimes B_{0}C_2[1]\qquad
\text{or}\\
\bigl(\text{\bf x}^{(2;2)}_{c_x},\text{\bf z}^{(2;2)}_{c_z}\bigr)
= 1 \otimes 1 \in B_{0}C_1[1] \otimes B_{0}C_2[1]
\end{gather*}
are included only for the first term of \eqref{form56}.
The contribution of the second term of \eqref{form56} in
those cases actually coincide with the
second term of the left-hand side of \eqref{form57}.
Therefore, the left-hand side of \eqref{form57}
coincides with \eqref{form58} including those cases.

We replace $\bf x$ by ${\rm OP}(\bf x)$ in \eqref{form57}.
We remark that
\[
{\rm OP}(\Delta_1(\text{\bf x}))
=
\sum
(-1)^{\deg'\text{\bf x}^{(2;1)}_{c_x}\deg'\text{\bf x}^{(2;2)}_{c_x}}{\rm OP}\bigl(\text{\bf x}^{(2;2)}_{c_x}\bigr) \otimes
{\rm OP}\bigl(\text{\bf x}^{(2;1)}_{c_x}\bigr).
\]
Therefore, \eqref{form57} and \eqref{form58}
becomes
\begin{gather}
\sum (-1)^{*_3}
\mathfrak n\bigl(\text{\bf x}^{(2;1)}_{c_x}
;\mathfrak n
\bigl(\text{\bf x}^{(2;2)}_{c_x};y;\text{\bf z}^{(2;1)}_{c_z}\bigr);
\text{\bf z}^{(2;2)}_{c_z}\bigr)\nonumber
\\
\qquad=\bigl(\mathscr F\bigl(\hat d_1({\rm OP}({\bf x})),\text{\bf z}\bigr)\bigr)(y)
+ (-1)^{\deg'\text{\bf x}}\bigl(\mathscr F\bigl({\rm OP}({\bf x}),\hat d_2(\text{\bf z})\bigr)\bigr)(y),\label{form510}
\end{gather}
where
\[
*_3 = \deg'\text{\bf x}^{(2;1)}_{c_x} + \deg'y\bigl(\deg'\text{\bf z}^{(2;1)}_{c_z}
+ \deg'\text{\bf z}^{(2;2)}_{c_z}\bigr)
=
\deg'\text{\bf x}^{(2;1)}_{c_x}
+ \deg'y\deg'\text{\bf z}.
\]
Using \eqref{form54new}, we can calculate the right-hand side of
\eqref{form510} to obtain
\[
- (-1)^{*_4}\mathfrak n\bigl(\hat d_1(\text{\bf x});y;\text{\bf z}\bigr)
+ (-1)^{*_5}\mathfrak n\bigl(\text{\bf x};y;\hat d_2(\text{\bf z})\bigr),
\]
where $*_4 = \deg'\text{\bf z}\deg 'y$
and $*_5 =\deg'\text{\bf x}
+ (\deg'\text{\bf z}+1)\deg 'y
= \deg'\text{\bf x} + \deg'\text{\bf z}\deg 'y
+ \deg y + 1$.
Therefore, \eqref{form510} becomes \eqref{form51},
as required.

\end{exm}
Note that the order of \smash{$\mathscr F\bigl(\text{\bf x}^{(2;1)}_{c_x},\text{\bf z}^{(2;1)}_{c_z}\bigr)$}
and
\smash{$\mathscr F\bigl(\text{\bf x}^{(2;2)}_{c_x},\text{\bf z}^{(2;2)}_{c_z}\bigr)$}
appearing in \eqref{form57}
is reversed in~\eqref{form58}.
This is because it is defined so in \eqref{form213}.
The sign
\[
(-1)^{\deg x(\deg y+1)}= (-1)^{\deg' x\deg' y+\deg' y}
\]
there is actually the Koszul sign, that is, associated to the exchange
of the symbols $\mathfrak m_2$, $x$, $y \mapsto y$, $\circ$, $x$.
This is an intuitive reason why rather complicate sign
calculation in Example~\ref{ex54} works.
\begin{rem}\label{rem55111}
If $D(1)$ has a structure of $A_{\infty}$ bi-module $\mathfrak n$, then
its degree shift $D$ has a $A_{\infty}$ bi-module $\mathfrak n'$
defined by
\begin{equation}\label{degzurashi}
\mathfrak n'({\bf x},ys,{\bf z}) = (-1)^{\deg'\bf z} \mathfrak n'({\bf x},y,{\bf z})s.
\end{equation}
Here $ys$ is an element $y \in D(1)$ regarded as an
element of $D$.
\end{rem}
\begin{exm}\label{ex53}
To any filtered $A_{\infty}$ category $\mathscr C$,
we can associate a filtered $A_{\infty}$ bi-module
$\mathscr C(1)$ as follows.
We put $
\mathscr C(1)^d(c_1,c_2) = \mathscr C^{d+1}(c_1,c_2)
$.
We define structure operation $\mathfrak n$ by
$
\mathfrak n({\bf x},y,{\bf z}) = \mathfrak m({\bf x},y,{\bf z})$.
It is easy to see that this satisfies the $A_{\infty}$
relation.

In view of Remark~\ref{rem55111}, it induces a
structure of filtered $A_{\infty}$ bi-module
on $\mathscr C$ (without degree shift) by
\[
\mathfrak n'({\bf x},y,{\bf z}) = (-1)^{\deg' {\bf z}}\mathfrak m({\bf x},y,{\bf z}).
\]
In case $\mathscr C$ is strict, the operator $\mathfrak n'$ induces
a strict filtered $A_{\infty}$-bi-functor
$
\mathscr F \colon \mathscr C^{\rm op} \times \mathscr C
\to \mathcal{CH}
$
as follows.
We put $
\mathscr F_{\rm ob}(c_1,c_2) = \mathscr C(c_1,c_2)
$.
We define a map
\[
\mathscr F'_{k_1,k_2} \colon\
B_{k_1}\mathscr C[1](c_0,c_1) \,\widehat{\otimes}\, B_{k_2}\mathscr C[1](c_2,c_3)
\to
\operatorname{Hom}(\mathscr C(c_1,c_2),\mathscr C(c_0,c_3))[1]
\]
by
\[
(\mathscr F_{k_1,k_2}(\text{\bf x};\text{\bf z}))(y)
=
(-1)^{*} \mathfrak m_{k_1+k_2+1}(\text{\bf x},y,\text{\bf z}),
\]
where
$* = \deg y\deg'\text{\bf z}$.
(Here $\deg y$ appears instead of $\deg' y$ because of the sign in
\eqref{degzurashi}.)
We then compose it with ${\rm OP}$,
and obtain the required map
\[
\mathscr F_{k_1,k_2} \colon\
B_{k_1}\mathscr C^{\rm op}[1](c_1,c_0)
 \,\widehat{\otimes}\, B_{k_2}\mathscr C[1](c_2,c_3)
\to
\operatorname{Hom}(\mathscr C(c_1,c_2),\mathscr C(c_0,c_3)).
\]

This construction is an analogue of the fact that an arbitrary algebra is a bi-module over itself.

\end{exm}
\begin{rem}\label{lem5757}
A reason why we shifted the degree in \cite{fooobook} is
Example~\ref{ex53}. Namely, we can put~${\mathfrak m = \mathfrak n}$ if
we shift the degree.
The reason why we do not shift the degree of bi-module will
be clear in Section~\ref{2-category formulation}.
There we will regard a left-$\mathscr C_1$ and right-$\mathscr C_2$
bi-module as a `morphism' from~$\mathscr C_1$ to $\mathscr C_2$.
In that case, the bi-module in Example~\ref{ex53} plays the role of
the identity morphism. However, if we shift the degree then it will not
behave as the identity morphism.
Until Section~\ref{2-category formulation}, we will use
the convention of \cite{fooobook}, that is, we shift the degree
of bi-module.
In the way explained in \eqref{degzurashi}, we can go
from one to the other.
\end{rem}

We next generalize the story of \cite[Section 5.2.2.1]{fooobook}
to our category case.
\begin{lem}\label{lem55}
Let $\mathscr{C}_1$, $\mathscr{C}_2$, $\mathscr{C}_3$ be
non-unital curved filtered $A_{\infty}$ categories
and $\mathscr{C}^s_1$, $\mathscr{C}^s_2$, $\mathscr{C}^s_3$
their associated strict categories.
Then
any filtered $A_{\infty}$
bi-functor
$
\mathscr{F} \colon \mathscr{C}_1 \times \mathscr{C}_2 \to \mathscr{C}_3
$
induces
a strict filtered $A_{\infty}$
bi-functor $
\mathscr{F}^s \colon \mathscr{C}^s_1 \times \mathscr{C}^s_2 \to \mathscr{C}^s_3$.
If $\mathscr{F}$ is unital or $G$-gapped, then
so is $\mathscr{F}^s$.\index[syindex]{Fs@$\mathscr{F}^s$}

\end{lem}
\begin{proof}
Let $c_i \in \mathfrak{OB}(\mathscr{C}_i)$, $(c_i,b_i) \in \mathfrak{OB}(\mathscr{C}'_i)$
for $i=1,2$.
We put
\[
b_3 = \sum_{k_1=0}^{\infty}\sum_{k_2=0}^{\infty}
\mathscr{F}_{k_1,k_2}\bigl(b_1^{k_1},b_2^{k_2}\bigr).
\]
We put $e^b = \sum_{k=0}^{\infty} b^k$ then
\smash{$
e^{b_3} = \widehat{\mathscr{F}}\big(e^{b_1}, e^{b_2}\big)$}.
Since $b_i$ are bounding cochains for $i=1,2$, we have
\smash{$\hat d_1\bigl(e^{b_1}\bigr) = \hat d_2\bigl(e^{b_2}\bigr) = 0$}.
(See \cite[Lemma 3.6.36]{fooobook}.)
Therefore, Definition~\ref{bi-functor}\,(4) implies~\smash{$\hat d_1\bigl(e^{b_3}\bigr) = 0$}. In other words, $b_3$ is a bounding
cochain. We define
\[
\mathscr{F}^s_{\rm ob}((c_1,b_1),(c_2,b_2)) = (c_3,b_3).\]
Let
${\bf x}_i \in B_{k_i}\mathscr C_i[1]\bigl(c^1_{i},c^2_{i}\bigr)
\cong B_{k_i}\mathscr C'_i[1]\bigl(\bigl(c^1_{i},b^1_{i}\bigr),\bigl(c^2_{i},b^2_{i}\bigr)\bigr)$,
$i=1,2$. We will define~$\mathscr{F}^s_{k_1,k_2}({\bf x}_1,{\bf x}_2)$.
For this purpose, we define
$
\mathfrak t^{b_i}\colon B\mathscr C_i[1]\bigl(c^1_{i},c^2_{i}\bigr)
\to B\mathscr C_i[1]\bigl(c^1_{i},c^2_{i}\bigr)
$
for $i=1,2,3$ as follows. Let
$
{\bf x}_i = x_{i,1} \otimes \dots \otimes x_{i,k_i}
$,
where $x_{i,j} \in \mathscr C_i[1](c_{i,j-1},c_{i,j})$,
$c_{i,j} \in \mathfrak{OB}(\mathscr C_i)$, with
$c_{i,0} = c^1_{i}$, $c_{i,k_i} = c^2_{i}$.
We define
\index[syindex]{tb@$\mathfrak t^{b}$}
\begin{equation}\label{defntttt}
\mathfrak t^{b_i}({\bf x}_i) = e^{b_{i,0}} \otimes x_{i,1}
\otimes e^{b_{i,1}}\otimes \dots \otimes e^{b_{i,k_i-1}}\otimes x_{i,k_i}
\otimes e^{b_{i,k_i}}.
\end{equation}
\begin{sublem}\label{sublem56}\quad
\begin{enumerate}\itemsep=0pt
\item[$(1)$]
$\mathfrak t^{b_i}$ is a $\Lambda_0$ module isomorphism.
\item[$(2)$]
$\Delta_i \circ \mathfrak t^{b_i} = \bigl(\mathfrak t^{b_1}\otimes \mathfrak t^{b_2}\bigr) \circ \Delta_i$.
\item[$(3)$]
$\hat d^{b_i}_i \circ \mathfrak t^{b_i} = \mathfrak t^{b_i} \circ \hat d_i$.
Here $\hat d^{b_i}_i$ is a coderivation induced from $\mathfrak m^{b_i}$.
\end{enumerate}

\end{sublem}
The proof is the same as the proof of \cite[Lemma 5.2.12]{fooobook}
and so is omitted.
By Sublemma~\ref{sublem56}\,(1) there exists uniquely a $\Lambda_0$ module
homomorphism
\begin{gather*}
\widehat{\mathscr{F}^s}\bigl(\bigl(c^1_1,b^1_1\bigr),\bigl(c^1_2,b^1_2\bigr);
\bigl(c^2_1,b^2_1\bigr),\bigl(c^2_2,b^2_2\bigr)\bigr)\colon\
B\mathscr C_1[1]\bigl(c^1_1,c^1_2\bigr)
\,\widehat{\otimes}\, B\mathscr C_2[1]\bigl(c^2_1,c^2_2\bigr) \\
\qquad\to
B\mathscr C_3[1]\bigl(\bigl(c^3_1,b^3_1\bigr),\bigl(c^3_2,b^3_2\bigr)\bigr)
\end{gather*}
\big(where $\bigl(c^3_i,b^3_i\bigr) = \mathscr{F}^s_{\rm ob}\bigl(\bigl(c^1_i,b^1_i\bigr),\bigl(c^2_i,b^2_i\bigr)\bigr)$\big)
such that
\smash{$
\mathfrak t^{b_3} \circ
\widehat{\mathscr{F}^s}
=
\widehat{\mathscr{F}}
\circ \bigl(\mathfrak t^{b_1} \otimes \mathfrak t^{b_2}\bigr)$}.
Here and hereafter, we write \smash{$\widehat{\mathscr{F}^s}$}
in place of
\smash{$\widehat{\mathscr{F}^s}\bigl(\bigl(c^1_1,b^1_1\bigr),\bigl(c^1_2,b^1_2\bigr);
\bigl(c^2_1,b^2_1\bigr),\bigl(c^2_2,b^2_2\bigr)\bigr)$},
for simplicity.
Sublem\-ma~\ref{sublem56}\,(3) implies
\begin{equation}\label{form511}
\hat d^{b_3} \circ
\widehat{\mathscr{F}^s}
=
\widehat{\mathscr{F}^s}
\circ \bigl(\hat d^{b_1} \,\widehat{\otimes}\, {\rm id} + {\rm id} \,\widehat{\otimes}\, \hat d^{b_2}\bigr).
\end{equation}
Sublemma~\ref{sublem56}\,(2) implies
\begin{equation}\label{form512}
\Delta_3 \circ
\widehat{\mathscr{F}^s}
=
\widehat{\mathscr{F}^s}
\circ \bigl(\Delta_1 \,\widehat{\otimes}\, {\rm id} + {\rm id}\, \widehat{\otimes}\, \Delta_2\bigr).
\end{equation}
\eqref{form512} implies that \smash{$\widehat{\mathscr{F}^s}$} is induced by
\smash{$\mathscr{F}^s_{k_1,k_2}$}.
In fact, \smash{$\mathscr{F}^s_{k_1,k_2}$} is a composition of the restriction of \smash{$\widehat{\mathscr{F}^s}$}
to
$B_{k_1}\mathscr C_1[1]\bigl(c^1_1,c^1_2\bigr)
\,\widehat{\otimes}\, B_{k_2}\mathscr C_2[1]\bigl(c^2_1,c^2_2\bigr)$
and the projection
$
B\mathscr C_3[1]\bigl(\bigl(c^3_1,b^3_1\bigr),\bigl(c^3_2,b^3_2\bigr)\bigr)
\to \mathscr C_3[1]\bigl(\bigl(c^3_1,b^3_1\bigr),\bigl(c^3_2,b^3_2\bigr)\bigr)
$.

Then \eqref{form511} implies that it satisfies the required property,
Definition~\ref{bi-functor}\,(4).

To show that $\mathscr{F}^s$ is strict, we observe
\[
\widehat{\mathscr{F}}
\circ \bigl(\mathfrak t^{b_1}(1) \otimes \mathfrak t^{b_2}(1)\bigr)
=
\widehat{\mathscr{F}}\bigl(e^{b_1}, e^{b_2}\bigr)
= e^{b_3} = \mathfrak t^{b_3}(1).
\]
Namely, \smash{$\widehat{\mathscr{F}}^s(1) = 1$}.
This implies ${\mathscr{F}}^s_{0,0}(1) = 0$.
\end{proof}

In the case when $\mathscr C$ is curved, we can not define
the filtered $A_{\infty}$ bi-functor in Example~\ref{ex53}.
However, we can still use the language of filtered $A_{\infty}$
bi-module to define a similar object.

Let $\mathscr C_1$, $\mathscr C_2$ be non-unital curved
filtered $A_{\infty}$ categories.
We define the notion of a left-$\mathscr C_1$ and right-$\mathscr C_2$
bi-module as follows.
\begin{defn}\label{bimodulecat}
A left-$\mathscr C_1$ and right-$\mathscr C_2$
filtered $A_{\infty}$ bi-module, is $\mathscr D = (\{D_{c_1,c_2}\},\{\mathfrak n_{c_1',c_1,c_2,c_2'}\})$, where
\begin{enumerate}\itemsep=0pt
\item[(1)]
The object $\{D_{c_1,c_2}\}$ assigns a completed free graded $\Lambda_0$ module
$D_{c_1,c_2}$ to each $c_1 \in \mathfrak{OB}(\mathscr C_1)$,
$c_2 \in \mathfrak{OB}(\mathscr C_2)$.
\item[(2)]
To each $c_1, c'_1 \in \mathfrak{OB}(\mathscr C_1)$,
$c_2, c'_2 \in \mathfrak{OB}(\mathscr C_2)$, we are given a
$\Lambda_0$ module homomorphism \index[syindex]{nc1'c1c2@$\mathfrak n_{c_1',c_1,c_2,c_2'}$}
\[
\mathfrak n_{c_1',c_1,c_2,c_2'} \colon\
B\mathscr C_1[1](c'_1,c_1) \,\widehat\otimes\, D_{c_1,c_2}
\,\widehat\otimes\, B\mathscr C_2[1](c_2,c'_2)
\to D_{c'_1,c'_2}
\]
of degree $+1$ which preserves the energy filtration.
\item[(3)]
The following $A_{\infty}$ relation is satisfied:
\begin{align}
0={}&\sum_{a_1}\sum_{a_2}(-1)^{*_1}\mathfrak n(\text{\bf x}_{1:a_1}, \mathfrak n(\text{\bf x}_{2:a_1},
z,\text{\bf y}_{1:a_2}),\text{\bf y}_{2:a_2})\nonumber
\\
& +
\mathfrak n\bigl(\hat d_1(\text{\bf x}),
z,\text{\bf y}\bigr)
+
(-1)^{*_2}\mathfrak n\bigl(\text{\bf x},
z,\hat d_2(\text{\bf y})\bigr).\label{form9300}
\end{align}
Here $*_1 = \deg'\text{\bf x}_{1:a_1}$,
$*_2 = \deg'\text{\bf x} + \deg' z$.
The notations are as follows:
${\bf x} \in B\mathscr C_1[1](c'_1,c_1)$,
$y \in D_{c_1,c_2}$,
${\bf z} \in B\mathscr C_2[1](c_2,c'_2)$.
$\Delta {\bf x} = \sum_{a_1}\text{\bf x}_{1:a_1}\otimes \text{\bf x}_{2:a_1}$.
$\Delta {\bf y} = \sum_{a_2}\text{\bf y}_{1:a_2}\otimes \text{\bf y}_{2:a_2}$.
The symbol \smash{$\hat d_i$} denotes the
coderivation induced by the structure operations of $\mathscr C_i$
and is defined in \eqref{hatddd}.
\end{enumerate}
A filtered $A_{\infty}$ bi-module over $G$-gapped unital curved filtered $A_{\infty}$
category is said to be $G$-gapped if all the structure operations are $G$-gapped.
It is said to be {\it unital} if the following holds:
\begin{enumerate}\itemsep=0pt
\item[(1)]
The equality
$
\mathfrak n_{1,0}({\bf e}_1,y) = (-1)^{\deg y}\mathfrak n_{1,0}(y,{\bf e}_2)= y$.
Here ${\bf e}_i$ is the unit of $\mathscr C_i$.
\item[(2)]
If ${\bf x}$ or ${\bf z}$ contains the unit, then $\mathfrak n({\bf x};y;{\bf z}) = 0$
except the cases appearing in item (1).
\end{enumerate}

\end{defn}
We define
\begin{gather*}
\widehat{\mathfrak n}_{c_1',c_2'} \colon\ \bigoplus_{c_1,c_2}
B\mathscr C_1[1](c'_1,c_1) \,\widehat\otimes\, D_{c_1,c_2}
\,\widehat\otimes\, B\mathscr C_2[1](c_2,c'_2) \\
\hphantom{\widehat{\mathfrak n}_{c_1',c_2'} \colon} \ \to \bigoplus_{c_1,c_2} B\mathscr C_1[1](c'_1,c_1) \,\widehat\otimes\,
D_{c'_1,c'_2}
\,\widehat\otimes\, B\mathscr C_2[1](c_2,c'_2),
\end{gather*}
by
\begin{align*}
\widehat{\mathfrak n}_{c_1',c_2'}(\text{\bf x},z,\text{\bf y})
={}
&\sum_{a_1}\sum_{a_2}(-1)^{*_1} \text{\bf x}_{1:a_1}\otimes \mathfrak n(\text{\bf x}_{2:a_1}
z,\text{\bf y}_{1:a_2})\otimes\text{\bf y}_{2:a_2}
\\
& +\hat d_1(\text{\bf x}) \otimes
z\otimes\text{\bf y}
+
(-1)^{*_2}\text{\bf x}\otimes
z\otimes\hat d_2(\text{\bf y}),
\end{align*}
where the notations are as in \eqref{form9300}.
Then the formula \eqref{form9300} is
equivalent to
\smash{$
\widehat{\mathfrak n}_{c_1',c_2'} \circ \widehat{\mathfrak n}_{c_1',c_2'} = 0$}.

\begin{defn}\label{defnbomodhomocat}
Let \smash{$\mathscr D^{(i)} = \big(\bigl\{D^{(i)}_{c_1,c_2}\bigr\},\bigl\{\mathfrak n^{(i)}_{c_1',c_1,c_2,c_2'}\bigr\}\big)$}
be a left-$\mathscr C_1$ and right-$\mathscr C_2$
filtered $A_{\infty}$ bi-module, for $i=1,2$.
A {\it pre-bi-module homomorphism} \index{pre-bi-module homomorphism} of degree $d$ from \smash{$\mathscr D^{(1)}$} to
\smash{$\mathscr D^{(2)}$}
is $\mathfrak f = \{\mathfrak f_{c_1',c_1,c_2,c_2'}\}$, where
\begin{enumerate}\itemsep=0pt\label{bomodhomocat}
\item[$(*)$]
To each $c_1, c'_1 \in \mathfrak{OB}(\mathscr C_1)$,
$c_2, c'_2 \in \mathfrak{OB}(\mathscr C_2)$, we are given a
$\Lambda_0$ module homomorphism
\[
\mathfrak f_{c_1',c_1,c_2,c_2'} \colon\
B\mathscr C_1[1](c'_1,c_1) \,\widehat\otimes\, D^{(1)}_{c_1,c_2}
\,\widehat\otimes\, B\mathscr C_2[1](c_2,c'_2)
\to D^{(2)}_{c'_1,c'_2},
\]
of degree $d$ which preserves the energy filtration.
\end{enumerate}
Let
\begin{gather}
\widehat{\mathfrak f}_{c_1',c_2'}\colon\
\bigoplus_{c_1,c_2}B\mathscr C_1[1](c'_1,c_1) \,\widehat\otimes\, D^{(1)}_{c_1,c_2}
\,\widehat\otimes\, B\mathscr C_2[1](c_2,c'_2) \nonumber\\
\hphantom{\widehat{\mathfrak f}_{c_1',c_2'}\colon} \ \to \bigoplus_{c_1,c_2} B\mathscr C_1[1](c'_1,c_1) \,\widehat\otimes\, D^{(2)}_{c_1,c_2}
\,\widehat\otimes\, B\mathscr C_2[1](c_2,c'_2),\label{defhatf}
\end{gather}
be the formal bi-comodule homomorphism
induced from $\mathfrak f_{c_1',c_1,c_2,c_2'}$.
More explicitly, the map \eqref{defhatf} is defined by
\[
\widehat{\mathfrak f}_{c_1',c_2'}(\text{\bf x},z,\text{\bf y})
:=
\sum_{a_1}\sum_{a_2}(-1)^{*}\text{\bf x}_{1:a_1}\otimes \mathfrak f(\text{\bf x}_{2:a_1}
z,\text{\bf y}_{1:a_2})\otimes\text{\bf y}_{2:a_2},
\]
where $* = \deg \mathfrak f\deg'\text{\bf x}_{1:a_1}
= \mathfrak{deg}' \mathfrak f \deg' {\bf x}_{1:a_1}$.
(Note that $\mathfrak{deg}' \mathfrak f = \deg \mathfrak f$, see Definition~\ref{defn215555}.)

We define
a pre-bi-module homomorphism $d(\mathfrak f)$
of degree $\deg \mathfrak f + 1$, so that
\[
\widehat{d(\mathfrak f)}
:=
\widehat{\mathfrak n}
\circ \widehat{\mathfrak f}
- (-1)^{\deg \mathfrak f}
\widehat{\mathfrak f}
\circ \widehat{\mathfrak n}
\]
holds.

We say a pre-bi-module homomorphism $\mathfrak f$
is a {\it bi-module homomorphism} \index{bi-module homomorphism} if its degree is $0$
and if $d(\mathfrak f)= 0$.
When $\mathfrak g = \{\mathfrak g_{c_1',c_1,c_2,c_2'}\}$ is
another pre-bimodule homomorphism, we define
$\mathfrak g \circ \mathfrak f$ so that
\smash{$
\widehat{\mathfrak g \circ \mathfrak f}
=
\widehat{\mathfrak g} \circ \widehat{\mathfrak f}$}.

\end{defn}
Note that $d(\mathfrak f)= 0$ is equivalent to
\begin{gather*}
\sum_{a_1}\sum_{a_2}\mathfrak n(\text{\bf x}_{1:a_1}, \mathfrak f(\text{\bf x}_{2:a_1},
z,\text{\bf y}_{1:a_2}),\text{\bf y}_{2:a_2})\\
\qquad
=
\mathfrak f\big(\hat d_1(\text{\bf x}),
z,\text{\bf y}\big)
+
(-1)^{\deg' {\bf x}+\deg z}\mathfrak f\big(\text{\bf x},
z,\hat d_2(\text{\bf y})\big)
\\
\phantom{\qquad
=}{}+
\sum_{a_1}\sum_{a_2}(-1)^{\deg' \text{\bf x}_{1:a_1}}\mathfrak f(\text{\bf x}_{1:a_1}, \mathfrak n(\text{\bf x}_{2:a_1},
z,\text{\bf y}_{1:a_2}),\text{\bf y}_{2:a_2}).
\end{gather*}
\begin{defn}
We define a DG-category $\mathcal{BIMOD}(\mathscr C_1,\mathscr C_2)$ \index[syindex]{BIMOD@$\mathcal{BIMOD}(\mathscr C_1,\mathscr C_2)$}
as follows:
\begin{enumerate}\itemsep=0pt
\item[(1)]
Its object is a left-$\mathscr C_1$, right-$\mathscr C_2$ filtered bi-module.
\item[(2)]
For two objects $\mathscr D_1$ and $\mathscr D_2$, a morphism
from $\mathscr D_1$ to $\mathscr D_2$ is a pre-filtered $A_{\infty}$-bimodule
homomorphism. We write it as $\mathcal{BIMOD}(\mathscr D_1,\mathscr D_2)$.
\item[(3)]
The composition and the differential of $\mathcal{BIMOD}(\mathscr C_1,\mathscr C_2)$
are defined as in Definition~\ref{bomodhomocat}.
\end{enumerate}

\end{defn}
It is obvious from definition that $\mathcal{BIMOD}(\mathscr C_1,\mathscr C_2)$
is a DG-category.
\begin{defn}
In Definitions \ref{bimodulecat} and \ref{bomodhomocat},
we can define $G$-gappedness and/or unitality
of bi-module and/or bi-module homomorphism
in an obvious way if $\mathscr C_i$ is
$G$-gapped and/or unital for $i=1,2$.

\end{defn}
We next explain the relation between a filtered $A_{\infty}$
bi-module and a bi-functor.
We need a~digression for this purpose.
\begin{defn}\label{lem56}
Let $\mathscr{C}_1$, $\mathscr{C}_2$, $\mathscr{C}_3$ be
strict non-unital curved filtered $A_{\infty}$ categories.
We will define bijections between the
following three objects:
\begin{enumerate}\itemsep=0pt
\item[(1)]
A filtered $A_{\infty}$
bi-functor
$
\mathscr{F} \colon \mathscr{C}_1 \times \mathscr{C}_2 \to \mathscr{C}_3.
$
\item[(2)]
A filtered $A_{\infty}$
bi-functor
$
\mathscr{F} \colon \mathscr{C}_2 \times \mathscr{C}_1 \to \mathscr{C}_3.
$
\item[(3)]
A filtered $A_{\infty}$ functor:
$
\mathscr{G} \colon \mathscr{C}_1 \to \mathcal{FUNC}(\mathscr{C}_2,\mathscr{C}_3).
$
\end{enumerate}
The bijection between (1) and (2) is constructed by using the
isomorphism
\[
\mathcal S \colon\ B_{k_1}\mathscr C_1[1]((c_{1,1},c_{1,2})
\,\widehat{\otimes}\, B_{k_2}\mathscr C_2[1]((c_{2,1},c_{2,2})
\to
B_{k_2}\mathscr C_2[1]((c_{2,1},c_{2,2})
\,\widehat{\otimes}\, B_{k_1}\mathscr C_1[1]((c_{1,1},c_{1,2}),
\]
which is
$
\mathcal S({\bf x} \otimes {\bf y})
=
(-1)^{\deg'{\bf x}\deg'{\bf y}}{\bf y} \otimes {\bf x}
$.

We next construct bijection between (1) and (3).

Suppose $
\mathscr{F}$ is given as in (1).
Let $c_1 \in \mathfrak{OB}(\mathscr C_1)$.
We first construct $\mathscr{G}_{\rm ob}(c_1)$
which is a~fil\-tered~$A_{\infty}$ functor: $
\mathscr C_2 \to \mathscr C_3$.
Let $c_2 \in \mathfrak{OB}(\mathscr C_2)$.
Then we put
$
(\mathscr{G}_{\rm ob}(c_1))_{\rm ob}(c_2) = \mathscr{F}_{\rm ob}(c_1,c_2)
\in \mathfrak{OB}(\mathscr C_3)$.
Let $c_{2,1}, c_{2,2} \in \mathfrak{OB}(\mathscr C_2)$.
We define
\[
(\mathscr{G}_{\rm ob}(c_1))_{k_2}(c_{2,1}, c_{2,2})
\colon\
B_{k_2}\mathscr C_2[1](c_{2,1},c_{2,2}) \to \mathscr C_3[1]
(\mathscr{F}_{\rm ob}(c_1,c_{2,1}),\mathscr{F}_{\rm ob}(c_1,c_{2,2}))
\]
by
$
(\mathscr{G}_{\rm ob}(c_1))_{k_2}(c_{2,1}, c_{2,2})({\bf y})
=
\mathscr{F}(c_1,c_1;c_{2,1},c_{2,2})_{0,k_2}
(1,{\bf y})$,
where
$
1 \in B_0\mathscr C_1[1](c_1,c_1) = \Lambda_0$,
$
{\bf y}
\in B_{k_2}\mathscr C_2[1](c_{2,1},c_{2,2})$.
We thus defined $\mathscr{G}_{\rm ob}(c_1)$, which is a~filtered
$A_{\infty}$ functor: $\mathscr C_2 \to \mathscr C_3$.

Let $c_{1,1},c_{1,2} \in \mathfrak{OB}(\mathscr C_1)$
and ${\bf x} \in B_{k_1}\mathscr C_1[1](c_{1,1},c_{1,2})$.

We will construct
$
\mathscr{G}_{k_1}(c_{1,1},c_{1,2})({\bf x})$,
which is a pre-natural transformation from
$\mathscr{G}_{\rm ob}(c_{1,1})$ to~$\mathscr{G}_{\rm ob}(c_{1,2})$.

Let $c_{2,1}, c_{2,2} \in \mathfrak{OB}(\mathscr C_2)$
and
${\bf y}
\in B_{k_2}\mathscr C_2[1](c_{2,1},c_{2,2})$.
Then
\[
(\mathscr{G}_{k_1}(c_{1,1},c_{1,2})({\bf x}))_{k_2}(c_{2,1},c_{2,2})
\colon\ B_{k_2}\mathscr C_2[1](c_{2,1},c_{2,2}) \to
\mathscr C_3[1]
(\mathscr{F}_{\rm ob}(c_{1,1},c_{2,1}),\mathscr{F}_{\rm ob}(c_{1,2},c_{2,2}))
\]
is defined by
\[
(\mathscr{G}_{k_1}(c_{1,1},c_{1,2})({\bf x}))_{k_2}(c_{2,1},c_{2,2})({\bf y})
=
\mathscr F_{k_1,k_2}(c_{1,1},c_{1,2};c_{2,1},c_{2,2}) ({\bf x},{\bf y}).
\]
It is straightforward to check that $\mathscr{G}$ is a filtered
$A_{\infty}$ functor.

The construction from (3) to (1) can be done by doing the same construction in the
opposite direction.

\end{defn}
\begin{exm}\label{ex58}
Let $
\mathscr F \colon \mathscr C^{\rm op} \times \mathscr C
\to \mathcal{CH}
$ be the filtered $A_{\infty}$
bi-functor in Example~\ref{ex53}.
Then by Definition~\ref{lem56}, we obtain a
filtered $A_{\infty}$ functor
$\mathscr C \to \mathcal{FUNC}(\mathscr C^{\rm op},\mathcal{CH})$.
This is nothing but the $A_{\infty}$ Yoneda functor.

\end{exm}
\begin{lem}\label{2122222333}
In the situation of Definition {\rm\ref{lem56}},
there exists an equivalence of $A_{\infty}$ categories
\[
\mathcal{FUNC}(\mathscr{C}_1,\mathcal{FUNC}(\mathscr{C}_2,\mathscr{C}_3))
\to \mathcal{FUNC}(\mathscr{C}_2,\mathcal{FUNC}(\mathscr{C}_1,\mathscr{C}_3)).
\]

\end{lem}
\begin{proof}
The proof is similar to the above construction and is a
straightforward calculation.
\end{proof}

Using Lemma~\ref{2122222333} and Definition~\ref{lem56},
we obtain a filtered $A_{\infty}$ category so
that its object is a filtered $A_{\infty}$ bi-functor:
$\mathscr{C}_1 \times \mathscr{C}_2 \to \mathscr{C}_3$.
This filtered $A_{\infty}$ category
is equivalent to
$\mathcal{FUNC}(\mathscr{C}_1,\mathcal{FUNC}(\mathscr{C}_2,\mathscr{C}_3))$
by definition.
We denote this filtered $A_{\infty}$ category
by
$
\mathcal{BIFUNC}(\mathscr{C}_1\allowbreak\times \mathscr{C}_2,\mathscr{C}_3).
$\index[syindex]{BIFUNC@$\mathcal{BIFUNC}(\mathscr{C}_1\times \mathscr{C}_2,\mathscr{C}_3)$}
A morphism between two filtered $A_{\infty}$
bi-functors in this category is called a
pre-natural transformation.
It is called a natural transformation
if its $\mathfrak m_1$ derivative is zero.

Note that during the discussion of Definition~\ref{lem56} and Lemma~\ref{2122222333},
we required the filtered~$A_{\infty}$ categories to be strict.
\begin{lem}\label{lem55214}
In the situation of Definitions {\rm\ref{bimodulecat}} and {\rm\ref{bomodhomocat}},
we assume $\mathscr C_i$ is strict for $i=1,2$.
Then there exists an equivalence of DG-categories
\[
\mathcal{BIFUNC}(\mathscr{C}^{\rm op}_1\times\mathscr{C}_2,\mathcal{CH})
\cong
\mathcal{BIMOD}(\mathscr C_1,\mathscr C_2).
\]

\end{lem}
\begin{proof}
In the same way as
Example~\ref{ex54}, we can find a bijection between the sets of
bi-modules and of bi-functors appearing as objects of the above
two categories.
The fact that morphisms and structure operations coincide
can be proved by a similar straightforward calculations.
\end{proof}

\begin{rem}
Note that in the case when $\mathscr C_1$, $\mathscr C_2$
are curved the author does not know
the way
to define a filtered $A_{\infty}$ category $\mathcal{BIFUNC}(\mathscr{C}_1\times\mathscr{C}_2,\mathscr{C}_3)$.
Only in the case when $\mathscr{C}_3 = \mathcal{CH}$,
we can use the notion of bi-module to define
DG-category equivalent to
$\mathcal{BIFUNC}(\mathscr{C}_1\times\mathscr{C}_2,\mathcal{CH})$
for curved categories $\mathscr C_1$, $\mathscr C_2$.
The functor category in the curved case is defined in \cite{DL},
which may be adapted to the bi-functor case.
\end{rem}
The next lemma is an analogue of Lemma~\ref{lem55}.
\begin{lem}\label{lem55revrev}
Let $\mathscr{C}_1$, $\mathscr{C}_2$ be
non-unital curved filtered $A_{\infty}$ categories
and $\mathscr{C}^s_1$, $\mathscr{C}^s_2$
their associated strict categories.
Then
a left-$\mathscr{C}_1$ and right-$\mathscr{C}_2$ filtered $A_{\infty}$
bi-module
$\mathscr D = (\{D_{c_1,c_2}\},\{\mathfrak n_{c_1',c_1,c_2,c_2'}\})$
induces
a~left-$\mathscr{C}^s_1$ and right-$\mathscr{C}^s_2$ filtered $A_{\infty}$
bi-functor
$\mathscr D^s$.
If $\mathscr D$ is unital or $G$-gapped, then
so is~$\mathscr D^s$.

If $\mathscr D_1$, $\mathscr D_2$ are left-$\mathscr{C}_1$ and right-$\mathscr{C}_2$ filtered $A_{\infty}$ bi-modules and $\mathfrak f$ is a pre-filtered $A_{\infty}$ bi-module homomorphism
from $\mathscr D_1$ to $\mathscr D_2$.
Then we can associate a pre-filtered bi-module homomorphism~$\mathfrak f^s$ from $\mathscr D^s_1$ to $\mathscr D^s_2$.
It induces a DG-functor from $\mathcal{BIMOD}(\mathscr C_1,\mathscr C_2)$
to $\mathcal{BIMOD}(\mathscr C^s_1,\mathscr C^s_2)$.
\end{lem}
\begin{proof}
The proof is the same as the proof of Lemma~\ref{lem55}.
In fact, we can take
$
D^s_{(c_1,b_1),(c_2,b_2)}: = D_{c_1,c_2}$,
and
\[
\mathfrak n^s(x_1,\dots,x_k;z;y_1,\dots,y_{\ell}):
=
\mathfrak n\bigl(e^b,x_1,e^b,\dots,e^b,x_k,e^b;z;e^b,y_1,e^b,\dots,e^b,y_{\ell},e^b\bigr).
\]
The proof of the statement on pre-filtered $A_{\infty}$ bi-module homomorphism
can be proved in the same way as \cite[Section~5.2.2.3]{fooobook}.
\end{proof}

We next discuss a composition of
$A_{\infty}$ (bi)-functors and pullback of bi-module by $A_{\infty}$ (bi)-functors.
To discuss them systematically we introduce the notion of a multi-$A_{\infty}$
functor.

\begin{defn}\label{tri-functor}
Let $m$ be a positive integer and let
$\mathscr{C}_i$, $i=1,\dots,m$, and $\mathscr{C}'$ be
non-unital curved filtered $A_{\infty}$ categories.
A {\it filtered $A_{\infty}$ multi-functor}\index{filtered $A_{\infty}$ multi-functor}
$
\mathscr{F} \colon \mathscr{C}_1 \times \cdots \times \mathscr{C}_m \to \mathscr{C}'
$
consists of~$\mathscr{F}_{\rm ob}$ and $\mathscr{F}_{k_1,\dots,k_m}$, $k_i = 0,1,2,3,\dots$, of degree $0$ such that
\begin{enumerate}\itemsep=0pt
\item[(1)]
A map:
$\mathscr{F}_{\rm ob} \colon \prod_{i=1}^m\mathfrak{OB}(\mathscr{C}_j)
\to\mathfrak{OB}(\mathscr{C}')$
is given.
\item[(2)]
Let $c_{i,1},c_{i,2} \in \mathfrak{OB}(\mathscr{C}_i)$, $i=1,\dots,m$.
We put $\vec c_j = (c_{1,j},\dots, c_{m,j})$, for $j=1,2$.
$\mathscr{F}_{k_1,\dots,k_m}$ associates a $\Lambda_0$
linear map
\[
\mathscr F_{k_1,\dots,k_m}(\vec c_{1};\vec c_{2})
\colon\ \bigotimes_{i=1}^m B_{k_i}\mathscr C_i[1](c_{i,1},c_{i,2})
\to
\mathscr C'[1](\mathscr F_{\rm ob}(\vec c_1),\mathscr F_{\rm ob}(\vec c_2))
\]
of degree $0$.
\item[(3)]
We require that $\mathscr F_{k_1,\dots,k_m}(\vec c_{1};\vec c_{2}) $
preserves the filtration in a similar sense as
Definition~\ref{defn22}\,(2).
\end{enumerate}
$\{\mathscr F_{k_1,\dots,k_m}\}$ induces uniquely a formal coalgebra
homomorphism
\[
\widehat{\mathscr F}(\vec c_{1};\vec c_{2})
\colon\ \bigotimes_{i=1}^3 B\mathscr C_i[1](c_{i,1},c_{i,2})
\to
B\mathscr C'[1](\mathscr F_{\rm ob}(\vec c_1),\mathscr F_{\rm ob}(\vec c_2)).
\]
Note that the structure operations of $\mathscr C_i$
induce coderivations
\[
\hat d_i \colon\ B\mathscr C_i[1](c_{i,1},c_{i,2})
\to B\mathscr C_i[1](c_{i,1},c_{i,2}).
\]
\begin{enumerate}\itemsep=0pt
\item[(4)] The homomorphism
$\widehat{\mathscr F}(\vec c_{1};\vec c_{2})$ is a cochain map.
\end{enumerate}
The unitality, strictness, $G$-gappedness of multi-functor are defined in the same way.

In the case when $m=3$, the multi-functor is called the {\it tri-functor}. \index{tri-functor}
\end{defn}
\begin{lem}\label{lem118ss}
A filtered $A_{\infty}$ multi-functor $\mathscr{F}$ induces a
strict filtered $A_{\infty}$ multi-functor $\mathscr{F}^s$ among the associated strict categories.
The unitality and/or $G$-gappedness is preserved.
\end{lem}

The proof is the same as Lemma~\ref{lem55} and is omitted.

Let
$\mathscr{C}_1,\ldots, \mathscr{C}_m$
and $\mathscr{C}'_1,\ldots, \mathscr{C}'_{m'}$ be
non-unial filtered $A_{\infty}$ categories
and
$
\mathscr F
\colon \mathscr{C}_1 \times \dots\times \mathscr{C}_m \to \mathscr{C}'_k$
and
$
\mathscr G
\colon \mathscr{C}'_1 \times \dots \times \mathscr{C}'_{\ell}
\to \mathscr{C}''$
be $A_{\infty}$ multi-functors.
We define its {\it composition}
\[
\mathscr G \circ \mathscr F\colon\
\mathscr{C}'_1 \times \dots \times\mathscr{C}'_{k-1} \times
\mathscr{C}_1 \times \dots\times \mathscr{C}_m
\times \mathscr{C}'_{k+1} \times \dots \times \mathscr{C}'_{\ell}
\to \mathscr{C}''
\]
by
\begin{gather}
(\mathscr G \circ \mathscr F)
({\bf x}_1 \otimes \dots \otimes {\bf x}_{k-1} \otimes {\bf y}_1 \otimes\dots \otimes {\bf y}_{m}
\otimes {\bf x}_{k+1} \otimes\dots \otimes {\bf x}_{\ell})\nonumber \\
\qquad=
\mathscr G
(
{\bf x}_1 \otimes \dots \otimes {\bf x}_{k-1} \otimes
\mathscr F({\bf y}_1 \otimes\dots \otimes {\bf y}_{m})
 \otimes {\bf x}_{k+1} \otimes\dots \otimes {\bf x}_{\ell}
).\label{henhen224}
\end{gather}
It is easy to check that \eqref{henhen224} defines a multi-functor.

\begin{lem}\label{lem256}
Suppose $\mathscr{C}_1,\dots,\mathscr{C}_m$
and $\mathscr{C}'$ are strict.
Then, we can define a filtered $A_{\infty}$
categories~$\mathcal{MULFUNC}(\mathscr{C}_1\times\cdots \times \mathscr{C}_m,\mathscr{C}')$
such that
\begin{enumerate}\itemsep=0pt
\item[$(1)$]
Its object is a filtered $A_{\infty}$ multi-functor
$
\mathscr{F} \colon \mathscr{C}_1 \times \cdots \times \mathscr{C}_m \to \mathscr{C}'
$.
\item[$(2)$]
There exists a filtered $A_{\infty}$ bi-functor
\begin{gather*}
\mathcal{MULFUNC}(\mathscr{C}_1\times\cdots \times \mathscr{C}_m,\mathscr{C}'')
\times
\mathcal{MULFUNC}(\mathscr{C}'_1\times\cdots \times \mathscr{C}'_{m'},\mathscr{C}_k) \\
\qquad\to
\mathcal{MULFUNC}(\mathscr{C}_1\times\cdots \times\mathscr{C}_{k-1} \times
\mathscr{C}'_1\times\cdots \times \mathscr{C}'_{m'}
\times \mathscr{C}_{k+1} \times \cdots\times \mathscr{C}_m,\mathscr{C}'')
\end{gather*}
such that $(\mathscr F,\mathscr G) \mapsto \mathscr G \circ \mathscr F$
is its object part.
\end{enumerate}

\end{lem}
The proof of (1) is straightforward.
(2) is a straightforward generalization of Theorem~\ref{thm94}.

Now it is rather obvious how to define the notion of multi-module over
(curved) filtered $A_{\infty}$ categories and define the notion of a pullback of a
multi-module structure by multi-functor.
We explain it below since we will use it.
\begin{defn}\label{bimodulecat2}
Let
$\mathscr{C}_{l,1},\ldots, \mathscr{C}_{l,m}$
and $\mathscr{C}_{r,1},\ldots, \mathscr{C}_{r,m'}$ be
non-unial filtered $A_{\infty}$ categories
A left-$\mathscr{C}_{l,1},\ldots, \mathscr{C}_{l,m}$ and right-$\mathscr{C}_{r,1},\ldots, \mathscr{C}_{r,m'}$
filtered $A_{\infty}$ multi-module, is $
(\{D_{\vec c_{l},\vec c_{r}}\},\{\mathfrak n_{\vec c_l,\vec c'_l,\vec c_r,\vec c'_r}\})$, where\index{filtered $A_{\infty}$ multi-module}
\begin{enumerate}\itemsep=0pt
\item[(1)]
To each \smash{$\vec c_{l} \in \prod_{i=1}^{m_l}\mathfrak{OB}(\mathscr{C}_{j,l})$},
\smash{$\vec c_{r} \in \prod_{i=1}^{m_r}\mathfrak{OB}(\mathscr{C}_{j,r})$},
a graded completed free $\Lambda_0$ module~$D_{\vec c_{l},\vec c_{r}}$
is assigned.
\item[(2)]
To each \smash{$\vec c_{l}, \vec c'_{l} \in \prod_{i=1}^{m_l}\mathfrak{OB}(\mathscr{C}_{j,l})$},
\smash{$\vec c_{r}, \vec c'_{r} \in \prod_{i=1}^{m_r}\mathfrak{OB}(\mathscr{C}_{j,r})$}, we are given a
$\Lambda_0$ module homomorphism
\[
\mathfrak n_{\vec c_{l},\vec c'_{l},\vec c_{r},\vec c'_{r}} \colon\
\bigotimes_{j=1}^{m_l}B\mathscr C_1(c'_{j,l},c_{j,l}) \,\widehat\otimes\, D_{\vec c_{l},\vec c_{r}}
\,\widehat\otimes\,\bigotimes_{j=1}^{m_r}B\mathscr C_r(c_{j,r},c'_{j,r})
\to D_{\vec c'_{l},\vec c'_{r}},
\]
of degree $+1$ which preserves the energy filtration.
\end{enumerate}
In case $m+m' = 3$, we call it a {\it tri-module}. \index{tri-module}

The unitality and/or gappedness of multi-module over unital and/or
gapped categories are defined in an obvious way.

When \smash{$\mathscr D^{\ell} =
\big(\bigl\{D^{\ell}_{\vec c_{l},\vec c_{r}}\bigr\},\{\mathfrak n_{\vec c_l,\vec c'_l,\vec c_r,\vec c'_r}\}\big)$}
is a left-$\mathscr{C}_{l,1},\ldots, \mathscr{C}_{l,m}$ and right-$\mathscr{C}_{r,1},\ldots, \mathscr{C}_{r,m'}$
filtered $A_{\infty}$ multi-module for $\ell =1,2$,
a multi-module pre-homomorphism \index{multi-module pre-homomorphism}
from $\mathscr D^{1}$ to $\mathscr D^{2}$ of degree $d$
is~${\mathfrak f = \{\mathfrak f_{\vec c_{l},\vec c_{r}}\}}$ where
the map $\mathfrak f_{\vec c_{l},\vec c_{r}}$,
\[
\mathfrak f_{\vec c_{l},\vec c'_{l},\vec c_{r},\vec c'_{r}} \colon\
\bigotimes_{j=1}^{m_l}B\mathscr C_1(c'_{j,l},c_{j,l}) \,\widehat\otimes\, D^1_{\vec c_{l},\vec c_{r}}
\,\widehat\otimes\,\bigotimes_{j=1}^{m_r}B\mathscr C_r(c_{j,r},c'_{j,r})
\to D^2_{\vec c'_{l},\vec c'_{r}},
\]
is a degree $d \Lambda_0$ module homomorphism which preserves filtration.

The maps
$\mathfrak f$ induces a formal bi-comodule homomorphism
\begin{gather*}
\widehat{\mathfrak f}_{\vec c_{l},\vec c_{r}} \colon\
\bigoplus_{c'_{j,l},c'_{j,r}} \bigotimes_{j=1}^{m_l}B\mathscr C_1(c'_{j,l},c_{j,l}) \,\widehat\otimes\, D^1_{\vec c_{l},\vec c_{r}}
\bigotimes_{j=1}^{m_r}B\mathscr C_r(c'_{j,r},c_{j,r})
\\
\hphantom{\widehat{\mathfrak f}_{\vec c_{l},\vec c_{r}} \colon}{} \
\to \bigotimes_{j=1}^{m_l}B\mathscr C_1(c'_{j,l},c_{j,l}) \,\widehat\otimes\, D^2_{\vec c_{l},\vec c_{r}}
\bigotimes_{j=1}^{m_r}B\mathscr C_r(c'_{j,r},c_{j,r}),
\end{gather*}
in the same way as Definition~\ref{tri-functor}\,(1).

We define $d\mathfrak f = \{(d\mathfrak f)_{{\vec c_{l},\vec c_{r}}}\}$ by
\[
(d\mathfrak f)_{{\vec c_{l},\vec c_{r}}} =
\hat{\mathfrak n}^2_{\vec c_{l},\vec c_{r}} \circ \widehat{\mathfrak f}_{\vec c_{l},\vec c_{r}}
-(-1)^{\deg \mathfrak f}
\widehat{\mathfrak f}_{\vec c_{l},\vec c_{r}} \circ
\hat{\mathfrak n}^1_{\vec c_{l},\vec c_{r}}.
\]
Here $\hat{\mathfrak n}^{\ell}_{\vec c_{l},\vec c_{r}}$ is the boundary
operator induced from the structure operations of $\mathscr D^{\ell}$
as in item (3) above.

When $\mathfrak f$ (resp.\ $\mathfrak g$) is a multi-module pre-homomorphism
from $\mathscr D^{1}$ to $\mathscr D^{2}$
\big(resp.\ $\mathscr D^{2}$ to $\mathscr D^{3}$\big),
we define a multi-module pre-homomorphism $\mathfrak f \circ \mathfrak g$
from $\mathscr D^{1}$ to $\mathscr D^{3}$
by
\smash{$\widehat{\mathfrak f \circ \mathfrak g}
= \widehat{\mathfrak f} \circ \widehat{\mathfrak g}$}.

{\samepage Thus we obtain the following filtered DG-category:
\begin{enumerate}\itemsep=0pt
\item[(1)] Its object is a left-$\mathscr{C}_{l,1},\ldots, \mathscr{C}_{l,m}$ and right-$\mathscr{C}_{r,1},\ldots, \mathscr{C}_{r,m'}$
filtered $A_{\infty}$ multi-module.
\item[(2)] The module of morphisms is the set of multi-module pre-homomorphisms.
\item[(3)] The differential $d$ and composition $\circ$ is defined as above.
\end{enumerate}}%
The unitality and/or gappedness of a multi-module homomorphism is defined in an obvious way.
\end{defn}

To a left-$\mathscr{C}_{l,1},\ldots, \mathscr{C}_{l,m}$ and right-$\mathscr{C}_{r,1},\ldots, \mathscr{C}_{r,m'}$
filtered $A_{\infty}$ multi-module $\mathscr D$
we can associate a~left-$\mathscr{C}^s_{l,1},\ldots, \mathscr{C}^s_{l,m}$ and right-$\mathscr{C}^s_{r,1},\ldots, \mathscr{C}^s_{r,m'}$
filtered $A_{\infty}$ multi-module $\mathscr D^s$
in the same way as Lemma~\ref{lem118ss}.

If $\mathscr{C}_{l,1},\ldots, \mathscr{C}_{l,m}$
and $\mathscr{C}_{r,1},\ldots, \mathscr{C}_{r,m'}$ are strict then
there exists a bijection between the set of all the
left-$\mathscr{C}_{l,1},\ldots, \mathscr{C}_{l,m}$ and right-$\mathscr{C}_{r,1},\ldots, \mathscr{C}_{r,m'}$
filtered $A_{\infty}$ multi-modules $\mathscr D$
and the set of all the filtered $A_{\infty}$ multi-functors
$
\mathscr F\colon
\mathscr{C}^{\rm op}_{l,1} \times \dots \times \mathscr{C}^{\rm op}_{l,m}
\times \mathscr{C}_{r,1} \times \dots \times \mathscr{C}_{r,m'}
\to \mathcal{CH}$.
Moreover, the set of multi-module homomorphisms can be identified with the
set of natural transformations in the category defined in Lemma~\ref{lem256}.

Let
$\mathscr{C}_{l,1},\ldots,\mathscr{C}_{l,m}$,
$\mathscr{C}_{r,1},\ldots,\mathscr{C}_{r,m_r}$
and $\mathscr{C}'_{1,l},\ldots,\mathscr{C}'_{1,m'}$ be
non-unial curved filtered $A_{\infty}$ categories
and
$
\mathscr F
\colon \mathscr{C}_1 \times \dots\times \mathscr{C}_m \to \mathscr{C}'_k$
be a multi-functor.
Let $\mathscr D$ be a
left $\mathscr{C}'_{1,l},\ldots ,\mathscr{C}'_{1,m'}$
and right~$\mathscr{C}_{r,1},\allowbreak\ldots, \mathscr{C}_{r,m_r}$
multi-module.
Then we can pull back $\mathscr D$ by
$\mathscr F$ and obtain
a left $\mathscr{C}'_{1,l},\dots, \mathscr{C}'_{1,k-1}$,
$\mathscr{C}_{l,1},\ldots,\allowbreak \mathscr{C}_{l,m} \mathscr{C}'_{1,k+1}, \ldots,\mathscr{C}'_{1,m'}$
and right $\mathscr{C}_{r,1},\ldots, \mathscr{C}_{r,m_r}$
bi-module,
which we denote $\mathscr F^*\mathscr D$.
We can perform a similar construction for $A_{\infty}$ categories which
act from right.
This construction commutes with the process to associate~$\mathscr D^s$ to~$\mathscr D$.

In the strict case, the above construction coincides
with the composition of multi-functors
via the identification between
a multi-functor to $\mathcal{CH}$ and
a multi-module.

\subsection[A geometric realization of an $A_\infty$ tri-module 1]{A geometric realization of an $\boldsymbol{A_{\infty}}$ tri-module 1}
\label{subsec:bi-functorgeo1}

\begin{situ}\label{situ14}
Let $(X_1,\omega_1)$, $(X_2,\omega_2)$ be
symplectic manifolds and $V_i$ an oriented real vector bundle on
the 3-skeleton $(X_i)_{[3]}$ of $X_i$, for $i=1,2$.
(Namely, $V_i$ is a background datum in the sense of Definition~\ref{relspin}.)

We consider a clean collection $\mathbb L_1$ (resp.\ $\mathbb L_2$)
of $V_1$ (resp.\ $V_2$)
relatively spin oriented and immersed Lagrangian submanifolds
of $X_i$.
(See Situation \ref{situ320}.)
We also take a finite set $\mathbb L_{12}$ of
$\pi_1^*(V_1\oplus TX_1) \oplus \pi_2^*V_2
$ relatively spin oriented and immersed Lagrangian submanifolds
of $-X_1 \times X_2$ that have clean intersection.
We also assume that $L_1 \times L_2$ has clean intersection with $L_{12}$
when $L_i \in \mathbb L_i$, $L_{12} \in \mathbb L_{12}$.
\end{situ}

In this subsection and the next, we will prove the following theorem.

\begin{thm}\label{trimain}
There exists a left-$\mathfrak{Fuk}(X_1,V_1,\mathbb L_1)$,
$\mathfrak{Fuk}(-X_1 \times X_2,\pi_1^*(V_1
\oplus TX_1)\oplus \pi_2^*V_2,\mathbb L_{12})
$
and right-$
\mathfrak{Fuk}(X_2,V_2,\mathbb L_2)
$
filtered $A_{\infty}$ tri-module
$\mathscr{CF}(\mathbb L_1,\mathbb L_{12};\mathbb L_2)$.\index[syindex]{CscrFL1L12@$\mathscr{CF}(\mathbb L_1,\mathbb L_{12};\mathbb L_2)$}
It is unital and gapped.

We call it the {\it correspondence tri-module}.
\index{correspondence tri-module}
\end{thm}
\begin{rem}
We consider associated tri-module\footnote{An
$A_{\infty}$ tri-module is a special case of
a multi-module (see Definition~\ref{bimodulecat2}), that is a multi-module over
three~$A_{\infty}$ categories.}
(see Lemma~\ref{lem55revrev}) over strict categories $\mathfrak{Fukst}(X_1,V_1,\mathbb L_1)$,
$\mathfrak{Fukst}(-X_1 \times X_2,\pi_1^*(V_1
\oplus TX_1)\oplus \pi_2^*V_2,\mathbb L_{12})$,
$\mathfrak{Fukst}(X_1,V_1,\mathbb L_1)$.
Then for objects $(L_1,b_1)$, $(L_{12},b_{12})$, $(L_{2},b_{2})$
of those categories, the tri-module of Theorem~\ref{trimain}
induces a chain complex
$
CF((L_1,b_1),(L_{12},b_{12});(L_{2},b_{2}))$.
Its cohomology is isomorphic to the Floer cohomology of
$HF((L_{12},b_{12});(L_1 \times L_2,b_1\times b_2))$.
This fact will be proved in Section~\ref{sec:kuneth} (see Theorem~\ref{thm164}).
The product $b_1\times b_2$ of bounding cochains is defined in Proposition~\ref{lem167}.

\end{rem}
\begin{proof}
The proof of Theorem~\ref{trimain} occupies this and the next subsections.
In the same way as Section~\ref{subsec:Opposite}, it suffices to consider the case when
$\mathbb L_1$, $\mathbb L_2$, $\mathbb L_{12}$ consist of single immersed Lagrangian
submanifolds $L_1$, $L_2$ and $L_{12}$, respectively, and construct
structure operations
\[
\mathfrak n\colon\ BCF[1](L_1) \otimes_{\Lambda_0} BCF[1](L_{12}) \otimes D[1]
\otimes_{\Lambda_0} BCF[1](L_{2})
\to D[1]
\]
for a certain graded completed free $\Lambda_0$ module $D$,
such that they satisfy $A_{\infty}$ relation.\footnote{Here we shifted
the degree of elements of $D$. This is because it is more consistent
with the discussion of sign in Section~\ref{sec:orient}.}

The construction of $\mathfrak n$ uses certain compactified moduli spaces
${\mathcal M}_{\rm QT}(\vec a_1,\vec a_{12},\vec a_2;a_-,a_+;E)$
of pseudo-holomorphic quilts,
which will be defined in Definition~\ref{defn528}.
We will define it in three steps.

We first define \smash{$\overset{\ \text{\tiny $\circ\circ$}}{\mathcal M}_{\rm QT}(\vec a_1,\vec a_{12},\vec a_2;a_-,a_+;E)$}
in Definition~\ref{def516}.
This moduli space is the set of pseudo-holomorphic quilts which do
not contain disk bubbles
and are not split into several quilts (see Figure~\ref{Figure51}). It contains objects with sphere bubbles.

We then include objects with disk bubbles and define
\smash{$\mathring{\mathcal M}_{\rm QT}(\vec a_1,\vec a_{12},\vec a_2;a_-,a_+;E)$}
in Definition~\ref{defn523000}.

Finally, we include the process where a sequence of a pseudo-holomorphic quilts
splits into several pseudo-holomorphic quilts in the limit
and define ${\mathcal M}_{\rm QT}(\vec a_1,\vec a_{12},\vec a_2;a_-,a_+;E)$ in Definition~\ref{defn528}.
The detail will follow.

We decompose the fiber products into connected components
\[
L_i(+) = \tilde L_i \times_{X_i} \tilde L_i
= \tilde L_i \sqcup \coprod_{a \in \mathcal A_{L_i}} L_i(a)
\]
for $i=1,2$, and
\[
L_{12}(+) = \tilde L_{12} \times_{X_1 \times X_2} \tilde L_{12}
= \tilde L_{12} \sqcup \coprod_{a \in \mathcal A_{L_{12}}} L_{12}(a).
\]
See Definition~\ref{def3131}\,(5).
We also decompose
\[
R =
\bigl(\tilde L_1 \times \tilde L_2\bigr) \times_{X_1 \times X_2} \tilde L_{12}
=
\coprod_{a \in \mathcal A_{R}} R(a).
\]
Let \smash{$\vec a_j = (a_{j,1},\dots,a_{j,k_j}) \in \bigl(\mathcal A_{L_{j}}^+\bigr)^{k_j}$},
\smash{$\vec a_{12} = (a_{12,1},\dots,a_{12,k}) \in \bigl(\mathcal A_{L_{12}}^+\bigr)^{k}$},
$a_+, a_- \in \mathcal A_{R}$
and $E \in \R_{\ge 0}$.
(Here \smash{$\mathcal A_{L_{12}}^+ := \mathcal A_{L_{12}} \cup \{o\}$} and $o$ is the
index of the diagonal component.)
Below, we identify~${\R \times \R \cong \C}$ by $(s,t) \mapsto s + \sqrt{-1}t$.
\begin{defn}\label{def516}
We consider
$(\Sigma;\vec z_1,\vec z_{12},\vec z_2;u_1,u_2;\gamma_1,\gamma_{12},\gamma_2)$
with the following properties (see Figure~\ref{Figure51}):
\begin{enumerate}\itemsep=0pt
\item[(1)]
The space
$\Sigma$ is a bordered Riemann surface with $\Sigma \supseteq ([-1,1]\times \R)$.
The closure of $\Sigma \setminus ([-1,1]\times \R)$ is a finite
union of (maximal) trees of spheres. We call its connected component a {\it tree of sphere components}.
\index{tree of sphere components}
We require that a tree of sphere components intersects
with
$[-1,1]\times \R$ at one point, which we call its {\it root}. \index{root}
All the roots are
points of~${((-1,0) \cup (0,1)) \times \R}$.\footnote{We require that the root is not on $\{0,\pm 1\} \times \R$.}
\item[(2)]
Let $\Omega_1$ (resp.\ $\Omega_2$) be the union of $[-1,0] \times \R$
(resp.\ $[0,+1] \times \R)$) and the trees of sphere components rooted on it.
We require the maps $u_1 \colon \Omega_1 \to (X_1,J_1)$ and $u_2 \colon \Omega_2 \to (X_2,J_2)$
to be pseudo-holomorphic.
\item[(3)]
We put
$\vec z_i = (z_{i,1},\dots,z_{i,k_i})$, $i=1,2$. Then
$z_{1,j} \in \{-1\} \times \R$, $z_{2,j} \in \{1\} \times \R$.
$\vec z_{12} = (z_{12,1},\dots,z_{12,k})$, $z_{12,j} \in \{0\} \times \R$.
If $j_1 < j_2$ then $\operatorname{Im} z_{1,j_1} > \operatorname{Im} z_{1,j_2}$,
$\operatorname{Im} z_{12,j_1} > \operatorname{Im} z_{12,j_2}$
and $\operatorname{Im} z_{2,j_1} < \operatorname{Im} z_{2,j_2}$.
See Remark~\ref{rem1521} for this enumeration.
We put $\vert \vec z_i\vert = \{z_{i,1},\dots,z_{i,k_i}\}$.
$\vert \vec z_{12}\vert$ is defined in the same way.
\item[(4)]
The maps
$\gamma_1 \colon (\{-1\}\times \R) \setminus \vert\vec z_1\vert \to \tilde L_1$,
$\gamma_2 \colon (\{1\}\times \R) \setminus \vert\vec z_2\vert \to \tilde L_2$,
$\gamma_{12} \colon (\{0\}\times \R) \setminus \vert\vec z_{12}\vert \to \tilde L_{12}$
are smooth and satisfy
\begin{gather*}
i_{L_1}(\gamma_1(z)) = u_1(z) \qquad \text{if}\quad z \in (\{-1\}\times \R) \setminus \vert\vec z_1\vert, \\
i_{L_2}(\gamma_2(z)) = u_2(z) \qquad \text{if}\quad z \in (\{1\}\times \R) \setminus \vert\vec z_2\vert, \\
i_{L_{12}}(\gamma_{12}(z)) = (u_1(z),u_2(z)) \qquad \text{if}\quad z \in
(\{0\}\times \R) \setminus \vert\vec z_{12}\vert.
\end{gather*}
\item[(5)]
At $\vec z_1$, $\vec z_2$, $\vec z_{12}$, the maps $\gamma_1$, $\gamma_2$,
$\gamma_{12}$ satisfy the switching condition, Condition~\ref{cond517} below.
\item[(6)]
When $z \in [-1,1]\times \R$, $\operatorname{Im}z \to \pm \infty$,
the maps $u_1(z)$ and $u_2(z)$
satisfy the asymptotic boundary condition, Condition~\ref{cond518} below.
\item[(7)]
The stability condition, Condition \ref{cond519} below, is satisfied.
\item[(8)]
$
\int_{\Omega_1}u_1^*\omega_1 + \int_{\Omega_2}u_1^*\omega_2 = E$.
\end{enumerate}

We will define an equivalence relation $\sim$ among the objects
$(\Sigma;\vec z_1,\vec z_{12},\vec z;u_1,u_2;\gamma_1,\gamma_{12},\gamma_2)$
satisfying (1)--(8), in Definition~\ref{defnevalama}. We denote by
\smash{$\overset{\ \text{\tiny $\circ\circ$}}{\mathcal M}_{\rm QT}(\vec a_1,\vec a_{12},\vec a_2;a_-,a_+;E)$}
\index[syindex]{M1QTa1a12@$\overset{\ \text{\tiny $\circ\circ$}}{\mathcal M}_{\rm QT}(\vec a_1,\vec a_{12},\vec a_2;a_-,a_+;E)$}
the set of all the equivalence classes of this equivalence relation.
We call its element a {\it pseudo-holomorphic quilt}.\index{pseudo-holomorphic quilt}

\end{defn}
\begin{figure}[ht]
\centering
\includegraphics[scale=0.4]{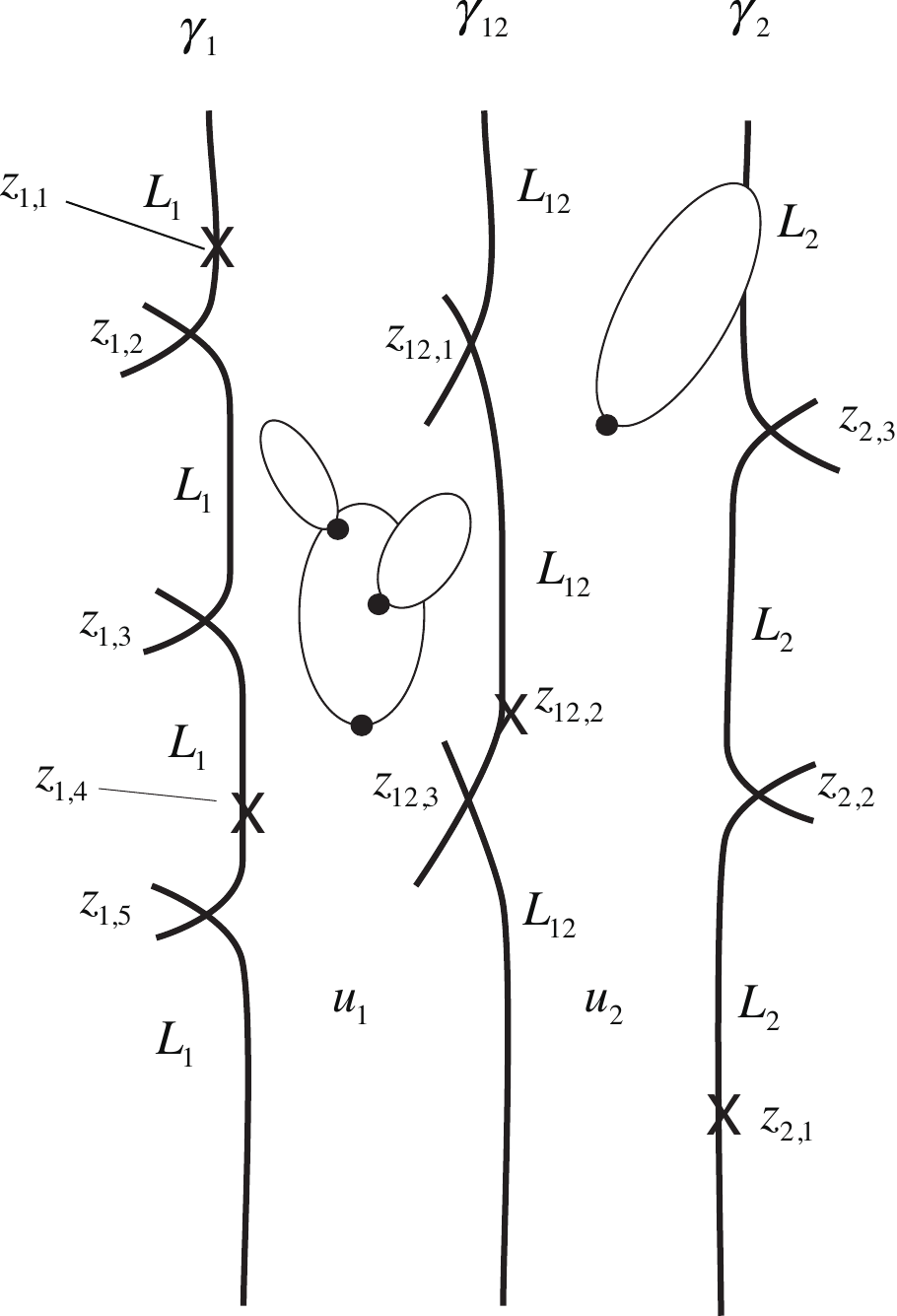}
\caption{An element $\overset{\ \text{\tiny $\circ\circ$}}{\mathcal M}_{\rm QT}(\vec a_1,\vec a_{12},\vec a_2;a_-,a_+;E)$.}
\label{Figure51}
\end{figure}
We next describe three of the conditions in Definition~\ref{def516}.
We put
\[
\partial_1\Sigma = \{-1\} \times \R,
\qquad
\partial_2\Sigma = \{1\} \times \R,
\qquad
\partial_{12} \Sigma = \{0\} \times \R.
\]
We call the line $\{0\} \times \R$ the {\it seam}. \index{seam}
We define the limit
$
p = \lim_{z \in \partial_1\Sigma, z \downarrow z_{1,j}} \gamma_1(z)
$
as follows. If $z_n, z \in (\{-1\}\times \R) \setminus \vert\vec z_1\vert$,
$\operatorname{Im}z_n > \operatorname{Im} z_{1,j}$ and $\lim_{n\to \infty} z_n = z_{1,j}$,
then
$
p = \lim_{n\to \infty} \gamma_1(z_n)
$.
The notations $\lim_{z \in \partial_1\Sigma, z \uparrow z_{1,j}}$
etc.\ are defined in the same way.

The two switching conditions below are
analogues of the switching conditions which
appeared in the immersed Lagrangian Floer theory.
See Definition~\ref{def3737}\,(5) and Figure~\ref{Figure35}.

\begin{conds}[switching condition 1]\label{cond517}
\quad
\index{switching condition}
\begin{enumerate}\itemsep=0pt
\item[(1)]
For each $j$,
\smash{$
(\lim_{z \in \partial_1\Sigma, z \downarrow z_{1,j}} \gamma_1(z),
\lim_{z \in \partial_1\Sigma, z \uparrow z_{1,j}} \gamma_1(z))
\in L_1(a_{1,j})$}.
\item[(2)]
For each $j$,
\smash{$
(\lim_{z \in \partial_2\Sigma, z \uparrow z_{2,j}} \gamma_2(z),
\lim_{z \in \partial_2\Sigma, z \downarrow z_{2,j}} \gamma_2(z))
\in L_2(a_{2,j})$}.
\item[(3)]
For each $j$,
\smash{$
(\lim_{z \in \partial_{12}\Sigma, z \downarrow z_{12,j}} \gamma_{12}(z),
\lim_{z \in \partial_{12}\Sigma, z \uparrow z_{12,j}} \gamma_{12}(z))
\in L_{12}(a_{12,j})$}.
\end{enumerate}

\end{conds}
\begin{conds}[switching condition 2]\label{cond518}
\quad
\begin{enumerate}\itemsep=0pt
\item[(1)]
There exists $(p_{+\infty,1},p_{+\infty,2})
\in R(a_+)$ such that
\[
\lim_{\tau \to + \infty} \bigl(\gamma_1\bigl(-1+\tau\sqrt{-1}\bigr),\gamma_2\bigl(+1+\tau\sqrt{-1}\bigr)
\bigr)= p_{+\infty,1},
\qquad \lim_{\tau \to + \infty}\gamma_{12}\bigl(\tau\sqrt{-1}\bigr) = p_{+\infty,2}.
\]
\item[(2)]
There exists $(p_{-\infty,1},p_{-\infty,2})
\in R(a_-)$ such that
\[
\lim_{\tau \to - \infty} \bigl(\gamma_1\bigl(-1-\tau\sqrt{-1}\bigr),\gamma_2\bigl(+1+\tau\sqrt{-1}\bigr)\bigr) = p_{-\infty,1},
\qquad \lim_{\tau \to - \infty}\gamma_{12}\bigl(\tau\sqrt{-1}\bigr) = p_{-\infty,2}.
\]
\end{enumerate}
See Figure~\ref{Figure53}.
\end{conds}
\begin{figure}[ht]
\centering
\includegraphics[scale=0.35]{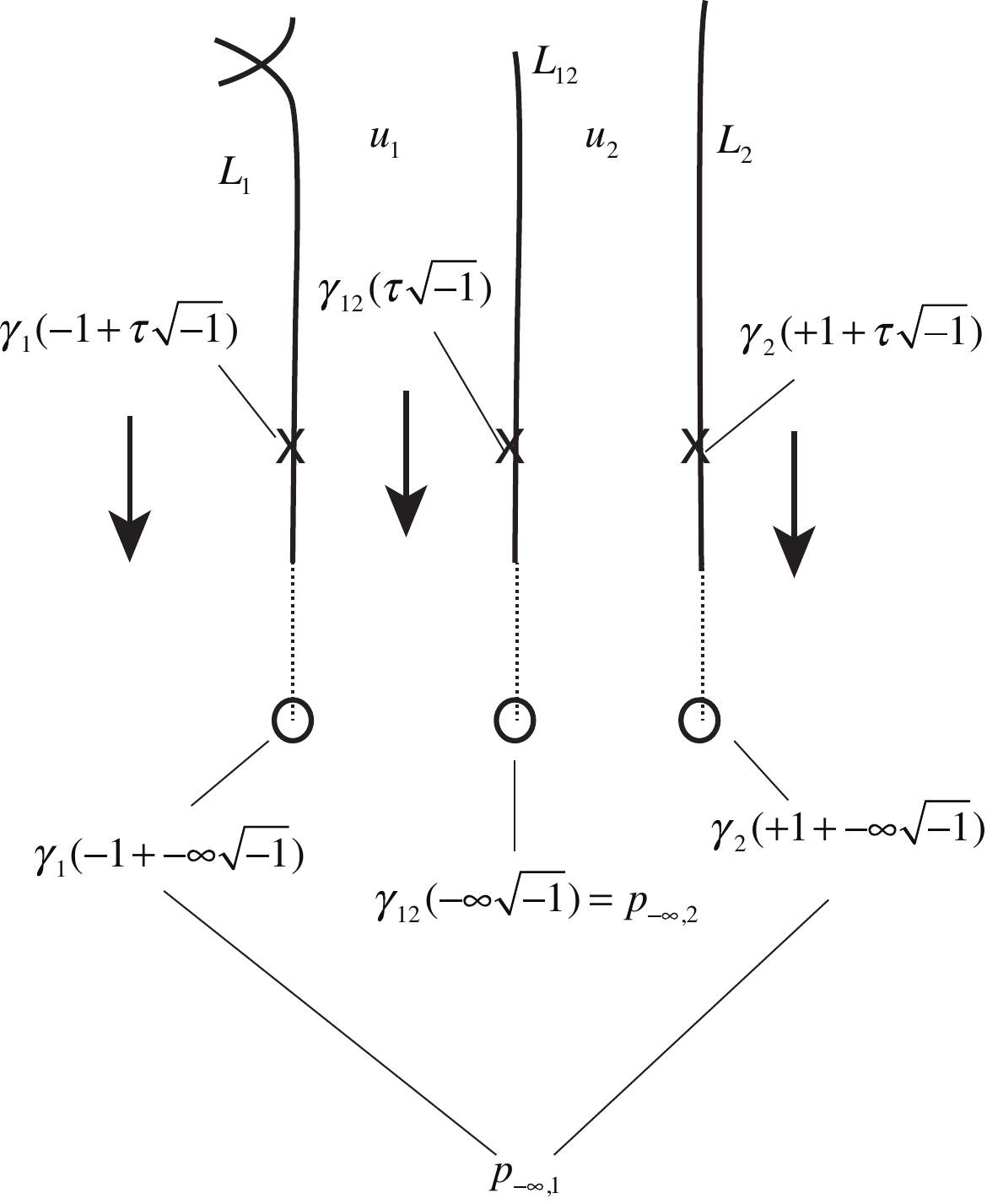}
\caption{Switching condition 2.}
\label{Figure53}
\end{figure}
\begin{rem}\label{rem1521}
Note that we enumerate the marked points on $\gamma_1$, $\gamma_{12}$
downward and the marked points on $\gamma_2$ upward.
This is related to the fact that we are constructing {\it left}
$\mathfrak{Fuk}(X_1,V_1,\mathbb L_1)$,
$\mathfrak{Fuk}(-X_1 \times X_2,\pi_1^*(V_1
\oplus TX_1)\oplus \pi_2^*V_2,\mathbb L_{12})$
and {\it right}
$\mathfrak{Fuk}(X_2,V_2,\mathbb L_2)$
filtered $A_{\infty}$ tri-module.
(We also remark the input $D$ corresponds
to the end $\tau \to -\infty$.)

In fact, we write the structure operation of this
filtered $A_{\infty}$ tri-module as
\[
\mathfrak n(x_1,\dots,x_{k_1};y_1,\dots,y_{k_{12}};z;w_1,\dots,w_{k_2}).
\]
Here $x_i$ corresponds to the evaluation map at the $i$-th marked point
of $\gamma_1$, $y_i$ corresponds to the evaluation map at the $i$-th marked point
of $\gamma_{12}$, $w_i$ corresponds to the evaluation map at the $i$-th marked point
of $\gamma_{2}$.
Thus the way we enumerate the marked points
is consistent with the way we write the
structure operation.
\end{rem}

\begin{conds}[stability condition]\label{cond519}
The set of all the maps $v\colon \Sigma \to \Sigma$ satisfying the
next conditions is finite.
\begin{enumerate}\itemsep=0pt
\item[(1)]
The map $v$ is a homeomorphism and is biholomorphic on each of the irreducible components.
\item[(2)]
$u_1 \circ v = u_1$, $u_2 \circ v = u_2$.
\end{enumerate}

\end{conds}
\begin{defn}\label{defnevalama}
We define {\it evaluation maps}
\index{evaluation map}
\begin{gather*}
{\rm ev} = \bigl({\rm ev}^1,{\rm ev}^{12},{\rm ev}^2\bigr) =
\bigl(\bigl({\rm ev}^1_1,\dots,{\rm ev}^1_{k_1}\bigr),\bigl({\rm ev}^{12}_1,\dots,{\rm ev}^{12}_{k_{12}}\bigr),
\bigl({\rm ev}^2_1,\dots,{\rm ev}^2_{k_2}\bigr)\bigr)\colon \\
\qquad
\overset{\ \text{\tiny $\circ\circ$}}{\mathcal M}_{\rm QT}(\vec a_1,\vec a_{12},\vec a_2;a_-,a_+;E)
\to \prod_{j=1}^{k_1} L_1(a_{1,j}) \times \prod_{j=1}^{k_{12}} L_{12}(a_{12,j})
\times \prod_{j=1}^{k_2} L_2(a_{2,j})
\end{gather*}
and
\[
{\rm ev}_{\infty} = ({\rm ev}_{\infty,-},{\rm ev}_{\infty,+})\colon \
\overset{\ \text{\tiny $\circ\circ$}}{\mathcal M}_{\rm QT}(\vec a_1,\vec a_{12},\vec a_2;a_-,a_+;E)
\to R(a_-) \times R(a_+)
\]
as follows.
\begin{enumerate}\itemsep=0pt
\item[(1)]
We use Condition \ref{cond517}\,(1) to define
\index[syindex]{ev1j@${\rm ev}^1_j$}
\[
{\rm ev}^1_j(\Sigma;\vec z_1,\vec z_{12},\vec z_2;u_1,u_2;\gamma_1,\gamma_{12},\gamma_2) =
\bigl(\lim_{z \in \partial_1\Sigma, z \downarrow z_{1,j}} \gamma_1(z),
\lim_{z \in \partial_1\Sigma, z \uparrow z_{1,j}} \gamma_1(z)\bigr)
\in L_1(a_{1,j}).
\]
\item[(2)]
We use Condition \ref{cond517}\,(3) to define
\index[syindex]{ev12j@${\rm ev}^{12}_j$}
\begin{gather*}
{\rm ev}^{12}_j(\Sigma;\vec z_1,\vec z_{12},\vec z_2;u_1,u_2;\gamma_1,\gamma_{12},\gamma_2)\\
\qquad
=
\bigl(\lim_{z \in \partial_{12}\Sigma, z \downarrow z_{12,j}} \gamma_{12}(z),
\lim_{z \in \partial_{12}\Sigma, z \uparrow z_{12,j}} \gamma_{12}(z)\bigr)
\in L_{12}(a_{12,j}).
\end{gather*}
\item[(3)]
The evaluation map ${\rm ev}^2_j$ is defined in the same way by using Condition \ref{cond517}\,(2).
\item[(4)]
We use Condition \ref{cond518}\,(1) to define
\index[syindex]{evinfty@${\rm ev}_{\infty,+}$}
\begin{gather*}
{\rm ev}_{\infty,+}
(\Sigma;\vec z_1,\vec z_{12},\vec z_2;u_1,u_2;\gamma_1,\gamma_{12},\gamma_2)\\
\qquad=
\lim_{\tau \to + \infty} \bigl(\bigl(\gamma_1\bigl(-1+\tau\sqrt{-1}\bigr),\gamma_2\bigl(+1+\tau\sqrt{-1}\bigr)
\bigl),\gamma_{12}\bigl(\tau\sqrt{-1}\bigr)\bigr).
\end{gather*}
The definition of
${\rm ev}_{\infty,-}$ is similar.
We call them {\it evaluation maps at infinity}.\index{evaluation map at infinity}
\end{enumerate}
\end{defn}
\begin{defn}\label{def521}
We say
$(\Sigma;\vec z_1,\vec z_{12},\vec z_2;u_1,u_2;\gamma_1,\gamma_{12},\gamma_2)$
as in Definition~\ref{def516}
is {\it equivalent} to~$(\Sigma';\vec z^{\,\prime}_1,\vec z^{\,\prime}_{12},\vec z^{\,\prime}_2;u'_1,u'_2;\gamma'_1,\gamma'_{12},\gamma'_2)$
if there exist $v \colon \Sigma \to \Sigma'$ satisfying the
next conditions.
\begin{enumerate}\itemsep=0pt
\item[(1)] The map $v$ is a homeomorphism and is biholomorphic on each connected component.
\item[(2)]
We require
$v(\Omega_1) = \Omega'_1$, $v(\Omega_2) = \Omega'_2$.
Here $\Omega'_1$ is the union of $[-1,0]\times \R \subset \Sigma'$
and the trees of sphere components rooted on it.
$\Omega'_2$ is defined in the same way.
\item[(3)]
$u'_1 \circ v = u_1$, $u'_2 \circ v = u_2$.
\item[(4)]
$v(z_{i,j}) = z'_{i,j}$, $v(z_{12,j}) = z'_{12,j}$, where $i=1,2$.
\item[(5)]
$\gamma'_1\circ v = \gamma_1$, $\gamma'_2\circ v = \gamma_2$,
$\gamma'_{12}\circ v = \gamma_{12}$.
\end{enumerate}
\end{defn}
\begin{rem}
In Floer theory, the moduli space which is used to define the boundary operator
is the quotient space by $\R$ action. (This $\R$ action is
induced by the translation of the $\R$ direction,
which is the second factor of $[-1,1]\times \R$ in our situation.)
The process to take the set of equivalence classes of the
equivalence relation in Definition~\ref{def521} includes the
process to take the quotient by this $\R$ action.
In other words, the object $\mathfrak x = (\Sigma;\vec z_1,\vec z_{12},\vec z_2;u_1,u_2;\gamma_1,\gamma_{12},\gamma_2)$
and $\tau\mathfrak x$ which is obtained from $\mathfrak x$ by shifting everything by $\tau \in \R$
are equivalent in the sense of Definition~\ref{def521}.

There is no mathematical difference between the way we take here and the usual
way to take quotient by $\R$ action. They are slightly different ways
to describe the same mathematical contents.
\end{rem}

In our situation, Condition \ref{cond518} is a consequence of the other
conditions. More precisely, we have the following.
\begin{lem}\label{lem526526}
Let $(\Sigma;\vec z_1,\vec z_{12},\vec z_2;u_1,u_2;\gamma_1,\gamma_2,\gamma_{12})$
be an object which satisfies conditions of Definition
{\rm\ref{def516}} except possibly $(6)$, for some
$\vec a_1$, $\vec a_2$, $\vec a$, $E$.
$($Note that $a_{\pm}$ appears only in $(6)$.$)$
Then there exists $a_{-}$, $a_+$ such that $(6)=$Condition {\rm\ref{cond518}} is satisfied.

Moreover, there exists $C_k, c_k > 0$ such that
\[
\bigl\Vert \nabla^k u_1 (z) \bigr\Vert \le C_ke^{-c_k \vert\operatorname{Im} z\vert},
\qquad
\bigl\Vert \nabla^k u_1 (z) \bigr\Vert \le C_ke^{-c_k \vert\operatorname{Im} z\vert}.
\]

\end{lem}
\begin{proof}
We use $(t,\tau)$ as a coordinate of $[0,1] \times [\tau_0,\infty)$
and the point
$(t,\tau) \in [0,1] \times [\tau_0,\infty)$
is identified with $z = t + \sqrt{-1}\tau \in \C$.

We may assume that there is no tree of sphere components
whose root is a point $z$ with~${\tau > \tau_0}$.
We may also assume that $\operatorname{Im}z_{i,j}, \operatorname{Im}z_{12,j} < -\tau_0$.
We define
$u \colon [0,1] \times [\tau_0,\infty) \to (X_1,-J_1) \times (X_2,J_2)$ by
$
u(z) = (u_1(\overline z),u_2(z))$.

The map $u$ is pseudo-holomorphic and
$
u(\{+1\} \times [\tau_0,\infty) ) \subset L_1 \times L_2$,
$
u(\{0\} \times [\tau_0,\infty) ) \subset L_{12}$.
Moreover,
\[
\int_{[0,1] \times [\tau_0,\infty)} u^*(-\pi_1^*(\omega_1)
+ \pi_2^*(\omega_2))
< \infty.
\]
Since $L_{12}$ and $L_1\times L_2$ have clean intersection
(see Situation \ref{situ14}),
there exists an element $p_{+\infty}
= (p_{+\infty,1},p_{+\infty,2}) \in L_{12} \cap (L_1\times L_2)$
such that
\[
d(u(z),p_{+\infty}) < C e^{-C \vert\operatorname{Im}z\vert},
\qquad
\big\Vert \nabla^k u(z) \big\Vert \le C_k^{-c_k \vert\operatorname{Im} z\vert}
\]
on $[0,1] \times (\tau_0,\infty)$.
(See \cite[Lemmas 2.4 and 2.5]{foooanalysis} for example.)
We can discuss in the same way for $\tau < -\tau_0$.
\end{proof}

We will next discuss the compactification of
\smash{$\overset{\ \text{\tiny $\circ\circ$}}{\mathcal M}_{\rm QT}(\vec a_1,\vec a_{12},\vec a_2;a_-,a_+;E)$}.
Note that we already included objects with sphere bubbles in
\smash{$\overset{\ \text{\tiny $\circ\circ$}}{\mathcal M}_{\rm QT}(\vec a_1,\vec a_{12},\vec a_2;a_-,a_+;E)$}.
We need to include disk bubbles and the process
where elements split into several pieces in the second factor of
$[-1,1]\times \R$.
Note that disk bubbles may occur at the boundaries
$\partial_{1}\Omega$, $\partial_2\Omega$ or the seam
$\partial_{12}\Omega$, where
pseudo-holomorphic disks in $X_1$, $X_2$, $-X_1 \times X_2$
with boundary in $L_1$, $L_2$, $L_{12}$ can bubble,
respectively.
The moduli spaces of such pseudo-holomorphic disks
are described by the moduli spaces we introduced
in Section~\ref{subsec:modpolygon} and hence
the moduli space of objects with disk bubbles is obtained by
an appropriate fiber product. We will describe it below.
\begin{defn}
Let ${\mathcal M}(L;\vec a;E)$ be the moduli space
introduced in \eqref{def33314}.
For the sake of simplicity of notations, we use the next
(slight abuse of) notations.
Let $\vec a = (a,a)$ ($a \in \mathcal A_L$).
We include ${\mathcal M}(L; \vec a;0) = {\mathcal M}(L; (a,a);0)$ and define it to be a
single point consisting of a constant map to $L(a)$.
In fact, this element is unstable. However, we include it as an
exception here. See Remark~\ref{rem5244}.

Let
\[
\vec a^{i,(j)} = \bigl(a^{i,(j)}_0,\dots,a^{i,(j)}_{m_{i,(j)}}\bigr)
\in \bigl(\mathcal A_{L_i}^+\bigr)^{m_{i,(j)}+1}\qquad \text{for $i=1,2$, $j=1,\dots,k_i$}
\]
and
\[
\vec a^{12,(j)} = \bigl(a^{12,(j)}_0, \dots,a^{12,(j)}_{m_{12,(j)}}\bigr)
\in \bigl(\mathcal A_{L_i}^+\bigr)^{m_{12,(j)}+1}\qquad \text{for $j=1,\dots,k_{12}$}.
\]
Here $m_{i,(j)}$ and $m_{12,(j)}$ are nonnegative integers.

We put
\begin{gather*}
\vec a'_i = (a'_{i,1},\dots,a'_{i,k_i}) = \bigl(a^{i,(1)}_0,\dots,a^{i,(k_i)}_0\bigr) \in \bigl(\mathcal A^+_{L_i}\bigr)^{k_i}, \qquad i=1,2, \\
\vec a'_{12} = (a'_{12,1},\dots,a'_{12,k_{12}})= \bigl(a^{12,(1)}_0,\dots,a^{12,(k_{12})}_0\bigr) \in \bigl(\mathcal A_{L_{12}}^+\bigr)^{k_{12}}.
\end{gather*}
We define
\[
\#_j\vec a^{i,(j)} := \bigl(a^{i,(1)}_1,\dots,a^{i,(j)}_{m_{i,(1)}},
a^{i,(2)}_1,\dots,a^{i,(j)}_{m_{i,(2)}},\dots,a^{i,(k_i)}_{1},
\dots,a^{i,(k_i)}_{m_{i,(k_i)}}\bigr)
\in \bigl(\mathcal A_{L_i}^+\bigr)^{m_{i}},
\]
where $i=1,2$ and \smash{$m_i = \sum_j m_{i,(j)}$}.
We moreover put
\[
\#_j\vec a^{12,(j)} := \bigl(a^{12,(1)}_1,\dots,a^{12,(j)}_{m_{12,(1)}},
\dots,a^{12,(k_{12})}_{1},\dots,a^{12,(k_{12})}_{m_{12,(k_{12})}}\bigr)
\in \bigl(\mathcal A_{L_{12}}^+\bigr)^{m_{12}},
\]
where $m_{12}= \sum_j m_{12,(j)}$.
(See Figure~\ref{Figure54}.)

\end{defn}
\begin{figure}[ht]
\centering
\includegraphics[scale=0.4]{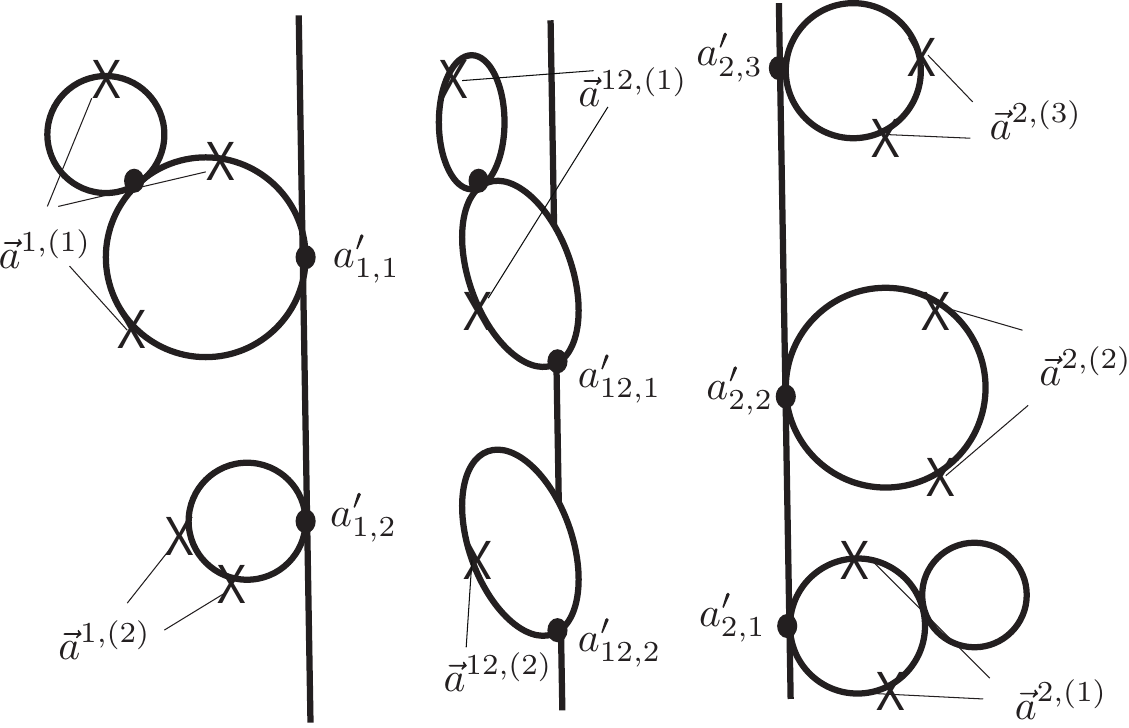}
\caption{Domain of an element of $\mathring{\mathcal M}_{\rm QT}(\vec a_1,\vec a_{12},\vec a_2;a_-,a_+;E)$.}
\label{Figure54}
\end{figure}
\begin{defn}\label{defn523000}
We define the set
\smash{$\mathring{\mathcal M}_{\rm QT}(\vec a_1,\vec a_{12},\vec a_2;a_-,a_+;E)$}
as the union of the fiber products
\begin{align}
\overset{\ \text{\tiny $\circ\circ$}}{\mathcal M}_{\rm QT}(\vec a'_1,\vec a'_{12},\vec a'_2;a_-,a_+;E')
&\times_{{\rm ev}_0,\dots,{\rm ev}_0}
\prod_{j=1}^{k_1}{\mathcal M}\bigl(L_1; \vec a^{1,(j)};E_{1,j}\bigr)\nonumber \\
&\times_{{\rm ev}_0,\dots,{\rm ev}_0}
\prod_{j=1}^{k_{12}}{\mathcal M}'\bigl(L_{12}; \vec a^{12,(j)};E_{12,j}\bigr)\nonumber\\
&\times_{{\rm ev}_0,\dots,{\rm ev}_0}
\prod_{j=1}^{k_2}{\mathcal M}\bigl(L_2; \vec a^{2,(j)};E_{2,j}\bigr),\label{eq518}
\end{align}
where \smash{$\#_j\vec a^{1,(j)} = \vec a_1$}, \smash{$\#_j\vec a^{12,(j)} = \vec a_{12}$},
\smash{$\#_j\vec a^{2,(j)} = \vec a_2$},
$E' + \sum_j E_{1,j} + \sum_j E_{12,j} + \sum_j E_{2,j} = E$.

\end{defn}
We remark that in the first line of \eqref{eq518}
the fiber product is taken over
\smash{$
\prod_{j=1}^{k_1} L_1\bigl(a^{1,(j)}_0\bigr)
$}
by the evaluation maps
\begin{gather*}
{\rm ev}^1_j \colon \ \overset{\ \text{\tiny $\circ\circ$}}{\mathcal M}_{\rm QT}(\vec a'_1,\vec a'_{12},\vec a'_2;a_-,a_+;E')
\to L_1\bigl(a^{1,(j)}_0\bigr), \\
{\rm ev}_0 \colon \ {\mathcal M}\bigl(L_{1}; \vec a^{1,(j)};E_{1,j}\bigr)
\to L_{1}\bigl(a^{1,(j)}_0\bigr).
\end{gather*}
The fiber product in the second line is taken over
\smash{$\prod_{j=1}^{k_{12}} L_{12}\bigl(a^{12,(j)}_0\bigr)
$}
by the evaluation maps
\begin{gather*}
{\rm ev}^{12}_j \colon \ \overset{\ \text{\tiny $\circ\circ$}}{\mathcal M}(\vec a'_1,\vec a'_{12},\vec a'_2;a_-,a_+;E')
\to L_{12}\bigl(a^{12,(j)}_0\bigr),\\
{\rm ev}_0 \colon \ {\mathcal M}'\bigl(L_{12}; \vec a^{12,(j)};E_{12,j}\bigr)
\to L_{12}\bigl(a^{12,(j)}_0\bigr).
\end{gather*}
The fiber product in the third line is taken in a similar way.

\begin{rem}\label{Remark524}
In the formula \eqref{eq518}, we used a compactification ${\mathcal M}'\bigl(L_{12}; \vec a^{12,(j)};E_{12,j}\bigr)$
\index[syindex]{M1primeL12a12@${\mathcal M}'(L_{12}; \vec a^{12,(j)};E_{12,j})$}
of the space \smash{$\overset{\ \text{\tiny $\circ\!\circ\!\circ$}}{\mathcal M}\bigl(L_{12}; \vec a^{12,(j)};E_{12,j}\bigr)$}.
Here \smash{$\overset{\ \text{\tiny $\circ\!\circ\!\circ$}}{\mathcal M}\bigl(L_{12}; \vec a^{12,(j)};E_{12,j}\bigr)$} is the moduli space of pseudo-holomorphic disks
whose source curve is $D^2$ without any disk or sphere bubbles.
This compactification is similar to the stable map compactification
${\mathcal M}\bigl(L_{12}; \vec a^{12,(j)};E_{12,j}\bigr)$
which we defined in Section~\ref{subsec:modpolygon}
but is slightly different from it.
It is necessary to use different
compactification for
our space ${\mathcal M}_{\rm QT}(\vec a_1,\vec a_{12},\vec a_2;a_-,a_+;E)$
to carry a Kuranishi structure.
We will explain
this point in detail in Section~\ref{sec:directcomp}.

\end{rem}

\begin{rem}\label{rem5244}
As we mentioned before, we include the case when a factor
${\mathcal M}\bigl(L_1; \vec a^{1,(j)};E_{1,j}\bigr)$ is ${\mathcal M}(L_1;(a,a);0)$.
This moduli space consists of one point and is a constant map to a point
in~$L_1(a)$. Note that this element actually is not a stable map since its
automorphism group is $\R$. This case corresponds to
the case when the corresponding marked point is on the
line~${\{-1\} \times \R}$ (and not on the disk bubble) and
is mapped to an element of $L_1(a)$.
We include this case in~\eqref{eq518} and etc.\
for the sake of simplicity of notation.
When we regard this element as a~`stable map' we shrink this
disk and regard the `root' as a marked point. (See Figure~\ref{Figure56}.)
We consider the case when ${\mathcal M}'(L_{12};(a,a);0)$ (resp.
${\mathcal M}(L_2;(a,a);0)$) appears in the
second (resp.\ third) line of~\eqref{eq518} in the same way.
\begin{figure}[ht]
\centering
\includegraphics[scale=0.4]{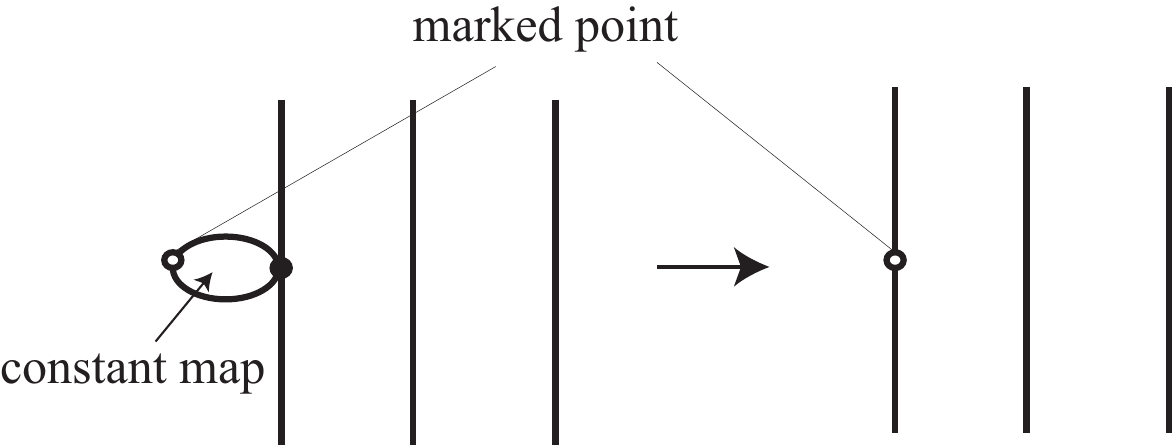}
\caption{Shrink an element of ${\mathcal M}(L_2;(a,a);0)$.}
\label{Figure56}
\end{figure}

\end{rem}
We have thus included the objects with disk bubbles.
We finally define our compactification as follows.

\begin{defn}\label{defn528}
We define the set
${\mathcal M}_{\rm QT}(\vec a_1,\vec a_{12},\vec a_2;a_-,a_+;E)$
\index[syindex]{M1a1a12a2a-@${\mathcal M}(\vec a_1,\vec a_{12},\vec a_2;a_-,a_+;E)$}
as the union of the fiber products
\begin{align}
\mathring{\mathcal M}_{\rm QT}(\vec a_{1,0},\vec a_{12,0},\vec a_{2,0};a_0,a_1;E_1)&\times_{R(a_1)} \mathring{\mathcal M}_{\rm QT}(\vec a_{1,1},\vec a_{12,1},\vec a_{2,1};a_1,a_2;E_2) \times_{R(a_2)} \cdots\nonumber
\\
& \times_{R(a_{\ell-1})} \mathring{\mathcal M}_{\rm QT}(\vec a_{1,\ell},\vec a_{12,\ell},
\vec a_{2,\ell};a_{\ell-1},a_\ell;E_{\ell}).\label{form519}
\end{align}
Here $\vec a_1 = \vec a_{1,0}, \vec a_{1,1},\dots, \vec a_{1,\ell}$,
$\vec a_{12} = \vec a_{12,0}, \vec a_{12,1},\dots, \vec a_{12,\ell}$,
$\vec a_2 = \vec a_{2,\ell}, \vec a_{2,\ell-1},\dots, \vec a_{2,0}$,
$E_1+\dots+E_{\ell} = E$ and
$
a_- = a_0, a_1,\dots,a_{\ell-1}, a_{\ell}= a_{+} \in \mathcal A_{R}$.
We use the maps ${\rm ev}_{\infty} = ({\rm ev}_{\infty,-},{\rm ev}_{\infty,+})$ to
define the fiber product.
\end{defn}
\begin{figure}[ht]
\centering
\includegraphics[scale=0.4]{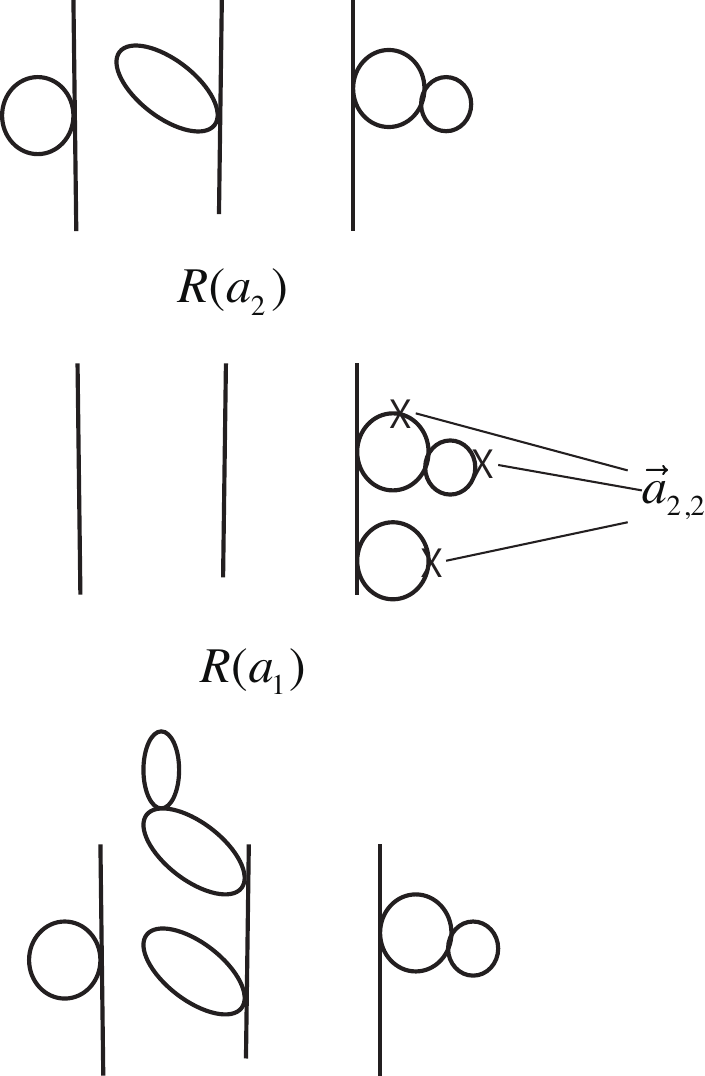}
\caption{Fiber product \eqref{form519}.}
\label{Figure57}
\end{figure}
\begin{defn}\label{defnevalama2}
We define the {\it evaluation maps}
\index{evaluation map}
\begin{gather*}
{\rm ev} = \bigl({\rm ev}^1,{\rm ev}^{12},{\rm ev}^2\bigr) =
\bigl(\bigl({\rm ev}^1_1,\dots,{\rm ev}^1_{k_1}\bigr),\bigl({\rm ev}^{12}_1,\dots,{\rm ev}^{12}_{k_{12}}\bigr),
\bigl({\rm ev}^2_1,\dots,{\rm ev}^2_{k_2}\bigr)\bigr) \colon\\
\qquad
{\mathcal M}_{\rm QT}(\vec a_1,\vec a_{12},\vec a_2;a_-,a_+;E)
\to \prod_{j=1}^{k_1} L_1(a_{1,j}) \times \prod_{j=1}^{k_{12}} L_{12}(a_{12,j})
\times \prod_{j=1}^{k_2} L_2(a_{2,j})
\end{gather*}
and
\[
{\rm ev}_{\infty} = ({\rm ev}_{\infty,-},{\rm ev}_{\infty,+}) \colon\
{\mathcal M}_{\rm QT}(\vec a_1,\vec a_{12},\vec a_2;a_-,a_+;E)
\to R(a_-) \times R(a_+)
\]
in the same way as Definition~\ref{defnevalama}.
\end{defn}
\begin{prop}\label{prop530}
We can define a topology on
${\mathcal M}_{\rm QT}(\vec a_1,\vec a_{12},\vec a_2;a_-,a_+;E)$
by which this space is compact and Hausdorff.
\end{prop}

The topology we use is the stable map topology which is
similar to \cite[Definitions 7.1.39 and~7.1.42]{fooobook2}
and \cite[Definition 10.3]{FO}.
The proof of the proposition is similar to one in
\cite[Definition~10.3]{FO}.
The only new point is the way how we handle disk bubbles
on the seam $\{0\} \times \R$
and more importantly the sphere bubbles
on such disk bubbles.
This is the point related to Remark~\ref{Remark524}.
We will discuss this point in detail in Section~\ref{sec:directcomp}.

\begin{thm}\label{therem530}
The space ${\mathcal M}_{\rm QT}(\vec a_1,\vec a_{12},\vec a_2;a_-,a_+;E)$
has a Kuranishi structure with corners with the following properties:
\begin{enumerate}\itemsep=0pt
\item[$(1)$]
We denote the codimension $d$ normalized corner
of the space with Kuranishi structure,
${\mathcal M}_{\rm QT}(\vec a_1,\vec a_{12},\vec a_2;a_-,a_+;E)$,
by $S_d{\mathcal M}_{\rm QT}(\vec a_1,\vec a_{12},\vec a_2;a_-,a_+;E)$.
Then it is a union of the fiber products
\begin{gather*}
S_{d_1}\mathring{\mathcal M}_{\rm QT}(\vec a_{1,0},\vec a_{12,0},\vec a_{2,0};a_0,a_1;E_1)\\
\qquad\times_{R(a_1)} S_{d_2}\mathring{\mathcal M}_{\rm QT}(\vec a_{1,1},\vec a_{12,1},\vec a_{2,1};a_1,a_2;E_2)
 \times_{R(a_2)} \cdots \\
\qquad\times_{R(a_{\ell-1})} S_{d_{\ell}}\mathring{\mathcal M}_{\rm QT}(\vec a_{1,\ell},\vec a_{12,\ell},
\vec a_{2,\ell};a_{\ell-1},a_\ell;E_{\ell})
\end{gather*}
of the form \eqref{form519}, where
$d_1+\dots+d_{\ell}+\ell -1 \ge d$
\item[$(2)$]
The codimension $d_j$ normalized corner
$S_{d_j}{\mathcal M}_{\rm QT}(\vec a_{1,j},\vec a_{12,j},\vec a_{2,j};a_j,a_{j+1};E)$
is the union of the closure of subsets
\begin{align*}
\overset{\ \text{\tiny $\circ\circ$}}{\mathcal M}_{\rm QT}(\vec a'_1,\vec a'_{12},\vec a'_2;a_-,a_+;E')
&\times_{{\rm ev}_0,\dots,{\rm ev}_0}
\prod_{j=1}^{k_1}S_{d'_{1,\ell_j}}{\mathcal M}(L_1; \vec a_{1,j};E_{1,j}) \\
&\times_{{\rm ev}_0,\dots,{\rm ev}_0}
\prod_{j=1}^{k_{12}}S_{d'_{12,\ell_j}}{\mathcal M}'(L_{12}; \vec a_{12,j};E_{12,j})\\
&\times_{{\rm ev}_0,\dots,{\rm ev}_0}
\prod_{j=1}^{k_2}S_{d'_{2,\ell_j}}{\mathcal M}(L_2; \vec a_{2,j};E_{2,j})
\end{align*}
of \eqref{eq518} such that there are $k'_1+k'_2+k'_3+1$ factors
other than those of the form of one of
${\mathcal M}(L_1; (a,a);0)$, ${\mathcal M}(L_{12}; (a,a);0)$, ${\mathcal M}(L_2; (a,a);0)$
and
\[
d_j = k'_1+k'_2+k'_3 + \sum_{j=1}^{k_1} d'_{1,\ell_j} + \sum_{j=1}^{k_{12}} d'_{12,\ell_j}
+ \sum_{j=1}^{k_2} d'_{2,\ell_j}.
\]
\item[$(3)$]
The evaluation maps defined in \eqref{defnevalama} are
the underlying continuous maps of strongly smooth maps.
\item[$(4)$]
The evaluation maps defined in \eqref{defnevalama2} are
the underlying continuous maps of strongly smooth maps. ${\rm ev}_{\infty,+}$ is weakly
submersive also.
\item[$(5)$]
The
fiber product description \eqref{eq518} and \eqref{form519} are compatible
with the Kuranishi structures.
Namely, there exists an isomorphism between
Kuranishi structures on the moduli space ${\mathcal M}_{\rm QT}(\vec a_1,\vec a_{12},\vec a_2;a_-,a_+;E)$
with ones obtained as the fiber product Kuranishi structures
of \eqref{eq518} or \eqref{form519}.
Here on the spaces appearing in the second, third and fourth factors of
\eqref{eq518} we take the Kuranishi structures
given in Theorem {\rm\ref{thekuraexist}}.
\item[$(6)$]
The isomorphisms of the Kuranishi structures in
item $(5)$ satisfies corner compatibility conditions
which are similar to Condition {\rm\ref{cccond}}.
\item[$(7)$] Given relative spin structures of
$L_1$, $L_{12}$, $L_2$ $($with
background data $V_1$, $\pi^*_1(V_1 \oplus TX_1) \oplus \pi_2^*(V_2)$, $V_2$, respectively$)$
we can define a principal ${\rm O}(1)$ bundle $\Theta^-_{12,a}$
on $R(a)$ such that the orientation
bundle of ${\mathcal M}_{\rm QT}(\vec a_1,\vec a_{12},\vec a_2;a_-,a_+;E)$
is canonically isomorphic to the tensor
product of the pullbacks of \smash{$\Theta^-_{a_{1,i}}$}, \smash{$\Theta^-_{12,a_{i}}$}
\smash{$\Theta^-_{a_{2,i}}$}, \smash{$\Theta^-_{a_{\pm}}$}. The isomorphism is
compatible with the description of the boundary which
is a part of item $(1)$.\footnote{In the case of moduli space of
holomorphic disks, a precise meaning of compatibility at boundary
with orientation is written as \cite[Condition 21.6\,(IX)]{fooonewbook},
when $L$ is embedded.
There is an explicit correction term of sign in
\cite[Condition 21.6\,(IX)]{fooonewbook} which coincides with one
in \cite{fooobook2} and \cite{ST}. However, the discussion of this
paper is not affected by the explicit form of correction terms.
See Remark~\ref{rem172}.
In the case $L$ is immersed with self-transversal intersection,
it is given in \cite[equation~(73)]{AJ}.
The way to generalize it to the self-clean case
is in Section~\ref{oriAinfMB} and in the paper \cite{ono2} by Kaoru Ono.
In the way we explain in Section~\ref{sec:orient}, the case of the moduli space
of quilt etc.\ can be reduced to the case of disks.}
\end{enumerate}

\end{thm}

Most of the proof of Theorem~\ref{therem530} is the same as the proof of
Theorem~\ref{thekuraexist} and is now becoming a routine,
in the study of pseudo-holomorphic curves based on the virtual fundamental
chain technique. (See \cite{fooo:const2}.)
The only point we need a discussion other than those in
Theorem~\ref{therem530} is the way how we handle
the point mentioned in Remark~\ref{Remark524}.
We will discuss it in Section~\ref{sec:directcomp}.\looseness=-1

See Sections~\ref{orisimpquilt}, \ref{oriAinfMB} and \cite{ono2} for item (7).

We finally mention the gappedness, which is related to
Gromov-compactness.
We define
\begin{align*}
\begin{split}
G_0(L_1,L_{12},L_2) :={}&
\bigl\{E \in \R_{\le 0} \mid\text{$\overset{\ \text{\tiny $\circ\circ$}}{\mathcal M}_{\rm QT}(\vec a_1,\vec a_{12},\vec a_2;a_-,a_+;E)$} \\
&\text{ \  is nonempty for some $\vec a_1,\vec a_{12},\vec a_2;a_-,a_+$}\bigr\}.
\end{split}
\end{align*}
Gromov compactness implies that $G_0(L_1,L_{12},L_2)$ is a discrete subset of
$\R_{\ge 0}$.
Let $G_0(L_1)$, $G_0(L_{12})$, $G_0(L_2)$ be as in \eqref{defnG0(L)}.
\begin{defn}
We define $G(L_1,L_{12},L_2)$ to be the discrete submonoid generated by the union of
$G_0(L_1,L_{12},L_2)$, $G_0(L_1)$, $G_0(L_{12})$, $G_0(L_2)$.

\end{defn}
The next lemma is obvious.
\begin{lem}
The set $G(L_1,L_{12},L_2)$ is a discrete submonoid.
If the moduli space ${\mathcal M}(\vec a_1,\vec a_{12},\vec a_2;\allowbreak a_-,a_+;E)$ is
non-empty, then $E \in G(L_1,L_{12},L_2)$.\index[syindex]{GL1L12@$G(L_1,L_{12},L_2)$}
\end{lem}
The filtered $A_{\infty}$ tri-module in Theorem~\ref{trimain}
will be $G(L_1,L_{12},L_2)$-gapped.

\subsection[A geometric realization of an $A_\infty$ tri-module 2]{A geometric realization of an $\boldsymbol{A_{\infty}}$ tri-module 2}
\label{subsec:bi-functorgeo2}

Using the system of Kuranishi structures given in Theorem~\ref{therem530},
we can define a system of CF-perturbations.
We will state it as Proposition~\ref{prop536} below.
We first describe the situation we work with precisely.

\begin{lem}
The conclusions of Theorem {\rm\ref{thekuraexist}} and Proposition {\rm\ref{prop330}} still hold when we replace
the compactification ${\mathcal M}(L_{12}; \vec a_{12};E_{12})$ by
the other compactification ${\mathcal M}'(L_{12}; \vec a_{12};E_{12})$.
\end{lem}

The proof is the same as the proof of Theorem~\ref{thekuraexist} and Proposition~\ref{prop330}
once the definition of ${\mathcal M}'(L_{12}; \vec a_{12};E_{12})$ is
understood. See Theorem~\ref{prop1417}.

\begin{situ}\label{situ535}
Let $E_0 > 0$.
We are given a system of CF-perturbations
of ${\mathcal M}(L_1; \vec a_{1};E_{1})$,
${\mathcal M}(L_2; \vec a_{2};E_{2})$,
${\mathcal M}'(L_{12}; \vec a_{12};E_{12})$
for $E_1, E_2, E_{12} < E_0$, so that they satisfy the
conclusions of Theorem~\ref{thekuraexist} and Proposition~\ref{prop330}.
\end{situ}
\begin{prop}\label{prop536}
Let $E_0 > 0$.
There exists a system of CF-perturbations
\smash{${\widehat{\mathfrak S}}$}
\index[syindex]{Sfrak@${\widehat{\mathfrak S}}$} on the moduli spaces
${\mathcal M}(\vec a_1,\vec a_{12},\vec a_2;a_-,a_+;E)$
with Kuranishi structures which are outer collarings of thickenings of those given in
Theorem {\rm\ref{therem530}}.
It enjoys the following properties:
\begin{enumerate}\itemsep=0pt
\item[$(1)$]
The CF-perturbations ${\widehat{\mathfrak S}}$ are transversal to zero.
\item[$(2)$]
The evaluation map ${\rm ev}_0$ is strongly submersive with respect to
this CF-perturbation
$($see {\rm\cite[\emph{Definition} 9.2]{foootech2}} for the definition of strong submersivity$)$.
\item[$(3)$]
They are compatible with the fiber product description
of their corners given in Theorem~{\rm\ref{therem530}}.
Here we use CF-perturbations in Situation {\rm\ref{situ535}} on those factors in the same sense as Proposition {\rm\ref{prop330}}.
\item[$(4)$]
They are compatible with the forgetful maps of the marked points which
corresponds to the diagonal component other than $0$-th one. The
precise definition of compatibility is written in {\rm\cite[Definition 5.1]{fooo091}}.
\end{enumerate}

\end{prop}
\begin{proof}
The proof is by the general theory of Kuranishi structures, such as those
developed in~\cite{foootech2,foootech22,fooonewbook}.
See \cite{fooo091} for item (4).
\end{proof}

\begin{defn}\label{defn53939}
\quad
\begin{enumerate}\itemsep=0pt
\item[(1)]
We put
\begin{equation}\label{form525}
\overline D = CF(L_1,L_{12},L_2;\R) \cong \Omega\bigl(\bigl(\tilde L_1 \times \tilde L_2\bigr)
\times_{X_1\times X_2}\tilde L_{12};\Theta^-\bigr),
\end{equation}
where $\Theta^-$ is a $\Z_2$ local system
defined on the fiber product $\bigl(\tilde L_1 \times \tilde L_2\bigr)
\times_{X_1\times X_2}\tilde L_{12}$
by Theorem~\ref{therem530}\,(7),
and
\smash{$D = CF(L_1,L_{12},L_2;\Lambda_0) = \overline D \,\widehat{\otimes}_{\R}\, \Lambda_0$}.\index[syindex]{CFL1C12@$CF(L_1,L_{12},L_2;\Lambda_0)$}
Then $D$ is a cochain complex with differential $\delta = d$.
\item[(2)]
We will define the structure operations\index[syindex]{nk1k12k2E@$\mathfrak{n}_{k_1,k_{12},k_2}^{E,\varepsilon}$}
\[
\mathfrak{n}_{k_1,k_{12},k_2}^{E,\varepsilon} \colon\
B_{k_1}CF(L_1;\R)[1] \otimes B_{k_{12}}CF(L_{12};\R)[1]
\otimes \overline D[1] \otimes
B_{k_2}CF(L_2;\R)[1]
\to \overline D[1]
\]
as follows.
Let
\begin{gather*}
{\bf x} =x_1 \otimes \dots \otimes x_{k_1} \in B_{k_1}CF(L_1;\R)[1],\\
{\bf y} =y_{12} \otimes \dots \otimes y_{k_{12}} \in B_{k_{12}}CF(L_{12};\R)[1],\\
{\bf z} =z_1 \otimes \dots\otimes z_{k_2} \in B_{k_2}CF(L_2;\R)[1],
\end{gather*}
and $w \in D$.
Then
\begin{align}
\mathfrak{n}_{k_1,k_{12},k_2}^{E,\varepsilon}({\bf x},{\bf y},w,{\bf z})
:={}&
{\rm ev}_{\infty,+}!\bigl({\rm ev}_{1,1}^*x_1\wedge \dots \wedge {\rm ev}_{1,k_1}^*x_{k_1} \wedge
{\rm ev}_{12,1}^*y_{1}\wedge \dots \wedge {\rm ev}_{12,k_{12}}^*y_{k_{12}}\nonumber
\\
& \wedge w \wedge {\rm ev}_{2,1}^*z_1\wedge \dots \wedge {\rm ev}_{2,k_2}^*z_{k_2}
;{\widehat{\mathfrak S^{\varepsilon}}}\bigr).\label{form526}
\end{align}
Here we use the integration along the fiber on the moduli spaces
${\mathcal M}(\vec a_1,\vec a_{12},\vec a_2;a_-,a_+;E)$
with Kuranishi structures and its CF-perturbations
\smash{$\widehat{\mathfrak S}$} in Proposition~\ref{prop536} to define the right-hand side
(see \cite[Definition 7.79]{fooonewbook}).\footnote{We remark that $\delta$ is the boundary operator of $D$.
The case $k_1,k_2,k_{12} = 0$, $E\ne 0$, the map
$\mathfrak{n}_{0,0,0}^{E,\varepsilon}$ may be nonzero and is a
deformation of the boundary operator of $D$ obtained by using
moduli spaces ${\mathcal M}(\varnothing,\varnothing,\varnothing;a_-,a_+;E)$.}
\item[(3)]
We finally put
\[
\mathfrak{n}_{k_1,k_{12},k_2}^{<E_0,\varepsilon}
=
\sum_{E<E_0,\, E \in G(L_1,L_{12},L_2)} T^{E}\mathfrak{n}_{k_1,k_{12},k_2}^{E,\varepsilon}.
\]
This is a map \index[syindex]{Fk1k12k2@$\mathscr{F}_{k_1,k_{12},k_2}^{<E_0,\varepsilon}$}
\[
\mathfrak{n}_{k_1,k_{12},k_2}^{<E_0,\varepsilon} \colon\
B_{k_1}CF(L_1)[1] \otimes B_{k_{12}}CF(L_{12})[1]\otimes D[1] \otimes
B_{k_2}CF(L_2)[1]
\to D[1].
\]
\end{enumerate}
\end{defn}

\begin{rem}
We remark that we need a certain sign $(-1)^*$ in \eqref{form526}.
We will prove in Section~\ref{sec:orient} that there {\it exists}
a choice of the sign so that $A_{\infty}$ relation~\eqref{form528} holds
with sign. The sign~$*$ is in principle calculable from the
discussion of Section~\ref{sec:orient} and the sign given in~\cite{AJ,fooobook2,fooonewbook}, Section~\ref{oriAinfMB} and~\cite{ono2}. Since all we need to prove the main results
of this paper are {\it existence} of sign~$*$ and not its explicit
formula we do not try to calculate it.
We do not repeat this remark in several other places.
\end{rem}
\begin{prop}\label{prop558}
\smash{$\mathfrak{n}_{k_1,k_{12},k_2}^{<E_0,\varepsilon}$}
defines a filtered $A_{\infty}$ tri-module modulo $T^{E_0}$.
Namely, it satisfies the congruence
\begin{align}
0 \equiv{}&\sum_{c_1,c_{12},c_2} (-1)^{*_1}
\mathfrak{n}^{<E_0,\varepsilon}
\bigl({\bf x}_{c_1;1},{\bf y}_{c_{12};1},
\mathfrak{n}^{<E_0,\varepsilon}({\bf x}_{c_1;2},{\bf y}_{c_{12};2},w,{\bf z}_{c_2;1}),{\bf z}_{c_2;2}\bigr)\nonumber\\
&+
(-1)^{*_2}\mathfrak{n}^{<E_0,\varepsilon}\bigl(\widehat d{\bf x},{\bf y},w,{\bf z}\bigr) +
(-1)^{*_3}\mathfrak{n}^{<E_0,\varepsilon}\bigl({\bf x},\widehat d{\bf y},w,{\bf z}\bigr)\nonumber\\
&+
(-1)^{*_4}\mathfrak{n}^{<E_0,\varepsilon}\bigl({\bf x},{\bf y},w,\widehat d{\bf z}\bigr) +
(-1)^{*_5}\delta \mathfrak{n}^{<E_0,\varepsilon}({\bf x},{\bf y},w,{\bf z})\nonumber\\
&+
(-1)^{*_6}\mathfrak{n}^{<E_0,\varepsilon}({\bf x},{\bf y},\delta w,{\bf z})
\mod T^{E_0}.\label{form528}
\end{align}
This filtered $A_{\infty}$ tri-module modulo $T^{E_0}$ is unital.
\end{prop}
The notations in \eqref{form528} is as follows.
We define ${\bf x}_{c_1;1}$, ${\bf x}_{c_1;2}$ by
$
\Delta ({\bf x}) = \sum_{c_1} {\bf x}_{c_1;1}\otimes {\bf x}_{c_1;2}
$.
Here~$c_1$ runs over a certain index set depending on ${\bf x}$.
The definitions of ${\bf y}_{c_{12};1}$, ${\bf y}_{c_{12};1}$, ${\bf z}_{c_2;1}$,
${\bf z}_{c_2;2}$ are similar.
The symbol $\widehat d$ in the second (resp.\ third, fourth) term of
\eqref{form528} is the derivation induced on
$BCF[1](L_1)$ (resp.\ $BCF[1](L_{12})$, $BCF[1](L_{2})$) by its filtered $A_{\infty}$
structure modulo~$T^{E_0}$.
$\delta$~is the operator induced from the de Rham differential in the same way as~\eqref{form3420000}.
We omit the indices $k_i$ etc.
of the operator $\mathfrak n$ since they are automatically determined by the variables
plugged in.
The signs $*_i$, $i=1,\dots,6$, are determined by Koszul rule.
We explain Koszul rule in detail in Section~\ref{Koszul}
\begin{proof}
The proof is a routine using Theorem~\ref{therem530}, Proposition~\ref{prop536}, Stokes' formula and the
composition formula and proceeds as follows.

By Stokes' formula (see \cite[Theorem 8.11]{fooonewbook}), we have
\begin{gather}
(-1)^{*_5}\delta \mathfrak{n}^{<E,\varepsilon}({\bf x},{\bf y},w,{\bf z})
+
(-1)^{*_7}\mathfrak{n}^{<E,\varepsilon}(\delta ({\bf x},{\bf y},w,{\bf z}))\nonumber
\\
\qquad=\sum_{E<E_0} T^E {\rm ev}_{\infty,+} !\bigl({\rm ev}_{1,1}^*x_1\wedge \dots \wedge {\rm ev}_{1,k_1}^*x_{k_1}\wedge
{\rm ev}_{12,1}^*y_{1}\wedge \dots \wedge {\rm ev}_{12,k_{12}}^*y_{k_{12}} \wedge w\nonumber
\\
\phantom{\qquad=}{}\wedge {\rm ev}_{2,1}^*z_1\wedge \dots \wedge {\rm ev}_{2,k_2}^*z_{k_2}
: \bigl(\partial {\mathcal M}_{\rm QT}(\vec a_1,\vec a_{12},\vec a_2;a_-,a_+;E),{\widehat{\mathfrak S^{\varepsilon}}}\bigr)\bigr).\label{form529}
\end{gather}
We include
the symbol $\partial {\mathcal M}_{\rm QT}(\vec a_1,\vec a_{12},\vec a_2;a_-,a_+;E)$
in the right-hand side to clarify the fact that we use this space to define
the integration along the fiber.
(We used ${\mathcal M}_{\rm QT}(\vec a_1,\vec a_{12},\vec a_2;a_-,a_+;E)$ in \eqref{form526}.)
There is actually a sign in the right-hand side of \eqref{form529}. We will explain it in Section~\ref{orisimpquilt}.

By Theorem~\ref{therem530} and \eqref{eq518},
the normalized boundary
$\partial {\mathcal M}_{\rm QT}(\vec a_1,\vec a_{12},\vec a_2;a_-,a_+;E)$
is the union of the following four types of fiber products.

The first type is
\begin{equation}\label{eq530}
{\mathcal M}_{\rm QT}(\vec a_{1,0},\vec a_{12,0},\vec a_{2,0};a_-,a;E_1)
\times_{R(a)} {\mathcal M}_{\rm QT}(\vec a_{1,1},\vec a_{12,1},\vec a_{2,1};a,a_+;E_2),
\end{equation}
where $\vec a_{1,0} \sqcup \vec a_{1,1} = \vec a_1$,
$\vec a_{12,0} \sqcup \vec a_{12,1} = \vec a_{12}$,
$\vec a_{2,0} \sqcup \vec a_{2,1} = \vec a_2$,
$E_1 + E_2 = E$.
See Figure~\ref{Figure58}.

\begin{figure}[ht]\centering
\begin{tabular}{cc}
\begin{minipage}[t]{0.45\hsize}
\centering
\includegraphics[scale=0.4]{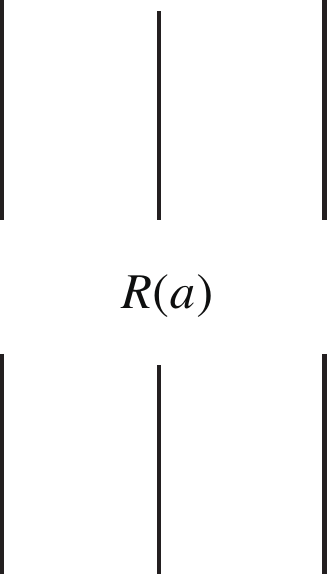}
\caption{Fiber product \eqref{eq530}.}
\label{Figure58}
\end{minipage} &
\begin{minipage}[t]{0.45\hsize}
\centering
\includegraphics[scale=0.4]{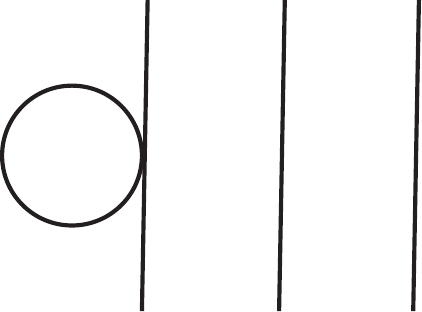}
\caption{Fiber product \eqref{eq531}.}
\label{Figure59}
\end{minipage}
\end{tabular}
\end{figure}

The second type is
\begin{equation}\label{eq531}
{\mathcal M}(L_1; \vec a^{ \prime\prime}_{1};E_{2})
\times_{{\rm ev}_0} {\mathcal M}_{\rm QT}(\vec a^{ \prime}_1,\vec a_{12},\vec a_2;a_-,a_+;E_1).
\end{equation}
Here $\vec a^{ \prime}_1 = (a_{1,1},\dots,a_{1,i-1},b,a_{1,j+1},\dots,a_{1,k_1})$,
$\vec a^{ \prime\prime}_{1} = (b,a_{1,i},\dots,a_{1,j})$ for some $1 \le i\le j \le k_1$ and $b \in \mathcal A_{L_1}$.
See Figure~\ref{Figure59}.

The third type is
\begin{equation}\label{eq532}
{\mathcal M}'(L_{12}; \vec a^{ \prime\prime}_{12};E_{2})
\times_{{\rm ev}_0} {\mathcal M}_{\rm QT}(\vec a_1,\vec a^{ \prime}_{12},\vec a_2;a_-,a_+;E_1).
\end{equation}
Here $\vec a^{ \prime}_{12} = (a_{12,1},\dots,a_{12,i-1},b,a_{12,j+1},\dots,a_{k_{12}})$,
$\vec a^{ \prime\prime}_{12} = (b,a_{12,i},\dots,a_{12,j})$ for some $1\le i\le j \le k_{12}$ and $b \in \mathcal A_{L_{12}}$.
See Figure~\ref{Figure510}.

\begin{figure}[ht]\centering
\begin{tabular}{cc}
\begin{minipage}[t]{0.45\hsize}
\centering
\includegraphics[scale=0.4]{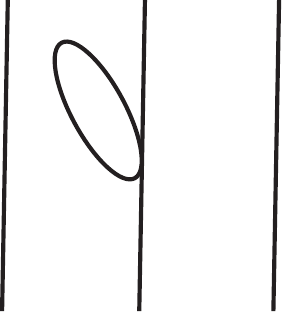}
\caption{Fiber product \eqref{eq532}.}
\label{Figure510}
\end{minipage} &
\begin{minipage}[t]{0.45\hsize}
\centering
\includegraphics[scale=0.4]{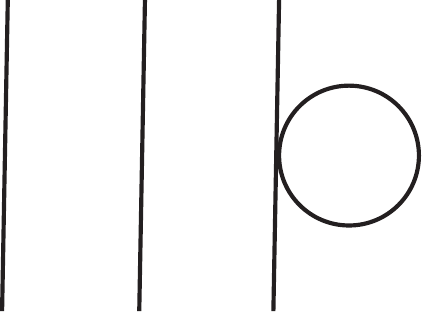}
\caption{Fiber product \eqref{eq533}.}
\label{Figure511}
\end{minipage}
\end{tabular}
\end{figure}

The fourth type is
\begin{equation}\label{eq533}
{\mathcal M}(L_1; \vec a^{ \prime\prime}_{2};E_{2})
\times_{{\rm ev}_0}
{\mathcal M}_{\rm QT}(\vec a_1,\vec a_{12},\vec a^{ \prime}_2;a_-,a_+;E_1).
\end{equation}
Here $\vec a^{ \prime}_2 = (a_{2,1},\dots,a_{2,i-1},b,a_{2,j+1},\dots,a_{2,k_2})$,
$\vec a^{ \prime\prime}_{2} = (b,a_{2,i},\dots,a_{2,j})$ for some $1 \le i\le j \le k_2$ and $b \in \mathcal A_{L_2}$.
See Figure~\ref{Figure511}.

By the composition formula
\cite[Theorem 10.21]{fooonewbook}, the integration along the fiber appearing in
\eqref{form529} on the spaces \eqref{eq530}
(resp.\ \eqref{eq531}, \eqref{eq532}, \eqref{eq533})
becomes the formula
\[
\sum_{c_1,c_{12},c_2} (-1)^{*_1}
\mathfrak{n}^{<E_0,\varepsilon}
({\bf x}_{c_1;1},{\bf y}_{c_{12};1},
\mathfrak{n}^{<E_0,\varepsilon}({\bf x}_{c_1;2},{\bf y}_{c_{12};2},w,{\bf z}_{c_2;1}),{\bf z}_{c_2;2}),
\]
\big(resp.\ the formula
\smash{$(-1)^{*_2}\mathfrak{n}^{<E_0,\varepsilon}\bigl(\widehat d{\bf x},{\bf y},w,{\bf z}\bigr)$}, the formula
\smash{$(-1)^{*_3}\mathfrak{n}^{<E_0,\varepsilon}\bigl({\bf x},\widehat d{\bf y},w,{\bf z}\bigr)$},
and the formula
\smash{$(-1)^{*_4}\mathfrak{n}^{<E_0,\varepsilon}\bigl({\bf x},{\bf y},w,\widehat d{\bf z}\bigr)$}\big).
This implies \eqref{form528}.
\end{proof}

Thus we defined a filtered $A_{\infty}$ tri-module modulo $T^{E_0}$.
The rest of the proof of
Theorem~\ref{trimain} is the same as the last step of the proof of
Theorem~\ref{AJtheorem}.
We first define the notion of a~pseudo-isotopy of $A_{\infty}$ tri-modules modulo $T^{E_0}$
in a similar way as Definition~\ref{pisotopydef}
(see Section~\ref{subsub:Pseudoisotopytri}).
We next show that for $E < E'$ the
$A_{\infty}$ tri-modulo modulo $T^{E'}$ we constructed in Proposition~\ref{prop558}
regarded as $A_{\infty}$ tri-module modulo $T^{E}$ is pseudo-isotopic
to the $A_{\infty}$ tri-module modulo $T^{E}$ we constructed in Proposition~\ref{prop558}.
We then prove a similar algebraic lemma as Lemma~\ref{lem339new}.
Using it, we complete the proof of Theorem~\ref{trimain}
in the same way as the last step of the proof of
Theorem~\ref{AJtheorem}.
Since this argument is now a routine, we omit the detail.
\end{proof}

\section[Unobstructedness is preserved by an unobstructed
Lagrangian\\ correspondence]{Unobstructedness is preserved by an unobstructed\\
Lagrangian correspondence}
\label{sec:Unobstructedness}

In this section, we prove Theorem~\ref{them4}.

\begin{situ}\label{situ61}
Suppose we are in Situation \ref{situ14}.
Moreover, we assume that, for $L_1 \in \mathbb L_1$ and $L_{12} \in \mathbb L_{12}$,
$L_1$ has clean transformation by $L_{12}$.
Let $(L_1,\sigma_1) \in \mathbb L_1$ and
$(L_{12},\sigma_{12}) \in \mathbb L_{12}$.
We consider the geometric transformation $(L_2,\sigma_2)
= L_1 \times_{X_1} L_{12}$ as in Definition~\ref{defn43}, where
the relative spin structure $\sigma_2$ is given later in Definition~\ref{def67}.
We assume $(L_2,\sigma_2)$ is contained in $\mathbb L_2$.

\end{situ}
\begin{situ}\label{situ62}
In Situation \ref{situ61},
we consider the filtered $A_{\infty}$ tri-module
$\mathscr{CF}(\mathbb L_1,\mathbb L_{12};\mathbb L_2)$ in
Theorem~\ref{trimain}.
We assume that $(L_1,\sigma_1) \in \mathbb L_1$ and
$(L_{12},\sigma_{12}) \in \mathbb L_{12}$ are unobstructed and
take their bounding cochains $b_1 \in CF(L_1)$, $b_{12}
 \in CF(L_{12})$.
\end{situ}
The main result of this section is as follows.
\begin{thm}\label{thm61}
In Situation {\rm\ref{situ62}},
we can choose a relative spin structure $\sigma_2$ such that
$(L_2,\sigma_2)$ is unobstructed.
Moreover, there exists a canonical choice of the gauge equivalence class of the bounding
cochain
$b_2$.
The gauge equivalence class of $b_2$ depends only on those
of $b_1$ and $b_{12}$.

\end{thm}

As we mentioned in Remark~\ref{rem16}\,(3), the bounding
cochain $b_2$ had been conjectured to be defined as the
virtual fundamental chain of a certain moduli space
(the moduli space of Figure~8 bubbles).
The author was trying to understand how we can use such a bounding
cochain to generalize the argument by Lekili--Lipyanskiy
beyond the monotone case using the $Y$-diagram.
Then he found that for this purpose we need an equality that
a certain element of the de Rham complex of a Lagrangian
submanifold becomes a cycle with respect to
the deformed Floer boundary operator via $b_1$, $b_{12}$, $b_2$.
The equality needed is \eqref{eq32}.
In fact, the homomorphism~\eqref{formjula711} becomes a
cochain map because of \eqref{eq32}.
As we will explain in Section~\ref{reltoBW},
the heuristic argument shows that the
bounding cochain obtained as the
virtual fundamental chain of the moduli space of
Figure 8 bubbles, after an appropriate gauge transformation,
satisfies \eqref{eq32}.
The author then found that the equality \eqref{eq32}
is strong enough to characterize $b_2$ (for given~$b_1$,~$b_{12}$)
and also \eqref{eq32} implies that $b_2$ is actually a bounding cochain.
Moreover, as we will see in Proposition~\ref{thm35}, we can solve \eqref{eq32}
uniquely. Thus we can use the algebraic
equation \eqref{eq32} in place of studying the moduli spaces,
to obtain the required bounding cochain.

Thus replacing the study of difficult moduli spaces
by a simple algebraic lemma (see Proposition~\ref{thm35})
is the key new idea of this paper.

\subsection[Right filtered $A_\infty$ modules and cyclic elements]{Right filtered $ \boldsymbol{A_{\infty}}$ modules and cyclic elements}
\label{subsec:cycliclgim}

The main idea of the proof of Theorem~\ref{thm61} is the same as
\cite[Section 3]{fu9} and is based on \cite[Proposition 3.5]{fu9}.
We repeat the argument here for the completeness sake and also
here we work over $\R$, while in \cite{fu9} we worked over
$\Z_2$.

\begin{defn}
Let $(C,\{\mathfrak m_k\})$ be a non-unital curved and filtered
$A_{\infty}$ algebra.
\begin{enumerate}\itemsep=0pt
\item[(1)]
A {\it filtered right $A_{\infty}$ module}
\index{filtered right $A_{\infty}$ module} over $(C,\{\mathfrak m_k\})$
is a left $\Lambda_0$ and right $(C,\{\mathfrak m_k\})$
filtered $A_{\infty}$ bi-module in the sense of
Definition~\ref{bimodulecat}.
\item[(2)]
We say a filtered right $A_{\infty}$ module is {\it $G$-gapped}
\index{$G$-gapped}
if its structure operations are all $G$-gapped.
\end{enumerate}
\end{defn}

More explicitly, a right filtered $A_{\infty}$ module over $(C,\{\mathfrak m_k\})$
is $(D,\{\mathfrak n_k \mid k=0,1,2,\dots\})$,
where
\begin{enumerate}\itemsep=0pt
\item[(1)] $D$ is a completed free $\Lambda_0$ module.
\item[(2)]
The operation $\mathfrak n_k$ is a $\Lambda_0$ moduli homomorphism\index[syindex]{nk@$\mathfrak n_k$}
\[
\mathfrak n_k \colon\ D[1] \,\widehat{\otimes}_{\Lambda_0}\, C[1]^{\otimes k} \to D[1]
\]
of degree $1$ which preserves filtration in the same sense as Definition
\ref{defn22}\,(2).\footnote{Here we shift the degree of elements of bi-module.}
\item[(3)]
The following holds for any $k$,
$y \in D$, $x_1,\dots,x_k \in C$:
\begin{align}
0 = {}&\sum_{k_1 + k_2 = k}
\mathfrak n_{k_1}(\mathfrak n_{k_2}(y;x_1,\dots,x_{k_2});x_{k_2+1},\dots,x_{k})\nonumber
\\
&+ \sum_{k_1+k_2 = k+1}\sum_{i=0}^{k_2}
(-1)^*
\mathfrak n_{k_1}(y;\dots,\mathfrak m_{k_1}(x_{i+1},\dots,x_{i+k_1}),\dots,x_k),\label{eq61}
\end{align}
where $* = \deg' y + \sum_{j=1}^{i-1}\deg' x_{j}$.
\end{enumerate}
\begin{defn}\label{defn33}
Let $(C,\{\mathfrak m_k\})$ be a $G$-gapped filtered $A_{\infty}$ algebra and
$(D,\{\mathfrak n_k\})$ a $G$-gapped right filtered $A_{\infty}$ module over $(C,\{\mathfrak m_k\})$.
We say an element ${\bf 1} \in D$ of degree $0$ a {\it cyclic element}\footnote{The word
cyclic element seems to be a standard one for an object satisfying a
condition such as (1). We remark that the notion of cyclic element has no
relation to the cyclic symmetry of the filtered $A_{\infty}$ algebra associated
to a Lagrangian submanifold.}\index{cyclic element} if the following holds:
\begin{enumerate}\itemsep=0pt\samepage
\item[(1)]
The map
$C \to D$ which sends $x$ to $\mathfrak n_1({\bf 1};x)$ is a
$\Lambda_{0}^{R}$ module isomorphism
$C \to D$.\footnote{Since $\deg'{\bf 1} = -1$,
$\deg'x = \deg'\mathfrak n_1({\bf 1};x)$.}
\item[(2)]
$\mathfrak n_0({\bf 1}) \equiv 0 \mod \Lambda_+^{R}$.
\end{enumerate}

\end{defn}

\begin{prop}\label{thm35}
Let $(C,\{\mathfrak m_k\})$ be a $G$-gapped filtered $A_{\infty}$ algebra and
$(D,\{\mathfrak n_k\})$ a $G$-gapped right filtered $A_{\infty}$ module over $(C,\{\mathfrak m_k\})$.
Suppose ${\bf 1} \in D$ is a cyclic element, which is $G$-gapped.
Then there exists a unique $G$-gapped bounding cochain $b$ of $(C,\{\mathfrak m_k\})$
such that
\begin{equation}\label{eq32}
\mathfrak n_0^b({\bf 1}) = 0,
\end{equation}
where we defined $\mathfrak n_0^b$ by
\begin{equation}\label{newfig63}
\mathfrak n_0^b(y) = \sum_{k=0}^{\infty} \mathfrak n_k(y;b,\dots,b).
\end{equation}

\end{prop}
\begin{proof}
We first prove the uniqueness.
Let $G = \{\lambda_i \mid i=0,1,2,\dots\}$,
where $0 = \lambda_0 < \lambda_1 < \lambda_2 < \cdots$.
We put
\[
{\bf 1} = \sum_{i=0}^{\infty} T^{\lambda_i}{\bf 1}_i,
\qquad
b = \sum_{i=1}^{\infty} T^{\lambda_i}b_i,\qquad
\mathfrak m_k = \sum_{i=0}^{\infty} T^{\lambda_i}\mathfrak m_{k,i},
\qquad
\mathfrak n_k = \sum_{i=0}^{\infty} T^{\lambda_i}\mathfrak n_{k,i}
\]
according to the definition of $G$-gappedness.
(Note that the coefficient of $T^{\lambda_0}$ ($\lambda_0 = 0$) of $b$ is~$0$ since
$b \in C \otimes \Lambda_{+,G}$.)

We calculate the coefficient of $T^{\lambda_n}$
of the equation \eqref{eq32} and obtain
\begin{equation}\label{cond33}
\mathfrak n_{1,0}({\bf 1}_0;b_{n})
+
\sum \mathfrak n_{k,m}({\bf 1}_{n_0};b_{n_1},\dots,b_{n_k}) = 0.
\end{equation}
Here the second term is the sum over all $k$, $m$, $n_0,n_1,\dots,n_k$ such that
\begin{equation}\label{cond34}
\lambda_n = \lambda_m + \lambda_{n_0} + \sum_{i=1}^k \lambda_{n_i}
\end{equation}
except the case $k=1$, $m=0$, $n_0 = 0$, $n_1 = n$.
(The case which we exclude here corresponds to the first term.)
Note that if $k$, $m$, $n_0,n_1,\dots,n_k$ satisfy
\eqref{cond34} then $n_i \le n$ for $i=0,\dots,k$.
Moreover, $n_i < n$ unless $k=1$, $m=0$, $n_0 = 0$, $n_1 = n$.
Therefore, we can solve \eqref{cond33} and obtain~$b_n$ uniquely
by induction on $n$.
(Here we use Definition~\ref{defn33}\,(1).)
Thus we proved that there exists a unique $G$-gapped element
$b \in C \otimes_{\Lambda_0^{R}} \Lambda_+^{R}$ satisfying \eqref{eq32}.
It remains to prove that this element $b$ satisfies the
Maurer--Cartan equation~\eqref{MCeq}.
We will prove
\begin{equation}\label{MCmod}
\sum_{k=0}^{\infty}\mathfrak m_k(b,\dots,b) \equiv 0 \mod T^{\lambda_c}
\end{equation}
by induction on $c \in \Z_+$.
We assume \eqref{MCmod} for $c \le n-1$ and will prove the case $c = n$ below.

We remark that the assumption implies that we have $\mathfrak n_{0,0} \circ \mathfrak n_{0,0} = 0$.
Using \eqref{eq61} and Definition~\ref{defn33}\,(2), we have
$
\mathfrak n_0(\mathfrak n_{1,0}({\bf 1}_0;x)) - \mathfrak n_{1,0}({\bf 1}_0;\mathfrak m_{1,0}(x)) = 0
$
for $x \in \overline C$.

We next consider $\mathfrak n_0 (\mathfrak n_{1,0}({\bf 1}_0;b_{n}))$.
Using \eqref{cond33}, we find
\[
\mathfrak n_0( \mathfrak n_{1,0}({\bf 1}_0;b_{n}))
=
- \sum \mathfrak n_0( \mathfrak n_{k,m}({\bf 1}_{n_0};b_{n_1},\dots,b_{n_k})).
\]
We calculate the right-hand side using \eqref{eq61} to obtain
\begin{gather}
 \sum \mathfrak n_{k_1,m_1}(\mathfrak n_{k_2,m_2}({\bf 1}_{n_0};b_{n_1},\dots,b_{n_{k_2}}),
\dots, b_{n_k}) \nonumber\\
\qquad- \sum \mathfrak n_{k_1,m_1}({\bf 1}_{n_0};b_{n_1},\dots,
\mathfrak m_{k_2,m_2}(b_{n_{i+1}},\dots,b_{n_{i+k_2}}),\dots,b_{n_k})\nonumber \\
\qquad - \sum \mathfrak n_{k,m}({\bf 1}_{n_0};b_{n_1},\dots,\mathfrak m_{1,0}(b_{n_j}),\dots,b_{n_k}).\label{eq3737}
\end{gather}
Here the sum in the first line is taken over
$k_1$, $k_2$, $m_1$, $m_2$, $n_0, \dots, n_k$ such that
$k_1 + k_2 = k$ and
$
\lambda_n =
\lambda_{m_1} +
 \lambda_{m_2} + \lambda_{n_0} + \sum_{i=1}^k \lambda_{n_i}
$, except $k_1 = 0$, $m_1 = 0$.

The sum in the second line is taken over
$k_1$, $k_2$, $m_1$, $m_2$, $n_0, \dots, n_k$ such that
$k_1 + k_2 = k+1$ and~${
\lambda_n =
\lambda_{m_1} +
 \lambda_{m_2} + \lambda_{n_0} + \sum_{i=1}^k \lambda_{n_i}
}$, except $m_2 = 0$, $k_2 =1$.
(The excluded case corresponds to the third line.)

The sum in the third line
is taken over
$k$, $m$, $j$, $n_0, \dots, n_k$ such that
$j=1,\dots, k$ and $\lambda_n = \lambda_m + \lambda_{n_0} + \sum_{i=1}^k \lambda_{n_i}$,
except $n_0 =0$, $k=1$, $m=0$.
We exclude this case since it is excluded in the second term of
\eqref{cond33}.

Note that the first line of \eqref{eq3737} vanishes because of the equality
\eqref{eq32}.

By using the induction hypothesis \eqref{MCmod} for $c \le n-1$,
the sum of the second and third lines cancel each other
except the sum
\[
-\sum \mathfrak n_{0,1}({\bf 1}_{0};
\mathfrak m_{k,m}(b_{n_1},\dots,b_{n_{k}})),
\]
which is taken over $k, m, n_1,\dots, n_k$
such that
$\lambda_n = \lambda_m + \sum_{i=1}^k \lambda_{n_i}$.
(In fact, this sum could be canceled with
$\mathfrak n_{0,1}({\bf 1}_{0};
\mathfrak m_{0,1}(b_{n}))$. However, this is the case excluded in the third line.)

Thus we have
\[
\mathfrak n_{1,0}({\bf 1}_0;\mathfrak m_{1,0} (b_{n}))
=
\mathfrak n_{1,0}(\mathfrak n_{1,0}({\bf 1}_0;b_{n}))
=
-\sum \mathfrak n_{0,1}({\bf 1}_{0};
\mathfrak m_{k,m}(b_{n_1},\dots,b_{n_{k}})).
\]
Using Definition~\ref{defn33}\,(1), it implies
\[
\mathfrak m_{1,0} (b_{n}) +
\sum \mathfrak m_{k,m}(b_{n_1},\dots,b_{n_{k}}) = 0.
\]
It implies \eqref{MCmod} for $c=n$.
The proof of Proposition~\ref{thm35} is now complete.
\end{proof}

\subsection{A geometric realization of a cyclic element}
\label{subsec:cycliclgeo}

In this section, we use Proposition~\ref{thm35}
to prove the existence part of Theorem~\ref{thm61}.

Suppose we are in Situation \ref{situ61}.
By definition (see \eqref{form525}),
\begin{equation}\label{form69}
CF((L_1,\sigma_1),(L_{12},\sigma_{12}),(L_2,\sigma_2)) \cong \Omega\bigl(\bigl(\tilde L_1 \times \tilde L_2\bigr)
\times_{X_1\times X_2}\tilde L_{12};\Theta^-\bigr) \,\widehat\otimes\, \Lambda_0.
\end{equation}
\begin{lem}\label{exirespi}
There exists a unique relative spin structure $\sigma_2$ such that
principal ${\rm O}(1)$ bundle~$\Theta^-$ in \eqref{form69} is trivial
on $\tilde L_2$.
\end{lem}
\begin{proof}
We have
\[
\bigl(\tilde L_1 \times \tilde L_2\bigr)
\times_{X_1\times X_2}\tilde L_{12}
\cong
\bigl(\tilde L_1 \times_{X_1} \tilde L_{12}\bigr) \times_{X_2} \tilde L_2
=
\tilde L_2 \times_{X_2} \tilde L_2
\]
(see Lemma~\ref{lem45}).
Therefore, the lemma follows from Lemmas~\ref{lem310} and \ref{lem44}.
\end{proof}

\begin{defn}\label{def67}
Let $\sigma_2$ be as in Lemma~\ref{exirespi}.
We call $(L_2,\sigma_2)$ the {\it geometric transformation}
\index{geometric transformation} of~$(L_1,\sigma_1)$ by $(L_{12},\sigma_{12})$.
\end{defn}
\begin{defn}\label{def672}
In the situation of Definition~\ref{def67},
let $b_1$ (resp.\ $b_{12}$) be a bounding cochain of~$CF(L_1,\sigma_1)$
(resp.\ $CF(L_{12},\sigma_{12})$).
We define\index[syindex]{nb1b12k@$\mathfrak n^{b_1,b_{12}}_k$}
\[
\mathfrak n^{b_1,b_{12}}_k \colon\ CF[1](L_1;L_{12};L_2) \otimes CF[1](L_2,\sigma_2)^{\otimes k}
\to CF[1](L_1;L_{12};L_2)
\]
by
\[
\mathfrak n^{b_1,b_{12}}_k(y;x_1,\dots,x_k)
=
\sum_{k_1=0}^{\infty}\sum_{k_{12}=0}^{\infty}
\mathfrak n_{k_1,k_{12},k}(b_1,\dots,b_1;b_{12},\dots,b_{12};y;x_1,\dots,x_k).
\]
The operation $\mathfrak n$ in the right-hand is defined by
Theorem~\ref{trimain}.
\end{defn}
Now we have the following.
\begin{lem}\label{lem6868}
In the situation of Definition {\rm\ref{def672}},
the operations $\mathfrak n^{b_1,b_{12}}_k$, $k=0,1,2,\dots$, define a~structure of right filtered $A_{\infty}$ module
on $CF(L_1;L_{12};L_2)$ over $CF(L_2,\sigma_2)$.
\end{lem}
The proof is a straightforward calculation and so is omitted.

In the simplest case $k = 0$,
Lemma~\ref{lem6868} becomes
\begin{equation}\label{4111}
\mathfrak n^{b_1,b_{12}}_0\bigl(\mathfrak n^{b_1,b_{12}}_0(h)\bigr) + \mathfrak n^{b_1,b_{12}}_1(h;\mathfrak m_0(1)) = 0.
\end{equation}
In a geometric language, its proof is roughly as follows.
We assume for simplicity that all the switching components of $L_2$ are zero-dimensional.
Let $(p_i,q_i,r_i) \in L_1 \times_{X_1} L_{12} \times_{L_2} L_2$ be in the switching
component $R(a_i)$ for $i=1,2$.
We consider the case $h = [p_1,q_1,r_1]$ and study
\[
\big\langle
\mathfrak n^{b_1,b_{12}}_0\bigl(\mathfrak n^{b_1,b_{12}}_0([p_1,q_1,r_1])\bigr),[p_2,q_2,r_2]\rangle.
\big\rangle
\]
As usual in various Floer theories, we consider the one-dimensional
moduli space ${\mathcal M}(a_1,a_2;E)$.
Its boundary contains the union of
$
{\mathcal M}(a_1,a;E_1) \times {\mathcal M}(a,a_2;E_2)
$
for various $a$ and $E_1$, $E_2$ with~${E_1 \!+\! E_2 = E}$.
The count of such boundary becomes
$
\langle
\mathfrak n_0(\mathfrak n_0([p_1,q_1,r_1])),[p_2,q_2,r_2]\rangle$.
(Here~$\mathfrak n_0$ is the boundary operator and we do not include bounding cochains
$b_1$, $b_{12}$.)
As usual in the Lagrangian Floer theory, the one-dimensional moduli space ${\mathcal M}(a_1,a_2;E)$
has other boundaries, which corresponds to various disk bubbles.
There are three kinds of disk bubbles, that are those on $L_1$, $L_{12}$, $L_{2}$.
By including bounding cochains $b_1$ and $b_{12}$, the effect of disk bubbles
on $L_1$, $L_{12}$ are cancelled. Therefore, only the disk bubble at $L_2$ remains.
It gives the term $\mathfrak n^{b_1,b_{12}}_1(h;\mathfrak m_0(1))$.
Thus \eqref{4111} follows.
Using the algebraic formalism, we have developed so far we can convert this
geometric argument to algebraic ones, which is the calculation to prove Lemma~\ref{lem6868}.
\begin{rem}
In Lemma~\ref{lem6868}, we do not need to assume that $(L_2,\sigma_2)$ is a geometric transform of
$(L_1,\sigma_1)$ by $(L_{12},\sigma_{12})$.
\end{rem}
\begin{prop}\label{prop610}
Let $(L_2,\sigma_2)$ be the geometric transformation of
$(L_1,\sigma_1)$ by $(L_{12},\sigma_{12})$.
Then we can choose our tri-module structure so that
\[
{\bf 1} \in \Omega^0\bigl(\bigl(\tilde L_1 \times \tilde L_2\bigr)
\times_{X_1\times X_2}\tilde L_{12};\R\bigr) \subset
CF((L_1,\sigma),(L_{12},\sigma_{12}),(L_2,\sigma_2))
\]
is a cyclic element of
$\bigl(CF((L_1,\sigma),(L_{12},\sigma_{12}),(L_2,\sigma_2)),\bigl\{\mathfrak n^{b_1,b_{12}}_k \bigr\}\bigr)$.

\end{prop}
Here ${\bf 1}$ is the zero form (function) 1 on the diagonal component $\tilde L_2 \subset
\bigl(\tilde L_1 \times_{X_1} \tilde L_{12}\bigr) \times_{X_2} \tilde L_2$.
\begin{proof}
Definition~\ref{defn33}\,(2) is the consequence of the fact that $d {\bf 1} =0$ and
$\mathfrak n_0^{b_1,b_{12}} \equiv \pm d \mod T^{\varepsilon}$.

We remark that \smash{$\mathfrak n_1^{b_1,b_{12}} \equiv \mathfrak n_{0,0,1} \mod T^{\varepsilon}$}.
We also remark that modulo $T^{\varepsilon}$, $\mathfrak n_{0,0,1}$ is
defined as the smooth correspondence via the moduli space
${\mathcal M}(\varnothing,\varnothing,a;o,b;0)$ of energy zero.
Namely,
\begin{equation}\label{formform41}
\mathfrak n_{0,0,1}(h)
\equiv
{\rm ev}_{\infty,+}!\bigl({\rm ev}_{\infty,-}^*(h); \widehat{\mathfrak S^{\varepsilon}};
{\mathcal M}(\varnothing,\varnothing,a;o,b;0)\bigr)
\mod T^{\varepsilon}.
\end{equation}
The notations are as follows.
In the notation ${\mathcal M}(\varnothing,\varnothing,a;o,b;0)$,
the symbol $\varnothing$ in the first component (resp.\ second component)
indicates that we do not put marked points on the line~${\operatorname{Re} z = -1}$
(resp.\ $\operatorname{Re} z = 0$). The symbol $a$ in the third component means that we put one marked
point on~${\operatorname{Re} z = 1}$ and require that this point goes to $L(a)$
in the sense of Condition~\ref{cond517}.
The symbol~$o$ in the fourth component means that we use the diagonal component
$\tilde L_2$ for the boundary condition (switching condition 2, Condition \ref{cond518}) when $\operatorname{Im} z \to - \infty$.
The symbol $b$ in the fourth component means that we use the component
$L_2(b)$ for the boundary condition (switching condition 2, Condition \ref{cond518}) when $\operatorname{Im} z \to + \infty$.
The symbol $0$ in the fifth component means that we consider the pseudo-holomorphic curve
with $0$ energy. (It is nothing but a constant map.)

In \eqref{formform41}, the maps ${\rm ev}_{\infty,+}$ and ${\rm ev}_{\infty,-}$
are evaluation maps defined on ${\mathcal M}(\varnothing,\varnothing,a;o,b;0)$
as in Definition~\ref{defnevalama2}. We pull back the differential form $h$ on
$\tilde L_2$ by ${\rm ev}_{\infty,-}$ and obtain a~differential
form on ${\mathcal M}(\varnothing,\varnothing,a;o,b;0)$,
a space with Kuranishi structure (see \cite[Definition~7.7.1]{fooonewbook}).
The symbol $\widehat{\mathfrak S^{\varepsilon}}$ denotes the CF-perturbation
defined on ${\mathcal M}(\varnothing,\varnothing,a;o,b;0)$ by Proposition~\ref{prop536}.
We use it to define the integration along the fiber ${\rm ev}_{\infty,+}!$
via the strongly submersive map~${\rm ev}_{\infty,+}$.
See Figure~\ref{Figure61}.

\begin{figure}[ht]\centering
\begin{tabular}{@{}cc@{}}
\begin{minipage}[t]{0.44\hsize}
\centering
\includegraphics[scale=0.35]{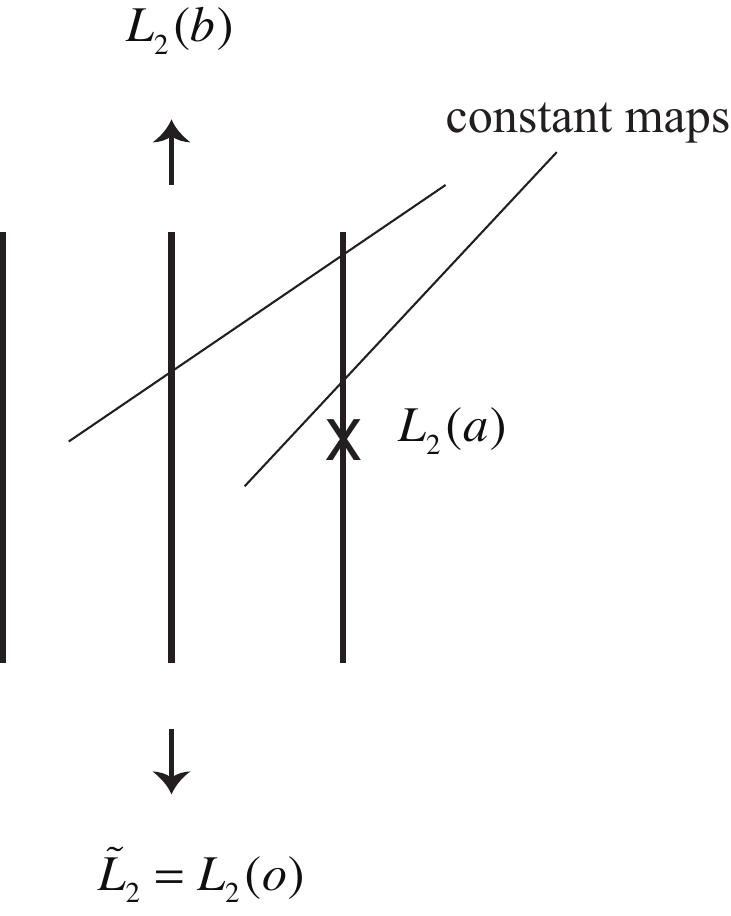}
\caption{An element of ${\mathcal M}(\varnothing,\varnothing,a;o,b;0)$.}
\label{Figure61}
\end{minipage} &
\begin{minipage}[t]{0.52\hsize}
\centering
\includegraphics[scale=0.35]{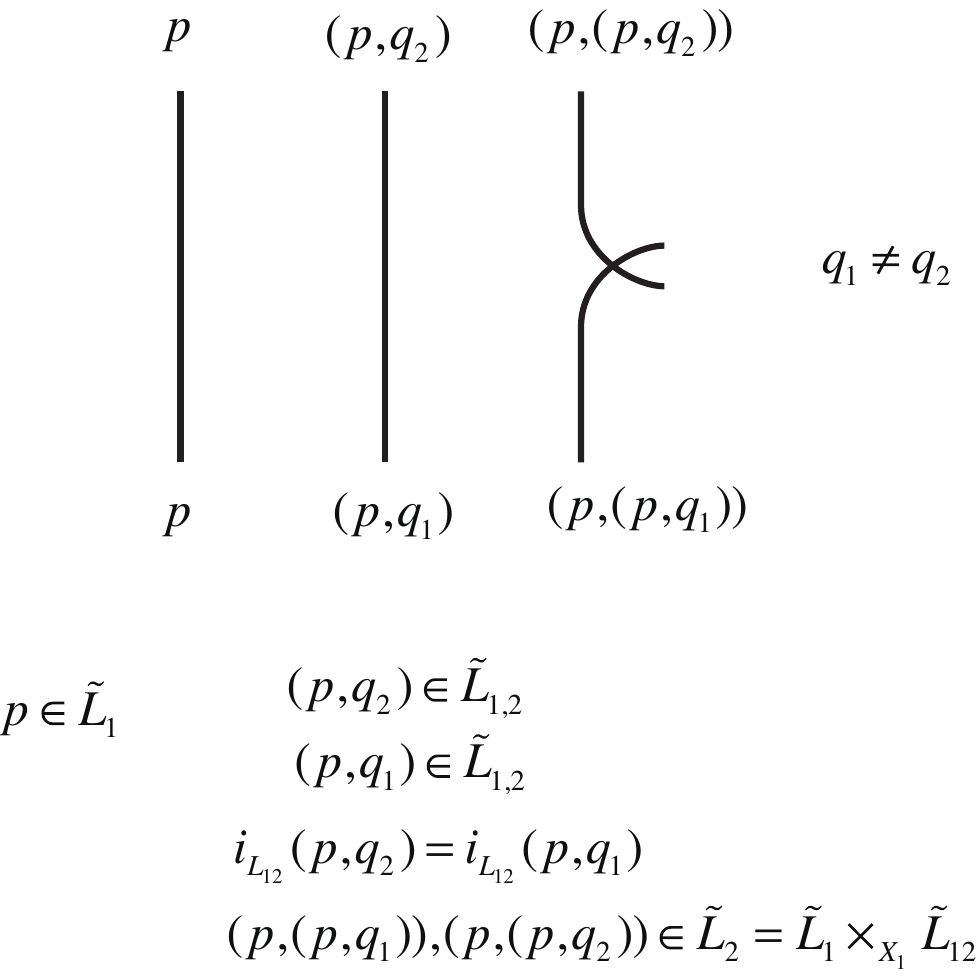}
\caption{An element of ${\mathcal M}(\varnothing,\varnothing,a;o,a;0)$, $a\ne o$.}
\label{Figure62}
\end{minipage}
\end{tabular}
\end{figure}

\begin{lem}\label{611}
${\mathcal M}(\varnothing,\varnothing,a;o,b;0)$ is an empty set if $a \ne b$.
If $a=b$ the space ${\mathcal M}(\varnothing,\varnothing,a;o,b;0)$
is diffeomorphic to $L_2(a)$ and evaluation map ${\rm ev}_2$ is a diffeomorphism.
Moreover, the moduli space ${\mathcal M}(\varnothing,\varnothing,a;o,b;0)$ is transversal.

\end{lem}
\begin{proof}
Since ${\mathcal M}(\varnothing,\varnothing,a;o,b;0)$ consists of constant maps,
the lemma is obvious except the statement about transversality.
See Figure~\ref{Figure62} in the case when $a=b$ is not diagonal component.

We show that ${\mathcal M}(\varnothing,\varnothing,a;o,a;0)$ is transversal.
We remark that this moduli space is identified with a connected component of the moduli space
of pseudo-holomorphic strip between $\tilde L_1 \times \tilde L_2$ and $\tilde L_{12}$.
Using the assumption that $\tilde L_1 \times \tilde L_2$ is of clean intersection with $\tilde L_{12}$, it is standard that this moduli space is transversal.
(In fact, the moduli space of pseudo-holomorphic strips with $0$ energy which bounds $L$ and $L'$ is
transversal if $L$ and $L'$ are of clean intersection.)
\end{proof}

Definition~\ref{defn33}\,(1) is an immediate consequence of Lemma~\ref{611}.
The proof of Proposition~\ref{prop610} is complete.
\end{proof}

Theorem~\ref{thm61} follows immediately from Propositions \ref{611} and \ref{thm35}.

\begin{defn}\label{def612}
In the situation of Theorem~\ref{thm61},
we call $(L_2,\sigma_2,b_2)$ the {\it geometric transformation} \index{geometric transformation}of $(L_1,\sigma_1,b_1)$
by $(L_{12},\sigma_{12},b_{12})$.
\end{defn}

\subsection{Well-definedness of bounding cochains up to gauge equivalence}

In this subsection, we prove that when we change the
bounding cochains $b_1$, $b_{12}$
by gauge equivalences the bounding
cochain $b_{2}$ in Definition~\ref{def612}
changes by a gauge equivalence.
Here we discuss only an algebraic part.
Namely, we fix the tri-module in Theorem~\ref{trimain}.
The independence of $b_2$
of the construction of the tri-module in Theorem~\ref{trimain}
will be proved in Section~\ref{sec:independence2},
Theorem~\ref{thm1466}.
The statement we prove is the next proposition.
\begin{situ}
Let $C_1$, $C_{12}$, $C_2$ be
curved filtered $A_{\infty}$ algebras and
Let $(D,\mathfrak n)$ be a left $C_1$, $C_{12}$
and right $C_2$ tri-module.
Let ${\bf 1} \in D$ be an element such that
\begin{enumerate}\itemsep=0pt
\item[(1)]
The map
$C_2 \to D$ which sends $x$ to $\mathfrak n_1({\bf 1};x)$ is an
$\Lambda_{0}^{R}$ module isomorphism
$C_2 \to D$.
\item[(2)]
$\mathfrak n_0({\bf 1}) \equiv 0 \mod \Lambda_+^{R}$.
\end{enumerate}

\end{situ}
A pair of bounding cochains $b_1$ and $b_{12}$ of $C_1$, $C_{12}$
defines a right filtered $A_{\infty}$ module structure on $D$ over $C_2$
by the next formula:
\begin{equation}\label{form61111}
\mathfrak n^{b_1,b_{12}}(y;x_1,\dots,x_k)
= \sum_{k_1,k_{12}} \mathfrak n_{k_1,k_{12},k}\bigl(b_1^{k_1},b_{12}^{k_{12}};y
; x_1,\dots,x_k\bigr).
\end{equation}
${\bf 1}$ is its cyclic element.
Therefore, by Proposition~\ref{thm35} there exists a unique bounding cochain $b_2$ such that
\[
\sum_k\mathfrak n^{b_1,b_{12}}\bigl({\bf 1};b_2^{k}\bigr) = 0.
\]
We write $b_2 = B(b_1,b_{12})$.
\begin{prop}\label{prop615615}
If $b_1$, $b_{12}$ are gauge equivalent to $b'_1$, $b'_{12}$,
then $B(b_1,b_{12})$ is gauge equivalent to $B(b'_1,b'_{12})$.

\end{prop}
\begin{proof}
We recall the definition of gauge equivalence in \cite[Section 4.3]{fooobook}.
For a completed free~$\Lambda_0$ module $C$, we define
$\operatorname{Poly}( [0,1],C)$ to be the set of all formal sums
\begin{equation}\label{polydefin}
\sum_{i=1}^{\infty} x_i(s) T^{\lambda_i}
+ \left(\sum_{i=1}^{\infty} y_i(s) T^{\lambda_i}\right)ds,
\end{equation}
where $x_i$, $y_i$ are polynomials (with variable $s$)
with coefficients in $\overline C$
and $\lambda_i \in \R_{\ge 0}$ with $\lim_{i\to \infty} \lambda_i = +\infty$.
$s$ and $ds$ are formal variables.

For $s_0 \in \R$, we define $\operatorname{Ev}(_{s_0} \colon \operatorname{Poly}( [0,1],C) \to C$
by sending the element \eqref{polydefin} to
\[
\sum_{i} x_i(s_0) T^{\lambda_i} \in C.
\]

In \cite[Definition 4.2.9]{fooobook}, we defined
filtered $A_{\infty}$ structures on the modules
$\operatorname{Poly}( [0,1],C_1)$, $\operatorname{Poly}( [0,1],C_{12})$, $\operatorname{Poly}( [0,1],C_2)$.\index[syindex]{Pzoly01C@$\operatorname{Poly}( [0,1],C_1)$}

During the proof of \cite[Theorem 5.2.3]{fooobook},
it is proved that if $D$ is a filtered $A_{\infty}$
bi-module over $C_1$, $C_2$ then
$\operatorname{Poly}( [0,1],D)$ is a filtered $A_{\infty}$
bi-module over $\operatorname{Poly}( [0,1],C_1)$, $\operatorname{Poly}( [0,1],C_2)$.
We can prove the same statement for
tri-module in the same way.
Thus
in our situation,
$\operatorname{Poly}( [0,1],D)$ is a filtered $A_{\infty}$
tri-module over $\operatorname{Poly}( [0,1],C_1)$, $\operatorname{Poly}( [0,1],C_{12})$,
$\operatorname{Poly}( [0,1],C_2)$.

Moreover, $\operatorname{Ev}_{s_0}$ defines a filtered $A_{\infty}$
algebra homomorphism or a filtered $A_{\infty}$
tri-module homomorphism.

The cyclic element ${\bf 1} \in D$ may be regarded as an element
of $\operatorname{Poly}( [0,1],D)$.

By assumption that $b_1$ (resp.\ $b_{12}$) is gauge equivalent to
$b'_1$ (resp.\ $b'_{12}$), there exists a
bounding cochain $\mathfrak b_1$ (resp.\ $\mathfrak b_{12}$)
of $\operatorname{Poly}( [0,1],C_1)$ (resp.\ $\operatorname{Poly}( [0,1],C_{12})$)
such that
\[
\operatorname{Ev}_{0}(\mathfrak b_1) = b_1, \qquad
\operatorname{Ev}_{1}(\mathfrak b_1) = b'_1,
\qquad
\operatorname{Ev}_{0}(\mathfrak b_{12}) = b_{12}, \qquad
\operatorname{Ev}_{1}(\mathfrak b_{12}) = b'_{12}.
\]
Using $\mathfrak b_1$ and $\mathfrak b_{12}$ in the same way as
\eqref{form61111}, we can define a
structure of right filtered $A_{\infty}$
module~$\bigl\{\mathfrak n^{\mathfrak b_1,\mathfrak b_{12}}_k\bigr\}$
on $\operatorname{Poly}( [0,1],D)$ over $\operatorname{Poly}( [0,1],C_2)$.

It is easy to see that ${\bf 1} \in \operatorname{Poly}( [0,1],D)$
is a cyclic element of $\bigl\{\mathfrak n^{\mathfrak b_1,\mathfrak b_{12}}_k\bigr\}$.
Therefore, by Proposition~\ref{thm35}
there exists a bounding cochain $\mathfrak b_2$ of
$\operatorname{Poly}( [0,1],C_2)$ such that
\[
\sum_k\mathfrak n^{\mathfrak b_1,\mathfrak b_{12}}_k\bigl({\bf 1};
\mathfrak b_2^{k}\bigr) = 0.
\]
It follows that
\[
\sum_k\mathfrak n^{b_1,b_{12}}_k\bigl({\bf 1};
\operatorname{Ev}_0(\mathfrak b_2)^{k}\bigr) = 0.
\]
Therefore, the uniqueness part of Proposition~\ref{thm35}
implies $\operatorname{Ev}_0(\mathfrak b_2) = b_2$.
In the same way, we can show
$\operatorname{Ev}_1(\mathfrak b_2) = b'_2$.
Thus $b_2$ is gauge equivalent to $b'_2$ as required.
\end{proof}

\section{Representability of correspondence functor}
\label{sec:represent}

\subsection{Statement}
\label{subsec:homotopyequiv}
Suppose we are in Situation \ref{situ61}.
We consider the correspondence tri-module $\mathscr{CF}(\mathbb L_1,\mathbb L_{12};\mathbb L_2)$
which is a left
$
\mathfrak{Fuk}(X_1,V_1,\mathbb L_1)
\times
\mathfrak{Fuk}(-X_1 \times X_2,\pi_1^*(V_1
\oplus TX_1)\oplus \pi_2^*V_2,\mathbb L_{12})
$
and right
$
\mathfrak{Fuk}(X_2,V_2,\mathbb L_2)
$
tri-module and which
we obtained in Theorem~\ref{trimain}.
\begin{notation}\label{not71}
Here and hereafter, we denote
\[
\mathfrak{Fuk}(-X_1 \times X_2) =
\mathfrak{Fuk}((X_1,-\omega_1) \times (X_2,\omega_2),\pi_1^*(V_1 \oplus TX_1) \oplus \pi_2^*(V_2),\mathbb L_{12})
\]
and
$
\mathfrak{Fuk}(X_1) = \mathfrak{Fuk}((X_1,\omega_1),V_1,\mathbb L_1)$,
$
\mathfrak{Fuk}(X_2)=\mathfrak{Fuk}((X_2,\omega_2),V_2,\mathbb L_2)$,
for simplicity of notations.
We also denote by $\mathfrak{Fukst}(-X_1 \times X_2)$,
$\mathfrak{Fukst}(X_1)$, $\mathfrak{Fukst}(X_2)$,
their associated strict categories
(see Definition~\ref{defn2333}\,(8)).\index[syindex]{Fukst@$\mathfrak{Fukst}$}
\end{notation}
By the tri-module analogue of Lemma~\ref{lem55revrev}, the tri-module
$\mathscr{CF}(\mathbb L_1,\mathbb L_{12};\mathbb L_2)$ induces a
left-$\mathfrak{Fukst}(X_1)$, $\mathfrak{Fukst}(-X_1 \times X_2)$
and right-$\mathfrak{Fukst}(X_2)$
filtered $A_{\infty}$ tri-module
$\mathscr{CF}^s(\mathbb L_1,\mathbb L_{12};\mathbb L_2)$.

It can be regarded as a tri-functor
\[
\mathfrak{Fukst}(X_1)^{\rm op}
\times \mathfrak{Fukst}(-X_1 \times X_2)^{\rm op}
\times \mathfrak{Fukst}(X_2)
\to \mathcal{CH}.
\]
By taking opposite functor and using Definition~\ref{lem56} and Lemma~\ref{lem256},\index[syindex]{MWW@$\widehat{\mathcal M\mathcal W \mathcal W}$}
we obtain\footnote{MWW stands for Ma'u--Wehrheim--Woodward.
As we mentioned in the introduction,
Ma'u--Wehrheim--Woodward proved Corollary \ref{cor73} in the case
all the Lagrangian submanifolds involved are embedded and monotone.}
\begin{gather}\label{form61}
\widehat{\mathcal M\mathcal W \mathcal W}\colon\ \mathfrak{Fukst}(-X_1 \times X_2) \to \mathcal{FUNC}(\mathfrak{Fukst}(X_1),\mathcal{FUNC}(\mathfrak{Fukst}(X_2)^{\rm op},\mathcal{CH}^{\rm op})).
\end{gather}

\begin{defn}\label{defn73}
Let
\[
(L_{12},b_{12},\sigma_{12}) = \mathcal L_{12} \in \mathfrak{OB}(\mathfrak{Fukst}(X_{12})),
\qquad
(L_1,b_1,\sigma_1) = \mathcal L_1 \in \mathfrak{OB}(\mathfrak{Fukst}(X_1)).
\]
By \eqref{form61}, we obtain a strict and unital filtered $A_{\infty}$
functor: $ \mathfrak{Fukst}(X_2)^{\rm op} \to
\mathcal{CH}^{\rm op}$.
We denote this functor by \smash{$\widehat{\mathcal W}_{\mathcal L_{12}}(\mathcal L_1)$}, where
$W$ stands for Wehrheim--Woodward. \index[syindex]{WL12@$\widehat{\mathcal W}_{\mathcal L_{12}}$}
We call \smash{$\widehat{\mathcal W}_{\mathcal L_{12}}$}
the {\it correspondence functor} associated to $\mathcal L_{12}$. \index{correspondence functor}

\end{defn}
Let
$\mathcal L_2 = (L_2,\sigma_2,b_2)
\in \mathfrak{OB}(\mathfrak{Fukst}(X_2))$
be the geometric transformation of $\mathcal L_1$ by $\mathcal L_{12}$
in the sense of Definition~\ref{def612}.

We defined
\begin{equation}\label{ophyonfff}
\mathfrak{OpYon}^{\rm op}\colon\ \mathfrak{Fukst}(X_2) \to \mathcal{FUNC}(\mathfrak{Fukst}(X_2)^{\rm op},\mathcal{CH}^{\rm op})
\end{equation}
in Section~\ref{subsec:Yoneda}.
The main result of this section is the following.
\begin{thm}\label{th72}
$\mathcal L_2$ represents \smash{$\widehat{\mathcal W}_{\mathcal L_{12}}(\mathcal L_1)$}
up to homotopy equivalence.
\end{thm}
We will prove Theorem~\ref{th72} in the next subsection.
Corollary \ref{cor73} below says that for each
pair~$(L_{12},b_{12})$ of a Lagrangian submanifold of $-X_1 \times X_2$
and its bounding cochain, we can associate a~filtered $A_{\infty}$ functor $\mathfrak{Fukst}(X_1) \to \mathfrak{Fukst}(X_2)$
in a canonical way.
\begin{cor}\label{cor73}
There exists a strict and unital filtered $A_{\infty}$ functor
\begin{equation}\label{form622}
{\mathcal M\mathcal W \mathcal W}\colon\ \mathfrak{Fukst}(-X_1 \times X_2)
\to \mathcal{FUNC} (\mathfrak{Fukst}(X_1), \mathfrak{Fukst}(X_2) )
\end{equation}
such that its composition with
\begin{gather*}
\mathfrak{OpYon}^{\rm op}_*\colon\
\mathcal{FUNC} (\mathfrak{Fukst}(X_1), \mathfrak{Fukst}(X_2) ) \\
\qquad\to \mathcal{FUNC} (\mathfrak{Fukst}(X_1), \mathcal{FUNC}(
\mathfrak{Fukst}(X_2)^{\rm op},
\mathcal{CH}^{\rm op}) )
\end{gather*}
 is homotopy equivalent to the functor
 \smash{$\widehat{\mathcal M\mathcal W \mathcal W}$} in \eqref{form61}.
 Here $\mathfrak{OpYon}^{\rm op}_*$ is induced by the functor $\mathfrak{OpYon}^{\rm op}$ in
{\rm(\eqref{ophyonfff})}.

\end{cor}
\begin{proof}
$A_{\infty}$-Yoneda lemma (see Theorem~\ref{Yoneda})
implies that there exists a homotopy inverse
\[
(\mathfrak{OpYon}^{\rm op})^{-1} \colon\ \mathfrak{Rep}(
\mathfrak{Fukst}(X_2)^{\rm op},
\mathcal{CH}^{\rm op}) \to \mathfrak{Fukst}(X_2)
\]
to the Yoneda functor $\mathfrak{Yon}$.
(Here $\mathfrak{Rep}$ denotes the full subcategory consisting of objects
which are homotopy equivalent to one in the image of Yoneda functor.
See Definition~\ref{defn241}.)\footnote{The filtered $A_{\infty}$ category,
functor, tri-module etc.\ which are defined by using the moduli space
of pseudo-holomorphic curves are always gapped because of Gromov compactness.}
It induces
\begin{gather*}
\bigl(\bigl(\mathfrak{OpYon}^{\rm op}\bigr)^{-1}\bigr)_* \colon\
 \mathcal{FUNC} (\mathfrak{Fukst}(X_1), \mathfrak{Rep}(
\mathfrak{Fukst}(X_2)^{\rm op},
\mathcal{CH}^{\rm op}) ) \\
\hphantom{\bigl(\bigl(\mathfrak{OpYon}^{\rm op}\bigr)^{-1}\bigr)_* \colon} \
\to
\mathcal{FUNC} (\mathfrak{Fukst}(X_1), \mathfrak{Fukst}(X_2) ).
\end{gather*}
 On the other hand, Theorem~\ref{th72} implies that
 the filtered $A_{\infty}$ functors \smash{$\widehat{\mathcal M\mathcal W \mathcal W}$}
 factor through
\begin{equation}\label{form733}
\mathfrak{Fukst}(-X_1 \times X_2)
\to \mathfrak{Rep} (\mathfrak{Fukst}(X_1), \mathfrak{Fukst}(X_2)) ).
\end{equation}
We compose \eqref{form733} with $\bigl((\mathfrak{OpYon}^{\rm op})^{-1}\bigr)_*$ to
obtain required filtered $A_{\infty}$ functor
${\mathcal M\mathcal W \mathcal W}$.
\end{proof}

\begin{defn}\label{defn7575}
We call the filtered $A_{\infty}$ functor ${\mathcal M\mathcal W \mathcal W}$ in
Corollary \ref{cor73} the {\it correspondence bi-functor}, when we regard it as a
bi-functor
\[
\mathfrak{Fukst}(-X_1 \times X_2) \times \mathfrak{Fukst}(X_1) \to \mathfrak{Fukst}(X_2).
\]

\index{correspondence bi-functor}

For a given unobstructed Lagrangian correspondence $\mathcal L_{12}$,
the correspondence bi-functor induces a filtered $A_{\infty}$ functor\index[syindex]{WL12@${\mathcal W}_{\mathcal L_{12}}$}
$
{\mathcal W}_{\mathcal L_{12}}
\colon \mathfrak{Fukst}(X_1) \to \mathfrak{Fukst}(X_2)$.
We call it
the \index{correspondence functor}{\it correspondence functor}
associated to the unobstructed immersed Lagrangian correspondence $\mathcal L_{12}$.
\end{defn}

\subsection{Proof}
\label{subsec:functorconst}

In this subsection, we prove Theorem~\ref{th72}.
\begin{proof}\label{prop77}
To prove Theorem~\ref{th72}, it suffices to show the next proposition.

\begin{prop}
There exists a natural transformation $\mathscr T$ from $\mathfrak{OpYon}^{\rm op}_{\rm ob}(\mathcal L_{2})$ to $\widehat{\mathcal W}_{\mathcal L_{12}}
(\mathcal L_1)$
which has a homotopy inverse.
\end{prop}
\begin{proof}
We remark
\[
\mathcal{FUNC}(\mathfrak{Fukst}(X_2)^{\rm op},\mathcal{CH}^{\rm op})
\cong \mathcal{FUNC}(\mathfrak{Fukst}(X_2),\mathcal{CH})^{\rm op}.
\]
We regard $\mathfrak{OpYon}_{\rm ob}^{\rm op}(\mathcal L_{2})$ and $\widehat{\mathcal W}_{\mathcal L_{12}}(\mathcal L_1)$ the objects of the right-hand side.

Let $\mathfrak c,\mathfrak c_{0}, \dots, \mathfrak c_{k}$ be objects of $\mathfrak{Fukst}(X_2)$.
We recall that the functor
$\mathfrak{OpYon}_{\rm ob}(\mathcal L_{2})$ for objects is defined by
$
\mathfrak c \mapsto CF(\mathcal L_{2},\mathfrak c)$.
The morphisms part of $\mathfrak{OpYon}^{\rm op}_{\rm ob}(\mathcal L_{2})$ is a map
\[
CF(\mathcal L_{2},\mathfrak c_0) \otimes B_k \mathfrak{Fukst}(X_2)[1](\mathfrak c_0,\mathfrak c_k)
\to CF(\mathcal L_{2},\mathfrak c_k)
\]
defined by
\[
z \otimes (y_1,\dots,y_{k})
\mapsto
\mathfrak m(z,y_1,\dots,y_{k})
\in CF(\mathcal L_2,\mathfrak c_{k}).
\]
Here $z \in CF(\mathcal L_{2},\mathfrak c_0)$,
$y_i \in CF(\mathfrak c_{i-1},\mathfrak c_i)$,
and $\mathfrak m$ is the structure operation
of the filtered $A_{\infty}$ category $\mathfrak{Fukst}(X_2)$.
(We remark that $\mathfrak m$ already includes the deformation by the bounding cochain.)
The Bar complex $B_k\dots$ of an $A_{\infty}$ category is defined in \eqref{form24new}.

On the other hand, the object part of \smash{$\widehat{\mathcal W}_{\mathcal L_{12}}
(\mathcal L_1)$} is
$
\mathfrak c \mapsto CF(L_1,L_{12};\mathfrak c)$.
Here, when the Lagrangian submanifold which is a part of the data in $\mathfrak c$
is $L'_2$, then we put
\[
CF(L_1,L_{12};\mathfrak c): = CF(L_1,L_{12};L'_2),
\]
where the right-hand side is defined in Definition~\ref{defn53939}\,(1).

The morphism part of \smash{$\widehat{\mathcal W}_{\mathcal L_{12}}
(\mathcal L_1)$} is a map
\[
CF(L_1,L_{12};\mathfrak c_0) \otimes B_k \mathfrak{Fukst}(X_2)[1](\mathfrak c_0,\mathfrak c_k)
\to CF(L_1,L_{12};\mathfrak c_k)
\]
and is defined by
\begin{equation}\label{form7107}
w \otimes (y_1,\dots,y_{k})
\mapsto
\mathfrak n(w;y_1,\dots,y_{k})
\in CF(L_1,L_{12};\mathfrak c_k).
\end{equation}
Here $w \in CF(L_1,L_{12};\mathfrak c_0)$,
$y_i \in CF(\mathfrak c_{i-1},\mathfrak c_i)$,
and $\mathfrak n$ is a filtered $A_{\infty}$ right module structure on~$CF(L_1,L_2;\mathfrak c_k)$.\footnote{We remark that
we take an opposite functor while defining \smash{$\widehat{\mathcal M\mathcal W \mathcal W}$}.}
Note that using the notation
$\mathfrak n^{b_1,b_{12}}$ appearing in Lemma~\ref{lem6868},
$\mathfrak n$ is defined~by
\begin{gather}\label{form71088}
\mathfrak n(w;y_1,\dots,y_{k})
=
\mathfrak n^{b_1,b_{12}}
\bigl(w;e^{b_{2,0}}y_1e^{b_{2,1}}\cdots e^{b_{2,k-1}}y_ke^{b_{2,k}}\bigr),
\end{gather}
where $b_{2,i}$ are bounding cochains for $i=1,2$.
Here we denote an object $\mathfrak c_i$ as a pair $(\mathcal L_{2,i},b_{2,i})$
of~${\mathcal L_{2,i} \in \mathbb L_2}$ and its bounding cochain
and $b_{2,i}$. Thus $b_{2,i}$ is a bounding cochain which is a part of data
consisting $\mathfrak c_i$.
The symbol $e^b$ is defined by
\[
e^b = \sum_{k=0}^{\infty} \underbrace{b \otimes \dots \otimes b}_{\text{$k$ times }}.
\]
The operation
\eqref{form7107} is a map
\[
CF(L_1,L_{12};\mathfrak c_0) \otimes B_k \mathfrak{Fukst}(X_2)[1](\mathfrak c_0,\mathfrak c_k)
\to CF(L_1,L_{12};\mathfrak c_k).
\]
See Figure~\ref{FigureSec7new1}.

Now the object part
$
\mathscr T_{\rm ob}(\mathfrak c) \colon CF(L_1,L_{12};\mathfrak c)
\to
CF(\mathcal L_2,\mathfrak c)
$
of $\mathscr T$ is defined by
\begin{equation}\label{formjula711}
\mathscr T_{\rm ob}(\mathfrak c)(z)
=
\mathfrak n({\bf 1};z),
\end{equation}
where $\mathfrak n$ is as in \eqref{form71088} and
${\bf 1} \in CF(L_1,L_{12};L_2)$ is the cyclic element
in Proposition~\ref{prop610}.

The morphism part
\[
\mathscr T_{k}(\mathfrak c_0,\mathfrak c_k)
\colon\
CF(L_1,L_{12};\mathfrak c_0) \otimes B_k\mathfrak{Fukst}[1](X_2)(\mathfrak c_0,\mathfrak c_k)
\to CF(\mathcal L_{2},\mathfrak c_k)
\]
is defined by
\begin{equation}\label{form7.11111}
\mathscr T_{k}(\mathfrak c)(z;y_1,\dots,y_k)
=
\mathfrak n({\bf 1};z,y_1,\dots,y_k).
\end{equation}
See Figure~\ref{Figurenew7-2}.

\begin{figure}[ht]\centering
\begin{tabular}{cc}
\begin{minipage}[t]{0.45\hsize}
\centering
\includegraphics[scale=0.4]{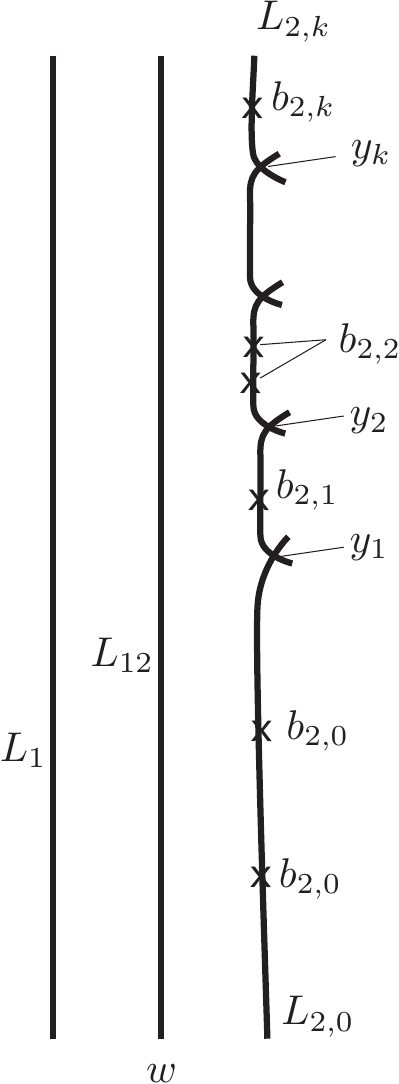}
\caption{$\mathfrak n(w;y_1,\dots,y_{k})$.}
\label{FigureSec7new1}
\end{minipage} &
\begin{minipage}[t]{0.45\hsize}
\centering
\includegraphics[scale=0.4]{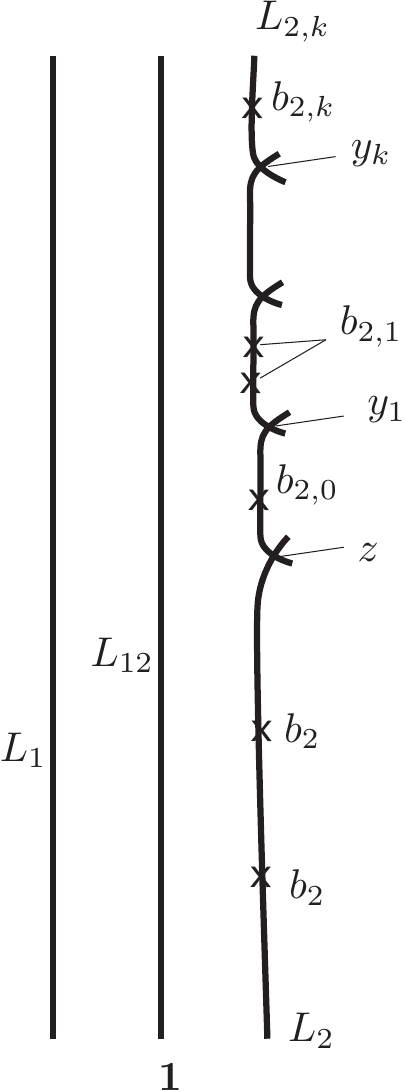}
\caption{$\mathscr T_{k}(\mathfrak c_0,\mathfrak c_k)(z;y_1,\dots,y_k)$.}
\label{Figurenew7-2}
\end{minipage}
\end{tabular}
\end{figure}
\begin{lem}

The maps $\mathscr T$ is a natural transformation.
$($Namely, its boundary in the functor category is $0$.$)$
In other words, it is a filtered right $A_{\infty}$
module homomorphism.
\end{lem}
\begin{proof}
\eqref{eq32} implies that \eqref{form7.11111} is a chain map.
Then the lemma follows from $A_{\infty}$ formula of
$\mathfrak n$.
(See \eqref{form51}. The element ${\bf x}$ there is empty here (that is, $1 \in B_0 CF(L_1,L_1)$).)
\end{proof}

\begin{lem}
$
\mathscr T_{\rm ob}(\mathfrak c)\colon CF(L_1,L_{12};\mathfrak c)
\to
CF(\mathcal L_2,\mathfrak c)
$
is an isomorphism of $\Lambda_0$ modules.
\end{lem}
\begin{proof}
Since ${\bf 1}$ is cyclic, the definition implies that
$\mathscr T_{\rm ob}(\mathfrak c) \mod \Lambda_+$ is an
isomorphism. The lemma then follows easily.
\end{proof}

Now Proposition~\ref{prop77} follows from the next Lemma
\ref{prop1015}.
\end{proof}

\begin{lem}\label{prop1015}
Let $\mathscr C_1$, $\mathscr C_2$ be unital and strict
filtered $A_{\infty}$ categories and $\mathscr{F}$, $\mathscr{G}$
unital and strict filtered $A_{\infty}$ functors from
$\mathscr C_1$ to $\mathscr C_2$.
Let $\mathcal T$ be a natural transformation from
$\mathscr{F}$ to $\mathscr{G}$.
We assume that, for each object $c$ of $\mathscr C_1$,
$\mathcal T_c \in \mathscr C_2(\mathscr{F}(c),\mathscr{G}(c))$
is a homotopy equivalence.
Then $\mathcal T$ is a homotopy equivalence in the functor category.
$($See Definition {\rm\ref{6.22}}.$)$
\end{lem}
This lemma seems to be well-known. For the sake of completeness, we will prove it below.
\begin{proof}
For simplicity of sign, we consider the case when the degree of $\mathcal T$
is $0$. (We use only such cases.)
We use the notation of Proposition~\ref{prop1015}.
We will construct natural transformations $\mathcal S\colon \mathscr{G} \to \mathscr{F}$ of degree $0$
and~${\mathcal H \colon \mathscr{F} \to \mathscr{F}}$ of degree $-1$
such that $\mathfrak M_1(\mathcal H) = \mathfrak M_2(\mathcal S,\mathcal T) -
\mathcal{ID}_{\mathscr F}$, where $\mathcal{ID}_{\mathscr F} \colon \mathscr F
\to \mathscr F$ is the identity natural transformation.
(Here $\mathfrak M_k$ is the structure operation of the functor category.)

We use the induction on the number filtration and will construct
\[
\mathcal S_k \colon\ B_k\mathscr C_1[1](c,c') \to \mathscr C_2(\mathscr{G}_0(c),\mathscr{F}(c')),
\qquad
\mathcal H_k \colon\ B_k\mathscr C_1[1](c,c') \to \mathscr C_2(\mathscr{F}_0(c),\mathscr{F}(c'))
\]
by induction on $k$ so that they satisfy the following conditions \eqref{eq1053}, \eqref{eq1054}
and \eqref{eq1055}.
Suppose~$\mathcal S_i$ is defined for $i\le k$ and
$\mathcal H_i$ is defined for $i\le k$.
We define
\begin{gather*}
\widehat{\mathcal S}_{(k)} \colon\
B\mathscr C_1[1](c,c') \to B\mathscr C_2[1](\mathscr{G}_0(c),\mathscr{F}(c')), \\
\widehat{\mathcal H}_{(k)} \colon\
B\mathscr C_1[1](c,c') \to B\mathscr C_2[1](\mathscr{F}_0(c),\mathscr{F}(c'))
\end{gather*}
by
\begin{gather*}
\widehat{\mathcal S}_{(k)}({\bf x})
=
\sum_c \widehat{\mathscr G}({\bf x}_{c;1}) \otimes \mathcal S_{\le k}({\bf x}_{c;2})
\otimes \widehat{\mathscr F}({\bf x}_{c;3}),
\\
\widehat{\mathcal H}_{(k)}({\bf x})
=
\sum_c (-1)^{\deg' {\bf x}_{c;1}} \widehat{\mathscr F}({\bf x}_{c;1}) \otimes \mathcal H_{\le k}({\bf x}_{c;2})
\otimes \widehat{\mathscr F}({\bf x}_{c;3}),
\end{gather*}
where $((\Delta \otimes {\rm id}) \circ \Delta)({\bf x})
= \sum_c {\bf x}_{c;1} \otimes {\bf x}_{c;2} \otimes {\bf x}_{c;3}$.
Here we define $\mathcal S_{\le k}$ such that it is $\mathcal S_i$
on $B_i\mathscr C_1(c,c')$ with $i\le k$ and is zero otherwise.
$\mathcal H_{\le k}$ is defined in a similar way.

We require
\begin{equation}
\label{eq1053}
\mathfrak m\bigl(\widehat{\mathcal S}_{\le k}({\bf x})\bigr) -
\widehat{\mathcal S}_{\le k}\bigl(\hat d{\bf x}\bigr) = 0
\qquad \text{for}\quad {\bf x} \in B_i\mathscr C_1[1](c,c')\  \text{with} \  i\le k.
\end{equation}
We also require
\begin{gather}
\sum_c \mathfrak m
\bigl(
\widehat{\mathscr F}({\bf x}_{c;1}) \otimes
\mathcal T_{\le k}({\bf x}_{c;2}) \otimes
\widehat{\mathscr G}({\bf x}_{c;3}) \otimes
\mathcal S_{\le k}({\bf x}_{c;4}) \otimes
\widehat{\mathscr F}({\bf x}_{c;5})\bigr)\nonumber
\\
\qquad=
\mathfrak m\bigl(\widehat{\mathcal H}_{\le k}({\bf x})\bigr) +
\widehat{\mathcal H}_{\le k}\bigl(\hat d{\bf x}\bigr)
\label{eq1054}
\end{gather}
for ${\bf x} \in B_i\mathscr C_1(c,c')$ with $0 < i\le k$.
Here
\begin{gather*}
((\Delta \otimes {\rm id}\otimes {\rm id}\otimes {\rm id})
\circ (\Delta \otimes {\rm id}\otimes {\rm id}) \circ (\Delta \otimes {\rm id}) \circ \Delta)({\bf x}) \\
\qquad=
\sum_c {\bf x}_{c;1}\otimes{\bf x}_{c;2}\otimes{\bf x}_{c;3}\otimes
{\bf x}_{c;4}\otimes{\bf x}_{c;5}.
\end{gather*}
Moreover, we require
\begin{equation}\label{eq1055}
\mathfrak m_2
(
\mathcal T_{0}(c) \otimes
\mathcal S_{\le 0}(c))
=
{\bf e}_{\mathscr F_{\rm ob}(c),\mathscr F_{\rm ob}(c)} + \mathfrak m_1(\mathcal H_0(c)).
\end{equation}

Let us start the construction of $\mathcal S_k$ and $\mathcal H_k$
by induction.
We first consider the case $k=0$.
By assumption, $\mathcal T_0(c) \in \mathscr C_2(\mathscr F(c),\mathscr G(c))$
is a homotopy equivalence.
Therefore, there exists
$\mathcal S_0(c) \in \mathscr C_2(\mathscr G(c),\mathscr F(c))$
and $\mathcal H_0(c) \in \mathscr C_2(\mathscr F(c),\mathscr F(c))$
such that
\[
\mathfrak m_1(\mathcal S_0(c)) = 0, \qquad
\mathfrak m_2(\mathcal T_0(c),\mathcal S_0(c)) = {\bf e}_{\mathscr F_{\rm ob}(c),\mathscr F_{\rm ob}(c)} +\mathfrak m_1(\mathcal H_0(c))
\]
We thus obtain required $\mathcal S_0(c)$ and $\mathcal H_0(c)$.

Suppose we have obtained $\mathcal S_i$ and
$\mathcal H_i$ for $i\le k$ such that \eqref{eq1053} and \eqref{eq1054} are satisfied.
We will construct $\mathcal S_{k+1}$ and
$\mathcal H_{k+1}$.

Let ${\bf x} \in B_{k+1}\mathscr C_1(c,c')$.
We put
\[
O({\bf x}) = \mathfrak m\bigl(\widehat{\mathcal S}_{(k)}({\bf x})\bigr) -
\widehat{\mathcal S}_{(k)}\bigl(\hat d{\bf x}\bigr)
\in \mathscr C_2(\mathscr G(c),\mathscr F(c')).
\]
Using \eqref{eq1053}, we can easily check that
\begin{equation}\label{form1056}
\mathfrak m_1\bigl(O({\bf x})\bigr) + O\bigl(\hat d_1({\bf x})\bigr) = 0.
\end{equation}
Here $\hat d_1$ is the coderivation induced by $\mathfrak m_1$.
In fact, \eqref{form1056} follows from $\mathfrak M_1(\mathfrak M_1(\mathcal S_{(k)})) = 0$
and~\eqref{eq1053}.

On the other hand, we use
$
\mathfrak M_1(\mathfrak M_2(\mathcal T,\mathcal S_{(k)}))
=
\mathfrak M_2(\mathcal T,\mathfrak M_1(\mathcal S_{(k)}))
$
together with \eqref{eq1054}, \eqref{eq1055}, and obtain
\begin{equation}\label{form1057}
\mathfrak m_2(\mathcal T_0(c),O({\bf x}))
=
\mathfrak m_1 (B({\bf x})) + B\bigl(\hat d_1{\bf x}\bigr),
\end{equation}
where
\[
B({\bf x})
=
-\sum_c
\mathfrak m
\bigl(
\widehat{\mathscr F}({\bf x}_{c;1}) \otimes
\mathcal T({\bf x}_{c;2}) \otimes
\widehat{\mathscr G}({\bf x}_{c;3}) \otimes
\mathcal S_{\le k}({\bf x}_{c;4}) \otimes
\widehat{\mathscr F}({\bf x}_{c;5})\bigr).
\]
\eqref{form1056}, \eqref{form1057} together with the fact that
$x \mapsto \mathfrak m_2(\mathcal T_0(c),x)$ is a chain homotopy
equivalence: $ \mathscr C_2(c,c') \to \mathscr C_2(c,c)$
imply that there exists
\[
{\mathcal S}'_{k+1} \colon\
B_{k+1}\mathscr C_1[1](c,c') \to \mathscr C_2(\mathscr{G}_0(c),\mathscr{F}(c'))
\]
such that when
we use this ${\mathcal S}'_{k+1}$ for ${\mathcal S}_{k+1}$
to define ${\mathcal S}'_{\le k+1}$, then
\eqref{eq1053} for $k+1$ replaced by $k$ holds.
We also use $0$ for ${\mathcal H}_{k+1}$
to define ${\mathcal H}'_{(k+1)}$.
We then consider
\begin{align*}
E_{(k+1)}({\bf x})
={}
&\sum_c \mathfrak m
\bigl(
\widehat{\mathscr F}({\bf x}_{c;1}) \otimes
\mathcal T_{\le k+1}({\bf x}_{c;2}) \otimes
\widehat{\mathscr G}({\bf x}_{c;3}) \otimes
\mathcal S'_{\le k+1}({\bf x}_{c;4}) \otimes
\widehat{\mathscr F}({\bf x}_{c;5})\bigr) \\
&
- \mathfrak m\bigl(\widehat{\mathcal H}'_{\le k+1}({\bf x})\bigr) -
\widehat{\mathcal H}'_{\le k+1}\bigl(\hat d{\bf x}\bigr)
\end{align*}
By induction hypothesis, $E_{(k+1)}({\bf x}) = 0$
for ${\bf x} \in B_{i}\mathscr C_1(c,c')$ with $0<i\le k$.
We use it and~\eqref{eq1055} to obtain
\smash{$
\mathfrak m_1(E_{(k+1)}({\bf x}))
- E_{(k+1)}\bigl(\widehat d({\bf x})\bigr) = 0
$}
by an easy calculation.
Then we again use the fact
$x \mapsto \mathfrak m_2(\mathcal T_0(c),x)$ is a chain homotopy
equivalence: $ \mathscr C_2(c,c') \to \mathscr C_2(c,c)$
to obtain
${\rm Corr} \colon B_{k+1}\mathscr C_1(c,c') \to \mathscr C_2(\mathscr{G}_0(c),\mathscr{F}(c'))$
and $\mathcal H_{k+1} \colon\allowbreak
B_{k+1}\mathscr C_1(c,c') \to \mathscr C_2(\mathscr{F}_0(c),\mathscr{F}(c'))$
such that
\begin{gather*}
E_{(k+1)}({\bf x}) +
\mathfrak m_2(\mathcal T_0(c),{\rm Corr}({\bf x}))
=
\mathfrak m(\mathcal H_{k+1}({\bf x}))
+
\mathcal H_{k+1}\bigl(\widehat d_1{\bf x}\bigr), \\
 \mathfrak m_1({\rm Corr}({\bf x}))
- {\rm Corr}\bigl(\widehat d_1{\bf x}\bigr) = 0.
\end{gather*}
Then ${\mathcal S}_{k+1} = {\mathcal S}'_{k+1} + E_{(k+1)}$
and the above $\mathcal H_{k+1}$ satisfy
\eqref{eq1053} and \eqref{eq1054} with $k$ replaced by~${k+1}$.

We thus obtained a natural transformation $\mathcal S \colon
\mathscr G \to \mathscr F$ such that
$\mathfrak M_2(\mathscr T,\mathcal S)$ is homotopic to the identity
natural transformation $\mathscr F \to \mathscr F$.

In the same way, we can find $\mathcal S' \colon
\mathscr G \to \mathscr F$ such that
$\mathfrak M_2(\mathcal S',\mathcal T)$ is homotopic to the identity
natural transformation $\mathscr G \to \mathscr G$.
Using associativity of $\mathfrak M_2$ up to homotopy,
it implies that $\mathcal S'$ is homotopic to $\mathcal S$.
Therefore, $\mathcal S'$ is a homotopy inverse to $\mathscr T$.
The proof of Lemma~\ref{prop1015} is now complete.
\end{proof}

The proof of Theorem~\ref{th72} is complete.
\end{proof}

\section{Compositions of Lagrangian correspondences}
\label{sec:comp}

\subsection{Unobstructedness of composed correspondences}

The main result of this subsection is Theorem~\ref{comp} below.

\begin{situ}\label{situ80}
Suppose that $\mathbb L_{1}$, $\mathbb L_{2}$ and $\mathbb L_{12}$ are
as in Situation \ref{situ61}.
We also assume that $\mathbb L_{2}$, $\mathbb L_{3}$ and $\mathbb L_{23}$ are
as in Situation \ref{situ61}.

For $(L_{12},\sigma_{12}) \in \mathbb L_{12}$ and $(L_{23},\sigma_{23}) \in
 \mathbb L_{23}$,
we assume that $\pi_{X_2}\circ i_{L_{12}} \colon \tilde L_{12} \to X_2$
is transversal to $\pi_{X_2}\circ i_{L_{23}} \colon \tilde L_{23} \to X_2$
and put
\begin{equation}\label{form81}
\tilde L_{13} = \tilde L_{12} \times_{X_2} \tilde L_{23}.
\end{equation}

\end{situ}
Together with $\tilde L_{13} \to -X_{1} \times X_3$ it becomes an
immersed Lagrangian submanifold $L_{13}$ of~${-X_1 \times X_3}$.
We assume that $L_{13}$ has clean self-intersection.
We remark that $L_{13}$ is
$(\pi_1^*(V_1 \oplus TX_1)\times \pi_3^*(V_3))$-relatively spin by
Definition--Lemma~\ref{defnlen47}.
\begin{thm}\label{comp}
There exists a $(\pi_1^*(V_1 \oplus TX_1)\times \pi_3^*(V_3))$-relatively spin
structure $\sigma_{13}$ of $L_{13}$ with the following properties.
Suppose that $b_{12}$ and $b_{23}$ are bounding cochains of $(L_{12},\sigma_{12})$
and~${(L_{23},\sigma_{23})}$, respectively.
Then there exists a bounding cochain $b_{13}$ of $(\mathcal L_{13},\sigma_{13})$.
Moreover, there is a canonical way to determine $b_{13}$ from $b_{12}$ and $b_{23}$
up to gauge equivalence.
\end{thm}

We can enhance the map $(L_{12},b_{12}), (L_{23},b_{23})
\mapsto (L_{13},b_{13})$ to an $A_{\infty}$ functor as in Theorem~\ref{comp2} below.

\begin{situ}\label{situ82}
\quad
\begin{enumerate}\itemsep=0pt
\item[(1)]
Suppose that $\mathbb L_{1}$, $\mathbb L_{2}$, $\mathbb L_{3}$ and $\mathbb L_{12}$, $\mathbb L_{23}$ are
as in Situation \ref{situ80}.
We also assume $\mathbb L_{1}$, $\mathbb L_{3}$ and $\mathbb L_{13}$ are
as in Situation \ref{situ61}.
\item[(2)]
Moreover, we assume the following.
Let $(L_{12},\sigma_{12}) \in \mathbb L_{12}$,
$(L_{23},\sigma_{23},b_{23}) \in \mathbb L_{23}$.
The fiber product $L_{13}$ as in \eqref{form81} together with $\sigma_{13}$
in Theorem~\ref{comp} gives a pair $(L_{13},\sigma_{13})$.
We require that $(L_{13},\sigma_{13})$ is an element of $\mathbb L_{13}$.
\end{enumerate}

\end{situ}
\begin{notation}
\quad
\begin{enumerate}\itemsep=0pt
\item[(1)]
In Situation \ref{situ82}, we write
$
(L_{13},\sigma_{13}) = (L_{23},\sigma_{23}) \circ (L_{12},\sigma_{12})
$
and call $(L_{13},\sigma_{13})$ the {\it geometric composition} of \index{geometric composition}
$(L_{23},\sigma_{23})$ and $(L_{12},\sigma_{12})$.
\item[(2)]
Suppose that $b_{12}$ and $b_{23}$ are bounding cochains of
$(L_{12},\sigma_{12})$ and $(L_{23},\sigma_{23})$, respectively.
Then by Theorem~\ref{comp}, we obtain a bounding cochain $b_{13}$ of $(L_{13},\sigma_{13})$.
We put
\begin{equation}\label{compobkj}
(L_{13},\sigma_{13},b_{13}) = (L_{23},\sigma_{23},b_{23}) \circ (L_{12},\sigma_{12},b_{12}).
\end{equation}
\item[(3)]
Let $\mathfrak{Fuk}(-X_1 \times X_2)$, $\mathfrak{Fuk}(-X_2 \times X_3)$, $\mathfrak{Fuk}(-X_1 \times X_3)$
be the filtered $A_{\infty}$ categories obtained in Theorem~\ref{AJtheorem},
the set of whose objects are $\mathbb {L}_{12}$, $\mathbb {L}_{23}$, $\mathbb {L}_{13}$, respectively.
We denote by $\mathfrak{Fukst}(-X_1 \times X_2)$, $\mathfrak{Fukst}(-X_2 \times X_3)$, $\mathfrak{Fukst}(-X_1 \times X_3)$
the associated strict categories.
\end{enumerate}

\end{notation}

\begin{thm}\label{comp2}
In Situation {\rm\ref{situ82}},
there exists a strict, unital and gapped filtered $A_{\infty}$ bi-functor
\begin{equation}\label{form91}
\mathfrak{Comp} \colon\ \mathfrak{Fukst}(-X_1 \times X_2)
\times \mathfrak{Fukst}(-X_2 \times X_3)
\to
\mathfrak{Fukst}(-X_1 \times X_3)
\end{equation}
such that its object part $\mathfrak{Comp}_{\rm ob}$ is the map given by
{\rm\eqref{compobkj}}.
\end{thm}
\begin{rem}
In the case when all the Lagrangian submanifolds involved are
embedded and monotone, Theorem~\ref{comp2} was proved by
Ma'u--Wehrheim--Woodwards in \cite{MWW}.
\end{rem}
\begin{proof}
The proofs of both Theorems~\ref{comp} and \ref{comp2} are similar to the proof of Theorem~\ref{thm61}, Corollary \ref{cor73} and
use tri-module and Proposition~\ref{thm35}.
Namely, we use the next result.
\begin{prop}\label{pro86}
In Situation {\rm\ref{situ82}}, there exists a
left-$\mathfrak{Fuk}(-X_1 \times X_3)$
and right-$\mathfrak{Fuk}(-X_1 \times X_2)$, $\mathfrak{Fuk}(-X_2 \times X_3)$
filtered $A_{\infty}$ tri-module
$\mathscr{CF}(\mathbb L_{13};\mathbb L_{12},\mathbb L_{23})$.\index[syindex]{CscrFL13L12L23@$\mathscr{CF}(\mathbb L_{13};\mathbb L_{12},\mathbb L_{23})$}

\end{prop}
The proof is similar to the proof of Theorem~\ref{trimain} and is given
in the next subsection.
We remark however that `left' and `right' appear in the opposite way
in Proposition~\ref{pro86} compared to Theorem~\ref{trimain}.
The reason will become clear when we discuss the $Y$-diagram in Section~\ref{sec:comptibility}.

We now prove Theorem~\ref{comp} assuming Proposition~\ref{pro86}.
Suppose we are in the situation of Theorem~\ref{comp}.
We define $L_{13}$ as in \eqref{form81}.
For each relative spin structure $\sigma_{13}$ of $L_{13}$,
the tri-module in Proposition~\ref{pro86} associates a
$\Lambda_0$
module
$
CF((L_{13},\sigma_{13});(L_{12},\sigma_{12}),(L_{23},\sigma_{23}))$.
We denote it by~$CF(L_{13};L_{12},L_{23})$ for simplicity.
\begin{lem}\label{lem87}
There exists a unique choice of $\sigma_{13}$ such that $CF(L_{13};L_{12},L_{23})$
is isomorphic to~\smash{$
\Omega\bigl(\tilde L_{13} \times_{X_1\times X_3} \tilde L_{13};\R\bigr) \,\widehat{\otimes}_{\R}\,
\Lambda_0
$}
on the diagonal component $\tilde L_{13}$.
\end{lem}
The proof is given at the end of Section~\ref{subsec:unobcomp}.

We define
\[
\mathfrak n_k \colon\ CF(L_{13})^{\otimes k} \otimes CF(L_{13};L_{12},L_{23})
\to CF(L_{13};L_{12},L_{23})
\]
by
\[
\mathfrak n_k(x_1,\dots,x_k;y) =
\sum_{k_{12}=0}^{\infty}\sum_{k_{23}=0}^{\infty}
\mathfrak n_{k;k_{12},k_{23}}(x_1,\dots,x_k;y;b_{12},\dots,b_{12};b_{23},\dots,b_{23}),
\]
where
$\mathfrak n_{k,k_{12},k_{23}}$ is a structure operation of the tri-module
of Proposition~\ref{pro86}.
\begin{lem}\label{lem88}
$\{\mathfrak n_k \mid k=0,1,2,\dots\}$ defines a structure of left filtered $A_{\infty}$ module on $CF(L_{13};L_{12},\allowbreak L_{23})$
over the filtered $A_{\infty}$ algebra $CF(L_{13})$.
\end{lem}
The proof is a straightforward calculation using Proposition~\ref{pro86}.

We remark that we can define the notion of a cyclic element
for a left filtered $A_{\infty}$ module and Proposition~\ref{thm35}
holds in the case of left filtered $A_{\infty}$ modules.
In fact, a left $\mathscr C$ module $D$ becomes a
right $\mathscr C^{\rm op}$ module, and
the Maurer--Cartan equation of $\mathscr C^{\rm op}$ is the same
as that of $\mathscr C$.
\begin{lem}\label{lem89}
We may take our tri-module structure so that the element
\[
{\bf 1} \in \Omega^0\bigl(\tilde L_{13}\bigr) \subset \Omega\bigl(\tilde L_{13} \times_{X_1\times X_3} \tilde L_{13};\R\bigr) \,\widehat{\otimes}_{\R}\,
\Lambda_0
\cong CF(L_{13};L_{12},L_{23})
\]
is a cyclic element of the left filtered $A_{\infty}$ module $CF(L_{13};L_{12},L_{23})$
in Lemma {\rm\ref{lem88}}.
\end{lem}
The proof is given at the end of Section~\ref{subsec:unobcomp}.

Now we use Proposition~\ref{thm35} to find uniquely
a bounding cochain $b_{13}$ of $L_{13}$ such that
\begin{equation}\label{newform87}
\mathfrak n^{b_{13}}({\bf 1}) = 0.
\end{equation}
By using Proposition~\ref{prop615615}, we can show that gauge equivalence class
of the bounding cochain~$b_{13}$ depends only on those of $b_{12}$ and $b_{23}$,
when the filtered $A_{\infty}$ tri-module $\mathscr{CF}(\mathbb L_{13};\mathbb L_{12},\mathbb L_{23})$ is given.
The independence of the choices to define $\mathscr{CF}(\mathbb L_{13};\mathbb L_{12},\mathbb L_{23})$
is Theorem~\ref{them1431} in Section~\ref{sec:independence2}.

We have proved Theorem~\ref{comp} assuming several results postponed
to later subsections.

We turn to the proof of Theorem~\ref{comp2}.
The proof is similar to Section~\ref{sec:represent}.
By Proposition~\ref{pro86}, we obtain a
strict and unital filtered $A_{\infty}$ bi-functor
\[
\mathscr F^{\rm bi} \colon\
\mathfrak{Fukst}(-X_1 \times X_2) \times \mathfrak{Fukst}(-X_2 \times X_3)
 \to
\mathcal{FUNC}(\mathfrak{Fukst}(-X_1 \times X_3)^{\rm op}, \mathcal{CH}).
\]

Let $\mathcal L_{12} = (L_{12},\sigma_{12},b_{12})$,
$\mathcal L_{23} = (L_{23},\sigma_{23},b_{23})$
be objects of $\mathfrak{Fukst}(-X_1 \times X_2) $ and
$ \mathfrak{Fukst}(-X_2 \times X_3)$, respectively.
By Lemma~\ref{lem87} and \eqref{newform87}, we obtain
$\mathcal L_{13} = (L_{12},\sigma_{13},b_{13})$
which is an object of~${\mathfrak{Fukst}(-X_1 \times X_3)}$.

Let $\mathscr C$ be a strict filtered $A_{\infty}$
category. Then there exists a
filtered $A_{\infty}$ functor
$\mathfrak{Yon}\colon \mathscr C \to \mathcal{FUNC}(\mathscr C^{\rm op},\mathcal{CH})
\cong \mathcal{BIMOD}(\mathscr C,\Lambda_0)$,
from $\mathscr C$ to the category of left-$\mathscr C$
modules such that its object part is
$c \mapsto (b \mapsto \mathscr C(b,c))$.

\begin{prop}\label{prop810}
$\mathscr F^{\rm bi}_{\rm ob}(\mathcal L_{12},\mathcal L_{23})$
is homotopy equivalent to $(\mathfrak{Yon})_{\rm ob}(\mathcal L_{13})$
as filtered $A_{\infty}$ functors:
$ \mathfrak{Fukst}(-X_1 \times X_3) \to \mathcal{CH}$.

\end{prop}
\begin{proof}
The proof is similar to the proof of Theorem~\ref{th72}.
We repeat the proof for completeness.

We denote by
\[
\mathscr F^{\rm tri}\colon\ \mathfrak{Fukst}(-X_1 \times X_3)^{\rm op}
\times \mathfrak{Fukst}(-X_1 \times X_2) \times \mathfrak{Fukst}(-X_2 \times X_3) \to
\mathcal{CH}
\]
the strict tri-functor
associated to $\mathscr F^{\rm bi}$.

Let $\mathcal L_{13}^{(i)}$, $i=0,\dots, m$, be objects of
$\mathfrak{Fukst}(-X_1 \times X_3)$.
We define
\[
\mathscr T_m \colon\
\bigotimes_{i=1}^m
CF\bigl(\mathcal L_{13}^{(i-1)},\mathcal L_{13}^{(i)}\bigr)
\otimes
CF\bigl(\mathcal L_{13}^{(m)},\mathcal L_{13}\bigr)
\to
CF\bigl(\mathcal L_{13}^{(0)};\mathcal L_{12},\mathcal L_{23}\bigr)
\]
by the next formula
\begin{equation}\label{form87}
\mathscr T_m(x_1,\dots,x_m;y)
=
\mathscr F^{\rm tri}_{0,0,m+1}(x_1,\dots,x_m;y;\varnothing,\varnothing;{\bf 1}).
\end{equation}
Note that
\[
x_1 \otimes \dots \otimes x_m \otimes y \in B_{m+1}\mathfrak{Fukst}(-X_1 \times X_3)\bigl(\mathcal L_{13},\mathcal L_{13}^{(m)}\bigr)
\]
and ${\bf 1} \in CF(L_{13};L_{12},L_{23})$. So
the right-hand side of \eqref{form87} is defined by Proposition~\ref{pro86}.
\begin{lem}\label{lem811}
\eqref{form87} defines a natural transformation $\mathscr T = \{\mathscr T_m \mid m=0,1,2,\dots\}$ from
$\mathscr F^{\rm bi}_{\rm ob}(\mathcal L_{12},\mathcal L_{23})$
to $\mathfrak{Yon}(\mathcal L_{13})$.
\end{lem}
\begin{proof}
Using the fact that ${\bf 1}$ is a cycle in $CF(\mathcal L_{13};\mathcal L_{12},
\mathcal L_{23})$,
the lemma is an immediate consequence of Proposition~\ref{pro86}.
\end{proof}

\begin{lem}\label{lem822}
$
\mathscr T_0 \colon
CF\bigl(\mathcal L_{13}^{(0)},\mathcal L_{13}\bigr)
\to
CF\bigl(\mathcal L_{13}^{(0)};\mathcal L_{12},\mathcal L_{23}\bigr)
$
is an isomorphism of $\Lambda_0$ module.
\end{lem}
\begin{proof}
Using the fact that ${\bf 1}$ is a cyclic element, we can easily show that
$\mathscr T_0$ becomes an isomorphism modulo $\Lambda_+$.
Therefore, $\mathscr T_0$ itself is also an isomorphism.
(We used $G$-gappedness here.
In fact, we construct the inverse by induction on energy filtration.
This induction works when the set of exponents of $T$ appearing
in the operations is discrete.)
\end{proof}

By Lemmas~\ref{lem811} and \ref{lem822}, we can use Lemma~\ref{prop1015}
to show that $\mathscr T$ is a homotopy equivalence.
The proof of Proposition~\ref{prop810} is complete.

Using Proposition~\ref{prop810} and $A_{\infty}$ Yoneda lemma,
we can prove Theorem~\ref{comp2} in the same way as
Corollary \ref{cor73}.
\end{proof}

\subsection{Construction of a tri-module}
\label{subsec:unobcomp}

In this subsection, we prove Proposition~\ref{pro86} and complete the proof of
Theorems~\ref{comp} and~\ref{comp2}.
The proof of Proposition~\ref{pro86} is based on a
moduli space of pseudo-holomorphic maps from a cylinder, which we describe below.

By the same trick as Section~\ref{subsec:Ainfcatim}, it suffices to consider the case when $\mathbb L_{12}$,
$\mathbb L_{23}$, $\mathbb L_{13}$ consist of single elements $\mathcal L_{12}
= (L_{12},\sigma_{12})$, $\mathcal L_{23}= (L_{23},\sigma_{23})$,
$\mathcal L_{13} = (L_{13},\sigma_{13})$, respectively.
We consider the cylinder
\begin{equation}\label{form89}
W = S^1 \times \R = [0,3]/{\sim} \times \R.
\end{equation}
Here $\sim$ identifies $0 \in [0,3]$ with $3 \in [0,3]$.
We define $W_1$, $W_2$, $W_3$ by
\begin{equation}\label{eq8100}
W_1 = [0,1] \times \R \subset W, \qquad
W_2 = [1,2] \times \R \subset W, \qquad
W_3 = [2,3] \times \R \subset W
\end{equation}
and also put
$
S_{(i-1) i} = \{i\}\times \R = W_{i-1} \cap W_i $, $i=1,2,3$.
(Here $S_{01} = S_{31}$, $W_0 = W_3$ by convention.)
Note that $\partial W_1 = S_{31} \cup S_{12}$ etc.
See Figure~\ref{Figure81}.
\begin{figure}[ht]
\centering
\includegraphics[scale=0.3]{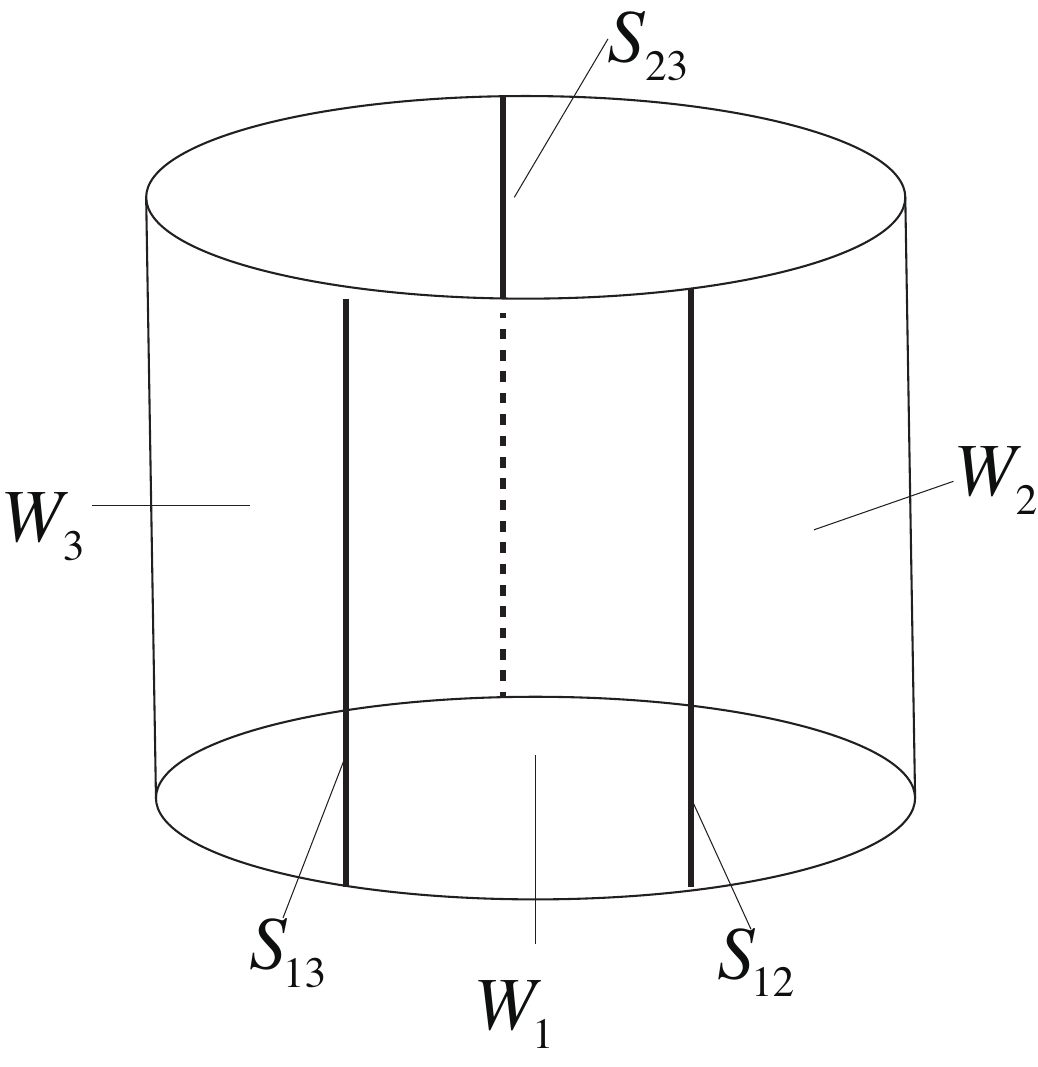}
\caption{Quilted drum $W$.}
\label{Figure81}
\end{figure}
We call $S_{12}$, $S_{23}$, $S_{31}$ the {\it seams}.\index{seam}
We decompose
\begin{gather}
\tilde L_{12} \times_{X_1 \times X_2} \tilde L_{12} = \bigcup_{a \in \mathcal A_{12}}L_{12}(a), \qquad
\tilde L_{23} \times_{X_2 \times X_3} \tilde L_{23} = \bigcup_{a \in \mathcal A_{23}}L_{23}(a), \nonumber\\
\tilde L_{13} \times_{X_1 \times X_3} \tilde L_{13} = \bigcup_{a \in \mathcal A_{13}}L_{13}(a),\nonumber \\
 \bigl(\tilde L_{12} \times \tilde L_{23} \times
\tilde L_{13}\bigr) \times_{X_1^2 \times X_2^2 \times X_3^2} \Delta= \bigcup_{a \in \mathcal A_{123}}R_{123}(a),\label{form812sss}
\end{gather}
where $\Delta$ in the fourth line is the diagonal
$X_1 \times X_2 \times X_3 \subset X_1^2 \times X_2^2 \times X_3^2$
(see Definition~\ref{def3131}\,(5)).\footnote{In \eqref{form812sss}, $\mathcal A_{12}$ etc.
contains the index of the diagonal component. So it corresponds to $\mathcal A_L^+$
in Definition~\ref{def3131}\,(5).}
Let $\vec a_{ii'} = (a_{ii',1},\dots,a_{ii',k_{ii'}}) \in (\mathcal A_{ii'})^{k_{ii'}}$
for $ii' = 12,23,13$.
We call $W$ the {\it quilted drum}. \index{quilted drum}

We define the moduli space
$\overset{\ \text{\tiny $\circ\circ$}}{\mathcal M}_{\rm DR}(\vec a_{12},\vec a_{23},\vec a_{13};a_-,a_+;E)$
for $a_-, a_+ \in \mathcal A_{123}$, $E \in [0,\infty)$ as follows.
\begin{rem}
In the case when $X_1$ is a point, this moduli space
is mostly the same as the one we used in Section~\ref{subsec:bi-functorgeo1}.
In this paper, the role of Lagrangian submanifolds of $X_i$ and of
$-X_i \times X_j$ are much different. The former gives an object of a
filtered $A_{\infty}$ category $\mathfrak{Fuk}(X_i)$,
the latter gives a filtered $A_{\infty}$ functor $\mathfrak{Fuk}(X_i)
\to \mathfrak{Fuk}(X_j)$. By this reason, we use different names and
notations to
those moduli spaces.

\end{rem}
\begin{defn}\label{def916}
We consider
$(\Sigma;\vec z_{12},\vec z_{23},\vec z_{13};u_1,u_2,u_3;\gamma_{1},\gamma_2,\gamma_{3})$
with the following properties (see Figure~\ref{Figure82}):
\begin{enumerate}\itemsep=0pt
\item[(1)]
The space $\Sigma$ is a bordered Riemann surface which is a union of $W$ and
trees of sphere components attached to $W$.
The roots of the trees of sphere components are not on
the seams $S_{12}$, $S_{23}$, $S_{13}$.
\item[(2)]
We denote by $\Sigma_{1}$ the union of $W_{1}$ together with
trees of sphere components rooted on~$W_{1}$.
We define $\Sigma_{2}$, $\Sigma_{3}$ in the same way.
The map
$u_i \colon \Sigma_i \to X_i$ is $-J_{X_i}$ holomorphic for $i=1,2,3$.\footnote
{The reason we consider $-J_{X_i}$ holomorphic maps and not
$J_{X_i}$ holomorphic maps will be explained in Remark~\ref{newrem94}.}
\item[(3)]
$\vec z_{ii'} = (z_{ii',1},\dots,z_{ii',k_{ii'}})$, $ii'=12,23,13$, and
$z_{ii',j} \in S_{ii'}$.
We put $\vert \vec z_{ii'}\vert = \{z_{ii',1},\dots,\allowbreak z_{ii',k_{ii'}}\}$.
\item[(4)]
The maps
$\gamma_{ii'} \colon S_{ii'} \setminus \vert\vec z_{ii'}\vert \to \tilde L_{ii'}$
are smooth and satisfies
$
i_{L_{ii'}} (\gamma_{ii'}(z)) = (u_i(z),u_{i'}(z))$.
When we identify $S_{i i'} \cong \R$ we require $z_{i,i';j} < z_{i,i';j'}$ for $j < j'$
and $(i,i') = (1,2)$ or $(2,3)$ and $z_{13;j} > z_{13;j'}$ for $j < j'$.\footnote{
See Remark~\ref{enumerationremarkkk} for this enumeration.}
\item[(5)]
At $\vec z_{ii'}$, the map $\gamma_{ii'}$ satisfies the switching condition
\begin{equation}\label{form814}
\bigl(\lim_{z \in S_{ii'} \uparrow z_{ii',j}}\gamma_{ii'}(z),\lim_{z \in S_{ii'} \downarrow z_{ii',j}}\gamma_{ii'}(z)
\bigr)
\in L_{ii'}(a_{ii',j})
\end{equation}
for $(i,i') = (1,2),(2,3)$
and
\[
\bigl(\lim_{z \in S_{ii'} \downarrow z_{ii',j}}\gamma_{ii'}(z),\lim_{z \in S_{ii'} \uparrow z_{ii',j}}\gamma_{ii'}(z)
\bigr)
\in L_{ii'}(a_{ii',j})
\]
for $(i,i') = (1,3)$.
Here we identify $S_{ii'} \cong \R$ and then $ \uparrow $, $ \downarrow$ have obvious meaning.
\item[(6)]
When $z \in S^1\times \R$ with $\pi_2(z) \to \pm \infty$,
the maps $u_1(z)$, $u_2(z)$, $u_3(z)$
satisfy the asymptotic boundary condition Condition \ref{cond814} below.
(Here $\pi_2 \colon S^1 \times \R \to \R$ is the projection to the second factor.)
\item[(7)]
The stability condition, Definition~\ref{defn815}\,(2) below, is satisfied.
\item[(8)]
$
\int_{\Omega_1}u_1^*\omega_1 + \int_{\Omega_2}u_1^*\omega_2
+ \int_{\Omega_3}u_3^*\omega_3 = - E$.
We remark that the left-hand side is non-positive since
$u_i$ is $-J_{X_i}$ holomorphic.
\end{enumerate}

 We will define an equivalence relation $\sim$ between objects
$(\Sigma;\vec z_{12},\vec z_{23},\vec z_{13};u_1,u_2,u_3;\gamma_{1},\gamma_2,\gamma_{3})$
which satisfy Conditions (1)--(8), in Definition~\ref{defn815}\,(3). We denote
the set of all the equivalence classes of this equivalence relation by
\smash{$\overset{\ \text{\tiny $\circ\circ$}}{\mathcal M}_{\rm DR}(\vec a_{12},\vec a_{23},\vec a_{13};a_-,a_+;E)$}.
\index[syindex]{M1Dra12a23@$\overset{\ \text{\tiny $\circ\circ$}}{\mathcal M}_{\rm DR}(\vec a_{12},\vec a_{23},\vec a_{13};a_-,a_+;E)$}
We call its element (or an element of its compactification) a {\it pseudo-holomorphic drum}.\index{pseudo-holomorphic drum}
\end{defn}

\begin{figure}[ht]
\centering
\includegraphics[scale=0.3]{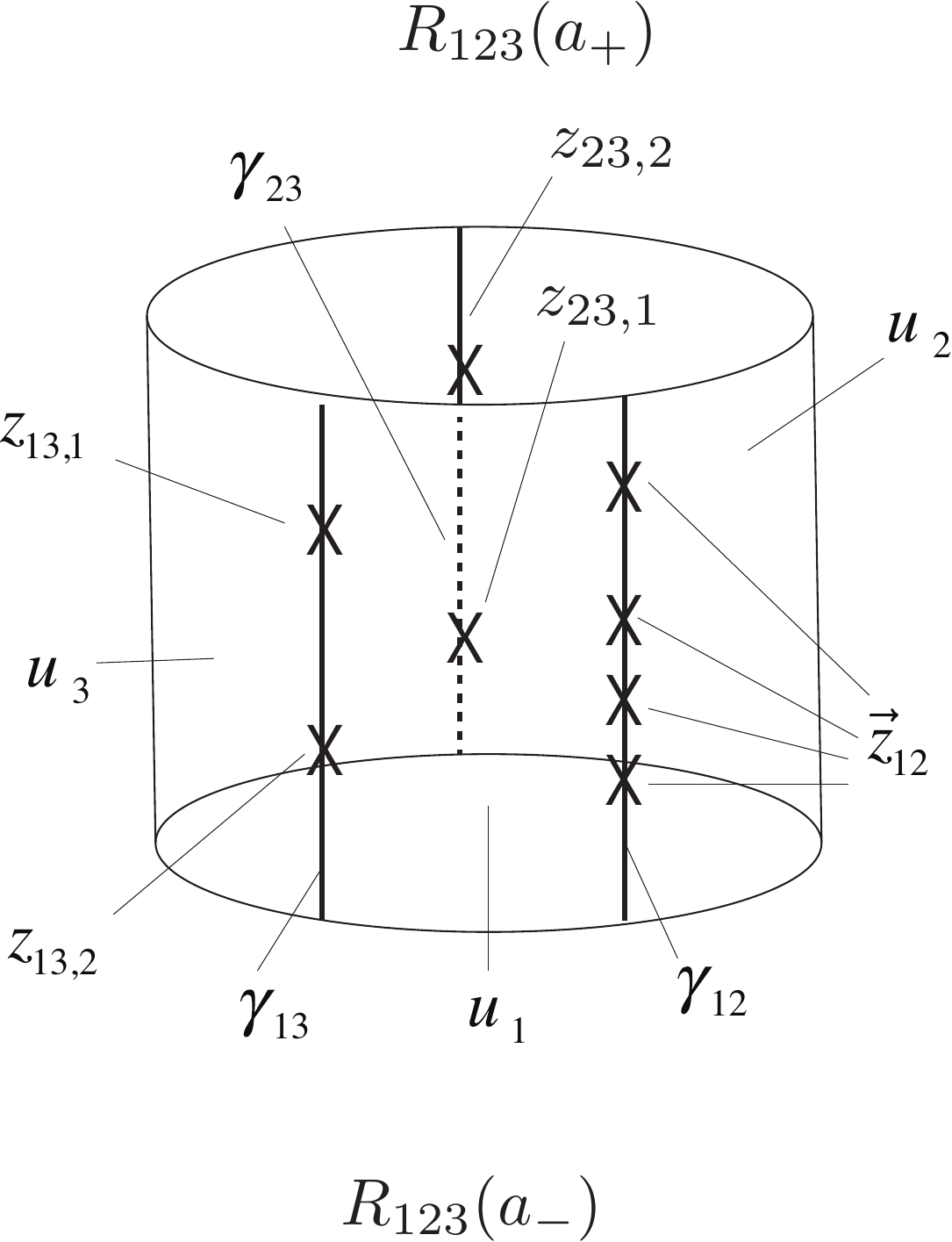}
\caption{An element of $\overset{\ \text{\tiny $\circ\circ$}}{\mathcal M}_{\rm DR}(\vec a_{12},\vec a_{23},\vec a_{13};a_-,a_+;E)$.}
\label{Figure82}
\end{figure}

\begin{rem}\label{enumerationremarkkk}
We enumerate
$z_{12,j}$ and $z_{23,j}$ upward and $z_{13,j}$ downward.
Therefore, we obtain a left-$\mathfrak{Fuk}(-X_1 \times X_3)$
and right-$\mathfrak{Fuk}(-X_1 \times X_2)$, $\mathfrak{Fuk}(-X_2 \times X_3)$
filtered $A_{\infty}$ tri-module by the same reason as explained
in Remark~\ref{rem1521}.

\end{rem}
\begin{conds}\label{cond814}
The asymptotic boundary condition for $\pi_2(z) \to -\infty$ is as follows.
\begin{enumerate}\itemsep=0pt
\item[(1)]
We require the limit $\lim_{\tau \to -\infty}u_1(t,\tau)$ exists and is independent of
$t \in [0,1]$. We write this limit
$\lim_{\pi_2(z) \to -\infty} u_1(z)$.
We require $\lim_{\pi_2(z) \to -\infty} u_2(z)$, $\lim_{\pi_2(z) \to -\infty} u_3(z)$
exist in a~similar sense.
\item[(2)]
\[
\bigl(\lim_{\pi_2(z) \to -\infty} u_1(z), \lim_{\pi_2(z) \to -\infty} u_2(z),
\lim_{\pi_2(z) \to -\infty} u_3(z)
\bigr) \in R_{123}(a_-).
\]
\end{enumerate}
The asymptotic boundary condition for $\pi_2(z) \to +\infty$ is defined in the same way
using $R_{123}(a_+)$.

\end{conds}
\begin{defn}\label{defn815}
Let
\[
\mathfrak x =(\Sigma;\vec z_{12},\vec z_{23},\vec z_{13};u_1,u_2,u_3;\gamma_{1},\gamma_2,\gamma_{3}),\qquad
\mathfrak x' = (\Sigma';\vec z^{\,\prime}_{12},\vec z^{\,\prime}_{23},\vec z^{\,\prime}_{13};u'_1,u'_2,u'_3;\gamma'_{1},\gamma'_2,\gamma'_{3})
\]
be objects satisfying Definition~\ref{def916}\,(1)--(6).
\begin{enumerate}\itemsep=0pt
\item[(1)]
An isomorphism from $\mathfrak x$ to $\mathfrak x'$ is a map $v \colon \Sigma \to \Sigma'$ such that
\begin{enumerate}\itemsep=0pt
\item It is biholomorphic.
\item It sends $\Sigma_i$ to $\Sigma'_i$.
\item It sends $\vec z_{ii'}$ to $\vec z^{\,\prime}_{ii'}$.
\item $u'_i \circ v = u_i$, $\gamma'_{ii'} \circ v = \gamma_{ii}$.
\end{enumerate}
\item[(2)]
$\mathfrak x$ is said to be {\it stable} \index{stable} if the set of
all isomorphisms from $\mathfrak x$ to $\mathfrak x$ is finite.
\item[(3)]
We say $\mathfrak x$ is {\it equivalent} to $\mathfrak x'$ if
there exists an isomorphism from $\mathfrak x$ to $\mathfrak x'$.
\end{enumerate}

\end{defn}
We define evaluation maps
\begin{equation}\label{form815}
{\rm ev}_{ii'} = ({\rm ev}_{ii',1},\dots,{\rm ev}_{ii',k_{ii'}})
\colon\
\overset{\ \text{\tiny $\circ\circ$}}{\mathcal M}_{\rm DR}(\vec a_{12},\vec a_{23},\vec a_{13};a_-,a_+;E)
\to \prod_{j=1}^{k_{ii'}} L_{ii'}(a_{ii',k})
\end{equation}
by the left-hand side of \eqref{form814}.

We also define
\begin{equation}\label{form816}
{\rm ev}_{\infty} = ({\rm ev}_{\infty,+},{\rm ev}_{\infty,-})\colon\
\overset{\ \text{\tiny $\circ\circ$}}{\mathcal M}_{\rm DR}(\vec a_{12},\vec a_{23},\vec a_{13};a_-,a_+;E)
\to R_{123}(a_+) \times R_{123}(a_-)
\end{equation}
by the left-hand side of Condition \ref{cond814}\,(2).
\begin{prop}\label{prop811}
We can define a topology on \smash{$\overset{\ \text{\tiny $\circ\circ$}}{\mathcal M}_{\rm DR}(\vec a_{12},\vec a_{23},\vec a_{13};a_-,a_+;E)$}
such that it has a~compactification
${\mathcal M}_{\rm DR}(\vec a_{12},\vec a_{23},\vec a_{13};a_-,a_+;E)$, \index[syindex]{M1DRa12a23@${\mathcal M}_{\rm DR}(\vec a_{12},\vec a_{23},\vec a_{13};a_-,a_+;E)$}
which is a compact metrizable space.
They have Kuranishi structures with corners and enjoy the following properties:
\begin{enumerate}\itemsep=0pt
\item[$(1)$]
The normalized boundary of ${\mathcal M}_{\rm DR}(\vec a_{12},\vec a_{23},\vec a_{13};a_-,a_+;E)$ is a disjoint union of $2$ types of
fiber products which we describe below.
\item[$(2)$]
The evaluation maps \eqref{form815} and \eqref{form816}
extend to strongly smooth maps with respect to this Kuranishi structure.
${\rm ev}_{\infty,+}$ is weakly submersive.
The extension is compatible with the description of the boundary
in item~{\rm (1)}.
\item[$(3)$]
The orientation bundle of ${\mathcal M}_{\rm DR}(\vec a_{12},\vec a_{23},\vec a_{13};a_-,a_+;E)$ is isomorphic to the tensor product of
the pullbacks of $\Theta^-$ by
the evaluation maps \eqref{form815} and \eqref{form816}. For the component~$R_{123}(a_{+})$, we take
$\Theta^+$ in place of $\Theta^-$.
\item[$(4)$]
It is compatible with the forgetful map of the marked points corresponding to the
diagonal components in the sense of {\rm\cite[\emph{Definition} 3.1]{fooo091}}.
\end{enumerate}
\end{prop}

We describe the boundary components:
\begin{enumerate}\itemsep=0pt
\item[(I)]
The first type of boundary corresponds to the bubble at one of the Lagrangian boundary conditions $L_{12}$,
$L_{23}$, $L_{13}$. We describe the case of $L_{12}$.
Let $b \in \mathcal A_{L_{12}}$ and $i \le j$. We put
$\vec a_{12}^1= (a_{12,0},\dots,a_{12,i},a_{12,j+1},\dots,a_{12,k_{12}})$,
$\vec a_{12}^2= (b,a_{12,i+1},\dots,a_{12,j})$.
This boundary corresponds to the fiber product
\begin{equation}\label{form817}
{\mathcal M}_{\rm DR}\bigl(\vec a^1_{12},\vec a_{23},\vec a_{13};a_-,a_+;E_1\bigr)
\times_{L_{12}(b)}
\mathcal M'\bigl(L_{12};\vec a^2_{12};E_2\bigr).
\end{equation}
Here $E_1 + E_2 = E$.
We remark that we use the compactification $\mathcal M'$ in the second factor.
(See Remark~\ref{Remark524} and Section~\ref{sec:directcomp} for
this compactification.)
The bubble at $L_{23}$ and $L_{13}$ are described by the following
fiber products:
\begin{gather}
{\mathcal M}_{\rm DR}\bigl(\vec a_{12},\vec a^1_{23},\vec a_{13};a_-,a_+;E_1\bigr)
\times_{L_{12}(b)}
\mathcal M'\bigl(L_{23};\vec a^2_{23};E_2\bigr),\label{form818}
\\
{\mathcal M}_{\rm DR}\bigl(\vec a_{12},\vec a_{23},\vec a^1_{13};a_-,a_+;E_1\bigr)
\times_{L_{12}(b)}
\mathcal M'\bigl(L_{13};\vec a^2_{13};E_2\bigr).\label{form819}
\end{gather}
Here $\vec a^1_{23}$, $\vec a^2_{23}$ and $\vec a^1_{13}$, $\vec a^2_{13}$
are defined in the same way as $\vec a^1_{12}$, $\vec a^2_{12}$.

\begin{figure}[ht]\centering
\begin{tabular}{cc}
\begin{minipage}[t]{0.45\hsize}
\centering
\includegraphics[scale=0.3]{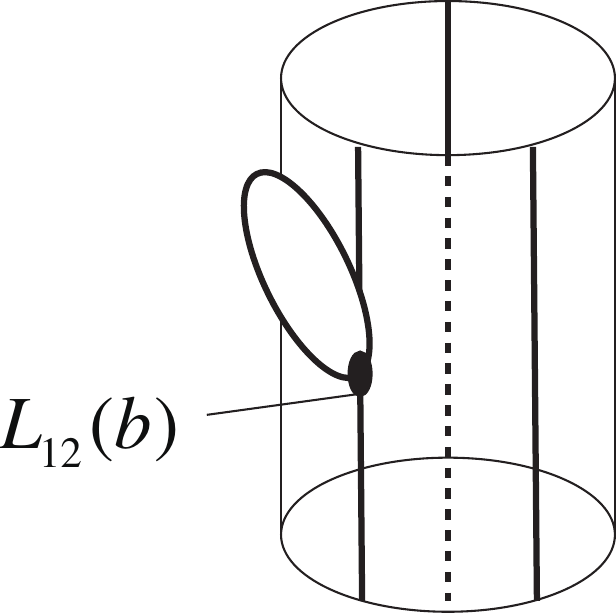}
\caption{An element of \eqref{form817}.}
\label{Figure83}
\end{minipage} &
\begin{minipage}[t]{0.45\hsize}
\centering
\includegraphics[scale=0.3]{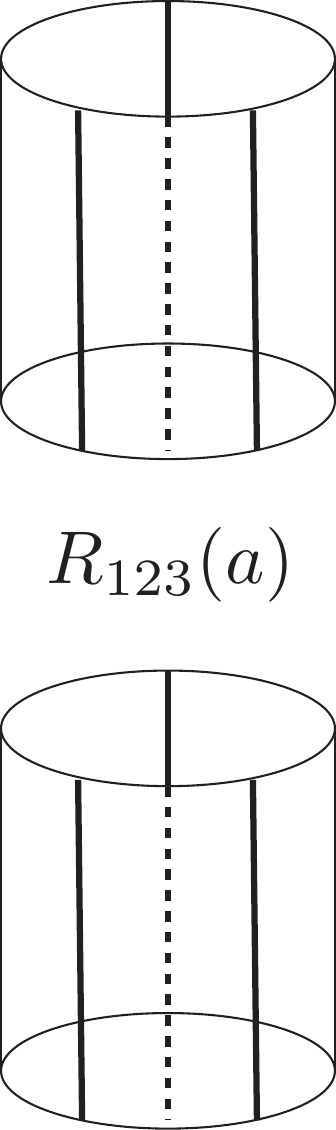}
\caption{An element of \eqref{form820}.}
\label{Figure84}
\end{minipage}
\end{tabular}
\end{figure}

\item[(II)]
The second type of boundary corresponds to the limit where the domain
will split into two parts along the second factor of $S^1 \times \R$.
It is described by the fiber product below.
Let~${j_{ii'} \in \{0,\dots,k_{ii'}\}}$.
We put
$\vec a_{ii'}^1= (a_{ii',1},\dots,a_{ii',j_{ii'}})$,
$\vec a_{ii'}^2= (a_{ii',j_{ii'}+1},\dots,a_{ii',k_{ii'}})$
if $ii' = 12$ or $23$ and
$\vec a_{ii'}^2= (a_{ii',1},\dots,a_{ii',j_{ii'}})$,
$\vec a_{ii'}^1= (a_{ii',j_{ii'}+1},\dots,a_{ii',k_{ii'}})$
if $ii' = 13$.

Note in case $j_{ii'} = 0$ (resp.\ $j_{ii'} = k_{ii'}$), $\vec a_{ii'}^1 = \varnothing$
\big(resp.\ $\vec a_{ii'}^2 = \varnothing$\big),
\begin{equation}\label{form820}
{\mathcal M}_{\rm DR}\bigl(\vec a^1_{12},\vec a^1_{23},\vec a^1_{13};a_-,a;E_1\bigr)
\times_{L_{123}(a)}
{\mathcal M}_{\rm DR}\bigl(\vec a^2_{12},\vec a^2_{23},\vec a^2_{13};a,a_+;E_2\bigr),
\end{equation}
where $E_1 + E_2 = E$ and $a \in \mathcal A_{123}$.
\end{enumerate}
We will discuss the orientation in Section~\ref{oridrum}.
The proof of the other parts of Proposition~\ref{prop811} is similar to the proof of
Theorem~\ref{therem530} and is now a routine. So we only
explain \eqref{form817}--\eqref{form819}.

We required that $u_i$ is $-J_{X_i}$ holomorphic.
Therefore, we may regard $(u_1,u_2)$ in a neighborhood of
$\gamma_{12}$ as a pseudo-holomorphic map from $(-\varepsilon,0] \times \R$ to
$-X_1 \times X_2$, by $(t,\tau) \mapsto (u_1(t,\tau),u_2(-t,\tau))$ where
$t=0$ is $S_{12}$.
See Figure~\ref{Figurebubbleseam1}.
Therefore, when a bubble on $\gamma_{12}$ occurs it corresponds
to a disk bubble as in Figure~\ref{Figurebubbleseam2}.
Note that the marked points on $\gamma_{12}$ is enumerated upward.
Therefore, the marked points on the boundary of the bubble
is enumerated according to the counter clockwise orientation
(see Figure~\ref{Figurebubbleseam2}). This implies that we can
describe such a bubble as in \eqref{form817}.
The explanation of \eqref{form818} is similar.

Let us discuss \eqref{form819}.
Note that the domain $\Omega_1$ (resp.\ $\Omega_3$) lies
right-hand side (resp.\ left-hand side) of the seam $\gamma_{13}$.
Therefore, $(u_1,u_3)$ in a neighborhood of $\gamma_{13}$
can be regarded as a pseudo-holomorphic map from $[0,\varepsilon) \times \R$ to
$-X_1 \times X_3$ by $(t,\tau) \mapsto (u_1(-t,\tau),u_3(t,\tau))$
where~${t=0}$ is $S_{13}$. See Figure~\ref{Figurebubbleseam3}.
Note that the marked points on $\gamma_{13}$ are enumerated downward.
Therefore, the marked points on the boundary of the bubble
are enumerated according to the counter clockwise orientation
(see Figure~\ref{Figurebubbleseam4}).
 This implies that we can
describe such a bubble as in \eqref{form819}.

\begin{figure}[ht]
\centering
\includegraphics[scale=0.4]{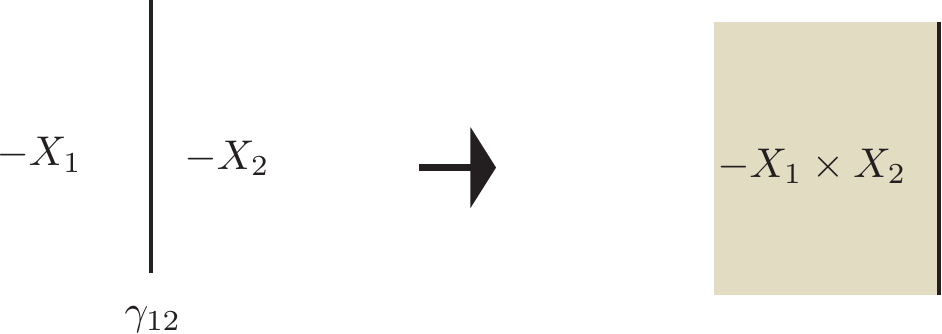}
\caption{Folding the pseudo-holomorphic map near the seam 1.}
\label{Figurebubbleseam1}
\end{figure}
\begin{figure}[ht]
\centering
\includegraphics[scale=0.4]{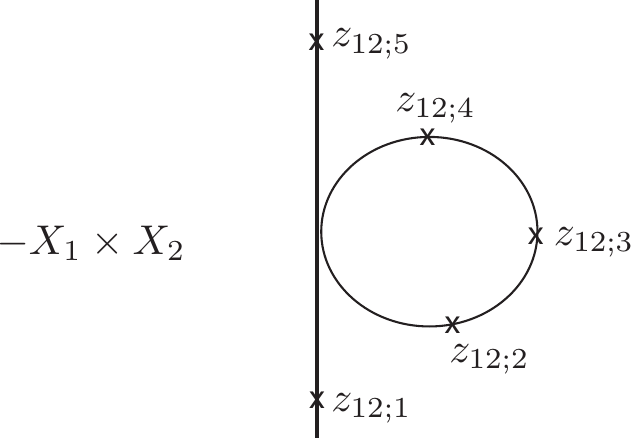}
\caption{Bubble on the seam 1.}
\label{Figurebubbleseam2}
\end{figure}
\begin{figure}[ht]
\centering
\includegraphics[scale=0.4]{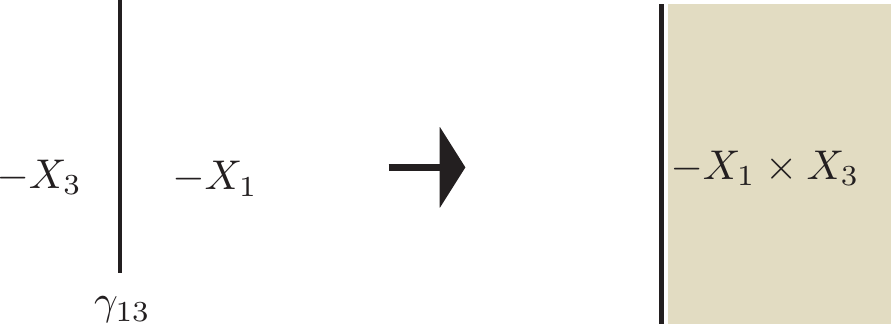}
\caption{Folding the pseudo-holomorphic map near the seam 2.}
\label{Figurebubbleseam3}
\end{figure}
\begin{figure}[ht]
\centering
\includegraphics[scale=0.4]{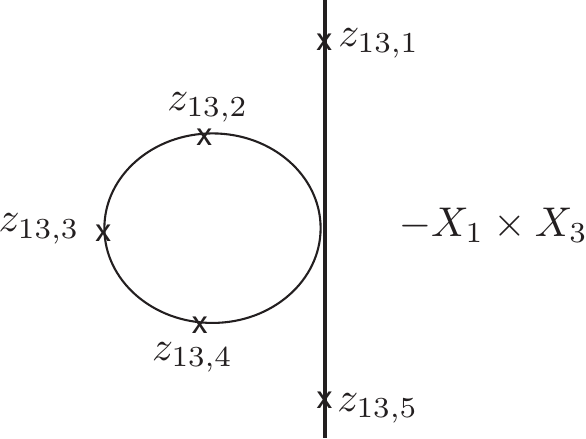}
\caption{Bubble on the seam 2.}
\label{Figurebubbleseam4}
\end{figure}

\begin{prop}\label{prop812}
For each $E_0$,
there exists a system of CF-perturbations $\widehat{\mathfrak S}$ on the spaces
${\mathcal M}(\vec a_{12},\vec a_{23},\vec a_{13};a_-,a_+;E)$ with
Kuranishi structures, which are outer collarings of thickenings of those in Proposition {\rm\ref{prop811}}, for $E < E_0$ and such that the following holds:
\begin{enumerate}\itemsep=0pt
\item[$(1)$]
The CF-perturbations \smash{$\widehat{\mathfrak S}$} are transversal to $0$.
\item[$(2)$]
The evaluation maps
${\rm ev}_{\infty,+}$, ${\rm ev}_{\infty,-}$ are strongly submersive with respect to these
CF-perturbations.\footnote{We do not require that the map $({\rm ev}_{\infty,+},{\rm ev}_{\infty,-})$ is strongly submersive.}
\item[$(3)$]
The CF-perturbations are compatible with the description of the
boundary. Namely, the restrictions of the CF-perturbations on the
boundaries coincide with the fiber product CF-perturbations
in the sense of {\rm\cite[\emph{Lemma--Definition} 10.6]{foootech2,fooonewbook}}.
\item[$(4)$]
The CF-perturbations are compatible with the forgetful maps of the boundary
marked points corresponding to the diagonal component in the same sense
as {\rm\cite[\emph{Definition} 5.1]{fooo091}}.
\end{enumerate}

\end{prop}
The proof is the same as Proposition~\ref{prop536} and is now a routine. We omit it.

We now use Propositions \ref{prop811} and \ref{prop812}
to define a filtered $A_{\infty}$ tri-module modulo $T^{E_0}$ as follows.
We put
\[
CF(L_{13};L_{12},L_{23})
=
\bigoplus_{a \in \mathcal A_{123}} \Omega(R(a)) \,\widehat{\otimes}\, \Lambda_0.
\]
We next define structure operations
\begin{align*}
\mathfrak n^{ < E_0,\varepsilon}_{k_{12},k_{23},k_{13}} \colon\
CF(L_{13})^{\otimes k_{13}
}
& \otimes
CF(L_{13};L_{12},L_{23})
\\
&\otimes CF(L_{12})^{\otimes k_{12}}\otimes
CF(L_{23})^{\otimes k_{23}}
\to
 CF(L_{12},L_{23},L_{13}).
 \end{align*}
Let
${\bf h}_{ii'} = (h_{ii',1} \otimes \dots \otimes h_{ii',k_{ii'}})
\in CF(L_{ii'})^{\otimes k_{ii'}}$.
We consider the case $h_{ii',j}$ is a differential form and
is in $\Omega(L_{ii'}(a_{ii',j}))$.
(See Definition~\ref{330}.)
Let $h_{-\infty} \in \Omega(R(a_-))$.

We define $\Omega(R(a_+))$ component of
\smash{$\mathfrak n^{{\rm tri}, E,\varepsilon}_{k_{12},k_{23},k_{13}}$ } by\index[syindex]{n1triEepsilon@$\mathfrak n^{{\rm tri}, E,\varepsilon}_{k_{12},k_{23},k_{13}}$}
\begin{equation}\label{form822}
{\rm ev}_{\infty,+}!
\bigl({\rm ev}_{13}^*{\bf h}_{13}
 \wedge {\rm ev}_{\infty,-}^*h_{-\infty} \wedge
 {\rm ev}_{12}^*{\bf h}_{12} \wedge {\rm ev}_{23}^*{\bf h}_{23}
; \widehat{\mathfrak S^{\varepsilon}}
\bigr).
\end{equation}
Here we use the space
${\mathcal M}(\vec a_{12},\vec a_{23},\vec a_{13};a_-,a_+;E)$
and its CF perturbation \smash{$ \widehat{\mathfrak S}$} to define the
integration along the fiber in \eqref{form822}.
We now put
\smash{$
\mathfrak n^{< E_0, \varepsilon}_{k_{13},k_{12},k_{23}}
:=
\sum_{E < E_0} T^E \mathfrak n^{E,\varepsilon}_{k_{13},k_{12},k_{23}}$}.
\begin{lem}
\smash{$\mathfrak n^{< E_0, \varepsilon}_{k_{13},k_{12},k_{23}}$} defines
a filtered $A_{\infty}$ tri-module modulo $T^{E_0}$.
Namely, it satisfies
\begin{align}
0 \equiv&\sum_{c_{13},c_{12},c_{23}} (-1)^{*_1}
\mathfrak n^{< E_0,\varepsilon}_{k_{c_{12};1},k_{c_{23};1},k_{c_{13};1}}
({\bf z}_{c_{13};1};\nonumber\\
&
\mathfrak n^{< E_0,\varepsilon}_{k_{c_{12},2};k_{c_{23};2},k_{c_{13};2}}
({\bf z}_{c_{13};2};w;{\bf x}_{c_{12};1},
{\bf y}_{c_{23};1});{\bf x}_{c_{12};2},{\bf y}_{c_{23};2})\nonumber\\
&+
(-1)^{*_2}\mathfrak n^{< E_0,\varepsilon}_{*,*,*}
\bigl({\bf z};w;\widehat d{\bf x},{\bf y}\bigr) +
(-1)^{*_3}\mathfrak n^{< E_0,\varepsilon}_{*,*,*}\bigl({\bf z};w;{\bf x},\widehat d{\bf y}\bigr)\nonumber \\
&+
(-1)^{*_4}\mathfrak n^{< E_0,\varepsilon}
_{*,*,*}\bigl(\hat d{\bf z};w;{\bf x},{\bf y}\bigr) +
(-1)^{*_5}\delta \bigl(\mathfrak n^{< E_0,\varepsilon}_{k_{12},k_{23},k_{13}}({\bf z};w;{\bf x},{\bf y})\bigr)\nonumber\\
&+
(-1)^{*_6}\mathfrak n^{< E_0,\varepsilon}_{k_{12},k_{23},k_{13}}({\bf z}
;\delta w;{\bf x},{\bf y})
\mod T^{E_0}.\label{form828}
\end{align}
Here $\Delta {\bf x} = \sum_{c_{12}} {\bf x}_{c_{12};1} \otimes {\bf x}_{c_{12};2}$.
We define ${\bf y}_{c_{23};1}$, ${\bf y}_{c_{23};2}$, ${\bf z}_{c_{13};1}$, ${\bf z}_{c_{13};2}$
in the same way.
The signs are by Koszul rule.
$\delta$ is the operator induced from the de Rham differential in the same way as~\eqref{form38}, \eqref{form3420000}.

\end{lem}
\begin{proof}
The proof is similar to the proof of Proposition~\ref{prop536} and is now a routine.
By Stokes' theorem (see \cite[Proposition 9.26]{foootech2} and \cite{fooonewbook}), the sum of fifth and six terms is obtained
by a similar formula as \eqref{form822} but using
the integration along the fiber on the boundary
$\partial{\mathcal M}(\vec a_{12},\vec a_{23},\vec a_{13};\allowbreak a_-,a_+;E)$.
This boundary is described by \eqref{form817}--\eqref{form820}.
By using the composition formula \cite[Theorem~10.20]{foootech2,fooonewbook},
we find that \eqref{form817}, \eqref{form818}, \eqref{form819} and \eqref{form820}
correspond to 2nd, 3rd, 4th and first term of \eqref{form828}, respectively.
\end{proof}

The rest of the proof of Proposition~\ref{pro86} is the same as the last step of the
proof of Theorem~\ref{trimain}. Namely, we show that $\mathfrak n^{< E',\varepsilon}$
is homotopic to $\mathfrak n^{ < E,\varepsilon}$ modulo $T^E$ if $E<E'$ and
also $\mathfrak n^{< E,\varepsilon}$ is homotopic to $\mathfrak n^{ < E,\varepsilon'}$.
We use this fact and homological algebra to find required filtered~$A_{\infty}$ tri-module.
\end{proof}

\begin{proof}[Proof of Lemma~\ref{lem87}]
The proof is the same as the proof of Lemma~\ref{exirespi}.
We first observe
\begin{equation}\label{eq824}
 \bigl(\tilde L_{12} \times \tilde L_{23} \times
\tilde L_{13}\bigr) \times_{X_1^2 \times X_2^2 \times X_3^2} \Delta
\cong
\tilde L_{13} \times_{X_1 \times X_3} \tilde L_{13}.
\end{equation}
Therefore,
$CF(\mathcal L_{12},\mathcal L_{23},\mathcal L_{13})$
is $\Omega\bigl(\tilde L_{13} \times_{X_1 \times X_3} \tilde L_{13},\Theta\bigr)
\,\widehat{\otimes}\,
\Lambda_0$
with some local system $\Theta$. By using Lemma~\ref{lem310},
we can uniquely choose relative spin structure $\sigma_{13}$
so that $\Theta$ is trivial.
\end{proof}

\begin{proof}[Proof of Proposition~\ref{lem89}]
The proof is the same as Proposition~\ref{prop610}.
It suffices to show that
$\mathfrak n_{0}$ is congruent to the identity map
modulo $\Lambda_+$.
By definition, $\mathfrak n_{0}$ is congruent
to the map determined by the moduli space
${\mathcal M}(\varnothing,\varnothing,a_{13};o,a_{13};0)$.
Here $o$ denotes the diagonal component
and we use the diffeomorphism \eqref{eq824}
to identify $\mathcal A_{123}$ with $\mathcal A_{L_{13}}$.
(Here $\mathcal A_{123}$ (resp.\ $\mathcal A_{L_{13}}$)
is the set of connected components of the left-hand side
(resp.\ right-hand side) of \eqref{eq824}.)
Using the fact that ${\mathcal M}(\varnothing,\varnothing,a_{13};o,a_{13};0)$
consists of constant maps, we can easily show that it induces the identity map.
\end{proof}

\section{Compatibility of compositions}
\label{sec:comptibility}

\subsection{Statement}
\label{subsec:compacompstate}

\begin{thm}\label{thm93}
Suppose we are in Situation {\rm\ref{situ82}}.
Let $\mathcal L_{12} \in \mathfrak {OB}(\mathfrak{Fukst}(-X_1 \times X_2))$,
$\mathcal L_{23} \in \mathfrak {OB}(\mathfrak{Fukst}(-X_2 \times X_3))$.
We put
$
\mathcal L_{13} = \mathcal L_{23} \circ \mathcal L_{12}
= \mathfrak{Comp}_{\rm ob}(\mathcal L_{12},\mathcal L_{23})$.
Then the correspondence functor $\mathcal W_{\mathcal L_{13}}$
associated to $\mathcal L_{13}$ is homotopy equivalent to
the composition $\mathcal W_{\mathcal L_{23}} \circ \mathcal W_{\mathcal L_{12}}$
of the correspondence functors associated to $\mathcal L_{12}$ and $\mathcal L_{23}$
respectively. Namely,
\begin{equation}\label{compocompfor}
\mathcal W_{\mathcal L_{23} \circ \mathcal L_{12}} \sim \mathcal W_{\mathcal L_{23}} \circ \mathcal W_{\mathcal L_{12}}.
\end{equation}
\end{thm}
Note that
\begin{gather*}
\mathcal W_{\mathcal L_{12}}\colon\ \mathfrak{Fukst}(X_1;\mathbb L_1) \to \mathfrak{Fukst}(X_2;\mathbb L_2), \qquad
\mathcal W_{\mathcal L_{23}}\colon\ \mathfrak{Fukst}(X_2;\mathbb L_2) \to \mathfrak{Fukst}(X_3;\mathbb L_3), \\
\mathcal W_{\mathcal L_{13}}\colon\ \mathfrak{Fukst}(X_1;\mathbb L_1) \to \mathfrak{Fukst}(X_3;\mathbb L_3).
\end{gather*}
\eqref{compocompfor} is a homotopy equivalence as strict, unital and gapped
filtered $A_{\infty}$ functors
from $\mathfrak{Fukst}(X_1;\allowbreak\mathbb L_1)$ to $\mathfrak{Fukst}(X_3;\mathbb L_3)$.

In this section, we prove the following weaker version of Theorem~\ref{thm93}.

\begin{prop}\label{prop912}
Suppose we are in the situation of Theorem {\rm\ref{thm93}}.
Let $\mathcal L_1 = (L_1,\sigma_1,b_1)$ be an object of $\mathfrak{Fukst}(X_1;\mathbb L_1)$.
We put
\[
(\mathcal W_{ \mathcal L_{23} \circ \mathcal L_{12}})_{\rm ob}(\mathcal L_1)
= \mathcal L_3^{(1)} = \bigl(L^{(1)}_3,\sigma_3^{(1)},b_3^{(1)}\bigr), \qquad
(\mathcal W_{ \mathcal L_{13}})_{\rm ob}(\mathcal L_1)
= \mathcal L_3^{(2)} = \bigl(L^{(2)}_3,\sigma_3^{(2)},b_3^{(2)}\bigr).
\]
Then we have the following:
\begin{enumerate}\itemsep=0pt
\item[$(1)$]
\smash{$\bigl(L^{(1)}_3,\sigma_3^{(1)}\bigr) = \bigl(L^{(2)}_3,\sigma_3^{(2)}\bigr)$}. Here the equality is as submanifolds
equipped with relative spin structures.
\item[$(2)$]
$b_3^{(1)}$ is gauge equivalent to $b_3^{(2)}$ in the sense of {\rm\cite[\emph{Definition} 4.3.1]{fooobook}}.
\end{enumerate}

\end{prop}
Proposition~\ref{prop912} is the object part of Theorem~\ref{thm93}.
The morphism part will be proved in the next section.
Proposition~\ref{prop912}\,(1) is proved in Section~\ref{oriYdiagarm}.

\subsection{Lekili--Lipyanskiy's Y-diagram}
\label{subsec:Ydiagram}

The proofs of Theorem~\ref{thm93} and Proposition~\ref{prop912} are based on
moduli spaces of configurations introduced by Lekili--Lipyanskiy in \cite{LL},
which they called a $Y$-diagram.
In this subsection, we define and study the moduli space of $Y$-diagrams.\index{$Y$-diagram}

\begin{figure}[ht]
\centering
\includegraphics[scale=0.28]{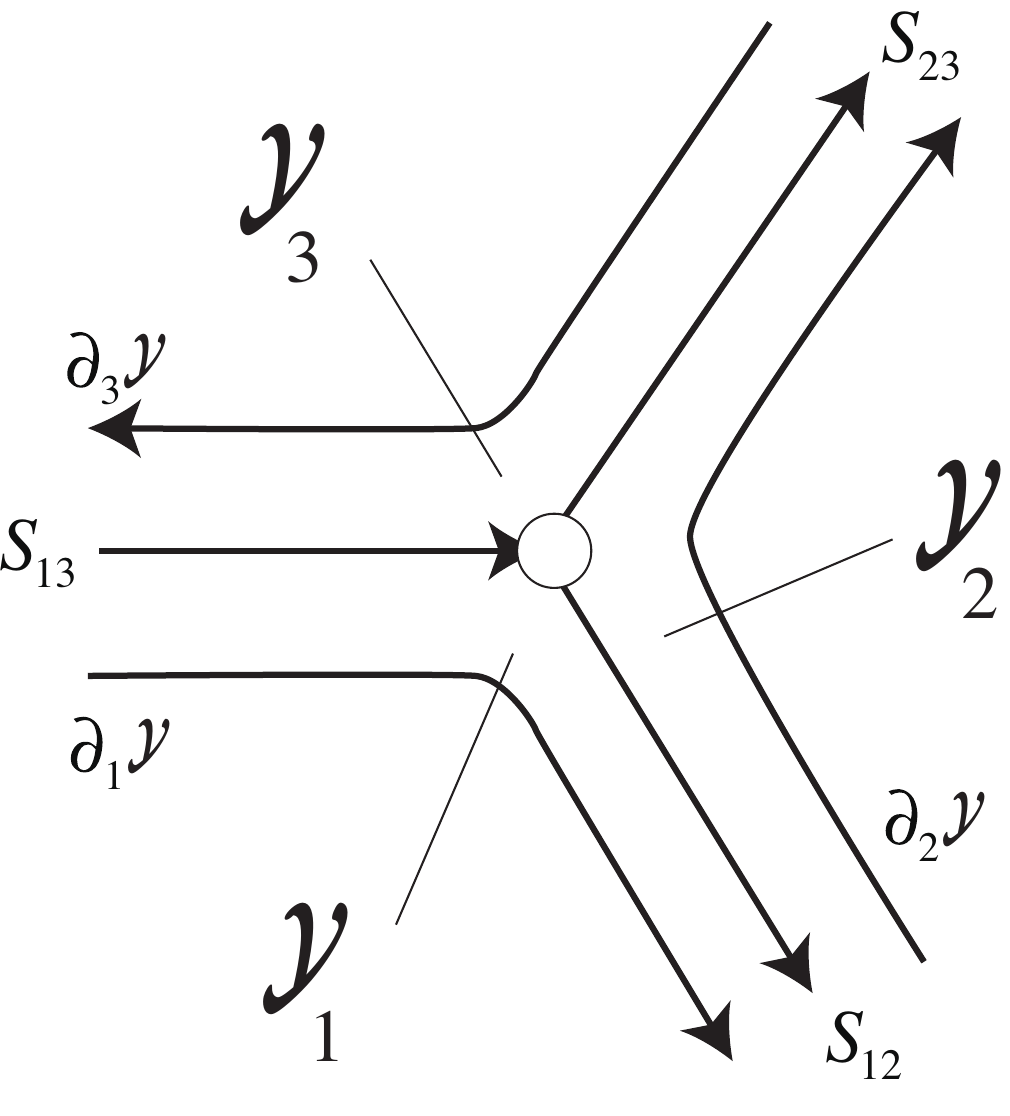}
\caption{Domain $\mathcal Y$.}
\label{Figure91}
\end{figure}

We consider the domain $\mathcal Y = \mathcal Y_1 \cup \mathcal Y_2 \cup
\mathcal Y_3 \subset \C$ as in Figure~\ref{Figure91}.
The boundary $\partial \mathcal Y$ has three
connected components $\partial_i \mathcal Y
= \partial \mathcal Y \cap \partial \mathcal Y_i$
($i=1,2,3$), which are diffeomorphic to $\R$.
We choose the diffeomorphism so that the direction of the arrow
in Figure~\ref{Figure91} coincides with the positive direction of $\R$.

The closure of the domain
$\mathcal Y$ minus a point $S_{12} \cap S_{23} \cap S_{13}$ has 4 ends.
We identify the end which is the neighborhood of the white circle in
Figure~\ref{Figure91} with
$S^1 \times (-\infty,0]$.
We take a~diffeomorphism
$
\phi_{123} \colon S^1 \times (-\infty,0] \to \mathcal Y$
to an open subset
such that
\begin{conds}
\quad
\begin{enumerate}\itemsep=0pt
\item[(1)]
$\phi_{123}$ is an anti-biholomorphic diffeomorphism to its image,
which is a~neighborhood of the point $S_{12} \cap S_{23} \cap S_{13}$ minus
$S_{12} \cap S_{23} \cap S_{13}$.
\item[(2)]
We identify $S^1 \times (-\infty,0] \subset S^1 \times (-\infty,\infty) = W$, where
$W$ is as in \eqref{form89}.
Then we require
\[
W_i \cap \bigl(S^1 \times (-\infty,0]\bigr) = \phi_{123}^{-1}(\mathcal Y_i)
\]
for $i=1,2,3$.
\end{enumerate}

\end{conds}
\begin{rem}\label{newrem94}
We emphasis that $\phi_{123}$ is an {\it anti}-biholomorphic map.
In fact, $\phi_{123}(t,\tau) = e^{2\pi (\tau+i t)}$ and
the complex structure of the domain is $j(\partial/\partial t) =
\partial/\partial \tau$.
We will identify the image of $\phi_{123}$ as a part of the domain
of the pseudo-holomorphic drum appearing in Section~\ref{sec:comp}.
Then a $J_{X_i}$ holomorphic map on $W_i$ will become
$-J_{X_i}$ holomorphic from an open set of the drum.
This is the reason why we required that the map $u_i$ is
$-J_{X_i}$ holomorphic in Definition~\ref{def916}\,(2).
\end{rem}
The other three ends intersect with $\mathcal Y_1$ and $\mathcal Y_2$
(resp.\ $\mathcal Y_2$ and $\mathcal Y_3$, $\mathcal Y_1$ and $\mathcal Y_3$).
We take a~diffeomorphism
$
\phi_{ii'} \colon [-1,1] \times (-\infty,0] \to \mathcal Y$
to an open subset for $(ii') = (12), (23)$, or $(13)$ such that the following conditions hold.
\begin{conds}\label{Cond95}
\quad
\begin{enumerate}\itemsep=0pt
\item[(1)]
The map $\phi_{ii'}$ is biholomorphic.
\item[(2)]
We require $
[-1,0] \times (-\infty,0] = \phi_{ii'}^{-1}(\mathcal Y_i)$,
$
[0,1] \times (-\infty,0] = \phi_{ii'}^{-1}(\mathcal Y_{i'})$
for $(ii') = (12)$ or $(23)$.
We also require
$
[-1,0] \times [0,+\infty) = \phi_{13}^{-1}(\mathcal Y_i)$,
$
[0,1] \times [0,+\infty) = \phi_{13}^{-1}(\mathcal Y_i)$
for $i = 1,3$.
\end{enumerate}
\end{conds}
\begin{figure}[ht]
\centering
\includegraphics[scale=0.3]{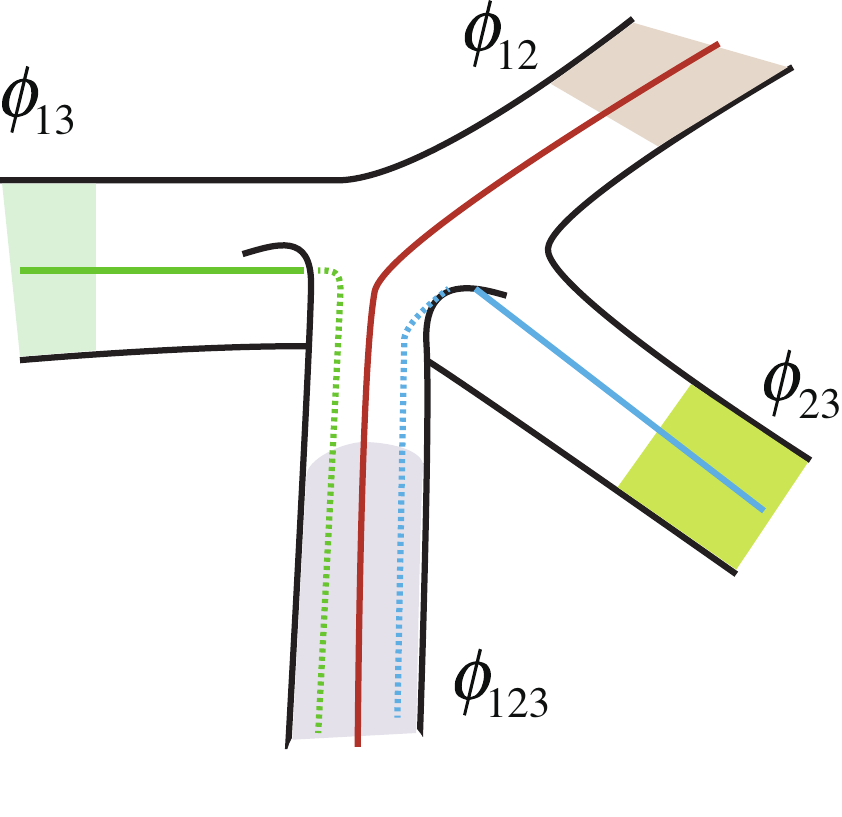}
\caption{$\phi_{123}$, $\phi_{ii'}$.}
\label{Figure92}
\end{figure}
We next put
$
S_{ii'} = \mathcal Y_i \cap \mathcal Y_{i'}
$
for $(ii') = (12), (23)$, or $(13)$.
$S_{ii'}$ is diffeomorphic to $\R$.
We call~$S_{ii'}$ a {\it seam} and the point
$S_{12} \cap S_{23} \cap S_{13}$ the {\it hole}. We take a \index{hole}
diffeomorphism between the seams and $\R$ as follows:
\begin{enumerate}\itemsep=0pt
\item[(so1)] Suppose that $(ii') = (12), (23)$.
Then for $-\tau$ which is sufficiently negative
the point of $S_{ii'}$ corresponding to $-\tau$ lies
in the image of $\phi_{ii'}$.
\item[(so2)]
Suppose that $(ii') = (13)$. Then for $\tau$ which is sufficiently positive
the point of $S_{13}$ corresponding to $\tau$ lies
in the image of $\phi_{13}$.
\end{enumerate}
See the arrows in Figure~\ref{Figure91} which show the orientation of
the seams.
Note that this orientation coincides with the way we enumerate the
marked points on the seams in the case of pseudo-holomorphic drums.

We orient the boundary of $\mathcal Y$ by the usual
counter clock-wise orientation of a boundary of a~domain of $\C$ (see the arrows in Figure~\ref{Figure91}).
Then on the images of $\phi_{ii'}$, the orientation
of the boundary and the seams coincide with the way we enumerate the
marked points in Definition~\ref{def516}\,(3).

We decompose fiber products to connected components
\begin{gather}
\tilde L_{ii'} \times_{X_i\times X_{i'}} \tilde L_{ii'}
= \bigcup_{a \in \mathcal A_{L_{ii'}}} L_{ii'}(a)
,\qquad
\tilde L_i \times_{X_i} \tilde L_i = \bigcup_{a \in \mathcal A_{L_{i}}} L_i(a),\nonumber
\\
\tilde L_i \times_{X_i} \tilde L_{ii'} \times_{X_{i'}} \tilde L_{i'}
=
\bigcup_{a \in \mathcal A_{R_{ii'}}} R_{ii'}(a),\label{decomp95}
\\
(L_{12} \times L_{23} \times L_{13}) \times_{(X_1\times X_2 \times X_3)^2} \Delta
=
\bigcup_{a \in \mathcal A_{123}} R_{123}(a),\label{decomp96}
\end{gather}
where $\Delta$ is the diagonal in $(X_1\times X_2 \times X_3)^2$.
See Definition~\ref{def3131}\,(5).

Let $\vec a_{ii'} = (a_{ii',1},\dots,a_{ii',k_{ii'}}) \in (\mathcal A_{L_{ii'}})^{k_{ii'}}$,
$\vec a_{i} = (a_{i,1},\dots,a_{i,k_{i}}) \in (\mathcal A_{L_{i}})^{k_{i}}$,
$a_{\infty,123} \in \mathcal A_{123}$.
Let $\vec a_{\infty} = (a_{\infty,12},a_{\infty,23},a_{\infty,13})$ with $a_{\infty,ii'} \in \mathcal A_{R_{ii'}}$.

\begin{figure}[ht]
\centering
\includegraphics[scale=0.3]{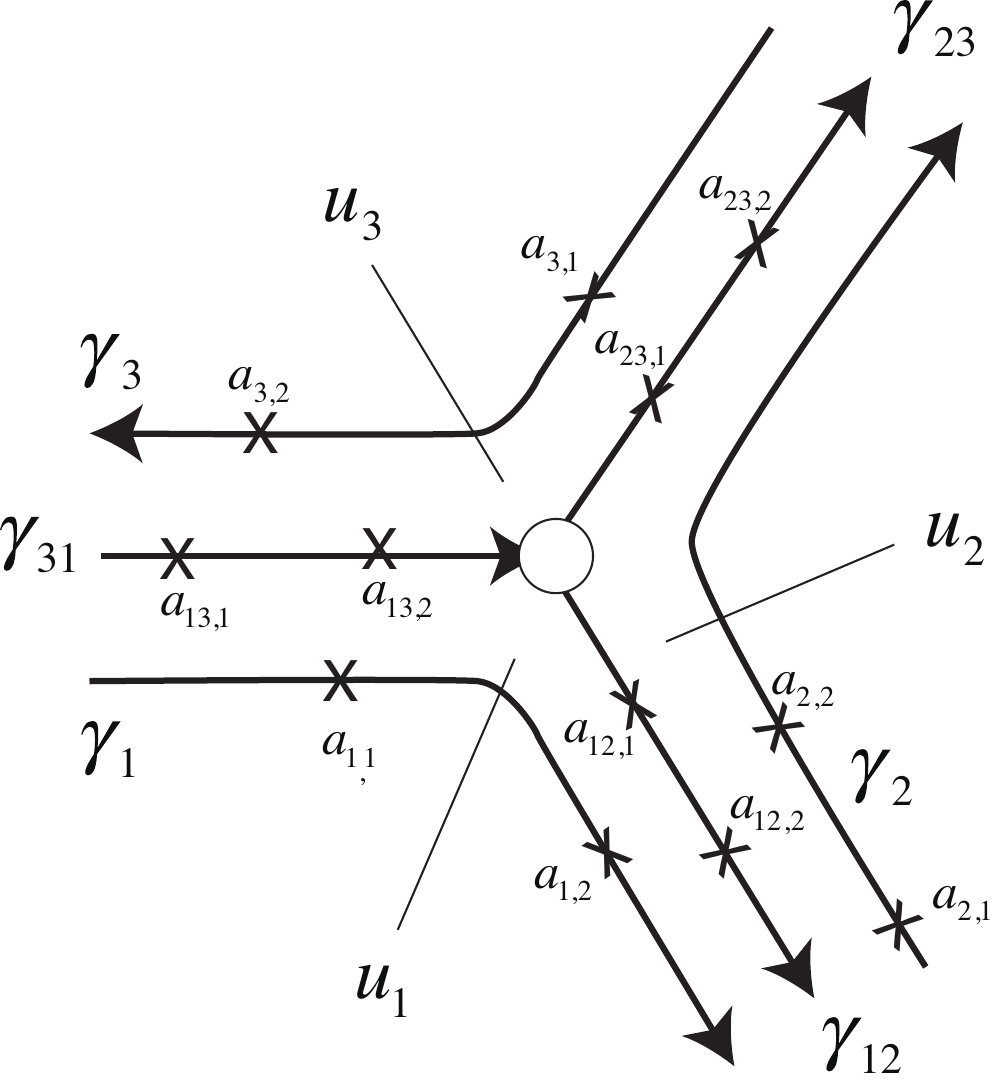}
\caption{An element of $\overset{\ \text{\tiny $\circ\circ$}}{\mathcal M}_{\rm Y}(\vec a_{12},\vec a_{23},\vec a_{13};\vec a_{1},\vec a_{2},\vec a_{3},a_{\infty,-},\vec a_{\infty,+};E)$.}
\label{Figure93}
\end{figure}

We next define the set
\smash{$\overset{\ \text{\tiny $\circ\circ$}}{\mathcal M}_{\rm Y}(\vec a_{12},\vec a_{23},\vec a_{13};\vec a_{1},\vec a_{2},\vec a_{3};a_{\infty,123},\vec a_{\infty};E)$}.
\begin{defn}\label{def1016}
We consider
\[
(\Sigma;
\vec z_1,\vec z_2,\vec z_3;\vec z_{12},\vec z_{23},\vec z_{13};u_1,u_2,u_3;\gamma_{1},\gamma_2,\gamma_{3};\gamma_{12},\gamma_{23},\gamma_{13})
\]
with the following properties (see Figure~\ref{Figure93}):
\begin{enumerate}\itemsep=0pt
\item[(1)]
The space $\Sigma$ is a bordered Riemann surface which is a union of $\mathcal Y$ and
trees of sphere components attached to $\mathcal Y$.
The roots of the trees of sphere components are neither on $S_{12}$, $S_{23}$, $S_{13}$
nor on $\partial \mathcal Y$.
\item[(2)]
We denote by $\Sigma_{1}$ the union of $\mathcal Y_{1}$ and the
trees of sphere components rooted on $\mathcal Y_{1}$.
We define $\Sigma_{2}$, $\Sigma_{3}$ in the same way.
The map $u_i \colon \Sigma_i \to X_i$ is $J_{X_i}$ holomorphic for $i=1,2,3$.
\item[(3)]
$\vec z_i = (z_{i,1},\dots, z_{i,k_i})$, $i=1,2,3$, and $z_{i,j} \in
\partial_i \mathcal Y$.
We require $z_{i,j} < z_{i,j'}$ for~${j<j'}$, where we identify
$\partial_i \mathcal Y \cong \R$ using the counter clockwise orientation.
\item[(4)]
$\vec z_{ii'} = (z_{ii',1},\dots,z_{ii',k_{ii'}})$, $ii'=12,23,13$, and
$z_{ii',j} \in S_{ii'}$.
We require $z_{ii',j} < z_{ii',j'}$ for $j<j'$, where we identify
$S_{ii'} \cong \R$ as in (so1),(so2).
We put $\vert \vec z_{ii'}\vert = \{z_{ii',1},\dots,z_{ii',k_{ii'}}\}$.
\item[(5)]
The maps $\gamma_i\colon \partial\Sigma \cap \mathcal Y_i \setminus \vert \vec z_i\vert\
\to \tilde L_i$ are smooth and satisfy
$
i_{L_{i}} (\gamma_{i}(z)) = u_i(z)$.

\item[(6)]
The maps $\gamma_{ii'} \colon S_{ii'} \setminus \vert\vec z_{ii'}\vert \to \tilde L_{ii'}$,
$(ii') = (12),(23),(13)$,
are smooth and satisfy
\[
i_{L_{ii'}} (\gamma_{ii'}(z)) = (u_i(z),u_{i'}(z)).
\]
\item[(7)]
On $\vec z_{i}$, the map $\gamma_{i}$ satisfies the switching condition
\begin{equation}\label{form9142}
\bigl(\lim_{z \in S_{i} \uparrow z_{i,j}}\gamma_{i}(z),\lim_{z \in \partial\Sigma \cap \mathcal Y_i \downarrow z_{i,j}}\gamma_{i}(z)
\bigr)
\in L_{i}(a_{i,j}).
\end{equation}
Here we identify $\partial\Sigma \cap \mathcal Y_i \cong \R$ by the counter clockwise orientation
and then $ \uparrow $, $ \downarrow$ have obvious meaning similar to Definition~\ref{def3737}\,(5).
\item[(8)]
On $\vec z_{ii'}$, the map $\gamma_{ii'}$ satisfies the switching condition
\begin{equation}\label{form914}
\bigl(\lim_{z \in S_{ii'} \uparrow z_{ii',j}}\gamma_{ii'}(z),\lim_{z \in S_{ii'} \downarrow z_{ii',j}}\gamma_{ii'}(z)
\bigr)
\in L_{ii'}(a_{ii',j}).
\end{equation}
Here we identify $S_{ii'} \cong \R$ by (so1), (so2) and then $ \uparrow $, $ \downarrow$ have obvious meaning similar to Definition~\ref{def3737}\,(5).
\item[(9)]
On the image of $\phi_{ii'}$, the map $\gamma_{ii'}$ satisfies the asymptotic boundary condition
\begin{gather}
\lim_{\tau \to + \infty}( (\gamma_i(-\tau),\gamma_{i'}(\tau)), \gamma_{ii'}(-\tau)
)
\in R_{ii'}(a_{\infty,ii'}) \qquad \text{if }\  (ii') = (12)\  \text{ or } \ (23),\nonumber \\
\lim_{\tau \to + \infty}( (\gamma_1(\tau),\gamma_{3}(-\tau)), \gamma_{13}(\tau)
)
\in R_{13}(a_{\infty,13}).\label{form915}
\end{gather}
\item[(10)]
On the image of $\phi_{123}$, the map $\gamma_{ii'}$
satisfies the asymptotic boundary condition
\begin{equation}\label{form916}
\lim_{\tau \to +\infty}\left(\gamma_{12}(-\tau),\gamma_{23}(-\tau),\gamma_{13}(\tau)
\right)
\in R_{123}(a_{\infty,123}).
\end{equation}
\item[(11)]
The stability condition, Definition~\ref{defn1015}\,(2) below, is satisfied.
\item[(12)]
$
\int_{\Sigma_1}u_1^*(\omega_1) + \int_{\Sigma_2}u_1^*(\omega_2)
+ \int_{\Sigma_3}u_3^*(\omega_3) = E$.
\end{enumerate}
In Definition~\ref{defn1015}\,(3), we will define an equivalence relation $\sim$ among the objects
\[
(\Sigma;\vec z_{12},\vec z_{23},\vec z_{13};u_1,u_2,u_3;\gamma_{1},\gamma_2,\gamma_{3};\gamma_{12},\gamma_{23},\gamma_{13})
\]
satisfying (1)--(12). We denote by
\smash{$\overset{\ \text{\tiny $\circ\circ$}}{\mathcal M}_{\rm Y}(\vec a_{12},\vec a_{23},\vec a_{13};\vec a_{1},\vec a_{2},\vec a_{3};a_{\infty,-},\vec a_{\infty,+};E)$}
\index[syindex]{M1a12a23a31a@$\overset{\ \text{\tiny $\circ\circ$}}{\mathcal M}(\vec a_{12},\vec a_{23},\vec a_{31};\vec a_{1},\vec a_{2},\vec a_{3};a_{\infty,-},\vec a_{\infty,+};E)$}
the set of all the equivalence classes of this equivalence relation.
\end{defn}

\begin{defn}\label{defn1015}
Let
\begin{gather*}
\mathfrak x =(\Sigma;\vec z_1,\vec z_2,\vec z_3;\vec z_{12},\vec z_{23},\vec z_{13};u_1,u_2,u_3;\gamma_{1},\gamma_2,\gamma_{3};\gamma_{12},\gamma_{23},\gamma_{13}), \\
\mathfrak x' = (\Sigma';\vec z^{\,\prime}_1,\vec z^{\,\prime}_2,\vec z^{\,\prime}_3;\vec z^{\,\prime}_{12},\vec z^{\,\prime}_{23},\vec z^{\,\prime}_{13};u'_1,u'_2,u'_3;\gamma'_{1},\gamma'_2,\gamma'_{3};\gamma'_{12},\gamma'_{23},\gamma'_{13})
\end{gather*}
be objects satisfying Definition~\ref{def1016}\,(1)--(10) and (12).
\begin{enumerate}\itemsep=0pt
\item[(1)]
An isomorphism from $\mathfrak x$ to $\mathfrak x'$ is a map $v\colon \Sigma \to \Sigma'$ such that
\begin{enumerate}\itemsep=0pt
\item It is biholomorphic.
\item It sends $\Sigma_i$ to $\Sigma'_i$.
\item It sends $\vec z_i$ to $\vec z^{\,\prime}_i$ and $\vec z_{ii'}$ to $\vec z^{\,\prime}_{ii'}$.
\item $u'_i \circ v = u_i$, $\gamma'_{i} \circ v = \gamma_{i}$, $\gamma'_{ii'} \circ v = \gamma_{ii'}$.
\end{enumerate}
\item[(2)]
$\mathfrak x$ is said to be {\it stable} \index{stable} if the set of
all isomorphisms from $\mathfrak x$ to $\mathfrak x$ is finite.
\item[(3)]
We say $\mathfrak x$ is {\it equivalent}
\index{equivalent} to $\mathfrak x'$ if
there exists an isomorphism from $\mathfrak x$ to $\mathfrak x'$.
\end{enumerate}

\end{defn}
We define the evaluation maps
\begin{gather}\label{form91550000}
{\rm ev}_{i} = ({\rm ev}_{i,1},\dots,{\rm ev}_{i,k_{i}}) \colon\ \overset{\ \text{\tiny $\circ\circ$}}{\mathcal M}_{\rm Y}(\vec a_{12},\vec a_{23},\vec a_{13};\vec a_{1},\vec a_{2},\vec a_{3},a_{\infty,123},\vec a_{\infty};E)
\to \prod_{j=1}^{k_{i}} L_{i}(a_{i,k})
\end{gather}
and
\begin{gather}
{\rm ev}_{ii'} = ({\rm ev}_{ii',1},\dots,{\rm ev}_{ii',k_{ii'}})\colon\nonumber \\
\qquad \overset{\ \text{\tiny $\circ\circ$}}{\mathcal M}_{\rm Y}(\vec a_{12},\vec a_{23},\vec a_{13};\vec a_{1},\vec a_{2},\vec a_{3},a_{\infty,123},\vec a_{\infty};E)
\to \prod_{j=1}^{k_{ii'}} L_{ii'}(a_{ii',k})\label{form915500}
\end{gather}
by the left-hand sides of \eqref{form9142} and \eqref{form914}, respectively.

We also define
\begin{gather}
\widehat{\rm ev}_{\infty} = ({\rm ev}_{\infty,123},{\rm ev}_{\infty}) = ({\rm ev}_{\infty,123},
({\rm ev}_{\infty,12},{\rm ev}_{\infty,23},{\rm ev}_{\infty,13}))\colon\nonumber \\
\qquad\overset{\ \text{\tiny $\circ\circ$}}{\mathcal M}_{\rm Y}(\vec a_{12},\vec a_{23},\vec a_{31};\vec a_{1},\vec a_{2},\vec a_{3},a_{\infty,123},\vec a_{\infty};E) \nonumber\\
\qquad\quad\to
R(a_{\infty,123}) \times L_{12}(a_{\infty,12}) \times L_{23}(a_{\infty,23}) \times L_{13}(a_{\infty,13})\label{form9166}
\end{gather}
by using the left-hand side of \eqref{form915} and \eqref{form916}.
\begin{prop}\label{prop9911}
We can define a topology, stable map topology, on
the moduli space
\[
\overset{\ \text{\tiny $\circ\circ$}}{\mathcal M}_{\rm Y}(\vec a_{12},\vec a_{23},\vec a_{13};\vec a_{1},\vec a_{2},\vec a_{3},a_{\infty,123},\vec a_{\infty};E)
\]
such that it has a compactification\index[syindex]{M1Ya12a23@${\mathcal M}_{\rm Y}(\vec a_{12},\vec a_{23},\vec a_{13};\vec a_{1},\vec a_{2},\vec a_{3},a_{\infty,123},\vec a_{\infty};E)$}
${\mathcal M}_{\rm Y}(\vec a_{12},\vec a_{23},\vec a_{13};\vec a_{1},\vec a_{2},\vec a_{3},a_{\infty,123},\vec a_{\infty};E)$, which is a compact metrizable space.
They have Kuranishi structures with corners which enjoy the following properties:
\begin{enumerate}\itemsep=0pt
\item[$(1)$]
The normalized boundary of ${\mathcal M}_{\rm Y}(\vec a_{12},
\vec a_{23},\vec a_{13};\vec a_{1},\vec a_{2},\vec a_{3},a_{\infty,123},\vec a_{\infty};E)$ is a disjoint union of $4$ types of
fiber products which we describe below.
\item[$(2)$]
The evaluation maps \eqref{form91550000}, \eqref{form915500} and \eqref{form9166}
extend to strongly smooth maps with respect to this Kuranishi structure.
The map ${\rm ev}_{\infty}$ in \eqref{form9166} is weakly submersive.
The extension is compatible with the description of the boundary
in item $(1)$.
\item[$(3)$]
The orientation local system of ${\mathcal M}_{\rm Y}(\vec a_{12},\vec a_{23},\vec a_{31};\vec a_{1},\vec a_{2},\vec a_{3},a_{\infty,123},\vec a_{\infty};E)$ is isomorphic to the tensor product of
the pullbacks of $\Theta^-$ by
the evaluation maps \eqref{form91550000}, \eqref{form915500} and \eqref{form9166}. For the component $L_{13}(a_{13})$ we take
$\Theta^+$ in place of $\Theta^-$.
\item[$(4)$]
The Kuranishi structures are compatible with the forgetful maps of the marked points corresponding to the
diagonal components.
\end{enumerate}
\end{prop}
We now describe the boundary components:
\begin{enumerate}\itemsep=0pt
\item[(I)]
The first type of boundary corresponds to a bubble at one of the Lagrangian boundary conditions $L_{12}$,
$L_{23}$, $L_{13}$. We describe the case of $L_{12}$.
Let $b \in \mathcal A_{L_{12}}$ and $i \le j$. We put
$\vec a_{12}^1= (a_{12,0},\dots,a_{12,i},b,a_{12,j+1},\dots,a_{12,k_{12}})$,
$\vec a_{12}^2= (b,a_{12,i+1},\dots,a_{12,j})$.
This boundary corresponds to the fiber product
\begin{equation}\label{form9817}
{\mathcal M}_{\rm Y}\bigl(\vec a^1_{12},\vec a_{23},\vec a_{31};\vec a_{1},\vec a_{2},\vec a_{3};
a_{\infty,123},\vec a_{\infty};E_1\bigr)
\times_{L_{12}(b)}
\mathcal M'\bigl(L_{12};\vec a^2_{12};E_2\bigr).
\end{equation}
Here $E_1 + E_2 = E$.
We remark that we use the compactification $\mathcal M'$ in the second factor.
The compactification $\mathcal M'$ is discussed in Remark~\ref{Remark524} and Section~\ref{sec:directcomp}.
See Figure~\ref{Figure94}.
The bubble at $L_{23}$ and $L_{13}$ are described by the following
fiber products:
\begin{gather}
{\mathcal M}_{\rm Y}\bigl(\vec a_{12},\vec a^1_{23},\vec a_{31};\vec a_{1},\vec a_{2},\vec a_{3};
a_{\infty,123},\vec a_{\infty};E_1\bigr)
\times_{L_{23}(b)}
\mathcal M'\bigl(L_{23};\vec a^2_{23};E_2\bigr),\label{form9818}
\\
{\mathcal M}_{\rm Y}\bigl(\vec a_{12},\vec a_{23},\vec a^1_{13};\vec a_{1},\vec a_{2},\vec a_{3};
a_{\infty,123},\vec a_{\infty};E_1\bigr)
\times_{L_{13}(b)}
\mathcal M'\bigl(L_{13};\vec a^2_{13};E_2\bigr).\label{form9819}
\end{gather}
Here $\vec a^1_{23}$, $\vec a^2_{23}$ and $\vec a^1_{13}$, $\vec a^2_{13}$
are defined in the same way as $\vec a^1_{12}$, $\vec a^2_{12}$.
\begin{figure}[ht]
\centering
\includegraphics[scale=0.3]{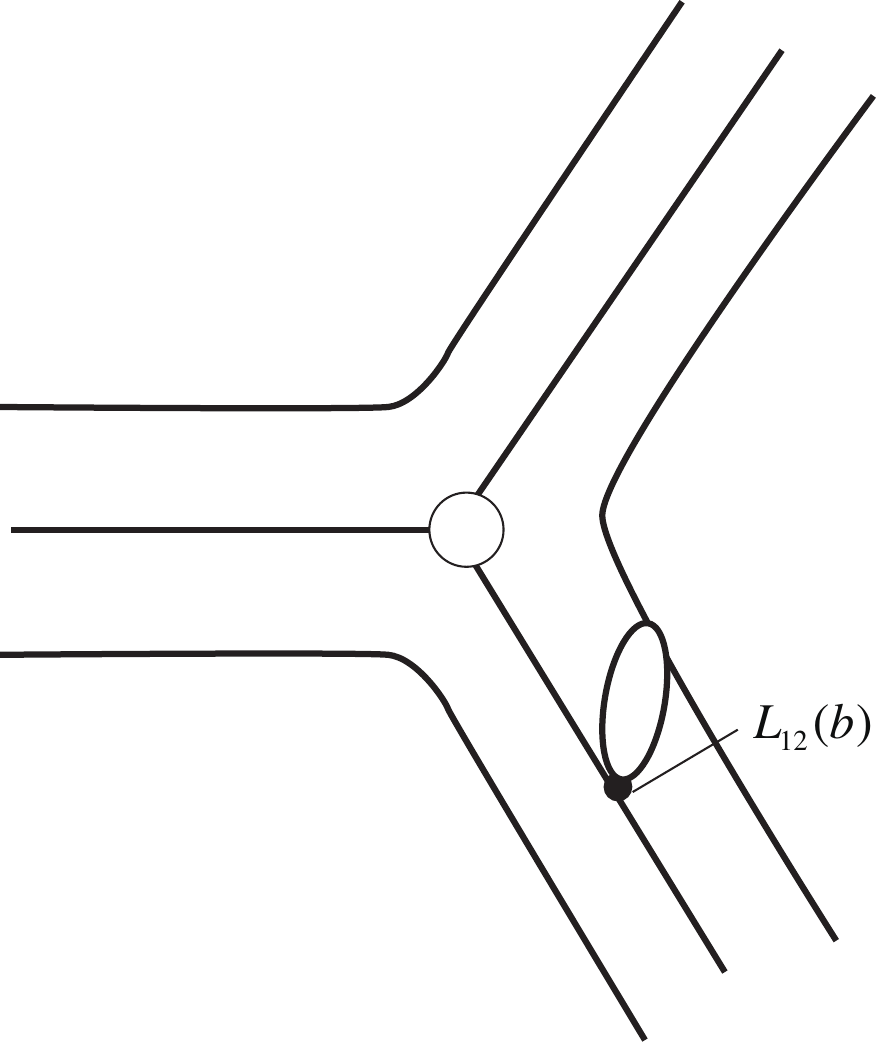}
\caption{Boundary of type (I).}
\label{Figure94}
\end{figure}
\item[(II)]
The second type of boundary corresponds to a bubble at one of the Lagrangian boundary conditions $L_{1}$,
$L_{2}$, $L_{3}$. We describe the case of $L_{1}$.
Let $b \in \mathcal A_{L_{1}}$ and $i \le j$. We put
$\vec a_{1}^1= (a_{1,0},\dots,a_{1,i},b,a_{1,j+1},\dots,a_{1,k_{1}})$,
$\vec a_{1}^2= (b,a_{1,i+1},\dots,a_{1,j})$.
This boundary corresponds to the fiber product
\begin{equation}\label{form98171}
{\mathcal M}_{\rm Y}\bigl(\vec a_{12},\vec a_{23},\vec a_{13};\vec a^1_{1},\vec a_{2},\vec a_{3};
a_{\infty,123},\vec a_{\infty};E_1\bigr)
\times_{L_{1}(b)}
\mathcal M\bigl(L_{1};\vec a^2_{1};E_2\bigr).
\end{equation}
Here $E_1 + E_2 = E$.
See Figure~\ref{Figure95}.

The bubble at $L_{2}$ and $L_{3}$ are described by the following
fiber products:
\begin{gather}
{\mathcal M}_{\rm Y}\bigl(\vec a_{12},\vec a_{23},\vec a_{13};\vec a_{1},\vec a^1_{2},\vec a_{3};
a_{\infty,123},\vec a_{\infty};E_1\bigr)
\times_{L_{2}(b)}
\mathcal M\bigl(L_{2};\vec a^2_{2};E_2\bigr),\label{form98181}
\\
{\mathcal M}_{\rm Y}\bigl(\vec a_{12},\vec a_{23},\vec a_{13};\vec a_{1},\vec a_{2},\vec a^1_{3};
a_{\infty,123},\vec a_{\infty};E_1\bigr)
\times_{L_{3}(b)}
\mathcal M\bigl(L_{3};\vec a^2_{3};E_2\bigr).\label{form98191}
\end{gather}
Here $\vec a^1_{2}$, $\vec a^2_{2}$ and $\vec a^1_{3}$, $\vec a^2_{3}$
are defined in the same way as $\vec a^1_{1}$, $\vec a^2_{1}$.
\begin{figure}[ht]
\centering
\includegraphics[scale=0.3]{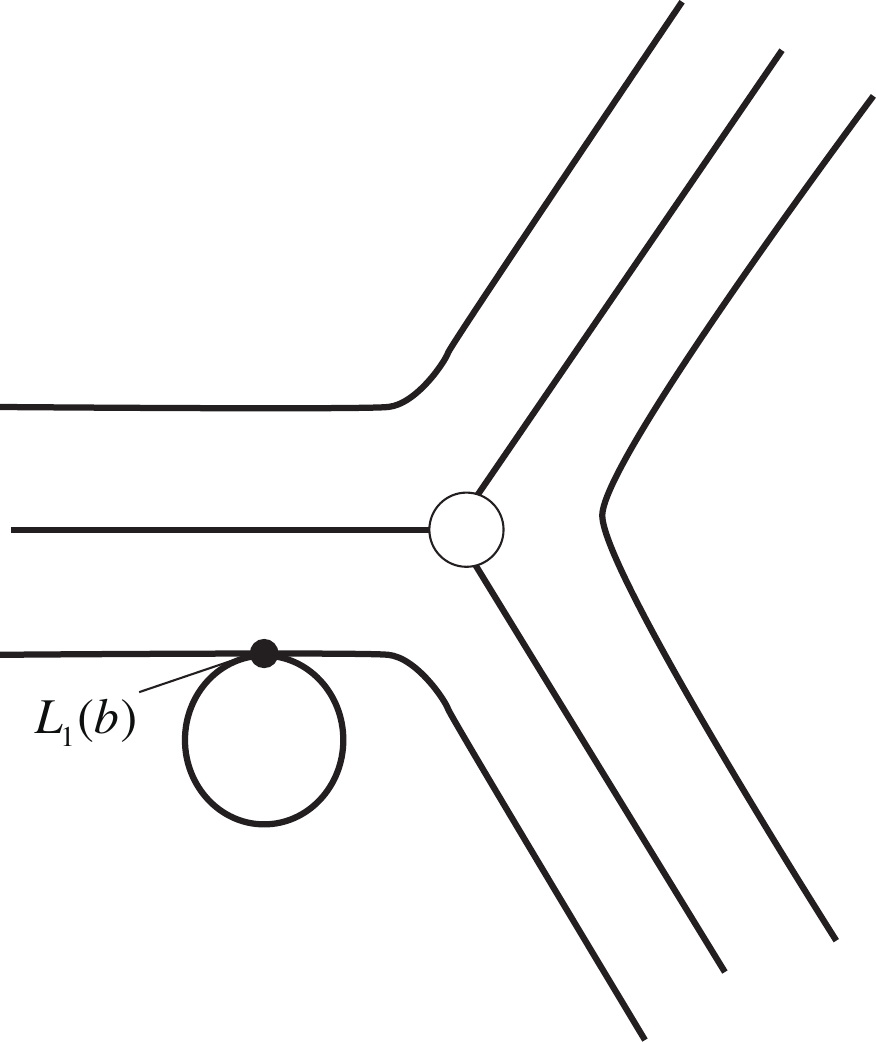}
\caption{Boundary of type (II).}
\label{Figure95}
\end{figure}
\item[(III)]
The third type of boundary corresponds to the limit where the domain
will split into two parts on the image of $\phi_{123}$.
It is described by the fiber product below.
Let $j_{ii'} \in \{0,\dots,k_{ii'}\}$.
We put
$\vec a_{ii'}^1= (a_{ii',1},\dots,a_{ii',j_{ii'}})$,
$\vec a_{ii'}^2= (a_{ii',j_{ii'}+1},\dots,a_{ii',k_{ii'}})$
for $(ii') = (12)$ or $(12)$.
We also put
$\vec a_{13}^2= (a_{13,1},\dots,a_{13,j_{13}})$,
$\vec a_{13}^1= (a_{13,j_{13}+1},\dots,a_{13,k_{13}})$.
Note that in case $j_{ii'} = 0$ (resp.\ $j_{ii'} = k_{ii'}$), $\vec a_{ii'}^1 = \varnothing$
\big(resp.\ $\vec a_{ii'}^2 = \varnothing$\big) for $(ii') = (12)$ or~$(13)$
(the case of $(ii') = (13)$ is similar):
\begin{gather*}
{\mathcal M}_{\rm Y}\bigl(\vec a^2_{12},\vec a^2_{23},\vec a^2_{13};\vec a_{1},\vec a_{2},\vec a_{3};a,\vec a_{\infty};E_2\bigr)\nonumber\\
\qquad
\times_{R_{123}(a)}
{\mathcal M}_{\rm DR}\bigl(\vec a^1_{12},\vec a^1_{23},\vec a^1_{13};a_{\infty,123},a;E_1\bigr), 
\end{gather*}
where $E_1 + E_2 = E$ and $a \in \mathcal A_{123}$.
See Figure~\ref{Figure96}.
\begin{figure}[ht]
\centering
\includegraphics[scale=0.3]{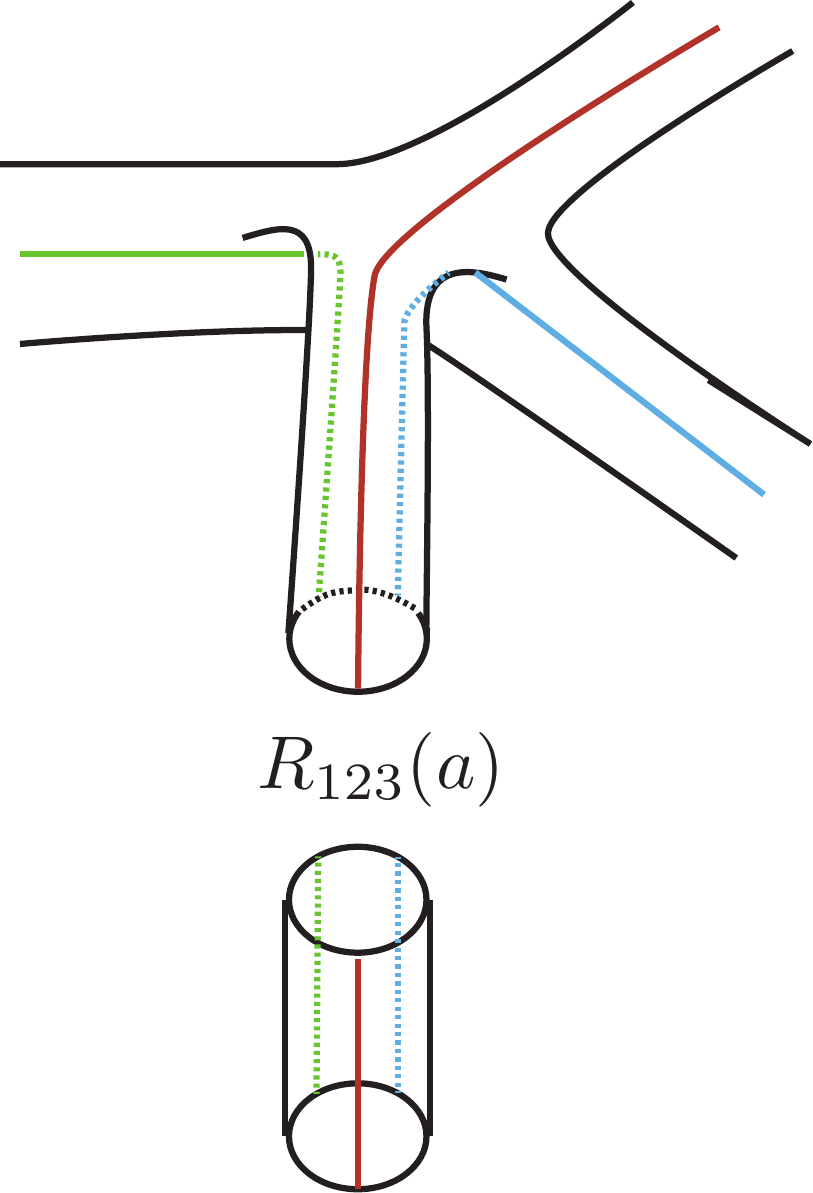}
\caption{Boundary of type (III).}
\label{Figure96}
\end{figure}
\item[(IV)]
The fourth type of boundary corresponds to the limit where the domain
will split into two parts on the image of $\phi_{ii'}$.
It is described by the fiber product below. We
consider the case of $\phi_{12}$.
Let $j \in \{0,\dots,k_{12}\}$, $j_1 \in \{0,\dots,k_{1}\}$, $j_2 \in \{0,\dots,k_{2}\}$.
We put
$\vec a_{12}^1= (a_{12,1},\dots,a_{12,j_{12}})$,
$\vec a_{12}^2= (a_{12,j_{12}+1},\dots,a_{12,k_{12}})$.
$\vec a_{i}^1= (a_{i,1},\dots,a_{i,j_{i}})$,
$\vec a_{i}^2= (a_{i,j_{i}+1},\dots,a_{i,k_{i}})$ for $i=1,2$:
\begin{gather}
{\mathcal M}_{\rm QT}\bigl(\vec a^1_{12},\vec a^1_{1},\vec a^1_{2};a_{\infty,12},a;E_1\bigr)\nonumber \\
\qquad\times_{L_{12}(a)}
{\mathcal M}_{\rm Y}\bigl(\vec a^2_{12},\vec a_{23},\vec a_{13};\vec a^2_{1},\vec a^2_{2},\vec a_{3};(a,a_{\infty,13},a_{\infty,31}),a_{\infty,123};E_2\bigr),\label{form99820}
\end{gather}
where $E_1 + E_2 = E$ and $a \in \mathcal A_{L_{12}}$.
See Figure~\ref{Figure97}.
The cases of $\phi_{23}$ and $\phi_{13}$ are described by the next fiber products:
\begin{gather}
{\mathcal M}_{\rm QT}\bigl(\vec a^1_{23},\vec a^1_{2},\vec a^1_{3};a_{\infty,23},a;E_1\bigr) \nonumber\\
\qquad\times_{L_{12}(a)}
{\mathcal M}_{\rm Y}\bigl(\vec a_{12},\vec a^2_{23},\vec a_{13};\vec a_{1},\vec a^2_{2},\vec a^2_{3};a_{\infty,123},(a_{\infty,12},a,a_{\infty,13});E_2\bigr),\label{form99821}
\\
 {\mathcal M}_{\rm Y}\bigl(\vec a_{12},\vec a_{23},\vec a^2_{13};\vec a^2_{1},\vec a_{2},\vec a^2_{3};
a_{\infty,123},(a_{\infty,12},a_{\infty,23},a);E_2\bigr)\nonumber\\
\qquad\times_{L_{13}(a)}
{\mathcal M}_{\rm QT}\bigl(\vec a^1_{13},\vec a^1_{1},\vec a^1_{3};a_{\infty,23},a;E_1\bigr).\label{form99822}
\end{gather}
Note that $\vec a^1_{23}$ and $\vec a^2_{23}$ is defined in the same way
as $\vec a^1_{12}$ and $\vec a^2_{12}$.
We define $\vec a_{13}^2= (a_{13,1},\dots,a_{13,j_{13}})$,
$\vec a_{13}^1= (a_{13,j_{13}+1},\dots,a_{13,k_{13}})$.
\begin{figure}[ht]
\centering
\includegraphics[scale=0.3]{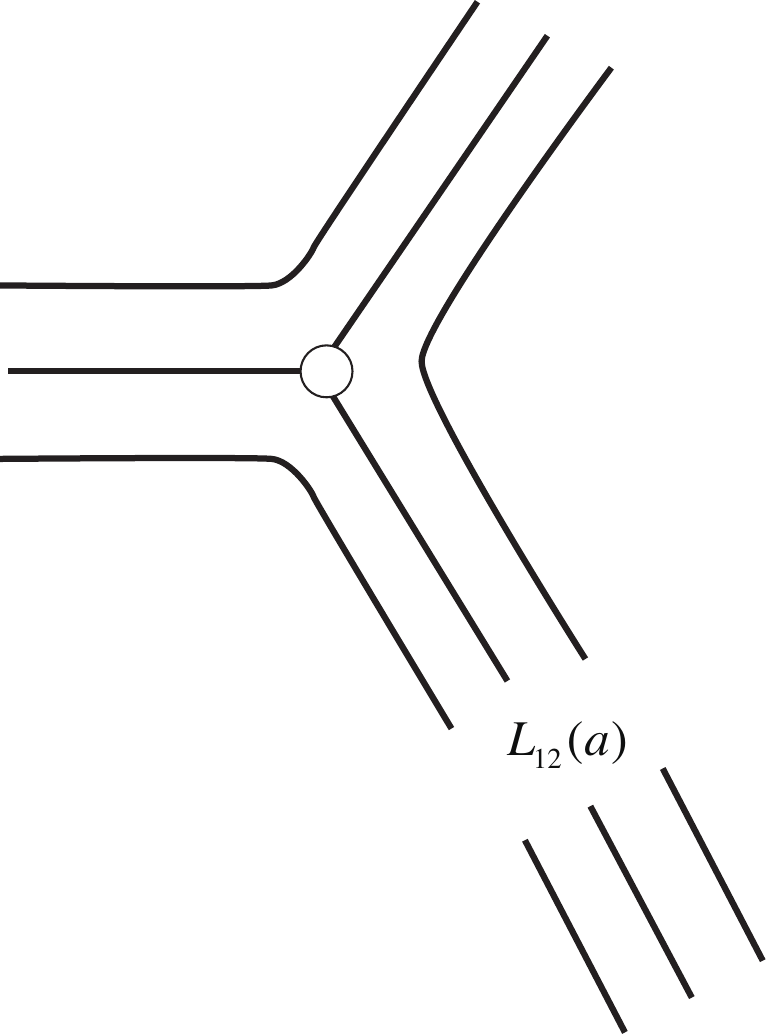}
\caption{Boundary of type (IV).}
\label{Figure97}
\end{figure}
\end{enumerate}
We will show item (3) of Proposition~\ref{prop9911} in Section~\ref{oriYdiagarm}.
We observe that the four types of the boundaries are described
by the fiber products explained above.
In the case of boundaries of types (I), (II), (IV), we only need to
check that the order of the marked points in the
moduli space of Y-diagrams coincides with those of previously defined
moduli spaces. We remark that the boundary of types (IV) with $(ii') = (13)$
the map $\phi_{13}$ identifies the domain of $Y$-diagram with $[-1,1] \times [0,\infty)$,
and for other $(ii')$
the map $\phi_{ii'}$ identifies the domain of $Y$-diagram with~${[-1,1] \times (-\infty,
0]}$
(see Condition \ref{Cond95}\,(2)).
Taking this fact into account the above mentioned coincidence of the
order of marked points is correct in this case also.

In the case of boundaries of type (III), we also remark that the
map $\phi_{123}$ is anti-holomorphic.
So the $J_{X_i}$ holomorphic map on the intersection of $\Omega_i$ with
the image of $\phi_{123}$ will become a~$-J_{X_i}$ holomorphic map on an open subset of the drum appearing
in Section~\ref{sec:comp}.

Once we observe these points, the proof of Proposition~\ref{prop9911} is now
a routine.

\begin{prop}\label{prop9812}
For each $E_0$,
there exists a system of CF-perturbations \smash{$\widehat{\mathfrak S}$} on
\[
{\mathcal M}_{\rm Y}(\vec a_{12},\vec a_{23},\vec a_{13};\vec a_{1},\vec a_{2},\vec a_{3},a_{\infty,-},\vec a_{\infty,+};E)
\]
$($with respect to Kuranishi structures which are outer collarings of thickenings of
those in Proposition~{\rm\ref{prop9911}}$)$
for $E < E_0$ such that the following holds:
\begin{enumerate}\itemsep=0pt
\item[$(1)$]
They are transversal to $0$.
\item[$(2)$]
The evaluation map
${\rm ev}_{\infty}$ is strongly submersive with respect to this
CF-perturbation.
\item[$(3)$]
The CF-perturbations are compatible with the description of the
boundary. Namely, the restriction of the CF-perturbation on the
boundary coincides with the fiber product CF-perturbation
in the sense of {\rm\cite[\emph{Lemma--Definition} 10.6]{foootech2} \emph{and} \cite{fooonewbook}}.
\item[$(4)$]
The CF-perturbations are compatible with the forgetful maps of the boundary
marked points corresponding to the diagonal component,
in the sense of {\rm\cite[\emph{Theorem} 5.1]{fooo091}}.
\end{enumerate}

\end{prop}
The proof is the same as Proposition~\ref{prop536} and is now a routine. We omit it.

The next step is to rewrite a geometric result
Proposition~\ref{prop9911} to an algebraic one.
This is a~process we have done in Sections \ref{subsec:Ainfalgim},
\ref{subsec:bi-functorgeo1} and \ref{subsec:functorconst}
as well as several other references especially in \cite[Part II]{fooonewbook}
and proceed as follows.
We regard the evaluation map
${\rm ev}_{13}$ as an `output' and other evaluation maps
as `input's.
In other words, we take differential forms on the targets
of the evaluation maps other than ${\rm ev}_{13}$,
we then pull them back to the moduli space in
Proposition~\ref{prop9911} and use the
CF-perturbation of Proposition~\ref{prop9812} to push it out to the target of
the evaluation map
${\rm ev}_{13}$.
We thus obtain a map
between de Rham complexes.
It will be the Y-diagram
transformation
below
\begin{gather}
\mathscr{YT}^{E,\varepsilon}_{k_{12},k_{23},k_{13};k_1,k_2,k_3}
\colon\
CF(L_{13};L_{12},L_{23}) \nonumber \\
\qquad{}\otimes B_{k_{12}}CF[1](L_{12}) \otimes B_{k_{23}}CF[1](L_{23})\nonumber
\otimes B_{k_{13}}CF[1](L_{13})
\\
\qquad{}\otimes
CF(L_1,L_{12};L_{2}) \otimes CF(L_2,L_{23};L_{3}) \nonumber\\
\qquad{}\otimes B_{k_1}CF[1](L_1) \otimes B_{k_2}CF[1](L_2)
\otimes B_{k_3}CF[1](L_3)
\to
CF(L_1,L_{13};L_{3}).\label{operation923}
\end{gather}
See \eqref{form923}.
Note that we can find the domain and codomain of the map \eqref{operation923}
by inspecting the targets of the evaluation maps of various kinds.

To obtain the basic property of the map \eqref{operation923}
we use Stokes' theorem and the composition formula as follows.
We consider the commutator of the map \eqref{operation923}
and the de Rham differential.
Stokes' theorem implies that the commutator is
equal to the map obtained from the
boundary of the moduli spaces of Proposition~\ref{prop9911}
in the same way as we obtain the map \eqref{operation923}.
We have described the boundary of the moduli space in
Proposition~\ref{prop9911} and
found that the boundary consists of four types of fiber products.
Actually each of types (I), (II), (IV) is a union of three
kinds of boundaries. In the case of type (I) it is a union
of components
corresponding to three kinds of disk bubbles, that are, those
at $L_{12}$, $L_{23}$, and $L_{13}$.
In the case of type (II) it is a union of components
corresponding to three kinds of disk bubbles, that are those at $L_{1}$, $L_{2}$,
and $L_{3}$.
In the case of type (IV) it is a union three different ends,
where strips escape at the image of~$\phi_{12}$,~$\phi_{23}$, or $\phi_{13}$.
Thus the formula \eqref{form925} contains ten terms corresponding to those different
kinds of boundaries.

Note that each boundary component is described as the fiber product
of a moduli space of Proposition~\ref{prop9911}
(whose energy is not greater than $E$) and another moduli space.
In the case of type (I), the another moduli space is one
we used to define the filtered $A_{\infty}$ category associated to $L_{ij}$.
In the case of type (II), the another moduli space is one
we used to define the filtered~$A_{\infty}$ category associated to $L_{i}$.
In the case of type (III), the another moduli space is
the moduli space of pseudo-holomorphic drums.
In the case of type (IV), the another moduli space is one
we used to define the filtered $A_{\infty}$ tri-module associated to $L_i$, $L_{ij}$, $L_{j}$.

Therefore, by the composition formula (see \cite[Theorem~10.21]{fooonewbook}),
the terms corresponding to those 4 types of boundary
components are obtained as compositions
of the map \eqref{operation923} (whose energy is smaller than $E$)
and one of the following:
a~map induced from the structure operations of the filtered $A_{\infty}$ category associated to $L_{ij}$;
a~map induced from the structure operations of the filtered $A_{\infty}$ category associated to $L_{i}$;
a~map induced from the structure operations of the filtered tri-module
$\mathscr{CF}(\mathbb L_{13};\mathbb L_{12},\mathbb L_{23})$;
a~map induced from the structure operations of the filtered tri-module
$\mathscr{CF}(\mathbb L_{i},\mathbb L_{ij};\mathbb L_{j})$.

The formula we obtain in this way is \eqref{form925} in Proposition~\ref{basiceqYdiagram}.

This process to go from geometry to algebra is straightforward and
is now becoming a routine.
Since the formula is long (contains many terms),
let us first describe it in a simple case and explain how it will be used
in this simple case.

We assume that $L_1$, $L_3$, $L_{ij}$ are embedded and monotone.
Suppose that $L_2$ is a union of embedded monotone Lagrangian
submanifolds $L_2^i$,
$i=0,\dots,k$,
which intersects transversally each other.
We consider the case when there is no marked points
which maps to $L_1$, $L_3$ or $L_{ij}$.
We use the cyclic element ${\bf 1}_{123}$ (that is
the function $1$ on the diagonal component)
and insert it at the hole in the middle of the Y-diagram.
The map \eqref{operation923} in this case becomes
\begin{gather}
CF\bigl(L_1,L_{12};L^0_{2}\bigr) \otimes
\bigotimes_{i=1}^k CF[1]\bigl(L_2^{i-1},L_2^{i}\bigr)
\otimes CF\bigl(L^k_2,L_{23};L_{3}\bigr)
\to
CF(L_1,L_{13};L_{3}).\label{form923923111}
\end{gather}
We recall that the tri-module $CF(L_i,L_{ij};L_{j})$
is used to define the filtered $A_{\infty}$
functor $\mathcal W_{\mathcal L_{ij}}$ via Yoneda functor.
In the simplified case we are discussing, we fix $L_{ij}$ and put no marked points
on the seam. So it is actually a bi-module.
Thus the right-hand side of \eqref{form923923111}
corresponds to the filtered $A_{\infty}$
functor $\mathcal W_{\mathcal L_{13}}$.

The direct sum of the left-hand side of \eqref{form923923111}
for various $L_2^0,\dots,L_2^{k}$
becomes the derived tensor product
$
\mathfrak{ten}(CF(L_1,L_{12};L^*_{2}),CF(L^*_2,L_{23};L_{3}))$.
See Lemma--Definition~\ref{defntensor}.
As we will discuss in Section~\ref{sec:compfunc} (see Proposition~\ref{lem9100}), the derived
tensor product of filtered $A_{\infty}$ bi-module corresponds to the
composition of the corresponding filtered $A_{\infty}$ functors.
Thus the left-hand side of
\eqref{form923923111} corresponds to the
composition $\mathcal W_{\mathcal L_{23}} \circ \mathcal W_{\mathcal L_{12}}$.

We will show that by taking the direct sum over various $L_2^0,\dots,L_2^{k}$
the map \eqref{form923923111} becomes a chain homotopy equivalence
and will use it to show \eqref{compocompfor}.

Actually, we need to include bounding cochains.
We also need to show that the map \eqref{form923923111}
becomes a left-$\mathfrak{Fukst}(X_1)$ and right-$\mathfrak{Fukst}(X_3)$ bi-module homomorphism.
Moreover, we need to show the functoriality when we have several components of $L_{ij}$
and morphisms (an element of Floer's chain complex) from $L_{ij}$ to $L'_{ij}$.
To work these out, we need \eqref{operation923}
and its basic property Proposition~\ref{basiceqYdiagram} in its
full generality.
(This part of the proof is carried out in Section~\ref{sec:compfuncmain} after preparing various
algebraic results.)

We go back to the general case and explain the way to define operations
\eqref{operation923} using Propositions \ref{prop9812} and \ref{prop9911}.

Let $h_{\infty,123} \in \Omega(R(a_{\infty,123}))$,
$h_{ii',j} \in \Omega(L_{ii'}(a_{ii',j})$,
${\bf h}_{ii'} = (h_{ii',1},\dots,h_{ii',k_{ii'}})$ ($ii' = 12, 23$ or $13$),
$h_{\infty,ii'} \in \Omega(L_{ii'}(a_{\infty,ii'}))$ ($ii' = 12$ or $23$),
$h_{i,j} \in \Omega(L_{i}(a_{i,j})$,
${\bf h}_{i} = (h_{i,1},\dots,h_{i,k_{i}})$ ($i = 1,2$ or $3$).
Then the $\Omega(a_{\infty,13})$ component of
\[
\mathscr{YT}^{E,\varepsilon}_{k_{12},k_{23},k_{13};k_1,k_2,k_3}
(h_{\infty,123};{\bf h}_{12},{\bf h}_{23},{\bf h}_{13};h_{\infty,12},h_{\infty,23};
{\bf h}_{1},{\bf h}_{2},{\bf h}_{3})
\]
is by definition
\begin{gather}
{\rm ev}_{\infty,13}!
\bigl(
 {\rm ev}_{\infty,123}^*h_{\infty,123}
\wedge
{\rm ev}_{12}^*{\bf h}_{12}
\wedge
{\rm ev}_{23}^*{\bf h}_{23}
\wedge
{\rm ev}_{13}^*{\bf h}_{13}\nonumber
\\
\qquad{}\wedge
{\rm ev}_{\infty,12}^*h_{\infty,12}
\wedge
{\rm ev}_{\infty,23}^*h_{\infty,23}
\wedge
{\rm ev}_{1}^*{\bf h}_{1}
\wedge
{\rm ev}_{2}^*{\bf h}_{2}
\wedge
{\rm ev}_{3}^*{\bf h}_{3}
; \widehat{\mathfrak S^{\varepsilon}}\bigr).\label{form923}
\end{gather}
Here the integration along the fiber
appearing in the formula \eqref{form923} is taken on the
moduli space~${{\mathcal M}_Y(\vec a_{12},\vec a_{23},\vec a_{13};\vec a_{1},\vec a_{2},\vec a_{3},a_{\infty,123},\vec a_{\infty};E)}$
using the CF-perturbation \smash{$\widehat{\mathfrak S^{\varepsilon}}$}.
\footnote{The sign is discussed in Section~\ref{oriYdiagarm}.}
We then put{\samepage
\[
\mathscr{YT}^{<E_0,\varepsilon}_{k_{12},k_{23},k_{13};k_1,k_2,k_3}
:=
\sum_{E<E_0}
T^E\mathscr{YT}^{E,\varepsilon}_{k_{12},k_{23},k_{13};k_1,k_2,k_3}.
\]
We call \smash{$\mathscr{YT}^{<E_0,\varepsilon}_{k_{12},k_{23},k_{13};k_1,k_2,k_3}$} the \index[syindex]{YTscrE0ep@$\mathscr{YT}^{<E_0,\varepsilon}_{k_{12},k_{23},k_{13};k_1,k_2,k_3}$}
{\it Y diagram transformation}.\index{Y diagram transformation}}

We usually omit the indices $k_{12}$, $k_{23}$, $k_{13}$; $k_1$, $k_2$, $k_3$
above since it is determined automatically from the
input.
\begin{rem}
The order how the variables appears in \eqref{form923}
does not coincide with the order of the tensor factors in \eqref{operation923}.
The former coincides with
\begin{gather*}
CF(L_{13};L_{12},L_{23}) \otimes BCF[1](L_{12})\otimes BCF[1](L_{23})\otimes BCF[1](L_{13})\\
\qquad{}\otimes CF(L_1,L_{12};L_2) \otimes CF(L_2,L_{23};L_3) \otimes
BCF[1](L_1) \otimes BCF[1](L_2) \otimes BCF[1](L_3).
\end{gather*}
The formula looks easier to read when written in this order.
The order of the tensor factors of~\eqref{form923} is
one the $Y$-diagram transformation will be applied
in Section~\ref{sec:compfuncmain}.
\end{rem}
\begin{prop}\label{basiceqYdiagram}
The Y diagram transformation
\smash{$\mathscr{YT}^{<E_0,\varepsilon}_{k_{12},k_{23},k_{13};k_1,k_2,k_3}$}
satisfies the following congruence:
\begin{gather}
(-1)^{*_1}\mathscr{YT}^{<E_0,\varepsilon}
\bigl(h_{\infty,123};\hat d({\bf h}_{12}),{\bf h}_{23},{\bf h}_{13};h_{\infty,12},h_{\infty,23};
{\bf h}_{1},{\bf h}_{2},{\bf h}_{3}\bigr)\nonumber
\\
\qquad{}
+ (-1)^{*_2}\mathscr{YT}^{<E_0,\varepsilon}
\bigl(h_{\infty,123};{\bf h}_{12},\hat d({\bf h}_{23}),{\bf h}_{13};h_{\infty,12},h_{\infty,23};
{\bf h}_{1},{\bf h}_{2},{\bf h}_{3}\bigr)\nonumber
\\
\qquad{}
+ (-1)^{*_3}\mathscr{YT}^{<E_0,\varepsilon}
\bigl(h_{\infty,123};{\bf h}_{12},{\bf h}_{23},\hat d({\bf h}_{13});h_{\infty,12},h_{\infty,23};
{\bf h}_{1},{\bf h}_{2},{\bf h}_{3}\bigr)\nonumber
\\
\qquad{}
+ (-1)^{*_4}\mathscr{YT}^{<E_0,\varepsilon}
\bigl(h_{\infty,123};{\bf h}_{12},{\bf h}_{23},{\bf h}_{13};h_{\infty,12},h_{\infty,23};
\hat d({\bf h}_{1}),{\bf h}_{2},{\bf h}_{3}\bigr)\nonumber
\\
\qquad{}
+ (-1)^{*_5}\mathscr{YT}^{<E_0,\varepsilon}
\bigl(h_{\infty,123};{\bf h}_{12},{\bf h}_{23},{\bf h}_{13};h_{\infty,12},h_{\infty,23};
{\bf h}_{1},\hat d({\bf h}_{2}),{\bf h}_{3}\bigr)\nonumber
\\
\qquad{}
+ (-1)^{*_6}\mathscr{YT}^{<E_0,\varepsilon}
\bigl(h_{\infty,123};{\bf h}_{12},{\bf h}_{23},{\bf h}_{13};h_{\infty,12},h_{\infty,23};
{\bf h}_{1},{\bf h}_{2},\hat d({\bf h}_{3})\bigr)\nonumber
\\
\qquad{}
+\sum_{c_{12},c_{23},c_{13}}
(-1)^{*_7}
\mathscr{YT}^{<E_0,\varepsilon}
\bigl(\mathfrak{n}^{<E_0,\varepsilon}
\bigl({\bf h}^{c_{13};2}_{13};h_{\infty,123};{\bf h}^{c_{12};1}_{12},{\bf h}^{c_{23};1}_{23}\bigr);\nonumber
\\
\qquad\qquad\quad\phantom{+}{}
{\bf h}^{c_{12};2}_{12},{\bf h}^{c_{23};2}_{23},{\bf h}^{c_{13};1}_{13};h_{\infty,12},h_{\infty,23};
;
{\bf h}_{1},{\bf h}_{2},{\bf h}_{3}\bigr)\nonumber
\\
\qquad{}
+\sum_{c_{1},c_{2},c_{12}}
(-1)^{*_8}
\mathscr{YT}^{<E_0,\varepsilon}
\bigl(h_{\infty,123};{\bf h}^{c_{12};1}_{12},{\bf h}_{23},{\bf h}_{13};\nonumber
\\
\qquad\qquad\quad\phantom{+}{}
\mathfrak{n}^{<E_0,\varepsilon}
\bigl({\bf h}^{c_{1};1}_{1},{\bf h}^{c_{12};2}_{12};h_{\infty,12};{\bf h}^{c_{2};1}_{2}\bigr),
h_{\infty,23};{\bf h}^{c_{1};2}_{1},{\bf h}^{c_{2};2}_{2},{\bf h}_{3}\bigr)\nonumber
\\
\qquad{}
+\sum_{c_{2},c_{3},c_{23}}
(-1)^{*_9}
\mathscr{YT}^{<E_0,\varepsilon}
\bigl(h_{\infty,123};{\bf h}_{12},{\bf h}^{c_{23};1}_{23},{\bf h}_{13};h_{\infty,12};\nonumber\\
\qquad\qquad\quad\phantom{+}{}
\mathfrak{n}^{<E_0,\varepsilon}
\bigl({\bf h}^{c_{2};1}_{2},{\bf h}^{c_{23};2}_{23};h_{\infty,23};{\bf h}^{c_{3};1}_{3}\bigr)
;{\bf h}_{1},{\bf h}^{c_{2};2}_{2},{\bf h}^{c_{3};2}_{3}\bigr)\nonumber \\
\qquad{}
-\sum_{c_{1},c_{3},c_{13}}
(-1)^{*_{10}}
\mathfrak{n}^{<E_0,\varepsilon}
\bigl({\bf h}^{c_{1};1}_{1},{\bf h}^{c_{13};1}_{13};
\mathscr{YT}^{<E_0,\varepsilon}
\bigl(h_{\infty,123};{\bf h}_{12},{\bf h}_{23},{\bf h}^{c_{13};2}_{13};\nonumber
\\
\qquad\qquad\quad\phantom{+}{}
h_{\infty,12},h_{\infty,23};
{\bf h}^{c_{1};2}_{1},{\bf h}_{2},{\bf h}^{c_{3};2}_{3}\bigr);{\bf h}^{c_{3};1}_{3}\bigr)\equiv 0 \mod T^{E_0}.\label{form925}
\end{gather}
Here \smash{$\Delta({\bf h}_i) = \sum_{c_i} {\bf h}^{c_i,1}_i \otimes {\bf h}^{c_i,2}_i$},
\smash{$\Delta({\bf h}_{ii'}) = \sum_{c_{ii'}} {\bf h}^{c_{ii'},1}_{ii'} \otimes {\bf h}^{c_{ii'},2}_{ii'}$}
and all the signs are by Koszul rule.\footnote{See Section~\ref{Koszul}
for the way the Koszul rule determines the sign.}

\end{prop}
\begin{proof}
The proof uses Propositions \ref{prop9911} and \ref{prop9812} together
with Stokes' theorem (see \cite[Proposition 9.26]{foootech2} and \cite{fooonewbook})
and the composition formula (see \cite[Theorem 10.20]{foootech2} and \cite{fooonewbook}).
It goes in the same way as the proofs of other similar statements
we proved before.
In fact, the first three terms correspond to the boundary
of type (I) and the fiber products
\eqref{form9817}, \eqref{form9818}, \eqref{form9819},
respectively.
The operator $\hat d$ in the first three terms
are induced by the structure operations of
$\mathfrak{Fuk}(-X_i \times X'_i)$.

The 4-th, 5-th and 6-th terms correspond to the
boundary of type (II) and the fiber products~\eqref{form98171}, \eqref{form98181} and \eqref{form98191},
respectively.
The operator $\hat d$ in the 4-th, 5-th and 6-th terms
are induced by the structure operations of
$\mathfrak{Fuk}(\mathbb L_{i})$.

The 7-th term corresponds to the boundary of type (III)
and the fiber product \eqref{form99820}.
Note that the structure map $\mathfrak n$ appearing
in the 7-th term is one of the tri-module
$\mathscr{CF}(\mathbb L_{13};\mathbb L_{12},\mathbb L_{23})$.

The 8-th, 9-th and 10-th terms
correspond to the boundary of type (IV)
and the fiber products
\eqref{form99820}, \eqref{form99821},
and \eqref{form99822}, respectively.
The structure map $\mathfrak n$ appearing in the 8-th, 9-th and 10-th terms
is structure operation of the tri-module
$\mathscr{CF}(\mathbb L_{i},\mathbb L_{ij};\mathbb L_{j})$.

The sign will be discussed in Section~\ref{oriYdiagarm}.
\end{proof}

In the same way as Definition~\ref{defn2333}\,(8), we can modify our operations
and change the congruence in \eqref{form925}
to the equality.
Namely, we have the following.
\begin{prop}\label{prop910}
There exists a map
\begin{gather}
\mathscr{YT} \colon\
CF(L_{12},L_{23},L_{13}) \otimes BCF[1](L_{12}) \otimes BCF[1](L_{23})
\otimes BCF[1](L_{13})\nonumber
\\
\qquad{}\otimes
CF(L_1,L_{12},L_{2}) \otimes CF(L_2,L_{23},L_{3}) \nonumber\\
\qquad{}\otimes BCF[1](L_1) \otimes \otimes BCF[1](L_2)
\otimes BCF[1](L_3)
\to
CF(L_1,L_{13},L_{3})\label{form926}
\end{gather}
such that if we replace $\mathscr{YT}^{<E_0,\varepsilon}$
by $\mathscr{YT}$ the formula \eqref{form925} holds
as an exact equality.
Namely,
\begin{gather}
(-1)^{*_1}\mathscr{YT}
\bigl(h_{\infty,123};\hat d({\bf h}_{12}),{\bf h}_{23},{\bf h}_{13};h_{\infty,12},h_{\infty,23};
{\bf h}_{1},{\bf h}_{2},{\bf h}_{3}\bigr)\nonumber
\\
\qquad{}
+ (-1)^{*_2}\mathscr{YT}
\bigl(h_{\infty,123};{\bf h}_{12},\hat d({\bf h}_{23}),{\bf h}_{13};h_{\infty,12},h_{\infty,23};
{\bf h}_{1},{\bf h}_{2},{\bf h}_{3}\bigr)\nonumber
\\
\qquad{}
+ (-1)^{*_3}\mathscr{YT}
\bigl(h_{\infty,123};{\bf h}_{12},{\bf h}_{23},\hat d({\bf h}_{13});h_{\infty,12},h_{\infty,23};
{\bf h}_{1},{\bf h}_{2},{\bf h}_{3}\bigr)\nonumber
\\
\qquad{}
+ (-1)^{*_4}\mathscr{YT}
\bigl(h_{\infty,123};{\bf h}_{12},{\bf h}_{23},{\bf h}_{13};h_{\infty,12},h_{\infty,23};
\hat d({\bf h}_{1}),{\bf h}_{2},{\bf h}_{3}\bigr)\nonumber
\\
\qquad{}
+ (-1)^{*_5}\mathscr{YT}
\bigl(h_{\infty,123};{\bf h}_{12},{\bf h}_{23},{\bf h}_{13};h_{\infty,12},h_{\infty,23};
{\bf h}_{1},\hat d({\bf h}_{2}),{\bf h}_{3}\bigr)\nonumber
\\
\qquad{}
+ (-1)^{*_6}\mathscr{YT}
\bigl(h_{\infty,123};{\bf h}_{12},{\bf h}_{23},{\bf h}_{13};h_{\infty,12},h_{\infty,23};
{\bf h}_{1},{\bf h}_{2},\hat d({\bf h}_{3})\bigr)\nonumber
\\
\qquad{}
+\sum_{c_{12},c_{23},c_{13}}
(-1)^{*_7}
\mathscr{YT}
\bigl(\mathfrak{n}
\bigl({\bf h}^{c_{13};2}_{13};h_{\infty,123};{\bf h}^{c_{12};1}_{12},{\bf h}^{c_{23};1}_{23}\bigr);\nonumber
\\
\qquad\qquad\quad\phantom{+}
{\bf h}^{c_{12};2}_{12},{\bf h}^{c_{23};2}_{23},{\bf h}^{c_{13};1}_{13};h_{\infty,12},h_{\infty,23};
;
{\bf h}_{1},{\bf h}_{2},{\bf h}_{3}\bigr)\nonumber
\\
\qquad{}
+\sum_{c_{1},c_{2},c_{12}}
(-1)^{*_8}
\mathscr{YT}
\bigl(h_{\infty,123};{\bf h}^{c_{12};1}_{12},{\bf h}_{23},{\bf h}_{13};
\mathfrak{n}
\bigl({\bf h}^{c_{1};1}_{1},{\bf h}^{c_{12};2}_{12};h_{\infty,12};{\bf h}^{c_{2};1}_{2}\bigr),\nonumber
\\
\qquad\qquad\quad\phantom{+}
h_{\infty,23};{\bf h}^{c_{1};2}_{1},{\bf h}^{c_{2};2}_{2},{\bf h}_{3}\bigr)\nonumber
\\
\qquad{}
+\sum_{c_{2},c_{3},c_{23}}
(-1)^{*_9}
\mathscr{YT}
\bigl(h_{\infty,123};{\bf h}_{12},{\bf h}^{c_{23};1}_{23},{\bf h}_{13};h_{\infty,12};\nonumber\\
\phantom{\qquad+}
\mathfrak{n}
\bigl({\bf h}^{c_{2};1}_{2},{\bf h}^{c_{23};2}_{23};h_{\infty,23};{\bf h}^{c_{3};1}_{3}\bigr)
;{\bf h}_{1},{\bf h}^{c_{2};2}_{2},{\bf h}^{c_{3};2}_{3}\bigr) \nonumber\\
\qquad{}
-\sum_{c_{1},c_{3},c_{13}}
(-1)^{*_{10}}
\mathfrak{n}
\bigl({\bf h}^{c_{1};1}_{1},{\bf h}^{c_{13};1}_{13};
\mathscr{YT}
\bigl(h_{\infty,123};{\bf h}_{12},{\bf h}_{23},{\bf h}^{c_{13};2}_{13};\nonumber
\\
\qquad\qquad\quad\phantom{+}
h_{\infty,12},h_{\infty,23};
{\bf h}^{c_{1};2}_{1},{\bf h}_{2},{\bf h}^{c_{3};2}_{3}\bigr);{\bf h}^{c_{3};1}_{3}\bigr)=0.\label{form9252}
\end{gather}
Moreover, $\mathscr{YT} \equiv \mathscr{YT}^{<E_0,\varepsilon} \mod T^{E_0}$.

\end{prop}
We call $\mathscr{YT}$, the {\it Y diagram transformation} also.\index{Y diagram transformation}

\subsection[Proof of Proposition~\ref{prop912}\,(2)]{Proof of Proposition~\ref{prop912}\,(2)}
\label{subsec:compatiobj}

In this subsection, we prove Proposition~\ref{prop912}\,(2).

Let $\mathcal L_{12} = (L_{12},\sigma_{12},b_{12})$
(resp.\ $\mathcal L_{23} = (L_{23},\sigma_{23},b_{23})$)
be an object of $\mathfrak{Fukst}(-X_1\times X_2)$
(resp.~${\mathfrak{Fukst}(-X_2\times X_3)}$).
Let $\mathcal L_{13} = (L_{13},\sigma_{13},b_{13})$
be the geometric composition
$\mathcal L_{23} \circ \mathcal L_{12}$.

Let $\mathcal L_1 = (L_{1},\sigma_{1},b_{1})$
and we put
\[
\mathcal L_2 = (L_{2},\sigma_{2},b_{2}) = \mathcal W_{\mathcal L_{12}}(\mathcal L_1),
\quad
\mathcal L_3^{(1)} = \bigl(L_{3},\sigma_{3},b^{(1)}_{3}\bigr) = \mathcal W_{\mathcal L_{23}}(\mathcal L_2)
\]
and
\[
\mathcal L_3^{(2)} = \bigl(L_{3},\sigma_{3},b^{(2)}_{3}\bigr) = \mathcal W_{\mathcal L_{13}}(\mathcal L_1).
\]
We remark that the underlying Lagrangian submanifold of
\smash{$\mathcal L_3^{(1)}$} is equal to the underlying Lagrangian submanifold
of \smash{$\mathcal L_3^{(2)}$}.
This is obvious since
\[
L_2 \times_{X_2} L_{23} =
L_1 \times_{X_1} L_{12} \times_{X_2} L_{23}
=
L_1 \times_{X_1} L_{13}.
\]
The coincidence of the relative spin structure is the main part of Proposition~\ref{prop912}\,(1)
which we will prove in Section~\ref{oriYdiagarm}.
We will prove in this subsection the next proposition.

\begin{prop}\label{prop911}
The bounding cochain
\smash{$b^{(1)}_{3}$} is gauge equivalent to \smash{$b^{(2)}_{3}$} in the sense of
{\rm\cite[\emph{Definition} 4.3.1]{fooobook}}.

\end{prop}
\begin{proof}
We use the next algebraic lemma to prove Proposition~\ref{prop911}.
\begin{lem}\label{lem912}
Let $(D,\{\mathfrak n_k\})$ be a $G$-gapped right filtered $A_{\infty}$ module over
$(C,\{\mathfrak m_k\})$.
Let \smash{${\bf 1}^{(1)}$}, \smash{${\bf 1}^{(2)}$} be cyclic elements of $D$
and \smash{$b^{(1)}$}, \smash{$b^{(2)}$} bounding cochains of $C$ such that
\[
\sum_{k=0}^{\infty} \mathfrak n_k\bigl({\bf 1}^{(i)};b^{(i)},\dots,b^{(i)}\bigr) = 0.
\]
We also assume
\begin{equation}\label{form929}
{\bf 1}^{(1)} \equiv {\bf 1}^{(2)} \mod \Lambda_+.
\end{equation}
Then $b^{(1)}$ is gauge equivalent to $b^{(2)}$.

\end{lem}
\begin{proof}
We use a certain result and notations of \cite{fooobook} in the proof.
Let $\mathfrak C$ be a model of $[0,1] \times C$ in the sense of \cite[Definition~4.2.1]{fooobook}.
Let $\mathfrak D$ be a model of $[0,1] \times D$ in the sense of \cite[Definition~5.2.21]{fooobook},
which is a right $\mathfrak C$ module.
Such $\mathfrak C$ and $\mathfrak D$ exists by \cite[Lemma 4.2.13 and Theorem 5.2.23]{fooobook}.
Since ${\rm Eval}_0 \oplus {\rm Eval}_1 \colon \mathfrak D \to D \oplus D$
is surjective (see \cite[Definition 5.2.23]{fooobook}),
we have
$\Delta {\bf 1} \in \mathfrak D$ such that
$
({\rm Eval}_0)(\Delta {\bf 1}) = 0$,
$
({\rm Eval}_1)(\Delta {\bf 1}) = {\bf 1}^{(2)} - {\bf 1}^{(1)}$.
Using \eqref{form929}, we may choose $\Delta {\bf 1}$ such that
\begin{equation}\label{form930}
\Delta {\bf 1} \equiv 0 \mod \Lambda_+.
\end{equation}
We put
\smash{$
\hat{\bf 1} = {\rm Incl}({\bf 1}^{(1)}) + \Delta {\bf 1}$}.
\eqref{form930} implies that $\hat{\bf 1}$ is a cyclic element of the
right $\mathfrak C$ module $\mathfrak D$.
Therefore, by Proposition~\ref{thm35} there exists a bounding cochain $\hat b$ of $\mathfrak c$ such that
\[
\sum_{k=0}^{\infty}\mathfrak n_k\bigl(\hat{\bf 1};\hat b,\dots,\hat b\bigr) = 0.
\]
We remark that
\smash{$
({\rm Eval}_0)\bigl(\hat{\bf 1}\bigr) = {\bf 1}^{(1)}$}, \smash{$
({\rm Eval}_1)\bigl(\hat{\bf 1}\bigr) = {\bf 1}^{(2)}$}.
Therefore, using the uniqueness part of Proposition~\ref{thm35}
we find that
\smash{$
({\rm Eval}_0)\bigl(\hat{b}\bigr) = b^{(1)}$}, \smash{$
({\rm Eval}_1)\bigl(\hat{b}\bigr) = b^{(2)}$}.
Hence \smash{$b^{(1)}$} is gauge equivalent to \smash{$b^{(2)}$}, as required.
\end{proof}

We go back to our geometric situation and use Y diagram transformation
$\mathscr{YT}$ to
define a~map
\begin{equation}\label{form531}
\mathscr{MY}\colon\
CF(L_1;L_{12};L_2) \otimes CF(L_2;L_{23};L_3) \to CF(L_1;L_{13};L_3)
\end{equation}
by
\[
\mathscr{MY}(h_{\infty,+,12},h_{\infty,+,23})
=
\mathscr{YT}
\bigl({\bf 1}_{123};e^{b_{12}},e^{b_{23}},e^{b_{13}};h_{\infty,+,12},h_{\infty,+,23};
e^{b_{1}},e^{b_{3}},e^{b_{3}^{(1)}}\bigr).
\]
Here ${\bf 1}_{123} \in CF(L_{12};L_{23};L_{13})$ is the cyclic element we used in Lemma~\ref{lem89}.
(See Figure~\ref{Figure107}.)
In other words, it is the function $1$ on the diagonal component of
\[
\bigl(\tilde L_{12} \times \tilde L_{23} \times \tilde L_{13}\bigr) \times_{(X_1\times X_2 \times X_3)^2} \Delta
\cong \tilde L_{13} \times_{X_1 \times X_3}\tilde L_{13}.
\]

\begin{figure}[ht]
\centering
\includegraphics[scale=0.4]{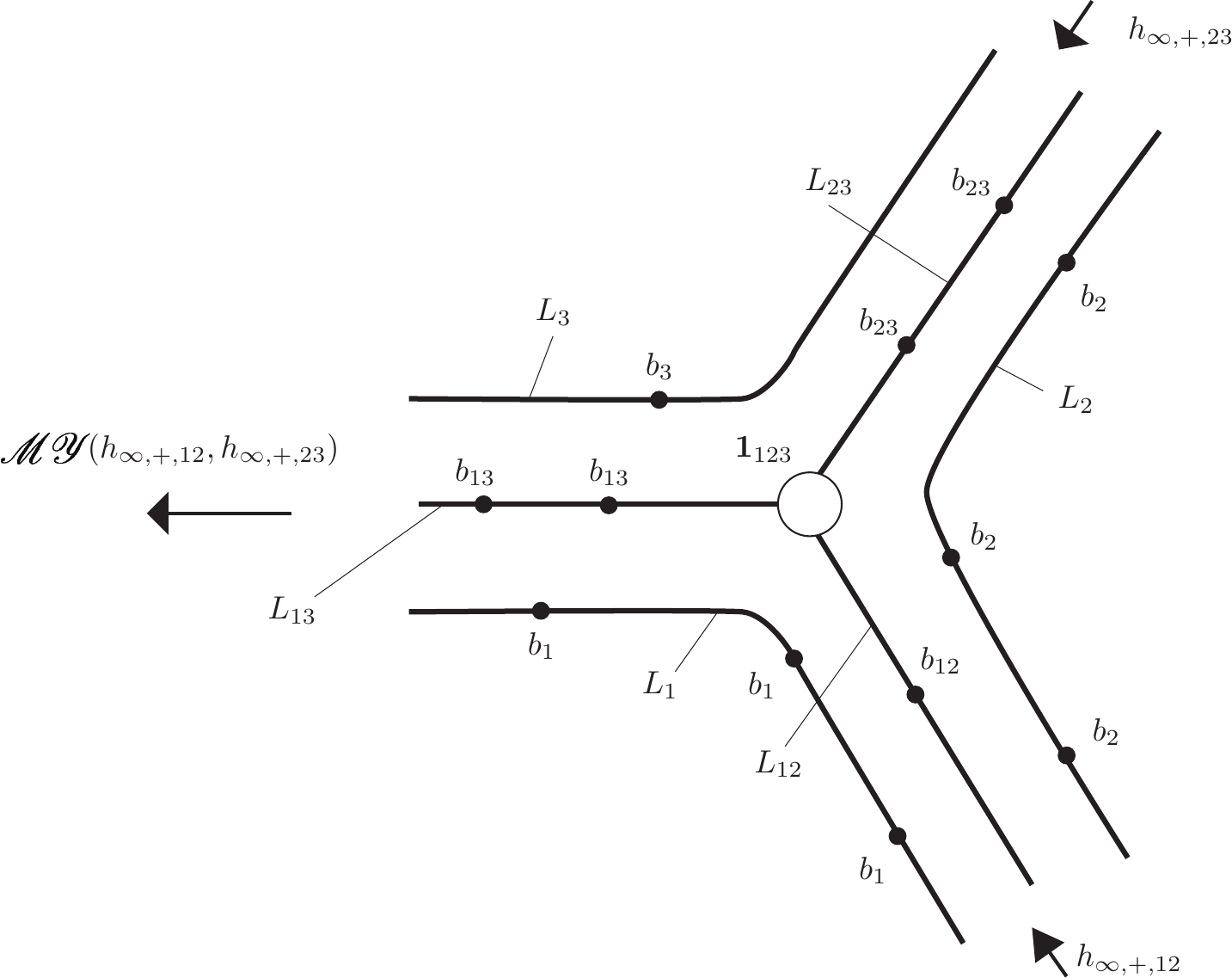}
\caption{The map $\mathscr{MY}$.}
\label{Figure107}
\end{figure}

Note that $CF(L_i,L_{ii'};L_{i'})$ for $ii' = 12, 23$ or $13$ is a filtered $A_{\infty}$
tri-module over $CF(L_i)$, $CF(L_{ii'})$, $CF(L_{i'})$.
Therefore, bounding cochains of $CF(L_i)$, $CF(L_{ii'})$, $CF(L_{i'})$
deform their `boundary operators' to obtain a boundary operator.
Namely, if $b_i b_{ii'}$, $b_{i'}$ are bounding cochains of
$CF(L_i)$, $CF(L_{ii'})$, $CF(L_{i'})$, we put
$
d^{b_i;b_{ii'};b_{i'}}(x) = \mathfrak n\bigl(e^{b_i},e^{b_{ii'}};x;e^{b_{i'}}\bigr)
$,
where $\mathfrak n$ is the structure operations of the tri-module
in Theorem~\ref{trimain}.
\begin{lem}\label{lem913}
The map $\mathscr{MY}$ in \eqref{form531} is a chain map
with respect to the boundary operators~$d^{b_1;b_{12};b_{2}}$, \smash{$d^{b_2;b_{23};b^{(1)}_{3}}$}, \smash{$d^{b_1;b_{13};b^{(1)}_{3}}$}.

\end{lem}
\begin{proof}
We put $h_{\infty,123} = {\bf 1}_{123}$, ${\bf h}_{12} = e^{b_{12}}$, ${\bf h}_{23} = e^{b_{23}}$, ${\bf h}_{13} = e^{b_{13}}$,
${\bf h}_{1} = e^{b_{1}}$, ${\bf h}_{2} = e^{b_{2}}$, \smash{${\bf h}_{3} = e^{b^{(1)}_{3}}$}
and apply Proposition~\ref{prop910}.
The first 6 terms of \eqref{form926} vanish because ${\bf h}_{ii'}$, ${\bf h}_{i}$ are exponentials
of the bounding cochains. The 7th term vanishes because ${\bf 1}_{123}$ is a cycle
with respect to the differential of ${CF}(L_{13};L_{12},L_{23})$ twisted by
$b_{12}$, $b_{23}$, $b_{13}$.
In fact, this is the definition of $b_{13}$. (See \eqref{newform87}.)
The 8th, 9th, 10th terms give the elements
$
\mathscr{YT}(d^{b_1;b_{12};b_{2}}(h_{\infty,12}),h_{\infty,23})
$,
\smash{$
\mathscr{YT}(h_{\infty,12},d^{b_1;b_{23};b^{(1)}_{3}}(h_{\infty,23}))
$},
and
\smash{$
d^{b_1;b_{13};b^{(1)}_{3}}(\mathscr{YT}(h_{\infty,12},h_{\infty,23}))
$}
respectively.
The lemma follows.
\end{proof}

For $ii' = 12, 23$ or $13$, we denote by
${\bf 1}_{ii'}$ the function $1$ on the diagonal component of
$\tilde L_i \times_{X_i \times X_{i'}} L_{i'}$.
\begin{lem}\label{lem914}
$
\mathscr{YT}({\bf 1}_{12},{\bf 1}_{23}) \equiv {\bf 1}_{13}
\mod \Lambda_+$.

\end{lem}
\begin{proof}
The operation $\mathscr{YT}$ are defined modulo $\Lambda_+$ by the integration along the
fiber of the~moduli space
${\mathcal M}_{\rm Y}(\vec a_{12},\vec a_{23},\vec a_{13};\vec a_{1},\vec a_{2},\vec a_{3},a_{\infty,123},\vec a_{\infty};0)$,
which consists of constant maps.
Using this fact and the definitions, we can prove the lemma easily in the same way as
Proposi\-tion~\ref{prop610}.
\end{proof}

We recall that on $CF(L_1,L_{13};L_3)$ we have a structure of
right $CF(L_3)$ module $\mathfrak n_k$. In fact, we put
$
\mathfrak n_k(y;x_1,\dots,x_k)
=
\mathfrak n\bigl(e^{b_1};e^{b_{13}};y;x_1,\dots,x_k\bigr)$
(see Lemma~\ref{lem6868}).

By the definition of \smash{$b_{3}^{(2)}$}, we have
\begin{equation}\label{form934}
\sum_{k=0}^{\infty}\mathfrak n_k\bigl({\bf 1}_{13};b_{3}^{(2)},\dots,b_{3}^{(2)}\bigr) = 0.
\end{equation}
We put
$
{\bf 1}'_{13} = \mathscr{YT}({\bf 1}_{12},{\bf 1}_{23})$.
Then by Lemma~\ref{lem913}, we have
\begin{equation}\label{form935}
\sum_{k=0}^{\infty}\mathfrak n_k\bigl({\bf 1}'_{13};b_{3}^{(1)},\dots,b_{3}^{(1)}\bigr) = 0.
\end{equation}
By \eqref{form934}, \eqref{form935} and Lemma~\ref{lem914}, we can apply Lemma~\ref{lem912}
to conclude that \smash{$b_{3}^{(1)}$} is gauge equivalent to \smash{$b_{3}^{(2)}$}.
The proof of Proposition~\ref{prop912}\,(2) is complete.
\end{proof}

\section{The compatibility as 2-functors}
\label{2-category formulation}

\subsection[The composition of $A_\infty$ functors defines a bi-functor]{The composition of $\boldsymbol{A_{\infty}}$ functors defines a bi-functor}
\label{sec:compfunc}

To obtain a more functorial version of Theorem~\ref{thm93},
we need the following algebraic result.

\begin{thm}\label{thm94}
Let $\mathscr C_i$ be a unital, strict and gapped filtered $A_{\infty}$ category for
$i=1,2,3$.
Then, there exists a filtered $A_{\infty}$ bi-functor\index[syindex]{comp@$\mathfrak{Comp}$}
\begin{equation}\label{form92}
\mathfrak{Comp} \colon\ \mathcal{FUNC}(\mathscr C_1,\mathscr C_2)
\times \mathcal{FUNC}(\mathscr C_2,\mathscr C_3)
\to \mathcal{FUNC}(\mathscr C_1,\mathscr C_3)
\end{equation}
such that
$
\mathfrak{Comp}_{\rm ob}(\mathscr F_{12},\mathscr F_{23})
= \mathscr F_{23}\circ \mathscr F_{12}$.

\end{thm}

We fix a discrete monoid $G \subset \R_{\ge 0}$. Here and hereafter the objects
$\mathcal{FUNC}(\mathscr C_1,\mathscr C_2)$
are strict, unital and $G$-gapped filtered $A_{\infty}$ functors.
\begin{rem}
Theorem~\ref{thm94} could be a part of the construction of an
($A_{\infty}$) 2-category whose object is a
filtered $A_{\infty}$ category. See Section~\ref{sec:catofAinfcat}.

\end{rem}

The unfiltered version of this statement is in \cite{Ly}.
We prove it here since we need the construction of
the functor $\mathfrak{Comp}$ for our application to geometry
in Sections~\ref{sec:compfunc21} and \ref{sec:compfuncalglem}.
Our proof below is different from the proof in \cite{Ly}.
\begin{proof}
Let $\mathscr C_1$, $\mathscr C_2$ be unital, strict and gapped
filtered $A_{\infty}$ categories.
\begin{lemdef}\label{lemdef97}
There exists a filtered $A_{\infty}$ functor
\[
\mathfrak{RYon} \colon\ \mathcal{FUNC}(\mathscr C_1,\mathscr C_2) \to \mathcal{BIMOD}(\mathscr C_1,\mathscr C_2)^{\rm op},
\]
which is a homotopy equivalence to its image.
We call this functor
the {\rm relative Yoneda functor}.\index{relative Yoneda functor}\index[syindex]{RYon@$\mathfrak{RYon}$}
\end{lemdef}

\begin{proof}
The functor $\mathfrak{OpYon}$ (for $\mathscr C_2$)
and the isomorphism in Lemma~\ref{opopoplemma} induces
\[
\mathcal{FUNC}(\mathscr C_1,\mathscr C_2)
\cong
\mathcal{FUNC}(\mathscr C_1^{\rm op},\mathscr C_2^{\rm op})^{\rm op}
\to \mathcal{FUNC}(\mathscr C^{\rm op}_1,\mathcal{FUNC}(\mathscr C_2,\mathcal{CH}))^{\rm op}.
\]
On the other hand, in Definition~\ref{lem56}, we defined an isomorphism
\[
\mathcal{FUNC}(\mathscr C^{\rm op}_1,\mathcal{FUNC}(\mathscr C_2,\mathcal{CH}))
\cong
\mathcal{BIMOD}(\mathscr C_1;\mathscr C_2).
\]
The lemma follows.
\end{proof}

Let us describe the functor $\mathfrak{RYon}$ more explicitly below.
Let $\mathfrak F \colon \mathscr C_1 \to \mathscr C_2$ be an object of~$\mathcal{FUNC}(\mathscr C_1,\mathscr C_2)$.
We first define left-$\mathscr C_1$ and right-$\mathscr C_2$ bi-module
$\mathfrak{RYon}_{\rm ob}(\mathfrak F)$.
Let $c_i$ be an object of~$\mathscr C_i$ for $i=1,2$.
We put
$
D_{c_1,c_2} = \mathscr C_2(\mathfrak F_{\rm ob}(c_1),c_2)$.
Let ${\bf x} \in B\mathscr C_1(c'_1,c_1)$,
$z \in D_{c_1,c_2}$,
${\bf y} \in B\mathscr C_2(c_2,c'_2)$.
We define
$\mathfrak n({\bf x},z,{\bf y})
\in D_{c'_1,c'_2} = \mathscr C_2(\mathfrak F_{\rm ob}(c'_1),c'_2)$
by
\begin{equation}\label{sec10form2}
\mathfrak n({\bf x},z,{\bf y}) =
(-1)^{\deg' {\bf y}}\mathfrak m\bigl(\widehat{\mathfrak F}({\bf x}),z,{\bf y}\bigr).
\end{equation}
\begin{rem}
In the definition of relative Yoneda functor, we use
$\mathfrak{OpYon}$ which is defined by using opposite category
$\mathscr C^{\rm op}$. Moreover, we use the isomorphism
$\mathcal{FUNC}(\mathscr C_1,\mathscr C_2)
\cong
\mathcal{FUNC}(\mathscr C_1^{\rm op},\allowbreak\mathscr C_2^{\rm op})^{\rm op}$.
Since we take the operation taking opposite twice
${\bf x}$ in the left-hand side becomes ${\bf x}$ in the right hand
side.
The $+1$ in Definition~\ref{opcate}\,(3) cancels with the minus
sign in Definition~\ref{defn2929}.

A rather complicated process to define $\mathfrak{OpYon}$
becomes a simple and natural formula \eqref{sec10form2},
when we write it explicitly.

The languages of functors and of bi-modules are mostly equivalent
when the target is $\mathcal{CH}$. However, the identification
includes the process taking opposite.
\end{rem}

We will check \eqref{form9300}.
Let
$\Delta {\bf x} = \sum_{a_1}\text{\bf x}_{a_1;1}\otimes \text{\bf x}_{a_1;2}$,
$\Delta {\bf y} = \sum_{a_2}\text{\bf y}_{a_2;1}\otimes \text{\bf y}_{a_2;2}$.
By definition, we have
\begin{gather}
\sum_{a_1}\sum_{a_2}(-1)^{*_2}\mathfrak n(\text{\bf x}_{a_1;1}, \mathfrak n(\text{\bf x}_{a_1;2},
z,\text{\bf y}_{a_2;1}),\text{\bf y}_{a_2;2})\nonumber
\\
\qquad=
\sum_{a_1}\sum_{a_2}(-1)^{*_3} \mathfrak m\bigl(\widehat{\mathfrak F}({\bf x}_{a_1;1}),
\mathfrak m(\widehat{\mathfrak F}({\bf x}_{a_2;1}),z,\text{\bf y}_{a_2;1}),\text{\bf y}_{a_2;2}\bigr),
\label{form95}
\end{gather}
where $*_2 = \deg'\text{\bf x}_{a_1;1}$
and $*_3 = \deg'{\bf x}_{a_2;1}+\deg' {\bf y}$.

Moreover,
\begin{equation}\label{form96}
\mathfrak n\bigl(\hat d(\text{\bf x}),z,\text{\bf y}\bigr)
=
(-1)^{\deg' {\bf y}}\mathfrak m\bigl(\widehat{\mathfrak F}(\hat d {\bf x}),z,{\bf y}\bigr)
=
(-1)^{\deg' {\bf y}}\mathfrak m\bigl(\hat d(\widehat{\mathfrak F}({\bf x})),z,{\bf y}\bigr).
\end{equation}
Here the second equality follows from the fact that $\mathfrak F$ is a filtered $A_{\infty}$ functor.
Furthermore,
\begin{equation}\label{form97}
(-1)^{\deg'{\bf x}+\deg z}\mathfrak n\bigl(\text{\bf x},
z,\hat d(\text{\bf y})\bigr)
=
\sum_{a_1}(-1)^{*_4}\mathfrak m\bigl(\widehat{\mathfrak F}({\bf x}),z,\hat d({\bf y})\bigr),
\end{equation}
where $*_4= \deg' {\bf x} + \deg z + \deg' {\bf y}$.

Formulas \eqref{form95}, \eqref{form96}, \eqref{form97} and
the $A_{\infty}$ formula for
$\mathfrak m$ imply \eqref{form9300} with sign modified (see Remark~\ref{sec10signerem}), using the fact that $\widehat{\mathfrak F}$ is a cohomomorphism.
Thus $D_{c_1,c_2}$ equipped with this bi-module structure is
$\mathfrak{RYon}_{\rm ob}(\mathcal F)(c_1,c_2)$.
\begin{rem}\label{sec10signerem}
In this and the next sections, we use the sign convention of filtered $A_{\infty}$
modules (multi-modules) so that the degree of elements of modules are not
shifted.
In other words, in~\eqref{form97}, $\deg z$ appears in place of $\deg' z$. The sign $(-1)^{\deg' {\bf y}}$ in \eqref{sec10form2}
appears by this reason.
See Remark~\ref{rem55111}.
\end{rem}
A natural transformation $\mathcal T$
from $\mathfrak F$ to $\mathfrak G$ gives
$
\mathcal T_0(c_1) \in \mathscr C_2(\mathfrak F_{\rm ob}(c_1),\mathfrak G_{\rm ob}(c_1))$.
It induces a~cochain map
$
\mathscr C_2(\mathfrak G_{\rm ob}(c_1),c_2)
\to \mathscr C_2(\mathfrak F_{\rm ob}(c_1),c_2).
$
This is a part of a bi-module homomorphism
from $\mathfrak{RYon}_{\rm ob}(\mathfrak G)$ to $\mathfrak{RYon}_{\rm ob}(\mathfrak F)$.
The direction of the arrows are opposite. This is the reason why
the opposite category appears in Lemma--Definition~\ref{lemdef97}.

The next lemma-definition, Propositions \ref{lem9100} and \ref{prop1019} are
closely related to the work~\cite{To} by To\"en.
\begin{lemdef}\label{defntensor}
Let $\mathscr C_1$, $\mathscr C_2$, $\mathscr C_3$ be filtered $A_{\infty}$
categories.
There exists a filtered $A_{\infty}$ bi-functor
\[
\mathfrak{ten} \colon\ \mathcal{BIMOD}(\mathscr C_1,\mathscr C_2)
\times \mathcal{BIMOD}(\mathscr C_2,\mathscr C_3)
\to
\mathcal{BIMOD}(\mathscr C_1,\mathscr C_3).
\]
We call it the {\em derived tensor product functor}.
\index{derived tensor product}
\index[syindex]{ten@$\mathfrak{ten}$}
\end{lemdef}
\begin{proof}
Let
\[
\mathfrak D^{12} = \bigl(\bigl\{D^{12}_{c_1,c_2}\bigr\},\bigl\{\mathfrak n^{12}_{c_1',c_1,c_2,c_2'}\bigr\}\bigr)
\]
be an object of
$\mathcal{BIMOD}(\mathscr C_1,\mathscr C_2)$
and let
\[
\mathfrak D^{23} = \bigl(\bigl\{D^{23}_{c_2,c_3}\bigr\},\bigl\{\mathfrak n^{23}_{c_2',c_2,c_3,c_3'}\bigr\}\bigr)
\]
be an object of
$\mathcal{BIMOD}(\mathscr C_1,\mathscr C_2)$.
We will define an object
\[
\mathfrak D^{13} =  \bigl(\bigl\{D^{13}_{c_1,c_3}\bigr\},\bigl\{\mathfrak n^{13}_{c_1',c_1,c_3,c_3'}\bigr\}\bigr)
\]
of $\mathcal{BIMOD}(\mathscr C_1,\mathscr C_3)$.

Let $c_1, c'_1 \in \mathfrak{OB}(\mathscr C_1)$,
$c_3, c'_3 \in \mathfrak{OB}(\mathscr C_3)$.
We put
\begin{equation}\label{form107}
D^{13}_{c_1,c_3}
=
\underset{c_2,c'_2}{\widehat{\bigoplus}}
D^{12}_{c_1,c_2} \,\widehat{\otimes}\, B\mathscr C_2[1](c_2,c'_2)
\,\widehat{\otimes}\, D^{23}_{c'_2,c_3}.
\end{equation}
We remark that $B\mathscr C_2[1](c_2,c'_2)$
contains $1 \in B_0\mathscr C_2[1](c_2,c_2) \cong \Lambda_0$
when $c_2 = c'_2$.

Let ${\bf x} \in B\mathscr C_1[1](c'_1,c_1)$,
${\bf y} \in B\mathscr C_3[1](c_3,c'_3)$
and
\[
z = u \otimes {\bf v} \otimes w\in
D^{12}_{c_1,c_2} \,\widehat{\otimes}\, B\mathscr C_2[1](c_2,c'_2)
\,\widehat{\otimes}\, D^{23}_{c'_2,c_3}
\subseteq D^{13}_{c_1,c_3}.
\]
We define
\smash{$\mathfrak n^{13}_{c_1',c_1,c_3,c_3'} \colon
B\mathscr C_1[1](c'_1,c_1)
\,\widehat\otimes\, D^{13}_{c_1,c_3}
\otimes B\mathscr C_3[1](c_3,c'_3) \to D^{13}_{c'_1,c'_3}$}
by
\begin{gather}
\mathfrak n^{13}_{c_1',c_1,c_3,c_3'}({\bf x},z,{\bf y})
:= \begin{cases}
\sum_{a}
\mathfrak n^{12}({\bf x},u,{\bf v}_{a;1}) \otimes {\bf v}_{a;2} \otimes w
& \text{if ${\bf y} =1 \in B_0\mathscr C_3(c_3,c'_3)$}, \\
\sum_{a} (-1)^*
u \otimes {\bf v}_{1:a} \otimes \mathfrak n^{23}({\bf v}_{a;2},w,{\bf y})
& \text{if ${\bf x} =1 \in B_0\mathscr C_1(c_1,c'_1)$}, \\
\sum_{a} \mathfrak n^{12}(u,{\bf v}_{a;1}) \otimes {\bf v}_{a;2} &\\
\quad+ \sum_{a} (-1)^* u \otimes {\bf v}_{a;1} \otimes \mathfrak n^{23}({\bf v}_{a;2},w) &\\
\quad+ (-1)^{\deg u} u \otimes \hat d({\bf v}) \otimes w
& \text{if ${\bf x} = {\bf y} =1$}, \\
0 & \text{otherwise},
\end{cases}\hspace{-10mm}\label{form108}
\end{gather}
where $* = \deg u + \deg' {\bf v}_{a;1}$
and $\Delta {\bf v} = \sum_{a}\text{\bf v}_{a;1}\otimes \text{\bf v}_{a;2}$.
It is straightforward to check \eqref{form9300} with sign modified
(see Remark~\ref{reem2109}).
We thus defined $\mathfrak{ten}_{\rm ob}$.
\begin{rem}\label{reem2109}
We remark that in the second, fourth and fifth lines of the right-hand side
we used $\deg u$ and not $\deg' u$.
\end{rem}
\begin{rem}
Note that in the case of $\mathfrak D^{13}
= \mathfrak{ten}\bigl(\mathfrak D^{12},\mathfrak D^{23}\bigr)$, the `left multiplication' and
the `right multiplication' exactly commute.
This is the reason why we take $0$ in the fourth case of~\eqref{form108}.
In fact, $\mathfrak n_{1,1}$ in the bi-module structure is a chain homotopy between
$\mathfrak n_{0,1}(\mathfrak n_{1,0}(x,z),y)$ and~${(-1)^{\deg' x}\mathfrak n_{1,0}(x,\mathfrak n_{0,1}(z,y))}$.
\end{rem}
We next define the morphism part of the bi-functor $\mathfrak{ten}$.
Let
\[
\mathfrak D^{(j),12} = \bigl(\bigl\{D^{(j),12}_{c_1,c_2}\bigr\},\bigl\{\mathfrak n^{12}_{c_1',c_1,c_2,c_2'}\bigr\}\bigr)
\]
be an object of
$\mathcal{BIMOD}(\mathscr C_1,\mathscr C_2)$ for $j=1,2$
and
\[
\mathfrak D^{(j),23} =\bigl(\bigl\{D^{(j),23}_{c_2,c_3}\bigr\},\bigl\{\mathfrak n^{23}_{c_2',c_2,c_3,c_3'}\bigr\}\bigr)
\]
an object of
$\mathcal{BIMOD}(\mathscr C_1,\mathscr C_2)$ for $j=1,2$.
A pre-bi-module homomorphism
\smash{$
\mathfrak f^{12} \colon \mathfrak D^{(1),12}\! \to \mathfrak D^{(2),12}$}
(resp.\ \smash{$
\mathfrak f^{23} \colon \mathfrak D^{(1),23} \to \mathfrak D^{(2),23}$} ) consists of
\begin{gather*}
\begin{split}
&\mathfrak f^{12}_{k_1,k_2} \colon\ B_{k_1}\mathscr C_1[1](c_1,c'_1) \otimes D^{(1),12}_{c'_1,c'_2} \otimes
B_{k_2}\mathscr C_2[1](c'_2,c_2)
\to D^{(2),12}_{c_1,c_2},\\
& \text{(resp.} \
\mathfrak f^{23}_{k_2,k_3} \colon\ B_{k_2}\mathscr C_2[1](c_2,c'_2) \otimes D^{(1),23}_{c'_2,c'_3} \otimes
B_{k_3}\mathscr C_3[1](c'_3,c_3)
\to D^{(2),23}_{c_2,c_3}).
\end{split}
\end{gather*}
See Definition~\ref{defnbomodhomocat}.
We define its tensor product
$\mathfrak f^{12} \otimes \mathfrak f^{23} = \mathfrak f^{13}$ as follows.
We define \smash{$D^{(j),13}_{c_1,c_3}$} in the same way as~\eqref{form107}.
$\mathfrak f^{13}$ consists of the maps
\[
\mathfrak f^{13}_{k_1,k_3}\colon\ B_{k_1}\mathscr C_1[1](c_1,c'_1) \otimes D^{(1),13}_{c'_1,c'_3} \otimes
B_{k_3}\mathscr C_3[1](c'_3,c_3)
\to D^{(2),13}_{c_2,c_3},
\]
which we define by the next formula.
Let ${\bf x} \in B\mathscr C_1[1](c'_1,c_1)$,
${\bf y} \in B\mathscr C_3[1](c_3,c'_3)$
and
\[
z = u \otimes {\bf v} \otimes w\in
D^{(1),12}_{c_1,c_2} \,\widehat{\otimes}\, B\mathscr C_2[1](c_2,c'_2)
\,\widehat{\otimes}\, D^{23}_{c'_2,c_3}
\subset D^{(1),13}_{c_1,c_3}.
\]
We put
\begin{align*}
&\mathfrak f^{13}_{k_1,k_3}({\bf x},z,{\bf y})\\
&\qquad =
\sum_{a}
(-1)^{\deg \mathfrak f^{23}_{*,k_3} (\deg'{\bf x} +\deg u +\deg' {\bf v}_{a;1}
+ \deg' {\bf v}_{a;2})}
\mathfrak f^{12}_{k_1,*}({\bf x},u,{\bf v}_{a;1}) \otimes {\bf v}_{a;2}
\otimes
\mathfrak f^{23}_{*,k_3}({\bf v}_{a;3},w,{\bf y}).
\end{align*}
Here $(1 \otimes \Delta) \circ \Delta {\bf v} = \sum_{a}\text{\bf v}_{a;1}\otimes \text{\bf v}_{a;2} \otimes \text{\bf v}_{a;3}$.
We can easily show that $\mathfrak f^{13}$ gives a chain map
\begin{gather*}
\mathcal{BIMOD}(\mathscr C_1,\mathscr C_2)\bigl(\mathfrak D^{(1),12},\mathfrak D^{(2),12}\bigr)
\otimes
\mathcal{BIMOD}(\mathscr C_2,\mathscr C_3)\bigl(\mathfrak D^{(1),23},\mathfrak D^{(2),23}\bigr) \\
\qquad\to
\mathcal{BIMOD}(\mathscr C_1,\mathscr C_3)\bigl(\mathfrak D^{(1),13},\mathfrak D^{(2),13}\bigr)
\end{gather*}
Moreover, this map $\bigl(\mathfrak f^{12},\mathfrak f^{23}\bigr) \mapsto \mathfrak f^{12} \otimes \mathfrak f^{23} = \mathfrak f^{13}$
is compatible with composition. Namely,
\begin{gather*}
\bigl(\mathfrak f^{(1),12} \circ \mathfrak f^{(2),12}\bigr) \otimes_s
\bigl(\mathfrak f^{(1),23} \circ \mathfrak f^{(2),23}\bigr) \\
\qquad=
(-1)^{\deg \mathfrak f^{(1),23}\deg \mathfrak f^{(2),12}}\bigl(\mathfrak f^{(1),12} \otimes_s \mathfrak f^{(1),23}\bigr) \circ \bigl(\mathfrak f^{(2),12} \otimes_s \mathfrak f^{(2),23}\bigr).
\end{gather*}
See \eqref{form2626} for $\otimes_s$.

Therefore, by putting other operations to be zero we obtain a required bi-functor.
The proof of Lemma--Definition~\ref{defntensor} is complete.
\end{proof}

The derived tensor product functor induces
\[
\mathcal{BIMOD}(\mathscr C_1,\mathscr C_2)^{\rm op}
\times \mathcal{BIMOD}(\mathscr C_2,\mathscr C_3)^{\rm op}
\to
\mathcal{BIMOD}(\mathscr C_1,\mathscr C_3)^{\rm op},
\]
which we denote also by $\mathfrak{ten}$.
\begin{rem}
The proof shows that $\mathfrak{ten}$ is actually a bi-DG-functor between DG-categories.
\end{rem}
The proof of the next proposition is the most nontrivial
part of the proof of Theorem~\ref{thm94}.

\begin{prop}\label{lem9100}
Assume that $\mathscr C_1$, $\mathscr C_2$, $\mathscr C_3$
are unital, strict and gapped.
Let $\mathscr F_{12} \colon \mathscr C_1 \to \mathscr C_2$
and $\mathscr F_{23} \colon \mathscr C_2 \to \mathscr C_3$
be filtered $A_{\infty}$ functors.
Then the object
$\mathfrak{ten}_{\rm ob}(\mathfrak{RYon}_{\rm ob}(\mathfrak F_{12}),\mathfrak{RYon}_{\rm ob}(\mathfrak F_{23}))$
of
$\mathcal{BIMOD}(\mathscr C_1,\mathscr C_3)^{\rm op}$ is homotopy equivalent to
$\mathfrak{RYon}_{\rm ob}(\mathfrak F_{23} \circ \mathfrak F_{12})$.
\end{prop}
The proof is given in Section~\ref{sec:compfuncalglem}.
\begin{rem}\label{rem1011}
Suppose that $\mathscr C_1$, $\mathscr C_2$ and $\mathscr C_3$
are associative rings with unity. They can be regarded as unital $A_{\infty}$ categories.
Let $\mathscr F_{12} \colon \mathscr C_1 \to \mathscr C_2$
and $\mathscr F_{23} \colon \mathscr C_2 \to \mathscr C_3$ be unital ring homomorphisms
which are special cases of $A_{\infty}$ functors.

The bi-module associated to $\mathscr F_{12}$ is $\mathscr C_2$
which is regarded as a right $\mathscr C_2$ module by right multiplication
and a left $\mathscr C_1$ module by $x\cdot y = \mathscr F_{12}(x)y$.
We write this bi-module as ${}_{\mathscr C_1} (\mathscr C_2)_{\mathscr C_2}$.
In the same way $\mathscr F_{23}$ corresponds to
${}_{\mathscr C_2} (\mathscr C_3)_{\mathscr C_3}$.
Their tensor product is
$
{}_{\mathscr C_1} (\mathscr C_2)_{\mathscr C_2} \otimes_{\mathscr C_2} {}_{\mathscr C_2} (\mathscr C_3)_{\mathscr C_3}
= {}_{\mathscr C_1} (\mathscr C_3)_{\mathscr C_3}$.
Here the left $\mathscr C_1$ module structure is induced by $\mathscr F_{23}\circ\mathscr F_{12}$.
This is Proposition~\ref{lem9100} in this case.
\end{rem}
Now we are in the position to complete the proof of
Theorem~\ref{thm94}.
We consider the bi-functor
\[
\mathfrak{ten}\circ (\mathfrak{RYon}\times \mathfrak{RYon}) \colon\
\mathcal{FUNC}(\mathscr C_1,\mathscr C_2)
\times \mathcal{FUNC}(\mathscr C_2,\mathscr C_3)
\to
\mathcal{BIMOD}(\mathscr C_1,\mathscr C_3)^{\rm op}.
\]
We consider the full subcategory
$\mathfrak{Rep}(\mathscr C_1,\mathscr C_3)$
of
$\mathcal{BIMOD}(\mathscr C_1,\mathscr C_3)^{\rm op}$
whose object is homotopy equivalent to an image
of
$\mathfrak{RYon} \colon
\mathcal{FUNC}(\mathscr C_1,\mathscr C_3)
\to \mathcal{BIMOD}(\mathscr C_1,\mathscr C_3)^{\rm op}$.
Proposition~\ref{lem9100} implies that the image of
$\mathfrak{ten}\circ (\mathfrak{RYon}\times \mathfrak{RYon})$
is contained in this full subcategory.

Moreover, by Lemma--Definition~\ref{lemdef97}
and Theorem~\ref{white},
there exists a filtered $A_{\infty}$ functor~${\mathfrak{Rep}(\mathscr C_1,\mathscr C_3)
\to \mathcal{FUNC}(\mathscr C_1,\mathscr C_3)}$
which is a homotopy inverse to $\mathfrak{RYon}$.
Therefore, there exists a filtered $A_{\infty}$
functor
$\mathfrak{Comp} \colon \mathcal{FUNC}(\mathscr C_1,\mathscr C_2)
\times \mathcal{FUNC}(\mathscr C_2,\mathscr C_3)
\to \mathcal{FUNC}(\mathscr C_1,\mathscr C_3)$
such that the next diagram commutes up to homotopy equivalence:
\begin{equation}\label{dia1010}
\begin{CD}
\mathcal{FUNC}(\mathscr C_1,\mathscr C_2)
\times \mathcal{FUNC}(\mathscr C_2,\mathscr C_3) @ >>>
\mathcal{FUNC}(\mathscr C_1,\mathscr C_3) \\
@ VVV @ VVV\\
\mathcal{BIMOD}(\mathscr C_1,\mathscr C_2)^{\rm op}
\times\mathcal{BIMOD}(\mathscr C_2,\mathscr C_3)^{\rm op} @ > >> \mathcal{BIMOD}(\mathscr C_1,\mathscr C_3)^{\rm op}.
\end{CD}
\end{equation}
This is the required functor.
\end{proof}

\begin{rem}
The construction of the composition functor we gave in this subsection is rather indirect.
In other words, we did not provide an explicit formula how the pre-natural transformations
are sent by this functor. This is because an explicit homotopy inverse to the Yoneda
functor is not given. We can certainly find some formula by following the proof.
In fact, the Yoneda functor is explicitly
given in \cite{fu4} and the proof of Theorem~\ref{white} in \cite{fu4}
is by induction each of whose steps is in principle can be made explicit.
However, the explicit formula which we may obtain in that way seems to be very complicated
and is not practical to use it.
\end{rem}

\subsection[Proof of Theorem~\ref{thm93}]{Proof of Theorem~\ref{thm93}}
\label{sec:compfunc21}

In this section, we prove Theorem~\ref{thm93}.
Before starting the proof, we twist the (category version of the) map $\mathscr{YT}$ in
 Proposition~\ref{prop910} by bounding cochains.
We denote by $\mathcal L_i$, $\mathcal L_{ii'}$
or \smash{$\mathcal L^{(j)}_i$}, \smash{$\mathcal L^{(j)}_{ii'}$} objects of
$\mathfrak{Fukst}(X_i)$, $\mathfrak{Fukst}(-X_i \times X_{i'})$.
We recall
\[
B_kCF[1](\mathcal L_i,\mathcal L'_i)
=
\bigoplus_{\mathcal L_i = \mathcal L^{(0)}_i,
\dots, \mathcal L^{(k)}_i =\mathcal L'_i}
\bigotimes_{j=1}^k CF\bigl(\mathcal L^{(j-1)}_i,\mathcal L^{(j)}_i\bigr)[1]
\]
and $BCF[1](\mathcal L_i,\mathcal L'_i)$
is their completed direct sum over $k$.
We define the modules
$B_kCF[1](\mathcal L_{ii'},\allowbreak\mathcal L'_{ii'})$,
$BCF[1](\mathcal L_{ii'},\mathcal L'_{ii'})$ in the same way.

We define a map
\smash{$
\mathfrak t_{\vec b} \colon
\bigotimes_{j=1}^k CF\bigl(\mathcal L^{(j-1)}_i,\mathcal L^{(j)}_i\bigr)[1]
\to
BCF[1](\mathcal L_i,\mathcal L'_i)
$}
by
\[
\mathfrak t_{\vec b}(x_1,\dots,x_k)
:= e^{b_0} x_1 e^{b_1} x_2 \cdots x_{k-1} e^{b_{k-1}} x_k e^{b_k}
\]
(see \eqref{defntttt}).
It induces
$
\mathfrak t_{\vec b} \colon BCF[1](\mathcal L_i,\mathcal L'_i)
\to BCF[1](\mathcal L_i,\mathcal L'_i)$.
We define
$
\mathfrak t_{\vec b} \colon BCF[1](\mathcal L_{ii'},\mathcal L'_{ii'})
\to BCF[1](\mathcal L_{ii'},\mathcal L'_{ii'})$ in the same way.

We now define the map
\begin{gather}
\mathscr{YT}^{\vec b} \colon \ BCF[1](\mathcal L_1,\mathcal L'_1)
\otimes
BCF[1](\mathcal L_{12},\mathcal L'_{12}) \otimes BCF[1](\mathcal L_{23},\mathcal L'_{23})\nonumber \\
\hphantom{\mathscr{YT}^{\vec b} \colon} \,
\otimes CF(\mathcal L'_1,\mathcal L'_{12};\mathcal L'_{2})
\otimes BCF[1](\mathcal L'_2,\mathcal L_2)
 \otimes
CF(\mathcal L_2,\mathcal L'_{23};\mathcal L'_{3})\otimes
BCF[1](\mathcal L'_{3},\mathcal L_{3})\nonumber
\\
\hphantom{\mathscr{YT}^{\vec b} \colon} \,  \otimes
BCF[1](\mathcal L'_{13},\mathcal L_{13})\otimes CF(\mathcal L_{13};\mathcal L_{12},\mathcal L_{23})
\to
CF(\mathcal L_1,\mathcal L_{13};\mathcal L_{3})\label{form10111}
\end{gather}
by composing $\mathfrak t_{\vec b}$ with $\mathscr{YT}$.
(We do not apply $\mathfrak t_{\vec b}$ to the factors
$CF(\mathcal L'_1,\mathcal L'_{12};\mathcal L'_{2})$, $CF(\mathcal L_2,\mathcal L'_{23};\mathcal L'_{3})$,
$CF(\mathcal L_1,\mathcal L_{13};\mathcal L_{3})$.)
\begin{lem}\label{lem1014}
\eqref{form9252} holds when we replace $\mathscr{YT}$,
$\hat d$, $\mathfrak{n}$ by
\smash{$\mathscr{YT}^{\vec b}$}, \smash{$\hat d^{b}$}, $\mathfrak{n}^b$,
respectively.
Here \smash{$\hat d^{b}$},~$\mathfrak{n}^b$ are defined by
\smash{$\mathfrak t_{\vec b} \circ \hat d^{b} = \hat d \circ \mathfrak t_{\vec b}$},
\smash{$\mathfrak t_{\vec b} \circ \mathfrak{n}^b = \mathfrak{n} \circ \mathfrak t_{\vec b}$}.
\end{lem}
This is immediate from Proposition~\ref{prop910}.
\begin{proof}[Proof of Theorem~\ref{thm93}]
Let $\mathcal L_{12}$, $\mathcal L_{23}$ be as in Theorem~\ref{thm93}
and $\mathcal L_{13} = \mathcal L_{23} \circ \mathcal L_{12}$.
We apply the relative Yoneda functor $\mathfrak{RYon}_{\rm ob}$ to
$\mathcal W_{\mathcal L_{13}}$. By definition we obtain
$
\mathscr{CF}(\mathbb L_1,\mathbb L_{13};\mathbb L_3).
$
We fixed~${\mathcal L_{13} \in \mathbb L_{13}}$
so
we regard $
\mathscr{CF}(\mathbb L_1,\mathbb L_{13};\mathbb L_3)
$ as a left-$\mathfrak{Fukst}(X_1)$
and right-$\mathfrak{Fukst}(X_3)$ bi-module.
It assigns
$
\mathcal W^{(1)}(\mathcal L_1,\mathcal L_3)
= CF(\mathcal L_1,\mathcal L_{13};\mathcal L_3)$
to $\mathcal L_i \in \mathfrak{Ob}(\mathfrak{Fukst}(X_i))$ for $i=1,3$.

We apply the relative Yoneda functor $\mathfrak{RYon}_{\rm ob}$ to
$\mathcal W_{\mathcal L_{12}}$ and $\mathcal W_{\mathcal L_{23}}$.
We then obtain tri-modules~${
\mathscr{CF}(\mathbb L_1,\mathbb L_{12};\mathbb L_2)}
$
and
$
\mathscr{CF}(\mathbb L_2,\mathbb L_{23};\mathbb L_3)
$,
respectively.
We fix $\mathcal L_{12}$ and $\mathcal L_{23}$
and regard them as left-$\mathfrak{Fukst}(X_1)$
and right-$\mathfrak{Fukst}(X_2)$
and
left-$\mathfrak{Fukst}(X_2)$
and right-$\mathfrak{Fukst}(X_3)$
modules respectively.
We consider $\mathcal W_{\mathcal L_{23}}\circ \mathcal W_{\mathcal L_{12}}$
and apply the relative Yoneda functor $\mathfrak{RYon}_{\rm ob}$ to it.
Then, by Proposition~\ref{lem9100},
we obtain
$
\mathfrak{ten}(\mathscr{CF}(\mathbb L_1,\mathbb L_{12};\mathbb L_2),
\mathscr{CF}(\mathbb L_2,\mathbb L_{23};\mathbb L_3))$.
We regard it as a left-$\mathfrak{Fukst}(X_1)$
and right-$\mathfrak{Fukst}(X_3)$ bi-module.
To $\mathcal L_i \in \mathfrak{Ob}(\mathfrak{Fukst}(X_i))$ for $i=1,3$,
it assigns
\begin{equation}\label{form101510}
\mathcal W^{(2)}(\mathcal L_1,\mathcal L_3) =
\bigoplus_{\mathcal L_2,\mathcal L'_2}
CF(\mathcal L_1,\mathcal L_{12};\mathcal L_2)
\otimes BCF[1](\mathcal L_2,\mathcal L'_2)
\otimes CF(\mathcal L'_2,\mathcal L_{23};\mathcal L_3).
\end{equation}
The pre-bi-module homomorphism we look for is a
system of maps
\[
BCF[1](\mathcal L_1,\mathcal L'_1)
\otimes
\mathcal W^{(2)}(\mathcal L'_1,\mathcal L'_3)
\otimes
BCF[1](\mathcal L'_3,\mathcal L_3)
\to
\mathcal W^{(1)}(\mathcal L_1,\mathcal L_3).
\]
Namely,
\begin{align*}
\mathcal T \colon\
BCF[1](\mathcal L_1,\mathcal L'_1)
&\otimes
CF(\mathcal L'_1,\mathcal L_{12};\mathcal L_2)
\otimes BCF[1](\mathcal L_2,\mathcal L'_2)
\\
&
\otimes CF(\mathcal L'_2,\mathcal L_{23};\mathcal L'_3)
\otimes
BCF[1](\mathcal L'_3,\mathcal L_3)
\to CF(\mathcal L_1,\mathcal L_{13};\mathcal L_3).
\end{align*}
We define
\begin{gather*}
\mathcal T({\bf h}_1,h_{\infty,+,12},{\bf h}_2,h_{\infty,+,23},{\bf h}_3) \\
\qquad=
(-1)^*
\mathscr{YT}^{\vec b}
({\bf 1}_{123};\varnothing_{12},\varnothing_{23},\varnothing_{13};h_{\infty,+,12},h_{\infty,+,23};
{\bf h}_{1},{\bf h}_{2},{\bf h}_{3}).
\end{gather*}
Here $\varnothing_{ii'} = 1 \in B_0CF(\mathcal L_{ii'})$ and
${\bf 1}_{123} \in CF(\mathcal L_{12},\mathcal L_{23};\mathcal L_{13})$ is the cyclic
element. The sign $(-1)^*$ is determined by the Koszul rule. We
count the way exchanging the order of the variables using the shifted degree $\deg'$
for elements of $BCF[1](\mathcal L'_i,\mathcal L_i)$
(or $BCF[1](\mathcal L_i,\mathcal L_i)$)
and $\deg$ for elements of $CF(\mathcal L'_1,\mathcal L_{12};\mathcal L_2)$ etc.
Then we put the sign according to whether the total count is even or odd.
\begin{rem}
The sign in \eqref{lem1014} is also by Koszul rule. However, it is different
from the one we describe above. Namely, $\deg'$ is used
also for elements of tri-modules, $CF(\mathcal L'_1,\mathcal L_{12};\mathcal L_2)$
etc.\ We change the sign of the maps in the same way as
\eqref{degzurashi} (see also \eqref{sec10form2})
to go from one to the other.
\end{rem}

The condition that $\mathcal T$ is a bi-module homomorphism
is a consequence of Lemma~\ref{lem1014} and the fact that ${\bf 1}_{123}$ becomes a cycle (after twisting
the boundary operator by the bounding cochains $b_{12}$, $b_{23}$, $b_{13}$).

We continue the proof of Theorem~\ref{thm93}
and prove that $\mathcal T$ is a homotopy equivalence.
In view of Lemma~\ref{prop1015}, the next step is to prove that
the chain map
\begin{equation}\label{1018888}
\mathcal T_{0,0;\mathcal L_1,\mathcal L_3} \colon\
\mathcal W^{(2)}(\mathcal L_1,\mathcal L_3)
\to
\mathcal W^{(1)}(\mathcal L_1,\mathcal L_3)
\end{equation}
is a chain homotopy equivalence for arbitrary
$\mathcal L_1$, $\mathcal L_3$.
By Proposition~\ref{lem9100},
the derived tensor product~\eqref{form101510}
is chain homotopy equivalent to
\begin{equation}\label{hattthhh}
CF(\mathcal W_{\mathcal L_{23}}(\mathcal W_{\mathcal L_{12}}(\mathcal L_1)),\mathcal L_3)
\cong CF(\mathcal W_{\mathcal L_{12}}(\mathcal L_1),\mathcal L_{23};\mathcal L_3).
\end{equation}
In fact, the chain homotopy equivalence from \eqref{hattthhh} to
\eqref{form101510} is given by
\begin{equation}\label{form101819}
x \mapsto {\bf 1}_{12}\otimes x,
\end{equation}
for $x \in CF(\mathcal W_{\mathcal L_{23}}(\mathcal W_{\mathcal L_{12}}(\mathcal L_{1})
),\mathcal L_3)$.
Here ${\bf 1}_{12} \in CF(\mathcal L_1,\mathcal L_{12};\mathcal W_{\mathcal L_{12}}
(\mathcal L_1))$
is the cyclic element which becomes the unity
in $CF(\mathcal W_{\mathcal L_{12}}
(\mathcal L_1))$ by the isomorphism
\[
 CF(\mathcal L_1,\mathcal L_{12}; \mathcal W_{\mathcal L_{12}}(\mathcal L_1))
\cong
CF(\mathcal \mathcal W_{\mathcal L_{12}}(\mathcal L_1)).
\]
Note that if we regard $x$ as an element of the right-hand side of \eqref{hattthhh},
then ${\bf 1}_{12}\otimes x$ is an element of
\[
CF(\mathcal L_1,\mathcal L_{12};\mathcal W_{\mathcal L_{12}}(\mathcal L_1))
\otimes CF(\mathcal W_{\mathcal L_{12}}(\mathcal L_1),\mathcal L_{23};\mathcal L_3),
\]
which is contained in \eqref{form101510} as the case
$\mathcal L_2 = \mathcal L'_2 = \mathcal W_{\mathcal L_{12}}(\mathcal L_1)$.

The map \eqref{form101819} is identified with the map
$\mathscr I_{12;0,0}$ in \eqref{form1024}, which we will use to
prove Proposition~\ref{lem9100}.

Thus to prove that \eqref{1018888} is a chain homotopy equivalence, it suffices to show that the composition
\begin{equation}\label{formula1020}
CF(\mathcal W_{\mathcal L_{12}}(\mathcal L_1),\mathcal L_{12};\mathcal L_3)
\to
\mathcal W^{(2)}(\mathcal L_1,\mathcal L_3)
\to
\mathcal W^{(1)}(\mathcal L_1,\mathcal L_3)
\end{equation}
is a chain homotopy equivalence. By definition,
\eqref{formula1020} is the map
\begin{gather}
h_{\infty,+,23} \mapsto \mathcal T(\varnothing_{1},{\bf 1}_{12},\varnothing_{2},h_{\infty,+,23},
\varnothing_{3})\nonumber\\
\hphantom{h_{\infty,+,23} \mapsto}{}
=
\mathscr{YT}^{\vec b}
({\bf 1}_{123};\varnothing_{12},\varnothing_{23},\varnothing_{13};{\bf 1}_{12},h_{\infty,+,23};
\varnothing_{1},\varnothing_{2},\varnothing_{3}).\label{form1021}
\end{gather}
Here $\varnothing_{i} = 1 \in B_0CF(\mathcal L_i)$, for $i=1,2,3$.
\begin{lem}
\quad
\begin{enumerate}\itemsep=0pt
\item[$(1)$]
$
CF(\mathcal W_{\mathcal L_{12}}(\mathcal L_1),\mathcal L_{23};\mathcal L_3)
\cong
\mathcal W^{(1)}(\mathcal L_1,\mathcal L_3)
\cong
\Omega({L_1\times_{X_1} L_{12}\times_{X_2} L_{23} \times_{X_3} L_3},\Theta_-)
\,\widehat\otimes\,\Lambda_0$.
\item[$(2)$]
The map \eqref{form1021} is congruent to the identity map modulo $\Lambda_+$
via the isomorphism of item~$(1)$.
\end{enumerate}

\end{lem}
\begin{proof}
(1) is immediate from the definition.
(2) then follows from the fact that energy $0$ part of the map
\smash{$\mathscr{FY}^{\vec b}$} is defined by the moduli space of constant maps.
\end{proof}

To complete the proof of Theorem~\ref{thm93},
we need to discuss the following point.
Note that while we proved Proposition
\ref{prop912} we showed that the two
bounding cochains, written as \smash{$b_3^{(1)}$} and \smash{$b_3^{(2)}$}
there, are gauge equivalent. However, they are not necessary
equal.
In the above argument, we used \smash{$b_3^{(2)}$}.
In fact, $CF(\mathcal L_2,\mathcal L_{23};\mathcal L_3)$
using \smash{$b_3^{(2)}$} gives
$\mathcal W_{\mathcal L_{23}}\colon\mathfrak{Fukst}(X_2) \to \mathfrak{Fukst}(X_3)$.

On the other hand, $CF(\mathcal L_1,\mathcal L_{13};\mathcal L_3)$
with \smash{$b_3^{(1)}$} gives
$\mathcal W_{\mathcal L_{13}}\colon\mathfrak{Fukst}(X_1) \to \mathfrak{Fukst}(X_3)$.
Therefore, to complete the proof of
Theorem~\ref{thm93}, we need
to compare
$CF(\mathcal L_1,\mathcal L_{13};\mathcal L_3)$
with two different choices of bounding cochains and show that they are
homotopy equivalent as left-$\mathfrak{Fukst}(X_1)$
and right-$\mathfrak{Fukst}(X_3)$ bi-modules.
We can prove it in the same way as the proof of
Proposition
\ref{prop912} as follows.

We consider
${\rm Poly}(CF(\mathcal L_1,\mathcal L_{13};\mathcal L_3))$
which is a left-$\mathfrak{Fukst}(X_1)$ and right-${\rm Poly}(\mathfrak{Fuk}(X_3))$ bi-module.
(See \cite[Section 5.2.3]{fooobook} and the proof of
Proposition~\ref{prop615615}.)
Here ${\rm Poly}(\mathfrak{Fuk}(X_3))$ is an $A_{\infty}$ category
obtained from $\mathfrak{Fuk}(X_3)$ replacing the morphism
modules $CF(\mathcal L_3,\mathcal L'_3)$ by
${\rm Poly}(CF(\mathcal L_3,\mathcal L'_3))$.

The $A_{\infty}$ category ${\rm Poly}(\mathfrak{Fuk}(X_3))$ is curved.
Note that each objects of $\mathfrak{Fukst}(X_3)$
which is in the image of the functors
$\mathcal W_{\mathcal L_{23}}\circ\mathcal W_{\mathcal L_{12}}$
(resp.\ $\mathcal W_{\mathcal L_{13}}$)
comes with a choice of bounding cochains \smash{$b_3^{(2)}$} \big(resp.\ \smash{$b_3^{(1)}$}\big).
We can lift those choices
to a bounding cochain $\hat b$ such that
\smash{${\rm Eval}_{s=0}\bigl(\hat b\bigr) = b_3^{(1)}$}
and~\smash{${\rm Eval}_{s=1}\bigl(\hat b\bigr) = b_3^{(2)}$}.
(See the proof of Lemma~\ref{lem912}.)
We use \smash{$\hat b$} to eliminate curvature
and obtain an object of associated strict category
of ${\rm Poly}(\mathfrak{Fuk}(X_3))$, which we denote by
${\rm Poly}_{st}(\mathfrak{Fuk}(X_3))$.

By the proof of Proposition~\ref{prop912},
there exists a commutative diagram of
$A_{\infty}$ functors
\begin{equation}\label{diagram1022}
\begin{CD}
\mathfrak{Fukst}(X_1) @ >>>
\mathcal{FUNC}(\mathfrak{Fukst}(X_3)^{\rm op},\mathcal{CH}) \\
@ V{=}VV @ VV{{\rm Eval}_{s=0}^*}V\\
\mathfrak{Fukst}(X_1) @ > >> \mathcal{FUNC}({\rm Poly}_{st}(\mathfrak{Fuk}(X_3))
^{\rm op},\mathcal{CH})\\
@ A{=}AA @ AA{{\rm Eval}_{s=1}^*}A \\
\mathfrak{Fukst}(X_1) @ >>>
\mathcal{FUNC}(\mathfrak{Fukst}(X_3)^{\rm op},\mathcal{CH}).
\end{CD}
\end{equation}
Here the first horizontal arrow is
obtained by using \smash{$b_3^{(1)}$} and the third horizontal arrow is
obtained by using \smash{$b_3^{(2)}$}.
Since the right vertical arrows are
homotopy equivalences, we obtained
the required homotopy equivalence.

The proof of Theorem~\ref{thm93} is complete.
\end{proof}

\subsection{The compatibility as bi-functors}
\label{sec:compfunc2}

We can strengthen Theorem~\ref{thm93} as follows.

\begin{thm}\label{thm109}
The next diagram commutes up to homotopy equivalence of
unital, strict and gapped filtered $A_{\infty}$ bi-functors:
\[
\begin{CD}
\mathfrak{Fukst}(-X_1 \times X_2)
\times \mathfrak{Fukst}(-X_2 \times X_3) @ >>>
\mathfrak{Fukst}(-X_1 \times X_3) \\
@ VVV @ VVV\\
\displaystyle{\mathcal{FUNC}(\mathfrak{Fukst}(X_1),\mathfrak{Fukst}(X_2))
\atop \times \mathcal{FUNC}(\mathfrak{Fukst}(X_2),\mathfrak{Fukst}(X_3))} @ > >> \mathcal{FUNC}(\mathfrak{Fukst}(X_1),\mathfrak{Fukst}(X_3)).
\end{CD}
\]
Here the first horizontal arrow is \eqref{form91}
and the second horizontal arrow is \eqref{form92}
in the case of~${\mathscr C_i = \mathfrak{Fukst}(X_i)}$.
The vertical arrows are correspondence bi-functors.
\end{thm}

The proof will be given in Section~\ref{sec:compfuncmain}.
\begin{rem}
Theorem~\ref{thm109} enhances Theorem~\ref{thm93},
and Theorem~\ref{thm93} enhances Proposition~\ref{prop912}.
Below we explain the difference between those three statements.
Theorem~\ref{thm109}
is a~coincidence between two
bi-functors
\begin{equation}\label{109defiefunc}
\mathfrak{Fukst}(-X_1 \times X_2)
\times \mathfrak{Fukst}(-X_2 \times X_3)
\to
\mathcal{FUNC}(\mathfrak{Fukst}(X_1),\mathfrak{Fukst}(X_3)).
\end{equation}
We first fix an object $\mathcal L_{12}$ (resp.\ $\mathcal L_{23}$)
of $\mathfrak{Fukst}(-X_1 \times X_2)$ (resp.\ $\mathfrak{Fukst}(-X_2 \times X_3)$).
Then the two bi-functors \eqref{109defiefunc} give two
objects of
$\mathcal{FUNC}(\mathfrak{Fukst}(X_1),\mathfrak{Fukst}(X_3))$.
The coincidence of those two objects, which are the
functors: $\mathfrak{Fukst}(X_1) \to \mathfrak{Fukst}(X_3)$,
is Theorem~\ref{thm93}.
Note that a functor: $\mathfrak{Fukst}(X_1) \to \mathfrak{Fukst}(X_3)$
gives a set theoretical map:
$\mathfrak{OB}(\mathfrak{Fukst}(X_1)) \to \mathfrak{OB}(\mathfrak{Fukst}(X_3))$.
The coincidence of two such set theoretical maps
is Proposition
\ref{prop912}.

Theorem~\ref{thm93} contains the coincidence of the morphism
parts of the functors: $\mathfrak{Fukst}(X_1) \to \mathfrak{Fukst}(X_3)$.
To prove Theorem~\ref{thm93}, we proved that \eqref{1018888}
is homotopy equivalence of left-$\mathfrak{Fukst}(X_1)$
and right-$\mathfrak{Fukst}(X_3)$ bi-modules.

Theorem~\ref{thm109} includes
statements on the coincidence of the way
the morphisms of $\mathfrak{Fukst}(-X_1\allowbreak \times X_2)$
and of $\mathfrak{Fukst}(-X_2 \times X_3)$
are mapped by
\eqref{109defiefunc}.
In the homology level, it implies the following.
Suppose $\mathcal L_{12}$, $\mathcal L'_{12}$
(resp.\ $\mathcal L_{23}$, $\mathcal L'_{23}$)
are objects of $\mathfrak{Fukst}(-X_1 \times X_2)$
(resp.\ $\mathfrak{Fukst}(-X_2 \times X_3)$)
and $\mathcal L_1$, $\mathcal L'_1$
are objects of $\mathfrak{Fukst}(X_1)$.

\eqref{109defiefunc} defines a homomorphism
\begin{equation}\label{form102325}
HF(\mathcal L_{12},\mathcal L'_{12})
\otimes
HF(\mathcal L_{23},\mathcal L'_{23})
\to
\operatorname{Hom}(HF(\mathcal L_1,\mathcal L'_1),
HF(\mathcal L_3,\mathcal L'_3)).
\end{equation}
Here $\mathcal L_3$ (resp.\ $\mathcal L'_3$) is obtained by
transforming $\mathcal L_1$ (resp.\ $\mathcal L'_1$)
via the composition of $\mathcal L_{12}$ and~$\mathcal L_{23}$
(resp.\ $\mathcal L'_{12}$ and $\mathcal L'_{23}$).
Theorem~\ref{thm109} implies that
the homomorphisms \eqref{form102325} obtained
by the following two different ways coincide.

The first way to obtain \eqref{form102325} is the following.
Let $\mathcal L_{13}$ (resp.\ $\mathcal L'_{13}$)
be the composition of~$\mathcal L_{12}$ (resp.\
$\mathcal L'_{12}$) and $\mathcal L_{23}$ (resp.\
$\mathcal L'_{23}$).
Then the composition bi-functor induces a homomorphism
\begin{equation}\label{form101414}
HF(\mathcal L_{12},\mathcal L'_{12})
\otimes
HF(\mathcal L_{23},\mathcal L'_{23})
\to
HF(\mathcal L_{13},\mathcal L'_{13}).
\end{equation}
On the other hand, by \eqref{form622}, we have
\begin{equation}\label{form101555}
HF(\mathcal L_{13},\mathcal L'_{13})
\to
\operatorname{Hom}(HF(\mathcal L_1,\mathcal L'_1),
HF(\mathcal L_3,\mathcal L'_3)).
\end{equation}
The composition of \eqref{form101414} and
\eqref{form101555} defines a homomorphism \eqref{form102325}.

The second way to obtain \eqref{form102325} is the following.
We have the following homomorphisms from \eqref{form622}:
\begin{align}
HF(\mathcal L_{12},\mathcal L'_{12})
&\to
\operatorname{Hom}(HF(\mathcal L_1,\mathcal L'_1),
HF(\mathcal L_2,\mathcal L'_2)), \nonumber\\
HF(\mathcal L_{23},\mathcal L'_{23})
&\to
\operatorname{Hom}(HF(\mathcal L_2,\mathcal L'_2),
HF(\mathcal L_3,\mathcal L'_3)).\label{form222201555}
\end{align}
Here $\mathcal L_{12}$ (resp.\ $\mathcal L'_{12}$) transforms
$\mathcal L_1$ (resp.\ $\mathcal L'_1$) to
$\mathcal L_2$ (resp.\ $\mathcal L'_2$).
On the other hand, the composition of homomorphisms
define a homomorphism
\begin{gather}
\operatorname{Hom}(HF(\mathcal L_1,\mathcal L'_1),
HF(\mathcal L_2,\mathcal L'_2))
\otimes
\operatorname{Hom}(HF(\mathcal L_2,\mathcal L'_2),
 HF(\mathcal L_3,\mathcal L'_3))\nonumber\\
\qquad\to
\operatorname{Hom}(HF(\mathcal L_1,\mathcal L'_1),
HF(\mathcal L_3,\mathcal L'_3)).\label{form233301555}
\end{gather}
The composition of \eqref{form222201555} and
\eqref{form233301555} is the second way to obtain
\eqref{form102325}.
\end{rem}

To prove Theorem~\ref{thm109}, we need more homological algebra.
In Section~\ref{sec:compfunc},
we used the derived tensor product to define the composition bi-functor
of functor categories.
In this subsection, we define the derived Hom functor.

\begin{defn}
Let $\mathscr C$ and $\mathscr C_{(i)}$, $i=1,2$
be strict, unital and gapped filtered $A_{\infty}$ categories and
$\mathfrak D_1$ a left-$\mathscr C$,$\mathscr C_{(1)}$
right-$\mathscr C_{(2)}$ tri-module. For $c \in \mathfrak{OB}(\mathscr C)$,
we define a left-$\mathscr C_{(1)}$ right-$\mathscr C_{(2)}$ bi-module
$\mathfrak D\vert_{c}$ as follows:
\begin{enumerate}\itemsep=0pt
\item[(1)]
If $c_i \in \mathfrak{OB}(\mathscr C_{(i)})$, then
$
\mathfrak D\vert_{c}(c_1,c_2) = \mathfrak D(c,c_1;c_2)$.
\item[(2)]
For ${\bf x} \in B_{k_1}\mathscr C_{(1)}[1](c_{1};c'_{1})$,
${\bf y} \in B_{k_2}\mathscr C_{(2)}[1](c'_{2},c_{2})$
and $v \in \mathfrak D\vert_{c}(c'_{1},c'_2) = \mathfrak D(c,c'_{1};c'_2)$,
we define~${
\mathfrak n_{k_1,k_2}({\bf x};v;{\bf y}) \in \mathfrak D\vert_{c}(c_1,c_2) = \mathfrak D(c,c_1;c_{2})}
$
by the tri-module structure on $\mathfrak D$.
This is the structure operation $\mathfrak n_{k_1,k_2}$ of $\mathfrak D\vert_{c}$.
\end{enumerate}

\end{defn}
\begin{defn}\label{defn104}
Let $\mathscr C$ and $\mathscr C_{(i)}$, $i=1,2,3,4$, be strict
filtered $A_{\infty}$ categories
and $\mathfrak D_{1}$ (resp.~$\mathfrak D_{2}$) be a left-$\mathscr C$,$\mathscr C_{(1)}$ right-$\mathscr C_{(2)}$ (resp.\ left-$\mathscr C$,$\mathscr C_{(3)}$ right-$\mathscr C_{(4)}$) filtered $A_{\infty}$ tri-module.

Let $c_i \in \mathfrak{OB}(\mathscr C_{(i)})$.
We define $\mathfrak{Hom}_{\mathscr C}(\mathfrak D_{1},\mathfrak D_{2})(c_2,c_3;c_1,c_4)$
\index[syindex]{Homc1d1@$\mathfrak{Hom}_{\mathscr C_1}(\mathfrak D_{1},\mathfrak D_{2})(c_1,c_3)$}
as the set of objects
\[
\mathfrak f = (\mathfrak f_{k_2;c,c'})_{c,c' \in \mathfrak{OB}(\mathscr C); k
= 0,1,2,\dots}
\]
such that
$
\mathfrak f_{k;c,c'} \colon B_{k}\mathscr C[1](c,c') \otimes \mathfrak D_{1}\vert_{c'}(c_1,c_2)
\to \mathfrak D_{2}\vert_{c}(c_3,c_4)
$
is a filtered $\Lambda_0$ module homomorphism.
\end{defn}
\begin{rem}
We remark that $\mathfrak{Hom}_{\mathscr C}(\mathfrak D_{1},\mathfrak D_{2})(c_2,c_3;c_1,c_4)$
is the direct {\it product}
\[
\prod_{c,c'}
\operatorname{Hom}(B_{k}\mathscr C[1](c,c') \otimes \mathfrak D_{1}\vert_{c'}(c_1,c_2),\mathfrak D_{2}\vert_{c'}(c_3,c_4)).
\]
In the definition of derived tensor product, we used
direct {\it sum}, see \eqref{form107}.
\end{rem}

\begin{lemdef}\label{lem1016}
There exists a left-$\mathscr C_{(2)}$,$\mathscr C_{(3)}$
right-$\mathscr C_{(1)}$,$\mathscr C_{(4)}$ multi-module, denoted by
$\mathfrak{Hom}_{\mathscr C}(\mathfrak D_{1},\mathfrak D_{2})$, so that
$
(c_2,c_3;c_1,c_4) \mapsto \mathfrak{Hom}_{\mathscr C}(\mathfrak D_{1},\mathfrak D_{2})(c_2,c_3;c_1,c_4)
$
in Definition {\rm\ref{defn104}} is its object part.
$($We define the boundary operator of the right-hand side during the proof.$)$

We write it $\mathfrak{Hom}_{\mathscr C}(\mathfrak D_{1},\mathfrak D_{2})$ and
call it the left $\mathscr C$ {\rm hom-module}. \index{hom-module}
\end{lemdef}
\begin{proof}
Let ${\bf x}_{(1)} \in B_{k_1}\mathscr C_{(1)}[1](c_1,c'_1)$,
${\bf y}_{(2)} \in B_{k_2}\mathscr C_{(2)}[1](c'_2,c_2)$,
${\bf x}_{(3)} \in B_{k_3}\mathscr C_{(3)}[1](c_3,c'_3)$,
${\bf y}_{(4)} \in B_{k_4}\mathscr C_{(4)}[1](c_4,c'_4)$,
and $\mathfrak f \in \mathfrak{Hom}_{\mathscr C}(\mathfrak D_{1},\mathfrak D_{2})(c_2,c_3;c_1,c_4)$.

We define
$
\mathfrak n_{k_1,k_2,k_3,k_4}({\bf y}_{(2)},{\bf x}_{(3)},\mathfrak f,{\bf x}_{(1)},{\bf y}_{(4)}) = \mathfrak g
\in \mathfrak{Hom}_{\mathscr C}(\mathfrak D_{1},\mathfrak D_{2})(c'_2,c_3;c'_1,c_4)$,
as follows.
We put $\mathfrak g = 0$ if $k_1+k_2 \ne 0$ and $k_3 +k_4 \ne 0$.

If $k_3 +k_4 = 0$ and $k_1+k_2 \ne 0$, we define
\begin{align}
\mathfrak n_{k_1,k_2,0,0}({\bf y}_{(2)},\mathfrak f,{\bf x}_{(1)})({\bf z},v)
={}&\mathfrak g_{k;c_2,c'_2}({\bf z},v) \nonumber\\
:={}&
- \sum_{c} (-1)^{*}\mathfrak f({\bf z}_{c;1},\mathfrak n({\bf z}_{c;z,2},{\bf x}_{(1)},v,{\bf y}_{(2)})),\label{form102200}
\end{align}
with
\[
* =\deg'{\bf z}_{c;1}+\deg \mathfrak f
+ \deg'{\bf y}_{(2)}
(\deg \mathfrak f+\deg v+\deg'{\bf x}+\deg'{\bf z})
+
\deg'{\bf x}_{(1)}\deg'{\bf z}.
\]
Here $v \in D_{1}(c',c'_2,c_3;c'_1,c_4)$, ${\bf z} \in B_k\mathscr C[1](c,c')$,
$\Delta {\bf z} = \sum_c {\bf z}_{c;1} \otimes {\bf z}_{c;2}$
and $\mathfrak n$ is the structure operation of $\mathfrak D_{1}$.

If $k_3 +k_4 \ne 0$ and $k_1+k_2 = 0$, we define
\begin{equation}\label{form1023koko}
\mathfrak n_{0,0,k_3,k_4}({\bf x}_{(3)},\mathfrak f,{\bf y}_{(4)})({\bf z},v)
=
\mathfrak g_{k;c_2,c'_2}({\bf z};v)
:=
\sum_{c} (-1)^*\mathfrak n'({\bf z}_{c;1},{\bf x}_{(3)},\mathfrak f({\bf z}_{c;2},v),{\bf y}_{(4)}),
\end{equation}
with
$
* = \deg{\bf y}_{(4)}(\deg'{\bf z} + \deg v)
+ \deg'{\bf x}_{(3)}\deg'{\bf z}_{c;1}
+ \deg\mathfrak f \deg'{\bf z}_{c;1}$.
Here
$v$, ${\bf z}$,
${\bf z}_{c;1},{\bf z}_{c;2}$ are as above
and $\mathfrak n'$ is the structure operation of $\mathfrak D_{32}$.

If $k_1 = k_2 = k_3 =k_4= 0$, we put
\begin{align}
\mathfrak n_{0,0}(1,\mathfrak f,1)({\bf z};v)
={}&
\mathfrak g_{k_2;c_2,c'_2}({\bf z};v)\nonumber \\
={}& \sum_{c} \mathfrak n'({\bf z}_{c;1},\mathfrak f({\bf z}_{c;2},v)) -
(-1)^{\deg \mathfrak f+\deg'{\bf z}_{c;1}}\sum_c\mathfrak f({\bf z}_{c;1},\mathfrak n({\bf z}_{c;2},v))\nonumber\\
&
- (-1)^{\deg \mathfrak f}\mathfrak f\bigl(\hat d{\bf z},v\bigr).\label{form1023+}
\end{align}
Note that all the signs in \eqref{form102200}, \eqref{form1023koko} and \eqref{form1023+} are by Koszul rule.

We can check $A_{\infty}$ relation as follows.
(Since the signs are by Koszul rule, the fact that the equality
holds {\it with signs} is in fact automatic.)
Let $\hat d$ be a map from
\begin{gather*}
\bigoplus_{c'_2,c'_1,c,c'}
B\mathscr{C}_{(2)}[1](c_2,c'_2) \\
\qquad{}\otimes {\operatorname{Hom}}(B\mathscr{C}[1](c,c')\otimes\mathfrak D_{1}(c',c'_2,c_3;c_1,c_4),
 \mathfrak D_{2}(c,c_2,c_3;c'_1,c_4))
\otimes B\mathscr{C}_{(1)}[1](c'_1,c_1)
\end{gather*}
to itself which is the coderivation induced by the structure operations.
We will prove
\begin{equation}\label{Ainfinhom}
\bigl(\mathfrak n\circ \hat d\bigr)({\bf y},\mathfrak f,{\bf x})({\bf z};v) = 0.
\end{equation}
Suppose $k_1 = k_2 = k_3 =k_4= 0$ for simplicity.
We have
\begin{align*}
\bigl(\mathfrak n'\circ \hat d\bigr)(\mathfrak f)({\bf z},v
={}&
\sum_c (-1)^{\deg'{\bf z}_{c;1}(\deg\mathfrak f + 1)}\mathfrak n'({\bf z}_{c;1},{\mathfrak n(\mathfrak f})({\bf z}_{c;2},v))\\
&+
\sum_c (-1)^{\deg' {\bf z}_{c;1}+\deg \mathfrak f}\mathfrak n(\mathfrak f)({\bf z}_{c;1},\mathfrak n({\bf z}_{c;2},v))
\\
&+ (-1)^{\deg \mathfrak f}
(\mathfrak n(\mathfrak f)) \bigl(\hat d{\bf z},v\bigr)
\\
={}& \sum_c (-1)^{\deg'{\bf z}_{c;1}+ \deg\mathfrak f(\deg'{\bf z}_{c;1}+\deg'{\bf z}_{c;2})} \mathfrak n'({\bf z}_{c;1},\mathfrak n'({\bf z}_{c;2},{\mathfrak f}({\bf z}_{c;3},v)))
\\
& + \sum_c (-1)^{\deg \mathfrak f + 1+\deg'{\bf z}_{c;2}
+ \deg'{\bf z}_{c;1}(\deg\mathfrak f + 1)}\mathfrak n'({\bf z}_{c;1},\mathfrak f({\bf z}_{c;2},\mathfrak n({\bf z}_{c;3},v))) \\
& + \sum_c (-1)^{\deg'{\bf z}_{c;1}(\deg\mathfrak f + 1)+(\deg\mathfrak f + 1)} \mathfrak n'\bigl({\bf z}_{c;1},\mathfrak f\bigl(\hat d {\bf z}_{c;2},v\bigr)\bigr)\\
& + \sum_c (-1)^{\deg' {\bf z}_{c;1}+\deg' {\bf z}_{c;2}+\deg \mathfrak f+
\deg'{\bf z}_{c;1}\deg\mathfrak f}\mathfrak n'({\bf z}_{c;1},\mathfrak f({\bf z}_{c;2},\mathfrak n({\bf z}_{c;3},v)))\\
& + \sum_c (-1)^{\deg'{\bf z}_{c;2}+1}\mathfrak f({\bf z}_{c;1},\mathfrak n({\bf z}_{c;2},\mathfrak n({\bf z}_{c;3},v)))
\\
&  + \sum_c (-1)^{1 + \deg' {\bf z}_{c;1}}\mathfrak f\bigl(\hat d {\bf z}_{c;1},\mathfrak n(
{\bf z}_{c;2},v)\bigr) \\
& +\sum_c (-1)^{\deg \mathfrak f+
\deg \mathfrak f\deg'{\bf z}_{c;1}+\deg'{\bf z}_{c;1}} \mathfrak n'\bigl({\bf z}_{c;1},\mathfrak f\bigl(\hat d{\bf z}_{c;2},v\bigr)\bigr)\\
& +\sum_c (-1)^{\deg \mathfrak f + \deg \mathfrak f(\deg'{\bf z}_{c;1}+1)} \mathfrak n'\bigl(
\hat d{\bf z}_{c;1},\mathfrak f({\bf z}_{c;2},v)\bigr)\\
& +\sum_c (-1)^{1} \mathfrak f\bigl({\bf z}_{c;1},\mathfrak n\bigl(\hat d{\bf z}_{c;2},v\bigr)\bigr)\\
& +\sum_c (-1)^{1+\deg'{\bf z}_{c;1}+1}\mathfrak f(\hat d{\bf z}_{c;1},\mathfrak n({\bf z}_{c;2},v)).
\end{align*}
The 1st and 8th terms of the right-hand side cancel by the $A_{\infty}$ relation of $\mathfrak n'$.
The 2nd and 4th terms cancel. The 3rd and 7th terms cancel.
The 5th and 9th terms cancel by the $A_{\infty}$ relation of~$\mathfrak n$.
The 6th and 10th terms cancel.
Thus we checked \eqref{Ainfinhom} in the case $k_1 = k_2 = k_3 =k_4= 0$.
The other cases are similar.
\end{proof}

\begin{lemdef}
There exists a filtered $A_{\infty}$ bi-functor
\[
\mathcal{MUMOD}(\mathscr C,\mathscr C_{(1)};\mathscr C_{(2)})
\times
\mathcal{MUMOD}(\mathscr C,\mathscr C_{(3)};\mathscr C_{(4)})
\to
\mathcal{MUMOD}(\mathscr C_{(2)},\mathscr C_{(3)}; \mathscr C_{(1)},\mathscr C_{(4)}),
\]
which is given by Lemma--Definition {\rm\ref{lem1016}} for the object part.

We call this bi-functor the {\rm derived hom-functor}
\index{derived hom-functor} and write its as $\mathfrak{Hom}$.
\index[syindex]{Hom@$\mathfrak{Hom}$}
\end{lemdef}
\begin{proof}
Let $\mathscr C$ and $\mathscr C_{(i)}$, $i=1,2,3,4$, be strict
filtered $A_{\infty}$ categories
and $\mathfrak D_{1}$, $\mathfrak D'_{1}$ (resp.\ $\mathfrak D_{2}$, $\mathfrak D'_{2}$) be left-$\mathscr C$,$\mathscr C_{(1)}$ right-$\mathscr C_{(2)}$ (resp.\ left-$\mathscr C$,$\mathscr C_{(3)}$ right-$\mathscr C_{(4)}$) filtered $A_{\infty}$ tri-module.

Suppose
$\mathfrak F_{1} \colon \mathfrak D'_{1} \to \mathfrak D_{1}$ and
$\mathfrak F_{2} \colon \mathfrak D_{2} \to \mathfrak D'_{2}$ are tri-module homomorphisms.
We will define
\[
(\mathfrak F^*_{1},\mathfrak F_{2 *})\colon\ \mathfrak{Hom}_{\mathscr C_2}(\mathfrak D_{1},\mathfrak D_{2})
\to \mathfrak{Hom}_{\mathscr C_2}(\mathfrak D'_{1},\mathfrak D'_{2}).
\]
Let $\widehat{\mathfrak f}
= (\mathfrak f^{c_1,c_2,c_3,c_4})_{c_i \in \mathfrak{Ob}(\mathscr C_{(i)})};
\mathfrak{Hom}_{\mathscr C}(\mathfrak D_{1},\mathfrak D_{2})$.
Here $\mathfrak f^{c_1,c_2,c_3,c_4} = \bigl(\mathfrak f^{c_1,c_2,c_3,c_4}_{c,c'}\bigr)_{
c, c' \in \mathfrak{Ob}(\mathscr C)}$ and
\[
\mathfrak f^{c_1,c_2,c_3,c_4}_{c,c'} \colon\ B\mathscr C[1](c,c') \otimes\mathfrak D_{1}\vert_{c'}(c_2,c_3;c_1,c_4)
\to \mathfrak D_{2}\vert_{c'}(c_2,c_3;c_1,c_4).
\]
We define
$\widehat{\mathfrak g} = (\mathfrak F^*_{1},\mathfrak F_{2 *})\bigl(\widehat{\mathfrak f}\bigr)$ as follows.
\smash{$\widehat{\mathfrak g}
= (\mathfrak g^{c_1,c_2,c_3,c_4})_{c_i \in \mathfrak{Ob}(\mathscr C_{(i)})}$},
$\smash{\mathfrak g^{c_1,c_2,c_3,c_4} = \bigl(\mathfrak g^{c_1,c_2,c_3,c_4}_{c,c'} };
c, c' \in \mathfrak{Ob}(\mathscr C)\bigr)$,
and \smash{$
\mathfrak g^{c_1,c_2,c_3,c_4}_{c,c'} \colon
B\mathscr C(c,c') \otimes \mathfrak D'_{1}\vert_{c'}(c_1,c_2;c_3,c_4)
\to \mathfrak D'_{2}\vert_{c}(c_1,c_2;c_3,c_4)
$} is
\[
\mathfrak g^{c_1,c_2,c_3,c_4}_{c,c'}({\bf z},v)
:=
\sum_c (-1)^{\deg\mathfrak f\deg' {\bf z}_{c:1}} \mathfrak F_{2}({\bf z}_{c;1},\mathfrak f({\bf z}_{2;1},\mathfrak F_{1}({\bf z}_{c;3},v))).
\]
Here $v \in \mathfrak D'_{1}\vert_{c}(c_1,c_2;c_3,c_4)$, ${\bf z} \in B\mathscr C(c,c')$
and $((\Delta \otimes 1)\circ \Delta)({\bf y})
= \sum_c {\bf z}_{c;1} \otimes {\bf z}_{c;2} \otimes {\bf z}_{c;3}$.

It is straightforward to check that $(\mathfrak F^*_{1},\mathfrak F_{2 *})$ is a chain map
and multi-module homomorphism.
Moreover, if $\mathfrak F_{1} \circ \mathfrak G_{1}= \mathfrak H_1$,
$\mathfrak G_{2} \circ \mathfrak F_{2}= \mathfrak H_2$,
then
$
(\mathfrak H^*_{1},
\mathfrak H_{2 *} ) = (\mathfrak G^*_{1},\mathfrak G_{2 *}) \circ
(\mathfrak F^*_{1},\mathfrak F_{2 *})
$.
Thus we obtain a required bi-functor.
(It is actually a DG-functor.)
\end{proof}

The next proposition is a Hom version of Proposition~\ref{lem9100}.
\begin{prop}\label{prop1019}
Let $\mathscr C$, $\mathscr C_{(i)}$, $i=1,2,3$, be strict unital and gapped
filtered $A_{\infty}$ categories, and
$\mathcal F\colon \mathscr C_{(1)} \to \mathscr C$ and
$\mathcal G\colon \mathscr C \times \mathscr C_{(2)} \to \mathscr C_{(3)}$
strict, unital and gapped
filtered $A_{\infty}$ $($bi-$)$ functors.
We consider
$
\mathfrak{Yon}\circ \mathcal F \colon \mathscr C_{(1)}
\to \mathcal{FUNC}(\mathscr C^{\rm op},\mathcal{CH})
$
and regard it as a bi-functor
$
\mathscr C^{\rm op} \times \mathscr C_{(1)} \to \mathcal{CH}$.
It can be regarded as a left-$\mathscr C$, right-$\mathscr C_{(1)}$,
bi-module, which we denote by $\mathfrak D_{(1)}$.
We apply $($bi-module analogue of$)$ the
relative Yoneda functor to $\mathcal G$
to obtain $\mathfrak{RYon}_{\rm ob}(\mathcal G)$, which becomes a~left-$\mathscr C$,$\mathscr C_{(2)}$ right $\mathscr C_{(3)}$ tri-module,
which we denote by $\mathfrak D_{(2)}$.

We next consider the composition $\mathcal G\circ \mathcal F\colon \mathscr C_{(1)}
\times \mathscr C_{(2)} \to \mathscr C_{(3)}$ and
apply $($the bi-module analogue of$)$ the
relative Yoneda functor.
We obtain a left-$\mathscr C_{(1)}$, $\mathscr C_{(2)}$ right-$\mathscr C_{(3)}$
bi-module and denote it by $\mathfrak D_{(3)}$.
Now we claim that $\mathfrak D_{(3)}$ is homotopy equivalent
to
$\mathfrak{Hom}_{\mathscr C}(\mathfrak D_{(1)},\mathfrak D_{(2)})$ as
a left-${\mathscr C}_{(1)}$,${\mathscr C}_{(2)}$ right-${\mathscr C}_{(3)}$
tri-module.
Here the left $\mathscr C_{(1)}$ module structure on
$\mathfrak{Hom}_{\mathscr C}(\mathfrak D_{(1)},\mathfrak D_{(2)})$
is induced from the right $\mathscr C_{(1)}$ module
structure on $\mathfrak D_{(1)}$.
$($We do not use left $\mathscr C$
module structure on $\mathfrak D_{(2)}$ to define this
left $\mathscr C_{(1)}$ module structure.$)$

\end{prop}

We remark that by definition
$\mathfrak D_{(3)}$ is induced from $\mathfrak D_{(2)}$
by $\mathcal F$.
The proof will be given in Section~\ref{sec:compfuncalglem}.
\begin{rem}
We consider the case when $\mathscr C$ and $\mathscr C_{(1)}$ are unital associative algebras,
$\mathscr C_{(2)}$ is trivial,
and $\mathcal F$ is a unital ring homomorphism.
We use the notation of Remark~\ref{rem1011}.
Then $\mathfrak D_{(1)}$ is the bi-module ${}_{\mathscr C}\mathscr C_{\mathscr C_{(1)}}$
and $\mathfrak D_{(2)}$ is given by a left $\mathscr C$ right ${\mathscr C_{(3)}}$ bimodule
${}_{\mathscr C}D_{\mathscr C_{(3)}}$.

Therefore, $\mathfrak{Hom}_{\mathscr C}(\mathfrak D_{(1)},\mathfrak D_{(2)})$
is $\operatorname{Hom}_{\mathscr C}({}_{\mathscr C}\mathscr C_{\mathscr C_{(1)}},
{}_{\mathscr C}D_{\mathscr C_{(3)}})$.
The map sending $\varphi$ to $\varphi({\bf e})$ gives an isomorphism
between $\operatorname{Hom}_{\mathscr C}({}_{\mathscr C}\mathscr C_{\mathscr C_{(1)}},
{}_{\mathscr C}D_{\mathscr C_{(3)}})$ and ${}_{\mathscr C_{(1)}}D_{\mathscr C_{(3)}}$
as left $\mathscr C_{(1)}$ right $\mathscr C_{(3)}$ modules.
Note that the left $\mathscr C_{(1)}$ action on ${}_{\mathscr C_{(1)}}D_{\mathscr C_{(3)}}$
is defined by $\mathcal F\colon \mathscr C_{(1)} \to \mathscr C$
and the left action of $\mathscr C$.

The bi-module ${}_{\mathscr C_{(1)}}D_{\mathscr C_{(3)}}$ corresponds to the
composition $\mathcal G\circ\mathcal F$.
We thus checked Proposition~\ref{prop1019} in this case.
\end{rem}

\subsection[Proof of Theorem~\ref{thm109}]{Proof of Theorem~\ref{thm109}}
\label{sec:compfuncmain}

\begin{proof}[Proof of Theorem~\ref{thm109}]
We first consider the composition
\begin{gather}
\mathfrak{Fukst}(-X_1 \times X_2)
\times \mathfrak{Fukst}(-X_2 \times X_3) \nonumber\\
\qquad\to
\mathcal{FUNC}(\mathfrak{Fukst}(X_1),\mathfrak{Fukst}(X_2))
\times \mathcal{FUNC}(\mathfrak{Fukst}(X_2),\mathfrak{Fukst}(X_3)) \nonumber\\
\qquad\to
\mathcal{FUNC}(\mathfrak{Fukst}(X_1),\mathfrak{Fukst}(X_3))\label{form1027}
\end{gather}
and compose it with the relative Yoneda functor.
By the commutativity of diagram \eqref{dia1010}
(see Propositions \ref{lem9100}), the composition \eqref{form1027}
is homotopy equivalent to
\begin{gather}
\mathfrak{Fukst}(-X_1 \times X_2)
\times \mathfrak{Fukst}(-X_2 \times X_3)\nonumber \\
\qquad\to
\mathcal{BIMOD}(\mathfrak{Fukst}(X_1),\mathfrak{Fukst}(X_2))^{\rm op} \times\mathcal{BIMOD}(\mathfrak{Fukst}(X_2),\mathfrak{Fukst}(X_3))^{\rm op}\nonumber
\\
\qquad\to
\mathcal{BIMOD}(\mathfrak{Fukst}(X_1),\mathfrak{Fukst}(X_3))^{\rm op},\label{1028form}
\end{gather}
where the first functor is the composition of the correspondence bi-functor and
the relative Yoneda functor and the
second functor is the derived tensor product.
Let $\mathcal L_{12}$ (resp.\ $\mathcal L_{23}$) be an object of
$\mathfrak{Fukst}(-X_1 \times X_2)$ (resp.\ $\mathfrak{Fukst}(-X_2 \times X_3)$).
By the definition of the correspondence bi-functor
the first functor is as follows.

Let $\mathcal L_i$ be an object of $\mathfrak{Fukst}(X_i)$ for $i=1,2,3$.
Then $\mathcal L_{12}$ (resp.\ $\mathcal L_{23}$) is sent to the
left-$\mathfrak{Fukst}(X_1)$ right-$\mathfrak{Fukst}(X_2)$ bi-module
$\mathscr{CF}(\mathbb L_1,\mathbb L_{12};\mathbb L_2)$
(resp.\ left-$\mathfrak{Fukst}(X_2)$ right-$\mathfrak{Fukst}(X_3)$ bi-module
$\mathscr{CF}(\mathbb L_2,\mathbb L_{23};\mathbb L_3)$),
which sends $\mathcal L_1$ and $\mathcal L_2$
(resp.\ $\mathcal L_2$ and $\mathcal L_3$) to
$
CF(\mathcal L_1,\mathcal L_{12};\mathcal L_{2})
$
(resp.~$
CF(\mathcal L_2,\mathcal L_{23};\mathcal L_{3})
$).
This is the object part of the functor.
The morphism part is determined by the
left-$\mathfrak{Fukst}(-X_1 \times X_2)$ module structure
of $
\mathscr{CF}(\mathbb L_1,\mathbb L_{12};\mathbb L_2)
$
(resp.
the left-$\mathfrak{Fukst}(-X_2\times X_3)$ module structure
of $
\mathscr{CF}(\mathbb L_2,\mathbb L_{23};\mathbb L_3)
$).

Therefore, by the definition of the derived tensor product,
the composition \eqref{1028form} sends
the pairs $(\mathcal L_{1},\mathcal L_{3})$, $(\mathcal L_{12},\mathcal L_{23})$
to
\begin{gather}
D_1(\mathcal L_{1},
\mathcal L_{12},\mathcal L_{23};\mathcal L_{3})\nonumber\\
\qquad =
\bigoplus_{\mathcal L_{2},\mathcal L'_{2}} CF(\mathcal L_1,\mathcal L_{12};\mathcal L_{2})
\otimes BCF[1](\mathcal L_{2},\mathcal L'_{2})
\otimes
CF(\mathcal L'_2,\mathcal L_{23};\mathcal L_{3}).\label{1029form}
\end{gather}
We consider \eqref{1029form} for various $\mathcal L_1$, $\mathcal L_3$,
$\mathcal L_{12}$, $\mathcal L_{23}$ and
obtain
the object part of the composition~\eqref{1028form}. The morphism part is
determined by the left $\mathfrak{Fukst}(-X_1 \times X_2)$, $\mathfrak{Fukst}(-X_2 \times X_3)$,
$\mathfrak{Fukst}(X_1)$, right $\mathfrak{Fukst}(X_3)$ quatro-module
structure of \eqref{1029form}.

We thus described the bi-functor \eqref{form1027}
composed with the relative Yoneda functor.

We next study the composition
\begin{align}
\mathfrak{Fukst}(-X_1 \times X_2)
\times \mathfrak{Fukst}(-X_2 \times X_3)
&\to
\mathfrak{Fukst}(-X_1 \times X_3) \nonumber\\
&\to
 \mathcal{FUNC}(\mathfrak{Fukst}(X_1),\mathfrak{Fukst}(X_3))\nonumber\\
 &\to
 \mathcal{BIMOD}(\mathfrak{Fukst}(X_1);\mathfrak{Fukst}(X_3))^{\rm op}.\label{form1030o}
\end{align}
By definition, the first functor
composed with
\begin{align*}
\mathfrak{Yon} \colon\ \mathfrak{Fukst}(-X_1 \times X_3)
&\to\mathcal{FUNC}(\mathfrak{Fukst}(-X_1 \times X_3)^{\rm op},\mathcal{CH})\\
&\cong \mathcal{BIMOD}(\mathfrak{Fukst}(-X_1 \times X_3),*)^{\rm op}
\end{align*}
is
given by $\mathscr{CF}(\mathbb L_{13};\mathbb L_{12},\mathbb L_{23})$ which is a left-$\mathfrak{Fukst}(-X_1 \times X_3)$ right-$\mathfrak{Fukst}(-X_1 \times X_2)$, $\mathfrak{Fukst}(-X_2\allowbreak \times X_3)$ tri-module.
(See Proposition~\ref{prop810}.)

We consider the composition of the second
and third functors in \eqref{form1030o}
and apply (the object part of) the relative
Yoneda functor $\mathfrak{YonR}_{\rm ob}$.
We then obtain a left-$\mathfrak{Fukst}(-X_1 \times X_3)$,
$\mathfrak{Fukst}(X_1)$ right-$\mathfrak{Fukst}(X_3)$
tri-module
$\mathscr{CF}(\mathbb L_{1},\mathbb L_{13};\mathbb L_{3})$.
(See Lemma--Definition~\ref{lemdef97}.)
Here left~${\mathfrak{Fukst}(-X_1 \times X_2)}$, $\mathfrak{Fukst}(-X_2 \times X_3)$
module structure on $\mathscr{CF}(\mathbb L_{1},\mathbb L_{13};\mathbb L_{3})$
is induced by its left $\mathfrak{Fukst}(-X_1 \times X_3)$ module structure
via the bi-functor $\mathfrak{Fukst}(-X_1 \times X_2)
\times \mathfrak{Fukst}(-X_2 \times X_3) \to \mathfrak{Fukst}(-X_1 \times X_3)$.

We next use Proposition~\ref{prop1019}.
We put
\begin{gather*}
\mathscr C_{(1)} = \mathfrak{Fukst}(-X_1 \times X_2)
\times \mathfrak{Fukst}(-X_2 \times X_3), \qquad
\mathscr C_{(2)} = \mathfrak{Fukst}(X_1), \\
\mathscr C_{(3)} = \mathfrak{Fukst}(X_3),\qquad
\mathscr C = \mathfrak{Fukst}(-X_1 \times X_3).
\end{gather*}
Then
\[
\mathfrak D_{(1)} = \mathscr{CF}(\mathbb L_{13};\mathbb L_{12},\mathbb L_{23}),\qquad
\mathfrak D_{(2)} = \mathscr{CF}(\mathbb L_{1},\mathbb L_{13};\mathbb L_{3}).
\]
$\mathfrak D_{(3)}$ is the pull-back of $\mathfrak D_{(2)}$
by
$\mathfrak{comp}\colon \mathfrak{Fukst}(-X_1 \times X_2)
\times \mathfrak{Fukst}(-X_2 \times X_3)
\to
\mathfrak{Fukst}(-X_1 \times X_3)$.

Proposition~\ref{prop1019} then implies that $\mathfrak D_{(3)}$ is homotopy equivalent to
\begin{equation}\label{form10352}
\mathfrak{Hom}_{\mathfrak{Fukst}(-X_1 \times X_3)}
(\mathscr{CF}(\mathbb L_{13};\mathbb L_{12},\mathbb L_{23}),
\mathscr{CF}(\mathbb L_{1},\mathbb L_{13};\mathbb L_{3}))
\end{equation}
as left-$\mathfrak{Fukst}(-X_1 \times X_2)$,$\mathfrak{Fukst}(-X_2 \times X_3)$,
$\mathfrak{Fukst}(X_1)$ and right-$\mathfrak{Fukst}(X_2)$
quatro-module.\footnote{Actually we use the variant of Proposition~\ref{prop1019}
where $\mathcal F\colon\mathcal C_{(1)} \to \mathcal C$ is replaced by
a bi-functor. The proof of the variant is the same as the proof of
 Proposition~\ref{prop1019}.}
Note that the quatro-module \eqref{form10352} associates
\begin{gather}
D_2(\mathcal L_{1},
\mathcal L_{12},\mathcal L_{23};\mathcal L_{3})\nonumber
\\
\qquad= \prod_{\mathcal L_{13},\mathcal L'_{13}} \operatorname{Hom}(BCF[1](\mathcal L_{13},\mathcal L'_{13})\otimes CF(\mathcal L'_{13};\mathcal L_{12},\mathcal L_{23}),CF(\mathcal L_{1},\mathcal L_{13};\mathcal L_{3}))\label{form1035}
\end{gather}
to $\mathcal L_{12}$, $\mathcal L_{23}$, $\mathcal L_{1}$, $\mathcal L_{3}$

We thus described two compositions
\[
\mathfrak{Fukst}(-X_1 \times X_2)\times \mathfrak{Fukst}(-X_2 \times X_3)
\to \mathcal{BIMOD}(\mathfrak{Fukst}(X_1),\mathfrak{Fukst}(X_3))^{\rm op},
\]
which are \eqref{1029form} and \eqref{form1035} together with
their quatro-module structures.
Theorem~\ref{thm109}
claims that they are homotopy equivalent as quatro-modules.
To prove it, we will construct a~quatro-module homomorphism
from \eqref{1029form} to \eqref{form1035}.

By definition, such a quatro-module homomorphism is a map
\begin{gather*}
\bigoplus_{\mathcal L'_1,\mathcal L'_3,\mathcal L'_{12},\mathcal L'_{23}}
BCF[1](\mathcal L_1,\mathcal L'_1) \otimes BCF[1](\mathcal L_{12},\mathcal L'_{12})
 \otimes BCF[1](\mathcal L_{23},\mathcal L'_{23})\\
\qquad{}\otimes D_1(\mathcal L'_{1},
\mathcal L'_{12},\mathcal L'_{23};\mathcal L'_{3}) \otimes BCF[1](\mathcal L'_{3},\mathcal L_{3}) \to D_2(\mathcal L_{1},
\mathcal L_{12},\mathcal L_{23};\mathcal L_{3}).
\end{gather*}
Therefore, it can be regarded as a homomorphism from
\begin{gather*}
BCF[1](\mathcal L_1,\mathcal L'_1)
\otimes
BCF[1](\mathcal L_{12},\mathcal L'_{12}) \otimes BCF[1](\mathcal L_{23},\mathcal L'_{23}) \\
\qquad
{}\otimes CF(\mathcal L'_1,\mathcal L'_{12};\mathcal L_{2})
\otimes BCF[1](\mathcal L_2,\mathcal L'_2)
\otimes
CF(\mathcal L'_2,\mathcal L'_{23};\mathcal L'_{3})\otimes
BCF[1](\mathcal L'_{3},\mathcal L_{3})
\\
\qquad{}\otimes
BCF[1](\mathcal L_{13},\mathcal L'_{13})\otimes CF(\mathcal L'_{13};\mathcal L_{12},\mathcal L_{23})
\end{gather*}
to
$CF(\mathcal L_{1},\mathcal L_{13};\mathcal L_{3})$.
The $Y$-diagram transformation \smash{$\mathscr{YT}^{\vec b}$} in \eqref{form10111}
is such a homomorphism and therefore defines a
pre-quatro-module homomorphism.
The condition that it becomes a~quatro-module homomorphism
is exactly the formula \eqref{form926}, which we proved
in Lemma~\ref{lem1014}.

To prove that this quatro-module homomorphism is a homotopy equivalence,
it suffices to show that the chain maps, which are parts of this
quatro-module homomorphism,
are chain homotopy equivalences
(see Proposition~\ref{prop1015}).
The chain map induced by \smash{$\mathscr{YT}^{\vec b}$} is
nothing but the chain homotopy equivalence \eqref{formula1020}
which we produced during the proof of Theorem~\ref{thm93}
in Section~\ref{sec:compfunc21}.

We can study the difference between two bounding cochains
\smash{$b_3^{(1)}$} and \smash{$b_3^{(2)}$} in the same way
as the last step of the proof of Theorem~\ref{thm93}
by enhancing diagram~\eqref{diagram1022}, so that it includes
left-$\mathfrak{Fukst}(-X_1 \times X_2)$, $\mathfrak{Fukst}(-X_2 \times X_3)$
structure.

The proof of Theorem~\ref{thm109} is now complete.
\end{proof}

\begin{rem}
To prove the commutativity of the diagram in Theorem~\ref{thm109}
for the object part, it suffices to show that
$D_1 = CF(\mathcal L_1;\mathcal L_{12};\mathcal L_{2})
\otimes BCF[1](\mathcal L_{2})
\otimes
CF(\mathcal L_2;\mathcal L_{23};\mathcal L_{3})$
is homotopy equivalent to
$D'_2 := CF(\mathcal L_{1},\mathcal L_{13};\mathcal L_{3})$
as left-$\mathfrak{Fukst}(X_1)$ right-$\mathfrak{Fukst}(X_3)$
bi-modules.

To prove the commutativity of the morphism part,
we need to include the compatibility of the
homotopy equivalence with the left $\mathfrak{Fukst}(-X_1 \times X_2) \times \mathfrak{Fukst}(-X_2 \times X_3)$
bi-module structures, as we have done above.
\end{rem}

\subsection[Proof of Propositions \ref{lem9100} and \ref{prop1019}]{Proof of Propositions \ref{lem9100} and \ref{prop1019}}
\label{sec:compfuncalglem}

In this subsection, we prove Propositions \ref{lem9100} and \ref{prop1019}.
We need certain calculations of the sign for the proof.
Note that in this paper the sign is almost always by the Koszul rule and
by this reason the cancellation with the sign is mostly automatic.
A certain nontrivial sign issue appears in this subsection
by the following reason.
We need to regard a filtered $A_{\infty}$ category~$\mathscr C$
itself as a left-$\mathscr C$ right-$\mathscr C$ bi-module.
In such a case an element $v$ of $\mathscr C(c,c')$ as an
element of bi-module appears with sign $(-1)^{\deg v}$ in the $A_{\infty}$
formula. In the case $v$ is regarded as an element of a morphism
complex of an $A_{\infty}$ category,
it appears with sign $(-1)^{\deg v+1}$ in the $A_{\infty}$
formula.

By several maps, which we will define in this subsection,
an element of $\mathscr C(c,c')$ as an
element of a bi-module in the domain becomes
an element of the morphism
complex in the co-domain or vice versa.
This process shifts the degree.
It is not obvious to understand the way how this process affects
the sign, since the Koszul rule does not tell it to us.
By this reason, we need to add a certain correction term
to the usual Koszul sign.
The author is unable to provide the general principle
on the way how the correction terms are determined.
Instead, he puts the correction terms `by hand'
(see, for example, \eqref{I12nocomplisign}) and check that
the sign works by a~calculation.\footnote{Actually a similar
problem occurs during the proof of Yoneda's lemma.}

Fortunately, this happens only in the purely algebraic situation
so that we do {\it not} need to understand the geometric
origin of the correction terms.
In fact, Propositions \ref{lem9100} and \ref{prop1019} are
algebraic statements and hold independent of the
origin of $A_{\infty}$ categories and functors in their statements.
For the construction of various operations
using moduli spaces, the fundamental formulas
among those operations are always with Koszul sign.
We will use this fact in Section~\ref{sec:orient}.

\begin{proof}[Proof of Proposition~\ref{lem9100}]
Let $\mathscr F_{i(i+1)} \colon \mathscr C_i \to \mathscr C_{i+1}$
be a filtered $A_{\infty}$ functor for $i=1,2$.
Let $c_1 \in \mathfrak{OB}(\mathscr C_1)$, $c_3 \in \mathfrak{OB}(\mathscr C_3)$.
We put
\begin{gather*}
 D^1(c_1,c_3) := \mathscr{C}_3((\mathscr F_{23})_{\rm ob}((\mathscr F_{12})_{\rm ob}(c_1)),c_3),
\\
D^2(c_1,c_3)
:= \bigoplus_{c_2,c'_2} \mathscr{C}_2((\mathscr F_{12})_{\rm ob}(c_1),c_2)
\,\widehat\otimes\, B\mathscr C_2(c_2,c'_2)
\,\widehat\otimes\, \mathscr{C}_3((\mathscr F_{23})_{\rm ob}(c'_2),c_3).
\end{gather*}
Note that $D^1$ is the object part of the bi-module
$\mathfrak{RYon}_{\rm ob}(\mathfrak F_{23} \circ \mathfrak F_{12})$ and
$D^2$ is the object part of the bi-module
$\mathfrak{ten}_{\rm ob}(\mathfrak{RYon}_{\rm ob}(\mathfrak F_{12}),\mathfrak{RYon}_{\rm ob}(\mathfrak F_{23}))$.

We define
$
\mathscr I_{12;0,0} \colon D^1(c_1,c_3) \to D^2(c_1,c_3)
$
by
\begin{equation}\label{form1024}
\mathscr I_{12;0,0}(z) = {\bf e}_{(\mathscr F_{12})_{\rm ob}(c_1)} \otimes 1 \otimes z.
\end{equation}
Here the symbol ${\bf e}_{(\mathscr F_{12})_{\rm ob}(c_1)}$ is the unity of the object
$(\mathscr F_{12})_{\rm ob}(c_1)$ and the symbol $1$ is an element of $B_0\mathscr C_2((\mathscr F_{12})_{\rm ob}(c_1),(\mathscr F_{12})_{\rm ob}(c_1))$, which is isomorphic to $\Lambda_0$. Hereafter, we omit $1$ from the notation.
It is obvious that $\mathscr I_{12;0,0}$ is a chain map.

We also define
$\mathscr I_{21} \colon D^2(c_1,c_3) \to D^1(c_1,c_3)$ by
\smash{$
\mathscr I_{21}(x,{\bf y},z) = (-1)^{\deg x}\mathfrak m_*\bigl(\widehat{\mathfrak F}_{23}(x,{\bf y}),z\bigr)$}.
Here $\mathfrak m$ is the structure operation of
$\mathscr C_3$.

Let $\mathfrak n$ is the $(0,0)$ part of the
left-$\mathscr C_1$ right-$\mathscr C_3$ bi-module structure
of $D^2(c_1,c_3)$.
Here we use the sign
convention so that degree of elements of bi-module is {\it not} shifted.
Namely,
\begin{align}
\mathfrak n(x\otimes{\bf y}\otimes z)
={}&
\sum_c (-1)^{\deg' {\bf y}_{1;c}} \mathfrak m(x,{\bf y}_{1;c}) \otimes
{\bf y}_{2;c} \otimes z\nonumber
\\
& +
\sum_c (-1)^{\deg x + \deg'{\bf y}_{1;c}}
x \otimes {\bf y}_{1;c} \otimes \mathfrak m\bigl(\widehat{\mathfrak F}_{23}({\bf y}_{2;c}),z\bigr)\nonumber
\\
& + (-1)^{\deg x} x \otimes \hat d({\bf y}) \otimes z.\label{defnnnnn}
\end{align}
This is a special case of \eqref{form108}.
We will check that $\mathscr I_{21;0,0}$ is a chain map.
We calculate
\begin{gather*}
(\mathscr I_{21}\circ \mathfrak n)(x,{\bf y},z) \\
\qquad= \sum_c (-1)^{\deg' {\bf y}_{c;1}}\mathscr I_{21}
(\mathfrak m(x,{\bf y}_{c;1}),{\bf y}_{c;2},z)
\\
\phantom{\qquad= }{}+ (-1)^{\deg x}\mathscr I_{21}\bigl(x,\hat d{\bf y},z\bigr)
+ \sum_c (-1)^{\deg x + \deg' {\bf y}_{c;1}} \mathscr I_{21}
\bigl(x,{\bf y}_{c;1},\mathfrak m\bigl(\widehat{\mathfrak F}_{23}({\bf y}_{c;2}),z\bigr)\bigr) \\
\qquad=
 \sum_c(-1)^{\deg x + 1}\mathfrak m\bigl(
\widehat{\mathfrak F}_{23}\bigl(\mathfrak m(x,{\bf y}_{c;1}),{\bf y}_{c;2}\bigr),z\bigr) \\
\phantom{\qquad= }{}+ \mathfrak m\bigl(\widehat{\mathfrak F}_{23}\bigl(x,\hat d {\bf y}\bigr),z\bigr)
+ (-1)^{\deg'{\bf y}_{c;1}}\sum_c\mathfrak m\bigl(\widehat{\mathfrak F}_{23}(x,{\bf y}_{c;1}),
\mathfrak m\bigl(\widehat{\mathfrak F}_{23}({\bf y}_{c;2}),z\bigr)\bigr).
\end{gather*}
By $A_{\infty}$ relation, this coincides with
\smash{$
(\mathfrak n\circ \mathscr I_{21})(x,{\bf y},z)
=
(-1)^{\deg x}\mathfrak m\bigl(\mathfrak m\bigl(\widehat{\mathfrak F}_{23}(x,{\bf y}),z\bigr)\bigr)$}.

\begin{lem}\label{lem1020}
$\mathscr I_{12;0,0}$ becomes a $(0,0)$ part of a
filtered bi-module homomorphism.
\end{lem}
\begin{proof}
We first define
\[
\mathscr I_{12;k_1,k_3} \colon\
B_{k_1}\mathscr C_1[1](c_1,c'_1) \,\widehat{\otimes}\, D^1(c'_1,c'_3)
\,\widehat{\otimes}\,B_{k_3}\mathscr C_3[1](c'_3,c'_3)
\to D^2(c_1,c_3)
\]
as follows.
If $k_3 \ne 0$, then $\mathscr I_{12;k_1,k_3}=0$. If $k_3 =0$, we put
$\smash{
\mathscr I_{12;k_1,0} ({\bf x}, z) =
{\bf e}_{(\mathscr F_{12})_{\rm ob}(c_1)} \otimes
\widehat{\mathscr F}_{12}({\bf x})}\allowbreak \otimes z$.
We will prove that they define an $A_{\infty}$ bi-module homomorphism.
Let
\begin{gather*}
\widehat{\mathscr I}_{12} \colon\
\bigoplus_{c'_1,c'_3}B\mathscr C_1[1](c_1,c'_1) \,\widehat{\otimes}\, D^1(c'_1,c_3)
\,\widehat{\otimes}\,B\mathscr C_3[1](c'_3,c_3) \\
\hphantom{\widehat{\mathscr I}_{12} \colon} \, \to
\bigoplus_{c'_1,c'_3}B\mathscr C_1[1](c_1,c'_1) \,\widehat{\otimes}\, D^2(c'_1,c_3)
\,\widehat{\otimes}\,B\mathscr C_3[1](c'_3,c_3)
\end{gather*}
be the formal bi-comodule homomorphism
induced by $\mathscr I_{12;k_1,0}$, $k_1 = 0,1,2,\dots$.
Let $\hat d$ be the
boundary operator on
$\bigoplus_{c'_1,c'_3}B\mathscr C_1[1](c_1,c'_1) \,\widehat{\otimes}\, D^i(c'_1,c'_3)
\,\widehat{\otimes}\,B\mathscr C_3[1](c'_3,c_3)$
induced by the bi-module structure
and
\[
\mathfrak n \colon\ B\mathscr C_1[1](c_1,c'_1) \,\widehat{\otimes}\, D^i(c'_1,c'_3)
\,\widehat{\otimes}\,B\mathscr C_3[1](c'_3,c_3)
\to D^i(c'_1,c_3),
\]
which is the structure operation of the bi-module structure as in \eqref{form108}.
Let ${\bf x} \in B\mathscr C_1[1](c_1,c'_1)$,
$z \in \mathscr D^1(c'_1,c'_3)$,
${\bf w} \in B_{k_3}\mathscr C_3[1](c'_3,c_3)$.
We calculate
\begin{gather*}
\bigl(\mathfrak n\circ \widehat{\mathscr I}_{12}\bigr)({\bf x},z,{\bf w}) \\
\qquad=
\sum_{c} \mathfrak n \bigl(\widehat{\mathscr F}_{12}({\bf x}_{c;1}) \otimes
\bigl(
{\bf e} \otimes \widehat{\mathscr F}_{12}({\bf x}_{c;2}) \otimes z
 \bigr)
 \otimes {\bf w}\bigr) \\
 \qquad= \begin{cases}
\sum_{c}(-1)^{\deg'{\bf x}_{c;1}}{\bf e} \otimes
\widehat{\mathscr F}_{12}({\bf x}_{c;1}) \otimes \mathfrak m\bigl(\bigl(\widehat
 {\mathscr F}_{23}\circ\widehat{\mathscr F}_{12}\bigr)({\bf x}_{c;2}),z,{\bf w}\bigr)
  &\text{if $k_3 \ne 0$}, \\
\sum_{c}(-1)^{\deg'{\bf x}_{c;1}}{\bf e} \otimes \widehat{\mathscr F}_{12}({\bf x}_{c;1}) \otimes \mathfrak m\bigl(\bigl(\widehat
 {\mathscr F}_{23}\circ\widehat{\mathscr F}_{12}\bigr)({\bf x}_{c;2}),z\bigr)&\\
 \qquad
 {}+
 {\bf e} \otimes \hat d\bigl(\widehat{\mathscr F}_{12}({\bf x})\bigr) \otimes z
&\text{if $k_3 = 0$}.
 \end{cases}
\end{gather*}
Note that in the case when $k_3 \ne 0$ the formula \eqref{form108}
implies that the summand in the second line
vanishes unless ${\bf x}_{c;1} = 1$.
In the case when $k_3 \ne 0$ and ${\bf x}_{c;1} = 1$, it becomes
the sum in the third line.

In the case when $k_3 = 0$ after a certain cancellation,
there remains another term, that is, the fifth line.
\begin{rem}
Note that $\deg' {\bf e} = -1$. However, as we remarked in Remark
\ref{reem2109} here the sign $\deg{\bf e} = 0$ is used.
\end{rem}
On the other hand,
\begin{align*}
\hat d({\bf x} \otimes z\otimes{\bf w})
={} &
\hat d({\bf x})\otimes z\otimes {\bf w}
+
(-1)^{\deg'{\bf x} + \deg z}{\bf x} \otimes z \otimes \hat d({\bf w}) \\
&+
\sum_{c_1,c_2}
(-1)^{\deg' {\bf x}_{c_1;1}}{\bf x}_{c_1;1}\otimes
\mathfrak m\bigl(\bigl(\widehat
 {\mathscr F}_{23}\circ\widehat{\mathscr F}_{12}\bigr)({\bf x}_{c_1;2}), z,
{\bf w}_{c_2;1}\bigr) \otimes {\bf w}_{c_2;2}.
\end{align*}
Therefore,
if $k_3 \ne 0$, we have
\[
\bigl({\mathscr I}_{12}\circ \hat d\bigr)({\bf x},z,{\bf w}) =
\sum_{c}(-1)^{\deg' {\bf x}_{c_1;1}}{\bf e}
\otimes \widehat{\mathscr F}_{12}({\bf x}_{c;1}) \otimes \mathfrak m\bigl(\bigl(\widehat
 {\mathscr F}_{23}\circ\widehat{\mathscr F}_{12}\bigr)({\bf x}_{c;2}),z,{\bf w}\bigr).
\]
If $k_3 = 0$, we have
\begin{gather*}
\bigl({\mathscr I}_{12}\circ \hat d\bigr)({\bf x},z) \\
\qquad=
\sum_{c}(-1)^{\deg' {\bf x}_{c_1;1}}{\bf e} \otimes \widehat{\mathscr F}_{12}({\bf x}_{c;1}) \otimes \mathfrak m\bigl(\bigl(\widehat
 {\mathscr F}_{23}\circ\widehat{\mathscr F}_{12}\bigr)({\bf x}_{c;2}),z\bigr)+
 {\bf e} \otimes \hat d({\mathscr F}_{12}({\bf x})) \otimes z.
\end{gather*}
Therefore, ${\mathscr I}_{12}$ is a filtered $A_{\infty}$ bi-module homomorphism.
\end{proof}

\begin{lem}\label{lem1021}
\quad
\begin{enumerate}\itemsep=0pt
\item[$(1)$]
The composition $\mathscr I_{21}\circ \mathscr I_{12,00}$ is equal to the
identity.
\item[$(2)$]
The composition $\mathscr I_{12,00}\circ \mathscr I_{21}$ is
chain homotopic to
the identity.
\end{enumerate}

\end{lem}
\begin{proof}
(1) follows by an easy and straightforward calculation.
We will prove (2).
Let $x \in \mathscr{C}_2((\mathscr F_{12})_{\rm ob}(c_1),c_2)$,
${\bf y} \in B\mathscr C_2[1](c_2,c'_2)$,
$z \in \mathscr{C}_3((\mathscr F_{23})_{\rm ob}(c'_2),c_3)$.
We observe
\[
(\mathscr I_{12,00}\circ \mathscr I_{21})(x,{\bf y},z)
= (-1)^{\deg x}{\bf e} \otimes \mathfrak m_*\bigl(\widehat{\mathfrak F}_{23}(x,{\bf y}),z\bigr).
\]
We define
$
\mathfrak H(x,{\bf y},z) := (-1)^{\deg' x}{\bf e} \otimes (x \otimes {\bf y}) \otimes z$.

Let $\mathfrak n$ be as in \eqref{defnnnnn}.
We calculate
\begin{align*}
(\mathfrak n \circ \mathfrak H)(x,{\bf y},z)
={}& (-1)^{\deg' x}\mathfrak n({\bf e} \otimes (x \otimes {\bf y}) \otimes z) \\
={}&x \otimes {\bf y}\otimes z + (-1)^{\deg' x}
\sum_c{\bf e} \otimes (\mathfrak m(x \otimes {\bf y}_{c;1})
\otimes {\bf y}_{c;2})\otimes z \\
&
+ {\bf e} \otimes \bigl(x \otimes \hat d{\bf y}\bigr) \otimes z \\
& + \sum_c (-1)^{\deg' {\bf y}_{c;1}}
{\bf e} \otimes (x \otimes {\bf y}_{c;1})
\otimes \mathfrak m\bigl(\widehat{\mathfrak F}_{23}({\bf y}_{c;2}),z\bigr) \\
&
- (-1)^{\deg x}{\bf e} \otimes \mathfrak m\bigl(\widehat{\mathfrak F}_{23}(x,{\bf y}),z\bigr).
\end{align*}
On the other hand,
\begin{align*}
(\mathfrak H\circ \mathfrak n)(x,{\bf y},z)
={}&
\sum_c (-1)^{\deg' {\bf y}_{c;1}}\mathfrak H(\mathfrak m(x,{\bf y}_{c;1})\otimes {\bf y}_{c;2} \otimes z)
\\
& + (-1)^{\deg x}\mathfrak H\bigl(x,\hat d {\bf y},z\bigr)
\\
& + \sum_c (-1)^{\deg x + \deg'{\bf y}_{c;1} }\mathfrak H(x,{\bf y}_{c;1},\mathfrak m({\bf y}_{c,2},z))\\
={}& \sum_c (-1)^{\deg' x+1}
{\bf e}\otimes (\mathfrak m(x \otimes {\bf y}_{c;1})
\otimes {\bf y}_{c;2})\otimes z \\
& - {\bf e} \otimes \bigl(x \otimes \hat d({\bf y})\bigr) \otimes z\\
& +\sum_c (-1)^{\deg' {\bf y}_{c;1}+1}
{\bf e} \otimes (x \otimes {\bf y}_{c;1})
\otimes \mathfrak m\bigl(\widehat{\mathfrak F}_{23}({\bf y}_{c;2}),z\bigr).
\end{align*}
Therefore,
\[
(\mathfrak n \circ \mathfrak H + \mathfrak H \circ \mathfrak n)(x,{\bf y},z)
=x \otimes {\bf y}\otimes z -
(-1)^{\deg x}{\bf e} \otimes \mathfrak m\bigr(\widehat{\mathfrak F}_{23}(x,{\bf y}),z\bigl),
\]
as required.
\end{proof}

Lemmas~\ref{lem1020} and \ref{lem1021}
together with Proposition~\ref{prop1015} imply that $\mathscr I_{12}$
is a homotopy equivalence. Proposition~\ref{lem9100} follows.
\end{proof}

\begin{proof}[Proof of Proposition~\ref{prop1019}]
We use the notation of Proposition~\ref{prop1019}.

Let $c_i \in \mathfrak{OB}(\mathscr C_i)$, $c,c' \in \mathfrak{OB}(\mathscr C)$.
By definition, we have
\begin{gather*}
\mathfrak{D}_{(1)}(c;c_1)
= \mathscr C(c,(\mathcal F)_{\rm ob}(c_1)), \qquad
\mathfrak{D}_{(2)}(c,c_2;c_3)
= \mathscr C_{(3)}((\mathcal G)_{\rm ob}(c,c_2);c_3), \\
\mathfrak{D}_{(3)}(c_1,c_2;c_3)
= \mathscr C_{(3)}((\mathcal G_{\rm ob}(\mathcal F_{\rm ob}(c_1),c_2);c_3).
\end{gather*}
We put
\begin{gather*}
 D^1(c_1,c_2;c_3) = \mathfrak{D}_{(3)}(c_1,c_2;c_3)
= \mathscr C_{(3)}((\mathcal G_{\rm ob}(\mathcal F_{\rm ob}(c_1),c_2);c_3),
\\
D^2(c_1,c_2;c_3) = \prod_{c,c',k}
\operatorname{Hom}\bigl(
B_k\mathscr C[1](c,c')\,\widehat\otimes\, \mathfrak{D}_{(1)}(c';c_1)
,
\mathfrak{D}_{(2)}(c,c_2;c_3)\bigr).
\end{gather*}
Note that $D^1$ is the object part of the tri-module associated to
$\mathcal G \circ \mathcal F$ and
$D^2$ is the object part of the tri-module
$(\mathfrak{Hom})_{\mathscr C}(\mathfrak D_{(1)},\mathfrak D_{(2)})$.

Note that the left module structure $\mathfrak n$ of
$\mathfrak{D}_{(1)}$
coincides with the $A_{\infty}$ operation $\mathfrak m$ of $\mathscr C$.
We define
\[
\mathscr I_{12;0,0;0} \colon\ D^1(c_1,c_2;c_3) \to D^2(c_1,c_2;c_3),
\qquad
\mathscr I_{21} \colon\ D^2(c_1,c_2;c_3) \to D^1(c_1,c_2;c_3)
\]
as follows.
Let $u \in D^1(c_1,c_2;c_3)$, ${\bf z} \in B\mathscr C(c,c')$,
$v \in \mathfrak{D}_{(1)}(c';c_1)$.
We put
\[
\mathscr I_{12;0,0;0}(u)({\bf z};v)
= (-1)^{(\deg u+1)(\deg v +\deg'{\bf z})}\mathfrak n({\bf z}\otimes v;u).
\]
Here $\mathfrak n$ is the left $\mathscr C$ module structure on $\mathfrak D_{(2)}$.
Note that the sign is different from Koszul sign and contains
the correction term $\deg v +\deg'{\bf z}$.
We will check that ${\mathscr I_{12;0,0;0}}$ is a chain map.
We have
\begin{equation}\label{10ten49}
\mathscr I_{12;0,0;0}(\mathfrak n(u))({\bf z};v)
= (-1)^{*_1}
\mathfrak n({\bf z};v;\mathfrak n(u)),
\end{equation}
where
\[
*_1 = (\deg u +1+1)(\deg v + \deg'{\bf z}) = \deg u(\deg' {\bf z} + \deg v).
\]
On the other hand,
\begin{align}
\mathfrak n({\mathscr I_{12;0,0;0}}(u))({\bf z};v)
={} &\sum_c (-1)^{\deg u\deg'{\bf z}_{c;1}}\mathfrak n({\bf z}_{c;1};{\mathscr I_{12;0,0;0}}(u)({\bf z}_{c;2};v))\nonumber \\
& +\sum_c (-1)^{\deg u+\deg'{\bf z}_{c;1}+1}
\mathscr I_{12;0,0;0}(u)({\bf z}_{c;1},\mathfrak n({\bf z}_{c;2},v))\nonumber
\\
&+ (-1)^{\deg u + 1} \mathscr I_{12;0,0;0}(u)\bigl(\hat d{\bf z},v\bigr)\nonumber\\
={} &\sum_c (-1)^{*_2} \sum_c \mathfrak n({\bf z}_{c;1},\mathfrak n({\bf z}_{c;2}\otimes v;u))\nonumber\\
&+ \sum_c (-1)^{*_3} \mathfrak n({\bf z}_{c;1}\otimes \mathfrak m({\bf z}_{c;2};v);u)
 + (-1)^{*_4} \mathfrak n\bigl(\hat d{\bf z}\otimes v;
u\bigr),\label{10ten50}
\end{align}
where
\begin{align*}
*_2={}& \deg u\deg'{\bf z}_{c;1} +
\deg u(\deg'{\bf z}_{c;2}+\deg v) + \deg v + \deg'{\bf z}_{c;2} \\
={}& \deg u(\deg'{\bf z}+\deg v)
+ \deg'{\bf z}_{c;2} + \deg v,
\\
*_3 ={}& \deg u+\deg'{\bf z}_{c;1} + 1
+ \deg u(\deg'{\bf z}_{c;1} + \deg'{\bf z}_{c;2} + \deg v+1)
\\
& + \deg'{\bf z}_{c;1}+ \deg'{\bf z}_{c;2} + \deg v + 1
\\
={}& \deg u(\deg v + \deg'{\bf z}) + \deg v + \deg'{\bf z}_{c;2}, \\
*_4 ={}& \deg u + 1 + \deg u(\deg'{\bf z} + \deg v + 1)
+ \deg'{\bf z} + \deg v + 1 \\
={}& \deg u(\deg'{\bf z} + \deg v) +\deg'{\bf z} + \deg v.
\end{align*}
Thus \eqref{10ten49} = \eqref{10ten50} is a consequence of the $A_{\infty}$
relation.
We remark that in the $A_{\infty}$ relation of $\mathfrak n$, the degree of
$v$ should be counted as $\deg'v$ (and not as $\deg v$),
since $v$ here appears as an element of the morphism complex of $A_{\infty}$
category (and not as an element of a bi-module).
We also remark that the operator $\mathfrak m$ appearing
in the second term of the right-hand side of \eqref{10ten50} coincides with $\mathfrak n$
in this case.

Let $\varphi = (\varphi_{c;c',k}) \in D^2(c_1,c_2;c_3)$.
We put
$
\mathscr I_{21}(\varphi) :=
\varphi_{c_1,c_2;c_3}
({\bf e}_{\mathcal F_{\rm ob}(c_1)})$.
Here ${\bf e}_{\mathcal F_{\rm ob}(c_1)} \in \mathscr C_2(\mathcal F_{\rm ob}(c_1),
\mathcal F_{\rm ob}(c_1))$ is the unity.
It is obvious that $\mathscr I_{21}$ is a chain map.
\begin{lem}\label{lem10202}
$\mathscr I_{12;0,0;0}$ becomes a $(0,0)$ part of
a filtered left $\mathscr C_{(1)}$ bi-module homomorphism.
\end{lem}
\begin{proof}
Let
${\bf x} \in B_{k_1}\mathscr{C}_{(1)}[1](c'_1,c_1)$,
$u \in D^1(c_1,c_2;c_3)$, ${\bf z} \in B\mathscr C(c,c')$, and
$v \in \mathfrak{D}_{(1)}(c';\mathcal F_{\rm ob}(c_1))
= \mathscr C(c',\mathcal F_{\rm ob}(c_1))$.
We put
\begin{equation}
\mathscr I_{12;k_1,0,0}({\bf x};u)({\bf z};v) :=
(-1)^{(\deg'{\bf x}+\deg u +1)(\deg'{\bf z}+\deg v)}\mathfrak n({\bf z} \otimes
v \otimes
{\bf x};u).\label{I12nocomplisign}
\end{equation}
We show
that this defines a left $\mathscr C_1$ module homomorphism.

We remark that the left $\mathscr C_1$ module structure on
$D^2$ is induced only from the right-$\mathscr C_1$ module structure
structure on $\mathfrak D_{(1)}$.
Namely,
\[
\mathfrak n({\bf x},\varphi)({\bf z};v)
= \sum_{c} (-1)^{\deg\varphi+\deg'{\bf z}_{c;1}+1
+ \deg'{\bf x}(\deg\varphi+\deg'{\bf z}+\deg v)}\varphi({\bf z}_{c;1},\mathfrak n({\bf z}_{c;2},v,{\bf x})).
\]
See \eqref{form102200}.
This is the case when ${\bf x} \notin B_0\mathscr C_1[1](c_1,c_2)$.
When ${\bf x} = 1$, we have
\begin{align*}
(\mathfrak n(\varphi))({\bf z};v)
={} &\sum_{c} (-1)^{1+\deg'{\bf z}_{c:1}+\deg\varphi}\varphi({\bf z}_{c;1},\mathfrak n({\bf z}_{c;2},v)) \\
&+ \sum_c (-1)^{\deg\varphi\deg'{\bf z}_{c:1}}\mathfrak n({\bf z}_{c;1},\varphi({\bf z}_{c;2},v))
+ (-1)^{\deg\varphi + 1}\varphi\bigl(\hat d{\bf z},v\bigr).
\end{align*}
See \eqref{form1023+}.
On the other hand, the left $\mathscr C_1$ module structure on
$D^1$ is induced from the left $\mathscr C$ module structure on
$\mathfrak D_{(1)}$ via $\mathcal F$.

We also remark that $\mathfrak m\bigl({\bf z},v,\widehat{\mathcal F}({\bf x})\bigr)
= (-1)^{\deg'{\bf x}}\mathfrak n({\bf z};v;{\bf x})$.

We now calculate
\begin{gather}
\mathscr I_{12}(\hat{\mathfrak n}({\bf x};u))({\bf z};v) \nonumber\\
\qquad= \sum_a (-1)^{\deg'{\bf x}_{a:1}}
\mathscr I_{12}({\bf x}_{a:1};{\mathfrak n}({\bf x}_{a:2};u))({\bf z};v)
+ \mathscr I_{12}\bigl(\hat d{\bf x};u\bigr)({\bf z};v)\nonumber \\
\qquad= \sum_a (-1)^{*_1}\mathfrak n\bigl({\bf z} \otimes
v \otimes \widehat{\mathcal F}({\bf x}_{a:1}) ; \mathfrak n\bigl(
\widehat{\mathcal F}({\bf x}_{a:2}); u\bigr)\bigr) + (-1)^{*_2}\mathfrak n\bigl({\bf z} \otimes
v \otimes \hat d
\widehat{\mathcal F}({\bf x})
; u\bigr),\label{Ainftyhomo1}
\end{gather}
with
\[
*_1 = \deg'{\bf x}_{a:1}+(\deg'{\bf x} + \deg u)(\deg'{\bf z} + \deg v),
\qquad
*_2 = (\deg'{\bf x} + \deg u)(\deg'{\bf z} + \deg v).
\]
Here $\Delta{\bf x} = \sum_a {\bf x}_{a:1} \otimes {\bf x}_{a:2}$.
On the other hand,
\begin{gather}
\bigl(\mathfrak n\bigl(\widehat{\mathscr I}_{12}({\bf x};u)\bigr)\bigr)({\bf z};v)\nonumber\\
\qquad= \sum_a \mathfrak n({\bf x}_{a:1};{\mathscr I}_{12}({\bf x}_{a:2};u))({\bf z};v)\nonumber\\
\qquad= \sum_{a,c} (-1)^{*_3}{\mathscr I}_{12}({\bf x}_{a:2};u)
({\bf z}_{c:1};\mathfrak n({\bf z}_{c:2};v
;{\bf x}_{a:1})) \nonumber\\
\phantom{\qquad=}{}+ \sum_{c} (-1)^{*_4}\mathfrak n({\bf z}_{c:1};\mathscr I_{12}({\bf x};u)({\bf z}_{c:2};v))
+ (-1)^{*_5}\mathscr I_{12}({\bf x};u)\bigl(\hat d{\bf z};v\bigr) \nonumber\\
\qquad=
 \sum_{a,c} (-1)^{*_6}\mathfrak n\bigl({\bf z}_{c:1}\otimes\mathfrak m({\bf z}_{c:2};v
;\widehat{\mathcal F}({\bf x}_{a:1}))\otimes {\bf x}_{a:2};u\bigr)\nonumber \\
\phantom{\qquad=}{}+
\sum_{c} (-1)^{*_7}\mathfrak n({\bf z}_{c:1};\mathfrak n({\bf z}_{c:2};v
;{\bf x};u))+ (-1)^{*_8}
\mathfrak n\bigl(\hat d{\bf z} \otimes
v \otimes {\bf x}
; u\bigr)\label{Ainftyhomo2}
\end{gather}
with
\begin{gather*}
*_3 = \deg'{\bf x}_{a:1}(\deg'{\bf x}_{a:2}+\deg u+\deg'{\bf z}+\deg v) + 1
+ \deg'{\bf x}_{a:2}+\deg u + \deg'{\bf z}_{c:1},
\\
*_4 = \deg'{\bf z}_{c:1}(\deg'{\bf x} + \deg u),\qquad
*_5 = \deg'{\bf x}+ \deg u + 1,
\end{gather*}
and\footnote{We remark that during the calculation of the
sign $*_6$ we use the fact that the operator $\mathfrak n$
appearing in the third line is related to the
operator $\mathfrak m$ in the sixth line by \eqref{sec10form2}.}
\begin{gather*}
*_6 = (\deg'{\bf x} + \deg u)(\deg'{\bf z} + \deg v)+\deg'{\bf z}_{c:2}+\deg v,
\\
*_7 = (\deg'{\bf x} + \deg u)(\deg'{\bf z} + \deg v)+\deg'{\bf z}_{c:2}+\deg v,
\\
*_8 = (\deg'{\bf x} + \deg u)(\deg'{\bf z} + \deg v)+\deg'{\bf z}+\deg v.
\end{gather*}
Therefore, the $A_{\infty}$ relation implies \eqref{Ainftyhomo1} = \eqref{Ainftyhomo2}.
We remark again that in the $A_{\infty}$ relation of $\mathfrak n$, the degree of
$v$ should be counted as $\deg'v$ (and not $\deg v$).

The proof that ${\mathscr I}_{12;0,0;0}$ extends to a tri-module homomorphism
is similar.
\end{proof}

\begin{lem}\label{lem10212}
\quad
\begin{enumerate}\itemsep=0pt
\item[$(1)$]
The composition $\mathscr I_{21}\circ \mathscr I_{12,00}$ is equal to the
identity.
\item[$(2)$]
The composition $\mathscr I_{12,00}\circ \mathscr I_{21}$ is
chain homotopic to
the identity.
\end{enumerate}
\end{lem}
\begin{proof}
(1) is easy to show.
We prove (2).
We remark that
\[
(\mathscr I_{12,00}\circ \mathscr I_{21})(\varphi)({\bf z};v)
=
(-1)^{(\deg \varphi+1) (\deg' {\bf z}+\deg v)}\mathfrak n({\bf z} \otimes v;\varphi({\bf e})),
\]
where the notations are as above.

We define
$\mathfrak H\colon D^2(c_1,c_2;c_3) \to D^2(c_1,c_2;c_3)$ by the next formula:
\[
\mathfrak H(\varphi)({\bf z};v) := (-1)^{\deg \varphi + \deg v+
\deg' {\bf z} +1 }\varphi({\bf z} \otimes v;{\bf e}).
\]
Let $\mathfrak n$ be the structure operations of $\mathfrak D_{(1)}$ and $\mathfrak D_{(2)}$,
which induce a boundary operator $\delta$ on $D^2(c_1,c_2; c_3)$.
See Lemma--Definition~\ref{lem1016}.
We calculate
\begin{align*}
(\delta (\mathfrak H(\varphi))({\bf z};v)
= {}&\sum_c (-1)^{\deg' {\bf z}_{c;1}(1+\deg \varphi)}
\mathfrak n({\bf z}_{c;1},\mathfrak H(\varphi)({\bf z}_{c;2};v))
\\
&
 + \sum_c (-1)^{\deg \varphi + \deg' {\bf z}_{c;1}}
\mathfrak H(\varphi)({\bf z}_{c;1};
\mathfrak n({\bf z}_{c;2};v))
 + (-1)^{\deg\varphi}\mathfrak H(\varphi)\bigl(\hat d{\bf z};v\bigr)
\\
={} &\sum_c (-1)^{*_1}
\mathfrak n({\bf z}_{c;1};\varphi({\bf z}_{c;2}\otimes v;{\bf e}))
\\
& + \sum_c (-1)^{*_2} \varphi({\bf z}_{c;1}\otimes \mathfrak n({\bf z}_{c;2},v);{\bf e})
+(-1)^{*3}\varphi\bigl(\hat d{\bf z}\otimes v ;{\bf e} \bigr),
\end{align*}
with
\begin{gather*}
*_1 = \deg' {\bf z}_{c;1}(1+\deg \varphi) + \deg\varphi
+ \deg'{\bf z}_{c;2}+ \deg v + 1, \\
*_2 = \deg \varphi + \deg'{\bf z}_{c;1}+\deg \varphi + \deg'{\bf z}+\deg v,\qquad
*_3 = \deg'{\bf z}+\deg v.
\end{gather*}
Here $\Delta({\bf z}) = \sum_c {\bf z}_{c;1} \otimes {\bf z}_{c;2}$.
On the other hand,
\begin{align*}
\bigl(\mathfrak H(\delta (\varphi))({\bf z};v)
={}&
 (-1)^{\deg\varphi + \deg'{\bf z}+\deg v} (\delta \varphi)({\bf z}\otimes v;{\bf e}) \\
={}&(-1)^{*_4}
\mathfrak n({\bf z} \otimes v ;\varphi({\bf e}))
+(-1)^{*_5}\sum_c
\mathfrak n({\bf z}_{c;1};\varphi(
{\bf z}_{c;2} \otimes v;{\bf e})) \\
& + (-1)^{*6}\varphi({\bf z} \otimes v) \\
 &+ \sum_c (-1)^{*_7}\varphi({\bf z}_{c;1};\mathfrak n({\bf z}_{c;2};v)
;{\bf e})
 + (-1)^{*_8}\varphi\bigl(\hat d{\bf z} \otimes v;{\bf e}\bigr)\bigr),
\end{align*}
with
\begin{gather*}
*_4 = \deg\varphi + \deg'{\bf z}+\deg v + \deg\varphi(\deg'{\bf z}+\deg v+1)= (\deg \varphi+1) (\deg' {\bf z}+\deg v),
\\
*_5 = \deg\varphi + \deg'{\bf z}+\deg v + \deg\varphi\deg'{\bf z}_{c;1}
= *_1 + 1, \qquad
*_6 = 1,\\
*_7 = \deg\varphi + \deg'{\bf z}+\deg v + \deg'{\bf z}_{c;1} + 1 + \deg\varphi
= *_2 + 1, \\
*_8 = \deg\varphi + \deg'{\bf z}+\deg v + 1 + \deg\varphi = *_3 + 1.
\end{gather*}
Therefore, $\mathfrak H$ is the chain homotopy we look for.
\end{proof}

Proposition~\ref{prop1019} follows from Lemmas~\ref{lem10202},
\ref{lem10212} and Proposition~\ref{prop1015}.
\end{proof}

\subsection[A note on 2-categories of $A_\infty$ categories]{A note on 2-categories of $\boldsymbol{A_{\infty}}$ categories}
\label{sec:catofAinfcat}

We remark that the diagram
\begin{equation}\label{Diagramctenten}
\begin{CD}
\displaystyle
\!\!\!\!\!\!\!\!\!\!\!\!
\!\!\!\!\!\!\!\!\!\!\!\!
\!\!\!\!\!\!\!\!\!\!\!\!
\!\!\!\!\!\!\!\!\!\!\!\!
\!\!\!\!\!\!\!\!
\mathcal{BIMOD}(\mathscr C_1;\mathscr C_2)
\atop\displaystyle
\times\mathcal{BIMOD}(\mathscr C_2;\mathscr C_3)
{\times \mathcal{BIMOD}(\mathscr C_3;\mathscr C_4)} @ >>>
\displaystyle\mathcal{BIMOD}(\mathscr C_1;\mathscr C_2)
\atop\displaystyle\times\mathcal{BIMOD}(\mathscr C_2;\mathscr C_4) \\
@ VVV @ VVV\\
\mathcal{BIMOD}(\mathscr C_1;\mathscr C_3)
\times\mathcal{BIMOD}(\mathscr C_3;\mathscr C_4) @ > >> \mathcal{BIMOD}(\mathscr C_1;\mathscr C_4)
\end{CD}
\end{equation}
{\it strictly} commutes. Here the arrows are the derived
tensor product functor $\mathfrak{ten}$.
The same holds if we replace $\mathcal{BIMOD}(*;*)$
by $\mathcal{BIMOD}(*;*)^{\rm op}$.
This implies that the diagram
\begin{equation}\label{Diagramcompss}
\begin{CD}
\displaystyle
\!\!\!\!\!\!\!\!\!\!\!\!
\!\!\!\!\!\!\!\!\!\!\!\!
\!\!\!\!\!\!\!\!\!\!\!\!
\!\!\!\!\!\!\!\!\!\!\!\!
\!\!\!\!
\mathcal{FUNC}(\mathscr C_1,\mathscr C_2)
\atop\displaystyle
\times\mathcal{FUNC}(\mathscr C_2,\mathscr C_3)
{\times \mathcal{FUNC}(\mathscr C_3,\mathscr C_4)} @ >>>
\displaystyle\mathcal{FUNC}(\mathscr C_1,\mathscr C_2)
\atop\displaystyle\times\mathcal{FUNC}(\mathscr C_2,\mathscr C_4) \\
@ VVV @ VVV\\
\mathcal{FUNC}(\mathscr C_1,\mathscr C_3)
\times\mathcal{FUNC}(\mathscr C_3,\mathscr C_4) @ > >> \mathcal{FUNC}(\mathscr C_1,\mathscr C_4)
\end{CD}
\end{equation}
commutes up to homotopy equivalence. Here the arrows are the composition functors
$\mathfrak{comp}$.
Since we take homotopy inverses to relative Yoneda functors
to obtain $\mathfrak{comp}$ from $\mathfrak{ten}$, the diagram \eqref{Diagramcompss}
does not commute strictly.
Using the version of Whitehead theorem with the notion `homotopic'
rather than `homotopy equivalent' (see Section~\ref{sec:homotopyfafunc}),
the `set of choices of homotopy inverse' seems to be `contractible'.
So we might be able to prove the associativity of $\mathfrak{comp}$ in a certain
$A_{\infty}$ sense. That might give a definition of
$A_{\infty}$ category of $A_{\infty}$ categories.
The author does not try to work it out here.
Instead, he points out the following.

Let $\mathscr A$ be a {\it set} whose elements are strict, unital and gapped
filtered $A_{\infty}$ categories.
We construct a DG-2-category $\mathfrak C(\mathscr A)$ whose object set is $\mathscr A$
and morphism category from $\mathscr C_1 \in \mathscr A$ to
$\mathscr C_2 \in \mathscr A$
is a full subcategory $\mathfrak C(\mathscr C_1,\mathscr C_2)$ of $\mathcal{BIMOD}(\mathscr C_1;\mathscr C_2)^{\rm op}$
such that the object set of~$\mathfrak C(\mathscr C_1,\mathscr C_2)$
consists of the bi-modules which are
homotopy equivalent to an element of the image of the relative
Yoneda functor $\mathfrak{RYon}_{\rm ob}\colon
\mathfrak{OB}(\mathcal{FUNC}(\mathscr C_1,\mathscr C_2))
\to \mathfrak{OB}(\mathcal{BIMOD}(\mathscr C_1,\mathscr C_2)^{\rm op})$.

The composition bi-functor of $\mathfrak C(\mathscr A)$ is $\mathfrak{ten}$.
By the strict commutativity of \eqref{Diagramctenten},
the composition bi-functors of $\mathfrak C(\mathscr A)$ are
{\it strictly} associative as DG-tri-functors.

Lemma--Definition~\ref{lemdef97} implies that
$\mathcal{FUNC}(\mathscr C_1,\mathscr C_2)$ is homotopy equivalent
to $\mathfrak C(\mathscr C_1,\mathscr C_2)$.
Moreover, this homotopy equivalence intertwines
composition bi-functors of $\mathfrak C(\mathscr A)$
with the composition bi-functors of $\mathcal{FUNC}(\mathscr C_1,\mathscr C_2)$
up to homotopy equivalence.

It is an opinion of the author that we can use $\mathfrak C(\mathscr A)$ as the
`2-category of $A_{\infty}$ categories' for most of the purposes.

\section{Associativity of compositions}
\label{sec:associ}

\subsection[Statement of the result of Section~\ref{sec:associ}]{Statement of the result of Section~\ref{sec:associ}}
\label{sateassoc}

In this section, we prove the associativity of the
composition functor defined in Theorem~\ref{comp2}.
\begin{situ}\label{situ11-1}
Let $(X_i,\omega_i,V_i)$ be a symplectic manifold $(X_i,\omega_i)$
equipped
with a background datum $V_i$.
Let $\mathbb L_{i (i+1)}$ for $i=1,2,3$ be a finite set of
$\pi_1^*(V_i \oplus TX_{i}) \oplus \pi_2^*(V_{i+1})$ relatively spin
Lagrangian submanifolds of $-X_i \times X_{i+1}$.
Let $\mathbb L_{i (i+2)}$, $i =1,2$, be a
finite set of
$\pi_1^*(V_i \oplus TX_{i+1}) \oplus \pi_2^*(V_{i+2})$ relatively spin
Lagrangian submanifolds of $-X_i\times X_{i+2}$.
Let $\mathbb L_{14}$ be a
finite set of~${\pi_1^*(V_1 \oplus TX_{1}) \oplus \pi_2^*(V_{4})}$ relatively spin
Lagrangian submanifolds of $-X_1\times X_{4}$.

We assume
\begin{enumerate}\itemsep=0pt
\item[(1)]
For $i=1,2,3$ and for any element $L_{i (i+1)}$ of
$\mathbb L_{i (i+1)}$ and $L_{(i+1) (i+2)}$ of $\mathbb L_{(i+1) (i+2)}$,
we assume that the fiber product $L_{i (i+1)}\times_{X_{i+1}}L_{(i+1) (i+2)}$
is transversal.
We also assume that its immersion to $X_{i (i+2)}$ has clean self-intersection.
\item[(2)]
For $i=1,2,3$,
the geometric composition of an element of
$\mathbb L_{i (i+1)}$ and of $\mathbb L_{(i+1) (i+2)}$ is contained in
$\mathbb L_{i (i+2)}$.
\item[(3)]
We assume the same condition as item (1) for the pairs
($\mathbb L_{12}$,$\mathbb L_{24}$),
($\mathbb L_{13}$,$\mathbb L_{34}$).
\item[(4)]
The geometric composition of an element of
$\mathbb L_{12}$ and of $\mathbb L_{24}$ is contained in $\mathbb L_{14}$.
The geometric composition of an element of
$\mathbb L_{13}$ and of $\mathbb L_{34}$ is contained in $\mathbb L_{14}$.
\end{enumerate}
For $1\le i<i'\le 4$, let $\mathfrak{Fuk}(-X_i \times X_{i'})$ be the filtered $A_{\infty}$ category defined in Theorem~\ref{prop333} whose objects
is an element of $\mathbb L_{i i'}$
and $\mathfrak{Fukst}(-X_i \times X_{i'})$ the strict category
associated to $\mathfrak{Fuk}(-X_i \times X_{i'})$.
\end{situ}

\begin{thm}\label{asscompmain}
Suppose we are in Situation {\rm\ref{situ11-1}}.
The next diagram commutes up to homotopy
equivalence:
\begin{equation}\label{diagram1111.1}
\begin{CD}
\displaystyle{\mathfrak{Fukst}(-X_1 \times X_2)
\times \mathfrak{Fukst}(-X_2 \times X_{3})
\atop\!\!\!\!\!\!\!\!\!\!\!\!\!\!\!\!\!
\!\!\!\!\!\!\!\!\!\!\!\!\!\!\!\!\!\!\!\!\!\!\!\!\times
\mathfrak{Fukst}(-X_3 \times X_{4})} @ >>>
\displaystyle{\mathfrak{Fukst}(-X_1 \times X_{3})
\atop
\times\mathfrak{Fukst}(-X_3 \times X_{4})} \\
@ VVV @ VVV\\
\displaystyle{\mathfrak{Fukst}(-X_1 \times X_{2})
\times \mathfrak{Fukst}(-X_2 \times X_{4})} @ > >> \mathfrak{Fukst}(-X_1 \times X_{4}),
\end{CD}
\end{equation}
where all the arrows are defined by the composition functor in
Theorem {\rm\ref{comp2}}. The homotopy equivalence is one of
unital, strict and gapped filtered $A_{\infty}$ tri-functors.

\end{thm}

The proof of Theorem~\ref{asscompmain} occupies the rest of this
section.
The proof is completed in Section~\ref{proofassoc}.
The argument of Section~\ref{proofassoc} is similar to
Section~\ref{sec:compfuncmain}.
The commutativity of~\eqref{diagram1111.1}
is homotopy equivalence between two
tri-functors.
Using relative Yoneda embedding, it is equivalent
to homotopy equivalence between certain
two quatro-modules.
For the proof, we will construct a quatro-module
homomorphism between them.
The quatro-module
homomorphism which we call Double-pants transformation
is defined by using a moduli space of objects
which we call Double-pants.
Double-pants in this section plays the role
Y-diagram played in Section~\ref{2-category formulation}.

\subsection{Opposite bi-modules and opposite drums}
\label{oppdrum}

For the proof of Theorem~\ref{asscompmain},
we need a certain digression.

\begin{defn}\label{definition113}
Let $\mathscr C_i$ be a filtered $A_{\infty}$
category for $i=1,2$ and $\mathfrak D = (D,\mathfrak n)$ a left-$\mathscr C_1$, right-$\mathscr C_2$
bi-module.
We define the {\it opposite bi-module} $\mathfrak D^{\rm op} = (D^{\rm op},\mathfrak n^{\rm op})$,\index{opposite bi-module}
which is a left-$\mathscr C^{\rm op}_2$, right-$\mathscr C_1^{\rm op}$ module\index[syindex]{Dop2@$\mathfrak D^{\rm op}$}
by the next formula. Let ${\bf x} \in B\mathscr C^{\rm op}_2(c_2,c'_2)$,
${\bf z} \in B\mathscr C^{\rm op}_1(c'_1,c_1)$, $y \in
\mathfrak D^{\rm op}(c'_2,c'_1): = D(c'_1;c'_2)$,
\begin{equation}\label{oppsitemoduledef}
\mathfrak n^{\rm op}({\bf x};y;{\bf z})
= (-1)^*\mathfrak n({\bf z}^{\rm op};y;{\bf x}^{\rm op}).
\end{equation}
Here the sign $*$ is by Kuszul rule $+1$. (See Definition~\ref{opcate}.)
We remark there are two convention of the degree of bi-module,
one is shifting the degree of an element of $D$ the other is not
shifting the degree of an element of $D$.
We put
$* = \varepsilon({\bf x}) + \deg'{\bf x}\deg'{\bf z} + \deg'y (\deg'{\bf x} + \deg'{\bf z})+1$
when we take the first convention and
$* = \varepsilon({\bf x}) +\deg'{\bf x}\deg'{\bf z} + \deg y (\deg'{\bf x} + \deg'{\bf z})+1$
when we take the second convention.
\end{defn}
It is easy to check \eqref{oppsitemoduledef} satisfies the $A_{\infty}$
relation.

\begin{exm}\label{sec1125}
In Section \ref{subsec:Yoneda}, we defined
the Yoneda functor $\mathfrak{Yon} \colon \mathscr C \to
\mathcal{FUNC}(\mathscr C^{\rm op},\mathcal{CH})$
and the opposite Yoneda functor $\mathfrak{OpYon} \colon \mathscr C^{\rm op} \to
\mathcal{FUNC}(\mathscr C,\mathcal{CH})$.
These two functors define
left-$\mathscr C$, right-$\mathscr C$ bi-module
structures on $\mathscr C(c,c')$.
It is easy to check that they are opposite bi-modules each other.

\end{exm}

We next define the opposite drum.
\begin{defn}\label{oppositedrm}
Suppose we are in the situation of Definition~\ref{def916}.
We consider the object~${(\Sigma;\vec z_{12},\vec z_{23},\vec z_{13};u_1,u_2,u_3;\gamma_{1},\gamma_2,\gamma_{3})}$
such that they enjoy the same properties as Definition~\ref{def916}
except the following:
\begin{enumerate}\itemsep=0pt
\item[(i)]
$u_1$ is a $J_{X_1}$-holomorphic map from $W_2$ and $u_2$ is a $J_{X_2}$-holomorphic map from $W_1$.
\item[(ii)]
We enumerate $\vec z_{12}$, $\vec z_{23}$ downward and $\vec z_{13}$ upward.
\end{enumerate}
We denote by
\smash{$\overset{\ \text{\tiny $\circ\circ$}}{\mathcal M}{}^{\rm op}_{\rm DR}(\vec a_{12},\vec a_{23},\vec a_{13};a_-,a_+;E)$}
the set of isomorphism classes of such objects.
We call an element of \smash{$\raisebox{-1pt}{$\overset{\ \text{\tiny $\circ\circ$}}{\mathcal M}{}^{\rm op}_{\rm DR}$}(\vec a_{12},\vec a_{23},\vec a_{13};a_-,a_+;E)$}
or its compactification an {\it opposite pseudo-holomorphic drum}.\index{opposite pseudo-holomorphic drum}

\end{defn}
\begin{prop}\label{propoppdrum}
The moduli space
\smash{$\overset{\ \text{\tiny $\circ\circ$}}{\mathcal M}{}^{\rm op}_{\rm DR}(\vec a_{12},\vec a_{23},\vec a_{13};a_-,a_+;E)$}\index[syindex]{M1opdra12@${\mathcal M^{\rm op}_{\rm DR}}(\vec a_{12},\vec a_{23},\vec a_{13};a_-,a_+;E)$}
has a compactification, abbreviated by ${\mathcal M^{\rm op}_{\rm DR}}(\vec a_{12},\vec a_{23},\vec a_{13};a_-,a_+;E)$, which is compact and Hausdorff. The compactifications have a system
Kuranishi structures and CF-perturbations.
They induce a left
$
\mathfrak{Fukst}(-X_1 \times X_2)
$,
$
\mathfrak{Fukst}(-X_2 \times X_3)
$ and right $\mathfrak{Fukst}(-X_1 \times X_3)$
tri-module.
\end{prop}

The proof is the same as the argument of Section~\ref{subsec:unobcomp}.
For example, Figures \ref{Figurebubbleseam1} and \ref{Figurebubbleseam2}
are replaced by the next Figures \ref{Figurebubbleopseam1} and
\ref{Figurebubbleopseam2}.

\begin{figure}[ht]
\centering
\includegraphics[scale=0.4]{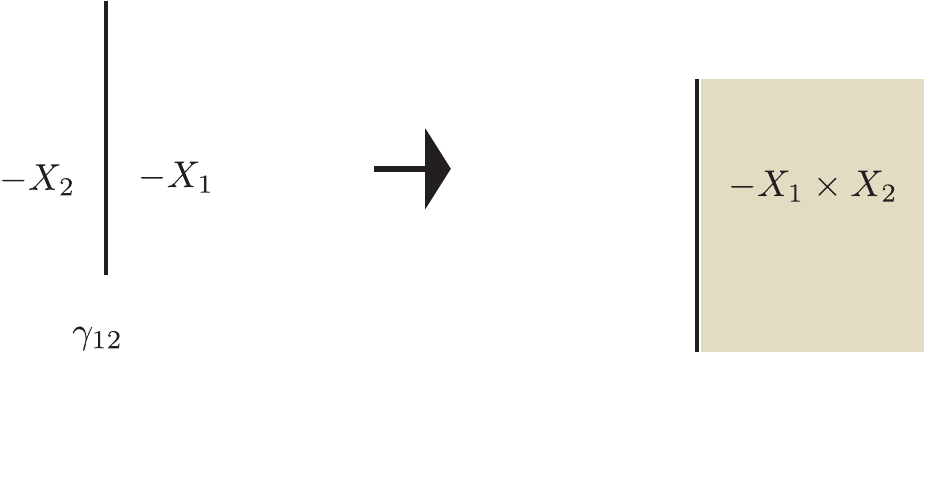}
\caption{Opposite version of Figure~\ref{Figurebubbleseam1}.}
\label{Figurebubbleopseam1}
\end{figure}

\begin{figure}[ht]
\centering
\includegraphics[scale=0.4]{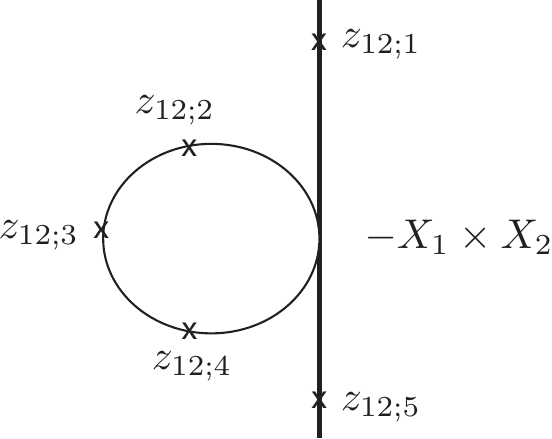}
\caption{Opposite version of Figure~\ref{Figurebubbleseam2}.}
\label{Figurebubbleopseam2}
\end{figure}

We denote the tri-module obtained in Proposition~\ref{propoppdrum}
by $\mathscr{CF}^{\rm op}(\mathbb L_{12},\mathbb L_{23};\mathbb L_{13})$.
We recall that in Section~\ref{sec:comp}
we defined the left-$\mathfrak{Fuk}(-X_1 \times X_3)$
right-$\mathfrak{Fuk}(-X_1 \times X_2)$, $\mathfrak{Fuk}(-X_2 \times X_3)$
tri-module $\mathscr{CF}(\mathbb L_{13};\mathbb L_{12},\mathbb L_{23})$.

\begin{lem}\label{lem118new}
$\mathscr{CF}^{\rm op}(\mathbb L_{12},\mathbb L_{23};\mathbb L_{13})$
is the opposite module to $\mathscr{CF}(\mathbb L_{13};\mathbb L_{12},\mathbb L_{23})$.

\end{lem}
\begin{proof}
We define
\[
\mathfrak I\colon\
\overset{\ \text{\tiny $\circ\circ$}}{\mathcal M}{}^{\rm op}_{\rm DR}(\vec a_{12},\vec a_{23},\vec a_{13};a_-,a_+;E)
\to
\overset{\ \text{\tiny $\circ\circ$}}{\mathcal M}_{\rm DR}(\vec a_{12},\vec a_{23},\vec a_{13};a_-,a_+;E)
\]
as follows.
We take
$
F\colon S^1 \times \R \to S^1 \times \R$
by
$
F(t,\tau) = (1-t,\tau)$.
This is an anti-holomorphic map.
In view of \eqref{eq8100},
this operation exchanges the domains $W_1$ and $W_2$.
Moreover, it revert the enumeration of the marked points on the seams.
Thus composing $F$ with the maps in the moduli space, we obtain a bijection $\mathfrak I$.
It is easy to see that the compactification
is preserved.
We~can take the Kuranishi structures and
CF-perturbations so that
they are preserved by $\mathfrak I$.
We~remark that the map $\mathfrak I$
reverse the enumeration of the marked points
on the seams.
This means that the operators obtained from
these two moduli spaces are
related by the operation taking the opposite category.
Therefore, in view of Example~\ref{sec1125},
the lemma holds up to sign.
In Section~\ref{oridrum}, we define
orientation of the moduli spaces of the drums and opposite drums
via appropriate doubling constructions. Therefore, Theorem~\ref{opthere} implies
that the sign becomes one of the opposite module.
\end{proof}

In Section~\ref{sec:comp}, we defined the functor
\begin{equation}\label{compsec11}
\mathfrak{comp} \colon\ \mathfrak{Fuk}(-X_1 \times X_2)
\times \mathfrak{Fuk}(-X_2 \times X_3)
\to \mathfrak{Fuk}(-X_1 \times X_3),
\end{equation}
so that
the composition
\[
\mathfrak{Yon} \circ \mathfrak{comp}
\colon\
\mathfrak{Fuk}(-X_1 \times X_2)
\times \mathfrak{Fuk}(-X_2 \times X_3)
\to
\mathcal{FUNC}(\mathfrak{Fuk}(-X_1 \times X_3)^{\rm op},\mathcal{CH})
\]
is the tri-module
$\mathscr{CF}(\mathbb L_{13};\mathbb L_{12},\mathbb L_{23})$.

On the other hand, the tri-module analogue of
Lemma--Definition~\ref{lemdef97} defines
\[
\mathfrak{RYon}\colon\
\mathcal{FUNC}(\mathscr C_1 \times \mathscr C_2,\mathscr C_3)
\to \mathcal{TRIMOD}(\mathscr C_1, \mathscr C_2;\mathscr C_3)^{\rm op}.
\]

\begin{cor}\label{Cor98new}
$\mathscr{CF}^{\rm op}(\mathbb L_{12},\mathbb L_{23};\mathbb L_{13})$
is homotopy equivalent to the tri-module
obtained by applying $\mathfrak{RYon}_{\rm ob}$
to {\rm\eqref{compsec11}}.

\end{cor}
This is a consequence of Lemma~\ref{lem118new} and
Example~\ref{oppositedrm}.

\subsection{Double pants}
\label{Surgery}

In this subsection, we work in Situation \ref{situ11-1}.
The proof of Theorem~\ref{asscompmain}
is based on a study of a moduli space of pseudo-holomorphic
maps from a space divided into pieces, which we explain now.
We consider the non-compact Riemann surface $W$ of
genus zero with 4 ends and its
division~${W = \bigcup_{i=1}^4 W_i}$ as in Figure~\ref{Figure11-1}
below.

\begin{figure}[ht]
\centering
\includegraphics[scale=0.4]{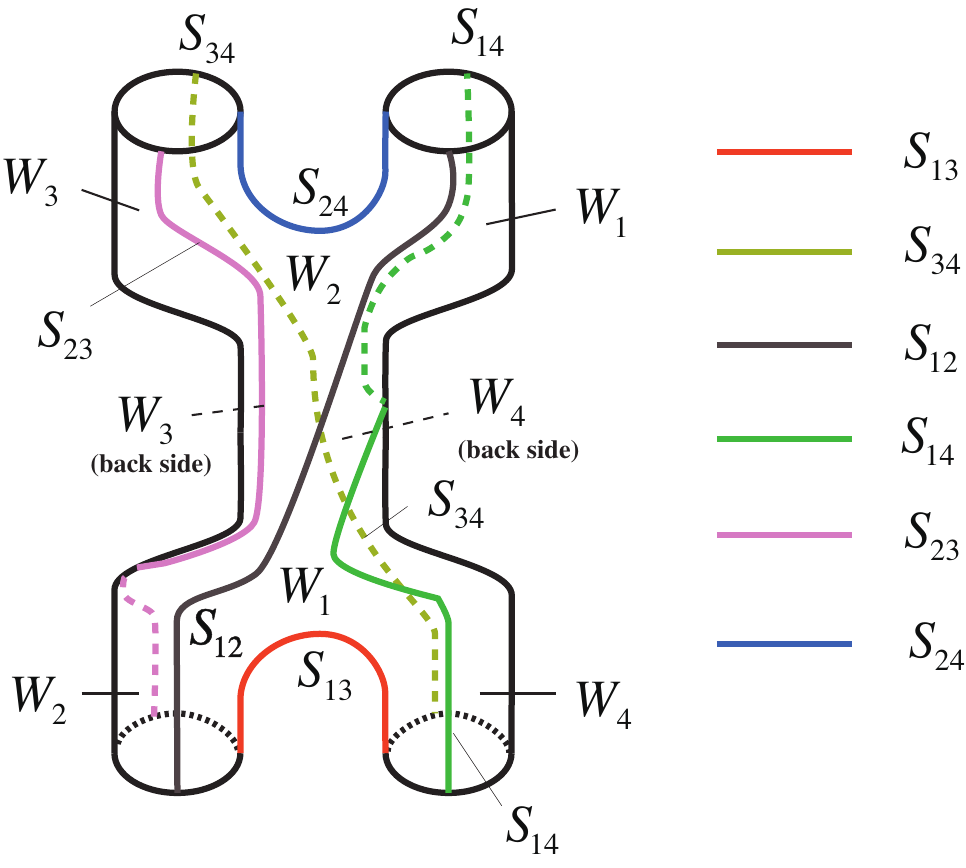}
\caption{Domain $W$.}
\label{Figure11-1}
\end{figure}

The domain $W$ is biholomorphic to $S^2$ minus 4 points.
It is divided into 4 domains $W_i$, $i=1,2,3,4.$
The intersection $S_{ii'} = W_i \cap W_{i'}$
is an arc for $ii' = 12,13,14,23,24,34$,
which we call a {\it seam}.
We call four points where three of the seams intersect
the {\it holes}.

\begin{figure}[ht]
\centering
\includegraphics[scale=0.4]{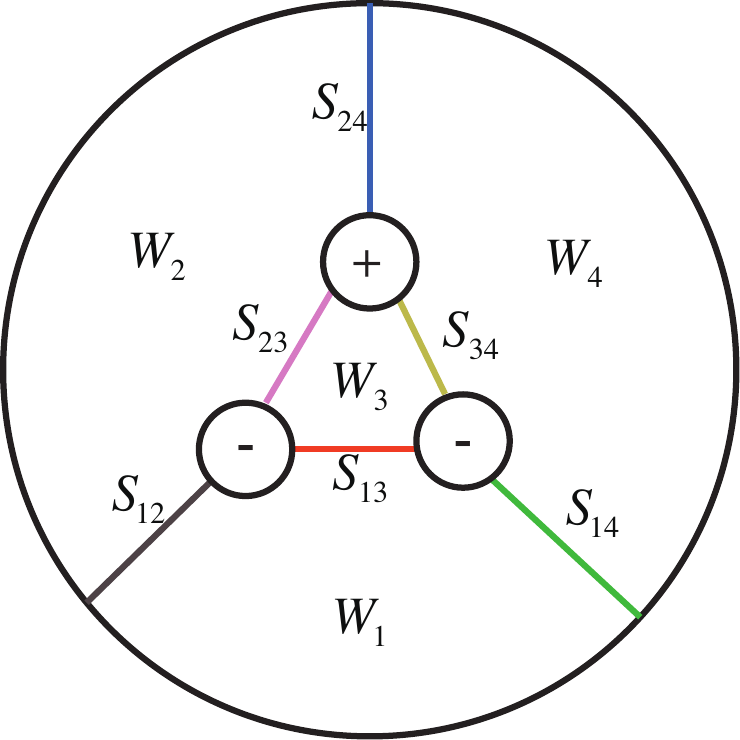}
\caption{Domain $W$ (alternative view).}
\label{Figure11-2}
\end{figure}
\begin{figure}[ht]
\centering
\includegraphics[scale=0.4]{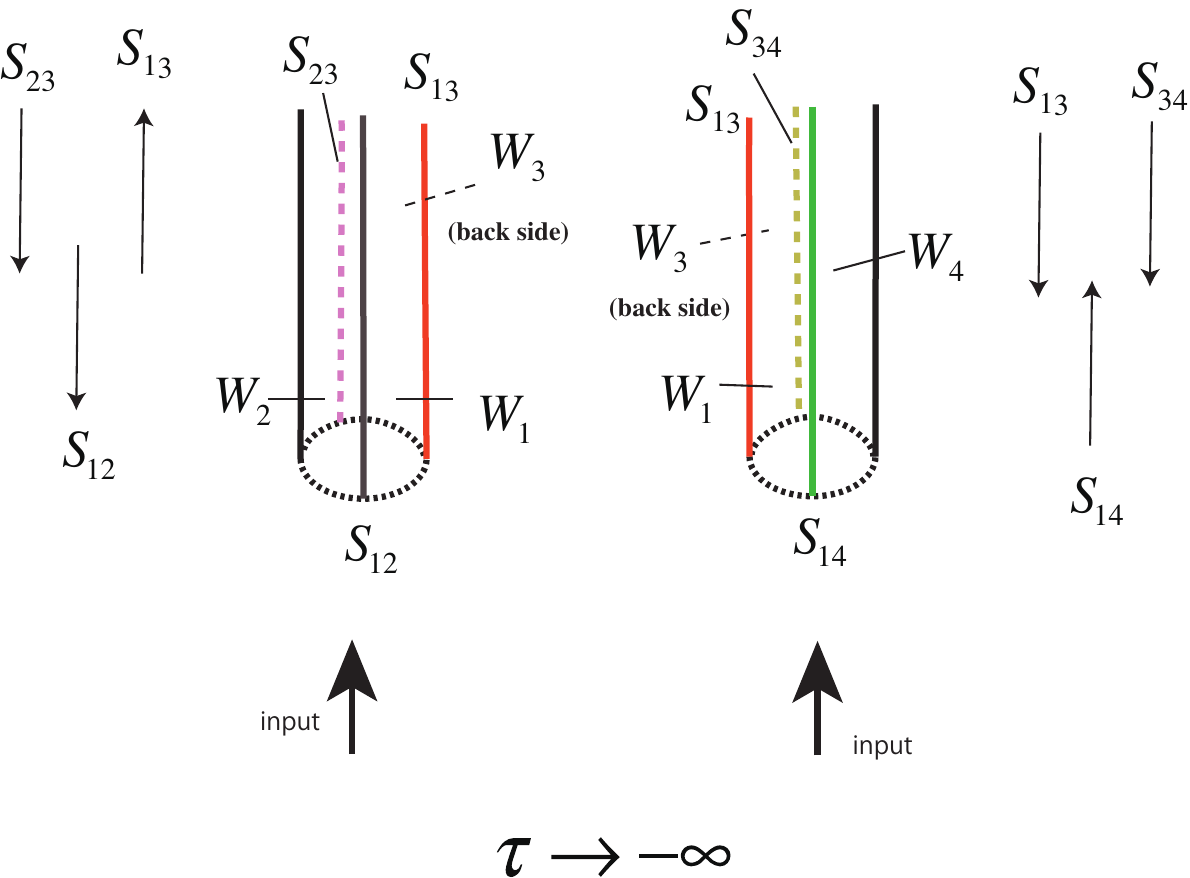}
\caption{Negative ends of the domain $W$.}
\label{Figure11-3}
\end{figure}
\begin{figure}[ht]
\centering
\includegraphics[scale=0.4]{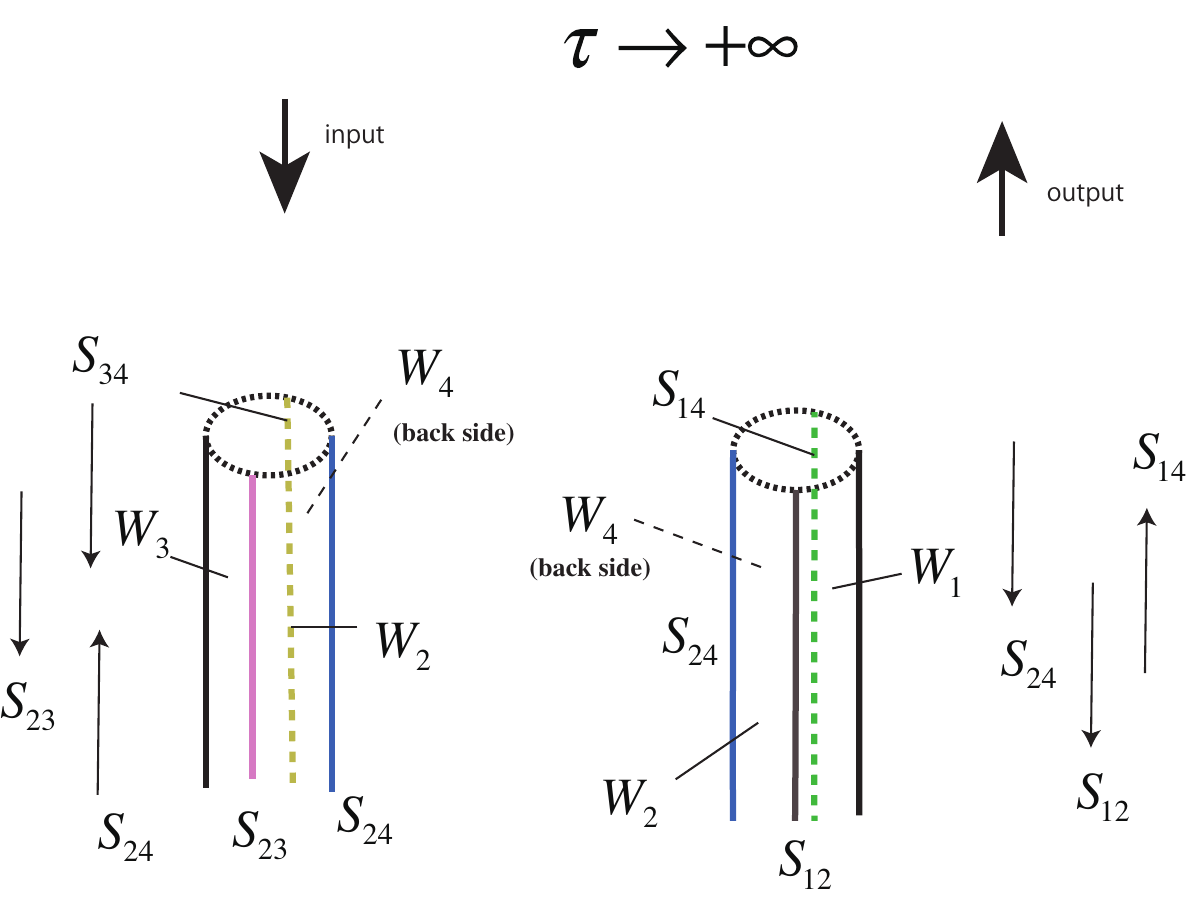}
\caption{Positive ends of the domain $W$.}
\label{Figure11-4}
\end{figure}

We consider the domain $W$ minus
holes and remove a relatively compact set from it.
Then the complement is biholomorphic to the disjoint union of the two copies
of $(-\infty,0] \times S^1$ and
the two copies of $[0,\infty) \times S^1$.
Each of those connected components are
divided into three pieces by seams.
In other words, each of them intersects with three
of $W_i$'s among four, as is shown in Figures
\ref{Figure11-3} and \ref{Figure11-4}.
We take and fix a bi-holomorphic map between
each of those ends and~${(-\infty,0] \times S^1}$
or $[0,\infty) \times S^1$.

We take the orientation of the seams $S_{ii'}$ as follows:
\begin{enumerate}\itemsep=0pt\setlength{\leftskip}{0.27cm}
\item[(seo1)]
For $ii' = 12,23,34$, we orient the seams so that it goes from the positive end
to the negative end.
\item[(seo2)]
For $ii'= 14$, we orient the seam so that it goes from the negative end
to the positive end.
\item[(seo3)]
For $i i' = 13$, we orient the seam so that it goes from the end
written in the left-hand side of Figure~\ref{Figure11-3}
to the end written in the
right-hand side of Figure~\ref{Figure11-3}.
For $i i' = 24$, we orient the seam so that it goes from the end
written in the right-hand side of
Figure~\ref{Figure11-4}
to the end written in the
left-hand side of Figure~\ref{Figure11-4}.
\end{enumerate}
See Figures \ref{Figure11-3} and \ref{Figure11-4}
for this orientation.

We observe that
Figure~\ref{Figure11-3} coincides with the negative end of the
opposite drum
used to define
$\mathscr{CF}^{\rm op}(\mathbb L_{12},\mathbb L_{23};\mathbb L_{13})$
and
$\mathscr{CF}^{\rm op}(\mathbb L_{13},\mathbb L_{34};\mathbb L_{14})$.
The right figure in Figure~\ref{Figure11-4}
coincides with the {\it positive} end of the
opposite drum
used to define $\mathscr{CF}^{\rm op}(\mathbb L_{12},\mathbb L_{24};\mathbb L_{14})$.
In the left side of the Figure~\ref{Figure11-4},
the positive end is actually an input.
So we rotate the figure by 180 degree so that
it becomes the negative end.
Then it coincides with the negative end of the
{\it drum}\footnote{We emphasise that this is not the opposite drum.}
used to define
$\mathscr{CF}(\mathbb L_{24};\mathbb L_{23},\mathbb L_{34})$.

We decompose the fiber product to the connected components as
\begin{equation}\label{eq114new}
\tilde L_{ii'} \times_{X_i\times X_{i'}} \tilde L_{ii'}
= \bigcup_{a \in \mathcal A_{L_{ii'}}} L_{ii'}(a).
\end{equation}
Situation \ref{situ11-1}\,(2) implies that the
fiber product in the left-hand side is
clean. Note that one of the components of \eqref{eq114new}
is the diagonal component.

For $i,i',i'' \in \{1,2,3,4\}$ with $i < i' < i''$, we decompose
\[
(L_{i i'} \times L_{i' i''} \times L_{i i''}) \times_{(X_i\times X_{i'} \times X_{i''})^2} \Delta
=
\bigcup_{a \in \mathcal A_{i i' i''}} R_{i i' i''}(a),
\]
where $\Delta$ is the diagonal in $(X_i\times X_{i i'} \times X_{i i' i''})^2$.
This is the decomposition to the connected components.
Situation \ref{situ11-1}\,(4) implies that the
fiber product in the left-hand side is
clean.

Let
\begin{gather*}
\vec a_{ii'} = (a_{ii',1},\dots,a_{ii',k_{ii'}}) \in (\mathcal A_{L_{ii'}})^{k_{ii'}} ,\qquad\!
\vec a_{i} = (a_{i,1},\dots,a_{i,k_{i}}) \in (\mathcal A_{L_{i}})^{k_{i}},\qquad\!
a_{ii'i''} \in \mathcal A_{i i' i''}.
\end{gather*}

In the next definition, we define
\smash{$\overset{\ \text{\tiny $\circ\circ$}}{\mathcal M}_{\rm DP}((\vec a_{ii'})_{ii'};(a_{ii'i''})_{ii'i''};E)$}.

\begin{defn}\label{def1116}
We consider
\smash{$(\Sigma;(\vec z_{ii'})_{1\le i<i' \le 4};(u_i;i=1,2,3,4);(\gamma_{ii'})_{1\le i<i' \le 4})$}
with the following properties.
\begin{enumerate}\itemsep=0pt
\item[(1)]
The bordered nodal curve $\Sigma$ is a union of $W$ and
trees of sphere components attached to~$W$.
The roots of the trees of sphere components are not on $\bigcup_{i,i'}S_{ii'}$.
\item[(2)]
For $i=1,2,3,4$, we denote by $\Sigma_{i}$ the union of $W_{i}$ together with
the trees of sphere components rooted on $W_{i}$.
The map $u_i \colon \Sigma_i \to X_i$ is $J_{X_i}$ holomorphic for $i=1,2,3,4$.
\item[(3)]
$\vec z_{ii'} = (z_{ii',1},\dots,z_{ii',k_{ii'}})$ and
$z_{ii',j} \in S_{ii'}$.
We require $z_{ii',j} < z_{ii',j'}$ for $j<j'$, where we identify
$S_{ii'} \cong \R$ by using the orientation defined in (seo1), (seo2), (seo3).
We put $\vert \vec z_{ii'}\vert = \{z_{ii',1},\dots,z_{ii',k_{ii'}}\}$.
\item[(4)]
The map
$\gamma_{ii'} \colon S_{ii'} \setminus \vert\vec z_{ii'}\vert \to \tilde L_{ii'}$
is smooth and satisfies
$
i_{L_{ii'}} (\gamma_{ii'}(z)) = (u_i(z),u_{i'}(z))$.
\item[(5)]
At $\vec z_{ii'}$, the map $\gamma_{ii'}$ satisfies the switching condition
\begin{equation}\label{form10914}
\bigl(\lim_{z \in S_{ii'} \uparrow z_{ii',j}}\gamma_{ii'}(z),\lim_{z \in S_{ii'} \downarrow z_{ii',j}}\gamma_{ii'}(z)
\bigr)
\in L_{ii'}(a_{ii',j}).
\end{equation}
Here we identify $S_{ii'} \cong \R$ and then $ \uparrow $, $ \downarrow$ have obvious meaning (see Definition~\ref{def3737}\,(5)) by using the orientation of $S_{ii'}$.
\item[(6)]
At the negative end of $W$, the following asymptotic boundary condition is satisfied:
\begin{gather}
\lim_{\tau \to -\infty} ( \gamma_{12}(\tau),\gamma_{23}(\tau), \gamma_{13}(-\tau)
 )
\in R_{123}(a_{123}),\nonumber
\\
\lim_{\tau \to -\infty} ( \gamma_{13}(\tau),\gamma_{34}(\tau), \gamma_{14}(-\tau)
 )
\in R_{134}(a_{134}).\label{form109152}
\end{gather}
\item[(7)]
At the positive end of $W$, the following asymptotic boundary condition is satisfied:
\begin{gather}
\lim_{\tau \to + \infty} ( \gamma_{23}(\tau),\gamma_{34}(\tau), \gamma_{14}(-\tau)
 )
\in R_{234}(a_{134}),\nonumber\\
\lim_{\tau \to + \infty} ( \gamma_{12}(\tau),\gamma_{24}(\tau), \gamma_{24}(-\tau)
 )
\in R_{124}(a_{124}).\label{form9152332}
\end{gather}
\item[(8)]
The stability condition, which is defined in the same way as Definition~\ref{defn1015}\,(2), is satisfied.
\item[(9)]
$
\sum_{i=1}^4\int_{\Sigma_i}u_i^*\omega_i = E$.
\end{enumerate}
In the same way as Definition~\ref{defn1015}\,(3), we define an equivalence relation $\sim$ among the objects~${(\Sigma;(\vec z_{ii'})_{1\le i<i' \le 4};(u_i)_{i=1,2,3,4};(\gamma_{ii'})
_{1\le i<i' \le 4})}$
satisfying (1)--(9). We denote
the set of all the equivalence classes of this equivalence relation by
\smash{$\overset{\ \text{\tiny $\circ\circ$}}{\mathcal M}_{\rm DP}((\vec a_{ii'})_{ii'};(a_{ii'i''})_{ii'i''};E)$}
\index[syindex]{M1aiiiiaiiiiE@$\overset{\ \text{\tiny $\circ\circ$}}{\mathcal M}((\vec a_{ii'})_{ii'};(a_{ii'i''})_{ii'i''};E)$}.
We call its element a pseudo-holomorphic double pants.\index{pseudo-holomorphic double pants}
\end{defn}
We define evaluation maps
\begin{equation}\label{form1111}
{\rm ev}_{ii',j} \colon\ \overset{\ \text{\tiny $\circ\circ$}}{\mathcal M}_{\rm DP}((\vec a_{ii'})_{ii'};(a_{ii'i''})_{ii'i''};E) \to L_{ii'}(a_{ii',j})
\end{equation}
by using \eqref{form10914}.
We define evaluation maps
\begin{equation}\label{form1112}
{\rm ev}_{ii'i''} \colon\ \overset{\ \text{\tiny $\circ\circ$}}{\mathcal M}_{\rm DP}((\vec a_{ii'})_{ii'};(a_{ii'i''})_{ii'i''};E) \to R_{ii'i''}(a_{ii'i''})
\end{equation}
by using one of \eqref{form109152}--\eqref{form9152332}.
\begin{prop}\label{prop10811}
We can define a topology on \smash{$\overset{\ \text{\tiny $\circ\circ$}}{\mathcal M}{}_{\rm DP}((\vec a_{ii'})_{ii'};(a_{ii'i''})_{ii'i''};E)$}
such that it has a compactification
${\mathcal M}_{\rm DP}((\vec a_{ii'})_{ii'};(a_{ii'i''})_{ii'i''};E)$, which is a compact metrizable space.
They have Kuranishi structures with corners which enjoy the following properties:\index[syindex]{M1DPaiiaiiii@${\mathcal M}_{\rm DP}((\vec a_{ii'})_{ii'};(a_{ii'i''})_{ii'i''};E)$}
\begin{enumerate}\itemsep=0pt
\item[$(1)$]
The normalized boundary of ${\mathcal M}_{\rm DP}((\vec a_{ii'})_{ii'};(a_{ii'i''})_{ii'i''};E)$ is a disjoint union of $2$ types of
fiber products, which we describe below.
\item[$(2)$]
The evaluation maps \eqref{form1111} and \eqref{form1112}
extend to strongly smooth maps with respect to this Kuranishi structure.
\eqref{form1112} is weakly submersive.
The extension is compatible with the description of the boundary
in item $(1)$.
\item[$(3)$]
The orientation local system of ${\mathcal M}_{\rm DP}((\vec a_{ii'})_{ii'};(a_{ii'i''})_{ii'i''};E)$ is isomorphic to the tensor product of
the pullbacks of $\Theta^-$ by
the evaluation maps \eqref{form1111} and \eqref{form1112}. For the component $R_{124}(a_{124})$, we take
$\Theta^+$ in place of $\Theta^-$.
\item[$(4)$]
The Kuranishi structure is compatible with the forgetful map of the marked points corresponding to the
diagonal components.
\end{enumerate}
\end{prop}
We describe the boundary components.

\begin{figure}[ht]\centering
\begin{tabular}{cc}
\begin{minipage}[t]{0.45\hsize}
\centering
\includegraphics[scale=0.4]{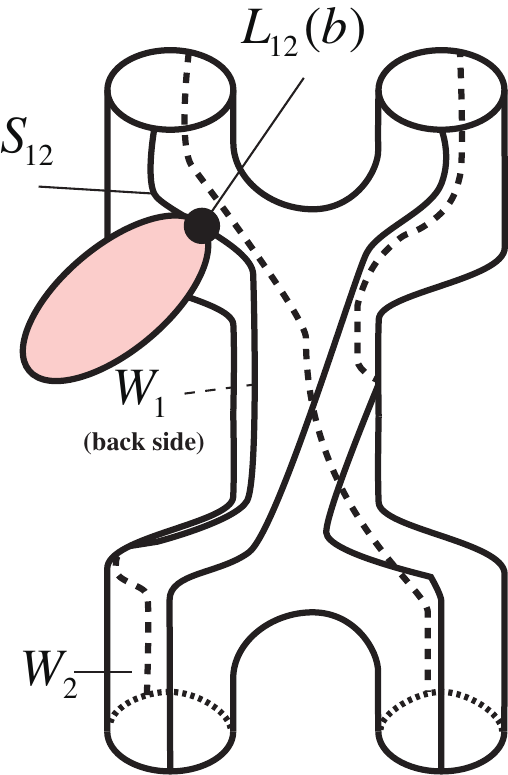}
\caption{Bubble of Type~I.}
\label{Figure11-5}
\end{minipage} &
\begin{minipage}[t]{0.45\hsize}
\centering
\includegraphics[scale=0.4]{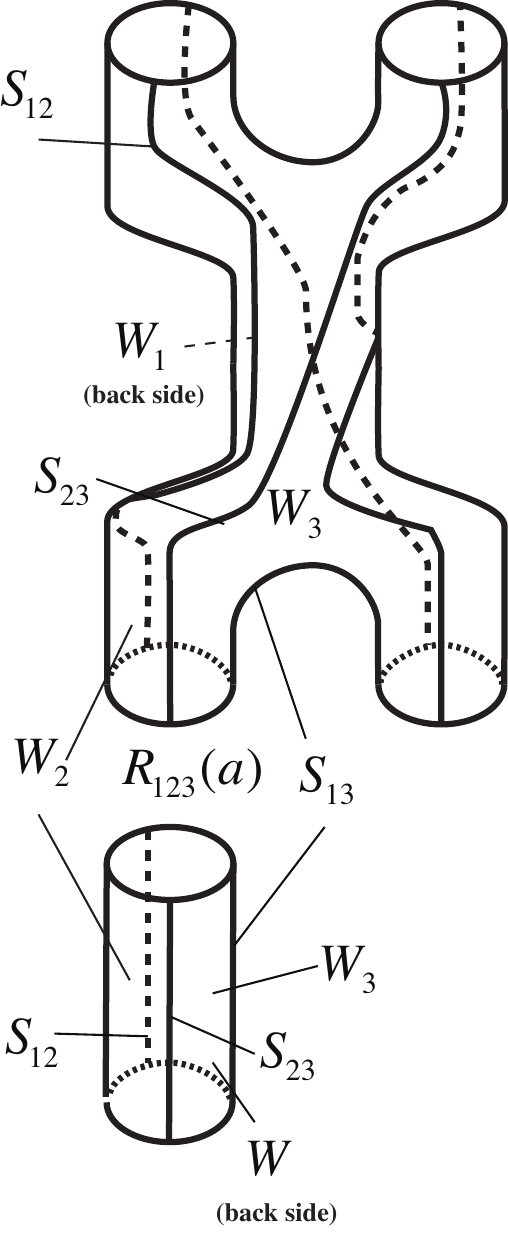}
\caption{Bubble of Type~II.}
\label{Figure11-6}
\end{minipage}
\end{tabular}
\end{figure}

(I)
The first type of boundary corresponds to the bubble at one of the Lagrangian boundary conditions $L_{ii'}$.
We describe the case of $L_{12}$.
Let $b \in \mathcal A_{L_{12}}$ and $i \le j$. We put
$\vec a_{12}^1= (a_{12,0},\dots,a_{12,i},b,a_{12,j+1},\dots,a_{12,k_{12}})$,
$\vec a_{12}^2= (b,a_{12,i+1},\dots,a_{12,j})$.
We put $\vec a'_{12} = \vec a_{12}^1$, $\vec a'_{ii'}= \vec a_{ii'}$ for $ii' \ne 12$.
This boundary corresponds to the fiber product
\[
{\mathcal M}_{\rm DP}((\vec a'_{ii'})_{i,i'};(a_{ii'i''})_{i,i',i''};E_1)
\times_{L_{12}(b)}
\mathcal M'\bigl(L_{12};\vec a^2_{12};E_2\bigr).
\]
Here $E_1 + E_2 = E$.
We remark that we use the compactification $\mathcal M'$ in the second factor,
which is a moduli space of pseudo-holomorphic disks
(see Remark~\ref{Remark524} and Section~\ref{sec:directcomp}).
See Figure~\ref{Figure11-5}.
The bubble at $L_{ii'}$ for $ii' \ne 12$ can be described in the same way.

(II)
The second type of boundary corresponds to the limit where the domain
will be divided into two parts at the ends.
There are 4 ends of our domain.
We first consider the case of the ends in the left-hand side of
Figure~\ref{Figure11-3}.

Let $j_{ii'} \in \{0,\dots,k_{ii'}\}$ for $ii' = 12, 23$ or $13$.
We put
$\vec a_{ii'}^1= (a_{ii',1},\dots,a_{ii',j_{ii'}})$,
$\vec a_{ii'}^2= (a_{ii',j_{ii'}+1},\dots,a_{ii',k_{ii'}})$
for $ii' = 12$ or $23$.
We also put
$\vec a_{ii'}^2= (a_{ii',1},\dots,a_{ii',j_{ii'}})$,
$\vec a_{ii'}^1= (a_{ii',j_{ii'}+1},\allowbreak \dots,a_{ii',k_{ii'}})$
for $ii' = 13$.
We then put $\vec a'_{ii'}= \vec a^2_{ii'}$ for $ii' = 12, 23$ or $13$
and $\vec a'_{ii'}= \vec a_{ii'}$ otherwise.

Let $a \in \mathcal A_{123}$.
We put $a'_{123} = a$ and $a'_{ii'i''} = a_{ii'i''}$ for $ii'i'' \ne 123$.

Now this boundary is described by the next fiber product
\[
{\mathcal M}_{\rm DP}\bigl(\bigl(\vec a_{ii'}^1\bigr)_{ii'};(a'_{ii'i''})_{ii'i''};E_1\bigr)
\times_{R_{123}(a)}
{\mathcal M}_{\rm DR}^{\rm op}((\vec a'_{ii'})_{ii'=12,23,13};a_{123},a;E_2),
\]
where $E_1 + E_2 = E$ and $a \in \mathcal A_{L_{12}}$.
See Figure~\ref{Figure11-6}.
Note that
${\mathcal M}_{\rm DR}^{\rm op}((\vec a'_{ii'})_{i,i'};a,a_{123};E_2)$
is the moduli space of opposite pseudo-holomorphic drums as in Definition
\ref{oppositedrm}.

In the case of the end in the right-hand side of
Figure~\ref{Figure11-3}, the end is described by
the fiber product{\samepage
\[
{\mathcal M}_{\rm DP}\bigl(\bigl(\vec a_{ii'}^1\bigr)_{ii'};(a'_{ii'i''})_{ii'i''};E_1\bigr)
\times_{R_{134}(a)}
{\mathcal M}_{\rm DR}^{\rm op}((\vec a'_{ii'})_{ii'=13,34,14};a_{134},a;E_2).
\]
Here $\vec a^1_{ii'}$, $\vec a'_{ii'i''}$ and $\vec a'_{ii'}$
are defined in a way similar to the first case.}

In the case of the end in the left-hand side of
Figure~\ref{Figure11-4}, the end is described by
the fiber product
\[
{\mathcal M}_{\rm DP}\bigl(\bigl(\vec a_{ii'}^1\bigr)_{ii'};(a'_{ii'i''})_{ii'i''};E_1\bigr)
\times_{R_{234}(a)}
{\mathcal M}_{\rm DR}((\vec a'_{ii'})_{ii'=23,34,24};a_{234},a;E_2).
\]
Here $\vec a^1_{ii'}$, $\vec a'_{ii'i''}$ and $\vec a'_{ii'}$
are defined in a way similar to the first case.
We remark that the second factor is the moduli space
of pseudo-holomorphic drums\footnote{Not opposite drum.} as in
Definition~\ref{def916}. The reason why pseudo-holomorphic drums appear
here is explained right after the orientation
of seams (seo1), (seo2), (seo3) are defined.

In the case of the end in the right-hand side of
Figure~\ref{Figure11-4}, the end is described by
the fiber product
\[
{\mathcal M}^{\rm op}_{\rm DR}((\vec a'_{ii'})_{ii'=12,24,14};a,a_{124};E_2)
\times_{R_{124}(a)}
{\mathcal M}_{\rm DP}\bigl(\bigl(\vec a_{ii'}^1\bigr)_{ii'};(a'_{ii'i''})_{ii'i''};E_1\bigr).
\]
Here the moduli space
of opposite pseudo-holomorphic drums appears.
Moreover, it appears as the first factor.
The reason is in the case of this end, $R_{124}(a)$ corresponds to the output of the second factor.

The proof of Proposition~\ref{prop10811} is similar to various other
propositions we discussed before in this and other papers
and so is omitted.
(See Section~\ref{oripants} for the proof of Proposition~\ref{prop10811}\,(3).)

\begin{prop}\label{prop109812}
For each $E_0$,
there exists a system of CF-perturbations $\widehat{\mathfrak S}$ on the space
${\mathcal M}_{\rm DP}((\vec a_{ii'})_{ii'};(a_{ii'i''})_{ii'i''};E)$
$($with respect to Kuranishi structures which are outer collarings of thickenings of those in Proposition
{\rm\ref{prop10811}}$)$
for $E < E_0$ such that the following holds:
\begin{enumerate}\itemsep=0pt
\item[$(1)$]
They are transversal to $0$.
\item[$(2)$]
The evaluation map
\eqref{form1112} is strongly submersive\footnote{See
\cite[Definition 9.2]{foootech2} and \cite{fooonewbook} for its definition.} with respect to this
CF-perturbation.
\item[$(3)$]
The CF-perturbations are compatible with the description of the
boundary. Namely, restriction of the CF-perturbation on the
boundary coincides with the fiber product CF-perturbation
in the sense of {\rm\cite[\emph{Lemma--Definition} 10.6]{foootech2}} and {\rm \cite{fooonewbook}}.
\item[$(4)$]
The CF-perturbations are compatible with the forgetful maps of the boundary
marked points corresponding to the diagonal component,
in the sense of {\rm\cite[\emph{Theorem}~5.1]{fooo091}}.
\end{enumerate}
\end{prop}

The proof is similar to the other similar statements we discussed already
and is now a routine. We omit it.

We now use Propositions \ref{prop10811} and \ref{prop109812}
to produce certain operations in a similar way as previous sections.
We need certain notations. For $1 \le i < i' \le 4$ and $1 \le j \le k_{ii'}$,
let $h_{ii',j} \in \Omega(L_{ii'}(a_{ii',j});\Theta^-)$.
We put
$
{\bf h}_{ii'} = (h_{ii',1},\dots,h_{ii',k_{ii'}})
\in B_{k_{ii'}}CF[1](\mathcal L_{ii'};\mathcal L'_{ii'})$.
For $ii'i'' = 123$ or $134$,
let
\[
h_{ii'i''} \in \Omega(R_{ii'i''}(a_{ii'i''});\Theta^-)
\subseteq CF^{\rm op}(\mathcal L_{ii'},\mathcal L_{i'i''};\mathcal L_{ii''}),
\]
and for $ii'i'' = 234$, let
\[
h_{234} \in \Omega(R_{234}(a_{234});\Theta^-)
\subseteq CF(\mathcal L_{24};\mathcal L_{23},\mathcal L_{34}).
\]
\begin{defn}\label{defnnew117}
We define
$
\mathscr{DPT}^{E,\varepsilon}(({\bf h}_{ii'})_{i,i'};h_{123},h_{134},h_{234})
\in \Omega(R_{124}(a_{124});\Theta^-)
$
by the next formula\index[syindex]{DPTscr@$\mathscr{DPT}^{E,\varepsilon}(({\bf h}_{ii'})_{i,i'};h_{123},h_{134},h_{234})$}
\begin{equation}\label{newe1117}
{\rm ev}_{124}!
\biggl(
\prod_{i<i'} {\rm ev}^*{\bf h}_{ii'} \wedge \prod_{ii'i'' = 123, 134,234}{\rm ev}_{ii'i''}^*h_{ii'i''}
; \widehat{\mathfrak S^{\varepsilon}}
\biggr).
\end{equation}
Here we use the moduli
space ${\mathcal M}((\vec a_{ii'})_{ii'};(a_{ii'i''})_{ii'i''};E)$
and its CF-perturbation $\widehat{\mathfrak S}$ to define~\eqref{newe1117}.
There is actually a sign in the right-hand side.
We will explain it in Section~\ref{oripants}.

We extend $\mathscr{DPT}^{E,\varepsilon}$ by $\Lambda_0$ linearly and
use it to define
\begin{gather}
\mathscr{DPT}^{<E_0,\varepsilon} \colon \
\prod_{i<i'} BCF[1](\mathcal L_{ii'};\mathcal L'_{ii'}) \otimes
CF^{\rm op}(\mathcal L'_{12},\mathcal L'_{23};\mathcal L_{13})\nonumber
\\
\hphantom{\mathscr{DPT}^{<E_0,\varepsilon} \colon}  \
{}\otimes
CF^{\rm op}(\mathcal L'_{13},\mathcal L'_{34};\mathcal L_{14})
\otimes
CF(\mathcal L'_{24};\mathcal L_{23},\mathcal L_{34})
\to
CF^{\rm op}(\mathcal L_{12},\mathcal L_{24};\mathcal L'_{14})\label{new1123}
\end{gather}
by the next formula
$
\mathscr{DPT}^{<E_0,\varepsilon} =
\sum_{E<E_0} T^E \mathscr{DPT}^{E,\varepsilon}$.
We call $\mathscr{DPT}^{<E_0,\varepsilon}$ the {\it double pants transformation}.
\index{double pants transformation}
\end{defn}
We next state the main property of the double pants transformation.
We need some notations.
Let
$
h_{ii'i''} \in CF^{\rm op}(\mathcal L_{ii'},\mathcal L_{i'i''};\mathcal L_{ii''})
$
for $ii'i'' = 123$ or $134$,
$
h_{234} \in CF(\mathcal L_{24};\mathcal L_{23};\mathcal L_{34})
$
and
$
{\bf h}_{ii'} = (h_{ii',1},\dots,h_{ii',k_{ii'}})
\in B_{k_{ii'}}CF[1](\mathcal L_{ii'},\mathcal L'_{ii'}).
$
We put
\smash{$
\Delta {\bf h}_{ii'}
=\sum_c {\bf h}_{ii'}^{c;1} \otimes {\bf h}_{ii'}^{c:2}$}.
For $\rho = jj'j'' = 123,\allowbreak 134$ or $234$
we define
${\bf h}^{\rho,c}_{ii'}$ as follows.
Let $W_{\rho}$ be one of the four ends corresponding to $\rho = jj'j''$:{\samepage
\begin{enumerate}\itemsep=0pt
\item[(1)]
${\bf h}^{\rho,c}_{ii'} = {\bf h}^{c,1}_{ii'}$
and ${\bf h}^{\rho,c;\prime}_{ii'} = {\bf h}^{c,2}_{ii'}$ if $S_{ii'}$ does not intersect with $W_{\rho}$.
\item[(2)]
${\bf h}^{\rho,c}_{ii'} = {\bf h}^{c;1}_{ii'}$
and ${\bf h}^{\rho,c;\prime}_{ii'} = {\bf h}^{c;2}_{ii'}$ if $S_{ii'} \cap W_{\rho}\ne \varnothing$ and
$W_{\rho}$ lies at the $-\infty$ side with respect to the orientation
of the seam $S_{ii'}$.
\item[(3)]
\smash{${\bf h}^{\rho,c}_{ii'} = {\bf h}^{c(\rho);2}_{ii'}$}
and ${\bf h}^{\rho,c;\prime}_{ii'} = {\bf h}^{c;1}_{ii'}$ if $S_{ii'} \cap W_{\rho}\ne \varnothing$ and
$W_{\rho}$ lies at the $+\infty$ side with respect to the orientation
of the seam $S_{ii'}$.
\end{enumerate}
In case $\rho = jj'j'' = 124$, we define ${\bf h}^{\rho,c}_{ii'}$ by exchanging the conditions~(2) and~(3).}

We also put
$
\hat d_{jj'}({\bf h}_{ii'})_{ii'} = ({\bf h}^*_{ii'})_{ii'}
$
where ${\bf h}^*_{ii'} = {\bf h}_{ii'}$ for $ii' \ne jj'$ and
${\bf h}^*_{jj'} = \hat d{\bf h}_{jj'}$.
We then put
\smash{$
\hat d ({\bf h}_{ii'})_{ii'} = \sum_{jj'} \hat d_{jj'}({\bf h}_{ii'})_{ii'}$}.

\begin{prop}\label{prop1118}
The double pants transformation $\mathscr{DPT}^{<E_0,\varepsilon}$
satisfies the next congruence modulo $T^{E_0}$:
\begin{gather}
\mathscr{DPT}^{<E_0,\varepsilon}
\bigl(\hat d(({\bf h}_{ii'})_{ii'});h_{123},h_{134},h_{234}\bigr)\nonumber \\
\qquad+
\sum_{c(12),c(23),c(13)} \mathscr{DPT}^{<E_0,\varepsilon}
\bigl(\bigl({\bf h}^{123,c(ii')}_{ii'}\bigr)_{ii'}
 ;\nonumber\\
 \phantom{\qquad\qquad-}{}
\mathfrak n^{\rm op}\bigl({\bf h}^{123,c(12);\prime}_{12},{\bf h}^{123,c(23);\prime}_{23},h_{123};{\bf h}^{123,c(13);\prime}_{13})
,h_{134},h_{234}\bigr)\nonumber
\\
\qquad+
\sum_{c(13),c(34),c(14)} \mathscr{DPT}^{<E_0,\varepsilon}\bigl(\bigl({\bf h}^{134,c(ii')}_{ii'}\bigr)_{ii'};\nonumber
\\
\phantom{\qquad\qquad-}{} h_{123},
\mathfrak n^{\rm op}\bigl({\bf h}^{134,c(13);\prime}_{13},{\bf h}^{134,c(34);\prime}_{34},h_{134};{\bf h}^{134,c(34);\prime}_{14}\bigr)
,h_{234}\bigr)\nonumber
\\
\qquad+
\sum_{c(23),c(34),c(24)} \mathscr{DPT}^{<E_0,\varepsilon}\bigl(\bigl({\bf h}^{234,c(ii')}_{ii'}\bigr)_{ii'}; h_{123},h_{134},\mathfrak n\bigl({\bf h}
^{234;\prime}_{24};h_{234};{\bf h}^{234;\prime}_{23},{\bf h}^{234;\prime}_{34}\bigr)\bigr)\nonumber
\\
\qquad-
\sum_{c(12),c(24),c(14)}
\mathfrak n^{\rm op}\bigl({\bf h}^{124,c(12)}_{12},{\bf h}^{124,c(24)}_{24};\nonumber
\\
\phantom{\qquad\qquad-}{}\mathscr{DPT}^{<E_0,\varepsilon}\bigl(\bigl({\bf h}^{124,c(ii')}_{ii'}\bigr)_{ii'};h_{123},h_{134},h_{234}\bigr);{\bf h}^{124,c(14);\prime}_{14}\bigr)
  \equiv 0 \mod T^{E_0}.\label{form1118}
\end{gather}
Here $\mathfrak n$ is the structure operation defined by the moduli space
of pseudo-holomorphic drums in Section {\rm\ref{sec:comp}}
and $\mathfrak n^{\rm op}$ is the structure operation defined by the moduli space
of opposite pseudo-holomorphic drums in Definition {\rm\ref{propoppdrum}}.
The signs $($which we omit from the above formula$)$ are by Koszul rule.
\end{prop}
\begin{proof}
Using Propositions \ref{prop10811}, \ref{prop109812},
 Stokes' formula (see \cite[Proposition 9.26]{foootech2} and \cite{fooonewbook}),
and the composition formula (see \cite[Theorem 10.20]{foootech2} and \cite{fooonewbook}),
the proof goes in the same way as the proof of Proposition~\ref{prop334}.
In fact, the first term of \eqref{form1118} corresponds to the
end of Type~I (see Figure~\ref{Figure11-5})
and the second-fifth terms of~\eqref{form1118} corresponds to the
end of Type~II (see Figure~\ref{Figure11-6}).

In fact, Type~I ends are described by the fiber products of
the moduli spaces of double pants diagrams and of pseudo-homomorphic
polygons. Type~II ends are described by the fiber products of
the moduli spaces of double pants diagrams and of (opposite) pseudo-homomorphic
drums.
\end{proof}

We can use Proposition~\ref{prop1118} to prove the next lemma
in the same way as we used
Propositions~\ref{prop330}, \ref{prop338} in Section~\ref{subsec:Ainfalgim}.
\begin{lem}\label{lem1199}
We can define $\mathscr{DPT}$
which is congruent to $\mathscr{DPT}^{<E_0,\varepsilon}$
modulo $T^{E_0}$ and which satisfies the same formula as
\eqref{form1118} except the congruence is replaced by
the equality.
\end{lem}
We call $\mathscr{DPT}$ in Lemma~\ref{lem1199} also
a double pants transformation.
We next twist the $\mathscr{DPT}$ by bounding cochains.
Let $b_{ii'}$ be bounding cochains of $\mathcal L_{ii'}$.
We define
\[
\mathfrak t^{\vec b} \colon\ \prod_{i<i'} BCF[1](\mathcal L_{ii'})
\to \prod_{i<i'} BCF[1](\mathcal L_{ii'})
\]
by the same formula as \eqref{defntttt}.
We then put
$
\mathscr{DPT}^{\vec b} = \mathscr{DPT} \circ \bigl(\mathfrak t^{\vec b} \otimes {\rm id}\bigr)$.
\begin{lem}\label{lem1110}
$\mathscr{DPT}^{\vec b}$ satisfies the same formula as \eqref{form1118}
except we twist $\hat d$ and $\mathfrak n$ by $\vec b$ and
the congruence is replaced by
the equality.
\end{lem}
The proof is easy and so is omitted.

\subsection{Proof of the associativity}
\label{proofassoc}

Now we use the double pants transformation to prove
Theorem~\ref{asscompmain}.
The proof is similar to the arguments of
Sections \ref{sec:comptibility} and
\ref{2-category formulation}.

We first prove the next proposition,
which is similar to Proposition~\ref{prop912}.

\begin{prop}\label{prop1212}
In Situation {\rm\ref{situ11-1}},
let $\mathcal L_{12} = (L_{12},\sigma_{12},b_{12})$
$($resp.\ $\mathcal L_{23} = (L_{23},\sigma_{23},b_{23})$,
$\mathcal L_{34} = (L_{34},\sigma_{34},b_{34})$
$)$ be an object of $\mathfrak{Fukst}(-X_1\times X_2)$
$($resp.\ $\mathfrak{Fukst}(-X_2\times X_3)$,
$\mathfrak{Fukst}(-X_3\times X_4))$.
We put
\[ \mathcal L_{13}= (L_{13},\sigma_{13},b_{13})
= \mathfrak{Comp}(\mathcal L_{12},\mathcal L_{23}), \qquad
\mathcal L^{(1)}_{14} = \bigl(L^{(1)}_{14},\sigma^{(1)}_{14},b^{(1)}_{14}\bigr)
= \mathfrak{Comp}(\mathcal L_{13},\mathcal L_{34}),
\]
and
\[ \mathcal L_{24}= (L_{24},\sigma_{24},b_{24})
= \mathfrak{Comp}(\mathcal L_{23},\mathcal L_{34}), \qquad
\mathcal L^{(2)}_{14}= \bigl(L^{(2)}_{14},\sigma^{(2)}_{14},b^{(2)}_{14}\bigr)
= \mathfrak{Comp}(\mathcal L_{12},\mathcal L_{24}).
\]
Then we have the following:
\begin{enumerate}\itemsep=0pt
\item[$(1)$]
\smash{$\bigl(L^{(1)}_{14},\sigma^{(1)}_{14}\bigr) = \bigl(L^{(2)}_{14},\sigma^{(2)}_{14}\bigr)$}. Here the equality is as submanifolds
equipped with relative spin structures.
\item[$(2)$]
\smash{$b^{(1)}_{14}$} is gauge equivalent to \smash{$b^{(2)}_{14}$} in the sense of {\rm\cite[\emph{Definition} 4.3.1]{fooobook}}.
\end{enumerate}

\end{prop}
\begin{proof}
(1) is proved in the same way as Proposition~\ref{prop912}\,(1),
which is proved in Section~\ref{oriYdiagarm}.
We prove (2) below.

We put
${\bf h}_{ii'} = e^{b_{ii'}}$ for $1 \le i<i' \le$ with $ii' \ne 14$
and ${\bf h}_{14} = e^{b^{(1)}_{14}}$.
Let $h_{ii'i''} = {\bf 1}_{ii'i''}$ for $ii'i'' = 123, 134,234$.
Here ${\bf 1}_{ii'i''}$ is the function $1$ on the diagonal component,
which is diffeomorphic to~$\tilde L_{ii''}$.
We define
\[
{\bf 1}^{(1)}_{124} = \mathscr{DPT}^{\vec b}(({\bf h}_{ii'})_{ii'};h_{123},h_{134},h_{234}).
\]
We consider the filtered $A_{\infty}$ tri-module~${\mathscr{CF}(\mathbb L_{14};\mathbb L_{12},\mathbb L_{24})}$ and twist
it by the bounding cochains $b_{12}$, $b_{24}$.
We then obtain a left filtered $A_{\infty}$
module $\mathscr{CF}(\mathbb L_{14};\mathbb L_{12},\mathbb L_{24})$ over
$\mathfrak{Fukst}(-X_{1} \times X_4)$. By Lemma~\ref{lem1110}, we have
\begin{equation}\label{form1120}
\mathfrak n\bigl(e^{b^{(1)}_{14}};{\bf 1}^{(1)}_{124}\bigr) = 0.
\end{equation}

Let \smash{${\bf 1}^{(2)}_{124}$} be the function $1 \in CF(L_{14};L_{12},L_{24})$
on the diagonal component, which is diffeomorphic to $\tilde L_{14}$.
By definition (see formulas \eqref{newform87} and \eqref{newfig63}), we have
\begin{equation}\label{form1121}
\mathfrak n(e^{b^{(2)}_{14}};{\bf 1}^{(2)}_{124}) = 0.
\end{equation}
By the definition of $\mathscr{DPT}^{\vec b}$ and \smash{${\bf 1}^{(1)}_{124}$},
we find that
\begin{equation}\label{form1122}
{\bf 1}^{(1)}_{124} \equiv {\bf 1}^{(2)}_{124} \mod \Lambda_+.
\end{equation}
Using \eqref{form1120}, \eqref{form1121}, \eqref{form1122},
we can apply (the left module analogue of) Lemma~\ref{lem912} to conclude
that \smash{$b^{(1)}_{14}$} is gauge equivalent to \smash{$b^{(2)}_{14}$}.
\end{proof}

\begin{proof}[Proof of Theorem~\ref{asscompmain}]
We first study the composition
\begin{gather}
\mathfrak{Fukst}(-X_1 \times X_2)
\times \mathfrak{Fukst}(-X_2 \times X_{3})
\times
\mathfrak{Fukst}(-X_3 \times X_{4})\nonumber \\
\qquad\to
\mathfrak{Fukst}(-X_1 \times X_{3})
\times\mathfrak{Fukst}(-X_3 \times X_{4})
\to \mathfrak{Fukst}(-X_1 \times X_{4}).\label{1123form}
\end{gather}
We apply the object part of the relative Yoneda functor to
\[
\mathfrak{Fukst}(-X_1 \times X_2)
\times \mathfrak{Fukst}(-X_2 \times X_{3})
\to
\mathfrak{Fukst}(-X_1 \times X_3).
\]
We then obtain the filtered $A_{\infty}$
tri-module $\mathscr{CF}^{\rm op}(\mathbb L_{12},\mathbb L_{23};\mathbb L_{13})$
by Corollary \ref{Cor98new}.

On the other hand, applying the object part of the relative Yoneda functor to the composition
\[
\mathfrak{Fukst}(-X_1 \times X_3)
\times \mathfrak{Fukst}(-X_3 \times X_{4})
\to
\mathfrak{Fukst}(-X_1 \times X_4),
\]
we obtain the filtered $A_{\infty}$
tri-module $\mathscr{CF}^{\rm op}(\mathbb L_{13},\mathbb L_{34};\mathbb L_{14})$
by Corollary \ref{Cor98new}.
Therefore, by Proposition~\ref{lem9100}, applying the
relative Yoneda functor to the composi\-tion \eqref{1123form} gives the derived tensor product
\[
D_1 = \mathfrak{ten}(\mathscr{CF}^{\rm op}(\mathbb L_{12},\mathbb L_{23};\mathbb L_{13}),
\mathscr{CF}^{\rm op}(\mathbb L_{13},\mathbb L_{34};\mathbb L_{14}))
\]
over $\mathfrak{Fukst}(-X_1 \times X_3)$.
The quatro-module structure on $D_1$ is defined in the same way as
one in the derived tensor product (see Lemma--Definition~\ref{defntensor}).
Here the quatro-module structure is
left-$\mathfrak{Fukst}(-X_1 \times X_2)$, $\mathfrak{Fukst}(-X_2 \times X_3)$, $\mathfrak{Fukst}(-X_3 \times X_4)$ and right $\mathfrak{Fukst}(-X_1 \times X_4)$ module structure.

We next consider the composition
\begin{gather}
\mathfrak{Fukst}(-X_1 \times X_2)
\times \mathfrak{Fukst}(-X_2 \times X_{3})
\times
\mathfrak{Fukst}(-X_3 \times X_{4}) \nonumber\\
\qquad\to
\mathfrak{Fukst}(-X_1 \times X_{2})
\times\mathfrak{Fukst}(-X_2 \times X_{4})
\to \mathfrak{Fukst}(-X_1 \times X_{4}).\label{1126form}
\end{gather}
By definition (see Proposition~\ref{prop810}), the composition functor
\[
\mathfrak{Fukst}(-X_2 \times X_3)
\times \mathfrak{Fukst}(-X_3 \times X_{4})
\to
\mathfrak{Fukst}(-X_2 \times X_4)
\]
composed with the Yoneda functor
\[
\mathfrak{Yon}\colon\ \mathfrak{Fukst}(-X_2 \times X_4)
\to \mathcal{FUNC}(\mathfrak{Fukst}(-X_2 \times X_4)^{\rm op};\mathcal{CH})
\]
gives a left-$\mathfrak{Fukst}(-X_2 \times X_4)$,
right-$\mathfrak{Fukst}(-X_2 \times X_3)$, $\mathfrak{Fukst}(-X_3 \times X_4)$
tri-module.
$
\mathscr{CF}(\mathbb L_{24};\allowbreak\mathbb L_{23},\mathbb L_{34}).
$
On the other hand, applying the relative Yoneda functor
to
\[
\mathfrak{Fukst}(-X_1 \times X_2)
\times \mathfrak{Fukst}(-X_2 \times X_{4})
\to
\mathfrak{Fukst}(-X_1 \times X_4)
\]
gives
the left-$\mathfrak{Fukst}(-X_1 \times X_2)$, $\mathfrak{Fukst}(-X_2 \times X_4)$
right-$\mathfrak{Fukst}(-X_1 \times X_4)$
tri-module $\mathscr{CF}^{\rm op}(L_{12},\allowbreak L_{24}; L_{14})$.

Therefore, we can apply Proposition~\ref{prop1019}
by putting
\begin{gather*}
\mathscr C_{(1)} = \mathfrak{Fukst}(-X_2 \times X_3)
\times \mathfrak{Fukst}(-X_3 \times X_4), \qquad
\mathscr C_{(2)} = \mathfrak{Fukst}(-X_1 \times X_2), \\
\mathscr C_{(3)} = \mathfrak{Fukst}(-X_1 \times X_4), \qquad
\mathscr C = \mathfrak{Fukst}(-X_2 \times X_4), \\
\mathfrak D_{(1)}= \mathscr{CF}^{\rm op}(\mathbb L_{12},\mathbb L_{24};\mathbb L_{14}), \qquad
\mathfrak D_{(2)}= \mathscr{CF}(\mathbb L_{24};\mathbb L_{23},\mathbb L_{34}),
\end{gather*}
to find that applying relative Yoneda functor to the map \eqref{1126form}
gives a
 left-$\mathfrak{Fukst}(-X_1 \times X_2)$, $\mathfrak{Fukst}(-X_2 \times X_3)$, $\mathfrak{Fukst}(-X_3 \times X_4)$ and right $\mathfrak{Fukst}(-X_1 \times X_4)$ quatro-module
\[
D_2 = \mathfrak{Hom}_{\mathfrak{Fukst}(-X_2 \times X_4)}(
\mathscr{CF}(\mathbb L_{24};\mathbb L_{23},\mathbb L_{34}),
\mathscr{CF}^{\rm op}(\mathbb L_{12},\mathbb L_{24};\mathbb L_{14}))
.
\]

Now a quatro-module homomorphism from $D_1$ to $D_2$ is a map
from{\samepage
\begin{gather}
BCF[1](\mathcal L_{12},\mathcal L'_{12}) \otimes BCF[1](\mathcal L_{23},\mathcal L'_{23})
\otimes BCF[1](\mathcal L_{34},\mathcal L'_{34}) \nonumber\\
\qquad{}\otimes CF^{\rm op}(\mathcal L'_{12},\mathcal L'_{23};\mathcal L_{13})
\otimes BCF[1](\mathcal L_{13},\mathcal L'_{13})
\otimes CF^{\rm op}(\mathcal L'_{13},\mathcal L'_{34};\mathcal L_{14})\nonumber \\
\qquad{}\otimes CF(\mathcal L_{24};\mathcal L'_{23},\mathcal L'_{34}) \otimes BCF[1](\mathcal L_{24},
\mathcal L'_{24})
\otimes CF(\mathcal L_{14},\mathcal L'_{14})\label{new1134}
\end{gather}
to $CF^{\rm op}(\mathcal L_{12},\mathcal L_{24};\mathcal L'_{14})$.}

The double pants transformation \smash{$\mathscr{DPT}^{\vec b}$} is such a map.
Note that in \eqref{new1123} $\mathcal L'_{24};\mathcal L_{23},\mathcal L_{34}$
appears in $CF(\dots)$ and all similar triple appears in
$CF^{\rm op}(\dots)$. This coincides with \eqref{new1134}.
The main property of double pants transformation, that is, Lemma~\ref{lem1110}, implies that \smash{$\mathscr{DPT}^{\vec b}$}
gives a quatro-module homomorphism.

Thus we obtained a natural transformation
from \eqref{1123form} to \eqref{1126form}. The fact that it induces an
isomorphism for objects
can be proved in the same way as the last part of the proof of Theorem~\ref{thm109} (see Section~\ref{sec:compfuncmain}).
We can combine it with the argument of the proof of Proposition~\ref{prop1212}
to complete the proof of Theorem~\ref{asscompmain} in the same way as the
last step of the proof of Theorem~\ref{thm93} (see Section~\ref{sec:compfunc21}).
\end{proof}

\section[Two different ways to compactify the moduli space\\
of pseudo-holomorphic disks in the direct product]{Two different ways to compactify the moduli space\\
of pseudo-holomorphic disks in the direct product}
\label{sec:directcomp}

\subsection{The reason why we need a different compactification}

Let $(L_{12},\sigma_{12})$ be a
$\pi_1^*(V_1 \oplus TX_1) \oplus \pi_2^*(V_2)$ relatively spin Lagrangian submanifold
of $-X_1 \times X_2$. (Here $V_i$ is a vector bundle
on $(X_i)_{[3]}$ for $i=1,2$.)
Let us consider the set
\smash{$\mathring{\mathcal M}(L_{12};\vec a_{12};E)$},
which we defined in Definition~\ref{defn314}.
\begin{defn}
The subset \smash{$\overset{\ \text{\tiny $\circ\circ$}}{\mathcal M}(L_{12};\vec a_{12};E)$}
of \smash{$\mathring{\mathcal M}(L_{12};\vec a_{12};E)$}
consists of the equivalence classes
$[(\Sigma;u;\vec z;\gamma)]$ such that
$\Sigma$ is a disk. In other words, it consists of
the stable maps with no sphere or disk bubbles.

\end{defn}
In Section~\ref{sec:HFIm} (see formula~\eqref{form317}), we compactified
\smash{$\overset{\ \text{\tiny $\circ\circ$}}{\mathcal M}(L_{12};\vec a_{12};E)$}
to ${\mathcal M}(L_{12};\vec a_{12};E)$.
In Definition~\ref{defn523000}, we did not
use this compactification but
mentioned that we use a
slightly different compactification
${\mathcal M}'(L_{12};\vec a_{12};E)$
to define the partial compactification
\smash{$\mathring{\mathcal M}(\vec a_1,\vec a_{12},\vec a_2;a_-,a_+;E)$}
of the space
\smash{$\overset{\ \text{\tiny $\circ\circ$}}{\mathcal M}(\vec a_1,\vec a_{12},\vec a_2;a_-,a_+;E)$}.
See Remark~\ref{Remark524}.
In this section, we define this compactification~${{\mathcal M}'(L_{12};\vec a_{12};E)}$ and its Kuranishi structure.

We first explain the reason why we need to use
different compactification from
${\mathcal M}(L_{12};\vec a_{12};E)$.
Actually,
the space ${\mathcal M}(\vec a_1,\vec a_{12},\vec a_2;a_-,a_+;E)$
will not carry Kuranishi structure if we use
the compactification ${\mathcal M}(L_{12};\vec a_{12,j};E)$
in \eqref{eq518}.
We explain its reason in the following example.
\begin{exm}\label{exa102}
We consider a neighborhood of
an element $(\xi,\eta)$ of the fiber product.
We define $\xi$, $\eta$ below.
Let
\[
\xi = ([-1,1]\times \R;\varnothing,(0,0),\varnothing;u_1,u_2;\gamma_1,\gamma_{12},\gamma_2)
\in \overset{\ \text{\tiny $\circ\circ$}}{\mathcal M}(\varnothing,(a_{12}),\varnothing;a_-,a_+;E_1).
\]
In other words, we consider the case when the source curve $\Sigma$
is $[-1,1]\times \R$ and has no sphere bubble,
and consider only one (boundary) marked point $(0,0)$
which is $(0,0) \in \{0\} \times \R$. See Figure~\ref{Figure14-1}.
\begin{figure}[ht]
\centering
\includegraphics[scale=0.25]{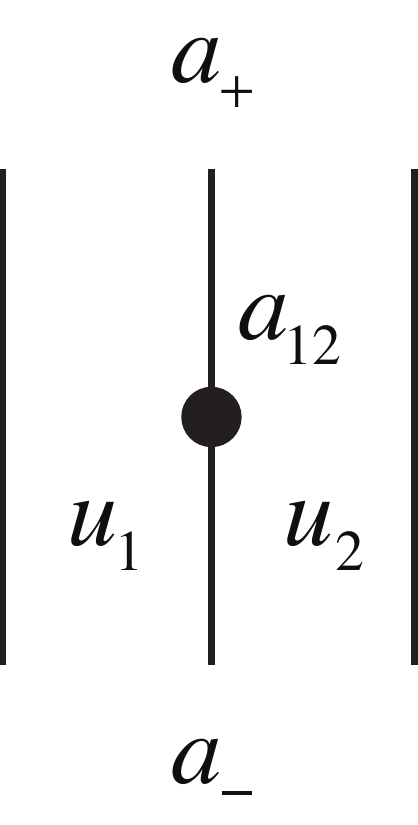}
\caption{Element $\xi$.}
\label{Figure14-1}
\end{figure}

We also consider
\begin{equation}\label{form101}
\eta = (\Sigma;u_3;(1);\gamma_3)
\in \mathring{\mathcal M}(L_{12};(a_{12});E_2).
\end{equation}
Here $\Sigma$ is the union of $D^2$ and $S^2$ glued at $0 \in D^2$
and $[\infty] \in S^2 = \C \cup \{\infty\}$,
$u_3 \colon \Sigma \to -X_1 \times X_2$ is a pseudo-holomorphic map
such that $u_3(\partial\Sigma_3) \subset L_{12}$ and
$1 \in \partial \Sigma$ is a boundary marked point.
See Figure~\ref{Figure14-2}.
\begin{figure}[ht]
\centering
\includegraphics[scale=0.25]{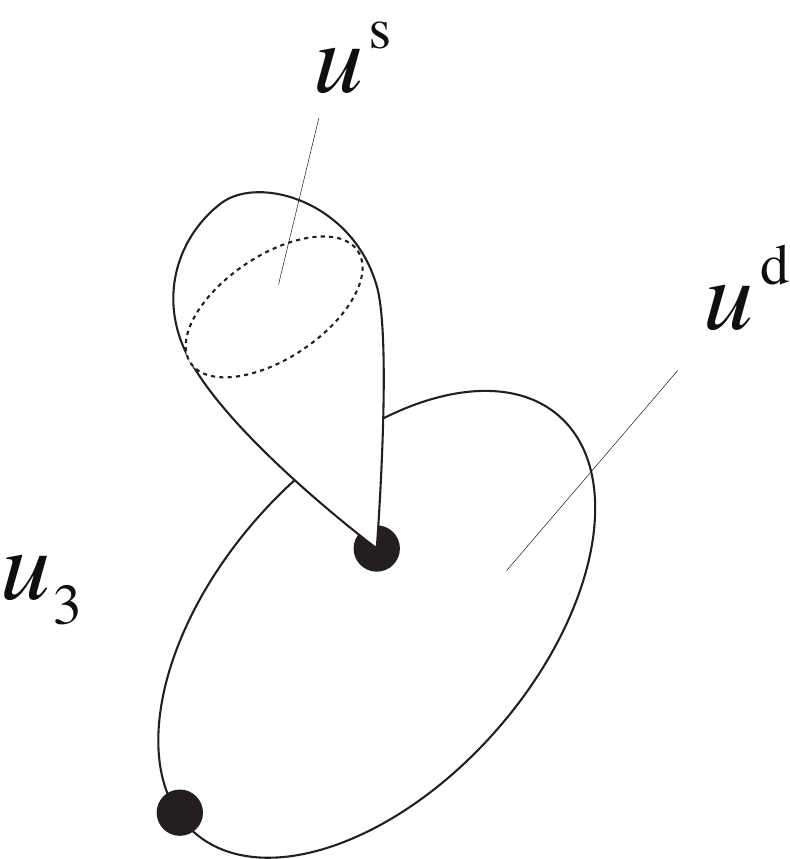}
\caption{Element $\eta$.}
\label{Figure14-2}
\end{figure}

We assume $(u_1,u_2)(0,0) = u_3(1)$ and
regard the pair $(\xi,\eta)$ as an
element of the fiber product
\[
\overset{\ \text{\tiny $\circ\circ$}}{\mathcal M}(\varnothing,(a_{12}),\varnothing;a_-,a_+;E_1)
\times_{L_{12}(a_{12})} {\mathcal M}(L_{12};(a_{12});E_2).
\]
This fiber product is similar to \eqref{eq518} but
we use
${\mathcal M}(L_{12};(a_{12});E_2)$
in place of ${\mathcal M}'(L_{12};(a_{12});E_2)$.
We assume $E = E_1 + E_2$.
See Figure~\ref{Figure14-3}.
\begin{figure}[ht]
\centering
\includegraphics[scale=0.25]{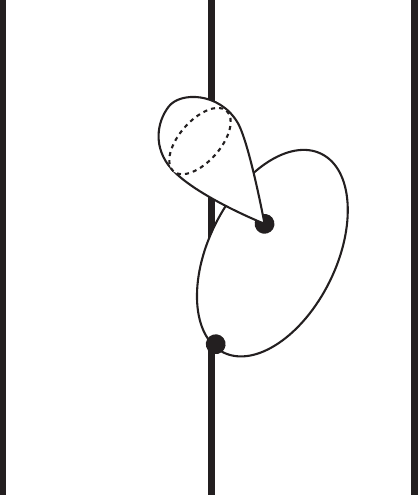}
\caption{Element $(\xi,\eta)$.}
\label{Figure14-3}
\end{figure}

Let us consider a neighborhood of this element
in the compactified moduli space.
For simplicity, we assume that the element $(\xi,\eta)$
is Fredholm regular in the fiber product.
We put
$u^{\rm d} = u_3\vert_{D^2}$ and
$u^{\rm s} = u_3\vert_{S^2}$.
We denote by
$V^{\rm quil}$, $V^{\rm d}$, $V^{\rm s}$
the parameter to deform $\xi$, $u^{\rm d}$,
$u^{\rm s}$, respectively.
We have two kinds of extra parameters
which resolve the singular point.
One is $[0,\varepsilon)$ which parametrizes the way
to resolve the boundary node and
the other is $D^2_{\varepsilon}$
which parametrizes the way
to resolve the interior node.
Therefore, we might imagine that the gluing analysis
implies that the neighborhood of
$(\xi,\eta)$ in
${\mathcal M}(\varnothing,\varnothing,\varnothing;a_-,a_+;E)$ is parametrized by
\begin{equation}\label{form102}
V^{\rm quil} \times_{L_{12}(a_{12})} V^{\rm s}
\times_{X_1 \times X_2} V^{\rm d}
\times [0,\varepsilon) \times D^2_{\varepsilon}.
\end{equation}
However, \eqref{form102} does {\it not} parametrize
a neighborhood of $(\xi,\eta)$ correctly.
To see this we examine the process of gluing more
carefully.
Actually it suffices to see the pre-gluing, which is the
process to obtain approximate solution of
the nonlinear Cauchy--Riemann equation by using partition of unity.
(The process to modify it to obtain an actual solution
is the same as other well-established cases.)

The maps $u^{\rm d}$ and $u^{\rm s}$ are
maps to the direct product
$-X_1 \times X_2$.
So we write $u^{\rm d} = \bigl(u^{\rm d}_1,u^{\rm d}_2\bigr)$
and $u^{\rm s} = (u^{\rm s}_1,u^{\rm s}_2)$.

We first glue $(u_1,u_2)$ with $u^{\rm d}$. This gluing is parametrized
by the parameter $\rho \in [0,\varepsilon)$
with~${\rho \ne 0}$. We put $\overline u_1(s,t) = u_1(-s,t)$ and regard
$(\overline u_1,u_2)$ as a map
from a neighborhood of~$(0,0)$ in $[0,1] \times \R$
to $-X_1 \times X_2$
and glue it with $u^{\rm d}$.
We obtain a map $(\overline u_1,u_2) \#_{\rho} u^{\rm d}$
from $[0,1] \times \R$
to~${-X_1 \times X_2}$.

We them regard $(\overline u_1,u_2) \#_{\rho} u^{\rm d}$
as a pair of maps $(u'_1,u'_2)$ where
$u'_1 \colon [0,1] \times \R \to X_1$,
$u'_2 \colon [0,1] \times \R \to X_2$
such that $(u'_1(0,t),u'_2(0,t)) \in L_{12}$.
By an abuse of notation, we may regard the pair~${(u'_1,u'_2)}$
as $\bigl(u_1\#_{\rho} u^{\rm d}_1,u_2\#_{\rho} u^{\rm d}_2\bigr)$.
See Figure~\ref{Figure14-4}.
\begin{figure}[ht]
\centering
\includegraphics[scale=0.25]{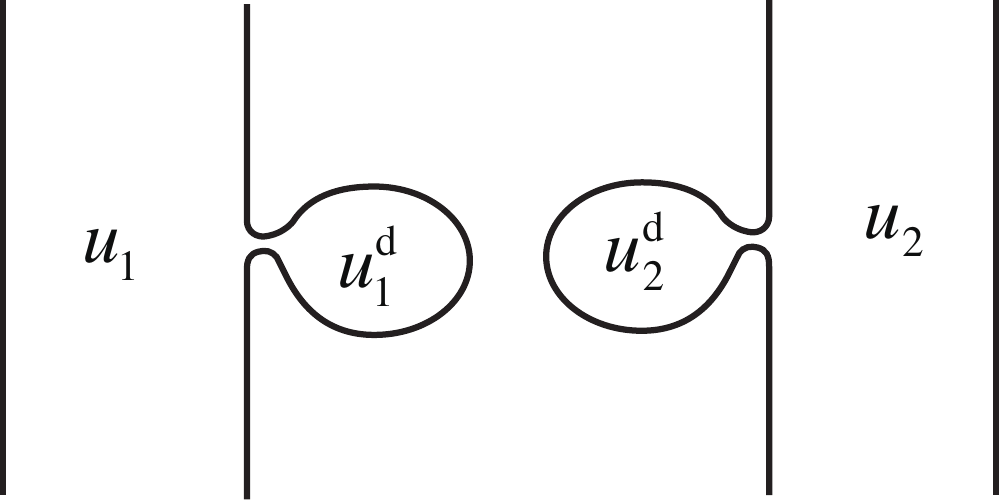}
\caption{$\bigl(u_1\#_{\rho} u^{\rm d}_1,u_2\#_{\rho} u^{\rm d}_2\bigr)$.}
\label{Figure14-4}
\end{figure}

We next glue $u^{\rm s}$ to this pair
$\bigl(u_1\#_{\rho} u^{\rm d}_1,u_2\#_{\rho} u^{\rm d}_2\bigr)$. This gluing is parametrized by the parameter~${\theta \in D^2_{\varepsilon}}$.
We assume $\theta \ne 0$.
We observe that the marked point $0 \in D^2$ at which we
glue the sphere bubble becomes a pair of points
$(-c(\rho),0) \in [-1,0] \times \R$ and $(c(\rho),0)
\in [0,1] \times \R$ after the first gluing.
Therefore, when we glue $u^{\rm s}$
we glue $\overline u^{\rm s}_1 \colon S^2 \to X_1$ to
$u_1\#_{\rho} u^{\rm d}_1$ at the point~${(-c(\rho),0)}$ and
glue $u^{\rm s}_2 \colon S^2 \to X_2$
to
$u_2\#_{\rho} u^{\rm d}_2$ at the point
$(0,c(\rho))$.
(Here $\overline u^{\rm s}_1$ is obtained from
$u^{\rm s}_1$ by using anti-holomorphic involution
of the source.)
We thus can write the element obtained by the
gluing as
$\bigl(u_1\#_{\rho} u^{\rm d}_1 \#_{\theta} \overline u^{\rm s}_1,u_2\#_{\rho} u^{\rm d}_2
\#_{\theta} u^{\rm s}_2\bigr)$.
See Figure~\ref{Figure14-5}.
\begin{figure}[ht]
\centering
\includegraphics[scale=0.2]{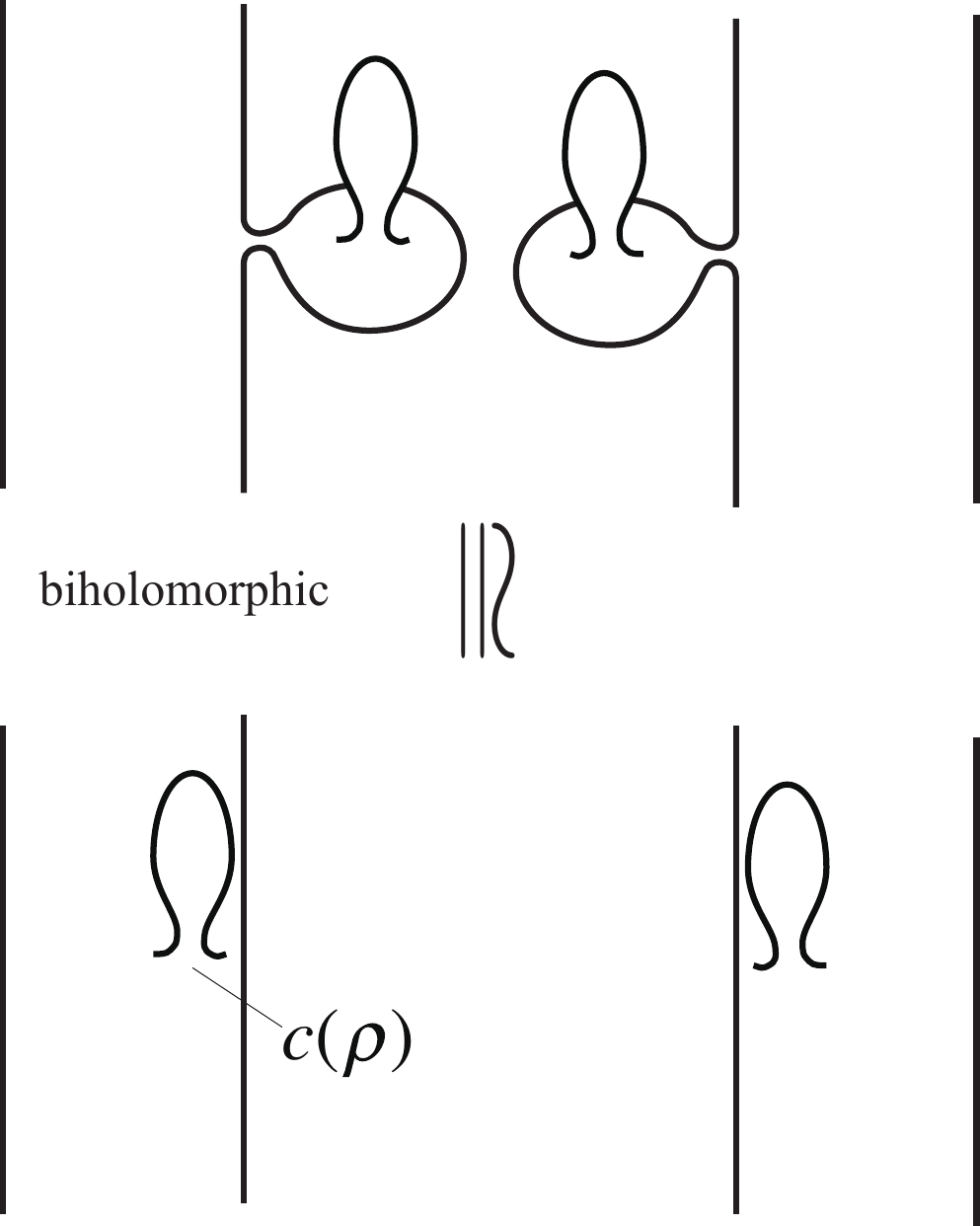}
\caption{$\bigl(u_1\#_{\rho} u^{\rm d}_1 \#_{\theta} \overline u^{\rm s}_1,u_2\#_{\rho} u^{\rm d}_2
\#_{\theta} u^{\rm s}_2\bigr)$.}
\label{Figure14-5}
\end{figure}

In this way, we obtain a family of approximate solutions
parametrized by \eqref{form102}.

Now the issue is that this family does not have
correct dimension. In fact, it has two more parameters
than the correct parameter.
Let us elaborate on this point below.

Let $v \colon S^2 \to S^2$ be a biholomorphic map which preserves
$\infty \in \C \cup \{\infty\} \cong S^2$.
We remark that $u_2^{\rm s} \circ v$ and $u_2^{\rm s} $ are the same
element of the moduli space of pseudo-holomorphic spheres with
one marked point in $X_1$.
However, $(u_1^{\rm s},u_2^{\rm s} \circ v)$ is a different element
from $(u_1^{\rm s},u_2^{\rm s} )$
in the moduli space of pseudo-holomorphic spheres
with one marked point in $-X_1 \times X_2$.
Thus~${\bigl(u_1\#_{\rho} u^{\rm d}_1 \#_{\theta} \overline u^{\rm s}_1,u_2\#_{\rho} u^{\rm d}_2
\#_{\theta} u^{\rm s}_2\circ v\bigr)}$
may become the same element as $\bigl(u_1\#_{\rho} u^{\rm d}_1 \#_{\theta} \overline u^{\rm s}_1,
u_2\#_{\rho} u^{\rm d}_2\allowbreak
\#_{\theta} u^{\rm s}_2\bigr)$ but
$(u_1,u_2\circ v) \ne (u_1,u_2)$.

Another point is that, using the notation
$\bigl(u_1\#_{\rho} u^{\rm d}_1 \#_{\theta} \overline u^{\rm s}_1,u_2\#_{\rho} u^{\rm d}_2
\#_{\theta} u^{\rm s}_2\bigr)$,
we can glue $\overline u^{\rm s}_1$ and~$u^{\rm s}_2$ by
different gluing parameter at interior nodes. Namely, we have a family
of elements of our moduli space
$\bigl(u_1\#_{\rho} u^{\rm d}_1 \#_{\theta_1} \overline u^{\rm s}_1,u_2\#_{\rho} u^{\rm d}_2
\#_{\theta_2} u^{\rm s}_2\bigr)$
where $\theta_1 \ne \theta_2$ may occur.

In fact, a part of the freedom to reparametrize the first (but not the second)
factor by $v$ corresponds to the freedom to choose $\theta_1 \ne \theta_2$.
We will elaborate on this point.
We identify~${\bigl(S^2,\infty\bigr) = (\C\cup \{\infty\},\infty)}$.
For $\mathfrak z \in \C$ in a neighborhood of $1$, we define
$v_{\mathfrak z} \colon (\C\cup \{\infty\},\infty)
\to (\C\cup \{\infty\},\infty)$ by $v_{\mathfrak z}(z) = \mathfrak z z$.
Then the element
$
\bigl(u_1\#_{\rho} u^{\rm d}_1 \#_{\mathfrak z\theta_1} \overline u^{\rm s}_1,u_2\#_{\rho} u^{\rm d}_2
\#_{\theta_2} u^{\rm s}_2\circ v_{\mathfrak z}\bigr)
$
represents the same element as
$
\bigl(u_1\#_{\rho} u^{\rm d}_1 \#_{\theta_1} \overline u^{\rm s}_1,u_2\#_{\rho} u^{\rm d}_2
\#_{\theta_2} u^{\rm s}_2\bigr)$.
See Remark~\ref{rem12333} below.
We now observe that the real dimension of the group of automorphisms of $(\C\cup \{\infty\},\infty)$
is $4$.
On the other hand, the extra parameter by allowing $\theta_1 \ne \theta_2$
is $2$. Thus we can conclude the dimension of~\eqref{form102} is~2 plus the correct dimension of our
moduli space.

In other words, we can not define Kuranishi structure
of our compactification if we use ${\mathcal M}(L_{12};\vec a_{12};E)$
in place of ${\mathcal M}'(L_{12};\vec a_{12};E)$
in \eqref{eq518}.
\end{exm}

\begin{rem}\label{rem12333}
To elaborate on the fact
\[
\bigl(u_1\#_{\rho} u^{\rm d}_1 \#_{\mathfrak z\theta_1} \overline u^{\rm s}_1,u_2\#_{\rho} u^{\rm d}_2
\#_{\theta_2} u^{\rm s}_2\circ v_{\mathfrak z}\bigr)
\sim
\bigl(u_1\#_{\rho} u^{\rm d}_1 \#_{\theta_1} \overline u^{\rm s}_1,u_2\#_{\rho} u^{\rm d}_2
\#_{\theta_2} u^{\rm s}_2\bigr),
\]
we consider the case when $\theta_1=\theta_2 = 0$, that is,
\begin{equation}\label{isotobeproved}
\bigl(u_1\#_{\rho} u^{\rm d}_1 \#_{0} \overline u^{\rm s}_1,u_2\#_{\rho} u^{\rm d}_2
\#_{0} u^{\rm s}_2\circ v_{\mathfrak z}\bigr)
\sim
\bigl(u_1\#_{\rho} u^{\rm d}_1 \#_{0} \overline u^{\rm s}_1,u_2\#_{\rho} u^{\rm d}_2
\#_{0} u^{\rm s}_2\bigr).
\end{equation}
In this case, the domain of those elements are depicted as in Figure~\ref{Figure56+1} below.
\begin{figure}[ht]
\centering
\includegraphics[scale=0.4]{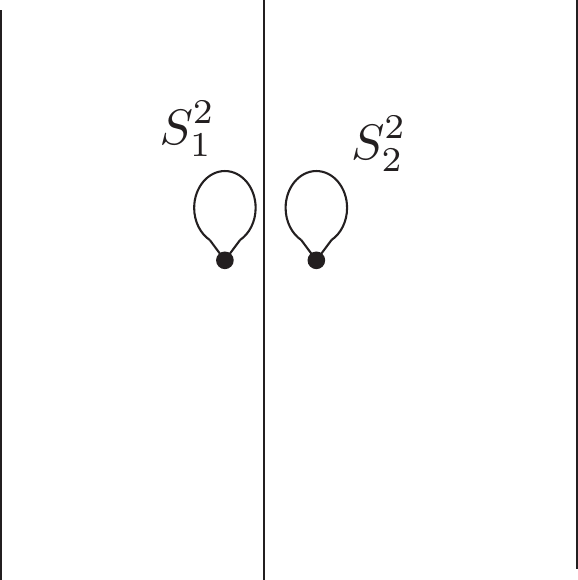}
\caption{The domain in the case when $\theta_1= \theta_2 = 0$.}
\label{Figure56+1}
\end{figure}

There are two sphere bubbles on the domain.
We denote by $S^2_1$ and $S^2_2$ the
sphere bubbles which lie in the left and the right of the
seam, respectively.
The maps on $S^2_1$ and $S^2_2$ are $\overline u^{\rm s}_1$
and~$u^{\rm s}_2$ for the right-hand side of
\eqref{isotobeproved}.
In the case of left-hand side of \eqref{isotobeproved},
the maps on $S^2_1$ and $S^2_2$ are $\overline u^{\rm s}_1$
and $u^{\rm s}_2\circ v_{\mathfrak z}$, respectively.
We define $\hat v_{\mathfrak z}$ to be an isomorphism from the configuration
as in Figure~\ref{Figure56+1} to itself
so that $\hat v_{\mathfrak z}$ is the identity map
outside $S^2_2$ and is $v_{\mathfrak z}$ on $S^2_2$.
Then it is easy to see that
\[
\bigl(u_1\#_{\rho} u^{\rm d}_1 \#_{0} \overline u^{\rm s}_1,u_2\#_{\rho} u^{\rm d}_2
\#_{0} u^{\rm s}_2\bigr)
\circ \hat v_{\mathfrak z}
=
\bigl(u_1\#_{\rho} u^{\rm d}_1 \#_{0} \overline u^{\rm s}_1,u_2\#_{\rho} u^{\rm d}_2
\#_{0} u^{\rm s}_2\circ v_{\mathfrak z}\bigr).
\]
This implies the equivalence \eqref{isotobeproved}.

We can choose the various additional data which we use
to perform the gluing process so that
the equivalence in the case $\theta_i = 0$ can be
extended to the case $\theta_i \ne 0$.
(We omit the detail of this part since the rigorous proof is
not necessary for the proof of our results. The discussion here
is a motivation to introduce new compactification.)

We also observe the following.
We take the limit as $\rho$ goes to $0$ in \eqref{isotobeproved}.
The domain depicted by Figure~\ref{Figure56+1} converges to
the domain depicted by Figure~\ref{Figure14-3}.
The automorphisms $\hat v_{\mathfrak z}$ however cannot be extended to this limit.
In fact, in the domain of Figure~\ref{Figure14-3} two
sphere bubbles become the one sphere bubble and so
we are not allowed to take two different biholomorphic maps
on the sphere bubble.
Therefore, `the limit' of left and right-hand sides (as $\rho$ goes to $0$)
are not equivalent.
By this reason, it seems likely that it is impossible to define an appropriate
topology which is Hausdorff, if we use ${\mathcal M}(L_{12}; \vec a_{12,j};E_{12,j})$
in place of
${\mathcal M}'(L_{12}; \vec a_{12,j};E_{12,j})$
in Theorem~\ref{therem530}\,(2).
\end{rem}
\begin{rem}
A similar problem already appeared in \cite{fooo:inv}
(see the proof of \cite[Lemma 6.62]{fooo:inv}).
In \cite[Lemma 6.62]{fooo:inv}, we compared
two moduli spaces. One is the space of pseudo-holomorphic
maps $u$ from a disk
to $-X \times X$ so that $u\bigl(\partial D^2\bigr)$ lies in the diagonal.
The other is the space of pseudo-holomorphic
maps $u'$ from a sphere to $X$.
We can use reflection principle to identify those two moduli spaces.
When we consider their stable map compactifications
the identification does not extend.
To explain this fact, we consider the case when $u$ is a map from $D^2$ with a~sphere bubble
to $-X \times X$ so that $u\bigl(\partial D^2\bigr)$ lies in the diagonal.
Suppose that $u$ is $(\overline u_1,u_2)$ on the bubble.
Then the corresponding element $u'$ is a map
from $S^2$ with two sphere bubbles and
the maps on those sphere bubbles are $u_1$ and $u_2$, respectively
(see Figure~\ref{Figure56+2}).
When we replace $u_2$ by~$u_2 \circ v_{\mathfrak z}$
the object in the compactification of the moduli space
of disks changes. However,
the corresponding objects in the compactification of the moduli space
of spheres are equivalent.
This is similar to the situation of Example~\ref{exa102}
and Remark~\ref{rem12333}.

\begin{figure}[ht]
\centering
\includegraphics[scale=0.4]{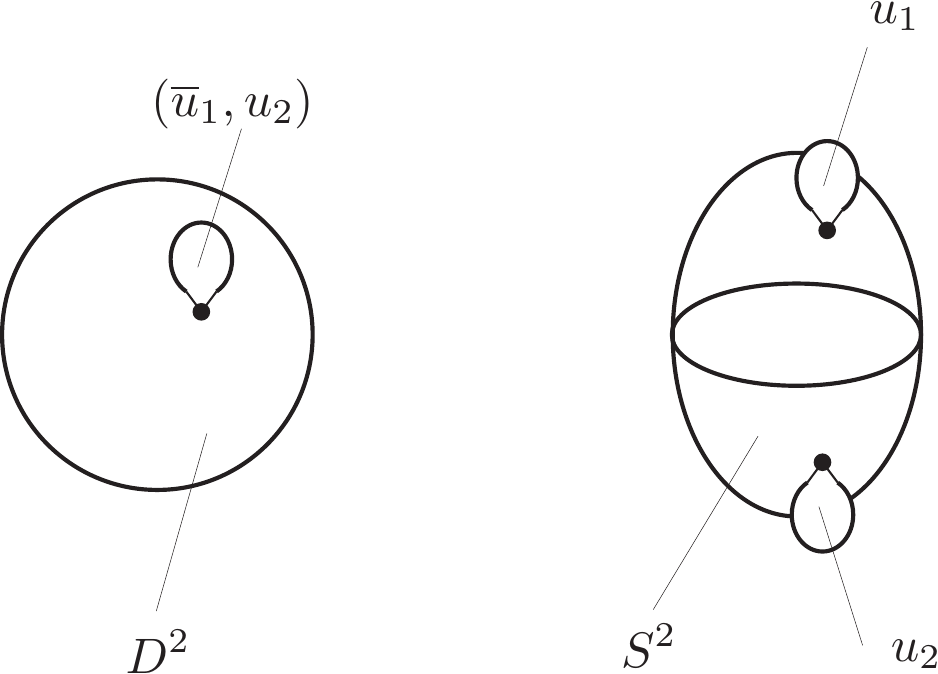}
\caption{Reflection principle at infinity.}
\label{Figure56+2}
\end{figure}
In \cite{fooo:inv}, the problem is slightly less serious since
there we need to show two well-defined numbers
to coincide. So we can use the fact that
the problem occurs only in codimension $\ge 2$ strata
and use dimension counting argument.
Here we need to work out the chain level argument.
So we describe the different compactification ${\mathcal M}'(L_{12};\vec a;E)$
in detail in this section.
\end{rem}
\begin{rem}
We remark that the problem of different reparametrizations
applied to the first and the second factors in the bubble,
which we described in Example~\ref{exa102}, does {\it not} occur
for the disk bubble but occurs only for the sphere bubble.
Let us elaborate on this point below.

Let us consider the same $\xi$ as Example~\ref{exa102}.
We replace $\eta$ as in \eqref{form101} by
\[
\eta = \bigl(D^2;u_3;(1);\gamma_3\bigr)
\in \mathring{\mathcal M}(L_{12};(a_{12});E_2).
\]
Namely, we assume the source curve of $\eta$ is a disk.
The group of automorphisms of the pair~$\bigl(D^2,1\bigr)$ of a
disk with one boundary marked point $1 \in \partial D^2$
is identified with the group of affine
transformations $z \mapsto \varphi_{a,b}(z) = az + b$ with
$a \in \R_+$ and $b \in \R$. Here we identify
$D^2 \setminus\{1\}$ with the upper half plane
$\{z \in \C \mid \operatorname{Im} z \ge 0\}$.

Let $u^{\rm d} = \bigl(u^{\rm d}_1,u^{\rm d}_2\bigr)$
be a representative of an element of
\smash{$\mathring{\mathcal M}(L_{12};(a_{12});E_2)$},
where $u^{\rm d}_1 \colon D^2 \to X_1$ and $u^{\rm d}_2 \colon D^2 \to X_2$.

Note that in this case $\bigl(u^{\rm d}_1,u^{\rm d}_2\circ\varphi_{a,b}\bigr)$
does not represent an element of \smash{$\mathring{\mathcal M}(L_{12};(a_{12});E_2)$}
in general,
since this element may not satisfy the boundary condition.
\end{rem}
\begin{rem}
The compactification ${\mathcal M}'_{\ell,\ell_1,\ell_2}(L_{12};\vec a;E)$
which we will define in the next subsection,
is `smaller' than
${\mathcal M}_{\ell,\ell_1,\ell_2}(L_{12};\vec a;E)$.
An intuitive reason why we need smaller compactification
lies in the fact that
\smash{$\overset{\ \text{\tiny $\circ\circ$}}{\mathcal M}_{\rm QT}(L_1,L_{12},L_2;p,q)$} is also `smaller' than
\smash{$\overset{\ \text{\tiny $\circ\circ$}}{\mathcal M}(L_{12},L_1 \times L_2;p,q)$}.
In fact, we can define
\[
\mathfrak{forget}_{\rm QT} \colon\
\overset{\ \text{\tiny $\circ\circ$}}{\mathcal M}(L_{12},L_1 \times L_2;p,q)
\to \overset{\ \text{\tiny $\circ\circ$}}{\mathcal M}_{\rm QT}(L_1,L_{12},L_2;p,q).
\]
Here $p = (p_1,p_2) \in (L_1\times L_2) \cap L_{12}$,
$q = (q_1,q_2) \in (L_1\times L_2) \cap L_{12}$.

The space \smash{$\overset{\ \text{\tiny $\circ\circ$}}{\mathcal M}(L_{12},L_1 \times L_2;p,q)$}
is a partially compactified moduli space of pseudo-holomor\-phic strips.
Its element
is an equivalence class of $(\Sigma;u)$ where $\Sigma$ is a strip
$[0,1]\times \R$ with
trees of sphere bubbles and $u \colon \Sigma \to -X_1 \times X_2$ is a
pseudo-holomorphic map. We require that
\smash{$u\vert_{\{0\} \times \R}$}
(resp.\ $u\vert_{\{1\} \times \R}$ ) lifts to a map to $\tilde L_{12}$
(resp.\ to $\tilde L_{1} \times \tilde L_{2}$)
and $u$ is asymptotic to $p$ (resp.\ $q$)
when the $\R$-factor of the domain goes to $-\infty$ (resp.\ $+\infty$).

The space \smash{$\overset{\ \text{\tiny $\circ\circ$}}{\mathcal M}_{\rm QT}(L_1,L_{12},L_2;p,q)$}
is a partially compactified moduli space of pseudo-holo\-mor\-phic quilt.
Its element is an equivalence class of $(\Sigma';u_1,u_2)$.
Here $\Sigma'$ is $[-1,1]\times \R$ with
trees of sphere bubbles, which is decomposed to
$\Sigma'_1 \cup \Sigma'_2$ such that $\Sigma'_1$ (resp.\ $\Sigma'_2$)
is $[-1,0]\times \R$ (resp.~${[0,1]\times \R}$) together with sphere bubbles.
$u_i \colon \Sigma'_i \to X_i$ is a pseudo-holomorphic map.
We require that~\smash{$u_1\vert_{\{-1\}\times \R}$} (rest.\ \smash{$u_2\vert_{\{1\}\times \R}$})
lifts to a~map to $\tilde L_1$ \big(resp.\ $\tilde L_2$\big).
We also require a~matching condition, that is, the map $\tau\mapsto (u_1(0,\tau),u_2(0,\tau))$
lifts to a map to \smash{$\tilde L_{12}$}.
Furthermore, we require asymptotic boundary condition given by $p$, $q$.

We define $\mathfrak{forget}_{\rm QT}$ as follows.
Let $(\Sigma;u)$ represent an element of
\smash{$\overset{\ \text{\tiny $\circ\circ$}}{\mathcal M}(L_{12},L_1 \times L_2;p,q)$}.
We write $u = (u_1,u_2)$ where $u_i$ is a map to $X_i$.
We consider $(\Sigma;u_1)$ and shrink all the unstable sphere components of
$\Sigma$ on which $u_1$ is constant to obtain $\Sigma''_1$.
Using the map $(t,\tau) \mapsto (-t,\tau)$ (which is a~map
$[-1,0] \times \R \to [0,1] \times \R$), we obtain $\Sigma'_1$ from
$\Sigma''_1$. The map $u_1$ induces a~map~${u'_1 \colon \Sigma'_1 \to X_1}$.
In a similar (and simpler) way we obtain $\Sigma'_2$ and $u'_2 \colon \Sigma'_2 \to X_2$.
We glue $\Sigma'_1$ and $\Sigma'_2$ on the line~${\{0\} \times \R}$ to
obtain $\Sigma'$.
It is easy to see that $(\Sigma';u'_1,u'_2)$ represents an element of
\smash{$\overset{\ \text{\tiny $\circ\circ$}}{\mathcal M}_{\rm QT}(L_1,L_{12},L_2;p,q)$}.

Using Lemma--Definition~\ref{lemdef14355}, we can extend the
map $\mathfrak{forget}_{\rm QT}$ so that it includes the case when
the objects have disk bubbles.

The map ${\rm Dob}$ in formula \eqref{dobe17111}
is an inverse of $\mathfrak{forget}_{\rm QT}$ on certain open dense subsets.
\end{rem}

\subsection[The definition of the compactification
${\mathcal M}'(L_{12};\vec a;E)$]{The definition of the compactification
$\boldsymbol{{\mathcal M}'(L_{12};\vec a;E)}$}
\label{sec:difcompex}

Based on the observation in the previous
subsection,
we define the compactification ${\mathcal M}'(L_{12};\vec a; \allowbreak E)$.\index[syindex]{M1primeL12aE@${\mathcal M}'(L_{12};\vec a;E)$}
For later use, we also include the case when
there are interior marked points and will define~\smash{${\mathcal M}'_{\ell,\ell_1,\ell_2}(L_{12};\vec a;E)$}.\index[syindex]{M1primeellell1@${\mathcal M}'_{\ell,\ell_1,\ell_2}(L_{12};\vec a;E)$}

\begin{defn}\label{defn145555}
We consider
objects
\[
\bigl(\bigl(\bigl(\Sigma_1,\vec z_1,\vec z_1^{\rm \,int},\vec w_1^{ \rm int}\bigr),u_1\bigr),\bigl(\bigl(\Sigma_2,\vec z_2,
,\vec z_2^{\rm \,int},\vec w_2^{\rm int}\bigr),u_2\bigr),\mathscr I,\gamma\bigr)
\]
with the following properties:
\begin{enumerate}\itemsep=0pt
\item[(1)]
The space
$\Sigma_i$, $i=1,2$, is a bordered curve of genus zero with one boundary
component. $\vec z_i = (z_{i,0},\dots,z_{i,k})$ are
mutually distinct boundary marked
points of $\Sigma_i$ such that the enumeration of the marked points respects
orientation of the boundary.
\smash{$\vec z_i^{\rm \,int}
= \bigl(z_{i,0}^{\rm int},\dots,z_{i,\ell}^{\rm int}\bigr)$}
and $\vec w_i^{ \rm int} = \bigl(w_{i,1},\dots,w_{i,\ell_i}\bigr)$
are mutually distinct interior marked points on $\Sigma_i$.
Marked points are not nodal points.
\item[(2)]
The maps $u_1 \colon \Sigma_1 \to -X_1$,
$u_2 \colon \Sigma_2 \to X_2$ are pseudo-holomorphic.
(We do not assume that $\bigl(\bigl(\Sigma_i,\vec z_i,\vec z_i^{\rm \,int},\vec w_i^{\rm int}\bigr),u_i\bigr)$
is stable. The stability condition we assume is Definition~\ref{148stable} below.)
\item[(3)]
We shrink all the unstable sphere components of
$\bigl(\Sigma_i,\vec z_i,\vec z_i^{\rm \,int}\bigr)$
(that is, the sphere components which have less than $3$ nodal or
marked points in $\vec z_i^{\rm \,int}$).
We denote by
$\bigl(\Sigma^0_i,\vec z_i,\vec z_i^{\rm \,int}\bigr)$ the marked bordered nodal curve
obtained by this shrinking.
(We use the same symbols~$\vec z_i$, $\vec z_i^{\rm \,int}$ for
marked points by an abuse of notation.)
(We remark that we forget $\vec w_i^{ \rm int}$ when we define $\Sigma^0_i$.)
Then,
$\mathscr I \colon \Sigma^0_1 \to \Sigma^0_2$ is a biholomorphic map
such that
$
\mathscr I(z_{1,j}) = z_{2,j}
$, $
\mathscr I(z^{\rm int}_{1,j}) = z^{\rm int}_{2,j}
$.
(See Figure~\ref{Figure139}.)
\item[(4)]
The map $\gamma \colon \partial \Sigma_1 \setminus \vec z_1 \to \tilde L_{12}$ is continuous and satisfies
\begin{equation}\label{form1444}
i_{L_{12}}(\gamma(z)) = (u_1(z),u_2(\mathscr I(z))).
\end{equation}
(Note that \eqref{form1444} implies $(u_1(z),u_2(\mathscr I(z))) \in L_{12}$ for $z \in \partial\Sigma_1$.)
(We also remark $\partial \Sigma_1 = \partial \Sigma^0_1$.)
\item[(5)]
We require the switching condition, Condition \ref{switching14}, below.
\item[(6)]
We require the stability condition, Definition~\ref{148stable}, below.
\item[(7)]
$
-\int_{\Sigma_1} u_1^*\omega_1 +
\int_{\Sigma_1} u_2^*\omega_2
= E$.
\end{enumerate}
\begin{figure}[ht]
\centering
\includegraphics[scale=0.5]{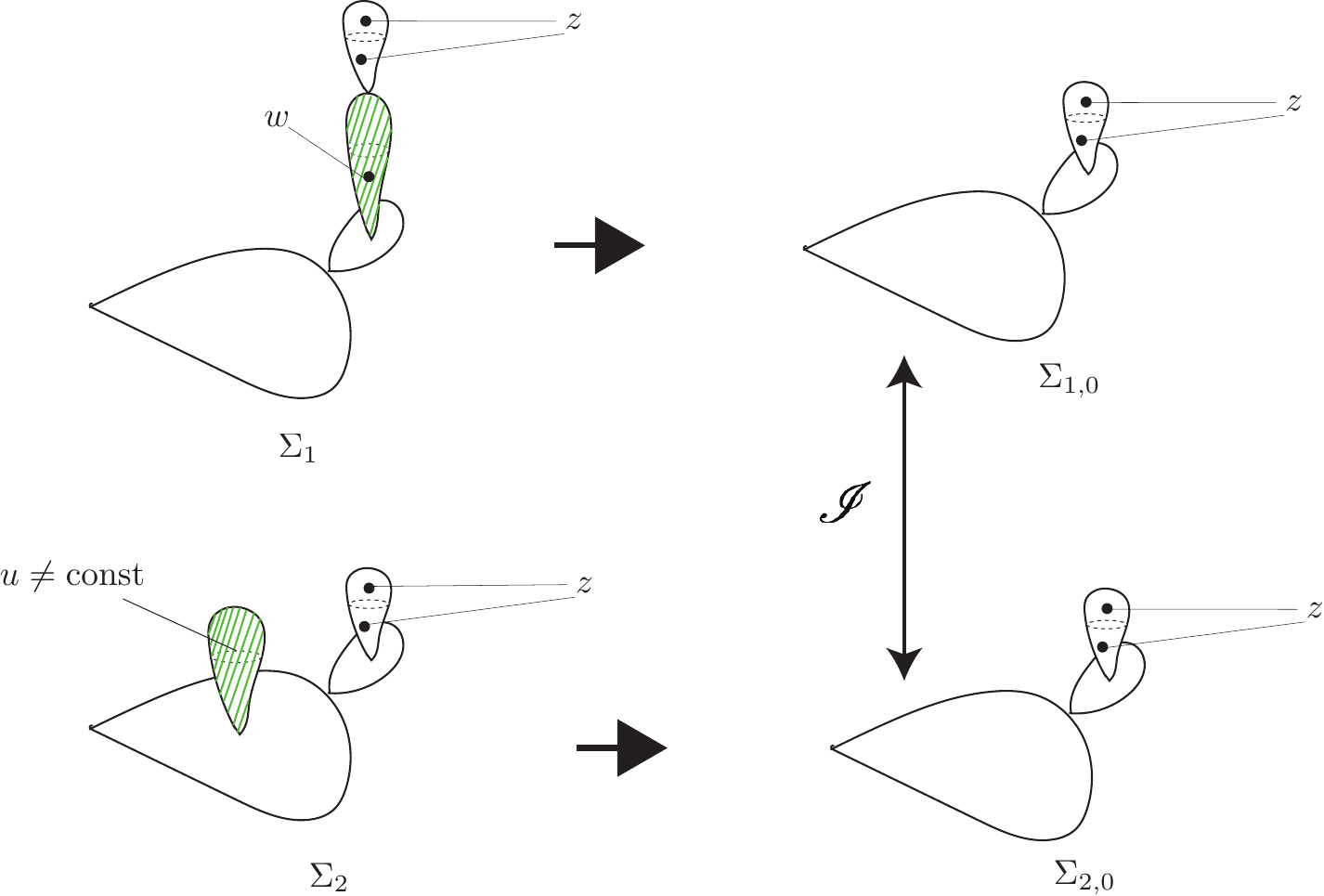}
\caption{Source curve of an element of ${\mathcal M}'_{\ell,\ell_1,\ell_2}(L_{12};\vec a;E)$.}
\label{Figure139}
\end{figure}
We call $z_{i,j}^{ \rm int}$ an {\it interior marked point of first kind}
\index{interior marked points of first kind}
and $w_{i,j}^{ \rm int}$ an {\it interior marked point of second kind}
\index{interior marked points of second kind}.

We denote by ${\mathcal M}'_{\ell,\ell_1,\ell_2}(L_{12};\vec a;E)$ the set of the
equivalence classes of such objects with respect to the
equivalence relation $\sim$ defined in Definition~\ref{defn1477}.

\end{defn}
\begin{conds}\label{switching14}
For each $j$,
$
(\lim_{z \in \partial\Sigma_1, z \uparrow z_{1,j}} \gamma(z),
\lim_{z \in \partial_1\Sigma, z \downarrow z_{1,j}} \gamma(z)
)
\in L_{12}(a_{1,j})$.

\end{conds}
\begin{defn}\label{defn1477}
Let
$\xi = \bigl(\bigl(\bigl(\Sigma_1,\vec z_1,\vec z_1^{\rm \, int},\vec w_1^{ \rm int}\bigr),u_1\bigr),\bigl(\bigl(\Sigma_2,\vec z_2
,\vec z_2^{\rm \, int},\vec w_2^{ \rm int}),u_2\bigr),\mathscr I,\gamma\bigr)$
and let
$\xi' = \bigl(\bigl(\bigl(\Sigma'_1,\vec z^{\,\prime}_1,\vec z_1^{\rm \, int \prime},\vec w_1^{\rm int \prime}\bigr),u'_1\bigr),\bigl(\bigl(\Sigma'_2,\vec z^{\,\prime}_2,\vec z_2^{\rm \,int \prime},
\vec w_2^{\rm int \prime}\bigr),u'_2\bigr),\mathscr I',\gamma'\bigr)$
be objects satisfying (1)--(5) of Definition~\ref{defn145555}.

A {\it weak isomorphism} \index{weak isomorphism} from $\xi$ to $\xi'$ is a pair of maps $(\psi_1,\psi_2)$
with the following properties:
\begin{enumerate}\itemsep=0pt
\item[(1)]
The map $\psi_i \colon \Sigma_i \to \Sigma'_i$ is biholomorphic.
\item[(2)]
$\psi_i(z_{i,j}) = z'_{i,j}$.
\item[(3)]
There exist permutations
$\sigma \colon \{1,\dots,\ell\} \to \{1,\dots,\ell\}$,
$\sigma_i \colon \{1,\dots,\ell_i\} \to \{1,\dots,\ell_i\}$ such that
$\psi_i\bigl(z^{\rm int}_{i,j}\bigr) = z^{\rm int \prime}_{i,\sigma(j)}$,
and that $\psi_i\bigl(w^{\rm int}_{i,j}\bigr) = w^{\rm int \prime}_{i,\sigma_i(j)}$,
for $i=1,2$.
\item[(4)]
$u'_i \circ \psi_i = u_i$ for $i=1,2$.
\item[(5)]
Note that (1)--(3) above implies that $\psi_i$ induces a map
$\overline \psi_i\colon \Sigma^0_i \to \Sigma^{\prime 0}_i$. We require:
$\mathscr I' \circ\overline\psi_1 = \overline\psi_2 \circ \mathscr I$ on $\Sigma^0_1$.
\end{enumerate}
A weak isomorphism $(\psi_1,\psi_2)$ is said to be an {\it isomorphism} \index{isomorphism}
if $\sigma$ and $\sigma_i$ in item (3) are the identity maps.

We say $\xi$
is {\it equivalent} \index{equivalent} to $\xi'$ and write $\xi \sim \xi'$ if there exists an isomorphism from $\xi$ to $\xi'$.
\end{defn}
\begin{defn}\label{148stable}
An object $\xi$
satisfying (1)--(5) of Definition~\ref{defn145555}
is said to be {\it stable} \index{stable} if the set of
isomorphisms from $\xi$ to $\xi$ is finite.

We say $\xi$ is {\it source stable}\index{source stable} if the set of $(\psi_1,\psi_2)$ which satisfies
(1), (2), (3), (5) of Definition~\ref{defn1477} (but not necessary (4)) is finite.

\end{defn}
\begin{rem}
We consider $ \bigl(\bigl(\bigl(\Sigma_1,\vec z_1,\vec z_1^{\rm \,int},\vec w_1^{ \rm int}\bigr),u_1\bigr),\bigl(\bigl(\Sigma_2,\vec z_2
,\vec z_2^{\rm \,int},\vec w_2^{ \rm int}),u_2),\mathscr I,\gamma\bigr)$
such that there exists an unstable disk component $\Sigma_1(a)$ of $\Sigma_1$ on which $u_1$ is constant.
Such an object can still be stable in the sense of Definition~\ref{148stable}.
In fact, if $u_2$ is non-constant on $\Sigma_2(a) = \mathscr I(\Sigma_1(a))$,
then by condition Definition~\ref{148stable}\,(4) there is no continuous family of
automorphisms supported on this component.
\end{rem}

\begin{figure}[ht]
\centering
\includegraphics[scale=0.3]{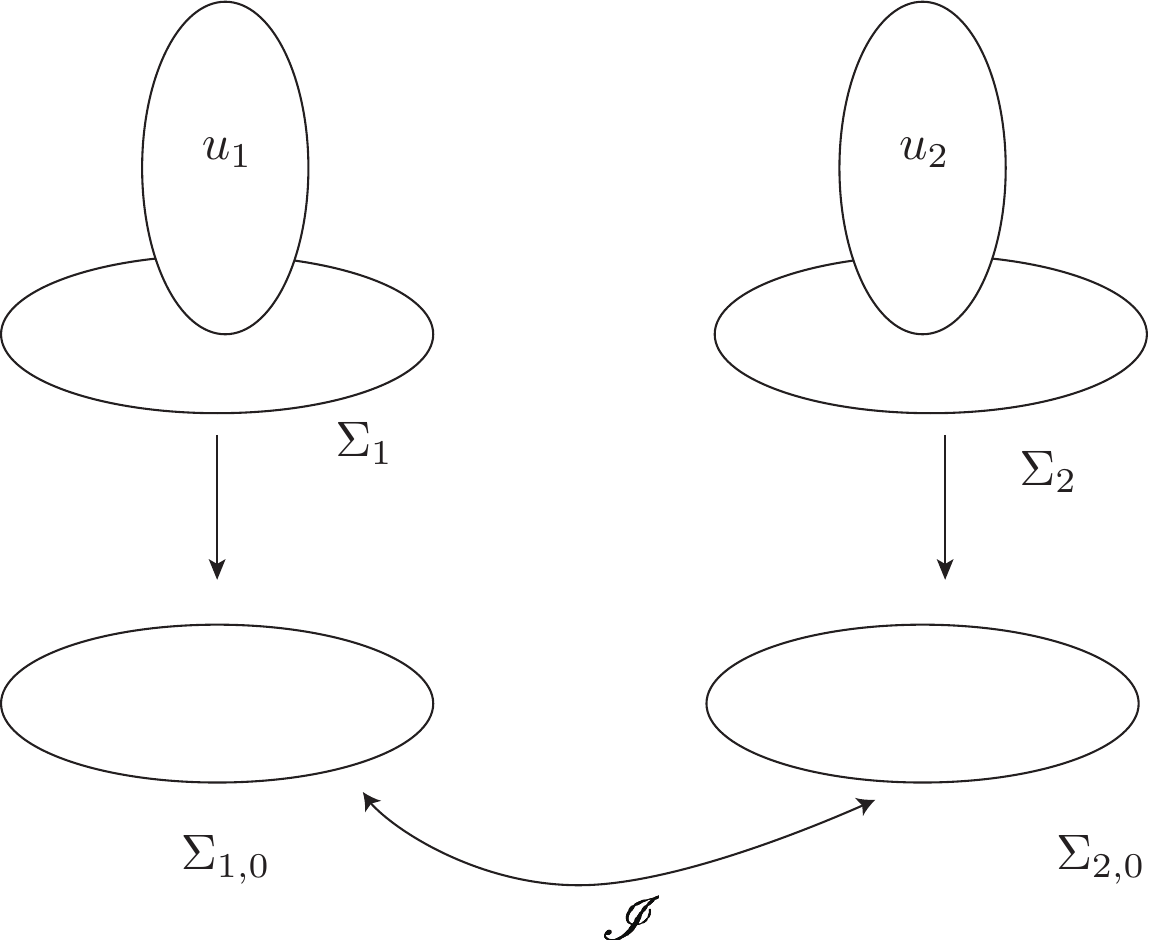}
\caption{Element $\eta'$.}
\label{Figurenewnew}
\end{figure}

\begin{exm}
We consider the situation of Example~\ref{exa102}.
The element $\eta$ corresponds in our
compactification to an element $\eta'$ from the domain
as in Figure~\ref{Figurenewnew}.
$\Sigma_{1,0} \cong \Sigma_{2,0}$ is a~disk in this case and
$\mathscr I$ is the identity map.
$u_1$, $u_2$ are defined on the sphere bubbles rooted on~$\Sigma_{1,0}$,~$\Sigma_{2,0}$, respectively.
By the definition of our equivalence relation,
the object is equivalent if we replace~$u_2$ by
$u_2\circ v$. Here $v \colon S^2 \to S^2$ is a biholomorphic
map which preserves the point $0$ where sphere bubble is
attached. Therefore, the problem mentioned in Example~\ref{exa102}
disappears.
\end{exm}

Let
\begin{equation}\label{form145555}
\mathfrak i = (\mathfrak i_0,\mathfrak i_1,\mathfrak i_2),
\qquad \mathfrak i_0 \colon\ \{1,\dots,\ell\} \to \{1,\dots,\ell'\},\qquad
 \mathfrak i_i \colon\ \{1,\dots,\ell_i\} \to \{1,\dots,\ell'_i\}
\end{equation}
be a triple of injective maps.
It induces a forgetful map
\begin{equation}\label{formass26}
\mathfrak i^* \colon\
{\mathcal M}'_{\ell',\ell'_1,\ell'_2}(L_{12};\vec a;E)
\to
{\mathcal M}'_{\ell,\ell_1,\ell_2}(L_{12};\vec a;E),
\end{equation}
as follows.

Let
\[
 \bigl(\bigl(\bigl(\Sigma_1,\vec z_1,\vec z_1^{\rm \,int},\vec w_1^{ \rm int}\bigr),u_1\bigr),\bigl(\bigl(\Sigma_2,\vec z_2
,\vec z_2^{\rm \,int},\vec w_2^{ \rm int}),u_2\bigr),\mathscr I,\gamma\bigr)
\in {\mathcal M}'_{\ell'}(L_{12};\vec a;E).
\]
We put
$
z^{\rm int \prime}_{i,j} = z^{\rm int}_{i,\mathfrak i(j)}$,
$w^{\rm int \prime}_{i,j} = w^{\rm int}_{i,\mathfrak i_i(j)}$.
We consider
\[
\xi =\bigl(\bigl(\bigl(\Sigma_1,\vec z_1,\vec z_1^{\rm \,int \prime},\vec w_1^{ \rm int \prime}\bigr),u_1\bigr),\bigl(\bigl(\Sigma_2,\vec z_2
,\vec z_2^{\rm \,int \prime},\vec w_2^{ \rm int \prime}\bigr),u_2\bigr),\mathscr I,\gamma\bigr).
\]
We shrink components such that there are infinitely many automorphisms
supported on it and obtain
\[
\xi' = \bigl(\bigl(\bigl(\Sigma'_1,\vec z^{\,\prime}_1,\vec z_1^{\rm \,int \prime},\vec w_1^{ \rm int \prime}),u_1\bigr),\bigl(\bigl(\Sigma'_2,\vec z^{\,\prime}_2
,\vec z_2^{\rm \,int \prime},\vec w_2^{ \rm int \prime}\bigr),u'_2\bigr),\mathscr I',\gamma'\bigr).
\]
We define
$\mathfrak i^*(\xi) = \xi'$.
\begin{defn}
Let
\[
\xi = \bigl(\bigl(\bigl(\Sigma_1,\vec z_1,\vec z_1^{\rm \,int},\vec w_1^{ \rm int}\bigr),u_1\bigr),\bigl(\bigl(\Sigma_2,\vec z_2
,\vec z_2^{\rm \,int},\vec w_2^{ \rm int}),u_2\bigr),\mathscr I,\gamma\bigr)
\]
be an element of \smash{${\mathcal M}'_{\ell',\ell'_1,\ell'_2}(L_{12};\vec a;E)$}.
We say an element
\[
\xi' = \bigl(\bigl(\bigl(\Sigma'_1,\vec z^{\,\prime}_1,\vec z_1^{\rm \,int \prime},\vec w_1^{ \rm int \prime}\bigr),u'_1\bigr),\bigl(\bigl(\Sigma'_2,\vec z^{\,\prime}_2
,\vec z_2^{\rm \,int \prime},\vec w_2^{ \rm int \prime}\bigr),u'_2\bigr),\mathscr I',\gamma'\bigr)
\]
of the space~\smash{${\mathcal M}'_{\ell,\ell_1,\ell_2}(L_{12};\vec a;E)$} is a {\it source stabilization}
\index{source stabilization}
of $\xi$ if the following holds:
\begin{enumerate}\itemsep=0pt
\item[(1)]
There exists $\mathfrak i$ as in \eqref{form145555} such that
$\mathfrak i^*(\xi' ) = \xi$.
\item[(2)]
For any isomorphism $(\psi_1,\psi_2) \colon \xi \to \xi$,
there exists an weak isomorphism $(\psi'_1,\psi'_2) \colon \xi' \to \xi'$
such that the next diagram commutes:
\[
\begin{CD}
\Sigma'_i
 @ >{\psi'_i}>>
\Sigma'_i
 \\
@ VVV @ VVV\\
\Sigma_i
@ >{\psi_i}>>
\Sigma_i,
\end{CD}
\]
where the vertical arrows are the maps shrinking unstable sphere components.
\item[(3)]
The element $\xi'$ is source stable.
\end{enumerate}

We call an interior marked point of $\xi$ an
{\it added marked point} \index{added marked point} if it does not correspond to a~marked point of $\xi'$.
(There are $\ell - \ell'+\ell_i-\ell'_i$ added marked points on each $\Sigma_i$ ($i=1,2$).)

\end{defn}
We next define a topology,
{\it stable map topology}, \index{stable map topology}
on ${\mathcal M}'_{\ell}(L_{12};\vec a;E)$,
in a similar way as \cite[Definition 10.3]{FO}, as follows.

We first consider the case of elements
\[
\xi = \bigl(\bigl(\bigl(\Sigma_1,\vec z_1,\vec z_1^{\rm \,int},\vec w_1^{ \rm int}\bigr),u_1\bigr),\bigl(\bigl(\Sigma_2,\vec z_2
,\vec z_2^{\rm \,int},\vec w_2^{ \rm int}),u_2\bigr),\mathscr I,\gamma\bigr)
\]
and
\begin{align*}
\xi(k)={}& \bigl(\bigl(\bigl(\Sigma_1(k),\vec z_1(k),\vec z_1^{\rm \,int}(k),\vec w_1^{ \rm int}(k)\bigr),u_1(k)\bigr),\\
&\bigl(\bigl(\Sigma_2(k),\vec z_2(k),
\vec z_2^{\rm \,int}(k),\vec w_2^{ \rm int}(k)),u_2(k)\bigr),\mathscr I(k),\gamma(k)\bigr)
\end{align*}
of ${\mathcal M}'_{\ell}(L_{12};\vec a;E)$
such that $\xi$ and $\xi(k)$ are all source stable. In such case, we define the following.
\begin{defn}\label{defn1211new}
We say
\smash{$\operatorname{lims}_{k\to\infty} \xi(k) = \xi$}
if the following holds:\index[syindex]{lims@${\rm lims}$}
\begin{enumerate}\itemsep=0pt
\item[(1)]
$\bigl(\Sigma_i(k),\vec z_i(k),\vec z_i^{\rm \,int}(k),\vec w_1^{ \rm int}(k)\bigr)$
converges to
$\bigl(\Sigma_i,\vec z_i,\vec z_i^{\rm \,int},\vec w_1^{ \rm int}\bigr)$
as $k\to \infty$
in the moduli space of bordered marked nodal curves, for $i=1,2$.
\item[(2)]
Let $\mathfrak N_i(\varepsilon)$ be the $\varepsilon$ neighborhood of the
set of the nodal points of $\Sigma_i$.
Using a universal family of nodal marked bordered curves
together with item~(1), we take a smooth embedding~${
\mathfrak I_{i,k} \colon \Sigma_i \setminus \mathfrak N_i(\varepsilon)
\to \Sigma_i(k)}$,
such that it converges to the identity map
as $k$ goes to infinity. (Here we regard $\Sigma_i(k)$
as a subset of the total space of the universal family.)
Moreover,
\begin{enumerate}\itemsep=0pt
\item
$\mathfrak I_{i,k}(z_{i,j}) = z_{i,j}(k)$.
\item
$\mathfrak I_{i,k}(z_{i,j}^{\rm int}) = z_{i,j}^{\rm int}(k)$.
\item
$\mathfrak I_{i,k}(w_{i,j}^{\rm int}) = w_{i,j}^{\rm int}(k)$.
\end{enumerate}
\item[(3)]
For any small $\varepsilon > 0$, we have
\[
\lim_{k\to\infty} \sup\{ d(u_i(k)(\mathfrak I_{i,k}(z)),u_i(z))
\mid z \in \Sigma_i \setminus \mathfrak N_i(\varepsilon)\}
= 0.
\]
\item[(4)]
There exist $\varepsilon_k \to 0$, $\delta_k \to 0$ such that
for each connected component
$S_{i,a}(k)$ of $\Sigma_i(k) \setminus \mathfrak I_{i,k}
(\Sigma_i \setminus \mathfrak N_i(\varepsilon_k))$
we have
$
\operatorname{Diam} S_{i,a}(k) \le \delta_k$.
\item[(5)]
We may choose $\mathfrak N_i(\varepsilon_k)$ and $\varepsilon_k \to 0$ with
the following property.
Let $\overline{\mathfrak N}_i(\varepsilon_k)$ be the image of $\mathfrak N_i(\varepsilon_k)$
in $\Sigma_i^0$. Then
\smash{$
\mathscr I\bigl(\overline{\mathfrak N}_1(\varepsilon_k) \cap \Sigma_1^0\bigr)
\subseteq
\overline{\mathfrak N}_2(\varepsilon_k) \cap \Sigma_2^0$}.
Now we require
\[
\lim_{k\to\infty} \sup\{ d(\mathscr I(k)(\mathfrak I_{1,k}(z)),\mathfrak I_{2,k}(\mathscr I(z)))
\mid z \in \Sigma_1 \setminus \mathfrak N_1(\varepsilon_k)\}
= 0.
\]
\end{enumerate}

\end{defn}
\begin{rem}
We require $C^0$ convergence in item (3). Since the maps are
pseudo-holomor\-phic, it implies $C^n$ convergence for any $n$.
\end{rem}
\begin{defn}\label{defn141111}
Let
$\xi, \xi(k) \in {\mathcal M}'_{\ell,\ell_1,\ell_2}(L_{12};\vec a;E)$.
We say
$
\lim_{k\to \infty} \xi(k) = \xi
$
if there exists $\ell'$, $\ell'_1$, $\ell'_2$, $\mathfrak i$ as in \eqref{form145555} and
\smash{$\xi^+, \xi(k)^+ \in {\mathcal M}'_{\ell',\ell'_1,\ell'_2}(L_{12};\vec a;E)$}
such that
\begin{enumerate}\itemsep=0pt
\item[(1)]
$\mathfrak i^* (\xi^+) = \xi$, $\mathfrak i^* (\xi(k)^+) = \xi(k)$.
Moreover, $\xi^+$, $\xi(k)^+$ are source stabilizations of
$\xi$, $\xi(k)$, respectively.
\item[(2)]
$\operatorname{lims}_{k\to\infty} \xi^+(k) = \xi^+$.
\end{enumerate}

\end{defn}

We will use the next lemma to show that Definition~\ref{defn141111}
determines a topology. (See Lemma~\ref{lem121515}.)
Lemma~\ref{lem1214140} is also used during the construction of the
Kuranishi structure, in Section~\ref{sec:kuradifcompex}.
\begin{lem}\label{lem1214140}
We consider
\begin{gather*}
\xi \in {\mathcal M}'_{\ell,\ell_1,\ell_2}(L_{12};\vec a;E) ,\qquad
 \xi^{(1)} \in {\mathcal M}'_{\ell^{(1)},\ell^{(1)}_1,\ell^{(1)}_2}(L_{12};\vec a;E),\\
\xi^{(2)} \in {\mathcal M}'_{\ell^{(2)},\ell^{(2)}_1,\ell^{(2)}_2}(L_{12};\vec a;E).
\end{gather*}
Suppose that
\smash{$\mathfrak i_{(1)}^*\xi^{(1)} = \mathfrak i_{(2)}^*\xi^{(2)} = \xi$} for
some forgetful maps \smash{$\mathfrak i_{(1)}$}, \smash{$\mathfrak i_{(2)}$}.
We assume that \smash{$\xi^{(1)}$}, \smash{$\xi^{(2)}$} are source stable.
Let \smash{$\xi^{(2)}(k)$} be a sequence of source stable objects such that
\smash{${\rm lims}_{k\to\infty}\xi^{(2)}(k) = \xi^{(2)}$}.
Then there exists a sequence of elements \smash{$\xi^{(1)}(k)$} which are source stable
and such that
\[
\mathfrak i_{(1)}^*\bigl(\xi^{(1)}(k)\bigr) = \mathfrak i_{(2)}^*\xi^{(2)}(k)
\]
and
$\operatorname{lims}_{k\to\infty}\xi^{(1)}(k) = \xi^{(1)}$.
\end{lem}

\begin{figure}[ht]\centering
\begin{tabular}{cc}
\begin{minipage}[t]{0.45\hsize}
\centering
\includegraphics[scale=0.4]{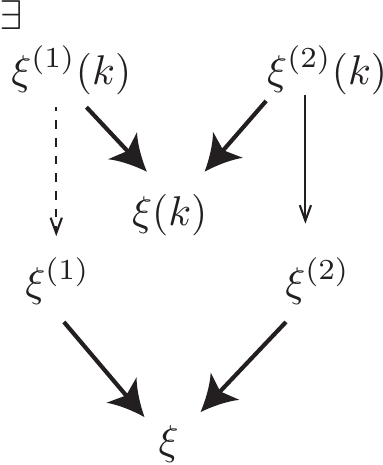}
\caption{$\xi^{(i)}$ and $\xi^{(i)}(k)$ in
Lemma~\ref{lem1214140}.}
\label{Figurelem12-14}
\end{minipage} &
\begin{minipage}[t]{0.45\hsize}
\centering
\includegraphics[scale=0.4]{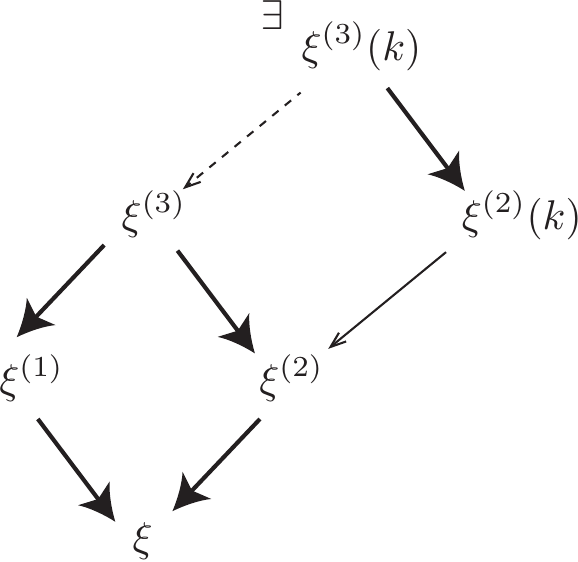}
\caption{$\xi^{(i)}$ and $\xi^{(i)}(k)$ in the claim.}
\label{FIgure12-14claim}
\end{minipage}
\end{tabular}
\end{figure}

\begin{proof}
We claim that there exist $\xi^{(3)}$, $\xi^{(3)}(k)$
such that
$
\operatorname{lims}_{k\to \infty}\xi^{(3)}(k) = \xi^{(3)}
$
and $\smash{\mathfrak i_{(32)}^*\xi^{(3)}} = \xi^{(2)}$,
\smash{$\mathfrak i_{(31)}^*\xi^{(3)} = \xi^{(1)}$},
\smash{$\mathfrak i_{(32)}^*\xi^{(3)}(k) = \xi^{(2)}(k)$}.
Here \smash{$\mathfrak i_{(32)}^*$}, \smash{$\mathfrak i_{(31)}^*$} are
appropriate forgetful maps.
(During the proof of this claim, we do not use the
assumption that $\xi^{(1)}$ is source stable.)

The proof of this claim is by an induction
on the number of added marked points of $\xi^{(1)}$.
We put
\[
\xi^{(i)} = \bigl(\bigl(\bigl(\Sigma^{(i)}_1,\vec z^{\,(i)}_1,\vec z_1^{\,(i),\rm int},\vec w_1^{ (i),\rm int}\bigr),u^{(i)}_1\bigr),\bigl(\bigl(\Sigma^{(i)}_2,\vec z^{\,(i)}_2
,\vec z_2^{\,(i),\rm int},\vec w_2^{ (i),\rm int}\bigr),u_2\bigr),\mathscr I^{(i)},\gamma^{(i)}\bigr)
\]
for $i=1,2$ and
\begin{gather*}
\xi^{(2)}(k) = \bigl(\bigl(\bigl(\Sigma^{(2)}_1(k),\vec z^{\,(2)}_1(k),
\vec z_1^{\,(2),\rm int}(k),\vec w_1^{ (2),\rm int}(k)),u^{(2)}_1(k)),
\\
\phantom{\xi^{(2)}(k) = }{}\bigl(\bigl(\Sigma^{(2)}_2(k),\vec z^{\,(i)}_2(k)
,\vec z_2^{\,(i),\rm int}(k),\vec w_2^{ (2),\rm int}(k)),u_2(k)\bigr),\mathscr I^{(2)}(k),\gamma^{(2)}(k)\bigr),
\\
\xi = \bigl(\bigl(\bigl(\Sigma_1,\vec z_1,\vec z_1^{\rm \,int}\bigr),u_1\bigr),
\bigl(\bigl(\Sigma_2,\vec z_2
,\vec z_2^{\rm \,int},\vec w_2^{ \rm int}\bigr),u_2\bigr),\mathscr I,\gamma\bigr).
\end{gather*}
Suppose the number of added marked point is one.
We consider the case when the added marked point is of type 2 and is \smash{$w^{(1),\rm int}_1$}
in $\Sigma_1$.
(The other cases are similar and so are omitted.)

Note that there are holomorphic maps \smash{$\pi^{(i)}_j \colon \Sigma^{(i)}_j \to \Sigma_j$},
which shrink certain irreducible components.

 Case 1: We assume that the irreducible
component containing \smash{$w^{(1),\rm int}_1$}
is not shrunk by $\pi^{(1)}_1 \colon \smash{\Sigma^{(1)}_1 \to \Sigma_1}$.

 Case 1-1:
Suppose \smash{$\pi^{(1)}_1\bigl(w^{(1),\rm int}_1\bigr)$} is not in the image of
a nodal or a marked point of
\smash{$\Sigma^{(2)}_1$}.~There exists a point $\hat w$ in \smash{$\Sigma^{(2)}_1$}
which goes to \smash{$\pi^{(1)}_1\bigl(w^{(1),\rm int}_1\bigr)$} by \smash{$\pi^{(2)}_1$}.
The point $\hat w$ is not nodal or marked.
We add $\hat w$ as an extra added marked point to \smash{$\Sigma^{(2)}_1$}
to obtain \smash{$\xi^{(3)}$}.
We then take one marked point~$\hat w(k)$ on \smash{$\Sigma^{(2)}_1(k)$} for each $k$,
which is `close' to $\hat w$ and
add $\hat w(k)$ to~\smash{$\xi^{(2)}(k)$} to obtain~\smash{$\xi^{(3)}(k)$}
such that $
{\rm lims}_{k\to \infty}\xi^{(3)}(k) \to \xi^{(3)}
$.
It is easy to see that $\xi^{(3)}$ and $\xi^{(3)}(k)$ have the required properties.

 Case 1-2:
Suppose \smash{$\pi^{(1)}_1\bigl(w^{(1),\rm int}_1\bigr)$} is the image of a marked point
$w'$ of
\smash{$\Sigma^{(2)}_1$}. We add a~sphere bubble~$S$ at $w'$ to \smash{$\Sigma^{(2)}_1$}
and add one marked point $\hat w$ on this bubble. (Then the sphere component~$S$ has one node and two
marked points. One of the two marked points
corresponds to $w'$ and the other is $\hat w$.)
We thus obtain $\xi^{(3)}$.
(See Figure~\ref{Figure12-1}.)
\begin{figure}[ht]
\centering
\includegraphics[scale=0.4]{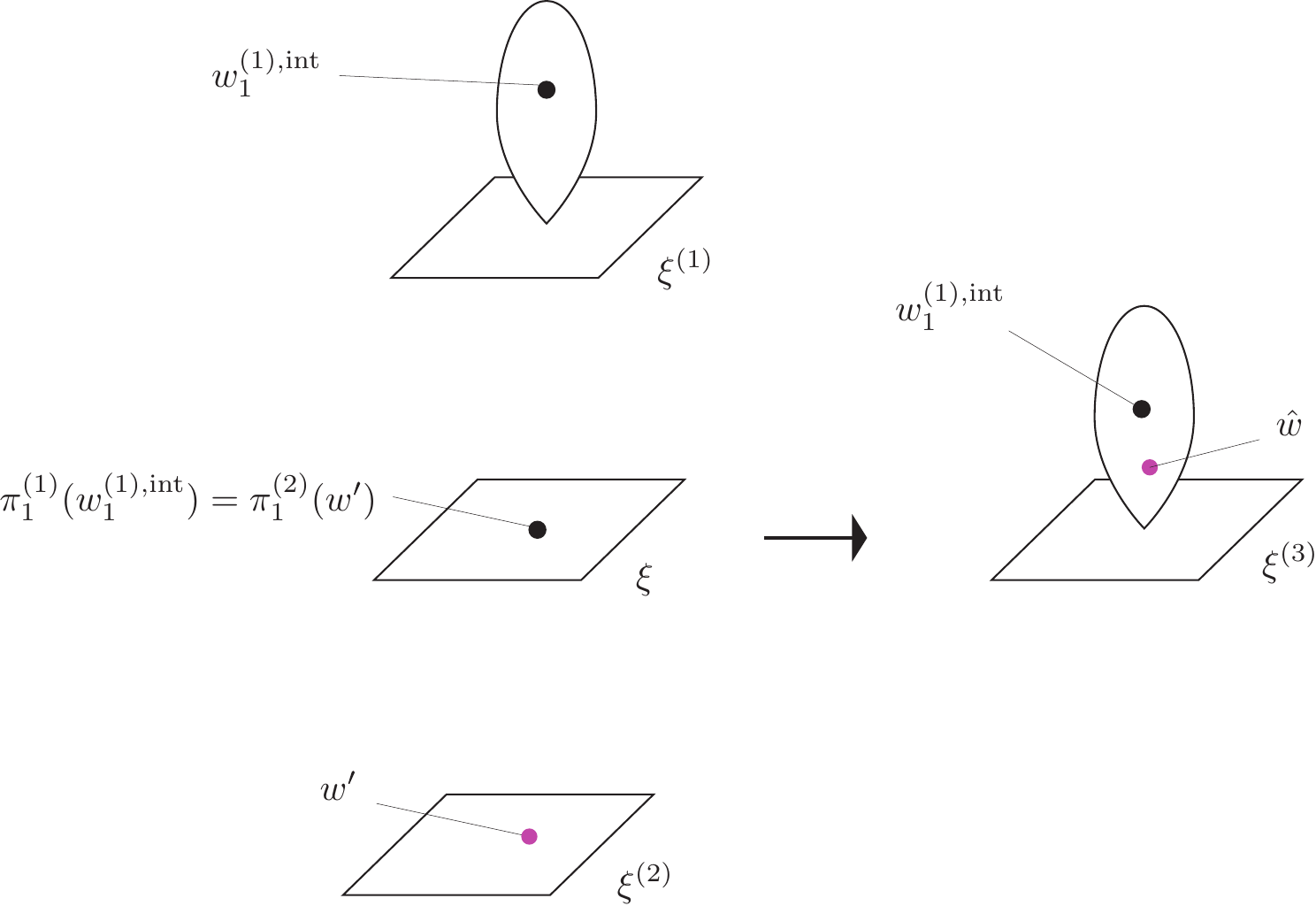}
\caption{Case 1--2.}
\label{Figure12-1}
\end{figure}

We consider \smash{$\Sigma^{(2)}_1(k)$}. We take the marked point $w'(k)$ corresponding to
$w'$. We add a sphere bubble $S(k)$ at $w'(k)$ and
a marked point $\hat w(k)$ on $S(k)$. We thus obtain $\xi^{(3)}(k)$
in the same way. It is easy to see that $\xi^{(3)}$, $\xi^{(3)}(k)$
have the required property.

 Case 1-3:
Suppose \smash{$\pi^{(1)}_1\bigl(w^{(1),\rm int}_1\bigr)$} is in the image of a node
$x$ of
\smash{$\Sigma^{(2)}_1$}.
We add a~sphere bubble at $x$ to \smash{$\Sigma^{(2)}_1$} and add one marked point on this bubble.
 (Then this sphere component
has two nodal points and one
marked point.)
We thus obtain $\xi^{(3)}$.
(See Figure~\ref{Figure12-2}.)
\begin{figure}[ht]
\centering
\includegraphics[scale=0.42]{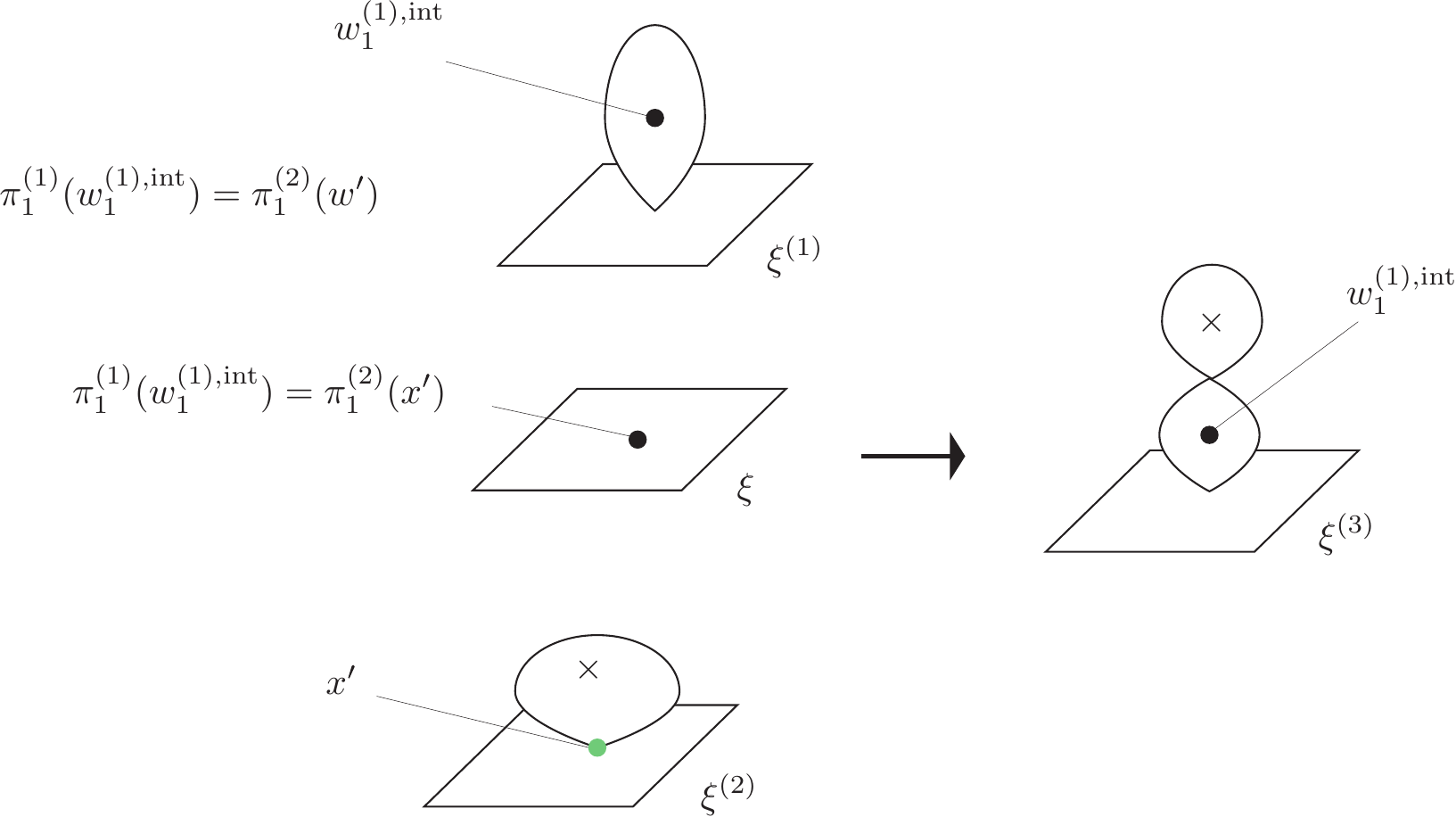}
\caption{Case 1--3.}
\label{Figure12-2}
\end{figure}

We consider \smash{$\Sigma^{(2)}_1(k)$}.
There are two cases. If there is a nodal point $x(k)$ corresponding to~$x$
in \smash{$\Sigma^{(2)}_1(k)$} then we add a~sphere bubble $S(k)$ at $x(k)$ and do the same construction
as above to obtain $\xi^{(3)}(k)$.
If there is no nodal point in \smash{$\Sigma^{(2)}_1(k)$} corresponding to $x$,
then there is a `neck region' corresponding to~$x$.
We add a marked point in this neck region to obtain
$\xi^{(3)}(k)$.
(See Figure~\ref{Figure12-3}.)
\begin{figure}[ht]
\centering
\includegraphics[scale=0.46]{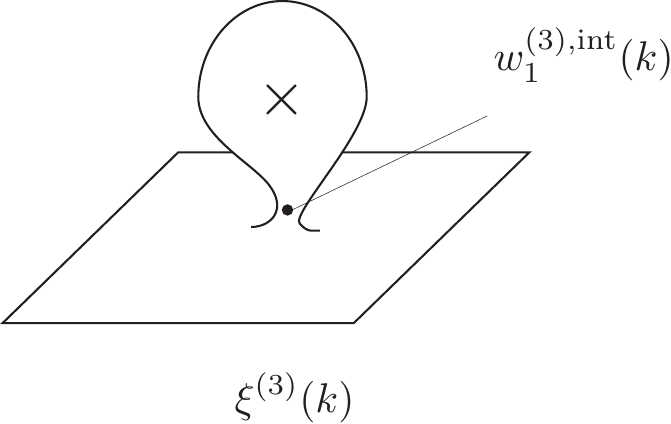}
\caption{Put a marked point on the neck region.}
\label{Figure12-3}
\end{figure}

It is easy to see that $\xi^{(3)}$, $\xi^{(3)}(k)$
have the required property.

We remark that Case 1-2 and Case 1-3 can occur at the same time.
Also the marked point~$w'$ in Case 1-2 or a node $x$ in Case 1-3
may not be unique. We can take any of such choices to prove
the claim in those cases.

 Case 2: We assume that the component containing $w^{(1),\rm int}_1$
is shrunk by \smash{$\pi^{(1)}_1 \colon \Sigma^{(1)}_1 \to \Sigma_1$}.
We~consider the point \smash{$\pi^{(1)}_1\bigl(w^{(1),\rm int}_1\bigr) \in \Sigma_1$}.
We~consider three subcases.

 Case 2-1:
\smash{$\pi^{(1)}_1\bigl(w^{(1),\rm int}_1\bigr) \in \Sigma_1$} is not
in the image of a nodal or a marked point of
\smash{$\Sigma^{(2)}_1$}.
The construction is the same as Case~1-1.

 Case 2-2:
\smash{$\pi^{(1)}_1\bigl(w^{(1),\rm int}_1\bigr) \in \Sigma_1$} is
in the image of a marked point of
\smash{$\Sigma^{(2)}_1$}.
The construction is the same as Case 1-2.

 Case 2-3:
\smash{$\pi^{(1)}_1\bigl(w^{(1),\rm int}_1\bigr) \in \Sigma_1$} is not
in the image of a nodal point of
\smash{$\Sigma^{(2)}_1$}.
The construction is the same as Case 1-3.

We thus proved the claim in the case when the number of
added marked points in $\xi^{(1)}$ is 1.

Now we prove the claim by the induction of
the number $n$ of added marked points in $\xi^{(1)}$.
(Such an induction is possible since during the proof of this claim we do not use the
assumption that $\xi^{(1)}$ is source stable.)

The case $n=1$ is already proved.
Suppose the claim is proved for $n-1$.
We remove one added marked point from $\xi^{(1)}$
and obtain $\xi^{(1),-}$.
We apply induction hypothesis to obtain~$\xi^{(3) -}$,~$\xi^{(3) -}(k)$.

Now we apply the case $n=1$ taking
$\xi^{(1)}$, $\xi^{(3) -}$, $\xi^{(3) -}(k)$
as $\xi^{(1)}$, $\xi^{(2)}$, $\xi^{(2)}(k)$.
It implies the claim in the case of $n$.

We have thus proved the claim.

We remark \smash{$\mathfrak i_{(31)}^*\bigl(\xi^{(3)}\bigr) = \xi^{(1)}$}.
Namely, $\xi^{(1)}$ is obtained by forgetting
certain marked points of $\xi^{(3)}$.
We forget the corresponding marked points of
$\xi^{(3)}(k)$ and obtain $\xi^{(1)}(k)$.
Since $\xi^{(1)}$ is source stable $\xi^{(1)}(k)$
is source stable for sufficiently large $k$.
Then
$
{\rm lims}_{k\to \infty}\xi^{(3)}(k) = \xi^{(3)}
$
implies~${
{\rm lims}_{k\to \infty}\xi^{(1)}(k) = \xi^{(1)}}$.
Since $\mathfrak i_{(32)}^*(\xi^{(3)}(k)) = \xi^{(2)}(k)$
we have $
\mathfrak i_{(1)}^*(\xi^{(1)}(k)) = \mathfrak i_{(2)}^*(\xi^{(2)}(k))$.
The proof of the lemma is complete.
\end{proof}

Note that we proved the next lemma also during the proof of
the claim in the proof of Lem\-ma~\ref{lem1214140}.
\begin{lem}\label{lemnew1215}
We consider
\begin{gather*}
\xi \in {\mathcal M}'_{\ell,\ell_1,\ell_2}(L_{12};\vec a;E) ,\qquad
\xi^{(1)} \in {\mathcal M}'_{\ell^{(1)},\ell^{(1)}_1,\ell^{(1)}_2}(L_{12};\vec a;E),\\
\xi^{(2)} \in {\mathcal M}'_{\ell^{(2)},\ell^{(2)}_1,\ell^{(2)}_2}(L_{12};\vec a;E).
\end{gather*}
Suppose that
\[
\mathfrak i_{(1)}^*\bigl(\xi^{(1)}\bigr) = \mathfrak i_{(2)}^*\bigl(\xi^{(2)}\bigr) = \xi
\]
 for
some $\mathfrak i_{(1)}$, $\mathfrak i_{(2)}$.
Then there exists
\smash{$\xi^{(3)} \in {\mathcal M}'_{\ell^{(1)},\ell^{(3)}_1,\ell^{(3)}_2}(L_{12};\vec a;E)$},
where \smash{$\ell^{(3)}_i = \ell^{(1)}_i + \ell^{(2)}_i - \ell_i$},
such that
\[
\mathfrak i_{(3,1)}^*\bigl(\xi^{(3)}\bigr) = \xi^{(1)},
\qquad
\mathfrak i_{(3,2)}^*\bigl(\xi^{(3)}\bigr) = \xi^{(2)}.
\]
Here $\mathfrak i_{(3,1)}^*$, $\mathfrak i_{(3,2)}^*$ are appropriate
forgetful maps.
\end{lem}

We now show that Definition~\ref{defn141111}
determines a topology on ${\mathcal M}'(L_{12};\vec a;E)$.
For a subset~${A \subset {\mathcal M}'(L_{12};\vec a;E)}$,
we define its closure $A^c$ as the set of all
elements $\xi$ such that there exists a~sequence $\xi(k) \in A$
which converges to $\xi$ in the sense of Definition~\ref{defn141111}.
Using Kuratowski's theorem (see, for example, \cite[Chapter~1, Theorem~8]{kelly}),
it suffices to show the next lemma to prove the
existence of the topology on ${\mathcal M}'(L_{12};\vec a;E)$
for which $A \mapsto A^c$ becomes the process taking the closure.

\begin{lem}\label{lem121515}
The following $4$ properties are satisfied:
{\rm (a)} $\varnothing^c = \varnothing$. {\rm (b)} $A \subseteq A^c$.
{\rm (c)} $A^{cc} = A^c$. {\rm (d)} $(A\cup B)^c = A^c \cup B^c$.
\end{lem}
\begin{proof}
(a), (b), (d) are trivial to check. We verify (c).
Let $\xi(i) \in A^c$ which converges to ${\xi \in A^{cc}}$.
We take $\xi(i,j) \in A$ such that $\lim_{j\to \infty} \xi(i,j) = \xi(i)$.
It suffices to find $j_i$ such that
$\lim_{i\to \infty} \xi(i,j_i) = \xi$.

Using Lemma~\ref{lem1214140},
we may assume that $\xi$, $\xi(i)$, $\xi(i,j)$ are all
source stable.
Let $\Sigma$, $\Sigma(i)$, $\Sigma(i,j)$ be the
source curves of $\xi$, $\xi(i)$, $\xi(i,j)$
and $u$, $u_{i}$, $u_{i,j}$ are maps on them, respectively.

Let $\varepsilon > 0$ be an arbitrary positive number.
We take sufficiently small neck of $\Sigma$ such that
the diameter of the image by $u$ of each of the neck is
smaller than $\varepsilon$. Let $\Sigma_0$ be the
complement of the neck.
We are given embedding of $\Sigma_0$
to $\Sigma(i)$ and to $\Sigma(i,j)$.

By Definition~\ref{defn1211new}, there exists
$I$ such that if $i\in I$ then
the diameter of each of the $u_i$ image of connected
component of $\Sigma(i) \setminus \Sigma_0$
is smaller than $2\varepsilon$.
Moreover, there exists $J_i$ such that if~${i > I}$, $j > J_i$, then the diameter
$u_{i,j}$ image of each of connected
component of $\Sigma(i,j) \setminus \Sigma_0$
is smaller than $3\varepsilon$.

By Definition~\ref{defn1211new} again, there exists
$I'$ such that if $i\in I'$ then
the $C^2$ distance between~$u_i\vert_{\Sigma_0}$ and
$u\vert_{\Sigma_0}$
is smaller than $\varepsilon$.
Moreover
there exists $I'_j$ such that
if $i > I'$, $j > J'_i$, then the $C^2$ distance between $u_{i,j}\vert_{\Sigma_0}$ and
$u\vert_{\Sigma_0}$
is smaller than $2\varepsilon$.

This implies that if $j_i > \max\{J_i,J'_i\}$
then $\xi(i,j_i)$ converges to $\xi$
in the sense of Definition~\ref{defn1211new}.
This proves~(c).
\end{proof}

In \eqref{def33314}, we defined a compactification ${\mathcal M}(L_{12};\vec a;E)$, whose element is
a bordered stable map with boundary marked points, switching specified by $\vec a$ and
with energy $E$.
We can include interior marked points and define ${\mathcal M}_{\ell}(L_{12};\vec a;E)$.
The way to include interior marked points is the same as
\cite[Definition 2.1.24]{fooobook} and so its detail is omitted.
\begin{lemdef}
We can define the forgetful map
\[
\mathfrak{fg} \colon\ {\mathcal M}_{\ell+\ell_1+\ell_2}(L_{12};\vec a;E)
\to
{\mathcal M}'_{\ell,\ell_1,\ell_2}(L_{12};\vec a;E),
\]
which is continuous.
\end{lemdef}
\begin{proof}
Let $\bigl(\bigl(\Sigma,\vec z,\vec z^{\,\rm int}\cup \vec w^{ \rm int}_1
\cup \vec w^{ \rm int}_2\bigr),u,\gamma\bigr)$ be an element of
${\mathcal M}_{\ell+\ell_1+\ell_2}\bigl(L_{12};\vec a;E,\gamma\bigr)$.
Here the object $\bigl(\Sigma,\vec z,\vec z^{\,\rm int}\cup \vec w^{ \rm int}_1
\cup \vec w^{ \rm int}_2\bigr)$ is a bordered nodal marked curve of
genus zero with one boundary component.
($\vec z$ are boundary marked points, $\vec z^{\,\rm int}$ are
first $\ell$ interior marked points,
$\vec w^{ \rm int}_1 = \bigl(w^{\rm int}_{1,1},\dots,w^{\rm int}_{1,\ell_1}\bigr)$
are next $\ell_1$ interior marked points and
$\vec w^{ \rm int}_2 = \bigl(w^{\rm int}_{2,1},\dots,w^{\rm int}_{2,\ell_2}\bigr)$
are last~$\ell_2$ interior marked points.)
The map
$u \colon (\Sigma,\partial \Sigma) \to (-X_1 \times X_2, L_{12})$
is pseudo-holomorphic and~${\gamma \colon \partial \Sigma \setminus \vec z \to \tilde L_{12}}$
is a lift of the restriction of $u$.

We put $u = (u_1,u_2)$, where $u_i$ is a map to $X_i$
from $\Sigma$.
We consider $\bigl(\bigl(\Sigma,\vec z,\vec z^{\,\rm int}\cup \vec w^{ \rm int}_i\bigr),u_i\bigr)$ for~${i=1,2}$.

We remark that, for $i=1$, we forget the marked points $\vec w^{ \rm int}_2$
and, for $i=2$, we forget the marked points $\vec w^{ \rm int}_1$.

We shrink unstable sphere components of $\bigl(\bigl(\Sigma,\vec z,\vec z^{\,\rm int}\cup \vec w^{ \rm int}_i\bigr),u_i\bigr)$.
Here an unstable sphere component of $\bigl(\bigl(\Sigma,\vec z,\vec z^{\,\rm int}\cup \vec w^{ \rm int}_i\bigr),u_i\bigr)$
is an unstable sphere component of the source curve
$\bigl(\Sigma,\vec z,\vec z^{\,\rm int}\cup \vec w^{ \rm int}_i\bigr)$ on which $u_i$ is constant.
We denote by
$\bigl(\bigl(\Sigma_i,\vec z_i,\vec z^{\,\rm int}_i\cup \vec w^{ \rm int}_i),u_i\bigr)$
the pair of a bordered marked curve and a map
obtained by this shrinking.

We next forget $\vec w^{ \rm int}_i$ and
let $\bigl(\Sigma^0_i,\vec z_i,\vec z^{\,\rm int}_i\bigr)$ be
the bordered marked curve obtained from
$\bigl(\Sigma_i,\vec z_i,\allowbreak\vec z^{\,\rm int}_i\bigr)$ by
shrinking all the unstable sphere components.

We remark that
$\bigl(\Sigma^0_1,\vec z_1,\vec z^{\,\rm int}_1\bigr)$
is canonically isomorphic to
$\bigl(\Sigma^0_2,\vec z_2,\vec z^{\,\rm int}_2\bigr)$.
In fact, they both are obtained by
shrinking all the unstable sphere components of
$(\Sigma,\vec z,\vec z^{\,\rm int}\bigr)$.
Therefore, we obtain a~biholomorphic map
$
\mathscr I \colon \bigl(\Sigma^0_1,\vec z_1,\vec z^{\,\rm int}_1\bigr)
\to
\bigl(\Sigma^0_2,\vec z_2,\vec z^{\,\rm int}_2\bigr)$.
We define
\begin{gather*}
\mathfrak{fg}((\Sigma,\vec z,\vec z^{\,\rm int}\cup \vec w^{ \rm int}_1
\cup \vec w^{ \rm int}_2\bigr),u,\gamma\bigr) \\
\qquad=
\bigl(\bigl(\bigl(\Sigma_1,\vec z_1,\vec z^{\,\rm int}_1
\cup \vec w^{ \rm int}_1),u_1\bigr),\bigl(\bigl(\Sigma_2,\vec z_2,\vec z^{\,\rm int}_2
\cup \vec w^{ \rm int}_2),u_2\bigr),\mathscr I,\gamma\bigr).
\end{gather*}
Note that we regard the interior marked points $\vec z_i^{\rm \,int}$ as interior
marked points of first kind and~$\vec w^{ \rm int}_i$
as interior
marked points of second kind.

The continuity of the map is easy to show from the definition.
\end{proof}

\begin{exm}
We consider the case when $L_{12}$ is embedded and $\vec a$ consists of one element
which corresponds to the diagonal component.
We define an element
\[
(((\Sigma_1,z_1),u_1),((\Sigma_2,z_2),u_2),\mathscr I,\gamma)
\]
of ${\mathcal M}'_{0,0,0}(L_{12};\vec a;E)$ as follows.

$\Sigma = \Sigma_1 = \Sigma_2$ is obtained by
gluing the disk $D^2$ with $S^2$ at
$0 \in D^2$ and $\infty \in S^2 \cong \C \cup \{\infty\}$.
We take $z_1=z_2 = 1 \in \partial D^2$ as the (boundary) marked point.
We take a holomorphic map
\[
u \colon\ (\Sigma,\partial \Sigma) \to (-X_1 \times X_2,L_{12}).
\]
We denote its restriction to $D^2$ by $u^{\rm d} = \bigl(u^{\rm d}_1,u^{\rm d}_2\bigr)$.
\big(Here $u^{\rm d}_i$ is a map to $X_i$.\big)
We denote its restriction to $S^2$ by $u^{\rm s} = (u^{\rm s}_1,u^{\rm s}_2)$.
$u_i \colon \Sigma_i \to X_i$ is a map which is $u^{\rm d}_i$ on $D^2$ and is
$u^{\rm s}_i$ on $S^2$.
Note $\Sigma_i^0$ (in the sense appearing in Definition~\ref{defn145555}\,(3))
is $D^2$ in this case.
Let $\mathscr I \colon \Sigma_1^0 \to \Sigma_2^0$ be the identity map.
We put $\gamma = u^{\rm s}\vert_{\partial \Sigma}$.
We thus obtain
$\xi = (((\Sigma_1,z_1),u_1),((\Sigma_2,z_2),u_2),\mathscr I,\gamma) \in {\mathcal M}'_{0}(L_{12};\vec a;E)$.
See Figure~\ref{Figure14-6}.
\begin{figure}[ht]
\centering
\includegraphics[scale=0.33]{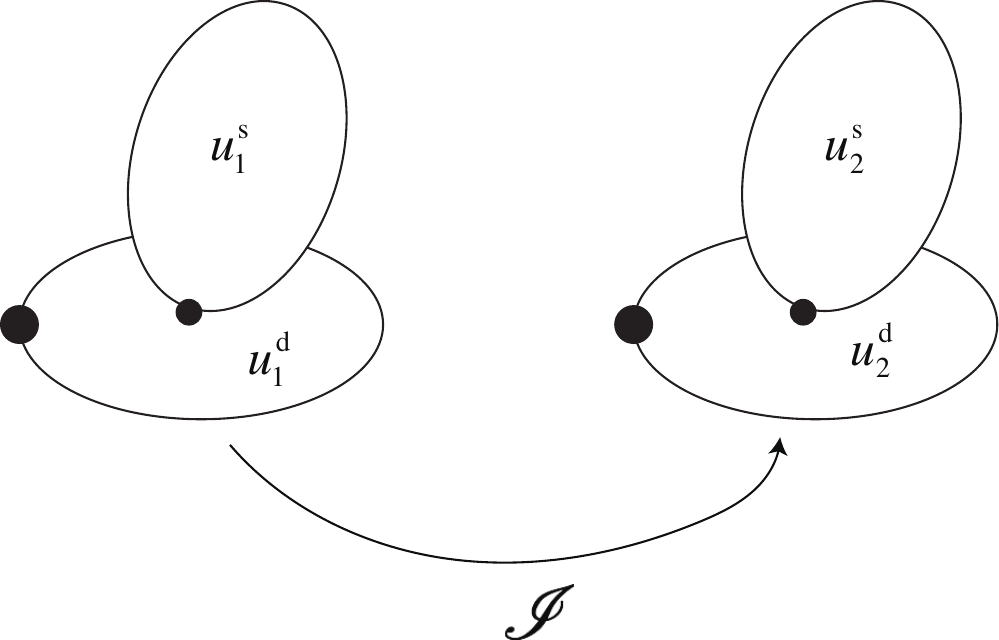}
\caption{$(((\Sigma_1,z_1),u_1),((\Sigma_2,z_2),u_2),\mathscr I,\gamma)$.}
\label{Figure14-6}
\end{figure}

We describe the fiber $\mathfrak{fg}^{-1}(\xi) \subset {\mathcal M}_{0}(L_{12};\vec a;E)$.
It is a real 4-dimensional compact space.
For~${a \in \C\setminus \{0\}}$ and $b\in \C$, we put
$
v_{a,b}(z) = az + b$,
and
\[
u^{\rm s}_{a,b} = (u^{\rm s}_1,u^{\rm s}_2 \circ v_{a,b}) \colon\ S^2 \to -X_1\times X_2.
\]
Since $\infty$ is a fixed point of $v_{a,b}$ we can glue it with $u^{\rm d}$
to obtain
$
u_{a,b} \colon (\Sigma,\partial \Sigma) \to (-X_1 \times X_2,L_{12})$.
Then
$
\xi_{a,b} = (((\Sigma,1),u_{a,b}),\gamma)
$
is an element of ${\mathcal M}_{0}(L_{12};\vec a;E)$ for any $a$, $b$ and
$
\mathfrak{fg}(\xi_{a,b}) = \xi$.
Those elements are parametrized by $(\C\setminus \{0\}) \times \C$ and consists
a non-compact space. See Figure~\ref{Figure14-92}.

The other elements of this fiber is described below in Figures \ref{Figure14-7}, \ref{Figure14-8}, \ref{Figure14-9}.

\begin{figure}[ht]\centering
\begin{tabular}{cc}
\begin{minipage}[t]{0.45\hsize}
\centering
\includegraphics[scale=0.48]{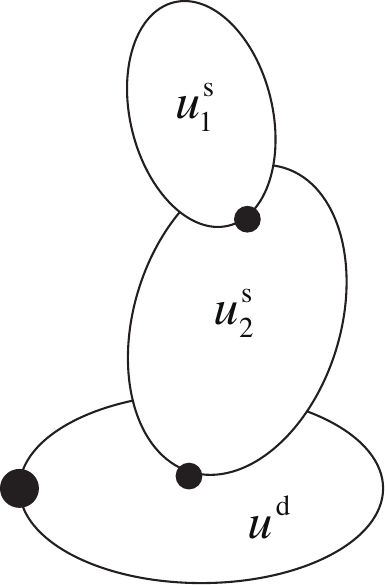}
\caption{First stratum.}
\label{Figure14-7}
\end{minipage} &
\begin{minipage}[t]{0.45\hsize}
\centering
\includegraphics[scale=0.48]{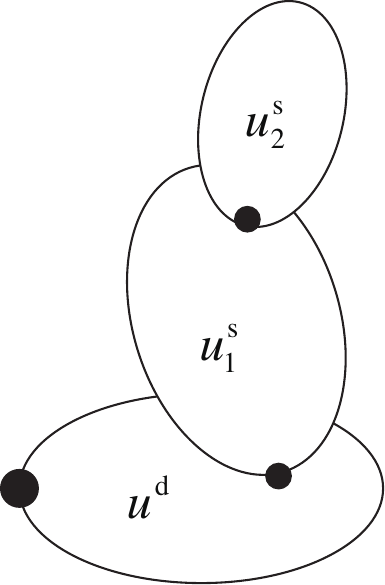}
\caption{Second stratum.}
\label{Figure14-8}
\end{minipage}
\end{tabular}
\end{figure}

\begin{figure}[ht]\centering
\begin{tabular}{cc}
\begin{minipage}[t]{0.45\hsize}
\centering
\includegraphics[scale=0.48]{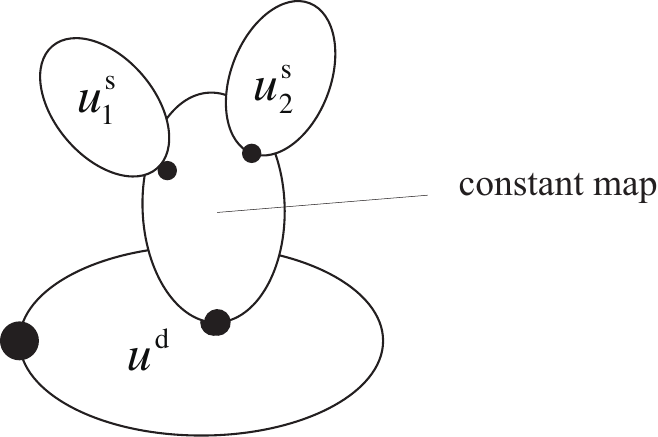}
\caption{Third stratum.}
\label{Figure14-9}
\end{minipage} &
\begin{minipage}[t]{0.45\hsize}
\centering
\includegraphics[scale=0.48]{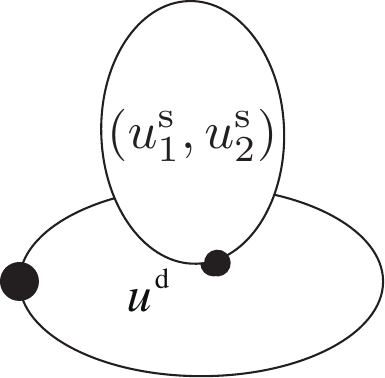}
\caption{Fourth stratum.}
\label{Figure14-92}
\end{minipage}
\end{tabular}
\end{figure}

Figure~\ref{Figure14-7} shows an element which has two sphere bubbles.
The map on the sphere component directly attached to a disk is constant in the
$X_1$ factor and the map on the other sphere component
is constant in the $X_2$ factor.
The element in the fiber $\mathfrak{fg}^{-1}(\xi)$ of the form Figure~\ref{Figure14-7}
is parametrized by the position of the nodal point
between two sphere components. So this part of the fiber is identified
with $\C$.

In Figure~\ref{Figure14-8}, the role of $X_1$ and $X_2$ is exchanged
from one in Figure~\ref{Figure14-7}. This part of fiber is also identified
with $\C$.

The closures of the parts Figures~\ref{Figure14-7} and~\ref{Figure14-8}
intersect at one point that is one depicted in Figure~\ref{Figure14-9}.
Here the map on one of the sphere components is a constant map.

The elements of the form depicted in Figures \ref{Figure14-7}, \ref{Figure14-8}, \ref{Figure14-9}
together with $\{\xi_{a,b} \mid (a,b) \in (\C\setminus \{0\}) \times \C\}$ consists
a compact 4-dimensional space, which is the fiber $\mathfrak{fg}^{-1}(\xi)$.
\end{exm}
\begin{prop}\label{prop141515}
The space
${\mathcal M}'_{\ell,\ell_1,\ell_2}(L_{12};\vec a;E)$ is compact and Hausdorff.

\end{prop}
\begin{proof}
The proof is similar to the proof of
\cite[Theorem 11.1]{FO} and \cite[Lemma 10.4]{FO}
and proceed as follows.

We first prove that the moduli space is sequentially compact.
Let $\xi_k$ be a sequence in
${\mathcal M}'_{\ell,\ell_1,\ell_2}(L_{12};\vec a;E)$.
We can add marked points to $\xi_k$ so that it becomes source stable.
Since the number of irreducible components of elements of
${\mathcal M}'_{\ell,\ell_1,\ell_2}(L_{12};\vec a;E)$
is bounded, we may assume that the number of marked points we add
to $\xi_k$ is independent of $k$.
Therefore, to prove the existence of convergent subsequence of
$\xi_k$ it suffices to assume that $\xi_k$ are source stable.
We assume so below.

Since the moduli space of stable marked curves is compact,
we may assume that the sequence of source (marked) curves of $\xi_k$ converges.
So using the local trivialization of the universal
family, we obtain a diffeomorphism between
source curves of $\xi_k$ and the limit, outside the neck region.
Therefore, the maps $u_{k,i}$, $i=1,2$, which is a part of $\xi_k$
can be regarded as a map $u_i$ from $\Sigma_i$, the limit curve.
If $u_{k,i}$ converges, there is nothing to show.

Suppose $u_{k,i}$ does not have a convergent subsequence. Then
the first derivative of $u_{k,i}$ diverges somewhere.

If it diverges on a disk component,
we can add two interior marked points of the first kinds
there in the same way as the proof of \cite[Theorem 11.1]{FO}
so that after we perform this replacement finitely many times
the sequence of maps $u_{k,i}$ does not diverge on the disk component.

Suppose $u_{k,i}$ diverges on a sphere component.
Then we can add two interior marked points of the second kind
around that point in the same way as the proof of \cite[Theorem 11.1]{FO}
so that after we perform this replacement finitely many times
the sequence of maps $u_{k,i}$ does not diverge on the sphere component either.

Thus by adding marked points the sequence of maps $u_{k,i}$ converges.
The proof of sequential compactness is complete.

We next prove the Hausdorffness.
It is easy to see from the definition and Lemma~\ref{lem1214140} that
${\mathcal M}'_{\ell,\ell_1,\ell_2}(L_{12};\vec a;E)$
satisfies the first axiom of countability.
Therefore, it suffices to show the following.
``For each sequence $\xi_k$ in ${\mathcal M}'_{\ell,\ell_1,\ell_2}(L_{12};\vec a;E)$
its limit is unique.''
We will prove it below.

Suppose $\lim_{k\to\infty} \xi_k = \xi$, $\lim_{k\to\infty} \xi_k = \xi'$.
By definition, there exists $\hat\xi_k$, $\hat\xi'_k$,
$\hat\xi$, $\hat\xi'$, such that
they are all source stable, \smash{$\mathfrak i^*(\hat\xi_k) = \xi_k$},
\smash{$\mathfrak i^*\bigl(\hat\xi'_k\bigr) = \xi_k$},
\smash{${\rm lims}_{k\to \infty} \hat\xi_k = \hat\xi$},
\smash{${\rm lims}_{k\to \infty} \hat\xi'_k = \hat\xi'$}
and \smash{$\mathfrak i^*\bigl(\hat\xi\bigr) = \xi$}, \smash{$\mathfrak i^*\bigl(\hat\xi'\bigr) = \xi'$}.
Here $\mathfrak i^*$ are forgetful maps.

By Lemma~\ref{lemnew1215},
we can find $\hat\xi''_k$ such that
\smash{$\mathfrak i_1^*\bigl(\hat\xi''_k\bigr) = \hat\xi_k$}, \smash{$\mathfrak i_2^*\bigl(\hat\xi''_k\bigr) = \hat\xi'_k$}
for certain forgetful maps~$\mathfrak i_1^*$ and $\mathfrak i_2^*$.

By taking a subsequence, we may assume that $\hat\xi''_k$
converges. Let $\hat\xi''$ be the limit.
Then by the continuity of forgetful map we have
$\mathfrak i_1^*\bigl(\hat\xi''\bigr) = \hat\xi$, $\mathfrak i_2^*\bigl(\hat\xi''\bigr) = \hat\xi'$.
$\xi = \xi'$ follows.

To complete the proof, it suffices to show the next lemma.
\begin{lem}\label{lem122000}
The space \smash{${\mathcal M}'_{\ell,\ell_1,\ell_2}(L_{12};\vec a;E)$}
satisfies the second axiom of countability.
\end{lem}
\begin{proof}
The proof is by induction on $E$.
In the case of smallest $E$ for which \smash{${\mathcal M}'_{\ell,\ell_1,\ell_2}(L_{12};\vec a;E)$}
is non-empty, we have
\smash{$
{\mathcal M}'_{\ell,\ell_1,\ell_2}(L_{12};\vec a;E)
=
\overset{\ \text{\tiny $\circ\circ$}}{\mathcal M}_{\ell,\ell_1,\ell_2}(L_{12};\vec a;E)$}.
It is easy to see that the right-hand side satisfies the second axiom of countability.

Suppose we have proved that \smash{${\mathcal M}'_{\ell',\ell'_1,\ell'_2}(L_{12};\vec a;E')$}
satisfies the second axiom of countabili\-ty for $E' < E$.
We consider the case of $E$.
Note that ${\mathcal M}'_{\ell,\ell_1,\ell_2}(L_{12};\vec a;E)$ has a stratification~${S_k{\mathcal M}'_{\ell,\ell_1,\ell_2}(L_{12};\vec a;E)}$
by its combinatorial types.
We will prove that $S_k{\mathcal M}'_{\ell,\ell_1,\ell_2}(L_{12};\vec a;E)$
satisfies the second axiom of countability by downward induction
on $k$.
For the stratum of smallest virtual dimension,
$S_k{\mathcal M}'_{\ell,\ell_1,\ell_2}(L_{12};\vec a;E)$
is a fiber product of various \smash{$\overset{\ \text{\tiny $\circ\circ$}}{\mathcal M}_{\ell',\ell'_1,\ell'_2}(L_{12};\vec a;E')$}
with~${E' \le E}$, and hence satisfies the second axiom of countability.
Suppose we have proved $S_{k+1}{\mathcal M}'_{\ell,\ell_1,\ell_2}(L_{12};\vec a;E)$
satisfies the second axiom of countability.
We will study the case of $k$.
As we will prove in Section~\ref{sec:glueglue} later, each point $p$ of $S_{k+1}{\mathcal M}'_{\ell,\ell_1,\ell_2}(L_{12};\vec a;E)$
has a Kuranishi neighborhood $(V_p,\mathcal E_p,s_p,\psi_p)$.
Therefore, $p$ has an open neighborhood $W_p$ in $S_k{\mathcal M}'_{\ell,\ell_1,\ell_2}(L_{12};\vec a;E)$
which satisfies the second axiom of countability.
In fact, $W_p$ is a closed subset of an orbifold.
Since $S_{k+1}{\mathcal M}'_{\ell,\ell_1,\ell_2}(L_{12};\vec a;E)$
satisfies the second axiom of countability by induction hypothesis
and since it is sequentially compact,
we can cover its open neighborhood by a finitely many $W_{p_i}$.
Note that
\[
S_k{\mathcal M}'_{\ell,\ell_1,\ell_2}(L_{12};\vec a;E)
\setminus \bigcup_i W_{p_i}
\]
is sequentially compact and is contained in a fiber product of
various \smash{$\overset{\ \text{\tiny $\circ\circ$}}{\mathcal M}_{\ell',\ell'_1,\ell'_2}(L_{12};\vec a;E')$}
for ${E' \le E}$. Therefore, it is contained
in an open subset which satisfies the second axiom of countability.

Thus $S_k{\mathcal M}'_{\ell,\ell_1,\ell_2}(L_{12};\vec a;E)$ is covered
by a finitely many open subsets
each of which satisfies the second axiom of countability.
This implies that $S_k{\mathcal M}'_{\ell,\ell_1,\ell_2}(L_{12};\vec a;E)$
satisfies the second axiom of countability.
The proof of Lemma~\ref{lem122000} is complete.
\end{proof}

The proof of Proposition~\ref{prop141515} is now complete.
\end{proof}

\subsection[Kuranishi structure of the compactification
${\mathcal M}'(L_{12};\vec a;E)$]{Kuranishi structure of the compactification
$\boldsymbol{{\mathcal M}'(L_{12};\vec a;E)}$}
\label{sec:kuradifcompex}

Let $\vec a = (a_0,\dots,a_k)$ and let
\smash{$
\xi = \bigl(\bigl(\bigl(\Sigma_1,\vec z_1,\vec z_1^{\rm \,int},\vec w_1^{ \rm int}\bigr),u_1\bigr),\bigl(\bigl(\Sigma_2,\vec z_2
,\vec z_2^{\rm \,int},\vec w_2^{ \rm int}\bigr),u_2\bigr),\mathscr I,\gamma\bigr)
$}
be an element of ${\mathcal M}'(L_{12};\vec a;E)$.
We define evaluation maps
\begin{gather}
{\rm ev} = \bigl({\rm ev}^{\partial},{\rm ev}^{{\rm int},(1)},{\rm ev}^{{\rm int},(2),1},{\rm ev}^{{\rm int},(2),2}\bigr)\colon\nonumber
\\
\hphantom{{\rm ev} =}{} {\mathcal M}'_{\ell,\ell_1,\ell_2}(L_{12};\vec a;E)
\to \prod_{j=0}^k L_{12}(a_j) \times (X_1 \times X_2)^{\ell} \times X_1^{\ell_1} \times X_2^{\ell_2},\label{newform1211}
\end{gather}
by
\begin{gather*}
{\rm ev}^{\partial}(\xi) := (\gamma_1(z_{1,0}),\dots,\gamma_1(z_{1,k})),\\
{\rm ev}^{{\rm int},(1)}(\xi) := \bigl(\bigl(u_1\bigl(z^{\rm int}_{1,1}\bigr),u_2\bigl(z^{\rm int}_{2,1}\bigr)\bigr),\dots,\bigl(u_1\bigl(z^{\rm int}_{1,\ell}\bigr),u_1\bigl(z^{\rm int}_{2,\ell}\bigr)\bigr)\bigr),
\\
{\rm ev}^{{\rm int},(2),i}(\xi) := \bigl(u_i\bigl(w^{\rm int}_{i,1}\bigr),\dots,u_i\bigl(w^{\rm int}_{i,\ell_i}\bigr)\bigr).
\end{gather*}
\begin{thm}\label{prop1417}
\smash{${\mathcal M}'_{\ell,\ell_1,\ell_2}(L_{12};\vec a;E)$} has a Kuranishi structure.
The evaluation map ${\rm ev}^{\partial}$
becomes an underlying continuous map of a strongly smooth map. The map, ${\rm ev}^{\partial}_0$,
the evaluation map at the $0$-th boundary marked point,
is weakly submersive.
They satisfy the same compatibility conditions as
Theorem {\rm\ref{thekuraexist}}.
\end{thm}
\begin{proof}
We prove the case of
${\mathcal M}'(L_{12};\vec a;E)
= {\mathcal M}'_{0,0,0}(L_{12};\vec a;E)$ below.
The general case is similar.
(We use the case ${\mathcal M}'(L_{12};\vec a;E)$ only
in this paper.) See Remark~\ref{re1234}.

Most of the proof is similar to the proof of Theorem~\ref{thekuraexist},
which was given in the reference quoted there.
We describe the place where the proof of
Theorem~\ref{prop1417} is different from
the proof of Theorem~\ref{thekuraexist}.
Especially we discuss the way how we include the maps $\mathscr I$, which is
a part of the data defining an element of ${\mathcal M}'(L_{12};\vec a;E)$
(see Definition~\ref{defn145555}\,(3)),
in the gluing analysis etc., which we use to construct a
Kuranishi chart.
The proof occupies this and the next subsections.

For this purpose, we review the construction of
the Kuranishi structure discussed in various literatures, explaining the places where the
construction here is to be modified.
Since the most detailed description of the gluing analysis
is given in \cite{foooanalysis}, we follow the
description of \cite[Section~8]{foooanalysis}.
(We follow \cite[Part 4]{foootech} on the discussion about stabilization
of the domain since that part is omitted in \cite{foooanalysis}.)

Let
$\xi = \bigl(\bigl(\bigl(\Sigma_1,\vec z_1,\vec z_1^{\rm \,int},\vec w_1^{ \rm int}\bigr),u_1\bigr),\bigl(\bigl(\Sigma_2,\vec z_2,\vec z_2^{\rm \,int},\vec w_2^{ \rm int}),u_2\bigr),\mathscr I,\gamma\bigr)$
be an element of the moduli space ${\mathcal M}'_{\ell,\ell_1,\ell_2}(L_{12};\vec a;E)$.
We first assume that it is source stable.

Let
$
\bigl\{\Sigma^{\rm d}_{i,\rm a} \mid {\rm a} \in \mathfrak{comp}^{\rm d}_i \bigr\}
$
(resp.\
$
\{\Sigma^{\rm s}_{i,\rm a} \mid {\rm a} \in \mathfrak{comp}^{\rm s}_i \}
$)
be the set of the disk (resp.\ sphere) components of $\Sigma_i$.
The (bordered) nodal curve \smash{$\Sigma^{\rm d}_{i,\rm a}$}
(resp.\ $\Sigma^{\rm s}_{i,\rm a}$)
together with marked or nodal points on it
determines an element of
\smash{${\mathcal M}^{\rm d}_{k_{i,\rm a},\ell_{i,\rm a}}$}
(resp.\ ${\mathcal M}^{\rm s}_{\ell_{i,\rm a}}$),
which we denote by \smash{$\xi^{\rm d}_{i,\rm a}$}
(resp.~$\xi^{\rm s}_{i,\rm a}$.)
Here~\smash{${\mathcal M}^{\rm d}_{k_{i,\rm a},\ell_{i,\rm a}}$} is the
moduli space of complex structures of disks with $k_{i,\rm a}$ boundary
and $\ell_{i,\rm a}$ interior marked points
and
${\mathcal M}^{\rm s}_{\ell_{i,\rm a}}$ is the
moduli space of complex structures of spheres with~$\ell_{i,\rm a}$ interior marked points.
(We require that the enumeration of the boundary marked points respects the orientation
of the boundary of the disk.)

Let ${\mathcal{CM}}^{\rm d}_{k,\ell}$,
${\mathcal{CM}}^{\rm s}_{\ell}$ be the Deligne--Mumford
type compactifications of ${\mathcal M}^{\rm d}_{k,\ell}$,
${\mathcal M}^{\rm s}_{\ell}$, respectively. Namely, we add stable nodal disks or spheres to compactify
them.
Let \smash{$\pi \colon {\mathcal C}^{\rm d}_{k,\ell} \to {\mathcal{CM}}^{\rm d}_{k,\ell}$}
be the universal family.
Namely, \smash{$\pi \colon {\mathcal C}^{\rm d}_{k,\ell} \to {\mathcal{CM}}^{\rm d}_{k,\ell}$}
comes with sections \smash{$\mathfrak s^{\rm d}_{j}$}, $j = 1,\dots,k$,
\smash{$\mathfrak s^{\rm s}_{j}$}, $j = 1,\dots,\ell$, such that for \smash{$\mathfrak x \in {\mathcal{CM}}^{\rm d}_{k,\ell}$}
the fiber $\pi^{-1}(\mathfrak x)$ together with the
marked points \smash{$\bigl(\bigl(\mathfrak s^{\rm d}_{j}(\mathfrak x)\bigr)_{j = 1,\dots,k}\bigr)$},
$(\mathfrak s^{\rm s}_{j}(\mathfrak x))_{j = 1,\dots,\ell}))$ becomes a representative of $\mathfrak x$.

Let $\pi \colon {\mathcal C}^{\rm s}_{\ell_{\rm a}} \to
{\mathcal{CM}}^{\rm s}_{\ell_{\rm a}}$ be the
universal family in a similar sense.
\begin{defn}[{compare \cite[Definition 16.2]{foootech} and \cite[Definition 8.6]{foooanalysis}}]\label{defn1222}
Suppose an element~$\xi$ of \smash{${\mathcal M}'_{\ell,\ell_1,\ell_2}(L_{12};\vec a;E)$} is source stable.
A {\it source gluing data} \index{source gluing data} $\mathscr{GL}$
\index[syindex]{GLscr@$\mathscr{GL}$} at $\xi$ is the following objects:
\begin{enumerate}\itemsep=0pt
\item[(1)]
A neighborhood $\mathcal V^{\rm d}_{i,\rm a}$
(resp.\ $\mathcal V^{\rm s}_{i,\rm a}$)
of $\xi^{\rm d}_{i,\rm a}$
(resp.\ $\xi^{\rm s}_{i,\rm a}$)
in ${\mathcal M}^{\rm d}_{k_{i,\rm a},\ell_{i,\rm a}}$
(resp.\ ${\mathcal M}^{\rm s}_{\ell_{i,\rm a}}$).
\item[(2)]
A trivialization of
$\pi \colon {\mathcal C}^{\rm d}_{k_{i,\rm a},\ell_{i,\rm a}} \to {\mathcal{CM}}^{\rm d}_{k_{i,\rm a},\ell_{i,\rm a}}$
(resp.\ $\pi \colon {\mathcal C}^{\rm s}_{\ell_{i,\rm a}} \to {\mathcal{CM}}^{\rm s}_{\ell_{i,\rm a}}$)
on $\mathcal V^{\rm d}_{i,\rm a}$
(resp.\ $\mathcal V^{\rm s}_{i,\rm a}$).
Here trivialization is one in $C^{\infty}$ category and
is required to be compatible with the sections~\smash{$((\mathfrak s^{\rm d}_{j})_{j = 1,\dots,k_{i,\rm a}})$},
\smash{$(\mathfrak s^{\rm s}_{j})_{j = 1,\dots,\ell_{i,\rm a}})$}.
(We remark that \smash{$\pi \colon {\mathcal C}^{\rm d}_{k_{i,\rm a},\ell_{i,\rm a}} \to {\mathcal{CM}}^{\rm d}_{k_{i,\rm a},\ell_{i,\rm a}}$}
is a~fiber bundle on~\smash{$\mathcal V^{\rm d}_{i,\rm a}$} since elements of
\smash{$\mathcal V^{\rm d}_{i,\rm a}$} are nonsingular.)
\item[(3)]
For each (boundary or interior) nodes of $\Sigma^{\rm d}_{i,\rm a}$ or $\Sigma^{\rm s}_{i,\rm a}$,
we take analytic families of coordinates~of the corresponding
marked points on \smash{$\mathcal V^{\rm d}_{i,\rm a}$} or $\mathcal V^{\rm s}_{i,\rm a}$.
(Note that one node is contained in two irreducible components.
We take an analytic family of coordinates at each of them.)
The notion of an analytic family of coordinates
is defined in
\cite[Definitions 8.1 and 8.5]{foooanalysis}.
\item[(4)]
The objects in (1), (2), (3) are preserved by all the weak isomorphisms
$(\psi_1,\psi_2) \colon \xi \to \xi$.
\end{enumerate}
The above conditions are mostly the same as one appearing in the
construction of Kuranishi structure on ${\mathcal M}(L_{12};\vec a;E)$,
for example. We need additional conditions to
include the map $\mathscr I$.
\begin{enumerate}\itemsep=0pt
\item[(5)]
All the interior marked points on the disk components are
of first kind.
All the marked points on the sphere components are
of second kind.
\item[(6)]
By (5) and Definition~\ref{defn145555}\,(3), for each of disk component
$\xi^{\rm d}_{1,a}$ of $\Sigma_1$ there exists corresponding
disk component of $\Sigma_2$, which we write \smash{$\xi^{\rm d}_{2,a}$}.
Namely, $\mathscr I$ gives an isomorphism between $\xi^{\rm d}_{1,a}$
and $\xi^{\rm d}_{2,a}$.
We require that $\mathcal V^{\rm d}_{1,a} = \mathcal V^{\rm d}_{2,a}$.
Moreover, we require the trivialization on \smash{$\mathcal V^{\rm d}_{1,a}$}
given in (2) is the same as the
trivialization on \smash{$\mathcal V^{\rm d}_{2,a}$}.
\item[(7)]
We require that the coordinate at nodal points given by (3) on
disk component $\xi^{\rm d}_{1,a}$ coincide with
those on \smash{$\xi^{\rm d}_{2,a}$}.
(We require this condition both for boundary and interior nodes.)
\item[(8)]
We will require all the analytic families of coordinates are extendable
in the sense we will define later in Definition~\ref{defn122812}.
\end{enumerate}

\end{defn}

We remark that for any element of
${\mathcal M}'(L_{12};\vec a;E)$
we can find its source stabilization such that
the conditions (5)--(8) are satisfied.

We next include the process to start with
$\xi \in {\mathcal M}'(L_{12};\vec a;E)$ which is
not necessary source stable and add marked points to
obtain an element of
\smash{${\mathcal M}'_{\ell,\ell_1,\ell_2}(L_{12};\vec a;E)$}
which is source stable.
\begin{defn}[{compare \cite[Definition 17.5]{foootech}}]\label{defn1223}
Let
\[
\xi = \bigl(\bigl(\bigl(\Sigma_1,\vec z_1,\vec z_1^{\rm \,int},\vec w_1^{ \rm int}\bigr),u_1\bigr),\bigl(\bigl(\Sigma_2,\vec z_2
,\vec z_2^{\rm \,int},\vec w_2^{ \rm int}\bigr),u_2\bigr),\mathscr I,\gamma\bigr)
\]
be an element of ${\mathcal M}'_{\ell,\ell_1,\ell_2}(L_{12};\vec a;E)$.

A {\it stabilization data} $\mathscr{ST}$ \index{stabilization data} \index[syindex]{STscr@$\mathscr{ST}$} at $\xi$ is the following objects:
\begin{enumerate}\itemsep=0pt
\item[(1)]
A source stabilization $\xi'$ of $\xi$ is given.
In particular, $\mathfrak i^*(\xi')= \xi$.
\item[(2)]
We require that the number of the irreducible components of
the source curve of $\xi'$ is the same as one of $\xi$.
\item[(3)]
A gluing data in the sense of Definition~\ref{defn1222} is given at $\xi'$.
\item[(4)]
We write
$\xi' = \bigl(\bigl(\bigl(\Sigma_1,\vec z_1,\vec z_1^{\rm \,int},\vec w_1^{ \rm int}\bigr),u_1\bigr),\bigl(\bigl(\Sigma_2,\vec z_2
,\vec z_2^{\rm \,int},\vec w_2^{ \rm int}\bigr),u_2\bigr),\mathscr I,\gamma\bigr)$.

Note that we use the same symbols $\Sigma_i$, $\vec z_i$, $u_i$, $\mathscr I$, $\gamma$
for $\xi'$ as $\xi$. In fact, item (2) implies that we can identify the source curves of
$\xi'$ and of $\xi$.
We do not put prime in the notation of interior marked points of $\xi'$.
Since $\xi$ has no interior marked points it does not cause confusion.
\item[(5)]
Let $z_{1,j}$ be an interior marked point of first kind,
which is necessary on the disk component by item (2) and Definition~\ref{defn1222}\,(5).
Suppose it is contained in \smash{$\Sigma^{\rm d}_{1,a_j}$}.
We put \smash{$\mathscr I(z_{1,j}) = z_{2,j}$} and
$\mathscr I(\Sigma^{\rm d}_{1,a_j}) = \Sigma^{\rm d}_{2,a_j}$.
We define $u_{j}^{\rm d} \colon \Sigma^{\rm d}_{1,a_j} \to -X_1 \times X_2$
by
$
u_{j}^{\rm d}(z) = (u_1(z),u_2(\mathscr I(z)))$.
\begin{enumerate}\itemsep=0pt
\item
If $u_{j}^{\rm d}$ is non-constant, we require that $u_{j}^{\rm d}$ is an
immersion at $z_{1,j}$.
\item
In the situation of (a), we take and fix a codimension $2$ submanifold
\smash{$\mathcal N^{(1)}_j$} of $-X_1 \times X_2$ which intersects transversally
with $u_{j}^{\rm d}$ at $u_{j}^{\rm d}(z_{1,j})$.
\end{enumerate}
\item[(6)]
Let $w_{i,j}$ be an interior marked point of second kind,
which is necessary on the sphere component by item (2) and Definition~\ref{defn1222}\,(5).
Suppose it is contained in $\Sigma^{\rm s}_{i,a_j}$.
\begin{enumerate}\itemsep=0pt
\item
If $u_{i}$ is non-constant on \smash{$\Sigma^{\rm s}_{i,a_j}$}, we require that $u_{i}$ is an
immersion at $w_{1,j}$.
\item
In the situation of (a), we take and fix a codimension $2$ submanifold
\smash{$\mathcal N^{(2)}_j$} of $X_i$ which intersects transversally
with $u_{i}$ at $u_{i}(w_{i,j})$.
\end{enumerate}
\item[(7)]
The data in item (6) are invariant under the action of the group of weak
isomorphisms of~$\xi$.
(Note that a weak isomorphism is the identity map on the disk components.)
\end{enumerate}

\end{defn}
It is easy to see that stabilization data always exist.
\begin{rem}
We need to add marked points of {\it second} kinds
to stabilize the source curve.
\end{rem}
We next describe the way how we use gluing data
to parametrize the deformation of the source objects.

Let $\xi$ be a source stable element of ${\mathcal M}'_{\ell,\ell_1,\ell_2}(L_{12};\vec a;E)$.
We take its gluing data as in Definition~\ref{defn1222}
and use the notation of Definition~\ref{defn1222}.\index[syindex]{Nodepluspar@${\rm Node}^+_{\partial}$}

Let
\smash{$\bigl\{ \bigl(\mathfrak z_{\rm b,\partial},\varphi^{(j)}_{\rm b,\partial}\bigr) \mid
{\rm b} \in {\rm Node}^+_{\partial},\, j=1,2\bigr\}
$}
be the set of pairs of boundary nodes and analytic families of coordinates
at those points. (Since each boundary node is contained in two irreducible
components, there are two choices $j=1,2$ of this pair for each
boundary node.)\index[syindex]{Nodeplusi@${\rm Node}^+_{i,\rm int}$}

Let
\smash{$
\bigl\{ \bigl(\mathfrak z_{i,\rm b,int},\varphi^{(j)}_{i,\rm b,int}\bigr) \mid
{\rm b} \in {\rm Node}^+_{i,\rm int}, \, j=1,2\bigr\}
$}
be the set of pairs of interior nodes of $\Sigma_i$ and analytic families of coordinates
at those points ($i=1,2$).

We denote by
$
\bigl\{\Sigma^{\rm d}_{i,\rm a} \mid {\rm a} \in \mathfrak{comp}^{\rm d}_i \bigr\}
$
(resp.\
$
\{\Sigma^{\rm s}_{i,\rm a} \mid {\rm a} \in \mathfrak{comp}^{\rm s}_i \}
$)
the set of disk (resp.\ sphere) components of $\Sigma_i$.
Together with nodal or marked points
the (bordered) Riemann surfaces~\smash{$\Sigma^{\rm d}_{i,\rm a}$},~$\Sigma^{\rm s}_{i,\rm a}$
determine~\smash{$\xi^{\rm d}_{i,\rm a}$},
$\xi^{\rm s}_{i,\rm a}$.
Its neighborhood $\mathcal V\bigl(\xi^{\rm d}_{i,\rm a}\bigr)$
and $\mathcal V\bigl(\xi^{\rm d}_{i,\rm s}\bigr)$ in Deligne--Mumford type
moduli spaces
are determined by Definition~\ref{defn1222}\,(1).

We consider the direct product
\begin{equation}\label{form1212}
\prod_{{\rm a} \in \mathfrak{comp}^{\rm d}_1} \mathcal V(\xi^{\rm d}_{1,\rm a})
\times
\prod_{i=1,2}\prod_{{\rm a} \in \mathfrak{comp}^{\rm s}_i} \mathcal V(\xi^{\rm s}_{i,\rm a})
\times
\prod_{{\rm b} \in {\rm Node}^+_{\partial}} [0,1)_{\rm b}
\times
\prod_{i=1,2}\prod_{{\rm b} \in {\rm Node}^+_{i,\rm int}} D^2_{\rm b}.
\end{equation}
Here $[0,1)_{\rm b}$ is a copy of $[0,1)$ taken for each
${\rm b} \in {\rm Node}^+_{\partial}$ and
$D^2_{\rm b}$ is a copy of $D^2 = \{z \in \C \mid \vert z\vert < 1\}$
taken for each ${\rm b} \in {\rm Node}^+_{\rm int}$.

The space \eqref{form1212} parametrizes the deformation of the source curve of
$\xi$.
We will define a~map~${{\rm Glue} = ({\rm Glue}_1,{\rm Glue}_2)}$ \index[syindex]{Glue@${\rm Glue}$}
\begin{equation}\label{form12113}
{\rm Glue}_i \colon\ \eqref{form1212}
\to \mathcal{CM}_{k,\ell+\ell_i}^{\rm d}
\end{equation}
to describe it.

Let $\sigma^{\rm d}_{\rm a} \in \mathcal V\bigl(\xi^{\rm d}_{1,\rm a}\bigr)
= \mathcal V(\xi^{\rm d}_{2,\rm a})$,
$\sigma^{\rm s}_{i,\rm a} \in \mathcal V(\xi^{\rm s}_{i,\rm a})$
and let
$
\xi^{\rm d}_{i,\rm a}\bigl(\sigma^{\rm d}_{\rm a}\bigr),
$
$\xi^{\rm s}_{i,\rm a}(\sigma^{\rm s}_{i,\rm a})$
be its representative. We denote by
\smash{$
\Sigma^{\rm d}_{i,\rm a}\bigl(\sigma^{\rm d}_{\rm a}\bigr)$},
$\Sigma^{\rm s}_{i,\rm a}(\sigma^{\rm s}_{i,\rm a})$
the underlying (bordered) Riemann surface.
(Actually it is either a disk or a sphere.)

We also denote
\begin{equation}\label{form1214}
\sigma
= \bigl(\bigl(\sigma^{\rm d}_{\rm a}\bigr)_{{\rm a}
\in \mathfrak{comp}^{\rm d}_1},(\sigma^{\rm s}_{1,\rm a})_{
{\rm a}
\in \mathfrak{comp}^{\rm s}_1},
(\sigma^{\rm s}_{2,\rm a})_{{\rm a}
\in \mathfrak{comp}^{\rm s}_2}\bigr).
\end{equation}
We call $\sigma$ the {\it source deformation parameter}
\index{source deformation parameter}.

Let $r_{\rm b} \in [0,1)_{\rm b}$ and
$\mathfrak r_{\rm b} \in D^2_{\rm b}$.
We write
\begin{equation}\label{form1215}
{\bf r} = ((r_{\rm b})_{
{\rm b} \in {\rm Node}^+_{\partial}},(\mathfrak r_{\rm b})_{{\rm b} \in {\rm Node}^+_{1,\rm int}},
(\mathfrak r_{\rm b})_{{\rm b} \in {\rm Node}^+_{2,\rm int}}).
\end{equation}
We call $\bf r$ the {\it gluing parameter}.
\index{gluing parameter}

We consider the disjoint union
\[
\hat \Sigma(\sigma) =\hat \Sigma_1(\sigma) \sqcup \hat \Sigma_2(\sigma)
=
\coprod_{i=1,2}\coprod_{a \in \mathfrak{comp}^{\rm d}_i}
\Sigma^{\rm d}_{i,\rm a}\bigl(\sigma^{\rm d}_{\rm a}\bigr)
\sqcup
\coprod_{i=1,2}\coprod_{a \in \mathfrak{comp}^{\rm s}_i}
\Sigma^{\rm s}_{i,\rm s}(\sigma^{\rm s}_{i,\rm a}).
\]

For each ${\rm b} \in {\rm Node}^+_{\partial}$
and ${\rm b} \in {\rm Node}^+_{i,\rm int}$,
the analytic families of coordinates we have taken in
Definition~\ref{defn1222}\,(3) induce holomorphic embeddings
\smash{$
\varphi^{(j),\partial}_{i,{\rm b},\sigma} \colon D^2_{\ge 0} \to \hat \Sigma(\sigma)
$}, \smash{$ \varphi^{(j),\rm int}_{{\rm b},\sigma} \colon D^2 \to \hat \Sigma(\sigma)$},
for $j=1,2$, $i=1,2$,
where \smash{$D^2_{\ge 0} = \bigl\{z \in D^2 \mid \operatorname{Im} z \ge 0\bigr\}$}.\index[syindex]{D20@$D^2_{\ge 0}$}
We put
\begin{align*}
\hat \Sigma(\sigma,{\bf r})
=
\hat \Sigma(\sigma)
\setminus
\bigcup_{j=1,2}\bigcup_{i=1,2}\bigcup_{{\rm b} \in {\rm Node}^+_{\partial}}
\varphi^{(j),\partial}_{i,{\rm b},\sigma}\bigl(\overline D_{\ge 0}^2(r_{\rm b})\bigr)
\setminus \bigcup_{j=1,2}\bigcup_{i=1,2}\bigcup_{{\rm b} \in {\rm Node}^+_{i,\rm int}}
\varphi^{(j),\rm int}_{{\rm b},\sigma}\bigl(\overline D^2(\vert \mathfrak r_{\rm b}\vert)\bigr),
\end{align*}
which we decompose to $\hat \Sigma(\sigma,{\bf r})
= \hat \Sigma_1(\sigma,{\bf r}) \sqcup \hat \Sigma_2(\sigma,{\bf r})$.
\begin{defn}\label{defn1225}
We define an equivalence relation $\sim$ on $\hat \Sigma(\sigma,{\bf r})$
as follows:
\begin{enumerate}\itemsep=0pt
\item[(1)]
If ${\rm b} \in {\rm Node}^+_{\partial}$ and \smash{$z,w \in
D^2_{\ge 0} \setminus \overline D_{\ge 0}^2(r_{\rm b})$},
$i=1,2$, with
$\vert zw \vert = r_{\rm b}$, $ \operatorname{Arg} z = -\operatorname{Arg} w$,
then
\smash{$
\varphi^{(1),\partial}_{i,{\rm b},\sigma}(z)
\sim
\varphi^{(2),\partial}_{i,{\rm b},\sigma}(w)
$}
for $i=1,2$. See Figure~\ref{Figure149-1}.
Note that $-\theta$ in the figure is $\operatorname{Arg} z$ and~$\theta'$ in the figure is $\operatorname{Arg} w$.
\item[(2)]
If \smash{${\rm b} \in {\rm Node}^+_{i,\rm int}$}, \smash{$z,w \in
D^2 \setminus \overline D^2(\vert \mathfrak r_{\rm b}\vert)$},
$i=1,2$, with
$
zw = \mathfrak r_{\rm b}$,
then
\smash{$
\varphi^{(1),{\rm int}}_{{\rm b},\sigma}(z)
\sim
\varphi^{(2),{\rm int}}_{{\rm b},\sigma}(w)$}.
See Figure~\ref{Figure149-2}.
\end{enumerate}
We put $\Sigma(\sigma,{\bf r}) = \hat \Sigma(\sigma,{\bf r})/\sim$
and decompose
$\Sigma(\sigma,{\bf r})
= \Sigma_1(\sigma,{\bf r}) \sqcup \Sigma_2(\sigma,{\bf r})$.

\end{defn}
\begin{figure}[ht]
\centering
\includegraphics[scale=0.48]{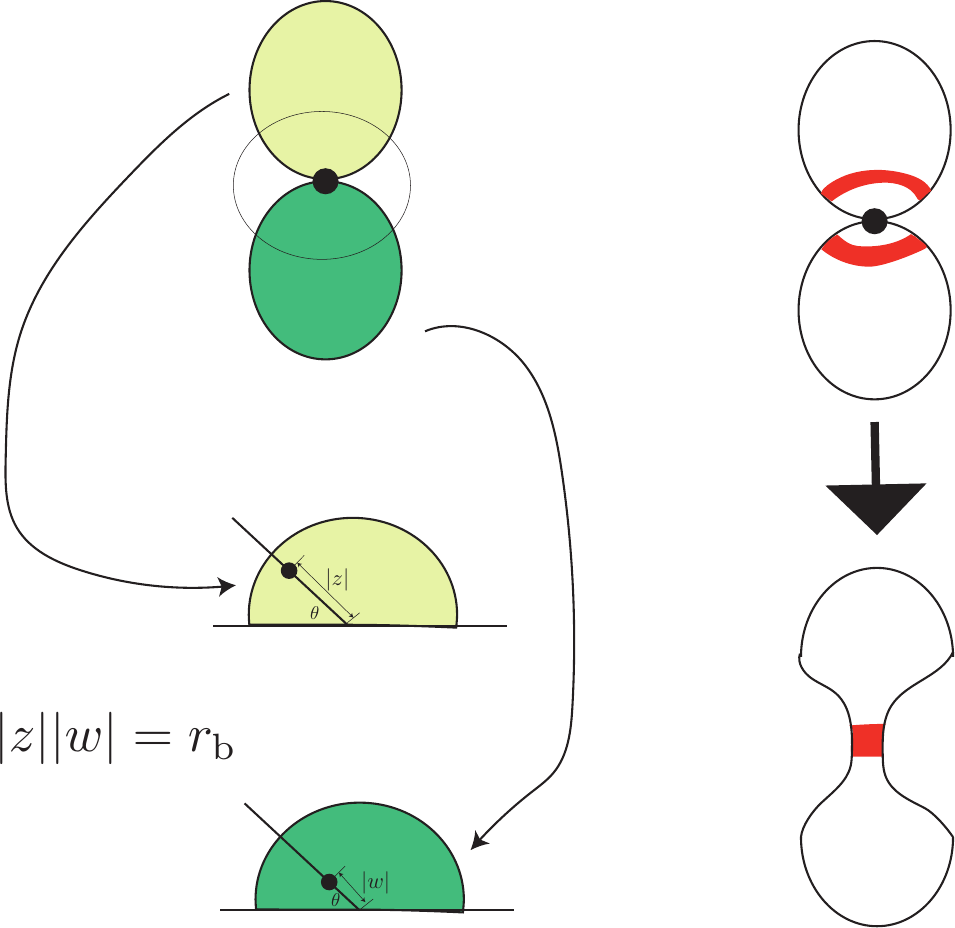}
\caption{Gluing at boundary node.}
\label{Figure149-1}
\end{figure}
\begin{figure}[ht]
\centering
\includegraphics[scale=0.36]{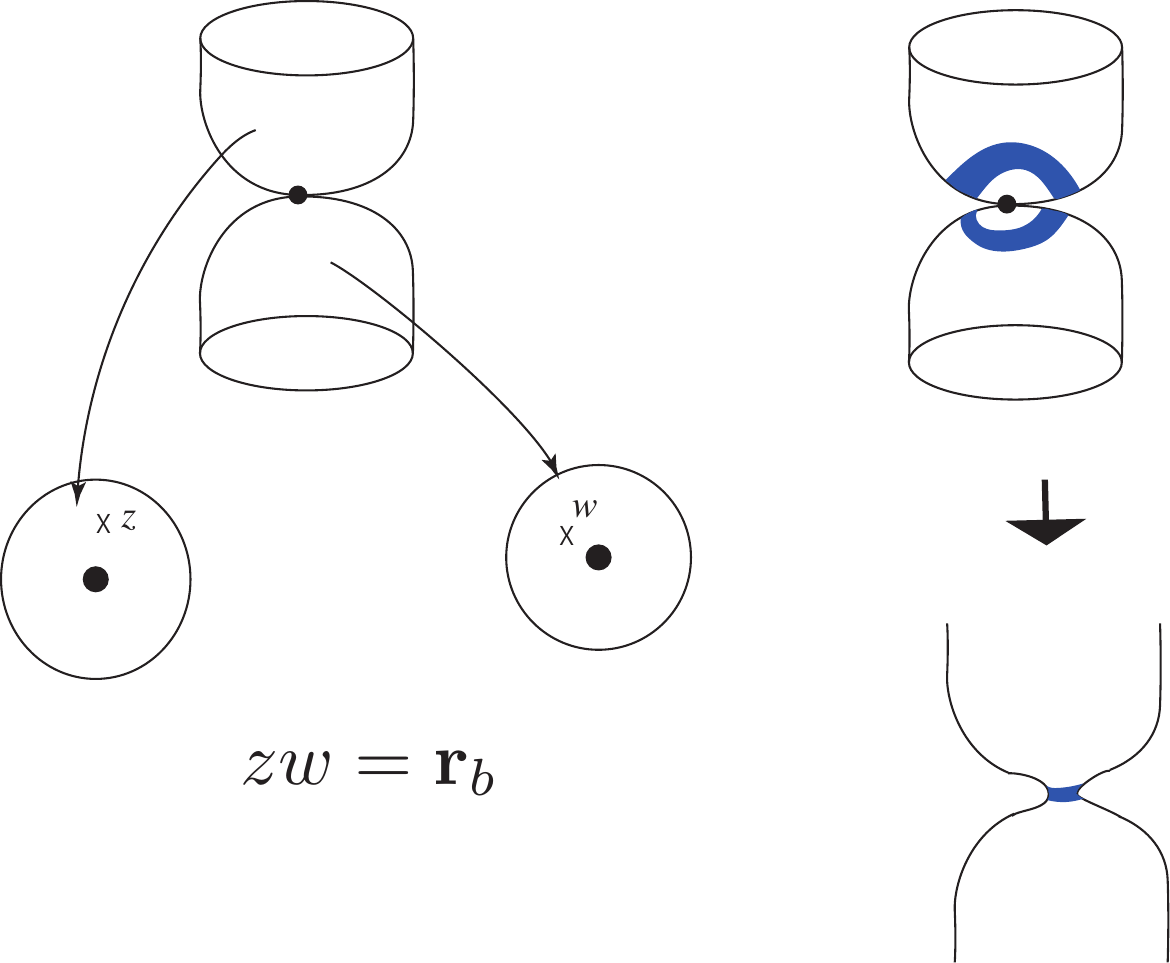}
\caption{Gluing at interior node.}
\label{Figure149-2}
\end{figure}
The marked points of $\xi$ determine marked points on
$\Sigma_i(\sigma,{\bf r})$ in an obvious way.
We denote them by
$\vec z_i(\sigma,{\bf r})$, \smash{$\vec z_i^{\rm \,int}(\sigma,{\bf r})$}, \smash{$\vec w
_i^{ \rm int}(\sigma,{\bf r})$}.
We put
$
\xi_i(\sigma,{\bf r}) = \bigl(\Sigma_i(\sigma,{\bf r}),
\vec z_i(\sigma,{\bf r}),\vec z_i^{\rm \,int}(\sigma,{\bf r}),\vec w_i^{ \rm int}(\sigma,{\bf r})\bigr)$.
\begin{defn}
We define
$
{\rm Glue}_i(\sigma,{\bf r}) = \xi_i(\sigma,{\bf r})$.
We call ${\rm Glue}_i$ and ${\rm Glue} := ({\rm Glue}_1,{\rm Glue}_2)$ the {\it source gluing maps}.
\index{source gluing map}
\end{defn}
For ${\rm a}
\in \mathfrak{comp}^{\rm d}_1
= \mathfrak{comp}^{\rm d}_2$, we put\index[syindex]{Kplusd@$K^{+, \rm d}_{i,\rm a}\bigl(\sigma^{\rm d}_{\rm a}\bigr)$}
\[
K^{+, \rm d}_{i,\rm a}\bigl(\sigma^{\rm d}_{\rm a}\bigr)
=
\Sigma^{\rm d}_{i,\rm a}\bigl(\sigma^{\rm d}_{\rm a}\bigr)
\setminus
\bigcup_{j=1,2}\bigcup_{{\rm b} \in {\rm Node}^+_{\partial}}
\varphi^{(j),\partial}_{i,{\rm b},\sigma}\bigl(\overline D_{\ge 0}^2(r_{\rm b})\bigr)
\setminus \bigcup_{j=1,2}\bigcup_{{\rm b} \in {\rm Node}^+_{i,\rm int}}
\varphi^{(j),\rm int}_{{\rm b},\sigma}\bigl(\overline D^2(\vert \mathfrak r_{\rm b}\vert)\bigr),
\]
and \index[syindex]{Kdia@$K^{\rm d}_{i,\rm a}\bigl(\sigma^{\rm d}_{\rm a}\bigr)$}
\begin{gather}\label{formdefcore}
K^{\rm d}_{i,\rm a}\bigl(\sigma^{\rm d}_{\rm a}\bigr)
=
\Sigma^{\rm d}_{i,\rm a}\bigl(\sigma^{\rm d}_{\rm a}\bigr)
\setminus
\bigcup_{j=1,2}\bigcup_{{\rm b} \in {\rm Node}^+_{\partial}}
\varphi^{(j),\partial}_{i,{\rm b},\sigma}\bigl(\overline D_{\ge 0}^2\bigr)
\setminus \bigcup_{j=1,2}\bigcup_{{\rm b} \in {\rm Node}^+_{i,\rm int}}
\varphi^{(j),\rm int}_{{\rm b},\sigma}\bigl(\overline D^2\bigr).
\end{gather}
For ${\rm a}
\in \mathfrak{comp}^{\rm s}_i$
we define \index[syindex]{Kplussi@$K^{+,s}_{i,\rm a}(\sigma^{\rm s}_{\rm a})$}
$K^{+,s}_{i,\rm a}(\sigma^{\rm s}_{\rm a})$
and $K^{\rm s}_{i,\rm a}(\sigma^{s}_{\rm a})$\index[syindex]{Ksia@$K^{\rm s}_{i,\rm a}(\sigma^{s}_{\rm a})$}
with
$
K^{\rm s}_{i,\rm a}(\sigma^{\rm s}_{\rm a})
\subset
K^{+, \rm s}_{i,\rm a}(\sigma^{\rm s}_{\rm a})
\subset
\Sigma^{\rm s}_{i,\rm a}(\sigma^{\rm s}_{\rm a})
$
in the same way.
We call $K^{\rm s}_{i,\rm a}(\sigma^{\rm s}_{\rm a})$ and
\smash{$K^{\rm d}_{i,\rm a}\bigl(\sigma^{\rm d}_{\rm a}\bigr)$} the {\it core}.
See Figures \ref{Figure159-3} and \ref{Figure159-4}.
\index{core}
\begin{figure}[ht]
\centering
\includegraphics[scale=0.3]{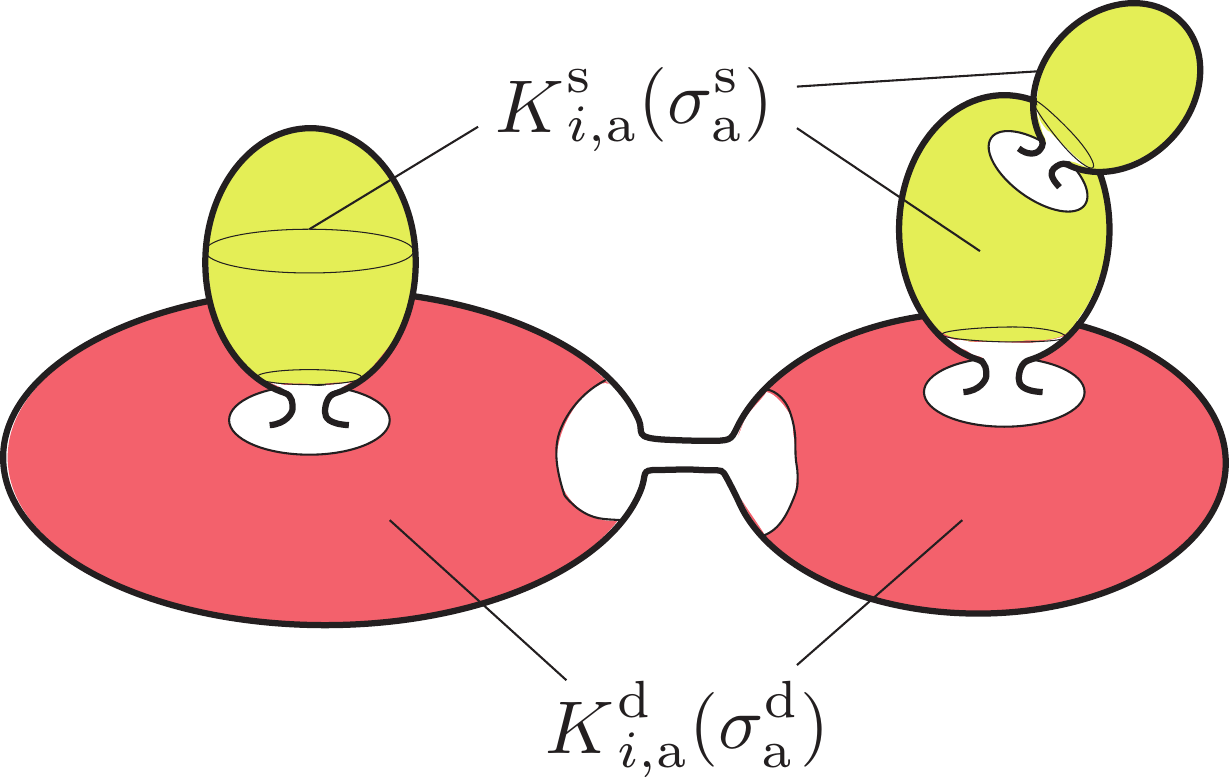}
\caption{Core.}
\label{Figure159-3}
\end{figure}
\index{core}
\begin{figure}[ht]
\centering
\includegraphics[scale=0.3]{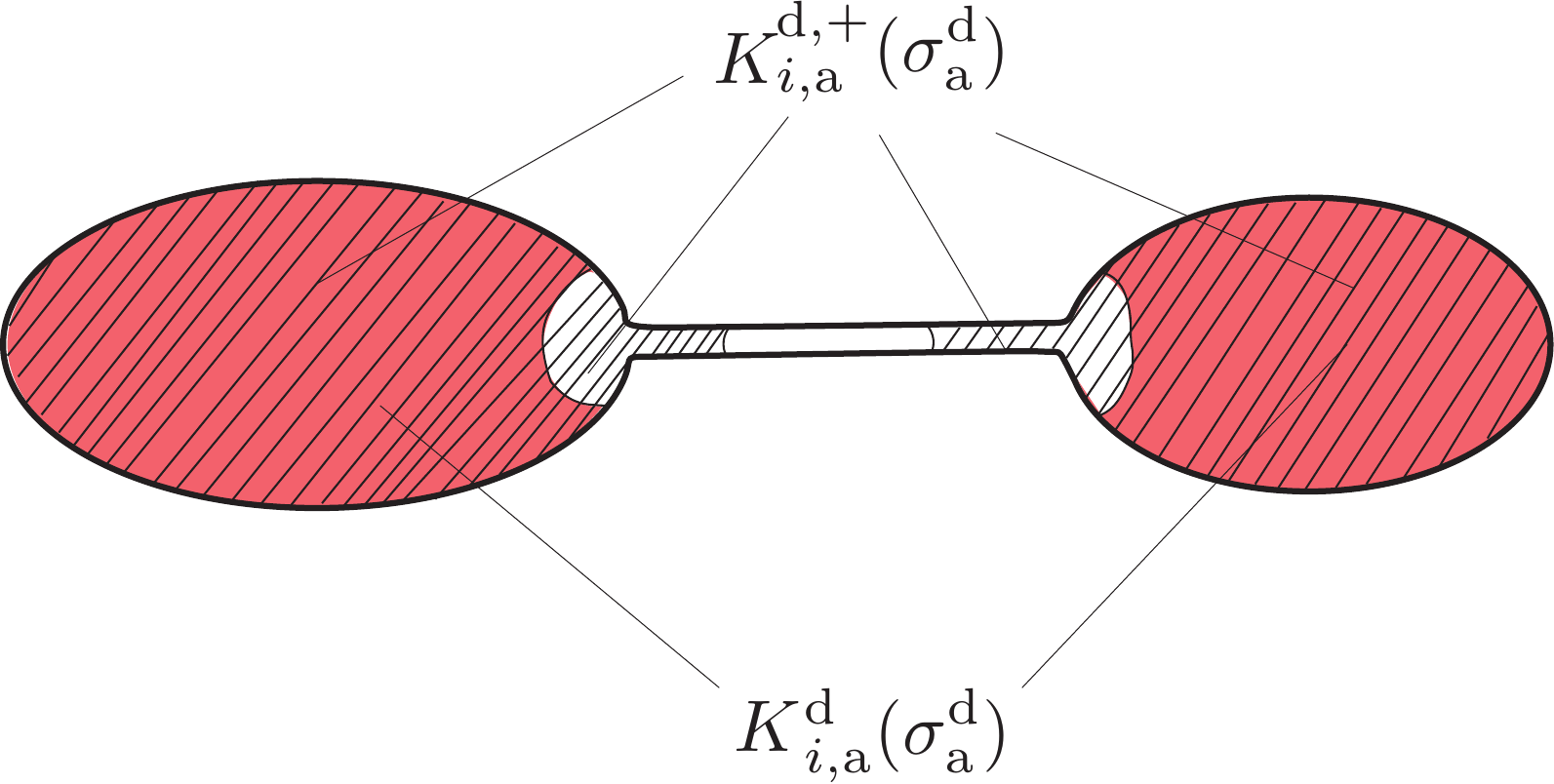}
\caption{$K^{\rm d}_{i,\rm a}\bigl(\sigma^{\rm d}_{\rm a}\bigr)$ and
$K^{\rm d,+}_{i,\rm a}\bigl(\sigma^{\rm d}_{\rm a}\bigr)$.}
\label{Figure159-4}
\end{figure}
\begin{defn}\label{defn1227}
By definition, we have holomorphic embeddings
\[
\mathfrak I^{+, \rm d}_{i,{\rm a},\sigma,{\bf r}} \colon\ K^{+, \rm d}_{i,\rm a}\bigl(\sigma^{\rm d}_{\rm a}\bigr)
\to \Sigma_i(\sigma,{\bf r}), \qquad
\mathfrak I^{+, \rm s}_{i,{\rm a},\sigma,{\bf r}} \colon\ K^{+, \rm s}_{i,\rm a}(\sigma^{\rm s}_{\rm a})
\to \Sigma_i(\sigma,{\bf r}).
\]
We call its restriction
\[
\mathfrak I^{\rm d}_{i,{\rm a},\sigma,{\bf r}} \colon\ K^{\rm d}_{i,\rm a}\bigl(\sigma^{\rm d}_{\rm a}\bigr)
\to \Sigma_i(\sigma,{\bf r}), \qquad
\mathfrak I^{\rm s}_{i,{\rm a},\sigma,{\bf r}} \colon\ K^{\rm s}_{i,\rm a}(\sigma^{\rm s}_{\rm a})
\to \Sigma_i(\sigma,{\bf r}),
\]
the {\it canonical holomorphic embedding}.
\index{canonical holomorphic embedding}
\end{defn}
\begin{lem}
All the weak isomorphisms $\psi = (\psi_1,\psi_2) \colon \xi \to \xi$
canonically induce biholomorphic maps
$
\psi_{i,\sigma,{\bf r}}
 \colon \Sigma_i(\sigma,{\bf r}) \to \Sigma_i((\psi_i)_*(\sigma,{\bf r}))$.

\end{lem}
\begin{proof}
The map $\psi_i$ permutes the interior nodes. We permute the components of the
gluing parameter ${\bf r}$ in the same way.
The map $\psi_i$ also permutes the sphere components.
We permute the components of $\sigma$ in the same way.
This is the definition of $(\psi_i)_*$.
The lemma is then an immediate consequence of
Definition~\ref{defn1222}\,(4) and the construction.
\end{proof}

Our next task is to define a biholomorphic map
$\mathscr I_{\sigma,{\bf r}} \colon \Sigma_1^0(\sigma,{\bf r})
\to \Sigma_2^0(\sigma,{\bf r})$.
Here $\Sigma_i^0(\sigma,{\bf r})$ is the union of disk components of
$\Sigma_i(\sigma,{\bf r})$.

Such an isomorphism is not canonically induced from $\mathscr I$,
since $\Sigma_1^0(\sigma,{\bf r})$ may contain a part of
the {\it sphere} components $K^{\rm s}_{i,\rm a}(\sigma^{\rm s}_{\rm a})$,
on which $\mathscr I$ is {\it not} defined.

We take a certain special choice of the coordinates around the
nodes which
we use to glue, so that we can define $\mathscr I_{\sigma,{\bf r}}$.
\begin{defn}\label{defn122812}
\quad
\begin{enumerate}\itemsep=0pt
\item[(1)]
A holomorphic embedding
$D^2 \to S^2$ is said to be {\it extendable}
\index{extendable} if it is a restriction
of a~biholomorphic map $S^2 \to S^2$.
\item[(2)]
A holomorphic embedding
$D^2 \to D^2$ is said to be {\it extendable}
\index{extendable} if it is a restriction
of biholomorphic map $D^2(R) \to D^2$
for some $R > 1$.
\item[(3)]
A holomorphic embedding
$\bigl(D_{\ge 0}^2,D^2\cap \R\bigr) \to \bigl(D^2,\partial D^2\bigr)$
is said to be {\it extendable} if its double is extendable
in the sense of $(1)$.
\item[(4)]
An analytic family of coordinates is said to be {\it extendable}
if its members are extendable in the sense of (1), (2) or (3).
\end{enumerate}

\end{defn}
We recall that we assumed that all the analytic families of coordinates
appearing as a part of gluing data are extendable.
(See Definition~\ref{defn1222}\,(8).)
\begin{lem}\label{loem1229}
We can canonically define a biholomorphic map $\mathscr I_{\sigma,{\bf r}} \colon \Sigma_1^0(\sigma,{\bf r})
\to \Sigma_2^0(\sigma,{\bf r})$
with the following properties:
\begin{enumerate}\itemsep=0pt
\item[$(1)$]
The next diagram commutes:
\begin{equation}\label{diag1218}
\begin{CD}
K^{+, \rm d}_{1,\rm a}\bigl(\sigma^{\rm d}_{\rm a}\bigr)
 @ >>>
K^{+, \rm d}_{2,\rm a}\bigl(\sigma^{\rm d}_{\rm a}\bigr)
 \\
@ V{\overline{\mathfrak I}^{+, \rm d}_{1,{\rm a},\sigma,{\bf r}}}VV @ VV{\overline{\mathfrak I}^{+, \rm d}_{2,{\rm a},\sigma,{\bf r}}}V\\
\Sigma_1^0(\sigma,{\bf r})
@ >{\mathscr I_{\sigma,{\bf r}}}>>
\Sigma_2^0(\sigma,{\bf r}),
\end{CD}
\end{equation}
where the first horizontal arrow is the isomorphism induced by $\mathscr I$.
The vertical arrows are maps induced by
\smash{$\mathfrak I^{+, \rm d}_{1,{\rm a},\sigma,{\bf r}}$}
and \smash{$\mathfrak I^{+, \rm d}_{2,{\rm a},\sigma,{\bf r}}$}.
\item[$(2)$]
If $\psi = (\psi_1,\psi_2)$ is a weak isomorphism: $ \xi \to \xi$,
then we have
$
\mathscr I_{\psi_*(\sigma,{\bf r})} \circ \psi_{1,\sigma,{\bf r}}
= \psi_{2,\sigma,{\bf r}} \circ \mathscr I_{\sigma,{\bf r}}$.
\item[$(3)$]
$\mathscr I_{\sigma,{\bf r}}\bigl(z_{1,j}^{ \rm int}(\sigma,{\bf r})\bigr)
= z_{2,j}^{ \rm int}(\sigma,{\bf r})$.
It also preserves boundary marked points.
\end{enumerate}

\end{lem}
\begin{proof}
We put
$
\Sigma_i(\sigma) := \Sigma_i(\sigma,{\bf 0})
$,
where the gluing parameter ${\bf 0}$ is
by definition $r_{\rm b} = 0$, $\mathfrak r_{\rm b} = 0$
for all $\rm b$.
Since we deform the disk components of $\Sigma_1$ and of $\Sigma_2$
in exactly the same way by definition,
we have a biholomorphic maps
$
\mathscr I_{\sigma,{\bf 0}} \colon \Sigma^0_1(\sigma) \to \Sigma^0_2(\sigma).
$
Therefore, to construct $\mathscr I_{\sigma,{\bf r}}$ it suffices to
find biholomorphic maps~${
\mathscr J_{\sigma,i} \colon \Sigma^0_i(\sigma,{\bf r}) \to \Sigma^0_i(\sigma)}
$
such that
\smash{$
\mathscr J_{\sigma,i}
\circ \mathfrak I^{+, \rm d}_{i,{\rm a},\sigma,{\bf r}}
=
\mathfrak I^{+, \rm d}_{i,{\rm a},\sigma,{\bf 0}}$}.

We describe the construction of $\mathscr J_{\sigma,i}$ in the following case.
$\Sigma_i(\sigma)$ is a union of $D^2$ and $S^2$ where we
glue them at $0 \in D^2$ and $0 \in S^2 = \C \cup \{\infty\}$.
By the definition of extendable coordinate, we take
our coordinate $\varphi^{\rm d}$ and $\varphi^{\rm s}$
by
$
\varphi^{\rm d}(z) = cz \in D^2$,
$
\varphi^{\rm s}(z) = c'z \in \C \cup \{\infty\}
$,
where $c \in \R_+$ is a small positive number and
$c' \in \C$ is a nonzero complex number with small absolute value.
We denote the gluing parameter by $\mathfrak r \in D^2$.

By definition, $\Sigma^0_i(\sigma) = D^2 \subset \Sigma_i(\sigma)$.
Let $\mathfrak r \ne 0$. Then
$\Sigma^0_i(\sigma,{\bf r})$
is obtained by gluing
\begin{equation}\label{form12221}
D^2 \setminus D^2(c \vert \mathfrak r\vert)
\end{equation}
and
\begin{equation}\label{form12222}
\C \cup \{\infty\} \setminus D^2(\vert c'\mathfrak r\vert)
\end{equation}
by the equivalence relation $\sim$.
The equivalence relation $\sim$ is defined in Definition~\ref{defn1225}.
In our case, it is described as follows.
Let $z \in \eqref{form12221}$ and
$w \in \eqref{form12222}$.
Then $z \sim w$ if and only if
$z/c \times w/c' = \mathfrak r$.
Namely, $z = cc' \mathfrak r/w$.
Therefore, we define $\mathscr J_{\sigma,i}$
such that
$\mathscr J_{\sigma,i}(z) = z$ if~${z \in \eqref{form12221}}$
and $\mathscr J_{\sigma,i}(w) = cc' \mathfrak r/w$ if
$w \in \eqref{form12222}$.
See Figure~\ref{Figure152-1}.
\begin{figure}[ht]
\centering
\includegraphics[scale=0.36]{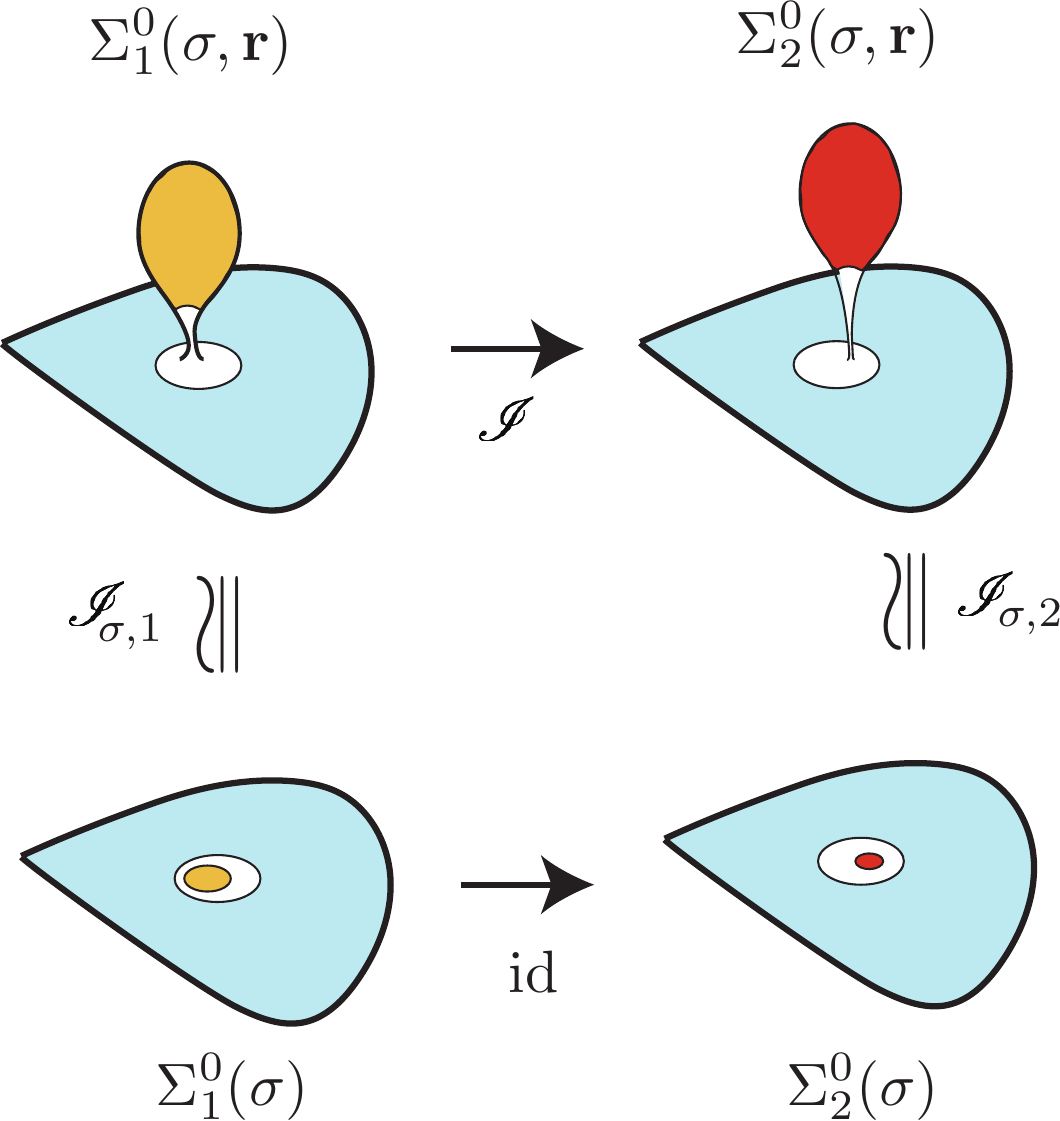}
\caption{Definition of $\mathscr I$.}
\label{Figure152-1}
\end{figure}

We thus defined $\mathscr J_{\sigma,i}$ in the above cases.
Its definition in the general case is similar.
See Figure~\ref{Figure152-2} below.
\begin{figure}[ht]
\centering
\includegraphics[scale=0.4]{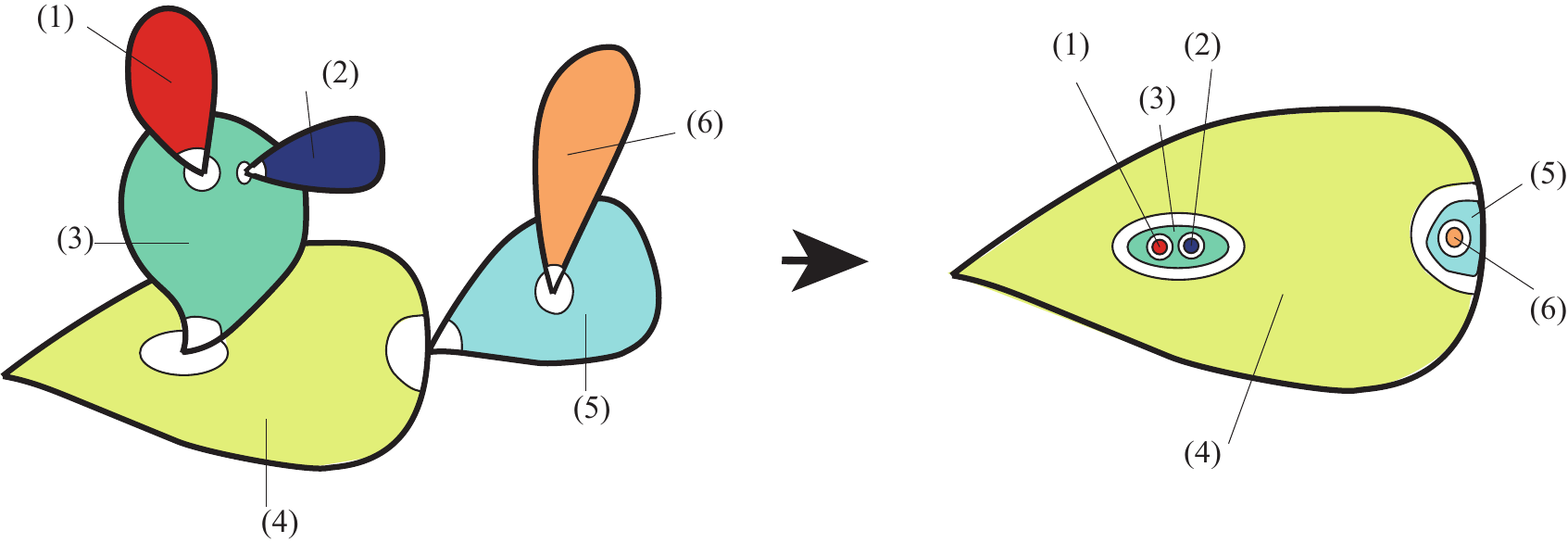}
\caption{Definition of $\mathscr I_{\sigma,i}$
in the general case.}
\label{Figure152-2}
\end{figure}

The properties (1), (2), (3) can be easily proved from the construction.
\end{proof}

\begin{rem}\label{re1234}
We remark that in our situation $\Sigma_i^0(\sigma,{\bf r})$
is a tree of disks {\it without} sphere bubbles.
This is because of Definition~\ref{defn1222}\,(5), that is,
all the marked points on the sphere components are of second kind.
Note that we forget all the marked points of second kind
to obtain $\Sigma_i^0(\sigma,{\bf r})$.
By this reason, we consider only core of
disk components $K^{+, \rm d}_{i,\rm a}\bigl(\sigma^{\rm d}_{\rm a}\bigr)$
in \eqref{diag1218}.
Since we are studying ${\mathcal M}'_{0,0,0}(L_{12};\vec a;E)$
as we mentioned at the beginning of the proof,
we can assume Definition~\ref{defn1222}\,(5).
When we generalize the construction to the case
of ${\mathcal M}'_{\ell,\ell_1,\ell_2}(L_{12};\vec a;E)$
then we need to study the case when there is a marked point
of the first kind on the sphere components.
So there may exist a core $K^{+, \rm s}_{i,\rm a}(\sigma^{\rm s}_{\rm a})$
in the sphere components contained in $\Sigma_i^0(\sigma,{\bf r})$.
In such cases, to define $\mathscr I_{\sigma,{\bf r}}$,
on such parts, we need to modify \eqref{form1212}.
Namely, for example, in place of
\smash{$\prod_{i=1,2}\prod_{{\rm a} \in \mathfrak{comp}^{\rm s}_i} \mathcal V(\xi^{\rm s}_{i,\rm a})$}
we need to consider its subset such that $\mathcal V(\xi^{\rm s}_{1,\rm a})$
factor and~$\mathcal V(\xi^{\rm s}_{2,\rm a'})$ factor
are the same for certain $a$, $a'$.
We do not discuss this point since we do not use it.
\end{rem}

We thus described the way to glue source curves.
To discuss the way to glue maps (that is the part where
nonlinear functional analysis enters),
we first describe the way to define obstruction spaces.
This part is mostly the same as the construction of the
Kuranishi structure on the moduli space of pseudo-holomorphic
disks.
(See \cite[Sections 17 and 18]{foootech}, \cite{fooo:const1,fooo:const2}.)
We include its discussion here for completeness.

Let
$\xi = \bigl(\bigl(\bigl(\Sigma_1,\vec z_1,\vec z_1^{\rm \,int},\vec w_1^{ \rm int}\bigr),u_1\bigr),\bigl(\bigl(\Sigma_2,\vec z_2
,\vec z_2^{\rm \,int},\vec w_2^{ \rm int}\bigr),u_2\bigr),\mathscr I,\gamma\bigr)$
be an element of the moduli space ${\mathcal M}'_{\ell,\ell_1,\ell_2}(L_{12};\vec a;E)$.
Using the above notations, we consider the source deformation parameter
$\sigma$ which corresponds to the source curve of $\xi$ itself.
We denote this $\sigma$ as ${\bf 0}$.
We put~${
K^{\rm s}_{i,\rm a} = K^{\rm s}_{i,\rm a}({\bf 0}) \subset \Sigma_i}$,
$
K^{\rm d}_{i,\rm a} = K^{\rm d}_{i,\rm a}({\bf 0}) \subset \Sigma_i$.
\begin{defn}\label{defn12310}
An {\it obstruction bundle data} \index{obstruction bundle data} $\mathscr{OB}$
\index[syindex]{OzzBscr@$\mathscr{OB}$} centered at $\xi$ is the following objects:
\begin{enumerate}\itemsep=0pt
\item[(1)]
We take a stabilization data at $\xi$.
\item[(2)]
We take a finite-dimensional linear subspace
$
\mathcal E^{\rm s}_{i,\rm a}
\subset C^{\infty}\bigl(K^{\rm s}_{i,\rm a},u_i^*TX_i \otimes \Lambda^{0,1}\bigr)
$
for each sphere component of $\Sigma_i$ and
$
\mathcal E^{\rm d}_{\rm a}
\subset C^{\infty}\bigl(K^{\rm d}_{1,\rm a},(u_1,u_2)^*\bigl(T(X_1 \times X_2) \otimes \Lambda^{0,1}\bigr)\bigr)
$
for each disk component of $\Sigma_i$.

Note that we regard \smash{$K^{\rm d}_{1,\rm a} \cong K^{\rm d}_{2,\rm a}$} by
$\mathscr I$.

We call them the {\it obstruction spaces}. \index{obstruction spaces}
We assume that the supports of the elements of obstruction spaces
are away from nodes and marked points.
We also assume that the supports of the elements of $\mathcal E^{\rm d}_{\rm a}$ is away from
boundary.
Furthermore, we assume the supports of the elements
of $\mathcal E^{\rm d}_{\rm a}$ (resp.\ $\mathcal E^{\rm s}_{\rm a}$)
are in a compact subset contained in the interior of
$K^{\rm d}_{1,\rm a}$ (resp.\ $K^{\rm s}_{i,\rm a}$).
\item[(3)]
We assume that the obstruction spaces satisfy the transversality conditions (see
Conditions~\ref{defn12332} and~\ref{conds1234} below).
\item[(4)]
We assume that $\{\mathcal E^{\rm s}_{i,\rm a}\}$ and
$\bigl\{\mathcal E^{\rm d}_{\rm a}\bigr\}$ are invariant
of the weak isomorphism $\xi' \to \xi'$, where
$\xi'$ is the source stabilization of $\xi$ which is a
part of the stabilization data given in (1).
\item[(5)]
We require that $\mathcal E^{\rm s}_{i,\rm a} = 0$ if $u_i$ is
constant on $\Sigma^{\rm s}_{i,{\rm a}}$ and
 $\mathcal E^{\rm d}_{\rm a} = 0$ if $u$ is
constant on $\Sigma^{\rm d}_{{\rm a}}$.
\item[(6)]
We require
\smash{$
\operatorname{Diam}
\bigl(u_i\circ \varphi^{(j),\rm int}_{{\rm b},\sigma}\bigr)\bigl(D^2\bigr)
\le \varepsilon_1
$}
for each ${\rm b} \in {\rm Node}^+_{i,\rm int}$ and
$\smash{
\operatorname{Diam}
\bigl(u_i\circ \varphi^{(j),\partial}_{i,{\rm b},\sigma}\bigr)\bigl(D_{\ge 0}^2\bigr)}\allowbreak
\le \varepsilon_1
$
for each \smash{${\rm b} \in {\rm Node}^+_{\partial}$}.
Here $\varepsilon_1$ is a sufficiently small number.
(It is smaller than the injectivity radius of $X_1 \times X_2$.
It is the constant appearing \cite[Condition 3.1]{foooanalysis}.)
\item[(7)]
We require all the marked points are in the core,
$K^{\rm s}_{i,\rm a}$, $K^{\rm d}_{i,\rm a}$.
\end{enumerate}

\end{defn}
Below we describe the transversality condition mentioned in item (3).
We review the linearization of the nonlinear Cauchy--Riemann equation for this purpose.
For each sphere component,
the linearization of the nonlinear Cauchy--Riemann equation
induces a linear differential operator of first order
\begin{equation}\label{form1228}
\bigl(D_{u_i}\overline \partial\bigr)_{i,\rm a}^{\rm s} \colon\
C^{\infty}(\Sigma_{i,\rm a}^{\rm s};u_i^*TX_i)
\to C^{\infty}\bigl(\Sigma_{i,\rm a}^{\rm s};u_i^*TX_i \otimes \Lambda^{0,1}\bigr).
\end{equation}
The definition of the function spaces appearing in \eqref{form1228} is obvious from
notation.
Let us discuss the case of disk component $\Sigma_{i,\rm a}^{\rm d}$.
We remark that $\Sigma_{1,\rm a}^{\rm d} \cong \Sigma_{2,\rm a}^{\rm d}$,
which we write $\Sigma_{\rm a}^{\rm d}$. The pair of maps
$u = (u_1,u_2)$ define a map \smash{$u \colon \Sigma_{\rm a}^{\rm d} \to -X_1 \times X_2$}.
Let $\vec z_{\rm a} = (z_{{\rm a},1},\dots,z_{{\rm a},k_{\rm a}})$ be the set of
all marked or nodal points on $\Sigma_{1,\rm a}^{\rm d}$.
$u(z_{{\rm a},j})$ lies on the image of $\tilde L_{12} \times_{X_1\times X_2} \tilde L_{12}
= \cup L_{12}(a)$.
We define $a_{{{\rm a},j}}$ such that $u(z_{{\rm a},j})$ lies in the image of $L_{12}(a_{{{\rm a},j}})$.
\begin{defn}\label{defn1231}
We define the function space
\[
C^{\infty}\bigl(\bigl(\Sigma_{\rm a}^{\rm d},\partial \Sigma_{\rm a}^{\rm d},\vec z_{\rm a}\bigr); (u^*TX,\gamma^* L_{12},L_{12}(\vec a_{\rm a}))\bigr)
\]
as the set of the pairs $(V,v)$ such that
\begin{enumerate}\itemsep=0pt
\item[(1)]
$V$ is a section of $u^{*}T(X_1 \times X_2)$ defined on
$\Sigma_{\rm a}^{\rm d}$.
\item[(2)]
$v$ is a section of $\gamma^* T\tilde L_{12}$ defined on
$\partial\Sigma_{\rm a}^{\rm d} \setminus \vec z_{\rm a}$.
\item[(3)]
If $z \in \partial\Sigma_{\rm a}^{\rm d} \setminus \vec z_{\rm a}$, then
$
V(z) := (d i_{L_{12}})(v(z))$.
Here $i_{L_{12}} \colon \tilde L_{12} \to L_{12}$ is the immersion.
\item[(4)]
Let $z_{{\rm a},j} \in \vec z_{\rm a}$. We then require
\[
(\lim_{z \in \partial\Sigma_{\rm a}^{\rm d} \uparrow z_{{\rm a},j}} v(z),
\lim_{z \in\partial\Sigma_{\rm a}^{\rm d} \downarrow z_{{\rm a},j}} v(z)
)
\in TL_{12}(a_{{{\rm a},j}}).
\]
\end{enumerate}
The operator
\[
\bigl(D_{u}\overline \partial\bigr)_{\rm a}^{\rm d} \colon\
C^{\infty}\bigl(\bigl(\Sigma_{\rm a}^{\rm d},\partial \Sigma_{\rm a}^{\rm d},\vec z_{\rm a}\bigr); (u^*TX,\gamma^* L_{12},L_{12}(\vec a_{\rm a}))\bigr) \to
C^{\infty}\bigl(\Sigma_{\rm a}^{\rm d}; u^*T(X_1\times X_2) \otimes \Lambda^{0,1}\bigr)
\]
is defined by
$
\bigl(D_{u}\overline \partial\bigr)_{\rm a}^{\rm d}(V,v) := \bigl(D_{u}\overline \partial\bigr)(V)$.

\end{defn}
\begin{conds}\label{defn12332}
We say that obstruction spaces $\mathcal E^{\rm s}_{i,\rm a}$, $\mathcal E^{\rm d}_{\rm a}$
satisfy {\it mapping transversality condition}
\index{mapping transversality condition}
if the following holds:
\begin{enumerate}\itemsep=0pt
\item[(1)]
For each sphere component $\Sigma_{i,\rm a}^{\rm s}$, we assume
\[
\operatorname{Im}\bigl(D_{u_i}\overline \partial\bigr)_{i,\rm a}^{\rm s}
+ \mathcal E^{\rm s}_{i,\rm a}
= C^{\infty}\bigl(\Sigma_{i,\rm a}^{\rm s};u_i^*TX_i \otimes \Lambda^{0,1}\bigr).
\]
\item[(2)]
For each disk component $\Sigma_{\rm a}^{\rm d}$, we assume
\[
\operatorname{Im} \bigl(D_{u}\overline \partial\bigr)_{\rm a}^{\rm d}
+ \mathcal E^{\rm d}_{\rm a}
= C^{\infty}\bigl(\Sigma_{\rm a}^{\rm d}; u^*T(X_1\times X_2) \otimes \Lambda^{0,1}\bigr).
\]
\end{enumerate}
\end{conds}
To describe another transversality condition, we define a linearized version $\mathcal{EV}$ of the
evaluation map.
The domain of this evaluation map is the direct sum
\begin{gather}
\bigoplus_{i=1,2}\bigoplus_{a} C^{\infty}(\Sigma_{i,\rm a}^{\rm s};u_i^*TX_i) \oplus
\bigoplus_{a} C^{\infty}\bigl(\bigl(\Sigma_{\rm a}^{\rm d},\partial \Sigma_{\rm a}^{\rm d},\vec z_{\rm a}\bigr); (u^*TX,\gamma^* L_{12},L_{12}(\vec a_{\rm a}))\bigr).\label{123000}
\end{gather}
Here the first direct sum is taken over all the sphere components $\Sigma_{i,\rm a}^{\rm s}$ and
the second direct sum is taken over all the disk components \smash{$\Sigma_{\rm a}^{\rm d}$}.

We next describe the target of $\mathcal{EV}$.
Let $\mathfrak z_{\rm b}$ be a boundary node.
There exists a component~$L_{12}(a_{\rm b})$ of $\tilde L_{12} \times_{X_1\times X_2} \tilde L_{12}$
such that it is mapped to $u(\mathfrak z_{\rm b}) = (u_1(\mathfrak z_{\rm b}),u_2(\mathfrak z_{\rm b}))$ by $i_{L_{12}}$.
The target space of $\mathcal{EV}$ is the direct sum
\begin{equation}
\bigoplus_{i=1,2}\bigoplus_{\rm b} T_{u_i(\mathfrak z_{\rm b})} X_i
\oplus
\bigoplus_{\rm b} T_{\gamma(\mathfrak z_{\rm b})} L_{12}(a_{\rm b}).\label{123131}
\end{equation}
Here the first direct sum is one over interior nodes $\mathfrak z_{\rm b}$.
The second direct sum is one over boundary nodes $\mathfrak z_{\rm b}$.
The point $\gamma(\mathfrak z_{\rm b}) \in L_{12}(a_{\rm b})$ is by definition
\[
\gamma(\mathfrak z_{\rm b}) =
(\lim_{z \in \partial\Sigma_{\rm a}^{\rm d} \uparrow \mathfrak z_{\rm b}} \gamma(z),
\lim_{z \in\partial\Sigma_{\rm a}^{\rm d} \downarrow \mathfrak z_{\rm b}} \gamma(z)
)
\in L_{12}(a_{{{\rm a},j}}).
\]
(See Definition~\ref{def3737}\,(5).)

Now we define
\begin{equation}\label{form1232}
\mathcal{EV}\colon\ \eqref{123000} \to \eqref{123131}.
\end{equation}
Let $\vec V = ((V_{a,1},(V_{a,2})),(V_a))$ be an element of domain \eqref{123000}.
Let $\mathfrak z_b$ be an interior node.
There are two components
\smash{$\Sigma^{c(1,b)}_{i,{\rm a}(1,b)}$}, \smash{$\Sigma^{c(2,b)}_{i,{\rm a}(2,b)}$}
containing it. Here $c(1,b)$, $c(2,b)$ are either ${\rm s}$ or ${\rm d}$.
Suppose~$c(1,b) = {\rm d}$, $c(2,b) = {\rm s}$. Then
we define
\begin{equation}\label{formula1234}
\text{$T_{u_i(\mathfrak z_{\rm b})} X_i$ component of
$\mathcal{EV}\bigl(\vec V \bigr)$}
=
\Pi_i(V_{{\rm a}(1,b)}(\mathfrak z_{\rm b}))
- V_{a(2,b),i}(\mathfrak z_{\rm b}),
\end{equation}
where $\Pi_i \colon T(X_1 \times X_2) \to T(X_i)$ is the
projection. The definitions in the other cases of $c(1,b)$, $c(2,b)$ are similar.

Let $\mathfrak z_{b}$ be a boundary node.
There are two disk components
$\Sigma^{\rm d}_{{\rm a}(1,b)}$, $\Sigma^{\rm d}_{{\rm a}(2,b)}$
containing it.
We define
\begin{equation}
\text{$T_{\gamma(\mathfrak z_{\rm b})}L_{12}({\rm a}_b)$ component of
$\mathcal{EV}\bigl(\vec V \bigr)$} =
(V_{{\rm a}(1,b)},v_{{\rm a}(1,b)})(\mathfrak z_{\rm b})
- (V_{{\rm a}(2,b)},v_{{\rm a}(2,b)})(\mathfrak z_{\rm b}).\label{formula1235}
\end{equation}
Here
\[
(V_{{\rm a}(i,b)},v_{{\rm a}(i,b)})(\mathfrak z_{\rm b})
=
(\lim_{z \in \partial\Sigma_{{\rm a}(i,b)}^{\rm d} \uparrow \mathfrak z_{\rm b}} v(z),
\lim_{z \in\partial\Sigma_{{\rm a}(i,b)}^{\rm d} \downarrow \mathfrak z_{\rm b}} v(z)
)
\in T_{\gamma(\mathfrak z_{\rm b})}L_{12}.
\]

\eqref{formula1234} and \eqref{formula1235}
define a map \eqref{form1232}.
\begin{conds}\label{conds1234}
We say that obstruction spaces $\mathcal E^{\rm s}_{i,\rm a}$, $\mathcal E^{\rm d}_{\rm a}$
satisfy {\it evaluation transversality condition}
\index{evaluation transversality condition}
if the restriction of $\mathcal{EV}$ to the direct sum
\[
\bigoplus_{i=1,2}\bigoplus_{{\rm a}} \bigl(\bigl(D_{u_i}\overline \partial\bigr)_{i,\rm a}^{\rm s}
\bigr)^{-1}(\mathcal E^{\rm s}_{i,\rm a})
\oplus
\bigoplus_{{\rm a}} \big(\bigl(D_{u}\overline \partial\bigr)_{\rm a}^{\rm d}
\big)^{-1}\bigl(\mathcal E^{\rm d}_{\rm a}\bigr)
\]
is surjective.

\end{conds}
We thus defined the notion of obstruction bundle data.
Our next task is to send obstruction spaces to a nearby object.
We make precise the meaning of `nearby object' below.
\begin{defn}\label{defn1233}
A {\it candidate of an element of the extended moduli space}
${\mathcal M}'_{\ell,\ell_1,\ell_2}(L_{12};\vec a;E)$
\index{candidate of an element of extended moduli space} is, by definition, an object
\[
\eta = \bigl(\bigl(\bigl(\Sigma^{\heartsuit}_1,\vec z^{ \heartsuit}_1,\vec z_1^{ \heartsuit,\rm int},\vec w_1^{ \heartsuit,\rm int}\bigr),u^{\heartsuit}_1\bigr),\bigl(\bigl(\Sigma^{\heartsuit}_2,\vec z^{ \heartsuit}_2
,\vec z_2^{ \heartsuit,\rm int},\vec w_2^{ \heartsuit,\rm int}\bigr),u^{\heartsuit}_2\bigr),\mathscr I^{\heartsuit},
\gamma^{\heartsuit}\bigr),
\]
which satisfies the same conditions as Condition \ref{defn145555}
except we do not assume
$u_i^{\heartsuit} \colon \Sigma^{\heartsuit}_i \to X_i$ is pseudo-holomorphic as in Condition \ref{defn145555}\,(2) but only
assume that it is of $C^{\infty}$ class.
\end{defn}
\begin{defn}\label{defn1236}
Let
$\xi = \bigl(\bigl(\bigl(\Sigma_1,\vec z_1,\vec z_1^{\rm \,int},\vec w_1^{ \rm int}\bigr),u_1\bigr),\bigl(\bigl(\Sigma_2,\vec z_2
,\vec z_2^{\rm \,int},\vec w_2^{ \rm int}\bigr),u_2\bigr),\mathscr I,\gamma\bigr)$
be an element of ${\mathcal M}'_{\ell,\ell_1,\ell_2}(L_{12};\vec a;E)$.
We assume that $\xi$ is source stable and fix a
source gluing data $\mathscr{GL}$ on it.

Let $\eta$ be a candidate of an element of the extended moduli space
of ${\mathcal M}'_{\ell,\ell_1,\ell_2}(L_{12};\vec a;E)$.

We say that $\eta$ is {\it $\varepsilon$ close} to $(\xi,\mathscr{GL})$,
if the following holds:
\begin{enumerate}\itemsep=0pt
\item[(1)]
There exists $\sigma$, $\bf r$ as in \eqref{form1214} and \eqref{form1215}
such that
\begin{equation}\label{firn1235}
\xi_i(\sigma,{\bf r})
= \bigl(\Sigma^{\heartsuit}_i,\vec z^{ \heartsuit}_i,\vec z_i^{ \heartsuit,\rm int},\vec w_i^{ \heartsuit,\rm int}\bigr).
\end{equation}
Moreover, via this isomorphism the biholomorphic map $\mathscr I_{\sigma,{\bf r}}$
in Lemma~\ref{loem1229} is
coincides with $\mathscr I^{\heartsuit}$.
\item[(2)]
The object
$(\sigma,{\bf r})$ is in the $\varepsilon$ neighborhood of $({\bf 0},{\bf 0})$.
\item[(3)]
The restriction of $u_i$ to each $K^{\rm d}_{i,\rm a}\bigl(\sigma^{\rm d}_{\rm a}\bigr)$
is $\varepsilon$ close to
the restriction of $u^{\heartsuit}_i$ to it in $C^2$ norm.
Here we use
$\mathfrak I^{\rm d}_{i,{\rm a},\sigma,{\bf r}}$
and the isomorphism \eqref{firn1235}
to regard the restrictions of $u_i$, $u^{\heartsuit}_i$ as a~map defined on
\smash{$K^{\rm d}_{i,\rm a}\bigl(\sigma^{\rm d}_{\rm a}\bigr)$}.
\item[(4)]
The restriction of $u_i$ to
each $K^{\rm s}_{i,\rm a}(\sigma^{\rm s}_{\rm a})$ is $\varepsilon$ close to
the restriction of $u^{\heartsuit}_i$ to it in $C^2$ norm.
Here we use
$\mathfrak I^{\rm s}_{i,{\rm a},\sigma,{\bf r}}$
and the isomorphism \eqref{firn1235}
to regard the restrictions of $u_i$, $u^{\heartsuit}_i$ as maps defined on
$K^{\rm s}_{i,\rm a}(\sigma^{\rm s}_{\rm a})$.
\item[(5)]
For any connected component $\mathcal S$ of
\[
\Sigma_i(\sigma,{\bf r})
\setminus
\bigcup_{a} K^{\rm d}_{i,\rm a}\bigl(\sigma^{\rm d}_{\rm a}\bigr)
\setminus
\bigcup_{a} K^{\rm s}_{i,\rm a}\bigl(\sigma^{\rm d}_{\rm a}\bigr),
\]
we require
$
\operatorname{Diam} u^{\heartsuit}_i(\mathcal S) < \varepsilon$.
(In other words, we require the diameter of
the images by $u^{\heartsuit}_i$ of the neck regions
are smaller than $\varepsilon$.)
\end{enumerate}

\end{defn}
Let $\eta'$ be a candidate of an element of the extended moduli space
of ${\mathcal M}'_{\ell,\ell_1,\ell_2}(L_{12};\vec a;E)$.
We forget all the interior marked points of $\eta'$
and shrink the components which become unstable.
We then obtain
a candidate of an element of the extended moduli space
of ${\mathcal M}'(L_{12};\vec a;E)$.
We denote it by $\eta = \mathfrak i^*(\eta')$.
Note this definition is a version of \eqref{formass26}.
Here $\mathfrak i$ is \eqref{formass26} with $\ell'=\ell'_1=\ell'_2 = 0$.
\begin{defn}\label{defn123534522}
Let $\xi$ be an element of
${\mathcal M}'(L_{12};\vec a;E)$.
We fix its stabilization data $\mathscr{ST}$.
Let $\eta$ be a candidate of an element of the extended moduli space
of ${\mathcal M}'(L_{12};\vec a;E)$.

We say that $\eta$ is {\it $\varepsilon$-close} \index{$\varepsilon$-close} to $(\xi,\mathscr{ST})$
if the following holds:
\begin{enumerate}\itemsep=0pt
\item[(1)] There exists a candidate of an element of the extended moduli space
of ${\mathcal M}'_{\ell,\ell_1,\ell_2}(L_{12};\vec a;E)$,
which we denote by $\eta'$ such that
$\mathfrak i^*(\eta') = \eta$.
\item[(2)]
Let $\xi'$ be the source stabilization of $\xi$ which is a part of
$\mathscr{ST}$. (See Definition~\ref{defn1223}\,(1).)
Let $\mathscr{GL}$ be the gluing data at $\xi'$
which is a part of $\mathscr{ST}$.
(See Definition~\ref{defn1223}\,(3).)

Then $\eta'$ is $\varepsilon$ close to $(\xi',\mathscr{ST})$
in the sense of Definition~\ref{defn1236}.
\item[(3)]
Let $z_{1,j}^{ \rm int}$ be an interior marked point of first
kind of $\xi'$.
Let $z_{1,j}^{\heartsuit \rm int}$ be the corresponding~interior marked point of first
kind of $\eta'$. (See \eqref{firn1235}.)
Let \smash{$\mathcal N^{(1)}_j$} be the codimension $2$ submanifold
of $-X_1 \times X_2$ which is a part of $\mathscr{ST}$.
(See Definition~\ref{defn1223}\,(5).)
We require
\begin{equation}\label{formtran1}
u^{\heartsuit}\bigl(z_{1,j}^{\heartsuit \rm int}\bigr) \in \mathcal N^{(1)}_j.
\end{equation}
Here $u^{\heartsuit}(z) = \bigl(u^{\heartsuit}_1(z),u^{\heartsuit}_2(\mathscr I'(z))\bigr)$.
\item[(4)]
Let $w_{i,j}^{ \rm int}$ be an interior marked point of second
kind of $\xi'$, and let $w_{i,j}^{\heartsuit \rm int}$ be the corresponding interior marked point of first
kind of $\eta'$ (see \eqref{firn1235}).
Let \smash{$\mathcal N^{(2)}_{j}$} be the codimension~$2$~sub\-manifold
of $X_i$ which is a part of $\mathscr{ST}$.
(See Definition~\ref{defn1223}\,(6).)
We require
\begin{equation}\label{formtran2}
u^{\heartsuit}_i\bigl(w_{i,j}^{\heartsuit \rm int}\bigr) \in \mathcal N^{(2)}_{j}.
\end{equation}
\end{enumerate}

\end{defn}
Let $\xi$ be an element of
${\mathcal M}'(L_{12};\vec a;E)$.
We fix an obstruction bundle data $\mathscr{OB}$ of it.
It includes a source stabilization data $\mathscr{ST}$.
Let $\eta$ be a candidate of an element of extended moduli space
of ${\mathcal M}'(L_{12};\vec a;E)$ which is $\varepsilon$ close
to $(\xi,\mathscr{ST})$.

Our next task is to send obstruction spaces (which is a part of
$\mathscr{OB}$) to a subspace of sections on the source curve
$\Sigma_{i}^{\heartsuit}$
of $\eta$.

Let $K^{\rm d}_{i,\rm a}$
(resp.\ $K^{\rm s}_{i,\rm a}$)
be a core of disk (resp.\ sphere) component of $\xi$.
We consider
\begin{equation}\label{form1237}
\mathcal I^{\rm d}_{i,\rm a} \colon\ K^{\rm d}_{i,\rm a} \cong
K^{\rm d}_{i,\rm a}\bigl(\sigma^{\rm d}_{\rm a}\bigr)
\to \Sigma_i(\sigma,{\bf r}) \cong \Sigma^{\heartsuit}_i.
\end{equation}
Here the first map is a diffeomorphism which is induced by the
trivialization given in Definition~\ref{defn1222}\,(2).
The second map is the map \smash{$\mathfrak I^{\rm d}_{i,{\rm a},\sigma,{\bf r}}$}
in Definition~\ref{defn1227}.
The third map is a~biholomorphic map \eqref{firn1235}.
Actually, the image of \eqref{form1237} lies in a certain
disk component of~\smash{$\Sigma^{\heartsuit}_i$} which we denote
\smash{$\Sigma^{\rm d,\heartsuit}_{i,{\rm a}^{\heartsuit}}$}. By $\mathscr I$ and $\mathscr I^{\heartsuit}$, we can
identify $\mathcal I^{\rm d}_{1,\rm a}$ with $\mathcal I^{\rm d}_{2,\rm a}$.
We write the composition~\eqref{form1237} as
\smash{$\mathcal I^{\rm d}_{\rm a} \colon K^{\rm d}_{\rm a}
\to \Sigma^{\rm d,\heartsuit}_{{\rm a}^{\heartsuit}}$}
\index[syindex]{Ida@$\mathcal I^{\rm d}_{\rm a}$}. It defines a complex linear map
\smash{$
C_0^{\infty}\bigl(K^{\rm d}_{\rm a};\Lambda^{0,1}\bigr)
\to C^{\infty}\bigl(\Sigma^{\rm d,\heartsuit}_{{\rm a}^{\heartsuit}};\Lambda^1\bigr)$}.
(Note that~$\mathcal I^{\rm d}_{i,\rm a}$ may not be holomorphic.)
Here $C_0^{\infty}$ denotes the set of smooth sections which have
compact support in the interior. We compose it with the projection to obtain{\samepage
\begin{equation}\label{map1238}
C_0^{\infty}\bigl(K^{\rm d}_{\rm a};\Lambda^{0,1}\bigr)
\to C^{\infty}\bigl(\Sigma^{\heartsuit}_{{\rm a}^{\heartsuit}};\Lambda^{0,1}\bigr).
\end{equation}
This map is complex linear.}

On the other hand, for each $z \in K^{\rm d}_{\rm a}$ we take
the minimal geodesic joining $u(z) \in X_1 \times X_2$ and
$u^{\heartsuit}\bigl(\mathcal I^{\rm d}_{\rm a}(z)\bigr)$. Then taking a
complex linear part of the parallel transport
(with respect to a~certain connection for which $L_{12}$ is
parallel), we obtain a complex linear map
\[
T_{u(z)}(-X_1 \times X_2) \to T_{u^{\heartsuit}
\bigl(\mathcal I^{\rm d}_{\rm a}(z)\bigr)}(-X_1 \times X_2).
\]
It induces
\begin{equation}\label{map1239}
C_0^{\infty}\bigl(K^{\rm d}_{\rm a};u^*(T(-X_1 \times X_2))\bigr)
\to C^{\infty}\bigl(\Sigma^{\rm d,\heartsuit}_{{\rm a}^{\heartsuit}};
\bigl(u^{\heartsuit}\bigr)^*T(-X_1 \times X_2)\bigr).
\end{equation}
\eqref{map1238} and \eqref{map1239} are induced by pointwise
complex linear maps. So we take pointwise tensor product
over $\C$ of \eqref{map1238} and \eqref{map1239} and obtain\index[syindex]{psiau@$\Psi_{{\rm a},u^{\heartsuit}}^{\rm d}$}
\[
\Psi_{{\rm a},u^{\heartsuit}}^{\rm d}\colon\
C_0^{\infty}\bigl(K^{\rm d}_{\rm a};u^*(T(-X_1 \times X_2)) \otimes \Lambda^{0,1}\bigr)
\to
C^{\infty}\bigl(\Sigma^{\rm d,\heartsuit}_{{\rm a}^{\heartsuit}};
\bigl(u^{\heartsuit}\bigr)^*T(-X_1 \times X_2)\otimes \Lambda^{0,1}\bigr).
\]
We consider the direct sum
\begin{gather}
\bigoplus_{{\rm a}^{\heartsuit}} C_0^{\infty}\bigl(\Sigma^{\rm d,\heartsuit}_{{\rm a}^{\heartsuit}};\bigl(u^{\heartsuit}\bigr)^*T(-X_1 \times X_2)\otimes \Lambda^{0,1}\bigr)\nonumber\\
\qquad
\oplus
\bigoplus_{i=1,2}\bigoplus_{{\rm a}^{\heartsuit}} C_0^{\infty}\bigl(\Sigma^{\rm s,\heartsuit}_{i,{\rm a}^{\heartsuit}};\bigl(u^{\heartsuit}_i\bigr)^*T(X_i)\otimes \Lambda^{0,1}\bigr).\label{fprm1241}
\end{gather}
Here the first direct sum is taken over disk components of $\Sigma^{\heartsuit}$ and
the second direct sum is taken over sphere components of $\Sigma^{\heartsuit}_i$.
The symbol $0$ in $C^{\infty}_0$ means sections with compact support away from nodal or
marked points and from boundary.

Taking direct sum of the maps $\Psi_{{\rm a},u^{\heartsuit}}^{\rm d} $, we obtain
\[
\Psi_{u^{\heartsuit}}^{\rm d} \colon\ \bigoplus_{\rm a} C_0^{\infty}\bigl(K^{\rm d}_{\rm a};u^*(T(-X_1 \times X_2)) \otimes \Lambda^{0,1}\bigr)
\to \eqref{fprm1241}.
\]
Here direct sum of the domain is taken over disk components.
In case we specify $\xi$ and $\mathscr{OB}$, we write
$\Psi_{\xi,u^{\heartsuit}}^{\rm d}$ or $\Psi_{\xi,\mathscr{OB},u^{\heartsuit}}^{\rm d}$.

We can perform a similar construction for sphere components
to obtain\index[syindex]{psiiu@$\Psi_{i,u^{\heartsuit}}^{\rm s}$}
\[
\Psi_{i,u^{\heartsuit}}^{\rm s} = \Psi_{\xi,\mathscr{OB},i,u^{\heartsuit}}^{\rm s} \colon\ \bigoplus_{\rm a} C_0^{\infty}\bigl(K^{\rm s}_{i,\rm a};u^*(TX_i) \otimes \Lambda^{0,1}\bigr)
\to \eqref{fprm1241}.
\]
\begin{defn}
We denote\index[syindex]{Exiob@$\mathcal E(\xi,\mathscr{OB})$}
\[
\mathcal E(\xi,\mathscr{OB})
:=
\bigoplus_{i=1,2}\bigoplus_{\rm a} \mathcal E^{\rm s}_{i,\rm a}
\oplus
\bigoplus_{\rm a} \mathcal E^{\rm d}_{\rm a}.
\]
We define the subspace
$
\mathcal E(\xi,\mathscr{OB};\eta)
\subset \eqref{fprm1241}
$
to be the image of $\mathcal E(\xi,\mathscr{OB})$
by the map $\Psi_{u^{\heartsuit}}^{\rm d} \oplus \Psi_{1,u^{\heartsuit}}^{\rm s} \oplus \Psi_{2,u^{\heartsuit}}^{\rm s}$.

We write $\mathcal E(\xi;\eta)$ in place of $\mathcal E(\xi,\mathscr{OB};\eta)$ in case
the choice of $\mathscr{OB}$ is obvious from the context.
\end{defn}
We remark that the choice of $\sigma,{\bf r}$ and the third isomorphism
in \eqref{form1237} is {\it not} unique.
The maps \smash{$\Psi_{u^{\heartsuit}}^{\rm d} \oplus \Psi_{1,u^{\heartsuit}}^{\rm s}
\oplus \Psi_{2,u^{\heartsuit}}^{\rm s}$}
depend on this choice.
However, two different choices are transformed each other by the
weak isomorphism of $\xi$. Therefore, by
Definitions \ref{defn1223}\,(7) and~\ref{defn12310}\,(4)
the image~$\mathcal E(\xi,\mathscr{OB})$ is independent of such
choices.

Roughly speaking, the underlying orbifold of Kuranishi chart consists of
$\eta$ such that
$\overline\partial u^{\heartsuit} \in \mathcal E(\xi;\eta)$.
To obtain Kuranishi chart so that we can define coordinate transformation
among them, we need one more steps to work out.

We remark that for each $(\xi,\mathscr{OB})$ there exists
$\varepsilon(\xi,\mathscr{OB})$ such that the construction of $\mathcal E(\xi,\mathscr{OB};\allowbreak\eta)$
works for $\eta$ which is $\varepsilon(\xi,\mathscr{OB})$ close to $(\xi,\mathscr{ST})$.

For each $\xi \in {\mathcal M}'(L_{12};\vec a;E)$, we choose and fix an
obstruction bundle data $\mathscr{OB}$. We also take a closed neighborhood
$\mathfrak N(\xi)$ of $\xi$ in ${\mathcal M}'(L_{12};\vec a;E)$ such that
each element of $\mathfrak N(\xi)$ is $\varepsilon(\xi,\mathscr{OB})/2$ close to
$(\xi,\mathscr{ST})$.

We take a finite subset
\begin{equation}\label{set1246}
\{\xi_{\mathfrak i} \mid \mathfrak i \in {\bf I}\} \subset {\mathcal M}'(L_{12};\vec a;E)
\end{equation}
such that
\begin{equation}\label{openset1247}
\bigcup_{\mathfrak i \in {\bf I}} {\rm Int}\,\mathfrak N(\xi_{\mathfrak i}) = {\mathcal M}'(L_{12};\vec a;E).
\end{equation}
Using the data we fixed above, we will construct a Kuranishi neighborhood of an
arbitrary element
$\xi$ of ${\mathcal M}'(L_{12};\vec a;E)$.
We put
\begin{equation}\label{set1248}
{\bf I}(\xi) := \{ \mathfrak i \in {\bf I} \mid \xi \in \mathfrak N(\xi_{\mathfrak i}) \}.
\end{equation}
By perturbing obstruction spaces of $(\xi_{\mathfrak i},\mathscr{OB})$ we may assume that
the sum
$
\bigoplus_{\mathfrak i \in {\bf I}(\xi)} \mathcal E(\xi_{\mathfrak i},\mathscr{OB};\xi)
$
is a direct sum.
See \cite[Lemma 18.8]{foootech}, which is proved in
\cite[Section 27]{foootech} and more detailed in~\cite[Section 11.4]{fooo:const1}.

We take stabilization data $\mathscr{ST}$ at $\xi$.
We assume that Definition~\ref{defn12310}\,(5) is satisfied.
($\mathscr{ST}$ may or may not coincide with one included in $\mathscr{OB}$ taken above.)
We take $\varepsilon_2(\xi)$ enough small so that if
$\eta$ is a candidate of an element of extended moduli space
of ${\mathcal M}'(L_{12};\vec a;E)$ which is $\varepsilon_2(\xi)$ close
to $(\xi,\mathscr{ST})$ then
$\eta$ is $\varepsilon(\xi_{\mathfrak i},\mathscr{OB})$ close to $(\xi_{\mathfrak i},\mathscr{OB})$
for each $\mathfrak i \in {\bf I}(\xi)$.
\begin{defn}\label{defn12383}
For $\varepsilon < \varepsilon_2(\xi)$, we define
$U(\xi;\varepsilon)$\index[syindex]{Uxiepsi@$U(\xi;\varepsilon)$} to be the isomorphism classes of $\eta$
with the following properties.
\begin{enumerate}\itemsep=0pt
\item[(1)]
$
\eta = \bigl(\bigl(\bigl(\Sigma^{\heartsuit}_1,\vec z^{ \heartsuit}_1\bigr),u^{\heartsuit}_1\bigr),\bigl(\bigl(\Sigma^{\heartsuit}_2,\vec z^{ \heartsuit}_2
\bigr),u^{\heartsuit}_2\bigr),\mathscr I^{\heartsuit},\gamma^{\heartsuit}\bigr)
$
is a candidate of an element of extended moduli space
of ${\mathcal M}'(L_{12};\vec a;E)$.
\item[(2)]
$\eta$ is $\varepsilon$ close to $(\xi,\mathscr{ST})$.
\item[(3)]
\begin{equation}\label{eq1246}
\overline\partial u^{\heartsuit}_i \in \bigoplus_{\mathfrak i \in {\bf I}(\xi)} \mathcal E(\xi_{\mathfrak i},\mathscr{OB};\eta)
\end{equation}
on the image of $\mathcal I^{\rm d}_{\rm a} \colon K^{\rm d}_{\rm a}
\to \Sigma^{\rm d,\heartsuit}_{{\rm a}^{\heartsuit}}$ and
\begin{equation}\label{eq1247}
\overline\partial u^{\heartsuit}_i \in \bigoplus_{\mathfrak i \in {\bf I}(\xi)} \mathcal E(\xi_{\mathfrak i},\mathscr{OB};\eta)
\end{equation}
on the image of $\mathcal I^{\rm s}_{i,\rm a}$.
Moreover, $u^{\heartsuit}_i$ is pseudo-holomorphic
outside the images of $\mathcal I^{\rm d}_{\rm a}$
and~$\mathcal I^{\rm s}_{i,\rm a}$.\index[syindex]{isia@$\mathcal I^{\rm s}_{i,\rm a}$}
\end{enumerate}

Let $\Gamma_{\xi}$ be the set of all automorphisms of $\xi$.
We denote \index[syindex]{Exi@$\mathcal E(\xi)$}
\smash{$
\mathcal E(\xi)
:=\bigoplus_{\mathfrak i \in {\bf I}(\xi)} \mathcal E(\xi_{\mathfrak i},\mathscr{OB})$}.

\end{defn}
The next proposition claims that we can construct a Kuranishi neighborhood of
$\xi$ using the above data.
\begin{prop}\label{prop1239}
For sufficiently small $\varepsilon > 0$, the following holds:
\begin{enumerate}\itemsep=0pt
\item[$(1)$]
There exists a
smooth manifold $V(\xi;\varepsilon)$ of finite dimension
on which $\Gamma_{\xi}$ acts smoothly
such that the quotient space
$V(\xi;\varepsilon)/\Gamma_{\xi}$ is homeomorphic
to $U(\xi;\varepsilon)$.
\item[$(2)$]
We can define a smooth $\Gamma_{\xi}$ equivalent map\index[syindex]{sxi@$\mathfrak s_{\xi}$}
$
\mathfrak s_{\xi} \colon V(\xi;\varepsilon) \to
\mathcal E(\xi)
$
as follows.
For $\hat\eta \in V(\xi;\varepsilon)$ whose equivalence class is mapped to $\eta \in U(\xi;\varepsilon)$
by the homeomorphism in item $(1)$,
we can take its representative and a choice of the map \eqref{form1237}
$($among finitely many possible choices$)$ such that
the components of $\mathfrak s_{\xi}(\hat\eta)$ is obtained by
applying \smash{$\bigl(\Psi_{\xi_{\mathfrak i},u'}^{\rm d}\bigr)^{-1}$} or \smash{$\bigl(\Psi_{\xi_{\mathfrak i},\mathscr{OB},u'}^{\rm d}\bigr)^{-1}$}
to $\overline\partial u$, $\overline\partial u_i$.
\item[$(3)$]
We define \smash{$\psi_{\xi} \colon \mathfrak s_{\xi}^{-1}(0) \to {\mathcal M}'(L_{12};\vec a;E)$}
by regarding an element \smash{$\hat\eta \in \mathfrak s_{\xi}^{-1}(0)$} as
an element of ${\mathcal M}'(L_{12};\vec a;E)$.
Then $\psi_{\xi}$ induces a homeomorphism from
$\mathfrak s_{\xi}^{-1}(0)/\Gamma_{\xi}$ to a neighborhood of $\xi$.
\item[$(4)$]
$
\mathcal U(\xi,\varepsilon) = (V(\xi;\varepsilon),\Gamma_{\xi},\mathcal E(\xi),
\mathfrak s_{\xi},\psi_{\xi})$
is a Kuranishi neighborhood of $\xi$ in the sense of {\rm\cite[\emph{Definition} A1.1]{fooobook2}}.
\end{enumerate}

\end{prop}
\begin{proof}
Below we provide the construction of
$(V(\xi;\varepsilon),\Gamma_{\xi},\mathcal E(\xi),
\mathfrak s_{\xi},\psi_{\xi})$ leaving the gluing analysis and
smoothness proof to the next subsection.

The construction of the manifold $V(\xi;\varepsilon)$ and
a homeomorphism $V(\xi;\varepsilon)/\Gamma_{\xi}
\cong U(\xi;\varepsilon)$ is the gluing construction of the
solution space of the equation \eqref{eq1246} and
\eqref{eq1247}.

The stabilization data of $\xi$ we take include a source stabilization
$\xi'$ and gluing data at it. It induces a source gluing map
whose domain is
\begin{equation}\label{form12122}
\prod_{{\rm a} \in \mathfrak{comp}^{ \rm d}_1} \mathcal V\bigl(\xi^{\prime,\rm d}_{1,\rm a}\bigr)
\times
\prod_{i=1,2}\prod_{{\rm a} \in \mathfrak{comp}^{\rm s}_i} \mathcal V(\xi^{
\prime, \rm s}_{i,\rm a})
\times
\prod_{{\rm b} \in {\rm Node}^+_{\partial}} [0,\varepsilon)_{\rm b}
\times
\prod_{i=1,2}\prod_{{\rm b} \in {\rm Node}^+_{i,\rm int}} D^2_{\rm b}(\varepsilon).
\end{equation}
(We restrict the gluing parameter so that the domain is smaller than \eqref{form1212}.)
We denote \eqref{form12122} by $\mathcal V(\xi',\mathscr{GL})$.

$V(\xi;\varepsilon)$ is a submanifold of the product of this space and
the other space which parametrizes the map. We define the latter space below.

For each disk component $\xi^{\prime,\rm d}_{\rm a}$ and
sphere component $\xi^{\prime,\rm s}_{i,\rm a}$, we consider
the set of maps
\[
u^{\heartsuit, \rm d}_{\rm a} \colon\
\bigl(\Sigma_{\rm a}^{\rm d},\partial \Sigma_{\rm a}^{\rm d},\vec z_{\rm a}\bigr)
\to (X_1 \times X_2,\gamma^* L_{12},L_{12}(\vec a_{\rm a})),\qquad
u^{\heartsuit, \rm s}_{i,\rm a} \colon\
\Sigma_{i,\rm a}^{\rm s} \to X_i
\]
(here the notation in the first line is as in Definition~\ref{defn1231}),
such that
$
\overline\partial u^{\heartsuit, \rm d}_{\rm a} \in \mathcal E(\xi')$,
\smash{$
\overline\partial u^{\heartsuit, \rm s}_{i,\rm a} \in \mathcal E(\xi')
$}
and that the $C^2$ distance between \smash{$u^{\heartsuit, \rm d}_{\rm a}$} \big(resp.\
\smash{$u^{\heartsuit, \rm s}_{i,\rm a}$}\big) and
\smash{$u^{\prime,\rm d}_{\rm a}$} \big(resp.\
\smash{$u^{\prime, \rm s}_{i,\rm a}$}\big) is smaller than $\varepsilon$.
Here~\smash{$u^{\prime, \rm d}_{\rm a}$},~\smash{$u^{\prime,\rm s}_{i,\rm a}$} are parts of
$\xi'$.

We denote the set of such maps $u^{\heartsuit, \rm d}_{\rm a}$ (resp.
$u^{\heartsuit, \rm s}_{i,\rm a}$)
by $\mathcal W^{\rm d}_{\rm a}(\xi';\varepsilon)$
\big(resp.\ $\mathcal W^{\rm s}_{i,\rm a}(\xi';\varepsilon)$\big)
and put
\[
\mathcal W^+(\xi';\varepsilon)
:=
\prod_{\rm a}\mathcal W^{\rm d}_{\rm a}(\xi';\varepsilon)
\times
\prod_{i=1,2}\prod_{\rm a}\mathcal W^{\rm s}_{i,\rm a}(\xi';\varepsilon).
\]
Here the first product is taken over disk components
and the second product is taken over sphere components.

We consider the direct product
\begin{equation}\label{form1254}
\prod_{i=1,2}\prod_{b \in {\rm node}^{+}_{i,\rm int}}(X_i)^2
\times
\prod_{b \in {\rm node}^{+}_{\partial}}(L_{12}(a_b))^2.
\end{equation}
Here $L_{12}(a_b)$ is as in \eqref{123131}.

Note each node is contained in exactly two irreducible
components. So using the evaluation maps on those components,
we define
$
{\rm ev}_{\rm node} \colon \mathcal W^+(\xi';\varepsilon) \to \eqref{form1254}$.
\begin{lem}
Let
\[
\Delta = \prod_{i=1,2}\prod_{b \in {\rm node}^{+}_{i,\rm int}}X_i
\times
\prod_{b \in {\rm node}^{+}_{\partial}}L_{12}(a_b)
\]
be the diagonal in \eqref{form1254}. Then the map
${\rm ev}_{\rm node}$ is transversal to $\Delta$ if
$\varepsilon$ is sufficiently small.
\end{lem}
\begin{proof}
This is a consequence of Condition \ref{defn12332},
which implies that ${\rm ev}_{\rm node}$ is transversal to $\Delta$
at $\bigl(\bigl(u^{\rm d}_{a}\bigr),(u^{\rm s}_{i,a})\bigr)$.
\end{proof}

\begin{defn}
We put
$
\mathcal W(\xi';\varepsilon) = {\rm ev}_{\rm node}^{-1}(\Delta)
\subset \mathcal W^+(\xi';\varepsilon)$.

\end{defn}
Note that $\Gamma_{\xi}$ the group of automorphisms of $\xi$ acts $\xi'$ as
a group of weak automorphisms. Then it acts on $\mathcal V(\xi',\mathscr{GL})$
and $\mathcal W^+(\xi';\varepsilon)$
by exchanging the factors. It then acts on $\mathcal W(\xi';\varepsilon)$.

The gluing construction proves the next proposition.

\begin{prop}\label{prop1242}
For each $\rho = \bigl(\bigl(\rho_a^{\rm d}\bigr),(\rho_{i,a}^{\rm s})\bigr) \in \mathcal W(\xi';\varepsilon)$
and $(\sigma,{\bf r}) \in \mathcal V(\xi',\mathscr{GL})$,
we obtain an object $\eta(\rho,\sigma,{\bf r})$
satisfying conditions in Definition {\rm\ref{defn12383}} except
\eqref{formtran1} and \eqref{formtran2} and whose source
object is $({\rm Glue}(\sigma,{\bf r}),\mathscr I_{\sigma,{\bf r}})$.

On the contrary, any object satisfying conditions in Definition {\rm\ref{defn12383}} except equation~\eqref{formtran1} and~\eqref{formtran2}
with sufficiently small $\varepsilon$ is equivalent to some $\eta(\rho,\sigma,{\bf r})$.

The isomorphism class of $\eta(\rho,\sigma,{\bf r})$ is the same as the
isomorphism class of $\eta(\gamma(\rho,\sigma,{\bf r}))$ for~${\gamma \in \Gamma_{\xi}}$.

\end{prop}
The proof is given in the next subsection.

We next cut down the space $\mathcal V(\xi',\mathscr{GL})
\times \mathcal W(\xi';\varepsilon)$ by conditions \eqref{formtran1} and \eqref{formtran2}.
We compose the map $(\rho,\sigma,{\bf r}) \mapsto \eta(\rho,\sigma,{\bf r})$ and
the evaluation maps at the interior marked points~\eqref{newform1211} to obtain
\begin{equation}\label{form1255form}
{\rm ev}_{\rm int} \colon\
\mathcal V(\xi',\mathscr{GL})
\times \mathcal W(\xi';\varepsilon)
\to (X_1 \times X_2)^{\ell} \times X_1^{\ell_1} \times X_2^{\ell_2}.
\end{equation}
The next proposition claims its smoothness.
We need carefully choose the smooth structure of~$\mathcal V(\xi',\mathscr{GL})$ so that ${\rm ev}_{\rm int}$ becomes a smooth map.

Let $r \in [0,\varepsilon)_{b}$ and $\mathfrak r \in D^2_{b}(\varepsilon)$
(see \eqref{form12122}).
We define $T$, $\theta$ by
\begin{equation}\label{form1256}
r = e^{-10\pi T},
\qquad
\mathfrak r = e^{-10\pi T + 2\pi\sqrt{-1}\theta}
\end{equation}
and put
$
s = 1/T$, \smash{$ \mathcal S = e^{2\pi\sqrt{-1}\theta}/T$}.
We use $s$ and $\mathcal S$ as coordinates in place of $r$ and $\mathfrak r$ to define a
smooth structure of $\mathcal V(\xi',\mathscr{GL})$.
\begin{prop}\label{prop1243}
When we put the above smooth structure on $\mathcal V(\xi',\mathscr{GL})$,
the map ${\rm ev}_{\rm int}$ in~\eqref{form1255form} is smooth.
\end{prop}
We will prove this proposition in the next subsection.
\begin{defn}
Let $V(\xi;\varepsilon)$ be the subset of
$\mathcal V(\xi',\mathscr{GL})
\times \mathcal W(\xi';\varepsilon)$ consisting of elements $(\rho,\sigma,{\bf r})$ such that
\[
{\rm ev}_{\rm int}(\rho,\sigma,{\bf r})
\in \prod_j \mathcal N^{(1)}_j \times \prod_j \mathcal N^{(2)}_j.
\]
Here $\mathcal N^{(1)}_j$ and $\mathcal N^{(2)}_j$ are as in
\eqref{formtran1} and \eqref{formtran2}.
The direct product in the first factor of the right-hand side is taken
over interior marked points on disk components and the direct product of
the second factor of the right-hand side is taken
over interior marked points on sphere components.
\end{defn}
\begin{cor}\label{cor1245}
If $\varepsilon$ is sufficiently small, then $V(\xi;\varepsilon)$
is a smooth submanifold of
$\mathcal V(\xi',\mathscr{GL})
\times \mathcal W(\xi';\varepsilon)$.
\end{cor}
\begin{proof}
In view of Proposition~\ref{prop1243}, it suffices to show that
${\rm ev}_{\rm int}$ is transversal to
$\smash{\prod_j \mathcal N^{(1)}_j} \times \smash{\prod_j \mathcal N^{(2)}_j}$
at $\bigl(\bigl(u^{\rm d}_{a}\bigr),(u^{\rm s}_{i,a})\bigr)$.
This is a consequence of Definition~\ref{defn1223}\,(5b) and (6b).
\end{proof}

From Corollary \ref{cor1245} and Proposition~\ref{prop1242}, it is easy
to see that there is a canonical isomorphism between
$V(\xi;\varepsilon)/\Gamma_{\xi}$ and $U(\xi;\varepsilon)$.
We thus have proved Proposition~\ref{prop1239}\,(1).

We next prove Proposition~\ref{prop1239}\,(2).
Let $\mathfrak i \in {\bf I}(\xi)$.
We take a source stabilization $\xi'_{\mathfrak i}$ of
 $\xi_{\mathfrak i}$ and a source stabilization $\xi'$ of $\xi$
 such that $\xi'$ is $\varepsilon$ close to $(\xi'_{\mathfrak i},\mathscr{GL})$.
 (Note $\xi'$ may depend on $\mathfrak i$.)
 We then fix a map
 \begin{gather}
\mathcal I^{\rm d}_{i,\rm a}({\bf 0}) \colon\ K^{\rm d}_{i,\rm a}(\xi'_{\mathfrak i}) \cong
K^{\rm d}_{i,\rm a}\bigl(\sigma^{\rm d}_{\rm a}({\bf 0})\bigr)
\to \Sigma_i(\sigma({\bf 0}),{\bf r}({\bf 0})) \cong \Sigma_i,\label{form12372}
\end{gather}
where $\Sigma_i$ is an irreducible component of $\xi$
(which is also an irreducible component of $\xi'$)
and~$(\sigma({\bf 0}),{\bf r}({\bf 0}))$ so that the source gluing map
at $\xi'_{\mathfrak i}$ sends $(\sigma({\bf 0}),{\bf r}({\bf 0}))$ to $\xi'$.
$K^{\rm d}_{i,\rm a}(\xi'_{\mathfrak i})$ is a~core of a disk component of the
source curve of $\xi'_{\mathfrak i}$.

Note that the image of \eqref{form12372} is in
a disk component $\Sigma^{\rm d}_{i,\rm a}$ of the source curve of $\xi'$.

Now let $(\sigma,{\bf r}) \in \mathcal V(\xi',\mathscr{GL})$.
The source curve of $\eta(\rho,\sigma,{\bf r})$ depends only on
$(\sigma,{\bf r})$ and is~$\varepsilon(\xi_{\mathfrak i},\mathscr{OB})$
close to $(\xi_{\mathfrak i},\mathscr{ST})$.
We write $\Sigma_i(\sigma,{\bf r})$ it.
Then we can uniquely choose
 \[
\mathcal I^{\rm d}_{i,\rm a}(\sigma,{\bf r}) \colon K^{\rm d}_{i,\rm a}(\xi'_{\mathfrak i})
\to \Sigma_i(\sigma,{\bf r})
\]
as \eqref{form1237}
which depends continuously on $(\sigma,{\bf r})$ and
becomes \eqref{form12372} when $(\sigma,{\bf r}) = (\bf 0,\bf 0)$.
We can choose a similar map for the sphere component
in the same way.

Using this choice, the map $
\mathfrak s_{\xi} \colon V(\xi;\varepsilon) \to
\mathcal E(\xi)
$ in Proposition~\ref{prop1239}\,(2) is continuous.

The $\Gamma_{\xi}$ equivalence of this map is
proved by $\Gamma_{\xi} \subseteq \Gamma_{\xi_{\mathfrak i}}$ and
$\Gamma_{\xi_{\mathfrak i}}$ invariance of various objects in the obstruction bundle
data.

The smoothness of $\mathfrak s_{\xi}$ follows from the exponential decay estimate
in the next subsection (see Proposition~\ref{propo1250}).

The proof of Proposition~\ref{prop1239}\,(3), (4) is now an immediate consequence
of the construction.
The proof of Proposition~\ref{prop1239} is complete modulo the points
postponed to the next subsection.\looseness=1
\end{proof}

We thus constructed a Kuranishi chart at each point of
${\mathcal M}'(L_{12};\vec a;E)$.
Let $\xi_1 \in {\mathcal M}'(L_{12};\vec a;E)$
and $\xi_2 \in {\mathcal M}'(L_{12};\vec a;E)$
is in the image of $\psi_{\xi_1}$.
Using the closedness of $\mathfrak N(\xi_1)$,
we may assume~${{\bf I}(\xi_2) \subseteq {\bf I}(\xi_1)}$,
by shrinking our Kuranishi neighborhood of $\xi_1$ if necessary.
Then by definition $U(\xi_2,\varepsilon_2) \subset U(\xi_1,\varepsilon_1)$
if we choose $\varepsilon_2$ sufficiently small.
We can use this fact and exponential decay estimate in the next subsection
to construct a smooth coordinate change from the Kuranishi neighborhood of
$\xi_2$ to one of $\xi_1$.
Thus we obtain the required Kuranishi structure.

\subsection{Gluing analysis for the construction of a Kuranishi
chart}
\label{sec:glueglue}
 In this subsection, we prove Propositions \ref{prop1242} and \ref{prop1243}.
 The proof is by gluing analysis similar to
\cite{fooobook2,foootech,foooanalysis}.
Since our compactification is slightly different from the stable
map compactification used
in those references, we explain the way we
modify the method of previous literatures so that
it works in our situation.
In \cite{fooobook2,foootech,foooanalysis},
a combination of the alternative method and the Newton's
iteration was used. We follow this method in this
subsection.
We follow \cite{foooanalysis} since the description is
the most detailed in this reference.
Below we provide the detail of the
formulation and the inductive scheme of the proof.
Once they are clarified the estimate, we need
on each step of the induction is entirely similar to \cite{foooanalysis}.

For the sake of simplicity of notation, we write the detail of our proof
in the following special case. This case contains all the points
we need to work out the general case.

We take $\Sigma_i^{\rm d} = D^2$ with one boundary marked point
$1 \in \partial D^2$ and two interior marked points~${0, \mathfrak z_i \in
{\rm Int}\,D^2}$, $0\ne \mathfrak z_i$.

We take $\Sigma_i^{\rm s} = S^2 = \C \cup \{\infty\}$
and consider three marked points $0$, $1$, $\infty$ on it.

We put $\Sigma^{\rm d} = \Sigma_1^{\rm d} = \Sigma_2^{\rm d}$.
We regard $\mathfrak z_1,\mathfrak z_2 \in \Sigma^{\rm d}$.
We glue $\Sigma_i^{\rm d}$ and $\Sigma_i^{\rm s}$ at $\mathfrak z_i
\in \Sigma_i^{\rm d}$ and $0 \in \Sigma_i^{\rm s}$.
The points $\mathfrak z_1$, $\mathfrak z_2$ may or may not coincide.
In case we use stable map compactification,
a~sphere component bubbles off when $\mathfrak z_1 = \mathfrak z_2$.
However, in our compactification,
the locus where~${\mathfrak z_1 = \mathfrak z_2}$ does {\it not} play
a special role.
\begin{figure}[ht]
\centering
\includegraphics[scale=0.36]{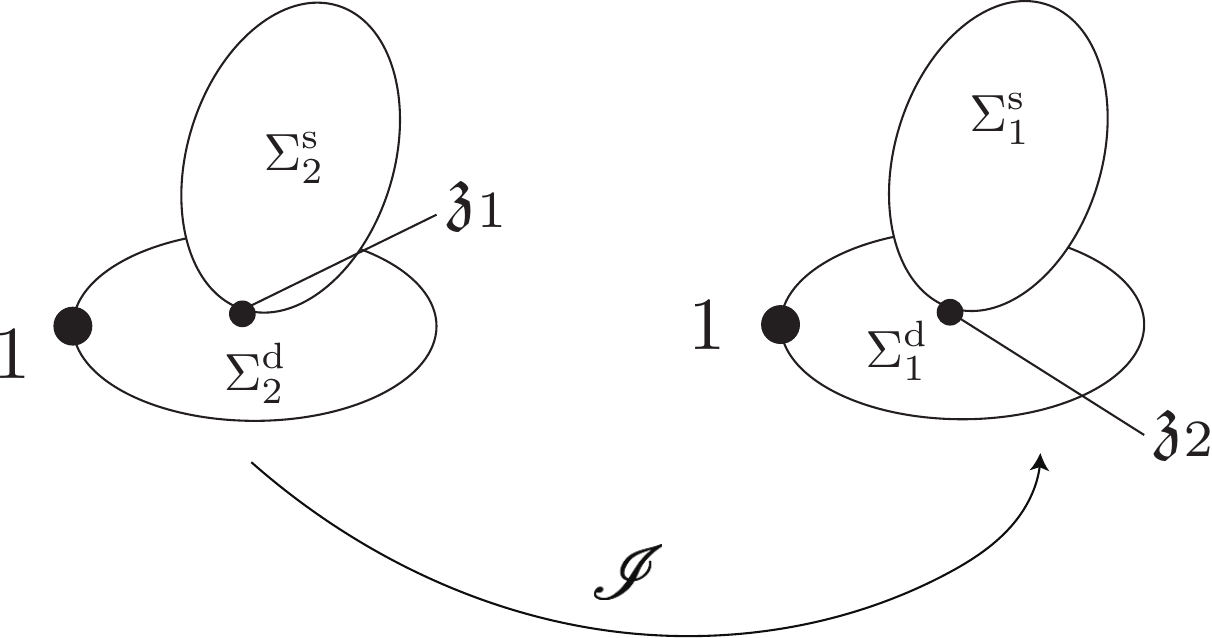}
\caption{The source domain we study. (We do not draw interior
marked points in the figure.)}
\label{Figure1241}
\end{figure}
We assume $L_{12} \subset -X_1 \times X_2$ is an {\it embedded}
Lagrangian submanifold. The immersed case can be worked out
in a similar way, given the formulation we have provided in the last
subsection. We assume $L_{12}$ is embedded for the sake of
simplicity of notation only.

We also remark that throughout this section,
we use almost complex structure $-J_{X_1}$ on $X_1$
and $J_{X_2}$ on $X_2$, unless otherwise mentioned explicitly.

We consider families of pseudoholomophic maps
\[
u^{{\rm d},\rho^{\rm d}}
= \bigl(u^{{\rm d},\rho^{\rm d}}_1,u^{{\rm d},\rho^{\rm d}}_2\bigr) \colon\ \bigl(\Sigma^{\rm d},\partial \Sigma^{\rm d}\bigr)
\to (-X_1 \times X_2,L_{12}),
\qquad
u^{{\rm s},\rho^{\rm s}_i}_i \colon\ \Sigma^{\rm s}_i \to X_i,
\]
parametrized by $\rho^{\rm d} \in \mathcal V^{\rm d}$,
$\rho^{\rm s}_i \in \mathcal V^{\rm s}_i$.
\begin{rem}\label{rem1248}
In the case we are studying here, there are two marked points and one nodal point
on the sphere bubble. We identify them as 0, 1, $\infty$. Therefore,
the domain coordinate is canonically determined.
In particular, the maps
\smash{$u^{{\rm d},\rho^{\rm d}}$}, \smash{$u^{{\rm s},\rho^{\rm s}_i}_i$}
are determined by $\rho^{\rm d}$ and $\rho^{\rm s}_i$ uniquely.
See Remark~\ref{rem1254}, for an explanation of the case when
the domain of the map
is not stable.

\end{rem}

Let $\mathfrak o \in {\rm Int}\,D^2 \setminus 0$ and
$\mathfrak O$ a small open neighborhood of it.
Let $\bf 0 \in \mathcal V^{\rm d}$ and
${\bf 0} \in \mathcal V^{\rm s}_i$.
We assume~\smash{$
u_i^{{\rm d},{\bf 0}}(\mathfrak o) = u_i^{{\rm s},{\bf 0}}(0),
$}
for $i=1,2$.

We consider the following element $\xi_0$ of
$\mathcal M'_{1,2,2}(L_{12};({\rm diag});E)$.
Here $\rm{diag}$ denotes the diagonal component of $L_{12} \times_{X_1\times X_2}L_{12}$
(which is actually the only component of it) and
\[
E = \sum_{i=1}^2 (-1)^{i}\int_{\Sigma^{\rm d}_i} \bigl(u^{{\rm d},{\bf 0}}_i\bigr)^* \omega_i
+ \sum_{i=1}^2 (-1)^{i}\int_{\Sigma^{\rm s}_i} \bigl(u^{{\rm s},{\bf 0}}_i\bigr)^* \omega_i.
\]
The source curve of $\xi_0$ is a pair of $\Sigma^{\rm d}_i$ and $\Sigma^{\rm s}_i$
glued at $\mathfrak o \in \Sigma^{\rm d}_i$ and $0 \in \Sigma^{\rm s}_i$.
The point $1 \in \Sigma^{\rm d}$ is its boundary marked point, the point $0 \in \Sigma^{\rm d}$ is an interior marked point of first kind, and the points~${1,\infty \in \Sigma^{\rm s}_i}$ are interior marked points of second kind.
$\mathscr I$ is the identity map $\Sigma^{\rm d}_1 = \Sigma^{\rm d} = \Sigma^{\rm d}_2$.
The maps $u_i$ are \smash{$u_i^{{\rm d},{\bf 0}}$}, \smash{$u^{{\rm s},{\bf 0}}_i$} on each of the components.

We study a neighborhood of $\xi_0$ in $\mathcal M'_{1,2,2}(L_{12};({\rm diag});E)$.

\begin{defn}
We consider the set $\mathcal V$ consisting of \smash{$\bigl(\rho^{\rm d},\rho^{\rm s}_1,
\rho^{\rm s}_2,\mathfrak z_1,\mathfrak z_2\bigr)$}
such that
\begin{enumerate}\itemsep=0pt
\item[(1)] $\rho^{\rm d} \in \mathcal V^{\rm d}$,
$\rho^{\rm s}_i \in \mathcal V^{\rm s}_i$.
\item[(2)]
$\mathfrak z_1, \mathfrak z_2 \in \mathfrak O$.
\item[(3)]
\smash{$
u_i^{{\rm d},\rho^{\rm d}}(\mathfrak z_i) = u_i^{{\rm s},\rho^{\rm s}_i }(0)
$}
for $i=1,2$.
\end{enumerate}

\end{defn}
We may regard $\mathcal V$ as a fiber product
\begin{equation}\label{form1261}
\mathcal V
= \bigl(\mathcal V^{\rm d} \times \mathfrak O \times \mathfrak O\bigr)
\times_{X_1 \times X_2}
(\mathcal V^{\rm s}_1
\times \mathcal V^{\rm s}_2).
\end{equation}
We work under the following assumptions.
\begin{assump}\label{ass1247}
\quad
\begin{enumerate}\itemsep=0pt
\item[(1)]
$\mathcal V^{\rm d}$, $\mathcal V^{\rm s}_1$
and $\mathcal V^{\rm s}_2$ are smooth manifolds.
Moreover, the linearizations of the nonlinear Cauchy--Riemann equations
\begin{gather*}
D_{u^{{\rm d},\rho^{\rm d}}}\overline\partial \colon\
C^{\infty}\bigl(\bigl(\Sigma^{\rm d};\partial \Sigma^{\rm d}\bigr);\bigl(
\bigl(u^{{\rm d},\rho^{\rm d}})^*T(X_1 \times X_2\bigr),\bigl(u^{{\rm d},\rho^{\rm d}}\bigr)^*TL_{12}\bigr)\bigr)
\\
\hphantom{D_{u^{{\rm d},\rho^{\rm d}}}\overline\partial \colon} \ \to
C^{\infty}\bigl(\Sigma^{\rm d},
\bigl(u^{{\rm d},\rho^{\rm d}}\bigr)^*T(X_1 \times X_2) \otimes \Lambda^{0,1}\bigr)
\end{gather*}
and
\[
D_{u_i^{{\rm s},\rho^{\rm s}_i}}\overline\partial \colon\
C^{\infty}\bigl(\Sigma^{{\rm s}}_i;
\bigl(u_i^{{\rm s},\rho^{{\rm s}}_i}\bigr)^*TX_i\bigr)
\to
C^{\infty}\bigl(\Sigma^{\rm d},
\bigl(u^{{\rm s},\rho^{{\rm s}}_i}\bigr)^*TX_i \otimes \Lambda^{0,1}\bigr)
\]
are surjective.
\item[(2)]
The fiber product \eqref{form1261} is transversal.
\end{enumerate}

\end{assump}
Note that Assumption \ref{ass1247} implies that $\mathcal V$ is a
smooth manifold.
\begin{rem}
In the general situation, we introduce obstruction bundles and
use the extended moduli space in place of
$\mathcal V^{\rm d}$, $\mathcal V^{\rm s}_i$,
so that a similar condition as Assumption~\ref{ass1247}
holds. The way to introduce obstruction bundles is
explained in detail in the last subsection.
Then the way to include the obstruction bundle
in the gluing analysis is the same as
\cite{foooanalysis} etc.
So, for the sake of simplicity of notation, we restrict ourselves to the case
when Assumption~\ref{ass1247} is satisfied in this subsection.
\end{rem}

We next recall the source gluing map
in our situation.

Let $c$ be a small positive number.
We define a map $\varphi^{\rm d}_{\mathfrak z_i}
\colon D^2 \to \Sigma^{\rm d}$ by
$
\varphi^{\rm d}_{\mathfrak z_i}(z) = cz + \mathfrak z_i$.
We also define
$\varphi^{\rm s}_{i}
\colon D^2 \to \Sigma^{\rm s}_i$ by
$
\varphi^{\rm s}_{i}(z) = cz$.
They are analytic families of coordinates
and are extendable. We use them to define
the source gluing map \eqref{form12113}.
Using also Lemma~\ref{loem1229},
we obtain
\begin{equation}\label{form1264}
(\Sigma_1(\mathfrak z_1,\mathfrak r_1),\Sigma_2(\mathfrak z_2,\mathfrak r_2),
\mathscr I_{\mathfrak z_1,\mathfrak z_2,\mathfrak r_1,\mathfrak r_2}),
\end{equation}
for $\mathfrak r_1,\mathfrak r_2 \in D^2$.
Here $\Sigma_i(\mathfrak z_i,\mathfrak r_i)$ is obtained
by gluing $\Sigma_i^{\rm d}$ and $\Sigma_i^{\rm s}$
by the gluing parameter $\mathfrak r_i$ using coordinates
$\varphi^{\rm d}_{\mathfrak z_i}$ and $\varphi^{\rm s}_{i}$.
$\mathscr I_{\mathfrak z_1,\mathfrak z_2,\mathfrak r_1,\mathfrak r_2}
 \colon\Sigma_1^0(\mathfrak z_1,\mathfrak r_1) \to \Sigma_2^0(\mathfrak z_2,\mathfrak r_2)$
 is a biholomorphic map obtained in Lemma~\ref{loem1229}
 by extending identity map.
\begin{prop}\label{prop1249}
We assume Assumption {\rm\ref{ass1247}}.
Then for sufficiently small $\varepsilon$ there exists a~map
\[
\mathscr G \colon\ \mathcal V \times D^2(\varepsilon) \times
D^2(\varepsilon)
\to \mathcal M'_{1,2,2}(L_{12};({\rm diag});E)
\]
with the following properties:
\begin{enumerate}\itemsep=0pt
\item[$(1)$]
The source object of $\mathscr G\bigl(\bigl(\rho^{\rm d},\rho^{\rm s}_1,
\rho^{\rm s}_2,\mathfrak z_1,\mathfrak z_2\bigr),\bigl(\mathfrak r_1,\mathfrak r_2\bigr)\bigr)$
is \eqref{form1264}.
\item[$(2)$]
$\mathscr G$ is a homeomorphism onto a
neighborhood of $\xi_0$.
\end{enumerate}

\end{prop}
Proposition~\ref{prop1249} is a special case of
Proposition~\ref{prop1239}\,(1).
To prove Propositions~\ref{prop1239}\,(2), (3),~(4),
\ref{prop1242}, \ref{prop1243} and the smoothness of
coordinate change, we use the next Proposition~\ref{propo1250}.
To state it we need notations.

We take a small open set $\mathfrak O^+ \subset \Sigma^{\rm d}$ which contains the
closure of $c$-neighborhood of $\mathfrak O$.
We put
$
K_{1,-}^{\rm d} = K_{2,-}^{\rm d} = K^{\rm d}_-
=
\Sigma^{\rm d} \setminus \mathfrak O^+
$
and
$
K_i^{\rm s}
=
\Sigma_i^{\rm s} \setminus \operatorname{Im}\,\varphi^{\rm s}_i$.
We may regard
$
K_{i,-}^{\rm d}, K_i^{\rm s} \subset \Sigma_i(\mathfrak r)
$
for all $\mathfrak r$.

Let
\smash{$
u_i^{(\rho^{\rm d},\rho^{\rm s}_1,
\rho^{\rm s}_2,\mathfrak z_1,\mathfrak z_2),(\mathfrak r_1,\mathfrak r_2)}
\colon \Sigma_i(\mathfrak r) \to X_i
$}
be the map part of
$\mathscr G\bigl(\bigl(\rho^{\rm d},\rho^{\rm s}_1,
\rho^{\rm s}_2,\mathfrak z_1,\mathfrak z_2\bigr),(\mathfrak r_1,\mathfrak r_2)\bigr)$.
We denote its restriction to
$K_{i,-}^{\rm d}$, $K_i^{\rm s}$ by
\begin{gather*}
\begin{split}
& \bigl({\rm Res}^{\rm d}_{i,-} \circ \mathscr G\bigr)\bigl(\bigl(\rho^{\rm d},\rho^{\rm s}_1,
\rho^{\rm s}_2,\mathfrak z_1,\mathfrak z_2\bigr),(\mathfrak r_1,\mathfrak r_2)\bigr)
\in C^{\infty}(K_{i,-}^{\rm d},X_i), \\
& ({\rm Res}^{\rm s}_i \circ \mathscr G)\bigl(\bigl(\rho^{\rm d},\rho^{\rm s}_1,
\rho^{\rm s}_2,\mathfrak z_1,\mathfrak z_2\bigr),(\mathfrak r_1,\mathfrak r_2)\bigr)
\in C^{\infty}(K_i^{\rm s},X_i)
\end{split}
\end{gather*}
and
\begin{gather*}
\bigl({\rm Res}^{\rm d}_- \circ \mathscr G\bigr)\bigl(\bigl(\rho^{\rm d},\rho^{\rm s}_1,
\rho^{\rm s}_2,\mathfrak z_1,\mathfrak z_2\bigr),(\mathfrak r_1,\mathfrak r_2)\bigr)\\
\qquad=
\bigl(\bigl({\rm Res}^{\rm d}_{1,-} \circ \mathscr G\bigr),
\bigl({\rm Res}^{\rm d}_{2,-} \circ \mathscr G\bigr)\bigr)
\bigl(\bigl(\rho^{\rm d},\rho^{\rm s}_1,
\rho^{\rm s}_2,\mathfrak z_1,\mathfrak z_2\bigr),(\mathfrak r_1,\mathfrak r_2)\bigr)
\in
C^{\infty}\bigl(K^{\rm d}_-,X_1 \times X_2\bigr).
\end{gather*}
We define $T_i \in \R_+$, $\theta_i \in [0,1]$, by
\begin{equation}\label{form1270}
\mathfrak r_i = \exp\bigl(-10\pi T_i - 2\pi\sqrt{-1}\theta_i\bigr).
\end{equation}
\begin{prop}\label{propo1250}
For $m > 10$,
there exists $\varepsilon_{m,n} > 0$ and $C_{m,n}, c_{m,n} > 0$ such that the following
holds if $\varepsilon < \varepsilon_{m,n}$
$($note $\varepsilon$ is the number in Proposition {\rm\ref{prop1249}}$)$:
\begin{enumerate}\itemsep=0pt
\item[$(1)$]
\[
\left\Vert
\nabla^n_{(\rho^{\rm d},\rho^{\rm s}_1,
\rho^{\rm s}_2,\mathfrak z_1,\mathfrak z_2)}
\frac{\partial^{\ell_1}}{\partial T_1^{\ell_1}}
\frac{\partial^{\ell'_1}}{\partial \theta_1^{\ell'_1}}
\frac{\partial^{\ell_2}}{\partial T_2^{\ell_2}}
\frac{\partial^{\ell'_2}}{\partial \theta_2^{\ell'_2}}
\bigl({\rm Res}^{\rm d}_- \circ \mathscr G\bigr)
\right\Vert_{L^2_{m-\ell}}
\le
C_{m,n} e^{-c_{m,n} T_1}
\]
if $\ell = \ell_1 +\ell'_1 + \ell_2 + \ell'_2 \le m-2$,
$\ell_1,\ell'_1,\ell_2, \ell'_2 \in \Z_{\ge 0}$
and $\ell_1 + \ell'_1 > 0$.
Here \smash{$\nabla^n_{(\rho^{\rm d},\rho^{\rm s}_1,
\rho^{\rm s}_2,\mathfrak z_1,\mathfrak z_2)}$} is $n$-th derivative
with respect to $\bigl(\rho^{\rm d},\rho^{\rm s}_1,
\rho^{\rm s}_2,\mathfrak z_1,\mathfrak z_2\bigr)$.
\item[$(2)$]
\[
\left\Vert
\nabla^n_{(\rho^{\rm d},\rho^{\rm s}_1,
\rho^{\rm s}_2,\mathfrak z_1,\mathfrak z_2)}
\frac{\partial^{\ell_1}}{\partial T_1^{\ell_1}}
\frac{\partial^{\ell'_1}}{\partial \theta_1^{\ell'_1}}
\frac{\partial^{\ell_2}}{\partial T_2^{\ell_2}}
\frac{\partial^{\ell'_2}}{\partial \theta_2^{\ell'_2}}
\bigl({\rm Res}^{\rm d}_- \circ \mathscr G\bigr)
\right\Vert_{L^2_{m-\ell}}
\le
C_{m,n} e^{-c_{m,n} T_2}
\]
if $\ell = \ell_1 +\ell'_1 + \ell_2 + \ell'_2 \le m-2$,
$\ell_1,\ell'_1,\ell_2, \ell'_2 \in \Z_{\ge 0}$
and $\ell_2 + \ell'_2 > 0$.
\item[$(3)$]
The same inequality as $(1)$, $(2)$ holds for
${\rm Res}^{\rm s}_i \circ \mathscr G$, $i=1,2$.
\end{enumerate}

\end{prop}
We can use the exponential decay estimate such as
Proposition~\ref{propo1250}
in the same way as~\cite[Chapter 8]{foooanalysis}
to prove Propositions \ref{prop1242}, \ref{prop1239}\,(2), (3), (4),
Proposition~\ref{prop1243} and the smoothness of
coordinate changes.
So to complete the proof of Theorem~\ref{prop1417},
it remains to prove Propositions~\ref{prop1249} and \ref{propo1250}.
The rest of this subsection is occupied by their proofs.
\end{proof}

\begin{rem}\label{rem1254}
As we mentioned in Remark~\ref{rem1248}, we study the case when there are
marked points on the sphere bubbles
so that the domain is stable in Propositions \ref{prop1249} and \ref{propo1250}.
In the general case, we follow the method of \cite[Appendix]{FO} and proceed as follows.
(This is a special case of the method we explained in the last subsection.)
Suppose we consider an element $\xi'_0$ of~$\mathcal M'_{1,0,0}(L_{12};({\rm diag});E)$ which is similar to the element
$\xi_0$ except we forget the 4 marked points on the sphere bubbles.
We consider the case when the maps \smash{$u_i^{s;{\bf 0}}$} on the sphere bubbles
which are parts of the data consisting $\xi'_0$ is non-constant.
We fix two points on each of the sphere bubbles such that
$u_i^{s;{\bf 0}}$ is an immersion at those points.
We change the objects by automorphisms so that the
marked points we add are $1,\infty \in S^2$.
We denote by $1_i$, $\infty_i$ ($i=1,2$) those added marked
points (of second kind) on the sphere bubbles $S^2_i$.
The nodal points on the sphere bubbles are identified with $0$.
We take codimension 2 submanifolds $\mathcal N_{i,1}$, $\mathcal N_{i,\infty}$
of $X_i$ which
intersects with the image of the map \smash{$u_i^{s;{\bf 0}}$} transversally
at $1_{i}$ and ${\infty}_i$.

We consider $\xi_0'$ with those extra four marked points added as an
element $\xi_0$ of the space $\mathcal M'_{1,2,2}(L_{12};({\rm diag});E)$.
We can then apply Propositions \ref{prop1249} and \ref{propo1250} to obtain a map
\[
\mathscr G \colon\ \mathcal V \times D^2(\varepsilon) \times
D^2(\varepsilon)
\to \mathcal M'_{1,2,2}(L_{12};({\rm diag});E).
\]
Then the Kuranishi neighborhood of
$\xi'_0$ of
$\mathcal M'_{1,0,0}(L_{12};({\rm diag});E)$
is the smooth submanifold of~$\mathcal V \times D^2(\varepsilon) \times
D^2(\varepsilon)$
which is cut out by the conditions
\begin{equation}\label{transversalconst}
({\rm ev}_{i,1}\circ \mathscr G)(x) \in \mathcal N_{i,1},
\qquad
({\rm ev}_{i,\infty}\circ \mathscr G)(x) \in \mathcal N_{i,\infty}.
\end{equation}
Here ${\rm ev}_{i,1}$, ${\rm ev}_{i,\infty}$ are the evaluation maps:
$ \mathcal M'_{1,2,2}(L_{12};({\rm diag});E) \to X_i$ at the
marked points corresponding to $1_{i}$ and ${\infty}_i$.

We remark that the Kuranishi neighborhood of $\xi'_0$ obtained in this way
depends on the choice of additional 4 marked points on the sphere bubbles
and also to the choice of transversals~$\mathcal N_{i,1}$,~$\mathcal N_{i,\infty}$.
However, using Proposition~\ref{propo1250} we can show the Kuranishi
neighborhood obtained is independent of such choices in a neighborhood of
$\xi_0'$ up to diffeomorphism. This independence is a special case of
the smoothness of the coordinate change, which is proved by using
Proposition~\ref{propo1250}. See \cite[Chapter 8]{foooanalysis}.

We will discuss this example more in Remark~\ref{rem1274}.
\end{rem}
\begin{proof}[Proof of Propositions \ref{prop1249} and \ref{propo1250}]
Proposition~\ref{prop1249} is similar to
\cite[Theorem 3.13]{foooanalysis} and
Proposition~\ref{propo1250} is similar
to \cite[Theorem 6.4]{foooanalysis}.
Their proofs are also similar.

We first modify the way to describe the disk component
of the source curve in a way convenient for our gluing analysis.
\begin{defn}
We take a $\mathfrak z \in \mathfrak O$ parametrized smooth family of
diffeomorphisms $h_{\mathfrak z} \colon D^2 \to D^2$ with the following
properties:
\begin{enumerate}\itemsep=0pt
\item[(1)]
$h_{\mathfrak z} =$ the identity map outside $\mathfrak O^+$.
Here $\mathfrak O^+$ is an open subset of $D^2$
which contains the closure of $\mathfrak O$ and
is disjoint from $\{0\} \cup \partial D^2$.
\item[(2)]
$h_{\mathfrak z} \circ \varphi^{\rm d}_{\mathfrak o} = \varphi^{\rm d}_{\mathfrak z}$.
In particular, $h_{\mathfrak z}(\mathfrak o) = \mathfrak z$.
\end{enumerate}
We pull back the standard complex structure $j$ of $D^2$ by $h_{\mathfrak z}$
to obtain
$j_{\mathfrak z} = h_{\mathfrak z}^* j$.
\end{defn}

We remark $\bigl(\bigl(D^2,j\bigr),(1,0,\mathfrak z)\bigr)$ is isomorphic to $\bigl(\bigl(D^2,j_{\mathfrak z}\bigr),(1,0,\mathfrak o)\bigr)$.
In other words, we move a~complex structure $j$ in place of moving a marked point $\mathfrak z$.
In this identification, the map $\mathscr I$ becomes
$
\mathscr I_{\mathfrak z_1,\mathfrak z_2} = h_{\mathfrak z_2} \circ (h_{\mathfrak z_1})^{-1}$.
We put
\smash{$
u^{{\rm d},\rho^{\rm d},\mathfrak z_i}_i = u^{{\rm d},\rho^{\rm d}}_i \circ h_{\mathfrak z_i}
$}
and \smash{$u^{{\rm d},\rho^{\rm d},\mathfrak z} =
\bigl(u^{{\rm d},\rho^{\rm d},\mathfrak z_1}_1,u^{{\rm d},\rho^{\rm d},\mathfrak z_2}_2\bigr)$}.
The map \smash{$u^{{\rm d},\rho^{\rm d},\mathfrak z_i}_i$} is holomorphic with respect to the complex
structure $j_{\mathfrak z_i}$ of the source.

Hereafter, to simplify the notation we write $\rho = \bigl(\rho^{\rm d},\rho^{\rm s},\mathfrak z_1,\mathfrak z_2\bigr)$
and write $u^{{\rm d},\rho}_i$ etc.\ in place of~\smash{$u^{{\rm d},\rho^{\rm d},\mathfrak z}_i$} etc.

We remark
$
\mathcal V \cong
\bigl\{\rho \mid u^{{\rm d},\rho}_i(\mathfrak o) =
u^{{\rm s},\rho}_i(0) \ \text{for $i=1,2$}\bigr\}$.

We use the cylindrical coordinate on neighborhoods of $\mathfrak o \in \Sigma^{\rm d}$
and of $0 \in \Sigma^{\rm s}_i$, which we describe below.

Hereafter, we write $\varphi^{\rm d}$ in place of $\varphi^{\rm d}_{\mathfrak o}$.
Let $z \in D^2$ and $p = \varphi^{\rm d}(z) \in \Sigma^{\rm d}$.
We then define $\tau'(p) \in [0,\infty)$ and $t'(p) \in [0,1)$ by
$
2\pi \bigl(\tau'(p) + \sqrt{-1}t'(p)\bigr) = - \log z$.
Let $q_i = \varphi_i^{\rm s}(w_i) \in \Sigma_i^{\rm s}$
We then define $\tau''_i(q_i) \in (-\infty,0]$ and $t''_i(q_i) \in [0,1)$ by~${
2\pi \bigl(\tau''_i(q_i) + \sqrt{-1}t''_i(q_i)\bigr) = \log w_i}$.
We glue $\Sigma^{\rm d}$ with $\Sigma^{\rm s}_i$
by the gluing parameter $\mathfrak r_i$ as follows.
If $p = \varphi^{\rm d}(z)$ and $q_i = \varphi_i^{\rm s}(w_i)$,
we identify $p$ and $q$ if and only if
\begin{equation}\label{form1277}
z w_i = \mathfrak r_i.
\end{equation}
See Definition~\ref{defn1225}.
In view of \eqref{form1270}, the condition \eqref{form1277} is
equivalent to
\begin{equation}
\tau''_i - \tau' = 10T_i, \qquad
t''_i - t' \equiv \theta_i \mod \Z.\label{form12798}
\end{equation}
Compare \cite[equations~(6.2) and (6.3)]{foooanalysis} and see Figure~\ref{Figure12-78}.
\begin{figure}[ht]
\centering
\includegraphics[scale=0.3]{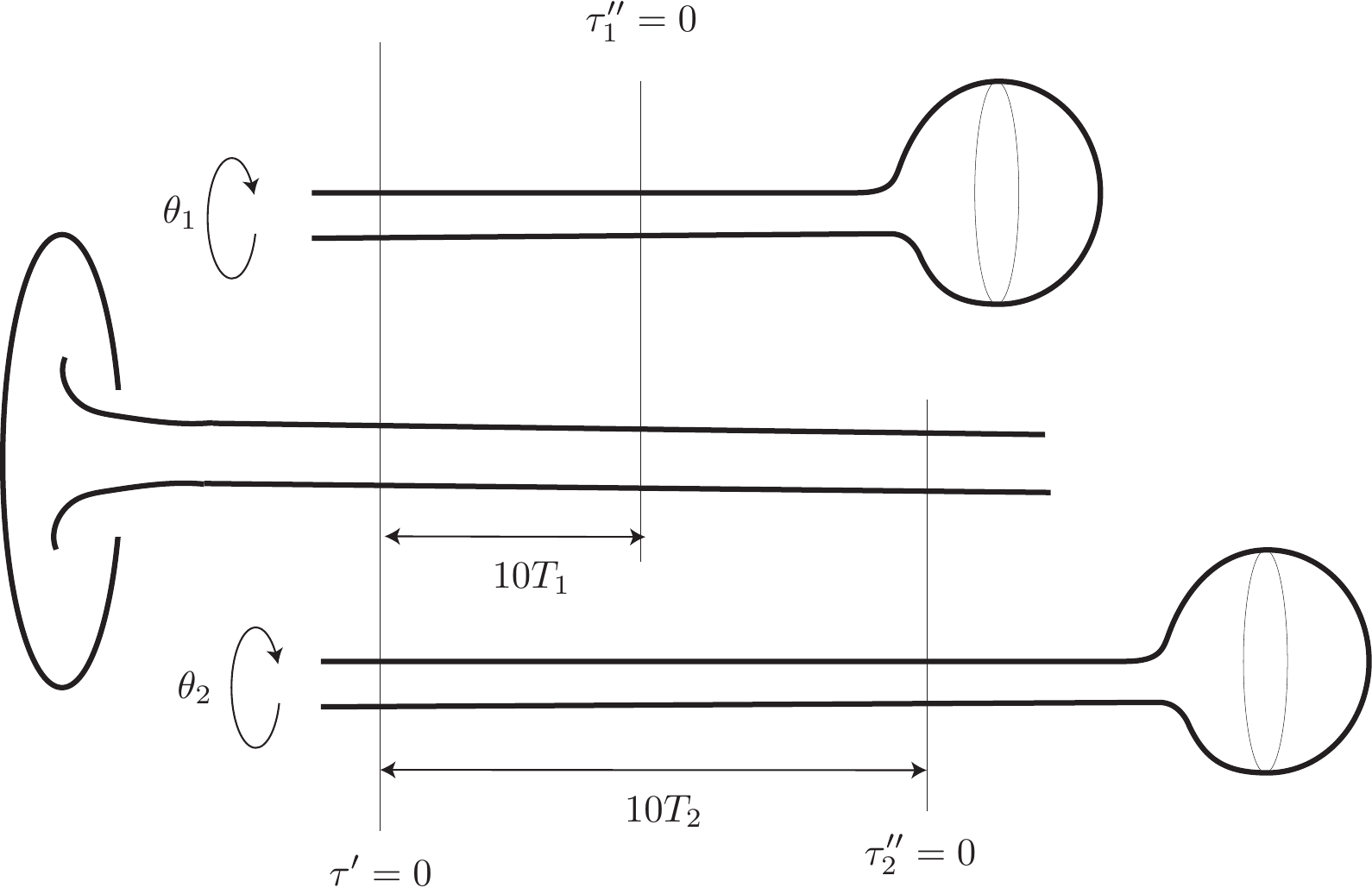}
\caption{Gluing disk and 2 spheres.}
\label{Figure12-78}
\end{figure}
We use Riemannian metric on $\Sigma^{\rm d} \setminus \{\mathfrak o\}$
(resp.~$\Sigma^{\rm s}_i \setminus \{0\}$) such that
on the image of $\varphi^{\rm d}$ (resp.~$\varphi^{\rm s}_i$)
it is isometric to $[0,\infty) \times S^1$
\big(resp.\ $(-\infty,0] \times S^1$\big) with
$(\tau',t')$ (resp.\ $(\tau''_i,t''_i)$) as coordinates.

We introduce the weighted Sobolev spaces which we use for our gluing analysis.
We follow~\cite[Section 3]{foooanalysis} here.
For $\rho \in \mathcal V$, we put $u^{\rm d,\rho}(\mathfrak o) = p^{\rho}$.
We take sufficiently small positive number~$\delta$ and fix it.
($\delta$ is taken to be small compared to the decay rate of the
pseudo-holomorphic curve at the neck. For example, $\delta < 1/100$.
See, for example, \cite[Section 2]{foooanalysis}.)

We take and fix connections of $X_1$ and of $X_2$ and
then direct product connection of~${X_1 \times X_2}$.
Let ${\rm Pal}_{\mathfrak o}$ be the parallel transport of the
tangent bundle of $X_1 \times X_2$ with respect to
this connection. We denote by the same symbol
the parallel transport of the
tangent bundle of $X_i$.
We denote by ${\rm Pal}^J_{\mathfrak o}$
the complex linear part of it.
(We remark that the almost complex
structure we use is $-J_{X_1} \oplus J_{X_2}$.)

\begin{defn}\label{defn1252}
We denote by $W^2_{m+1,\delta}\bigl(\bigl(\Sigma^{\rm d};\partial \Sigma^{\rm d}\bigr),
\bigl(u^{\rm d,\rho}\bigr)^*T(X_1\times X_2);\bigl(u^{\rm d,\rho}\bigr)^*T(L_{12})\bigr)$
the set of all pairs $(s,v)$ such that
\begin{enumerate}\itemsep=0pt
\item[(1)] $s$ is a section of $\bigl(u^{\rm d,\rho}\bigr)^*T(X_1\times X_2)$ on $\Sigma^{\rm d}
\setminus \{\mathfrak o\}$ which is locally of $L^2_{m+1}$ class.
\item[(2)]
$v \in T_{p^{\rho}}(X_1 \times X_2)$.
\item[(3)]
$s(z) \in T_{u^{\rm d,\rho}(z)}L_{12}$ if $z \in \partial \Sigma^{\rm d}$.
\item[(4)]
\begin{equation}\label{form1279}
\sum_{k=0}^{m+1}\int_{0}^{\infty}{\rm d}\tau' \int_{S^1}
e^{2\delta \tau'}
\vert
\nabla^k\bigl(s - v^{\rm pal}\bigr)
\vert^2 {\rm d}t'
< \infty.
\end{equation}
Here \smash{$v^{\rm pal}(\tau',t') = ({\rm Pal}_{\mathfrak o})_{\mathfrak o}^{u^{\rm d,\rho}(\tau',t')}(v)$}.
\end{enumerate}
The $W^{2}_{m+1,\delta}$ norm of $(s,v)$ is by definition
\[
\Vert (s,v) \Vert^2_{W^{2}_{m+1,\delta}}
:=
\sum_{k=0}^{m+1}\int_{\Sigma^{\rm d}\setminus \varphi^{\rm d}_{\mathfrak o}(D^2)}
\vert
\nabla^k(s)
\vert^2
+
\eqref{form1279}
+
\vert v\vert^2.
\]
We define the $L^2$ inner product between two elements $(s_1,v_1)$ and $(s_2,v_2)$
of the function space~${W^2_{m+1,\delta}\bigl(\bigl(\Sigma^{\rm d};\partial \Sigma^{\rm d}\bigr),
\bigl(u^{\rm d,\rho}\bigr)^*T(X_1\times X_2);\bigl(u^{\rm d,\rho}\bigr)^*T(L_{12})\bigr)}$
by
\[
\langle\!\langle (s_1,v_1),(s_2,v_2)\rangle\!\rangle_{L^2}
=
 \int_{[0,\infty)\times S^1} \bigl(s_1-v_1^{\rm Pal},s_2-v^{\rm Pal}_2\bigr)
+\int_{\Sigma^{\rm d}\setminus \varphi^{\rm d}_{\mathfrak o}(D^2)} (s_1,s_2)
+
( v_1,v_2).
\]

We denote by
$W^2_{m+1,\delta}\bigl((\Sigma^{\rm s}),
\bigl(u^{\rm s,\rho}_i\bigr)^*TX_i\bigr)$
the set of all pairs $(s,v)$ such that
\begin{enumerate}\itemsep=0pt
\item[(1)] $s$ is a section of $(u^{\rm s,\rho})^*TX_i$ on $\Sigma^{\rm s}_i
\setminus \{0\}$ which is locally of $L^2_{m+1}$ class.
\item[(2)]
$v \in T_{p^{\rho}}X_i$.
\item[(3)]
\[
\sum_{k=0}^{m+1}\int_{-\infty}^{0} {\rm d}\tau'' \int_{t'' \in S^1}
e^{-2\delta \tau'}
\vert
\nabla^k\bigl(s - v^{\rm pal}\bigr)
\vert^2 {\rm d}t''
< \infty.
\]
Here \smash{$v^{\rm pal}(\tau'',t'') = ({\rm Pal}_{\mathfrak o})_{\mathfrak o}^{u^{\rm s,\rho}_i(\tau'',t'')}(v)$}.
\end{enumerate}

The $W^{2}_{m+1,\delta}$ norm and the $L^2$ inner product is defined in a similar way.

\end{defn}
\begin{defn}\label{defn1253}
We denote by $L^2_{m,\delta}\bigl(\Sigma^{\rm d}_i,
\bigl(u^{\rm d,\rho}_i\bigr)^*TX_i \otimes \Lambda^{0,1}_{\rho}\bigr)$
the set of all $s$ such that
\begin{enumerate}\itemsep=0pt
\item[(1)]
 $s$ is a section of $\bigl(u^{\rm d,\rho}_i\bigr)^*TX_i \otimes \Lambda^{0,1}\bigl(\Sigma^{\rm d},j_{\mathfrak z_i}\bigr)$ on $\Sigma^{\rm d}
\setminus \{\mathfrak o\}$ which is locally of $L^2_{m}$ class.
Note that we use the complex structure $j_{\mathfrak z_i}$ to define the notion of $(0,1)$ forms on $\Sigma^{\rm d}$.
\item[(2)]
\begin{equation}\label{form1283}
\sum_{k=0}^{m}\int_{0}^{\infty}{\rm d}\tau' \int_{S^1}
e^{2\delta \tau'}
\bigl\vert
\nabla^k s
\bigr\vert^2 {\rm d}t'
< \infty.
\end{equation}
\end{enumerate}
The square of the $L^2_m$ norm of $s$ is by definition the sum of \eqref{form1283}
and the square of $L^2_m$ norm of the restriction of $s$ to
$\Sigma^{\rm d}\setminus \varphi^{\rm d}\bigl(D^2\bigr)$.

The weighted Sobolev space
$L^2_{m,\delta}\bigl(\Sigma^{\rm s}_i,
\bigl(u^{\rm s,\rho}_i\bigr)^*TX_i \otimes \Lambda^{0,1}\bigr)$
of sections of
$\bigl(u^{\rm s,\rho}_i\bigr)^*TX_i \otimes \Lambda^{0,1}(\Sigma^{\rm s},j)$
and its \smash{$L^2_{m,\delta}$} norm is defined in a similar way.

The direct sum
\[
\bigoplus_{i=1,2}
L^2_{m,\delta}\bigl(\Sigma^{\rm d}_i,
\bigl(u^{\rm d,\rho}_i\bigr)^*TX_i \otimes \Lambda^{0,1}_{\rho}\bigr)
\]
is denoted by
$L^2_{m,\delta}\bigl(\Sigma^{\rm d},
\bigl(u^{\rm d,\rho}\bigr)^*T(X_1 \oplus X_2) \otimes \Lambda^{0,1}_{\rho}\bigr)$,
by a slight abuse of notation.
(Note that the complex structure we use for \smash{$\Sigma^{\rm d}$}
is different between $X_1$ factor and $X_2$ factor.)
\end{defn}
We next define the linearization operator of
the nonlinear Cauchy--Riemann equation.

We use the parallel transport and the exponential map for this purpose.
\begin{defn}\label{defn1254}
We take a $z \in \Sigma^{\rm d}$ depending family of
connections $\nabla^{z}$ of the tangent bundle of $X_1 \times X_2$ such that
\begin{enumerate}\itemsep=0pt
\item[(1)] If $z \in \mathfrak O^+$, then $\nabla^{z}$ coincides with direct product
connection mentioned right above Definition~\ref{defn1252}.
\item[(2)]
There exists a neighborhood of $\partial \Sigma^{\rm d}$ such that if $z$ is in this neighborhood
then $\nabla^{z}$ coincides with a connection $\nabla^0$ for which
$L_{12}$ is totally geodesic.
\end{enumerate}
Let \index[syindex]{Expz@${\rm Exp}^z$}
\begin{equation}\label{form1284}
{\rm Exp}^z\colon\ T(X_1 \times X_2) \to (X_1 \times X_2)^2
\end{equation}
be the exponential map defined by $\nabla^{z}$.

 \index[syindex]{Expi@${\rm Exp}_i$}

If $z \in \mathfrak O^+$, item (1) implies that \eqref{form1284} becomes a direct product of two exponential
maps $
{\rm Exp}_i \colon TX_i \to (X_i)^2$.
The restriction of the exponential maps are diffeomorphisms onto a neighborhood of the diagonal,
which contain $U(\Delta_{X_i})$ etc.
We denote their inverses by\index[syindex]{Ez@${\rm E}^z$}\index[syindex]{Ei@${\rm E}_i$}
\begin{gather*}
{\rm E}^z \colon\ U(\Delta_{X_1 \times X_2)}) \to T(X_1 \times X_2),\qquad
{\rm E}_i \colon\ U(\Delta_{X_i)}) \to TX_i.
\end{gather*}

Let $x,y \in X_1 \times X_2$ which is sufficiently close each other.
Then we can use the $\nabla^z$-parallel transport with respect to the $\nabla^z$ geodesic
to define\index[syindex]{Pal@$({\rm Pal}_z)_x^y$}
$
({\rm Pal}_z)_x^y \colon T_x(X_1\times X_2) \to T_y(X_1\times X_2)$.
We denote by $\bigl({\rm Pal}^J_z\bigr)_x^y$ its complex linear part.

(1) implies that it splits into direct product if $z \in \mathfrak O^+$.
(2) implies that if $x,y \in L_{12}$ and~${z \in \partial \Sigma^{\rm d}}$, then
\begin{equation*}%
({\rm Pal}_z)_x^y(T_xL_{12}) \subset T_yL_{12}.
\end{equation*}
\end{defn}
\begin{rem}
Note that there may not exist a connection satisfying both of
Defini\-tion~\ref{defn1254} (1),~(2). This is the reason
why we use $z$ dependent family of connections.
\end{rem}
\begin{defn}
We define an operator
\begin{gather}
D^{\rho}_{u^{\rm d,\rho}} \overline\partial
:=
\bigl(D^{\rho}_{u_1^{\rm d,\rho}} \overline\partial,D^{\rho}_{u_2^{\rm d,\rho}} \overline\partial\bigr)
\colon\ W^2_{m+1,\delta}\bigl(\bigl(\Sigma^{\rm d};\partial \Sigma^{\rm d}\bigr),
\bigl(u^{\rm d,\rho}\bigr)^*T(X_1\times X_2);\bigl(u^{\rm d,\rho}\bigr)^*T(L_{12})\bigr)\nonumber \\
\hphantom{D^{\rho}_{u^{\rm d,\rho}} \overline\partial :=
\bigl(D^{\rho}_{u_1^{\rm d,\rho}} \overline\partial,D^{\rho}_{u_2^{\rm d,\rho}} \overline\partial\bigr)\colon }{} \
\to
L^2_{m,\delta}\bigl(\Sigma^{\rm d},
\bigl(u^{\rm d,\rho}\bigr)^*T(X_1 \oplus X_2) \otimes \Lambda^{0,1}_{\rho}\bigr) \label{eq128787}
\end{gather}
as follows. Let \smash{$(s,v) \in W^2_{m+1,\delta}\bigl(\bigl(\Sigma^{\rm d};\partial \Sigma^{\rm d}\bigr),
\bigl(u^{\rm d,\rho}\bigr)^*T(X_1\times X_2);\bigl(u^{\rm d,\rho}\bigr)^*T(L_{12})\bigr)$}.

Let $z \in \mathfrak O^+$. We put $s = (s_1,s_2)$, where $s_i$ is a section of
$\bigl(u_i^{\rm d,\rho}\bigr)^*TX_i$.
Then we define
\begin{gather}
\bigl(D^{\rho}_{u_i^{\rm d,\rho}} \overline\partial\bigr)(s,v)
:= \frac{\rm d}{{\rm d}t}\mathcal P^{-1} \bigl(\overline\partial_{j_{\mathfrak z_i}} {\rm Exp}_i\bigl(u_i^{\rm d,\rho}, ts_i\bigr)\bigr)\Big\vert_{t=0}
\label{eq1288}
\end{gather}
in a neighborhood of $z$.
Here \smash{${\rm Exp}_i\bigl(u_i^{\rm d,\rho},ts_i\bigr)$} is a map
\smash{$z \mapsto {\rm Exp}_i\bigl(u_i^{\rm d,\rho}(z),ts_i(z)\bigr)$}.
Then \linebreak $\overline\partial_{j_{\mathfrak z_i}} \smash{\bigl({\rm Exp}_i\bigl(u_i^{\rm d,\rho},ts_i\bigr)\bigr)}$ at $z$ is an
element of \smash{$T_{y(t)}X_i \otimes \Lambda_x^{0,1}\bigl(\Sigma^{\rm d}_i,j_{\mathfrak z_i}\bigr)$},
where
\[
y(t) = {\rm Exp}_i\bigl(u_i^{\rm d,\rho}(z),ts_i(z)\bigr).
\]
$\mathcal P$ is induced by \smash{$({\rm Pal}_z)_x^{y(t)}$} where $x = u_i^{\rm d,\rho}(z)$.
(We remark that ${\rm Pal}_z = {\rm Pal}_{\mathfrak o}$ in our case.)

Let $z \notin \mathfrak O^+$. Then the complex structure $j_{\mathfrak z_i}$ is the standard complex structure
$j$ in a~neighborhood of $z$. We define
\begin{equation}\label{eq1289}
\bigl(D^{\rho}_{u^{\rm d,\rho}} \overline\partial\bigr)(s,v)
:=\frac{\rm d}{{\rm d}t}\mathcal P_z^{-1} \bigl(\overline\partial {\rm Exp}_z\bigl(u^{\rm d,\rho}, ts\bigr)\bigr)\Big\vert_{t=0}
\end{equation}
in a neighborhood of $z$.
The notation in \eqref{eq1289} is similar to \eqref{eq1288}.
We however remark that in \eqref{eq1289} we work on the product space $X_1 \times X_2$
and use $z$ parametrized family of connections to define the exponential map and
the parallel transport.

By Definition~\ref{defn1254}\,(1), it is easy to see that \eqref{eq1289} coincides with \eqref{eq1288}
on the overlapped part and define \eqref{eq128787}.

The definition of the linearization map
\begin{equation}\label{eq1290}
D^{\rho}_{u^{\rm s,\rho}_i} \overline\partial
\colon\ W^2_{m+1,\delta}\bigl(\Sigma^{\rm s}_i,
\bigl(u_i^{\rm s,\rho}\bigr)^*TX_i\bigr) \to
L^2_{m,\delta}\bigl(\Sigma^{\rm s}_i,
\bigl(u^{\rm s,\rho}_i\bigr)^*TX_i \otimes \Lambda^{0,1}\bigr)
\end{equation}
is similar to and easier than \eqref{eq128787}.
\end{defn}
\begin{defn}\label{defn1257}
We denote by $W\bigl(m;u^{\rm d,\rho},u^{\rm s,\rho}_1,u^{\rm s,\rho}_2\bigr)$
the subspace of direct sum
\[
W^2_{m+1,\delta}\bigl(\bigl(\Sigma^{\rm d};\partial \Sigma^{\rm d}\bigr),
\bigl(u^{\rm d,\rho}\bigr)^*T(X_1\times X_2);\bigl(u^{\rm d,\rho}\bigr)^*T(L_{12})\bigr) \oplus
\bigoplus_{i=1,2}W^2_{m+1,\delta}\bigl(\Sigma^{\rm s}_i,
\bigl(u_i^{\rm s,\rho}\bigr)^*TX_i\bigr)
\]
consisting $((s,v),(s_1,v_1),(s_2,v_2))$ with
$
v = (v_1,v_2)$.
We consider the direct sum
\begin{gather}
L^2_{m,\delta}\bigl(\Sigma^{\rm d},
\bigl(u^{\rm d,\rho}\bigr)^*T(X_1 \oplus X_2) \otimes \Lambda^{0,1}_{\rho}\bigr)
\oplus
\bigoplus_{i=1,2}
L^2_{m,\delta}\bigl(\Sigma^{\rm s}_i,
\bigl(u^{\rm s,\rho}_i\bigr)^*TX_i \otimes \Lambda^{0,1}\bigr).\label{eq1292}
\end{gather}
We define
\begin{equation}\label{form1293}
D^{\rho}_{u^{\rm d,\rho},u^{\rm s,\rho}_1,u^{\rm s,\rho}_2} \overline\partial \colon\
W\bigl(m;u^{\rm d,\rho},u^{\rm s,\rho}_1,u^{\rm s,\rho}_2\bigr)
\to \eqref{eq1292}
\end{equation}
as the restriction of the direct sum of \eqref{eq128787} and \eqref{eq1290}.

\end{defn}
\begin{lem}\label{lema1258}
The map
\smash{$D^{\rho}_{u^{\rm d,\rho},u^{\rm s,\rho}_1,u^{\rm s,\rho}_2} \overline\partial$}
in \eqref{form1293} is surjective if $\delta$ is sufficiently small.

\end{lem}
\begin{proof}
Since we assumed Conditions \ref{defn12332} and \ref{conds1234} with trivial
obstruction bundle (that is, Assumption \ref{ass1247})
this is a consequence of the standard exponential decay estimate and regularity of
linear operators.
\end{proof}

\begin{defn}
We denote by $\mathfrak H\bigl(m;u^{\rm d,\rho},u^{\rm s,\rho}_1,u^{\rm s,\rho}_2\bigr)$
the $L^2$ orthonormal complement of the kernel of
\smash{$D^{\rho}_{u^{\rm d,\rho},u^{\rm s,\rho}_1,u^{\rm s,\rho}_2} \overline\partial$}
in $W\bigl(m;u^{\rm d,\rho},u^{\rm s,\rho}_1,u^{\rm s,\rho}_2\bigr)$.
\end{defn}
We next introduce bump functions we use in our gluing analysis.
(This part is similar to \cite[Section 3.1]{foooanalysis}.)
\begin{notation}
Hereafter, we use $[a,b]_{\tau'}$, $[a,b]_{\tau''_1}$, $[a,b]_{\tau''_2}$
for the interval $[a,b]$ to specify the coordinates $\tau'$, $\tau''_1$ or $\tau''_2$
we use.
\end{notation}
\begin{defn}
We define $\mathcal A^i_{T_i}$, $\mathcal X^i_{T_i}$, $\mathcal B^i_{T_i}$ for $i=1,2$ as follows:
\begin{gather*}
\mathcal A^i_{T_i} = [4T_i-1,4T_i+1]_{\tau'} \times S^1
= [-6T_i-1,-6T_i+1]_{\tau''_i} \times S^1, \\
\mathcal X^i_{T_i} = [5T_i-1,5T_i+1]_{\tau'} \times S^1
= [-5T_i-1,-5T_i+1]_{\tau''_i} \times S^1, \\
\mathcal B^i_{T_i} = [6T_i-1,6T_i+1]_{\tau'} \times S^1
= [-4T_i-1,-4T_i+1]_{\tau''_i} \times S^1.
\end{gather*}
They may be regarded as subsets of $\Sigma^{\rm d}_i$ or of $\Sigma^{\rm s}_i$ or of
$\Sigma_i(\bf r)$. (Here ${\bf r} = (\mathfrak r_1,\mathfrak r_2)$ and $T_i$ is related to $\mathfrak r_i$ by
\eqref{form1270}.)

\begin{figure}[ht]
\centering
\includegraphics[scale=0.34]{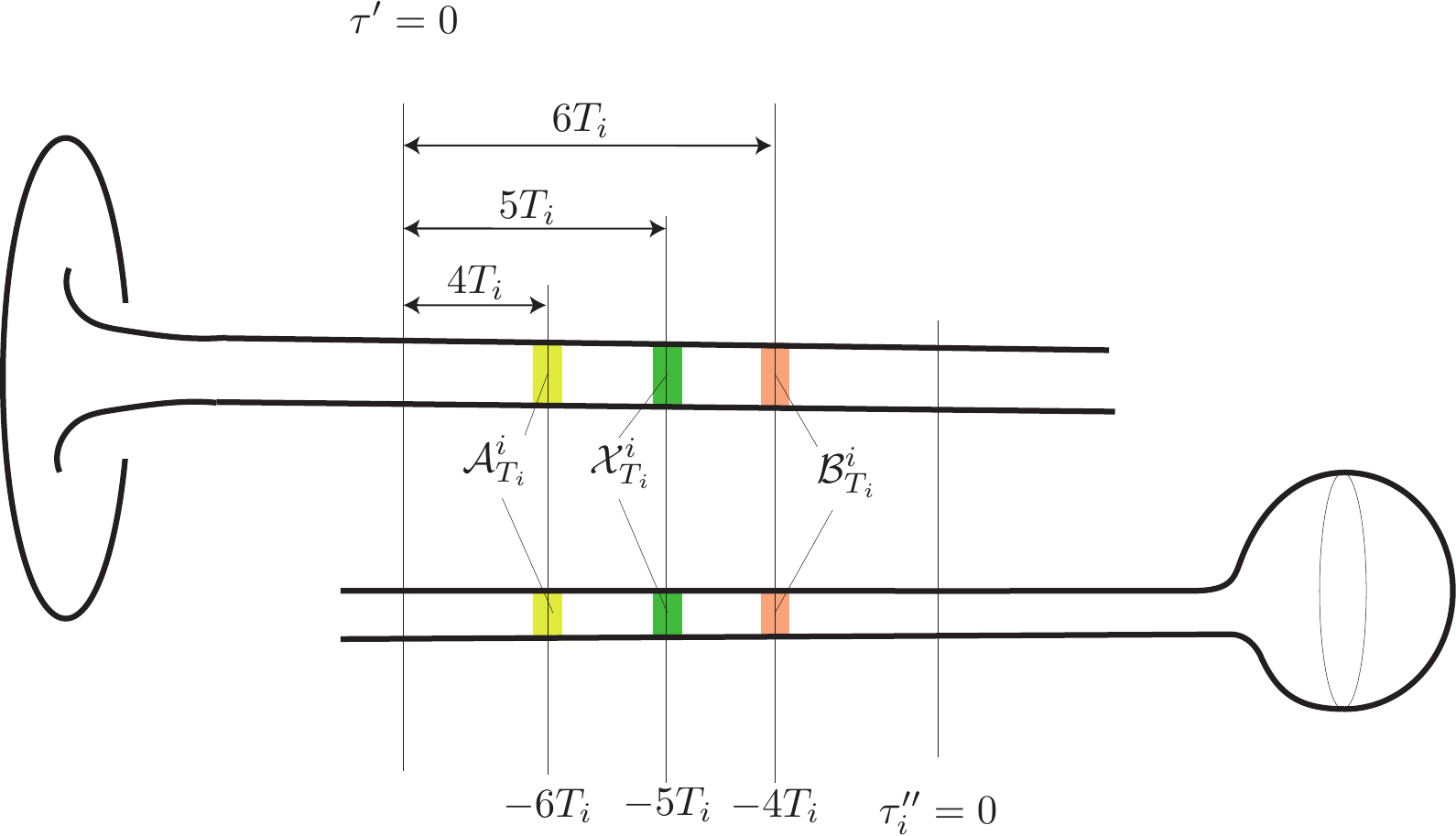}
\caption{$\mathcal A^i_{T_i}$, $\mathcal X^i_{T_i}$, $\mathcal B^i_{T_i}$.}
\label{Figure12-174}
\end{figure}

Let $\chi \colon \R \to [0,1]$ be a nondecreasing smooth function such that
\[
\chi(\tau)
=
\begin{cases}
0 &\text{if $\tau \le 1$}, \\
1 &\text{if $\tau \ge 1$}.
\end{cases}
\]
We use it to define functions on $[0,10T_i]_{\tau'} \times S^1 \cong
[-10T_i,0]_{\tau''_i} \times S^1$ as follows:
\begin{gather*}
\chi^{\rightarrow}_{\mathcal A^i_{T_i}}(\tau',t') := \chi(\tau'-4T_i),
\qquad \chi^{\rightarrow}_{\mathcal A^i_{T_i}}(\tau''_i,t'') := \chi(\tau''_i+6T_i),\\
\chi^{\rightarrow}_{\mathcal X^i_{T_i}}(\tau',t') := \chi(\tau'-5T_i),
\qquad \chi^{\rightarrow}_{\mathcal X^i_{T_i}}(\tau''_i,t'') := \chi(\tau''_i+5T_i),
\\
\chi^{\rightarrow}_{\mathcal B^i_{T_i}}(\tau',t') := \chi(\tau'-6T_i),
\qquad \chi^{\rightarrow}_{\mathcal B^i_{T_i}}(\tau''_i,t'') := \chi(\tau''_i+4T_i)
\end{gather*}
and
\smash{$
\chi^{\leftarrow}_{\mathcal A^i_{T_i}} := 1 - \chi^{\rightarrow}_{\mathcal A^i_{T_i}}$},
\smash{$ \chi^{\leftarrow}_{\mathcal X^i_{T_i}} := 1 - \chi^{\rightarrow}_{\mathcal X^i_{T_i}}$},
\smash{$ \chi^{\leftarrow}_{\mathcal B^i_{T_i}} := 1 - \chi^{\rightarrow}_{\mathcal B^i_{T_i}}$}.
We can extend them outside of
\[
[0,10T_i]_{\tau'} \times S^1 \cong
[-10T_i,0]_{\tau''_i} \times S^1
\]
 as locally constant functions
and regard them as functions on $\Sigma^{\rm d}_i$ or on $\Sigma^{\rm s}_i$ or on
$\Sigma_i(\bf r)$.
See Figure~\ref{Figure12vump}.
\begin{figure}[ht]
\centering
\includegraphics[scale=0.5]{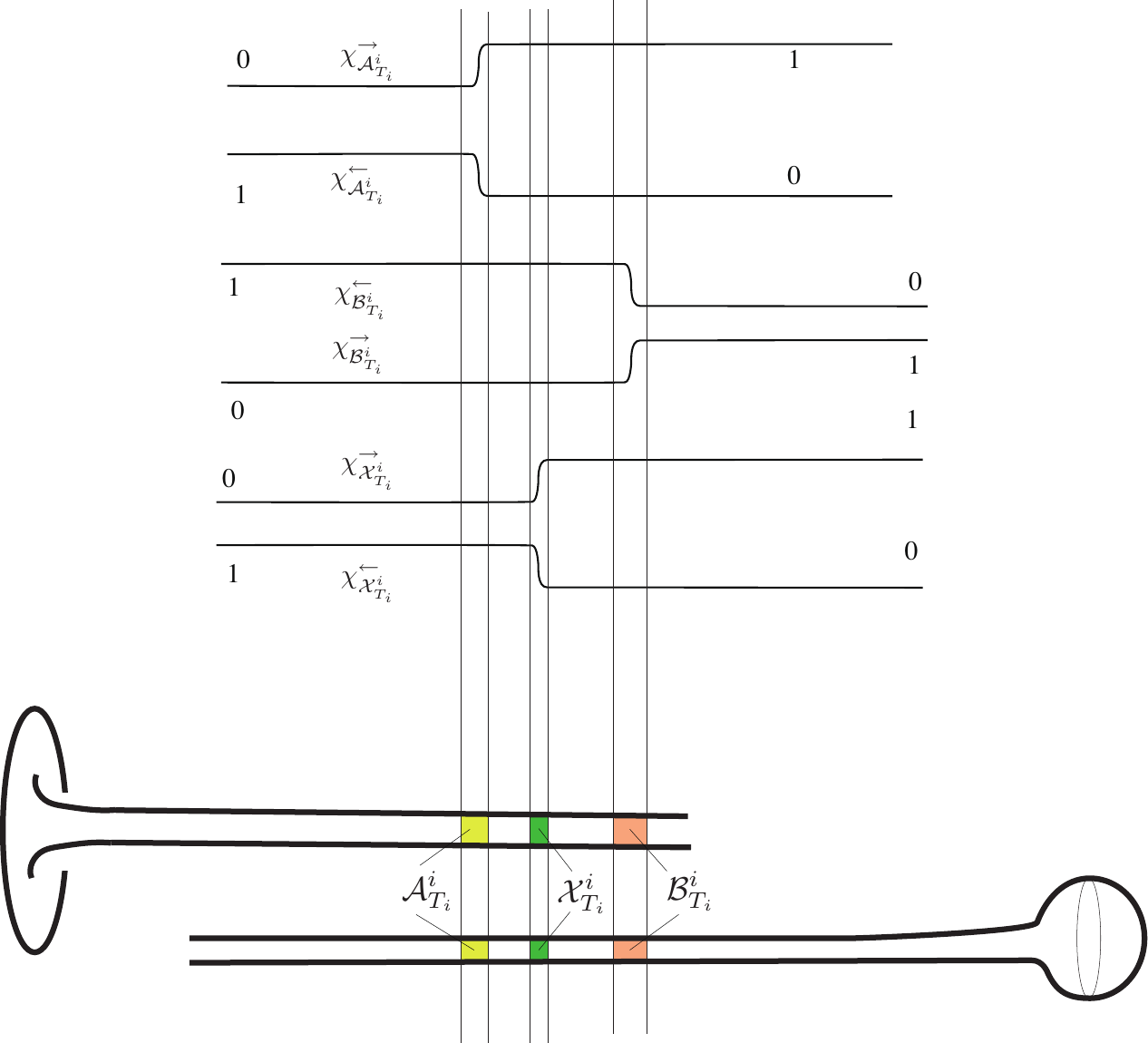}
\caption{Bump functions.}
\label{Figure12vump}
\end{figure}

\end{defn}
Now we are ready to start our inductive construction of gluing.
We discuss the case when the gluing parameters $\mathfrak r_1$ and $\mathfrak r_2$ are
both nonzero. (In the case when $\mathfrak r_1=\mathfrak r_2 = 0$, there is nothing
to do. The discussion when one of $\mathfrak r_1$, $\mathfrak r_2$ is zero
is similar and is omitted.)

 {\bf Pregluing.}
We put $u^{\rho}(\mathfrak o) = p^{\rho} = \bigl(p^{\rho}_1,p^{\rho}_2\bigr)$.
We recall ${\bf r} = (\mathfrak r_1,\mathfrak r_2)$.
A pair of complex numbers~${\bf r}$ corresponds to $T_1$, $\theta_1$, $T_2$, $\theta_2$ by
\eqref{form1256}.
We define
\smash{$
u_{{\bf r},(0),i}^{\rho} \colon \Sigma^{\rho}_i({\mathfrak r}_i) \to
X_i
$}
as follows.
We put $K^{\rm d}_i = \Sigma^{\rm d}_i \setminus \varphi^{\rm d}_{\mathfrak o_i}\bigl(D^2\bigr)$.
Then
$
\Sigma^{\rho}_i({\mathfrak r}_i)
= K^{\rm d}_i \cup K^{\rm s}_i \cup [0,10T_i]_{\tau'}$.

On $[0,10T_i]_{\tau'} \times S^1 \cong [-10T_i,0]_{\tau''_i} \times S^1$
we put
\[
u_{i,{\bf r},(0)}^{\rho}(\tau',t') =
 {\rm Exp}_i\bigl(p^\rho_i,
\chi_{\mathcal B^i_{T_i}}^{\leftarrow}(\tau',t') {\rm E}_i\bigl(p_i^{\rho},u_1^{\rho_1}(\tau',t')\bigr)
+ \chi_{\mathcal A^i_{T_i}}^{\rightarrow}(\tau''_i,t''_i) {\rm E}_i\bigl(p_i^\rho, u_2^{\rho_2}
(\tau''_i,t''_i)\bigr)\bigr).
\]
Here $(\tau''_i,t''_i)$ is related to $(\tau',t')$ by
\eqref{form12798}.
We also put
\smash{$u_{{\bf r},(0)}^{\rho} = (u_{1,{\bf r},(0)}^{\rho},u_{2,{\bf r},(0)}^{\rho})$}
with
$
\smash{u_{{\bf r},(0)}^{\rho} }:= u^{\rm d,\rho}
$
on $K^{\rm d}_i$ and
\smash{$
u_{i,{\bf r},(0)}^{\rho} := u^{\rm s,\rho}_i
$}
on $K^{\rm s}_i$.
$u_{{\bf r},(0)}^{\rho}$ is an approximate solution of
our pseudo-holomorphic curve equation.

 {\bf Step 0-(3+4) (Separating error terms into two parts).}
\begin{defn}
We define
\[
\hat u^{\rm d,\rho}_{{\bf r},(0)}
= \bigl(\hat u^{\rm d,\rho}_{1,{\bf r},(0)},\hat u^{\rm d,\rho}_{2,{\bf r},(0)}\bigr)
\colon\
\bigl(\Sigma^{\rm d},\partial \Sigma^{\rm d}\bigr)
\to (-X_1 \times X_2,L_{12})
\]
as follows:
\begin{gather*}
\hat u^{{\rm d},\rho}_{i,{\bf r},(0)}(z)
  :=
\begin{cases}
{\rm Exp}_i\bigl(p_i^\rho,\chi_{\mathcal B^i_{T_i}}^{\leftarrow}(\tau'-T_i,t')&\\
\quad {}\times
{\rm E}_i\bigl(p_i^\rho,u_{i,{\bf r},(0)}^{\rho}(\tau',t')\bigr)\bigr)&\text{if $z = (\tau',t') \in [0,10T_i]_{\tau'} \times S^1$}, \\
u_{i,{\bf r},(0)}^{\rho}(z)
&\text{if $z \in K^{\rm d}_i$}, \\
p_i^{\rho}
&\text{if $z \in [10T_i,\infty)_{\tau'}\times S^1$}.
\end{cases}
\end{gather*}
We also define
$
\hat u^{\rm s,\rho}_{i,{\bf r},(0)}
\colon
\Sigma_i^{\rm s}
\to X_i
$
as follows:
\begin{gather*}
\hat u^{{\rm s},\rho}_{i,{\bf r},(0)}(z)
:=
\begin{cases} {\rm Exp}_i\bigl(p_i^\rho,\chi_{\mathcal A^i_{T_i}}^{\rightarrow}(\tau''_i+T_i,t''_i)&\\
\quad{}\times
{\rm E}_i\bigl(p_i^\rho,u_{i,{\bf r},(0)}^{\rho}(\tau'_i,t'_i)\bigr)\bigr)&
\text{if $z = (\tau''_i,t''_i) \in [-10T_i,0]_{\tau''_i}$}\times S^1,\\
u_{i,{\bf r},(0)}^{\rho}(z)
&\text{if $z \in K^{\rm s}_i$}, \\
p_i^{\rho}
&\text{if $z \in (-\infty,-10T_i]_{\tau''_i}\times S^1$}.
\end{cases}
\end{gather*}

\end{defn}
\begin{defn}\label{defn523}
We put
\[
{\rm Err}^{\rho,\rm d}_{i,{\bf r},(1)}
= \chi_{\mathcal X^i_{T_i}}^{\leftarrow} \overline\partial_{j_{\mathfrak z_i}} u^{\rho} _{T,(0),i},\qquad
{\rm Err}^{\rho,\rm s}_{i,{\bf r},(1)}
= \chi_{\mathcal X^i_{T_i}}^{\rightarrow} \overline\partial_{j} u^{\rho} _{T,(0),i}.
\]
We regard them as elements of the weighted Sobolev spaces
\[
L^2_{m,\delta}\bigl(\Sigma^{\rm d}_i,
\bigl(\hat u^{\rm d,\rho}_{i,{\bf r},(0)}\bigr)^*TX_i \otimes \Lambda^{0,1}_{\rho}\bigr)\qquad
\text{and}\qquad L^2_{m,\delta}\bigl(\Sigma^{\rm s}_i,
\bigl(\hat u^{\rm s,\rho}_{i,{\bf r},(0)}\bigr)^*TX_i \otimes \Lambda^{0,1}\bigr)
\]
by extending them to be $0$ outside the support of~\smash{$\chi_{\mathcal X^i_{T_i}}^{\leftarrow}$}
and~\smash{$\chi_{\mathcal X^i_{T_i}}^{\rightarrow}$}, respectively.
Note that
\[
L^2_{m,\delta}\bigl(\Sigma^{\rm d}_i,
\bigl(\hat u^{\rm d,\rho}_{i,{\bf r},(0)}\bigr)^*TX_i \otimes \Lambda^{0,1}_{\rho}\bigr)
\qquad \text{and}\qquad
L^2_{m,\delta}\bigl(\Sigma^{\rm s}_i,
\bigl(\hat u^{\rm s,\rho}_{i,{\bf r},(0)}\bigr)^*TX_i \otimes \Lambda^{0,1}\bigr)
\]
are defined in the same way as Definition~\ref{defn1253}.
\end{defn}

 {\bf Step 1-1 (approximate solution for linearization).}
We define
\[W\bigl(m;\hat u^{\rm d,\rho}_{{\bf r},(0)},\hat u^{\rm s,\rho}_{1,{\bf r},(0)},
\hat u^{\rm s,\rho}_{2,{\bf r},(0)}\bigr)\]
in the same way as Definition~\ref{defn1257}.
Using \smash{$\bigl({\rm Pal}^J_z\bigr)^{\hat u^{\rm d,\rho}_{\bf r}(z)}_{u^{\rm d,\rho}(z)}$}
and \smash{$\bigl({\rm Pal}^J_{\mathfrak o}\bigr)^{\hat u^{\rm s,\rho}_{i,\bf r}(z)}_{u^{\rm s,\rho}_i(z)}$},
we obtain a linear map
\[
\Phi_{\rho,(0)}\colon\
W\bigl(m;u^{\rm d,\rho},u^{\rm s,\rho}_1,u^{\rm s,\rho}_2\bigr)
\to
W\bigl(m;\hat u^{\rm d,\rho}_{{\bf r},(0)},\hat u^{\rm s,\rho}_{1,{\bf r},(0)},
\hat u^{\rm s,\rho}_{2,{\bf r},(0)}\bigr).
\]
(See \cite[Definition 5.10 and Lemma 5.11]{foooanalysis}.)
We consider the direct sum
\begin{equation}\label{eq129222}
L^2_{m,\delta}\bigl(\Sigma^{\rm d},
\bigl(\hat u^{\rm d,\rho}_{{\bf r},(0)}\bigr)^*T(X_1 \oplus X_2) \otimes \Lambda^{0,1}_{\rho}\bigr)
\oplus
\bigoplus_{i=1,2}
L^2_{m,\delta}\bigl(\Sigma^{\rm s}_i,
\bigl(\hat u^{\rm s,\rho}_{i,\bf r}\bigr)^*TX_i \otimes \Lambda^{0,1}\bigr).
\end{equation}
We define
\[
D^{\rho}_{\hat u^{\rm d,\rho}_{{\bf r},(0)},\hat u^{\rm s,\rho}_{1,{\bf r},(0)},
\hat u^{\rm s,\rho}_{2,{\bf r},(0)}} \overline\partial \colon \
W\bigl(m;\hat u^{\rm d,\rho}_{{\bf r},(0)},\hat u^{\rm s,\rho}_{1,{\bf r},(0)},
\hat u^{\rm s,\rho}_{2,{\bf r},(0)}\bigr)
\to \eqref{eq129222}
\]
in the same way as \eqref{form1293}.
\begin{lem}
There exists a unique
element
\[
{\bf V}_{\rho,(1)} =
\bigl(\bigl(V^{\rm d}_{\rho,(1)},\Delta p_{\rho,(1)}\bigr),(V^{\rm s}_{1,\rho,(1)},\Delta p_{1,\rho,(1)}),(V^{\rm s}_{2,\rho,(1)},\Delta p_{2,\rho,(1)})\bigr)
\]
which is contained
in the image of the restriction of $\Phi_{\rho,(0)}$
to the $L^2$ orthogonal complement of the kernel
\[
\operatorname{Ker}D^{\rho}_{\hat u^{\rm d,\rho}_{{\bf r},(0)},\hat u^{\rm s,\rho}_{1,{\bf r},(0)},
\hat u^{\rm s,\rho}_{2,{\bf r},(0)}} \overline\partial \qquad \text{in}\quad
W\bigl(m;\hat u^{\rm d,\rho}_{{\bf r},(0)},\hat u^{\rm s,\rho}_{1,{\bf r},(0)},
\hat u^{\rm s,\rho}_{2,{\bf r},(0)}\bigr)
\]
such that
\[
\bigl(D^{\rho}_{\hat u^{\rm d,\rho}_{{\bf r},(0)},\hat u^{\rm s,\rho}_{1,{\bf r},(0)},
\hat u^{\rm s,\rho}_{2,{\bf r},(0)}} \overline\partial\bigr)({\bf V}_{\rho,(1)})
=
\bigl({\rm Err}^{\rho,\rm d}_{{\bf r},(1)},
{\rm Err}^{\rho,\rm s}_{1,{\bf r},(1)},{\rm Err}^{\rho,\rm s}_{2,{\bf r},(1)}\bigr).
\]
\end{lem}

This is a consequence of Lemma~\ref{lema1258} and can be proved in the same
way as
\cite[Lem\-ma~5.13]{foooanalysis}.\looseness=1

{\bf Step 1-2 (gluing solutions).}
\begin{defn}\label{defn518}
We define \smash{$u^{\rho}_{{\bf r},(1)} = (u^{\rho}_{1,{\bf r},(1)},
u^{\rho}_{2,{\bf r},(1)})$} as follows:
\begin{enumerate}\itemsep=0pt
\item[(1)] If $z \in K^{\rm d}$, we put
$
u^{\rho}_{{\bf r},(1)}(z)
=
{\rm Exp}^z\bigl(\hat u^{\rm d,\rho}_{{\bf r},(0)},V^{\rm d}_{\rho,(1)}(z)\bigr)$.
\item[(2)] If $z \in K^{\rm s}_i$, we put
$
u^{\rho}_{i,{\bf r},(1)}(z)
=
{\rm Exp}^z\bigl(\hat u^{\rm s,\rho}_{i,{\bf r},(0)},V^{\rm s}_{i,\rho,(1)}(z)\bigr)$.
\item[(3)]
If $z = (\tau',t') \in [0,10T_i]_{\tau'}\times S^1$,
we put
\begin{align*}
u^{\rho}_{i,{\bf r},(1)}(\tau',t')
= {}& {\rm Exp}^z\bigl(u_{{\bf r},(0)}^{\rho}(\tau',t'),
\chi_{\mathcal B^{i}_{T_i}}^{\leftarrow}(\tau,t) \bigl(V^{\rm d}_{i,\rho,(1)}(\tau',t') - (\Delta p_{i,{\rho},(1)})^{\rm Pal}\bigr)\\
& + \chi_{\mathcal A^{i}_{T_i}}^{\rightarrow}(\tau',t')\bigl(V^{\rm}_{i,\rho,(1)}(\tau''_i,t'')-
(\Delta p_{i,{\rho},(1)})^{\rm Pal}\bigr)
+ (\Delta p_{i,{\rho},(1)})^{\rm Pal}\bigr).
\end{align*}
\end{enumerate}
We also define
\smash{$
p_{(1)}^{\rho} = \bigl(p_{(1)}^{1,\rho},p_{(1)}^{2,\rho}\bigr)$}
by
\smash{$
p_{i,(1)}^{\rho} = {\rm Exp}_i(p_{i}^{\rho},\Delta p_{i,{\rho},(1)})$}.
\end{defn}

\smash{$u^{\rho}_{{\bf r},(1)}(z)$} is an improved approximate
solution.
Note that by Definition~\ref{defn1254}\,(1), \smash{$u^{\rho}_{{\bf r},(1)}(z)$}
satisfies the boundary condition at $\partial \Sigma(\mathfrak r)$.

 {\bf Step 1-4 (separating error terms into two parts).}
\begin{defn}
We define
\[
\hat u^{\rm d,\rho}_{{\bf r},(1)}
= \bigl(\hat u^{\rm d,\rho}_{1,{\bf r},(1)},\hat u^{\rm d,\rho}_{2,{\bf r},(1)}\bigr)
\colon\
\bigl(\Sigma^{\rm d},\partial \Sigma^{\rm d}\bigr)
\to (-X_1 \times X_2,L_{12})
\]
as follows:
\begin{gather*}
\hat u^{{\rm d},\rho}_{i,{\bf r},(1)}(z)
:=
\begin{cases} {\rm Exp}_i\bigl(p_{i,(1)}^\rho,\chi_{\mathcal B^i_{T_i}}^{\leftarrow}(\tau'-T_i,t')&\\
\quad{}\times
{\rm E}_i\bigl(p_i^\rho,u_{i,{\bf r},(0)}^{\rho}(\tau',t')\bigr)\bigr)
&\text{if $z = (\tau',t') \in [0,10T_i]_{\tau'}$}\times S^1, \\
u_{i,{\bf r},(1)}^{\rho}(z)
&\text{if $z \in K^{\rm d}_i$}, \\
p_{i,(1)}^{\rho}
&\text{if $z \in [10T_i,\infty)_{\tau'}\times S^1$}.
\end{cases}
\end{gather*}
We also define
\smash{$\hat u^{\rm s,\rho}_{i,{\bf r},(1)}
\colon
\Sigma_i^{\rm s}
\to X_i
$}
as follows:
\begin{gather*}
\hat u^{{\rm s},\rho}_{i,{\bf r},(1)}(z)
 :=
\begin{cases} {\rm Exp}_i\bigl(p_{i,(1)}^\rho,\chi_{\mathcal A^i_{T_i}}^{\rightarrow}(\tau''_i+T_i,t''_i)&\\
\quad{}\times
{\rm E}_i\bigl(p_i^\rho,u_{i,{\bf r},(0)}^{\rho}(\tau'_i,t'_i)\bigr)\bigr)
&\text{if $z = (\tau''_i,t''_i) \in [-10T_i,0]_{\tau''_i} \times S^1$},
\\
u_{i,{\bf r},(1)}^{\rho}(z)
&\text{if $z \in K^{\rm s}_i$}, \\
p_{i,(1)}^{\rho}
&\text{if $z \in (-\infty,-10T_i]_{\tau''_i}\times S^1$}.
\end{cases}
\end{gather*}
\end{defn}

\begin{defn}\label{defn52300}
We put
\[
{\rm Err}^{\rho,\rm d}_{i,{\bf r},(2)}
= \chi_{\mathcal X^i_{T_i}}^{\leftarrow} \overline\partial_{j_{\mathfrak z_i}} u^{\rho} _{T,(1),i}, \qquad
{\rm Err}^{\rho,\rm s}_{i,{\bf r},(2)}
= \chi_{\mathcal X^i_{T_i}}^{\rightarrow} \overline\partial_{j} u^{\rho} _{T,(1),i}.
\]
We regard them as elements of the weighted Sobolev spaces
\smash{$L^2_{m,\delta}\bigl(\Sigma^{\rm d}_i,
\bigl(\hat u^{\rm d,\rho}_{i,{\bf r},(1)}\bigr)^*TX_i \otimes \Lambda^{0,1}_{\rho}\bigr)$}
and
\smash{$L^2_{m,\delta}(\Sigma^{\rm s}_i,
(\hat u^{\rm s,\rho}_{i,{\bf r},(1)})^*TX_i \otimes \Lambda^{0,1})$},
respectively, by extending them to be $0$ outside the support of~\smash{$\chi_{\mathcal X^i_{T_i}}^{\leftarrow}$}
and \smash{$\chi_{\mathcal X^i_{T_i}}^{\rightarrow}$}, respectively.

\end{defn}
We now come back to Step 2-1 and continue.
We thus obtain a sequence of maps
$\smash{u^{\rho}_{{\bf r},(\kappa)} =}\allowbreak (u^{\rho}_{1,{\bf r},(\kappa)},
u^{\rho}_{2,{\bf r},(\kappa)})$
for $\kappa =0,1,2,\dots$ inductively on $\kappa$.

In the same way as \cite[Section 5]{foooanalysis},
we can show that it converges to
$u^{\rho}_{{\bf r}} = (u^{\rho}_{1,{\bf r}},
u^{\rho}_{2,{\bf r}})$ in $L^2_{m+1}$ norm as $\kappa$ goes to infinity.

Together with source object \eqref{form1264},
$u^{\rho}_{{\bf r}}$ defines an element of
$\mathcal M'_{1,2,2}(L_{12};({\rm diag});E)$.
This element is by definition
$\mathscr G\bigl(\bigl(\rho^{\rm d},\rho^{\rm s}_1,
\rho^{\rm s}_2,\mathfrak z_1,\mathfrak z_2\bigr),(\mathfrak r_1,\mathfrak r_2)\bigr)$.
The proof that it is injective and its image is an
open neighborhood of $\xi_0$ is the same as
\cite[Section 7]{foooanalysis}.

Let us elaborate on the latter proof now.
Below we discuss the case when the source object is unstable.
(So it is slightly more involved than the case of Proposition~\ref{prop1249}.)
We consider the situation at the beginning of Section~\ref{sec:glueglue}, depicted in Figure~\ref{Figure1241}.
We write the element of
$\mathcal M'_{1,2,2}(L_{12};{\rm (diag)};E)$ described there as
${\bf x}^+$.
Denote the four interior marked points of ${\bf x}^+$
by $w_{i,j}$, $i=1,2$, $j=1,2$, where $w_{i,j}$ is on $\Sigma_{i}^{\rm s}$.
Note that $0$ on the disk is an interior marked point of
first kind.
This element ${\bf x}^+$ comes with one boundary marked point $1$
(the symbol ${\rm (diag)}$ shows existence of one boundary marked point).
We forget the four interior marked points of second kind and one
interior marked point of first kind
to obtain~${{\bf x} \in \mathcal M'_{0,0,0}(L_{12};{\rm (diag)};E)}$.

We add 4 codimension two transversals $\mathcal N_{i,j} \subset X_i$
which intersect with the image of $u_i^{{\rm s},{\bf 0}}$ transversally
at $w_{i,j}$.
We also add 1 codimension two transversal $\mathcal N \subset X_1 \times X_2$ which intersect
 with \smash{$\bigl(u_1^{{\rm d},{\bf 0}},u_2^{{\rm d},{\bf 0}}\bigr)$}
transversally at $0$.
We will prove that the set of the image of $\mathscr G$ satisfying
the transversal constraint contains a
neighborhood of
${\bf x}$ in $\mathcal M'_{0,0,0}(L_{12};{\rm (diag)};E)$.
Suppose that ${\bf y}_n$ is a sequence of elements~${\mathcal M'_{0,0,0}(L_{12};{\rm (diag)};E)}$
converging to ${\bf x}$. By the definition of topology and
Lemma~\ref{lem1214140}, there exists ${\bf y}^+_n$
such that $i^*({\bf y}^+_n) = {\bf y}_n$
and ${\rm lims}_{n\to \infty} {\bf y}^+_n = {\bf x}^+$.
Here $i^*$ is the forgetful map of the marked points.
In particular, the source curves of
${\bf y}^+_n$ converge to the source curve of ${\bf x}^+$.
Let $w_{i,j}^n$ be the four interior marked points of second kind of
${\bf y}^+_n$ and $z_i^n$ be the interior marked points of first kind of ${\bf y}^+_n$.
Then, the source curve of~${\bf y}^+_n$
is a pair $((\Sigma^n_1,(1;z_1^n,w_{1,1}^n,;w_{1,2}^n)),
(\Sigma^n_2,(1;z_2^n,w_{2,1}^n,;\allowbreak w_{2,2}^n)))$
together with an isomorphism of the disk part of
$\Sigma^n_1$ to the disk part of $\Sigma^n_2$,
which sends~$1$ and $z_1^n$ to $1$ and $z_2^n$,
respectively.
(We denote the boundary marked point by $1$.)
Therefore, ${(\Sigma^n_i,(1;z_1^n,w_{i,1}^n,;w_{i,2}^n))}$
converges to $(\Sigma,(1;0,w_{i,1},w_{i,2}))$ in
the moduli space of marked stable bordered curves.
Here $\Sigma$ is a disk with one sphere bubble on it.
We change the representative and assume that~$z_1^n$ is~$0$.\looseness=-1

So we obtain a gluing parameter $\mathfrak r^n_1$ and the parameter of the position of the node
$\mathfrak z^n_1$ uniquely
such that $(\Sigma^n_1,(1;0;w_{1,1}^n,;w_{1,2}^n))$ is conformal
to $(\Sigma_1(\mathfrak z^n_1,\mathfrak r^n_1),(1,0,w_{1,1},w_{1,2}))$.
We obtain $\mathfrak r^n_2$ and~$\mathfrak z^n_2$ in a similar way.
We remark that $\lim_{n\to \infty} \mathfrak r^n_{i} = 0$.

We can identify $\Sigma^n_i \cong \Sigma_i(\mathfrak z^n_i,\mathfrak r^n_i)$
and $\Sigma_i(\mathfrak z^n_i,{\bf 0})$
outside the neck region\footnote{There exists a one neck for each $i \in \{1,2\}$.}
using the local trivialization of the universal family (outside the node).
Via this identification, the map $u^n_{i}$ which is a part of
${\bf y}^+_n$ converges to $u_{i}$ which is a part of
${\bf x}^+$ outside the neck region in the compact $C^{\infty}$
topology as $n$ goes to infinity.
On the neck region, we can use the exponential decay estimate
such as \cite[Proposition~7.1]{foooanalysis}.
Therefore, we can take $\rho^n = \bigl(\rho^{n,\rm d},\rho^{n,\rm s}_1,\rho^{n,\rm s}_2\bigr)$
such that
the difference of two elements
${\bf y}^+_n$ and $\mathscr G(\rho^n,\mathfrak z^n_1,\mathfrak z^n_2,\mathfrak r^n_1,\mathfrak r^n_2)$
goes to zero.

Hence we can interpolate $u^n_{i}$ which is a part of
${\bf y}^+_n$ and $u^{n,\prime}_{i}$
which is a map part of $\mathscr G(\rho^n,\mathfrak z^n_1,\mathfrak z^n_2,\allowbreak\mathfrak r^n_1,\mathfrak r^n_2)$ to obtain a
one parameter family of maps $u^{n,\mathfrak s}_{i}\colon \Sigma^n_i \to X_i
$ for $\mathfrak s \in [0,1]$
such that it becomes~$u^n_{i}$ and $u^{n,\prime}_{i}$ at $\mathfrak s = 0,1$.
(Note that the domain curves of them are isomorphic each other.)
We may also assume that the transversality constraints are satisfied.

For a sufficiently large $n$, this path $\mathfrak s \mapsto u^{n,\mathfrak s}_{i}$
can be arbitrary short. (The shortness is taken in the sense of
the weighted Sobolev norm we used in the gluing analysis.)

Now we run the Newton's iteration in the one parameter family
and modify $u^{n,\mathfrak s}_{i}$ so that it is pseudo-holomorphic.
Since $u^n_{i}$ and $u^{n,\prime}_{i}$ are both
pseudo-holomorphic, Newton's iteration does not change them.
Hence the path $\mathfrak s \mapsto u^{n,\mathfrak s}_{i}$ still joins them.
We may still assume that the transversality constraint are satisfied
by using implicit function theorem.

By index calculation, the image of the map $\mathscr G$ has the same
dimension as the moduli space for each fixed domain.
Therefore, we can lift our path to the domain of $\mathscr G$
for sufficiently large~$n$.
This is the proof of openness of the image.

We thus proved Proposition~\ref{prop1249}.

Then the proof of Proposition~\ref{propo1250}
is entirely the same as \cite[Section 6]{foooanalysis}.

The proof of Theorem~\ref{prop1417} is now complete.
\end{proof}

\begin{exm}
We consider the sequence
$\rho(n) = \bigl(\rho^{\rm d}_n,\rho^{\rm s}_n\bigr),(\mathfrak z_1(n),\mathfrak z_2(n))$ and a sequence
${\bf r}(n) = ((\mathfrak r_1(n),\mathfrak r_2(n)),(\mathfrak z_1(n),\mathfrak z_2(n)))$
which converges to $({\bf 0},\mathfrak z_0,\mathfrak z_0)$ and to $({\bf 0},{\bf 0})$ as
$n$ goes to infinity.
The limit of $\mathscr G(\rho(n),{\bf r}(n))$
in our compactification $\mathcal M'_{1,2,2}(L_{12};({\rm diag});E)$
is the object~${\xi_0 = \mathscr G(({\bf 0},({\bf 0},{\bf 0})),{\bf r}(n))}$ and
is independent of the choice of such sequences $\rho(n)$, ${\bf r}(n)$.

On the other hand, the limit of the sequence $\mathscr G(\rho(n),{\bf r}(n))$
in the stable map compactification~${\mathcal M_{1,2,2}(L_{12};({\rm diag});E)}$
depend on the choice of $\rho(n)$, ${\bf r}(n)$ as follows.

We put $d(n) = \vert \mathfrak z_1(n) - \mathfrak z_2(n)\vert$.

Case 1: If $d(n)/\vert \mathfrak r_1(n)\vert \to 0$, $d(n)/\vert \mathfrak r_2(n)\vert \to 0$.
Then the source curve of the limit $\mathscr G(\rho(n),{\bf r}(n))$
in the stable map compactification is as in Figure~\ref{Figure14-9}.

Case 2: $d(n)/\vert \mathfrak r_1(n)\vert > c > 0$, $\vert \mathfrak r_2(n)\vert/\vert \mathfrak r_1(n)\vert \to 0$.
Then the source curve of the limit $\mathscr G(\rho(n),{\bf r}(n))$
in the stable map compactification is as in Figure~\ref{Figure14-8}.

Case 3: $d(n)/\vert \mathfrak r_2(n)\vert > c > 0$, $\vert \mathfrak r_1(n)\vert/\vert \mathfrak r_2(n)\vert \to 0$.
Then the source curve of the limit $\mathscr G(\rho(n),{\bf r}(n))$
in the stable map compactification is as in Figure~\ref{Figure14-7}.

Case 4: $d(n)/\vert \mathfrak r_2(n)\vert > c > 0$, $c_2 > \vert \mathfrak r_1(n)\vert/\vert \mathfrak r_2(n)\vert > c_1$.
Then the source curve of the limit $\mathscr G(\rho(n),{\bf r}(n))$
in the stable map compactification is as in Figure~\ref{Figure14-92}.

We can prove these facts by looking the proof of Lemma~\ref{loem1229}.

Thus the stable map compactification $\mathcal M_{1,2,2}(L_{12};({\rm diag});E)$
is a kind of blow up of the space $\mathcal M'_{1,2,2}(L_{12};({\rm diag});E)$.
Note that the fact that the blow up of a variety $Z$ is smooth does not imply
the smoothness of~$Z$
in algebraic geometry. By the same reason, the fact that~$\mathcal M_{1,2,2}(L_{12};({\rm diag});E)$
has Kuranishi structure, which was proved in previous literatures, does not imply
$\mathcal M'_{1,2,2}(L_{12};({\rm diag});E)$ has Kuranishi structure.
This is the reason why we provide the detail of the proof of Theorem~\ref{prop1417}
in this subsection.
\end{exm}

\begin{rem}\label{rem1274}
We discuss the example in Remark~\ref{rem1254} and how the
gluing analysis works in that case.
Moreover, we will compare it to the gluing analysis in the case of
stable map compactification.

In the situation of Remark~\ref{rem1254}, we consider
the case when the configuration $\xi_0'$ which is defined by
$\bigl(u^d_1,u^d_2,u^s_1,u^s_2\bigr)$
is isolated among the objects in this
combinatorial type,\footnote{Here `this combinatorial type'
contains the condition $\mathfrak z_1 = \mathfrak z_2,$ that is, the roots
$\mathfrak z_1,\mathfrak z_2 \in D^2$ of the two sphere
bubbles coincide.} up to an automorphism on the sphere bubbles
which preserves the point $D_i^2 \cap S^2_i$ that is $0_i \in S^2_i$.
Here
\[
u^d_1 \colon\ D^2 \to -X_1, \qquad u^d_2 \colon\ D^2 \to X_2, \qquad
u^s_1 \colon\ S^2 \to -X_1, \qquad u^s_2 \colon\ S^2 \to X_2.
\]
with the constraint
$
u^d_i(\mathfrak o) = u^s_i(0)$, $ i=1,2$
and
$
\bigl(u^d_1(z),u^d_2(z)\bigr) \in L_{12}
$
for $z \in \partial D^2$.
Since we assume this configuration is isolated,
$u^d_i$ is uniquely determined by these conditions.\footnote{This follows from the fact that the group of automorphisms of the source curve
acts as an identity map on the disk component. This is because
the disk component has one boundary marked point and
one interior node.}
The automorphisms on the sphere bubbles which
preserves $0$ consist a complex two-dimensional
group $G_i$. So $u^s_i$ has a freedom
$u^s_i \mapsto u^s_i \circ g_i$, $g_i \in G_i$.
The space $\mathcal V$ has complex dimension~6,
parametrized $g_1 $, $g_2$ and $\mathfrak z_1$, $\mathfrak z_2$.
Here $\mathfrak z_i \in D^2_i$ which parametrizes
the `root' of the sphere bubble.
Note that we assumed that $\xi_0'$ is isolated among the object in this
combinatorial type. Therefore, such element is unique if $\mathfrak z_1 = \mathfrak z_2$.
So there is one constraint and the dimension is~${\dim_{\C} G_1 + \dim_{\C} G_2 + 2 - 1 = 5}$.

Together with two gluing parameters $\rho_1$, $\rho_2$
the domain of the map $\mathscr G$
has 7 complex dimension.

The Kuranishi neighborhood of $\xi_0'$ in $\mathcal M'_{1,0,0}(L_{12};({\rm diag});E)$
is a submanifold of this real 14-dimensional
manifold cutting out by the constraint \eqref{transversalconst}.
The constraint is by codimension 2 submanifolds at 4 added marked points.
Therefore, it decreases dimension by 8.

Thus Kuranishi neighborhood of $\xi_0'$ in $\mathcal M'_{1,0,0}(L_{12};({\rm diag});E)$
is a real 6-dimensional manifold.
It can be depicted schematically in the Figure~\ref{FigSec1710}\,(a).
\begin{figure}[ht]
\centering
\includegraphics[scale=0.45]{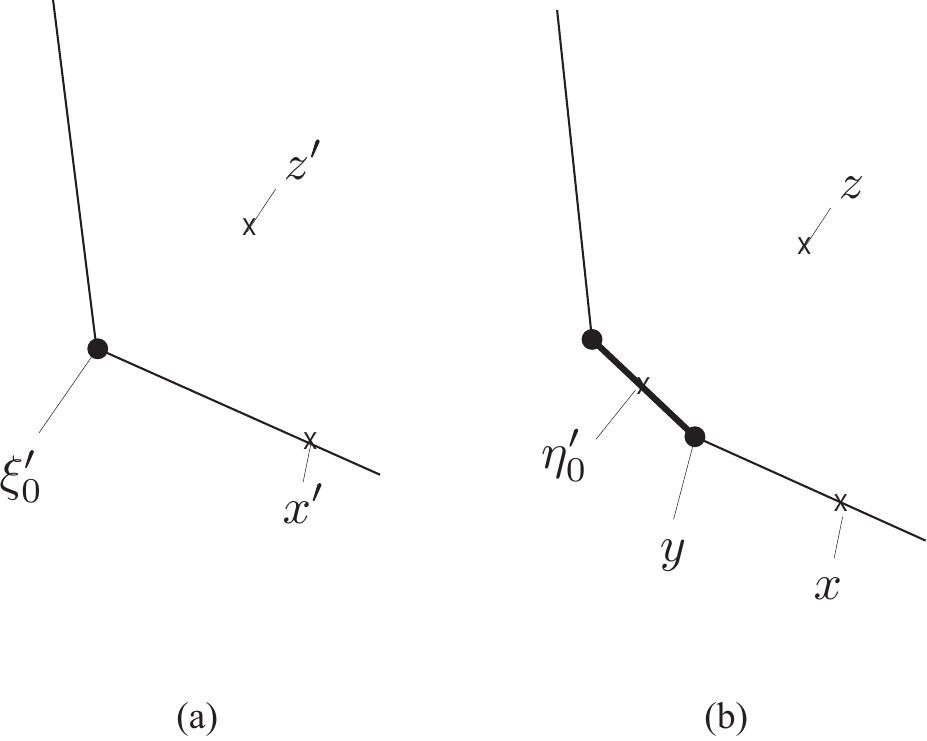}
\caption{Schematic pictures of Kuranishi neighborhoods.}
\label{FigSec1710}
\end{figure}
\begin{figure}[ht]
\centering
\includegraphics[scale=0.46]{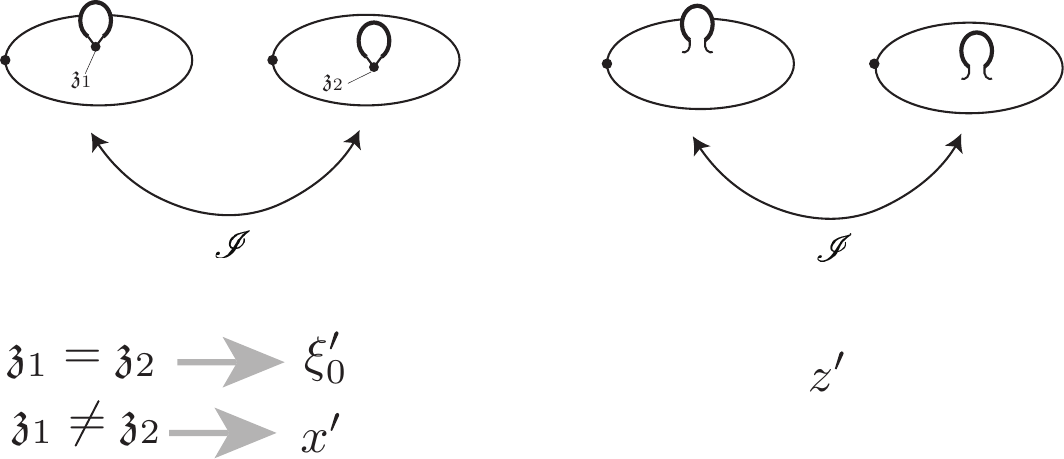}
\caption{Source curves of objects in Figure~\ref{FigSec1710}\,(a).}
\label{FigSec1711}
\end{figure}
\begin{figure}[ht]
\centering
\includegraphics[scale=0.5]{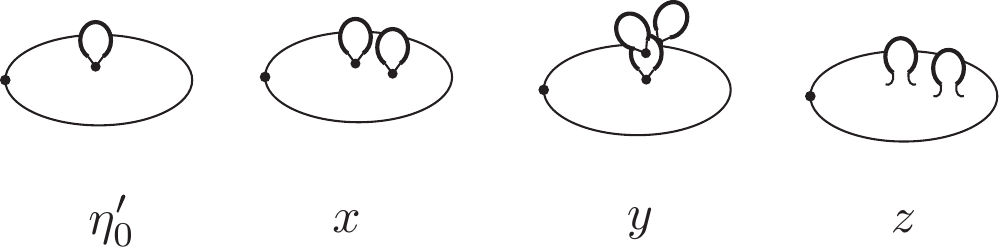}
\caption{Source curves of objects in Figure~\ref{FigSec1710}\,(b).}
\label{FigSec1712}
\end{figure}
\begin{figure}[ht]
\centering
\includegraphics[scale=0.5]{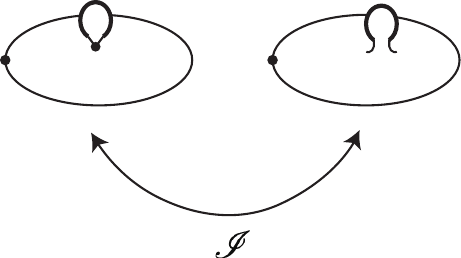}
\caption{Source curves of objects mentioned in the footnote.}
\label{FigSec1713}
\end{figure}

We next compare it with the Kuranishi neighborhood of the
corresponding object in the stable map compactification.
The element $\eta_0'
\in \mathcal M_{1,0}(L_{12};({\rm diag});E)$ corresponding to $\xi_0'$
in the stable map compactification,\footnote{Note that there is no prime in $\mathcal M_{1,0,0}(L_{12};({\rm diag});E)$.} has a source curve consisting of $D^2$ and one sphere bubble.
The map on the disk is $\bigl(\overline u^d_1, u^d_2\bigr)$ and
the map on the sphere is~$(\overline u^s_1, u^s_2)$.
Note that the element~$\eta_0'$ is {\it not} isolated
in this combinatorial type.
Namely, if we change $(\overline u^s_1, u^s_2)$ to
$(\overline u^s_1, u^s_2\circ g_2)$ with $g_2 \in G_2$
it represents a different element in
$\mathcal M_{1,0}(L_{12};({\rm diag});E)$.
Thus this stratum is a real 4-dimensional manifold.
(This family is depicted in Figure~\ref{FigSec1710}\,(b) by a thick line containing $\eta_0'$.)
Together with real 2-dimensional gluing parameter
it gives 6-dimensional family of objects.
The dimension coincides to the dimension
of the Kuranishi neighborhood of $\xi_0'$ in $\mathcal M'_{1,0,0}(L_{12};({\rm diag});E)$.
Those two moduli spaces are identical before compactification
and so the dimension must coincide.\footnote{When we work out the gluing process starting from $\eta_0'$,
we proceed as follows. We take two marked points, say $1$, $\infty$, on the sphere bubble
and take $\mathcal N_1$, $\mathcal N_{\infty}$ codimension 2 submanifold of $X_1 \times X_2$
which intersects with~${u^s: =(\overline u^s_1, u^s_2)}$ transversally at
$u^s(1)$ and $u^{s}(\infty)$. We thus obtain $\eta_0$. Note that
the neighborhood of $\eta_0$ has~4 extra (real) parameter corresponding to the group $G$
of automorphisms of $S^2$ preserving $0$.
After gluing we cutout by using the constraint defined by $\mathcal N_1, \mathcal N_{\infty}$.
Which decrease the dimension by 4.}
The Kuranishi neighborhood of $\eta_0'$ in $\mathcal M_{1,0,0}(L_{12};E)$ can be depicted schematically as in Figure~\ref{FigSec1710}\,(b).\footnote{This is an oversimplified picture.
In fact, there are other kinds of strata such as those depicted in Figure~\ref{FigSec1713}.}
The points $x$ and $z$ in Figure~\ref{FigSec1710}\,(b) are objects whose domain are depicted in Figure~\ref{FigSec1712}.

One can observe that Figure~\ref{FigSec1710}\,(b) is a kind of blow up of Figure~\ref{FigSec1710}\,(a) at the
stratum containing $\xi_0'$.
\end{rem}

\begin{rem}
Note that we did not care about the compatibility
of the Kuranishi structure with the forgetful map
of the boundary marked points in this section.
We actually use such a~compatibility to prove that the operations we obtain from our moduli spaces
is unital.
We can prove the compatibility with the forgetful map in the same way as,
for example,
\cite[Sections~3 and~5]{fooo091}.
Note that we only need to consider the forgetful map at the
diagonal component to study unitality.
Let $\vec a$ be as in Theorem~\ref{prop1417}.
We remove all the diagonal components from it except possibly $a_0$.
We denote it by $\vec a'$.
We use Theorem~\ref{prop1417} to obtain a Kuranishi structure on
$\mathcal M'(L_{12};\vec a';E)$.
Now for each element of $\mathcal M'(L_{12};\vec a;E)$
we use the obstruction space used in the construction of $\mathcal M'(L_{12};\vec a';E)$.
We then perform the gluing analysis in the same way to obtain a
Kuranishi structure on $\mathcal M'(L_{12};\vec a;E)$.
This Kuranishi structure is obviously compatible with the forgetful map.

We omit the detail since there is nothing new in this construction
compared to those which have already appeared in the literature.
\end{rem}

\section[Homotopy equivalence and homotopy between filtered $A_{\infty}$
functors]{Homotopy equivalence and homotopy between\\ filtered $\boldsymbol{A_{\infty}}$
functors}
\label{sec:homotopyfafunc}

In Section~\ref{subsec:Ainfcat} (see Definition~\ref{defn216}),
we defined the notion of two filtered $A_{\infty}$ functors
being homotopy equivalent and built homotopy theory
of filtered $A_{\infty}$ categories based on this notion.
This is the way taken in \cite{fu4}.
The way taken in \cite{fooobook} (in the case
of filtered $A_{\infty}$ algebras) is slightly
different. We describe the method of \cite{fooobook}
in the filtered $A_{\infty}$ category case
and discuss its relation to the method of Section~\ref{subsec:Ainfcat}.

There are certain issues to state it correctly
because the
category of categories is
rather a~2-category than a 1-category
and so claiming two morphisms of the category of categories
to be `the same' is a nontrivial issue.
A certain part of the discussion of this section is related to this point.

We say a filtered $A_{\infty}$ functor
$\mathscr F \colon \mathscr C_1 \to \mathscr C_2$ is
{\it linear} \index{linear} if
 $\mathscr F_k \colon B_k\mathscr C_1[1] \to \mathscr C_2[1]$
 is $0$ for $k\ne 1$.

\begin{defn}[{compare \cite[Definition 4.2.1]{fooobook}}]\label{defn131}
Let $\mathscr C$ be a non-unital curved filtered
$A_{\infty}$ category.
A {\it model of $\mathscr C \times [0,1]$}
\index{model of $\mathscr C \times [0,1]$}
consists of
$(\mathfrak C,{\rm Incl}, {\rm Eval}_0,{\rm Eval}_1)$
with the following properties:
\begin{enumerate}\itemsep=0pt
\item[(1)] $\mathfrak C$ is a curved non-unital filtered
$A_{\infty}$ category.\index[syindex]{incl@${\rm Incl}$}.\index[syindex]{eval@${\rm Eval}$}
\item[(2)]
${\rm Incl} \colon \mathscr C \to \mathfrak C$,
$ {\rm Eval}_0 \colon \mathfrak C \to \mathscr C$,
$ {\rm Eval}_1 \colon \mathfrak C \to \mathscr C$
are linear filtered $A_{\infty}$ functors,
such that
${\rm Incl}_{\rm ob} \colon \mathfrak{OB}(\mathscr C) \to \mathfrak{OB}(\mathfrak C)$,
$({\rm Eval}_0)_{\rm ob} \colon \mathfrak{OB}(\mathfrak C) \to \mathfrak{OB}(\mathscr C)$,
$({\rm Eval}_1)_{\rm ob} \colon \mathfrak{OB}(\mathfrak C) \to \mathfrak{OB}(\mathscr C)$
are bijections.
\item[(3)]
${\rm Eval}_0 \circ {\rm Incl} = {\rm Eval}_1 \circ {\rm Incl} =$ the identity functor:
$ \mathscr C \to \mathscr C$.
\item[(4)]
For $c,c' \in \mathscr C_{\rm ob}$, the map
$\overline{{\rm Incl}}_1(c,c') \colon \overline{\mathscr C}(c,c') \to \overline{\mathfrak C}({\rm Incl}_{\rm ob}(c),{\rm Incl}_{\rm ob}(c'))$
is a chain homotopy equivalence of the chain complexes, where $\overline{\mathfrak m}_1$ is the boundary
operators.
(Here~$\overline{{\rm Incl}}$ etc.\ denotes the $R$-reduction.)
For $c,c' \in \mathfrak C_{\rm ob}$ and $j=0,1$, the map
$\bigl(\overline{{\rm Eval}}_j\bigr)_1(c,c') \colon \overline{\mathfrak C}(c,c')
\to \overline{\mathscr C}(({\rm Eval}_j)_{\rm ob}(c),({\rm Eval}_j)_{\rm ob}(c'))$
is a chain homotopy equivalence of the chain complexes, where $\overline{\mathfrak m}_1$ is the boundary
operators.
\item[(5)]
For $c,c' \in \mathfrak C_{\rm ob}$,
the $\Lambda_0$ module homomorphism
\begin{gather*}
({{\rm Eval}}_0)_1(c,c') \oplus ({{\rm Eval}}_1)_1(c,c')\colon\\
\qquad {\mathfrak C}(c,c')
\to
{\mathscr C}(({\rm Eval}_1)_{\rm ob}(c),({\rm Eval}_1)_{\rm ob}(c'))
\oplus
{\mathscr C}(({\rm Eval}_2)_{\rm ob}(c),({\rm Eval}_2)_{\rm ob}(c'))
\end{gather*}
is split surjective.
\end{enumerate}
In the case when $\mathscr C$ is strict (resp.\ unital, $G$-gapped),
the model of $\mathscr C \times [0,1]$ is said strict (resp.\ unital, $G$-gapped) if
$\mathfrak C$, ${\rm Incl}$, ${\rm Eval}_0$, ${\rm Eval}_1$ are
all strict (resp.\ unital, $G$-gapped).

Sometimes, we say $\mathfrak C$ is a model of $\mathscr C \times [0,1]$
(and do not specify ${\rm Eval}_j$ and ${\rm Incl}$) by an
abuse of notation.
\end{defn}
By (2) and (3), we can identify $\mathfrak{OB}(\mathscr C)$ and $\mathfrak{OB}(\mathfrak C)$.
So we identify these two sets from now on.

\begin{prop}
For any curved non-unital filtered $A_{\infty}$ category $\mathscr C$, a model of $\mathscr C \times [0,1]$
exists.
If $\mathscr C$ is strict $($resp.\ unital, $G$-gapped$)$, then
we take the model so that it is strict $($resp.\ unital, $G$-gapped$)$.
\end{prop}

The proof is the same as the proof of
\cite[Lemma 4.2.13]{fooobook} (if $R$ contains $\Q$)
\cite[Lemma 4.2.25]{fooobook} (in general).
Those are the cases of a filtered $A_{\infty}$ algebra
but the proof of the category case is the same.

\begin{prop}\label{prop133}
Let $\mathscr C_j$ $(j=1,2)$ be non-unital curved filtered
$A_{\infty}$ categories and
$\mathscr F \colon \mathscr C_1 \to \mathscr C_2$
a filtered $A_{\infty}$ functor.
Let $\mathfrak C_j$ be a model of $\mathscr C_j\times [0,1]$
for $j=1,2$.
Then there exists a~filtered~$A_{\infty}$ functor
$\mathfrak F \colon \mathfrak C_1 \to \mathfrak C_2$
such that
$
{\rm Eval}_j \circ \mathfrak F = \mathscr F \circ {\rm Eval}_j
$
for~${j=0,1}$.
If $\mathscr C_j$ and $\mathscr F$ are strict $($resp.\ unital, $G$-gapped$)$,
we may choose $\mathfrak F$ to be strict $($resp.\ unital, $G$-gapped$)$.
\end{prop}

The proof is the same as the proof of \cite[Theorem~4.2.34]{fooobook}
and so is omitted.
Note that in Proposition~\ref{prop133} the case $\mathscr C_1 = \mathscr C_2 = \mathscr C$
and $\mathscr F$ is the identity functor is included.
In that case Proposition~\ref{prop133} implies the following.
\begin{cor}\label{cor134}
Let $\mathfrak C_j$ be a model of $\mathscr C \times [0,1]$ for $j=1,2$.
Then there exists a filtered $A_{\infty}$ functor
$\mathfrak F \colon \mathfrak C_1 \to \mathfrak C_2$
such that
${\rm Eval}_j \circ \mathfrak F = {\rm Eval}_j$
for $j=0,1$.
If $\mathscr C$ is strict $($resp.\ unital, $G$-gapped$)$,
we may choose $\mathfrak F$ to be strict $($resp.\ unital, $G$-gapped$)$.
\end{cor}

\begin{defn}\label{defn135}
Let $\mathscr C_j$ be a non-unital curved filtered
$A_{\infty}$ category for $j=1,2$ and
$\mathscr F, \mathscr G\colon \mathscr C_1 \allowbreak\to \mathscr C_2$
filtered $A_{\infty}$ functors.
Let $\mathfrak C_2$ be a model of $\mathscr C_2\times [0,1]$.

We say $\mathscr F$ is {\it homotopic}\index{homotopic} to $\mathscr G$
and write $\mathscr F \approx \mathscr G$
if there exists a filtered $A_{\infty}$ functor
$\mathscr H \colon \mathscr C_1 \to \mathfrak C_2$
such that
${\rm Eval}_0 \circ \mathscr H = \mathscr F$,
${\rm Eval}_1 \circ \mathscr H = \mathscr G$.
We call $\mathscr H$ the {\it homotopy functor}.
\index{homotopy functor}

We can define a strict (resp.\ unital, $G$-gapped) version in an obvious way.
\end{defn}

\begin{rem}\label{rem136}
If $\mathscr F$ is homotopic to $\mathscr G$, then
$\mathscr F_{\rm ob} = \mathscr G_{\rm ob}$.
This is a consequence of Definition~\ref{defn131}\,(2)(3).
\end{rem}
\begin{lem}\label{lem1347}
\quad
\begin{enumerate}\itemsep=0pt
\item[$(1)$] The notion `homotopic' is independent of the
choice of the model of $\mathscr C_2\times [0,1]$.
\item[$(2)$] `homotopic' is an equivalence relation.
\item[$(3)$]
If $\mathscr F \approx \mathscr F'$,
then
$\mathscr F \circ \mathscr G \approx \mathscr F' \circ \mathscr G$,
$\mathscr G \circ \mathscr F \approx \mathscr G \circ \mathscr F'$.
\end{enumerate}
The strict $($resp.\ unital, $G$-gapped$)$ version of these statements also hold.
\end{lem}
\begin{proof}
(1) follows from Corollary~\ref{cor134}
(see \cite[Lemma 4.2.36]{fooobook}).
(2) can be proved in the same way as \cite[Proposition~4.2.37]{fooobook}.
The proof of (3) is the same as \cite[Lemma~4.2.43]{fooobook}.
\end{proof}

\begin{defn}\label{def138}
Let $\mathscr F \colon \mathscr C_1 \to \mathscr C_2$
be a filtered $A_{\infty}$ functor between
non-unital curved filtered~$A_{\infty}$ categories.
We say that $\mathscr F$ is a {\it strong homotopy equivalence}\index{strong homotopy equivalence}
if there exists a
filtered~$A_{\infty}$ functor
$\mathscr G \colon \mathscr C_2 \to \mathscr C_1$
such that
$\mathscr F \circ \mathscr G \colon
\mathscr C_2 \to \mathscr C_2$
and
$\mathscr G \circ \mathscr F \colon
\mathscr C_1 \to \mathscr C_1$
are homotopic to the identity functor.

We call $\mathscr G$ the {\it strong homotopy inverse} to $\mathscr F$.\index{strong homotopy inverse}
We say two non-unital curved filtered
$A_{\infty}$ categories are {\it strongly homotopy equivalent}
\index{strongly homotopy equivalent} to each other
if there exists a strong homotopy equivalence between them.

We can define a strict (resp.\ unital, $G$-gapped) version in an obvious way.
\end{defn}

\begin{rem}
If $\mathscr F \colon \mathscr C_1 \to \mathscr C_2$ is a strong homotopy equivalence,
then it induces a bijection~${\mathfrak{OB}(\mathscr C_1) \to \mathfrak{OB}(\mathscr C_2)}$.
This is a consequence of Definition~\ref{rem136}.

This is a rather restrictive requirement. To define an appropriate notion
of equivalence between ($A_{\infty}$) categories, it is {\it not} a
correct idea to require that the set of objects are equal.
This point is related to the basic concept of category,
where an equality should be replaced by an equivalence.
This is a point where the notion of a homotopy equivalence
which we introduced in Definition~\ref{defn225}
is more natural from the point of view of category
theory than the notion of a strong homotopy equivalence
we defined above.

We will further discuss the relation between these two notions
later in this section.
\end{rem}

\begin{lem}\label{lem1310}
Let $\mathscr F \colon \mathscr C_1 \to \mathscr C_2$ be a
strong homotopy equivalence.
\begin{enumerate}\itemsep=0pt
\item[$(1)$]
Let $\mathscr G,\mathscr G'\colon \mathscr C \to \mathscr C_1$
be filtered $A_{\infty}$ functors. Then
$\mathscr G$ is homotopic to $\mathscr G'$
if and only if~$\mathscr F\circ\mathscr G$ is homotopic to $\mathscr F\circ\mathscr G'$.
\item[$(2)$]
Let $\mathscr G,\mathscr G' \colon \mathscr C_2 \to \mathscr C$
be filtered $A_{\infty}$ functors. Then
$\mathscr G$ is homotopic to $\mathscr G'$
if and only if~$\mathscr G\circ\mathscr F$ is homotopic to $\mathscr G'\circ\mathscr F$.
\item[$(3)$]
Composition of strong homotopy equivalences is a strong homotopy equivalence.
\end{enumerate}
The strict $($resp.\ unital, $G$-gapped$)$ version of these statements also hold.

\end{lem}
The proof is easy and is omitted.

Now a strong homotopy equivalence version of Theorem~\ref{white}
is the following.
We assume that the ground ring $R$ is a field.
\begin{thm}\label{white2}
Let $\mathscr C_1$, $\mathscr C_2$ be $G$-gapped filtered $A_{\infty}$ categories
and $\mathscr F \colon \mathscr C_1 \to \mathscr C_2$ a $G$-gapped
filtered $A_{\infty}$
functor such that
\begin{enumerate}\itemsep=0pt
\item[$(1)$]
For any $c,c' \in \mathfrak{OB}(\mathscr C_1)$,
the map
$\overline{\mathscr F}_1 \colon \overline{\mathscr C}_1(c_1,c'_1)
\to \overline{\mathscr C}_2(\mathscr F_{\rm ob}(c_1),
\mathscr F_{\rm ob}(c'_1))$ induces an isomorphism on
$\overline{\mathfrak m}_1$ homology.
\item[$(2)$] The map $\mathscr F_{\rm ob}
\colon \mathfrak{OB}(\mathscr C_1) \to \mathfrak{OB}(\mathscr C_2)$
is a bijection.
\end{enumerate}
Then $\mathscr F$ is a strong homotopy equivalence.
The strong homotopy inverse can be taken to be $G$-gapped.
If $\mathscr C_1$, $\mathscr C_2$, $\mathscr F$ are strict $($resp.\ unital$)$, then
we may take the strong homotopy inverse to be strict $($resp.\ unital$)$.
\end{thm}

The proof is the same as the proof of \cite[Theorem~4.2.45]{fooobook}.

We next discuss a relation between strong homotopy equivalence
and homotopy equivalence.
\begin{lem}\label{lem1312}
Suppose $\mathfrak C$ is a model of $\mathscr C \times [0,1]$
and assume that $\mathfrak C$ is $G$-gapped.
Then ${\rm Incl}$ is a strong homotopy inverse of
${\rm Eval}_0$. It is a strong homotopy inverse of
${\rm Eval}_1$ also.
\end{lem}

\begin{proof}
By Theorem~\ref{white2}, ${\rm Incl}$ is strong homotopy
equivalence.
The lemma then follows from ${\rm Eval}_j \circ {\rm Incl} =$
identity and Lemma~\ref{lem1310}.
\end{proof}

\begin{prop}\label{prop1313}
In the situation of Definition {\rm\ref{defn135}}
we assume that $\mathscr C_1$, $\mathscr C_2$,
$\mathscr F$ and $\mathscr G$ are strict and $G$-gapped.
We also assume that $\mathscr C_2$ is unital.
Then the following holds.
If $\mathscr F$ is homotopic to $\mathscr G$ in
the sense of Definition {\rm\ref{defn135}}, then
$\mathscr F$ is homotopy equivalent to $\mathscr G$
in the sense of Definition {\rm\ref{defn216}}.
\end{prop}
\begin{proof}
We first prove the following analogue of
Lemma~\ref{lem1312}.
\begin{lem}\label{lem1314}
Suppose $\mathfrak C$ is a model of $\mathscr C \times [0,1]$
and assume that $\mathfrak C$ is strict unital and $G$-gapped.
Then ${\rm Incl}$ is a homotopy inverse of
${\rm Eval}_0$. It is a homotopy inverse of
${\rm Eval}_1$ also.
\end{lem}
\begin{proof}
Using Theorem~\ref{white} in place of Theorem~\ref{white2},
the proof is the same as the proof of Lemma~\ref{lem1312}.
\end{proof}

We also remark that Lemma~\ref{lem1310}
still holds when we replace strong homotopy equivalence
by homotopy equivalence.

Now we prove Proposition~\ref{prop1313}.
We assume that $\mathscr F$ is homotopic to $\mathscr G$ in
the sense of Definition~\ref{defn135}
and let $\mathscr H \colon \mathscr C_1 \to \mathfrak C_2$
be the homotopy.
Since ${\rm Eval}_0\circ\mathscr H = \mathscr F$,
Lemma~\ref{lem1314} implies that
${\rm Incl} \circ \mathscr F$ is homotopy equivalent
to~$\mathscr H$.
In the same way, we can show that ${\rm Incl} \circ \mathscr G$
to~$\mathscr H$.
Therefore, ${\rm Incl} \circ \mathscr F$
is homotopy equivalent to ${\rm Incl} \circ \mathscr G$.
Since ${\rm Incl}$ is a homotopy equivalence,
the analogue of Lemma~\ref{lem1310} we mentioned above
implies $\mathscr F$ is homotopy equivalent to~$\mathscr G$.\looseness=-1
\end{proof}

We remark that the converse to Lemma~\ref{lem1314}
is false. Namely, there is
a pair of strict, unital and $G$-gapped filtered
$A_{\infty}$ functors $\mathscr F$,
$\mathscr G$ such that
they are homotopy equivalent,
$\mathscr F_{\rm ob} = \mathscr G_{\rm ob}$,
but $\mathscr F$ is not homotopic to
$\mathscr G$.
A counterexample is the following.
\begin{exm}\label{example1315}
Let $C$ be an associative ring with unit.
We regard it as a differential graded algebra
with trivial boundary operator and grading.
We then regard it as a (filtered) $A_{\infty}$ category
$\mathscr C$
(with trivial filtration) as in Definition~\ref{defnDGcate}, Remark~\ref{remDGcat}.
Let $f_1,f_2 \colon C\to C$ be ring homomorphisms.
We regard them as (filtered) $A_{\infty}$ functors $\mathscr C \to \mathscr C$.
We remark that
$f_1$ is homotopic to $f_2$ in the sense of
Definition~\ref{def138} if and only if
$f_1 = f_2$.
On the other hand, $f_1$ is homotopy equivalent to $f_2$
in the sense of Definition~\ref{defn216} if and only if
there exists an invertible element $g \in C$
\big(that is, an element such that there exists $g^{-1} \in C$
with $g \cdot g^{-1} = g^{-1} \cdot g = 1$\big)
such that $f_1(x) = g^{-1}f_2(x)g$.
Thus they are different notion in this case.
\end{exm}

\begin{cor}
Let $\mathscr C_i$ be a $G$-gapped filtered
$A_{\infty}$ category for $i=1,2$. We assume that they are
strict and $\mathscr C_2$ is unital.
Let $\mathscr F \colon \mathscr C_1 \to \mathscr C_2$ be a
filtered $A_{\infty}$ functor, which is strict and $G$-gapped.
Assume $\mathscr F_{\rm ob} \colon \mathfrak{OB}(\mathscr C_1) \to
\mathfrak{OB}(\mathscr C_2)$ is a bijection.
Then the next two conditions are equivalent:
\begin{enumerate}\itemsep=0pt
\item[$(1)$]
$\mathscr F$ is a strong homotopy equivalence in the sense of
Definition {\rm\ref{def138}}.
\item[$(2)$]
$\mathscr F$ is a homotopy equivalence in the sense of
Definition {\rm\ref{defn225}}.
\end{enumerate}
\end{cor}
\begin{proof}
This is immediate from Theorems~\ref{white} and \ref{white2}.
\end{proof}

\begin{rem}
Let $\mathscr F, \mathscr G \colon \mathscr C_1 \to \mathscr C_2$ be
two strict filtered $A_{\infty}$ functors between strict filtered~$A_{\infty}$
categories. We assume $\mathscr F$ is homotopy equivalent to
$\mathscr G$. It means that there exists
natural transformations $\mathcal T \colon \mathscr F \to \mathscr G$,
$\mathcal S \colon \mathscr G \to \mathscr F$ of degree $0$
and pre-natural transformations
$\mathcal U \colon \mathscr F \to \mathscr F$, $\mathcal V \colon \mathscr G \to \mathscr G$
such that
\begin{equation}\label{form131}
\mathfrak m_2(\mathcal S,\mathcal T) = \mathcal{ID} + \delta \mathcal U,
\qquad
\mathfrak m_2(\mathcal T,\mathcal S) = \mathcal{ID} + \delta \mathcal V.
\end{equation}
Let us elaborate on these equalities.
For $c_1,c_2 \in \mathfrak{OB}(\mathscr C_1)$, the functors $\mathscr F$ and $\mathscr G$
induce homomorphisms
\begin{gather}
(\mathscr F_1)_* \colon\ H(\mathscr C_1(c_1,c_2))
\to H(\mathscr C_2(\mathscr F_{\rm ob}(c_1),\mathscr F_{\rm ob}(c_2)),\nonumber \\
(\mathscr G_1)_* \colon\ H(\mathscr C_1(c_1,c_2))
\to H(\mathscr C_2(\mathscr G_{\rm ob}(c_1),\mathscr G_{\rm ob}(c_2)).\label{form1311}
\end{gather}
Here $H$ in the right and left-hand sides are $\mathfrak m_1$-homologies.
We show that the two maps in \eqref{form1311} coincide
as follows.
We observe that $\mathcal T$ and $\mathcal S$ induce
\[
\mathcal T(c_i) \in
H(\mathscr C_2(\mathscr F_{\rm ob}(c_i),\mathscr G_{\rm ob}(c_i))),
\qquad
\mathcal S(c_i) \in
H(\mathscr C_2(\mathscr G_{\rm ob}(c_i),\mathscr T_{\rm ob}(c_i))).
\]
We define
\begin{gather*}
\varphi \colon\ H(\mathscr C_2(\mathscr F_{\rm ob}(c_1),\mathscr F_{\rm ob}(c_2))
\to H(\mathscr C_2(\mathscr G_{\rm ob}(c_1),\mathscr G_{\rm ob}(c_2)), \\
\psi \colon\ H(\mathscr C_2(\mathscr G_{\rm ob}(c_1),\mathscr G_{\rm ob}(c_2))
\to H(\mathscr C_2(\mathscr F_{\rm ob}(c_1),\mathscr F_{\rm ob}(c_2))
\end{gather*}
by
$\varphi([x]) = [\mathfrak m_2(\mathfrak m_2(\mathcal S(c_1),x),\mathcal T(c_2))]$, $ \psi([y]) = [\mathfrak m_2(\mathfrak m_2(\mathcal T(c_1),y),\mathcal S(c_2))]$.
Using \eqref{form131} and definitions, we can show that
$
\varphi \circ \psi = {\rm id}$, $
\psi \circ \phi = {\rm id}$, $
\varphi \circ (\mathscr F_1)_* = (\mathscr G_1)_*$.
In other words, two maps~$(\mathscr F_1)_*$ and $(\mathscr G_1)_*$ are identified
by the isomorphism $\varphi$, $\psi$.

\end{rem}
We also can show the proposition on
a relation between associated strict functors
and homotopies.
\begin{prop}\label{prop1318}
Let $\mathscr F, \mathscr G \colon \mathscr C_1 \to \mathscr C_2$ be
two filtered $G$-gapped $A_{\infty}$ functors between $G$-gapped non-unital curved filtered $A_{\infty}$
categories. We assume $\mathscr C_2$ is unital.
Let $\mathscr F^s, \mathscr G^s \colon \mathscr C^s_1 \to \mathscr C^s_2$
be associated strict functors between associated strict categories.
If $\mathscr F$ is homotopic to $\mathscr G$, then~$\mathscr F^s$ is homotopy equivalent to $\mathscr G^s$.
\end{prop}

\begin{proof}
Let $\mathfrak C_2$ be a model of $\mathscr C_2 \times [0,1]$ and
$\mathcal H \colon \mathscr C_1 \to \mathfrak C_2$
a homotopy between $\mathcal F$ and $\mathcal G$.
It induces a strict filtered $A_{\infty}$ functor
$\mathcal H^s \colon \mathscr C^s_1 \to \mathfrak C^s_2$.
The linear filtered $A_{\infty}$ functors
${\rm Incl}$, ${\rm Eval}_0$, ${\rm Eval}_1$
induce
${\rm Incl}^s \colon \mathscr C^s_1 \to \mathfrak C^s_2$,
${\rm Eval}^s_0, {\rm Eval}^s_1 \colon \mathfrak C^s_2 \to \mathscr C^s_1$, respectively.
We obtain equalities
\begin{gather}
{\rm Eval}^s_0\circ {\rm Incl}^s = {\rm Eval}^s_1\circ {\rm Incl}^s
= \mathscr{ID},\qquad
 {\rm Eval}^s_0\circ \mathcal H^s = \mathcal F^s,
{\rm Eval}^s_1\circ \mathcal H^s = \mathcal G^s\label{form133}
\end{gather}
from the corresponding equalities between $\mathcal F$,
$\mathcal G$ and etc.

Moreover, by Theorem~\ref{white}, the first line of \eqref{form133} and
Definition~\ref{defn131}\,(4) imply that ${\rm Incl}^s$, ${\rm Eval}^s_0$
and ${\rm Eval}^s_1$ are homotopy equivalences
and ${\rm Incl}^s$ is a homotopy inverse to ${\rm Eval}^s_i$, $i=0,1$.
The second line of \eqref{form133} then implies
\[
{\rm Incl}^s\circ \mathcal F^s \approx {\rm Incl}^s\circ {\rm Eval}^s_0\circ \mathcal H^s
\approx \mathcal H^s
\approx {\rm Incl}^s\circ {\rm Eval}^s_1\circ \mathcal H^s
\approx {\rm Incl}^s\circ \mathcal G^s.
\]
Then using Proposition~\ref{defn2929}, we conclude
$\mathcal F^s \approx \mathcal G^s$.
\end{proof}

\begin{rem}
In the situation of Proposition~\ref{prop1318}, we can not expect $\mathscr F^s$ is homotopic to~$\mathscr G^s$.
In fact, the object $\mathscr F^s_{\rm ob}(c,b)$ is
$(\mathcal F_{\rm ob}(c),\mathcal F_*(b))$
and the object $\mathscr G^s_{\rm ob}(c,b)$ is
$(\mathcal G_{\rm ob}(c),\mathcal G_*(b))$.
They are in general different objects.
Note that $\mathcal F_{\rm ob}(c) = \mathcal G_{\rm ob}(c)$
but $\mathcal F_*(b) \ne \mathcal G_*(b)$
in general.
We can show that
$\mathcal F_*(b)$ is gauge equivalent to
$\mathcal G_*(b)$, in the sense of \cite[Definition 4.3.1]{fooobook}.
So they are not so far away from being `equal'.
However, because of well-known problem to
distinguish saying equal and equivalent
this small difference should be taken seriously.
\end{rem}
We remark that a filtered $A_{\infty}$
bi-functor
$
\mathscr F \colon \mathscr C_1 \times \mathscr C_2 \to \mathscr C_3
$
is identified with a filtered $A_{\infty}$ functor
$
\mathscr C_1 \to \mathcal{FUNC}(\mathscr C_2,\mathscr C_3)
$
by Definition~\ref{lem56}.
We can use this fact to define the notion
that two filtered $A_{\infty}$
bi-functors to be homotopic each other,
in an obvious way.
The case of a~tri-functor etc.\ is similar.

\section[Independence of the filtered $A_\infty$ functors of the
choices]{Independence of the filtered $\boldsymbol{A_{\infty}}$ functors of the
choices}
\label{sec:independence2}

\subsection{Statement}
\label{sec:independence1state}

In this section, we prove that the
correspondence functor and
correspondence bi-functor are
independent of the choices involved in the
construction.
In this subsection, we state the main result
of this section.

\begin{choice}\label{chicesec151}
Suppose we are in Situation \ref{situ61}.
We choose a compatible almost complex structure
$J_{X_i}$ on $X_i$.
We also choose Kuranishi structures
and a system of their CF-perturbations
on the moduli spaces of the pseudo-holomorphic
disks which appear in the definition of
$\mathfrak{Fuk}(X_i;\mathbb L_i)$.
(See Theorem~\ref{thekuraexist}
and Proposition~\ref{prop330}.)

\end{choice}

\begin{choice}\label{chicesec152}
Suppose $(X_i,\omega_i)$, $\mathbb L_i$, $\mathbb L_{12}$
etc.\ are as in Situation \ref{situ61}.
We take $-J_{X_1} \times J_{X_2}$ as the
compatible almost complex structure of
$-X_1\times X_2$.
\begin{enumerate}\itemsep=0pt
\item[(1)]
We choose Kuranishi structures and
their CF-perturbations of the moduli spaces
used to define filtered $A_{\infty}$ category
$\mathfrak{Fuk}(-X_1\times X_2,\mathbb L_{12})$.
This construction is the same as Theorem~\ref{thekuraexist}
and Proposition~\ref{prop330},
except we use the compactification ${\mathcal M}'(L_{12};\vec a;E)$
etc.\ which we discussed in Section~\ref{sec:directcomp} instead of
${\mathcal M}(L_{12};\vec a;E)$.
\item[(2)]
Suppose we made Choices \ref{chicesec151} and~\ref{chicesec152}\,(1).
Finally, we take Kuranishi structures and their CF-perturbations
of the moduli spaces of pseudo-holomorphic quilts
appearing in the construction of the filtered $A_{\infty}$
tri-functor in Theorem~\ref{trimain}.
See Theorem~\ref{therem530} and Proposition~\ref{prop536}.
These Kuranishi structures and their CF-perturbations
should be compatible with those we took
already in Choice \ref{chicesec151} and
item (1).
\end{enumerate}

\end{choice}
\begin{rem}
In Choice \ref{chicesec152}, we take the compactification ${\mathcal M}'(L_{12};\vec a;E)$
in Section~\ref{sec:directcomp}
to define a filtered $A_{\infty}$ category $\mathfrak{Fuk}(-X_1\times X_2,\mathbb L_{12})$.
We can use the stable map compactification~${\mathcal M}(L_{12};\vec a;E)$ also to define
a filtered $A_{\infty}$ category whose objects are identified with
elements of $\mathbb L_{12}$. We will show in Section~\ref{sec:proofindepen}
that those two categories
are homotopy equivalent.
\end{rem}
\begin{thm}\label{thm1444}
Suppose we take two different ways of
Choice {\rm\ref{chicesec151}}, which we denote by $\Xi_{i,1}$
and $\Xi_{i,2}$, respectively.
We denote by $\mathfrak{Fuk}(X_i;\mathbb L_i;\Xi_{i,1})$,
$\mathfrak{Fuk}(X_i;\mathbb L_i;\Xi_{i,2})$,
the filtered $A_{\infty}$ categories
obtained by these two different choices, respectively.
\begin{enumerate}\itemsep=0pt
\item[$(1)$]
The filtered $A_{\infty}$ category
$\mathfrak{Fuk}(X_i;\mathbb L_i;\Xi_{i,1})$ is
strongly homotopy equivalent to
$\mathfrak{Fuk}(X_i;\mathbb L_i;\allowbreak\Xi_{i,2})$.
\item[$(2)$]
There is a choice of the strong homotopy equivalence in item $(1)$
which is canonical up to homotopy.
\end{enumerate}

When we take two different ways of Choice {\rm\ref{chicesec152}\,(1)},
$\Xi_{12,1}$
and $\Xi_{12,2}$, then for two filtered $A_{\infty}$
categories
$\mathfrak{Fuk}(-X_1,\times X_2,\mathbb L_{12};\Xi_{12,1})$,
$\mathfrak{Fuk}(-X_1,\times X_2,\mathbb L_{12};\Xi_{12,2})$
the same conclusion as above {\rm (1)}, {\rm (2)} holds.
\end{thm}

The proof is given in Section~\ref{sec:welldefAinfcat}.

We denote by
\begin{gather*}
\begin{split}
& \mathscr G^i \colon\
\mathfrak{Fuk}(X_i;\mathbb L_i;\Xi_{i,1})
\to \mathfrak{Fuk}(X_i;\mathbb L_i;\Xi_{i,2}),\nonumber
\\
& \mathscr G^{12} \colon\
\mathfrak{Fuk}(-X_1,\times X_2,\mathbb L_{12};\Xi_{12,1})
\to \mathfrak{Fuk}(-X_1,\times X_2,\mathbb L_{12};\Xi_{12,2}),
\end{split}
\end{gather*}
the strong homotopy equivalences given
in Theorem~\ref{thm1444}.

\begin{situ}\label{situ1415}
Suppose we are in the situation of Theorem~\ref{thm1444}.
In particular, we made a~choice of $\Xi_{i,j}$ for $i=1,2$, $j=1,2$
and of $\Xi_{12,j}$ for $j=1,2$.

For each $j=1,2$, we make Choice \ref{chicesec152}\,(2)
so that this choice is compatible with $\Xi_{1,j}$, $\Xi_{2,j}$,
$\Xi_{12,j}$ at the boundaries.
We denote this choice by \smash{$\Xi^{\rm quilt}_{12,j}$}.
By Corollary \ref{cor73}, those choices determine a filtered $A_{\infty}$
functor
\[
\mathfrak{Fukst}(-X_1,\times X_2,\mathbb L_{12};\Xi_{12,j})
\to
\mathcal{FUNC}
(\mathfrak{Fukst}(X_1;\mathbb L_1;\Xi_{1,j}),
\mathfrak{Fukst}(X_2;\mathbb L_2;\Xi_{2,j})).
\]
(Here we put $\mathfrak{st}$ to denote the associated
strict category.)
We denote~by \smash{$\mathcal{MWW}^{\Xi^{\rm quilt}_{12,j}}$}
this filtered~$A_{\infty}$
functor.
\end{situ}

\begin{thm}\label{thm1466}
In Situation {\rm\ref{situ1415}}, the next diagram commutes up to homotopy equivalence:
\begin{equation*}
\begin{CD}
\displaystyle
\!\!\!\!\!\!\!\!\!\!\!\!\!\!\!\!\!\!\!\!\!\!\mathfrak{Fukst}(X_1;\mathbb L_{1};\Xi_{1,1})
\atop\displaystyle\times
\mathfrak{Fuk}(-X_1 \times X_2;\mathbb L_{12};\Xi_{12,1})
@ >{\mathcal{MWW}^{\Xi^{\rm quilt}_{12,1}}}>>
\mathfrak{Fukst}(X_2;\mathbb L_{2};\Xi_{2,1})
\\
@ V{\mathscr G^1 \times \mathscr G^{12}}VV @ V
{\mathscr G^2}VV\\
\displaystyle
\!\!\!\!\!\!\!\!\!\!\!\!\!\!\!\!\!\!\!\!\!\!
\mathfrak{Fukst}(X_1;\mathbb L_{1};\Xi_{1,2})
\atop\displaystyle\times
\mathfrak{Fuk}(-X_1 \times X_2;\mathbb L_{12};\Xi_{12,2})
@ >{\mathcal{MWW}^{\Xi^{\rm quilt}_{12,2}}}>>
\mathfrak{Fukst}(X_2;\mathbb L_{2};\Xi_{2,2}).
\end{CD}
\end{equation*}

\end{thm}
The proof is in Section~\ref{sec:proofindepen}.

\subsection{Higher pseudo-isotopy}
\label{sec:highpiso}

We will prove Theorem~\ref{thm1444}\,(1) by constructing pseudo-isotopy
between two filtered $A_{\infty}$ algebras.
As we will see in the next subsection, pseudo-isotopy
induces a homotopy equivalence.
To prove Theorem~\ref{thm1444}\,(2),
we need to show that the homotopy equivalence is independent
of the choice of pseudo-isotopy up to homotopy.
To prove it, we use pseudo-isotopy of pseudo-isotopies.

As we explained in \cite[Section 7.2.3]{fooobook2}
and Section~\ref{subsec:Ainfalgim},
during the construction of various structures,
to obtain structure operations directly from geometry
(moduli spaces),
we need to fix an arbitrary but finite $E_0$ and
define structure operations up to energy level $E_0$
only.
We then take homotopy inductive limit
as $E_0 \to \infty$.
To work out homotopy inductive limit argument, we need one extra parameter.
Namely, to obtain a pseudo-isotopy of pseudo-isotopies\index{pseudo-isotopy of pseudo-isotopies}
we need to define a pseudo-isotopy of pseudo-isotopies
up to energy level $E_0$
and a pseudo-isotopy between two
pseudo-isotopies of pseudo-isotopies,
one up to energy level $E_0$
and the other up to energy level $E_1$.
In other words, we need
pseudo-isotopy of
pseudo-isotopies of pseudo-isotopies.
To define such objects,
it seems simpler to define
a family of filtered $A_{\infty}$
structures parametrized by a cornered manifold.
Such a construction is worked out in detail
in \cite[Section 21]{foootech22}, \cite[Chapter~22]{fooonewbook} and \cite{AFOOO}.
In this subsection, we provide its summary.

Let $P$ be an $n$-dimensional manifold with corners.
We consider only the case when $P \subset \R^n$.
We consider
$\tilde L \times_X \tilde L$,
where $L = \bigl(\tilde L,i_L\bigr)$ is an immersed Lagrangian submanifold of
a symplectic manifold $X$, which has
clean self intersection.

Let $\Theta$ is a principal ${\rm O}(1)$ bundle on $\tilde L \times_X \tilde L$
and we put
\begin{gather*}
CF(P \times L;\Theta;\Lambda_0) =
\Omega\bigl(P \times \bigl(\tilde L \times_X \tilde L\bigr);\Theta\bigr) \,\widehat{\otimes}\, \Lambda_0,
\\
CF(P \times L;\Theta;\R) =
\Omega\bigl(P \times \bigl(\tilde L \times_X \tilde L\bigr);\Theta\bigr).
\end{gather*}
(Compare \eqref{form315}.)
\begin{defn}[{\cite[Definition 21.27]{fooonewbook}}]\label{defn1926}
A multilinear map
$
F \colon B_k(CF(P \times L;\Theta;\R)[1]) \allowbreak\to CF(P \times L;\Theta;\R)[1]
$
is said to be
{\it pointwise in $P$ direction}
\index{pointwise in $P$ direction}
if the following holds.
For each $I,J_1,\dots,J_k \subseteq \{1,\dots,d\}$
and ${\bf t} \in P$, there exists a continuous map
\[
F^{\bf t}_{I;J_1,\dots,J_k} \colon\  B_k(CF(L;\Theta;\R)[1]) \to CF(L;\Theta;\R)[1]
\]
such that
\begin{gather*}
F( d t_{J_1} \wedge h_1, \dots,  d t_{J_k} \wedge h_k)\vert_{\{{\bf t}\} \times L} =
\sum_{I}  d t_I \wedge  d t_{J_1} \wedge \dots \wedge  d t_{J_k}
\wedge F^{\bf t}_{I;J_1,\dots,J_k}\bigl(h^{\bf t}_1,\dots,h^{\bf t}_k\bigr),
\end{gather*}
where \smash{$\vert_{\{{\bf t}\} \times L}$}
means the restriction to \smash{$\{{\bf t}\} \times
\bigl(\tilde L\times_X \tilde L\bigr)$}. Moreover, \smash{$F^{\bf t}_{I;J_1,\dots,J_k}$}
depends smoothly on~${\bf t}$ with respect to the operator topology.
Here $h^{\bf t}_i$ is the restriction of $h_i$ to $\{{\bf t}\} \times \bigl(\tilde L\times_X \tilde L\bigr)$
and $t_I = t_{i_1,\dots,i_{\vert I\vert}}$ if $I = \{i_1,\dots,i_{\vert I\vert}\}$
with $i_1 < \dots < i_{\vert I\vert}$.
\end{defn}
Here the continuity of \smash{$F^{\bf t}_{I;J_1,\dots,J_k}$} mentioned above
is one in $C^{\infty}$ topology.
\begin{defn}\label{PparaAinfdef}
A {\it $P$-parametrized family of $G$-gapped filtered $A_{\infty}$
structures on $CF(P \times L;\Theta;\Lambda_0)$}
\index{parametrized family of $G$-gapped filtered $A_{\infty}$ algebra}
is $\bigl\{\mathfrak m^P_{k,\beta}\bigr\}$
for $\beta \in G$ and $k = 0,1,2,\dots$,
that satisfies the following:
\begin{enumerate}\itemsep=0pt
\item[(1)]
$
\mathfrak m^P_{k,\beta} \colon B_k(\Omega(P\times L)[1]) \to \Omega(P\times L)[1]
$
is a multilinear map of degree $1$.
\item[(2)]
$\mathfrak m^P_{k,\beta}$ is
pointwise in $P$ direction if $\beta \ne \beta_0$.
\item[(3)]
$\mathfrak m^P_{k,\beta_0} = 0$ for $k \ne 1,2$.
\item[(4)]
$\mathfrak m^P_{1,\beta_0}(h) = (-1)^*dh$. Here $d$ is the de Rham differential
and $*$ is as in \eqref{form3420000}.
\item[(5)]
$\mathfrak m^P_{k,\beta}$ satisfies the following $A_{\infty}$ relation:
\begin{gather}
\sum_{k_1{+}k_2=k{+}1}\sum_{\beta_1{+}\beta_2{=}\beta}\!\!\!\sum_{i=1}^{k-k_2+1}\!\!\! (-1)^*{\mathfrak m}^P_{k_1,\beta_1}\!\bigl(h_1,\ldots,{\mathfrak m}^P_{k_2,\beta_2}(h_i,\ldots,h_{i+k_2-1}),\ldots,h_{k}\bigr)
  = 0,\!\!\!\label{Ainfinityrelbeta2}
\end{gather}
where $* = \deg' h_1 + \dots + \deg'h_{i-1}$.
\end{enumerate}
\end{defn}
We put
$
\mathfrak m^P_k = \sum_{\beta \in G} T^{\beta} \mathfrak m^P_{k,\beta}$.
\eqref{Ainfinityrelbeta2} the implies $A_{\infty}$ relation for
$\mathfrak m^P_k$.
\begin{rem}
We may choose $\mathfrak m^P_{2,\beta_0}$ such that
$\mathfrak m^P_{2,\beta_0}(h_1 \wedge h_2) = (-1)^* h_1 \wedge h_2$
holds if $h_1$ or $h_2$ are supported on the diagonal component.
Here $\wedge$ is the wedge product and $* = \deg h_1(\deg h_2+1)$.
See Remark~\ref{rem134141}.
\end{rem}
\begin{defn}
A {\it partial $P$-parametrized family of $G$-gapped filtered $A_{\infty}$ algebra structures on $CF(P \times L;\Theta;\Lambda_0)$ of energy cut level $E$ and of minimal energy $e_0$} \index{energy cut level}
\index{minimal energy}
\index{partial $P$-parametrized family of $G$-gapped filtered $A_{\infty}$ algebras}
is $\bigl\{\mathfrak m^P_{k,\beta}\bigr\}$
that satisfies the same properties as above except the following points:
\begin{enumerate}\itemsep=0pt
\item[(a)]
$\mathfrak m^P_{k,\beta}$ is defined only for
$\beta \in G, k = 0,1,2,\dots$ with $\beta + ke_0\le E$.
\item[(b)]
We require the $A_{\infty}$ relation \eqref{Ainfinityrelbeta2}
only for $\beta$, $k$ with $\beta + ke_0 \le E$.
\item[(c)]
$\mathfrak m^P_{k,\beta} = 0$ if $0 < \beta < e_0$.
\end{enumerate}

\end{defn}
We can restrict a $P$-parametrized family of $G$-gapped filtered $A_{\infty}$ algebra
structures to the normalized boundary of $P$ and corners of $P$
in an obvious way.
\begin{exm}
The $[0,1]$-parametrized family of $G$-gapped filtered $A_{\infty}$ algebra
structures of energy cut level $E$ is
nothing but a pseudo-isotopy modulo $T^E$ of
$G$-gapped filtered $A_{\infty}$ algebra as in Definition~\ref{pisotopydef}.
\end{exm}

We next define the notion of collared structure.
We define the case
$
P = [0,1]^n,
$
only. See~\cite{foootech22,fooonewbook} for the general case.
\begin{defn}\label{Pparamet}
Let $\bigl\{ \mathfrak m^{P}_{k,\beta}\bigr\}$
be a $P$-parametrized family of $G$-gapped partial
$A_{\infty}$-structure of energy cut level $E_0$ and minimal energy $e_0$.
We say it is {\it $\tau$-collared}\index{$\tau$-collared}
if the following holds for~${(t_1,\dots,t_n) \in [0,1]^n}$.

We consider the case $t_1 \le t_2 \le \dots \le t_{n-1} \le t_n$,
only. The general case is similar.

Let $t_i \in [0,\tau]$, $t_{j+1} \in [1-\tau,1]$
and $t_{i+1},\dots,t_{j} \in (\tau,1-\tau)$.
Let $(s_1,\dots,s_{i+n-j}) = (t_1,\dots,t_i,\allowbreak t_{j+1},\dots,t_n)$
and
$(s'_1,\dots,s'_{j-i}) = (t_{i+1},\dots,t_{j})$.
A differential form of $P$ in a neighborhood is written as
$\sum f_{II'}  d s_{I} \wedge  d s^{\prime}_{I'}$,
where $ d s_{I}$ are wedge products of~$ d s_i$'s and
$ d s^{\prime}_{I'}$ are wedge products of $ d s'_i$'s.

By definition, $\mathfrak m^{P}_{k,\beta}$ is written
on this neighborhood as the form
\[
\mathfrak m^{P}_{k,\beta}(h_1,\dots,h_k)
=
\sum_{I,I'}
 d s_{I} \wedge  d s^{\prime}_{I'}
\wedge \mathfrak m^{s,s'}_{k,\beta;I,I'}(h_1,\dots,h_k),
\]
where $h_i$ are smooth differential forms
$P \times L$ twisted by $\Theta$ which does not have
$ d t_i$ components.
We now require:
\begin{enumerate}\itemsep=0pt
\item[(1)]
\smash{$\mathfrak m^{s,s'}_{k,\beta;I,I'}(h_1,\dots,h_k) = 0$}
unless $I = \varnothing$.
\item[(2)]
If $I=\varnothing$,
\smash{$\mathfrak m^{s,s'}_{k,\beta;I,\varnothing}(h_1,\dots,h_k)$}
is independent of $s \in [0,\tau]^i \times [1-\tau,1]^{n-j}$.
\end{enumerate}
We say $\bigl\{ \mathfrak m^{P}_{k,\beta}\bigr\}$ is collared if it is
$\tau$-collared for some $\tau>0$.
\end{defn}
The main lemma we will use to prove
Theorems~\ref{thm1444} and \ref{thm1466}
is Proposition~\ref{prop2013} below.
\begin{situ}\label{situ2011}
Let $P$ be a manifold with corner and $E_1 > E_0 \ge 0$, $e_0 > 0$.
We assume that we are given the following objects:
\begin{enumerate}\itemsep=0pt
\item[(1)]
A $P\times [0,1]$-parametrized collared
partial
$A_{\infty}$ structure \smash{$\bigl\{ \mathfrak m^{P\times [0,1]}_{k,\beta}\bigr\}$}
of energy cut level $E_0$ and of minimal energy $e_0$ on $CF(P \times L;\Theta;\Lambda_0)$.
\item[(2)]
A collared partial
$A_{\infty}$ structure \smash{$\bigl\{ \mathfrak m^{P\times \{1\}}_{k,\beta}\bigr\}$}
on $CF(P\times \{1\} \times L;\Theta;\Lambda_0)$
 of energy cut level~$E_1$~is given.
We require that it coincides with the restriction of
 \smash{$\bigl\{ \mathfrak m^{P\times [0,1]}_{k,\beta}\bigr\}$}
to $P \times \{1\}$ as the partial
$A_{\infty}$ structures of energy cut level $E_0$.
\item[(3)]
Assume $\partial P = \coprod \partial_i P$ is the decomposition of the normalized
boundary of~$P$ into the connected components.
Then for each $i$, we are given a collared filtered
$A_{\infty}$ structure \smash{$\bigl\{ \mathfrak m^{\partial_i P\times [0,1]}_{k,\beta}\bigr\}$}
of energy cut level $E_1$.
We  require that it coincides with the restriction of structure \smash{$\bigl\{ \mathfrak m^{P\times [0,1]}_{k,\beta}\bigr\}$} to ${\partial_i P \times [0,1]}$
as the partial
$A_{\infty}$ structures of energy cut level $E_0$.
\item[(4)]
We assume that the restriction of the
structure \smash{$\bigl\{ \mathfrak m^{P\times \{1\}}_{k,\beta}\bigr\}$} in item (2)
coincides with the structure
\smash{$\bigl\{ \mathfrak m^{\partial_i P\times [0,1]}_{k,\beta}\bigr\}$} in item (3) on
 $\partial_i P \times \{1\}$.
\item[(5)]
Suppose that the images of $\partial_i P$ and $\partial_j P$ intersect
each other in
$P$ at the component $\partial_{ij} P$ of the codimension 2 corner of $P$.
(Note that the case $i=j$ is included.
In this case, $\partial_{ii} P$ is the `self intersection' of
$\partial_i P$.)
We then assume that the restriction of
\smash{$\bigl\{ \mathfrak m^{\partial_i P\times [0,1]}_{k,\beta}\bigr\}$}
to~${\partial_{ij} P \times [0,1]}$
coincides with
the restriction of
\smash{$\bigl\{ \mathfrak m^{\partial_j P\times [0,1]}_{k,\beta}\bigr\}$}.
\end{enumerate}
\end{situ}

\begin{figure}[ht]
\centering
\includegraphics[scale=0.4]{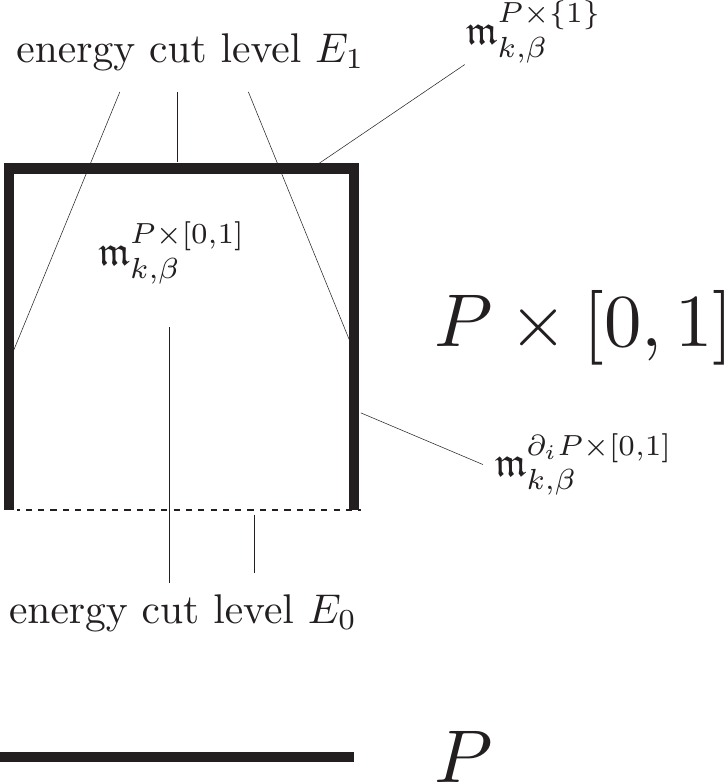}
\caption{$P\times [0,1]$-parametrized family.}
\label{Figuresec14-2}
\end{figure}

See Figure~\ref{Figuresec14-2}.

\begin{prop}\label{prop2013}
In Situation {\rm\ref{situ2011}}, there exists
a collared partial $P$-parametrized family of $G$-gapped filtered $A_{\infty}$ algebra structures on $CF(P \times L;\Theta;\Lambda_0)$ of energy cut level $E_1$ and of minimal energy $e_0$, which we denote by
\smash{$\bigl\{ \mathfrak m^{P\times [0,1]}_{+,k,\beta}\bigr\}$}.
It has the following properties:
\begin{enumerate}\itemsep=0pt
\item[$(1)$] If we regard \smash{$\bigl\{ \mathfrak m^{P\times [0,1]}_{+,k,\beta}\bigr\}$}
as a partial structure of energy cut level $E_0$,
then it coincides with~\smash{$\bigl\{ \mathfrak m^{P\times [0,1]}_{k,\beta}\bigr\}$}.
\item[$(2)$]
If we restrict \smash{$\bigl\{ \mathfrak m^{P\times [0,1]}_{+,k,\beta}\bigr\}$}
to
$P \times \{1\}$, then it coincides with
the structure \smash{$\bigl\{ \mathfrak m^{P\times \{1\}}_{k,\beta}\bigr\}$} in
Situation {\rm\ref{situ2011}\,(2)}.
\item[$(3)$]
If we restrict \smash{$\bigl\{ \mathfrak m^{P\times [0,1]}_{+,k,\beta}\bigr\}$}
to
$\partial_i P \times [0,1]$, then it coincides with
the structure \smash{$\bigl\{ \mathfrak m^{\partial_i P\times [0,1]}_{k,\beta}\bigr\}$} in
Situation {\rm\ref{situ2011}\,(3)}.
\end{enumerate}

\end{prop}
In the case when $P$ is a one point, this is nothing but
Lemma~\ref{lem339}. The proposition in this generality
is proved in \cite{foootech22}, \cite[Proposition 22.13]{fooonewbook}.
See also \cite[Section 14]{fooo091}.\footnote{The
singular homology version (of the case $P=[0,1]$)
is \cite[Theorem 7.2.212]{fooobook2}.
Actually singular homology version is
harder to state and prove.}

\subsection[Well definedness of a filtered $A_\infty$
category up to strong homotopy equivalence]{Well definedness of a filtered $\boldsymbol{A_{\infty}}$
category\\ up to strong homotopy equivalence}
\label{sec:welldefAinfcat}

\begin{proof}[Proof of Theorem~\ref{thm1444}\,(1)]
We prove the case of $(X_1,\omega_1,\mathbb L_1)$.
The other cases are the same.
By the trick we used in Section~\ref{subsec:Ainfcatim},
it suffices to consider the case when $\mathbb L_1$ consists
of a single immersed Lagrangian submanifold $L_1$.
Let $J_{1,j}$ ($j=1,2$) be the almost complex structure
on~$X_1$ chosen as a part of Choice~$\Xi_{1,j}$.
We take one parameter family of
compatible almost complex structures $J_{1,s}$
parametrized by $s \in [0,1]$ such that
$J_{1,s} = J_1$ for~${s \in [0,\tau]}$
and $J_{1,s} = J_2$ for~${s \in [1-\tau,1]}$.
We use the notations of Section~\ref{subsec:Ainfalgim}.
The moduli space ${\mathcal M}_{k+1}(L_1;E)$ is as in~\eqref{form317}.
We write ${\mathcal M}_{k+1}((L_1,J);E)$ to specify the almost
complex structure $J$ we use.

Hereafter, in this subsection we omit the suffix $1$ and
write $L$, $X$, $J_j$, $\Xi_j$ etc.
in place of $L_1$, $X_1$, $J_{1,j}$, $\Xi_{1,j}$ etc.
\begin{defn}
We define
\[
{\mathcal M}_{k+1}(L;E;[0,1]_s)
= \bigcup_{s \in [0,1]_s}
{\mathcal M}_{k+1}((L,J_{s});E) \times \{s\}.
\]
The evaluation map
\begin{equation}\label{146666}
{\rm ev} = ({\rm ev}_0,\dots,{\rm ev}_k) \colon\
{\mathcal M}_{k+1}(L;E;[0,1]_s)
\to L^{k+1}
\end{equation}
is defined by Definition~\ref{defn32222}.
The other evaluation map
\begin{equation}\label{form147777}
{\rm ev}_{[0,1]_s} \colon\ {\mathcal M}_{k+1}(L;E;[0,1]_s)
\to [0,1]_s
\end{equation}
is defined by sending
${\mathcal M}_{k+1}((L,J_{s});E) \times \{s\}$ to $s$.\footnote{
We use the notation $[0,1]_s$ here to distinguish
it from the interval which we use for different parameter.}
\end{defn}
\begin{prop}\label{prop14555}
We can define a topology on ${\mathcal M}_{k+1}(L;E;[0,1]_s)$,
the {\rm stable map topology}, \index{stable map topology}
which is Hausdorff and compact.
There exists a system of Kuranishi structures with
a boundary and corners on
${\mathcal M}_{k+1}(L;E;[0,1]_s)$ for various
$k$ and $E$ with the following properties:
\begin{enumerate}\itemsep=0pt
\item[$(1)$]
The evaluation maps \eqref{146666} and
\eqref{form147777} extend to strongly smooth maps.
The map $({\rm ev}_0,\allowbreak{\rm ev}_{[0,1]_s})$ is
weakly submersive.\footnote{Since $[0,1]_s$ has boundary,
we need to define weakly submersivity a bit carefully.
See \cite[Section 25]{foootech22}, \cite[Chapter 26]{fooonewbook}.}
\item[$(2)$]
The normalized boundary of ${\mathcal M}_{k+1}(L;E;[0,1]_s)$
is a disjoint union of the following two types of spaces:
\begin{enumerate}\itemsep=0pt
\item[$(I)$] The fiber product
\[
{\mathcal M}_{k_1+1}(L;E_1;[0,1]_s)
 {}_{({\rm ev}_i,{\rm ev}_{[0,1]_s})}
\times_{({\rm ev}_0,{\rm ev}_{[0,1]_s})}
{\mathcal M}_{k_2+1}(L;E_2;[0,1]_s),
\]
over $L \times [0,1]_s$.
Here $k_1 + k_2 = k$, $i=1,\dots,k_2$, $E_1 + E_2 = E$.
The fiber product carries a Kuranishi structure
because of the weak submersivity of $({\rm ev}_0,{\rm ev}_{[0,1]})$.
See {\rm\cite[\emph{Definition} 4.9]{foootech2}, \cite[\emph{Chapter} 26]{fooonewbook}}.
\item[$(II)$]
The inverse image
$
{\rm ev}_{[0,1]}^{-1}(\partial [0,1]_s)
\subset {\mathcal M}_{k_1+1}(L;E;[0,1]_s)$.
\end{enumerate}
\item[$(3)$]
For sufficiently small $\tau$, the following holds.
The restriction of the Kuranishi structure to
\smash{${\rm ev}_{[0,1]}^{-1}([0,\tau])$}
coincides with the direct product of the
trivial Kuranishi structure on $[0,\tau]$
and the Kuranishi structure of~${\mathcal M}_{k+1}((L,J_{1});E)$, which is a part of the data
$\Xi_{1}$.
The restriction of the Kuranishi structure to
\smash{${\rm ev}_{[0,1]}^{-1}([1-\tau,1])$}
coincides with the direct product of the
trivial Kuranishi structure on $[1-\tau,1]$
and the
Kuranishi structure of~${\mathcal M}_{k+1}((L,\allowbreak J_{2});E)$ which is a part of the data
$\Xi_{2}$.
\item[$(4)$]
The orientation bundles are compatible with the description of the
boundary in item {\rm(2)}.
\end{enumerate}

\end{prop}
The proof is a one parameter version of the
proof of Theorem~\ref{thekuraexist}.
Note that the fact that our Kuranishi structure (and the moduli space)
is of product type near the boundary of Type~(II), which is stated
as item (3), is a consequence of our choice of
family of almost complex structures.

Let $G(L;\Xi_{j})$ be the discrete submonoid
defined in Definition~\ref{defn314}.
Note that it depends on the almost complex structure and so
on $\Xi_{j}$.
However, we may choose a discrete submonoid
$G(L)$ with the following properties:
\begin{enumerate}\itemsep=0pt\setlength{\leftskip}{0.40cm}
\item[(Mo.1)]
The submonoid
$G(L)$ contains both $G(L;\Xi_{1})$
and $G(L;\Xi_{2})$.
\item[(Mo.2)]
If ${\mathcal M}_{k+1}(L;E;[0,1])$
is nonempty then $E \in G(L)$.
\end{enumerate}
We put $G(L) = \{E_1,E_2,\dots,E_n,\dots\}$,
where $E_1 < E_2 < \cdots$.
Let $E_i \in G(L)$.
Proposition~\ref{prop330} assigns a CF-perturbation
of ${\mathcal M}_{k+1}((L,J_{j});E)$ with $E \le E_i$
($j=1,2$).
These CF-perturbations and the Kuranishi structures
on which they are defined are parts of the data
$\Xi_{j}$.
We denote this CF-perturbation by \smash{$\widehat{\mathfrak S}(\Xi_{j};E_i)$}.
\begin{prop}\label{prop1416}
There exists a system of CF-perturbations
\smash{$\widehat{\mathfrak S}([0,1]_s;E_i)$} on outer collarings of thickenings of
${\mathcal M}_{k+1}(L;E;[0,1]_s)$ with $E \le E_i$
with the following properties:
\begin{enumerate}\itemsep=0pt
\item[$(1)$] The CF-perturbation \smash{$\widehat{\mathfrak S}([0,1]_s;E)$} is transversal to $0$.
\item[$(2)$]
The map $({\rm ev}_0,{\rm ev}_{[0,1]_s})$ is
strongly submersive with respect to
\smash{$\widehat{\mathfrak S}([0,1]_s;E)$}.
\item[$(3)$]
The restriction of \smash{$\widehat{\mathfrak S}([0,1]_s;E)$}
to the boundary in Proposition {\rm\ref{prop14555}\,(2),~(I)}
coincides with the fiber product CF-perturbation
$($see {\rm\cite[\emph{Definition} 10.13]{foootech22}}$)$,
which is well-defined by item~$(2)$.
\item[$(4)$]
For sufficiently small $\tau$, the following holds.
The restriction of \smash{$\widehat{\mathfrak S}([0,1]_s;E)$} to
\smash{${\rm ev}_{[0,1]_s}^{-1}([0,\tau])$}
coincides with the pullback of \smash{$\widehat{\mathfrak S}(\Xi_{1};E)$}.
The restriction of \smash{$\widehat{\mathfrak S}([0,1]_s;E)$} to
\smash{${\rm ev}_{[0,1]_s}^{-1}([1-\tau,1])$}
coincides with the pullback of
\smash{$\widehat{\mathfrak S}(\Xi_{2};E)$}.
\end{enumerate}
\end{prop}

\begin{proof}
We define the CF-perturbation on the neighborhood of
the boundary component in Proposition~\ref{prop14555}\,(2),~(II)
by item (4). Then we can extend it using the
(relative version) of the existence of CF-perturbation.
(See \cite[Chapter 17]{foootech22,fooonewbook}.)
\end{proof}

\begin{rem}
During the proof of Proposition~\ref{prop1416}, we
construct Kuranishi structures on which our
CF-perturbations are defined at the same time.
See \cite{foootech22,fooonewbook} for this point.
\end{rem}
Now for $\beta = E \in G(L)$ with $E \le E_i$, we define
\[
\mathfrak m_{k,E}^{E_i;[0,1]_s} \colon\
CF([0,1]_s \times L;\Theta;\R)^{\otimes k}
\to
CF([0,1]_s \times L;\Theta;\R)
\]
by
\begin{align}
\mathfrak m_{k,E}^{E_i;[0,1]_s}(h_1,\dots,h_k) ={}&
({\rm ev}_0,{\rm ev}_{[0,1]_s})!
(({\rm ev}_1,{\rm ev}_{[0,1]_s})^*(h_1) \wedge\nonumber \\
&
\dots \wedge
({\rm ev}_k,{\rm ev}_{[0,1]_s})^*(h_k);
\widehat{\mathfrak S}([0,1]_s;E_i)).\label{form148}
\end{align}
Here the integration by parts is taken on the space
${\mathcal M}_{k+1}(L;E;[0,1]_s)$ using the
CF-perturba\-tion~\smash{$\widehat{\mathfrak S}([0,1]_s;E_i)$}.
See \cite[Section 2.2.4]{fooonewbook} and \cite[Section 4.1]{ST} for the sign.

\begin{lem}\label{lem1418}
The operations $\bigl\{\mathfrak m_{k,E}^{E_i;[0,1]_s} ; E \le E_i\bigr\}$
define a collared partial $[0,1]_s$-parametrized family of $G(L)$-gapped filtered $A_{\infty}$ algebra structures on $CF([0,1]_s \times L;\Theta;\Lambda_0)$ of energy cut level~$E_i$ and of minimal energy $e_0$.\footnote{The minimal energy is
always $e_0$ in this section. So we omit it from now on.}

\end{lem}
This is a consequence of Proposition~\ref{prop1416}.
Point-wiseness in $[0,1]_s$ direction follows from~\cite[Proposition 22.17]{fooonewbook}.
Moreover, the restrictions of the
structure operations $\smash{\bigl\{\mathfrak m_{k,E}^{E_i;[0,1]_s}} ; E \le E_i\bigr\}$
to $\{0\} \in [0,1]_s$
(resp.\ $\{1\} \in [0,1]_s$) coincide with the
partial $[0,1]_s$-parametrized family of $G(L)$-gapped filtered $A_{\infty}$ algebra structures on $CF(L;\Theta;\Lambda_0)$ of energy cut level $E_i$,
which we used during the construction of
$\mathfrak{Fuk}(X;\mathbb L;\Xi_{1})$
(resp.\ $\mathfrak{Fuk}(X;\mathbb L;\Xi_{2})$).

We remark however that \smash{$\mathfrak m_{k,E}^{E_i;[0,1]_s}$} itself is not
the structure operation of the pseudo-isotopy between
$\mathfrak{Fuk}(X;\mathbb L;\Xi_{1})$ and $\mathfrak{Fuk}(X;\mathbb L;\Xi_{2})$,
which we look for. This is because this structure is yet a partial structure
where \smash{$\mathfrak m_{k,E}^{E_i;[0,1]_s}$} is defined for $E \le E_i$ only.
We will combine the process of taking homotopy limit
with the construction of pseudo-isotopy as follows.

During the construction of the
structure operations of $\mathfrak{Fuk}(X;\mathbb L;\Xi_{j})$,
we used a Kuranishi structure
on ${\mathcal M}_{k+1}((L,J_{j});E) \times [0,1]_t$
and its CF-perturbation such that the
restriction of this CF-perturbation
to ${\mathcal M}_{k+1}((L,J_{j});E) \times \{0\}$
is \smash{$\widehat{\mathfrak S}(\Xi_{j};E_i)$}
and that the
restriction of this CF-perturbation
to ${\mathcal M}_{k+1}((L,J_{j});E) \times \{1\}$
is \smash{$\widehat{\mathfrak S}(\Xi_{j};E_{i+1})$}
(see Lemma~\ref{lem337}).
We denote this CF-perturbation
by \smash{$\widehat{\mathfrak S}([0,1]_t,\Xi_{j};E_i,E_{i+1})$}.
Note that we can take this CF-perturbation so that it is
constant in $t$ direction for $t \in [0,\mu] \cup [1-\mu,1]$.

During the proof of Proposition~\ref{prop336},
we used \smash{$\widehat{\mathfrak S}([0,1]_t,\Xi_{j};E_i,E_{i+1})$}
in the same way as~\eqref{form148} to define
a collared partial $[0,1]_t$-parametrized family of $G(L)$-gapped filtered $A_{\infty}$ algebra structures on $CF([0,1]_t \times L;\Theta;\Lambda_0)$ of energy cut level $E_{i+1}$ and of minimal energy~$e_0$  (see~\eqref{form3833}).
We denote the structure operation of this structure
by \smash{$\{\mathfrak m_{k,E}^{E_i,E_{i+1};[0,1]_t} ; E \le E_i\}$}.

We then construct a pseudo-isotopy of pseudo-isotopies
using the next proposition.
\begin{prop}\label{prop1419}
There exists a system of CF-perturbations, which we denote by $\smash{\widehat{\mathfrak S}}([0,1]_s \times [0,1]_t;E_i,E_{i+1})$, on outer collarings of thickenings of
${\mathcal M}_{k+1}(L;E;[0,1]_s)
\times [0,1]_t$ for $E \le E_{i+1}$
with the following properties:
\begin{enumerate}\itemsep=0pt
\item[$(1)$] The CF-perturbation
\smash{$\widehat{\mathfrak S}([0,1]_s \times [0,1]_t;E_i,E_{i+1})$} is transversal to $0$.
\item[$(2)$]
The map
\[
({\rm ev}_0,{\rm ev}_{[0,1]_s},{\rm ev}_{[0,1]_t})
\colon\
{\mathcal M}_{k+1}(L;E;[0,1]_s)
\times [0,1]_t
\to \bigl(\tilde L \times_X \tilde L\bigr) \times [0,1]_s \times [0,1]_t
\]
is
strongly submersive with respect to
\smash{$\widehat{\mathfrak S}([0,1]_s \times [0,1]_t,\Xi_{j};E_i,E_{i+1})$}.
\item[$(3)$]
We consider the restriction of \smash{$\widehat{\mathfrak S}([0,1]_s \times [0,1]_t;E_i,E_{i+1})$}
to the boundary component, which is a product of $[0,1]_t$ and the
boundary component of ${\mathcal M}_{k+1}(L;E;[0,1]_s)$ in Proposition~{\rm\ref{prop14555}\,(2) (I)}.
It then coincides with the fiber product CF-perturbation,
which is well-defined by item $(2)$.
\item[$(4)$]
For sufficiently small $\tau$, the following holds.
The restriction of \smash{$\widehat{\mathfrak S}([0,1]_s \times [0,1]_t;E_i,E_{i+1})$} to
\smash{${\rm ev}_{[0,1]_s}^{-1}([0,\tau])$}
coincides with the pullback of \smash{$\widehat{\mathfrak S}([0,1]_t,\Xi_{1};E_i,E_{i+1})$}.
The restriction of \smash{$\widehat{\mathfrak S}([0,1]_s \times [0,1]_t;E_i,E_{i+1})$} to
\smash{${\rm ev}_{[0,1]_s}^{-1}([1-\tau,1])$}
coincides with the pullback of $\smash{\widehat{\mathfrak S}([0,1]_t,\Xi_{2}};\allowbreak E_i,E_{i+1})$.
\item[$(5)$]
For sufficiently small $\tau$, the following holds.
The restriction of \smash{$\widehat{\mathfrak S}([0,1]_s \times [0,1]_t;E_i,E_{i+1})$} to
\smash{${\rm ev}_{[0,1]_t}^{-1}([0,\tau])$}
coincides with the pullback of \smash{$\widehat{\mathfrak S}([0,1]_s;E_{i})$}.
The restriction of $\smash{\widehat{\mathfrak S}([0,1]_s} \times [0,1]_t;E_i,E_{i+1})$ to
\smash{${\rm ev}_{[0,1]_s}^{-1}([1-\tau,1])$}
coincides with the pullback of \smash{$\widehat{\mathfrak S}([0,1]_s;E_{i+1})$}.
\end{enumerate}

\end{prop}
\begin{figure}[ht]
\centering
\includegraphics[scale=0.44]{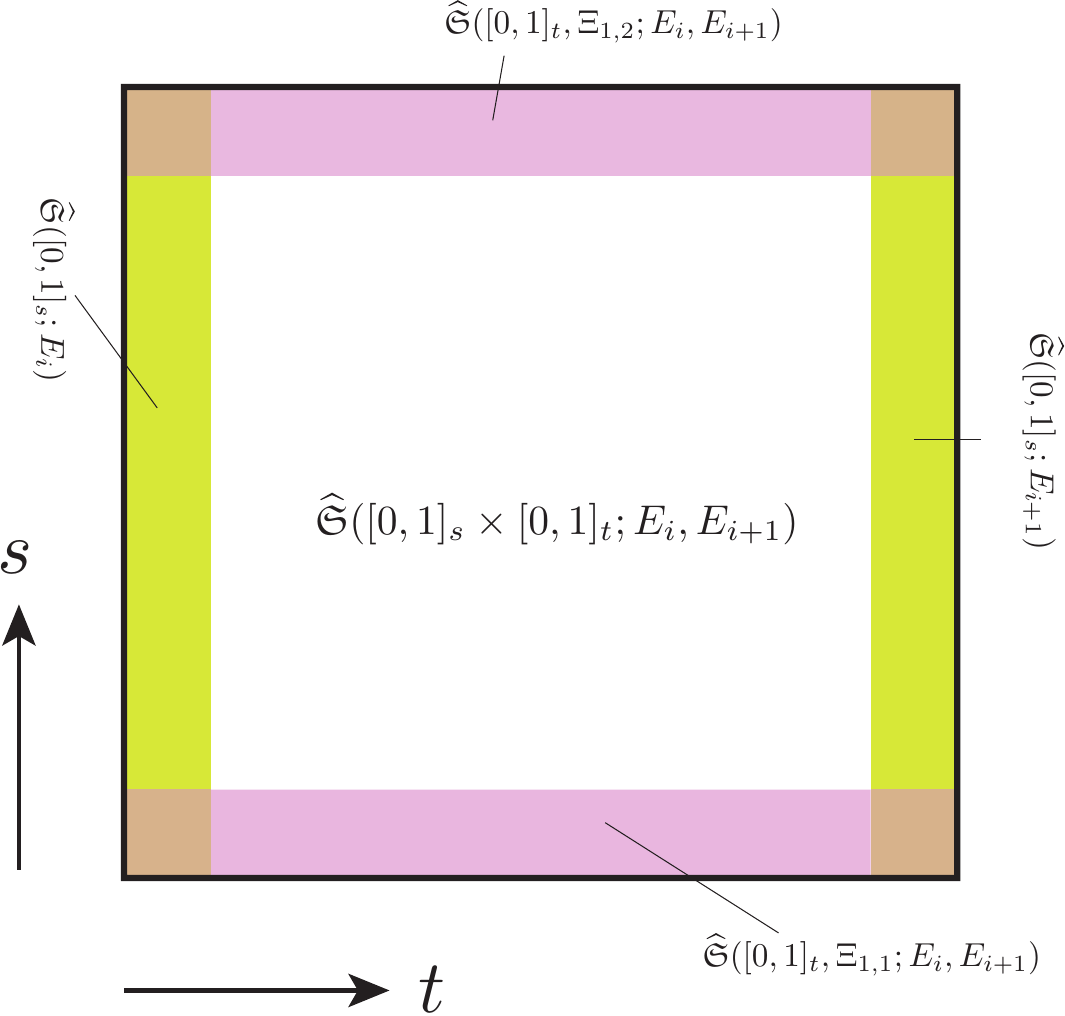}
\caption{$\widehat{\mathfrak S}([0,1]_s \times [0,1]_t;E_i,E_{i+1})$.}
\label{Figuresec14222}
\end{figure}
The proof of Proposition~\ref{prop1419} is by induction
on $E$.
On each step of the induction, the CF-perturbation
on the boundary is determined by the statement we are proving.
So we can extend it. (See \cite[Chapter 17]{foootech22,fooonewbook}.)

We now recall the construction at the end of
Section~\ref{subsec:Ainfalgim}.
We consider the restriction to~${s=0}$.
We use $\Xi_{1}$
to obtain a system of partial filtered $A_{\infty}$
structures
\smash{$\bigl\{\mathfrak m_{k,\beta}^{\Xi_{1},i,E\le E_i}\bigr\}$} and pseudo-isotopies~\smash{$\bigl\{\mathfrak m_{k,\beta}^{[0,1]_t,\Xi_{1},i,E\le E_i}\bigr\}$} among them.
Then we used Lemma~\ref{lem339}, which is
nothing but the case of $P = [0,1]$ of
Proposition~\ref{prop2013}.
We then obtain a sequence of filtered $A_{\infty}$
structures~\smash{$\bigl\{\mathfrak m_{k,\beta}^{\Xi_{1},i}\bigr\}$} on
$CF(L;\Theta;\Lambda_0)$
such that it coincides with
\smash{$\bigl\{\mathfrak m_{k,\beta}^{\Xi_{1},i,E\le E_i}\bigr\}$}
as partial structures with energy cut level $E_i$,
for each $i$.
Moreover, there exists a
pseudo-isotopy
\smash{$\bigl\{\mathfrak m_{k,\beta}^{[0,1]_t,\Xi_{1},i}\bigr\}$}
between~\smash{$\bigl\{\mathfrak m_{k,\beta}^{\Xi_{1},i,E\le E_i}\bigr\}$}
and
\smash{$\bigl\{\mathfrak m_{k,\beta}^{\Xi_{1},i+1,E\le E_{i+1}}\bigr\}$}
which coincides with
\smash{$\bigl\{\mathfrak m_{k,\beta}^{[0,1]_s,\Xi_{1},i,E\le E_i}\bigr\}$}
as a pseudo-isotopy with energy cut level $E_i$.
See Figure~\ref{FigureSec1473}.
\begin{figure}[ht]
\centering
\includegraphics[scale=0.4]{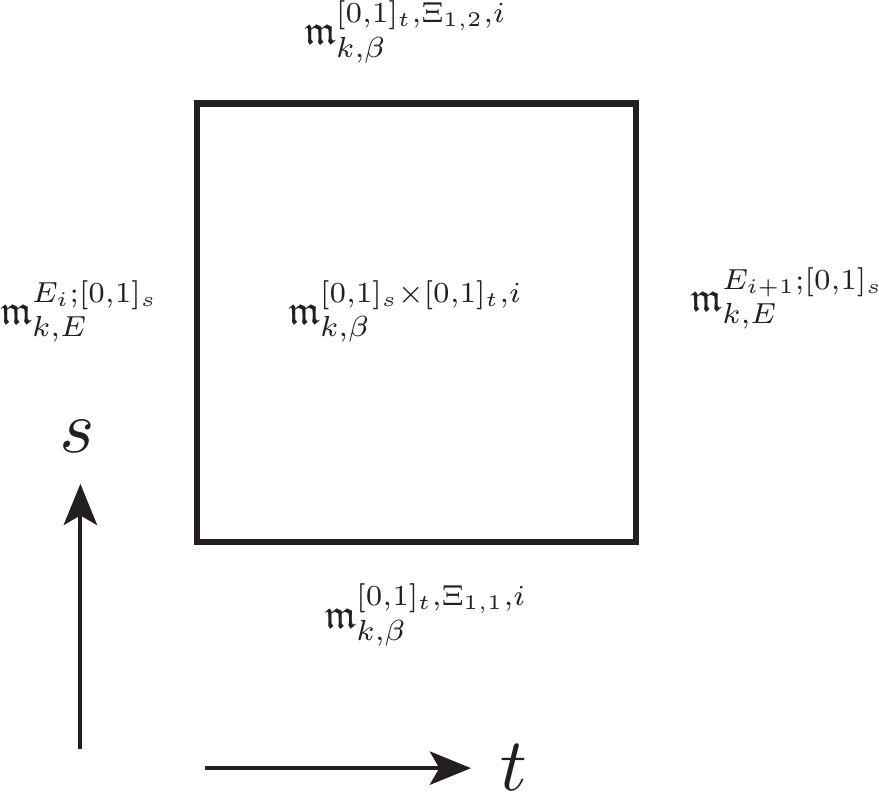}
\caption{Pseudo-isotopy of pseudo-isotopies.}
\label{FigureSec1473}
\end{figure}

We can perform the same construction
for $s=1$ using $\Xi_{1}$ and obtain operations
\smash{$\bigl\{\mathfrak m_{k,\beta}^{[0,1]_t,\Xi_{2},i}\bigr\}$},
\smash{$\bigl\{\mathfrak m_{k,\beta}^{\Xi_{2},i+1,E\le E_{i+1}}\bigr\}$}.

Now we apply Proposition~\ref{prop2013} inductively
and obtain the following.
\begin{lem}\label{lem1420}
There exists a sequence of
$P = [0,1]_s \times [0,1]_t$ parametrized
family of filtered~$A_{\infty}$ algebra
\smash{$\bigl\{\mathfrak m_{k,\beta}^{[0,1]_s \times[0,1]_t,i}\bigr\}$}
on
$CF([0,1]_s\times [0,1]_t \times L;\Theta;\R)$
with the following properties.
\begin{enumerate}\itemsep=0pt
\item[$(1)$]
It coincides with one obtained by
the CF-perturbation
\smash{$\widehat{\mathfrak S}([0,1]_s \times [0,1]_t;E_i,E_{i+1})$}
as partial structures with energy cut level $E_i$.
\item[$(2)$]
Its restriction to
$s=0$ coincides with
\smash{$\bigl\{\mathfrak m_{k,\beta}^{[0,1]_t,\Xi_{1},i}\bigr\}$}.
\item[$(3)$]
Its restriction to
$s=1$ coincides with
\smash{$\bigl\{\mathfrak m_{k,\beta}^{[0,1]_t,\Xi_{2},i}\bigr\}$}.
\item[$(4)$]
Its restriction to $t=1$ coincides with
\smash{$\bigl\{\mathfrak m_{k,E}^{E_{i+1};[0,1]_s} ; E \le E_{i+1}\bigr\}$}
in Lemma~{\rm\ref{lem1418}} as partial structures
with energy cut level $E_{i+1}$.
\item[$(5)$]
Its restriction to $t=0$ coincides with
\smash{$\bigl\{\mathfrak m_{k,E}^{E_{i};[0,1]_s} ; E \le E_i\bigr\}$}
in Lemma~{\rm\ref{lem1418}} as partial structures
with energy cut level $E_{i}$.
\end{enumerate}

\end{lem}
See Figure~\ref{FigureSec1474}.
\begin{figure}[ht]
\centering
\includegraphics[scale=0.48]{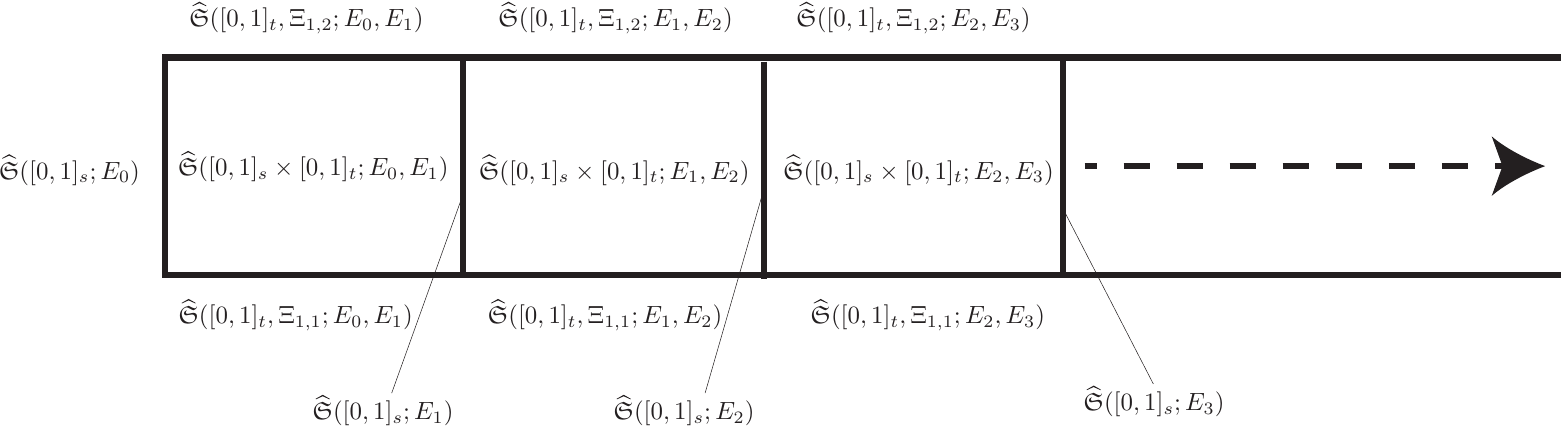}
\caption{Inductive limit construction of pseudo-isotopy.}
\label{FigureSec1474}
\end{figure}

We restrict
\smash{$\bigl\{\mathfrak m_{k,\beta}^{[0,1]_s \times[0,1]_t,0}\bigr\}$}
to $t=0$ and obtain the following.

\begin{cor}\label{cor1421}
There exists a pseudo-isotopy of filtered $A_{\infty}$
structures
\smash{$\bigl\{\mathfrak m_{k,E}^{[0,1]_s}\bigr\}$}
on $CF([0,1]_s\allowbreak \times L;\Theta;\R)$
with the following properties:
\begin{enumerate}\itemsep=0pt
\item[$(1)$]
The structure \smash{$\bigl\{\mathfrak m_{k,E}^{[0,1]_s}\bigr\}$} coincides with
\smash{$\bigl\{\mathfrak m_{k,E}^{E_{0};[0,1]_s} ; E \le E_0\bigr\}$}
in Lemma {\rm\ref{lem1418}} as partial structures
of energy cut level $E_{0}$.
\item[$(2)$]
The restriction of \smash{$\bigl\{\mathfrak m_{k,E}^{[0,1]_s}\bigr\}$}
to $s=0$ and $s=1$ coincide with
$\bigl\{\mathfrak m_{k,E}^{\Xi_{1}}\bigr\}$ and
$\bigl\{\mathfrak m_{k,E}^{\Xi_{2}}\bigr\}$, respectively.
\end{enumerate}

\end{cor}
Corollary \ref{cor1421} implies that
$\bigl\{\mathfrak m_{k,E}^{\Xi_{1}}\bigr\}$ is pseudo-isotopic to
$\bigl\{\mathfrak m_{k,E}^{\Xi_{2}}\bigr\}$.
It particular they are strongly homotopy equivalent.
The proof of Theorem~\ref{thm1444}\,(1) is complete.
\end{proof}

\begin{proof}[Proof of Theorem~\ref{thm1444}\,(2)]
The proof of Theorem~\ref{thm1444}\,(2)
is similar to the proof of (1)
but we need to iterate once more the process to take
higher homotopy as the following.

During the construction of pseudo-isotopy
in Corollary \ref{cor1421} we made
various choices.
Especially we made a choice of CF-perturbations
$\widehat{\mathfrak S}([0,1]_s \times [0,1]_t;E_i,E_{i+1})$
in Proposition~\ref{prop1419}.
We will prove the
homotopy equivalence we obtained
in the proof of Theorem~\ref{thm1444}\,(1)
is independent of such choices up to homotopy.

Suppose
$\widehat{\mathfrak S}([0,1]_s \times [0,1]_t;j;E_i,E_{i+1})$,
$j=1,2$, are two choices.
We denote by
$\bigl\{\mathfrak m_{k,E}^{[0,1]_s,j=1}\bigr\}$,
\smash{$\bigl\{\mathfrak m_{k,E}^{[0,1]_s,j=2}\bigr\}$}
the pseudo-isotopies obtained by these two choices,
respectively.
\begin{lem}\label{lemma14222}
There exists a system of CF-perturbations, which we denote by $\widehat{\mathfrak S}([0,1]_s \times [0,1]_t\times [0,1]_u;E_i,E_{i+1})$, on outer collarings of thickenings of
${\mathcal M}_{k+1}(L_1;E;[0,1]_s)
\times [0,1]_t\times [0,1]_u$ for~${E \le E_{i+1}}$
with the following properties.
\begin{enumerate}\itemsep=0pt
\item[$(1)$] The CF-perturbation
$\widehat{\mathfrak S}([0,1]_s \times [0,1]_t\times [0,1]_u;E_i,E_{i+1})$ is transversal to $0$.
\item[$(2)$]
The map
\begin{gather*}
({\rm ev}_0,{\rm ev}_{[0,1]_s},{\rm ev}_{[0,1]_t},{\rm ev}_{[0,1]_u})
\colon
\\
\qquad{\mathcal M}_{k+1}(L_1;E;[0,1]_s)
\times [0,1]_t\times [0,1]_u \to R \times [0,1]_s \times [0,1]_t\times [0,1]_u
\end{gather*}
is
strongly submersive with respect to
$\widehat{\mathfrak S}([0,1]_s \times [0,1]_t\times [0,1]_u;E_i,E_{i+1})$.
\item[$(3)$]
We consider the restriction of
$\widehat{\mathfrak S}([0,1]_s \times [0,1]_t\times [0,1]_u;E_i,E_{i+1})$
to the boundary components, which are products of $[0,1]_t
\times [0,1]_u$ and the
boundary components of the space ${\mathcal M}_{k+1}(L_1;E;[0,1]_s)$ in Proposition {\rm\ref{prop14555}\,(2),~(I)}.
It then coincides with the fiber product CF-perturbation,
which is well-defined by item $(2)$.
\item[$(4)$]
For sufficiently small $\tau$, the following holds.\!
The restriction of the CF-perturbation
$\smash{\widehat{\mathfrak S}([0,1]_s}\allowbreak \times [0,1]_t\times [0,1]_u;E_i, E_{i+1})$ to
\smash{${\rm ev}_{[0,1]_u}^{-1}([0,\tau])$}
coincides with the pullback of $\widehat{\mathfrak S}([0,1]_s \times [0,1]_t;1;\allowbreak E_i, E_{i+1})$. The restriction of
$\widehat{\mathfrak S}([0,1]_s \times [0,1]_t\times [0,1]_u;E_i,E_{i+1})$ to
\smash{${\rm ev}_{[0,1]_u}^{-1}([1-\tau,1])$}
coincides with the pullback of $\widehat{\mathfrak S}([0,1]_s \times [0,1]_t,2;E_i,E_{i+1})$.
\item[$(5)$]
For sufficiently small $\tau$, the following holds.
The restriction of $\widehat{\mathfrak S}([0,1]_s \times [0,1]_t\times [0,1]_u;E_i,\allowbreak E_{i+1})$ to
\smash{${\rm ev}_{[0,1]_s}^{-1}([0,\tau])$}
coincides with the pullback of $\widehat{\mathfrak S}([0,1]_t,\Xi_{1};E_i,E_{i+1})$.
$($In particular, this restriction is constant in $u$ direction.$)$
The restriction of
$\widehat{\mathfrak S}([0,1]_s \times [0,1]_t\times [0,1]_u;E_i,E_{i+1})$ to
\smash{${\rm ev}_{[0,1]_s}^{-1}([1-\tau,1])$}
coincides with the pullback of $\widehat{\mathfrak S}([0,1]_t,\Xi_{2};E_i,E_{i+1})$.
\item[$(6)$]
For sufficiently small $\tau$, the following holds.
The restriction of
$\widehat{\mathfrak S}([0,1]_s \times [0,1]_t\times [0,1]_u;\allowbreak E_i,E_{i+1})$ to
\smash{${\rm ev}_{[0,1]_t}^{-1}([0,\tau])$}
coincides with the pullback of $\widehat{\mathfrak S}([0,1]_s;\allowbreak E_{i})$.
The~restriction of
$\widehat{\mathfrak S}([0,1]_s \times [0,1]_t\times [0,1]_u;E_i,E_{i+1})$ to
\smash{${\rm ev}_{[0,1]_s}^{-1}([1-\tau,1])$}
coincides with the pullback of $\widehat{\mathfrak S}([0,1]_s;E_{i+1})$.
\end{enumerate}
\end{lem}

See Figure~\ref{FigureSec1476}.
The proof of Lemma~\ref{lemma14222} is the same as other similar
results such as Proposition~\ref{prop1419}.
\begin{figure}[ht]
\centering
\includegraphics[scale=0.44]{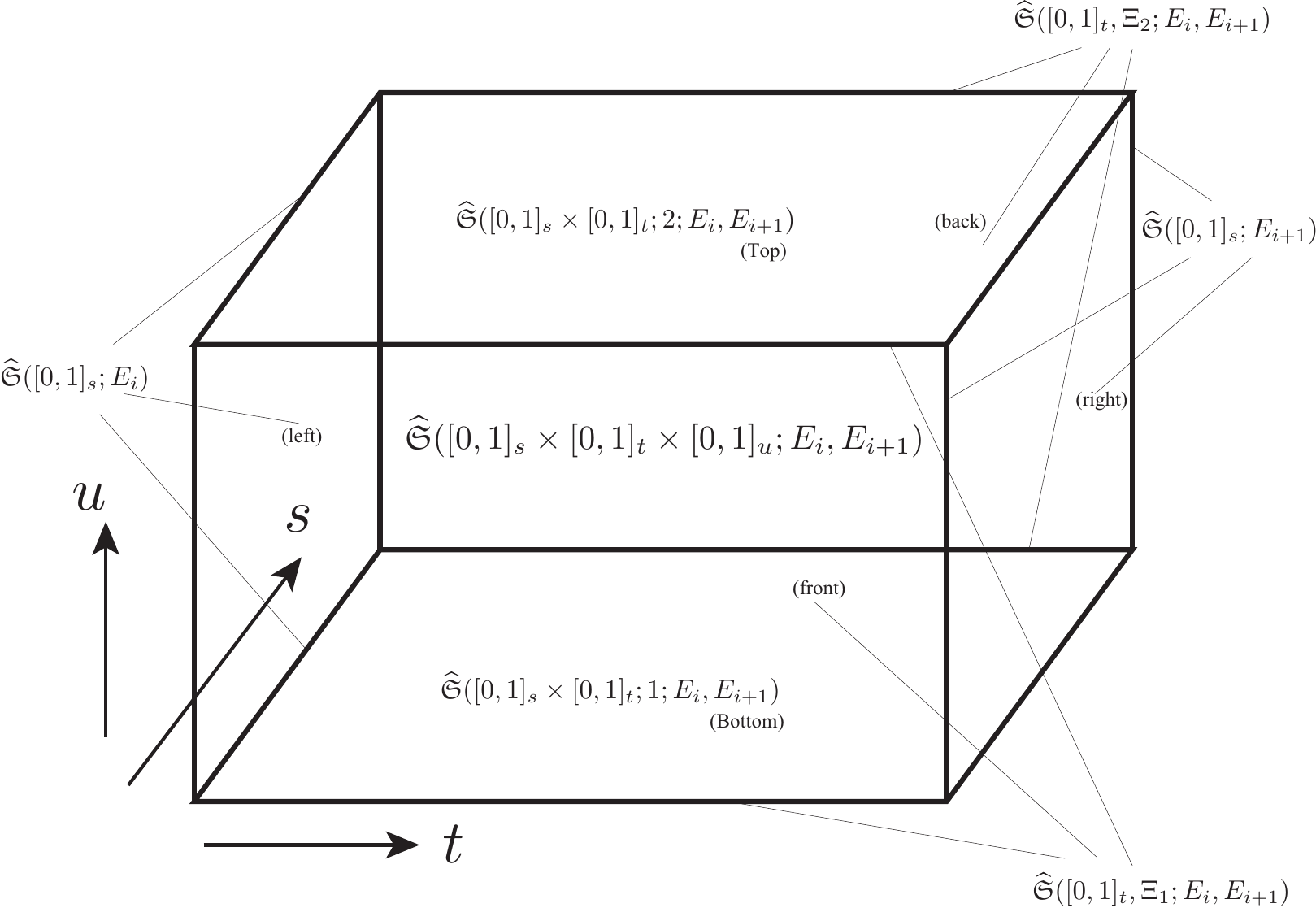}
\caption{pseudo-isotopy of pseudo-isotopies of pseudo-isotopies.}
\label{FigureSec1476}
\end{figure}

Now we discuss in the same way as
Lemma~\ref{lem1420} and Corollary \ref{cor1421}
using Proposition~\ref{prop2013} and obtain:
\begin{lem}\label{lem1423}
There exists a $P=[0,1]_s \times [0,1]_u$
parametrized family of filtered $A_{\infty}$
structures~\smash{$\bigl\{\mathfrak m_{k,E}^{[0,1]_s\times [0,1]_u}\bigr\}$}
on $CF([0,1]_s\times [0,1]_u \times L;\Theta;\R)$
with the following properties:
\begin{enumerate}\itemsep=0pt
\item[$(1)$]
The restriction of the structure \smash{$\bigl\{\mathfrak m_{k,E}^{[0,1]_s\times [0,1]_u}\bigr\}$}
to $u=0$ $($resp.\ $u=1)$ coincides with the
pseudo-homotopy of Corollary {\rm\ref{cor1421}}
obtained by using
\smash{$\widehat{\mathfrak S}([0,1]_s \times [0,1]_t;1;E_i,E_{i+1})$}
$($resp.\
$\widehat{\mathfrak S}([0,1]_s \times [0,1]_t;2;E_i,E_{i+1}))$.
\item[$(2)$]
The restriction of $\bigl\{\mathfrak m_{k,E}^{[0,1]_s}\bigr\}$
to $s=0$ and $s=1$ coincide with the pullback of
$\bigl\{\mathfrak m_{k,E}^{\Xi_{1}}\bigr\}$ and
\smash{$\bigl\{\mathfrak m_{k,E}^{\Xi_{2}}\bigr\}$}, respectively.
In particular, they are trivial in $[0,1]_u$ factor.
\end{enumerate}
\end{lem}

In other words, we have the following commutative diagram:
\[
\begin{CD}
\bigl(CF(L),\bigl\{\mathfrak m_{k,E}^{\Xi_{2}}\bigr\}\bigr)
@ <{{\rm Eval}_{u=0}}<< CF([0,1]_u\times L)
@>{{\rm Eval}_{u=1}}>> \bigl(CF(L),\bigl\{\mathfrak m_{k,E}^{\Xi_{2}}\bigr\}\bigr)
\\
@ A{{\rm Eval}_{s=1}}AA @ AA{{\rm Eval}_{s=1}}A
@ AA{{\rm Eval}_{s=1}}A\\
\displaystyle\bigl(CF([0,1]_s\times L),
\atop
\displaystyle\!\!\bigl\{\mathfrak m_{k,E}^{[0,1]_s,j=1}\bigr\}\bigr)
@ <{{\rm Eval}_{u=0}}<<
\displaystyle
\bigl(CF([0,1]_s\times[0,1]_u\times L),
\atop\displaystyle\!\!\!\!\!\!\!\!\!\!\!\!\!\!\!\!\!\!\!\!\!\!\!\!\!\!\!\bigl\{\mathfrak m_{k,E}^{[0,1]_s,j=1}\bigr\}\bigr)
@>{{\rm Eval}_{u=1}}>>
\displaystyle\bigl(CF([0,1]_s\times L),
\atop\displaystyle\!\!\!\!\!\!\bigl\{\mathfrak m_{k,E}^{[0,1]_s,j=2}\bigr\}\bigr)
\\
@ V{{\rm Eval}_{s=0}}VV @ VV{{\rm Eval}_{s=0}}V
@ VV{{\rm Eval}_{s=0}}V\\
\bigl(CF(L),\bigl\{\mathfrak m_{k,E}^{\Xi_{1}}\bigr\}\bigr)
@ <{{\rm Eval}_{u=0}}<< CF([0,1]_t\times L)
@>{{\rm Eval}_{u=1}}>> \bigl(CF(L),\bigl\{\mathfrak m_{k,E}^{\Xi_{1}}\bigr\}\bigr).
\end{CD}
\]
All the arrows in the diagram are strong homotopy equivalences.
By Lemma~\ref{lem1423}\,(2), we find that
${\rm Eval}_{u=1}$ in the first horizontal line
is homotopic to ${\rm Eval}_{u=0}$
in the first horizontal line.
The same holds for the third horizontal line.
The composition
\[
{\rm Eval}_{s=1} \circ ({\rm Eval}_{s=0})^{-1}
\colon\ \bigl(CF(L),\bigl\{\mathfrak m_{k,E}^{\Xi_{1}}\bigr\}\bigr)
\to \bigl(CF(L),\bigl\{\mathfrak m_{k,E}^{\Xi_{2}}\bigr\}\bigr)
\]
of maps in the first vertical line
is the strong homotopy equivalence
obtained from the choice \smash{$\widehat{\mathfrak S}([0,1]_s \times [0,1]_t;1;E_i,E_{i+1})$}.
In the same way,
the third vertical line gives the strong homotopy equivalence
obtained from the choice $\widehat{\mathfrak S}([0,1]_s \times [0,1]_t;2;E_i,E_{i+1})$.
Thus those two homotopy equivalences are homotopic each
other. The proof of Theorem~\ref{thm1444} is complete.
\end{proof}

\begin{rem}
The above diagram is similar to \cite[Figure 10.1]{AJ},
which is used for a similar purpose.
\end{rem}

\subsection[Proof of Theorem~\ref{thm1466}]{Proof of Theorem~\ref{thm1466}}
\label{sec:proofindepen}

\subsubsection{Pseudo-isotopy of tri-modules}
\label{subsub:Pseudoisotopytri}

\begin{situ}
Let $R_m$, $m=1,2,3$, be a compact smooth manifold
without boundary and~$\Theta_m$ a principal ${\rm O}(1)$
bundles on it. Let $R$ be a compact smooth manifold
without boundary and~$\Theta$ a principal ${\rm O}(1)$
bundle on it.

Let $P$ be a manifold with corners and
$C(P \times R_m;\R) = C^{\infty}(\Omega(P \times R_m))$
(resp.\ $
C(P \times R;\R) = C^{\infty}(\Omega(P \times R))
$)
the set of smooth forms on $P \times R_m$
(resp.\ $R$)
twisted by $\Theta_m$ (resp.\ $\Theta$).

We define $C(P \times R_m;\Lambda_0)$
(resp.\ $C(P \times R;\Lambda_0)$)
as a completion of
the tensor product $C(P \times R_m;\R) \otimes \Lambda_0$
(resp.\ $C(P \times R;\R) \otimes \Lambda_0$).

Suppose that, for each $m$, we are given a
$P$-parametrized family of $G$-gapped filtered $A_{\infty}$
structures on $C(P \times R;\Lambda_0)$,
which we denote by $\bigl\{\mathfrak m^{P,m}_{k,\beta}\bigr\}$.
We put $\mathscr C_m^P = (C(P \times R;\Lambda_0),\bigl\{\mathfrak m^{P,m}_{k,\beta}\bigr\})$.

\end{situ}
\begin{defn}\label{PparaAinfmokddef}
A {\it $P$-parametrized family of $G$-gapped filtered $A_{\infty}$
tri-module structures on $CF(P \times R;\Theta;\Lambda_0)$} over $C\bigl(P \times R_m;\Lambda_0,\bigl\{\mathfrak m^{P,m}_{k,\beta}\bigr\}\bigr)$,
$m=1,2,3$,
is $\bigl\{\mathfrak n^P_{k_1,k_2,k_3;\beta}\bigr\}$
for $\beta \in G$ and $k_i = 0,1,2,\dots$,
that satisfies the following:\index{parametrized family of tri-module structures}
\begin{enumerate}\itemsep=0pt
\item[(1)]
\[
\mathfrak n^P_{k_1,k_2,k_3;\beta} \colon
\bigotimes_{i=1}^3B_{k_i}(\Omega(P\times R_i)[1])
\otimes \Omega(P\times R)[1] \to \Omega(P\times R)[1]
\]
is a multilinear map of degree $1$.
\item[(2)]
The maps $\mathfrak n^P_{k_1,k_2,k_3;\beta}$ is
pointwise in $P$ direction if $\beta \ne \beta_0$
or $k_1+k_2+k_3 \ge 1$.
\item[(3)]
$\mathfrak n^P_{0,0,0;\beta_0}(h) = (-1)^*dh$. Here $d$ is the de Rham differential
and $*$ is as in \eqref{form3420000}.
\item[(4)]
The operations $\bigl\{\mathfrak n^P_{k_1,k_2,k_3;\beta}\bigr\}$
define a filtered $A_{\infty}$ tri-module
over $\mathscr C_m(P)$
($m=1,2,3$).
\end{enumerate}
In the case when $P = [0,1]$ we call $CF([0,1] \times R;\Theta;\Lambda_0)$
together with its $P$-parametrized family of $G$-gapped filtered $A_{\infty}$
tri-module structures, a {\it pseudo-isotopy of
$G$-gapped filtered $A_{\infty}$
tri-modules} over the pseudo-isotopies
\smash{$\mathscr C_m^{[0,1]}$}, $m=1,2,3$,
of filtered $A_{\infty}$ categories.
\end{defn}

\subsubsection{Existence of a pseudo-isotopy of tri-modules}
\label{subsub:exiPseudoisotopytri}

We go back to our geometric situation of Theorem~\ref{thm1466}.
We consider the case when the sets $\mathbb L_1$, $\mathbb L_2$ and
$\mathbb L_{12}$ consist of single immersed Lagrangian
submanifolds.

We put
$
R_1 = \tilde L_1 \times_{X_1} \tilde L_1$,
$
R_2 = \tilde L_2 \times_{X_2} \tilde L_2$,
$
R_3 = \tilde L_{12} \times_{X_1 \times X_2} \tilde L_{12}$.
The pseudo-isotopy $\mathscr C_m^{[0,1]}$ are given by
Corollary \ref{cor1421}.
In particular, we make Choices
$\Xi_{1,j}$, $\Xi_{2,j}$, $\Xi_{12,j}$ in Situation \ref{situ1415}.
They give filtered $A_{\infty}$ structures of $\mathscr C_m^{s=0}$
and of $\mathscr C_m^{s=1}$.
We also take
\begin{equation}\label{choiceR}
R = \tilde L_1 \times_{X_1} \tilde L_{12} \times_{X_2} \tilde L_2.
\end{equation}
We make a choice of \smash{$\Xi_{12,j}^{\rm quilt}$} in Situation \ref{situ1415}.
It determines a filtered $A_{\infty}$ tri-module structure
on $CF(R;\Theta;\Lambda_0)$ over $\mathscr C_m^{s=0}$
or $\mathscr C_m^{s=1}$ for $j=1,2$.
\begin{prop}\label{prop1430}
There exists a pseudo-isotopy of filtered $A_{\infty}$
tri-module on
$CF([0,1]\times R;\Theta;\Lambda_0)$ over
\smash{$\mathscr C_m^{[0,1]}$} for $m=1,2,3$.
We may choose it so that the restriction to $s=0,1$
coincides with the tri-module structure induced by Choices
$\Xi_{12,j}^{\rm quilt}$ for $j=1,2$.
\end{prop}
\begin{proof}
The proof of this proposition is mostly the same as the proof of
Lemma~\ref{lem1420} in Section~\ref{sec:welldefAinfcat}.
We first define the notion of a partial pseudo-isotopy of
tri-module structures with energy cut level $E$.
We then can show the existence of
a partial pseudo-isotopy of
tri-module structures with energy cut level $E$
for any $E$. Then we proceed in the same way to
define the notion of
a partial pseudo-isotopy of pseudo-isotopies of
tri-module structures and use it
to work out the homotopy inductive limit construction.
The way to modify the proof of Lemma~\ref{lem1420}
is thus a routine, which we omit.
\end{proof}

\begin{situ}\label{situ1431}
Let $\mathbb L_1$, $\mathbb L_{12}$, $\mathbb L_2$ be as in
Situation \ref{situ61}.
We consider the (disjoint) union of all the elements of
$\mathbb L_1$ (resp.\ $\mathbb L_{12}$, $\mathbb L_2$)
and denote them by $L_1$, $L_{12}$, $L_2$.
We consider $R$ as in \eqref{choiceR}.
\end{situ}
We remark that since we are in Situation \ref{situ61}
the fiber product
$
\tilde L_2^0 = \tilde L_1 \times_{X_1} \tilde L_{12}
$
is an open subset of $L_2$.
We put
$
R_0 = \tilde L_1 \times_{X_1} \tilde L_{12} \times_{X_2} \tilde L_2^0
\subseteq R$.
We remark also
$
R_0 = \tilde L_2^0 \times_{X_2} \tilde L_2^0 \subseteq R_2$.
We remark that the tri-module structure
we used in Theorem~\ref{thm61} satisfies the following properties.
If $h_2,h \in C^{\infty}(R_0,\Omega(R)\otimes \Theta)$,
then
$
\mathfrak n_{0,0,1;\beta_0}(1,1,h_2;h) = (-1)^{\deg h_2} h_2 \wedge h$.
This fact is used during the proof of Proposition~\ref{prop610}.
\begin{lem}\label{lem1431}
We can take the pseudo-isotopy in Proposition {\rm\ref{prop1430}}
such that the following holds in addition.
If $h_2,h \in C^{\infty}([0,1] \times R_0,\Omega(R)\otimes \Theta)$,
then
\[
\mathfrak n^{[0,1]}_{0,0,1;\beta_0}(1,1,h_2;h) = (-1)^{\deg h_2} h_2 \wedge h.
\]

\end{lem}
Using the fact that the moduli space defining $\mathfrak n^{[0,1]}_{0,0,1;\beta_0}$
on $[0,1] \times R_0$ consists of constant maps, and has the
required transversality and submersivity properties
without perturbation, the proof of the lemma is similar
to an argument during the proof of Proposition~\ref{prop610} and so is omitted.

\subsubsection[Completion of the proof of Theorem~\ref{thm1466}]{Completion of the proof of Theorem~\ref{thm1466}}
\label{subsub:laststepofuniquess}

Now we are in the position to complete the proof of Theorem~\ref{thm1466}.

Suppose we are in Situation \ref{situ1431}.
We use the same trick as Section~\ref{subsec:Ainfcatim} to obtain a
filtered~$A_{\infty}$
category from a filtered $A_{\infty}$ algebra
\smash{$\mathscr C_m^{[0,1]}$}, $m=1,2,3$. Here \smash{$\mathscr C_m^{[0,1]}$}
is obtained in Proposition~\ref{prop1430}.
We denote them by
$\mathfrak{Fukst}(X_1;\mathbb L_1)^{[0,1]}$,
$\mathfrak{Fukst}(-X_1\times X_2;\mathbb L_{12})^{[0,1]}$
and
$\mathfrak{Fukst}(X_2;\mathbb L_{2})^{[0,1]}$.
The sets of their objects are the same as the sets of objects of
$\mathfrak{Fukst}(X_1;\mathbb L_1)$,
$\mathfrak{Fukst}(-X_1\times X_2;\mathbb L_{12})$
and
$\mathfrak{Fukst}(X_2;\mathbb L_{2})$,
respectively.

Hereafter, we omit $\mathbb L_1$, $\mathbb L_{12}$, $\mathbb L_2$
from the notation for simplicity.
The pseudo-isotopy of tri-modules we produced in Proposition~\ref{prop1430}
induces a tri-module structure
over the strict categories $\mathfrak{Fukst}(X_1)^{[0,1]}$,
$\mathfrak{Fukst}(-X_1 \times X_2)^{[0,1]}$ and
$\mathcal{FUNC}\bigl(\mathfrak{Fukst}(X_2)^{[0,1]}\bigr)$.
It induces a strict filtered $A_{\infty}$ bi-functor
\begin{equation}\label{form141313}
\mathfrak{Fukst}(X_1)^{[0,1]} \times \mathfrak{Fukst}(-X_1 \times X_2)^{[0,1]}
\to \mathcal{FUNC}\bigl(\bigl(\mathfrak{Fukst}(X_2)^{[0,1]}\bigr)^{\rm op},\mathcal{CH}\bigr).
\end{equation}
We denote by
$\mathcal{REP}\bigl(\mathfrak{Fukst}(X_2)^{[0,1]}\bigr)$
the full subcategory of the filtered $A_{\infty}$ category
\[
\mathcal{FUNC}\bigl(\bigl(\mathfrak{Fukst}(X_2)^{[0,1]}\bigr)^{\rm op},\mathcal{CH}\bigr)
\]
whose object is homotopy equivalent to the image of
the Yoneda-functor
\[
\mathfrak{Fukst}(X_2)^{[0,1]} \to
\mathcal{FUNC}\bigl(\bigl(\mathfrak{Fukst}(X_2)^{[0,1]}\bigr)^{\rm op},\mathcal{CH}).
\]
We define $\mathcal{REP}(\mathfrak{Fukst}(X_2))$ in the same way.

Lemma~\ref{lem1431} implies that the image of the
functor \eqref{form141313} lies in $\mathcal{REP}\bigl(\mathfrak{Fukst}(X_2)^{[0,1]}\bigr)$.
Thus we obtain the next diagram, which commutes up to
homotopy equivalence:
\[
\begin{CD}
\displaystyle\!\!\!\!\!\!\!\!\!\!\!\!\mathfrak{Fukst}(X_1;\Xi_{1,2})
\atop \displaystyle\times \mathfrak{Fukst}(-X_1 \times X_2;\Xi_{12,2})
@ >{\mathscr F}>> \mathcal{REP}(\mathfrak{Fukst}(X_2;\Xi_{2,2}))
@<<< \mathfrak{Fukst}(X_2;\Xi_{2,2})
\\
@ A{{\rm Eval}_{s=1}}AA @ VV{{\rm Eval}^*_{s=1}}V
@ AA{{\rm Eval}_{s=1}}A\\
\displaystyle\!\!\!\!\!\!\!\!\!\!\!\!\mathfrak{Fukst}(X_1)^{[0,1]}
\atop \displaystyle\times \mathfrak{Fukst}(-X_1 \times X_2)^{[0,1]}
@ >{\mathscr F}>>
\mathcal{REP}\bigl(\mathfrak{Fukst}(X_2)^{[0,1]}\bigr)
@<<<
\mathfrak{Fukst}(X_2)^{[0,1]}
\\
@ V{{\rm Eval}_{s=0}}VV @ AA{{\rm Eval}^*_{s=0}}A
@ VV{{\rm Eval}_{s=0}}V\\
\displaystyle\!\!\!\!\!\!\!\!\!\!\!\!\mathfrak{Fukst}(X_1;\Xi_{1,1})
\atop \displaystyle\times \mathfrak{Fukst}(-X_1 \times X_2;\Xi_{12,1})
@ >{\mathscr F}>> \mathcal{REP}(\mathfrak{Fukst}(X_2;\Xi_{2,1}))
@<<< \mathfrak{Fukst}(X_2;\Xi_{2,1}).
\end{CD}
\]
All arrows in the diagram are homotopy equivalences except the three horizontal
arrows in the left-hand side, which are written as $\mathscr F$.
By definition, the composition of the arrows of the first~line
is the filtered $A_{\infty}$ functor \smash{$\mathcal{MWW}^{\Xi^{\rm quilt}_{12,2}}$}.
The composition of the arrows of the third line
is the filtered $A_{\infty}$ functor
\smash{$\mathcal{MWW}^{\Xi^{\rm quilt}_{12,1}}$}.

The composition of two arrows in the first
column is the functor $\mathscr G^1 \times \mathscr G^{12}$.
The composition of the two arrows in the
third column is the functor $\mathscr G^2$.
Thus Theorem~\ref{thm1466} follows from the
commutativity of the diagram.

Note that we can prove the next theorem in the same way.
\begin{thm}\label{them1431}
The composition functor $\mathfrak{Comp}$ in Theorem {\rm\ref{comp2}}
is independent of the choices up to homotopy
equivalence.
\end{thm}

We omit the proof.

\subsection[Coincidence of $A_\infty$ structures defined by the
two compactifications]{Coincidence of $\boldsymbol{A_{\infty}}$ structures defined by the
two compactifications}
\label{sec:comptwocomp}

Let $\mathbb L_{12}$ be a clean collection of
$\pi_1^*(TX_1 \oplus V_1) \oplus \pi_2^*(V_2)$
relatively spin immersed Lagrangian submanifolds
of $-X_1 \times X_2$.

In Section~\ref{sec:HFIm}, we used the stable map compactification
${\mathcal M}(L_{12};\vec a;E)$ of the moduli space of pseudo-holomorphic disks
to define a filtered $A_{\infty}$ category
the set of whose objects is $\mathbb L_{12}$.
We denote it by $\mathfrak{Fuk}(-X_1\times X_2,\mathbb L_{12})$.
In Section~\ref{sec:directcomp}, we
introduced a different compactification~${\mathcal M}'(L_{12};\vec a;E)$.
We use it also to define a filtered $A_{\infty}$ category
the set of whose objects is~$\mathbb L_{12}$.
We denote it by $\mathfrak{Fuk}'(-X_1\times X_2,\mathbb L_{12})$.
In this subsection, we prove the following.

\begin{prop}\label{prop143434}
$\mathfrak{Fuk}(-X_1\times X_2,\mathbb L_{12})$
is pseudo-isotopic to $\mathfrak{Fuk}'(-X_1\times X_2,\mathbb L_{12})$.
\end{prop}
\begin{proof}
By the same trick as Section~\ref{subsec:Ainfcatim}, it suffices
to consider the case when $\mathbb L_{12}$ consists of a single
immersed Lagrangian submanifold $L_{12}$ and construct a pseudo-isotopy of filtered
$A_{\infty}$ algebras.
\begin{lemdef}\label{lemdef14355}
We can define the forgetful map\index[syindex]{fg@$\mathfrak{fg}$}
\[
\mathfrak{fg} \colon\ {\mathcal M}_{\ell}(L_{12};\vec a;E)
\to
{\mathcal M}'_{\ell,0,0}(L_{12};\vec a;E),
\]
which is continuous.
\end{lemdef}
\begin{proof}
Let $\bigl(\bigl(\Sigma,\vec z,\vec z^{\,\rm int}\bigr),u,\gamma\bigr)$ be an element of
${\mathcal M}_{\ell}(L_{12};\vec a;E,\gamma)$.
Here $\bigl(\Sigma,\vec z,\vec z^{\,\rm int}\bigr)$ is a bordered nodal marked curve of
genus zero with one boundary component.
($\vec z$ are boundary marked points and $\vec z^{\,\rm int}$ are
interior marked points.)
The map
$u \colon (\Sigma,\partial \Sigma) \to (-X_1 \times X_2, L_{12})$
is pseudo-holomorphic and the map
$\gamma \colon \partial \Sigma \setminus \vec z \to \tilde L_{12}$
is a lift of the restriction of $u$.

We put $(u_1,u_2):=u$, where $u_i$ is a map to $X_i$
from $\Sigma$.
We consider $\bigl(\bigl(\Sigma,\vec z,\vec z^{\,\rm int}\bigr),u_i\bigr)$ for $i=1,2$
and shrink unstable sphere components.
Here an unstable sphere component of $\bigl(\bigl(\Sigma,\vec z,\vec z^{\,\rm int}\bigr),u_i\bigr)$
is an unstable sphere component of the source curve
$\bigl(\Sigma,\vec z,\vec z^{\,\rm int}\bigr)$ on which $u_i$ is constant.
We denote by
$\bigl(\bigl(\Sigma_i,\vec z_i,\vec z^{\,\rm int}_i\bigr),u_i\bigr)$
the pair of a bordered marked curve and a map
obtained by this shrinking.

Let $\bigl(\Sigma^0_i,\vec z_i,\vec z^{\,\rm int}_i\bigr)$ be
the bordered marked curve obtained from
$\bigl(\Sigma_i,\vec z_i,\vec z^{\,\rm int}_i\bigr)$ by
shrinking all the unstable sphere components.

We remark that
$\bigl(\Sigma^0_1,\vec z_1,\vec z^{\,\rm int}_1\bigr)$
is canonically isomorphic to
$\bigl(\Sigma^0_2,\vec z_2,\vec z^{\,\rm int}_2\bigr)$.
In fact, they both are obtained by
shrinking all the unstable sphere components of
$\bigl(\Sigma,\vec z,\vec z^{\,\rm int}\bigr)$.
Therefore, we obtain a biholomorphic map
$
\mathscr I \colon \bigl(\Sigma^0_1,\vec z_1,\vec z^{\,\rm int}_1\bigr)
\to
\bigl(\Sigma^0_2,\vec z_2,\vec z^{\,\rm int}_2\bigr)$.
We define
\[
\mathfrak{fg}\bigl(\bigl(\Sigma,\vec z,\vec z^{\,\rm int}\bigr),u,\gamma\bigr)
=
\bigl(\bigl(\bigl(\Sigma_1,\vec z_1,\vec z^{\,\rm int}_1\bigr),u_1\bigr),\bigl(\bigl(\Sigma_2,\vec z_2,\vec z^{\,\rm int}_2),\mathscr I,\gamma\bigr)\bigr).
\]
Note that we regard the interior marked points $\vec z^{\,\rm int}$ as interior
marked points of first kind in the sense of Definition~\ref{defn145555}.

The continuity of the map $\mathfrak{fg}$ is easy to show from the definition.
\end{proof}

We consider the case when $\ell = 0$ to obtain a map
$
\mathfrak{fg} \colon {\mathcal M}(L_{12};\vec a;E)
\to
{\mathcal M}'(L_{12};\vec a;E)$.
We start with a Kuranishi structure which we defined on
${\mathcal M}'(L_{12};\vec a;E)$ and will
pull it back to one on
${\mathcal M}(L_{12};\vec a;E)$.
We describe the detail of this pullback construction now.

Let $\tilde\xi \in {\mathcal M}(L_{12};\vec a;E)$.
We take
$\hat\xi = \bigl(\bigl(\Sigma,\vec z,\vec z^{\,\rm int}\bigr),u,\gamma\bigr) \in {\mathcal M}_{\ell}(L_{12};\vec a;E)$ such that $\tilde\xi = [(\Sigma,\vec z),u,\gamma]$ and
$\bigl(\Sigma,\vec z,\vec z^{\,\rm int}\bigr)$ is stable.
We use it to define a notion that
$((\Sigma',\vec z^{\,\prime}),u',\gamma')$
is $\varepsilon$-close to \smash{$\hat\xi$} in a similar way as
Definition~\ref{defn1236}
as follows.
\begin{defn}\label{defn14360}
Let $\bigl(\bigl(\Sigma^{\heartsuit},\vec z^{ \heartsuit}\bigr),u^{\heartsuit},\gamma^{\heartsuit}\bigr)$
be an object which has the same properties as an element of
${\mathcal M}(L_{12};\vec a;E)$
except we do not require $u^{\heartsuit}$ to be pseudo-holomorphic.
We call such an object a \index{candidate of an element of extended moduli space}
{\it candidate of an element of the extended moduli space}.
\end{defn}
\begin{defn}\label{defn1436}
We say
$\bigl(\bigl(\Sigma^{\heartsuit},\vec z^{ \heartsuit}\bigr),u^{\heartsuit},\gamma^{\heartsuit}\bigr)$
is {\it $\varepsilon$-close to $\bigl(\tilde\xi,\hat\xi\bigr)$} \index{$\varepsilon$-close} if there exists
$\vec z^{ \,\rm int,\heartsuit}$ with the following properties:
\begin{enumerate}\itemsep=0pt
\item[(1)]
$\bigl(\Sigma^{\heartsuit},\vec z^{ \heartsuit},\vec z^{ \rm int,\heartsuit}\bigr)$
is $\varepsilon$-close to $\bigl(\Sigma,\vec z,\vec z^{\,\rm int}\bigr)$ in the
moduli space of marked stable disks.\footnote{We take and
fix a metric of the moduli space of marked stable disks to
define this $\varepsilon$-closeness.}
\item[(2)]
We define the {\it core} $K^{\rm s}_{a}$ and $K^{\rm d}_{a}$ in the same way as
\eqref{formdefcore}. Here $a$ is an index of the irreducible component
of $\Sigma$. $K^{\rm s}_{a}$ lies in a sphere component and
$K^{\rm d}_{a}$ lies in a disk component.
Then we obtain
smooth embeddings
$
\mathcal I^{\rm d}_{\heartsuit} \colon\ K^{\rm d}_{a}
\to \Sigma^{\heartsuit}$, $
\mathcal I^{\rm s}_{\heartsuit} \colon\ K^{\rm s}_{\rm a}
\to \Sigma^{\heartsuit}$,
in the same way as~\eqref{form1237} and Definition~\ref{defn1227}.
(We use analytic family of coordinates at the nodal points of
$\Sigma$ and also a trivialization of the universal family of
marked stable disks on the $\varepsilon$ neighborhood of
$\bigl(\Sigma,\vec z,\vec z^{\,\rm int}\bigr)$ to define them.
See \cite[Section 18]{foootech}, \cite[Section 3]{foootech2} etc.)

We now require
\begin{enumerate}\itemsep=0pt
\item
The restriction of $u$ to each $K^{\rm d}_{\rm a}$,
is $\varepsilon$ close to $u^{\heartsuit}\circ \mathcal I^{\rm d}_{\heartsuit}$ in $C^2$ norm.
\item
The restriction of $u$ to
each $K^{\rm s}_{\rm a}$ is $\varepsilon$ close to
$u^{\heartsuit}\circ \mathcal I^{\rm s}_{\heartsuit}$ in $C^2$ norm.
\end{enumerate}
We remark that these conditions are similar to Definition~\ref{defn1236}\,(3),~(4), respectively.
\item[(3)]
For any connected component $\mathcal S$ of
\[
\Sigma^{\heartsuit}
\setminus
\bigcup_{a} \mathcal I^{\rm d}_{\heartsuit}\bigl(K^{\rm d}_{\rm a}\bigr)
\setminus
\bigcup_{a} \mathcal I^{\rm s}_{\heartsuit}(K^{\rm s}_{\rm a}),
\]
we require
$
\operatorname{Diam} u^{\heartsuit}(\mathcal S) < \varepsilon$.
(In other words, we require the diameters of
the images by~$u^{\heartsuit}_i$ of the neck regions
are smaller than $\varepsilon$.)
We remark that these conditions are similar to Definition~\ref{defn1236}\,(6).
\end{enumerate}

\end{defn}
\begin{lem}\label{lema1437}
Let $\hat\xi^{\prime} = \bigl(\bigl(\Sigma,\vec z,\vec z^{ \rm int,\prime}\bigr),u,\gamma\bigr) \in {\mathcal M}_{\ell'}(L_{12};\vec a;E)$ such that $\tilde\xi = [(\Sigma,\vec z),u,\gamma]$ and
$(\Sigma,\vec z,\vec z^{\, \rm int,\prime})$ is stable.
Then for each $\varepsilon$ there exists $\delta$ with the following properties.
If $\bigl(\bigl(\Sigma^{\heartsuit},\vec z^{ \heartsuit}\bigr),\allowbreak u^{\heartsuit},\gamma^{\heartsuit}\bigr)$
is $\delta$-close to $(\hat\xi^{\prime},\tilde\xi)$,
then $\bigl(\bigl(\Sigma^{\heartsuit},\vec z^{ \heartsuit}\bigr),u^{\heartsuit},\gamma^{\heartsuit}\bigr)$ is
$\varepsilon$-close to $(\hat\xi,\tilde\xi)$.
\end{lem}
The proof of Lemma~\ref{lema1437} and the next Lemma~\ref{lema1438} are
similar, for example, to the proof of \cite[Lemma 7.26]{equikura}.
So we omit it.

Let
$\tilde\eta = \bigl(\bigl(\Sigma^{\heartsuit},\vec z^{ \heartsuit}\bigr),u^{\heartsuit},\gamma^{\heartsuit}\bigr)$
be a candidate of an element of
${\mathcal M}(L_{12};\vec a;E)$.
We define
${\eta = \mathfrak{fg}(\tilde\eta)}$, which is a candidate of an element of
${\mathcal M}'(L_{12};\vec a;E)$ in the sense of
Definition~\ref{defn1233} in the same way
as Lemma--Definition~\ref{lemdef14355}.

\begin{lem}\label{lema1438}
Let $\tilde\xi \in {\mathcal M}(L_{12};\vec a;E)$ and
$\hat\xi \in {\mathcal M}_{\ell}(L_{12};\vec a;E)$
as in Definition {\rm\ref{defn1436}}.
We fix a stabilization data $\mathscr{ST}$ $($see Definition {\rm\ref{defn1223})}
for $\xi = \mathfrak{fg}\bigl(\tilde\xi\bigr) \in {\mathcal M}(L_{12};\vec a;E)$.
Then for each $\varepsilon>0$ there exists $\delta>0$
with the following properties. If $\eta$ is $\delta$-close to $\bigl(\tilde\xi,\hat\xi\bigr)$
in the sense of Definition {\rm\ref{defn1436}},
then $\eta = \mathfrak{fg}(\tilde\eta)$ is
$\varepsilon$-close to $(\xi,\mathscr{ST})$
in the sense of Definition~{\rm\ref{defn123534522}}.
\end{lem}
Now we will describe the process to pullback Kuranishi structure on
${\mathcal M}'(L_{12};\vec a;E)$ to one on
${\mathcal M}(L_{12};\vec a;E)$.
Let $\xi = \mathfrak{fg}\bigl(\tilde\xi\bigr) \in {\mathcal M}'(L_{12};\vec a;E)$.
We take an obstruction bundle data $\mathscr{OB}$ at
$\xi$ in the sense of Definition~\ref{defn12310}.

Let $\tilde\eta$ be a candidate of an element of
${\mathcal M}(L_{12};\vec a;E)$ which is $\varepsilon$-close to
$\bigl(\tilde\xi,\hat\xi\bigr)$.
(Here we fix $\hat\xi$. Lemma~\ref{lema1437}
shows that the Kuranishi chart we obtain
below is independent of this choice in a~neighborhood of the origin.)

We put $\tilde\eta = \bigl(\bigl(\Sigma^{\heartsuit},\vec z^{ \heartsuit}\bigr),u^{\heartsuit},\gamma^{\heartsuit}\bigr)$
and
\[
{\eta} = \mathfrak{fg}(\tilde\eta)
=\bigl(\bigl(\bigl(\Sigma^{\heartsuit}_1,\vec z^{ \heartsuit}_1\bigr),u^{\heartsuit}_1\bigr),\bigl(\bigl(\Sigma^{\heartsuit}_2,\vec z^{ \heartsuit}_2
\bigr),u^{\heartsuit}_2\bigr),\mathscr I^{\heartsuit},\gamma^{\heartsuit}\bigr).
\]
By Definition~\ref{defn1236}, we have a finite-dimensional linear subspace
$\mathcal E(\xi,\mathscr{OB};\eta)$
of
\[
\bigoplus_{{\rm a}^{\heartsuit}} C_0^{\infty}\bigl(\Sigma^{\rm d,\heartsuit}_{{\rm a}^{\heartsuit}};\bigl(u^{\heartsuit}_1,u^{\heartsuit}_2\bigr)^*T(-X_1 \times X_2)\otimes \Lambda^{0,1}\bigr)
\oplus
\bigoplus_{i=1,2}\bigoplus_{{\rm a}^{\heartsuit}} C_0^{\infty}\bigl(\Sigma^{\rm s,\heartsuit}_{i,{\rm a}^{\heartsuit}};\bigl(u^{\heartsuit}_i\bigr)^*T(X_i)\otimes \Lambda^{0,1}\bigr).
\]
We observe that there exists a map
\smash{$\mathfrak I^{\heartsuit}_i \colon \Sigma^{\heartsuit} \to \Sigma^{\heartsuit}_i$}
which is either bi-holomorphic or a constant map, on each irreducible
component.
We can pull back the subspace $\mathcal E(\xi,\mathscr{OB};{\eta})$
by $\mathfrak I^{\heartsuit}_1$, $\mathfrak I^{\heartsuit}_2$ and
obtain a finite-dimensional linear subspace
of
\[
\bigoplus_{{\rm a}} C_0^{\infty}\bigl(\Sigma^{\heartsuit}_{{\rm a}};
\bigl(u^{\heartsuit}_a\bigr)^*T(-X_1 \times X_2)\otimes \Lambda^{0,1}\bigr).
\]
Here the index ${\rm a}$ runs in the set of irreducible components of
$\Sigma^{\heartsuit}$ and
$\Sigma^{\heartsuit}_{{\rm a}}$ is the irreducible component
corresponding to ${\rm a}$.
We denote this subset by $\mathcal E(\xi,\mathscr{OB};\tilde\eta)$.
(Note that this subspace depends only on $\xi$, $\tilde\eta$ but is
independent of the lift $\tilde\xi$. This is because
$\mathcal E(\xi,\mathscr{OB};{\eta})$ is zero
on the part where we shrink $\Sigma$ to define $\mathfrak{fg}$.)

While we defined a Kuranishi structure on
${\mathcal M}'(L_{12};\vec a;E)$ we made choices of
a finite set
$\{\xi_{\mathfrak i} \mid \mathfrak i \in {\bf I}\} \subset {\mathcal M}'(L_{12};\vec a;E)
$
\eqref{set1246} and
a closed set $\mathfrak N(\xi_{\mathfrak i}) \subset {\mathcal M}'(L_{12};\vec a;E)$ satisfying
\eqref{openset1247}.
We defined a subset
${\bf I}(\xi)$ in \eqref{set1248}.
We use them to define a Kuranishi chart at $\tilde\xi \in {\mathcal M}(L_{12};\vec a;E)$
in the same way as Definition~\ref{defn12383} as follows.

\begin{defn}\label{defn1238322}
We fix $\hat\xi$ and take a sufficiently small
positive number $\varepsilon$ and define
$U\bigl(\tilde\xi;\varepsilon\bigr)$ to be the isomorphism classes of $\tilde\eta
= \bigl(\bigl(\Sigma^{\heartsuit},\vec z^{ \heartsuit}\bigr),u^{\heartsuit},\gamma^{\heartsuit}\bigr)$
with the following properties:
\begin{enumerate}\itemsep=0pt
\item[(1)]
$\tilde\eta$ is a candidate of an element of extended moduli space
 ${\mathcal M}(L_{12};\vec a;E)$.
\item[(2)]
$\eta$ is $\varepsilon$ close to $\bigl(\tilde\xi,\hat\xi\bigr)$.
\item[(3)]
$
\overline\partial u^{\heartsuit} \in \bigoplus_{\mathfrak i \in {\bf I}(\xi)} \mathcal E(\xi_{\mathfrak i},\mathscr{OB};\tilde\eta)$.
\end{enumerate}

Let \smash{$\Gamma_{\tilde\xi}$} be the set of all automorphisms of $\tilde\xi$.
It acts on $U\bigl(\tilde\xi;\varepsilon\bigr)$
and the quotient space is an orbifold $V\bigl(\tilde\xi;\varepsilon\bigr)$.
\end{defn}
We can define $\mathcal E\bigl(\tilde\xi\bigr)$ \big(an orbibundle on $V\bigl(\tilde\xi;\varepsilon\bigr)$\big),
its section \smash{$s_{\tilde\xi}$}, and a map
$\smash{\psi_{\tilde\xi} \colon s_{\tilde\xi}^{-1}(0) }\to {\mathcal M}(L_{12};\vec a;E)$
which is a homeomorphism onto an open neighborhood of \smash{$\tilde\xi$}.
We can show that~$\bigl(V\bigl(\tilde\xi;\varepsilon\bigr),\mathcal E\bigl(\tilde\xi\bigr),s_{\tilde\xi},
\psi_{\tilde\xi}\bigr)$ is a Kuranishi chart at $\tilde\xi$ of
${\mathcal M}(L_{12};\vec a;E)$ in the same way as
the proof of Proposition~\ref{prop1239}.
We thus defined a Kuranishi structure on ${\mathcal M}(L_{12};\vec a;E)$.
We call it the induced Kuranishi structure.
\begin{lem}
For a given system of CF-perturbations on
${\mathcal M}'(L_{12};\vec a;E)$,
which induces a~filtered $A_{\infty}$ algebra
structure $\mathfrak m'_k$ on $CF(L_{12},\Lambda_0)$,
we can define a system of CF-perturbations on the induced
Kuranishi structures of ${\mathcal M}(L_{12};\vec a;E)$,
so that the filtered $A_{\infty}$ algebra
structure~$\mathfrak m''_k$ induced by it on $CF(L_{12},\Lambda_0)$
is exactly the same as $\mathfrak m'_k$.
\end{lem}
\begin{proof}
There exists a group homomorphism \smash{$\Gamma_{\tilde\xi}\to \Gamma_{\xi}$}
and an equivariant map $U\bigl(\tilde\xi;\varepsilon\bigr) \to U(\xi;\varepsilon)$.
Moreover, there exists \smash{$\mathcal E\bigl(\tilde\xi\bigr) \to \mathcal E(\xi)$}
which can be identified with an equivariant bundle map
which covers \smash{$U\bigl(\tilde\xi;\varepsilon\bigr) \to U(\xi;\varepsilon)$}.
Thus the given CF-perturbation on ${\mathcal M}'(L_{12};\vec a;E)$
can be lifted to a~CF-perturbation on the
induced Kuranishi structure.
Since evaluation maps are compatible with
$U\bigl(\tilde\xi;\varepsilon\bigr) \to U(\xi;\varepsilon)$,
and this map is an isomorphism outside a set of
codimension 2,
the operations $\mathfrak m'_k$ obtained by the
CF-perturbation is the same as the
operations $\mathfrak m'_k$ obtained by the
pull-backed CF-perturbation.
The lemma follows.
\end{proof}

Now we have two systems of Kuranishi structures and its
CF-perturbations. One (the induced Kuranishi structures
and its induced CF-perturbations) gives $\mathfrak m'_k = \mathfrak m''_k$.
The other gives $\mathfrak m_k$. In other words, $\mathfrak m_k$ is obtained
from the Kuranishi structures and the CF-perturbations,
which we described in
Section~\ref{sec:HFIm}.
We can find a system of Kuranishi structures of
${\mathcal M}(L_{12};\vec a;E) \times [0,1]$
and their CF-perturbations which interpolates the two systems of
Kuranishi structures and CF-perturbations,
in the same way as Propositions \ref{prop14555} and
\ref{prop1416}. We use it in the same way as
Lemma~\ref{lem1418} to obtain the required pseudo-isotopy.
(We need to use a pseudo-isotopy of pseudo-isotopies
to take homotopy inductive limit. This step again is the same as Section~\ref{sec:highpiso} and so is omitted.)
The proof of Proposition~\ref{prop143434} is complete.
\end{proof}

\section[Independence of the filtered $A_\infty$ functors of the
Hamiltonian isotopy]{Independence of the filtered $\boldsymbol{A_{\infty}}$ functors\\ of the
Hamiltonian isotopy}
\label{sec:independence3}

\subsection{Algebraic preliminary}
\label{sec:alghammilton}

In this section, we prove that
if $(L_1,b_1)$ is Hamiltonian equivalent to
$(L'_1,b'_1)$ and $(L_{12},b_{12})$
is Hamiltonian equivalent to $(L'_{12},b'_{12})$
then the functor
$\mathcal W_{(L_{12},b_{12})}(L_1,b_1)$
is homotopy equivalent to the functor $\mathcal W_{(L'_{12},b'_{12})}(L'_1,b'_1)$
in the category $\mathfrak{Fukst}(X_2)$ over $\Lambda$ coefficient.

To state and prove this result, we start with an algebraic preliminary.
Let $\mathscr C$ be a filtered~$A_{\infty}$ category.
We consider its associated strict category $\mathscr C^s$.
\begin{defn}\label{defn151}
We define an $A_{\infty}$ category $\mathscr C^{\Lambda}$\index[syindex]{Cscrlambda@$\mathscr C^{\Lambda}$}
as follows:
\begin{enumerate}\itemsep=0pt
\item[(1)]
$\mathfrak{OB}\bigl(\mathscr C^{\Lambda}\bigr) =
\mathfrak{OB}(\mathscr C^s)$.
\item[(2)]
For $(c_1,b_1), (c_2,b_2) \in \mathfrak{OB}\bigl(\mathscr C^{\Lambda}\bigr) =
\mathfrak{OB}(\mathscr C^s)$, we put
\[
\mathscr C^{\Lambda}((c_1,b_1), (c_2,b_2))
=
\mathscr C^s((c_1,b_1), (c_2,b_2)) \otimes_{\Lambda_0}\Lambda.
\]
\item[(3)]
The structure operations of $\mathscr C^{\Lambda}$
is obtained by extending the structure operations of
$\mathscr C^s$ by $\Lambda$ linearity.\footnote{We remark
that structure operations of
$\mathscr C^s$ are $\Lambda_0$ linear.\index{$\Lambda_0$ linear}}
\end{enumerate}

\end{defn}
\begin{defn}\label{defn15no2}
In the situation of Definition~\ref{defn151},
let $(c_1,b_1), (c_2,b_2) \in \mathfrak{OB}\bigl(\mathscr C^{\Lambda}\bigr) =
\mathfrak{OB}(\mathscr C^s)$.
We assume $\mathscr C$ is unital.
We say $(c_1,b_1)$ is {\it homotopy equivalent to $(c_2,b_2)$ over $\Lambda$}
\index{homotopy equivalent over $\Lambda$}
and write $(c_1,b_1) \sim_{\Lambda} (c_2,b_2)$\index[syindex]{szimlambda@$\sim_{\Lambda} $}
if they are homotopy equivalent as objects of $\mathscr C^{\Lambda}$.
Suppose $(c_1,b_1) \sim_{\Lambda} (c_2,b_2)$.
We define the {\it Hofer distance}\index{Hofer distance} $d_{\rm Hof}((c_1,b_1),(c_2,b_2))$
between them as the infimum of the positive numbers $\varepsilon$ \index[syindex]{dHof@$d_{\rm Hof}$}
such that the following holds:
\begin{enumerate}\itemsep=0pt
\item[(1)]
There exists $x_{12} \in \mathscr C^{\Lambda}((c_1,b_1), (c_2,b_2))$
$x_{21} \in \mathscr C^{\Lambda}((c_2,b_2), (c_1,b_1))$,
$y_{1} \in \mathscr C^{\Lambda}((c_1,b_1),\allowbreak (c_1,b_1))$,
$y_{2} \in \mathscr C^{\Lambda}((c_2,b_2), (c_2,b_2))$,
such that
\begin{enumerate}\itemsep=0pt\samepage
\item $\mathfrak m_2(x_{21},x_{12}) = {\bf e}_{c_2} + \mathfrak m_1(y_2)$.
\item $\mathfrak m_2(x_{12},x_{21}) = {\bf e}_{c_1} + \mathfrak m_1(y_1)$.
\item $\mathfrak m_1(x_{21}) = 0$. $\mathfrak m_1(x_{12}) = 0$.
\end{enumerate}
\item[(2)]
We require
$T^{\varepsilon_1}x_{21} \in \mathscr C^s((c_1,b_1), (c_2,b_2))$,
$T^{\varepsilon_2}x_{12} \in \mathscr C^s((c_2,b_2), (c_1,b_1))$,
where
$\varepsilon_1$, $\varepsilon_2$ are positive numbers with
$\varepsilon_1 + \varepsilon_2 \le \varepsilon$.
We also require
$T^{\varepsilon}y_{1} \in \mathscr C^s((c_1,b_1), (c_1,b_1))$,
$T^{\varepsilon}y_{2} \in \mathscr C^s((c_2,b_2), (c_2,b_2))$.
\end{enumerate}

\end{defn}
It is easy to see that $\sim_{\Lambda}$ is an equivalence relation.
It is also easy to see that
\begin{equation}\label{form151}
d_{\rm Hof}((c_1,b_1),(c_2,b_2))
+
d_{\rm Hof}((c_2,b_2),(c_3,b_3))
\ge
d_{\rm Hof}((c_1,b_1),(c_3,b_3)).
\end{equation}
We also remark that if $(c_1,b_1)$ is homotopy equivalent
to $(c_2,b_2)$ as objects of $\mathscr C^s$, then
\[
d_{\rm Hof}((c_1,b_1),(c_2,b_2)) = 0.
\]
The next lemma is also easy to show.
\begin{lem}\label{lema153333}
Let $\mathscr F \colon \mathscr C_1 \to \mathscr C_2$
be a strict and unital homotopy equivalence of filtered
$A_{\infty}$ categories.
Then the following holds for $c_1,c_2 \in \mathfrak{OB}(\mathscr C_1)$:
\begin{enumerate}\itemsep=0pt
\item[$(1)$]
$c_1 \sim_{\Lambda} c_2$ if and only if
$\mathscr F(c_1) \sim_{\Lambda} \mathscr F(c_2)$.
\item[$(2)$]
$
d_{\rm Hof}(c_1,c_2)
=
d_{\rm Hof}(\mathscr F(c_1),\mathscr F(c_2))$.
\end{enumerate}

\end{lem}

\subsection[Homotopy equivalence over Lambda in the
geometric situation]{Homotopy equivalence over $\boldsymbol{\Lambda}$ in the
geometric situation}
\label{sec:equivalencehammilton}

\begin{situ}\label{situ153}
Let $(X,\omega)$ be a symplectic manifold which is compact or
tame and $V$ a background datum.
Suppose that a map
$\Phi \colon X \to X$ is a Hamiltonian diffeomorphism
generated by a compactly supported time dependent
Hamiltonian $H \colon X \times [0,1] \to \R$.
We take a finite set
$\mathbb L$ of $V$-relatively spin compact Lagrangian
submanifolds of $X$.
We assume that it is a clean collection.
We assume $L \in \mathbb L$ and $\Phi(L) \in \mathbb L$.
\end{situ}

\begin{thm}\label{thm154}
In Situation {\rm\ref{situ153}},
let $b \in CF(L)$ be a bounding cochain.
\begin{enumerate}\itemsep=0pt
\item[$(1)$]
There exists a bounding cochain $\Phi_*(b) \in CF(\Phi(L))$.
\item[$(2)$]
$(L,b)$ is equivalent to $(\Phi(L),\Phi_*(b))$ over $\Lambda$.
$($Note that they are objects of $\mathfrak{Fukst}(X,\mathbb L)$.$)$
\item[$(3)$]
The Hofer distance between
$(L,b)$ and $(\Phi(L),\Phi_*(b))$
is not greater than the Hofer distance~{\rm\cite{Ho}} between $\Phi$ and the identity
map.
\end{enumerate}

\end{thm}
Theorem~\ref{thm154}\,(1) is a slightly stronger version of \cite[Theorem G\,(G4)]{fooobook}.
Theorem~\ref{thm154}\,(2) is a~slightly stronger version of \cite[Theorem 6.1.25]{fooobook}
(see also \cite{fooo:polydisk}).
We explain how Theorem~\ref{thm154} follows from the
argument of the above quoted papers
\cite{fooobook,fooo:polydisk} in Section~\ref{sec:equivalencehammilton}.

We also remark the following.
\begin{prop}\label{prop156}
If $(L,b),(L',b') \in \mathfrak{OB}(\mathfrak{Fukst}(X;\mathbb L))$
and $L \ne L'$, then
\[
d_{\rm Hof}((L,b), (L',b')) > 0.
\]

\end{prop}
\begin{proof}
Let $L$ be a relatively spin (immersed) Lagrangian submanifold of $(X,\omega)$.
In \cite[Definition 6.5.42]{fooobook}, we defined the notion of a bounding cochain modulo $T^{E}$ as an element $b$ of~$CF(L;\Lambda_+)$ such that
\[
\sum_{k=0}^{\infty}\mathfrak m_k(b,\dots,b) \equiv 0 \mod T^{E}.
\]
When $(L_1,b_1)$, $(L_2,b_2)$ are pairs of Lagrangian submanifolds
with bounding cochains modulo~$T^{E}$, we can define
Floer homology over $\Lambda_0/T^E$ as follows
(see \cite[Definition 6.5.45]{fooobook}).
Let~$CF(L_1,L_2)$ be the left $CF(L_1;\Lambda_0)$
and right $CF(L_2;\Lambda_0)$ bi-module,
which is nothing but the morphism space from
$(L_1,b_1)$ to $(L_2,b_2)$ in the curved $A_{\infty}$
category of $X$.
Let
\[
\mathfrak n_{k_1,k_2}\colon\ B_{k_1}CF(L_1)[1] \otimes CF(L_1,L_2)
\otimes B_{k_2}CF(L_2)[1] \to CF(L_2;\Lambda_0)
\]
be the structure operations.
We put
\[
\delta_{b_1,b_2}(x) = \sum_{k_1,k_2=0}^{\infty}\mathfrak n_{k_1,k_2}\bigl(b_1^{k_1},x,b_2^{k_2}\bigr).
\]
The $A_{\infty}$ relations imply $\delta_{b_1,b_2} \circ \delta_{b_1,b_2} \equiv 0 \mod T^E$.
Therefore, $\delta_{b_1,b_2}$ becomes a boundary operator
on $CF(L_1,L_2) \otimes_{\Lambda_0} \Lambda_0/T^E$.
Its cohomology is by definition $HF\bigl((L_1,b_1),(L_2,b_2);\Lambda_0/T^E\bigr)$.
It is independent of the choices of perturbations and almost complex structures.
\begin{lem}\label{lem27rev}
Let $(L_1,b_1)$, $(L_2,b_2)$ be objects of $\mathfrak{Fukst}(X,\mathbb L)$
and $(L,b)$ a pair of an element of~$\mathbb L$ and its bounding cochain modulo $T^E$.
Suppose $d_{\rm Hof}((L_1,b_1), (L_2,b_2)) = 0$.
Then
\[
HF\bigl((L_1,b_1),(L,b);\Lambda_0/T^E\bigr)
\cong
HF\bigl((L_2,b_2),(L,b);\Lambda_0/T^E\bigr).
\]

\end{lem}
\begin{proof}
By perturbing a bit and using \cite[Theorem 6.5.47]{fooobook} we may assume
that $L_1$ and~$L_2$ are transversal to $L$.
Then, for an arbitrary small $\varepsilon$, there exist
$x_{12} \in CF(L_1,L_2;\Lambda)$ and $x_{21} \in CF(L_2,L_1;\Lambda)$
as in Definition~\ref{defn15no2}.
Multiplications with $T^{\varepsilon}x_{12}$ and with $T^{\varepsilon}x_{21}$
define chain maps
\begin{gather*}
\varphi_{12} \colon\ CF(L_1,L) \otimes_{\Lambda_0} \Lambda_0/T^{E'} \to CF(L_2,L) \otimes_{\Lambda_0} \Lambda_0/T^{E'},
\\
\varphi_{21} \colon\ CF(L_2,L) \otimes_{\Lambda_0} \Lambda_0/T^{E'} \to CF(L_1,L) \otimes_{\Lambda_0} \Lambda_0/T^{E'}
\end{gather*}
for any $E' \le E$. Moreover, using Definition~\ref{defn15no2}\,(2)
we can show that
\begin{gather*}
\varphi_{12} \circ \varphi_{21} \colon\
CF(L_1,L) \otimes_{\Lambda_0} \Lambda_0/T^{E} \to CF(L_1,L) \otimes_{\Lambda_0} \Lambda_0/T^{E}
\\
\varphi_{21} \circ \varphi_{12} \colon\
CF(L_2,L) \otimes_{\Lambda_0} \Lambda_0/T^{E} \to CF(L_2,L) \otimes_{\Lambda_0} \Lambda_0/T^{E}
\end{gather*}
are chain homotopic to $T^{2\varepsilon}$ times the identity map.
We write
\[
HF\bigl((L_i,b_i),(L,b);\Lambda_0/T^E\bigr) = \sum_{j=1}^{N_i} \Lambda_0/T^{a_{i,j}},
\]
where $a_{i,j} \le E$ and $a_{i,j} > a_{i,j+1}$.
Then using $\varphi_{12}$, $\varphi_{21}$
and their properties explained above, we have the following: if $a_{1,j} > 4\varepsilon$, we have $\vert a_{1,j} - a_{2,j} \vert \le 2\varepsilon$.
(See \cite[pp 391--392]{fooobook}.)
Since $\varepsilon$ is arbitrary small, we obtain the lemma by taking the limit $\varepsilon \to 0$.
\end{proof}

Now we are in the position to prove Proposition~\ref{prop156}.
Suppose $L \ne L'$. We may assume that there exists $p \in L \setminus L'$.
Let $d = d(p,L')$. Let $\rho$ be a positive number sufficiently small
compared to $d$. We can take a small Clifford type torus $T_{\rho}$
such that
$T_{\rho} \cap L' = \varnothing$,
$T_{\rho} \cap L \ne \varnothing$ and $T_{\rho}$ intersects transversally with $L$.
We may also assume that $T_{\rho}$ admits a bounding cochain~$b_{\rho}$ modulo
$T^{\rho}$.
Since $T_{\rho} \cap L' = \varnothing$,
$
HF((L',b'),(T_{\rho},b_{\rho});\Lambda_0/T^{\rho}) = 0$.
On the other hand, using the fact $T_{\rho} \cap L \ne \varnothing$ and
all the non-constant holomorphic strips have positive energy we can show~${
HF((L,b),(T_{\rho},b_{\rho});\Lambda_0/T^{\rho}) \ne 0}$.
(This is a classical fact going back to Chekanov \cite{chekanov}.)
This contradicts Lemma~\ref{lem27rev}.
\end{proof}

\subsection{The main theorem}
\label{sec:mainhammilton}
\begin{situ}\label{situ155}
Let $X_1$, $V_1$, $\mathbb L_1$, $\Phi_1$, $L_1$ be as in
Situation \ref{situ153}.
Let $(X_2,\omega_2)$ be a compact symplectic manifold and $V_2$ a background datum.
Let $-X_1 \times X_2, \pi_1^*(V_1\oplus TX_1) \oplus \pi_2^*(V_2),
\mathbb L_{12}, L_{12}, \Phi_{12}$ be
also as in Situation \ref{situ153}.
Let $\mathbb L_2$ be
a finite set of $V_2$-relatively spin compact Lagrangian
submanifolds of $X_2$.
We assume that it is a clean collection.
We assume also that for any $L'_1 \in \mathbb L_1$ and $L'_{12} \in \mathbb L_{12}$
the geometric transformation of $L'_1$ by~$L'_{12}$
is contained in $\mathbb L_2$.
\end{situ}

\begin{thm}\label{them1566}
In Situation {\rm\ref{situ155}}, let $b_1$ be a bounding cochain of
$L_1$ and $b_{12}$ a bounding cochain of $L_{12}$.
Then
\begin{enumerate}\itemsep=0pt
\item[$(1)$]
$\mathcal W_{(L_{12},b_{12})}(L_1,b_1)$
is equivalent to $\mathcal W_{(\Phi_{12}L_{12},(\Phi_{12})_*b_{12})}
(\Phi_1(L_1),(\Phi_1)_*(b_1))$ over $\Lambda$.
\item[$(2)$]
The Hofer distance
\[
d_{\rm Hof}
(\mathcal W_{(L_{12},b_{12})}(L_1,b_1),
\mathcal W_{(\Phi_{12}L_{12},(\Phi_{12})_*b_{12})}
(\Phi_1(L_1),(\Phi_1)_*b_1))
\]
is not greater than the sum of the Hofer distance
{\rm\cite{Ho}} between $\Phi_1$ and the identity
map and the Hofer distance
between $\Phi_{12}$ and the identity
map.
\end{enumerate}

\end{thm}
The proof is given in the next subsection.

The next result is a more functorial version of Theorem~\ref{them1566}.
\begin{situ}\label{situ1588}
Let $(X_i,\omega_i)$ be a compact symplectic manifold,
$V_i$ a background datum of~$X_i$,
and $\mathbb L_i$ a finite set of $V_i$ relatively spin
immersed Lagrangian submanifolds, for $i=1,2$.
We assume $\mathbb L_i$ are clean collections.
Let $L_{12}$ be a $\pi_1^*(V_1 \oplus TX_1) \oplus \pi_2^*(V_2)$ relatively spin
Lagrangian submanifold of $-X_1 \times X_2$
and~${\Phi \colon -X_1 \times X_2 \to -X_1 \times X_2}$ a
Hamiltonian diffeomorphism.
We assume that for each $L_1 \in \mathbb L_1$
the geometric transformations
$
L_1 \times_{X_1} L_{12}$,
$
L_1 \times_{X_1} \Phi(L_{12})
$
are both elements of $\mathbb L_2$.
We assume that $L_{12}$ is unobstructed and $b_{12}$
is its bounding cochain.
By Theorem~\ref{thm154}, we obtain a bounding cochain
$\Phi_*(b_{12})$ of $\Phi(L_{12})$.
\end{situ}

\begin{thm}\label{thm15999}
In Situation {\rm\ref{situ1588}}, we consider two filtered $A_{\infty}$ functors
\begin{gather*}
\mathcal W_{(L_{12},b_{12})} \colon\ \mathfrak{Fukst}(X_1;\mathbb L_1)
\to \mathfrak{Fukst}(X_2;\mathbb L_2), \\
\mathcal W_{(\Phi(L_{12}),\Phi_*(b_{12}))} \colon\ \mathfrak{Fukst}(X_1;\mathbb L_1)
\to \mathfrak{Fukst}(X_2;\mathbb L_2).
\end{gather*}
They induce the following filtered $A_{\infty}$ functors
of $\Lambda$ linear categories in an obvious way:
\begin{gather*}
\mathcal W_{(L_{12},b_{12})}^{\Lambda} \colon\ \mathfrak{Fukst}(X_1;\mathbb L_1)^{\Lambda}
\to \mathfrak{Fukst}(X_2;\mathbb L_2)^{\Lambda}, \\
\mathcal W_{(\Phi(L_{12}),\Phi_*(b_{12}))}^{\Lambda} \colon\ \mathfrak{Fukst}(X_1;\mathbb L_1)^{\Lambda}
\to \mathfrak{Fukst}(X_2;\mathbb L_2)^{\Lambda}.
\end{gather*}
\begin{enumerate}\itemsep=0pt
\item[$(1)$]
$\mathcal W_{(L_{12},b_{12})}^{\Lambda}$ is homotopy
equivalent to
$\mathcal W_{(\Phi(L_{12}),\Phi_*(b_{12}))}^{\Lambda}$.
\item[$(2)$]
The Hofer distance between $\mathcal W_{(L_{12},b_{12})}$
and $\mathcal W_{(\Phi(L_{12}),\Phi_*(b_{12}))}$
in the filtered $A_{\infty}$ category
$\mathcal{FUNC}(\mathfrak{Fukst}(X_1;\mathbb L_1),\mathfrak{Fukst}(X_2;\mathbb L_2))$
is not greater than the Hofer distance between $\Phi$ and the identity map.
\end{enumerate}

\end{thm}
The proof is given in the next subsection.
\begin{rem}\label{rem151212}
\quad
\begin{enumerate}\itemsep=0pt
\item[(1)]
The two immersed Lagrangian submanifolds $L_1 \times_{X_1} L_{12}$
and $L_1 \times_{X_1} \Phi(L_{12})$ may not be isotopic each other in general.
So Theorem~\ref{thm15999}
provides a lot of examples of a pair of Lagrangian submanifolds
which are not isotopic but are Floer theoretically equivalent.
\item[(2)]
K. Ono \cite{ono} studied a Lagrangian intersection between $L$ and $L'$
where the lifts of $L$ and $L'$ to the prequantum bundle are Hamiltonian
isotopic each other. Theorem~\ref{thm15999} is related to his study.
\item[(3)]
We recall that two Lagrangian submanifolds $L,L' \in X$ are said
to be Lagrangian cobordant if there exists
a Lagrangian submanifold $\tilde L$ in $\C \times X$
and a sufficiently large ball $D(R)$ of~$\C$ centered at $0$ such that
\[
\tilde L \cap ((\C\setminus D(R)) \times X)
=
(((-\infty,0) \times L) \cup ((0,\infty) \times L'))
\cap ((\C\setminus D(R)) \times X).
\]
We can show $L_1 \times_{X_1} L_{12}$
is Lagrangian cobordant to $L_1 \times_{X_1} \Phi(L_{12})$ in this sense.
\item[(4)]
In the situation of item (3),
assuming $L,L',\tilde L$ are monotone and $L'' \subset X$ is also monotone
Biran--Cornea \cite{BC} proved
$
HF(L,L'') \cong HF(L',L'')$.
It seems likely that we can generalize it as follows.
Suppose $L$, $L'$ have bounding cochains $b$, $b'$, respectively.
Moreover, we assume that there exists a bounding cochain $\tilde b$ of
$\tilde L$ such that on $ ((\C\setminus D(R)) \times X)$
it coincides with the pullbacks of $b$ and $b'$.
Then
\[
HF((L,b),(L'',b'');\Lambda) \cong HF((L',b'),(L'',b'');\Lambda).
\]
We say $(L,b)$ is unobstructed-Lagrangian cobordant to
$(L',b')$ in this situation.
We can then try to use the argument of the proof of Theorem~\ref{thm61}
to prove the following.
Let $(L_1,b_1)$, $(L_{12},b_{12})$ be objects of
$\mathfrak{Fukst}(X_1;\mathbb L_1)$, $\mathfrak{Fukst}(X_1;\mathbb L_1)$,
$\mathfrak{Fukst}((X_1,\omega_1),\mathbb L_1)
\times
\mathfrak{Fukst}((X_1 \times X_2,-\pi_1^*(\omega_1)
+\pi_2^* (\omega_2)),\mathbb L_{12})$,
respectively.
Let $\Phi \colon -X_1 \times X_2 \to -X_1 \times X_2$ be a
Hamiltonian diffeomorphism and
$(\Phi(L_{12}),\Phi_*(b_{12})$) be also
an object of $\mathfrak{Fukst}((X_1 \times X_2,-\pi_1^*(\omega_1)
+\pi_2^* (\omega_2)),\mathbb L_{12})$.
We put $L_2 = L_1 \times_{X_1} L_{12}$ and
$L'_2 = L_1 \times_{X_1} \Phi(L_{12})$.
We obtain their bounding cochains by Theorem~\ref{thm61},
which we denote by $b_2$, $b'_2$.
Then $(L_2,b_2)$ is unobstructed-Lagrangian cobordant to
$(L'_2,b'_2)$.

This argument can be an alternative proof of Theorem~\ref{thm15999}\,(1).
\item[(5)]
Cornea--Shelukhin \cite{CS} study the area of the image $\pi\bigl(\tilde L\bigr)$
of the Lagrangian cobordism $\tilde L$ by the projection $\pi\colon
\C \times X \to \C$. Including the bounding cochain,
their argument may imply that
if $(L,b)$ is unobstructed-Lagrangian cobordant to
$(L',b')$ by a pair $\bigl(\tilde L,\tilde b\bigr)$,
then the Hofer distance (in the sense of Definition~\ref{defn15no2})
is not greater than the area of \smash{$\pi\bigl(\tilde L\bigr)$}.
This statement is related to Theorem~\ref{thm15999}\,(2).
\end{enumerate}

\end{rem}

\subsection{Proof of the main theorem}
\label{sec:proofmainhammilton}

Theorem~\ref{them1566} is an immediate consequence of
Theorem~\ref{thm154} and the following purely algebraic
result.

\begin{prop}
Let $\mathscr C_1$, $\mathscr C_2$, $\mathscr C_3$
be strict and unital filtered $A_{\infty}$ categories
and
$
\mathscr F \colon \mathscr C_1 \times \mathscr C_2\allowbreak
\to \mathscr C_3
$
a strict and unital filtered $A_{\infty}$ bi-functor.
Suppose $c_1,c'_1 \in \mathfrak{OB}(\mathscr C_1)$,
$c_2,c'_2 \in \mathfrak{OB}(\mathscr C_2)$.
\begin{enumerate}\itemsep=0pt
\item[$(1)$]
If $c_1 \sim_{\Lambda} c'_1$, $c_2 \sim_{\Lambda} c'_2$,
then
$
\mathscr F(c_1,c_2) \sim_{\Lambda} \mathscr F(c'_1,c'_2)$.
\item[$(2)$]
$
d_{\rm Hof}(\mathscr F(c_1,c_2),\mathscr F(c'_1,c'_2))
\le
d_{\rm Hof}(c_1,c'_1)
+
d_{\rm Hof}(c_2,c'_2)$.
\end{enumerate}

\end{prop}
\begin{proof}
By \eqref{form151}, it suffices to show the case $c_1 = c'_1$
and the case $c_2 = c'_2$. By symmetry, it suffices to prove
the case $c_2 = c'_2$.
Let $\varepsilon > d_{\rm Hof}(c_1,c'_1)$ and we take
$x_1 \in \mathscr C_1^{\Lambda}(c_1,c'_1)$,
$x_2 \in \mathscr C_1^{\Lambda}(c'_1,c_1)$,
$y_1 \in \mathscr C_1^{\Lambda}(c_1,c_1)$,
$y_2 \in \mathscr C_1^{\Lambda}(c'_1,c'_1)$
such that
\[\mathfrak m_2(x_1,x_2) = {\bf e}_{c_1} + \mathfrak m_1(y_1),
\qquad
\mathfrak m_2(x_2,x_1) = {\bf e}_{c'_1} + \mathfrak m_1(y_2),
\qquad
 \mathfrak m_1(x_1) = \mathfrak m_1(x_2) = 0.
\]
We also assume that
$T^{\varepsilon_1}x_1 \in C_1(c_1,c'_1)$,
$T^{\varepsilon_2}x_2 \in C_1(c'_1,c_1)$
with $\varepsilon_1 + \varepsilon_2 \le \varepsilon$
and
$T^{\varepsilon}y_1 \in C_1(c_1,c_1)$,
$T^{\varepsilon}y_2 \in C_1(c'_1,c'_1)$.

We put
\begin{gather*}
\mathfrak x_1 := \mathscr F_{1,1}(x_1,{\bf e}_{c_2}),
\qquad
\mathfrak x_2 := \mathscr F_{1,1}(x_2,{\bf e}_{c_2}),
\qquad
\mathfrak y_1 := \mathscr F_{1,1}(y_1,{\bf e}_{c_2}),
\\
\mathfrak y_2 := \mathscr F_{1,1}(y_2,{\bf e}_{c_2}).
\end{gather*}
(Note that we extend $\mathscr F_{1,1}$ by $\Lambda$
linearity to define the right-hand sides.)
Since $\mathscr F$ is strict, we have
\begin{align*}
\mathfrak m_2(\mathfrak x_1,\mathfrak x_2)
&= \mathscr F_{1,1}(\mathfrak m_2(x_1,x_2),{\bf e}_{c_2})
= \mathscr F_{1,1}({\bf e}_{c_1} + \mathfrak m_1(y_1),{\bf e}_{c_2})
\\
&= {\bf e}_{\mathscr F(c_1,c_2)} + \mathfrak m_1(\mathscr F_{1,1}(y_1,{\bf e}_{c_2}))
= {\bf e}_{\mathscr F(c_1,c_2)} +\mathfrak m_1(\mathfrak y_1).
\end{align*}
Similarly, we have
$
\mathfrak m_2(\mathfrak x_2,\mathfrak x_1)
=
{\bf e}_{\mathscr F(c'_1,c_2)} +\mathfrak m_1(\mathfrak y_2)$,
and $\mathfrak m_1(\mathfrak x_1) = \mathfrak m_1(\mathfrak x_2) = 0$.
Therefore,
$\mathscr F(c_1,c_2) \allowbreak\sim_{\Lambda} \mathscr F(c'_1,c_2)$.
(1) follows.
(2) follows from
\begin{gather*}
T^{\varepsilon_1}\mathfrak x_1 \in \mathcal C_3(\mathscr F(c_1,c_2),\mathscr F(c'_1,c_2)),
\qquad
T^{\varepsilon_2}\mathfrak x_2 \in \mathcal C_3(\mathscr F(c'_1,c_2),\mathscr F(c_1,c_2)),
\\
T^{\varepsilon}\mathfrak y_1 \in \mathcal C_3(\mathscr F(c_1,c_2),\mathscr F(c_1,c_2)),
\qquad
T^{\varepsilon}\mathfrak y_2 \in \mathcal C_3(\mathscr F(c'_1,c_2),\mathscr F(c'_1,c_2)).\tag*{\qed}
\end{gather*}\renewcommand{\qed}{}
\end{proof}

Theorem~\ref{thm15999} is an immediate consequence of
Theorem~\ref{thm154} and the following purely algebraic
result.
\begin{lem}
Let $\mathscr C_1$, $\mathscr C_2$, $\mathscr C_3$
be strict and unital filtered $A_{\infty}$ categories
and
$
\mathscr F \colon \mathscr C_1 \times \mathscr C_2
\to \mathscr C_3
$
a strict and unital filtered $A_{\infty}$ bi-functor.
It induces a strict and unital filtered $A_{\infty}$ functor~${
\mathscr F_* \colon \mathscr C_2
\to \mathcal{FUNC}(\mathscr C_1,\mathscr C_3)}
$
by Lemma {\rm\ref{lem56}} $($and its unital and strict analogue$)$.
Suppose
$c_2,c'_2 \in \mathfrak{OB}(\mathscr C_2)$.
\begin{enumerate}\itemsep=0pt
\item[$(1)$] If $c_2 \sim_{\Lambda} c'_2$, then
the two $(\Lambda$ linear$)$ filtered $A_{\infty}$ functors
$
\mathscr F_*(c_2)^{\Lambda},
\mathscr F_*(c'_2)^{\Lambda} \colon
\mathscr C^{\Lambda}_1 \to \mathscr C^{\Lambda}_3
$
are homotopy equivalent.
\item[$(2)$]
The inequality
$
d_{\rm Hof}(\mathscr F_*(c_2),\mathscr F_*(c'_2))
\le
d_{\rm Hof}(c_2,c'_2)
$
holds.
\end{enumerate}
\end{lem}

The proof is easy and so is omitted.

\subsection[Proof of Theorem~\ref{thm154}]{Proof of Theorem~\ref{thm154}}

In this subsection, we explain how Theorem~\ref{thm154}
follows from (the proof of) \cite[Theorem G\,(G4) and Theorem 6.1.25]{fooobook}
and \cite{fooo:polydisk}.
Suppose we are in Situation \ref{situ153}.

We put $\Phi(L) = L'$.
We take a compatible almost complex structure $J$
and consider filtered~$A_{\infty}$ structures
\[
\mathfrak m_k^{J,L} \colon\ \Omega\bigl(\tilde L \times_X \tilde L\bigr)^{\otimes k} \to \Omega\bigl(
\tilde L \times_X \tilde L\bigr) \,\widehat{\otimes}\, \Lambda_0
,
\qquad
\mathfrak m_k^{J,L'} \colon\ \Omega\bigl(\tilde L' \times_X \tilde L'\bigr)^{\otimes k} \to
\Omega\bigl(\tilde L' \times_X \tilde L'\bigr)
\,\widehat{\otimes}\, \Lambda_0.
\]
Note that we can decompose \smash{$\mathfrak m_k^{J,L}$}, \smash{$\mathfrak m_k^{J,L'}$}
to a sum
\[
\mathfrak m_k^{J,L} = \sum_{E}T^{E} \mathfrak m_{k,E}^{J,L},
\qquad
\mathfrak m_k^{J,L'} = \sum_{E}T^{E} \mathfrak m_{k,E}^{J,L'},
\]
where \smash{$\mathfrak m_{k,E}^{J,L}$} and \smash{$\mathfrak m_{k,E}^{J,L'}$} are $\R$ linear.
\begin{rem}
The right-hand side is an infinite sum. However, for each $E_0$
the set of $E<E_0$ such that
\smash{$\mathfrak m_{k,E}^{J,L}$}, \smash{$\mathfrak m_{k,E}^{J,L'}$} is nonzero is a finite set.
This is a consequence of Gromov compactness.
\end{rem}
We denote $CF(L) = \Omega\bigl(
\tilde L \times_X \tilde L\bigr) \,\widehat{\otimes}\, \Lambda_0$
and
$CF(L') = \Omega\bigl(
\tilde L' \times_X \tilde L'\bigr) \,\widehat{\otimes}\, \Lambda_0$.
The next theorem is the de Rham version of \cite[Corollary 4.6.3]{fooobook}.

\begin{thm}\label{them1513}
There exists a $($curved$)$ filtered $A_{\infty}$ homomorphism
\begin{gather*}
\hat{\mathfrak f} = \{\mathfrak f_k \mid k=0,1,2,\dots\}
\colon\\
 \qquad \bigl(CF(L),\bigl\{\mathfrak m_k^{J,L};k=0,1,2,\dots\bigr\}\bigr) \to \bigl(CF(L'),\bigl\{\mathfrak m_k^{J,L'};k=0,1,2,\dots\bigr\}\bigr),
\end{gather*}
$
\mathfrak f_k \colon CF(L)^{\otimes k} \to CF(L')
$
such that
$
\mathfrak f_1(h) = \bigl(\Phi^{-1}\bigr)^*(h) \mod \Lambda_+$.

\end{thm}
\begin{proof}
Let $J' = \bigl(\Phi^{-1}\bigr)_*J$.
The moduli space of
$J'$ holomorphic disks with the boundary conditions
given by $L$ is canonically identified with
the moduli space of
$J$ holomorphic disks with
the boundary condition
given by $L'$. Therefore,
the following diagram commutes:
\[
\begin{CD}
B_kCF[1](L') @ >{\mathfrak m_k^{J,L'}}>>
CF(L') \\
@ A{(\Phi^{-1})^*}AA @ AA{(\Phi^{-1})^*}A\\
B_kCF[1](L) @ >{\mathfrak m_k^{J',L}}>>
CF(L).
\end{CD}
\]
Therefore, it suffices to construct a filtered $A_{\infty}$
homomorphism $\mathfrak g =
\{\mathfrak g_k\}$ from $\bigl(CF(L),\allowbreak\bigl\{\mathfrak m_k^{J,L}\bigr\}\bigr)$ to \smash{$\bigl(CF(L),\bigl\{\mathfrak m_k^{J',L}\bigr\}\bigr)$}
such that $\mathfrak g_1 \equiv {\rm id}\mod \Lambda_+$.

We take a one parameter family of compatible
almost complex structures $\mathcal J = \bigl\{J^{(\rho)}\bigr\}$
such that
\begin{enumerate}\itemsep=0pt
\item[(1)] $J^{(0)} = J$,
\item[(2)] $J^{(1)} = J'$.
\end{enumerate}
For the proof of Theorem~\ref{them1513},
we can take any such $\mathcal J$.
We will specify $\mathcal J$ later during the proof of
Proposition~\ref{prop1518}.

We use the `time ordered product'
\index{time ordered product} moduli spaces
$
\mathcal M_{k+1}(L;\mathcal J;E;{\rm top}({\rho}))
$
introduced in \cite[Section 4.6.1]{fooobook},
which have the properties spelled out in Proposition~\ref{prop1514} below.
We use the following notation in Proposition~\ref{prop1514}.

The moduli space $\mathcal M_{k+1}(L;E)$ is defined in
\eqref{form317}. To specify the almost complex structure
we use, we write $\mathcal M_{k+1}(L;E;J)$.
It comes with evaluation maps
\[
{\rm ev} = ({\rm ev}_0,{\rm ev}_1,\dots,{\rm ev}_k) \colon\
\mathcal M_{k+1}(L;E;J)
\to
\bigl(\tilde L \times_X \tilde L\bigr)^{k+1},
\]
which is strongly smooth
\index{strongly smooth} and such that ${\rm ev}_0$ is
weakly submersive.
\index{weakly submersive}
(See \cite{foootech2,fooonewbook} and \cite[Part 7]{fooospectr}
for the definition of strong smoothness and weak submersivity.)
\begin{prop}\label{prop1514}
There exists a compact Hausdorff space
$
\mathcal M_{k+1}(L;\mathcal J;E;{\rm top}({\rho}))
$\index[syindex]{M1k+1LJ@$\mathcal M_{k+1}(L;\mathcal J;E;{\rm top}({\rho}))$}
equipped with a
Kuranishi structure with corners, which enjoys
the following properties:
\begin{enumerate}\itemsep=0pt
\item[$(1)$]
There exists an evaluation map
\[
{\rm ev} = ({\rm ev}_0,{\rm ev}_1,\dots,{\rm ev}_k) \colon\
\mathcal M_{k+1}(L;\mathcal J;E;{\rm top}({\rho}))
\to
\bigl(\tilde L \times_X \tilde L\bigr)^{k+1}
\]
which is strongly smooth. Moreover, ${\rm ev}_0$ is weakly submersive.
\item[$(2)$]
The normalized boundary of $
\mathcal M_{k+1}(L;\mathcal J;E;{\rm top}({\rho}))
$
is the union of two types of the fiber products:
\begin{itemize}\itemsep=0pt
\item[$(I)$]
The fiber product
\begin{equation}\label{form151414}
\mathcal M_{k_1+1}(L;E_1;J)
{}_{{\rm ev}_0}\times_{{\rm ev}_i}
\mathcal M_{k_2+1}(L;\mathcal J;E_2;{\rm top}({\rho})),
\end{equation}
where $k_1 + k_2 = k$, $E_1 + E_2 = E$ and
$i=1,\dots,k_2$.

\item[$(II)$]
The
fiber product
\begin{equation}\label{form151515}
\prod_{i=1}^m
\mathcal M_{k_i+1}(L;\mathcal J;E_i;{\rm top}({\rho}))
{}_{({\rm ev}_0,\dots,{\rm ev}_0)}\times_{({\rm ev}_1,\dots,{\rm ev}_m)}
\mathcal M_{m+1}(L;E_0;J'),
\end{equation}
where $k_1+\dots+k_m =k$ and $E_1+\dots+E_m + E_0 = E$.
\end{itemize}
\item[$(3)$]
In the case when $E=0$ and $k=0$,
$
\mathcal M_{1}(L;\mathcal J;0;{\rm top}({\rho}))
=
\tilde L \times_X \tilde L$,
and ${\rm ev}_0$ is the identity map.
\item[$(4)$]
The set of $E$ such that $
\mathcal M_{k+1}(L;\mathcal J;E;{\rm top}({\rho}))
\ne \varnothing $
is discrete.
\end{enumerate}

\end{prop}
\begin{proof}[A sketch of the proof.]
The construction of the moduli spaces
$\mathcal M_{k+1}(L;\!\mathcal J;\!E;\!{\rm top}({\rho}))$
is worked out in detail in \cite[Definition 4.6.1]{fooobook}.\footnote{Actually, we need a slight modification since our Lagrangian submanifolds
are immersed. This modification is the same as the argument of
Section~\ref{sec:HFIm}.}
Its element is an object $((\Sigma,\vec z),u,\gamma,\{\rho_{\alpha_i}\})$
as depicted in Figure~\ref{fig19FOOO}.
Here $(\Sigma,\vec z)$ is a bordered marked curve of genus zero with
one boundary component and $k+1$ boundary marked points, and
$u \colon (\Sigma,\partial \Sigma) \to (X,L)$ is a smooth map.
The restriction of $u$ to $\partial \Sigma \setminus (\vec z
\cup \{\text{boundary node}\})$ is lifted to a map
$\gamma \colon \partial \Sigma \setminus (\vec z
\cup \{\text{boundary node}\})
\to \tilde L$. The map $\alpha_i \mapsto \rho_{\alpha_i}$ assigns a number $\rho_{\alpha_i}
\in [0,1]$ to each irreducible component $\Sigma_{\alpha_i}$
of $\Sigma$.
We require the next Condition \ref{conds1515} for $\rho_{\alpha_i}$.
We also require that the restriction of $u$ to $\Sigma_{\alpha_i}$ is $J^{(\rho_{\alpha_i})}$-holomorphic.
At boundary marked points and boundary nodes, we
require switching conditions similar to those appeared in Section~\ref{sec:HFIm}.
(See Definition~\ref{def3737}\,(5).)
\begin{figure}[ht]
\centering
\includegraphics[scale=0.8]{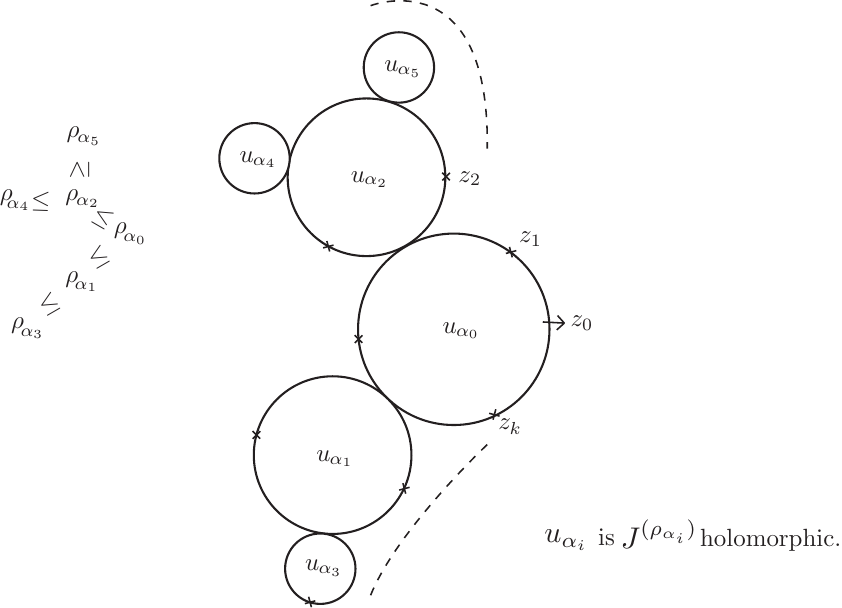}
\caption{Time ordered product moduli space.}
\label{fig19FOOO}
\end{figure}

\begin{conds}\label{conds1515}
Let $p \in \Sigma$ be a boundary node and
$p \in \Sigma_{\alpha_i} \cap \Sigma_{\alpha_j}$,
$\Sigma_{\alpha_i} \ne \Sigma_{\alpha_j}$.
We suppose~$\Sigma_{\alpha_i}$ is contained
in the connected component of $\Sigma \setminus \{p\}$
which contains the zero-th marked point~$z_0$.
Then we require
$
\rho_{\alpha_i} \ge \rho_{\alpha_j}$.
\end{conds}

The definition of the topology of this moduli space
and proof of its compactness and Haus\-dorffness
are similar to those of Theorem~\ref{thekuraexist}.
The construction of the Kuranishi structure is similar
to the proof of Theorem~\ref{thekuraexist}.

We next describe the boundary.
We observe that
\eqref{form151414} corresponds to the case when one of $\rho_{\alpha_i}$
becomes $0$
and that
\eqref{form151515} corresponds to the case when one of $\rho_{\alpha_i}$
becomes $1$.
Actually, such $\Sigma_{\alpha_i}$ is necessary the
irreducible component containing $z_0$, the zero-th marked point.
(This is a consequence of Condition \ref{conds1515}.)

The other possible boundary components of $\mathcal M_{k+1}(L;\mathcal J;E;{\rm top}({\rho}))$
cancel out each other as is explained in \cite[p.~246]{fooobook}.
The key observation is the cancellation between
two types of potential boundaries.
One is depicted in Figure~\ref{fig19.2}
 below and the other is
depicted in Figure~\ref{fig19.3} below.
(Those two figures are \cite[Figure 4.6.2]{fooobook} and
\cite[Figure 4.6.3]{fooobook}, respectively.)
\begin{figure}[ht]
\centering
\includegraphics[scale=0.8]{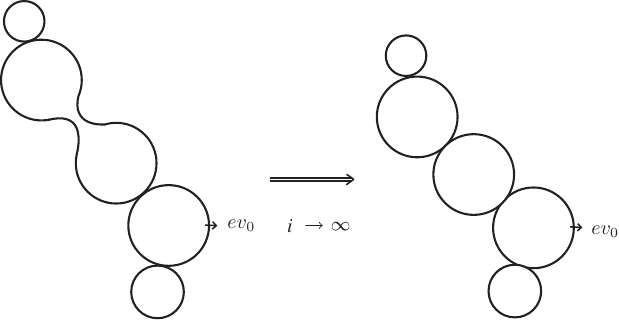}
\caption{Cancellation in \cite[Section 4.6]{fooobook} : 1.}
\label{fig19.2}
\end{figure}
\begin{figure}[ht]
\centering
\includegraphics[scale=0.8]{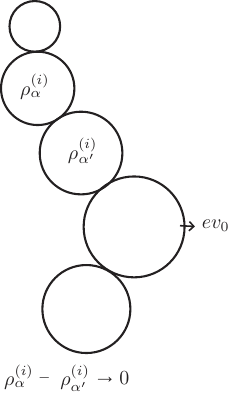}
\caption{Cancellation in \cite[Section 4.6]{fooobook} : 2.}
\label{fig19.3}
\end{figure}

Item (3) of Proposition~\ref{prop1514} is a consequence of the
fact that left-hand side is the moduli space of
constant maps, which is transversal.
Item (4) follows from Gromov compactness.
\end{proof}

For later use, we define
$
\rho_0\colon \mathcal M_{k+1}(L;\mathcal J;E_2;{\rm top}({\rho}))
\to [0,1]
$
\index[syindex]{rzho0@$\rho_0$}
as follows:
\begin{equation}\label{defrho0}
\rho_0((\Sigma,\vec z),u,\gamma,\{\rho_{\alpha_i}\})
= \rho_{\alpha_0},
\end{equation}
where $\Sigma_{\alpha_0}$ is the irreducible component
which contains the $0$-th marked point.

Now using Proposition~\ref{prop1514}, we define $\mathfrak g_{k,E}$ by
the next formula
\begin{equation}\label{form1516}
\mathfrak g_{k,E}(h_1,\dots,h_k) =
{\rm ev}_0 !
\bigl(
{\rm ev}_1^*h_1 \wedge\dots \wedge {\rm ev}_k^*h_k;
(\mathcal M_{k+1}(L;\mathcal J;E;{\rm top}({\rho}))
;\widehat{\mathfrak S}_{\varepsilon})
\bigr).
\end{equation}
Here we take a system of CF-perturbations
$\widehat{\mathfrak S}_{\varepsilon}$
on $\mathcal M_{k+1}(L;\mathcal J;E;{\rm top}({\rho}))$\footnote{More precisely, the outer collaring of its thickening.}
such that
\begin{enumerate}\itemsep=0pt
\item[(1)]
${\rm ev}_0$ is strongly submersive with respect to this
CF-perturbation.
\item[(2)]
Those CF-perturbations are compatible with the identification
of the boundary as \eqref{form151414}, \eqref{form151515}.
\end{enumerate}
We use this system of CF-perturbations to define the integration along
the fiber ${\rm ev}_0 !$
in \eqref{form1516}.
We now define
\[
\mathfrak g_{k}(h_1,\dots,h_k) := \sum_E T^E \mathfrak g_{k,E}(h_1,\dots,h_k).
\]
This is well-defined by Proposition~\ref{prop1514}\,(4).

Stokes' formula (see \cite[Proposition 9.26]{foootech2} and \cite{fooonewbook}) and the
composition formula
(see \cite[Theorem~10.20]{foootech2} and \cite{fooonewbook})
imply that $\mathfrak g_{k}$ defines a filtered $A_{\infty}$
homomorphism.
In fact, \eqref{form151414} corresponds~to
\[
\mathfrak g_{k_1,E_1}\bigl(h_1,\dots,\mathfrak m^{J,L}_{k_1,E_2}(h_{i+1},\dots,h_{i+k_1}),\dots,h_k\bigr)
\]
and \eqref{form151515} corresponds to
\smash{$
\mathfrak m^{J',L}_{m,E_0}
\bigl(\mathfrak g_{k_1,E_1}\bigl(\vec h_1\bigr),\dots,\mathfrak g_{k_m,E_m}\bigl(\vec h_m\bigr)\bigr)$}.
Here $\vec h_1 = (h_1,\dots,h_{k_1})$,
$\vec h_2 = (h_{k_1+1},\dots,h_{k_1+k_2})$, etc.

The congruence $\mathfrak g_1 \equiv {\rm id}\mod \Lambda_+$
follows from Proposition~\ref{prop1514}\,(3).
The proof of Theorem~\ref{them1513} is complete.
\end{proof}

\begin{rem}
We omit the argument needed to take the homotopy inductive limit
$E \to \infty$, since it is similar to the other cases.
(This process is necessary since we work with only finitely
many moduli spaces consisting of moduli spaces of
objects with energy $< E_0$, to construct a~system
of Kuranishi structures and its CF-perturbations.)
\end{rem}
Theorem~\ref{thm154}\,(1) follows from
Theorem~\ref{them1513}.
We turn to the proof of Theorem~\ref{thm154}\,(2),~(3).

We use Yoneda embedding for the proof.
The objects $(L,b)$ and $(\Phi(L),\Phi_*(b))$
define filtered~$A_{\infty}$ right modules
$\mathfrak{Yon}(L,b)$ and $\mathfrak{Yon}(\Phi(L),\Phi_*(b))$
over $\mathfrak{Fukst}(\mathbb L)$, respectively.
By Lem\-ma~\ref{lema153333} and $A_{\infty}$ Yoneda lemma
(see Theorem~\ref{Yoneda}),
it suffices to prove the following.
\begin{enumerate}\itemsep=0pt
\item[(2)$'$]
The equivalence
$\mathfrak{Yon}(L,b)\sim_{\Lambda}\mathfrak{Yon}(\Phi(L),\Phi_*(b))$ holds as objects
of the functor category~$\mathcal{FUNC}(\mathfrak{Fukst}(\mathbb L)^{\rm op},\mathcal{CH})$.
\item[(3)$'$]
The Hofer distance
$d_{\rm Hof}(\mathfrak{Yon}(L,b),\mathfrak{Yon}(\Phi(L),\Phi_*(b)))$
is not greater than the Hofer distance between $\Phi$ and
the identity map.
\end{enumerate}

The proof of (2)$'$, (3)$'$ occupies the rest of this subsection.
We put $L' = \Phi(L)$
and
$M = \text{the disjoint union of elements of $\mathbb L$}$.
Note that $M = \bigl(\tilde M,i_M\bigr)$ is an immersed Lagrangian
submanifold of $X$.
We put
\begin{equation}\label{flor1588}
R := \tilde L \times_X \tilde M \qquad
R' := \tilde L' \times_X \tilde M.
\end{equation}
They are submanifolds
of $\tilde L \times \tilde M$, $\tilde L' \times \tilde M$, respectively.
We define
\begin{gather*}
CF(L,M) = \Omega(R)
\,\widehat{\otimes}\, \Lambda_0, \qquad\!
CF(L',M) = \Omega(R')
\,\widehat{\otimes}\, \Lambda_0, \qquad\!
CF(M) = \Omega\bigl(\tilde M \times_X \tilde M\bigr)
\,\widehat{\otimes}\, \Lambda_0.
\end{gather*}
We take a bounding cochain $b_M$ of $M$.
Then together with the bounding cochain $b$ of $L$ and~$\Phi_*(b)$
of $L'$ we obtain a right $\bigl(CF(M),\bigl\{\mathfrak m_k^{b_M}\bigr\}\bigr)$
module structures
\begin{gather}
\mathfrak n^L_k \colon\ CF(L,M) \otimes B_kCF(M) \to CF(L,M),\nonumber \\
\mathfrak n^{L'}_k\colon\ CF(L',M) \otimes B_kCF(M) \to CF(L',M).\label{form1510}
\end{gather}
They are nothing but $\mathfrak{Yon}(L,b)$ and $\mathfrak{Yon}(\Phi(L),\Phi_*(b))$.
We set
\[
CF(L,M)^{\Lambda} = CF(L,M) \otimes_{\Lambda_0} \Lambda,
\qquad
CF(L',M)^{\Lambda} = CF(L',M) \otimes_{\Lambda_0} \Lambda.
\]

We first review the moduli spaces we use to define
the filtered bi-module structure \eqref{form1510} on
$CF(L,M)$ (resp.\ $CF(L',M)$) over $CF(L)$-$CF(M)$
(resp.\ $CF(L')$-$CF(M)$).
See \cite[Sections~3.7.4 and 3.7.5]{fooobook} for detail.
We
consider the equation
\begin{equation}\label{eq1511}
\frac{\partial u}{\partial \tau} + J \frac{\partial u}{\partial t} = 0
\end{equation}
for a map $u \colon \R \times [0,1] \to X$
with boundary conditions:
\begin{enumerate}\itemsep=0pt
\item[(a)]
$u(\tau,0) \in M$.
\item[(b)]
$u(\tau,1) \in L$ (resp.\ $u(\tau,1) \in L'$).
\end{enumerate}
We consider $\vec z_0$, $\vec z_1$ such that
$\vec z_0 = (z_{0,1},\dots,z_{0,k_0})$, where
$z_{0,i} = (\tau_{0,i},0)$ with
$\tau_{0,1} < \dots < \tau_{0,k_0}$, and
$\vec z_1 = (z_{1,1},\dots,z_{1,k_1})$,
where
$z_{1,i} = (\tau_{1,i},0)$ with
$\tau_{1,1} > \dots > \tau_{1,k_1}$.\footnote{Note that we use
the counter clock-wise
ordering to enumerate the marked points.}
We also consider $\gamma_0 \colon \R \times \{0\} \setminus \vec z_0
\to \tilde M$, $\gamma_1 \colon \R \times \{1\} \setminus \vec z_1
\to \tilde L$, lifts of the restriction of $u$.
We assume an appropriate switching condition
similar to those appeared in Section~\ref{sec:HFIm}.
(See Definition~\ref{def3737}\,(5).)
We finally require\looseness=-1
\begin{enumerate}\itemsep=0pt
\item[(c)]
$\lim_{\tau \to \pm \infty} (\gamma_0(\tau),\gamma_1(\tau)) \in R$
(resp.
$\lim_{\tau \to \pm \infty} (\gamma_0(\tau),\gamma_1(\tau)) \in R'$).
\item[(d)]
$
\int_{\R \times [0,1]} u^* \omega
= E$.
\end{enumerate}
We consider such $(\vec z_0,\vec z_1;u;\gamma_0,\gamma_1)$
satisfying the above conditions and the moduli space
of such objects. We then take its quotient
by the $\R$ action induced by the translation of
the first factor of the source $\R \times [0,1]$.
We denote this space by \smash{$\mathring{\mathcal M}_{k_1,k_0}(L,M;E;J)$}
(resp.\ \smash{$\mathring{\mathcal M}_{k_1,k_0}(L',M;E;J)$}).
\begin{rem}
In equation \eqref{eq1511} (and
in other places of this subsection), we take $\R \times [0,1]$ as a~strip, while
in Section~\ref{sec:Kunneth} (and in other places of this paper) we took $[0,1] \times \R$.
In this subsection, we use $\R \times [0,1]$ for the sake of consistency with
\cite[Sections 3.7.4 and~3.7.5]{fooobook}.
In Section~\ref{sec:Kunneth}, we identified
$(t,\tau) \in [0,1] \times \R$ with $t + \sqrt{-1}\tau \in \C$ to define
complex structure.
(See the proof of Lemma~\ref{lem526526}.) Here we identify $(\tau,t)$ with
$\tau + \sqrt{-1}t \in \C$.
(Note that in Section~\ref{sec:Kunneth} the equation
corresponding to \eqref{eq1511} is
$\frac{\partial u}{\partial \tau} = J \frac{\partial u}{\partial t}$.)

In both cases, if we regard the first coordinate ($t$ in case of Section~\ref{sec:Kunneth}
and $\tau$ in case of this subsection) as the $x$-axis and
the second coordinate as the $y$-axis, then the above choice is consistent
with the standard conformal structure of the $xy$-plane.

We also remark that in Section~\ref{sec:Kunneth} we construct
{\it right} $CF(L_1)$ module and $L_1$ is assigned at the right, that is, $t=1$.
In this subsection, we construct {\it right} $CF(M)$ module
and $M$ is assigned at the bottom, that is, $t=0$.
This is consistent with our choice of orientation
and conformal structure of the domain.

\end{rem}
\begin{figure}[ht]
\centering
\includegraphics[scale=0.3]{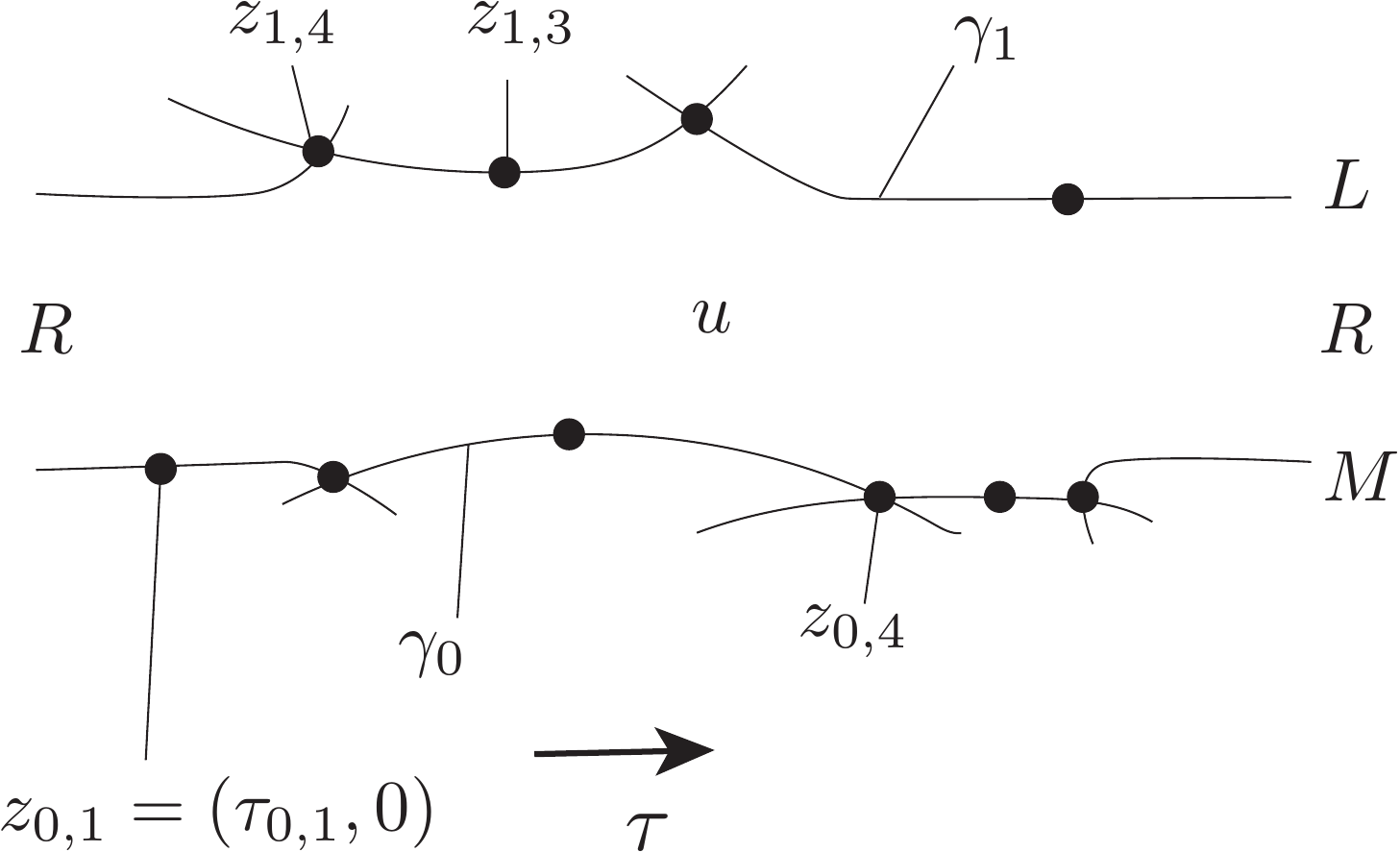}
\caption{Elements of ${\mathcal M}_{k_1,k_0}(L',M;E;J)$.}
\label{Figure15-1}
\end{figure}
\begin{prop}\label{prop1517}
The space
\smash{$\mathring{\mathcal M}_{k_1,k_0}(L,M;E;J)$}
$($resp.\ \smash{$\mathring{\mathcal M}_{k_1,k_0}(L',M;E;J)$}$)$
has a compactification
${\mathcal M}_{k_1,k_0}(L,M;E;J)$
$($resp.\ ${\mathcal M}_{k_1,k_0}(L',M;E;J))$
which is compact and Hausdorff, with respect to the
stable map topology.
They have Kuranishi structures with corners, and enjoy the following
properties:
\begin{enumerate}\itemsep=0pt
\item[$(1)$]
There exist evaluation maps
\begin{gather*}
\begin{split}
& {\rm ev} = \bigl({\rm ev}^{(1)},{\rm ev}^{(0)}\bigr)
=
\bigl(\bigl({\rm ev}^{(1)}_1,\dots,{\rm ev}^{(1)}_{k_1}\bigr),
\bigl({\rm ev}^{(0)}_1,\dots,{\rm ev}^{(0)}_{k_0}\bigr)\bigr)\colon \\
& \qquad
{\mathcal M}_{k_1,k_0}(L,M;E;J)
\to \bigl(\tilde L \times_X \tilde L\bigr)^{k_1}
\times \bigl(\tilde M \times_X \tilde M\bigr)^{k_0}
\end{split}
\end{gather*}
{\rm(}resp.
\begin{gather*}
{\rm ev} = \bigl({\rm ev}^{(1)},{\rm ev}^{(0)}\bigr)
=
\bigl(\bigl({\rm ev}^{(1)}_1,\dots,{\rm ev}^{(1)}_{k_1}\bigr),
\bigl({\rm ev}^{(0)}_1,\dots,{\rm ev}^{(0)}_{k_0}\bigr)\bigr)\colon\\
\qquad
{\mathcal M}_{k_1,k_0}(L',M;E;J)
\to \bigl(\tilde L' \times_X \tilde L'\bigr)^{k_1}
\times \bigl(\tilde M \times_X \tilde M\bigr)^{k_0}.)
\end{gather*}
These maps are evaluation maps at the marked points $\vec z_1$,
$\vec z_0$ and are underlying continuous maps
of strongly smooth maps.
\item[$(2)$]
There exists also evaluation maps at infinity
$
({\rm ev}_{-\infty},{\rm ev}_{+\infty})
\colon {\mathcal M}_{k_1,k_0}(L,M;E;J)
\to R\times R$,
{\rm(}resp.\
$
({\rm ev}_{-\infty},{\rm ev}_{+\infty})
\colon {\mathcal M}_{k_1,k_0}(L',M;E;J)
\to R' \times R'$.$)$
These maps are defined by the limit in item $(c)$ above and are underlying continuous maps
of strongly smooth maps.
The map~${\rm ev}_{+\infty}$ is weakly submersive.
\item[$(3)$]
The normalized boundary of
${\mathcal M}_{k_1,k_0}(L,M;E;J)$
$($resp.\ ${\mathcal M}_{k_1,k_0}(L',M;E;J))$
is the disjoint union of the following
three types of fiber products:
\begin{itemize}\itemsep=0pt
\item[$(I)$]
The fiber product
\begin{equation*}
\mathcal M_{k_{1,1}+1}(L;E_1;J)
{}_{{\rm ev}_0}\times_{{\rm ev}^{(1)}_i}
\mathcal M_{k_{1,2},k_0}(L,M;E_2;J),
\end{equation*}
where $k_{1,1} + k_{1,2} = k_1$, $E_1 + E_2 = E$ and
$i=1,\dots,k_{1,2}$
$($resp.\ the same except we replace $L$ by $L')$.
$($See Figure {\rm\ref{Figure15-2}.)}

\item[$(II)$]
The fiber product
\[
\mathcal M_{k_{0,1}+1}(M;E_1;J)
{}_{{\rm ev}_0}\times_{{\rm ev}^{(0)}_i}
\mathcal M_{k_1,k_{0,2}}(L,M;E_2;J),
\]
where $k_{0,1} + k_{0,2} = k_0$, $E_1 + E_2 = E$ and
$i=1,\dots,k_{0,2}$
$($resp.\ the same except we replace $L$ by $L')$.
$($See Figure {\rm\ref{Figure15-3}.)}

\item[$(III)$]
The fiber product
\begin{equation*}
\mathcal M_{k_{1,1},k_{0,1}}(L,M;E_1;J)
{}_{{\rm ev}_{+\infty}}\times_{{\rm ev}_{-\infty}}
\mathcal M_{k_{1,2},k_{0,2}}(L,M;E_2;J),
\end{equation*}
where $k_{0,1} + k_{0,2} = k_0$, $k_{1,1} + k_{1,2} = k_1$, $E_1 + E_2 = E$
$($resp.\ the same except we replace~$L$ by~$L')$.
$($See Figure {\rm\ref{Figure15-4}.)}
\end{itemize}
The evaluation maps are compatible with these identifications
of the boundary with fiber product.
\item[$(4)$]
There exists a principal ${\rm O}(1)$ bundle on $\tilde L \times_X \tilde L$,
$\tilde L \times_X \tilde L$, $\tilde L' \times_X \tilde L'$, $R$ and $R'$
and the trivialization of the orientation
bundle of ${\mathcal M}_{k_1,k_0}(L,M;E;J)$ tensored
with the pullbacks of those principal ${\rm O}(1)$ bundles.
These trivializations are compatible with the above
identification of the boundary.
\item[$(5)$]
The set of $E$ for which
${\mathcal M}_{k_1,k_0}(L,M;E;J)$
$($resp.\ ${\mathcal M}_{k_1,k_0}(L',M;E;J))$
is nonempty is discrete.
\end{enumerate}
\end{prop}

\begin{figure}[ht]\centering
\begin{tabular}{cc}
\begin{minipage}[t]{0.45\hsize}
\centering
\includegraphics[scale=0.33]{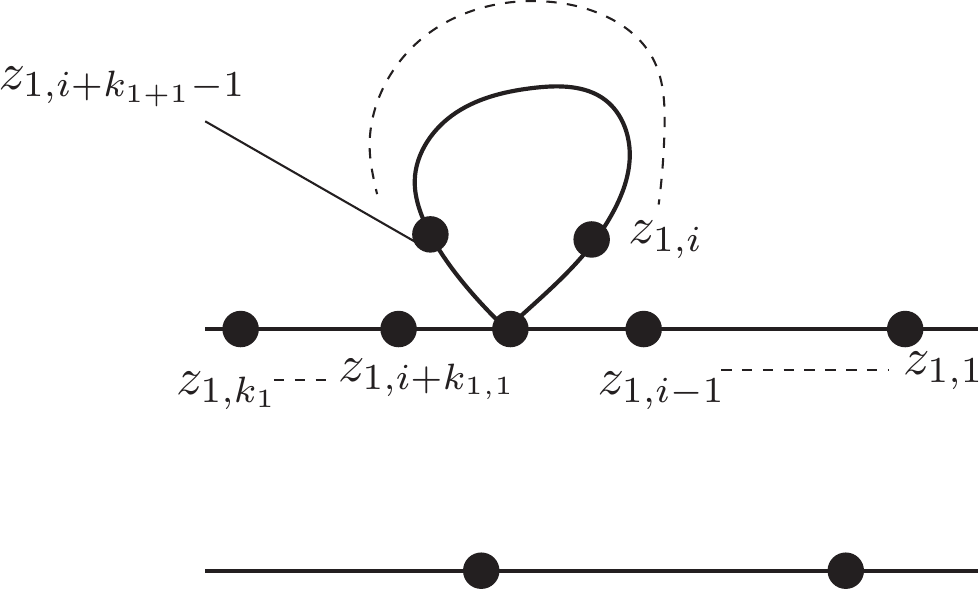}
\caption{Boundary of type I.}
\label{Figure15-2}
\end{minipage} &
\begin{minipage}[t]{0.45\hsize}
\centering
\includegraphics[scale=0.33]{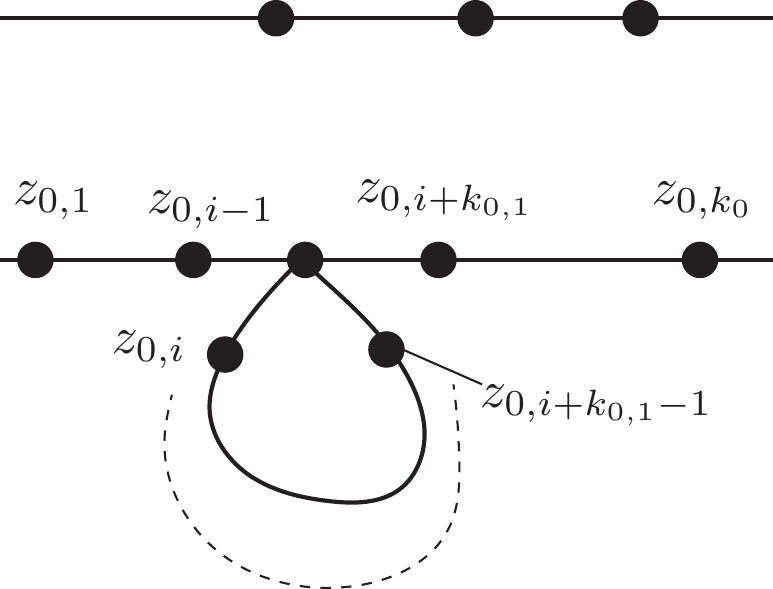}
\caption{Boundary of type II.}
\label{Figure15-3}
\end{minipage}
\end{tabular}
\end{figure}

\begin{figure}[ht]
\centering
\includegraphics[scale=0.25]{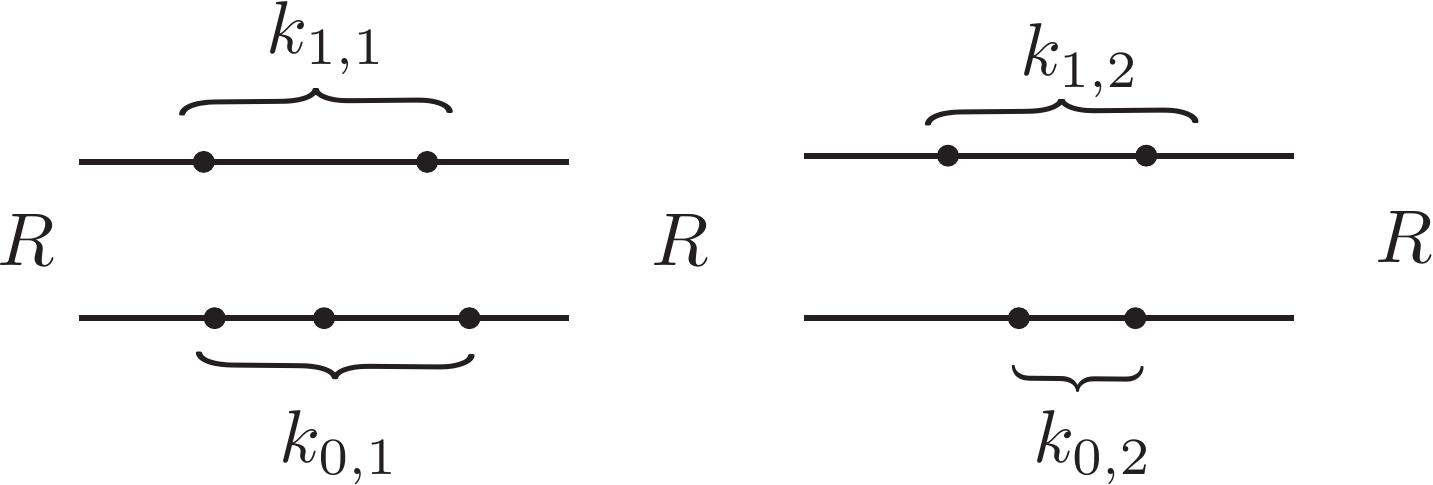}
\caption{Boundary of type III.}
\label{Figure15-4}
\end{figure}
The proof is now a routine. (See also \cite[Sections 3.7.4 and 3.7.5]{fooobook}, \cite{fooo:const2}
and Section~\ref{subsec:modpolygon} of this paper.)
We now define
\begin{gather*}
\mathfrak n^L_{k_1,k_0}\colon\ B_{k_0}CF(L)[1] \otimes CF(L,M) \otimes B_{k_1} CF(M)[1] \to CF(L,M), \\
\mathfrak n^{L'}_{k_1,k_0} \colon\ B_{k_0}CF(L')[1] \otimes CF(L',M) \otimes B_{k_1}CF(M)[1] \to CF(L',M),
\end{gather*}
by
\begin{gather}
\mathfrak n^L_{k_1,k_0}\bigl(h^{(1)}_1,\dots,h^{(1)}_{k_1};h;
h^{(0)}_1,\dots,h^{(0)}_{k_0}\bigr) \nonumber\\
\qquad:=
\sum_E T^E
{\rm ev}_{+\infty}!
\bigl(
\bigl({\rm ev}^{(1)}_1\bigr)^*h^{(1)}_1 \wedge
\bigl({\rm ev}^{(1)}_{k_1}\bigr)^*h^{(1)}_{k_1}
\wedge {\rm ev}_{-\infty}^*h\nonumber
 \\
\phantom{\qquad:=}{}
\wedge \bigl({\rm ev}^{(0)}_1\bigr)^*h^{(0)}_1 \wedge
\bigl({\rm ev}^{(0)}_{k_0}\bigr)^*h^{(0)}_{k_0}
; {\mathcal M}_{k_1,k_0}(L,M;E;J);\widehat{\mathfrak S}_{\varepsilon}
\bigr).\label{form151616}
\end{gather}
Here we take a system of CF-perturbations $\widehat{\mathfrak S}_{\varepsilon}$ on
${\mathcal M}_{k_1,k_0}(L,M;E;J)$ such that
${\rm ev}_{+\infty}$ is strongly submersive with respect to \smash{$\widehat{\mathfrak S}_{\varepsilon}$} and
that the CF-perturbations \smash{$\widehat{\mathfrak S}_{\varepsilon}$} are compatible with the
fiber product description of the boundaries in
Proposition~\ref{prop1517}\,(3). We use the CF-perturbation to define
the integration along the fiber ${\rm ev}_{+\infty}!$ in \eqref{form151616}.
(See \cite[Definitions~7.78 and~9.13]{foootech2} and \cite{fooonewbook}.)
The definition of $\mathfrak n^{L'}_{k_1,k_0}$ is similar.
We can show that these maps define structures of
filtered~$A_{\infty}$ bi-module by using
Stokes' formula (see \cite[Proposition 9.26]{foootech2} and \cite{fooonewbook}) and
the composition formula~(see \cite[Theorem 10.20]{foootech2} and \cite{fooonewbook}) together
with Proposition~\ref{prop1517}\,(3).

We now define the map \eqref{form1510} by
\[
\mathfrak n^L_{k}(y;x_1,\dots,x_k)
:=
\sum_{\ell,m_0,\dots,m_k} \mathfrak n^L_{\ell,k+\sum m_i}
\bigl(b^{\ell};y;b_M^{m_0}x_1 b_M^{m_1} \cdots b_M^{m_{k-1}} x_k b_M^{m_{k}}\bigr).
\]
Here and hereafter, for example, $b^2_M x b_M$ means $b_M \otimes b_M \otimes x \otimes b_M$.
The $A_{\infty}$ relation of $\mathfrak n^L_{k_1,k_0}$ and
the fact that $b$, $b_M$ are bounding cochains
imply that \smash{$\mathfrak n^L_{k}$} defines a (strict and unital) filtered right $A_{\infty}$
$\bigl(CF(M);\bigl\{\mathfrak m_k^{b_M}\bigr\}\bigr)$ module structure on $CF(L,M)$.

We can define $\mathfrak n^{L'}_{k}$ in the same way.

We next describe the moduli spaces which we use to define
a filtered right $A_{\infty}$ module homomorphism $CF(L,M)^{\Lambda} \to CF(L',M)^{\Lambda}$.
We follow \cite[Section 5.3.1]{fooobook} with modification given in \cite{fooo:polydisk}.
We will use a two parameter family of almost complex structures $\mathcal{JJ}$,
which is defined in Definition~\ref{defn1519}, to define the moduli space
appearing in
Proposition~\ref{prop1518}.

\begin{figure}[ht]\centering
\begin{tabular}{cc}
\begin{minipage}[t]{0.45\hsize}
\centering
\includegraphics[scale=0.3]{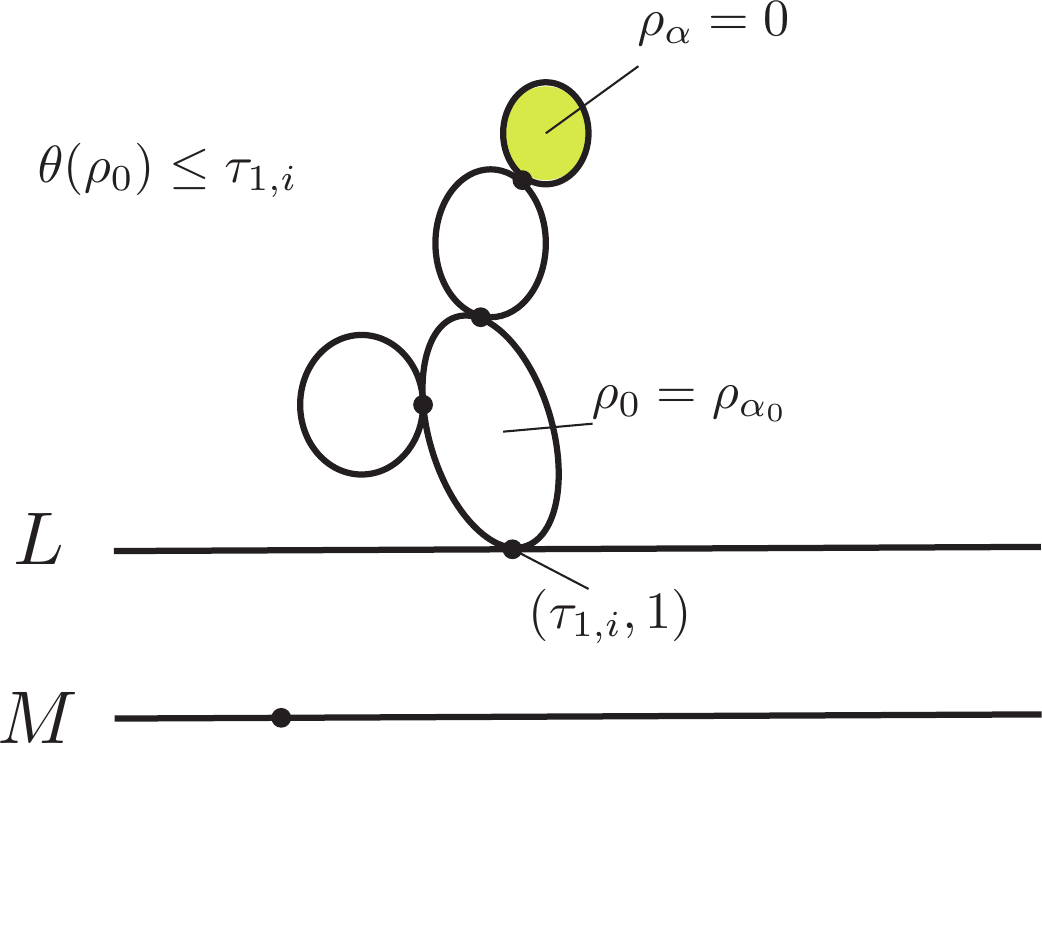}
\caption{Boundary of type I.}
\label{Figure15-5}
\end{minipage} &
\begin{minipage}[t]{0.45\hsize}
\centering
\includegraphics[scale=0.3]{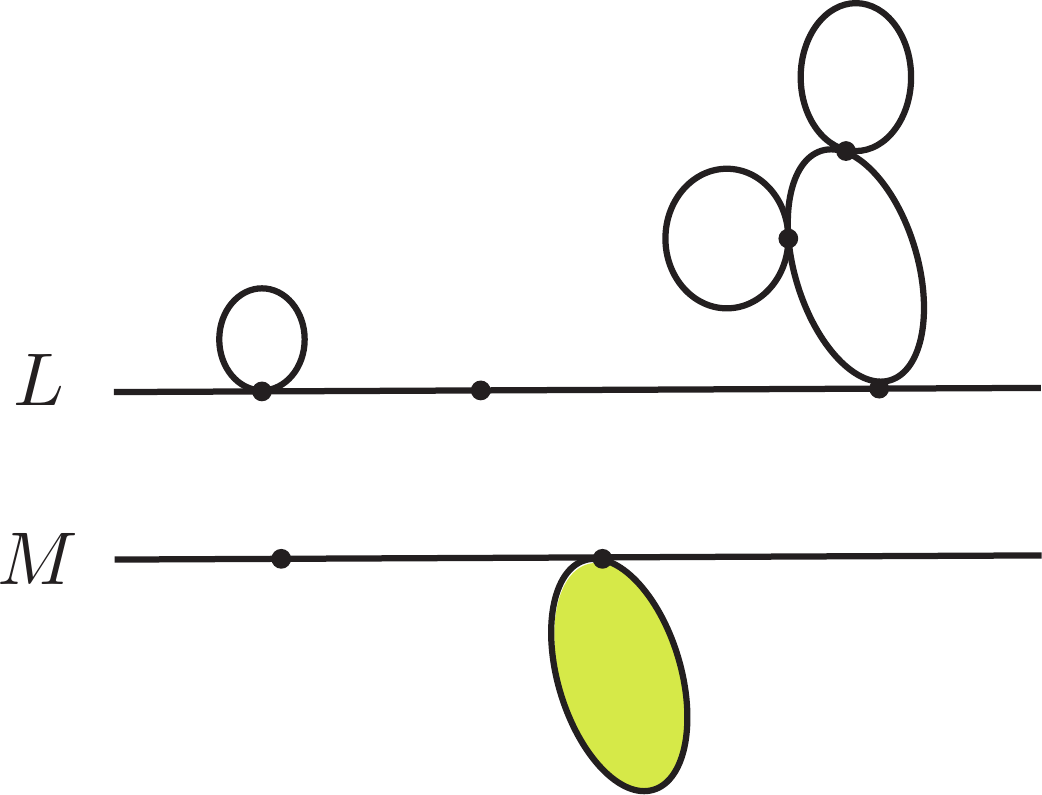}
\caption{Boundary of type II.}
\label{Figure15-6}
\end{minipage}
\end{tabular}
\end{figure}

\begin{figure}[ht]
\centering
\includegraphics[scale=0.4]{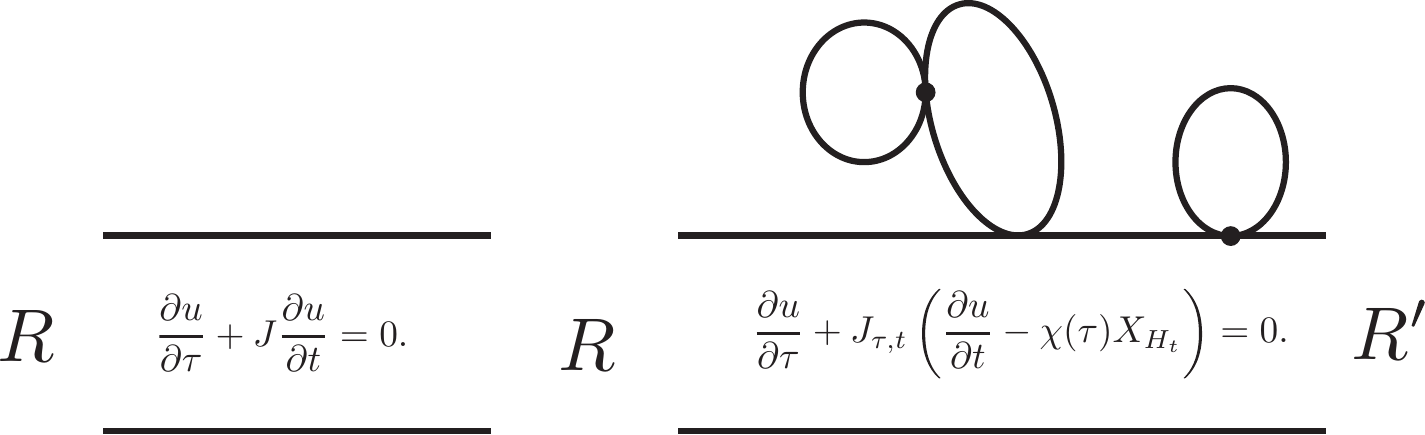}
\caption{Boundary of type III.}
\label{Figure15-7}
\end{figure}

\begin{figure}[ht]
\centering
\includegraphics[scale=0.4]{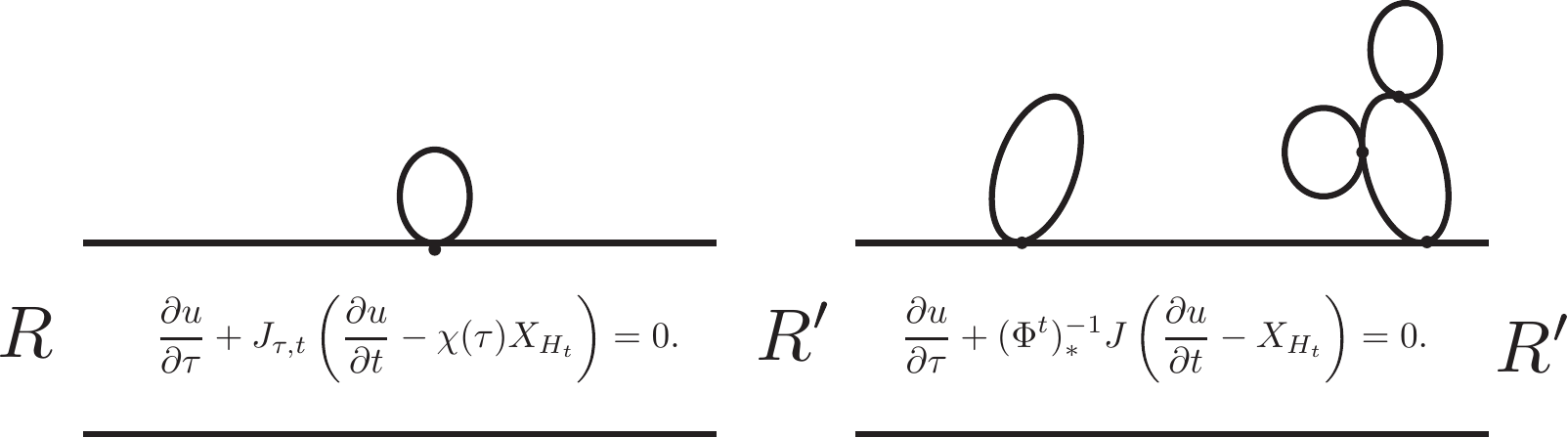}
\caption{Boundary of type IV.}
\label{Figure15-8}
\end{figure}

\begin{prop}\label{prop1518}
There exists a system of compact Hausdorff spaces\index[syindex]{M1k1k0LLprime@${\mathcal M}_{k_1,k_0}(L,L';M;E;\mathcal{JJ};{\rm top}(\rho))$}
\[
{\mathcal M}_{k_1,k_0}(L,L';M;E;\mathcal{JJ};{\rm top}(\rho))
\]
with the
following properties.
The spaces ${\mathcal M}_{k_1,k_0}(L,L';M;E;\mathcal{JJ};{\rm top}(\rho))$ carry
Kuranishi structures with corners.
\begin{enumerate}\itemsep=0pt
\item[$(1)$]
There exist evaluation maps
\begin{gather*}
{\rm ev} = \bigl({\rm ev}^{(1)},{\rm ev}^{(0)}\bigr)
=
\bigl(\bigl({\rm ev}^{(1)}_1,\dots,{\rm ev}^{(1)}_{k_1}\bigr),
\bigl({\rm ev}^{(0)}_1,\dots,{\rm ev}^{(0)}_{k_0}\bigr)\bigr)\colon \\
\qquad
{\mathcal M}_{k_1,k_0}(L,L';M;E;\mathcal{JJ};{\rm top}(\rho))
\to \bigl(\tilde L \times_X \tilde L\bigr)^{k_1}
\times \bigl(\tilde M \times_X \tilde M\bigr)^{k_0}.
\end{gather*}
These maps are underlying continuous maps
of strongly smooth maps.
\item[$(2)$]
There exist also evaluation maps at infinity
\[
({\rm ev}_{-\infty},{\rm ev}_{+\infty})
\colon\ {\mathcal M}_{k_1,k_0}(L,L';M;E;\mathcal{JJ};{\rm top}(\rho))
\to R \times R'.
\]
These maps are underlying continuous maps
of strongly smooth maps.
${\rm ev}_{+\infty}$ is weakly submersive.
$R$ and $R'$
are defined in \eqref{flor1588}.
\item[$(3)$]
The normalized boundary of
${\mathcal M}_{k_1,k_0}(L,L';M;E;\mathcal{JJ};{\rm top}(\rho))$
is the disjoint union of the following
four types of fiber products:
 \begin{itemize}\itemsep=0pt
\item[$(I)$]
The fiber product
\begin{equation}\label{form15141411}
\mathcal M_{k_{1,1}+1}(L;E_1;J)
{}_{{\rm ev}_0}\times_{{\rm ev}^{(1)}_i}
{\mathcal M}_{k_{1,2},k_0}(L,L';M;E_2;\mathcal{JJ};{\rm top}(\rho)),
\end{equation}
where $k_{1,1} + k_{1,2} = k_1$, $E_1 + E_2 = E$ and
$i=1,\dots,k_{1,2}$
$($see Figure {\rm\ref{Figure15-5})}.
\item[$(II)$]
The fiber product
\[
\mathcal M_{k_{0,1}+1}(M;E_1;J)
{}_{{\rm ev}_0}\times_{{\rm ev}^{(0)}_i}
{\mathcal M}_{k_1,k_{0,2}}(L,L';M;E_2;\mathcal{JJ};{\rm top}(\rho)),
\]
where $k_{0,1} + k_{0,2} = k_0$, $E_1 + E_2 = E$ and
$i=1,\dots,k_{0,2}$
$($see Figure {\rm\ref{Figure15-6})}.
\item[$(III)$]
The fiber product
\begin{equation}\label{form1514142222}
\mathcal M_{k_{0,1},k_{1,1}}(L,M;E_1;J)
{}_{{\rm ev}_{+\infty}}\times_{{\rm ev}_{-\infty}}
{\mathcal M}_{k_{1,2},k_{0,2}}(L,L';M;E_2;\mathcal{JJ};{\rm top}(\rho)),
\end{equation}
where $k_{0,1} + k_{0,2} = k_0$, $k_{1,1} + k_{1,2} = k_1$, $E_1 + E_2 = E$
$($see Figure {\rm\ref{Figure15-7})}.
\item[$(IV)$]
The fiber product of
\begin{equation}\label{form1520}
{\mathcal M}_{k_{1,1},k_{0,1}}(L,L';M;E;\mathcal{JJ};{\rm top}(\rho))
\end{equation}
and
\begin{gather}
\prod_{j=1}^{\ell}
\mathcal M_{m_j+1}(L;\mathcal J;E_{2,j};{\rm top}({\rho}))\nonumber\\
{}_{({\rm ev}_0,\dots,{\rm ev}_0)}\times_{({\rm ev}^{(1)}_1,\dots,{\rm ev}^{(1)}_{1,2})}
\mathcal M_{k_{0,2},k_{1,2}}(L',M;E_{2,0};J),\label{form1521}
\end{gather}
where $k_{1,1} + \sum_{j=1}^{\ell}m_j = k_1$, $k_{0,1} + k_{0,2} = k_0$,
$E_1 + \sum_{j=0}^{k_{1,2}}E_{2,j} = E$.
We use ${\rm ev}_{+\infty} \colon \eqref{form1520} \to R'$
and ${\rm ev}_{-\infty} \colon \eqref{form1521} \to R'$
to take fiber product between \eqref{form1520} and \eqref{form1521}
$($see Figure {\rm\ref{Figure15-8})}.
\end{itemize}

The evaluation maps are compatible with this identifications.
\item[$(4)$]
There exists a principal ${\rm O}(1)$ bundle on $\tilde L \times_X \tilde L$,
$\tilde L \times_X \tilde L$, $\tilde L' \times_X \tilde L'$, $R$ and $R'$
and the trivializations of the orientation
bundle of ${\mathcal M}_{k_1,k_0}(L,L';M;E;\mathcal{JJ};{\rm top}(\rho))$ tensored
with the pullbacks of those principal ${\rm O}(1)$ bundles.
These trivializations are compatible with the above
identification of the boundary.
\item[$(5)$]
The set of $E$ for which
${\mathcal M}_{k_1,k_0}(L,L';M;E;\mathcal{JJ};{\rm top}(\rho))$
is nonempty is discrete.
\end{enumerate}

\end{prop}
The elements of the moduli spaces corresponding to the
boundaries of types (I)(II)(III)(IV) are depicted in the
Figures \ref{Figure15-5}--\ref{Figure15-8} below. The explanation of the figures will be given
during the proof of Proposition~\ref{prop1518}.

\begin{proof}
We take one parameter family of Hamiltonian diffeomorphisms $\Phi^{\rho}$
such that $\Phi^0 = {\rm id}$ and~${\Phi^1 = \Phi}$.
When $H \colon X \times [0,1] \to \R$ is the time dependent family of Hamiltonians
generating~$\Phi$, we take $\Phi^{\rho}$ so that{\samepage
\begin{equation}\label{phirho}
\frac{d\Phi^{\rho}}{d{\rho}} = X_{H_{\rho}} \circ \Phi^{\rho},
\qquad \Phi^0 = {\rm id},
\end{equation}
where $H_{\rho}(x) = H(x,{\rho})$ and $X_{H_{\rho}}$ is
the Hamiltonian vector field associated to $H_{\rho}$.}

We replace $H$ by $cH$ and obtain
one parameter family of Hamiltonian diffeomorphisms, which
we denote by $\Phi_{cH}^{\rho}$.

We take a non-decreasing function
$
\chi \colon \R \to [0,1]
$
such that
\begin{enumerate}\itemsep=0pt
\item[(1)]
$\chi(\tau) = 0$ for sufficiently small $\tau$.
\item[(2)]
$\chi(\tau) = 1$ for sufficiently large $\tau$.
\end{enumerate}
\begin{defn}\label{defn1519}
We take a two parameter family of
complex structures $\mathcal{JJ} = \{J_{\tau,t}\}$
with the following properties:
\begin{enumerate}\itemsep=0pt
\item[(1)]
There exists $A>0$ such that $J_{\tau,t} =
J$
if $\tau< -A$.
\item[(2)]
$
J_{\tau,t} = \bigl(\Phi^{t}\bigr)^{-1}_* J
$
if $\tau > + A$.
\item[(3)]
We denote by $\Phi_{c H}^{\rho}$ the one parameter family
of Hamilton diffeomorphisms generated by
the time dependent Hamiltonian $cH \colon X \times [0,1] \to \R$.
Then
\smash{$
J_{\tau,1} = \bigl(\Phi_{\chi(\tau) H}^{1}\bigr)^{-1}_* J
$}
if $\tau > 0$.
\item[(4)]
$J_{\tau,0} = J$ for any $\tau$.
\end{enumerate}
\end{defn}
We take the one parameter family of almost
complex structures $\mathcal J = \bigl\{J^{(\rho)}\bigr\}$
which we used to prove Proposition~\ref{prop1514}
as follows.
We take and fix an order preserving {\it diffeomorphism} $
\theta \colon (0,1) \to \R$.
We then put
\begin{equation}\label{form15261526}
J^{(\rho)} =\bigl (\Phi_{\chi(\theta(\rho)) H}^{1}\bigr)^{-1}_*J.
\end{equation}

We consider maps
\begin{equation}\label{form152424}
u \colon\ \R \times [0,1] \to X,
\end{equation}
which satisfy the following conditions.
\begin{conds}\label{conds152020}
\quad
\begin{enumerate}\itemsep=0pt
\item[(1)]
$u$ satisfies the equation
\begin{equation}\label{eq1520}
\frac{\partial u}{\partial \tau} + J_{\tau,t}
\left(\frac{\partial u}{\partial t}
- \chi(\tau) X_{H_t}\right) = 0.
\end{equation}
Here $H$ is the time dependent Hamiltonian as in \eqref{phirho}.
\item[(2)]
$u(\tau,0) \in M$.
\item[(3)]
$
u(\tau,1) \in L$.
\end{enumerate}

\end{conds}
\begin{rem}
In \cite[Section 5.3.1]{fooobook}, we used
pseudo-holomorphic curve equation
(without Hamiltonian term) with a moving boundary condition,
(which becomes the condition $u(\tau,1) \in \Phi_{\chi_+(\tau)}(L)$
in our situation).
(See \cite[equations~(5.3.18.1) and (5.3.18.2)]{fooobook}.)
Here we use the equation \eqref{eq1520} (which has a Hamiltonian
term) and the boundary conditions are given by fixed Lagrangian
submanifolds $M$ and $L$.
The way taken here is the same as \cite{fooo:polydisk}.
(See \cite[equations~(3.3) and (3.4)]{fooo:polydisk}.)
The relation between these two formulations are explained in \cite[Section 4]{fooo:polydisk}.
We use the current formulation since then we can obtain energy
estimate (see Lemma~\ref{lem524})
easier.
\end{rem}
\begin{defn}\label{definition152666}
We define \smash{$\overset{\ \text{\tiny $\circ\circ$}}{\mathcal M}_{k_1,k_0}(L,L';M;E;\mathcal{JJ};{\rm top}(\rho))$}
as the set of
objects
\[((\R \times [0,1];\vec z_0,\vec z_1);u;\gamma;\vec \rho)
\]
 such that
\begin{enumerate}\itemsep=0pt
\item[(1)]
The map $u$ is as in \eqref{form152424} and satisfying Conditions \ref{conds152020}.
\item[(2)]
$\vec z_0$ (resp.\ $\vec z_1$) is a $k_0$ (resp.\ $k_1$) tuple of points, that is,
$\vec z_0 = (z_{0,1},\dots,z_{0,k_0})$, where
$z_{0,i} = (\tau_{0,i},0)$ with
$\tau_{0,1} < \dots < \tau_{0,k_0}$
(resp.\ $\vec z_1 = (z_{1,1},\dots,z_{1,k_1})$,
where
$z_{1,i} = (\tau_{1,i},0)$ with
$\tau_{1,1} > \dots > \tau_{1,k_1}$).
\item[(3)]
The maps $\gamma_0 \colon (\R \times \{0\}) \setminus \vec z_0
\to \tilde M$, $\gamma_1 \colon (\R \times \{1\}) \setminus \vec z_1
\to \tilde L$ are lifts of the restrictions of~$u$.
Namely,
$
u(\tau,1)
= i_{L}(\gamma_1(\tau)))$, $
u(\tau,0) = i_{M}(\gamma_0(\tau))$.
We assume an appropriate switching condition
similar to those appeared in Section~\ref{sec:HFIm}.
(See Definition~\ref{def3737}\,(5).)

\item[(4)]
$\vec \rho = (\rho_{1},\dots,\rho_{k_1})$,
where $\rho_{i}$ are real numbers.
We require
\begin{equation}\label{form1526}
\theta(\rho_{i}) \le \tau_{1,i}.
\end{equation}
\item[(5)]
We require
\[
E = \int_{\R \times [0,1]} u^*\omega
+ \lim_{\tau\to+\infty}\int_{[0,1]} H(t,u(\tau,t))  d t.
\]
\end{enumerate}
\end{defn}

\begin{figure}[ht]
\centering
\includegraphics[scale=0.4]{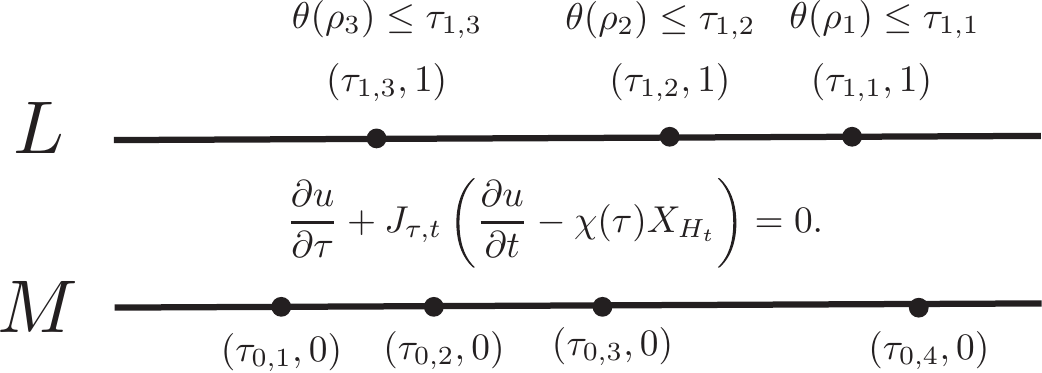}
\caption{An element of $\overset{\ \text{\tiny $\circ\circ$}}{\mathcal M}_{k_1,k_0}(L,L';M;E;\mathcal{JJ};{\rm top}(\rho))$.}
\label{Figure15-22}
\end{figure}

We define evaluation maps
\[
{\rm ev}_{1,i} \colon\ \overset{\ \text{\tiny $\circ\circ$}}{\mathcal M}_{k_1,k_0}(L,L';M;E;\mathcal{JJ};{\rm top}(\rho))
\to \tilde L\times_X \times \tilde L
\]
(resp.
\[
{\rm ev}_{0,i}\colon\ \overset{\ \text{\tiny $\circ\circ$}}{\mathcal M}_{k_1,k_0}(L,L';M;E;\mathcal{JJ};{\rm top}(\rho))
\to \tilde M\times_X \times \tilde M),
\]
as the evaluation maps at the marked points $z_{0,i}$ (resp.\ $z_{1,i}$)
using the switching condition in the same way as \eqref{form3188}.
Namely,
\[
{\rm ev}_{1,i}((\R \times [0,1];\vec z_0,\vec z_1);u;\gamma;
\vec \rho)
=(\lim_{\tau \downarrow \tau_{1,i}}\gamma_1(\tau),\lim_{\tau \uparrow \tau_{1,i}}\gamma_1(\tau)),
\]
and
\[
{\rm ev}_{0,i}((\R \times [0,1];\vec z_0,\vec z_1);u;\gamma;
\vec \rho)
=(\lim_{\tau \uparrow \tau_{0,i}}\gamma_0(\tau),\lim_{\tau \downarrow \tau_{0,i}}\gamma_0(\tau)),
\]
respectively.

We also define the evaluation maps at infinity
\[
{\rm ev}_{-\infty} \colon\ \overset{\ \text{\tiny $\circ\circ$}}{\mathcal M}_{k_1,k_0}(L,L';M;E;\mathcal{JJ};{\rm top}(\rho))
\to R
\]
(resp.
\[
{\rm ev}_{+\infty} \colon\ \overset{\ \text{\tiny $\circ\circ$}}{\mathcal M}_{k_1,k_0}(L,L';M;E;\mathcal{JJ};{\rm top}(\rho))
\to R').
\]
Here we define ${\rm ev}_{+\infty}$ by
\begin{equation}\label{form1534}
{\rm ev}_{+\infty}((\R \times [0,1];\vec z_1;u;\gamma;
\vec \rho)
=
\lim_{\tau \to +\infty}(\gamma_0(\tau),\Phi(\gamma_1(\tau))).
\end{equation}
Note that the limit $\ell(t) = \lim_{\tau\to + \infty}u(\tau,t)$
satisfies
$
\frac{ d \ell}{ d t} = X_{H_t}\circ \ell$.
This is a consequence of \eqref{conds152020} and
$\lim_{\tau\to +\infty} \frac{\partial u }{\partial \tau} = 0$.
Therefore,
$
\lim_{\tau\to + \infty}\Phi(u(\tau,0))
=
\lim_{\tau\to + \infty}u(\tau,1)$.
Hence the right-hand side of \eqref{form1534} is an element of
$R'$.

Finally, we define
\[
{\rm ev}_{i}^{\rm deti} \colon\ \overset{\ \text{\tiny $\circ\circ$}}{\mathcal M}_{k_1,k_0}(L,L';M;E;\mathcal{JJ};{\rm top}(\rho))
\to [0,1]
\]
\index[syindex]{evdetii@${\rm ev}_{i}^{\rm deti}$}
by
$
{\rm ev}_{i}^{\rm deti}
((\R \times [0,1];\vec z_0,\vec z_1);u;\gamma;
\vec \rho)
=
\rho_{i}$.
(Here ${\rm deti}$ stands for `time with delay'.)
\begin{defn}\label{1522222}
\smash{$\mathring{\mathcal M}_{k_1,k_0}(L,L';M;E;\mathcal{JJ};{\rm top}(\rho))$}
is a union of fiber products of
\begin{equation}\label{153434}
\prod_{j=1}^{k'_1}
\mathcal M_{m_j+1}(L;\mathcal J;E_{1,j};{\rm top}({\rho}))
\end{equation}
and
\begin{equation}\label{153435}
\overset{\ \text{\tiny $\circ\circ$}}{\mathcal M}_{k'_1,k_0}(L,L';M;E_2;\mathcal{JJ};{\rm top}(\rho)),
\end{equation}
where the union is taken over $k'_1$, $\{m_j\}$, $\{E_{1,j}\}$, $E_2$ with
$\sum_{j=1}^{k'_1}m_j = k_1$, $\sum_{j=1}^{k'_1}E_{1,j} + E_2= E$.

The fiber product is taken over
\smash{$\prod_{j=1}^{k'_1} \bigl(\bigl(\tilde L \times_X \tilde L\bigr) \times \R\bigr)$}.
We use the map
$\eqref{153434} \to\smash{ \prod_{j=1}^{k'_1} }\bigl(\bigl(\tilde L\allowbreak \times_X \tilde L\bigr) \times \R\bigr)$
which is $(({\rm ev}_0,\rho_0),\dots,({\rm ev}_0,\rho_0))$
and the map
\smash{$\eqref{153435} \to \prod_{j=1}^{k'_1} \bigl(\bigl(\tilde L \times_X \tilde L\bigr) \times \R\bigr)$}
which is $\bigl(\bigl({\rm ev}_{1,1},{\rm ev}_{1,1}^{\rm deti}\bigr),\dots,\bigl({\rm ev}_{1,k'_1},{\rm ev}_{1,k'_1}^{\rm deti}\bigr)\bigr)$ to define the fiber product.
(Note $\rho_0$ is defined by~\eqref{defrho0}.)
\end{defn}
Figure~\ref{Figure15-23} below depicts an element of \smash{$\mathring{\mathcal M}_{k_1,k_0}(L,L';M;E;\mathcal{JJ};{\rm top}(\rho))$}. It is a map to~$X$ from the domain which is a union of a strip $\R \times [0,1]$
and trees of disks attached at $t=1$. It is pseudo-holomorphic with respect to the almost complex structure of $X$
which depends on the components of the domain.
The almost complex structure we use is $J^{(\rho_i)}$ on the disk components
depicted in Figure~\ref{Figure15-23}. Note that $J^{(\rho)}$ is defined in~\eqref{form15261526}.
\begin{figure}[ht]
\centering
\includegraphics[scale=0.48]{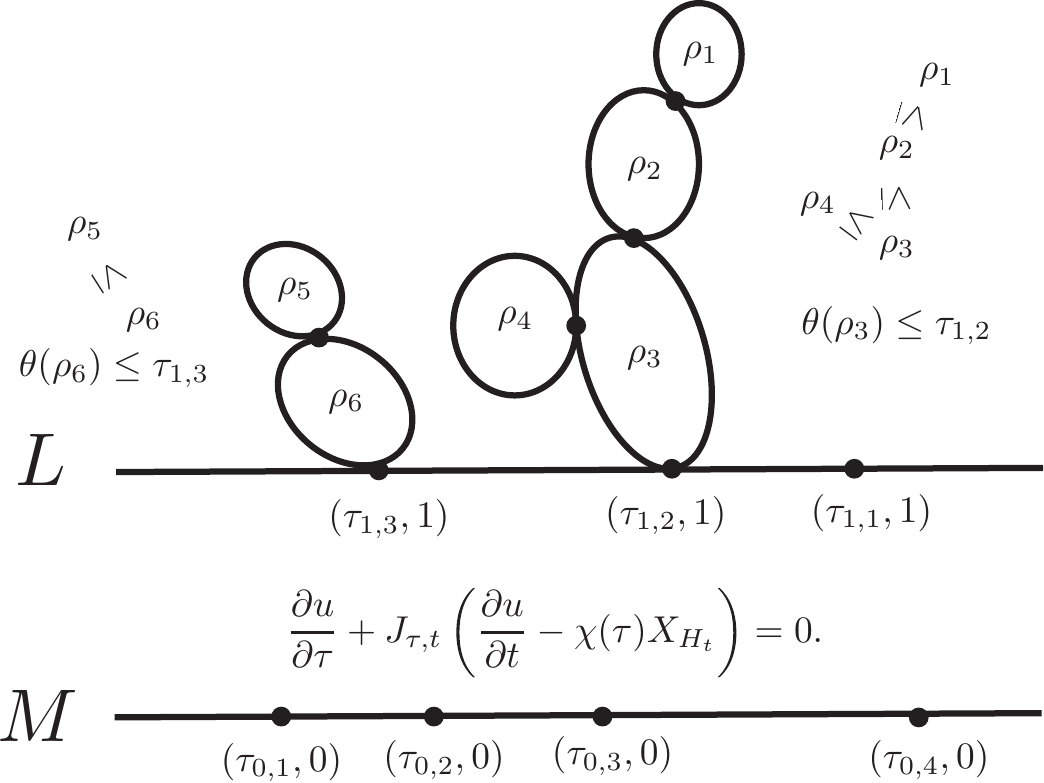}
\caption{An element of $\mathring{\mathcal M}_{k_1,k_0}(L,L';M;E;\mathcal{JJ};{\rm top}(\rho))$.}
\label{Figure15-23}
\end{figure}

We remark that all the fiber products of $\eqref{153434}$ and $\eqref{153435}$ have the same virtual dimension
that is independent of $k'_1$, $\{m_j\}$, $\{E_{1,j}\}$, $E_2$
but depends only on the total homology class of the map, and $k_1$, $k_2$.

We also remark that the union includes the case when
$m_j = 1$ and $E_{1,j} = 0$ for all $j$.
In this case, the fiber product of $\eqref{153434}$ and $\eqref{153435}$
is nothing but \smash{$\overset{\ \text{\tiny $\circ\circ$}}{\mathcal M}_{k_1,k_0}(L,L';M;E;\mathcal{JJ};{\rm top}(\rho))$}.

\begin{rem}
We remark that the space \smash{$\mathring{\mathcal M}_{k_1,k_0}(L,L';M;E;\mathcal{JJ};{\rm top}(\rho))$}
contains several components with the same virtual
dimension as the space \smash{$\raisebox{-0.5pt}{$\overset{\ \text{\tiny $\circ\circ$}}{\mathcal M}_{k_1,k_0}$}(L,L';M;E;\mathcal{JJ};{\rm top}(\rho))$}. So,
even in the case when all the elements are Fredholm regular, the subset
$\smash{\raisebox{-1.0pt}{$\overset{\ \text{\tiny $\circ\circ$}}{\mathcal M}_{k_1,k_0}$}(L,L';M;E};\allowbreak \mathcal{JJ};{\rm top}(\rho))$ may not
be a dense subset of the moduli space
\smash{$\mathring{\mathcal M}_{k_1,k_0}(L,L';M;E;\mathcal{JJ};{\rm top}(\rho))$}.
\end{rem}

The moduli space
${\mathcal M}_{k_1,k_0}(L,L';M;E;\mathcal{JJ};{\rm top}(\rho))$
is the stable map compactification of \smash{$\mathring{\mathcal M}_{k_1,k_0}(L,L';M;E;\mathcal{JJ};{\rm top}(\rho))$}.
The compactification is obtained by
adding the following:
\begin{enumerate}\itemsep=0pt\setlength{\leftskip}{0.50cm}
\item[(Bub.1)]
We include the case when the source curve has a sphere bubble.
\item[(Bub.2)]
We include the case when the source curve has a disk bubble
at $t=0$.
(The disk bubble at $t=1$ is already included when we take
the fiber product in Definition~\ref{1522222}.)
\item[(Bub.3)]
We include the case when the source curve splits into several
pieces in the $\tau$-direction.
The cases when it splits into two pieces are depicted in
Figures \ref{Figure15-7}, \ref{Figure15-8}.
\end{enumerate}
The detail of this stable map compactification is written in
\cite[Section 5.3.1]{fooobook} and is now a~routine. So we omit it here.
(We remark that all the components corresponding to one of
(Bub.1), (Bub.2), (Bub.3) have codimension $\ge 1$.)

We now study the boundary of the moduli space
${\mathcal M}_{k_1,k_0}(L,L';M;E;\mathcal{JJ};{\rm top}(\rho))$.

Case (I) in Proposition~\ref{prop1518}\,(3) occurs at the point
$((\xi_j),\xi_0)$
where one of the factors $\xi_j$ of~\eqref{153434} lies in
the boundary point of that factor. This corresponds to the case
when some $\rho_{\alpha_i}$ is $0$, where~${\xi_j = ((\Sigma,\vec z),u,\gamma,\{\rho_{\alpha_i}\})}$.
This boundary corresponds to Proposition~\ref{prop1514}\,(2), \eqref{form151414}.
Therefore, this case is described by the fiber product
\eqref{form15141411}. See Figure~\ref{Figure15-5}.

Note that the boundary which corresponds to
Proposition~\ref{prop1514}\,(2), \eqref{form151515} does not appear here.
In fact, \eqref{form1526} implies
$\theta(\rho_{i}) \le \tau_{1,i}$.
Here $\rho_i, \tau_{1,i}$ are parts of the data of $\xi_0$.
Moreover, by the definition of the fiber product appearing in Definition~\ref{1522222},
we find $\rho_0(\xi_i) = \rho_{i}$.
Therefore, since $\theta$ is a diffeomorphism
$\rho_0(\xi_i) = 1$ occurs only in the limit
which we discuss in Case (IV).

Case (II) in Proposition~\ref{prop1518}\,(3) occurs
when a disk bubble occurs at $t=0$, that is, (Bub.2).
See Figure~\ref{Figure15-6}.

We remark that the situations of the bubbles at $t=0$ and $t=1$ are
different. This is because the boundary conditions are different.

Case (III) in Proposition~\ref{prop1518}\,(3) occurs
when the domain splits into two pieces one of which
moves to the direction $\tau \to -\infty$.

Note that in this limit some of the trees of disk bubbles
at $t=1$ may be attached to the piece which
moves to the direction $\tau \to -\infty$.
If such a tree of disk bubbles corresponds to $\xi_i$
(that is, one of the factors of the fiber product
\eqref{153434} and
the root of such piece is $(\tau_i,1)$, then $\tau \to -\infty$.
Therefore, ${\rho_0(\xi_i) = 0}$.
(In fact, $\rho_0(\xi_i) \le {\rm ev}_{1,i}^{\rm deti}(\xi)$,
where $\xi$ is an element of the factor~\eqref{153435}.)

Therefore, by definition this case is
described by the fiber product \eqref{form1514142222}.
See Figure~\ref{Figure15-7}.

Case (IV) in Proposition~\ref{prop1518}\,(3) occurs
when
the domain splits into two pieces one of which
moves to the direction $\tau \to +\infty$.

Note that the piece which moves to the direction $\tau \to +\infty$
consists of a map from a~strip $\R \times [0,1]$ plus a union of trees
of disk bubbles.
The map $u_{\infty} \colon \R \times [0,1] \to X$,
which is a part of this map, satisfies the equation
\begin{equation}\label{153915}
\frac{\partial u_{\infty}}{\partial \tau} + \bigl(\bigl(\Phi^t\bigr)^{-1}_*J\bigr)
\left(\frac{\partial u_{\infty}}{\partial t}
- X_{H_t}\right) = 0,
\end{equation}
together with the boundary condition $u_{\infty}(\tau,0) \in M$
and $u_{\infty}(\tau,1) \in L$.
We put
$
v(\tau,t) = \Phi^t(u_{\infty}(\tau,t))$.
Then \eqref{153915} is equivalent to
$
\frac{\partial v}{\partial \tau} + J
\frac{\partial v}{\partial t} = 0$.
The corresponding boundary condition for $v$ is
$v(\tau,0) \in M$
and $v(\tau,1) \in L'$.
This is the equation and the boundary condition which we used in the definition of
$\mathcal M_{k_{0},k_{1}}(L',M;E;J)$.

The trees of disk bubbles
at $t=1$ may be attached to the piece which
moves to the direction~${\tau \to +\infty}$.
Such a tree of disk bubbles corresponds to $\xi_i$,
that is, one of the factors of the fiber product
\eqref{153434}.
Since ${\rm ev}_{1,i}^{\rm deti}$ can take any value between $0$ and $1$,
there is no constraint on~$\rho_0$ for $\xi_i$.

Therefore, by definition this case is
described by the fiber product of \eqref{form1520}
and \eqref{form1521}.
See Figure~\ref{Figure15-8}.

We have thus checked that all the fiber products described by
(I), (II), (III), (IV) in Proposition~\ref{prop1518}\,(3)
appear as boundary components of ${\mathcal M}_{k_1,k_0}(L,L';M;E;\mathcal{JJ};{\rm top}(\rho))$.
To show that all other potential boundary components cancel out each other,
the most important point to observe is the following.
We consider the case when a disk bubble occurs at $t=1$
for a limit of a sequence of elements of
\smash{$\overset{\ \text{\tiny $\circ\circ$}}{\mathcal M}_{0,k_0}(L,L';M;E;\mathcal{JJ};{\rm top}(\rho))$}.
Let $E_1$ be the energy of the disk \mbox{bubble.}
The set of elements of the compactification
\smash{$\overset{\ \text{\tiny $\circ\circ$}}{\mathcal M}_{0,k_0}(L,L';M;E;\mathcal{JJ};{\rm top}(\rho))$}
corresponding to such a~disk bubble is
described by the pair $(\xi,\mathfrak x,\rho_1)$,
where
\begin{enumerate}\itemsep=0pt\setlength{\leftskip}{0.40cm}
\item[(DB.1)]
$\xi \in \mathcal M_{0+1}\bigl(L;E_1;J^{(\rho_1)}\bigr)$.
\item[(DB.2)]
\smash{$\mathfrak x \in \overset{\ \text{\tiny $\circ\circ$}}{\mathcal M}_{1,k_{0}}(L,L';M;E_2;\mathcal{JJ};{\rm top}(\rho))$}.
\item [(DB.3)]
$(\theta(\rho_1),1)$ is the (unique) boundary marked point of the element $\mathfrak x$.
\item[(DB.4)]
$\rho_0(\xi) = \rho_1$.
Here $\rho_0 \colon \mathcal M_{0+1}\bigl(L;E_1;J^{(\rho_1)}\bigr) \to [0,1]$
is as in \eqref{defrho0}.
\item [(DB.5)]
${\rm ev}_0(\xi) = {\rm ev}_{1,i}(\mathfrak x)$.
\end{enumerate}
See Figure~\ref{FiguretoDB1} below.
We remark that the disk bubble at $(\tau,1)$ is
identified with an element of
$\mathcal M_{0+1}\bigl(L;E_1;J^{(\rho(\tau))}\bigr)$.
\begin{figure}[ht]
\centering
\includegraphics[scale=0.40]{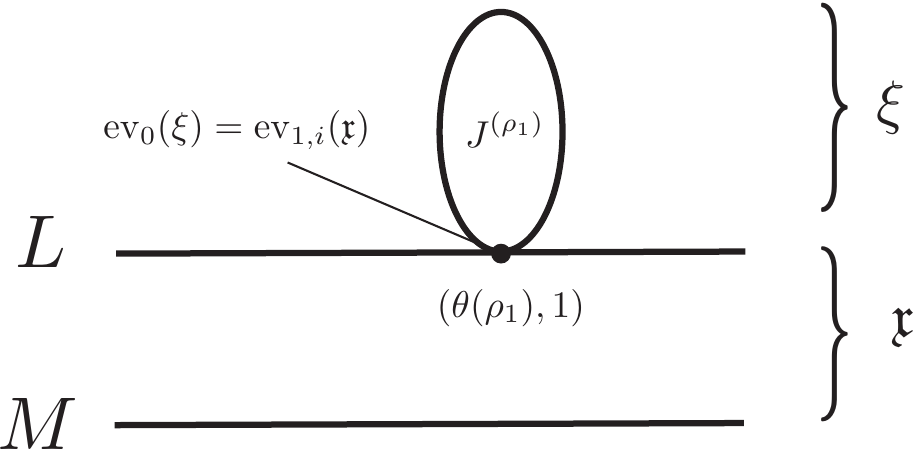}
\caption{An element $(\xi,\mathfrak x,\rho_1)$.}
\label{FiguretoDB1}
\end{figure}

We next consider the fiber product
\begin{gather}
\mathcal M_{0+1}(L;\mathcal J;E_{1};{\rm top}({\rho}))\nonumber
\\
\qquad {}_{({\rm ev}_0,\rho_0)}\times_{({\rm ev}_{1,1},{\rm ev}_{1,1}^{\rm deti})}
\overset{\ \text{\tiny $\circ\circ$}}{\mathcal M}_{1,k_0}(L,L';M;E_2;\mathcal{JJ};{\rm top}(\rho)).\label{form1536}
\end{gather}
This is a part of \smash{$\mathring{\mathcal M}_{0,k_0}(L,L';M;E;\mathcal{JJ};{\rm top}(\rho))$}
defined in Definition~\ref{1522222}.
We consider a~part~of the boundary of
\smash{$\raisebox{-1.3pt}{$\overset{\ \text{\tiny $\circ\circ$}}{\mathcal M}_{1,k_0}$}(L,L';M;E_2;\mathcal{JJ};{\rm top}(\rho))$}
which consists of
$((\R \times [0,1];z_{0,1},\vec z_1);u;\gamma;\vec \rho)$
such that $z_{0,1} = \rho_1$.
This is the case when the equality holds in the inequality \eqref{form1526}.
(See Figure~\ref{FiguredegDB1}.)
\begin{figure}[ht]
\centering
\includegraphics[scale=0.36]{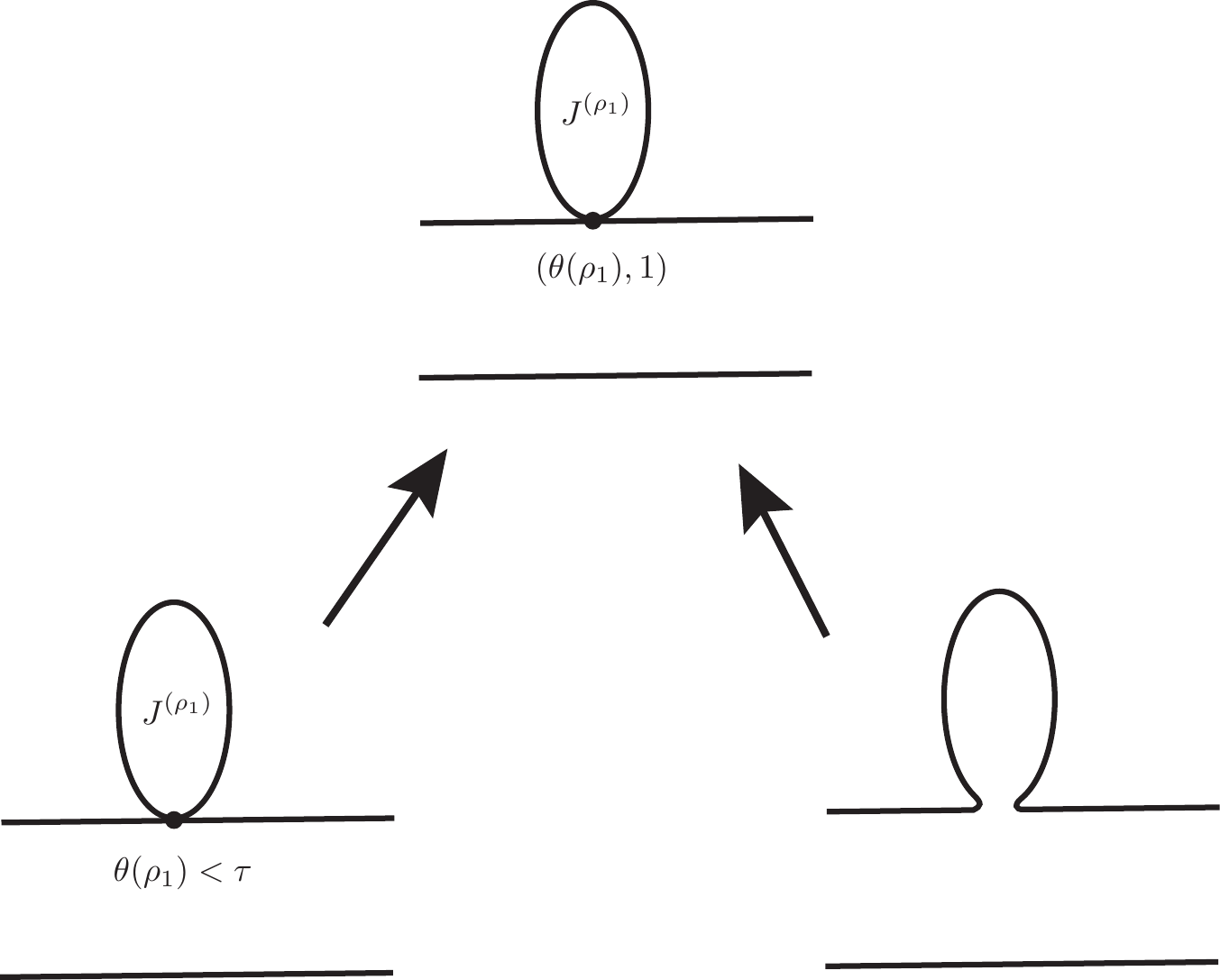}
\caption{Cancellation at $(\xi,\mathfrak x,\rho_1)$.}
\label{FiguredegDB1}
\end{figure}

Now it is easy to see that the part of the boundary of \eqref{form1536} which
we describe above cancels with the part of the boundary
corresponding to the disk bubble at $t=1$, which we describe by
(DB.1), (DB.2), (DB.3), (DB.4), (DB.5).

In a similar way as above, we can show that all the potential boundaries of
the moduli space~${\mathcal M}_{1,k_0}(L,L';M;E_2;\mathcal{JJ};{\rm top}(\rho))$
other than those spelled out in Proposition~\ref{prop1518}\,(3), (I)--(IV)
cancel each other.
The proof of Proposition~\ref{prop1518} is complete.
\end{proof}

We use the next energy estimate which is due to Chekanov \cite{chekanov}.
(See also \cite[Section 5]{fooo:polydisk}.)
From now on, we will assume\index[syindex]{H1plusnorm@$\Vert H\Vert_+$}\index[syindex]{H1minusnorm@$\Vert H\Vert_-$}
$
\int_X H_t \omega^n = 0
$
and $\int_X \omega^n = 1$.
We put
\[
\Vert H\Vert_+ =
\int_{0}^1 \sup(H_t)  d t,
\qquad
\Vert H\Vert_- =
-\int_{0}^1 \inf(H_t)  d t.
\]
They are non-negative numbers.
We remark that the Hofer distance \cite{Ho}
from $\Phi$ to identity is the infimum of
$\Vert H\Vert_- + \Vert H\Vert_+$ for all
$H$ with $\Phi_H^1 = \Phi$.
\begin{lem}\label{lem524}
If ${\mathcal M}_{k_1,k_0}(L,L';M;E;\mathcal{JJ};{\rm top}(\rho))$
is nonempty, then
$E
\ge
-\Vert H\Vert_-$.
\end{lem}
\begin{proof}
We remark
$u^*\omega = \omega\bigl(\frac{\partial u}{\partial \tau},\frac{\partial u}{\partial t}\bigr) d \tau \wedge  d t$.
By equation \eqref{eq1520}, we have
\[
\omega\left(\frac{\partial u}{\partial \tau},\frac{\partial u}{\partial t}\right)
=
-\omega\left(J_{\tau,t}\left(\frac{\partial u}{\partial \tau}\right) + \chi(\tau) X_{H_t},
\frac{\partial u}{\partial \tau}\right)
=
g\left(\frac{\partial u}{\partial \tau},\frac{\partial u}{\partial \tau}\right)
- \chi(\tau) \frac{\partial (H \circ u)}{\partial \tau}.
\]
Therefore,
\begin{align*}
\int_{\R \times [0,1]} u^* \omega
&\ge
-\int_{\R \times [0,1]} \chi(\tau)\frac{\partial (H \circ u)}{\partial \tau}
 d \tau  d t
\\
&\ge
+\int_{\R \times [0,1]} \frac{\partial\chi}{\partial \tau} (H \circ u)
 d \tau  d t
- \lim_{\tau\to+\infty}\int_{[0,1]} H(t,u(\tau,t))  d t
\\
&\ge
- \Vert H\Vert_{-}
- \lim_{\tau\to+\infty}\int_{[0,1]} H(t,u(\tau,t))  d t.
\end{align*}
Here the first inequality is a consequence of positivity of the Riemannian metric
$g$, the second equality is proved by integration by parts,
and the third inequality follows from the definition of $\Vert H\Vert_-$.

We remark that the energy $E$ is defined in Definition~\ref{definition152666}\,(5).

The lemma follows.
\end{proof}

\begin{rem}
When we identified the solution space of the equation \eqref{153915}
with the moduli space $\mathcal M_{k_{0},k_{1}}(L',M;E;J)$,
we identify $u \colon \R \times [0,1] \to X$ with $v\colon \R \times [0,1] \to X$ by
$
v(\tau,t) = \Phi^t(u_{\infty}(\tau,t))$.
The term
\smash{$ \lim_{\tau\to+\infty}\int_{[0,1]} H(t,u(\tau,t))  d t$}
which appear in the definition of the energy~$E$ and the above calculation
is related to this point.
In fact, for a solution $u$ of equation \eqref{153915} we define its energy
by
\[
\int_{\R \times [0,1]} u^*\omega
+ \lim_{\tau\to+\infty}\int_{[0,1]} H(t,u(\tau,t))  d t
-
\lim_{\tau\to-\infty}\int_{[0,1]} H(t,u(\tau,t))  d t.
\]
The third term is $0$ in our case.

Note that equation \eqref{153915} is regarded as a gradient flow equation
of certain action functional and the above energy is the difference
between values of action functional at $\tau=+\infty$ and~${\tau=-\infty}$.

\end{rem}

We define
\[
\varphi_{k_1,k_0} \colon\ B_{k_1}CF[1](L) \otimes CF(L;M) \otimes B_{k_0}CF[1](M) \to
CF(L';M)\otimes_{\Lambda_0} \Lambda
\]
by
\begin{gather*}
\varphi_{k_1,k_0}(h_{1,1},\dots,h_{1,k_1};h;h_{0,1},\dots,h_{0,k_0})
\\
\qquad=
\sum_E T^E
{\rm ev}_{+\infty}!
\bigl(
{\rm ev}_{1,1}^*h_{1,1}
\wedge
\dots \wedge {\rm ev}_{1,k_1}^*h_{1,k_1}
\wedge {\rm ev}_{-\infty}^* h
 \\
\phantom{\qquad=}{} \wedge
{\rm ev}_{0,1}^*h_{0,1}
\wedge
\dots\wedge {\rm ev}_{0,k_0}^*h_{1,k_0}
;
{\mathcal M}_{k_1,k_0}(L,L';M;E;\mathcal{JJ};{\rm top}(\rho)),
\widehat{\mathfrak S^{\varepsilon}}
\bigr).
\end{gather*}
\begin{lem}\label{lem1525}
$\{\varphi_{k_1,k_0}\}$ is a filtered $A_{\infty}$
$CF(L)$-$CF(M)$ bi-module homomorphism
over the filtered $A_{\infty}$ homomorphisms
$\mathfrak g\colon CF(L) \to CF(L')$ and ${\rm id}\colon CF(M) \to CF(M)$.\footnote{See \cite[Definition 3.7.7]{fooobook}
for the definition of $A_{\infty}$ bi-module homomorphism
over a pair of $A_{\infty}$ homomorphisms.}
\end{lem}
\begin{proof}
This is a consequence of Proposition~\ref{prop1518}
together with Stokes' formula (see \cite[Proposition 9.26]{foootech2} and \cite{fooonewbook}) and the composition formula
(see \cite[Theorem 10.20]{foootech2} and \cite{fooonewbook}).
In fact, the boundaries of type (I), (II), (III), (IV) corresponds to
\eqref{compI}, \eqref{compI2}, \eqref{compI23} and \eqref{compI24} below,
respectively,
\begin{gather}
\varphi_{k_1,k_{0,2}}(h_{1,1},\dots,h_{1,k_1};h;h_{0,1},\dots,
\mathfrak m_{k_{0,1}}(h_{0,i},\dots,h_{0,i+k_{0,1}-1}), \dots,h_{0,k_0}),\label{compI}
\\
\varphi_{k_{1,2},k_0}(h_{1,1},\dots,
\mathfrak m_{k_{1,1}}(h_{1,i},\dots,h_{1,i+k_{1,1}-1}),\dots, h_{1,k_0};h;h_{0,1},\dots,h_{0,k_0}),\label{compI2}
\\
\varphi_{k_{1,2},k_{0,2}}(h_{1,1},\dots, \mathfrak n_{k_{1,1},k_{0,1}}(h_{1,k_{1,2}-1},\dots,h_{1,k_1};h;h_{0,1},\dots,
h_{0,k_{0,1}}),\dots,
h_{0,k_0}),\label{compI23}
\\
\mathfrak n_{k_{1,2},k_{0,2}}
(\mathfrak g_{m_1}(h_{1,1},\dots,h_{1,m_1}), \dots,
\mathfrak g_{m_{\ell}}(h_{1,k_{1,2}-m_{\ell}+1},
\dots,h_{1,k_{1,2}}),\nonumber \\
\hphantom{\mathfrak n_{k_{1,2},k_{0,2}}(}
\varphi_{k_{1,2},k_{0,1}}(h_{1,k_{1,2}+1}\dots,h_{1,k_1}; h; h_{0,1}
,\dots,h_{0,k_{0,2}}),\dots,h_{0,k_2}).\label{compI24}
\end{gather}
Here in \eqref{compI24} we put $k_{1,2} = \sum_{j=1}^{\ell} m_j$.
(Note $k_{1,1} + k_{1,2} = k_1$.)
\end{proof}

Note that the filtered $A_{\infty}$ bi-module structure etc.\ appearing
in Lemma~\ref{lem1525} are curved.

We define
\[
\psi_m \colon\ CF(L;M) \otimes B_{m}CF[1](M) \to
CF(L';M)\otimes_{\Lambda_0} \Lambda
\]
by
\[
\psi_m(h;h_1,\dots,h_m) =
\sum_{k_0=0}^{\infty}\sum_{k_{1,0}=0}^{\infty}
\dots \sum_{k_{1,m}=0}^{\infty}\varphi_{k_0,k_1+\sum_{i=0}^m k_{1,i}}
\bigl(b_L^{k_0};h;b_M^{k_{1,0}}h_1 \cdots h_k b_M^{k_{1,m}}\bigr).
\]
Lemma~\ref{lem1525} now implies the following.
We use $b_M$ and $b_L$ to define a
strict and unital filtered~$A_{\infty} \bigl(CF(M),\bigl\{\mathfrak m^{b_M}_k\bigr\}\bigr)$
right module structure on $CF(L;M)$.
We use $b_M$ and $b_{L'}$ to define a
strict and unital filtered \smash{$A_{\infty} \bigl(CF(M),\bigl\{\mathfrak m^{b_M}_k\bigr\}\bigr)$}
right module structure on $CF(L';M)$.
\begin{lem}
$\{\psi_m \mid m=0,1,2,\dots\}$ define a
strict and unital right filtered $A_{\infty}$
module homomorphism: $CF(L;M)\otimes_{\Lambda_0}{\Lambda} \to CF(L';M)\otimes_{\Lambda_0}{\Lambda}$.
\end{lem}
Lemma~\ref{lem524} implies the next lemma.
\begin{lem}
$\bigl\{T^{\Vert H\Vert_+}\psi_m \mid m=0,1,2,\dots\bigr\}$
define a
strict and unital right filtered $A_{\infty}$
module homomorphism: $CF(L;M) \to CF(L';M)$.
\end{lem}
By exchanging the role of $L'$ and $L$, we obtain
the following.
\begin{lem}
There exists $\{\psi'_m \mid m=0,1,2,\dots\}$
which define a
strict and unital left filtered~$A_{\infty}$
module homomorphism: $CF(L';M)\otimes_{\Lambda_0}{\Lambda} \to CF(L;M)\otimes_{\Lambda_0}{\Lambda}$.
Moreover, $\bigl\{T^{\Vert \Phi\Vert_-}\psi'_m \mid m=0,1,2,\dots\bigr\}$
define a
strict and unital left filtered $A_{\infty}$
module homomorphism: $CF(L';M) \to CF(L;M)$.
\end{lem}
We put $\psi = \{\psi_m \mid m=0,1,2,\dots\}$
and $\psi' = \{\psi'_m \mid m=0,1,2,\dots\}$.
We can use a similar argument (one parameter version)
to show that there exists a strict and unital
filtered $A_{\infty}$ pre-natural transformations
$\phi$ and $\phi'$ such that
\[\psi' \circ \psi - \mathfrak m_1(\phi) = {\rm identity},\qquad \psi \circ \psi' - \mathfrak m_1(\phi') = {\rm identity}.
\]
Moreover, $T^{\Vert H\Vert_++ \Vert H\Vert_-}\phi$ determines pre-natural transformation
$CF(L;M) \to CF(L;M)$ and
 $T^{\Vert H\Vert_++ \Vert H\Vert_-}\phi'$ determines pre-natural transformation
$CF(L';M) \to CF(L';M)$.
We omit the detail of the proof of this statement.
See \cite[Sections 5.3.3 and 5.3.4]{fooobook}
and \cite[Lemma 6.4]{fooo:polydisk} for the proof of this part.
(The way to adapt the argument there to the current situation is the
same as the way we do so for $\psi$ which we explained in detail above.)

We thus proved that $CF(L';M)$ is equivalent to $CF(L;M)$ over $\Lambda$
in the category of right $CF(M)$ module.
This proves Theorem~\ref{thm154}\,(2).

To prove Theorem~\ref{thm154}\,(3),
it suffices to recall that
the infimum of $\Vert H\Vert_+ + \Vert H\Vert_-$
over all $H$ which generates the Hamiltonian diffeomorphism~$\Phi$
is nothing but the Hofer distance between~$\Phi$ and the
identity map.

The proof of Theorem~\ref{thm154} is now complete.

\subsection{Completion by Hofer distance}
\label{sec:comphofer}

In \cite[Section 6.5.4]{fooobook},
we proved that if $\mathfrak c_i$, $\mathfrak c_{i'}$
are Cauchy sequences of objects of a strict fil\-tered~$A_{\infty}$
category $\mathscr C$ with respect to the Hofer distance
$d_H$, then we can define an inductive limit
\begin{equation}\label{1533333}
\lim_{i\to\infty} HF(\mathfrak c_i, \mathfrak c_{i'})
\end{equation}
as $\Lambda_0$ modules.
Namely, we consider the $\mathfrak m_1$ cohomology
$HF(\mathfrak c_i, \mathfrak c_{i'})$ and write it as
\[
HF(\mathfrak c_i, \mathfrak c_{i'})
= \Lambda_0^{n}\oplus
\bigoplus_{j=1}^{k_i} \frac{\Lambda_0}{T^{\lambda_{i,j}}\Lambda_0}.
\]
Here $\lambda_{i,j}$ are positive numbers such that $\lambda_{i,j}\ge \lambda_{i,j+1}$.\footnote{See \cite[Theorem 6.1.20]{fooobook}.}
In fact, $d_H(\mathfrak c_i, \mathfrak c_{j}) < \infty$
the rank $n$ is independent of $i$.
The torsion exponent $\lambda_{i,j}$ is one Lipschitz
by \cite[Theorem 6.1.25]{fooobook} and \cite[Theorem 6.2]{fooo:polydisk}.
We use it to define the inductive limit \eqref{1533333}.
The inductive limit \eqref{1533333} has a~form
\begin{equation}\label{15333332}
\lim_{i\to\infty} HF(\mathfrak c_i, \mathfrak c_{i'})
=
\Lambda_0^{n}\oplus
\bigoplus \frac{\Lambda_0}{T^{\lambda_{j}}\Lambda_0}.
\end{equation}
Here the direct sum in the second factor
may be infinite sum with $\lim_{j\to \infty}\lambda_j \to 0$.
See \cite[Proposition 6.5.38]{fooobook} and \cite[Example 6.5.40]{fooobook}.
It seems likely that we can prove the next conjecture purely
algebraically.
\begin{conj}\label{cong1653}
The $A_{\infty}$ operations $\mathfrak m_k$ extend `continuously'
to the limit \eqref{15333332} and define a
`filtered $A_{\infty}$ category' whose object
is a Cauchy sequence of $\mathfrak{OB}(\mathscr{C})$.
\end{conj}
\begin{rem}
Conjecture \ref{cong1653} appeared in the preprint version of this
paper in 2017. A~version of its positive answer is now given in \cite{Hausdoredd}.
\end{rem}
One reason why proving Conjecture \ref{cong1653}, taking completion of
$\mathfrak{OB}(\mathscr{C})$ and trying to find a~filtered $A_{\infty}$ category whose object is an element
of such a completion, could be interesting is as follows.
In this paper, we consider only a set of Lagrangian submanifolds
$L_1$, $L_{12}$ etc.\ which satisfy certain `clean intersection'
properties. If $L_1$ and $L_{12}$ do not necessary have clean
intersection, the geometric transformation of $L_1$ by $L_{12}$
may not exist. However, Theorem~\ref{thm154} implies that
it exists as an object of a
certain completion of $\mathfrak{Fukst}(X_2)$.
So by taking the completion with respect to the Hofer distance,
we may take the geometric transformation and the composition of Lagrangian
correspondences without assuming any kinds of transversality or
cleanness of the Lagrangian submanifolds involved.

\section{K\"unneth bi-functor revisited}\label{sec:kuneth}

\subsection[Tensor product of filtered $A_\infty$ categories]{Tensor product of filtered $\boldsymbol{A_{\infty}}$ categories}
\label{sec:TensorAinfini}

We begin with defining the tensor product of
filtered $A_{\infty}$ categories.
There are various works such as \cite{Limo1,KS,marsh,SU} etc.\ on this subject.
We describe it using the notion of filtered $A_{\infty}$ bi-functor.
We remark that in this section, we use the sign convention of the
filtered $A_{\infty}$ multi-module so that its element $v$
contributes $\deg' v$ to the sign. This convention is different
from one we used in Section~\ref{2-category formulation}, where
the contribution is $\deg v$.
\begin{rem}
In this section, we always assume the ground ring $R$ is a field.
We also assume that filtered $A_{\infty}$ categories are always
gapped. Moreover, for two objects $c$, $c'$ of $\mathscr C$ we assume
that the $\overline{\mathfrak m}_1$ cohomology $H(\mathscr C(c,c');\overline{\mathfrak m}_1)$
is finitely generated. (It is then a finite direct sum of~$\Lambda_0$.)

Under this assumption, the cohomology $H(\mathscr C(c,c');\mathfrak m_1)$
is isomorphic to a direct sum of finitely many copies of
$\Lambda_0$ or $\Lambda_0/T^a\Lambda_0$
(see \cite[Proposition 6.3.14]{fooobook}).
Using this fact, cohomology of completed tensor product behaves
in the same way as the case of usual tensor product over Dedekind ring.
In fact,
\[
\frac{\Lambda_0}{T^a\Lambda_0} \,\widehat{\otimes}\, \frac{\Lambda_0}{T^b\Lambda_0}
=
\frac{\Lambda_0}{T^a\Lambda_0}, \qquad
{\rm Tor}\left(\frac{\Lambda_0}{T^a\Lambda_0},\frac{\Lambda_0}{T^b\Lambda_0}\right)
=
\frac{\Lambda_0}{T^a\Lambda_0},
\]
if $a \le b$.
\end{rem}

It seems that the construction of the tensor product below is
a category version of one
suggested by Kontsevich and Soibelman \cite[p.~174, line~6]{KS}.

Let $\mathscr C_i$ be a unital filtered $A_{\infty}$
category for $i=1,2$.
There are 2 versions of the story of tensor products of filtered $A_{\infty}$
categories, that are, strict and $G$-gapped versions.
Let
\[
\mathcal{BIFUNC}(\mathscr C_1^{\rm op} \times \mathscr C_2^{\rm op}
;\mathcal{CH})
\]
be the filtered $A_{\infty}$ category whose objects are filtered $A_{\infty}$ bi-functors
$\mathscr C_1^{\rm op} \times \mathscr C_2^{\rm op}
\to \mathcal{CH}$.
We require objects of $\mathcal{BIFUNC}$ to be strict (resp.\ $G$-gapped).
In other words, it is a category of left~$\mathscr C_1$,~$\mathscr C_2$ bimodules.\footnote{In the gapped case, we use the language of left filtered $A_{\infty}$ modules.}
\begin{lem}\label{lem1610}
There exists a strict $($resp.\ $G$-gapped$)$ and unital filtered $A_{\infty}$
bi-functor
\[
\mathfrak{BIYON}\colon\
\mathscr C_1 \times \mathscr C_2
\to
\mathcal{BIFUNC}(\mathscr C_1^{\rm op} \times \mathscr C_2^{\rm op}
;\mathcal{CH}).
\]
\end{lem}

\begin{defn}
We call $\mathfrak{BIYON}$
\index[syindex]{BIyon@$\mathfrak{BIYON}$} the $A_{\infty}$ {\it bi-Yoneda functor}. \index{bi-Yoneda functor}
\end{defn}
\begin{proof}
We discuss the strict case. The construction of $G$-gapped case is similar.
Let $\mathfrak c = (c_1,c_2)$ be an objects of $\mathscr C_1 \times \mathscr C_2$
(namely, $c_i \in \mathfrak{OB}(\mathscr C_i)$).
We construct
$\mathfrak{BiYon}({\mathfrak c}) \colon \mathscr C_1^{\rm op} \times \mathscr C_2^{\rm op}
\to \mathcal{CH}$.
Let $\mathfrak b = (b_1,b_2)$ be an object of
$\mathscr C_1^{\rm op} \times \mathscr C_2^{\rm op}$.

We define
$
\mathfrak{BiYon}_{\rm ob}({\mathfrak c})(\mathfrak b)
=
\mathscr C_1(b_1,c_1) \,\widehat{\otimes}\, \mathscr C_2(b_2,c_2)
$
which is a chain complex.
This is the object part.
We next define
\begin{gather*}
\mathfrak{BiYon}_{k_1,k_2}({\mathfrak c})\colon\
B_{k_1}\mathscr C_1^{\rm op}[1](b_{1,1},b_{1,2})
\,\widehat\otimes\, B^{\rm op}_{k_2}\mathscr C_2[1](b_{2,1},b_{2,2}) \\
\hphantom{\mathfrak{BiYon}_{k_1,k_2}({\mathfrak c})\colon} \
\to
\operatorname{Hom}
(\mathscr C_1(b_{1,1},c_1) \,\widehat{\otimes}\, \mathscr C_2(b_{1,2},c_2),
\mathscr C_1(b_{2,1},c_1) \,\widehat{\otimes}\, \mathscr C_2(b_{2,1},c_2)).
\end{gather*}
If $k_1 \ne 0$ and $k_2 \ne 0$, we put
$\mathfrak{BiYon}_{k_1,k_2}({\mathfrak c}) = 0$.
Otherwise, we define
\begin{gather*}
\mathfrak{BiYon}_{k_1,0}({\mathfrak c})({\bf x}_1 \otimes 1)(z_1 \otimes
z_2)
=
(-1)^{*_1}\mathfrak m({\bf x}_1^{\rm op},z_1) \otimes z_2,
\\
\mathfrak{BiYon}_{0,k_2}({\mathfrak c})(1 \otimes {\bf x}_2)(z_1 \otimes
z_2)
=
(-1)^{*_2} z_1 \otimes \mathfrak m({\bf x}_2^{\rm op},z_2).
\end{gather*}
Here
$
*_1 = \varepsilon({\bf x}_1)$, $
*_2 =
\varepsilon({\bf x}_2)
+ (1 + \deg' {\bf x}_2) \deg' y_1
$
are Koszul sign.
(The symbol $\varepsilon({\bf x})$ is defined in \eqref{form215}.)
It is easy to check that
$\mathfrak{BiYon}({\mathfrak c}) =
(\mathfrak{BiYon}_{\rm ob}({\mathfrak c}),\{\mathfrak{BiYon}_{k_1,k_2}({\mathfrak c})\})$
becomes a~filtered~$A_{\infty}$ bi-functor.
(The calculation is similar to \cite[p.~93]{fu4}, which is the case of usual Yoneda functor.)

We thus constructed the object part of bi-Yoneda functor.
We next construct its morphism part.
Let $\mathfrak c_i = (c_{i,1},c_{i,2})$ be an element of $\mathfrak{OB}(\mathscr C_1) \times \mathfrak{OB}(\mathscr C_2)$
for $i=1,2$.
We denote by $
({\bf C}({\mathfrak c}_1,{\mathfrak c}_2),d)
$
the complex of all pre-natural transformations from
$\mathfrak{BiYon}_{\rm ob}({\mathfrak c}_1)$ to $\mathfrak{BiYon}_{\rm ob}({\mathfrak c}_2)$.
We define the product $
\circ \colon {\bf C}({\mathfrak c}_1,{\mathfrak c}_2)
\otimes {\bf C}({\mathfrak c}_2,{\mathfrak c}_3) \to {\bf C}({\mathfrak c}_1,{\mathfrak c}_3)$
by the composition of pre-natural transformations.
We thus obtain a DG-category
$
({\bf C},d,\circ)$.
(We use the fact $\mathcal{CH}$ is not only an
$A_{\infty}$ category but also DG-category to obtain
this DG category.)
We regard it as an $A_{\infty}$ category.\looseness=-1

We now define a filtered map
\[
\mathfrak{BiYon}_{\ell_1,\ell_2} \colon\
B_{\ell_1}\mathscr C_1[1](c_{1,1},c_{1,2})
\otimes
B_{\ell_2}\mathscr C_2[1](c_{2,1},c_{2,2})
\to {\bf C}({\mathfrak c}_1,{\mathfrak c}_2).
\]
Let
$\mathfrak b_i = (b_{i,1},b_{i,2})$, $\mathfrak c_i = (c_{i,1},c_{i,2})$
an element of $\mathfrak{OB}(\mathscr C_1) \times \mathfrak{OB}(\mathscr C_2)$
$i=1,2$.
Let
\begin{gather}
({\bf x}_1 \otimes {\bf x}_2)
\in B_{k_1}\mathscr C_1^{\rm op}[1](b_{1,1},b_{2,1})
\,\widehat\otimes\, B_{k_2}\mathscr C^{\rm op}_2[1](b_{1,2},b_{2,2}),\nonumber \\
({\bf y}_1 \otimes {\bf y}_2)
\in B_{\ell_1}\mathscr C_1[1](c_{1,1},c_{2,1})
\,\widehat\otimes\, B_{\ell_2}\mathscr C_2[1](c_{1,2},c_{2,2})\label{formdef1688}
\end{gather}
and $(z_1,z_2) \in \mathscr F(\mathfrak c_1)_{\rm ob}(\mathfrak b_2)$.
We define
\[
((\mathfrak{BiYon}_{\ell_1,\ell_2}({\bf y}_1\otimes {\bf y}_2)_{k_1,k_2})({\bf x}_1 \otimes {\bf x}_2))(z_1,z_2)
=0
\]
if
$(k_1,\ell_1) \ne (0,0)$ and $(k_2,\ell_2) \ne (0,0)$.
In case either $(k_1,\ell_1) = (0,0)$ or $(k_2,\ell_2) = (0,0)$,
we define
\begin{gather}
((\mathfrak{BiYon}_{\ell_1,0}({\bf y}_1\otimes 1)_{k_1,0})({\bf x}_1 \otimes 1))(z_1,z_2)\nonumber \\
\qquad=(-1)^{*_1} \mathfrak m({\bf x}_1^{\rm op},z_1,{\bf y}_1) \otimes z_2 \in
\mathscr C_1(b_{1,1},c_{2,1}) \,\widehat{\otimes}\, \mathscr C_2(b_{1,2},c_{2,2})\label{form161999}
\end{gather}
(note that $b_{1,2} = b_{2,2}$, $c_{1,2} = c_{2,2}$ in this case)
and
\begin{gather}
((\mathfrak{BiYon}_{0,\ell_2}(1 \otimes {\bf y}_2)_{0,k_2})(1 \otimes{\bf x}_2))(z_1,z_2) \nonumber\\
\qquad=(-1)^{*_2} z_1 \otimes \mathfrak m({\bf x}_2^{\rm op},z_2,{\bf y}_2) \in
\mathscr C_1(b_{1,1},c_{2,1}) \,\widehat{\otimes}\, \mathscr C_2(b_{1,2},c_{2,2}).\label{form1619992}
\end{gather}
(Note that $b_{1,1} = b_{2,1}$, $c_{1,1} = c_{2,1}$ in this case.)
Here $*_i$ is the Koszul sign, that is,
\begin{gather*}
*_1 = \varepsilon({\bf x}_1) + (\deg' z_1+\deg' {\bf x}_1)\deg'{\bf y}_1,
\\
*_2 = \varepsilon({\bf x}_2) + \deg' {\bf y}_2(\deg' {\bf x}_2+ \deg' z_1+ \deg' z_2)
+ \deg' z_1(\deg' {\bf x}_2 + 1).
\end{gather*}
\begin{sublem}
\eqref{form161999} and \eqref{form1619992} define filtered $A_{\infty}$ bi-functor.
\end{sublem}
The proof is a straightforward calculation similar to the proof of \cite[Lemma 9.8]{fu4}
and so is omitted.
The proof of Lemma~\ref{lem1610} is complete.
\end{proof}

\begin{defn}\label{defn16111}
Let $\mathscr C_1$ and $\mathscr C_2$ be unital
filtered $A_{\infty}$ categories.
We define the full subcategory
of $\mathcal{BIFUNC}(\mathscr C_1^{\rm op} \times \mathscr C_2^{\rm op}
;\mathcal{CH})$
whose objects are image of the bi-Yoneda functor
the {\it tensor product} \index{tensor product of filtered $A_{\infty}$ category} of $\mathscr C_1$ and $\mathscr C_2$
and write $\mathscr C_1 \otimes \mathscr C_2$.\index[syindex]{Cscr1timesC2@$\mathscr C_1 \otimes \mathscr C_2$}
By definition, there exists a strict and unital filtered $A_{\infty}$
bi-functor $\mathscr C_1 \times \mathscr C_2 \to \mathscr C_1 \otimes \mathscr C_2$.
\end{defn}

It is easy to show that
$\mathscr C_1 \otimes \mathscr C_2$ is homotopy
equivalent to $\mathscr C'_1 \otimes \mathscr C'_2$
if $\mathscr C_i$ is homotopy equivalent to $\mathscr C'_i$.

\begin{lem}\label{DGcaselem}
Suppose $\mathscr C_1$, $\mathscr C_2$ are DG-categories.
Then the tensor product as filtered $A_{\infty}$ category $\mathscr C_1 \otimes \mathscr C_2$
is homotopy equivalent to the $($DG-category$)$ tensor product $\mathscr C_1 \otimes \mathscr C_2$
as filtered~$A_{\infty}$ categories.
\end{lem}

We prove Lemma~\ref{DGcaselem} in Section~\ref{prooflemmata}.
Lemma~\ref{DGcaselem} implies that the tensor product
defined in \cite{Limo1,DL} etc.\ is the tensor product in the sense of
Definition~\ref{defn16111}.

We put ${\bf C} = \mathscr C_1 \otimes \mathscr C_2$.
Note that by construction there exists a
left $\mathscr C_1$, $\mathscr C_2$ and
right ${\bf C}$, filtered~$A_{\infty}$ tri-module
$
{\bf M}(\mathscr C_1,\mathscr C_2;{\bf C}),
$
and left ${\bf C}$ right $\mathscr C_1$, $\mathscr C_2$ tri-module
$
{\bf M}({\bf C};\mathscr C_1,\mathscr C_2),
$
as follows.
Let~$\mathfrak c = (c_1,c_2) \in
\mathfrak{OB}(\mathscr C_1) \times \mathfrak{OB}(\mathscr C_2)$,
$\mathfrak b = (b_1,b_2) \in \mathfrak{OB}({\bf C}) =
\mathfrak{OB}(\mathscr C_1) \times \mathfrak{OB}(\mathscr C_2)$.
Then we put
\begin{gather}
{\bf M}(\mathscr C_1,\mathscr C_2;{\bf C})(c_1,c_2;\mathfrak b)
= \mathscr C_1[1](c_1,b_1) \otimes \mathscr C_2[1](c_2,b_2), \nonumber\\
{\bf M}({\bf C};\mathscr C_1,\mathscr C_2)(\mathfrak b;c_1,c_2)
= \mathscr C_1[1](b_1,c_1) \otimes \mathscr C_2[1](b_2,c_2).\label{form161212}
\end{gather}
This is the object part of our tri-module.
The morphism part is defined as follows.
Let $\mathcal T \in {\bf C}({\mathfrak b}_1,{\mathfrak b}_2)$
and ${\bf T} \in B_{k}{\bf C}({\mathfrak b}_1,{\mathfrak b}_2)$.
Let ${\bf x}_1 \otimes {\bf x}_2$ be as in
\eqref{formdef1688}.
Let $(z_1,z_2) \in \mathcal C_1[1](c_{2,1},b_{1,1}) \otimes \mathcal C_2[1](c_{2,2},b_{1,2})
= {\bf M}(\mathscr C_1,\mathscr C_2;{\bf C})(c_{2,1},c_{2,2};\mathfrak b_1)$.
Now we define the bi-module structure $\mathfrak n_{\ell_1,\ell_2;k}$
as
\begin{equation}\label{form16130}
\mathfrak n_{\ell_1,\ell_2;k}({\bf x}_1 \otimes {\bf x}_2;(z_1,z_2);{\bf T})
\in {\bf M}(\mathscr C_1,\mathscr C_2;{\bf C})(c_{1,1},c_{1,2};\mathfrak b_2),
\end{equation}
where
$
\eqref{form16130} = 0
$
unless
$(\ell_1,\ell_2;k) = (\ell_1,0;0)$,
$(\ell_1,\ell_2;k) = (0,\ell_2;0)$,
or $k=1$.
In case $k=1$, we define
\begin{equation}\label{form16132}
\eqref{form16130} = (-1)^{\mathfrak {deg}'\mathcal T_{\ell_1,\ell_2}(\deg' {\bf x}_1
+ \deg' {\bf x}_2 +\deg'z_1 + \deg'z_2)}\mathcal T_{\ell_1,\ell_2}({\bf x}_1\otimes {\bf x}_2)
(z_1,z_2).
\end{equation}
Here ${\bf T} = \mathcal T$ is a pre-natural transformation.
In case $k=0$, the structure $\mathfrak n_{\ell_1,\ell_2;k}$
is nothing but the left filtered $A_{\infty}$ module structure
over $\mathscr C_1$, $\mathscr C_2$,
which is nothing but the filtered $A_{\infty}$ bi-functor
$\mathscr C^{\rm op}_1 \times \mathscr C^{\rm op}_2 \to \mathcal{CH}$.
More explicitly, it is
\begin{gather}
\mathfrak n_{\ell_1,0;0}({\bf x}_1\otimes 1;(z_1,z_2);1)
= \mathfrak m({\bf x}_1,z_1) \otimes z_2,\nonumber\\
\mathfrak n_{0,\ell_2;0}(1\otimes {\bf x}_2;(z_1,z_2);1)
= (-1)^{*}z_1 \otimes \mathfrak m({\bf x}_2,z_2),\label{form1615}
\end{gather}
with $* = (\deg' {\bf x}_2 + 1) \deg' z_1$.

The definition of tri-module structure
on ${\bf M}({\bf C};\mathscr C_1,\mathscr C_2)$ is similar.

By the definition of a pre-natural transformation,
the composition and the differential,
it is straightforward to check that
\eqref{form161212}--\eqref{form1615} define
filtered $A_{\infty}$ tri-module structure.
(We remark that ${\bf C}$ is a DG-category
because $\mathcal{CH}$ is a DG-category.)

We next define a filtered $A_{\infty}$
bi-module ${\bf M}({\bf C};{\bf C})$
over ${\bf C} \times {\bf C}$ as follows:
$
{\bf M}({\bf C};{\bf C})(\mathfrak c,\mathfrak b)
= {\bf C}(\mathfrak c,\mathfrak b)$,
and the structure operations of the bi-moduli structure
are given by the structure operation of the
filtered $A_{\infty}$ category ${\bf C}$.
(This is actually a DG bi-module.)
Using $A_{\infty}$ bi-functor $\mathscr C_1 \times \mathscr C_2 \to {\bf C}$, we regard
${\bf M}({\bf C};{\bf C})$ as a left $\mathscr C_1$, $\mathscr C_2$ and right ${\bf C}$
$A_{\infty}$ tri-modules or left ${\bf C}$ and right $\mathscr C_1$, $\mathscr C_2$
$A_{\infty}$ tri-module.
\begin{lem}\label{lem1613}
There exists a homotopy equivalence
$
{\bf M}({\bf C};{\bf C}) \sim {\bf M}(\mathscr C_1,\mathscr C_2;{\bf C})
$
as left $\mathscr C_1$, $\mathscr C_2$
and right ${\bf C}$ tri-modules.
There exists also a homotopy equivalence
$
{\bf M}({\bf C};{\bf C}) \sim {\bf M}({\bf C};\mathscr C_1,\mathscr C_2)
$
as left ${\bf C}$
and right $\mathscr C_1$, $\mathscr C_2$ tri-modules.

\end{lem}
The proof is given in Section~\ref{prooflemmata}.

\subsection{K\"unneth functor in Lagrangian Floer theory}
\label{sec:Kunnethfunctor}

\begin{situ}
Let $(X_i,\omega_i)$ be a compact symplectic manifold,
$V_i$ a background datum of $X_i$,
and $\mathbb L_i$ a finite set of $V_i$-relatively spin
immersed Lagrangian submanifolds for $i=1,2$.
We assume $\mathbb L_i$ is a clean collections for $i=1,2$.

We obtain a curved filtered $A_{\infty}$ category
$\mathfrak{Fuk}(X_i;\mathbb L_i)$ for $i=1,2$.
We denote by $\mathfrak{Fukst}(X_i;\mathbb L_i)$ the strict
category associated to $\mathfrak{Fuk}(X_i;\mathbb L_i)$.

We consider the direct product
$(X_1\times X_2,\omega_1\oplus \omega_2)$
and the background datum $\pi_1^*V_1 \oplus \pi_1^*V_2$ on it.
We put
$
\mathbb L_1 \times \mathbb L_2
:=
\{
L_1 \times L_2 \mid L_1 \in \mathbb L_1, L_2 \in \mathbb L_2
\}$.
This is a clean collection of $\pi_1^*V_1 \oplus \pi_1^*V_2$ relatively
spin immersed Lagrangian submanifolds.
We then obtain a curved filtered~$A_{\infty}$ category
$\mathfrak{Fuk}(X_1 \times X_2;\mathbb L_1 \times \mathbb L_2)$.
We denote by $\mathfrak{Fukst}(X_1 \times X_2;\mathbb L_1 \times \mathbb L_2)$ the strict
category associated to $\mathfrak{Fuk}(X_1 \times X_2;\mathbb L_1 \times \mathbb L_2)$.
\end{situ}
The next theorem is the main result of this section.
\begin{thm}\label{thm1615}
There exists a strict and unital filtered $A_{\infty}$ functor
\[
\mathfrak{Fukst}(X_1;\mathbb L_1) \otimes \mathfrak{Fukst}(X_2;\mathbb L_2)
\to
\mathfrak{Fukst}(X_1 \times X_2;\mathbb L_1 \times \mathbb L_2),
\]
which is a homotopy equivalence to the image.
\end{thm}
Theorem~\ref{thm1615} was obtained previously by L.\ Amorin \cite{Limo}
by a different method.

We call the functor in Theorem~\ref{thm1615} the K\"unneth bi-functor
and denote it by $\mathscr K$.\index{K\"unneth bi-functor}
\begin{cor}
Let $L_i \subset X_i$ be a $V_i$-relatively spin immersed Lagrangian
submanifold
for ${i=1,2}$. Suppose $L_1$, $L_2$ are unobstructed.
Then $L_1 \times L_2$ is also unobstructed.
Moreover, bounding cochains $b_1$ and $b_2$ of $L_1$ and $L_2$
determine a bounding cochain $b_1 \times b_2$ of~${L_1 \times L_2}$
canonically up to gauge equivalence and we have an exact sequence
\begin{align*}
0
&\to
{\rm Tor}(HF((L_1,b_1),(L'_1,b'_1)), HF((L_2,b_2),(L'_2,b'_2)) \\
&\to
HF((L_1\times L_2,b_1\times b_2),(L'_1\times L'_2,b'_1\times b'_2)) \\
&\to
HF((L_1,b_1),(L'_1,b'_1)) \otimes_{\Lambda_0} HF((L_2,b_2),(L'_2,b'_2))
\to 0.
\end{align*}

\end{cor}
\begin{proof}
Let $\mathscr C_1$, $\mathscr C_2$ be strict and unital filtered $A_{\infty}$ categories
and $c_i,c'_i \in \mathfrak{OB}(\mathscr C_i)$.
It suffices to show the existence of the next exact sequence
\begin{align*}
0
&\to
{\rm Tor}(H(\mathscr C_1(c_1,c'_1),\mathfrak m_1), H(\mathscr C_2(c_2,c'_2),\mathfrak m_1)) \\
&\to
H((\mathscr C_1\otimes \mathscr C_2)((c_1,c_2),(c'_1,c'_2)),\mathfrak m_1)
\to
H(\mathscr C_1(c_1,c'_1),\mathfrak m_1) \otimes_{\Lambda_0} H(\mathscr C_2(c_2,c'_2),\mathfrak m_1)
\to 0.
\end{align*}
This is immediate in case $\mathscr C_1$, $\mathscr C_2$ are DG-categories
by Lemma~\ref{DGcaselem}.
The corollary then follows from the fact that any filtered $A_{\infty}$ category is homotopy equivalent to
a DG-category.
\end{proof}

\begin{proof}[Proof of Theorem~\ref{thm1615}]
We apply Theorems~\ref{trimain} and \ref{opthere}. We then obtain a
left $\mathfrak{Fuk}((X_1 \times X_2,\pi_1^*(V_1)\oplus \pi_2^*(V_2);\mathbb L_1
\times \mathbb L_2)$, right $\mathfrak{Fuk}((X_1,V_1);\mathbb L_1)$,$\mathfrak{Fuk}((X_2,V_2);\mathbb L_2)$ filtered $A_{\infty}$ tri-module
\begin{equation}\label{trifindd1610}
\mathscr{CF}(\mathbb L_{12};\mathbb L_1,\mathbb L_2)
\end{equation}
as follows.
We replace $(X_1,\omega_1), V_1$ in Theorems~\ref{trimain}
by $(X_1,-\omega_1), V_1\oplus TX_1$.
Since
\[
\mathfrak{Fuk}((-X_1,V_1\oplus TX_1);\mathbb L_2)
\cong \mathfrak{Fuk}((X_1,V_1);\mathbb L_1)^{\rm op}
\]
by Theorem~\ref{opthere}, we obtain \eqref{trifindd1610}.

The tri-module \eqref{trifindd1610} induces a filtered $A_{\infty}$ bi-functor
\[
\mathscr F \colon\
\mathfrak{Fukst}(X_1;\mathbb L_1)
\times
\mathfrak{Fukst}(X_2;\mathbb L_2)
\to
\mathcal{FUNC}(\mathfrak{Fukst}(X_1 \times X_2;\mathbb L_1 \times \mathbb L_2)^{\rm op},
\mathcal{CH}).
\]
\begin{prop}\label{lem167}
Let $(L_i,b_i) \in \mathfrak{OB}(\mathfrak{Fukst}(X_i;\mathbb L_i))$.
\begin{enumerate}\itemsep=0pt
\item[$(1)$]
$L_1 \times L_2$ is unobstructed.
Moreover, $b_1$, $b_2$ determine a bounding cochain $b_1\times b_2$
of $L_1 \times L_2$ up to gauge equivalence canonically.\index[syindex]{b1timesb2@$b_1\times b_2$}
\item[$(2)$]
The object $(L_1 \times L_2,b_1\times b_2)$
of the category $\mathfrak{Fukst}(X_1 \times X_2;\mathbb L_1 \times \mathbb L_2)$
represents the functor~${\mathscr F((L_1,b_1),(L_2,b_2)) \colon
\mathfrak{Fukst}(X_1 \times X_2;\mathbb L_1 \times \mathbb L_2)^{\rm op} \to
\mathcal{CH}}$.
\end{enumerate}
\end{prop}
\begin{proof}
The proof of (1) is similar to the proof of Theorem~\ref{thm61}
and the proof of (2) is similar to the proof of Theorem~\ref{th72}.

We start with (1).
The tri-module
applied to $L_1$, $L_2$ and $L_1 \times L_2$
induces a left $CF(L_1\times L_2)$, right $CF(L_1)$, $CF(L_2)$
tri-module $D$.
Its structure operations are
\[
\mathfrak n_{k_{12},k_1,k_2} \colon\
B_{k_{12}}CF(L_1 \times L_2)[1] \otimes D \otimes
 B_{k_1}CF(L_1)[1] \otimes B_{k_2}CF(L_2)[1]
\to D.
\]
We use bounding cochains $b_1$, $b_2$ to deform it to obtain
$
\mathfrak n^{b_1,b_2}_{k_{12}} \colon
B_{k_{12}}CF(L_1 \times L_2) \otimes D
\to D$.
Namely, we put
\[
\mathfrak n^{b_1,b_2}_{k} (x_1,\dots,x_{k};y)
=
\sum_{k_1,k_2} \mathfrak n_{k,k_1,k_2} \bigl(x_1,\dots,x_{k};y;b_1^{k_1},b_2^{k_2}\bigr).
\]
The maps $\bigl\{\mathfrak n^{b_1,b_2}_{k} \mid k=0,1,2,\dots\bigr\}$ define a structure of
left $CF(L_1 \times L_2)$ module over $D$.

By definition, there exists an isomorphism
$
D = \Omega\bigl(\bigl(\tilde L_1 \times_{X_1} \bigl(\tilde L_1 \times \tilde L_2\bigr) \times_{X_2} \tilde L_2\bigr)\bigr)
\otimes \hat \Lambda_0
$
as $\Lambda_0$ modules.
So $D$ is actually isomorphic to $CF(L_1 \times L_2)$,
as a $\Lambda_0$ module.
The differential $0$-form (function) $1$ on the diagonal component
$
\tilde L_1 \times \tilde L_2
\subset
\tilde L_1 \times_{X_1} \bigl(\tilde L_1 \times \tilde L_2\bigr) \times_{X_2} \tilde L_2
$
is a cyclic element of the left $CF(L_1\times L_2)$
module $D$.
(We can prove it in the same way as Proposition~\ref{prop610}.)

Proposition~\ref{lem167}\,(1) now follows from (the left module analogue of) Proposition~\ref{thm35}.
We remark that the bounding cochain $b_1 \times b_2$ is characterized
by the formula
\[
0 =
\sum_{k_1,k_2,k_{12}} \mathfrak n_{k_{12},k_1,k_2} \bigl(\bigl(b_1 \times b_2\bigr)^{k_{12}};1;b_1^{k_1},b_2^{k_2}\bigr).
\]
We turn to the proof of (2).
Let $K_1,\dots,K_k$ be elements of $\mathbb L_1 \times \mathbb L_2$
and $b_{K_i}$ a bounding cochain of $K_i$.
We consider
$
D(K_i) = \mathscr F(L_1,L_2)(K_i)$.
It is a left $CF(K_i)$, right $CF(L_1)$, $CF(L_2)$
tri-module.
Note that we use $b_1$, $b_2$, $b_{K_i}$ to obtain a~left $CF(K_i,b_{K_i})$,
right $CF(L_1,b_1)$, $CF(L_2,b_2)$
tri-module structure on $D(K_i)$, which we write
\[
\mathfrak n_{k_{12},k_1,k_2}^{b_{K_i},b_1,b_2}\colon\
B_{k_{12}}CF(K_i)[1]
 \otimes D(K_i) \otimes B_{k_1}CF(L_1) \otimes B_{k_2}CF(L_2) \to D(K_i).
\]
We also have{\samepage
\begin{equation}\label{form16333}
\mathfrak n_m \colon\
\bigotimes_{i=1}^m CF(K_{i-1},K_i) \otimes D(K_m)
\to D(K_0).
\end{equation}
In \eqref{form16333}, we include corrections by bounding cochains of
$L_i$ and $K_i$.}

We consider two left $\mathfrak{Fukst}(X_1 \times X_2;\mathbb L_1 \times \mathbb L_2)$
modules $\mathfrak D_1$ and $\mathfrak D_2$ as follows.
We write $\mathcal K = (K,b_K), \mathcal K_i = (K_i,b_{K_i})$.
They are elements of $\mathfrak{OB}(\mathfrak{Fukst}(X_1 \times X_2;\mathbb L_1 \times \mathbb L_2))$.
\begin{enumerate}\itemsep=0pt\setlength{\leftskip}{0.40cm}
\item[(M1-1)]
As a $\Lambda_0$ module $\mathfrak D_1(\mathcal K)$ is $D(K)$.
\item[(M1-2)]
The structure operations of the left $\mathfrak{Fukst}(X_1 \times X_2;\mathbb L_1 \times \mathbb L_2)$
module structure are \eqref{form16333}.
\item[(M2-1)]
$
\mathfrak D_2(\mathcal K)
:= CF((K,b_K),(L_1\times L_2,b_1\times b_2))$.
Here the right-hand side is the module of morphisms in
$\mathfrak{Fukst}(X_1 \times X_2;\mathbb L_1 \times \mathbb L_2)$.
\item[(M2-2)]
The structure operations of the left $\mathfrak{Fukst}(X_1 \times X_2;\mathbb L_1 \times \mathbb L_2)$
module structure are the structure operations of the filtered $A_{\infty}$ structure of
$\mathfrak{Fukst}(X_1 \times X_2;\mathbb L_1 \times \mathbb L_2)$.
\end{enumerate}
Note that $\mathfrak D_1$ is nothing but the left filtered $A_{\infty}$ module $\mathscr F((L_1,b_1),(L_2,b_2))$
and $\mathfrak D_2$ is nothing but the left filtered $A_{\infty}$ module $\mathfrak{Yon}(X_1 \times X_2;\mathbb L_1 \times \mathbb L_2)$.
Therefore, the next lemma completes the proof of Proposition~\ref{lem167}.
\begin{lem}\label{lem1688}
$\mathfrak D_1(\mathcal K)$ is homotopy equivalent to $\mathfrak D_2(\mathcal K)$ as left
$\mathfrak{Fukst}(X_1 \times X_2;\mathbb L_1 \times \mathbb L_2)$
modules.
\end{lem}
\begin{proof}
Let
$
{\bf 1} \in D(L_1\times L_2,b_1 \times b_2)
$
be the cyclic element.
We define
\begin{gather*}
\mathfrak T_0(\mathcal K) \colon\ \mathfrak D_2(\mathcal K) \to \mathfrak D_1(\mathcal K), \\
\mathfrak T_{m-1}(\mathcal K_1,\dots,\mathcal K_{m})\colon\ \bigotimes_{i=1}^{m-1} CF(\mathcal K_{i},\mathcal K_{i+1})\mathfrak \otimes D_2(\mathcal K_m)
\to \mathfrak D_1(\mathcal K_1)
\end{gather*}
as follows. We put $\mathcal K_0 = \mathcal K$, $\mathcal K_1 = (L_1\times L_2,b_1 \times b_2)$
in \eqref{form16333} and define
$
\mathfrak T_0(\mathcal K)(z) = \mathfrak n_1(z;{\bf 1})$.
Here $z \in CF(K,L_1\times L_2)$ and ${\bf 1} \in CF(L_1 \times L_2;L_1,L_2)$ is the
cyclic element.

We put $\mathcal K_m = (L_1\times L_2,b_1 \times b_2)$ and take $\mathcal K_i$
for $i=1,\dots,m-1$ in \eqref{form16333} and define
\[
 \mathfrak T_{m-1}(\mathcal K)(x_1,\dots,x_{m-1})(z) = \mathfrak n_{m-1}(x_1,\dots,x_{m-1};z;{\bf 1})
\]
for $z \in D_2(K_m)\in CF(K_m,L_1\times L_2)$, $x_i \in CF(K_i,K_{i+1})$.
See Figure~\ref{Newfiguresec16}.

\begin{figure}[ht]
\centering
\includegraphics[scale=0.4]{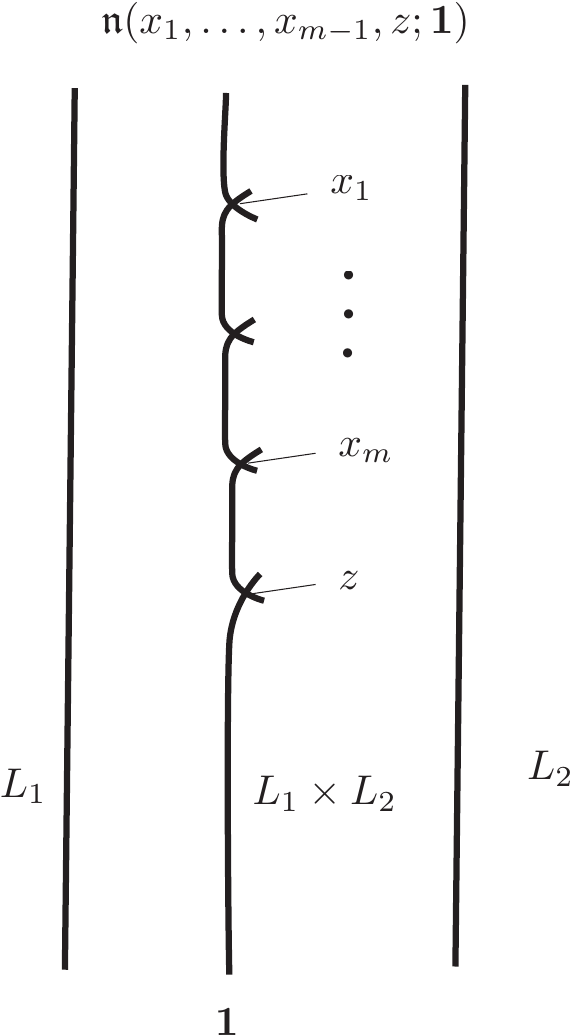}
\caption{$\mathfrak T_{m-1}$.}
\label{Newfiguresec16}
\end{figure}

The $A_{\infty}$ relations for $\{ \mathfrak n_m \}$ imply that
$\{\mathfrak T_i \mid i=0,1,\dots\}$ is a left $\mathfrak{Fukst}(X_1 \times X_2;\mathbb L_1 \times \mathbb L_2)$
module homomorphism.

To show that it is a homotopy equivalence, we first observe that
$\mathfrak D_1(\mathcal K)$, $\mathfrak D_2(\mathcal K)$ both are isomorphic to
$
\Omega\bigl(\tilde K \times_{X_1\times X_2}\bigl(\tilde L_1 \times \tilde L_2\bigr)\bigr)
\,\widehat{\otimes}\, \Lambda_0
$
as $\Lambda_0$ modules.
\big(Here $\tilde K \to X_1 \times X_2$ is the immersed Lagrangian submanifold
which is a part of $\mathcal K$.\big)

We next observe that $\mathfrak T_0(\mathcal K)$ is congruent to the identity map
(via the above identification~${\mathfrak D_1(\mathcal K) = \mathfrak D_2(\mathcal K)}$).
In fact, $\mathfrak T_0 = \mathfrak n_1$ is defined by using the moduli space of pseudo-holomorphic
disks and we put $T^{E}$ as a part of the weight when we use
the moduli space of pseudo-holomorphic disks with symplectic area $E$. The disk
with symplectic area $0$ is nothing but a constant map,
whose contribution to $\mathfrak T_0$ is the identity map.

The proof of Lemma~\ref{lem1688} is complete.
\end{proof}

The proof of Proposition~\ref{lem167} is complete.
\end{proof}

We consider the full subcategory of
$\mathfrak{Fukst}(X_1 \times X_2;\mathbb L_1 \times \mathbb L_2)$ the set of
whose objects consists of $(L_1 \times L_2,b_1 \times b_2)$ where $L_i \in \mathbb L_i$
and $b_1 \times b_2$ is a bounding cochain of $L_1 \times L_2$ obtained by
Lemma~\ref{lem167} from bounding cochains $b_1$ and $b_2$ of $L_1$ and $L_2$.
We denote it by $\mathscr L$.

We put $\mathscr C_i = \mathfrak{Fukst}(X_i,\mathbb L_i)$ and ${\bf C} = \mathscr C_1 \otimes \mathscr C_2$.
The formulas \eqref{form161212} and \eqref{form16132} define a
left $\mathscr C_1$, $\mathscr C_2$ and right ${\bf C}$ filtered
$A_{\infty}$ tri-module ${\bf M}(\mathscr C_1,\mathscr C_2;{\bf C})$.
(See also \eqref{form1613}.)

By \eqref{trifindd1610}, we obtain a left $\mathscr L$ right $\mathscr C_1$, $\mathscr C_2$
tri-module, which we denote by
$
{\bf M}(\mathscr L;\mathscr C_1,\mathscr C_2)
$.

Note that the set of objects of ${\bf C}$ is canonically identified with
the set of objects of $\mathscr L$.
(In fact, they both are $\mathfrak{Ob}\mathfrak{Fukst}(X_1,\mathbb L_1) \times \mathfrak{Ob}\mathfrak{Fukst}(X_2,\mathbb L_2)$.)
For any objects $\mathfrak c$ of ${\bf C}$, the tri-module
${\bf M}({\bf C};\mathscr C_1,\mathscr C_2)$ determines
a right $\mathscr C_1$, $\mathscr C_2$ filtered $A_{\infty}$ bi-module,
which we write
${\bf M}({\bf C};\mathscr C_1,\mathscr C_2)(\mathfrak c;*,*)$.
We define a right $\mathscr C_1$, $\mathscr C_2$ filtered $A_{\infty}$ bi-module
${\bf M}(\mathscr L;\mathscr C_1,\mathscr C_2)(\mathfrak c;*,*)$
in the same way.

\begin{lem}\label{lem1611222211}
For any object $\mathfrak c$ of ${\bf C}$,
the module ${\bf M}({\bf C};\mathscr C_1,\mathscr C_2)(\mathfrak c;*,*)$
 is isomorphic\footnote{Here two right filtered $A_{\infty}$ bi-modules are said
 to be isomorphic each other if there exist a bi-module homomorphism between them which has an inverse.} to
${\bf M}(\mathscr L;\mathscr C_1,\mathscr C_2)(\mathfrak c;*,*)$ as right $\mathscr C_1$, $\mathscr C_2$ filtered $A_{\infty}$
 bi-modules.
\end{lem}
\begin{rem}
Let $\mathfrak c_i = (L_i,b_i), \mathfrak c'_i = (L'_i,b'_i)\in \mathbb L_i$ and
we put $\mathfrak c = (\mathfrak c_1,\mathfrak c_2)$, $\mathfrak c' = (\mathfrak c'_1,\mathfrak c'_2)$.
The chain complex
${\bf M}(\mathscr L;\mathscr C_1,\mathscr C_2)(\mathfrak c;\mathfrak c'_1,\mathfrak c'_2)$
is chain homotopy equivalent to
\[
CF((L_1 \times L_2,b_1 \times b_2),(L'_1 \times L'_2,b'_1 \times b'_2))
\]
by Proposition~\ref{lem167}\,(2).

On the other hand, the chain complex ${\bf M}({\bf C};\mathscr C_1,\mathscr C_2)(\mathfrak c;\mathfrak c'_1,\mathfrak c'_2)$
is chain homotopy equivalent to
$CF((L_1,b_1),(L'_1,b'_1)) \otimes CF((L_2,b_2),(L'_2,b'_2))$, by definition.
Therefore, Lemma~\ref{lem1611222211} implies
\[
CF((L_1 \times L_2,b_1 \times b_2),(L'_1 \times L'_2,b'_1 \times b'_2))
\sim
CF((L_1,b_1),(L'_1,b'_1)) \otimes CF((L_2,b_2),(L'_2,b'_2)),
\]
where $\sim$ means chain homotopy equivalence. This is K\"unneth formula.
The proof of Theorem~\ref{thm1615} shows that this
chain homotopy equivalence is functorial.
\end{rem}
\begin{proof}[Proof of Lemma~\ref{lem1611222211}]
We put $\mathfrak c = (\mathfrak c_1,\mathfrak c_2)$ and $\mathfrak c'_j = (\mathfrak c'_{j,1},\mathfrak c'_{j,2})$ for $j=1,2$.
The structure operations of the right
bi-module structure on ${\bf M}({\bf C};\mathscr C_1,\mathscr C_2)(\mathfrak c;*,*)$,
which we denote by $\mathfrak m'$ is defined by
\begin{gather}
\mathfrak m'(((z_1,z_2),{\bf x}, {\bf y}) :=
\begin{cases}
0 &\text{if ${\bf x} \ne 1$ or ${\bf y} \ne 1$}, \\
(-1)^{\deg'{\bf x}\deg'z_2}( \mathfrak m(z_1,{\bf x}),z_2) &\text{if ${\bf y} = 1$}, \\
(-1)^{\deg' z_1}(z_1, \mathfrak m(z_2,{\bf y})) &\text{if ${\bf x} = 1$}.
\end{cases}\label{form1613}
\end{gather}
Here $\mathfrak m$ is the structure operations of $\mathscr C_i$, which is obtained
by `counting' pseudo-holomorphic polygons,
${\bf x} \in B\mathscr C_1[1](\mathfrak c'_{1,1},\mathfrak c'_{2,1})$,
${\bf y} \in B\mathscr C_2[1](\mathfrak c'_{1,2},\mathfrak c'_{2,2})$
and $z_i \in \mathscr C_i[1](\mathfrak c_i,\mathfrak c'_{1,i})$ for $i=1,2$.

On the other hand,
the structure operations of the bi-module structure on ${\bf M}(\mathscr L;\mathscr C_1,\mathscr C_2)(\mathfrak c;\allowbreak*,*)$
is the operations $\{\mathfrak n\}$ which are structure operations of the tri-module \eqref{trifindd1610}
and obtained by `counting' pseudo-holomorphic quilts.

For $\mathfrak b_i, \mathfrak b'_i\in \mathfrak{Ob}\mathfrak{Fukst}(X_i,\mathbb L_i)$, we define
\begin{gather*}
\mathscr F_{k_1,k_2} \colon\ {\bf M}({\bf C};\mathscr C_1,\mathscr C_2)(\mathfrak c;\mathfrak b_1,\mathfrak b_2)
\otimes
B_{k_1}\mathscr C_1[1](\mathfrak b_1,\mathfrak b'_1) \otimes B_{k_2}\mathscr C_2[1](\mathfrak b_2,\mathfrak b'_2)
\\
\hphantom{\mathscr F_{k_1,k_2} \colon} \ \to {\bf M}(\mathscr L;\mathscr C_1,\mathscr C_2)(\mathfrak c;\mathfrak b'_1,\mathfrak b'_2)
\end{gather*}
by
\begin{gather}
\mathscr F_{k_1,k_2}((z_1,z_2),({\bf x},{\bf y})):=
\sum_c(-1)^{*}\mathfrak n(\mathfrak n({\bf 1};1\otimes (z_2\otimes {\bf y}_c));(z_1 \otimes {\bf x})\otimes {\bf y}'_c
),
\label{formlast1614}
\end{gather}
where
$
* =\deg' z_1(\deg'z_2 + \deg'{\bf y}_c) + \deg'{\bf x}\deg'{\bf y}_c
$
is the Koszul sign.
Here ${\bf x} \in B_{k_1}\mathscr C_1[1](\mathfrak b_1,\mathfrak b'_1)$,
${\bf y} \in B_{k_2}\mathscr C_2[1](\mathfrak b_2,\mathfrak b'_2)$,
and $z_i \in \mathscr C_i[1](\mathfrak c_{i},\mathfrak b_i)$.
The symbol ${\bf 1}$ is the fundamental class of ${\bf M}(\mathscr L;\mathscr C_1,\allowbreak\mathscr C_2)(\mathfrak c;(\mathfrak c_1,\mathfrak c_2))$,
that is, the cyclic element
and $
\Delta({\bf y}) = \sum_c {\bf y}_{c} \otimes {\bf y}'_c
$.
The idea behind this definition can been seen from Figure~\ref{FigureSec16} below.
\begin{figure}[ht]
\centering
\includegraphics[scale=0.4]{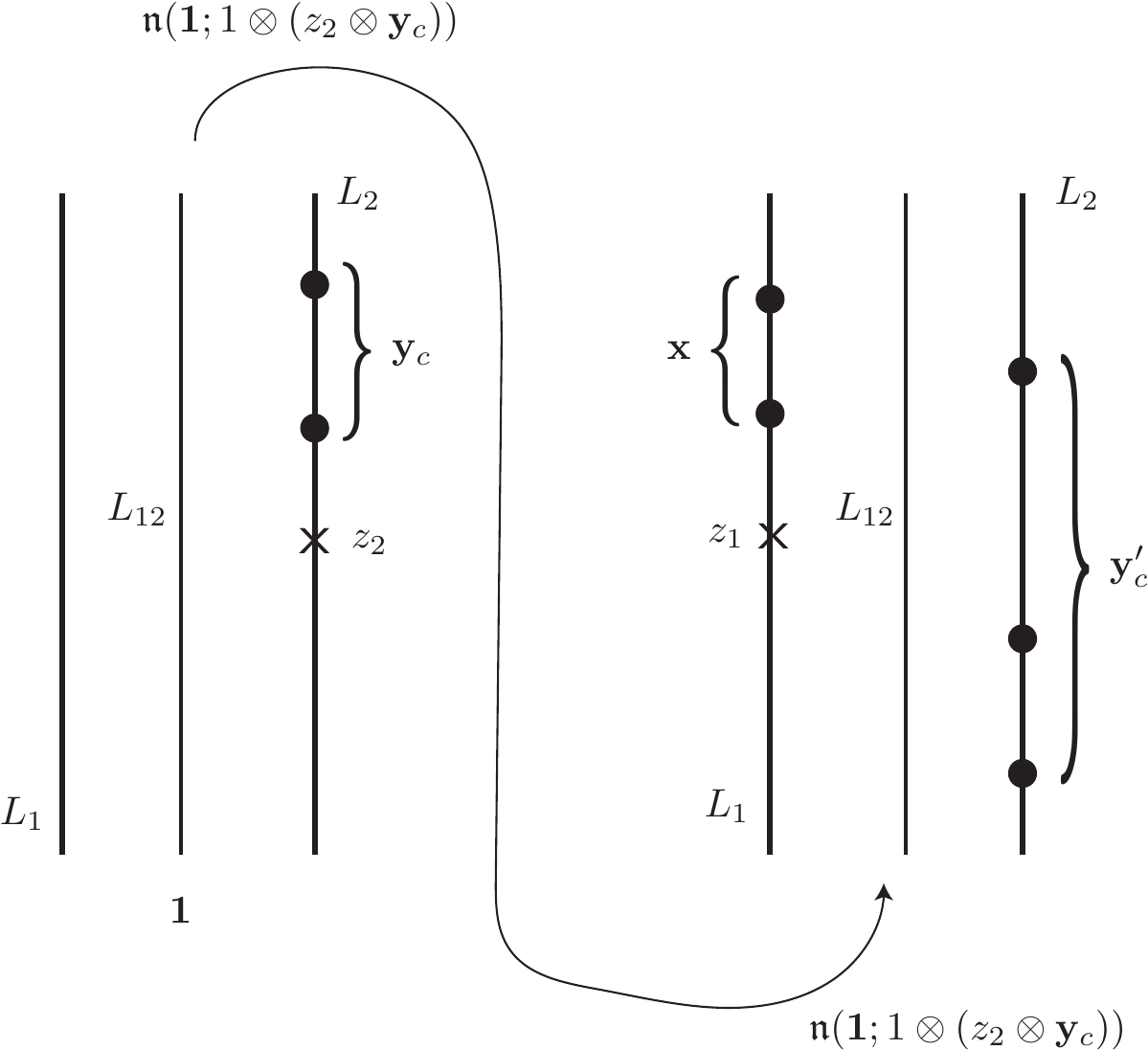}
\caption{$\mathscr F_{k_1,k_2}((z_1,z_2),({\bf x},{\bf y}))$.}
\label{FigureSec16}
\end{figure}

We will prove that $\{\mathscr F_{k_1,k_2}\}$ is a right filtered $A_{\infty}$ bi-module homomorphism.
We denote by \smash{$\widehat{\mathfrak m'}$} the maps
\begin{gather*}
\bigoplus_{\mathfrak b_1,\mathfrak b_2}{\bf M}({\bf C};\mathscr C_1,\mathscr C_2)(\mathfrak c;\mathfrak b_1,\mathfrak b_2)
\otimes
B\mathscr C_1[1](\mathfrak b_1,\mathfrak c'_1) \otimes B\mathscr C_2[1](\mathfrak b_2,\mathfrak c'_2)
\\
\qquad\to
\bigoplus_{\mathfrak b_1,\mathfrak b_2}{\bf M}({\bf C};\mathscr C_1,\mathscr C_2)(\mathfrak c;\mathfrak b_1,\mathfrak b_2)
\otimes
B\mathscr C_1[1](\mathfrak b_1,\mathfrak c'_1) \otimes B\mathscr C_2[1](\mathfrak b_2,\mathfrak c'_2)
\end{gather*}
induced by $\mathfrak m'$ (see \eqref{form1613}) and denote by $\mathfrak m$
the structure operations of $\mathscr C_i$.
We also denote~by
\begin{gather*}
\widehat{\mathscr F} \colon\
\bigoplus_{\mathfrak b_1,\mathfrak b_2}{\bf M}({\bf C};\mathscr C_1,\mathscr C_2)(\mathfrak c;\mathfrak b_1,\mathfrak b_2)
\otimes
B\mathscr C_1[1](\mathfrak b_1,\mathfrak c'_1) \otimes B\mathscr C_2[1](\mathfrak b_2,\mathfrak c'_2) \\
\hphantom{\widehat{\mathscr F} \colon} \
\to
\bigoplus_{\mathfrak b_1,\mathfrak b_2}{\bf M}(\mathscr L;\mathscr C_1,\mathscr C_2)(\mathfrak c;\mathfrak b_1,\mathfrak b_2)
\otimes
B\mathscr C_1[1](\mathfrak b_1,\mathfrak c'_1) \otimes B\mathscr C_2[1](\mathfrak b_2,\mathfrak c'_2)
\end{gather*}
the map induced by $\mathscr F_{k_1,k_2}$.

We will check
\begin{equation}\label{formula1614}
\bigl({\mathfrak n} \circ \widehat{\mathscr F}\bigr)((z_1,z_2),({\bf x},{\bf y}))
= \bigl(\mathscr F \circ \widehat{\mathfrak m'}\bigr)((z_1,z_2),({\bf x},{\bf y})).
\end{equation}
\begin{rem}
Intuitively \eqref{formula1614} can be proved easily by studying the boundary of the moduli space
depicted by Figure~\ref{FigureSec16}.
The proof below is an algebraic analogue of such a geometric argument.
\end{rem}
Let $\widehat{\mathfrak m_i} \colon B\mathscr C_i \to B\mathscr C_i $ is the map induced by
the structure operations of $\mathscr C_i$.
We denote $\widehat{\mathfrak m} = \widehat{\mathfrak m_1} \,\widehat{\otimes}\, {\rm id} +{\rm id} \,\widehat{\otimes}\, \widehat{\mathfrak m_2}$.
Let $\mathfrak m$ be the $\mathscr C_1 \otimes \mathscr C_2$ component of
$\widehat{\mathfrak m}$.

We put $
\Delta({\bf x}) = \sum_{b} {\bf x}_{b} \otimes {\bf x}'_{b}
$
,
$((\Delta\otimes {\rm id})\circ \Delta)({\bf y})
= \sum_{c'} {\bf y}_{c'} \otimes {\bf y}'_{c'} \otimes {\bf y}''_{c'}$.

Now the right-hand side of \eqref{formula1614} is
\begin{gather}
\sum_c(-1)^{*_1}\mathfrak n\bigl(\mathfrak n({\bf 1};1\otimes (z_2\otimes {\bf y}_c));\bigl(z_1 \otimes \widehat{\mathfrak m}(
{\bf x})\bigr)\otimes {\bf y}'_c
\bigr)\nonumber\\
\qquad+
\sum_c(-1)^{*_2}\mathfrak n\bigl(\mathfrak n({\bf 1};1\otimes (z_2\otimes {\bf y}_c));(z_1\otimes {\bf x})\otimes \widehat{\mathfrak m}({\bf y}'_c)\bigr)\nonumber
\\
\qquad+
\sum_c(-1)^{*_3}\mathfrak n\bigl(\mathfrak n({\bf 1};1\otimes \bigl(z_2\otimes \widehat{\mathfrak m}
({\bf y}_{c}))\bigr);(z_1\otimes {\bf x})\otimes {\bf y}'_{c}\bigr)\nonumber
\\
\qquad+
\sum_{c,c'}(-1)^{*_4}\mathfrak n(\mathfrak n({\bf 1};1\otimes (
\mathfrak m(z_2\otimes {\bf y}_{c'})\otimes {\bf y}'_{c'});(z_1\otimes {\bf x})\otimes {\bf y}''_{c'})\nonumber
\\
\qquad+\sum_{b,c}
(-1)^{*_5} \mathfrak n({\bf 1};\mathfrak n(1\otimes (z_2\otimes {\bf y}_c));(\mathfrak m(z_1
\otimes {\bf x}_b)
\otimes {\bf x}'_b )\otimes {\bf y}'_c
),\label{form16151615}
\end{gather}
where
the signs $*_i$ are by Koszul rule.
Note that the sign here is always by Koszul rule and so it is actually not necessary
to calculate the sign as we will explain at the end of the proof of~\eqref{formula1614}.

On the other hand, the left-hand side of \eqref{formula1614} is
\begin{equation}\label{form161616}
\sum_{b,c'}
\mathfrak n(\mathfrak n(\mathfrak n({\bf 1};1\otimes (z_2\otimes {\bf y}_{c'}));(z_1 \otimes{\bf x}_b) \otimes {\bf y}'_{c'});{\bf x}'_{b}\otimes {\bf y}''_{c'})
\end{equation}
up to sign.
We put
$
z_{c} := \pm \mathfrak n({\bf 1};1\otimes (z_2 \otimes {\bf y}_{c}))$,
where $\Delta({\bf y}) = \sum_c {\bf y}_{c} \otimes {\bf y}'_{c}$.

\eqref{form161616}
plus the next formula \eqref{form16171617222}
 is obtained by applying $\mathfrak n$ twice
to the element $\sum_{c} z_{c}\otimes ((z_1 \otimes {\bf x}) \otimes {\bf y}'_{c} ) )$ up to sign,
\begin{equation}
\sum_{c''}
\mathfrak n(\mathfrak n(\mathfrak n({\bf 1};1 \otimes (z_2 \otimes {\bf y}_{c})); 1\otimes {\bf y}'_{c'});
(z_1\otimes {\bf x}) \otimes {\bf y}''_{c'}).\label{form16171617222}
\end{equation}

Therefore, applying $A_{\infty}$ formula for $\mathfrak n$, the formula
\eqref{form161616} is equal to
\begin{gather}
+\sum_{c}
\mathfrak n\bigl(z_{c};\widehat{\mathfrak m}(z_1\otimes {\bf x})\otimes {\bf y}'_{c}\big)
+\sum_{c}
\mathfrak n(z_{c};(z_1\otimes {\bf x}) \otimes
\widehat{\mathfrak m}({\bf y}'_{c}))
+ \eqref{form16171617222}\label{form16181811}
\end{gather}
up so sign.
This cancels with \eqref{form16151615} up to sign.
In fact, the first term of \eqref{form16181811}
cancel with the first and fifth terms of \eqref{form16151615},
the second term of \eqref{form16181811}
cancel with the second term of~\eqref{form16151615},
and \eqref{form16171617222} cancel with
the third and fourth terms of \eqref{form16151615}, using $A_{\infty}$ relation
(of left bi-module structure)
applied to ${\bf 1} \otimes (1 \otimes (z_2 \otimes {\bf y}_c))$.

The calculation of sign looks complicated.
However, we actually do {\it not} need to check
the sign by calculation to see the corresponding terms cancel out with sign.
In fact, all the signs are caused by changing the order of operators
or elements (which are graded), that is by Koszul rule. So except
the minus sing which comes from exchanging the order of $\mathfrak m$ and
$\mathfrak n$ (both of which have degree $\pm 1$) the
sign of the corresponding terms coincide.
Therefore, the cancellation occurs with signs.\footnote{This is the standard magic
of Koszul sign.
}

Thus $\{\mathscr F_{k_1,k_2}\}$ defines a filtered bi-module homomorphism.

We remark that ${\bf M}({\bf C};\mathscr C_1,\mathscr C_2)(\mathfrak c;\mathfrak c')$
as $\Lambda_0$ module is a $T$-adic completion of the tensor product
of de Rham complex $\Omega\bigl(\tilde L_1 \times_{X_1} \tilde L'_1\bigr) \,\widehat{\otimes}\,
\Lambda_0$ and $\Omega\bigl(\tilde L_2 \times_{X_2} \tilde L'_2\bigr) \,\widehat{\otimes}\,
\Lambda_0$.
On the other hand, ${\bf M}(\mathscr L;\mathscr C_1,\mathscr C_2)(\mathfrak c;\mathfrak c')$
is \smash{$\Omega\bigl(\bigl(\tilde L_1 \times \tilde L_2\bigr) \times_{X_1
\times X_2} \bigl(\tilde L'_1 \times \tilde L'_2\bigr)\bigr) \,\widehat{\otimes}\,
\Lambda_0$}.
Therefore, their reductions to the ground ring is isomorphic each other.
It is easy to see that $R$ reduction of $\mathscr F_{0,0}$
is this isomorphism.
Therefore, $\mathscr F_{0,0}$ is an isomorphism.
Hence $\mathscr F$ is an isomorphism.
\end{proof}

We recall that
${\bf C} = \mathscr C_1 \otimes \mathscr C_2$ can be regarded as the category of
right bi-module homomorphisms~${{\bf M}({\bf C};\mathscr C_1,\mathscr C_2) \to
{\bf M}({\bf C};\mathscr C_1,\mathscr C_2)}$ in the following sense.
An object of ${\bf C}$ is identified with
a pair of objects $(\mathfrak c_1,\mathfrak c_2)$ of $\mathscr C_1$ and of $\mathscr C_2$.
For a fix $(\mathfrak c_1,\mathfrak c_2)$, by moving $(\mathfrak c'_1,\mathfrak c'_2)$
this defines a right~$\mathscr C_1$,~$\mathscr C_2$ module, which is nothing but
${\bf M}({\bf C};\mathscr C_1,\mathscr C_2)(*,*,(\mathfrak c_1,\mathfrak c_2))$.
The morphisms and operations in ${\bf C}$ are defined to be the
right $\mathscr C_1$, $\mathscr C_2$ bi-module homomorphisms and their compositions.\footnote{
Here we use the operation taking opposite module (see Definition~\ref{definition113}) to go from
left module to right module and vice versa.}

By Lemma~\ref{lem1611222211}, ${\bf M}({\bf C};\mathscr C_1,\mathscr C_2)$ is
isomorphic to ${\bf M}(\mathscr L;\mathscr C_1,\mathscr C_2)$.
Now using the fact that
${\bf M}(\mathscr L;\mathscr C_1,\mathscr C_2)$ is a
left $\mathscr L$,
right $\mathscr C_1$, $\mathscr C_2$ tri-module, we obtain a filtered $A_{\infty}$
functor
$
\mathscr G \colon \mathscr L \to \mathscr C_1 \otimes \mathscr C_2$.
Note that the object part of $\mathscr G$ is the identity map.

\begin{lem}\label{G1chainlem}
The linear part $\mathscr G_1$ of $\mathscr G$ is a chain homotopy equivalence from $\mathscr L((\mathfrak c'_1,\mathfrak c'_2),(\mathfrak c_1,\mathfrak c_2))$ to
${\bf C}((\mathfrak c'_1,\mathfrak c'_2),(\mathfrak c_1,\mathfrak c_2))$.
\end{lem}
We prove Lemma~\ref{G1chainlem} at the end of Section~\ref{reltotrimod}.
Lemma~\ref{G1chainlem} implies $\mathscr G$ is a homotopy equivalence.
The proof of Theorem~\ref{thm1615} is complete.
\end{proof}

\subsection{The K\"unneth functor and the correspondence tri-module}
\label{reltotrimod}

Suppose that we are in Situation \ref{situ14}.
We consider the set $\mathbb L_1 \times \mathbb L_2$
of Lagrangian submanifolds of $-X_1 \times X_2$ which
consists of direct products $L_1 \times L_2$ of elements
$L_1 \in \mathbb L_1$ and $L_2 \in \mathbb L_2$.
The K\"unneth functor defines
\begin{gather}
\mathscr K\colon\ \mathfrak{Fukst}((-X_1,V_1\oplus TX_1);\mathbb L_1)
\times \mathfrak{Fukst}((X_2,V_2);\mathbb L_2)\nonumber\\
\hphantom{\mathscr K\colon} \
\to
\mathfrak{Fukst}((-X_1\times X_2,\pi_1^*(V_1\oplus TX_1)\oplus \pi_2^*(V_2));\mathbb L_1
\times \mathbb L_2).\label{formunew1620}
\end{gather}
Note that we replace $X_1$, $V_1$ by $-X_1, V_1\oplus TX_1$ when we apply Theorem~\ref{thm1615}
to obtain \eqref{formunew1620}.

\begin{thm}\label{thm164}
Let $L_1 \in \mathbb L_1$, $L_{12} \in \mathbb L_{12}$
and $b_1$, $b_{12}$ their bounding cochains.
We put $\mathcal L_1 = (L_1,b_1)$, $\mathcal L_{12} = (L_{12},b_{12})$.
Using correspondence bi-functor, we obtain
$(L_2,b_2) = \mathcal W_{\mathcal L_{12}}(\mathcal L_1)$.
Then we have the following isomorphism for any
$L'_2 \in \mathbb L_2$ and its bounding cochain $b'_2$:
\[
HF((L_2,b_2);(L'_2,b'_2);\Lambda_0)
\cong
HF((L_{12},b_{12});(L_1\times L'_2,b_1\times b'_2);\Lambda_0).
\]
Here $(L_1\times L'_2,b_1\times b'_2): = \mathscr K_{\rm ob}((L_1,b_1),(L'_2,b'_2))$.
\end{thm}
\begin{proof}
In our situation, where $X_1$, $V_1$ are replaced by $-X_1$, $V_1\oplus TX_1$,
the tri-module \eqref{trifindd1610} is the correspondence tri-module.
Therefore, the theorem is an immediate consequence of
Theorems~\ref{th72} and~\ref{thm1615}.
\end{proof}

\subsection[Proof of Lemmas \ref{DGcaselem}, \ref{lem1613}
and \ref{G1chainlem}]{Proof of Lemmas~\ref{DGcaselem}, \ref{lem1613}
and \ref{G1chainlem}}
\label{prooflemmata}

In this subsection, we prove Lemmas~\ref{DGcaselem}, \ref{lem1613}
and \ref{G1chainlem}.
It suffices to consider the case when $\mathscr C_1$, $\mathscr C_2$ are DG-categories.

\begin{proof}[Proof of Lemma~\ref{lem1613}]
We prove the first half. The proof of the second half is similar.
We define a tri-module homomorphism
$
\mathfrak I \colon {\bf M}(\mathscr C_1,\mathscr C_2;{\bf C}) \to {\bf M}({\bf C};{\bf C})
$
as follows.

We first define
\begin{align*}
& \mathfrak n \colon\ (B\mathscr C_1[1](b_{1,1},b_{2,1})  \otimes B\mathscr C_2[1](b_{1,2},b_{2,2}))
\otimes
\mathscr C_1[1](b_{2,1},c_{1,1}) \otimes \mathscr C_2[1](b_{2,2},c_{1,2})\otimes {\bf C}(\mathfrak c_1,\mathfrak c_2)
\\
&\hphantom{\mathfrak n \colon} \
\to
\mathscr C_1[1](c_{1,1},c_{2,1}) \otimes \mathscr C_2[1](c_{1,2},c_{2,2}),
\end{align*}
as follows.
Note an element $\mathcal T \in {\bf C}(\mathfrak c_1,\mathfrak c_2)$ is a
pre-natural transformation from ${\bf C}(\mathfrak c_1)$ to ${\bf C}(\mathfrak c_2)$.
Such pre-natural transformation assigns to each $\mathfrak b_1= (b_{1,1},b_{1,2})$,
$\mathfrak b_2 = (b_{2,1},b_{2,2})$ a map
\begin{gather*}
(B\mathscr C_1[1](b_{1,1},b_{2,1}) \otimes B\mathscr C_2[1](b_{1,2},b_{2,2}))
\otimes
\mathscr C_1[1](b_{2,1},c_{1,1}) \otimes \mathscr C_2[1](b_{2,2},c_{1,2})\\
\qquad
\to
\mathscr C_1[1](c_{1,1},c_{2,1}) \otimes \mathscr C_2[1](c_{1,2},c_{2,2}).
\end{gather*}
For ${\bf x} \otimes {\bf y} \otimes z \in (B\mathscr C_1[1](b_{1,1},b_{2,1}) \otimes B\mathscr C_2[1](b_{1,2},b_{2,2}))
\otimes \mathscr C_1[1](b_{2,1},c_{1,1}) \otimes \mathscr C_2[1](b_{2,2},c_{1,2})$,
we denote by
$
\mathfrak n({\bf x} \otimes {\bf y},z,\mathcal T)
$
the image of ${\bf x} \otimes {\bf y}\otimes z$ by this map.

We next define
$
\mathfrak I_{0,0,0}(\mathfrak c_1,\mathfrak c_2) \colon \mathscr C_1[1](c_{1,1},c_{1,2}) \otimes \mathscr C_2[1](c_{2,1},c_{2,2})
\to {\bf C}(\mathfrak c_1,\mathfrak c_2)
$
by the formula
\begin{gather}
\mathfrak n({\bf x} \otimes {\bf y}, z, \mathfrak I_{0,0,0}(\mathfrak c_1,\mathfrak c_2)(a_1,a_2))\nonumber \\
\qquad=
\begin{cases}
(-1)^{\deg' a_1\deg'z_2} (\mathfrak m_2(z_1,a_1),\mathfrak m_2(z_2,a_2)) & \text{if ${\bf x} \otimes {\bf y} = 1\otimes 1$}, \\
0 & \text{otherwise.}
\end{cases}\label{1618formoo}
\end{gather}
Here $\mathfrak n$ is defined as above\footnote{We remark
that the pre-natural transformation $\mathcal T$ is determined if
$\mathfrak n({\bf x} \otimes {\bf y}, z,\mathcal T)$ are given for all ${\bf x}$, ${\bf y}$, $z$.} and
$z = (z_1,z_2) \in \mathscr C_1[1](b_{2,1},c_{1,1}) \otimes \mathscr C_2[1](b_{2,2},c_{1,2})$,
${\bf x} \otimes {\bf y} \in B\mathscr C_1[1](b_{1,1},b_{2,1}) \otimes B\mathscr C_2[1](b_{1,2},b_{2,2})$,
$(a_1,a_2) \in \mathscr C_1[1](c_{1,1},c_{1,2}) \otimes \mathscr C_2[1](c_{2,1},c_{2,2})$.

Hereafter, we write $\mathfrak I_{0,0,0}$ in place of $\mathfrak I_{0,0,0}(\mathfrak c_1,\mathfrak c_2)$.
We define all other $\mathfrak I_{k_1,k_2;\ell}$ to be $0$.
\begin{sublem}
$\mathfrak I\colon {\bf M}(\mathscr C_1,\mathscr C_2;{\bf C}) \to {\bf M}({\bf C};{\bf C}) $ is a tri-module homomorphism.
\end{sublem}
\begin{proof}
Since $\mathscr C_i$ and ${\bf C}$ are DG categories, \eqref{1618formoo} implies that
$
\mathfrak n({\bf x} \otimes {\bf y};z;\delta(\mathfrak I_{0,0,0}(a_1,a_2))) = 0
$
unless ${\bf x}\otimes {\bf y} \in B_{k_1}\mathscr C_1[1]
\otimes B_{k_2}\mathscr C_2[1]$ with $(k_1,k_2) = (0,0),(1,0),(1,1)$.
Here $\delta$ is the boundary operator of ${\bf C}$.
In case $(k_1,k_2) = (1,0)$, we calculate\footnote{In the calculation below,
we omit the Koszul sign unless otherwise mentioned.}
\begin{align*}
\mathfrak n(x \otimes 1,z,\delta(\mathfrak I_{0,0,0}(a_1,a_2))
={}&
(\mathfrak m_2(\mathfrak m_2(x,z_1),a_1),\mathfrak m_2(z_2,a_2))\\
&
+
(\mathfrak m_2(x,\mathfrak m_2(z_1,a_1)),\mathfrak m_2(z_2,a_2))\\
={}&0
=
\mathfrak n(x \otimes 1,z,\mathfrak I_{0,0,0}(\mathfrak m_1(a_1,a_2))).
\end{align*}
Note that the second equality follows from the fact that the
product structures on $\mathscr C_i$ are strictly associative.

The case $(k_1,k_2) = (0,1)$ is similar.
In case $(k_1,k_2) = (0,0)$, we calculate
\begin{align*}
\mathfrak n(1\otimes 1,z,\delta(\mathfrak I_{0,0,0}(a_1,a_2)))
={}&
\mathfrak m_1(\mathfrak n(1\otimes 1,z,\mathfrak I_{0,0,0}(a_1,a_2)))
- \mathfrak n(1\otimes 1,\mathfrak m_1(z),\mathfrak I_{0,0,0}(a_1,a_2)))\\
={}& (\mathfrak m_1(\mathfrak m_2(z_1,a_1)),\mathfrak m_2(z_2,a_2)) +
(\mathfrak m_2(z_1,a_1),\mathfrak m_1(\mathfrak m_2(z_2,a_2)))\\
&
-(\mathfrak m_2(\mathfrak m_1(z_1),a_1),\mathfrak m_2(z_2,a_2))
- (\mathfrak m_2(z_1,a_1),\mathfrak m_2(\mathfrak m_1(z_2),a_2))\\
={}&
\mathfrak n(1\otimes 1;z;\mathfrak I_{0,0,0}(\mathfrak m_1(a_1,a_2))).
\end{align*}
We thus proved that $\mathfrak I_{0,0,0}$ is a chain map.
Note that the fact the equality holds with the sign since
we always use Koszul sign here, as we mentioned during the proof of
\eqref{formula1614}.

The calculation to show that $\mathfrak I$ is a tri-module homomorphism
is similar. We omit it.
\end{proof}

We remark that $\mathfrak I$ is a bijection on objects.
So to prove Lemma~\ref{lem1613}, it suffices
to show that~$\mathfrak I_{0,0,0}$ is a chain homotopy equivalence.
We will prove it below.\footnote{The proof below is
similar to a proof of $A_{\infty}$ Yoneda's lemma.
The key idea of the proof of Yoneda's lemma is plugging in the identity
morphisms to obtain an inverse to the Yoneda embedding. We follow this idea. In fact~\eqref{defnJJJJ}
is nothing but plugging in the identity morphism.
In the case of usual Yoneda's lemma then it defines an inverse. In the
$A_{\infty}$ case, it is only a chain homotopy inverse. So we need to
find a chain homotopy. The chain homotopy~\eqref{form16223}
is similar to one given in \cite{fu4} for the $A_{\infty}$ Yoneda's lemma.
} We write $\mathfrak I$ in place of $\mathfrak I_{0,0,0,}$
for simplicity.

We define a map
$
\mathfrak J \colon {\bf M}({\bf C};{\bf C})(\mathfrak c,\mathfrak b) \to {\bf M}(\mathscr C_1,\mathscr C_2;{\bf C})(\mathfrak c,\mathfrak b)
$
by
\begin{equation}\label{defnJJJJ}
\mathfrak J(\mathcal T) = \mathfrak n(1;({\bf e},{\bf e});\mathcal T).
\end{equation}
Here ${\bf e}$ is the unit of $\mathscr C_i$ and $\mathscr T \in {\bf C}$.
\begin{sublem}
$\mathfrak J$ is a chain map.
\end{sublem}
\begin{proof}
$
\mathfrak m_1(\mathfrak J(\mathcal T))
=
\mathfrak m_1(\mathfrak n(1;\mathfrak m_1({\bf e},{\bf e});\mathcal T))
+
\mathfrak n(1;({\bf e},{\bf e});\delta(\mathcal T))
=
\mathfrak n(1;({\bf e},{\bf e});\delta(\mathcal T))$.
\end{proof}

We calculate
$
\mathfrak J(\mathfrak I(a_1,a_2))
=
\mathfrak n(1;({\bf e},{\bf e});\mathfrak I(a_1,a_2))
= (a_1,a_2)$.
Therefore, $\mathfrak J\circ \mathfrak I = {\rm id}$.
We finally prove $\mathfrak I\circ \mathfrak J$
is chain homotopic to the identity map.
We define the maps
$
\mathcal H_i \colon {\bf C}(\mathfrak c_1,\mathfrak c_2) \to {\bf C}(\mathfrak c_1,\mathfrak c_2)
$
by the next formula
\begin{gather}
\mathfrak n(({\bf x}\otimes {\bf y});z;\mathcal H_1(\mathcal T))
=
(-1)^{*_1}\mathfrak n((({\bf x} \otimes z_1)\otimes {\bf y});({\bf e}\otimes z_2);\mathcal T), \nonumber\\
\mathfrak n(({\bf x}\otimes {\bf y});z;\mathcal H_2(\mathcal T))
=
(-1)^{*_2}\mathfrak n(({\bf x}\otimes ({\bf y} \otimes z_2));(z_1\otimes{\bf e});\mathcal T),\label{form16223}
\end{gather}
where $*_1 = \deg'z_2 + \deg' {\bf x} \deg'z_1$,
$*_2 = \deg'z_1 \deg' z_2$.
We will calculate $\delta \circ \mathcal H_i + \mathcal H_i\circ \delta$,
where $\delta$ is the boundary operator of ${\bf C}$.
We define maps
$
\Phi_i\colon {\bf C}(\mathfrak c_1,\mathfrak c_2) \to {\bf C}(\mathfrak c_1,\mathfrak c_2)
$
by the next formula
\begin{gather*}
\mathfrak n(({\bf x}\otimes {\bf y});z;\Phi_1(\mathcal T)) =
\begin{cases}
(-1)^{*_3}\mathfrak n(z_1 \otimes 1;\mathfrak n(1 \otimes {\bf y};({\bf e}\otimes z_2);\mathcal T))
&\text{if ${\bf x}=1$}, \\
0 & \text{otherwise},
\end{cases}\\
\mathfrak n(({\bf x}\otimes {\bf y});z;\Phi_2(\mathcal T)) =
\begin{cases}
(-1)^{*_4}\mathfrak n(1 \otimes z_2;\mathfrak n({\bf x} \otimes 1;(z_1 \otimes {\bf e});\mathcal T))
&\text{if ${\bf y}=1$}, \\
0 & \text{otherwise},
\end{cases}
\end{gather*}
where $*_3 = \deg'z_1 \deg' {\bf y} + (\deg'{\bf y} + \deg'z_2) + \deg'z_2$,
$*_4 = \deg'z_2 \deg' {\bf x} + \deg'{\bf x}\deg'z_1$.
\begin{sublem}\label{subem1618}
$
\delta \circ \mathcal H_i + \mathcal H_i\circ \delta
= {\rm id} + \Phi_i$.

\end{sublem}
\begin{proof}
We write ${\bf x} = x_f \otimes {\bf x}_R = {\bf x}_L\otimes x_l$.
and ${\bf y} = y_f \otimes {\bf y}_R = {\bf y}_L\otimes y_l$.
We first calculate omitting all the signs
\begin{gather*}
\mathfrak n({\bf x}\otimes {\bf y};z;(\mathcal H_1 \circ \delta)(\mathcal T))\\
\qquad= \mathfrak n(({\bf x}\otimes z_1) \otimes {\bf y};({\bf e}\otimes z_2);\delta(\mathcal T))
= \mathfrak m(\mathfrak n(({\bf x}\otimes z_1) \otimes {\bf y};({\bf e}\otimes z_2);\mathcal T)) \\
\phantom{\qquad= }{}+
\mathfrak n^{\mathscr C}
(x_f \otimes 1;\mathfrak n(({\bf x}_R\otimes z_1) \otimes {\bf y};({\bf e}\otimes z_2);\mathcal T)) +\mathfrak n(({\bf x}\otimes {\bf y});z;\Phi_1(\mathcal T))\\
\phantom{\qquad= }{}+\mathfrak n^{\mathscr C}(1 \otimes y_f;\mathfrak n(({\bf x} \otimes z_1) \otimes {\bf y}_R;({\bf e}\otimes z_2);\mathcal T)) + \mathfrak n((\widehat{\mathfrak m}({\bf x}) \otimes z_1)
\otimes {\bf y}; ({\bf e}\otimes z_2);\mathcal T) \\
\phantom{\qquad= }{}+ \mathfrak n(({\bf x} \otimes {\mathfrak m}(z_1))
\otimes {\bf y}; ({\bf e}\otimes z_2);\mathcal T)
+ \mathfrak n(({\bf x}_L \otimes \mathfrak m_2(x_l,z_1))
\otimes {\bf y}; ({\bf e}\otimes z_2);\mathcal T) \\
\phantom{\qquad= }{}+ \mathfrak n\bigl(({\bf x} \otimes z_1)
\otimes \widehat{\mathfrak m}({\bf y}); ({\bf e}\otimes z_2);\mathcal T\bigr)+ \mathfrak n(({\bf x} \otimes z_1)
\otimes {\bf y}; ({\bf e}\otimes \mathfrak m_1(z_2));\mathcal T)\\
\phantom{\qquad= }{}+ \mathfrak n({\bf x}
\otimes {\bf y}; (z_1 \otimes z_2);\mathcal T)
.
\end{gather*}
Here $\mathfrak n^{\mathscr C}$ is the structure operation of left ${\mathscr C}_1$, ${\mathscr C}_2$
bimodule structure on ${\mathscr C}_1 \otimes {\mathscr C}_2$.
($\mathfrak n$ is defined at the beginning of the proof of Lemma~\ref{lem1613}.)

We also calculate
\begin{gather*}
\mathfrak n({\bf x}\otimes {\bf y};z;(\delta \circ \mathcal H_1)(\mathcal T)) \\
\qquad= \mathfrak m(\mathfrak n(({\bf x}\otimes z_2) \otimes {\bf y};({\bf e}\otimes z_1);\mathcal T)) +
\mathfrak n^{\mathscr C}(x_f \otimes 1;\mathfrak n(({\bf x}\otimes z_1) \otimes {\bf y};({\bf e}\otimes z_2);\mathcal T)) \\
\phantom{\qquad= }{}+\mathfrak n^{\mathscr C}(1 \otimes y_f;\mathfrak n(({\bf x} \otimes z_1) \otimes {\bf y}_R;({\bf e}\otimes z_2);\mathcal T))
+ \mathfrak n\bigl(\widehat{\mathfrak m}({\bf x} \otimes {\bf y}) \otimes (z_1 \otimes 1)
; ({\bf e}\otimes z_2);\mathcal T\bigr) \\
\phantom{\qquad= }{}+ \mathfrak n(({\bf x} \otimes {\mathfrak m}(z_1))
\otimes {\bf y}; ({\bf e}\otimes z_2);\mathcal T)
+ \mathfrak n(({\bf x} \otimes z_1)
\otimes {\bf y}; ({\bf e}\otimes \mathfrak m_1(z_2));\mathcal T)\\
\phantom{\qquad= }{}+ \mathfrak n(({\bf x}_L \otimes \mathfrak m_2(x_l,z_1))
\otimes {\bf y}; ({\bf e}\otimes z_2);\mathcal T) + \mathfrak n(({\bf x} \otimes z_1) \otimes {\bf y}_L;
({\bf e} \otimes \mathfrak m_2(y_l\otimes z_2));\mathcal T)
.
\end{gather*}
We remark that all the terms of the first formula except the
3rd and 10th ones appear in the second formula.
The formula for $\mathcal H_1$ thus follows up to sign.

Remark again that all the signs are caused by changing the order of operators
or elements (which are graded), that is by Koszul rule. So except
the minus sing which comes from exchanging the order of $\delta$ and
$\mathcal H_i$ (both of which are degree $\pm 1$) the
sign of the corresponding terms coincide.

Therefore, the cancellation occurs with signs.
Thus the formula for $\mathcal H_1$ holds with sign.

The proof of the formula for $\mathcal H_2$ is similar.
\end{proof}

Sublemma~\ref{subem1618} implies that
$\Phi_1 \circ \Phi_2$ is chain homotopic to the identity.
It is easy to see that~${\Phi_1 \circ \Phi_2 = \mathfrak I\circ \mathfrak J}$.
The proof of Lemma~\ref{lem1613} is complete.
\end{proof}

\begin{proof}[Proof of Lemma~\ref{DGcaselem}]
We define
\[
\mathscr I_{1,1} \colon\ \mathscr C_1[1](c_{1,1},c_{1,2}) \otimes \mathscr C_2[1](c_{2,1},c_{2,2})
\to {\bf C}(\mathfrak c_1,\mathfrak c_2)
\]
by \eqref{1618formoo}
and define all the other $\mathscr I_{k,\ell}$ to be zero.
It is easy to see that it defines a DG-functor.
(We use the assumption that $\mathscr C_1$ and $\mathscr C_2$ are
DG-categories here.)
We proved, during the proof of Lemma~\ref{lem1613}, that
$\mathscr I_{1,1}$ is a chain homotopy equivalence.
The lemma now follows from Theorem~\ref{white}.
\end{proof}

\begin{proof}[Proof of Lemma~\ref{G1chainlem}]
Let $\mathfrak c_i = (L_i,b_i)$, $\mathfrak c'_i = (L'_i,b'_i)$.
Note that \[
\mathscr L((\mathfrak c'_1,\mathfrak c'_2),(\mathfrak c_1,\mathfrak c_2))
\cong CF(L'_1,L_1) \otimes CF(L'_2,L_2)
\cong \mathscr C_1(c'_1,c_1) \otimes \mathscr C_2(c'_2,c_2).
\]

On the other hand,
we use the fact that the filtered $A_{\infty}$ category
obtained by Lagrangian Floer theory becomes a DG-category,
after reduction of coefficient to the ground ring, the reduction of the map $\mathscr I_{1,1}$
is given by formula \eqref{1618formoo}.
Here we use the fact that $\mathscr F_{k_1,k_2}$ in \eqref{formlast1614} is $0$
for~${(k_1,k_2) \ne (0,0)}$ and is an isomorphism for $(k_1,k_2) = (0,0)$.

It is easy to see that the reduction of $\mathscr G_1$
is the same map.
Therefore, the reduction of $\mathscr G_1$ to the ground ring is a chain homotopy equivalence.
It implies that $\mathscr G_1$ is a chain homotopy equivalence.
\end{proof}

\section{Orientation and sign}
\label{sec:orient}

In this section, we discuss the orientation and the sign.
The orientation of the moduli spaces of pseudo-holomorphic quilts is
studied by \cite{WW4}.
The orientation of the moduli spaces of pseudo-holomorphic disks
(polygons) and
its relation to $A_{\infty}$ structures is
studied in detail in \cite[Chapter~8]{fooobook}.
In this section, we will prove that orientation and sign appearing in various
moduli spaces and operations in this paper can be reduced to ones
of the moduli spaces of pseudo-holomorphic disks and operations defined by it.

\subsection[Koszul rule in $A_\infty$ structures]{Koszul rule in $\boldsymbol{A_{\infty}}$ structures}
\label{Koszul}

As we mentioned several times, the sign
in various formulas in this paper is by Koszul rule\index{Koszul rule}
(except a few cases which appear in purely algebraic situations,
see the beginning of Section~\ref{sec:compfuncalglem}).
By this reason, we do not write the explicit sign in many of
those formulas. In principle, it is possible (and not so difficult)
to calculate and put the explicit sign to those formulas.
However, actually it is unnecessary to calculate the sign
for the purpose of this paper.
This is because the check of the signs in the equalities
needed in this paper is carried out based on the fact
that the sign is always by Koszul rule
and {\it not} by an explicit calculation of the signs.
Since some of such formulas are complicated,
checking the signs by an explicit calculation could
be cumbersome and lengthy.
Fortunately, we never need it in this paper.

In this subsection, we describe what we mean by
Koszul rule precisely and demonstrate how it works
in certain examples.

We first consider the $A_{\infty}$ formula
\begin{equation}\label{formula2522}
0=\sum_{k_1+k_2=k+1}\sum_{i=0}^{k_1-1}
(-1)^* \mathfrak m_{k_1}(x_1,\dots,x_i,\mathfrak m_{k_2}(x_{i+1},\dots,x_{k_2}),
\dots,x_k).
\end{equation}
Here the sign is given by
\begin{equation}\label{formula2522sign}
* = i +\sum_{j=1}^i \deg x_j.
\end{equation}
We explain how this sign is determined by the Koszul rule.
We order variables and operations appearing in the formula
as follows\footnote{The particular choice
of orders in \eqref{form1733} is not important.
If we take another choice, then the sign in the formulas
changes in exactly the same way for all the terms of \eqref{formula2522}.}
\begin{equation}\label{form1733}
\mathfrak m, \mathfrak m, x_1, x_2, \dots, x_k
\end{equation}
In one of the terms of \eqref{formula2522}, it appears
in the following order
\begin{equation}\label{form1744}
\mathfrak m, x_1, x_2, \dots, \mathfrak m,
 x_{i+1}, \dots, x_k.
\end{equation}
The permutation of operators and variables we need to
go from \eqref{form1733} to \eqref{form1744}
is a composition of permutations of $\mathfrak m$ and $x_j$
for $j=1,\dots,i$.
We remark that the degree of $\mathfrak m$ is $1$ and
the (shifted) degree of $x_j$ is $1+\deg x_j$.
So the sign we pick up by exchanging them
is $1 + \deg x_j$.
Summing them up for $j=1,\dots,i$ we obtain
the sign \eqref{formula2522sign}.\footnote{There are several other sign conventions of
$A_{\infty}$ structure in the literature.
For example, Stasheff's original convention \cite{St-I,St-II}
and Seidel's convention~\cite{Se} are different
from our convention, which is introduced in \cite{fooobook}.
An advantage of our convention lies in the fact that it is
{\it exactly} by Koszul rule.
Therefore, we can automatically determine all the
signs appearing in various formulas by requiring them to be
the Koszul convention. Since there are many operators
which are related to but are slightly different each other
in this paper, putting the sign `by hand' and checking the
consistency by a calculation becomes much cumbersome and lengthy. We can avoid it
by using Koszul convention in all the places.}

In this way, we can obtain the signs appearing in
various formulas systematically.
The author emphasis that this is not only an idea to define a sign but is also a
logical and rigorous {\it definition} of the sign.

To elaborate on this rule, let us describe one more example.
We consider formula \eqref{form925} in Proposition~\ref{basiceqYdiagram}.
In a similar way to \eqref{form1733}, we start with
\begin{equation}\label{175new}
\mathscr{YT}, \mathfrak m, h_{\infty,123},{\bf h}_{12},{\bf h}_{23},{\bf h}_{13},h_{\infty,12},h_{\infty,23}, {\bf h}_{1},{\bf h}_{2},{\bf h}_{3}.
\end{equation}
(Here the symbol $\mathfrak m$ appears. It is identified with $\mathfrak n$
while studying certain terms of \eqref{form925}.)
In the third term, where $(-1)^{*_3}$ appears,
the operations and variables appear in the order
\begin{gather}\label{176new}
\mathscr{YT}, h_{\infty,123},{\bf h}_{12},{\bf h}_{23},{\bf h}_{13}^{c;1}, \mathfrak m, {\bf h}_{13}^{c;2},
{\bf h}_{13}^{c;3},h_{\infty,12},h_{\infty,23},
{\bf h}_{1},{\bf h}_{2},{\bf h}_{3}.
\end{gather}
Here we put
\[
((\Delta\otimes {\rm id}) \circ \Delta)({\bf h}_{13})
= \sum_c {\bf h}_{13}^{c;1} \otimes {\bf h}_{13}^{c;2} \otimes {\bf h}_{13}^{c;3}
\]
and remark that
\begin{equation}\label{177new}
\hat d({\bf h}_{13})
= \sum_c (-1)^{\deg' {\bf h}_{13}^{c;1}}{\bf h}_{13}^{c;1} \otimes \mathfrak m\bigl({\bf h}_{13}^{c;2}\bigr) \otimes {\bf h}_{13}^{c;3}.
\end{equation}
The sign we pick up to go from \eqref{175new} to \eqref{176new}
is $(-1)^*$ with
\[
* = \deg'h_{\infty,123} + \deg'{\bf h}_{12} + \deg'{\bf h}_{23} + \deg'{\bf h}_{13}^{c;1}.
\]
Since $\deg'{\bf h}_{13}^{c;1}$ cancels with the corresponding sign in
\eqref{177new}, we have
\[
*_3 = \deg' h_{\infty,123} + \deg'{\bf h}_{12} + \deg'{\bf h}_{23}.
\]
We next consider the 9th term where $(-1)^{*_9}$ appears.
The order of the operations and variables appearing in this term is
\begin{gather}
\mathscr{YT},h_{\infty,123}, {\bf h}_{12},{\bf h}_{23}^{c_{23};1},{\bf h}_{13},
h_{\infty,12},\mathfrak n,{\bf h}_{2}^{c_2;1},
{\bf h}_{3}^{c_3;1}, {\bf h}_{23}^{c_{23};2}, h_{\infty,23},
{\bf h}_{1},{\bf h}_{2}^{c_2;2},{\bf h}_{3}^{c_3;2}.\label{178new}
\end{gather}
The sign we pick up to go from \eqref{175new} to \eqref{178new}
is
\begin{gather*}
\deg' h_{\infty,123} + \deg'{\bf h}_{12} + \deg' {\bf h}_{23}^{c_{23};1}
+ \deg' {\bf h}_{13} + \deg' h_{\infty,12} \\
\qquad+ \deg' {\bf h}_{23}^{c_{23};2}
\bigl(\deg' {\bf h}_{13} +
\deg' h_{\infty,12} + \deg'{\bf h}_{2}^{c_2;1}+ \deg' {\bf h}_{3}^{c_3;1}\bigr)
\\
\qquad+
\deg' h_{\infty,23} \bigl(\deg'{\bf h}_{2}^{c_2;1}+ \deg' {\bf h}_{3}^{c_3;1}\bigr)
+
\deg' {\bf h}_{1} \bigl(\deg'{\bf h}_{2}^{c_2;1} + \deg'{\bf h}_{3}^{c_3;1}\bigr)
\\
\qquad+
\deg' {\bf h}_{2}^{c_2;2}\deg'{\bf h}_{3}^{c_3;1}.
\end{gather*}
This is by definition $*_9$.
The other $*_k$ is defined in the same way.
Note that there is a minus sign in front of the 10th term.
This minus sign is caused by the fact that
the order of $\mathfrak n$ and~$\mathscr{YT}$ is
exchanged (only) in this term.

The formula we gave for $*_9$ above is rather complicated and actually
it is not useful to write it down explicitly.
On the other hand, it is important that there is a well-defined and canonical way
to determine the signs.

The latter fact is used, for example, in the following way.
During the proof of
Theorem~\ref{thm109} in Section~\ref{sec:compfuncmain},
we claimed that the Y-diagram transformation
is a quatro-module homomorphism.
In other words, the formula which implies that
a pre-quatro-module homomorphism
is a quatro-module homomorphism
coincides with formula \eqref{form925} in Proposition~\ref{basiceqYdiagram}.
It is easy to see that the terms appearing
in those two formulas
are the same except possibly the sign.
We also need to check the signs appearing
in those two formulas coincide.
Since there are many terms
to be checked and since the signs (such as $*_9$ above)
are rather complicated to write down explicitly,
verifying this coincidence by calculating the signs in
those formulas
could be cumbersome and lengthy.
Fortunately, we do not need to
carry out any calculation to check the
coincidence of the signs,
since this fact is an {\it immediate} consequence
of the fact that both signs are by Koszul rule.

The author also remarks that the way we use Koszul
rule here is equivalent to a certain point in the study of $A_{\infty}$,
$L_{\infty}$ structures and their cousins by using
the language of supermanifolds and super-vector-fields
on it. See, for example, \cite{ASKZ97} for such a method.
In those methods, calculation of the
explicit sign is avoided by saying several objects are
functions, vector fields, etc.\ in the sense of supermanifolds.

Another point where we use the fact that
all the signs are by Koszul rule is the proof of the fact
that after adding appropriate correction terms
and putting appropriate orientations,
the operations obtained from moduli spaces
satisfy the basic formulas with {\it Koszul sign}.
The way we will prove it in this section is as follows.
We describe the way how
various moduli spaces such as those
used to define tri-module structures,
Y-diagram transformations,
Double pants transformations, and etc.
can be identified to
moduli spaces of holomorphic
disks (polygons) outside a certain subspace lying in strata of
positive codimension.
Then we use the conclusion of the papers
on the construction of $A_{\infty}$ operations
in Lagrangian Floer theory {\it with sign}
(such as \cite{fooobook2,fooonewbook,ST})
so that there {\it exists} a way to define orientations
of those moduli spaces and correction terms of the signs,
by importing ones of the corresponding moduli space of pseudo-holomorphic disks
via the identification we will give in this section.
Then the $A_{\infty}$ formula of operations in Lagrangian Floer theory
implies the basic formulas on tri-module structures,
Y-diagram transformations,
Double pants transformations, and etc.
{\it with signs}.
This is because they both are by Koszul rule.
We will explain this process more in a concrete situation
in Example~\ref{exm1711}.
We emphasis that this proof does not need to use
the proof of the signs for the~$A_{\infty}$ formula of
Lagrangian Floer theory, in the literature.
It uses only the {\it conclusion} of those papers.
In fact, there could be several different ways to
put orientations and correction terms of the signs
so that the $A_{\infty}$ formula of
Lagrangian Floer theory can be proved.
The argument of this section is independent on such
choices.
For each choice of system of orientations and correction terms
in Lagrangian Floer theory,
we can expand it to the case of
tri-module structures,
Y-diagram transformations, Double pants transformations, and etc.
We also remark that we will {\it not}
provide explicit correction terms
to define tri-module structures,
Y-diagram transformations,
Double pants transformations, and etc.
In principle it is possible to
find it by going back to the
corresponding discussions in the
case of Lagrangian Floer theory
and modify it by Koszul rule. See Section~\ref{orisimpquilt}
and Example~\ref{exm1711}.\footnote{If we want to do so, we would need to see the detail of the
proof of the signs in Lagrangian Floer theory. For
example, it occupies more than 70 pages in \cite{fooobook2}.
So, it seems likely that many of the readers do {\it not}
want to go back and see the proof in the literature and try to understand
how it is adapted to our situation. The way we take in this
section is written in such a way that it is unnecessary for the readers
to do so.}
However, doing so in many places are rather cumbersome
and lengthy process.
Fortunately, we do not need to do so, since we only claim the {\it existence}
of the correction terms of the signs.
Existence of such correction terms is certainly
enough to prove all the results in this paper.

\subsection{Orientation of the moduli space of the simplest quilt}
\label{orisimpquilt}

In this subsection, we consider the case of
the moduli space
\smash{$\overset{\ \text{\tiny $\circ\circ$}}{\mathcal M}_{\rm QT}(\vec a_1,\vec a_{12},\vec a_2;a_-,a_+;E)$}
which is defined in Definition~\ref{def516}.
For simplicity we begin with the case when
$\vec a_1$, $\vec a_{12}$, $\vec a_2$ are empty sets, that is, the
case we do not consider marked points.
(We will discuss the case when there are marked points
later in this subsection.)
We write this moduli space as
$\smash{\overset{\ \text{\tiny $\circ\circ$}}{\mathcal M}_{\rm QT}(L_1,L_{12},L_2;a_-,a_+;}\allowbreak E)$,
where $a_{\pm}$ are connected components of
$\tilde L_1 \times_{X_1} \tilde L_{12} \times_{X_2} \tilde L_2$.
When we are interested in defining orientation only,
it suffices to consider its subset consisting of
a map from a strip $\Sigma = [-1,1] \times \R$.
We write this subset as
${\mathcal M}^{\rm reg}_{\rm QT}(L_1,L_{12},L_2;a_-,a_+;E)$.
It is an equivalence class of maps $((u_1,u_2),(\gamma_1,\gamma_{12},\gamma_2))$,
where
$u_1 \colon [-1,0] \times \R \to X_1$, $u_2\colon [0,1] \times \R \to X_2$ and $\gamma_{i} \colon \R \to L_i$
($i=1,2$), $\gamma_{12} \colon \R \to L_{12}$
and they have the following properties:
\begin{enumerate}\itemsep=0pt\setlength{\leftskip}{0.15cm}
\item[(A.1)]
 $u_1(-1,\tau) = i_{L_1}(\gamma_1(\tau))$,
$u_2(1,\tau) = i_{L_2}(\gamma_2(\tau))$ and
$
(u_1(0,\tau),u_2(0,\tau)) = i_{L_{12}}(\gamma_{12}(\tau))$.
\item[(A.2)]
We require asymptotic boundary condition Condition \ref{cond518}.
\item[(A.3)]
$u_1$, $u_2$ are assumed to be $J_1$, $J_2$ holomorphic, respectively.
\item[(A.4)]
\[
\int_{[-1,0] \times \R} u_1^* \omega_1 + \int_{[0,1] \times \R} u_2^* \omega_2 = E.
\]
\end{enumerate}
We define ${\rm Dub}((u_1,u_2),(\gamma_1,\gamma_{12},\gamma_2))$
as $(u;\gamma_l,\gamma_r)$ such that
\begin{enumerate}\itemsep=0pt\setlength{\leftskip}{0.15cm}
\item[(B.1)]
$u \colon [0,1] \times \R \to -X_1 \times X_2$ is defined by
$
u(\tau,t) = (u_1(-t,\tau),u_2(t,\tau))$.
\item[(B.2)]
$\gamma_r = \gamma_{12} \colon \R \to \tilde L_{12}$.
$\gamma_l \colon \R \to \tilde L_{1} \times \tilde L_2$ is defined by
$
\gamma_l(\tau) = (\gamma_1(\tau),\gamma_2(\tau))$.
\end{enumerate}
By definition, $u$ is $-J_1 \times J_2$ holomorphic.

We consider the disjoint union $L = (L_1 \times L_2) \cup L_{12}$.
Then $(u,\gamma_l,\gamma_+)$ becomes an element of
\smash{$
\mathring{\mathcal M}(L,(a_-,a_+);E)
$}
which is defined in Definition~\ref{defn314}.
Here we write it as
$\smash{
\mathring{\mathcal M}}(L_{12},L_1 \times L_2;(a_-,a_+);E)
$.
This is the moduli space used in \cite[Section 3.7.4]{fooobook} to define
the boundary operator on $CF(L_{12},L_1 \times L_2)$.

We thus obtain an open embedding\index[syindex]{Dob@${\rm Dob}$}
\begin{equation}\label{dobe17111}
{\rm Dob} \colon\ {\mathcal M}^{\rm reg}_{\rm QT}(L_1,L_{12},L_2;a_-,a_+;E) \to \mathring{\mathcal M}(L_{12},L_1 \times L_2;(a_-,a_+);E).
\end{equation}

We assumed that $L_{12}$ is $\pi_1^*(V_1 \oplus TX_1) \oplus \pi_2^*V_2$
relatively spin.
We also assumed $L_1$ is $V_1$ relatively spin and so is $V_1 \oplus TX_1$
relatively spin. In fact, since~${TX_1\vert_{L_1} = TL_1 \otimes_{\R} \C}$ and
$TL_1$ is oriented. So $TX_1\vert_{L_1}$ has canonical
spin structure.
We assumed~$L_2$ is $V_2$ relatively spin.
Therefore, $L_1 \times L_2$ is also $\pi_1^*(V_1 \oplus TX_1) \oplus \pi_2^*V_2$
relatively spin.

Thus by Proposition~\ref{prop329}, we have an isomorphism of
principal ${\rm O}(1)$ bundle
\begin{equation}\label{oriaso172}
O_{{\mathcal M}^{\rm reg}(a_-,a_+;E)} \cong
O_{\mathring{\mathcal M}(L,(a_-,a_+);E)}
\cong {\rm ev}_-^* \Theta_{a_-}^- \otimes {\rm ev}_+^* \Theta_{a_+}^+.
\end{equation}
We can use the isomorphism \eqref{oriaso172} to define
$
{\rm ev}_+! \circ {\rm ev}_-^* \colon
\Omega(R_{a_-};\Theta_{a_-}^-)
\to
\Omega(R_{a_-};\Theta_{a_+}^-)
$
by smooth correspondence.
This is \eqref{form526} in case we do not have boundary marked points.

We next show the consistency of orientations at the boundary.
\begin{rem}\label{rem171}
Before doing so, we explain what we mean by
`consistency of orientations at the boundary'
more precisely.
We consider the `open inclusion'
\[
(-1)^{*_1}\mathcal M_{k_1+1}(L;\beta_1) {}_{{\rm ev}_i}\times_{{\rm ev}_0}
\mathcal M_{k_2+1}(L;\beta_2)
\subseteq
\partial \mathcal M_{k+1}(L;\beta),
\]
where $k_1+k_2 = k$, $\beta_1+\beta_2=\beta$.
Here $*_1$ is a certain correction term of the sign.\footnote
{$\mathcal M_{k+1}(L;\beta)$ is the compactified
moduli space of pseudo-holomorphic disks with $k+1$ boundary marked
points and of homology class $\beta \in \pi_2(X,L)$.}
This is an example of consistency of orientations at the
boundary.
Namely, the orientations of the moduli spaces
appearing in the left and right-hand sides of the
formula coincide.
To give a rigorous meaning to its
coincidence, we
also need to fix a convention of the orientation
of the fiber product (as well as the boundary).

The `consistency of orientations at the boundary'
are supposed to imply the fundamental equation
(in this case the $A_{\infty}$ relation)
{\it with sign}, which is the Koszul sign in this paper.
Let us elaborate on this point.
Let $C(L)$ be a certain chain model of the cohomology of the space~${\tilde L\times_X \tilde L}$.
(In this paper we take de Rham model.)
The moduli space $\mathcal M_{k+1}(L;\beta)$
regarded as a correspondence from $L^k$ to $L$
gives an operation,
which is
$
\mathfrak m_{k,\beta}\colon C(L)^{\otimes k} \to C(L)$.
In the case of de Rham model, it is
$
(h_1,\dots,h_k) \mapsto
(-1)^{*_2} {\rm ev}_{0 !} ({\rm ev}^*_{1} h_1\wedge \dots\wedge{\rm ev}^*_kh_k)$.
(More precisely, we need CF-perturbations.)
To make sense of this formula, we need to fix a~convention of sign for integration along fiber ${\rm ev}_{0 !}$.
(Provably, the sign convention for the pullback ${\rm ev}_{i}^*$
is mostly obvious.)
Here $*_2$ is a certain correction term of the sign.
Thus we have to fix all the conventions mentioned above
together with correction terms $*_1$, $*_2$
so that the operator~$\mathfrak m_{k,\beta}$
satisfies the~$A_{\infty}$ relation with Koszul sign.

In the case of this $A_{\infty}$ relation,
this is worked out in singular chain complex model
in \cite{fooobook2} and in de Rham model in \cite{fooonewbook}
and \cite{ST}, in the case when our Lagrangian submanifold
is embedded. In the case of an immersed Lagrangian
submanifold which has transversal self-intersection,
it is worked out in \cite{AJ} in singular
chain complex model.
In the case of an immersed Lagrangian
submanifold which has self-clean intersection, it is written in Section~\ref{oriAinfMB}
in singular
chain complex model and in the paper \cite{ono2} by Kaoru Ono
in de Rham model.
We use the conclusion of those results (but not the proof of them).

The sign convention of \cite{fooobook2} and of \cite{fooonewbook}
are different\footnote{In fact, the sign of $\mathfrak m_1$ is different.}
but they both satisfy the same $A_{\infty}$ formula with the same sign.
(Note that if we regard smooth singular chains as currents and
approximate them by smooth differential forms,
then we can ask whether various conventions
(the convention of the sign of pushout (=~integration along the fiber)
or pullback), together with correction terms $*_2$\footnote{Which is not the same in two books.}
gives the same operator $\mathfrak m_{k,\beta}$ with sign or not.)\footnote{The correction term $*_1$ coincide in those two books.}
The author did not check whether the convention of~\cite{fooonewbook}
coincides with \cite{ST} or not.
The convention of \cite{AJ} is slightly different from \cite{fooobook2}
at the point which we mention in Proposition~\ref{form3472}.

The sign part of the works \cite{AJ,fooobook2,fooonewbook,ST},
is computational.\footnote{There are geometric
ideas behind those computations in many cases. However,
such ideas are not used during the proof.} The conventions and correction terms are
defined by `hand' and the
sign part of the $A_{\infty}$ formula is checked
by computation.
It might be possible to give more conceptional proof.
So far no such proof is written in the literature.
Since the check of sign in many cases are complicated and
pains taking such a proof would be desired.
However, it is not a theme of this paper.

The discussion of this section,
which reduces the
sign issue of this paper
to one of $A_{\infty}$ relation among $\mathfrak m$,
is not computational.

\end{rem}
\begin{rem}\label{rem172}
In several places, we write explicit correction terms
(written as $(-1)^{*_1}$, $(-1)^{*_2}$ in Remark~\ref{rem171}),
following the convention of \cite{fooonewbook,ST}.
However, actually we never used these particular
choices or the choices of other conventions.
We use only the fact that there {\it exist} such
choices which induce $A_{\infty}$ formula with Koszul sign.
\end{rem}
We go back to the discussion of consistency of orientations at the boundary.
We first~ob\-ser\-ve that the boundary of the
compactification
${\mathcal M}_{\rm QT}(L_1,L_{12},L_2;a_-,a_+;E)$
of the moduli space ${{\mathcal M}^{\rm reg}_{\rm QT}(L_1,L_{12},L_2;a_-,a_+;E)}$
consists of four kinds of components depicted in
Figures~\ref{Figure58}--\ref{Figure511}.

On the other hand, the codimension one boundary
component of the compactification
${\mathcal M}(L_{12},\allowbreak L_1 \times L_2;(a_-,a_+);E)$
of
\smash{$\mathring{\mathcal M}(L_{12},L_1 \times L_2;(a_-,a_+);E)$}
is described by one of the configurations (1), (2), (3) depicted in
Figures \ref{FigureSec17-1} below.

\begin{figure}[ht]
\centering
\includegraphics[scale=0.38]{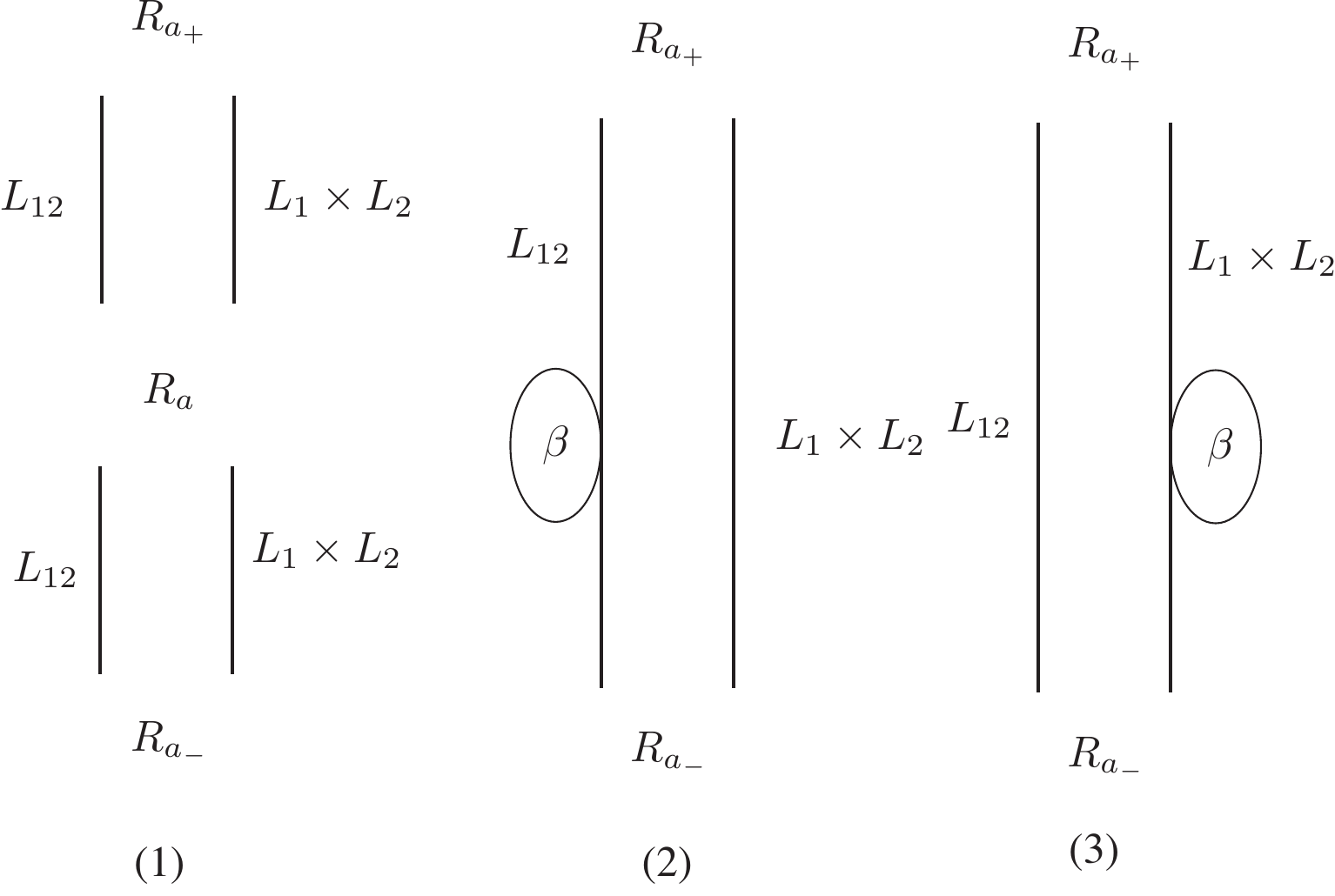}
\caption{Boundary components of ${\mathcal M}(L_1,L_{12},L_2;a_-,a_+;E)$.}
\label{FigureSec17-1}
\end{figure}

We observe that Figure~\ref{FigureSec17-1}\,(1) and (2)
correspond to Figures \ref{Figure58} and \ref{Figure510},
respectively.
Since the orientation is defined so that
\eqref{dobe17111} lifts to the isomorphism
of orientation bundles (principal~${\rm O}(1)$ bundles)
the compatibility of the orientation of
the moduli space, the compactification of
${\mathcal M}^{\rm reg}_{\rm QT}(L_1,L_{12},L_2;a_-,a_+;E)$,
at the boundary described by Figures~\ref{Figure58} and \ref{Figure510}
follows from the corresponding compatibility of
${\mathcal M}(L_{12},L_1\times L_2;(a_-,a_+);E)$
at the boundary described by
Figure~\ref{FigureSec17-1}\,(1) and (2).
The latter is established in
\cite[Chapter 8, Theorem 8.8.10 etc.]{fooobook}.

We finally consider
the boundary described by
Figure~\ref{FigureSec17-1}\,(3).
The homology class of the bubbled disk
is $\beta \in \pi_2(X_1\times X_2;L_1 \times L_2;\Z)
= \pi_2(X_1;L_1;\Z) \times \pi_2(X_2;L_2;\Z)$.
We write it~$(\beta_1,\beta_2)$.
We consider the following three cases separately.

Case 1: $\beta_1\ne 0 \ne\beta_2$.
The configuration which corresponds to an element of the space ${\mathcal M}^{\rm reg}(L_1, L_{12},L_2;a_-,a_+;E)$
of this case is depicted
in Figure~\ref{Figure17-2} below.
We remark that this component has codimension
greater than $1$. Therefore, we do not need to study
this case to show the consistency of the orientation
at the boundary.
\begin{figure}[ht]
\centering
\includegraphics[scale=0.4]{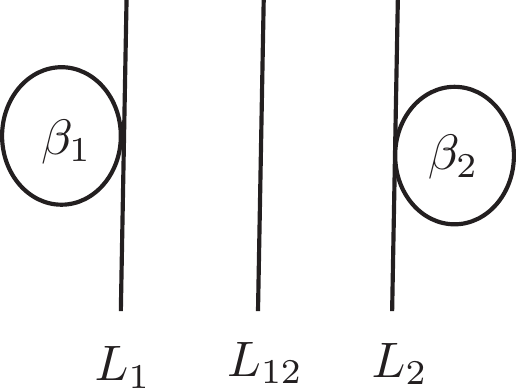}
\caption{Case 1.}
\label{Figure17-2}
\end{figure}

Case 2: $\beta_1 = 0 \ne\beta_2$.
This case corresponds to Figure~\ref{Figure511}.
Therefore, the consistency of orientation at this
boundary component follows from
the corresponding discussion for Figure~\ref{Figure511}, which is in
\cite[Chapter 8, Theorem 8.8.10 etc.]{fooobook}.

Case 3: $\beta_1 \ne 0 = \beta_2$.
This case corresponds to Figure~\ref{Figure59}.
Therefore, the consistency of orientation at this
boundary component follows from \cite[Chapter~8]{fooobook},
except the following points.

For an element $((u_1,u_2),(\gamma_1,\gamma_{12},\gamma_2))$ of
${\mathcal M}^{\rm reg}(L_1,L_{12},L_2;a_-,a_+;E)$,
we require $u_1$ to be $-J_1$ holomorphic.
Since in Case 3 bubble occurs at the line $t = -1$, the
map in the bubble is $-J_1$ holomorphic.
We also consider the map $(t,\tau) \mapsto u_1(-t,\tau)$.
Note we use $V_1 \oplus TX_1$ relative spin structure
of $L_1$ for our orientation.
As was shown as Theorem~\ref{opthere}, the orientation of the
moduli space of $-J_1$ holomorphic
map $(t,\tau) \mapsto u_1(-t,\tau)$ using $V_1 \oplus TX_1$ relative spin
structure, coincides with one of $u_1$ using
$J_1$ holomorphic
map moduli space using $V_1$ relative spin structure,
after reversing the enumeration of the boundary marked
points. This is consistent with the fact that
we study opposite category for $L_1$. Moreover, we use
$V_1 \oplus TX_1$ relative spin structure
in place of $V_1$ relative spin structure
for orientation.
As is shown in Section~\ref{subsec:Opposite},
this is equivalent to use $-J_{X_1}$ instead of
$J_{X_1}$.

Note that the above discussion proves Theorem~\ref{therem530}\,(7).

We next include the case when there are marked points
and explain the way to fix the sign of the operations
$\mathfrak{n}_{k_1,k_{12},k_2}^{E,\varepsilon}$ in \eqref{form526}.
We use the notations in \eqref{form526}. We put $m = k_1+k_2$
We denote by ${\rm Shuf}(k_1,k_2)$ the set of pairs of maps
$(I_1,I_2)$, where $I_j \colon \{1,\dots,k_j\} \to \{1,\dots,m\}$
such that
\begin{enumerate}\itemsep=0pt
\item[(1)] The image of $I_1$ and $I_2$ are disjoint.
\item[(2)] $I_1$ reverses the order.
\item[(3)] $I_2$ preserves the order.
\end{enumerate}
For $I=(I_1,I_2) \in {\rm Shuf}(k_1,k_2)$, we write
\[
\mathfrak x^I_j =
\begin{cases}
x_i &\text{if $I_1(i) = j$}, \\
z_i &\text{if $I_2(i) = j$}.
\end{cases}
\]
Let
\[
\mathfrak n^E_{m,k_{12}} \colon\ B_{m}CF(L_1\times L_2) \otimes CF(L_1\times L_2;L_{12})
\otimes B_{k_{12}}CF(L_{12})
\to CF(L_1\times L_2;L_{12})
\]
be the filtered $A_{\infty}$ bimodule structure
for the pair of Lagrangian submanifolds $L_1 \times L_2$,
$L_{12}$ of $X_1 \times X_2$. \big(More precisely, its coefficient of $T^E$.\big)
This is defined in \cite[Definition 3.7.41]{fooobook}
in the singular homology version.
The de Rham version is a part of the structure operation
of the filtered $A_{\infty}$ category associated to
the symplectic manifold $X_1 \times X_2$, which we
defined in Theorem~\ref{prop333}.
See also \cite{fooonewbook,ST}.

The discussion in the case without marked point,
implies that in our case the moduli space
we use to define
the filtered tri-module structure \smash{$\mathfrak{n}_{k_1,k_{12},k_2}^{E,\varepsilon}$}
coincides with the closure of union of
the moduli spaces defining $\mathfrak n^E_{m,k_{12}}$ for
various $I=(I_1,I_2) \in {\rm Shuf}(k_1,k_2)$,
outside codimension 1 set.
Namely, we have
\begin{equation}\label{form171422}
\mathfrak{n}_{k_1,k_{12},k_2}^{E,\varepsilon}({\bf x},{\bf y},w,{\bf z})
=
\sum_{I\in {\rm Shuf}(k_1,k_2)} (-1)^{*_I}
\mathfrak n_{m,k_{12}}\bigl({\bf y},w,{\mathfrak x}^I\bigr).
\end{equation}
Here the sign $(-1)^{*_I}$ is the Koszul sign,
which is determined as follows.
We remark that $\mathfrak x^I$, $w$,~${\bf y}$ coincide with
${\bf x}$, ${\bf y}$, ${\bf z}$, $w$ up to exchanging the order.
So we shift the degree of them by one and put the sign
which arises when we exchange the order of those
variables via Koszul rule.

\begin{figure}[ht]
\centering
\includegraphics[scale=0.3]{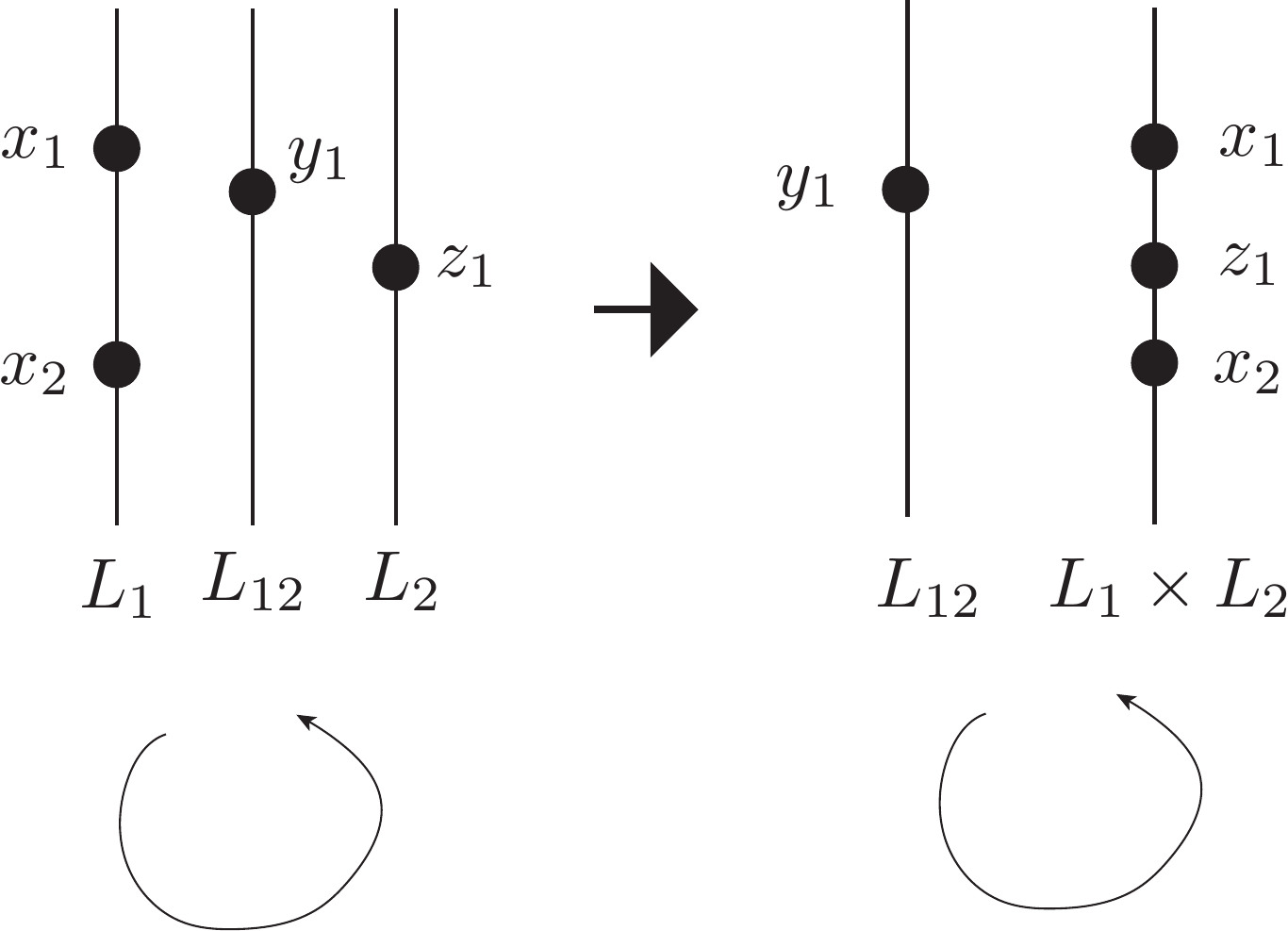}
\caption{Enumeration of marked points assigned to $L_1\times L_2$.}
\label{Figure17-3}
\end{figure}

For example, if $k_1 = 2$, $k_2 = 1$ and $k_{12} = 1$
and $\operatorname{Im}(I_{1}) = \{1,3\}$, then
the corresponding term is
$
(-1)^{*}\mathfrak n_{1,3}(y_1,w,x_2,z_1,x_1)
$,
where
\[
* = \deg'x_1(\deg' y_1+ \deg w' + \deg' x_2+ \deg' z_1)
+ \deg' x_2(\deg' y_1+ \deg w')
+ \deg'w \deg'y_1
\]
is the sign which we get to exchange
$
y_1,w,x_2,z_1,x_1
\mapsto
x_1,x_2,y_1,w,z_1$.
See Figure~\ref{Figure17-3} and Example~\ref{exm1711}.

The bimodule structure $\mathfrak n$ satisfies
the relation
\begin{align}
\sum_{a_1,a_2}(-1)^*{\mathfrak n}\bigl({\bf y}_{a_1}^{(1)},
\mathfrak n\bigl({\bf y}_{a_1}^{(2)},w,{\mathfrak x}_{a_2}^{(1)}\bigr),{\mathfrak x}_{a_2}^{(2)}\bigr)
&+
\sum_a (-1)^* \mathfrak n\bigl({\bf y}_{a}^{(1)}\otimes\mathfrak m\bigl({\bf y}_{a}^{(2)}\bigr)
\otimes{\bf y}_{a}^{(3)},w, {\mathfrak x}\bigr) \nonumber\\
&+ \sum_a (-1)^* \mathfrak n\bigl({\bf y},w, {\mathfrak x}_{a}^{(1)}
\otimes \mathfrak m\bigl({\mathfrak x}_{a}^{(2)}\bigr)\otimes {\mathfrak x}_{a}^{(3)}\bigr).\label{form175175}
\end{align}
Here $\Delta({\bf y}) = \sum_{a_1} {\bf y}_{a_1}^{(1)} \otimes {\bf y}_{a_1}^{(2)}$,
$(\Delta\otimes 1)(\Delta({\bf y})) = \sum_{a} {\bf y}_{a}^{(1)} \otimes {\bf y}_{a}^{(2)}
\otimes {\bf y}_{a}^{(3)}$ etc.

The sign $*$ in \eqref{form175175} is the Koszul sign.
The relation \eqref{form175175} is \cite[Theorem 3.7.72]{fooobook}
in singular homology version and is a part of Theorem~\ref{prop333}
in de Rham version. The Koszul sign rule is a consequence of
\cite[Chapter 8, Theorem 8.8.10]{fooobook}.

Since the sign is always by Koszul rule in this paper, the
(tri-module analogue of) the formula~\eqref{form9300}, where the sign is also by Koszul rule, is
a consequence of \eqref{form175175},
once we take the next two points (dif.1) (dif.2) into account.

 (dif.1)
When we put ${\mathfrak x} = {\mathfrak x}^I$ in formula \eqref{form175175},
the third sum contains a term where $\mathfrak m$ is applied to both
$x_i$'s and $z_i$'s.
For example, if $k_1 = 2$, $k_2 = 1$ and $k_{12} = 1$
and $\operatorname{Im}(I_{1}) = \{1,3\}$, a~term such as
\begin{equation}\label{form1766}
\pm \mathfrak n_{1,2}(y_1,w,x_2,\mathfrak m_2(z_1,x_1))
\end{equation}
appears.
There is no corresponding term in \eqref{form9300}.
The reason is as follows. The compactification we take for
the moduli space ${\mathcal M}_{\rm QT}(L_1,L_{12},L_2;\vec a_1,\vec a_{12},\vec a_2;a_-,a_+;E)$
which we used to define~\eqref{form9300} is different
from the compactification of
\smash{$\mathring{\mathcal M}(L_{12},L_1 \times L_2;(a_-,a_+);E)$}
which we use to define
$\mathfrak n_{m,k_{12}}$.
More specifically, the configuration such as
Figure~\ref{Figure17-2} (with marked points included)
do {\it not} appear in the
codimension one boundary of ${\mathcal M}_{\rm QT}(L_1,L_{12},L_2;\vec a_1,\vec a_{12},\vec a_2;a_-,\allowbreak a_+; E)$.
As we observed before, this is a special case of
the codimension one boundary described in Figure~\ref{FigureSec17-1}\,(2),
which gives a term such as \eqref{form1766}.

Note that this fact does not affect the discussion of sign of the
other components which both appear in \eqref{form9300} and \eqref{form175175}.

 {(def.2)}
We remark the following three points:
\begin{itemize}\itemsep=0pt
\item[$\bullet$] For $(I_1,I_2) \in {\rm Shuf}(k_1,k_2)$, we require $I_1$ to be
order {\it reversing}.

\item[$\bullet$]
When we consider the bubble which occurs at $L_1$ for the compactified
moduli space ${\mathcal M}_{\rm QT}(L_1,L_{12},L_2; \vec a_1,\vec a_{12},\vec a_2;a_-,a_+;E)$,
the map on this bubble is regarded as a $J_1$ holomorphic map.
On the other hand, when we consider the corresponding object
as a bubble in an element of the compactification of
\smash{$\mathring{\mathcal M}(L_{12},L_1 \times L_2;(a_-,a_+);E)$}, the map on the bubble is regarded as a $-J_1$ holomorphic
map by using appropriate anti-holomorphic map~${D^2 \to D^2}$.

\item[$\bullet$]
We regard $L_1$ as $V_1 \oplus TX_1$ relatively spin
when we consider the moduli space ${\mathcal M}_{\rm QT}(L_1,L_{12},\allowbreak L_2;\vec a_1,\vec a_{12},\vec a_2;a_-,a_+;E)$
but as $V_1$ relatively spin when we consider
$\smash{\mathring{\mathcal M}}(L_{12},L_1 \times L_2;\allowbreak(a_-,a_+);E)$.
\end{itemize}
By Theorem~\ref{opthere}, these three points cancel out and
we obtain the correct sign.

\begin{exm}\label{exm1711}
Let us elaborate on this fact more by an explicit example.
Let us consider the moduli space depicted in Figure~\ref{Figure17-3}.
We first study the boundary component depicted in Figure~\ref{Figure17-3-2} below.

\begin{figure}[ht]
\centering
\includegraphics[scale=0.6]{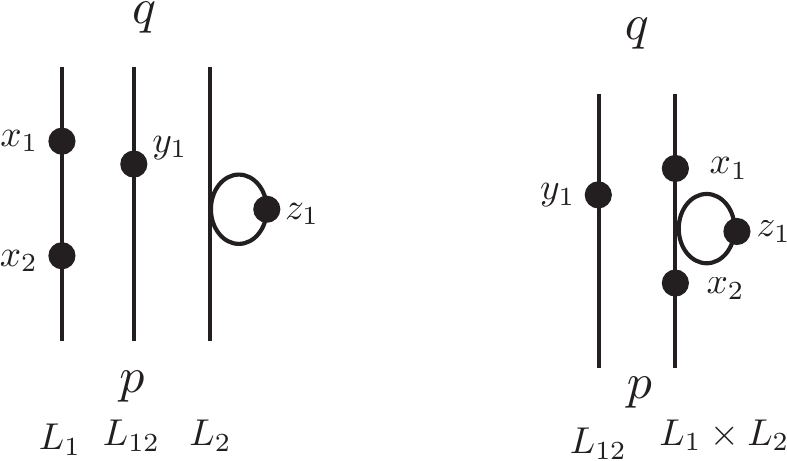}
\caption{A boundary component of Figure~\ref{Figure17-3}.}
\label{Figure17-3-2}
\end{figure}

The left figure corresponds to
\begin{equation}\label{form17440}
\langle \mathfrak n_{2,1;1}(x_1,x_2;y_1;p;\mathfrak m_1(z_1)),q \rangle.
\end{equation}
Here $p$ and $q$ are chains
of $(L_1 \times L_2)\times_{X_1\times X_2} L_{12}$ appearing at $\tau \to \pm \infty$.
The right figure corresponds to
\begin{equation}\label{form1755}
\langle \mathfrak n_{1,3}(y_1;p;x_2,\mathfrak m_1(z_1),x_1),q\rangle.
\end{equation}
Note that in the $A_{\infty}$ relation of tri-module
\eqref{form17440} appears with sign $(-1)^{*_1}$, where
$
*_1 = \deg' x_1 + \deg' x_2 + \deg' y_1 + \deg' p$.
In the $A_{\infty}$ relation of Lagrangian Floer theory (the
$A_{\infty}$ bi-module structure on $CF(L_{12},L_1\times L_2)$),
\eqref{form1755} appears with sign $(-1)^{*_2}$, where
$
*_2 = \deg' y_1 + \deg' p + \deg' x_2$.
We next study the boundary component depicted in Figure~\ref{Figure17-3-3} below.
\begin{figure}[ht]
\centering
\includegraphics[scale=0.6]{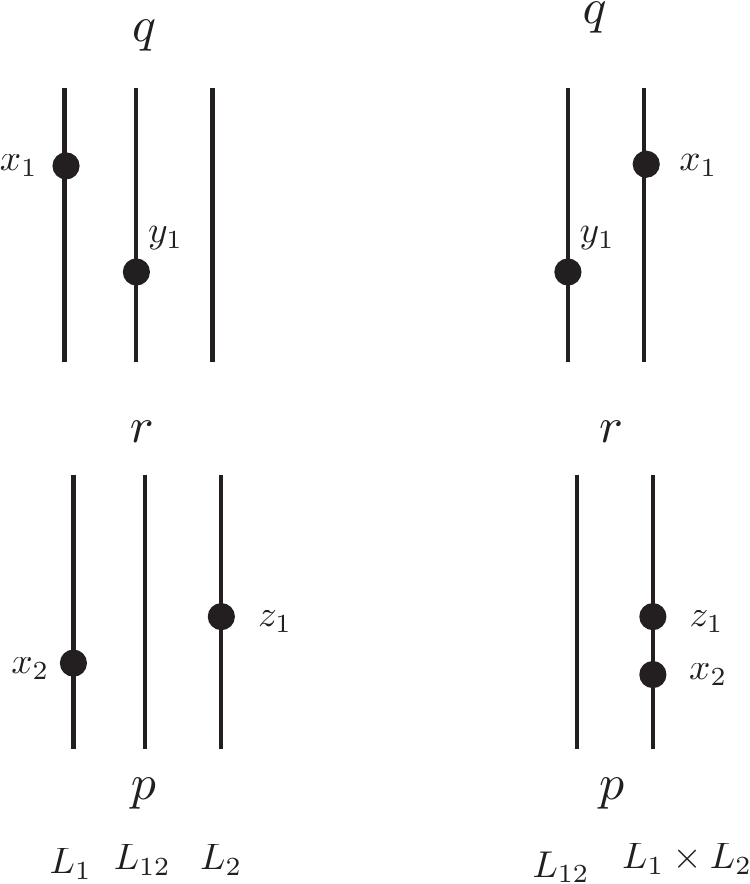}
\caption{Another boundary component of Figure~\ref{Figure17-3}.}
\label{Figure17-3-3}
\end{figure}

The left figure corresponds to
\begin{equation}\label{form1746}
\langle \mathfrak n_{1,1;0}(x_1,y_1;\mathfrak n_{1,0;1}(x_2;p;z_1)),q \rangle
\end{equation}
and the right figure corresponds to
\begin{equation}\label{form1747}
\langle \mathfrak n_{1,1}(y_1,\mathfrak n_{0,2}(p,x_2,z_1),x_1),q\rangle.
\end{equation}
\eqref{form1746} comes with sign $(-1)^{*_3}$ where
$
*_3 = \deg' x_1 + \deg' y_1 + \deg' x_2 \deg' y_1$.
\eqref{form1747} comes with sign $(-1)^{*_4}$ where
$
*_4 = \deg' y_1 $.
Thus
the $A_{\infty}$ formula of the tri-module structure $\mathfrak n$
which we want to prove
is of the form
\begin{align}
0 ={} &\dots
+ (-1)^{*_1}\mathfrak n_{2,1;1}(x_1,x_2;y_1;p;\mathfrak m_1(z_1)) + \cdots\nonumber\\
&
 + (-1)^{*_3}\mathfrak n_{1,1;0}(x_1,y_1;\mathfrak n(x_2;p;z_1))
 + \cdots,\label{form1719}
\end{align}
and the $A_{\infty}$ formula for operations $\mathfrak n$, $\mathfrak m$
which was proved in the literature is
\begin{align}
0 ={} &\dots
+ (-1)^{*_2}\mathfrak n_{1,3}(y_1,p,x_2,\mathfrak m_1(z_1),x_1)
+ \cdots\nonumber\\
& + (-1)^{*_4}\mathfrak n_{1,1}(y_1,\mathfrak n_{0,2}(p,x_2,z_1),x_1) + \cdots.\label{form1720}
\end{align}
We claim that \eqref{form1719} is a consequence of
\eqref{form1720}. To see this, we calculate
the difference of signs between $\mathfrak n$'s and $\mathfrak m$'s.

We note that for $v = \mathfrak m_1(z_1)$ we have
$
\mathfrak n_{2,1;1}(x_1,x_2;y_1;p;v)
= (-1)^{*_5}\mathfrak n_{1,3}(y_1,p,x_2,v,x_1)
$,
where $*_5$ is the Koszul sign induced by the change of the
order of variables
$
x_1,x_2,y_1,p,v
\longrightarrow
y_1,p,x_2,v,x_1$.
Therefore,
\begin{align*}
*_5 ={} &\deg' x_1 (\deg' x_2 + \deg'y_1 + \deg' p + \deg' v)
+ \deg' x_2 (\deg'y_1 + \deg' p) \\
={} &\deg' x_1 (\deg' x_2 + \deg'y_1 + \deg' p + \deg' z_1 +1)
+ \deg' x_2 (\deg'y_1 + \deg' p).
\end{align*}
We have also
$
\mathfrak n_{1,1}(x_2;p;z_1)
= (-1)^{*_6}\mathfrak n_{0,2}(p,x_2,z_1)
$,
where $*_6$ is the Koszul sign induced by the change of the
order of variables
$
x_2,p,z_1
\longrightarrow
p,x_2,z_1$.
Therefore,
$
*_6 = \deg' x_2 \deg' p$.
For $w = \pm \mathfrak n_{0,2}(p,x_2,z_1)$, we have
$
\mathfrak n_{1,1;0}(x_1,y_1;w)
= (-1)^{*_7}\mathfrak n_{1,1}(y_1,w,x_1)
$,
where $*_7$ is the Koszul sign induced by the change of the
order of variables
$
x_1,y_1,w
\longrightarrow
y_1,w,x_1$.
Therefore,
\[
*_7 = \deg' x_1 (\deg' y_1 + \deg'w)
= \deg' x_1 (\deg' y_1 + \deg' p + \deg' x_2 + \deg'z_1 + 1).
\]
The claim that \eqref{form1719} is a consequence of
\eqref{form1720}
follows from the congruence
\begin{equation}\label{1721}
*_1 + *_2 + *_3 + *_4 + *_5 + *_6 + *_7 \equiv 0 \mod 2.
\end{equation}
One can check \eqref{1721} by calculating the formula
of $*_i$ given explicitly above.
However, actually we do not need any calculation to
show \eqref{1721},
since \eqref{1721} is an immediate consequence of the fact
that the map from permutation group to $\{\pm 1\}$ which
associates the Koszul sign to each permutation is a
group homomorphism.

In fact, both $*_1+*_2+*_5$ and
$*_3+*_4+*_6$ are the Koszul sign
associated to the permutation~${
\mathfrak n, \mathfrak n, x_1, x_2, y_1, p, z_1
\to
\mathfrak n, \mathfrak n, y_1, p, x_2, z_1, x_1}$.

By this reason, the discussion of this example can be easily
generalized to other cases, as far as the sign of the
formulas we want to prove is by Koszul rule and
we are given an identification of the moduli spaces we use to
moduli spaces of pseudo-holomorphic disks (polygons).
\end{exm}

\subsection{Orientation of the moduli space of pseudo-holomorphic drums}
\label{oridrum}

\begin{figure}[ht]
\centering
\includegraphics[scale=0.3]{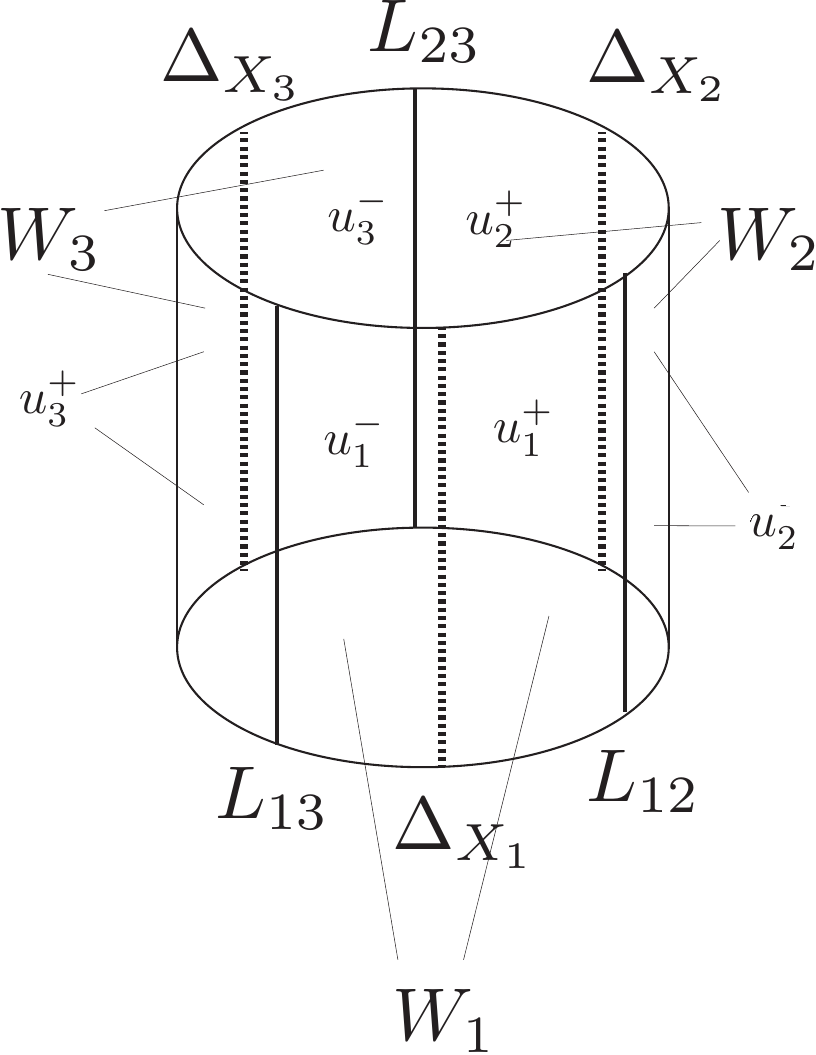}
\caption{Adding diagonal to a drum.}
\label{Drumcut}
\end{figure}

\begin{figure}[ht]
\centering
\includegraphics[scale=0.3]{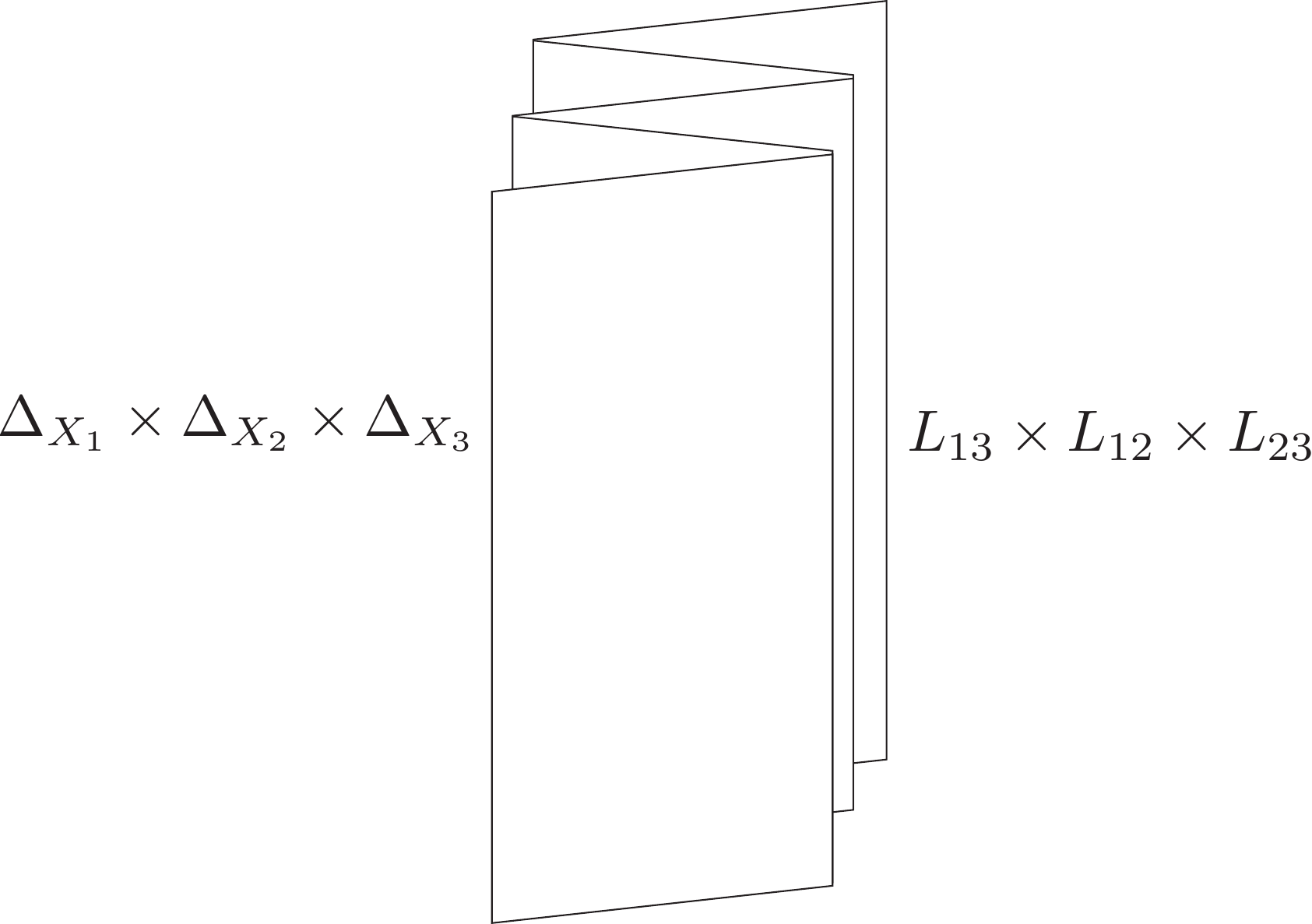}
\caption{Regard a drum as a strip.}
\label{reglueddrum}
\end{figure}

In this subsection, we study the orientation of the moduli space
of pseudo-holomorphic drums, Definition~\ref{def916}.
The quilted domain $W$ there is divided into
three pieces $W_1$, $W_2$ and $W_3$ and
$u_i\colon W_i \to X_i$ is $-J_{X_i}$ holomorphic.
We identify~${
W_i = [-1,1] \times \R}
$
and put
$
W_i^- = [-1,0] \times \R$, $
W_i^+ = [0,1] \times \R$.
Let $u_i^+$, $u_i^-$ be the restriction of $u_i$
to $W_i^+$, $W_i^-$. We define
$
\hat u_i = (\hat u_i^-,\hat u_i^+) \colon [0,1] \times \R \to X_i^2
$
by~${
\hat u_i^-(t,\tau) = u_i^-(t,\tau)}$,
$
\hat u_i^+(t,\tau) = u_i^+(-t,\tau)
$.
Then
\[
(\hat u_1,\hat u_2,\hat u_3)
\colon\ [-1,0] \times \R \to (-X_1 \times X_1)
\times (-X_2 \times X_2)\times (-X_3 \times X_3)
\]
is pseudo-holomorphic.
See Figures \ref{Drumcut} and \ref{reglueddrum}.
In Figure~\ref{Drumcut}, we add 3 extra seams that are
depicted by dotted lines in Figure~\ref{Drumcut}.
The boundary condition becomes
the product of diagonals $\prod_{i=1}^3 \Delta_{X_i}$
at the boundary $\{0\} \times \R$
and is $L_{13} \times L_{12} \times L_{23}$
at the boundary $\{1\} \times \R$.

Let $L$ be the disjoint union of
$\prod_{i=1}^3 \Delta_{X_i}$
and $L_{13} \times L_{12} \times L_{23}$.
It is an immersed Lagrangian submanifold
of $\prod_{i=1}^3 (-X_i \times X_i)$.

We decompose
\[
\prod_{i=1}^3 \Delta_{X_i}
\times_{\prod_{i=1}^3 (-X_i \times X_i)}
\bigl(\tilde L_{13} \times \tilde L_{12} \times \tilde L_{23}\bigr)
\]
into components $R_{123}(a)$, $a \in \mathcal A$.

Thus the above construction defines a map
\begin{gather}
{\rm Dob} \colon\ {\mathcal M}^{\rm reg}_{\rm DR}(L_{13},L_{12},L_{23};a_-,a_+;E) \nonumber\\
\hphantom{{\rm Dob} \colon} \ \to
\overset{\ \text{\tiny $\circ\circ$}}{\mathcal M}(\Delta_{X_1} \times \Delta_{X_2} \times \Delta_{X_3},L_{13} \times L_{12} \times L_{23};a_-,a_+;E).\label{dobe17111rev}
\end{gather}
Here the moduli space
\smash{$\overset{\ \text{\tiny $\circ\circ$}}{\mathcal M}_{\rm DR}(L_{13},L_{12},L_{23};a_-,a_+;E)$}
is a special case of the moduli space
\smash{$\overset{\ \text{\tiny $\circ\circ$}}{\mathcal M}_{\rm DR}(\vec a_{13},\vec a_{12},\vec a_{23};a_-,a_+;E)$}
in Definition~\ref{def916}, where $\vec a_{13}$, $\vec a_{12}$, $\vec a_{23}$
are empty sets.\footnote{In other words, we do not put marked points on
the seams.}
$\smash{\overset{\ \text{\tiny $\circ\circ$}}{\mathcal M}}(\Delta_{X_1} \times\allowbreak \Delta_{X_2} \times \Delta_{X_3},
L_{12} \times L_{23} \times L_{13};a_-,a_+;E)$
is the part of
\smash{$\overset{\ \text{\tiny $\circ\circ$}}{\mathcal M}(L,a_-,a_+;E)$}
which is used to define the boundary operator
\begin{gather*}
\mathfrak n_{0,0} \colon\
CF(\Delta_{X_1} \times \Delta_{X_2} \times \Delta_{X_3},
L_{12} \times L_{23} \times L_{13})\\
\hphantom{\mathfrak n_{0,0} \colon} \
\to CF(\Delta_{X_1} \times \Delta_{X_2} \times \Delta_{X_3},
L_{12} \times L_{23} \times L_{13}).
\end{gather*}

\eqref{dobe17111rev} is an isomorphism of Kuranishi structure
and so we can use orientation of
the right-hand side to define orientation of the left-hand side.
This implies Proposition~\ref{prop811}\,(3).
We remark that once Proposition~\ref{prop811}\,(3) is proved
then the choice of $\sigma_{13}$, the relative spin
structure of the geometric composition
$L_{13} = L_{12} \times_{X_2} L_{23}$ is obtained in the same way as the proof of
Lemma~\ref{exirespi}.
Namely, we apply Proposition~\ref{prop811} in the case $L_{13} = L_{12}
\times_{X_2} L_{23}$.
Then the triple fiber product
\begin{equation}\label{171888-}
\Delta\times_{X_1^2 \times X_2^2 \times X_3^2}
\bigl(\tilde L_{12} \times \tilde L_{23} \times
\tilde L_{13}\bigr) = \bigcup_{a \in \mathcal A_{123}}R_{123}(a)
\end{equation}
contains a `diagonal component' which is diffeomorphic to
$\tilde L_{13}$.
We can use Lemma~\ref{lem310} to prove the unique existence of
the relative spin structure $\sigma_{13}$ on $\tilde L_{13}$
so that the orientation bundle $\Theta_-$ induced on the
diagonal component $\tilde L_{13}$ is trivial.

Now we include boundary marked points.
In other words, we use the structure operation of the
bimodule structure
\begin{gather}
\mathfrak n_{k,m} \colon\
B_kCF[1](\Delta_{X_1} \times \Delta_{X_2} \times \Delta_{X_3})
\otimes CF(\Delta_{X_1} \times \Delta_{X_2} \times \Delta_{X_3},L_{12} \times L_{23} \times L_{13})\nonumber \\
\hphantom{\mathfrak n_{k,m} \colon} \ {} \otimes
B_mCF[1](L_{12} \times L_{23} \times L_{13})
\to
CF(\Delta_{X_1} \times \Delta_{X_2} \times \Delta_{X_3},L_{12} \times L_{23} \times L_{13})\label{171888}
\end{gather}
to define the structure operations
\begin{gather}
\mathfrak n^{\rm tri, < E_0,\varepsilon}_{k_{13},k_{12},k_{23}} \colon\
CF(L_{13})^{\otimes k_{13}}
\otimes
CF(L_{12})^{\otimes k_{12}}\nonumber
\\
\hphantom{\mathfrak n^{\rm tri, < E_0,\varepsilon}_{k_{13},k_{12},k_{23}} \colon} \ {} \otimes
 CF(L_{13},L_{12},L_{23})
\otimes
CF(L_{23})^{\otimes k_{23}}
\to
 CF(L_{13},L_{12},L_{23})\label{17188899}
\end{gather}
of the tri-module structure.
We use appropriate triples
$I_{13}$, $I_{12}$, $I_{23}$ which
splits $\{1,\dots,m\}$ $(m= k_{13} + k_{12} + k_{23})$,
in the same way as \eqref{form171422}.
We use again the Koszul sign rule.
Namely, \eqref{17188899} is a sum with Koszul sign
of \eqref{171888} over the choice of $I_{13}$, $I_{12}$, $I_{23}$.
(In \eqref{171888}, we put $k=0$, see (def.3).)

Now taking into account similar points as (dif.1) (dif.2) and
the next point (dif.3), the bi-module property
of \eqref{171888} implies the tri-module property \eqref{form828}
of
\eqref{17188899}, {\it with sign}.

 {(def.3)}
The bubble at the Lagrangian submanifold
$\Delta_{X_1} \times \Delta_{X_2} \times \Delta_{X_3}$
appears in
the compactification
${\mathcal M}(\Delta_{X_1} \times \Delta_{X_2} \times \Delta_{X_3},L_{12} \times L_{23} \times L_{13};a_-,a_+;E)$
but there is no corresponding boundary component in the
compactification of the moduli space
${\mathcal M}^{\rm reg}_{\rm DR}(L_{13},L_{12},L_{23};a_-,a_+;E)$.

In fact, the disk bubble at $\Delta_{X_1} \times \Delta_{X_2} \times \Delta_{X_3}$
corresponds to the {\it sphere} bubble
of an element of ${\mathcal M}^{\rm reg}_{\rm DR}(L_{13},L_{12},L_{23};a_-,a_+;E)$
at the seams depicted by the dotted lines in Figure~\ref{Drumcut}.
Since they are sphere bubbles and occurs in codimension $\ge 2$,
they do not contribute the formula.
In other words, we can consider only $k=0$ case of
\eqref{171888} and obtain a left
$CF(L_{13})$, $CF(L_{12})$ and right $CF(L_{23})$ tri-module
structure.

\subsection[Orientation of the moduli space of $Y$-diagrams]{Orientation of the moduli space of $\boldsymbol{Y}$-diagrams}
\label{oriYdiagarm}

In this subsection, we study orientation of the
moduli space of $Y$-diagrams.
We consider $Y$-diagram as in Figure~\ref{Figure91}
and Definition~\ref{def1016}.
We put 3 extra seams which are depicted by dotted
lines in Figure~\ref{Ydiagramori1} below.
The domain $\mathcal Y$ in Figure~\ref{Figure91}
is divided into three pieces $\mathcal Y_i$ ($i=1,2,3$).
The added seams divide each of $\mathcal Y_i$ into
two pieces $\mathcal Y_{i,+}$ and $\mathcal Y_{i,-}$
as depicted in Figure~\ref{Ydiagramori1}.

Definition~\ref{def1016}
defines a moduli space
\smash{$\overset{\ \text{\tiny $\circ\circ$}}{\mathcal M}_{\rm Y}(\vec a_{12},\vec a_{23},\vec a_{13};\vec a_{1},\vec a_{2},\vec a_{3},a_{\infty,-},\vec a_{\infty,+};E)$}.
We consider the case
when $\vec a_{12},\vec a_{23},\vec a_{13};\vec a_{1},\vec a_{2},\vec a_{3}$
are all empty sets and write it as
$\smash{\overset{\ \text{\tiny $\circ\circ$}}{\mathcal M}}_{\rm Y}(L_1,L_2,L_3;L_{12},L_{23},\allowbreak
L_{13};a_{\infty,-},\vec a_{\infty,+};E)$.

We consider an element
$(\Sigma;\vec z_{1},\vec z_{2},\vec z_{3};\vec z_{12},\vec z_{23},\vec z_{13};u_1,u_2,u_3;\gamma_{1},\gamma_2,\gamma_{3};\gamma_{12},\gamma_{23},\gamma_{13})$
of the moduli space
\smash{$\overset{\ \text{\tiny $\circ\circ$}}{\mathcal M}_{\rm Y}(L_1,L_2,L_3;L_{12},L_{23},
L_{13};a_{\infty,-},\vec a_{\infty,+};E)$}.
We restrict $u_i$ to
$\mathcal Y_{i,+}$ and $\mathcal Y_{i,-}$
and obtain~$u_i^+$ and $u_i^-$.

\begin{figure}[ht]
\centering
\includegraphics[scale=0.28]{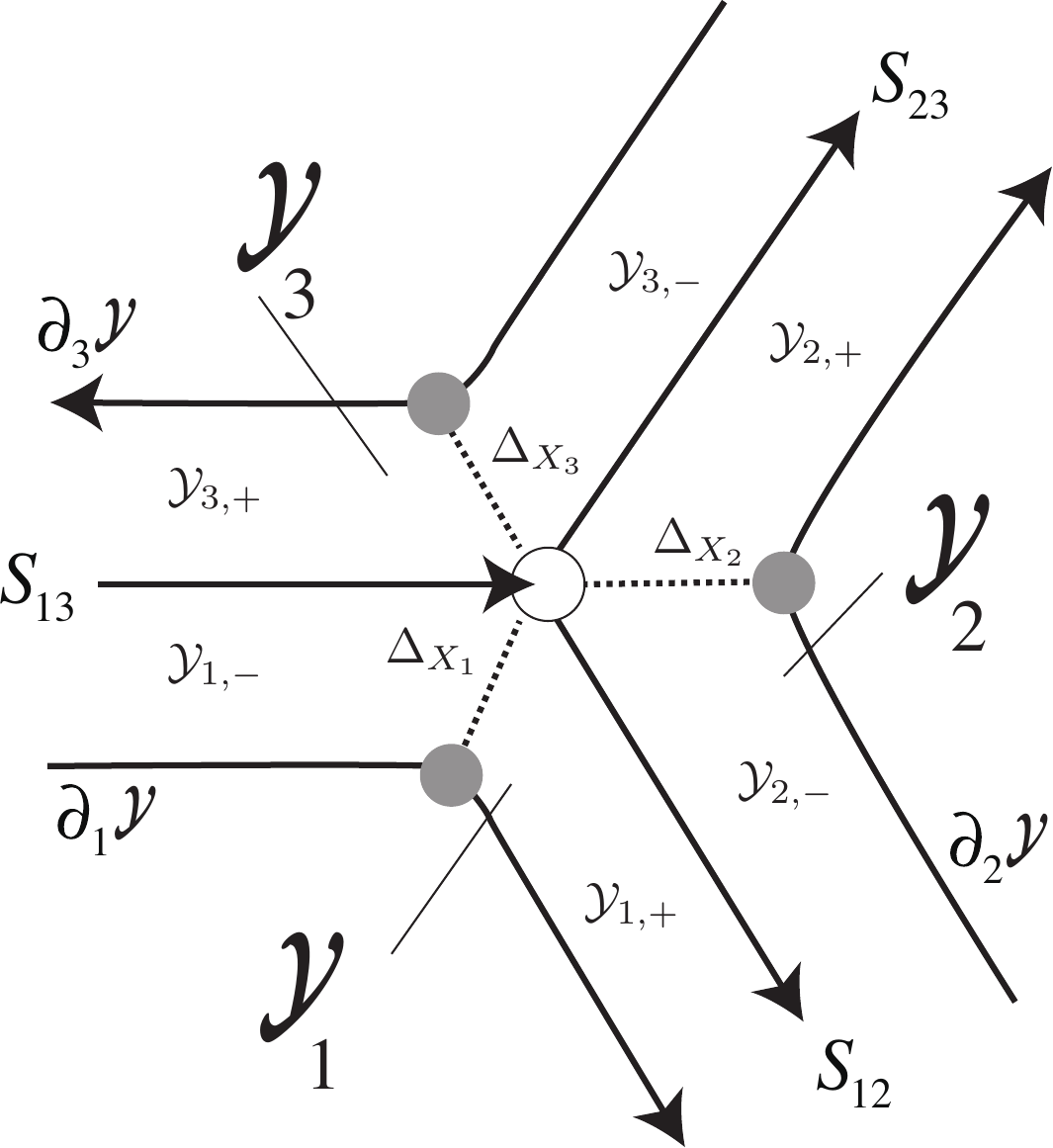}
\caption{Split domains in the Y-diagram.}
\label{Ydiagramori1}
\end{figure}

\begin{figure}[ht]
\centering
\includegraphics[scale=0.3]{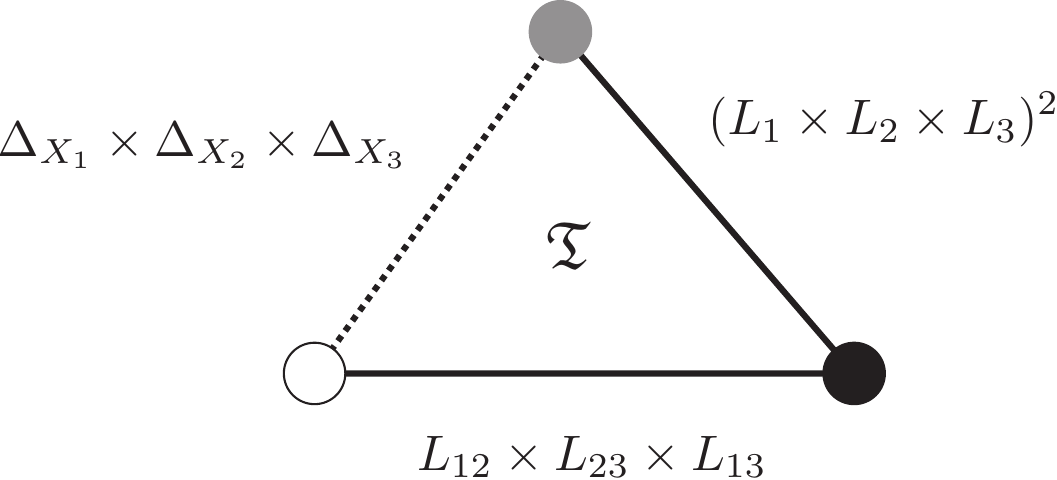}
\caption{Reglue maps from the Y-diagram.}
\label{FigureYdiagram2}
\end{figure}

We identify $\mathcal Y_{i,+}$ and $\mathcal Y_{i,-}$
with a triangle $\mathfrak T$ in Figure~\ref{FigureYdiagram2}.
We use a holomorphic map to identify $\mathcal Y_{i,+}$ with $\mathfrak T$
and an anti-holomorphic map to identify
$\mathcal Y_{i,-}$ with $\mathfrak T$.
By this identification, the point depicted by the white circle
(resp.\ the gray circles) in Figure~\ref{Ydiagramori1}
is sent to the point depicted by the white circle
(resp.\ the gray circle) in Figure~\ref{FigureYdiagram2}.
Three ends of the domain $\mathcal Y$ in Figure~\ref{Ydiagramori1}
is sent to the black circle in Figure~\ref{FigureYdiagram2}.
Thus $u_i^+$ and $u_i^-$, $i=1,2,3$, altogether induce a
pseudo-holomorphic map
$\hat u \colon \mathfrak T \to \prod_{i=1}^3 (-X_i \times X_i)$.
At the three boundary components, the map
$\hat u$ satisfies the boundary condition
given by the Lagrangian submanifolds
$\prod_{i=1}^3 \Delta_{X_i}$,
$L_{12} \times L_{23} \times L_{13}$,
$L_1^2 \times L_2^2 \times L_3^2$,
respectively.

The boundary conditions at the three vertices are obtained as follows.
$a_{\infty,-}$ assigns a
component $R_{123}(a_{\infty,-})$ of
the fiber product \eqref{171888-}
(see \eqref{decomp96}).
This boundary condition is used at the vertex
depicted by white circles in Figure~\ref{FigureYdiagram2}.
We next use $\vec a_{\infty,+} = (a_{\infty,+,12},a_{\infty,+,23},a_{\infty,+,13})$.
Then determine components $R_{ii'}(a_{\infty,+,ii'})$ of
$
\tilde L_i \times_{X_i} \tilde L_{ii'} \times_{X_{i'}} \tilde L_{i'}$
(see \eqref{decomp95}). Here $(ii') = (12),(23),(13)$.
Then the boundary condition at the black circles is given by
\[
R(\vec a_{\infty,+}): = R_{12}(a_{\infty,+,12}) \times R_{23}(a_{\infty,+,23})
\times R_{13}(a_{\infty,+,13}).
\]
We finally describe the boundary condition at
the vertex drawn by gray circle in
Figure~\ref{FigureYdiagram2}.

It should be a component of the fiber product
\[
\left(\prod_{i=1}^3 \Delta_{X_i}\right)
\times_{\prod_{i=1}^3 (-X_i \times X_i)}
L_1^2 \times L_2^2 \times L_3^2.
\]
We take the diagonal component
$\cong \tilde L_1 \times \tilde L_2 \times \tilde L_3$
and its fundamental class as the boundary condition.
We denote by
\[
\mathcal M^{\rm reg}_3\Bigl(\prod\Delta_{X_i},
L_{12} \times L_{23} \times L_{13},
L_1^2 \times L_2^2 \times L_3^2 ;\Delta,R_{123}(a_{\infty,-}),R(\vec a_{\infty,+});E\Bigr)
\]
the moduli space of such holomorphic triangles.

Thus the above construction defines a map
\begin{gather}
{\rm Dob} \colon\ \overset{\ \text{\tiny $\circ\circ$}}{\mathcal M}_{\rm Y}(L_1,L_2,L_3;L_{12},L_{23},
L_{13};a_{\infty,-},\vec a_{\infty,+};E)\nonumber \\
\hphantom{{\rm Dob} \colon} \ \to
\mathcal M^{\rm reg}_3\Bigl(\prod\Delta_{X_i},
L_{12} \times L_{23} \times L_{13},
L_1^2 \times L_2^2 \times L_3^2 ;
\nonumber \\
\hphantom{{\rm Dob} \colon \ \to{} } \Delta,R_{123}(a_{\infty,-}),R(\vec a_{\infty,+});E\Bigr).\label{form17121712}
\end{gather}
The orientation of the right-hand side of
\eqref{form17121712}
is defined by Proposition~\ref{prop329}.
We thus define the orientation
of
\smash{$\overset{\ \text{\tiny $\circ\circ$}}{\mathcal M}_{\rm Y}(L_1,L_2,L_3;L_{12},L_{23},
L_{13};a_{\infty,-},\vec a_{\infty,+};E)$}
so that \eqref{form17121712} preserves orientation.
This proves Proposition~\ref{prop9911}
(3).

We show the compatibility of the orientation
at the boundary below.
We consider the codimension one boundary component
of the target of \eqref{form17121712}.
We divide it into various cases.

Case 1.
Disk bubble at $\prod\Delta_{X_i}$.
There is no corresponding codimension one
boundary component in the source of \eqref{form17121712}.
In fact, this corresponds to the sphere bubble at the seams
depicted by the dotted lines in Figure~\ref{Ydiagramori1}.
This occurs in codimension $\ge 2$.

Case 2.
Disk bubble at
$L_{12} \times L_{23} \times L_{13}$.
This corresponds to the disk bubble
in Figure~\ref{Figure94}.
The homotopy class $\beta$ of such disk is determined by
\[
\pi_2\left(\prod_{i=1}^3 (X_i \times X_i),L_{12} \times L_{23} \times L_{13}
\right)
\cong \prod_{(ii') = (12),(23),(13)}
\pi_2(X_i\times X_{ii'};L_{ii'}).
\]
If $\beta = (\beta_{12},\beta_{23},\beta_{13})$ and
at least two of $\beta_{12}$, $\beta_{23}$, $\beta_{13}$ are
nonzero, then there is no corresponding component
in the source of \eqref{form17121712}.
The reason is the same as the reason why Figure~\ref{Figure17-2}
appears in codimension $\ge 2$.
Therefore, it suffices to consider the case
when only one of $\beta_{12}$, $\beta_{23}$, $\beta_{13}$
is nonzero.
The boundary component corresponding to such cases corresponds to the boundary component
described by Figure~\ref{Figure94}.
Therefore, the orientation is consistent at this boundary component.

Case 3.
Disk bubble at $L_1^2 \times L_2^2 \times L_3^2$.
The homotopy class of such bubble is given
by $(\prod(\pi_2(X_i,L_i)))^2$.
By the same reason as above it suffices to consider the
case only one of those~6 factors is nonzero.
Then it corresponds to the boundary
component depicted by Figure~\ref{Figure95} in the
right-hand side.
Thus the orientation is consistent at this boundary.

We next consider the boundary component corresponding to the
three vertices of $\mathfrak T$.

Case 4.
The boundary component corresponding to the white vertex.
This is described by the fiber product of
\begin{equation}\label{form171313}
\mathcal M\Bigl(\prod\Delta_{X_i};
L_{12} \times L_{23} \times L_{13};R_{123}(a)\Bigr)
\end{equation}
with
\[
\mathcal M^{\rm reg}_3\Bigl(\prod\Delta_{X_i},
L_{12} \times L_{23} \times L_{13},
L_1^2 \times L_2^2 \times L_3^2
;\Delta,R_{123}(a),R(\vec a_{\infty,+});E\Bigr)
\]
over $R_{123}(a)$.
We apply the identification of \eqref{dobe17111rev}
and \eqref{form171313}.
Then this boundary component corresponds to one
in Figure~\ref{Figure96}.
Thus the orientation is consistent at this boundary.

Case 5.
The boundary component corresponding to the black vertex.
This is described by the fiber product of
\[
\mathcal M^{\rm reg}_3\Bigl(\prod\Delta_{X_i},
L_{12} \times L_{23} \times L_{13},
L_1^2 \times L_2^2 \times L_3^2
;\Delta,(a_{\infty,-}),R(\vec a'_{\infty,+});E\Bigr)
\]
with
$
\mathcal M\bigl(
L_{12} \times L_{23} \times L_{13},
L_1^2 \times L_2^2 \times L_3^2;R(\vec a'_{\infty,+}),R(\vec a_{\infty,+})\bigr)
$
taken over $R(\vec a'_{\infty,+})$.
We put $\vec a'_{\infty,+} = (a'_{\infty,+,12},a'_{\infty,+,23},
a'_{\infty,+,13})$.
In the case when two among the three inequalities
$a'_{\infty,+,12} \ne a_{\infty,+,12}$,
$a'_{\infty,+,23} \ne a_{\infty,+,23}$,
$a'_{\infty,+,13} \ne a_{\infty,+,13}$
hold, the corresponding component in the
source of~\eqref{form17121712}
has codimension $\ge 2$.
(The reason is the same as Case 2.)
Therefore, it suffices to consider
the case when exactly one of
the three inequalities~${a'_{\infty,+,12} \ne a_{\infty,+,12}}$,
$a'_{\infty,+,23} \ne a_{\infty,+,23}$,
$a'_{\infty,+,13} \ne a_{\infty,+,13}$
hold.
This case corresponds to one
depicted in Figure~\ref{Figure97}.
Thus the orientation is consistent at this boundary.

Case 6.
The boundary component corresponding to the gray vertex.
The corresponding component in the
source of \eqref{form17121712}
has codimension $\ge 2$.
In fact, it corresponds to the case
when there is a disk bubble exactly
at the gray vertex of Figure~\ref{Ydiagramori1}.
This is a codimension 2 phenomenon.

We thus checked the consistency of
the orientation at the codimension
one component.
We proved Proposition~\ref{prop9911}\,(3).

\begin{proof}
[Proof of Proposition~\ref{prop912}\,(1)]
Given $L_1$, $L_2$, $L_{12}$, $L_{23}$
we consider the case when
$\tilde L_{13} = \tilde L_{12} \times_{X_2}
\tilde L_{23}$ and
$\tilde L_3 = \tilde L_2 \times_{X_2}
\tilde L_{23} = \tilde L_1 \times_{X_1}
\tilde L_{23}$.

We take the diagonal component as $a$ for $R_{123}(a)$
and $a_{\infty,+,12}$, $a_{\infty,+,23}$, $a_{\infty,+,13}$
for $R(\vec a_{\infty,+})$.

Given relative spin structure $\sigma_1$,
$\sigma_{12}$ of $\tilde L_{1}$, $\tilde L_{12}$,
we have chosen the relative spin structure
$\sigma_2$ of $\tilde L_2$ so that the
local system associated to
$R_{12}(a_{\infty,+,12}) \cong \tilde L_2$
is trivial.
We also have chosen the relative
spin structure $\sigma_{13}$ of $\tilde L_{13}$ so that the
local system associated to
$R_{123}(a) \cong \tilde L_{13}$ is trivial.

We consider the moduli space
\smash{$\overset{\ \text{\tiny $\circ\circ$}}{\mathcal M}(L_1,L_2,L_3;L_{12},L_{23},
L_{13};a_{\infty,-},\vec a_{\infty,+};0)$}
consisting of constant map.
It corresponds to the
moduli space of constant maps
\[
\mathcal M^{\rm reg}_3\Bigl(\prod\Delta_{X_i},
L_{12} \times L_{23} \times L_{13},
L_1^2 \times L_2^2 \times L_3^2
;\Delta,(a_{\infty,-}),R(\vec a'_{\infty,+});0\Bigr).
\]
This space is diffeomorphic to $\tilde L_3$ and is oriented.

Therefore, by Proposition~\ref{prop329},
for any choice of relative
spin structure $\sigma_3$ of $\tilde L_3$,
the local system induced on
$R_{13}(a_{\infty,+,13}) \cong \tilde L_3$
is isomorphic to one on
$R_{23}(a_{\infty,+,23}) \cong \tilde L_3$.

If \smash{$\sigma_3 = \sigma_3^{(1)}$}, then
it is trivial for $R_{23}(a_{\infty,+,23}) \cong \tilde L_3$.
It \smash{$\sigma_3 = \sigma_3^{(2)}$} then
it is trivial for $R_{13}(a_{\infty,+,13}) \cong \tilde L_3$.
Therefore, \smash{$\sigma_3^{(1)} = \sigma_3^{(2)}$}.
Proposition~\ref{prop912}\,(1) is proved.
\end{proof}

We next include marked points on the boundary of
$Y$-diagram and will prove the equality~\eqref{form925}
{\it with sign}.

Including marked points on the boundary
the target of \eqref{form17121712}
becomes the moduli space
which is used to define
structure operation
\begin{gather}\label{form171515}
\mathfrak m_{m_1+m_2+m_3+3}\colon\
BCF_{m_1}\left(\prod\Delta_{X_i}\right)
\otimes
CF(\Delta) \otimes
BCF_{m_2}
\left(
L_{12} \times L_{23} \times L_{13}
\right) \otimes CF(R_{123}(a))\nonumber \\
\hphantom{\mathfrak m_{m_1+m_2+m_3+3}\colon} \ {} \otimes
BCF_{m_3}
\bigl(L_1^2 \times L_2^2 \times L_3^2\bigr)
\to CF(R(\vec a_{\infty,+});E)),
\end{gather}
of a filtered $A_{\infty}$ category
assigned to $\prod(-X_i\times X_i)$
and its Lagrangian submanifolds
$\bigl\{\prod\Delta_{X_i},
L_1^2 \times L_2^2 \times L_3^2,L_{12} \times L_{23} \times L_{13}\bigr\}$.
It satisfies the $A_{\infty}$ relation.
We convert
first and second factor of the output
$
CF(R(\vec a_{\infty,+}))
=
CF(R(a_{\infty,+,12})) \otimes
CF(R(a_{\infty,+,23})) \otimes
CF(R(a_{\infty,+,13}))
$
to the input by duality.

Then
the operation \eqref{form171515} with an appropriate sign
becomes the operator $\mathscr{Y\!\!T}^{E,\varepsilon}_{k_{12},k_{23},k_{13};k_1,k_2,k_3}$
in \eqref{operation923}.
Here $m_1 = 0$, $m_2 = k_{12} + k_{23} + k_{13}$,
$m_3 = k_1 + k_2 + k_3$.
We use the Koszul rule in the same way as
Sections~\ref{orisimpquilt}, \ref{oridrum}, to define the sign.
Then taking into account (def.2),
the~$A_{\infty}$ relation for \eqref{form171515}
becomes the equality \eqref{form925}
with sign.
(We use the fact that some of the terms of
the~$A_{\infty}$ relation of \eqref{form171515}
is absent in \eqref{form925}.
The reason is explained in Cases 1--6 above.)

\subsection{Orientation of the moduli space of double pants
diagrams}
\label{oripants}

In this subsection, we study orientation of the
moduli space of double pants
diagrams.
We draw double pants as in Figure~\ref{Figure11-2} and
put 12 seams as in Figure~\ref{oridpants1} below.
In Figure~\ref{oridpants1}, the new seams are depicted
by dotted lines.
We have new vertices also. There are 4 black vertices
which are new.
The circles of Figure~\ref{Figure11-2} are
now depicted by white vertices in Figure~\ref{oridpants1}.
There are 4 white vertices in Figure~\ref{oridpants1}.
Here the outer circle in Figure~\ref{oridpants1} should
be regarded as a vertex. \big(In other words, the
domain should be regarded as $S^2$.\big)

\begin{figure}[ht]
\centering
\includegraphics[scale=0.4]{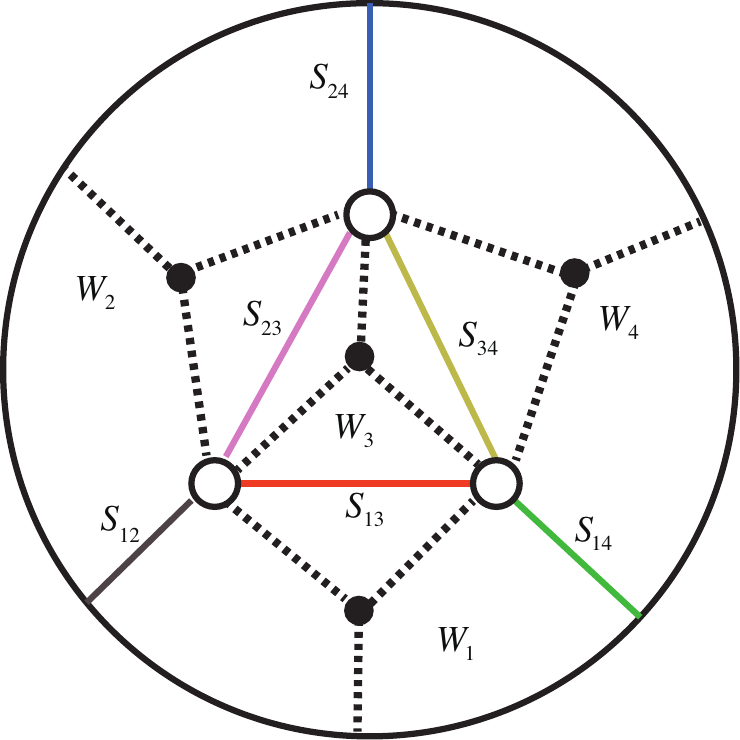}
\caption{Adding seams to double pants.}
\label{oridpants1}
\end{figure}

We cut the domain in Figure~\ref{oridpants1}
and obtain the triangle $\mathfrak T$ in the Figure~\ref{Figurepants2} below.
The maps $u_i$ ($i=1,2,3,4$) in Definition~\ref{def1116}
induces a map
$
\hat u \colon\mathfrak T \to \prod_{i=1}^4 (-X_i \times X_i)$.
Its boundary condition is given by
$\Delta_{X_1} \times \Delta_{X_2} \times \Delta_{X_3} \times \Delta_{X_4}$
for the (two) dotted edges and $
\prod_{(ij)=(12),(13),(14),(23),(24),(34)}L_{ij}
$
 for the
other edge.

\begin{figure}[ht]
\centering
\includegraphics[scale=0.4]{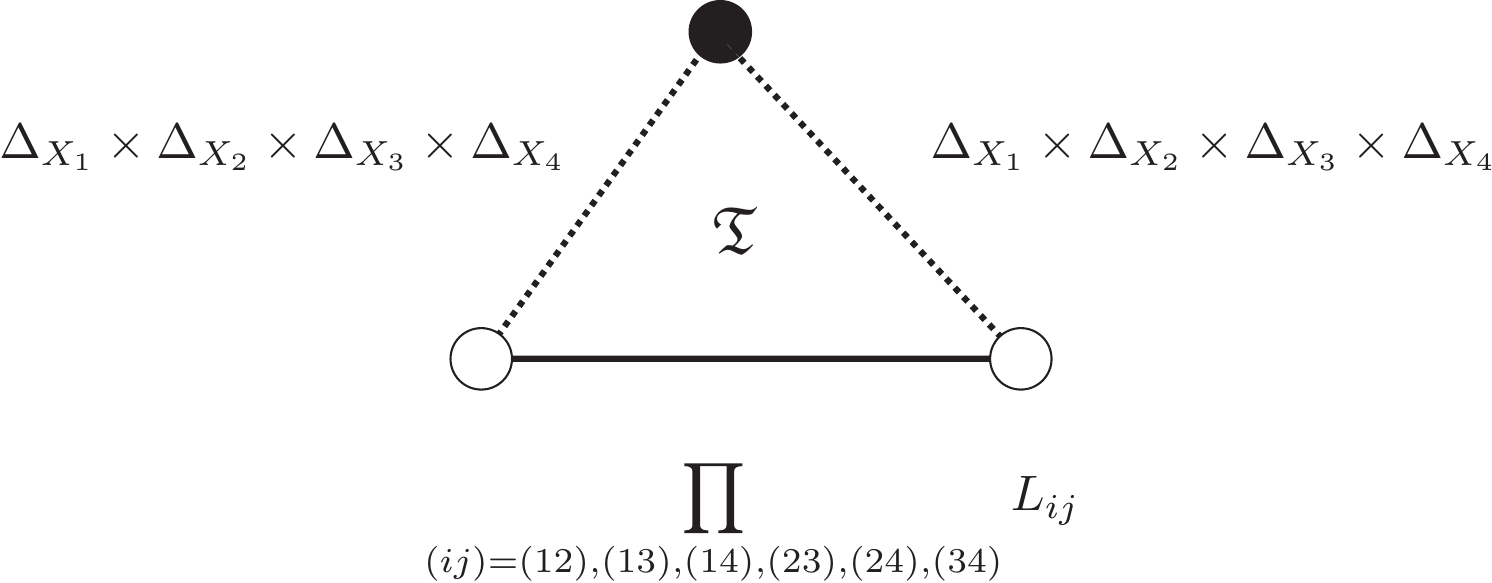}
\caption{Regluing double pants.}
\label{Figurepants2}
\end{figure}

We thus obtain an identification of
the moduli space
\smash{$\overset{\ \text{\tiny $\circ\circ$}}{\mathcal M}((\vec a_{ii'};i,i');(a_{ii'i''};i,i',i'');E)$}
of Definition~\ref{def1116} with the
moduli space of pseudo-holomorphic triangles
depicted in Figure~\ref{Figurepants2}.
Therefore, applying Proposition~\ref{prop329}
to the moduli space of pseudo-holomorphic triangles
depicted in Figure~\ref{Figurepants2}, we obtain
an orientation of the moduli space
\smash{$\overset{\ \text{\tiny $\circ\circ$}}{\mathcal M}_{\rm DP}((\vec a_{ii'};i,i');(a_{ii'i''};i,i',i'');E)$}.
The compatibility at the boundary
can be proved in the same way as Section~\ref{oriYdiagarm}.
It implies Proposition~\ref{prop10811}\,(3).

The proof of Proposition~\ref{prop1212}\,(1) is the same as the proof of Proposition~\ref{prop912}\,(1)
given in Section~\ref{oriYdiagarm}.

We consider marked points on the seam in Figure~\ref{oridpants1}.
Then the double pants transformation~$\mathscr{DPT}^{E,\varepsilon}$
in Definition~\ref{defnnew117}
is defined by using appropriate $A_{\infty}$
operation (which is defined from \eqref{Figurepants2}
with marked points added).
The sign is defined by Koszul rule.
So again taking into account (def.2),
$A_{\infty}$ formula implies formula \eqref{form1118}
{\it with Koszul sign}.

\subsection[Orientation and sign for $A_\infty$-structure in the Morse--Bott case]{Orientation and sign for $\boldsymbol{ A_{\infty}}$-structure in the Morse--Bott case}
\label{oriAinfMB}

The proof of $A_{\infty}$-formula {\it with sign} is written
in detail in the case of a single embedded Lagrangian submanifold
in \cite{fooobook2,fooonewbook,ST} etc.
For an immersed Lagrangian submanifold which has transversal self-intersection,
it is written in detail in \cite{AJ}. The latter implies the $A_{\infty}$
formula with sign in the
case when we have finitely many immersed Lagrangian submanifolds
which have transversal self-intersection and are mutually transversal.
We can prove it in the case of immersed Lagrangian submanifold
which has clean self-intersection (Morse--Bott type) in a~similar way.
Since it is not easy to find a reference which describes Morse--Bott case in detail,
we below explain the way to obtain the orientation which gives
$A_{\infty}$-formula with Koszul sign in such a~case.

In this subsection, we follow Akaho--Joyce's method in \cite{AJ} and
will explain how we modify it to generalize to the Morse--Bott case.
In the paper \cite{ono2}, written by Kaoru Ono, the way to generalize
\cite{fooonewbook} to the Morse--Bott situation will be written.

We consider the moduli space
${\mathcal M}(L;\vec a;E)$
defined in equation~\eqref{def33314}.
It comes with evaluation maps~\eqref{form3666}
\begin{equation}\label{form36662}
{\rm ev} = ({\rm ev}_0,\dots,{\rm ev}_k)\colon\
{\mathcal M}(L;\vec a;E) \to L(\vec a).
\end{equation}
Here
$L(\vec a) = L(a_0) \times \dots \times L(a_k)$
is a direct products of connected components
of $\tilde L \times_X \tilde L$. In
\cite{AJ}, $\tilde L \times_X \tilde L$ minus diagonal components
is written as $R$. (In their case, $R$ is a finite set. In our
case, it is a disjoint union of smooth compact manifolds.)
In \cite[p.~425, equation~(50)]{AJ}, Akaho and Joyce
take a product of ${\mathcal M}(L;\vec a;E)$
with vector spaces associated to each point $R$
to obtain \smash{$\widehat{\mathcal M}(L;\vec a;E)$}.
(They use the notation \smash{$\widetilde{\mathcal M}$}.
The author changes it to \smash{$\widehat{\mathcal M}$} since
\smash{$\widetilde{\mathcal M}$} is used in Definition~\ref{defn314}.)

Let $x \in L(a)$ which is not in a diagonal component.
We consider operators
\begin{gather*}
\overline\partial_{Z_-,\lambda_x} \colon\ L^2_k(Z_-;T_xX;\lambda_a;\delta)
\to L^2_{k-1}(Z_-;T_xX;\delta),
\\
\overline\partial_{Z_+,\lambda_x} \colon\ L^2_k(Z_+;T_xX;\lambda_a;\delta)
\to L^2_{k-1}(Z_+;T_xX;\delta)
\end{gather*}
as in \eqref{3535}.
Here we fix a choice of $\lambda_x$.
As is proved in \cite[Proposition 5.15]{AJ}, the
definition of orientation and sign which we describe below is independent
of such a choice.
To study orientation problem we can work locally on $L(\vec a)$.
So when $x$ is in a (small) neighborhood of given $x_0$ we can and will take a choice of
$\lambda_x$ depending continuously on $x$.
We can also perturb appropriately so that \smash{$\overline\partial_{Z_{\pm},\lambda_x}$}
are surjective.
Then \smash{$\operatorname{Ker}\overline\partial_{Z_{-},\lambda_x}$}
defines a vector bundle on~$L(a)$. (More precisely, on a neighborhood of $x_0$ of $L(a)$.)
We denote its total space by $\hat L(a)$.
In case~$L(a)$ is a diagonal component, we define $\hat L(a) = L(a)$.
We put $\hat L(\vec a) = \hat L(a_0) \times \dots \times \hat L(a_k)$.
Following \cite[equation~(59)]{AJ} we put
\begin{equation}\label{1733}
\widehat{\mathcal M}(L;\vec a;E) := {\mathcal M}(L;\vec a;E) \times_{L(\vec a)} \hat L(\vec a).
\end{equation}
Note that the line bundle \smash{$\Theta^-_{a_i}$} appearing in
Proposition~\ref{prop329} is the determinant line bundle of~${\hat L(a_i) \to L(a_i)}$.
Therefore, Proposition~\ref{prop329} implies
that $\widehat{\mathcal M}(L;\vec a;E)$ has a canonical orientation.

For $i=0$, the convention of $\tilde{\rm ev}_0$ in \cite[equation~(62)]{AJ}
is slightly inconsistent with our convention of ${\rm ev}_0$ at the point
which we explain below.
Let $\sigma$ be the involution $\tilde L \times_X \tilde L \to \tilde L \times_X \tilde L$
defined by $\sigma(x,y) = (y,x)$. We denote by $L(\sigma(a))$ the component
to which $L(a)$ is
sent by~$\sigma$. (If $L(a)$ is a diagonal component, then $\sigma(a) = a$.)
The $0$-th evaluation map used by \cite{AJ} is the composition $\sigma\circ{\rm ev}_0$,
where ${\rm ev}_0$ is the evaluation map \eqref{form36662}.
From now on, {\it in this subsection} we use \cite{AJ}'s convention.
Namely, we change the definition of ${\rm ev}_0$ to those
by Akaho--Joyce.\footnote{The map ${\rm ev}_0$ used to define \eqref{1733} is the one in \eqref{form36662}.
With this choice our $\widehat{\mathcal M}$ coincides with \cite{AJ}'s $\widetilde{\mathcal M}$.}
This is only a matter of notation and there is no mathematical
difference.
Note that then~${L(\vec a) = L(\sigma(a_0)) \times L(a_1) \times \dots \times L(a_k)}$
in this convention.

We next describe the evaluation map following \cite[equations~(61) and (62)]{AJ}.
We first define~$\tilde L(a)$ by modifying a bit \cite{AJ}'s $\tilde R$. (We need slight
modification since our $L(a)$ may not be discrete.)
Let $x = (p,q) \in L(a)$. We take $\lambda_x \in \mathcal P^a_x$. (See Definition~\ref{lem34}.)
Then as we proved in~\eqref{form382}, we have
\begin{equation}\label{split333}
T_x\tilde L \cong \operatorname{Ker}\overline\partial_{Z_-,\lambda_x} \oplus \operatorname{Ker}\overline\partial_{Z_+,\lambda_x} \oplus T_xL(a).
\end{equation}
We take a vector bundle on $L(a)$ (more precisely on a neighborhood of $x_0$ in $L(a)$)
whose fiber at $x$ is
$\operatorname{Ker}\overline\partial_{Z_-,\lambda_x} \oplus \operatorname{Ker}\overline\partial_{Z_+,\lambda_x}$ and
define $\tilde L(a)$ to be the total space of this vector bundle.
We remark that $L(a)$ may not be orientable. However, since we assume $L$ to be
oriented \smash{$\tilde L(a)$} is oriented.

We put
$\tilde L(\vec a) = \tilde L(\sigma(a_0)) \times \tilde L(a_1) \times \dots \times \tilde L(a_k)$.
For $\sigma(x) \in L(\sigma(a))$, we take $\lambda_{\sigma(x)}$ to be the
opposite path to $\lambda_x$.
Then we have canonical isomorphisms
\[
\operatorname{Ker}\overline\partial_{Z_-,\lambda_{\sigma(x)}}
\cong \operatorname{Ker}\overline\partial_{Z_+,\lambda_x},\qquad
\operatorname{Ker}\overline\partial_{Z_+,\lambda_{\sigma(x)}}
\cong \operatorname{Ker}\overline\partial_{Z_-,\lambda_x}.
\]
In particular, the involution $\sigma \colon L(a) \to L(\sigma(a))$ lifts
to an involution $\sigma \colon \tilde L(a) \to \tilde L(\sigma(a))$.

We remark that \eqref{split333} implies that
$
\dim \tilde L(a) = n
$
for any $a$. (Here $n = \dim L$.)
\big(The right-hand side is independent of $a$. This is an advantage to replace $L(a)$ by $\tilde L(a)$.\big)

Now we define
\smash{$
\tilde{\rm ev} = (\tilde{\rm ev}_0,\dots,\tilde{\rm ev}_k)\colon
\widehat{\mathcal M}(L;\vec a;E) \to \tilde L(\vec a)
$}
in a similar way as \cite[equations~(61) and (62)]{AJ} as follows.
We remark that an element of~$\widehat{\mathcal M}(L;\vec a;E)$ is~$(\mathfrak x;(\xi_i)_{i=0}^k)$, where
$\xi_i \in \operatorname{Ker}\overline\partial_{Z_-,\lambda_{{\rm ev}_i(\mathfrak x)}}$
for $i\ne 0$ and $\xi_0 \in \operatorname{Ker}\overline\partial_{Z_-,\lambda_{\sigma({\rm ev}_0(\mathfrak x))}}$.
We put
\[
\tilde{\rm ev}_i(\mathfrak x) = ({\rm ev}_i(\mathfrak x),\xi_i) \in \tilde L(a_i), \quad {i\ne 0}, \qquad
\tilde{\rm ev}_i(\mathfrak x) = ({\rm ev}_0(\mathfrak x),\sigma(\xi_0)) \in \tilde L(\sigma(a_0)), \quad {i = 0}.
\]
We remark that
\smash{$
\sigma(\xi_0) \in \operatorname{Ker}\overline\partial_{Z_+,\lambda_{{\rm ev}_0(\mathfrak x)}}$}.
We use this fact and Theorem~\ref{thekuraexist} to prove the following:
\begin{prop}
$\widehat{\mathcal M}(L;\vec a;E)$
has Kuranishi structure with corners whose normalized boundary is the disjoint union of the
fiber products as follows
\begin{align}
\partial \widehat{\mathcal M}(L;\vec a;E)
={}&
\coprod_{b,i,j
\atop E_1 + E_2 = E} (-1)^*
\widehat{\mathcal M}(L;\vec a(b,i,j,2);E_2) {}_{\tilde{\rm ev}_{0}}\nonumber\\
& \times_{\tilde{\rm ev}_{i+1}}
\widehat{\mathcal M}(L;\vec a(b,i,j,1);E_1).\label{form3472}
\end{align}
\end{prop}
\eqref{form3472} is mostly the same as \eqref{form347},\footnote{\eqref{form347} is the case of
Theorem~\ref{thekuraexist}\,(3) when the graph $\Gamma$ has one interior vertex.} but we replace
$\mathcal M$ by $\widehat{\mathcal M}$ and ${\rm ev}$ by $\tilde{\rm ev}$.
We also remark that the order of first and second factors in the right-hand side
of \eqref{form3472} is the same as \cite[equation~(73)]{AJ}
but is opposite to \cite[equation~(20.11)]{fooonewbook}.
(See \cite[the last part of Section~4.2]{AJ}.)
In this subsection, we follow \cite{AJ}.
Note that in~\eqref{form3472} $\tilde{\rm ev}_{i+1}$ is used
at the place where $\tilde{\rm ev}_{i}$ is used in \cite{AJ}.
This is because of the convention used in $\vec a(b,i,j,2)$ and $\vec a(b,i,j,1)$
and is not related to the mathematical contents of the formula.

Now we state the compatibility of the orientations
with the isomorphism \eqref{form3472}, that is, the sign $*$ in \eqref{form3472}.
As we mentioned already, Proposition~\ref{prop329} can be restated that
$\widehat{\mathcal M}(L;\vec a;E)$ is oriented.
Also in \eqref{form3472} we take the fiber product over $\tilde L(a)$, which is always
$n$-dimensional and oriented.

Therefore, the situation is the same as the case of $A_{\infty}$ algebra
associated to a single embedded Lagrangian submanifold.
We require
$* = n+(i+1)+(i+1)k_2$.
This is the same as \cite[equation~(73)]{AJ} except $i$ is replaced by $i+1$.
(It is different from \cite[equation~(21.7)]{fooonewbook} by the above mentioned reason.)

The rest of the argument is mostly the same as \cite{AJ}.
Let $P_i = (P_i,f_i)$ be a chain in $L(a_i)$. More precisely, it is a smooth singular chain with
an orientation of $\Theta_{a_i,-} \otimes \operatorname{Det} N_{P_i}L(a_i)$ given.\footnote
{Here we denote by $N_{P_i}L(a_i)$ the normal bundle. The determinant line bundle of the
normal bundle is defined even in the case $f_i \colon P_i \to L(a_i)$ is not an immersion.}
(Such a chain can be used to calculate the cohomology of $L(a_i)$ with $\Theta_{a_i,-}$
coefficient. Note that in case $\Theta_{a_i,-}$ is trivial, the chain is co-oriented and
so is related to cohomology rather than homology.
We remark that $L(a_i)$ may not be orientable.
Even in such a case the set of singular chains with co-orientation
is a model of its cohomology.)
We put
\begin{equation}\label{741form}
\tilde P_i = P_i \times \operatorname{Ker}\overline\partial_{Z_-,\lambda_{\sigma(x_0)}}.
\end{equation}
Here we take $\lambda_{\sigma(x_0)}$ for $x_0 \in L(a_i)$ and assume the image of $f_i$ is in a small
neighborhood of $x_0$.
Compare \cite[equation~(68)]{AJ}.
Note that
\begin{equation}\label{sigmaplusminus}
\operatorname{Ker}\overline\partial_{Z_-,\lambda_{\sigma(x_0)}}
\cong
\operatorname{Ker}\overline\partial_{Z_+,\lambda_{x_0}}.
\end{equation}
Then by the definition of \smash{$\Theta_{a_i,-}$} and
\eqref{split333} the chain $\tilde P_i$ is oriented.\footnote
{It might be more natural to say that it is co-oriented.
However, in our situation the ambient manifold $\tilde L(a_i)$
is oriented. So `oriented' and `co-oriented' are equivalent.}
Using
\eqref{sigmaplusminus}, we obtain~${\tilde f_i \colon \tilde P_i \to \tilde L(a_i)}$ in an obvious way.

Now we define
\begin{equation}\label{17403form}
\widehat{\mathcal M}(L;\vec a;E;\vec P)
:= (-1)^* \widehat{\mathcal M}(L;\vec a;E) \times_{\tilde L(a_1) \times \dots \times \tilde L(a_k)}
\tilde P_1 \times \dots \times \tilde P_k.
\end{equation}
Here we use $\tilde{\rm ev}_i$ and $\tilde f_i$ to define the fiber product.
The sign is
\begin{equation}\label{1743form}
* = (n+1)\sum_{\ell=1}^{k} (k-\ell)\deg P_{\ell}.
\end{equation}
Since $\deg$ in \cite{AJ} is the shifted degree $\deg'$ in FOOO's notation (see \cite[p.~418]{AJ}),
\eqref{1743form} exactly coincides with the sign in \cite[equation~(79)]{AJ}.
(In \eqref{1743form}, $\deg P_{\ell}$ is one in FOOO's convention.)
Note that the degree of the chain in $L(a_i)$ as an element of $CF(L(a_{i-1}),L(a_{i}))$
is shifted from its codimension in $L(a_i)$
by the dimension of $\operatorname{Ker}\overline\partial_{Z_-,\lambda_{x_0}}$.
(It is the Morse index in the related context of Morse--Bott
theory.)
Therefore, the degree of $P_i$ as an element of $CF(L(a_{i-1}),L(a_{i}))$
is equal to the codimension of $\tilde P_i$ in $\tilde L(a_i)$.

It is easy to see that
$\widehat{\mathcal M}(L;\vec a;E;\vec P)$ coincides with
\begin{equation}\label{17403form2}
{\mathcal M}(L;\vec a;E;\vec P) =
{\mathcal M}(L;\vec a;E) \times_{L(a_1) \times \dots \times L(a_k)}
P_1 \times \dots \times P_k
\end{equation}
as spaces with Kuranishi structure (if we forget the orientation).
The reason we rewrite \eqref{17403form2} to~\eqref{17403form} is
then the correction term to orientation is easier to write down.
The map $\tilde{\rm ev}_0 \colon \smash{\widehat{\mathcal M}}(L;\allowbreak \vec a;E) \to \tilde L(a_0)$
induces
\smash{$
\tilde{\rm ev}_0 \colon \widehat{\mathcal M}(L;\vec a;E;\vec P) \to \tilde L(a_0)$}.
If we triangulate the domain, it gives singular chains of $\tilde L(a_0)$.
It is easy to see that those singular chains are related to the
singular chains obtained from
$
{\rm ev}_0\colon {\mathcal M}(L;\vec a;E;\vec P) \to L(a_0)
$
by the formula \eqref{741form}.

Now the rest of the construction is entirely the same as
\cite[pp.~434--444]{AJ}
and we obtain operations which satisfy $A_{\infty}$ relations with Koszul sign,
in the singular chain complex model.
As is explained in \cite{fooo010} (see the discussion around
formula \cite[pp.~190--191]{fooo010}), the sign and orientation in the
singular homology model induces one in the
de Rham model.\footnote{Since we converted the situation to the case
when all the fiber products involved are taken over $n$-dimensional
oriented manifolds, we can also import the method of \cite{fooonewbook}.}

In the paper \cite{ono2} by Kaoru Ono, the direct discussion
based on de Rham model is given.

\section{Concluding remarks}
\label{sec:remark}

\subsection[What we need to convert informal Definition \ref{def11} / Informal Summary \ref{thm02} into formal ones]{What we need to convert informal Definition~\ref{def11} /\\ Informal Summary \ref{thm02} into formal ones}
\label{issueformal}

In this subsection, we explain certain issues which will appear when
one tries to give a rigorous versions of
Definition~\ref{def11} or Informal Summary \ref{thm02}.
We explain the following three points:
\begin{enumerate}\itemsep=0pt
\item[(A)]
In Theorem~\ref{thm03}, we take a finite set of
Lagrangian submanifolds (not all the Lagrangian
submanifolds) and the object set of the curved
filtered $A_{\infty}$ category is this finite set.
\item[(B)]
The geometric transformation of a Lagrangian submanifold
$L_1$ by a Lagrangian correspondence $L_{12}$ is defined
under certain transversality assumptions.
The composition of Lagrangian correspondences is defined
under certain transversality assumptions.
\item[(C)]
The commutativity of several diagrams such as \eqref{Diagram1} or
\eqref{thm18} is up to homotopy equivalence and is not strict.
\end{enumerate}

We elaborate on those points below.
\begin{enumerate}\itemsep=0pt
\item[(A)]
In Theorem~\ref{thm03},
we take a {\it finite} set of (immersed and spin) Lagrangian
submanifolds $\mathbb L$ of $(X,\omega)$ and
an object of our filtered $A_{\infty}$ category is a pair
$(L,b)$ of an element $L$ of $\mathbb L$ and its
bounding cochain $b$.
This category of course depends on the choice of $\mathbb L$
and so is not canonically associated to $(X,\omega)$.
A natural way to make it more canonical is
taking all the Lagrangian submanifolds.
There is an issue for such a construction.

First we use a trick in Section~\ref{subsec:Ainfcatim}
to reduce the construction of a filtered $A_{\infty}$ category
to one of a filtered $A_{\infty}$ algebra, by taking disjoint
union of all the elements of $\mathbb L$ and regarding it as a
single immersed Lagrangian submanifold.
This trick does not work if $\mathbb L$ has infinite order.
However, this point itself does not seem to be so serious
since we used this trick mainly to shorten the paper.

The other and more essential issue is gappedness.
Our construction in Section~\ref{sec:HFIm} is
based on the induction on energy filtration.
We took and fix a discrete submonoid $G =
\{0=E_0,E_1,\dots\}$ of $\R_{\ge 0}$
and construct an $A_{\infty}$ category modulo $T^{E_i}$
by an induction on $i$.
The monoid $G$ is generated by the set of symplectic
areas of the all pseudo-holomorphic maps (polygons, strips and etc.)
which appear during the construction. Such $G$ is discrete
by Gromov compactness when $\mathbb L$ is finite.
In the case $\mathbb L$ is infinite we cannot take
such a {\it discrete} submonoid $G$.

In a certain situation, we can overcome this problem
by using `homotopy inductive limit' as follows.
Suppose we have a countable set of spin
immersed Lagrangian submanifolds $\mathbb L$
of~$(X,\omega)$.
We take finite subsets $\mathbb L^{(j)}$ of $\mathbb L$
for each $j =1,2,3,\dots$ such that $\mathbb L^{(j)} \subset \mathbb L^{(j+1)}$
and the union of all $\mathbb L^{(j)}$ is $\mathbb L$.
For each $j$, we can take a discrete submonoid $G_j$ such that
we can construct a $G_j$-gapped filtered $A_{\infty}$ category $\mathfrak{Fukst}\bigl((X,\omega),\mathbb L^{(j)}\bigr)$
from the finite set $\mathbb L^{(j)}$.
We may assume~${G_j \subset G_{j+1}}$.
We next regard \smash{$\mathfrak{Fukst}\bigl((X,\omega),\mathbb L^{(j)}\bigr)$}
as a $G_{j+1}$-gapped filtered $A_{\infty}$ category.
Then we can construct a
$G_{j+1}$-gapped filtered $A_{\infty}$ functor
$\mathfrak{Fukst}\bigl( (X,\omega),\mathbb L^{(j)}\bigr)
\to \mathfrak{Fukst}\bigl( (X,\omega),\mathbb L^{(j+1)}\bigr)$
which is a homotopy equivalence to the image.
In this way we can construct an inductive system of
filtered $A_{\infty}$ categories.

In the case when the completion (with respect to the Hofer--Chekanov distance \cite{chekanov2})
of the set of Lagrangian submanifolds we study is separable, we can use
the above sequence and construct
the inductive limit
$
{\varinjlim} \mathfrak{Fukst}\bigl((X,\omega),\mathbb L^{(j)}\bigr) = \mathfrak{Fukst}((X,\omega),\mathbb L)$,
see \cite{Hausdoredd}.
The author does not know how much the separability assumption is essential.

\item[(B)]
Let $L_{12} \subset -X_1 \times X_2$ and $L_{23} \subset -X_2 \times X_3$
be immersed Lagrangian correspondences.
If the fiber product $L_{12} \times_{X_2} L_{23}$ is transversal,
then it becomes a Lagrangian correspondence $\subset -X_1 \times X_3$.
If those Lagrangian correspondences are self-clean,
then assuming $L_{12}$, $L_{23}$ are unobstructed,
we proved that the composition of correspondence functors
$\mathcal W_{(L_{12},b_{12})}$ and
$\mathcal W_{(L_{23},b_{23})}$
are represented by an unobstructed Lagrangian correspondence
$(L_{13},b_{13})$.
(We need to restrict ourselves to a~finite set of Lagrangian submanifolds
because of point~(A).)

However, if the fiber product $L_{12} \times_{X_2} L_{23}$ is not transversal,
there is no good candidate of a~Lagrangian correspondence
representing the composition of correspondence functors.

A possible way to resolve this issue is using the result of \cite{Hausdoredd}
as follows.
We perturb $L_{23}$ to~$L^{\varepsilon}_{23}$ by a small Hamiltonian isotopy.
Then we obtain
$L_{13}^{\varepsilon} = L_{12} \times_{X_2} L^{\varepsilon}_{23}$
which is an immersed Lagrangian correspondence.
If $b_{12}$ and $b_{23}$ are bounding cochains of
$L_{12}$ and $L_{23}$ respectively, then
we obtain a bounding cochain $b_{13}^{\varepsilon}$ of $L_{13}^{\varepsilon}$.
We can show that for $\varepsilon_n \to 0$, the sequence
$(L_{13}^{\varepsilon_n},b_{13}^{\varepsilon_n})$ becomes a
Cauchy sequence with respect to the Hofer distance as objects of
$\mathfrak{Fukst}(-X_1 \times X_3)$ (see Definition~\ref{defn151}).\footnote{Moreover, it is a Cauchy sequence with respect to the Hofer infinite distance,
which is introduced in \cite{Hausdoredd}.}
Generalizing various constructions of this paper to the
completion of filtered $A_{\infty}$ category
via Gromov--Hausdorff distance (which is introduced in \cite{Hausdoredd}),
it seems likely that we can define
the composed functor as the limit of \smash{$\mathcal W_{(L_{13}^{\varepsilon_n},b_{13}^{\varepsilon_n})}$}.

We remark that if we change the coefficient from Novikov ring $\Lambda_0$ to
its field of fractions $\Lambda$, then the problem becomes easier to
handle.
In fact, over $\Lambda$ two objects $(L_{23},b_{23})$ and $(L^{\varepsilon}_{23},b_{23}^\varepsilon)$
are equivalent. So we do not need to take the limit as above.
On the other hand, the Lagrangian Floer theory over $\Lambda_0$ is much richer and
contains much more information than the Lagrangian Floer theory over $\Lambda$.

\item[(C)]
By inspecting the proofs of the commutativity of diagrams \eqref{Diagram1} and
\eqref{thm18} given in this paper, we find that they actually do not strictly commute
but commute only up to homotopy equivalence.
It seems likely that there is a certain pseudo-isotopy which
interpolates two compositions appearing in the diagram.
Those pseudo-isotopies are well-defined up to pseudo-isotopy of pseudo-isotopies.
For the composition, we can also try to understand the
`higher associativity', as follows.
In the case when we consider four unobstructed immersed Lagrangian correspondences $(L_{i(i+1)},b_{i(i+1)})$,
$i=1,2,3,4$, from $X_{i}$ to $X_{i+1}$,
the correspondence functors
$\mathcal W_{(L_{12},b_{12})}$, $\mathcal W_{(L_{23},b_{23})}$, $\mathcal W_{(L_{34},b_{34})}$,
$\mathcal W_{(L_{45},b_{45})}$ can be composed in various different orders.
For example,
\begin{gather*}
\mathcal W_{(L_{45},b_{45})} \circ (\mathcal W_{(L_{34},b_{34})} \circ (\mathcal W_{(L_{23},b_{23})} \circ \mathcal W_{(L_{12},b_{12})})), \\
((\mathcal W_{(L_{45},b_{45})} \circ \mathcal W_{(L_{34},b_{34})}) \circ \mathcal W_{(L_{23},b_{23})}) \circ \mathcal W_{(L_{12},b_{12})}
\end{gather*}
and etc.
There exist pseudo-isotopies between the compositions with different orders.
Moreover, it seems likely that one can construct a pseudo-isotopies of pseudo-isotopies
parametrized by the Stasheff 2-gon.
It seems likely that one can continue and obtain a certain infinite category type
construction, if the issues (A), (B) are resolved.
\end{enumerate}

The above discussions sketch a possible way to proceed to
overcome (A), (B), (C) and actually prove `Informal Summary \ref{thm02}'.
However,
the actual works needed to carry out those plans are extremely heavy and
likely become extremely lengthy.
So I think taking a break at the point where we proved the results in this paper
before going further is a reasonable choice.

\subsection{Relations to the works by Bottman--Wehrheim}
\label{reltoBW}

In this subsection, we mention relations of this paper with several papers by
Bottman \cite{bott2,bott} and Bottman--Wehrheim \cite{bww}.
First we review briefly the method of a strip shrinking, introduced by
Wehrheim--Woodward, in the simplest case.
In Section~\ref {subsec:bi-functorgeo1}, we consider
a moduli space consisting of $u_1$, $u_2$ where
$u_1 \colon [-1,0]\times \R \to X_1$ and $u_2 \colon [0,1] \to X_2$
are pseudo-holomorphic maps,
which satisfy a certain matching (boundary) condition
at $\{0\} \times \R$.
One can generalize this moduli space
so that $u_1$ is a map
from $[-S,0]\times \R$ for a certain $S > 0$.
We denote by~${\mathcal M}_{\rm QT}(L_1,L_{12},L_{2};S)$
the moduli space obtained in this way.\footnote{It is an analogue of the
moduli space ${\mathcal M}_{\rm QT}(\vec a_1,\vec a_{12},\vec a_2;a_-,a_+;E)$
introduced in Section~\ref{subsec:bi-functorgeo1}.
Since the discussion here is heuristic, I do not
include the marked points or energy in the notation.}
We can then proceed in the same way to obtain a~tri-module,
which we denote by
$
\mathscr{CF}(\mathbb L_1,\mathbb L_{12};\mathbb L_2;S)$.
We can use it instead of $\mathscr{CF}(\mathbb L_1,\mathbb L_{12};\mathbb L_2)$
to obtain a filtered $A_{\infty}$ functor,
$
\mathcal W_{L_{12}}^S\colon \mathfrak{Fukst}(X_1;\mathbb L_1) \to \mathfrak{Fukst}(X_2;\mathbb L_2)
$
in the same way as Sections~\ref{sec:Unobstructedness} and~\ref{sec:represent}.

\begin{figure}[ht]
\centering
\includegraphics[scale=0.4]{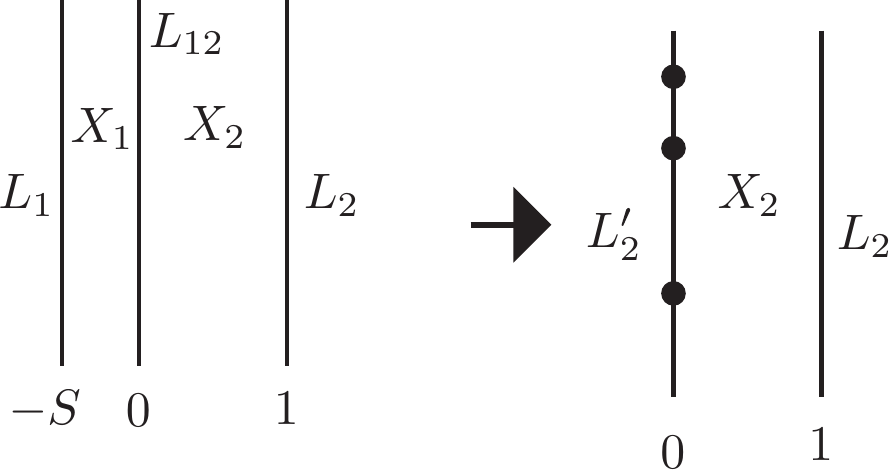}
\caption{Strip shrinking 1.}
\label{shrink1}
\end{figure}

Wehrheim--Woodward--Ma'u--Bottman studied the limit when $S$ goes to zero.
It is believed that the limit
$\lim_{S\to 0}{\mathcal M}_{\rm QT}(L_1,L_{12},L_{2};S)$
becomes a
moduli space $\mathcal M(L'_2,L_2)$ together with bubbles
on the boundary $\{0\} \times \R$.
Here $\mathcal M(L'_2,L_2)$ is a moduli space of pseudo-holomorphic
maps~${u\colon [0,1] \to X_2}$ such that
$u(0,\tau) \in L'_2$ and $u(1,\tau) \in L_2$
and $L'_2$ is the geometric transformation~${L_1 \times_{X_1} {L_{12}}}$.
See Figure~\ref{shrink1}.

\begin{figure}[ht]
\centering
\includegraphics[scale=0.4]{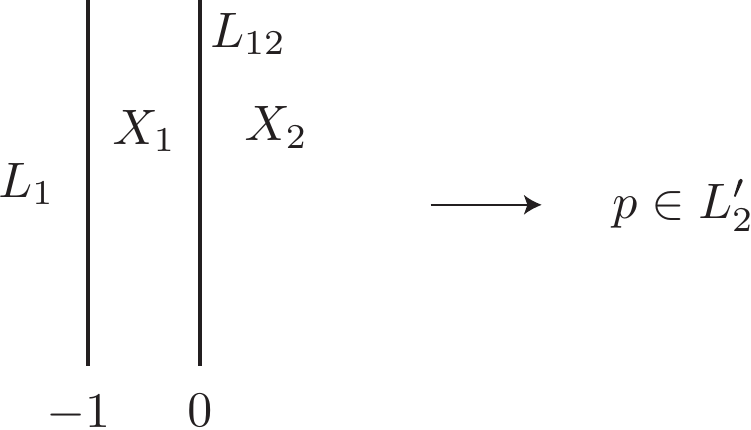}
\caption{Figure 8 bubble 1.}
\label{fig8bubble1}
\end{figure}

The bubble on the line $\{0\} \times \R$ is called a Figure 8 bubble\index{Figure 8 bubble}
and in this case it is expected to be described by a
moduli space of $(u_1,u_2)$ which are pseudo-holomorphic maps
\[
u_1 \colon\ [-1,0] \times \R \to X_1, \qquad
u_2 \colon\ [0,\infty) \times \R \to X_2
\]
with boundary conditions
$
u_1(-1,\tau) \in L_1$, $(u_1(0,\tau),u_2(0,\tau)) \in L_{12}$,
$ \lim_{t \to \infty} u_2(t,\tau) = p$,
where~${p \in L'_2}$ is independent of $\tau$.
See Figure~\ref{fig8bubble1}.

The conjecture mentioned in Remark~\ref{rem16} claims
that the virtual fundamental chain of the moduli space of Figure 8 bubbles
becomes a bounding cochain $b'_2$ of $L_2'$ and the homology
of the tri-module
$
\mathscr{CF}((L_1,b_1),(L_{12},b_{12});(L_2,b_2))
$
becomes isomorphic to $HF((L'_2,b'_2),(L_2,b_2))$.

We conjecture also that $b_2'$ is
gauge equivalent to the bounding cochain we obtained
in Theorem~\ref{them4} as follows.

We consider the bounding cochains
$b'_2(S)$ such that $(L'_2,b'_2(S)) = \mathcal W_{L_{12}}^S(L_1,b_1)$.
Using the fact that
$\mathscr{CF}(\mathbb L_1,\mathbb L_{12};\mathbb L_2;S)$
is pseudo-isotopic to
$\mathscr{CF}(\mathbb L_1,\mathbb L_{12};\mathbb L_2;S')$
for $S,S' > 0$,
we can show that $b'_2(S)$ is independent of $S$
up to gauge equivalence.
Note that $b'_2(S)$ is characterized by the condition that
\begin{equation}\label{eq1855}
\mathfrak n_S\bigl(e^{b_1},e^{b_{12}};{\bf 1};e^{b'_2(S)}\bigr) = 0,
\end{equation}
where $\mathfrak n_S$ is the structure operation of the
tri-module $\mathscr{CF}(L_1,L_{12};L'_2;S)$ and ${\bf 1}$ is the
cyclic element, which is the 0-form $1$ on
the diagonal component of the fiber product $L_1 \times_{X_1} L_{12} \times_{X_2} L'_2
\cong L'_2 \times_{X_2} L'_2$.

The tri-module $\mathscr{CF}(L_1,L_{12};L'_2;S)$ is expected
to `converge' to the Floer chain complex
$CF(L'_2,L'_2)$ with the boundary operator corrected by $b'_2$,
which is
\begin{equation}\label{eq18552}
d(x) = \mathfrak m\bigl(e^{b'_2},x\bigr).
\end{equation}
Here $b'_2$ is the conjectured bounding cochain
obtained from Figure 8 bubbles.\footnote{Note that $d \circ d = 0$
may not hold for the operator \eqref{eq18552}. We need to add bounding cochain of $CF(L'_2,L'_2)$
which acts from the {\it right} also to obtain
$d'$ such that $d' \circ d' = 0$.}
It is easy to see that
\begin{equation}\label{eq1865}
\mathfrak m\bigl(e^{b'_2},{\bf 1},e^{b'_2}\bigr) = 0.
\end{equation}
Here $\mathfrak m$ is the structure operation of the $A_{\infty}$ algebra
associated to $L'_2$ and ${\bf 1}$ is the fundamental class, which is
the $0$ form $1$ of $L'_2$.
Comparing \eqref{eq1855} and \eqref{eq1865},
we expect
$
\lim_{S \to 0} b'_2(S) = b'_2$.
Namely, the bounding cochain obtained from the moduli space of Figure 8 bubbles
is gauge equivalent to one in Theorem~\ref{them4}.

We mention a reason\footnote{Which was known to
various researchers before the year 2010.} why the
virtual fundamental chain of the moduli space
of Figure 8 bubbles is not yet rigorously constructed.
We draw Figure~\ref{fig8bubble1} on the 2 sphere as in Figure~\ref{fig8bubble2} below.
\begin{figure}[ht]
\centering
\includegraphics[scale=0.4]{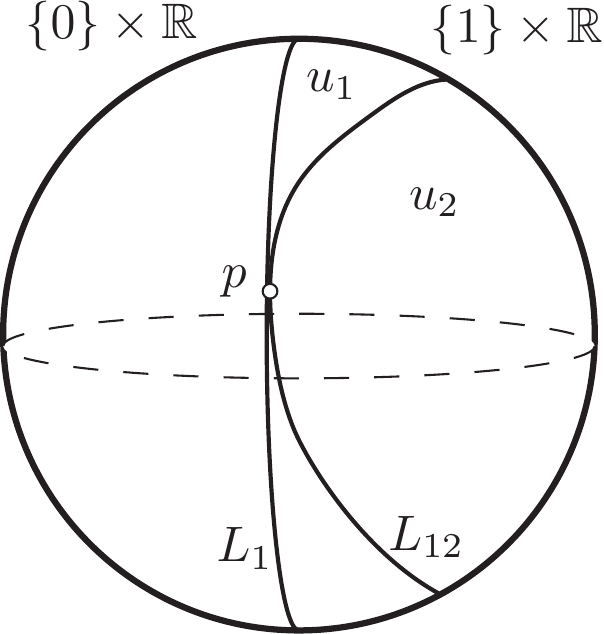}
\caption{Figure 8 bubble 2.}
\label{fig8bubble2}
\end{figure}

Two lines (seams) $\{0\} \times \R$ and $\{1\} \times \R$
on which we require boundary conditions
are tangent at the point $\infty$ (which is required to be sent to $p$).
This is different from the situation of the Y-diagram,
where 3 seams intersect {\it transversally} at the hole.
The existence of tangency between seems is a new phenomenon
and Fredholm theory for such boundary valued problem is not
yet established.
We like to mention Bottman \cite{bott2,bott} and Bottman--Wehrheim~\cite{bww}
established compactness and removable singularity, which is a
very important step toward constructing the virtual
fundamental chain of the moduli space of Figure 8 bubbles.

\subsection{Relation to the works by Ma'u--Wehrheim--Woodward}
\label{relWWM}

As we mentioned in the introduction, Weinstein \cite{Wi} proposed to regard a Lagrangian
submanifold of the product $-X \times Y$
as a morphism $X \to Y$ between symplectic manifolds.
Since Weinstein's proposal looks so natural,
there had been attempts to associate a functor~${\mathfrak F_{\mathcal L} \colon \mathfrak{Fukst}(X) \to \mathfrak{Fukst}(Y)}$
to an unobstructed immersed Lagrangian correspondence $\mathcal L = (L,b)$.
A possible naive idea to do so is the following.
Let $L_1$ be a Lagrangian submanifold of~$X$.
Instead of associating an object of~$\mathfrak{Fukst}(Y)$ to $L_1$,
we try to define a right~$\mathfrak{Fukst}(Y)$ module~$\mathfrak F_{\mathcal L}(L_1)$.
In the cohomology level, $\mathfrak F_{\mathcal L}(L_1)$ can be defined
by associating the Floer homology~$HF(L;L_1\times L_2)$ in the product $-X \times Y$
to a Lagrangian submanifold $L_2$ of $Y$.
Actually we can construct an $A_{\infty}$ functor
\begin{equation}\label{corr71}
\mathfrak F_{\mathcal L}\colon \ \mathfrak{Fukst}(X) \to \mathcal{RMOD}(\mathfrak{Fukst}(Y))
\end{equation}
in this way, as we did in Section~\ref{sec:Kunneth}.
Here $\mathcal{RMOD}(\mathfrak{Fukst}(Y))$ is the DG-category
of right $\mathfrak{Fukst}(Y)$ modules.
Because of Yoneda's lemma, an object of
$\mathcal{RMOD}(\mathfrak{Fukst}(Y))$ can be
regarded as an `extended object' of $\mathfrak{Fukst}(Y)$. Thus \eqref{corr71}
could be regarded as a version of
$\mathfrak F_{\mathcal L} \colon \mathfrak{Fukst}(X) \to \mathfrak{Fukst}(Y)$.

However, the problem is in this formulation it is difficult
to compose $\mathfrak F_{\mathcal L_{12}}$ and $\mathfrak F_{\mathcal L_{23}}$
where $\mathcal L_{i(i+1)} = (L_{i(i+1)},b_{i(i+1)})$ is an unobstructed immersed Lagrangian submanifold of
$-X_i \times X_{i+1}$, for $i=1,2$. This point is mentioned also in the first page of
\cite{MWW}.
In the early 2000's, the author tried to resolve this problem
by a purely algebraic method of homological algebra of
$A_{\infty}$ categories, but he was not successful.\footnote{Theorem~\ref{them4} resolves
this problem by using more geometric input.}
\begin{rem}

The above naive idea can be regarded as a `finite-dimensional analogue'
of the proposal \cite{fu2} to construct instanton Floer homology of 3-manifolds
with boundary as an $A_{\infty}$ module.
In the moduli space introduced in \cite{gafa} during the attempt to realize the
proposal, a line in the domain of $\C$ where the equation changes from the ASD-equation
(on a 4-manifold) to the pseudo-holomorphic curve equation, appears.
This line plays the same role as seams play in the study of
Lagrangian correspondences.
The moduli space introduced by Lipyanskiy \cite{Li} is more directly
an infinite-dimensional analogue of the moduli space
of pseudo-holomorphic quilts.
\end{rem}

As mentioned in Remark~\ref{rem16}, Wehrheim--Woodward--Ma'u
used the following idea to go around this problem.
For a given symplectic manifold $X$,
they consider a series of Lagrangian correspondences
$
L_{i} \subset -X_i \times X_{i+1}
$
such that $X_0$ is a~point and $X_n = X$.
They regard such a~system $(L_0,\dots,L_n)$ as
an object of expanded category $\mathfrak{Fuk}^{\#}(X)$.
Then, if $L'' \subset -X \times Y$ is a~Lagrangian correspondence,
one can define
$
(\mathcal{W}_{{\mathcal L}})_{\rm ob}\colon
\mathfrak{OB}\bigl(\mathfrak{Fuk}^{\#}(X)\bigr) \to \mathfrak{OB}\bigl(\mathfrak{Fuk}^{\#}(Y)\bigr),
$
by~${(L_0,\dots,L_n) \mapsto (L_0,\dots,L_n,L'')}$.

To define the $A_{\infty}$ category $\mathfrak{Fuk}^{\#}(X)$, one needs to define the
Floer homology
 between extended objects
$(L_0,\dots,L_n)$, $(L'_0,\dots,L'_{n'})$,
where $L_{i} \subset -X_i \times X_{i+1}$
and $L'_{i} \subset -X'_i \times X'_{i+1}$,
$X_0$, $X'_0$ are points and $X_n = X'_{n'} = X$.
They denote this Floer homology by $HF(L_0,\dots,L_n,L'_{n'},\dots,\allowbreak L'_0)$.
Wehrheim--Woodward--Ma'u used the notion of
a pseudo-holomorphic quilt to define it.
The pseudo-holomorphic quilt used to define
$HF(L_0,\dots,L_n,L'_{n'},\dots,L'_0)$
is as in Figure~\ref{Figure7} below. Here $u_i$ (resp.\ $u'_i$) is a pseudo-holomorphic map to $X_i$ (resp.\ $X'_i$)
and $u$ is a pseudo-holomorphic map to~$X$.

\begin{figure}[ht]
\centering
\includegraphics[scale=0.45]{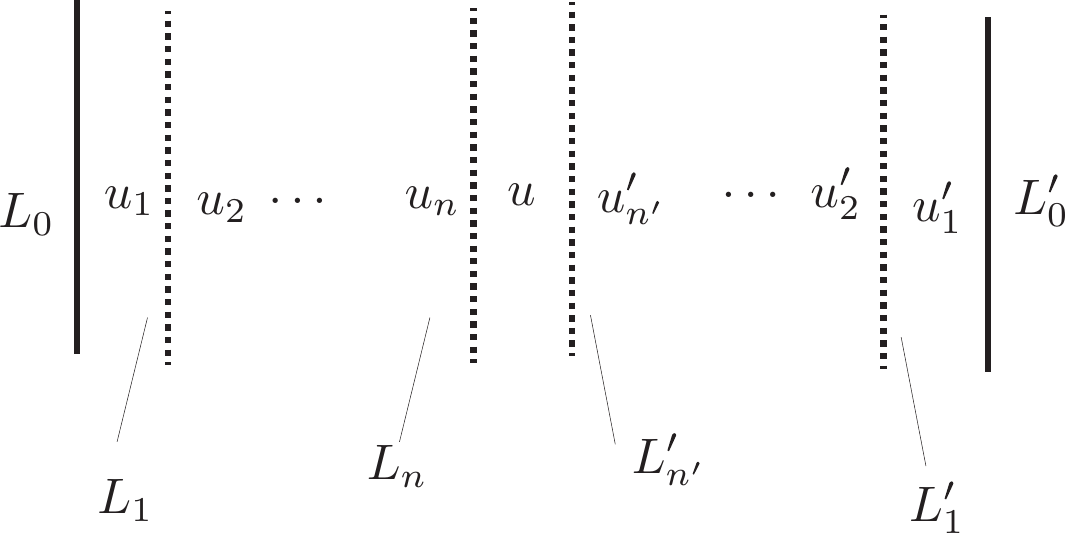}
\caption{A pseudo-holomorphic quilt.}
\label{Figure7}
\end{figure}

Wehrheim--Woodward--Ma'u went further to define a version of the correspondence bi-func\-tor~${
\mathcal{MWW}\colon \mathfrak{Fuk}^{\#}(-X\times Y) \times \mathfrak{Fuk}^{\#}(X)
\to \mathfrak{Fuk}^{\#}(Y)}$.
Their works are very important contributions to the study of
Lagrangian correspondence and Lagrangian Floer homology.

We remark that
in a way similar to Theorem~\ref{thm164}
(and using reflection principle in a similar way as we used
in Section~\ref{sec:orient}),
we can show the next isomorphism.
\begin{equation}\label{formula72}
HF(L_0,\dots,L_n,L'_{n'},\dots,L'_0)
\cong
HF(L_0\times\dots\times L_n\times L'_0\times \dots\times L'_{n'};
\Delta).
\end{equation}
Here
\begin{equation}\label{formula73}
\Delta \subset \left(\prod_{i=1}^{n-1} (-X_i\times X_i)\right) \times
\left(\prod_{i=1}^{n'-1} (-X'_{i'}\times X'_{i'})\right) \times (-X \times X)
\end{equation}
is the product of diagonals. The right-hand side of \eqref{formula72}
is the Floer homology of two Lagrangian submanifolds
in the symplectic manifold given in \eqref{formula73}.

\begin{figure}[ht]
\centering
\includegraphics[scale=0.4]{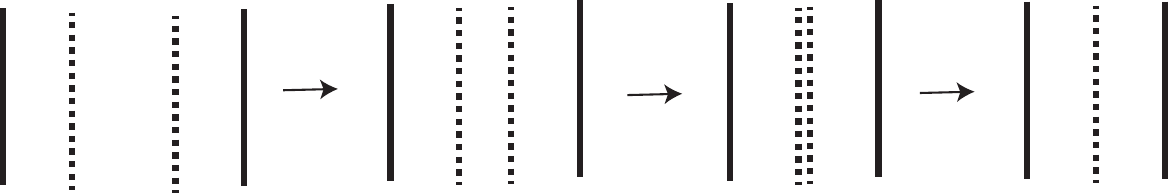}
\caption{Strip shrinking 2.}
\label{Figure8}
\end{figure}

The advantage to use a pseudo-holomorphic quilt rather than
Floer homology in the direct product (as in \eqref{formula72})
lies in the fact that, then, one can use a strip shrinking to prove the
next important isomorphism
\begin{equation}\label{form74}
HF(L_0,\dots,L_n,L'_{n'},\dots,L'_0)
\cong
HF(L_0,\dots,L_{n-1},L_n\times_XL'_{n'},L'_{n'-1},\dots,L'_0).
\end{equation}
As mentioned in the last subsection, a strip shrinking is a process to change the width between two seams
until it becomes $0$
(see Figure~\ref{Figure8}).
Note that the method of using reflection principle to replace
Wehrheim--Woodward's definition by~\eqref{formula72} works only
in the case when all the strips have the same width.
Therefore, it is not consistent with strip shrinking.

Wehrheim--Woodward proved the isomorphism \eqref{form74}
under the assumption that all the Lagrangian submanifolds
involved (including the fiber product $L_n\times_XL'_{n'}$)
are embedded and monotone.
The isomorphism \eqref{form74} is a version of
composability of filtered $A_{\infty}$ functors
associated to the composition of Lagrangian correspondences.

The reason why one does not need to study Figure 8 bubbles
in the case when all the Lagrangian submanifolds
involved are embedded and monotone is as follows.
One can show that
if the Figure 8 bubble occurs then
it carries a strictly positive energy and so
the virtual dimension of the moduli space
of the configuration drops at least 2
in the monotone case.
By this dimension counting argument,
one can avoid Figure 8 bubbles
in the monotone situation. Later
Lekili and Lipyanskiy \cite{LL} gave an alternative proof of \eqref{form74}
using Y-diagram. (They assume embeddedness and monotonicity.)

This is somewhat similar to the
usual Floer theory or Gromov--Witten theory.
In the semi-positive case,
one can avoid sphere bubbles by the
dimension counting argument.
Therefore, one does not need to
find a Kuranishi chart at such `infinity'.
When we study a symplectic manifold which is not semi-positive then
we need an abstract perturbation and so we
need a~chart centered at a point of infinity, which
corresponds to a stable map with sphere bubbles.

\printindex[syindex]
\printindex

\subsection*{Acknowledgements}

The author would like to thank Simons Center for Geometry and Physics where
most of the research written in this paper is performed.
The author would like to thank J.~Evans, Y.~Lekili,
Y.-G.~Oh, H.~Ohta, K.~Ono for helpful discussions
while he is working on the contents of this paper.
He also would like to thank anonymous referees for careful and serious reading,
which is a heavy and pains taking work, and for huge number of important comments,
which improve this paper much compared to its earlier version.
Special thanks are to K.~Ono who agreed to write an article \cite{ono2}
on sign and orientation which we need in this paper.


\LastPageEnding

\end{document}